\documentclass[a4paper,10pt]{preprint}
\usepackage[a4paper,top=1.6in,bottom=1.6in,left=1.3in,right=1.3in]{geometry}
\usepackage[full]{textcomp}
\usepackage[osf]{newtxtext}

\usepackage{xstring}
\usepackage{amsmath}
\usepackage{amssymb}
\usepackage{amsthm}

\usepackage{mathtools}
\usepackage[scr=boondoxo]{mathalfa}
\usepackage{bbm} 
\usepackage{esint}
\usepackage{halloweenmath}
\usepackage{graphicx}
\usepackage{afterpage}

\usepackage[utf8]{inputenc}
\usepackage{xcolor}

\usepackage{margincomment}

\usepackage{cprotect} 
\usepackage{tikz} 

\usepackage{SYM-diagrams}
\usepackage[colorlinks=true,linkcolor=darkred,citecolor=darkblue]{hyperref}
\usepackage{mhequ}
\usepackage{orcidlink}

\usepackage{microtype}
\usepackage{tabularx}
\usepackage{colortbl} 
\usepackage{booktabs} 
\usepackage{longtable} 
\usepackage{ifthen}

\usepackage{enumitem}
\usepackage{subcaption}
\usepackage{array}

\newcommand{\hlmath}[1]{\!\colorbox{black!7!white}{$\displaystyle #1$}\!}

\let\oldskull\skull
\def\skull{\mathord{\oldskull}}

\usetikzlibrary{external}
\tikzsetexternalprefix{tikz/}
\ifnum\pdfstrcmp{\jobname}{output}=0
\else
\fi

\DeclareSymbolFont{timesoperators}{T1}{ptm}{m}{n}
\SetSymbolFont{timesoperators}{bold}{T1}{ptm}{b}{n}
\DeclareMathAlphabet{\mathbb}{U}{jkpsyb}{m}{n}
\SetMathAlphabet{\mathbb}{bold}{U}{jkpsyb}{bx}{n}
\SetMathAlphabet{\mathrm}{normal}{T1}{ptm}{m}{n}  
\SetMathAlphabet{\mathbf}{normal}{T1}{ptm}{b}{n}

\makeatletter
\renewcommand{\operator@font}{\mathgroup\symtimesoperators}
\makeatother

\newcolumntype{P}[1]{>{\centering\arraybackslash}p{#1}}		
\renewcommand\arraystretch{1.6} 

\theoremstyle{definition}
\newtheorem{definition}{Definition}
\newtheorem{lemma}[definition]{Lemma}
\newtheorem{proposition}[definition]{Proposition}
\newtheorem{corollary}[definition]{Corollary}
\newtheorem{notation}[definition]{Notation}

\newtheorem{assumption}[definition]{Assumption}
\newtheorem{setting}[definition]{Setting}

\theoremstyle{remark}
\newtheorem{remark}[definition]{Remark}
\newtheorem{example}[definition]{Example}

\theoremstyle{plain}
\newtheorem{theorem}[definition]{Theorem}

\numberwithin{equation}{section}
\numberwithin{definition}{section}

\newcommand{\lo}{\textup{lo}}

\newcommand{\covd}{\textup{D}}

\renewcommand{\Im}{\textup{Im}}
\renewcommand{\Re}{\textup{Re}}

\newcommand{\R}{\ensuremath{\mathbb{R}}}
\newcommand{\C}{\ensuremath{\mathbb{C}}} 
\newcommand{\Z}{\mathbb{Z}}
\newcommand{\T}{\mathbb{T}}
\newcommand{\N}{\mathbb{N}}

\newcommand{\initialclass}{\mathfrak{I}}

\newcommand{\Cs}{\mathscr{C}}
\let\CC\Cs
\let\cC\Cs

\newcommand{\nl}{\mathrm{nl}}

\newcommand{\md}{\mathrm{d}}
\newcommand{\du}{\mathrm{d}u}

\newcommand{\dr}{\mathrm{d}r}
\newcommand{\ds}{\mathrm{d}s}
\newcommand{\dt}{\mathrm{d}t}
\newcommand{\dx}{\mathrm{d}x} 
\newcommand{\dy}{\mathrm{d}y}
\newcommand{\dz}{\mathrm{d}z}

\newcommand*{\medcap}{
\mathbin{\raisebox{1pt}{\scalebox{0.8}{\ensuremath{\bigcap}}}}}

\DeclareMathOperator{\Kil}{Kil}
\DeclareMathOperator{\Duh}{Duh}
\DeclareMathOperator{\hol}{hol}

\newcommand{\extraC}{\mathsf{C}}
\newcommand{\Crho}{\mathsf{C}_\rho}

\newcommand{\hCrho}{\hat{\mathsf{C}}_\rho}
\newcommand{\Czero}{\mathsf{C}_0}
\newcommand{\Cone}{\mathsf{C}_1}
\newcommand{\Ctwo}{\mathsf{C}_2}

\newcommand{\mass}{\mathsf{M}}
\newcommand{\massA}{\mathsf{N}}

\def\Wick#1{\mathopen{\scalebox{1.2}{:}}#1\mathclose{\scalebox{1.2}{:}}}

\newcommand{\beq}{\begin{equation}}
\newcommand{\eeq}{\end{equation}}

\newcommand{\mbf}[1]{\mathbf{#1}}
\newcommand{\mrm}[1]{\mathrm{#1}}
\newcommand{\mc}[1]{\mathcal{#1}}

\newcommand{\ostar}{\star}
\newcommand{\odiamond}{\diamond}

\newcommand{\ra}{\rightarrow}
\newcommand{\Tr}{\mrm{Tr}}
\newcommand{\sse}{\subseteq}
\newcommand{\ptl}{\partial}

\newcommand{\ind}{\mathbbm{1}}
\newcommand{\fa}{\mathfrak{a}}

\newcommand{\frkh}{\mathfrak{h}}
\newcommand{\dyadic}{2^{\mathbb{N}_0}}
\newcommand{\sdyadic}{2^{-\mathbb{N}_0}}

\newcommand{\op}{\operatorname{op}}

\newcommand{\toinf}{\ra \infty}

\newcommand{\bP}{\mathbb{P}}
\DeclareMathOperator{\Law}{Law}
\newcommand{\eps}{\varepsilon}
\newcommand{\varep}{\varepsilon}
\newcommand{\Wrem}{W}

\newcommand{\reg}{\delta}      
\newcommand{\lossdel}{\delta} 

\newcommand{\icomplex}{\textup{\textbf{i}}}
\newcommand{\Ad}{\mrm{Ad}}
\newcommand{\ad}{\mrm{ad}}

\def\lsup#1{{}^{#1}\kern-0.3ex}
\def\TV{{\mathrm{TV}}}

\newcommand{\E}{\mathbb{E}} 

\newcommand{\frkg}{\mathfrak{g}}

\newcommand{\End}{\textup{End}}
\DeclareMathOperator{\Cov}{Cov}
\def\Cas{\Cov_{\frkg}^{\frkg}}

\makeatletter 
\DeclareRobustCommand*{\act}{%
  \mathbin{\mathpalette\bigcdot@{}}%
}
\DeclareRobustCommand*{\place}{%
  \mathord{\kern0.1em\mathpalette\bigcdot@{}\kern0.1em}%
}
\newcommand*{\bigcdot@scalefactor}{.5}
\newcommand*{\bigcdot@widthfactor}{1.15}
\newcommand*{\bigcdot@}[2]{%
  \sbox0{$#1\vcenter{}$}
  \sbox2{$#1\cdot\m@th$}%
  \hbox to \bigcdot@widthfactor\wd2{%
    \hfil
    \raise\ht0\hbox{%
      \scalebox{\bigcdot@scalefactor}{%
        \lower\ht0\hbox{$#1\bullet\m@th$}%
      }%
    }%
    \hfil
  }%
}
\makeatother

\def\dash{\leavevmode\unskip\kern0.15em\scalebox{0.85}[1]{---}\penalty\exhyphenpenalty\kern0.15em}
\def\slash{\leavevmode\unskip\kern0.15em/\penalty\exhyphenpenalty\kern0.15em}

\colorlet{symbols}{blue!90!black}

\makeatletter
\def\DeclareSymbol#1#2#3{%
	\expandafter\gdef\csname MH@symb@#1b\endcsname{{\tikzsetnextfilename{symbol#1b}\colorlet{symbols}{black}%
	\tikz[baseline=#2,scale=0.15,draw=symbols,line join=round,line cap=round]{#3}}}%
	\expandafter\gdef\csname MH@symb@#1bs\endcsname{\scalebox{0.75}{\tikzsetnextfilename{symbol#1b}\colorlet{symbols}{black}%
	\tikz[baseline=#2,scale=0.15,draw=symbols,line join=round,line cap=round]{#3}}}%
	\expandafter\gdef\csname MH@symb@#1bss\endcsname{\scalebox{0.65}{\tikzsetnextfilename{symbol#1b}\colorlet{symbols}{black}%
	\tikz[baseline=#2,scale=0.15,draw=symbols,line join=round,line cap=round]{#3}}}%
	\expandafter\gdef\csname MH@symb@#1r\endcsname{{\tikzsetnextfilename{symbol#1r}\colorlet{symbols}{darkred}%
	\tikz[baseline=#2,scale=0.15,draw=symbols,line join=round,line cap=round]{#3}}}%
	\expandafter\gdef\csname MH@symb@#1rs\endcsname{\scalebox{0.75}{\tikzsetnextfilename{symbol#1r}\colorlet{symbols}{darkred}%
	\tikz[baseline=#2,scale=0.15,draw=symbols,line join=round,line cap=round]{#3}}}%
	\expandafter\gdef\csname MH@symb@#1rss\endcsname{\scalebox{0.65}{\tikzsetnextfilename{symbol#1r}\colorlet{symbols}{darkred}%
	\tikz[baseline=#2,scale=0.15,draw=symbols,line join=round,line cap=round]{#3}}}%
	\expandafter\gdef\csname MH@symb@#1\endcsname{\tikzsetnextfilename{symbol#1}%
	\tikz[baseline=#2,scale=0.15,draw=symbols,line join=round,line cap=round]{#3}}%
	\expandafter\gdef\csname MH@symb@#1s\endcsname{\scalebox{0.75}{\tikzsetnextfilename{symbol#1}%
	\tikz[baseline=#2,scale=0.15,draw=symbols,line join=round,line cap=round]{#3}}}%
	\expandafter\gdef\csname MH@symb@#1ss\endcsname{\scalebox{0.65}{\tikzsetnextfilename{symbol#1}%
	\tikz[baseline=#2,scale=0.15,draw=symbols,line join=round,line cap=round]{#3}}}%
	}
\def\<#1>{\ifthenelse{\boolean{mmode}}{\mathchoice{\csname MH@symb@#1\endcsname}{\csname MH@symb@#1\endcsname}{\csname MH@symb@#1s\endcsname}{\csname MH@symb@#1ss\endcsname}}{\csname MH@symb@#1\endcsname}}
\makeatother

\makeatletter
\pgfdeclareshape{crosscircle}
{
  \inheritsavedanchors[from=circle] 
  \inheritanchorborder[from=circle]
  \inheritanchor[from=circle]{north}
  \inheritanchor[from=circle]{north west}
  \inheritanchor[from=circle]{north east}
  \inheritanchor[from=circle]{center}
  \inheritanchor[from=circle]{west}
  \inheritanchor[from=circle]{east}
  \inheritanchor[from=circle]{mid}
  \inheritanchor[from=circle]{mid west}
  \inheritanchor[from=circle]{mid east}
  \inheritanchor[from=circle]{base}
  \inheritanchor[from=circle]{base west}
  \inheritanchor[from=circle]{base east}
  \inheritanchor[from=circle]{south}
  \inheritanchor[from=circle]{south west}
  \inheritanchor[from=circle]{south east}
  \inheritbackgroundpath[from=circle]
  \foregroundpath{
    \centerpoint%
    \pgf@xc=\pgf@x%
    \pgf@yc=\pgf@y%
    \pgfutil@tempdima=\radius%
    \pgfmathsetlength{\pgf@xb}{\pgfkeysvalueof{/pgf/outer xsep}}%
    \pgfmathsetlength{\pgf@yb}{\pgfkeysvalueof{/pgf/outer ysep}}%
    \ifdim\pgf@xb<\pgf@yb%
      \advance\pgfutil@tempdima by-\pgf@yb%
    \else%
      \advance\pgfutil@tempdima by-\pgf@xb%
    \fi%
    \pgfpathmoveto{\pgfpointadd{\pgfqpoint{\pgf@xc}{\pgf@yc}}{\pgfqpoint{-0.707107\pgfutil@tempdima}{-0.707107\pgfutil@tempdima}}}
    \pgfpathlineto{\pgfpointadd{\pgfqpoint{\pgf@xc}{\pgf@yc}}{\pgfqpoint{0.707107\pgfutil@tempdima}{0.707107\pgfutil@tempdima}}}
    \pgfpathmoveto{\pgfpointadd{\pgfqpoint{\pgf@xc}{\pgf@yc}}{\pgfqpoint{-0.707107\pgfutil@tempdima}{0.707107\pgfutil@tempdima}}}
    \pgfpathlineto{\pgfpointadd{\pgfqpoint{\pgf@xc}{\pgf@yc}}{\pgfqpoint{0.707107\pgfutil@tempdima}{-0.707107\pgfutil@tempdima}}}
  }
}
\makeatother

\tikzset{
root/.style={circle,fill=black!50,inner sep=0pt, minimum size=3mm},
        var/.style={circle,fill=black!10,draw=black,inner sep=0pt, minimum size=3mm},
        kernel/.style={semithick,shorten >=2pt,shorten <=2pt},
        kernel1/.style={thick},
        kernels/.style={snake=zigzag,shorten >=2pt,shorten <=2pt,segment amplitude=1pt,segment length=4pt,line before snake=2pt,line after snake=5pt,},
	kernels2/.style={ultra thick,draw=symbols,segment length=12pt},
        rho/.style={densely dashed,semithick,shorten >=2pt,shorten <=2pt},
           testfcn/.style={dotted,semithick,shorten >=2pt,shorten <=2pt},
           tau/.style={circle,inner sep=1pt,draw=black,fill=white,text=black,thin},
        renorm/.style={shape=circle,fill=white,inner sep=1pt},
        labl/.style={shape=rectangle,fill=white,inner sep=1pt},
        xibig/.style={very thin,circle,fill=symbols!25,draw=symbols,inner sep=0pt,minimum size=1.8mm},
        xi2big/.style={very thin,fill=symbols!25,draw=symbols,inner sep=0pt,minimum size=1.6mm},
        xi/.style={very thin,circle,fill=symbols!25,draw=symbols,inner sep=0pt,minimum size=1.2mm},
        xi2/.style={very thin,fill=symbols!25,draw=symbols,inner sep=0pt,minimum size=1mm},
        xix/.style={very thin,crosscircle,fill=symbols!25,draw=symbols,inner sep=0pt,minimum size=1.2mm},
	not/.style={thin,regular polygon, regular polygon sides=3,draw=connection,fill=connection,inner sep=0pt,minimum size=1.2mm},
	dot/.style={thin,circle,fill=symbols,draw=symbols,inner sep=0pt,minimum size=0.5mm},
	>=stealth,
  }

\let\oldprop\propto
\def\propto{\mathbin{\oldprop}}

\DeclareSymbol{Xi}{-2.5}{\node[xibig] {} (0,0);}
\DeclareSymbol{Xi2}{-2.5}{\node[xi2big] {} (0,0);}

\DeclareSymbol{IXi}{-1}{
\draw (0,1) node[xi] {} -- (0,-.4);}

\DeclareSymbol{I'Xi}{-1}{
\draw (0,1)[kernels2] node[xi] {} -- (0,-.4);}

\DeclareSymbol{I'Xi2}{-1}{
\draw (0,1)[kernels2] node[xi2] {} -- (0,-.4);}

\DeclareSymbol{I'XiX}{-1}{
\draw (0,1)[kernels2] node[xix] {} -- (0,-.4);}

\DeclareSymbol{IXiX}{-1}{
\draw (0,1) node[xix] {} -- (0,-.4);}

\DeclareSymbol{IXi2}{-1}{
\draw (0,1) node[xi2] {} -- (0,-.4);}

\DeclareSymbol{IXiI'Xi}{-1}{\draw (-1,1) node[xi] {} -- (0,0);
\draw (1,1)[kernels2] node[xi] {} -- (0,0) {};}

\DeclareSymbol{I'XiIXi2}{-1}{\draw (-1,1)[kernels2] node[xi] {} -- (0,0);
\draw (1,1) node[xi2] {} -- (0,0) {};}

\DeclareSymbol{IXi2I'Xi}{-1}{\draw (-1,1) node[xi2] {} -- (0,0);
\draw (1,1)[kernels2] node[xi] {} -- (0,0) {};}

\DeclareSymbol{IXiI'Xi2}{-1}{\draw (-1,1) node[xi] {} -- (0,0);
\draw (1,1)[kernels2] node[xi2] {} -- (0,0) {};}

\DeclareSymbol{IXi^3}{0}{\draw (-1.3,1) node[xi] {} -- (0,0);
\draw (1.3,1) node[xi] {} -- (0,0);
\draw (0,1.5) node[xi] {} -- (0,0) {};}

\DeclareSymbol{IXiI'[IXiI'Xi]}{1}{\draw[kernels2] (2,2) node[xi] {} -- (1,1);
\draw (0,2) node[xi] {} -- (1,1) {};
\draw[kernels2] (1,1) -- (0,0) {};
\draw (-1,1) node[xi] {} -- (0,0);}

\DeclareSymbol{I'[IXiI'Xi]}{1}{\draw[kernels2] (2,2) node[xi] {} -- (1,1);
\draw (0,2) node[xi] {} -- (1,1) {};
\draw[kernels2] (1,1) -- (0,0) {};}

\DeclareSymbol{I[IXiI'Xi]I'Xi}{1}{\draw (-2,2) node[xi] {} -- (-1,1);
\draw[kernels2] (0,2) node[xi] {} -- (-1,1) {};
\draw (-1,1) -- (0,0) {};
\draw[kernels2] (1,1) node[xi] {} -- (0,0);}

\DeclareSymbol{I[I[Xi]]Xi}{1}{
\draw (0,2) node[xi] {} -- (-1,1) {};
\draw[snake=snake , segment length=1pt, segment amplitude=1pt] (-1,1) -- (0,0) node[xi] {};}

\DeclareSymbol{I[I'Xi]I'Xi}{1}{\draw[kernels2] (0,2) node[xi] {} -- (-1,1);
\draw (-1,1) {} -- (0,0);
\draw[kernels2] (1,1) node[xi] {} -- (0,0);}

\DeclareSymbol{IXiI'[I'Xi]}{1}{\draw[kernels2] (0,2) node[xi] {} -- (-1,1);
\draw[kernels2] (-1,1)  {} -- (0,0);
\draw (1,1) node[xi] {} -- (0,0);}

\DeclareSymbol{I'[I'Xi]}{1}{\draw[kernels2] (0,2) node[xi] {} -- (-1,1);
\draw[kernels2] (-1,1)  {} -- (0,0);}

\DeclareSymbol{I[I'Xi]}{1}{\draw[kernels2] (0,2) node[xi] {} -- (-1,1);
\draw (-1,1)  {} -- (0,0);}

\DeclareSymbol{IXiI'XXi}{-1}{\draw (-1,1) node[xi] {} -- (0,0);
\draw[kernels2] (1,1) node[xi] {} -- (0,0) {};
\draw (1,1) node[xix] {};}

\DeclareSymbol{IXXiI'Xi}{-1}{\draw (1,1) node[xix] {} -- (0,0);
\draw[kernels2] (-1,1) node[xi] {} -- (0,0) {};}

\DeclareSymbol{IXi^2}{-1}{\draw (-1,1) node[xi] {} -- (0,-.4);
\draw (1,1) node[xi] {} -- (0,-.4);}

\DeclareSymbol{XiIIXi}{0.5}{\draw (0,0) node[xi] {} -- (-1,1) node[dot] {} -- (0,2) node[xi] {};}

\DeclarePairedDelimiter{\sscal}{\langle\kern-.27em\langle}{\rangle\kern-.27em\rangle}
\DeclarePairedDelimiter{\scal}{\langle}{\rangle}
\DeclarePairedDelimiter{\abs}{\lvert}{\rvert}
\DeclarePairedDelimiter{\norm}{\|}{\|}
\DeclarePairedDelimiter{\tripn}{\lvert\!|\!\lvert}{\rvert\!|\!\rvert}

\DeclarePairedDelimiter{\Fnorm}{\|}{\|^{\mrm{pert}}}

\DeclarePairedDelimiter{\revFnorm}{\|}{\|^{\mrm{trep}}}

\DeclarePairedDelimiter{\YMFnorm}{\|}{\|^{\mrm{YM}}}

\definecolor{darkgreen}{rgb}{0.0, 0.55, 0.0}
\definecolor{darkred}{rgb}{0.7, 0.0, 0.0}
\definecolor{darkblue}{rgb}{0, 0.0, 0.7}
\definecolor{orange}{rgb}{0.8, 0.3, 0.0}
\definecolor{purple}{rgb}{0.7, 0.0, 0.7}

\newcommand{\eqeps}{\stackrel{\mbox{\tiny$\eps$}}{=}}
\newcommand{\eqdef}{\stackrel{\mbox{\tiny\rm def}}{=}}
\newcommand{\eqlaw}{\stackrel{\mbox{\tiny\rm law}}{=}}

\newcommand{\renormY}{\mathsf{C}^Y_{\rho,\eps}}
\newcommand{\renormYb}{\bar{\mathsf{C}}^Y_{\rho,\eps}}
\newcommand{\renormA}{\mathsf{C}^A_\rho}
\newcommand{\renormphi}{\mathsf{C}^\higgs_{\rho,\eps}}
\newcommand{\solutionmap}{\mathbb{S}}
\newcommand{\pt}{\tau}
\newcommand{\tstop}{t_{\mathsf{st}}}
\newcommand{\xstop}{x_{\mathsf{st}}}

\newcommand{\liegroup}{{\mathbf{G}}}
\newcommand{\Hliegroup}{{\mathbf{H}}}
\newcommand{\groupid}{\id_\liegroup}
\newcommand{\qres}{\mrm{Res}}
\newcommand{\supp}{\operatorname{supp}}
\newcommand{\Hess}{\operatorname{Hess}}
\newcommand{\Dec}{\operatorname{Dec}}
\newcommand{\Decmod}{\widetilde{\operatorname{Dec}}}
\newcommand{\driftmod}{\widetilde{\operatorname{Drift}}}
\newcommand{\fluctmod}{\widetilde{\operatorname{Fluct}}}
\newcommand{\roughdataweight}{2\kappa}
\newcommand{\roughdataweightmultiple}[1][1]{\ensuremath{\number\numexpr 2*#1\relax\kappa}}
\newcommand{\Blin}{B_{\mrm{lin}}}
\newcommand{\Brem}{B_{\mrm{rem}}}
\newcommand{\higgs}{\Phi}
\newcommand{\higgsdet}{\varphi}
\newcommand{\higgsdetmod}{\widetilde{\higgsdet}}
\newcommand{\rhiggsdet}{\higgsdet\minus}

\newcommand{\rlinhiggsdet}{\rhiggsdet_{\mrm{lin}}}
\newcommand{\rlinhiggsdetmplusone}{\rhiggsdet_{\mrm{lin},m+1}} 
\newcommand{\slinhiggsdet}{\shiggsdet_{\mrm{lin}}}

\newcommand{\invnorm}[2]{\big\|#1 \big\|_{\cov, #2}}
\newcommand{\shiggsdet}{\higgsdet^{+}}

\newcommand{\reghiggsdet}{\higgsdet_{\mrm{reg}}}
\newcommand{\remhiggsdet}{\higgsdet_{\mrm{rem}}}

\newcommand{\higgsdetcombined}{\| \higgsdet_0^+\|_{L_x^r}+\bigl\|\rhiggsdet_0\bigr\|_{\covHolder{B(0), T}^{-\kappa}}}
\newcommand{\Asmoothpart}[1][]{A^{+ #1}}
\newcommand{\Aroughpart}[1][]{A^{- #1}}
\newcommand{\Arem}[1][]{R_{#1} }
\newcommand{\Aremmod}[1][]{\widetilde{R}_{#1}}
\newcommand{\Alinmod}[1][]{\widetilde{A}_{\mrm{lin} #1}}
\newcommand{\higgsroughpart}{\higgsdet^{\scriptstyle -}}
\newcommand{\phismoothpart}{\higgs^{+}}
\newcommand{\phiroughpart}{\higgs^{-}}
\newcommand{\phirem}{\phir}
\newcommand{\phiremmod}{\widetilde{\phir}}
\newcommand{\maglapF}[1]{\operatorname{Rot}(#1)}
\newcommand{\YMHFsmooth}[1][]{B^{+ #1}}
\newcommand{\YMHFsmoothmod}[1][]{\widetilde{B}^{+ #1}}
\newcommand{\Bmod}{\widetilde{B}}
\newcommand{\YMHFrough}[1][]{B^{- #1}}

\newcommand{\Sdet}{S}
\newcommand{\Sdetmod}{\widetilde{S}}
\newcommand{\itobjects}{\Upsilon}
\newcommand{\itobjectsmod}{\widetilde{\Upsilon}}
\newcommand{\stscaling}[1]{|#1|_{\mathfrak{s}}}
\newcommand{\ptcutoff}{\sigma}
\newcommand{\bil}{\mathbf{B}}
\newcommand{\tril}{\hat{\mathbf{T}}}
\newcommand{\linTerm}{\kern-0.25em\hat{\kern0.25em\mathbf{L}}}
\newcommand{\HH}{\mathbf{H}}
\newcommand{\centre}{\mathrm{z}}

\newcommand{\gronwall}{\mathsf{Gr}}
\newcommand{\lossfactor}{\mathsf{Ls}}
\newcommand{\perturblemnorm}{\mathfrak{X}}

\newcommand{\gronwallconstant}{C_0}
\newcommand{\stepupperbound}{\delta_{\mathsf{step}}}
\newcommand{\logtimestep}{\eps_0}
\newcommand{\roughpartsize}{K}

\newcommand{\flow}{\mathcal{P}}
\newcommand{\shift}{\mathrm{sh}}
\newcommand{\wt}{\widetilde}
\newcommand{\tmod}{\widetilde{t}}
\newcommand{\tdet}{t}
\newcommand{\taumod}{\widetilde{\tau}}
\newcommand{\Lmod}{\widetilde{L}}
\newcommand{\Asmoothpartmod}{\widetilde{A}^{+}}
\newcommand{\Aroughpartmod}{\widetilde{A}^{-}}
\newcommand{\phismoothpartmod}{\widetilde{\higgs}^{+}}
\newcommand{\phiroughpartmod}{\widetilde{\higgs}^{-}}
\newcommand{\event}{\mathfrak{A}}
\newcommand{\eventmod}{\widetilde{\mathfrak{A}}}
\newcommand{\covariantenergy}{\mathsf{E}}

\newcommand{\ONBg}{u}
\newcommand{\ONBV}{v}
\newcommand{\ONBY}{\mathbf{u}}

\DeclareMathOperator{\trace}{tr}
\DeclareMathOperator{\Area}{Area}

\def\discrep{\mathsf{K}}

\DeclareMathOperator{\curl}{curl}
\let\d\partial
\let\restr\upharpoonright

\let\div\undefined
\DeclareMathOperator{\div}{div}
\def\cov{\mathrm{inv}}

\def\sol{{\mathrm{sol}}}
\def\modsol{{\mathrm{msol}}}
\def\Euc{\mrm{Euc}}
\def\bUpsilon{\boldsymbol{\Upsilon}}

\newcommand{\covHolder}[1]{\mathcal{C}_{#1}}

\def\id{\mathrm{id}} 
\let\Id\id 

\def\moll{\hat \rho_\eps}

\def\X{\mathbf{X}}
\def\narrowdots{\makebox[.9em][c]{$.$\hfil$.$\hfil$.$}}

\newcommand{\one}{\mathbf{1}}
\newcommand{\h}{\mathbf{h}}
\newcommand{\Qc}{\mathcal{Q}}

\newcommand{\bigO}{O}  
\newcommand{\schO}{\mathcal{O}} 

\newcommand{\fF}{\mathcal{F}}

\newcommand{\rR}{\mathcal{R}}

\newcommand{\CR}{\mathcal{R}}

\newcommand{\CV}{\mathcal{V}}
\newcommand{\CT}{\mathcal{T}}
\newcommand{\CI}{\mathcal{I}}
\newcommand{\CE}{\mathcal{E}}
\newcommand{\CY}{\mathcal{Y}}
\newcommand{\CU}{\mathcal{U}}
\newcommand{\CB}{\mathcal{B}}
\newcommand{\CH}{\mathcal{H}}
\newcommand{\CP}{\mathcal{P}}

\newcommand{\CJ}{\mathcal{J}}
\newcommand{\CK}{\mathcal{K}}
\newcommand{\CG}{\mathcal{G}}
\newcommand{\CA}{\mathcal{A}}

\newcommand{\LL}{\mathscr{L}}

\newcommand{\CM}{\mathcal{M}} 
\renewcommand{\models}{\mathscr{M}} 
\newcommand{\CX}{\mathcal{X}} 
\newcommand{\CS}{\mathcal{S}} 
\newcommand{\CZ}{\mathcal{Z}} 
\newcommand{\DD}{\mathscr{D}} 
\newcommand{\CD}{\mathcal{D}} 
\newcommand{\CN}{\mathcal{N}} 
\newcommand{\RR}{\mathfrak{R}} 
\def\Picard{\CD^{\gamma,\eta}} 
\def\hPicard{\hat\CD^{\gamma,\eta}} 
\def\YMH{{\mathrm{YMH}}}
\def\YMHs{{\scalebox{0.4}{$\mathrm{YMH}$}}}

\newcommand{\fL}{\mathfrak{L}}
\newcommand{\fl}{\mathfrak{l}}
\newcommand{\ft}{\mathfrak{t}}
\newcommand{\fb}{\mathfrak{b}}

\newcommand{\fx}{\mathfrak{x}}
\newcommand{\fX}{\mathfrak{X}}
\newcommand{\fy}{\mathfrak{y}}
\newcommand{\fh}{\mathfrak{h}}
\newcommand{\fu}{\mathfrak{u}}
\newcommand{\fg}{\mathfrak{g}}
\newcommand{\fp}{\mathfrak{p}}
\newcommand{\fq}{\mathfrak{q}}
\newcommand{\fr}{\mathfrak{r}}
\newcommand{\fs}{\mathfrak{s}}

\newcommand{\moralfp}{\mathfrak{r}}
\newcommand{\moralind}{\mathfrak{I}}

\newcommand{\mehler}{\mathfrak{m}}
\newcommand{\fO}{\mathfrak{O}}

\newcommand{\tauYB}{\tau_{\, \Y,B}}

\newcommand{\phiw}{\mathsf{w}}
\newcommand{\phix}{\mathsf{x}}
\newcommand{\phiy}{\mathsf{y}}
\newcommand{\phiz}{\mathsf{z}}
\newcommand{\phir}{\mathsf{r}}

\def\minus{^{-}}

\def\Lip{\mathrm{Lip}}

\def\P{\mathbb{P}}
\def\PPi{\boldsymbol{\Pi}}

\def\V{\mathbf{V}}
\def\genV{\mathbf{W}}       
\def\Y{\mathbf{Y}}
\def\boldZ{\mathbf{Z}}
\def\boldU{\mathbf{U}}
\def\Deltam{\Delta^{\!-}}

\def\SU{\textup{SU}}
\def\Gl{\textup{GL}}
\def\su{\mathfrak{su}}

\def\lin{\textup{lin}}
\def\ans{\textup{ans}}
\def\der{\textup{der}}
\def\pow{\textup{pow}}
\def\para{\textup{para}}

\def\ss{\mathrm{ss}}

\def\kernel{\mathscr{k}}

\newcommand{\F}{\mathscr{F}}

\def\${|\kern-0.18em|\kern-0.18em|}

\newcommand{\Xblack}{X^{\bullet}}
\newcommand{\Xwhite}{X^{\circ}}

\newcommand{\Lb}{\boldsymbol{\mathrm{L}\hspace{-0.1ex}\mathrm{b}}}

\DeclareFontFamily{U}{bbm}{}
\DeclareFontShape{U}{bbm}{m}{n}{<-> bbm10}{}
\newcommand{\bbx}{\mathord{\text{\usefont{U}{bbm}{m}{n}x}}}

\newcommand{\Xop}{\mathbb{X}}
\newcommand{\Xblackop}{\mathbb{X}^{\bullet}}
\newcommand{\Xwhiteop}{\mathbb{X}^{\circ}}
\newcommand{\Xblackindexop}[1]{\mathbb{X}^{\bullet, #1}}
\newcommand{\Xwhiteindexop}[1]{\mathbb{X}^{\circ, #1}}
\newcommand{\phixop}{\bbx}

\newcommand{\boldxi}{\boldsymbol{\xi}}

\def\VCas{\Cov_\V^\fg}
\def\PhiCas{\Cov_{\fg}^{\V}}

\def\Repr{\pi}
\def\repr{\pi_*}
\newcommand{\rep}{{\boldsymbol{\varrho}}}
\let\f\frac

\def\step#1{\medskip\par\noindent\textbf{#1}}
\def\proofstep#1{\medskip\par\noindent\textit{#1}}

\makeatletter

\DeclareRobustCommand{\TitleEquation}[2]{\texorpdfstring{\StrLeft{\f@series}{1}[\@firstchar]$\if%
		b\@firstchar\boldsymbol{#1}\else#1\fi$}{#2}}

\makeatother

\def\YM{\mathrm{YM}}

\begin{document}

\title{Construction of the 2D non-Abelian\\ Yang--Mills--Higgs measure via stochastic quantisation}
\author{Bjoern Bringmann$^1$\orcidlink{0000-0001-7709-623X}, Sky Cao$^2$\orcidlink{0000-0003-0345-0608}, 
Martin Hairer$^{3}$\orcidlink{0000-0002-2141-6561}, Wenhao Zhao$^{4}$\orcidlink{0000-0002-7489-1884}}

\institute{
California Institute of Technology, USA. \email{bjoern@caltech.edu}
\and
UPenn, USA. \email{skycao@wharton.upenn.edu}
\and
EPFL, Switzerland and Imperial College London, UK.
\email{martin.hairer@epfl.ch}
\and
EPFL, Switzerland. \email{wenhao.zhao@epfl.ch}
}


\linespread{1.1}

\maketitle

\begin{abstract}
We provide the first mathematical construction of the non-Abelian Yang--Mills--Higgs measure
on a two-dimensional torus by implementing Parisi and Wu's ``stochastic quantisation''
procedure. We use a specific gauge covariant (but non-Markovian) approximation
to the previously constructed associated Langevin dynamics and develop a covariant
paracontrolled calculus for it. This allows us to obtain short-time control over times that
are much longer than what the cubic nonlinearity naïvely suggests.  
When combined with the decay properties of the Yang--Mills--Higgs action under the deterministic
gradient flow, as well as a quantitative measure of how the noise pushes the dynamic away from Yang--Mills connections, 
this is good enough to yield uniform-in-time bounds. 

We finally show that
there exists a renormalisation under which our gauge covariant approximation converges
to the Langevin dynamics (modulo gauge transformation) and that the latter is topologically irreducible and strong Feller. 
It therefore admits a unique invariant measure, providing the sought-after construction.
Our construction is non-perturbative and works in particular for any fixed values of the 
coupling constants.
\vspace{1em}

\end{abstract}

\tableofcontents
\thispagestyle{empty}

\section{Introduction}\label{section:intro}

There has been significant recent interest in the construction of Euclidean quantum field theories using Parisi and Wu's stochastic quantisation procedure \cite{ParisiWu}, which is also the main focus of this article. 
At a formal level, Euclidean quantum field theories (EQFTs) can be interpreted as probability measures of the form
\begin{equation}\label{intro:eq-EQFT} 
\mathrm{d}\mu(\varphi)= \text{``}\mathcal{Z}^{-1} \exp\bigl( - 2\CE(\varphi) \bigr) \mathrm{d}\varphi\text{''}. 
\end{equation} 
Here, $\varphi$ is an element of an infinite-dimensional state space, typically a space of distributions, $\mathcal{Z}$ is a normalisation constant, $\CE$ is an energy functional, and $\mathrm{d}\varphi$ represents a Lebesgue (or volume) measure on the state space. At a rigorous level, however, one cannot make sense of the individual expressions on the right-hand side of \eqref{intro:eq-EQFT}. Indeed, the Lebesgue measure $\mathrm{d}\varphi$ only exists in finite dimensions, the energy $\CE$ is expected to be infinite on the support of the probability measure $\mu$, and the normalisation constant $\mathcal{Z}$ is not defined. However, even though the constitutive components of such expressions cannot be given a rigorous definition individually, it may still be possible to give their combination, in the form of $\mu$, a rigorous meaning. In the literature, the rigorous definition of $\mu$ is referred to as the construction of the EQFT. Often, but not always, this construction relies on proving the convergence of lattice approximations of $\mu$. In the classical constructive quantum field theory literature, several different EQFTs have been constructed. The $\Phi^4_d$-models, which are perhaps the simplest models, have been constructed in two and three dimensions in \cite{Nelson1,Nelson2} and \cite{Glimm,GJ,FO,GJBook}, respectively.

The aforementioned stochastic quantisation method of Parisi and Wu is a specific approach towards the construction of \eqref{intro:eq-EQFT} which has seen growing interest in the singular SPDE community in recent years. Rather than starting directly from \eqref{intro:eq-EQFT}, it relies on the Langevin dynamics
\begin{equation}\label{intro:eq-Langevin-general}
\partial_t \varphi = - \nabla \CE(\varphi) + \xi\;, 
\end{equation}
where both the gradient $\nabla$ and the law of the stochastic forcing term $\xi$ are uniquely determined by a chosen inner product 
on the state space (typically the $L^2$ inner product, leading to local equations if the
energy functional $\CE$ is local). The idea is to prove that 
\eqref{intro:eq-Langevin-general} is globally well-posed and admits a unique invariant measure, which is then taken as 
the \textit{definition} of the EQFT from \eqref{intro:eq-EQFT}. The reason why this is sensible is that, at 
least in finite dimensions, it is easily seen that the right-hand side of \eqref{intro:eq-EQFT} is an invariant measure 
of \eqref{intro:eq-Langevin-general}. The stochastic quantisation method has now been implemented for many different models, 
including $\Phi^4_d$-models \cite{MW17,AK2017,MW2020,GubHof,EW24}, the sine-Gordon model \cite{GHOZ25,CdLFW26,BC24sineGordon}, 
and the Abelian Higgs model \cite{BC24}. Closely related works on the global well-posedness of the stochastic 
Navier--Stokes equations \cite{HR23,HZ25} and the parabolic Anderson model \cite{CdLFW26,SZZ26} encounter similar difficulties, although these equations 
are not of the gradient form \eqref{intro:eq-Langevin-general}. 
One point that we would like to emphasise is that the original motivation for Parisi and Wu 
to introduce stochastic quantisation was precisely the hope to study \slash construct 
gauge theories without having to fix a specific gauge, thanks to the fact that the dynamical
perspective allows one to adjust the gauge ``on the fly''. This aspect of the stochastic quantisation method 
will be crucially exploited here in order to obtain global estimates, so the present work can be viewed as 
vindicating their original vision.

In this article, our primary interest lies in the two-dimensional non-Abelian Yang--Mills--Higgs model 
(as well as the pure Yang--Mills model as a special case) originally introduced 
in \cite{PhysRevLett.13.321,Higgs,GHK} in order to explain the mechanism giving mass to 
elementary particles. 
In order to state (an informal version of) our main theorem, we first formulate our setting where we use the
notation $\End(\genV)$ for the space of linear maps $\genV\to \genV$.

\begin{setting}\label{intro:assumption-setup}
Throughout this paper, we work in the following setting.
\begin{enumerate}[label=(\roman*),itemsep=0ex,topsep=0.5ex]
\item We fix a connected compact Lie group $\liegroup$ with Lie algebra $\frkg$. We write $\Ad\colon \liegroup \rightarrow \Gl(\frkg)$ for its adjoint representation defined by $\Ad_g(X)\coloneq g X g^{-1}$ for $g\in \liegroup$ and $X\in \frkg$. Its derivative $\ad\colon \frkg \rightarrow \End(\frkg)$ is given by $\ad_X(Y)=[X,Y]$, where $[\place,\place]$ denotes the Lie bracket in $\fg$.
\item We equip the Lie algebra $\frkg$ with an $\Ad$-invariant inner product $\langle\place,\place \rangle_\frkg$, i.e., an inner product satisfying $\langle \Ad_g X, \Ad_g Y\rangle_\frkg=\langle X, Y \rangle_\frkg$ for all $g\in \liegroup$ and $X,Y\in \frkg$. 
\item We let $\V$ be a finite-dimensional inner product space which is equipped with an orthogonal representation $\rep\colon \liegroup \rightarrow O(\V)$. Its derivative is denoted by $\rep_\ast \colon \frkg \rightarrow \End(\V)$.
\item We fix a symmetric $\mass \in \End(\V)$ which is $\liegroup$-invariant in the sense that
$\scal{\rep(g)\Phi,\mass\rep(g)\Phi}_\V = \scal{\Phi,\mass\Phi}_\V$ for all $g \in \liegroup$ and $\Phi \in \V$.
\end{enumerate}
\end{setting}

\begin{remark}
The $\Ad$-invariance of the inner product on $\frkg$ is equivalent to the fact that it extends to a bi-invariant (i.e.\ invariant under both left and right multiplication) Riemannian metric on $\liegroup$, which always exists
on a compact Lie group. In the case when $\frkg$ is semisimple, its Killing form is non-degenerate and we 
obtain an $\Ad$-invariant inner product by $\langle X,Y \rangle_\frkg=-\Kil(X,Y)=-\Tr(\ad_X \circ \ad_Y)$. 
When $\frkg$ is simple, this is the only possible choice up to scalar multiples.
\end{remark}

In order to formulate Yang--Mills (YM) theory in the above setting, one furthermore needs to fix a principal
$\liegroup$-bundle $P$. In our case, we will always consider the simplest choice with two-dimensional base manifold, namely
$P = \liegroup \times \T^2$, where the base manifold is a torus $\T^2 \coloneq (\R/ 2\pi \Z)^2$.
The main object of interest in YM theory is a $\liegroup$-equivariant connection on the principal bundle $P$
\dash the so-called \textit{gauge field} \dash
which in our simplified setting can always be identified with a $\frkg$-valued one-form 
$A$ on $\T^2$. In coordinates, we write it as $A(x)=A_j(x)\,\mathrm{d}x^j$ with $j \in\{1,2\}$ and $A_j \colon \T^2 \rightarrow \frkg$. 
In Yang--Mills--Higgs theory, one
furthermore considers a Higgs field $\Phi$, which in our simplified setting can be viewed 
as a $\V$-valued map $\Phi\colon \T^2 \rightarrow \V$.

The Yang--Mills--Higgs (YMH) energy, which depends on the gauge and Higgs fields, is defined as
\begin{equation}\label{intro:eq-YMH-energy}
\CE_{\YMH}(A,\Phi) = \int_{\T^2} \Bigl(\frac12 |F_A|_\frkg^2 + \frac12 |\covd_A \Phi|_\V^2 + \frac14 |\Phi|_\V^4 + \frac12 \scal{\Phi,\mass\Phi}_\V \Bigr)\, \mathrm{d}x\;.
\end{equation}
Here, $F_A$ and $\covd_A$ are the curvature tensor and covariant derivative induced by $A$, respectively. We postpone their definitions in our setting, as well as a discussion of the gauge symmetries of $\CE_{\YMH}$, until after the theorem statement.

\begin{remark}
The special case $\V=\{0\}$ covers pure Yang--Mills theory.
\end{remark}

\begin{remark}
One could in principle consider a torus of arbitrary size $L$ and add coupling constants in \eqref{intro:eq-YMH-energy}
and the definition \eqref{e:defDA} of $\covd_A$. It turns out, however, that this can always be reduced
to a scalar multiple of \eqref{intro:eq-YMH-energy} by rescaling the fields as well as the norms $|\place|_\fg$ and
$|\place|_\V$.
\end{remark}

The Langevin equation associated to \eqref{intro:eq-YMH-energy} is then formally given by the system of SPDEs
\begin{equation}\label{intro:eq-YMH-Langevin-generic}
\partial_t A = - \nabla_A \CE_{\YMH}(A,\higgs) + \xi, \qquad \partial_t \higgs = - \nabla_\higgs \CE_{\YMH}(A,\higgs) + \zeta, 
\end{equation}
where $\xi$ and $\zeta$ are independent space-time white noises on the torus, taking values in $\frkg^2$ and $\V$, respectively. We refer to the solutions to \eqref{intro:eq-YMH-Langevin-generic} as the stochastic YMH flow. Over both two- and three-dimensional tori, its\footnote{To be precise, \cite{CCHS22,CCHS22+} considered a gauge transformation of \eqref{intro:eq-YMH-Langevin-generic}, which includes an additional DeTurck term \cite{DeTurck,D85}. Furthermore, \cite{CCHS22} considered only the Yang--Mills model, but its arguments carry over to the Yang--Mills--Higgs model.} local well-posedness modulo time-dependent gauge transformations
was recently established in the seminal works \cite{CCHS22,CCHS22+}. However, these articles did not rule out the existence of finite-time blowup, which was left as an important open problem. The goal of this article is to rule out finite-time blowup, at least in two dimensions and, more ambitiously, to show the existence of a unique invariant measure.

\begin{theorem}[Informal version]
\label{thm:main_informal}
In Setting~\ref{intro:assumption-setup} with base manifold $\T^2$, the stochastic Yang--Mills--Higgs 
flow is globally well-posed as a Markov process on the space of gauge orbits, obeys uniform-in-time 
estimates and admits a unique invariant probability measure.
\end{theorem}

\begin{remark}
The invariant measure from Theorem~\ref{thm:main_informal} is a natural candidate for the 
non-Abelian YMH measure, which had not 
been constructed previously. To the best of our knowledge, Theorem~\ref{thm:main_informal} yields 
the first example of a 
classical Euclidean quantum field theory that has only been constructed via stochastic 
quantisation, i.e.\ using stochastic PDE methods. 
\end{remark}

\begin{remark}\label{rem:GaugeOrbits}
The space of gauge orbits $\fO_{\kappa}^\YMHs$ underlying Theorem~\ref{thm:main_informal} 
is analogous to the space $\fO_{1-\kappa}$ 
appearing in \cite[Def.~3.37]{CCHS22}, with the addition of the Higgs field in $\cC^{-\kappa}_x$. It is, however, \textit{not} given
by $\fO_{1-\kappa} \times \cC^{-\kappa}_x$ since the gauge group acts on both fields. Instead, just like one has
$\fO_{\alpha} = \Omega_\alpha^1 / \CG_\alpha$ with $\Omega_\alpha^1$ a suitable space of distributional one-forms and 
$\CG_\alpha$ the space of little $\alpha$-Hölder gauge transformations, one has $\fO_{\kappa}^\YMHs = (\Omega_{1-\kappa}^1 \times \cC^{-\kappa}_x)/ \CG_{1-\kappa}$, see Section~\ref{sec:stateSpace} below.

In particular, it admits Wilson loop observables as jointly continuous functions 
$\CC^2(S^1,\T^2) \times \fO_{\kappa} \to \liegroup/\Ad(\liegroup)$
given by $(A,\gamma) \mapsto [\hol(A,\gamma)]$. It also admits continuous Wilson line observables, but their endpoints need to
be smoothed out; see Corollary~\ref{cor:WilsonLines} below. Regarding the possible values for $\kappa$,
we will always assume that $\kappa$ is sufficiently small for our arguments, which leads to a finite number of upper bounds.
\end{remark}

In order to define covariant derivatives, we let $\genV$ be a general finite-dimensional vector space endowed with an orthogonal representation $\Repr$ of $\liegroup$ and let $\repr$ be its derivative at the identity. For us, the most important cases are when the pair $(\genV,\Repr)$ is given by $(\frkg,\Ad)$ or $(\V,\rep)$. Any connection one-form $A$ then defines a covariant derivative $\covd_A$ on $\genV$-valued functions
by
\begin{equ}[e:defDA]
\covd_{A,j} \coloneq \partial_j + \repr(A_j).
\end{equ}
The curvature tensor $F_A$ captures the lack of commutativity of the covariant derivatives, being 
defined as the unique $\frkg$-valued two-form satisfying, for all tuples $(\genV,\Repr)$, 
\begin{equation}\label{intro:eq-curvature-commutator-identity}
[\covd_{A,j}, \covd_{A,k}] = \repr (F_{A,jk}). 
\end{equation}
A direct computation shows that it is explicitly given by 
\begin{equation*}
F_{A,jk}= \partial_j A_k - \partial_k A_j + [A_j, A_k]. 
\end{equation*}
We write $\CG$\label{intro:def-CG} for the group of smooth automorphisms of our principal bundle $P$ that preserve the fibres. Each element of $\CG$ can be identified with a smooth function $g\colon \T^2 \rightarrow \liegroup$ which acts on the fibres of $P$ by pointwise multiplication. For $\alpha \in (\frac{1}{2},1]$, we write $\CG_\alpha$ for the closure
of $\CG$ in the space of $\alpha$-Hölder continuous functions. The group $\CG_\alpha$ then acts on 
gauge fields $A$ and Higgs fields $\higgs$ by 
\begin{equation}\label{intro:eq-gauge-transform}
g \act A = \Ad_g A - \mathrm{d}g\,g^{-1} 
\qquad \text{and} \qquad g \act \higgs = \rep(g)\higgs. 
\end{equation}
Here, the assumption $\alpha>\frac{1}{2}$ is crucial to make sure that the product $\mathrm{d}g\,g^{-1} $ 
is well-defined thanks to the fact that multiplication then extends to a map
$\cC^\alpha_x \times \cC^{\alpha-1}_x \to \cC^{\alpha-1}_x$. As for the Higgs field, the gauge transformation of a $\genV$-valued map $w$ is defined as $g \act w = \Repr(g)w$. 
We say that two connection one-forms $A$ and $B$ are gauge-equivalent if there exists an element $g\in \CG_{1-\kappa}$ such that 
$B=g\act A$. A similar definition can also be made for two Higgs fields $\higgs$ and $\Psi$, as well as for pairs $(A,\higgs)$ 
and $(B,\Psi)$. 

\begin{remark}
Geometrically, these definitions make sense if we think of all functions as being defined on the whole principal 
bundle $P$, in such a way
that they are covariant under the action of $\liegroup$ on each fibre. Then elements of $\CG$ represent simply coordinate 
changes on $P$ that preserve the fibres.
\end{remark}

\begin{remark}\label{intro:rem-genuine-one-form}
While we refer to $A$ as a connection one-form, it is important to keep in mind that gauge transformations act differently on $A$ than on a genuine, $\frkg$-valued one-form $\omega$, which transforms like $g\act \omega = \Ad_g \omega$. Thus, only the difference between two connection one-forms should be treated like a genuine one-form. 
\end{remark}

The gauge covariance of the covariant derivatives $\covd_A$ is encoded by the fact that, for any sufficiently regular 
$\genV$-valued map $w$ on the torus, it holds that $\Repr(g) \covd_A w = \covd_{g\act A}(\Repr(g)w)$, which we also write as $g\act \covd_A w= \covd_{g\act A}(g\act w)$. The curvature tensor $F_A$ transforms under gauge transformations as $g\act F_A \coloneq \Ad_g F_A = F_{g\act A}$. Thus, unlike the connection one-form $A$ itself, the curvature tensor $F_A$ transforms like a genuine differential form. An important observation (in fact the whole \textit{raison d'être} of gauge theories) is that, since both $\Ad$ and $\rep$ are orthogonal representations, the above readily implies the gauge invariance of the YMH energy,
namely $\CE_{\YMH}(g\act A,g\act \Phi)=\CE_{\YMH}(A,\Phi)$. Since $\CE_{\YMH}$ is constant along the (non-compact) infinite-dimensional gauge orbits, 
it is in particular non-coercive. The consequences of this will be felt throughout the article and will be discussed further below. 

Since we discussed the covariant derivatives, curvature tensor, and gauge transformations, we can now turn our attention to the dynamic. While our main result concerns the Langevin dynamics of the YMH energy \eqref{intro:eq-YMH-energy}, the Yang--Mills gradient flow also plays a crucial role in our proof. It is the gradient flow associated with the half squared $L^2$-norm of the curvature tensor, i.e., the first summand in \eqref{intro:eq-YMH-energy}, and it is given by solutions of
\begin{equ}
\label{e:YMHF}
\ptl_t B = -\covd_B^\ast F_B\;.
\end{equ}
The local and global well-posedness of \eqref{e:YMHF} in dimension $d=2$ (as well as $d=3$) was first obtained in \cite{Rad92}. A different proof of local well-posedness, which is closer to the arguments used in this article, was obtained in \cite{Oh2014}. The Langevin equation of the Yang--Mills--Higgs energy from \eqref{intro:eq-YMH-energy} is then (at least formally) given by\footnote{Unlike in \cite[Eq.~1.4]{CCHS22+}, the covariant derivative $\covd_A \higgs$ appears in the second argument of $\bil$, and moreover the $\bil$ term appears with a positive sign. Since $\bil$ is skew-symmetric, this results in the same equation as in \cite[Eq.~1.4]{CCHS22+}.}
\begin{equation}\label{intro:eq-YMH-Langevin}
\begin{aligned}
\partial_t A &= - \covd_A^\ast F_A + \bil( \higgs \otimes \covd_A \higgs) + \xi, \\ 
\partial_t \higgs &= \Delta_A \higgs - \mass\higgs - |\higgs|^2 \higgs + \zeta. 
\end{aligned}
\end{equation}
Here, $\xi$ and $\zeta$ are independent space-time white noises taking values in $\frkg^2$ and $\V$, respectively. Furthermore,  $\Delta_A = \delta^{k\ell} \covd_{A,k} \covd_{A,\ell} = \covd_A^\ell  \covd_{A,\ell}$ is the covariant Laplacian and $\bil \colon \V \otimes \V \rightarrow \frkg$ is the linear map uniquely determined by the relation
\begin{equ}\label{intro:eq-bil}
\scal{\bil(u\otimes v),h}_\fg = \scal{u,\rep_*(h)v}_\V\;,
\end{equ}
see \cite[Eq.~1.5]{CCHS22+}. Since the representation $\rep$ of $\liegroup$ on $\V$ is orthogonal, it follows that $\rep_*(h)$ is skew-symmetric and therefore $\bil(u\otimes v) = -\bil(v \otimes u)$. For a formulation (of a smooth approximation) of \eqref{intro:eq-YMH-Langevin} in coordinates, we also refer the reader to \eqref{local:eq-SYM-coordinates-eps} below.

We recall from above that the Yang--Mills--Higgs energy $\CE_{\YMH}$ is constant along gauge orbits. As a result, the $L^2$-inner product between its gradient $\nabla \CE_{\YMH}(A,\higgs)$ and directions in the tangent space of the gauge orbit vanishes, which implies that \eqref{intro:eq-YMH-Langevin} cannot possibly be parabolic since elliptic operators on compact manifolds have compact resolvent and therefore a finite-dimensional kernel. Indeed, a direct calculation shows that the principal part of $-(\covd_A^\ast F_A)_k$ is given by $\Delta A_k - \partial_k (\partial^\ell A_\ell)$. Due to the low regularity of the space-time white noises $\xi$ and $\zeta$, this lack of parabolicity makes a direct analysis of \eqref{intro:eq-YMH-Langevin} infeasible. However, one needs to keep in mind that we are not interested in the solution $(A,\higgs)$ of \eqref{intro:eq-YMH-Langevin} itself, but only its equivalence class under gauge transformations, i.e., its gauge orbit. By using a time-dependent gauge transformation, we can\footnote{For an explanation of the behaviour of \eqref{intro:eq-YMH-Langevin} under time-dependent gauge transformations, we refer the reader to Section~\ref{section:identities-YM-gradient-flow}.
Note that the term $\covd_A^\ast(A-B)$ appearing in both equations could have been replaced by any
$\fg$-valued process, but this specific choice turns out to be very useful.} then convert \eqref{intro:eq-YMH-Langevin} into the stochastic Yang--Mills--Higgs equations
\begin{equs}[intro:eq-SYMH]
\partial_t A &= - \covd_A^\ast F_A - \covd_A \covd_A^\ast ( A-B) + \bil( \higgs \otimes \covd_A \higgs) + \xi, \\ 
\partial_t \higgs &= \Delta_A \higgs + (\covd_A^\ast(A-B)) \higgs - \mass\higgs - |\higgs|^2 \higgs + \zeta. 
\end{equs}
(This is all still formal at this stage.)
Here, $B$ denotes the solution to the Yang--Mills gradient flow \eqref{e:YMHF} with the same
initial data as $A$. This is a variant of the Donaldson--DeTurck trick \cite{DeTurck,D85}
which, following \cite{OT22} where a similar gauge was used, we call the ``infinitesimal DeTurck gauge''.

\begin{remark}\label{rem:DeTurck}
We will show that, when renormalising \eqref{intro:eq-SYMH} in a natural way, the process one obtains is indeed gauge equivalent to the one considered in \cite{CCHS22,CCHS22+} (see Proposition \ref{main:prop-gauge-equivalence}). In particular, when quotienting by the action
of the gauge group, we do obtain a representation for the Markov process constructed in \cite{CCHS22} for any choice of smooth $B$, and the DeTurck gauge used in \cite{CCHS22, BC23} can be interpreted as the special choice $B=0$ (but see the discussion before Theorem~2.9 in \cite{CCHS22}). 
\end{remark}

\begin{remark}
We do not prove the global well-posedness for \eqref{intro:eq-SYMH} with $B$ solving \eqref{e:YMHF}, nor do we expect it to hold. For this reason, we only work with \eqref{intro:eq-SYMH} for short times and then restart the equation, updating the initial data $A_0$ and hence modifying the choice of $B$. 
\end{remark}

\subsection{Related works}

At the beginning of this introduction, we briefly mentioned works related to Theorem~\ref{thm:main_informal}, but we now discuss the literature in more detail. In the following, we focus our discussion on the quantum Yang--Mills--Higgs (or Yang--Mills) model. The literature concerning the classical YM(H) model, such as work on geometric elliptic equations \cite{U82,TT04}, geometric evolution equations \cite{Rad92,Struwe,KM95,KT99,Oh2015,OT19,Wal19}, and topology \cite{AB83,D83,H87}, will not be discussed here. The reader interested only in the main ideas of this article, and not in its relationship to the earlier literature, is encouraged to skip ahead to 
Section~\ref{section:main-ideas-brief}.

\step{Stochastic quantisation and global bounds on singular SPDEs:} 
As already mentioned, this article makes progress on the rigorous implementation of Parisi 
and Wu's stochastic quantisation method. This method is best understood in the context of the 
$\Phi^4_d$-model in two or three dimensions \cite{MW17,MW17plane,MW2020,GubHof} (see also 
\cite{JLM,MR1113223,DPDPhi4} for earlier attempts that do, however, exploit a priori knowledge of the
$\Phi^4_2$ measure). For a real-valued unknown $\varphi\colon [0,\infty) \times \T^d \rightarrow \R$, the 
parabolic evolution associated to this model can formally be written as 
\begin{equation}\label{intro:eq-Phi4d}
(\partial_t - \Delta) \varphi = - \varphi^3 + \infty\cdot \varphi + \xi.  
\end{equation}
(Here, given a regularisation procedure depending on a small parameter $\eps$, 
`$\infty$' should be interpreted as an $\eps$-dependent constant that diverges as
$\eps\to0$ in such a way that the sequence of solutions admits a non-trivial limit.)
The proof of global well-posedness and uniform-in-time estimates for \eqref{intro:eq-Phi4d} relies heavily on the strong damping term $-\varphi^3$. To see the effect of this damping term, consider the deterministic model
\begin{equation}\label{intro:eq-deterministic-Phi4d} 
(\partial_t - \Delta) \psi = - \psi^3. 
\end{equation}
By comparison with the ODE $y^\prime = - y^3$ and the maximum principle, 
the solution to \eqref{intro:eq-deterministic-Phi4d} satisfies 
$\| \psi (t) \|_{L^\infty} \leq (\| \psi(0) \|_{L^\infty}^{-2} + 2 t)^{-1/2}$ so that, no matter the size of the initial data, the $L^\infty$-norm of $\psi(t)$ is bounded by $(2t)^{-1/2}$. 
This property is known as ``coming down from infinity'' and has proven itself instrumental in the analysis of \eqref{intro:eq-Phi4d}. In fact, by treating the $\Phi^4_d$-model as a perturbation of \eqref{intro:eq-deterministic-Phi4d}, similar estimates (depending also on a suitable measure of the ``size'' of $\xi$) can be established for solutions of \eqref{intro:eq-Phi4d} (see e.g.\ \cite[Thm~2.1 and Lem.~2.7]{MW2020}).

More recently, there has been growing interest in proving global well-posedness and uniform-in-time estimates for singular SPDEs without such strong damping terms. One model of particular interest is the stochastic Navier--Stokes equation in two dimensions, which was studied in \cite{HR23,HZ25}. For a divergence-free vector-valued unknown $u\colon [0,\infty) \times \T^2 \rightarrow \R^2$, it is given by 
\begin{equation}\label{intro:eq-SNS}
(\partial_t - \Delta) u + \mathbf{P} \operatorname{div} ( u\otimes u) = \eta. 
\end{equation}
Here, $\mathbf{P}$ is the Leray projection and $\eta$ is a divergence-free stochastic forcing term. In \cite{HR23,HZ25}, $\eta$ is more general than the Leray projection of space-time white noise, but has the same regularity. Unlike in the $\Phi^4_d$-model, the nonlinearity in \eqref{intro:eq-SNS} does not have a strong damping effect since $\langle u, \mathbf{P} \operatorname{div} ( u\otimes u)  \rangle_{L_x^2}=0$ and hence, at least at the level of the $L_x^2$-norm, it does not lead to any decrease in the size of the solution. In \cite{HR23,HZ25}, the global well-posedness of \eqref{intro:eq-SNS} was shown. To be more precise, the solution of \eqref{intro:eq-SNS} is first decomposed into a large, low-frequency component and a small, high-frequency component. The low- and high-frequency components are then estimated using a monotonicity formula for the energy and the stochastic local theory, respectively. The most delicate parts of the argument concern interactions between the low- and high-frequency components, which witness both the large size of the low-frequency component and the roughness of the high-frequency component. In the later article \cite{HZ25}, the third and fourth authors not only proved the global well-posedness of \eqref{intro:eq-SNS}, but also uniform-in-time estimates for the solution. Compared to \cite{HR23}, the main ingredient was a new trick involving It\^{o}'s formula, which led to more refined estimates of the low-frequency component. 

In \cite{BC24}, the first two authors studied the massless two-dimensional stochastic Abelian Higgs model. In this model, one considers the Abelian group\footnote{In this case, the Lie algebra $\mathfrak{u}(1)$ is often identified with $\R$.} $\liegroup = U(1)$, the vector space $\V=\C$, which is equipped with the canonical representation $\rep(e^{i\theta})z=e^{i\theta}z$, and the mass $\mass=0$. After restricting to the Coulomb slice $\partial_j A^j=0$, the Langevin equation associated to \eqref{intro:eq-YMH-energy} can then be written as 
\begin{equation}\label{intro:eq-AH}
\partial_t A = \Delta A - \mathbf{P} \Im \bigl( \overline{\higgs} \covd_A \higgs \bigr) + \mathbf{P} \xi, \qquad
\partial_t \higgs = \Delta_A \higgs - |\higgs|^2 \higgs+\zeta,
\end{equation}
where $\mathbf{P}$ denotes again the Leray projection. While the $\higgs$-equation in \eqref{intro:eq-AH} contains the strong damping term $-|\higgs|^2 \higgs$, no such term is contained in the $A$-equation. In particular, since $A$ enters into the $\higgs$-equation via the term $\Delta_A \higgs$, and then enters back into its own equation via the term $\Im \bigl( \overline{\higgs} \covd_A \higgs \bigr)$, it is quite difficult to control. The main novelty in \cite{BC24} lies in the so-called covariant stochastic objects. For a suitable time $T>0$ and background connection one-form $B\colon (0,T]\times \T^2 \rightarrow \R^2$, the corresponding covariant stochastic object $\philinear[B]$ is defined as the solution of
\begin{equation}\label{intro:eq-AH-covariant-stochastic}
(\partial_t - \Delta_B ) \philinear[B] = \zeta. 
\end{equation}
The estimates of the covariant stochastic object itself, as well as the estimates of the derivative nonlinearity $\Im \bigl( \overline{\philinear[B]} \covd_B \philinear[B] \bigr)$, form one of the main parts of \cite{BC24}. They are obtained using covariant heat kernel estimates which, in turn, were proven using covariant monotonicity formulas. As will be explained in Section~\ref{section:main_results} and the body of this article, our argument will significantly expand on the ideas in \cite{BC24} in many different ways. 

Lastly, we mention that the rigorous construction of the EQFTs from \eqref{intro:eq-EQFT} has also been tackled using the Bou\'{e}--Dupuis formula \cite{BD98}. This approach was pioneered by Barashkov and Gubinelli in \cite{VarPhi} and, since then, has had several further applications \cite{B22,GOTW22,OST22,OOT24,BG25,TW26}. While it is not based on the Langevin equation \eqref{intro:eq-Langevin-general}, its implementation often uses many of the same tools as the stochastic quantisation method. Very recently, the variational approach has also been used to construct the two-dimensional Abelian Higgs measure in the Coulomb gauge~\cite{BCCKM26}.

\step{The pure Yang--Mills model using exact solvability:} There is a rich literature on the pure YM model on two-dimensional Riemannian manifolds $\mathcal{M}$, much of which focuses on the YM holonomy process. To describe it, let $A$ be a $\frkg$-valued connection one-form on $\mathcal{M}$ and $\gamma \colon [0,1]\rightarrow \mathcal{M}$ be a path, both satisfying suitable regularity conditions (see e.g.\ Section~\ref{sec:stateSpace}). Then, one can define the holonomy $\hol(A,\gamma)\eqdef h(1)\in \liegroup$, where $h\colon [0,1]\rightarrow \liegroup$ is the solution to the ODE
\begin{equation}\label{intro:eq-holonomy-ODE}
dh(t) h^{-1}(t) = - A(d\gamma(t)), \qquad h(0)=\groupid. 
\end{equation}
The holonomy $\hol(A,\gamma)$ is gauge covariant, i.e., it holds that
\begin{equation}\label{intro:eq-holonomy-covariance}
\hol(g\act A,\gamma)=g(\gamma(1)) \hol(A,\gamma) g(\gamma(0))^{-1}. 
\end{equation}
To obtain \eqref{intro:eq-holonomy-covariance}, one only needs to verify that if $h(t)$ is the solution of \eqref{intro:eq-holonomy-ODE}, then $g(\gamma(t))h(t)g(\gamma(0))^{-1}$ is the solution of \eqref{intro:eq-holonomy-ODE} with $A$ replaced by $g\act A$.  
From \eqref{intro:eq-holonomy-covariance}, it follows that if $x_0 \in \mathcal{M}$ is a given base point and $\gamma$ is a loop starting from and ending at $x_0$, then $\hol(g\act A,\gamma)=g(x_0) \hol(A,\gamma) g(x_0)^{-1}$. From the perspective of this article, it is therefore natural to first construct the YM measure in two dimensions, and then consider the holonomy process 
\begin{equation}\label{intro:eq-holonomy-process}
\bigl( \hol(A,\gamma) \bigr)_{\textup{loops }\gamma},
\end{equation}
which is defined modulo simultaneous conjugation by elements of the Lie group $\liegroup$. However, the perspective taken in much of the literature on the YM holonomy process, such as \cite{L03,L10,L20}, is opposite. The YM holonomy process is constructed directly as a stochastic process, without explicit reference to a continuum YM measure. This is achieved using invariance under subdivision (see e.g.\ \cite[Section 1.6]{L03}). The invariance under subdivision, which is very specific to pure YM theory in two dimensions, implies that the behaviour of the discrete YM model at large scales does not change under further refinements of the discretisation. For properties of the YM holonomy process, as well as different perspectives on two-dimensional YM theory, we refer the reader to \cite{DM,Driver,Sengupta1,Sengupta2,S97,L03,L10,L20, DL2023, DL2025, LemoineSurvey2025} and the references therein. 

More recently, there has been interest in going beyond the construction of \eqref{intro:eq-holonomy-process} as a stochastic process. The goal is then to construct a state space (consisting of gauge orbits of connection one-forms) on which the holonomies along smooth loops are well-defined (modulo the action of $\liegroup$) as continuous observables, as well as a probability measure on that state space. If the joint law of all the holonomies under this probability measure agrees with that of the YM holonomy process constructed in \cite{L03,L10}, the probability measure can be rightfully called ``the'' YM measure. For $\mathcal{M}=\T^2$, the YM measure was first constructed by Chevyrev in~\cite{C19}. An important feature of the construction in \cite{C19} is that, using a suitable gauge, the state space can be embedded into $\Cs_x^{-\kappa}(\T^2, \frkg^2)$, where $\kappa>0$. The regularity of the random connection one-forms therefore agrees with the regularity of the Gaussian free field, which is both needed for and expected from the analysis of the corresponding Langevin dynamics given in \cite{CCHS22}. The argument in \cite{C19} goes through lattice approximations of the YM measure. To obtain uniform estimates in the lattice spacing, Chevyrev first works with an axial gauge on large scales \cite[Prop.~4.15]{C19}, and then with a Coulomb gauge on small scales \cite[Thm~4.12]{C19}. The axial gauge is used here because, at the price of giving up smoothing in one spatial variable, it eliminates nonlinear interactions, which reflects the exact solvability of the two-dimensional, pure YM model. However, once the Higgs field is included (as in Theorem~\ref{thm:main_informal}), the axial gauge loses much of its utility. The reason is that while the axial gauge eliminates nonlinear interactions involving only the connection one-form $A$, it does not eliminate nonlinear interactions involving the Higgs field $\higgs$. And without smoothing in all spatial variables, it is rather difficult to control the impact of the Higgs field $\higgs$ on the connection one-form $A$. 
We also mention that the YM measure on $\T^2$ has been constructed using a singular elliptic PDE in \cite{CKM26}, and that the construction of the YM measure has recently been extended from $\T^2$ to general compact surfaces $\mathcal{M}$ in \cite{NN26}.

In \cite{IlyaHaoYM}, Chevyrev and Shen study the Langevin equation of the pure YM model on the two-dimensional torus $\T^2$, the so-called stochastic YM equation. In their main theorem, Chevyrev and Shen prove that the stochastic YM equation is globally well-posed as a Markov process on the space of gauge orbits, and that the YM measure is its unique invariant measure. Similar to \cite{C19}, the article \cite{IlyaHaoYM} relies on lattice approximations of the stochastic YM equation
which are shown in \cite[Thm~2.12]{IlyaHaoYM} to converge to the continuum dynamic previously constructed in \cite{CCHS22}. In order to identify the continuum limit, the authors prove a strong form of uniqueness for the renormalisation constant leading to gauge covariant continuum dynamic \cite[Thm~2.15]{IlyaHaoYM}, which is one of the main novelties in their article. The global well-posedness of the stochastic YM equations, as well as the invariance of the YM measure, is then proven using Bourgain's globalisation argument \cite{B94} which relies on new moment estimates for the YM measure \cite[Thm~9.1]{IlyaHaoYM}, refining earlier estimates from \cite{C19}. This approach is of course not available to us since the 
construction of the Yang--Mills--Higgs measure on $\T^2$ was an open problem before the present work, which is what motivated it in the first place.

\step{Lattice gauge theories:}
A natural approach towards constructing continuum YM theories involves taking scaling limits of lattice discretisations. Such discretisations are known as lattice gauge theories (alternatively, lattice YM theories), and were first proposed by Wilson \cite{Wilson1974}. For a more thorough introduction to these models, see the survey \cite{Ch19}. Brydges, Fröhlich and Seiler \cite{BFS1, BFS2, BFS3} constructed the two-dimensional Abelian Higgs model via this approach. Soon after, Balaban \cite{Balaban82I, Balaban82II, 
Balaban82III, Balaban83} proved ultraviolet stability for the two- and three-dimensional Abelian Higgs models. 
(In the 2D case, this was recently strengthened in \cite{AjayIlyaAbelian}.)
Later on, King \cite{King-I, King-II} built on Balaban's Abelian Higgs work by constructing the continuum theory in three dimensions. In all these cited works on the Abelian Higgs model, the main technique was a renormalisation group analysis to understand how the various theories behave under changes of length scales. Then, by inductively applying this analysis, the authors were able to control the unit length scale behaviour of theories which were originally defined on microscopic length scales. 

In higher dimensions, Gross \cite{Gross1983} established convergence of the three-dimensional $U(1)$ lattice gauge theory to a Gaussian continuum limit, and later Driver \cite{Driver1987} extended the analysis to four dimensions. 

Next, we turn to the non-Abelian setting, i.e.\ the YM model. In the series of works
\cite{Balaban84I,Balaban84II,Balaban85I,Balaban85II,Balaban85III,Balaban85IV,Balaban85V}, Balaban studied this model in three dimensions (again using renormalisation group methods). To discuss Balaban's work, we let $\T^3_\eps$ be the lattice approximation of $\T^3$ with lattice spacing $\eps>0$ and, for simplicity, assume that $\liegroup \in \{ \textup{U}(N),\textup{SU}(N),\textup{O}(N)\}$. The state space then consists of maps $U_\eps\colon E_{\eps} \rightarrow \liegroup$, with $ E_{\eps}$ the set of directed edges of $\T^3_\eps$, which satisfy $U_\eps(x,y)=U_\eps(y,x)^{-1}$ for all $(x,y)\in E_{\eps}$. The Wilson approximation of the YM energy is given by
\begin{equation*}
\CE_{\eps}(U_\eps) = \CE_{\YM,\eps}(U_\eps)  \eqdef \frac{1}{g^2 \eps} \sum_{p \in P_\eps} \Bigl(1- \Re \trace \bigl(U_\eps(\partial p)\bigr)\Bigr). 
\end{equation*}
Here, $g$ is a coupling constant, $P_\eps$ is the set of plaquettes of $\T^3_\eps$, $U_\eps(\partial p)$ is the holonomy of $U_\eps$ along the boundary of the plaquette $p$, and $\trace$ is the normalised trace. 
Balaban then iteratively constructs effective energies $\CE_\eps^{(k)}$ on $\T_{L^k \eps}^3$, where $L\geq 2$ is an integer, through a local, gauge covariant averaging procedure (cf. \cite{Balaban85I}). This yields effective lattice gauge theories at coarse scales and, in
the main theorem \cite[Thm~1]{Balaban85IV} of the series of papers \cite{Balaban84I}--\cite{Balaban85V}, Balaban proves uniform estimates on these effective theories as $\eps \to 0$, a property known as ultraviolet stability.
Note, though, that even in the case of 2D pure YM where a good continuum theory is known to exist, 
it is not clear how the observables of the effective coarse theory relate to those of the continuum theory.
Due to this, one cannot directly compare the results in \cite{Balaban85IV} with the results in \cite{IlyaHaoYM} or in this article. 
In the later works \cite{Balaban87,Balaban88I,Balaban88II,Balaban89I,Balaban89II}, Balaban also made progress\footnote{OpenAI recently announced a rigorous construction of the continuum limit of a four-dimensional Yang--Mills theory, building in part on Balaban's work \cite{OpenAI26}. At the time of writing, the proof has not yet been made publicly available.} towards proving ultraviolet stability in four dimensions.
It is believable that a better understanding of \cite{Balaban84I}--\cite{Balaban85V} may lead to further 
progress, also in the stochastic quantisation approach to YM theory, but this remains to be seen.

In more recent years, there has been a renewed interest in the mathematical analysis of lattice gauge theory. In \cite{Chatterjee2016}, Chatterjee computed the leading term for the free energy of the three- and four-dimensional pure $U(N)$ lattice gauge theories in an appropriate scaling limit, and in particular showed that this term is given by the free energy of a corresponding Gaussian theory, which itself was explicitly calculated by Brennecke \cite{Brennecke2025}. In another vein, Chatterjee \cite{Chatterjee2026} showed that under a certain scaling (where the inverse coupling is taken to be much larger than the lattice mesh size), a simplified version of the $SU(2)$ lattice YMH theory (where the Higgs field takes values in the unit sphere) converges to a continuum Gaussian theory for all $d \ge 2$. This was soon generalised to arbitrary compact connected matrix Lie groups by Rajasekaran, Yakir, and Zhou \cite{RYZ2026}. Chandra and Chevyrev \cite{AjayIlyaAbelian} revisited the two-dimensional Abelian Higgs model and obtained estimates on the gauge field marginal which are uniform in the lattice spacing.

We also mention various works which primarily focus on the long-distance behaviour of lattice gauge theories. The primary questions here are whether there is ``quark confinement'' and a ``mass gap''; see \cite[Sections~4, 5]{Ch19} for an introduction to these properties. Generally speaking, a theory that possesses these properties looks very disordered at long distances, while a theory without these properties exhibits a form of long-range order. Osterwalder and Seiler \cite{OsterwwalderSeiler78} showed that lattice gauge theories with arbitrary gauge group $\liegroup$ have a mass gap in the ``strongly coupled'' phase, and moreover $SU(N)$ lattice gauge theories are confining in the same phase. Göpfert and Mack \cite{GopfertMack82} showed confinement and a mass gap for the three-dimensional pure $U(1)$ lattice gauge theory, while Guth \cite{Guth80} and separately Fröhlich and Spencer \cite{FrohlichSpencer82} showed that the four-dimensional $U(1)$ theory has a deconfining phase (which is also massless). It is a (long-standing) open problem to show that in four dimensions, non-Abelian lattice gauge theories are confining and have a mass gap at all values of the coupling constant. (Actually, one needs the gauge group to have a non-trivial centre in order for the corresponding theory to be confining, as discussed in \cite{GlimmJaffe1979, Seiler1982, Chatterjee2021, Nissim2026}.)

One approach towards this problem is to relate lattice gauge theories to certain ``string theories'', i.e.\ surface models, via a ``gauge\slash string duality''. See e.g.\ Polyakov's retrospective \cite{Polyakov2008} on this subject from a physics point of view. In mathematics, Chatterjee \cite{Chatterjee2019a} established a gauge\slash string duality for the large-$N$ limit of $SO(N)$ lattice gauge theories, which was further analysed by Basu and Ganguly \cite{BasuGanguly2018} (primarily in the case of two dimensions). The gauge\slash string duality for $SU(N)$ was later established by Jafarov \cite{Jafarov2016}, while Chatterjee and Jafarov \cite{ChatterjeeJafarov2016} gave a large-$N$ expansion in terms of strings. Later, Cao, Park, and Sheffield \cite{CaoParkSheffield2025} derived a gauge\slash string duality at finite-$N$ (i.e.\ without taking the large-$N$ limit), with follow-up work by various authors \cite{BCSK2024, BCSK2025, CNS2025}. See also Lemoine \cite{Lemoine2026} for a recent generalisation of \cite{CaoParkSheffield2025}; we also mention the companion work \cite{Lemoine2026b} which considers the large-$N$ limit of lattice gauge theories equipped with the Villain (i.e.\ heat kernel) action. 

Finally, there has also been a recent line of work initiated by Shen, Smith, R. Zhu, and X. Zhu \cite{ShenZhuZhu2023, ShenSmithZhu2024, ShenSmithZhu2025, ShenZhuZhu2026}, which uses the lattice Langevin dynamics to prove results such as the Makeenko--Migdal equations and a mass gap (in certain strong coupling regimes) for lattice gauge theories. See also \cite{CNS2025b}, which proves confinement at strong coupling. 

\subsection{Brief discussion of main ideas and techniques}\label{section:main-ideas-brief}

We conclude the introduction with a brief discussion of the main ideas behind the proof of Theorem~\ref{thm:main_informal}, but refer the reader to Section~\ref{section:main_results} (as well as the 
body of this article) for more detailed discussions. For most of this brief discussion, we focus on the 
stochastic Yang--Mills (rather than Yang--Mills--Higgs) equations in the infinitesimal DeTurck gauge 
(see \eqref{intro:eq-SYMH}), which are (formally) given by 
\begin{equation}\label{intro:eq-pure-SYM}
\partial_t A = - \covd_A^\ast F_A - \covd_A \covd_A^\ast (A-B) + \xi. 
\end{equation}

\step{Infinitesimal DeTurck gauge:} 
The idea behind this gauge is to use the DeTurck trick (which would correspond to taking $B=0$ in 
\eqref{intro:eq-pure-SYM}) only at the level of the linearised dynamic, not at the level of the full Langevin dynamics. To make this more precise, take $B$ a solution to the
YM gradient flow and consider the difference $C\coloneq A-B$. Using Lemma~\ref{prelim:lem-linearization} below, it can then be shown that $C$ solves
\begin{equation}\label{intro:eq-linearized}
\begin{aligned}
( \partial_t - \Delta_B ) C^k - 2 \big[ F_B^{kj},C_j \big] &=  \big[ C_j , 2 \covd_B^j C^k - \covd_B^k C^j \big]
+ \big[ C_j, \big[ C^j , C^k \big]\big] + \xi^k. 
\end{aligned}
\end{equation}
We emphasise that the connection one-form $B$ only enters into \eqref{intro:eq-linearized} through $\covd_B$ and $F_B$, but never by itself, which will be important in our covariant local theory. We also remind the reader that our background $B$ solves the Yang--Mills gradient flow without any DeTurck term \eqref{e:YMHF}. This is important because the Yang--Mills heat flow (with DeTurck term) does not obey suitable long-time estimates. The reader may think of $B$ as the smooth (or low-frequency) part of $A$ and of $C$ as the rough (or high-frequency) part of $A$. By using the infinitesimal DeTurck gauge, rather than the original DeTurck gauge, we therefore use the DeTurck trick only for the rough part of $A$.

\step{Covariant heat kernel estimates:} Motivated by \eqref{intro:eq-linearized}, we examine the covariant heat kernel $\fq_B(s,y;t,x)$ corresponding to the covariant heat operator $\partial_t - \Delta_B - 2 \ad(F_B)$
appearing on the left-hand side of \eqref{intro:eq-linearized}. We emphasise that, while the curvature term $2\ad(F_B)$ is only a zeroth-order differential operator, it has a significant impact on the behaviour of $\fq_B(s,y;t,x)$. In particular, the curvature term can be responsible for exponential growth. 
Our main estimates of $\fq_B(s,y;t,x)$ are pointwise and, for all $0\leq s < t \leq T$ and $x,y\in \T^2$, show that 
\begin{equ}\label{intro:eq-kernel-pointwise}
\big|\covd_B^{j_1} \hdots \covd_B^{j_k} \fq_B(s,y;t,x) \big| 
\lesssim \exp\bigl( C \| F_B \|_{L_t^1L_x^\infty([0,T])} \bigr)  \frac{(1+\Fnorm{F_B}_{T,k})^k}{(t-s)^{1+\frac{k}{2}}} \exp\Bigl( - c \tfrac{|x-y|^2}{t-s} \Bigr). 
\end{equ}
Here, $k\geq 1$, $1\leq j_1,\hdots,j_k \leq 2$, and $\Fnorm{F_B}_{T,k}$ is as in Definition~\ref{def:Fnorm}. These pointwise estimates form the foundation of our analytical machinery and are used heavily in almost all later parts of our argument. The proof of \eqref{intro:eq-kernel-pointwise} relies on covariant monotonicity formulas, as well as delicate estimates of curvature terms. 

\step{Gauge-covariant smooth approximations:} In singular SPDEs, local and global well-posedness are studied using limits of smooth approximations. In the earlier works \cite{CCHS22,CCHS22+,BC23}, the smooth approximations broke the gauge invariance of the law of the solution, which could only be restored in the limit. In this work, we introduce a gauge covariant smooth approximation of \eqref{intro:eq-pure-SYM}, which is defined as 
\begin{equation}\label{intro:eq-gauge covariant-smooth-approximation}
\partial_t A_\eps = - \covd_{A_{\eps}}^\ast F_{A_{\eps}} - \covd_{A_{\eps}} \covd_{A_{\eps}}^\ast (A_\eps-B) + \renormA ( A_\eps - B ) +\xi_{\Beps}. 
\end{equation}
Here, $\renormA (A_\eps-B)$ is a renormalisation term and the smooth approximation of the space-time white noise is defined as
\begin{equ}[e:mollifNoise]
\xi_{\Beps}(t,x) =   \int_0^t \ds \int_{\T^2} \dy \, \rho_\eps^\prime(t-s) \fq_B(s,y;t,x) \xi(s,y), 
\end{equ}
where $\rho_\eps\colon [0,\infty) \rightarrow [0,1]$ is a cutoff function as in \eqref{main:eq-mcP}--\eqref{main:eq-rho-eps} below. The gauge covariance of \eqref{intro:eq-gauge covariant-smooth-approximation}, however, comes at a price. Unlike for the non-covariant smooth approximations from \cite{CCHS22,CCHS22+,BC23}, it is not clear that the resulting stochastic process (at the level of the gauge orbits) is Markov. To establish the Markov property, we will show the gauge equivalence of the limiting solutions from this article and \cite{CCHS22}, see Proposition~\ref{main:prop-gauge-equivalence}. 

\step{Covariant local theory:} We develop a covariant local theory for \eqref{intro:eq-gauge covariant-smooth-approximation}. The resulting local time of existence only depends on the size of the curvature $F_B$ (rather than the connection one-form $B$ itself), which is gauge invariant. Since monotonicity formulas for \eqref{e:YMHF} only control the size of $F_B$, this will be crucial in our global arguments. 
As part of our covariant local theory, we examine the covariant stochastic objects 
\begin{equation}\label{intro:eq-linearB}
(\partial_t - \Delta_B - 2 \ad(F_B)) \linear[B] = \xi, \qquad \linear[B] (0)=0. 
\end{equation}
Similar covariant stochastic objects were previously used in \cite{BC24} for the Abelian Higgs model, but our non-Abelian setting prevents us from using many of the estimates in \cite{BC24}. In addition, we introduce a covariant paracontrolled calculus, which comes with covariant commutator estimates, Littlewood--Paley decompositions, paraproduct estimates, and Schauder estimates.

\step{Global iteration argument:} 
In order to prove the global well-posedness in Theorem~\ref{thm:main_informal}, we iterate over well-chosen times $(t_m)_{m=0}^\infty$. The increments $\tau_{m+1} \eqdef t_{m+1}-t_m$ will ultimately be chosen such that our covariant local theory can be used for \eqref{intro:eq-pure-SYM} on each interval $[t_m,t_{m+1}]$. For the global theory, the most important term in our solution $A$ is a background $B$. On each time interval, this background $B$ is a Yang--Mills gradient flow. To bound its curvature at the time $t=t_{m+1}$, we use a geometric monotonicity formula, which implies that\footnote{In our global iteration scheme (see Section \ref{section:global-theory}), the background connection one-form $B$ will be a gauge transformation of~$B^++B^-$, where $B^+$ and $B^-$ are smooth and rough parts. The monotonicity formula \eqref{intro:eq-monotonicity-curvature} will then not be used for $B$ itself, but for $B^+$. For the purposes of this introduction, we believe it is best to ignore this.}
\begin{equation}\label{intro:eq-monotonicity-curvature}
\| F_B(t_{m+1}) \|_{L_x^2}^2 - \| F_B(t_m) \|_{L_x^2}^2 
= - 2 \| \covd_B^\ast F_B \|_{L_t^2 L_x^2([t_m,t_{m+1}])}^2. 
\end{equation}
Due to the space-time white noise $\xi$, our solution of \eqref{intro:eq-pure-SYM} also contains several rough terms. On a single time interval, the rough terms can be controlled using our covariant local theory. However, due to the exponential growth in \eqref{intro:eq-kernel-pointwise}, the rough terms have to enter into the background $B$ after finitely many steps of our iteration. In the next step of the iteration, a covariant mollification of the covariant stochastic object $\linear[B]$ is absorbed into the background $B$, which increases the squared  $L^2$-norm of its curvature by roughly $\tau_{m+1}^{-1}$. Our global iteration argument, which is the most delicate part of this article, then carefully balances the increase by $\tau_{m+1}^{-1}$ with the decrease coming from \eqref{intro:eq-monotonicity-curvature}.

\step{Decay mechanism near Yang--Mills connections:}
If the background $B$ is close to a Yang--Mills connection, i.e., if $\covd_B^\ast F_B$ is small, then \eqref{intro:eq-monotonicity-curvature} does not provide us with any appreciable decay. In order to obtain the uniform-in-time estimates of Theorem~\ref{thm:main_informal}, we therefore need to uncover a decay mechanism near Yang--Mills connections. It is important to realise that Yang--Mills connections with large curvature are rather unstable (see e.g.\ Section~\ref{section:perturbing-around-YM-connection}). While the unstable directions lead to growth of the linearised dynamic, this growth should lead to a decrease in energy (see Figure~\ref{fig:Morse_Langevin_saddle} as well as the identity \eqref{main:eq-hessian-toy} below). To be more precise, if $B$ is close to a Yang--Mills connection, we expect a decrease of the form
\begin{equation}\label{intro:eq-decrease-near-YM-connection}
(\nabla^2 \CE_\YM)(v,v)\;,\qquad v =  e^{(t_{m+1} - t_{m+1/2}) (\Delta_{B}+2\ad(F_B))} \linear[B](t_{m+1/2})\;.
\end{equation}
Here, $\CE_\YM$ is the Yang--Mills energy, $t_{m+1/2}$ is a suitable time in $(t_{m},t_{m+1})$, and the operator
acting on $\linear[B](t_{m+1/2})$ in \eqref{intro:eq-decrease-near-YM-connection} is the covariant heat flow induced by the time-dependent connection one-form~$B$. To make this precise, we construct a parametrix for the covariant heat kernel around Yang--Mills connections (see Proposition~\ref{main:prop-kernel-YM-connection}). 

The mechanism here is strongly reminiscent of the situation arising when one considers the
projective dynamic for a linear stochastic PDE, see \cite{Tommaso,Lyapunov}. There, each eigenvector 
of the linear operator yields a saddle point for the deterministic dynamic and one also relies
on the presence of the noise to show that the dynamic get ``unstuck''.

An even simpler model, illustrated in Figure \ref{fig:Morse_Langevin_saddle}, is the stochastic gradient flow with respect to a Morse function on a surface. In the figure, the Morse function $f$ attains its maximum on the right and its minimum on the left. 
Pairs of deterministic gradient flow trajectories form loops representing different homotopy classes.  
The deterministic trajectories can be stuck at the saddle points drawn in the figure, whose presence is due to the high genus of the surface. Adding noise allows the trajectories to escape these saddle points and eventually reach a neighbourhood of the minimum on the left. The 2D stochastic Yang--Mills flow can be viewed as an infinite-dimensional analogue of this model, and the Yang--Mills energy can be viewed as a Morse function on the space of connections, with infinitely many saddle points. 

\begin{figure}[htbp]
\centering
\includegraphics[width=\textwidth]{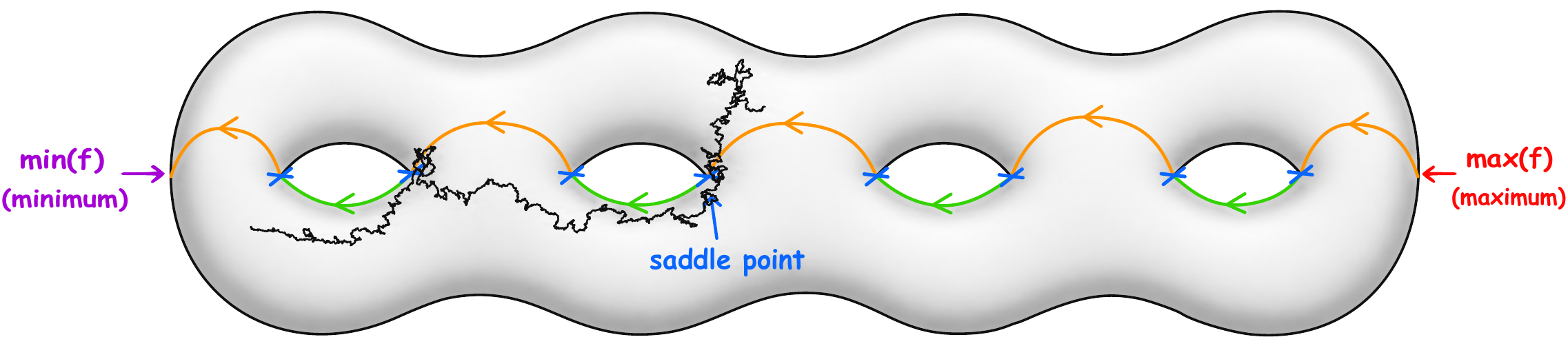}
\caption{Illustration of Langevin dynamics on a surface equipped with a Morse function.}
\label{fig:Morse_Langevin_saddle}
\end{figure}

\step{Control of the Higgs field:} While we have so far focused our discussion on the stochastic Yang--Mills 
equations \eqref{intro:eq-pure-SYM}, including the Higgs field requires at least one new ingredient. The
 counterpart of the background connection one-form $B$ is a field $\higgsdet$, which is defined 
 using a deterministic evolution equation, see \eqref{eq:Sdet-equation}--\eqref{eq:higgsdet-equation}. In 
 our global iteration argument, it is then tempting to replace the Yang--Mills energy of $B$ with the 
Yang--Mills--Higgs energy of $(B,\higgsdet)$, but it is better to only use norms of the Higgs field at regularity zero. For this reason, we work with the quantity
\begin{equation}\label{intro:eq-higgs-quantity}
\| F_{B} \|_{L_x^2}^2 + \| \higgsdet \|_{L_x^4}^{\gamma_1} + \| \higgsdet \|_{L_x^r}^{\gamma_2}. 
\end{equation}
Here, $\gamma_1$ is slightly larger than $4$, $\gamma_2$ is equal to $3$, and $r$ is chosen very large. 
Due to geometric monotonicity formulas, the $L_x^4$-term in \eqref{intro:eq-higgs-quantity} is useful 
in our estimates of the $\bil(\higgsdet \otimes \covd_B \higgsdet)$-term appearing in the $A$ equation (recall \eqref{intro:eq-YMH-Langevin}). In contrast, the 
$L_x^r$-term in  \eqref{intro:eq-higgs-quantity} is useful in closing our estimates for the Higgs field itself. 

\step{Remnants of exact solvability:}
We previously emphasised that, unlike the two-dimensional pure Yang--Mills model, the Yang--Mills--Higgs model is not exactly solvable. Nevertheless, since we take as our reference connection in \eqref{intro:eq-SYMH} 
the solution to the (pure) Yang--Mills gradient flow,
some remnants of the exact solvability of the pure Yang--Mills model can be seen in parts of our argument. In two dimensions, the four components of the curvature tensor $F_A$ satisfy $F_A^{11}=F_A^{22}=0$ and $F_A^{21}=-F_A^{12}$, so that $[F_A^{ij},F_A^{k\ell}]=0$. Furthermore, the components of the divergence $-\covd_A^\ast F_A$ are given by $(-\covd_A^\ast F_A)^1 = -\covd_{A,2} F_A^{12}$ and $(-\covd_A^\ast F_A)^2 = \covd_{A,1} F_A^{12}$. Thus, unlike in higher dimensions, the divergence $-\covd_A^\ast F_A$ controls the full covariant derivative of the curvature tensor. 

For any solution $B$ of the Yang--Mills gradient flow \eqref{e:YMHF}, it can then be shown that the component $F_B^{12}$ 
satisfies the covariant heat equation $(\partial_t - \Delta_B)F_B^{12}=0$. Strictly speaking, it is still a nonlinear evolution 
equation, since $B$ and $F_B^{12}$ are related. However, it still shares many features of the linear covariant heat equation. In 
particular, $F_B^{12}(t)$ can be written as an integral of a covariant heat kernel against the initial data $F_B^{12}(0)$. In 
addition, for any Yang--Mills connection $B$, it can be shown that the curvature tensor $F_B$ is covariantly constant. 
These two facts will be used heavily in our covariant heat kernel estimates and decay mechanism.

\step{Going for the kill:}
Once we have global-in-time curvature bounds thanks to our control of 
\eqref{intro:eq-higgs-quantity}, the remainder of the argument
is by now quite standard. There are, however, a few subtleties that need verifying. First, we need to make 
sure that the process constructed as the limit $\eps \to 0$ of \eqref{intro:eq-gauge covariant-smooth-approximation}
coincides with the process constructed in \cite{CCHS22}, modulo gauge transformation (and with the Higgs field
added on both sides). 
Since the approximation used there is
quite different, this is actually far from obvious! However, the problems are virtually
identical to those arising when trying to show gauge covariance in \cite{CCHS22}. (This is not too
surprising since the process \eqref{intro:eq-gauge covariant-smooth-approximation} is already
gauge covariant for fixed $\eps > 0$, so this would be another way of proving gauge covariance.)
The main difference is that the presence of the Higgs field as well as the somewhat non-standard
regularisation of the noise in \eqref{intro:eq-gauge covariant-smooth-approximation} create quite a 
few more counterterms, so that the computation of the counterterms is more involved. 

Once we know that we are looking at the correct process, it remains to show that it admits a unique
invariant probability measure. Existence is immediate from the usual tightness argument, while
uniqueness follows if we can show that the process is strong Feller and topologically irreducible.
The strong Feller property follows as in \cite{IlyaHaoYM}. Regarding irreducibility, it unfortunately appears
quite difficult to apply the results of \cite{SupportSimple,SupportGeneral}. Instead, we develop a new
(and simpler) argument which is more widely applicable and should therefore be of independent interest.

\subsection{AI disclosure statement}

The main mathematical ideas contained in this article as well as the vast majority of
its write-up were generated without the assistance 
of AI. The authors started working on this project in April~2025 and arrived at the main ideas by 
February~2026, which predates the most recent AI-powered breakthroughs.
These ideas and results were presented 
by the authors at various conferences and seminar talks in early 2026, including at MPI 
Leipzig (March 10), Stanford University (April 6), and a Simons Collaboration Meeting 
(April 30). 

While AI assistance was not used to generate new ideas, AI tools were used for other research tasks. 
Different versions of ChatGPT Pro and Claude were used to tweak the exposition, search the literature, and 
proofread parts of the article. During the final stages of this project, GPT-6 Astra and GPT-5.6 Sol were also used to prove a small number of technical lemmas and propositions, namely, Proposition~\ref{prelim:prop-X-supermartingale}, Lemma~\ref{prelim:lem-Z-supermartingale-bound}, Proposition~\ref{prelim:prop-stopping-time}, Lemma~\ref{le:B_Sdet_heat_comparison}, Lemma~\ref{le:Decmod_rare_event}, Lemma~\ref{le:taumod_rare_event} and Lemma~\ref{lemma:inverse-time-tail-bound}.  
GPT-6 Astra was also used by WZ to assist with typing some of the formulas in Section~\ref{section:decay-estimates-in-expectation}. They were generated one at a time and carefully reviewed before being incorporated into the manuscript. 

In addition, during the development of the main ideas in 2025, conversations (for lack of a better word) with ChatGPT about existing tools in differential geometry and geometric analysis were extremely helpful.

\subsection*{Acknowledgements}

{\small BB's research was partially supported by the National Science Foundation under Grant No. DMS-2453126; and he thanks Princeton University for providing a stimulating environment during most stages of this work. SC's research was partially supported by the National Science Foundation under Grant No. DMS-2303165; he would also like to thank MIT for providing excellent working conditions under which much of this work was completed. MH's research was partially funded by the Simons Collaboration ``Probabilistic Paths to Quantum Field Theory'', as well as by the Swiss National Science Foundation, grant number 10006872. We thank the Courant Institute of Mathematical Sciences, EPFL, the Institute for Mathematical Sciences at NUS, the Hausdorff Research Institute for Mathematics, Peking University, Princeton University, the Simons Foundation, and SLMath for hosting various subsets of the authors during the completion of this work.}

\section{Main results}\label{section:main_results}

In the previous section, we stated our main theorem (Theorem \ref{thm:main_informal}) and briefly discussed our main ideas, but this was all done at an informal level. In this section, we will be detailed and rigorous. In particular, we introduce rigorous definitions and state several precise theorems capturing the main results from different parts of this article. In Section \ref{section:main-heat}, we discuss our covariant heat kernels and their estimates. In Section \ref{section:main-local}, we state our covariant local well-posedness theorem. In Section \ref{section:main-global}, we discuss our global iteration scheme and uniform-in-time estimates. Finally, in Section \ref{section:Markov-and-Measure}, we discuss the Markov process on the space of gauge orbits and its unique invariant measure, namely, the Yang--Mills--Higgs measure. 

\subsection{Covariant heat kernel estimates}\label{section:main-heat}

Throughout most of this article (except when comparing the process constructed here to the one 
constructed in \cite{CCHS22,CCHS22+}, which is required to show that it has the Markov property), we adopt 
the following central principle: we always work covariantly with respect to a (large) background connection, 
because we want to treat the background non-perturbatively. This requires us to develop an entire covariant Schauder theory which constitutes the analytic backbone of this article and is the subject of Section~\ref{section:covariant-heat-kernel-estimates}. In particular, we obtain estimates on covariant heat kernels, which we heavily use in the ensuing covariant local theory (see Section~\ref{section:main-local} for more discussion). We now begin to describe these estimates. Unless otherwise specified, the following discussion applies to a general dimension $d \geq 2$.

Throughout, we fix a time $T \in (0, 1]$ and a smooth solution $B \colon [0, T] \times \T^d \ra \frkg^d$ to the Yang--Mills gradient flow. Given this, we may consider the following covariant heat equations:
\begin{equs}
(\ptl_t - \Delta_B) \phi &= G, \label{eq:intro-covariant-heat-equation} \\
(\ptl_t - \Delta_B) \psi_i &= 2 (F_B)_i^{~\ell} \psi_\ell + H_i, ~~ i \in [d]. \label{eq:intro-covariant-heat-equation-with-curvature}
\end{equs}
In the first equation, $\phi, G$ are $\genV$-valued, while in the second equation, $\psi, H$ are 
$\genV^d$-valued.\footnote{Recall that $\genV$ is a finite-dimensional vector space endowed with an orthogonal representation of $\liegroup$.} We will often write the second equation in the component-free form
\begin{equs}
(\ptl_t - \Delta_B) \psi = 2 F_B \psi + H,
\end{equs}
which is in line with the interpretation of $F$ as an operator on $\genV^d$, which acts by $(F X)_i = F_{i}^{~j} X_j$ for given $X = (X_i, i \in [d]) \in \genV^d$.

\begin{definition}[Covariant heat kernels]\label{kernel:def-kernels}
Let $\fp_B(s, y; t, x; \genV)$ be the fundamental solution to \eqref{eq:intro-covariant-heat-equation} and $\fq_B(s, y; t, x; \genV^d)$ be the fundamental solution to \eqref{eq:intro-covariant-heat-equation-with-curvature}, i.e., for $\phi_s \in C^\infty(\T^d, \genV)$ we have that
\begin{equs}\label{eq:intro-pB-characterization}
\phi(t, x) \coloneq \int \dy \, \fp_B(s, y; t, x; \genV) \phi_s(y)
\end{equs}
is the solution to \eqref{eq:intro-covariant-heat-equation} on $[s,T]$ with $G = 0$ and initial data $\phi_s$, and for $\psi_s \in C^\infty(\T^d, \genV^d)$, 
\begin{equs}\label{eq:intro-qB-characterization}
\psi(t, x) \coloneq \int \dy \, \fq_B(s, y; t, x; \genV^d) \psi_s(y) 
\end{equs}
is the solution to \eqref{eq:intro-covariant-heat-equation-with-curvature} with $H = 0$ and initial data $\psi_s$. Note that, for any points $(s, y), (t, x)\in [0,T]\times \T^d$ satisfying $s<t$, $\fp_B(s, y; t, x; \genV) \in \End(\genV)$, while $\fq_B(s, y; t, x; \genV^d) \in \End(\genV^{d})$. We will also often view $\fp_B \in \End(\genV^{d})$ as acting componentwise (i.e.\ in the 
natural block-diagonal fashion). We refer to $\fp_B, \fq_B$ as covariant heat kernels. For notational simplicity, we will typically drop the dependence on $\genV, \genV^d$ and simply write $\fp_B(s, y; t, x)$, $\fq_B(s, y; t, x)$.
For notational convenience, we also write $\fp_B(t, x) = \fp_B(0, 0; t, x)$, $\fq_B(t, x) = \fq_B(0, 0; t, x)$.
\end{definition}

\begin{remark}\label{rem:geometric}
The correct geometric setting consists of a vector bundle $E \ra M$ (where $M = \T^d$) whose fibres are modelled on $\genV$, from which we also obtain the tensor bundle $T^* M \otimes E$. Then one should really view $\fp_B(s, y; t, x) \in L(E_y, E_x)$ and $\fq_B(s, y; t, x) \in L((T_y^* M) \otimes E_y, (T_x^* M) \otimes E_x)$, where $L(X, Y)$ denotes the space of linear maps from a vector space $X$ to a vector space $Y$. This is consistent with the characterisations of $\fp_B$ and $\fq_B$ given by \eqref{eq:intro-pB-characterization} and \eqref{eq:intro-qB-characterization}, where we view $\phi \in \Gamma(E)$ and $\psi \in \Gamma(T^* M \otimes E)$ as sections of the appropriate bundles, so that $\phi(t, x) \in E_x$ and $\psi(t, x) \in T_x^*M \otimes E_x$.

The more explicit viewpoints of $\fp_B \in \End(\genV)$ and $\fq_B \in \End(\genV^d)$ follow if 
$E$ is trivial (noting that the tangent bundle $TM = T \T^d$ admits a global trivialisation), so that there is a global identification of the fibres $E_x \cong \genV$, $T_x^* M \otimes E_x \cong \genV^d$.
\end{remark}

In order to obtain covariant Schauder estimates, we need to estimate the solutions $\phi, \psi$ to \eqref{eq:intro-covariant-heat-equation} and \eqref{eq:intro-covariant-heat-equation-with-curvature}. This will require us to assume control on certain quantities involving the curvature $F_B$, which we next define.

\begin{definition}[Curvature norm]\label{def:Fnorm}
Let $I = [s_0, t_0] \sse [0, T]$ be an interval and let $k \geq 0$. Define the following norm:
\begin{equs}
\Fnorm{F_B}_{I, k} \coloneq \sum_{j=0}^k  \Bigl(\big\|(t - s_0)^{(j+1)/2} \covd_B^{(j)} F_B\big\|_{L_t^2 L_x^\infty(I)} + \big\|(t-s_0)^{j/2} \covd_B^{(j)} F_B \big\|_{L_t^1 L_x^\infty(I)}\Bigr).
\end{equs}
In particular, when $k = 0$, we have that
\begin{equs}
\Fnorm{F_B}_{I, 0} = \big\|(t - s_0)^{1/2}F_B\big\|_{L_t^2 L_x^\infty(I)} + \big\|F_B\big\|_{L_t^1 L_x^\infty(I)}.
\end{equs}
If $I = [0, t]$ is an interval starting at $0$, we will often write $\Fnorm{F_B}_{t, k}$ as a shorthand for $\Fnorm{F_B}_{[0, t], k}$.
\end{definition}

We now present one of the main results of Section~\ref{section:covariant-heat-kernel-estimates} (more precisely, this is a combination of two results, Propositions~\ref{prop:pB-ptwise-bounds} and~\ref{prop:qB-ptwise-bds}). This result says that the covariant heat kernels $\fp_B, \fq_B$ (and their higher derivatives) satisfy bounds similar to those of the usual heat kernel, as long as $\Fnorm{F_B}$ is controlled. In the following, recall our convention that $\fp_B(t, x) = \fp_B(0, 0; t, x)$ and similarly for $\fq_B(t, x)$.

\begin{theorem}[Pointwise estimates of covariant heat kernels]\label{main:thm-pointwise-kernel}
There is a constant $c > 0$ such that the following holds. For any $(t_0, x_0) \in (0, T] \times \T^d$ and $k \geq 0$, we have the pointwise bounds:
\begin{equ}
\bigl| \covd_B^{(k)} \fp_B(t_0, x_0)\bigr| \lesssim_k \bigl(1 + \bigl(\Fnorm{F_B}_{t_0, (k-1)\vee 0}\bigr)^{k} \bigr) t_0^{-\frac{k+d}{2}} \exp(-c |x_0|^2 / t_0).
\end{equ}
Additionally, there is a constant $C_0 > 0$ such that the following holds. For any $(t_0, x_0) \in (0, T] \times \T^d$ and $k \geq 0$, we have the pointwise estimates on the difference between $\fq_B$ and $\fp_B$:
\begin{equs}
\bigl|&\covd_B^{(k)} \bigl(\fq_B - \fp_B\bigr)(t_0, x_0) \bigr| \\
&\lesssim_k \exp\bigl(C_0 \bigl\|F_B\bigr\|_{L_t^1 L_x^\infty([0, t_0])} \bigr) \Fnorm{F_B}_{t_0, k} \bigl(1 + (\Fnorm{F_B}_{t_0, k})^{(k-1)\vee 0}\bigr) t_0^{-\frac{k+d}2} \exp(-c |x_0|^2/t_0).
\end{equs}
\end{theorem}

\begin{remark}\label{remark:main-thm-ptwise-bds}
We make the following remarks regarding Theorem~\ref{main:thm-pointwise-kernel}.
\begin{enumerate}[label=(\alph*)]
\item For simplicity, we state the theorem with the specific base point $(0, 0)$, but the result directly extends to provide estimates on $\fp_B(s_0, y_0; t_0, x_0)$, $\fq_B(s_0, y_0; t_0, x_0)$ by translation.
\item In \cite[Section 3]{BC24}, various covariant heat kernel estimates (in the case $\liegroup = U(1)$) were also proven. Therein, the estimates were mostly stated in terms of operator norms, whereas here we obtain more explicit pointwise estimates. These pointwise estimates make the ensuing analysis of the covariant stochastic objects much simpler than in \cite{BC24}.
\end{enumerate}
\end{remark}

Next, we discuss some ideas behind the proof of Theorem~\ref{main:thm-pointwise-kernel}. In particular, we discuss how to estimate a solution $\phi$ to \eqref{eq:intro-covariant-heat-equation}. In one sentence, the main technique is to use monotonicity formulas. 
We begin by computing the evolution of $|\phi|^2$:
\begin{equ}
\frac{1}{2}(\ptl_t - \Delta) \bigl(|\phi|^2\bigr) + \bigl|\covd_B \phi\bigr|^2  = \langle \phi, G \rangle_{\genV}.
\end{equ}
Now let $K(t, x)$ be the backward heat kernel  centred at the point $(t_0, x_0)$ (i.e.\ so that $(\ptl_t + \Delta)K = 0$ for $t < t_0$ and $K(t_0,\cdot) = \delta_{x_0}$). Upon multiplying the above equation by $K$, we may obtain
\begin{equ}
\frac{1}{2}\ptl_t(K |\phi|^2) + K \bigl|\covd_B \phi\bigr|^2  = \frac{1}{2} K \Delta (|\phi|^2) - \frac{1}{2} (\Delta K) |\phi|^2 +  K \langle \phi, G \rangle_{\genV}.
\end{equ}
Thus upon integrating in space-time (and applying integration by parts in space to cancel the first two terms on the right-hand side above) we obtain
\begin{equation}\label{eq:intro-monotonicity-integrated}
\begin{aligned}
\bigl\|& (K^{1/2} \phi)(t_0)\bigr\|_{L_x^2}^2 + 2 \bigl\|K^{1/2} \covd_B \phi\bigr\|_{L_{t,x}^2([0, t_0])}^2 \\ 
&= \bigl\|(K^{1/2} \phi)(0)\bigr\|_{L_x^2}^2 + 2 \int_0^{t_0} \dt \int_{\T^d} \dx \, \bigl(K \langle \phi, G \rangle_{\genV}\bigr)(t, x).
\end{aligned}
\end{equation}
In particular, if $G = 0$ (so that $\phi$ exactly solves a covariant heat equation), then we obtain
\begin{equ}
\bigl\|(K^{1/2} \phi)(t_0)\bigr\|_{L_x^2} \leq \bigl\|(K^{1/2} \phi)(0)\bigr\|_{L_x^2}.
\end{equ}
The heat kernel localisation $K$ allows us to upgrade the above $L_x^2$ estimate into an $L_x^\infty$ estimate for $\phi$. In particular, recalling that $K(t_0, \cdot) = \delta_{x_0}$, we have the identity $\bigl\| (K^{1/2} \phi)(t_0)\bigr\|_{L_x^2} = |\phi(t_0, x_0)|$ (of course, this is slightly formal since $K(t_0,\cdot)$ is 
only a distribution, but we ignore this for the purposes of the discussion). Thus, we may for instance estimate
\begin{equ}
|\phi(t_0, x_0)| \leq \bigl\|(K^{1/2} \phi)(0)\bigr\|_{L_x^2} \leq \bigl\|K^{1/2}(0)\bigr\|_{L_x^\infty} \bigl\|\phi(0)\bigr\|_{L_x^2} \lesssim t_0^{-\frac{d}{4}} \bigl\|\phi(0)\bigr\|_{L_x^2}.
\end{equ}
By varying $x_0$, we may thus obtain the $L_x^\infty$ estimate
\begin{equ}
\bigl\|\phi(t_0)\bigr\|_{L_x^\infty} \lesssim t_0^{-\frac{d}{4}} \bigl\|\phi(0)\bigr\|_{L_x^2}.
\end{equ}
Note that this is the same estimate one would obtain if $\phi(t) = e^{t \Delta} \phi(0)$ were a solution to the usual heat equation. 

\begin{remark}
The heat kernel localisation $K$ was also used in \cite{BC24} to obtain operator bounds on covariant heat kernels.
\end{remark}

For general $G$, we can estimate the additional term in \eqref{eq:intro-monotonicity-integrated} by
\begin{equs}
\biggl| \int_0^{t_0} \dt \int_{\T^d} \dx \, \bigl(K \langle \phi, G \rangle_{\genV}\bigr)(t, x) \biggr| &\leq \bigl\|K^{1/2} \phi\bigr\|_{L_t^\infty L_x^2([0, t_0])} \bigl\|K^{1/2} G \bigr\|_{L_t^1 L_x^2([0, t_0])} \\
& \leq \delta \bigl\|K^{1/2} \phi\bigr\|_{L_t^\infty L_x^2([0, t_0])}^2 + C_\delta \bigl\|K^{1/2} G \bigr\|_{L_t^1 L_x^2([0, t_0])}^2,
\end{equs}
for any $\delta > 0$ (and $C_\delta$ a constant depending on $\delta$). Combining this with \eqref{eq:intro-monotonicity-integrated}, we may obtain the further estimate
\begin{equ}\label{eq:intro-integrated-heat-kernel-monotonicity-G}
\bigl\|K^{1/2} \phi\bigr\|_{L_t^\infty L_x^2([0, t_0])} + \bigl\|K^{1/2} \covd_B \phi\bigr\|_{L_{t,x}^2([0, t_0])} \lesssim \bigl\|(K^{1/2} \phi)(0)\bigr\|_{L_x^2} + \bigl\|K^{1/2} G \bigr\|_{L_t^1 L_x^2([0, t_0])}.
\end{equ}
A suitably time-weighted version of this estimate will allow us to control higher derivatives of $\phi$, because we may compute
\begin{equ}
(\ptl_t - \Delta_B) (\covd_{B, i} \phi) = 2 (F_B)_i^{~\ell} \covd_{B, \ell} \phi + \covd_{B, i} G, \quad i \in [d].
\end{equ}
In particular, $\covd_B \phi$ itself satisfies a forced heat equation (in fact, the equation is precisely of the form \eqref{eq:intro-covariant-heat-equation-with-curvature}). Applying the estimate \eqref{eq:intro-integrated-heat-kernel-monotonicity-G} with $\phi$ replaced
by $\covd_B \phi$ (and $G$ taken to be $2 F_B \covd_B \phi + \covd_B G$), we see that we will be able to obtain an $L_t^\infty L_x^2$ estimate of $K^{1/2} \covd_B \phi$, provided that we can estimate the $L_t^1 L_x^2$ norm of $K^{1/2} F_B \covd_B \phi$ (and we assume suitable control on the forcing term $\covd_B G$). This is possible because \eqref{eq:intro-integrated-heat-kernel-monotonicity-G} also gives an $L_{t,x}^2$ estimate of $K^{1/2} \covd_B \phi$, and thus we may simply apply Hölder. Even higher derivatives of $\phi$ may be estimated in a similar inductive manner: we start by computing the equation for $\covd_B^{(k)} \phi$, which schematically looks like
\begin{equ}
(\ptl_t - \Delta_B) (\covd_B^{(k)} \phi) = \sum_{\substack{j_1 + j_2 = k \\ j_2 \geq 1}} \schO((\covd_B^{(j_1)} F_B) \covd_B^{(j_2)} \phi) + \covd_B^{(k)} G. 
\end{equ}
All terms involving $\phi$ on the right-hand side will have already been estimated in previous steps of the induction, and therefore we may close the inductive step.

Besides the pointwise estimates of Theorem~\ref{main:thm-pointwise-kernel}, Section~\ref{section:covariant-heat-kernel-estimates} collects various other Schauder estimates that are needed later on. Additionally, the same section contains two different short-time expansions for the covariant heat kernels, which we next discuss. The first is a heat kernel parametrix, which approximates $\fp_B$ by $\fp \pt_B$, where $\fp$ is the usual heat kernel and $\pt_B \in \End(\genV)$ is the parallel transport associated with $B$ (which we discuss in more detail in Section~\ref{section:parallel-transport}). We state this result as follows.

\begin{proposition}[Covariant heat kernel parametrix]\label{main:prop-parametrix}
There is a constant $c > 0$ such that the following holds. Let $t_0 \in (0, T]$. Let
\begin{equs}
\mc{N}(B) \eqdef \Big(t_0^{1/2} \big(\|\ptl_t B\|_{L_t^1 L_x^\infty([0, t_0])} + \|F_B\|_{L_t^2 L_x^\infty([0, t_0])} \big) \Big)^2  & + t_0^{1/2}\|\ptl_t B\|_{L_t^1 L_x^\infty([0, t_0])}  \\
&+ t_0^{10} \bigl(1 + \big\|F_B\big\|_{L_t^2 L_x^\infty([0, t_0])}\bigr) .
\end{equs}
For any point $x_0$ within distance $1/10$ (say) of $0$, we have the pointwise bounds
\begin{equs}
|\fp_B - \fp  \pt_B |(t_0, x_0) &\lesssim \mc{N}(B) t_0^{-\tfrac{d}{2}} \exp(-c|x_0|^2/t_0), \\
\bigl|\covd_B (\fp_B - \fp  \pt_B)\bigr|(t_0, x_0) &\lesssim \mc{N}(B) (1 + \Fnorm{F_B}_{t_0, 0})  t_0^{-\tfrac{d+1}{2}} \exp(- c|x_0|^2/t_0). 
\end{equs}
\end{proposition}

We prove this proposition in Section~\ref{section:heat-kernel-parametrix}. In this paper, we will generally work on time scales such that the quantity $\mathcal{N}(B)$ appearing in Proposition~\ref{main:prop-parametrix} is small, so that $\fp \pt_B$ is indeed a good approximation of $\fp_B$. This parametrix plays a critical role in computing resonances of covariant stochastic objects, which is needed to arrive at the correct renormalisation counterterms in our covariant local theory (to be discussed in Section~\ref{section:main-local}). To prove Proposition~\ref{main:prop-parametrix}, we begin by deriving the evolution equation of the difference $\fp_B - \fp \pt_B$:
\begin{equ}
(\ptl_t - \Delta_B) (\fp_B - \fp \pt_B) = \text{perturbative terms},
\end{equ}
and then we proceed to estimate all perturbative terms appearing on the right-hand side. These terms involve covariant derivatives of parallel transport, which we estimate (in Section~\ref{section:parallel-transport-estimates}) using the quantity $\mathcal{N}(B)$ appearing in the proposition.

By combining Theorem~\ref{main:thm-pointwise-kernel} and Proposition~\ref{main:prop-parametrix}, we also obtain a parametrix for $\fq_B$. However, it is apparent from these estimates that once the curvature $F_B$ becomes too large or the time scale $t_0$ too long, the parametrix $\fp \pt_B$ ceases to be an accurate approximation of $\fq_B$. Our second short-time expansion enters precisely at this point \dash it becomes particularly effective when $B$ is near a Yang--Mills connection with large curvature, and in this case it moreover (crucially) works on slightly larger time scales. This alternative expansion provides an approximation of $\fq_B$ by an explicit kernel which \emph{exhibits exponential growth}. This exponential growth arises due to instabilities caused by the additional curvature term in the equation \eqref{eq:intro-covariant-heat-equation-with-curvature} and is one of the main inputs needed for our decay mechanism (which ultimately leads to uniform-in-time estimates), to be discussed in more detail in Section~\ref{section:main-global}. We now state an informal version of the second short-time expansion. We avoid stating the precise version here, as it takes a fair amount of notational buildup. See the discussion at the beginning of Section~\ref{section:perturbing-around-YM-connection} as well as Proposition~\ref{prop:qB-perturb-near-YM}.

\begin{proposition}[Informal. Covariant heat kernel near Yang--Mills connection]\label{main:prop-kernel-YM-connection}
Let $d = 2$. Suppose that $B$ is sufficiently close to a Yang--Mills connection, in the sense that $\covd_B F_B$ and higher derivatives of $F_B$ are suitably controlled. Then $\fq_B$ is suitably close to $\pt_B \fs$, where the kernel $\fs(t, x) \in \End(\genV^2)$ satisfies
\begin{equ} 
\Bigl(\ptl_t - \Delta - 2 F_B(0, 0) - \frac{1}{4} |x|^2 F_B^{12}(0, 0)^2\Bigr)  \fs = 0, \quad \fs(0) = \delta_0 \Id_{\genV^2}.
\end{equ}
Here $F_B^{12}(0, 0) \in \frkg$ and, as usual, we implicitly apply the representation $\repr$ to view $F_B^{12}(0, 0) \in \End(\genV)$. We then extend $F_B^{12}(0, 0)$ to an operator on $\genV^2$ in the componentwise fashion (this is equivalent to just directly applying $\repr \oplus \repr$). This provides the interpretation of the term $F_B^{12}(0, 0)^2 \fs$ on the left-hand side. The term $F_B(0, 0) \in \End(\genV^{ 2})$ acts by $(F_B(0, 0)X)_i = (F_B(0, 0))_{i}^{~j} X_j$ for $X \in \genV^{ 2}$, $i \in [2]$.

Additionally, $\fs(t, x) \in \End(\genV^2)$ has eigenvalues as follows. Note that $F_B^{12}(0, 0) \in \End(\genV)$ is skew-symmetric, and thus has purely imaginary eigenvalues. For each such eigenvalue $\icomplex \lambda$, $\lambda \in \R$, there is an eigenvalue of $\fs(t, x)$ given by:
\begin{equ}
e^{2 \lambda t} \frac{\lambda}{4\pi \sinh(\lambda t)} \exp\biggl(-\frac{\lambda |x|^2}{4 \tanh(\lambda t)}\biggr).
\end{equ}
(When $\lambda = 0$, the above expression should be interpreted as the limit as $\lambda \ra 0$, which is readily seen to be the usual heat kernel.) In particular, when $\lambda > 0$ and $x$ is sufficiently close to $0$, the above exhibits exponential growth with rate $e^{\lambda t}$. (By the skew-symmetry of $F_B^{12}(0, 0)$, the nonzero eigenvalues of $F_B^{12}(0, 0)$ come in pairs $\pm \icomplex \lambda$, so as long as $F_B^{12}(0, 0)$ is nonzero, there is an eigenvalue $\icomplex \lambda$ with $\lambda > 0$.)
\end{proposition}

\begin{remark}
Unlike the other results of this section, we must take $d = 2$ in Proposition~\ref{main:prop-kernel-YM-connection}. This reflects the fact that two-dimensional Yang--Mills connections exhibit special structure not seen in higher dimensions, as briefly discussed at the end of Section~\ref{section:main-ideas-brief}. As one particular consequence, when working on the whole-plane $\R^2$, if $B$ is Yang--Mills, the kernel $\pt_B^{-1} \fq_B$ may be exactly computed \dash it is precisely the (whole-plane analogue of the) kernel $\fs$ in Proposition~\ref{main:prop-kernel-YM-connection} whose eigenvalues we completely understand. Thus Proposition~\ref{main:prop-kernel-YM-connection} (and our decay mechanism which rests on this proposition) is one of the places in this article where we most heavily use the assumption that $d = 2$. 

We also remark that Avramidi \cite{Avramidi1993, Avramidi1995, Avramidi1995b} previously studied similar questions. He showed that in general dimensions, if one assumes that the curvature is covariantly constant (when $d \geq 3$ this seems to be strictly stronger than assuming that the connection is Yang--Mills), then one can obtain a similar approximation of the heat kernel as in the case $d = 2$. Put another way, it appears that if one wants a result that works in general dimensions, then one must assume more than just that the underlying connection is Yang--Mills. 
\end{remark}

We prove Proposition~\ref{main:prop-kernel-YM-connection} (or more precisely, Proposition~\ref{prop:qB-perturb-near-YM}) in Section~\ref{section:perturbing-around-YM-connection}. To briefly discuss the proof, a natural approach would be to proceed as in the parametrix construction, by deriving the following form for the evolution equation of $\fq_B - \pt_B\fs$:
\begin{equ}
(\ptl_t - \Delta_B - 2 F_B) (\fq_B - \pt_B\fs) = \text{perturbative terms}.
\end{equ}
However, the above is not exactly true, as there is still a single term on the right-hand side which is non-perturbative, because it involves the curvature $F_B$, rather than covariant derivatives of $F_B$. Indeed, by ``perturbative term'' in this case, we mean a term which is small when $B$ is close to Yang--Mills, because that is the assumption of Proposition~\ref{main:prop-kernel-YM-connection}. Therefore, perturbative terms should only involve covariant derivatives of $F_B$ and not $F_B$ itself. To handle the extra term, our only choice is to group it with the operator on the left-hand side, to obtain a new operator $(\ptl_t -  \mathcal{L}_B)$ for some suitable $\mathcal{L}_B$. Fortunately, this new operator satisfies the same Schauder estimates as the original operator $(\ptl_t - \Delta_B - 2F_B)$, which is ultimately due to the fact that the single additional non-perturbative term only contributes a phase (essentially, the new operator behaves as a magnetic Laplacian). This argument is carried out in Section~\ref{section:perturbing-around-YM-connection}. To be precise, in that section we find it more convenient to compare $\pt_B^{-1} \fq_B - \fs$ instead of $\fq_B - \pt_B \fs$, which results in the evolution equation looking slightly different than in the present proof sketch.

\subsection{Covariant local well-posedness}\label{section:main-local}

One of the major goals of this article is to establish a covariant local theory of \eqref{intro:eq-SYMH}. This covariant local theory will heavily rely on the covariant heat kernel estimates discussed in Section~\ref{section:main-heat} and will be an important ingredient in our global theory (see Section~\ref{section:main-global} below). As a first step towards our covariant local theory, we discuss the Yang--Mills gradient flow \eqref{e:YMHF} with rough initial data. This is not only needed to formulate the infinitesimal DeTurck gauge appearing in \eqref{intro:eq-SYMH}, but will later serve as a warm-up for our full covariant local theory (cf. Section~\ref{section:deterministic-rough} and Sections~\ref{section:covariant-Ansatz}--\ref{section:objects}).

\begin{proposition}[Yang--Mills gradient flow with rough initial data]\label{main:prop-YMHF-rough-data}
Let $\kappa>0$ be as in \eqref{e:small_constants}. Let $C_0 >1$ and $0< c_0 \leq 1$ be sufficiently large and small constants, respectively. 
Let $B_0^\sharp, \YMHFsmooth_0\colon \T^2 \rightarrow \frkg^2$ be connection one-forms and let $\YMHFrough_0 \colon \T^2 \rightarrow \frkg^2$ be a one-form such that 
$\YMHFsmooth_0$ is smooth, $\YMHFrough_0$ is in $\Cs_x^{-\kappa}$, and $B_0^\sharp=\YMHFsmooth_0+\YMHFrough_0$. Let $\YMHFsmooth$ be the Yang--Mills gradient flow with initial data $\YMHFsmooth(0)=\YMHFsmooth_0$. Finally, let $0<T\leq 1$ and $K\geq 1$ be such that
\begin{align}
\exp\bigl(C_0 \| F_{\YMHFsmooth} \|_{L_t^1 L_x^\infty([0,T])}\bigr) \bigl(1 + \Fnorm{F_{\YMHFsmooth}}_{T, 2}\bigr)^2 T^{\frac{1-\kappa}{2}} K &\leq c_0, \label{main:eq-bound-YMHFsmooth} \\ 
\|\YMHFrough_0\|_{\covHolder{\YMHFsmooth_0, F, T, 2}^{-\kappa}}  &\leq \roughpartsize , \label{main:eq-bound-YMHFrough}
\end{align}
where the last condition involves the covariant Hölder norm at regularity $-\kappa$ from Definition \ref{defn:covariant_Holder} below. Then, there exist a solution $B\in C_t^0 \Cs_x^{-\kappa}([0,T]\times \T^2, \frkg^2)$ of \eqref{e:YMHF} and a gauge transformation $g\in C_t^0 \CG_{1-\kappa}([0, T])$ (recall the definition of $\CG_{1-\kappa}$ just above \eqref{intro:eq-gauge-transform}), both of which are smooth on $(0,T]\times \T^2$, such that $B(0)=g(0)\act B_0^\sharp$. Furthermore, for all $t\in (0,T]$, we have the curvature estimates
\begin{equation}\label{main:eq-YMHF-curvature-estimate}
\begin{aligned}
\bigl\|F_B(t)\bigr\|_{L_x^\infty} &\leq \bigl\|F_{\YMHFsmooth}(t)\bigr\|_{L_x^\infty} + C t^{-\frac{1+\kappa}{2}} K, \\ 
\bigl\|\covd_B^\ast F_B(t)\bigr\|_{L_x^\infty} &\leq \bigl\| \covd_{\YMHFsmooth}^\ast F_{\YMHFsmooth}(t)\bigr\|_{L_x^\infty} + C t^{-\frac{\kappa}{2}} \bigl( t^{-1} + \| F_{\YMHFsmooth}(t) \|_{L_x^\infty}\bigr) K, 
\end{aligned}
\end{equation}
where $C$ is a constant depending on the parameters in \eqref{e:small_constants}. 
\end{proposition}
For more detailed statements concerning the Yang--Mills gradient flow with low-regularity initial data, we refer the reader to Proposition~\ref{pr:admissible_timescale_general} and Lemma~\ref{lemma:ymgf-low-regularity-uniqueness} below. The former result contains additional estimates of the solution, and the latter yields the uniqueness of $B$ (for given initial data $B(0)$). We emphasise that Proposition~\ref{main:prop-YMHF-rough-data} does not construct a Yang--Mills gradient flow for the initial data $B_0^\sharp$ itself, but only for the gauge-transformed initial data $g(0)\act B_0^\sharp$. The reason we are allowed to gauge-transform the initial data is that we are ultimately only interested in the Markov process on the space of gauge orbits (see Theorem~\ref{thm:main_informal}) and not in the solution of \eqref{intro:eq-SYMH} itself. In any case, the aspect of Proposition~\ref{main:prop-YMHF-rough-data} we use the most is the curvature estimate \eqref{main:eq-YMHF-curvature-estimate}, which is gauge invariant. In applications of Proposition~\ref{main:prop-YMHF-rough-data}, $K$ will not be much larger than one, and hence \eqref{main:eq-YMHF-curvature-estimate} implies that the size of $F_B(t)$ is dictated by the size of $F_{\YMHFsmooth}(t)$. In particular, the $L_t^1L_x^\infty$-norm of $F_B$, which plays a crucial role in our covariant heat kernel estimates, is often comparable to the $L_t^1L_x^\infty$-norm of $F_{\YMHFsmooth}$.  

The proof of Proposition~\ref{main:prop-YMHF-rough-data}, which is the subject of Section~\ref{sec:rough_YMHF}, is rather simple. The main idea is to first show the existence of a unique solution $B^\sharp\colon [0,T]\times \T^2 \rightarrow \frkg^2$ of 
\begin{equation}\label{main:eq-YMHF-rough-data-DeTurck}
\partial_t B^\sharp = - \covd_{B^\sharp}^\ast F_{B^\sharp} - \covd_{B^\sharp} \covd_{B^\sharp}^\ast (B^\sharp-\YMHFsmooth), \qquad \qquad 
B^\sharp(0) = B_0^\sharp. 
\end{equation}
That is, we first solve the Yang--Mills gradient flow from \eqref{e:YMHF} in the infinitesimal DeTurck gauge. Since the $- \covd_{B^\sharp} \covd_{B^\sharp}^\ast (B^\sharp-\YMHFsmooth)$-term makes \eqref{main:eq-YMHF-rough-data-DeTurck} parabolic, this can be done using our covariant heat kernel estimates from Section~\ref{section:main-heat} and a standard contraction-mapping argument. Once \eqref{main:eq-YMHF-rough-data-DeTurck} has been solved, we then construct a gauge transformation $g\colon [0,T]\times \T^2 \rightarrow \liegroup$ by solving the ODE
\begin{equation}\label{main:eq-YMHF-rough-data-gauge-transformation}
(\partial_t g) g^{-1}= - \Ad_g \bigl( \covd_{B^\sharp}^\ast \bigl( B^\sharp - \YMHFsmooth \bigr)\bigr), 
\qquad \qquad 
g(T) = \groupid. 
\end{equation}
The reason for imposing a condition on $g$ at the terminal time $t=T$, rather than at the initial time $t=0$, is that $B^\sharp$ is only smooth away from the initial time. Once \eqref{main:eq-YMHF-rough-data-gauge-transformation} has been solved, the Yang--Mills gradient flow $B$ can then be defined as $g\act B^\sharp$.  
Since the assumptions in Proposition~\ref{main:prop-YMHF-rough-data} are rather long and our covariant local theory only relies on bounds of the curvature tensor $F_B$, it is convenient to introduce a simpler class of initial data. To treat \eqref{intro:eq-SYMH}, this class also includes initial data for the Higgs field. It is defined as
\begin{equation}\label{main:eq-initial-class}
\begin{aligned}
\initialclass_T \coloneq \Bigl \{& (A_0, \higgs_0) \in \Cs_x^{-\kappa}(\T^2 , \frkg^2 \times \V ) \colon  \, \exists\,\textup{YM gradient flow }\\ 
& B \in C_t^0 \Cs_x^{-\kappa}([0,T]\times \T^2 , \frkg^2) \medcap C_{t,x}^\infty((0,T]\times \T^2 , \frkg^2) \textup{ with initial data $A_0$} \\ 
&\textup{satisfying } \| t^{\frac{1+\kappa}{2}} F_B \|_{L_{t,x}^\infty},  
\| t^{\frac{2+\kappa}{2}} \covd_B^\ast F_B \|_{L_{t,x}^\infty}<\infty, \,  
\|\higgs_0\|_{\covHolder{A_0, T}^{-\kappa}} < \infty \Bigr\}. 
\end{aligned}
\end{equation}
Here, $\|\place \|_{\covHolder{A_0, T}^{-\kappa}}$ is the gauge invariant Hölder norm from Definition \ref{defn:covariant_Holder}.

\begin{remark}\label{remark:initial-class-time}
The initial-data class $\initialclass_T$ defined in
\eqref{main:eq-initial-class} is independent of $T\in(0,1]$. Indeed, it holds that $\initialclass_T=\initialclass_1$ for all $T\in (0,1]$. This follows from the global well-posedness result of \cite{Rad92}, which implies that the
Yang--Mills gradient flow $B$ appearing in
\eqref{main:eq-initial-class} extends from $[0,T]$ to $[0,1]$, and the extension is smooth on $(0,1]$. Moreover, global
well-posedness of the linear covariant heat equation implies
that
$\|\higgs_0\|_{\covHolder{A_0,T}^{-\kappa}}<\infty$
if and only if
$\|\higgs_0\|_{\covHolder{A_0,1}^{-\kappa}}<\infty$.

We nevertheless retain $T$ in the notation to emphasize that,
in Theorem~\ref{main:thm-covariant-local} below, we only work
with the Yang--Mills gradient flow on $[0,T]$.
\end{remark}

After the above discussion of the Yang--Mills gradient flow, we now turn to our gauge covariant smooth approximation of \eqref{intro:eq-SYMH}. To this end, we first let $C_{\mc{P}}$ be a sufficiently large constant (say $C_{\mc{P}}=1000!$), which will be fixed throughout the rest of the article. Then, we introduce the set of cutoff functions 
\begin{equation}\label{main:eq-mcP}
\begin{aligned}
\mc{P} \eqdef \Big\{ \rho \in C^\infty([0,\infty)\ra [0,1])\Big|& \, \, \rho|_{[0,1]}=0, \, \rho^\prime|_{[1,2]}\geq 0, \, \rho|_{[2,\infty)}=1, \, \, \\
&\mbox{and} \, \, \| \rho^{(k)}\|_{L^\infty} \leq C_{\mc{P}} \, \, \mbox{for all}\, \, 0 \leq k \leq 20 \Big\}. 
\end{aligned}
\end{equation}
For all $\rho\in \mc{P}$ and $\eps \in \sdyadic$, we further set
\begin{equation}\label{main:eq-rho-eps}
\rho_\eps(t) \eqdef \rho\bigl(\eps^{-2} t\bigr). 
\end{equation}
From \eqref{main:eq-mcP} and \eqref{main:eq-rho-eps}, it follows that $\rho_\eps$ is a smooth cutoff function to times $t\gtrsim \eps^2$. Furthermore, it follows that its derivative $\rho_\eps^\prime$ is an $L^1$-normalised smooth cutoff function to times $t\sim \eps^2$ and hence converges to the Dirac delta at $t=0$. Using our smooth cutoff functions, we can now define our gauge covariant, smooth approximations of the space-time white noises.

\begin{definition}[\protect{Gauge-covariant approximations of $\xi$ and $\zeta$}]\label{main:def-gauge covariant-approximation-space-time-white-noise}
Let $T>0$, let $B$ be a Yang--Mills gradient flow as in \eqref{main:eq-initial-class}, let $\rho\in \mc{P}$, and let $\eps \in \sdyadic$. On the time interval $[0,T]$, we then define 
\begin{equation}\label{main:eq-xi-zeta-smooth-approximation}
\begin{aligned}
\xi_{\Beps}(t,x) &\eqdef   \int_0^t \ds \int_{\T^2} \dy \, \rho_\eps^\prime(t-s) \fq_B(s,y;t,x) \xi(s,y), \\ 
\zeta_{\Beps}(t,x) &\eqdef   \int_0^t \ds \int_{\T^2} \dy \, \rho_\eps^\prime(t-s) \fp_B(s,y;t,x) \zeta(s,y),
\end{aligned}
\end{equation}
where $\fp_B$ and $\fq_B$ are the covariant heat kernels from Definition~\ref{kernel:def-kernels}.
\end{definition}

In the earlier works \cite{CCHS22,CCHS22+,BC23}, the smooth approximations of space-time white noises were defined using standard mollifiers. In Definition~\ref{main:def-gauge covariant-approximation-space-time-white-noise}, however, we use covariant heat kernels, which makes both $\xi_{\Beps}$ and $\zeta_{\Beps}$ gauge covariant. To be more precise, let $g\in \CG$ be a smooth, spatial gauge transformation. The smoothness of $g$ is assumed for simplicity. We believe that, as long as $g$ is probabilistically independent of $\xi$ and $\zeta$, similar results should hold under the weaker assumption $g\in \CG_{1-\kappa}$.  For any $\xi \in \Cs_{t,x}^{-2-\kappa}([0,T]\times \T^2 , \frkg^2)$ and $\zeta \in \Cs_{t,x}^{-2-\kappa}([0,T]\times \T^2 , \V)$, one can show that $g\act \xi_{\Beps} = (g\act \xi)_{g\act B,\eps}$ and $g\act \zeta_{\Beps}  = (g\act \zeta)_{g\act B,\eps}$. Here, the action of $g$ on $\xi$ or $\zeta$ is defined using the adjoint representation or the representation $\rep$, respectively.

In our gauge covariant, smooth approximation of \eqref{intro:eq-SYMH}, we also need renormalisation constants
 $\renormA \in \End_{\liegroup}(\frkg)$ and $\renormphi\in \End_{\liegroup}(\V)$, where $\End_{\liegroup}$ denotes the space of $\liegroup$-covariant endomorphisms. Since $\frkg$ is equipped with the adjoint representation and $\V$ is equipped with the representation $\rep$, this means that $\renormA \Ad_g = \Ad_g \renormA$ and $\renormphi \rep(g) = \rep(g) \renormphi$ for all $g\in \liegroup$. For now, both $\renormA$ and $\renormphi$ are entirely arbitrary, but for our covariant local and global well-posedness results, $\renormA$ and $\renormphi$ will need to be chosen as in Definition~\ref{main:def-renormalization-constant} below. Equipped with the renormalisation constants, we then consider  
\begin{equation}\label{main:eq-SYMH-eps}
\begin{aligned}
\partial_t A_\eps &= - \covd_{A_{\eps}}^\ast F_{A_{\eps}} - \covd_{A_{\eps}} \covd_{A_{\eps}}^\ast (A_\eps-B) + \bil\bigl( \higgs_\eps \otimes \covd_{A_\eps}\higgs_\eps \bigr) + \renormA ( A_\eps - B ) +\xi_{\Beps}, \\ 
\partial_t \higgs_\eps &= \Delta_{A_\eps} \higgs_\eps + (\covd_{A_\eps}^\ast(A_\eps-B)) \higgs_\eps 
-  |\higgs_\eps|^2 \higgs_\eps + \renormphi (\higgs_\eps)   + \zeta_{\Beps},  \\
A_\eps(0)&=A_{0}, \qquad \higgs_\eps(0)= \higgs_0.   
\end{aligned}
\end{equation}
(We have absorbed the mass term $\mass\higgs_\eps$ into the term $\renormphi \higgs_\eps$.)
For notational reasons, it is convenient to denote the unique solution $(A_\eps,\higgs_\eps)$ of \eqref{main:eq-SYMH-eps} on the time interval $[0,T]$, provided that it exists, by 
\begin{equation}\label{main:eq-data-to-solution}
\solutionmap_{T,\rho,\eps}(A_0,\higgs_0,\xi,\zeta;\renormA,\renormphi). 
\end{equation} 
In the following, we refer to \eqref{main:eq-data-to-solution} as the data-to-solution map.

\begin{lemma}[Gauge-covariant approximations]\label{main:lem-gauge-covariance}
Let $T>0$, let $\rho \in \mc{P}$, and let $\eps\in \sdyadic$. In addition, let $\renormA\in \End_{\liegroup}(\frkg)$ and $\renormphi\in \End_{\liegroup}(\V)$, i.e., let $\renormA$ and $\renormphi$ be $\liegroup$-covariant. Furthermore, let $(A_0,\higgs_0)\in \initialclass_T$,  let $\xi \in \Cs_{t,x}^{-2-\kappa}([0,T]\times \T^2 , \frkg^2)$, and let $\zeta \in \Cs_{t,x}^{-2-\kappa}([0,T]\times \T^2 , \V)$. Finally, assume that the unique solution of \eqref{main:eq-SYMH-eps} exists on the time interval $[0,T]$. Then, for all smooth spatial gauge transformations $g\in \CG$, it holds that
\begin{equation}\label{main:eq-gauge covariance}
g \act \solutionmap_{T,\rho,\eps}(A_0,\higgs_0,\xi,\zeta;\renormA,\renormphi) = 
\solutionmap_{T,\rho,\eps}(g \act A_0,g \act\higgs_0,g \act\xi,g \act\zeta;\renormA,\renormphi).
\end{equation}
In other words, the data-to-solution map is gauge covariant. 
\end{lemma} 

\begin{remark}
The gauge covariance statement in Lemma~\ref{main:lem-gauge-covariance} is very different from the gauge covariance statements proven in \cite{CCHS22,CCHS22+,BC23}. In  \cite{CCHS22,CCHS22+,BC23}, the gauge covariance is only shown after taking the limit $\eps \rightarrow 0$ and only for suitable renormalisation constants $\renormA$ and $\renormphi$ (see also \cite{CS26}).  In contrast, Lemma~\ref{main:lem-gauge-covariance} holds for all values of $\eps$ and all $\liegroup$-covariant renormalisation constants. We believe that the gauge covariant smooth approximation \eqref{main:eq-SYMH-eps} may be useful in situations other than just the proof of Theorem~\ref{thm:main_informal}. In particular, it is possible that the gauge covariant smooth approximation \eqref{main:eq-SYMH-eps} can be compared more easily with lattice dynamics (such as in \cite{S21,IlyaHaoYM}) than the non-covariant smooth approximations from \cite{CCHS22,CCHS22+,BC23}. 
\end{remark}

Since the proof of Lemma~\ref{main:lem-gauge-covariance} is rather short, we present it here. 
\begin{proof}
Let $(A_\eps,\higgs_\eps)$ be the solution of \eqref{main:eq-SYMH-eps} and, to simplify the notation, let $(\widetilde{A}_\eps,\widetilde{\higgs}_\eps)\eqdef g \act (A_\eps,\higgs_\eps)$ be the gauge-transformed solution. Furthermore, let $(\widetilde{A}_{0},\widetilde{\higgs}_0)\eqdef g\act (A_0,\higgs_0)$ be the gauge-transformed initial data. Due to the gauge covariance of the Yang--Mills gradient flow, it is clear that $(\widetilde{A}_{0},\widetilde{\higgs}_0)$ is still contained in $\initialclass_T$ and that $\widetilde{B}=g\act B$ is the Yang--Mills gradient flow with initial data $\widetilde{A}_0$. Since both $\covd_{A_\eps}$ and $\bil$ are gauge covariant and since the difference $A_\eps-B$ transforms like a genuine one-form (see Remark \ref{intro:rem-genuine-one-form}), it then follows that 
\begin{align*}
\partial_t \widetilde{A}_\eps &= - \covd_{\widetilde{A}_{\eps}}^\ast F_{\widetilde{A}_{\eps}} - \covd_{\widetilde{A}_{\eps}} \covd_{\widetilde{A}_{\eps}}^\ast (\widetilde{A}_\eps-\widetilde{B}) + \bil\bigl( \widetilde{\higgs}_\eps \otimes \covd_{\widetilde{A}_\eps}\widetilde{\higgs}_\eps \bigr) + g \act \renormA ( A_\eps - B ) + g\act \xi_{\Beps}, \\ 
\partial_t \widetilde{\higgs}_\eps &= \Delta_{\widetilde{A}_\eps} \widetilde{\higgs}_\eps + (\covd_{\widetilde{A}_\eps}^\ast(\widetilde{A}_\eps-\widetilde{B})) \widetilde{\higgs}_\eps 
-  |\widetilde{\higgs}_\eps|^2 \widetilde{\higgs}_\eps + g\act \renormphi \higgs_\eps   + g\act \zeta_{\Beps},  \\
\widetilde{A}_\eps(0)&=\widetilde{A}_{0}, \qquad \widetilde{\higgs}_\eps(0)= \widetilde{\higgs}_0.   
\end{align*}
Since both $\renormA$ and $\renormphi$ commute with the action of $\liegroup$, we have that 
\begin{equation*}
g \act \renormA ( A_\eps - B ) = \renormA (\widetilde{A}_\eps - \widetilde{B}) \quad \text{and} \quad g\act \renormphi \higgs_\eps = \renormphi \widetilde{\higgs}_\eps. 
\end{equation*}
Furthermore, since the heat kernels $\fp_B$ and $\fq_B$ from Definition \ref{kernel:def-kernels} are gauge covariant, it holds that 
\begin{equation*}
g\act \xi_{\Beps} = (g\act \xi)_{\widetilde{B},\eps} \quad \text{and} \quad g\act \zeta_{\Beps} = (g\act \zeta)_{\widetilde{B},\eps}.
\end{equation*}
From this, we obtain \eqref{main:eq-gauge covariance}.
\end{proof}

As previously mentioned in the introduction, the gauge covariance of \eqref{main:eq-SYMH-eps} comes at a price. Unlike for the non-covariant smooth approximations from \cite{CCHS22,CCHS22+,BC23}, the process on the space of gauge orbits induced by \eqref{main:eq-SYMH-eps} is not necessarily Markov, and establishing the Markov property will require us to choose suitable renormalisation constants (see Section~\ref{section:Markov-and-Measure}). Of course, aside from the Markov property, we will also need to choose suitable renormalisation constants to obtain the local convergence of $(A_\eps,\higgs_\eps)$ as $\eps \rightarrow 0$. Our next definition makes a choice of renormalisation constants that will ensure both properties. 

\begin{definition}[Renormalisation constants]\label{main:def-renormalization-constant}
Let $\rho \in \mc{P}$ (recall \eqref{main:eq-mcP}) and let $\eps \in \sdyadic$. First, we define the real-valued constants
\begin{equation}\label{main:eq-sigma-eps-Crho}
\sigma_{\rho,\eps}^2 \eqdef \int_0^1 \dt  \frac{\rho_\eps^2(t)}{8\pi t} 
\qquad \text{and} \qquad \Crho \eqdef   \int_0^\infty \dt \Big( \int_0^t \ds \rho_\eps(s) \Big)  \frac{\rho^\prime_\eps(t)}{8\pi t}. 
\end{equation}
Second, we define $\renormA\in \End(\frkg)$ and $\renormphi\in \End(\V)$ as 
\begin{equs}\label{main:eq-renormalization-operators}
\renormA &\eqdef \Crho \bigl( \Cov^{\frkg}_{\frkg} + \Cov^{\frkg}_{\V}\bigr) -\tfrac{1}{8\pi} \Cov^{\frkg}_{\V} \;,\\
\renormphi &\eqdef (-\sigma_{\rho,\eps}^2 + \Crho ) \Cov^{\V}_{\frkg} + (\dim(\V)+2) \sigma_{\rho,\eps}^2 \Id_{\V} - \mass\;,
\end{equs}
where $\Cov^{\frkg}_{\frkg}$, $\Cov^{\frkg}_{\V}$, and $\Cov^{\V}_{\frkg}$  are as in Definition~\ref{prelim:def-Casimir-induced} below. 
\end{definition}

In Lemma \ref{prelim:lem-G-covariance} below, it will be shown that $\renormA \in \End_{\liegroup}(\frkg)$ and $\renormphi \in \End_{\liegroup}(\V)$, i.e., that both $\renormA$ and $\renormphi$ are $\liegroup$-covariant. 

\begin{remark}
While $\sigma_{\rho,\eps}^2$ diverges logarithmically as $\eps \rightarrow 0$, $\Crho$ does not depend on $\eps$. Indeed, a simple change of variables shows that the two instances of $\rho_\eps$
in its definition can be replaced by $\rho$. At this point in the article, it is not possible to motivate the precise definitions in \eqref{main:eq-sigma-eps-Crho}--\eqref{main:eq-renormalization-operators}. For the explicit calculations that led us to \eqref{main:eq-sigma-eps-Crho}--\eqref{main:eq-renormalization-operators}, we refer the reader to Sections~\ref{section:covariant-Ansatz},~\ref{section:objects}, and~\ref{section:langevin} below.
In particular, this specific choice for $\renormA$ is the only one guaranteeing that \eqref{main:eq-SYMH-eps}
yields a Markov process in the limit $\eps \to 0$. We will also see that in the notation of Section~\ref{section:langevin},
$\Crho$ is equivalent to $\<XiIIXir>$ in the sense that the latter depends on $\eps$ but converges to $\Crho$
as $\eps \to 0$.
\end{remark}

\begin{remark}
We emphasise that, while most of our estimates will be performed with respect to the background connection one-form $B$, the renormalisation constants in Definition~\ref{main:def-renormalization-constant} do not depend on $B$. To make this work, we need to prove several estimates of resonant parts of covariant stochastic objects, and these estimates heavily rely on our heat kernel parametrix $\fp \pt_B$ from Proposition \ref{main:prop-parametrix}. For more details, we refer the reader to Section~\ref{section:non-paracontrolled-resonant}. 
\end{remark}

\begin{remark}
The choice of upper limit $1$ in the definition of $\sigma_{\rho,\eps}^2$ is of course completely arbitrary,
but changing it is equivalent to a change of the mass $\mass$, so one can fix it without any loss of generality.
\end{remark}

In the following, the renormalisation constants in \eqref{main:eq-SYMH-eps} will always be chosen as in Definition~\ref{main:def-renormalization-constant}. Since the renormalisation constants were the last ingredient we needed, we can now state the main result of this subsection, namely, the covariant local well-posedness of the stochastic Yang--Mills--Higgs equations.

\begin{theorem}[Covariant local well-posedness]\label{main:thm-covariant-local}
Let $C_0>1$ and $c_0>0$ be sufficiently large and small constants, respectively. Let $0<T\leq c_0$, let $\rho\in \mc{P}$, and let $(A_0,\higgs_0) \in \initialclass_T$ be initial data satisfying
\begin{align}
T^\lossdel \exp\bigl( C_0 \| F_B \|_{L_t^1 L_x^\infty([0,T])}\bigr) \bigl( 1 + \Fnorm{F_B}_{T,100} + T^{(1-\roughdataweight)/2} \Fnorm{t^{-(1-\roughdataweight)/2} F_B}_{T,1} \bigr)^{100} &\leq 1, \label{main:eq-covariant-local-condition-1}\\ 
T  \| \higgs_0\|_{\covHolder{A_0, T}^{-\kappa}}^{\gamma_0} &\leq 1. \label{main:eq-covariant-local-condition-2}
\end{align}
Here, $B$ is as in \eqref{main:eq-initial-class}, $\lossdel$ and $\kappa$ are as in \eqref{e:small_constants}, and $\gamma_0$ is as in \eqref{e:constants_global_iteration}. Then, there exists an event $E(T,A_0,\higgs_0)$, depending on both the time $T$ and the initial data, such that the following properties hold: 
\begin{enumerate}[label=(\roman*),itemsep=0ex,topsep=0.1ex]
\item (High probability) It holds that 
\begin{equation*}
\bP \bigl( E(T,A_0,\higgs_0)\bigr) \geq 1 - c^{-1} \exp\bigl( - c T^{-c} \bigr), 
\end{equation*}
where $c$ is a constant depending on the parameters in \eqref{e:small_constants}. 
\item (Local convergence) For each $\eps\in \sdyadic$, let $(A_\eps,\higgs_\eps)$ be the unique solution of \eqref{main:eq-SYMH-eps} with the renormalisation constants from Definition~\ref{main:def-renormalization-constant}. Then, on the event $E(T,A_0,\higgs_0)$, the limit 
\begin{equation*}
(A,\higgs)= \lim_{\eps\rightarrow 0} (A_\eps,\higgs_\eps)
\end{equation*}
exists in $C_t^0 \Cs_x^{-\kappa}([0,T]\times \T^2 , \frkg^2 \times \V)$.
\end{enumerate}
\end{theorem}

\begin{remark}\label{main:rem-Cminuskappa}
Since the $C_t^0 \Cs_x^{-\kappa}$-norm of $(A,\higgs)$ is gauge dependent, the reader may be surprised to see it used in the convergence result of Theorem \ref{main:thm-covariant-local}. Our reasons for this are twofold: First, the convergence of the smooth approximations is a qualitative statement. And while we need all of our quantitative estimates, like the conditions in \eqref{main:eq-covariant-local-condition-1}--\eqref{main:eq-covariant-local-condition-2}, to only involve gauge invariant norms, we do not require this in qualitative statements. Second, we could use our solution norms (from Definition \ref{local:def-solution-space}) to state a gauge covariant convergence result. However, this would require notation that has not yet been introduced.
\end{remark}

\begin{remark}\label{main:rem-local-condition}
The conditions \eqref{main:eq-covariant-local-condition-1}--\eqref{main:eq-covariant-local-condition-2} will later be included in Assumptions \ref{local:assumption-B} and \ref{local:assumption-higgs-initial}, which will be used throughout all of our local arguments.
Since the condition \eqref{main:eq-covariant-local-condition-1} is important for our global iteration scheme, we discuss its main aspects. The polynomial factor in $\Fnorm{F_B}_{T,100}$ and $T^{(1-\roughdataweight)/2} \Fnorm{t^{-(1-\roughdataweight)/2} F_B}_{T,1}$ in \eqref{main:eq-covariant-local-condition-1} is less important than the factors $T^\lossdel$ and $\exp( C_0 \| F_B \|_{L_t^1 L_x^\infty([0,T])})$ and should be ignored for now. Since $\lossdel>0$ is small, we can then think of \eqref{main:eq-covariant-local-condition-1} as 
\begin{equation*}
\| F_B \|_{L_t^1 L_x^\infty([0,T])} \ll \abs{\log(T)}. 
\end{equation*}
If the $L_x^\infty$-norm\footnote{Since it appears in \eqref{main:eq-covariant-local-condition-1}, we use the $L_x^\infty$-norm of $F_B$ for this discussion, but our global iteration scheme will mostly rely on the $L_x^2$-norm of $F_{B^+}$ (see Section~\ref{section:main-global}), where the relation between $B$ and $B^+$ is as in Proposition \ref{main:prop-YMHF-rough-data}.} of $F_B$ is roughly constant on the time interval $[0,T]$, one can further think of \eqref{main:eq-covariant-local-condition-1} as 
\begin{equation*}
T \| F_B(0) \|_{L_x^\infty} \ll \abs{\log(T)}. 
\end{equation*}
Hence, $T$ only has to be as small as a single inverse power of the size of $F_B$. In fact, there is even a logarithmic amount of room, and this will later be crucial for our decay mechanism near Yang--Mills connections (recall Section~\ref{section:main-ideas-brief}). The quantitative condition on the size of the local time of existence from \eqref{main:eq-covariant-local-condition-1} is one of the main benefits of our covariant approach over the non-covariant approach from \cite{CCHS22,CCHS22+,BC23}. If we only used the existing non-covariant local theory, we would still be able to work only with gauge invariant norms of the initial data $A_0$, since we allow gauge transformations of the initial data (see Proposition~\ref{main:prop-YMHF-rough-data}). However, the cubic term in \cite[(1.7)]{CCHS22} means that the resulting local time of existence $T$ can at best be of the order of the inverse square of the size of $A_0$, which is insufficient for our global iteration scheme. 
\end{remark}

We emphasise that the precise formulation of Theorem~\ref{main:thm-covariant-local} was chosen for its relative simplicity. While Theorem~\ref{main:thm-covariant-local} itself contains little information on the structure of the solution, its proof (as well as Proposition~\ref{local:prop-main}) shows that the solution can be decomposed into covariant stochastic objects, paracontrolled terms, and smooth, nonlinear remainders. The structure of the solution will not only be needed in the proof, but will also be used in our global iteration argument. 

In the following, we describe the main ideas behind the proof of Theorem~\ref{main:thm-covariant-local}. Our discussion will be more detailed than in Section~\ref{section:main-ideas-brief}, but less detailed than in the body of the article (Sections~\ref{section:covariant-Ansatz}--\ref{section:objects}). In our covariant local theory for \eqref{main:eq-SYMH-eps}, we expand around the background connection one-form $B$. That is, we set $C_\eps\coloneq A_\eps- B$ and consider the system
\begin{equation}\label{main:local-eq-C-system}
\begin{aligned}
( \partial_t - \Delta_B ) C_\eps^k - 2 \big[ F_B^{kj},C_{\eps,j} \big] &=  \big[ C_{\eps,j} , 2 \covd_B^j C_\eps^k - \covd_B^k C_\eps^j \big]
+ \big[ C_{\eps,j}, \big[ C_\eps^j , C_\eps^k \big]\big] \\
&\quad +\bil\bigl( \higgs_\eps \otimes \covd_B^k \higgs_\eps \bigr) + \bil \bigl( \higgs_\eps \otimes C_\eps^k \higgs_\eps\bigr) + \renormA ( C_\eps^k ) + \xi_{\Beps}^k, \\ 
(\partial_t - \Delta_B ) \higgs_\eps &= 2 C_{\eps,j} \covd_B^j \higgs_\eps + C_\eps^j C_{\eps,j} \higgs_\eps -   |\higgs_\eps|^2 \higgs_\eps + \renormphi(\higgs_\eps) + \zeta_{\Beps}. 
\end{aligned}
\end{equation}
For the derivation of \eqref{main:local-eq-C-system}, we refer the reader to Lemma~\ref{prelim:lem-linearization} and Section~\ref{section:covariant-Ansatz}. The first step in our covariant ansatz for $(C_\eps,\higgs_\eps)$ is to introduce the deterministic terms $\Sdet$ and $\higgsdet$, which do not depend on the stochastic forcing terms $\xi_{\Beps}$ and $\zeta_{\Beps}$. They are defined as solutions of
\begin{align}
(\partial_t - \Delta_B)  \Sdet^k - 2 [F_B^{kj},\Sdet_j]&= [ \Sdet_j, 2 \covd_B^j \Sdet^k - \covd_B^k \Sdet^j] + [ \Sdet_j , [\Sdet^j , \Sdet^k ]] \label{main:local-eq-Sdet} \\
&\qquad + \bil\bigl( \higgsdet \otimes \covd_B^k \higgsdet\bigr) + \bil \bigl( \higgsdet \otimes \Sdet^k \higgsdet \bigr), \notag \\ 
(\partial_t - \Delta_B) \higgsdet &= 2 \Sdet_j \covd_B^j \higgsdet  + \Sdet^j \Sdet_j \higgsdet - |\higgsdet|^2 \higgsdet, \label{main:local-eq-higgsdet} \\ 
\Sdet(0)&=0, \qquad \higgsdet(0)=\higgs_0. \notag 
\end{align}
The equations~\eqref{main:local-eq-Sdet}--\eqref{main:local-eq-higgsdet}
are chosen so that $(B+\Sdet,\higgsdet)$ solves the deterministic
Yang--Mills--Higgs equations (with zero mass) in the infinitesimal DeTurck gauge
with background $B$. Since both $\Sdet$ and $\higgsdet$ will be rather smooth (see e.g.\ Proposition~\ref{prop:Sdet-higgsdet-estimates}), $\Sdet$ and $\higgsdet$ will not be too important in our covariant local theory. Their real significance will only become clear in our global iteration argument, where we will use that $\Sdet$ and $\higgsdet$ capture the largest parts of our solution in terms of the size of $\higgs_0$. Instead of $\Sdet$ and $\higgsdet$, the covariant local theory is mostly dictated by the covariant stochastic objects $\linear[\Beps]$ and $\philinear[\Beps]$. As briefly mentioned in Section~\ref{section:main-ideas-brief}, they are defined as solutions of 
\begin{align}
(\partial_t - \Delta_B) \linear[\Beps][k][] - 2 [F_B^{kj},\linear[\Beps][][j] ]  &= \xi_{\Beps}^k, \qquad \linear[\Beps][k][] (0)=0, \label{main:eq-linear-stochastic-A} \\ 
(\partial_t - \Delta_B )  \philinear[\Beps]&= \zeta_{\Beps}, \qquad \philinear[\Beps] (0)=0. \label{main:eq-linear-stochastic-phi}
\end{align}
We also denote the limits of $\linear[\Beps]$ and $ \philinear[\Beps]$ as $\eps\rightarrow 0$ by $\linear[B]$ and $ \philinear[B]$.  
In addition to the linear stochastic objects in \eqref{main:eq-linear-stochastic-A}--\eqref{main:eq-linear-stochastic-phi}, we will also need quadratic stochastic objects, which are defined as the solutions with zero initial data of 
\begin{align}
(\partial_t - \Delta_B) \quadratic[\Beps][k] - 2 [F_B^{kj},\quadratic[\Beps][][j] ]  
    &= \big[ \linear[\Beps][][j], 2 \covd_B^j \linear[\Beps][k][] - \covd_B^k \linear[\Beps][j]\big], \label{main:eq-A-quadratic-stochastic} \\ 
(\partial_t - \Delta_B) \phiquadratic[\Beps][k] - 2 [F_B^{kj},\phiquadratic[\Beps][][j] ]   
    &= \bil \bigl(  \philinear[\Beps] \otimes \covd_B^k \philinear[\Beps]\bigr), \label{main:eq-phi-quadratic-stochastic}\\ 
(\partial_t - \Delta_B)  \mixedquadratic[\Beps] 
&= 2 \linear[\Beps][][j] \covd_{B}^j \philinear[\Beps]. \label{main:eq-mixed-quadratic-stochastic}
\end{align}
We emphasise that the covariant objects from \eqref{main:eq-linear-stochastic-A}--\eqref{main:eq-mixed-quadratic-stochastic} can also be represented using the covariant heat kernels $\fp_B$ and $\fq_B$. For instance, we can write 
\begin{equation}\label{main:eq-linear-integral-representation}
 \linear[\Beps](t,x) = \int_0^t \ds \int_{\T^2} \dy \rho_\eps(t-s) \fq_B(s,y;t,x) \xi(s,y).
\end{equation}
In our analysis of the covariant stochastic objects, we then heavily rely on covariant heat kernel estimates. While a similar approach was previously used in work of the first and second authors on the Abelian Higgs model \cite{BC24}, our implementation of this approach differs significantly from \cite{BC24}. As already mentioned in Remark \ref{remark:main-thm-ptwise-bds}, \cite{BC24} relies on operator estimates of the covariant heat kernel and its covariant derivatives, whereas we rely on the pointwise estimates from Section~\ref{section:main-heat}. This drastically simplifies one of the hardest parts of the proof in \cite{BC24}, namely, the estimates of the resonant part of the quadratic stochastic object (cf. \cite[Proposition 5.29]{BC24} and Lemmas~\ref{lemma:qres-estimate-sol-space} and~\ref{lemma:covd-philinear-times-philinear-resonance} below). 

Another important difference between \cite{BC24} and the present article is that the covariant methods in our local theory are not only used for the covariant stochastic objects. Instead, we build an entire covariant functional framework, which encompasses covariant Littlewood--Paley decompositions, function spaces, Schauder estimates, and paraproducts. To motivate this framework, we first briefly recall the standard Littlewood--Paley decompositions. Throughout this discussion only, we let $\varsigma\colon \R^2 \rightarrow [0,\infty)$ be a smooth function satisfying $\varsigma|_{B(0,1)}\equiv 1$ and $\varsigma|_{\R^2\backslash B(0,2)}\equiv 0$, where $B(0,r) \subseteq \R^2$ is the open Euclidean ball centred at the origin with radius $r$. For any $\phi \colon \T^2 \rightarrow \C$, one can then define the standard Littlewood--Paley operators $(P_{\leq N}\phi)_{N\in \dyadic}$ and $(P_N \phi)_{N\in \dyadic}$ using the Fourier transform as 
\begin{equation}\label{main:eq-standard-Littlewood-Paley}
\begin{aligned}
\widehat{P_{\leq N}\phi}(n) \eqdef  \varsigma(n/N) \widehat{\phi}(n) \qquad \text{and} \qquad  P_{N} \phi \eqdef (P_{\leq N} - P_{\leq N/2}) \phi, 
\end{aligned}
\end{equation}
where we use the convention $P_{\leq 1/2}\eqdef 0$. 
Starting from the Littlewood--Paley projections in \eqref{main:eq-standard-Littlewood-Paley}, one can then build up Hölder spaces, Schauder estimates, and paraproducts (cf. \cite{GIP15}). 

Unfortunately, since they do not interact nicely with gauge transformations, the standard Littlewood--Paley projections from \eqref{main:eq-standard-Littlewood-Paley} will not be too useful in this article. In order to introduce our covariant version, recall from \eqref{main:eq-linear-integral-representation} that the covariant stochastic object $\linear[\Beps]$ can be written as an integral of the kernel $\rho_\eps \fq_B$ against the space-time white noise $\xi$. For any $N\in \dyadic$, we can then define
\begin{equation}\label{main:eq-covariant-Littlewood-Paley-block}
\linear[\Beps][][][N] \eqdef \iint \ds \dy \,  \rho_\eps(t-s) \chi_N(t-s) \fq_B(s,y;t,x) \xi(s,y), 
\end{equation}
where $\chi_N(t-s)$ is the smooth cutoff to times $(t-s) \sim N^{-2}$  from \eqref{prelim:eq-chi-N}.  Due to parabolic scaling and the uncertainty principle, the restriction to times $N^2 (t-s) \sim 1$ in \eqref{main:eq-covariant-Littlewood-Paley-block} is morally similar to a projection to frequencies of size $\sim N$. For a more detailed discussion of \eqref{main:eq-covariant-Littlewood-Paley-block}, as well as potential alternatives, we refer the reader to Section~\ref{section:Littlewood--Paley} below. More generally, we introduce localised versions of the Duhamel integral operators of \eqref{main:local-eq-C-system}. For example, for any $H \colon [0,T] \times \T^2 \rightarrow \frkg^2$, we introduce
\begin{equation}\label{main:eq-covariant-Duhamel-localized}
\Duh_{B,N}^F(H)(t,x) \eqdef \iint \ds \dy \, \chi_N(t-s) \fq_B(s,y;t,x) H(s,y). 
\end{equation}
Definitions similar to those in \eqref{main:eq-covariant-Littlewood-Paley-block} and \eqref{main:eq-covariant-Duhamel-localized} can also be made for the Higgs component of our solution by replacing the kernel $\fq_B$ with $\fp_B$. Starting from \eqref{main:eq-covariant-Littlewood-Paley-block} and \eqref{main:eq-covariant-Duhamel-localized}, we then build an entire covariant functional framework. Compared to the non-covariant paracontrolled calculus of \cite{GIP15}, the main difficulty is that our covariant Littlewood--Paley decomposition does not obey the Littlewood--Paley trichotomy; see Section \ref{section:Littlewood--Paley}. Since the full definitions and estimates included in the functional framework are well beyond the scope of this section, we do not include them here, and instead only provide the following brief descriptions and examples: 
\begin{enumerate}[label=(\roman*)]
\item (Solution spaces) Our solution spaces will be defined in terms of the covariant Littlewood--Paley blocks. The covariant linear stochastic object will be controlled in our solution space at regularity $-\kappa$, which requires estimates of 
\begin{equation*}
\sup_{N\in \dyadic} N^{-m-\kappa} \| \covd_{B,j_1} \hdots \covd_{B,j_m} \linear[\Beps][][][N] \|_{L_{t,x}^\infty}. 
\end{equation*}
Here, $0 \leq m \leq 10$ and $1 \leq j_1,\hdots, j_m \leq 2$. 
\item (Schauder estimates) Our Schauder estimates will be phrased in terms of the localised Duhamel integrals from \eqref{main:eq-covariant-Duhamel-localized}. For example, we will prove that 
\begin{equation*}
\| \covd_{B,j_1} \hdots \covd_{B,j_m} \Duh_{B,N}^F ( H) \|_{L_{t,x}^\infty} 
\lesssim \lossfactor_T(F_B) N^{m-2} \| H \|_{L_{t,x}^\infty}. 
\end{equation*}
Here, as before, $0\leq m \leq 10$ and $1\leq j_1,\hdots,j_m \leq 2$. Furthermore, the loss-factor $\lossfactor_T(F_B)$ from \eqref{local:eq-lossfactor} captures the $F_B$-terms in our covariant heat kernel estimates (Theorem \ref{main:thm-pointwise-kernel}). 
\item (Paraproducts) We will define several covariant paraproducts. For example, if $(W_M)_M$ is a sequence of $\frkg^2$-valued Littlewood--Paley blocks and $(\linear[\Beps][][][N])_N$ is as in \eqref{main:eq-covariant-Littlewood-Paley-block}, we define the low$\times$high-product 
\begin{equation*}
\bigl[ W_\star^i \parall \covd_{B}^k \linear[\Beps][j][][\ostar] \bigr] 
= \sum_{\substack{M,N\in \dyadic\colon \\ M \leq N^{1-\theta}}} \bigl[ W_M^i ,  \covd_B^k \linear[\Beps][j][][N] \bigr].
\end{equation*}
Here, $1\leq i,j,k \leq 2$ and the symbol $\star$ is used to denote sequences of Littlewood--Paley blocks. Furthermore, the parameter $\theta$ is as in \eqref{e:small_constants}, and it is needed due to the absence of the Littlewood--Paley trichotomy. 
\item (Commutator estimates) We prove several commutator estimates that are needed to understand the structure of the paracontrolled term. For example, in the same setting as above, we will estimate\footnote{Here, $\Duh_{B,N}$ denotes the same localised Duhamel integral operator as in \eqref{main:eq-covariant-Duhamel-localized}, except that the kernel $\fq_B$ has been replaced with $\fp_B$.} 
\begin{equation*}
\Duh_{B,K} \bigl(\bigl[ W_M^i ,  \covd_{B}^k \linear[\Beps][j][][N] \bigr]  \bigr) 
- \bigl[ W^i_M, \covd_{B}^k \Duh_{B,K} \bigl( \linear[\Beps][j][][N] \bigr) \bigr].
\end{equation*}
\end{enumerate}

\begin{remark} As described above, this article relies on a covariant version of the paracontrolled calculus of \cite{GIP15}. In order to make progress towards the global well-posedness of the stochastic Yang--Mills--Higgs model in dimension $d=3$ (rather than $d=2$), it may be useful to also develop a covariant version of the theory of regularity structures \cite{H14}. The theory of regularity structures has already been studied in variable-coefficient settings \cite{HairerSingh23,Singh25,BSS25}, which include \eqref{main:local-eq-C-system}. However, the focus of \cite{HairerSingh23,Singh25,BSS25} and this article is quite different. In \cite{HairerSingh23,Singh25,BSS25}, the geometry is fixed once and for all, and the main task is to show that the counterterms only depend on local aspects of the geometry. In this article, the geometry is determined by the background connection one-form $B$, which changes throughout the global iteration scheme (see Section \ref{section:main-global}). In particular, one of the main tasks in the proof of Theorem \ref{main:thm-covariant-local} is to determine how the size of the covariant stochastic objects depends on the size of $F_B$. In addition to \cite{HairerSingh23,Singh25,BSS25}, we also highlight the works \cite{Hos25, HT25}. The main similarity between \cite{Hos25, HT25} and this article is that both start from a heat kernel estimate (compare Definition 2.6.(iii) in \cite{Hos25} and Theorem \ref{main:thm-pointwise-kernel}). 
\end{remark}

\subsection{Global well-posedness and uniform-in-time estimates}\label{section:main-global}

In this section, we state both soft and hard global well-posedness theorems (Theorems \ref{thm:global_soft} and \ref{thm:global_hard}). The statement of the soft result (Theorem \ref{thm:global_soft}) will be relatively simple and, after combining it with the Krylov--Bogoliubov theorem, it will be sufficient to conclude the existence of the invariant measure from Theorem \ref{thm:main_informal}. However, the soft result contains little information on the structure or the decay of the solution to the stochastic Yang--Mills--Higgs equations. The  $\invnorm{\place}{-\kappa}$-norm used in its statement does not separate the smooth and rough terms in the solution and, as a result, the $\invnorm{\place}{-\kappa}$-norm will not be used directly in our global iteration scheme. 
Our hard result (Theorem \ref{thm:global_hard}) has a more complicated statement, but it makes up for it by providing more information on the structure of the solution. In particular, the statement of our hard result will be closer to the estimates used in our proof. 

In order to state our soft global well-posedness result, we first make the following definition.\ 

\begin{definition}\label{main:def-invnorm}
Let $A_0\in \Cs_x^{-\kappa}(\T^2, \frkg^2)$ and $\higgs_0 \in \Cs_x^{-\kappa}(\T^2, \V)$. We define a gauge invariant Hölder norm as follows. First, let $T > 0$ and $g \in \CG_{1-\kappa}$ be a gauge transformation such that $g \act (A_0, \higgs_0) \in \initialclass_T$, where $\initialclass_T$ is as defined in \eqref{main:eq-initial-class}. Recall from Remark~\ref{remark:initial-class-time} that we then also have $g \act (A_0, \higgs_0) \in \initialclass_1$. Thus by definition, there exists $B \in C_t^0 \Cs_x^{-\kappa}([0, 1] \times \T^2, \frkg^2) \cap C_{t,x}^\infty((0, 1] \times \T^2, \frkg^2)$ which is a smooth classical solution to the Yang--Mills gradient flow at positive times, and such that $B(0) = g \act A_0$. Let $\Psi(s)$ be the solution to the covariant heat equation $(\ptl_s - \Delta_{B(s)}) \Psi(s) = 0$ with initial datum $g \act \higgs_0$. We then let
\begin{equ}\label{main:eq-invnorm}
\invnorm{(A_0, \higgs_0)}{-\kappa} \eqdef \sup_{s \in (0, 1]} s^{\frac{\kappa}{2}} \Bigl( s^{\frac{1}{2}} \big\| F_B(s) \big\|_{L^\infty_x} + s \big\| \covd_B F_B(s) \big\|_{L_x^\infty} + \big\| \Psi(s) \big\|_{L^\infty_x} \Bigr).
\end{equ}
\end{definition}

The fact that this gauge invariant Hölder norm is well-defined follows by Lemma~\ref{lemma:invnorm-well-def}. We define the norm in this way, because it is not necessarily possible to start a Yang--Mills gradient flow from arbitrary $\Cs_x^{-\kappa}(\T^2, \frkg^2)$ initial data. On the other hand, it is in fact always possible to find a $\CG_{1-\kappa}$ gauge transformation such that after applying this transformation, we can start a (unique) Yang--Mills gradient flow. Lemma~\ref{lemma:invnorm-well-def} shows that the resulting quantity in \eqref{main:eq-invnorm} does not depend on this choice of gauge transformation (and moreover that such a gauge transformation exists).

We emphasise that the gauge invariant Hölder norm in \eqref{main:eq-invnorm} does not require the choice of a background connection one-form, which is not the case for the norms used in our covariant local theory (or the global iteration scheme described below). This allows us to measure the size of the solution to the stochastic Yang--Mills--Higgs equations without any reference to its random structure.

\begin{theorem}[Soft global well-posedness and uniform-in-time estimates]\label{thm:global_soft}
Let $A_0\in \Cs_x^{-\kappa}(\T^2 , \frkg^2)$ and $\higgs_0\in \Cs_x^{-\kappa}(\T^2 , \V)$. Then, there exist a connection one-form $A\in L_{t,\textup{loc}}^\infty \Cs_x^{-\kappa}([0,\infty)\times \T^2 , \frkg^2)$, a Higgs field $\higgs\in L_{t,\textup{loc}}^\infty \Cs_x^{-\kappa}([0,\infty)\times \T^2 , \V)$, and stopping times $(t_m)_{m=0}^\infty$ starting at $t_0=0$ and a.s.\ converging to infinity, which satisfy the following properties:
\begin{enumerate}[label=(\roman*),itemsep=0ex]
\item\label{main:item-soft-SYMH} (Stochastic YMH equations) On each time interval $(t_m,t_{m+1}]$, where $m\geq 0$, $(A,\higgs)$ solves \eqref{intro:eq-SYMH}. To be more precise, it is the limit of the smooth, gauge covariant approximations from \eqref{main:eq-SYMH-eps} in $C_t^0 \Cs_x^{-\kappa}((t_m,t_{m+1}]\times \T^2, \frkg^2 \times \V)$.
\item\label{main:item-soft-initial} (Initial conditions and right-limits) The initial value $(A(0),\higgs(0))$ is equal to $(A_0,\higgs_0)$. Furthermore, for each $m\geq 0$, the right-limit $(A(t_{m}+),\higgs(t_{m}+))$ exists in $\Cs_x^{-\kappa}(\T^2,\frkg^2\times \V)$. This right-limit belongs to the initial class  $\initialclass_{t_{m+1}-t_{m}}$ defined in \eqref{main:eq-initial-class} and is gauge-equivalent to $(A(t_{m}),\higgs(t_{m}))$, i.e.,  there exists a $g_{m+1}\in \CG_{1-\kappa}$ such that $(A(t_{m}+),\higgs(t_{m}+))=g_{m+1}\act (A(t_{m}),\higgs(t_{m}))$. 
\item \label{item:thm:global_soft_unif_in_time} (Uniform-in-time estimate) For all $p\geq 1$, we have the gauge invariant uniform-in-time estimate
\begin{equation*}
\sup_{t\in [0,\infty)} \E \Bigl[ \invnorm{(A(t),\higgs(t))}{-\kappa}^p \Bigr] \leq C(A_0,\higgs_0,\kappa,p). 
\end{equation*}
Moreover, as $t\rightarrow \infty$, the bound on $(A,\higgs)$ no longer depends on the initial data $(A_0,\higgs_0)$, i.e., 
\begin{equation*}
\limsup_{t\rightarrow \infty} \E \Bigl[ \invnorm{(A(t),\higgs(t))}{-\kappa}^p \Bigr] \leq C(\kappa,p).
\end{equation*}
\end{enumerate}
\end{theorem}

We emphasise that, as in Theorem \ref{main:thm-covariant-local}, we only use the $C_t^0 \Cs_x^{-\kappa}$- and $\Cs_x^{-\kappa}$-norms for qualitative statements (see Remark \ref{main:rem-Cminuskappa}). Note also that the stochastic Yang--Mills--Higgs equations \eqref{intro:eq-SYMH} themselves are not solved globally in time. Instead, we only solve \eqref{intro:eq-SYMH} on a short time interval, then apply a gauge transformation, and then solve \eqref{intro:eq-SYMH} again on the next short time interval. The reason is that, as previously discussed below \eqref{intro:eq-linearized}, the Yang--Mills--Higgs gradient flow with a DeTurck term does not obey suitable long-time estimates. Since we are ultimately interested in the global dynamics of the Markov process on the space of gauge orbits, this use of gauge transformations is allowed.

\begin{remark}
By \ref{main:item-soft-SYMH} and \ref{main:item-soft-initial}, the process $(A,\higgs)$ in Theorem~\ref{thm:global_soft} is left-continuous with right limits. This is purely a matter of convention: by redefining the values of $(A,\higgs)$ at the stopping times $(t_m)_{m=0}^\infty$ and replacing the interval $(t_m,t_{m+1}]$ in \ref{main:item-soft-SYMH} with $[t_m,t_{m+1})$, one obtains an analogous result for a right-continuous process with left limits. This distinction is not important for our purposes since, as noted above, we are ultimately interested only in the projection of $(A,\higgs)$ onto the space of gauge orbits, and this projection will be continuous.
\end{remark}

As stated at the beginning of this subsection, the soft result from Theorem~\ref{thm:global_soft} will be the main ingredient needed for the existence of the invariant measure from Theorem~\ref{thm:main_informal}. For this reason, Theorem~\ref{thm:global_soft} can perhaps be regarded as the most important statement in the article. However, our soft result will not be proven on its own, but together with our hard result (Theorem \ref{thm:global_hard} below). In the following, we first discuss the main definitions and ideas behind the global iteration scheme that leads to our hard estimates. The statement of our hard global well-posedness result will be included after this discussion.

To start off, we assume that our solution $(A,\higgs)$ has already been constructed up to a stopping time $t_m$, where $m\geq 0$. We further assume that $(A,\higgs)(t_m)$ has been decomposed into a smooth part $(\Asmoothpart,\phismoothpart)(t_m)$ and a rough part $(\Aroughpart,\phiroughpart)(t_m)$. The size of the smooth part will be measured using 
\begin{equation}\label{main:eq-covariant-energy}
\covariantenergy(\Asmoothpart(t_m),\phismoothpart(t_m)) \eqdef \| F_{\Asmoothpart}(t_m) \|_{L^2_x}^2 
+ \| \phismoothpart(t_m) \|_{L_x^4}^{\gamma_1} + \| \phismoothpart(t_m) \|_{L_x^r}^{\gamma_2} + 10, 
\end{equation}
where $r\gg 1$ is as in \eqref{e:small_constants} and $\gamma_1 >4 > \gamma_2 = 3$ are as in \eqref{e:constants_global_iteration}. We emphasise that, since $\gamma_1>\gamma_2$ and $4<r$, neither the $L^4$-term nor the $L^r$-term is stronger than the other. Later on, we will see that the $L^4$-term is needed to control the $\bil(\higgs \otimes \covd_A \higgs)$-term in \eqref{intro:eq-SYMH}, and that the $L^r$-term is needed to control all other interactions stemming from the smooth part of the Higgs field. We refer to \eqref{main:eq-covariant-energy} as the monotone quantity. The reason is that it decreases under the deterministic dynamics\footnote{Unlike the Yang--Mills--Higgs energy from \eqref{intro:eq-YMH-energy}, the monotone quantity \eqref{main:eq-covariant-energy} is not necessarily monotone under the Yang--Mills--Higgs gradient flow. However, it is monotone if the connection one-form is a Yang--Mills gradient flow and the Higgs field solves a cubic covariant heat equation. Since the $\bil(\higgs \otimes \covd_A \higgs)$-term in the Yang--Mills--Higgs gradient flow can be controlled separately, this will be sufficient for our purposes. For more details, we refer the reader to Section \ref{section:global-iteration-monotone-quantities-deterministic}.} and, under suitable conditions on the time scale $\tau_{m+1}$, the deterministic decrease on the time interval $(t_m,t_{m+1}]$ can be bounded below by
\begin{equs}\label{main:eq-decrease}
\Dec_{m+1} &\eqdef 2 \| \covd_{B^+_{m+1}}^\ast F_{B^+_{m+1}} \|_{L_{t,x}^2((t_m,t_{m+1}])}^2 \\
&+ C_{\gamma_2,r} \one_{\| \phismoothpart(t_m) \|_{L_x^r} \geq 1} \| \phismoothpart(t_m) \|_{L_x^r}^{\gamma_2 + 2} \tau_{m+1} \\
&+ \one_{\| \phismoothpart(t_m) \|_{L_x^4} \geq 1} \Bigl( \gamma_1 \| \shiggsdet_{m+1}(t_{m+1}) \|_{L_x^4}^{\gamma_1 - 4} \| |\covd_{B_{m+1} + \Sdet_{m+1}} \shiggsdet_{m+1}| |\shiggsdet_{m+1}| \|_{L_{t,x}^2((t_m,t_{m+1}])}^2 \\
&\qquad \qquad \qquad \qquad \qquad \qquad \qquad \qquad \qquad \qquad \qquad  \vee C_{\gamma_1,4} \| \phismoothpart(t_m) \|_{L_x^4}^{\gamma_1 + 2} \tau_{m+1} \Bigr) \;.
\end{equs}
Here, $B^+_{m+1}$ is the Yang--Mills gradient flow with initial data $\Asmoothpart(t_m)$. Furthermore, after replacing the time interval with $[t_m,t_{m+1}]$, $\Sdet_{m+1}$ is as in \eqref{main:local-eq-Sdet} and $\shiggsdet_{m+1}$ is the smooth part of the solution of \eqref{main:local-eq-higgsdet}. Finally, $C_{\gamma_1,4},C_{\gamma_2,r}>0$ are suitable constants. For more details, see the discussion below \eqref{main:eq-list-iterates} and Section \ref{section:deterministic-rough}. The size of the rough part $(\Aroughpart,\phiroughpart)(t_m)$ will be measured using the covariant Hölder norms from Definition \ref{defn:covariant_Holder_global}, which use the Yang--Mills gradient flow $B_{m+1}^+$ starting from $\Asmoothpart(t_m)$ as a background connection one-form. For the purposes of the current discussion, we encourage the reader to think of the rough part as being of unit size, i.e., that 
\begin{equation}\label{main:eq-rough-part-heuristic-size}
\|\Aroughpart (t_m) \|_{\covHolder{\Asmoothpart(t_m), F}^{-\kappa}} \lesssim 1 
\qquad \text{and} \qquad 
\|\phiroughpart (t_m) \|_{\covHolder{\Asmoothpart(t_m)}^{-\kappa}} \lesssim 1.
\end{equation}
For a more rigorous treatment of the rough part, we refer the reader to Definition \ref{defn:L_admissible} (as well as the rest of Section \ref{section:global-theory}). After the above discussion of the sizes of the smooth and rough parts, we now turn to the choice of the time scale $\tau_{m+1}$, which will determine the next stopping time $t_{m+1}\eqdef t_m + \tau_{m+1}$. The most important (but not only) question is how the time scale $\tau_{m+1}$ depends on the smooth part $(\Asmoothpart,\phismoothpart)(t_m)$. To capture this, we make the following definition.

\begin{definition}[Time scales]\label{main:def-time-scales}
Let  $\Asmoothpart\colon \T^2 \rightarrow \frkg^2$ be a smooth connection one-form and $\phismoothpart \colon \T^2 \rightarrow \V$ be a smooth Higgs field. Furthermore, let $T_0$ be as in \eqref{prelim:eq-T0} below and let $0<\logtimestep \leq 10^{-10}$ and $0<\stepupperbound \leq T_0$ be small constants that remain to be chosen. Then, we define 
\begin{equation}\label{main:eq-time-step-1}
T^{(1)}(\Asmoothpart,\phismoothpart) \eqdef  
\min \Bigl\{ \logtimestep \| F_{\Asmoothpart} \|_{L^2_x}^{-1} 
\log\bigl( \| F_{\Asmoothpart} \|_{L^2_x} \vee 2\bigr), \| \phismoothpart \|_{L_x^r}^{-\gamma_0}, \stepupperbound, \logtimestep \Bigr\} 
\end{equation}
and
\begin{equation}\label{main:eq-time-step-2}
T^{(2)}(\Asmoothpart,\phismoothpart) = T^{(2)}(\Asmoothpart)  \eqdef \inf \Bigl\{ t>0 \colon  2 \| \covd_{B^+}^\ast F_{B^+} \|_{L^2_{t,x}([0,t])}^2 \geq \logtimestep t^{-1} \log(t^{-1}) \Bigr\}. 
\end{equation}
Here, $\gamma_0>2$ is as in \eqref{e:constants_global_iteration} and $B^+$ is the Yang--Mills gradient flow starting from $\Asmoothpart$. 
\end{definition}
For the rest of this discussion, the reader may assume that the time scale $\tau_{m+1}$ is given by 
\begin{equation*}
\tau_{m+1} = (T^{(1)} \wedge T^{(2)})(\Asmoothpart(t_m),\phismoothpart(t_m)). 
\end{equation*}
In the full definition of our time scale $\tau_{m+1}$ (see Definition \ref{defn:tau_L}), this will be true except in the rare cases in which either the rough part or one of the covariant stochastic objects is much larger than its expected size. With Definition \ref{main:def-time-scales} at hand, it is now important to ask whether our covariant local theory (from Section \ref{section:main-local}) can be used on the time interval $(t_m,t_{m+1}]$. Due to Proposition~\ref{main:prop-YMHF-rough-data} and Theorem~\ref{main:thm-covariant-local}, the most important question for this is whether we can control $\| F_{\smash{B_{m+1}^+}} \|_{L_t^1 L_x^\infty((t_m,t_{m+1}])}$, where $B^+_{m+1}$ is the Yang--Mills gradient flow with initial data $B_{m+1}^+(t_m)=\Asmoothpart(t_m)$. And indeed, by using Lemmas \ref{lemma:curvature-equation}, \ref{prelim:lem-diamagnetic}, and \ref{lemma:strichartz-estimate-using-decrease} from the body of this article, one has that 
\begin{equation}\label{main:eq-L1Linf-curvature}
\begin{aligned}
\|& F_{B_{m+1}^+} \|_{L_t^1L_x^\infty((t_m,t_{m+1}])} \\ 
&\lesssim \tau_{m+1}^{1/2} \| F_{B_{m+1}^+} \|_{L_t^2 L_x^\infty((t_m,t_{m+1}])} \\ 
&\lesssim \tau_{m+1} \| F_{B_{m+1}^+}(t_m) \|_{L_x^2} + \tau_{m+1}^{1/2} \abs{\log(\tau_{m+1})}^{1/2} \| \covd_{B_{m+1}^+}^\ast F_{B_{m+1}^+}\|_{L_t^2 L_x^2((t_m,t_{m+1}])} \\ 
&\lesssim \sqrt{\logtimestep} \abs{\log(\tau_{m+1})}. 
\end{aligned}
\end{equation}

\begin{remark}
The $\log$-factors in \eqref{main:eq-time-step-1} and \eqref{main:eq-time-step-2} will be needed in our global iteration scheme. The reason is that the decay mechanism near Yang--Mills connections (see Section \ref{section:main-ideas-brief}) relies on their instability, and only the $\log$-factors allow us to venture slightly into the unstable regime.
\end{remark}

\begin{remark}
Since the $L_t^1 L_x^\infty$-norm of the curvature plays a central role in our approach, it is natural to ask why it does not directly appear in the definition of the time scales (Definition \ref{main:def-time-scales}). The reason for this is that it appears neither in the monotone quantity \eqref{main:eq-covariant-energy} nor in its decrease \eqref{main:eq-decrease}, and including it in the time scales would therefore not fit nicely into our global iteration scheme. In particular, the fact that $T^{(2)}$ controls $\covd_{B^+}^\ast F_{B^+}$ will be needed to control several cross-terms in our global iteration scheme; see, e.g. the estimates of \eqref{main:eq-curvature-linear-S} and \eqref{main:eq-curvature-linear-lollipop}. 
\end{remark}

Assuming now that Proposition \ref{main:prop-YMHF-rough-data} and Theorem \ref{main:thm-covariant-local} (as well as their more detailed versions from Propositions~\ref{pr:admissible_timescale_general} and~\ref{local:prop-main}) can be applied\footnote{In Theorem \ref{main:thm-covariant-local}, the conclusions hold on an event with high probability, but not almost surely. In our global iteration scheme, however, most estimates will hold almost surely. This can be achieved by modifying the time scale $\tau_{m+1}$ on an event with small probability.} on the time interval $[t_m,t_{m+1}]$ with the initial data $(\Asmoothpart,\Aroughpart,\phismoothpart,\phiroughpart)(t_m)$, we then obtain the maps
\begin{equation}\label{main:eq-list-iterates}
g_{m+1}, A_{m+1}, \higgs_{m+1}, B_{m+1}, \Sdet_{m+1}, \higgsdet_{m+1}, \linear[B_{m+1}][][][], \philinear[B_{m+1}][][][], R_{m+1}, \quad \text{and} \quad \phir_{m+1}, 
\end{equation}
which are all defined on $[t_m,t_{m+1}]\times \T^2$ as follows:
\begin{enumerate}[label=(\alph*),itemsep=0ex]
\item\label{main:item-gauge} $g_{m+1}$ is the gauge transformation from Proposition \ref{main:prop-YMHF-rough-data}. 
\item $(A_{m+1},\higgs_{m+1})$ is the solution of the stochastic Yang--Mills--Higgs equations \eqref{intro:eq-SYMH} with initial data $g_{m+1}(t_m) \act (A,\higgs)(t_m)$ from Theorem \ref{main:thm-covariant-local}. It can be decomposed into 
\begin{equation*}
A_{m+1} = B_{m+1} + \Sdet_{m+1} + \linear[B_{m+1}][][][] + R_{m+1} \qquad \text{and} \qquad 
\higgs_{m+1} = \higgsdet_{m+1} + \philinear[B_{m+1}][][][] + \phir_{m+1}. 
\end{equation*}
\item $B_{m+1}$ is the Yang--Mills gradient flow with initial data $g_{m+1}(t_m) \act A(t_m)$ from Proposition~\ref{main:prop-YMHF-rough-data}. Furthermore, it can be written as $B_{m+1} = g_{m+1} \act (B_{m+1}^+ + B_{m+1}^-)$, where the smooth part $B_{m+1}^+$ is the Yang--Mills gradient flow with initial data $A^+(t_m)$ and the rough part $B_{m+1}^-$ is not too large. 
\item $(\Sdet_{m+1},\higgsdet_{m+1})$ is the solution of \eqref{main:local-eq-Sdet}--\eqref{main:local-eq-higgsdet} with initial data $(0, g_{m+1}(t_m) \act \higgs(t_m))$. The Higgs component can be further decomposed into $\higgsdet_{m+1}=\shiggsdet_{m+1} + \rhiggsdet_{m+1}$, where $\shiggsdet_{m+1}$ solves \eqref{main:local-eq-higgsdet} with initial data $g_{m+1}(t_m) \act \phismoothpart(t_m)$ (and $\Sdet$ replaced by $\Sdet_{m+1}$) and $\rhiggsdet_{m+1}$ is a rough term. 
\item $\linear[B_{m+1}][][][]$ and $\philinear[B_{m+1}][][][]$ are covariant linear stochastic objects, which are defined using $B_{m+1}$.
\item\label{main:item-remainder}$R_{m+1}$ and $\phir_{m+1}$ are nonlinear remainders which are relatively small and more regular; see Proposition~\ref{local:prop-main}.   
\end{enumerate}
The properties in \ref{main:item-gauge}--\ref{main:item-remainder} provide us with detailed information on the structure of our solution. Equipped with this information, we can now decide on the definition of the smooth and rough parts $(\Asmoothpart,\phismoothpart)$ and $(\Aroughpart,\phiroughpart)$. Our initial attempt, which will be revised once in \eqref{main:eq-modified-iteration-scheme} below, is to define 
\begin{equation}\label{main:eq-iteration-smooth-rough}
\begin{aligned}
\Asmoothpart(t) &\eqdef (B_{m+1}+\Sdet_{m+1})(t),  \quad & \quad
\Aroughpart(t) &\eqdef (\linear[B_{m+1}][][][]+R_{m+1})(t), \\ 
\phismoothpart(t) &\eqdef \higgsdet_{m+1}(t), \quad & \quad 
\phiroughpart(t) &\eqdef (\philinear[B_{m+1}][][][]+\phir_{m+1})(t) 
\end{aligned}
\end{equation}
for all $t\in (t_m,t_{m+1}]$. One important aspect of \eqref{main:eq-list-iterates} and \eqref{main:eq-iteration-smooth-rough} is that many of the terms are essentially conditionally independent. To illustrate this in one example, note that the initial data $A(t_m)$ is measurable with respect to $\fF_{t_m}$, i.e., the $\sigma$-algebra generated by the noises $\xi$ and $\zeta$ up to the stopping time $t_m$. Since $B_{m+1}$ is a deterministic Yang--Mills gradient flow with initial data $g_{m+1}(t_m) \act A(t_m)$, this implies that $B_{m+1}$ is essentially\footnote{There is one caveat: The domain of $B_{m+1}$ is given by $[t_m,t_{m+1}]\times \T^2$. Since the stopping time $t_{m+1}$ is not $\fF_{t_m}$-measurable, this means that technically $B_{m+1}$ is not $\fF_{t_m}$-measurable either. However, the stopping time $t_{m+1}$ will be equal to the time $t_{m+1}^\star$ from Definition \ref{defn:time_steps} with high probability, and $t_{m+1}^\star$ will be $\fF_{t_m}$-measurable.} $\fF_{t_m}$-measurable. By the strong Markov property, the noise $\xi(\place+t_m)|_{[0,\infty)}$ is independent of $\fF_{t_m}$, and this implies that $B_{m+1}$ is essentially conditionally independent of $\linear[B_{m+1}][][][]$. 

Using the estimates from our covariant local theory, it is relatively easy to control the rough part $(\Aroughpart,\phiroughpart)(t_{m+1})$ in our covariant Hölder norms.  The main task is then to estimate the monotone quantity 
\begin{equation}\label{main:eq-iteration-monotone-quantity}
\covariantenergy( \Asmoothpart(t_{m+1}) , \phismoothpart(t_{m+1})) 
= \covariantenergy( (B_{m+1}+\Sdet_{m+1})(t_{m+1}), \higgsdet_{m+1}(t_{m+1})), 
\end{equation}
and to hopefully see it decrease from its earlier value  $\covariantenergy( \Asmoothpart(t_{m}) , \phismoothpart(t_{m}))$. This task is difficult primarily because the initial datum of $(B_{m+1},\higgsdet_{m+1})$ not only contains the gauge transformation of the smooth part $(\Asmoothpart,\phismoothpart)(t_m)$, but also the gauge transformation of the rough part $(\Aroughpart,\phiroughpart)(t_m)$. That is, the later smooth part $(\Asmoothpart,\phismoothpart)(t_{m+1})$ contains contributions from the earlier rough part $(\Aroughpart,\phiroughpart)(t_m)$. Of course, it is natural to ask whether this is necessary, but there is a good reason for it: Due to the exponential growth in Theorem \ref{main:thm-pointwise-kernel}, each term in our solution can only be treated perturbatively for one step of our iteration (or perhaps a bounded number of steps). But all of our non-perturbative arguments can only be used for the smooth part, and hence the rough part has to enter into the smooth part after one step. In order to estimate \eqref{main:eq-iteration-monotone-quantity}, it is insufficient to rely only on the covariant Hölder norm of $(\Aroughpart,\phiroughpart)(t_m)$, and we instead need to use information on its structure. After noting that \eqref{main:eq-iteration-smooth-rough} also holds at the earlier step in our iteration, this structure is given by $(\Aroughpart,\phiroughpart)(t_m)= (\linear[B_{m}][][][], \philinear[B_m][][][])(t_m) + (R_m,\phir_m)(t_m)$. After several approximations, which will be detailed in Sections \ref{section:deterministic-rough} and \ref{section:global-theory}, we then obtain the following expansions of the terms in \eqref{main:eq-iteration-monotone-quantity}. The curvature term can be expanded\footnote{Since $B_{m+1}=g_{m+1}\act (\smash{B_{m+1}^{+}}+\smash{B_{m+1}^{-}})$, the main term in $B_{m+1}$ is a gauge transform of $\smash{B_{m+1}^+}$, but not necessarily $\smash{B_{m+1}^+}$ itself. The reader may therefore be surprised that the gauge transform $g_{m+1}$ does not appear more frequently in \eqref{main:eq-curvature-YMGF}--\eqref{main:eq-curvature-quadratic}, but the reason for this is simple: The $L^2$-norm of the curvature is gauge invariant.} as
\begin{align}
\|& F_{B_{m+1}+\Sdet_{m+1}}(t_{m+1}) \|_{L_x^2}^2 \notag \\ 
&\approx \| F_{B_{m+1}^+}(t_{m+1}) \|_{L_x^2}^2 \label{main:eq-curvature-YMGF} \\ 
&+ 2 \int_{\T^2} \Bigl \langle F_{B_{m+1}^+}^{12}(t_{m+1}), \curl_{B_{m+1}^+(t_{m+1})} \bigl( \Ad_{g_{m+1}^{-1}} \Sdet_{m+1}\bigr)(t_{m+1}) \Bigr \rangle_{\frkg} \dx \label{main:eq-curvature-linear-S} \allowdisplaybreaks[4] \\ 
&+ 2 \int_{\T^2} \Bigl \langle F_{B_{m+1}^+}^{12}(t_{m+1}), \curl_{B_{m+1}^+(t_{m+1})} e^{\tau_{m+1} \Delta_{\smash{B_{m+1}^+(t_m)}}^F} \linear[B_m][][][](t_m) \Bigr \rangle_{\frkg} \dx \label{main:eq-curvature-linear-lollipop} \allowdisplaybreaks[4]\\ 
&+ \Bigl \|  \curl_{B_{m+1}^+(t_{m+1})} e^{\tau_{m+1} \Delta_{\smash{B_{m+1}^+(t_m)}}^F} \linear[B_m][][][](t_m) \Bigr\|_{L_x^2}^2  \label{main:eq-curvature-curl}\allowdisplaybreaks[4] \\ 
&+ 2 \int_{\T^2} \Bigl \langle F_{B_{m+1}^+}^{12}(t_{m+1}), \Bigl[ \Bigl( e^{\tau_{m+1} \Delta_{\smash{B_{m+1}^+(t_m)}}^F} \linear[B_m][][][](t_m) \Bigr)^1 , \Bigl( e^{\tau_{m+1} \Delta_{\smash{B_{m+1}^+(t_m)}}^F} \linear[B_m][][][](t_m) \Bigr)^2 \Bigr] \Bigr\rangle_{\frkg} \dx. \label{main:eq-curvature-quadratic}
\end{align}
Here, $\approx$ is used to indicate that only the leading terms are kept, the covariant curl is as in Notation~\ref{notation:covariant-curl}, and the covariant heat propagator $e^{\tau_{m+1} \Delta_{\smash{B_{m+1}^+(t_m)}}^F}$ is as in  Notation \ref{notation:linear-propagation}. Furthermore, the $L_x^4$-term involving the Higgs field can be expanded as
\begin{align}
\|& \higgsdet_{m+1}(t_{m+1}) \|_{L_x^4}^{\gamma_1}  \notag \\
&\approx \| \shiggsdet_{m+1}(t_{m+1}) \|_{L_x^4}^{\gamma_1}  
\label{main:eq-higgs-l4-smooth}  \\ 
&+ \gamma_1 \| \shiggsdet_{m+1}(t_{m+1}) \|_{L_x^4}^{\gamma_1-4} 
\Bigl( \int_{\T^2} \bigl| \shiggsdet_{m+1}(t_{m+1}) \bigr|^2
\bigl\langle \shiggsdet_{m+1}(t_{m+1}), 
\rlinhiggsdetmplusone (t_{m+1}) \bigr\rangle_{\V} \dx  
\Bigr).  \label{main:eq-higgs-l4-first-order}
\end{align}
Here, $\rlinhiggsdetmplusone (t_{m+1})$ is defined as 
$g_{m+1}(t_{m+1}) \act (e^{\tau_{m+1}  \Delta_{\smash{B_{m+1}^+(t_m)}}}  \philinear[B_{m}](t_m))$, where $e^{\tau_{m+1}  \Delta_{\smash{B_{m+1}^+(t_m)}}}$ is as in Notation \ref{notation:linear-propagation}. That is, it is the gauge transformation of a covariant linear evolution of $\philinear[B_{m}](t_m)$. Finally, while the $L_x^r$-term involving the Higgs field can be expanded similarly, we omit this expansion here. The reason is that, unlike the parameter $\gamma_1$, the parameter $\gamma_2$ is smaller than four, and this turns out to make it easier to control the $L_x^r$-term. 

In the following, we discuss the main ideas behind the estimates of \eqref{main:eq-curvature-YMGF}--\eqref{main:eq-higgs-l4-first-order}. 
We first focus on estimates that are sufficient for global well-posedness, but \textit{not} for the uniform-in-time estimates. The necessary ingredients for the latter will only be discussed at the end of this subsection. We start with the $\curl$-term from \eqref{main:eq-curvature-curl}, since it will serve as a benchmark for all remaining terms. To determine the expected size of \eqref{main:eq-curvature-curl}, let us briefly turn back to the non-covariant setting. For the non-covariant linear stochastic object, a direct calculation (in either frequency space or physical space) yields that $\E[ \| \nabla e^{\tau \Delta} \linear[][][][](t)\|_{L^2}^2] \lesssim \tau^{-1}$ for all $t\in [0,1]$ and $\tau \in (0,1]$. We emphasise that there is no $\log$-loss, which only occurs in $L^\infty$-based spaces or at the borderline regularity zero. The probabilistic independencies discussed just below \eqref{main:eq-iteration-smooth-rough} suggest that this second-moment estimate may not only hold in the non-covariant setting, but also in the covariant setting. The best estimate one can hope for is therefore
\begin{equation}\label{main:eq-curl-estimate}
\E \Bigl[ \Bigl \|  \curl_{B_{m+1}^+(t_{m+1})} e^{\tau_{m+1} \Delta_{\smash{B_{m+1}^+(t_m)}}^F} \linear[B_m][][][](t_m) \Bigr\|_{L_x^2}^2 \Bigr] 
\lesssim \E \bigl[ \tau_{m+1}^{-1} \bigr]. 
\end{equation}
The estimate \eqref{main:eq-curl-estimate} turns out to be true, and it is a consequence of Lemma \ref{le:curl_square_expectation} below. However, its proof is more difficult than one may expect. The reason is that we cannot directly use our covariant heat kernel estimates of $\fq_{\smash{B_{m}}}$ and $\fq_{\smash{B_{m+1}^+}}$ from Theorem \ref{main:thm-pointwise-kernel}. Indeed, due to the $\log$-factor on the right-hand side of \eqref{main:eq-L1Linf-curvature} and the $\exp$-factor in Theorem \ref{main:thm-pointwise-kernel}, this would only yield an upper bound of the form $\E[ \tau_{m+1}^{-(1+C\sqrt{\logtimestep})}]$. And, as can be seen from the discussion in the next paragraph, this upper bound cannot prevent finite-time blowup. In order to prove \eqref{main:eq-curl-estimate}, we need to rely on the fact that \eqref{main:eq-curl-estimate} only involves the covariant curl, and not a general covariant derivative. As will be shown in Lemma \ref{lemma:curl-equation} below, the geometric evolution equation for the covariant curl primarily involves the covariant Laplacian and the covariant derivatives of the curvature, but not the curvature term itself as in \eqref{main:eq-linear-stochastic-A}. 
Together with the control provided by the decrease \eqref{main:eq-decrease}, this allows us to treat the covariant derivatives of the curvature perturbatively and replace the pointwise estimates for $\fq_{\smash{B_m}}$ and $\fq_{\smash{B_{m+1}^+}}$ with the diamagnetic inequalities for $\fp_{\smash{B_m}}$ and $\fp_{\smash{B_{m+1}^+}}$, which ultimately leads to the desired estimate \eqref{main:eq-curl-estimate}.

We now discuss why the exponent $-1$ in the factor $\tau_{m+1}^{-1}$ appearing in \eqref{main:eq-curl-estimate} is significant for our global iteration scheme. The discussion of this exponent will be entirely heuristic, and making this discussion more precise will occupy much of Section \ref{section:global-theory}. First, we assume\footnote{In our global iteration scheme, this assumption may fail frequently, but it is still useful for building intuition.} that the time increment $\tau_{m+1}$ is always determined by the first argument in \eqref{main:eq-time-step-1}, i.e., that 
\begin{equation*}
\tau_{m+1} = \logtimestep \| F_{\Asmoothpart}(t_m) \|_{L^2_x}^{-1} 
\log\bigl( \| F_{\Asmoothpart}(t_m) \|_{L^2_x} \vee 2\bigr)
\end{equation*}
for all $m\geq 0$. In that case, after bounding \eqref{main:eq-curvature-YMGF} by the squared $L_x^2$-norm of $F_{B_{m+1}^+}(t_m)$, omitting the contributions of \eqref{main:eq-curvature-linear-S}, \eqref{main:eq-curvature-linear-lollipop}, and \eqref{main:eq-curvature-quadratic}, and omitting all expectations, one obtains an estimate of the form 
\begin{equation}\label{main:eq-exponent-motivation-e1}
\begin{aligned}
\| F_{\Asmoothpart}(t_{m+1}) \|_{L_x^2}^2 
&\leq \| F_{\Asmoothpart}(t_m) \|_{L_x^2}^2 
+ C \tau_{m+1}^{-1} \\ 
&\leq \| F_{\Asmoothpart}(t_m) \|_{L_x^2}^2 
+ C \logtimestep^{-1} \| F_{\Asmoothpart}(t_m) \|_{L_x^2} \log\bigl( \| F_{\Asmoothpart}(t_m) \|_{L_x^2} \vee 2 \bigr)^{-1}.   
\end{aligned}
\end{equation}
By iterating \eqref{main:eq-exponent-motivation-e1} until the time $t=1$ is reached\footnote{Strictly speaking, one should only expect bounds of the form \eqref{main:eq-exponent-motivation-e1} exactly at the times $(t_m)_{m=0}^\infty$, which need not include $t=1$. For this heuristic discussion, however, we ignore this technical problem.}, we then obtain 
\begin{equation}\label{main:eq-exponent-motivation-e2} 
\| F_{\Asmoothpart}(1) \|_{L_x^2}^2 
\leq \| F_{\Asmoothpart}(0) \|_{L_x^2}^2 
+ C \Bigl( \logtimestep^{-1} \| F_{\Asmoothpart}(0) \|_{L_x^2} \log\bigl( \| F_{\Asmoothpart}(0) \|_{L_x^2} \vee 2 \bigr)^{-1} \Bigr)^2.  
\end{equation}
The estimate \eqref{main:eq-exponent-motivation-e2} is sufficient for global well-posedness and, upon further iteration, even yields a growth rate of the form $\exp(Ct^{1/3})$. We emphasise, however, that a worse exponent of $\| F_{\Asmoothpart}(0) \|_{L_x^2}$ in \eqref{main:eq-exponent-motivation-e2} would not allow us to rule out finite-time blowup. As a result, we cannot allow for a worse exponent on the right-hand side of \eqref{main:eq-curl-estimate}, and thus $\tau_{m+1}^{-1}$ is a natural benchmark in our analysis. After this extensive discussion of \eqref{main:eq-curvature-curl}, we now turn to a brief discussion of the remaining terms in \eqref{main:eq-curvature-YMGF}--\eqref{main:eq-higgs-l4-first-order}. In this discussion, we only try to convey the main ideas, but will often neglect various technical difficulties.

\step{Estimate of \eqref{main:eq-curvature-YMGF}:} In order to control the $L^2$-norm of $F_{B^+_{m+1}}(t_{m+1})$, we use the monotonicity formula for the Yang--Mills gradient flow. It implies that 
\begin{equation}\label{main:eq-curvature-YMGF-estimate}
\| F_{B_{m+1}^+}(t_{m+1})\|_{L_x^2}^2 
= \| F_{B_{m+1}^+}(t_{m})\|_{L_x^2}^2  - 2 \| \covd_{B_{m+1}^+}^\ast F_{B_{m+1}^+} \|_{L_t^2 L_x^2([t_m,t_{m+1}])}^2. 
\end{equation}
The first summand in \eqref{main:eq-curvature-YMGF-estimate} appears in our monotone quantity \eqref{main:eq-covariant-energy} at the time $t=t_{m}$. The second summand in \eqref{main:eq-curvature-YMGF-estimate} appears in our decrease \eqref{main:eq-decrease}, and it will be used to absorb several different contributions. 

\step{Estimate of \eqref{main:eq-curvature-linear-S}:} 
By using integration by parts and the $\Ad$-invariance of the inner product on $\frkg$, one obtains that 
\begin{equation*}
\bigl| \eqref{main:eq-curvature-linear-S} \bigr| 
\lesssim \| \covd_{\smash{B_{m+1}^+}}^\ast F_{B_{m+1}^+} (t_{m+1}) \|_{L_x^2} \| \Sdet_{m+1}(t_{m+1}) \|_{L_x^2}. 
\end{equation*}
The first factor can essentially be estimated by $\tau_{m+1}^{-\frac{1}{2}} \| \covd_{\smash{B_{m+1}^+}}^\ast F_{\smash{B_{m+1}^+}} \|_{L_t^2 L_x^2([t_m,t_{m+1}])}$ (see Lemma \ref{le:DF_bound} below), which involves the curvature term in the decrease \eqref{main:eq-decrease}. In the equation for $\Sdet_{m+1}$ (see \eqref{main:local-eq-Sdet}), the largest term is $\bil(\shiggsdet_{m+1} \otimes \covd_{B_{m+1}+\Sdet_{m+1}} \shiggsdet_{m+1})$. Using our covariant Schauder estimates, the resulting contribution to the $L^2$-norm of $\Sdet_{m+1}(t_{m+1})$ can then essentially be bounded by 
\begin{equation*}
 \tau_{m+1}^{\frac{1}{2}} \bigl\| \bil(\shiggsdet_{m+1} \otimes \covd_{B_{m+1}+\Sdet_{m+1}} \shiggsdet_{m+1}) \bigr\|_{L_t^2 L_x^2([t_m,t_{m+1}])}, 
\end{equation*}
which involves the decrease of the $L^4$-term appearing in \eqref{main:eq-decrease}. Provided\footnote{Since we used the word ``essentially'' to hide small losses, we actually need $\gamma_1-4$ to be bounded below by certain small parameters appearing in our argument.} that $\gamma_1>4$, we can then use Young's inequality to absorb the contribution of \eqref{main:eq-curvature-linear-S} in the decrease. 

\step{Estimate of \eqref{main:eq-curvature-linear-lollipop}:}
As in our estimate of \eqref{main:eq-curvature-linear-S}, we first use integration by parts. Then, we use that $B_{m+1}^+(t_{m+1})$ and $\linear[B_m][][][](t_m)$ are essentially\footnote{To be more precise, $B_{m+1}^+(t_{m+1})$ is equal to $B_{m+1}^{+, \star \star}(t_{m+1}^{\star\star})$ from Definition \ref{defn:star_objects} with high probability and, after conditioning on  $\fF_{t_{m-1}}$, $B_{m+1}^{+,\star\star}(t_{m+1}^{\star\star})$ and $\linear[B_m][][][](t_m)$ are probabilistically independent.} probabilistically independent, which allows us to control \eqref{main:eq-curvature-linear-lollipop} using moment and covariant heat kernel estimates, as well as our previous estimates of $\covd_{\smash{B_{m+1}^+}}^\ast F_{\smash{B_{m+1}^+}} (t_{m+1})$.

\step{Estimate of \eqref{main:eq-curvature-quadratic}:} We now turn to the interaction involving $F_{\smash{B_{m+1}^+}} (t_{m+1})$ and two factors of $\linear[B_{m}][][][](t_m)$, which is perhaps the most important term among \eqref{main:eq-curvature-YMGF}--\eqref{main:eq-higgs-l4-first-order}. A naïve estimate yields 
\begin{equation*}
\bigl| \eqref{main:eq-curvature-quadratic} \bigr| 
\lesssim \bigl\|  F_{\smash{B_{m+1}^+}} (t_{m+1}) \bigr\|_{L_x^2} 
\bigl\|  e^{\tau_{m+1} \Delta_{\smash{B_{m+1}^+(t_m)}}^F} \linear[B_m][][][](t_m) \bigr\|_{L_x^4}^2. 
\end{equation*} 
By using the definition of $T^{(1)}$ from \eqref{main:eq-time-step-1}, the first factor can be estimated by $\tau_{m+1}^{-1} \abs{\log(\tau_{m+1})}$. Thus, even if any logarithmic losses from $\linear[B_m][][][](t_m)$ are ignored, the naïve estimate falls short of our $\tau_{m+1}^{-1}$-benchmark. Rather than estimating \eqref{main:eq-curvature-quadratic} directly, we rely on the geometric identity
\begin{equation}\label{main:eq-curvature-quadratic-identity}
\begin{aligned}
2 \int_{\T^2} \bigl \langle& F_{B_{m+1}^+}^{12}(t_{m+1}), \bigl[ X_1, X_2 \bigr] \bigr\rangle_{\frkg} \dx \\ 
&= \bigl\|  \curl_{B_{m+1}^+(t_{m+1})} X \bigr\|_{L_x^2}^2 
+ \bigl\|  \div_{B_{m+1}^+(t_{m+1})} X \bigr\|_{L_x^2}^2  
- \bigl\|  \covd_{B_{m+1}^+(t_{m+1})} X \bigr\|_{L_x^2}^2 
\end{aligned}
\end{equation}
with $X=e^{\tau_{m+1} \Delta_{\smash{B_{m+1}^+(t_m)}}^F} \linear[B_m][][][](t_m)$, which is proven in Lemma \ref{prelim:lem-Hessian-identity} below. The $\curl$-term has already been estimated in \eqref{main:eq-curl-estimate}. In two dimensions, the curl and divergence of the covariant stochastic object have the same law (see Lemma \ref{lemma:linear-object-rotation-invariance}), and hence the $\div$-term satisfies the same estimate as the $\curl$-term. Finally, due to the negative sign, the last summand in \eqref{main:eq-curvature-quadratic-identity} can be bounded from above by zero. All in all, we obtain an upper bound of the expectation of \eqref{main:eq-curvature-quadratic} by $\tau_{m+1}^{-1}$.

Since this will be important in our proof of uniform-in-time estimates, which will be discussed momentarily, we emphasise that the expectation of \eqref{main:eq-curvature-quadratic} cannot be bounded from below by $-\tau_{m+1}^{-1}$. That is, \eqref{main:eq-curvature-quadratic} can be rather negative. 

\step{Estimate of \eqref{main:eq-higgs-l4-smooth}:}
For any $p\geq 2$, the $L^p$-norm of $\shiggsdet_{m+1}$ can be estimated using a covariant monotonicity formula (see Proposition \ref{prop:monotonicity}). The covariant Laplacian in \eqref{main:local-eq-higgsdet} leads to smoothing, and the cubic nonlinearity in \eqref{main:local-eq-higgsdet} leads to decay. Altogether, we obtain
\begin{align*}
\|\shiggsdet_{m+1}&(t_{m+1})\|_{L_x^4}^{\gamma_1}  \\ 
&\leq \| \shiggsdet_{m+1}(t_m)\|_{L_x^4}^{\gamma_1}  
-c\| \shiggsdet_{m+1}(t_m) \|_{L_x^4}^{\gamma_1+2} \tau_{m+1}
\\
&-c \| \shiggsdet_{m+1}(t_{m+1}) \|_{L_x^4}^{\gamma_1-4} \| \abs{\shiggsdet_{m+1}} \abs{\covd_{B_{m+1}+\Sdet_{m+1}} \shiggsdet_{m+1}} \|_{L_{t,x}^2 ([t_m,t_{m+1}])}^2. 
\end{align*}
We emphasise that the $\abs{\shiggsdet_{m+1}} \abs{\covd_{B_{m+1}+\Sdet_{m+1}} \shiggsdet_{m+1}}$-term is the reason for including the $\| \place \|_{L^4}^{\gamma_1}$-term in the monotone quantity \eqref{main:eq-covariant-energy}, and this term was already used in our estimate of \eqref{main:eq-curvature-linear-S}. 

\step{Estimate of \eqref{main:eq-higgs-l4-first-order}:} 
Recall from below \eqref{main:eq-higgs-l4-first-order} that $\rlinhiggsdetmplusone$ is defined using $\philinear[B_m](t_m)$. Since $\philinear[B_m](t_m)$ is essentially\footnote{As in our estimate of \eqref{main:eq-curvature-linear-lollipop}, the probabilistic independence is not actually between $\shiggsdet_{m+1}(t_{m+1})$ and  $\philinear[B_m](t_m)$. For the rigorous details, see Sections \ref{section:deterministic-rough} and \ref{section:global-theory}.} probabilistically independent of $\shiggsdet_{m+1}(t_{m+1})$, we can then control the integral in \eqref{main:eq-higgs-l4-first-order} using moment and covariant heat kernel estimates. The first moment of this contribution is then bounded, up to the preceding approximations, by 
\begin{equation}\label{main:eq-higgs-l4-first-order-treatment}
\| \shiggsdet_{m+1}(t_m) \|_{L_x^4}^{\gamma_1-4} \times \tau_{m}^{\frac{1}{2}} \| (\abs{\shiggsdet_{m+1}}^2 \shiggsdet_{m+1})(t_{m+1}) \|_{L_x^2} 
\lesssim \tau_{m}^{\frac{1}{2}} \| \shiggsdet_{m+1}(t_m) \|_{L_x^4}^{\gamma_1-4} \| \shiggsdet_{m+1}(t_m) \|_{L_x^6}^{3}.
\end{equation}
In passing from $\shiggsdet_{m+1}(t_{m+1})$ to $\shiggsdet_{m+1}(t_{m})$, we also used the covariant monotonicity formula. The reason behind the $\tau_{m}^{1/2}$-factor is that $\philinear[B_m][][][](t_m)$ should behave like a covariant Gaussian free field, and we therefore should see a covariant $H^{-1}$-norm of the cubic nonlinearity. But, due to parabolic scaling, we can gain a factor of $\tau_{m}^{1/2}$ by using the $L^2$-norm. The fact that the gain in \eqref{main:eq-higgs-l4-first-order-treatment} involves the previous time increment $\tau_m$, and not the current time increment $\tau_{m+1}$, makes the final estimate of \eqref{main:eq-higgs-l4-first-order-treatment} rather technical (see Lemma \ref{le:iteration_basic_bounds}). However, the desired upper bound by $\tau_{m+1}^{-1}$ can eventually be deduced using the upper bound on $T^{(1)}$ from \eqref{main:eq-time-step-1} involving the Higgs field. 

Altogether, one obtains an estimate of the schematic form\footnote{Some of the details behind \eqref{main:eq-global-scheme-altogether} will only be shown for the modification of our global iteration scheme indicated in \eqref{main:eq-modified-iteration-scheme}, which will be needed for uniform-in-time estimates. Also, since $\tau_{m+1}$ will not actually be $\fF_{t_m}$-measurable, we will later need to replace $\tau_{m+1}^{-1}$ by $(\tau_{m+1}^{\star})^{-1}$, see Definition~\ref{defn:time_steps}.}
\begin{equation}
\begin{aligned}	
\label{main:eq-global-scheme-altogether}
 \covariantenergy \bigl(\Asmoothpart(t_{m+1}),\phismoothpart(t_{m+1})\bigr) 
\leq  \covariantenergy \bigl(\Asmoothpart(t_{m}),\phismoothpart(t_{m})\bigr)  - \tfrac{1}{4} \Dec_{m+1}  + C  \tau_{m+1}^{-1} + C. 
\end{aligned}
\end{equation}
For the precise formulation of our estimate of $\covariantenergy \bigl(\Asmoothpart(t_{m+1}),\phismoothpart(t_{m+1})\bigr)$, we refer the reader to Section~\ref{section:global-theory}; see Propositions~\ref{pr:expansion_iteration},~\ref{pr:expectation_one_step}, and~\ref{pr:main_iteration}. To see whether the $\tau_{m+1}^{-1}$-term can be absorbed into the $\Dec_{m+1}$-term, we first recall the time scales from Definition~\ref{main:def-time-scales}. The $\stepupperbound$- and $\logtimestep$-terms in \eqref{main:eq-time-step-1} can be ignored, since they can only be the minimum in \eqref{main:eq-time-step-1} when the solution is bounded. When our stopping criterion is triggered by the $\| \phismoothpart(t_m)\|_{L^r}^{-\gamma_0}$-term from \eqref{main:eq-time-step-1}, the $\| \phismoothpart(t_m)\|_{L^r}^{\gamma_2+2}\tau_{m+1}$-term in the decrease will allow us to beat $\tau_{m+1}^{-1}$; see Lemma \ref{le:close_YM_connection}. Furthermore, when our stopping criterion is triggered by the $\| \covd_{\smash{B_{m+1}^+}}^\ast F_{\smash{B_{m+1}^+}} \|_{L^2}$-term from \eqref{main:eq-time-step-2}, the curvature term in the decrease will allow us to beat $\tau_{m+1}^{-1}$; see also Lemma \ref{le:close_YM_connection}. The only case in which the decrease is unable to beat $\tau_{m+1}^{-1}$ is when our stopping criterion is triggered by the $\|F_{\Asmoothpart(t_m)}\|_{L_x^2}$-term from \eqref{main:eq-time-step-1},  and this is for a good reason: If the initial connection one-form $\Asmoothpart(t_m)$ is close to a Yang--Mills connection and $\phismoothpart(t_m)$ is close to zero, then $(\Asmoothpart,\phismoothpart)(t_m)$ is close to a stationary point of the Yang--Mills--Higgs gradient flow. In that case, the decrease from the deterministic dynamics may be too small to absorb the $\tau_{m+1}^{-1}$-term. 
We emphasize that the failure to absorb the $\tau_{m+1}^{-1}$-term into the $\Dec_{m+1}$-term does not rule out a proof of global well-posedness based on \eqref{main:eq-global-scheme-altogether}. Indeed, this obstruction only arises when the time step is determined by the curvature term in \eqref{main:eq-time-step-1}, placing us in the regime discussed around \eqref{main:eq-exponent-motivation-e1}--\eqref{main:eq-exponent-motivation-e2}. However, the obstruction shows that \eqref{main:eq-global-scheme-altogether}, by itself, is insufficient for uniform-in-time estimates.

To obtain uniform-in-time estimates, we need to exploit our decay mechanism near Yang--Mills connections, which was briefly discussed in Section \ref{section:main-ideas-brief}. To motivate this decay mechanism further, we briefly examine a finite-dimensional toy model. In this toy model, $M\in \R^{n\times n}$ is a symmetric matrix with eigenvalues $(\lambda_j)_{j=1}^{n}$, $W$ is an $\R^n$-valued standard Brownian motion, and $X$ is a solution to the SDE
\begin{equation}\label{main:eq-SDE-toy}
\mathrm{d}X= - MX \, \mathrm{d}t + \,  \mathrm{d}W, \qquad X(0)=0. 
\end{equation}
We note that \eqref{main:eq-SDE-toy} is the Langevin equation corresponding to the quadratic energy $\frac{1}{2}\langle M X, X \rangle$. Using It\^{o}'s isometry, one then easily obtains that
\begin{equation}\label{main:eq-hessian-toy}
\E \langle M X(t), X(t) \rangle = \tfrac{1}{2} \sum_{j=1}^{n} (1-e^{-2\lambda_j t}). 
\end{equation}
For any negative eigenvalue $\lambda_j<0$, the summand $1-e^{-2\lambda_j t}$ is negative and tends to $-\infty$. In particular, provided that $M$ is indefinite, \eqref{main:eq-hessian-toy} eventually becomes negative. 

To connect the toy problem back to the Yang--Mills--Higgs model, we note that Yang--Mills connections with large curvature are unstable (see Section \ref{section:perturbing-around-YM-connection}). Based on our analysis of \eqref{main:eq-hessian-toy}, we then hope that the combined contribution of \eqref{main:eq-curvature-curl}--\eqref{main:eq-curvature-quadratic}, which comes from the Hessian of the Yang--Mills energy, is negative. Unfortunately, it is difficult to exhibit our decay mechanism directly in \eqref{main:eq-curvature-curl}--\eqref{main:eq-curvature-quadratic}. To see why, recall that we are currently only thinking of $B_{m+1}^+$ as being close to a Yang--Mills connection, which does not imply the same for $B_m^+$. Since $\linear[B_m][][][]$ is defined using $B_m$, and not $B_{m+1}$, it is then hard to exhibit decay in \eqref{main:eq-curvature-curl}--\eqref{main:eq-curvature-quadratic}. For this reason, we exhibit our decay mechanism by using a part of the stochastic object $\linear[B_{m+1}][][][]$, i.e., a part of the stochastic object on our current time interval. To be more precise, we introduce a modified version of the global iteration scheme discussed above. In the case when the decrease cannot absorb the $\tau_{m+1}^{-1}$-term, and only in this case, the iterates $\Asmoothpart(t_{m+1})$ and $\Aroughpart(t_{m+1})$  are replaced by\footnote{Of course, by replacing the iterate $\Asmoothpart(t_{m+1}) $ with $\Asmoothpartmod(t_{m+1})$, we change all the terms appearing in our global iteration scheme. For this reason, \eqref{main:eq-modified-iteration-scheme} should actually involve the modifications $\tmod_{m+1/2}$, $\tmod_{m+1}$, and $\widetilde{B}_{m+1}$, but we omit this here.}
\begin{equation}\label{main:eq-modified-iteration-scheme}
\begin{aligned}
\Asmoothpartmod(t_{m+1}) 
&\eqdef \Asmoothpart(t_{m+1})
+ e^{(t_{m+1}-t_{m+1/2}) \smash{\Delta_{B_{m+1}(t_{m+1/2})}^F}} \linear[B_{m+1}](t_{m+1/2}), \\ 
\Aroughpartmod(t_{m+1}) &\eqdef \Aroughpart(t_{m+1}) - e^{(t_{m+1}-t_{m+1/2}) \smash{\Delta_{B_{m+1}(t_{m+1/2})}^F}} \linear[B_{m+1}](t_{m+1/2}). 
\end{aligned}
\end{equation}
Here, the intermediate time $t_{m+1/2}$ is defined as the midpoint between $t_m$ and $t_{m+1}$. The reason for using $e^{(t_{m+1}-t_{m+1/2}) \smash{\Delta_{B_{m+1}(t_{m+1/2})}^F}} \linear[B_{m+1}](t_{m+1/2})$, and not just $\linear[B_{m+1}](t_{m+1})$, is that the latter object is not smooth. For a complete and rigorous definition of our modified global iteration scheme, we refer the reader to Section \ref{section:decay-estimates-in-expectation}. Once the modification \eqref{main:eq-modified-iteration-scheme} has been included, the desired decay can be obtained using the following proposition.

\begin{proposition}[Decay mechanism near Yang--Mills connections]\label{main:prop-YM-connection-decay-mechanism}
Let $C_1,C_2>0$ be any positive constants. Furthermore, let the time $T>0$ and the constants $C,c>0$ be sufficiently small depending on $C_1$, $C_2$, and the parameters in \eqref{e:small_constants}. Finally, let $B\colon [0,T]\times \T^2 \rightarrow \frkg^2$ be a smooth Yang--Mills gradient flow satisfying the same assumption as in our covariant local theory, i.e., \eqref{main:eq-covariant-local-condition-1}, as well as the bounds 
\begin{equ}
\label{e:close_to_YM_average}
		\| \covd_B F_B \|_{L^2_{t, x}((0, T])} \leq \frac{C_1}{\sqrt{T}}\;, 
\end{equ}
and 
\begin{equ}
\label{e:decay_mechanism_F_lower_bound}
		\inf_{t \in (0, T]} \|F_B(t) \|_{L^2} = \|F_B(T) \|_{L^2} \geq C_2 \frac{\abs{\log T}}{T} \;. 
\end{equ}
Then, we have that
\begin{equs}
\label{e:decay_mechanism_interval}
		\,&\E  \int_{\T^2}  \Big \langle F_B^{12}(T), \Big [\bigl(e^{(T/2)\Delta^F_{B(T/2)}} \linear[B](T/2) \bigr)^1, \bigl(e^{(T/2)\Delta^F_{B(T/2)}} \linear[B](T/2) \bigr)^2 \Big] \Big \rangle_{\frkg} \dx \\
		&\leq - C \frac{\abs{\log T}^{1 + c }}{T}  \;. 
\end{equs}
\end{proposition}

Our assumptions in \eqref{e:close_to_YM_average} and \eqref{e:decay_mechanism_F_lower_bound} imply that the average size of $\covd_B^\ast F_B$ on the time interval $[0,T]$ is much smaller than the minimal size of $F_B$. Since Yang--Mills connections are solutions of $\covd_B^\ast F_B=0$, one may then interpret \eqref{e:close_to_YM_average} and \eqref{e:decay_mechanism_F_lower_bound} as saying that $B$ is almost a Yang--Mills connection. The idea behind the proof of Proposition \ref{main:prop-YM-connection-decay-mechanism} is to first use Proposition \ref{main:prop-kernel-YM-connection} or, more precisely, Proposition \ref{prop:qB-perturb-near-YM}. Using this, we can approximate the covariant heat kernel $\fq_B$ using the explicit kernel $\fs$ and parallel transport. After this approximation, the estimate \eqref{e:decay_mechanism_interval} can be obtained using a rather explicit calculation; see Section \ref{section:refined_decay_mechanism}. 

As we will see in Section \ref{section:global-theory}, a combination of our earlier estimates of \eqref{main:eq-curvature-YMGF}--\eqref{main:eq-higgs-l4-first-order} and the decay mechanism from Proposition \ref{main:prop-YM-connection-decay-mechanism} can be used to prove a hard global well-posedness result, which is the subject of the following theorem. In its statement, all explicit and implicit constants may depend on the parameters from \eqref{e:small_constants} and \eqref{e:constants_global_iteration}.

\begin{theorem}[Hard global well-posedness and uniform-in-time estimates]\label{thm:global_hard}
Let $A_0\in \Cs_x^{-\kappa}(\T^2, \frkg^2)$ and $\higgs_0\in \Cs_x^{-\kappa}(\T^2, \V)$. Then, there exist connection one-forms $A$ and $\Asmoothpartmod$, a genuine one-form $\Aroughpartmod$, and Higgs fields $\higgs$, $\phismoothpartmod$, and $\phiroughpartmod$, all of which belong to $L_{t,\textup{loc}}^\infty \Cs_x^{-\kappa}([0,\infty)\times \T^2, \genV)$ with $\genV=\frkg^2$ or $\genV=\V$,  and stopping times $(\tmod_m)_{m=0}^\infty$ starting at $\tmod_0=0$ and a.s.\ converging to infinity, which satisfy the following properties:
\begin{enumerate}[label=(\roman*),itemsep=0ex]
\item\label{main:item-global-SYMH} (Stochastic YMH equations) On each time interval $(\tmod_m,\tmod_{m+1}]$, where $m\geq 0$, $(A,\higgs)$ solves \eqref{intro:eq-SYMH}. To be more precise, it is the limit of the smooth, gauge covariant approximations from \eqref{main:eq-SYMH-eps} in $C_t^0 \Cs_x^{-\kappa}((\tmod_m,\tmod_{m+1}]\times \T^2, \frkg^2 \times \V)$.
\item\label{main:item-global-initial} (Initial conditions and right-limits) The initial value $(A(0),\higgs(0))$ is equal to $(A_0,\higgs_0)$. Furthermore, for each $m\geq 0$, the right-limit $(A(\tmod_{m}+),\higgs(\tmod_{m}+))$ exists in $\Cs_x^{-\kappa}(\T^2,\frkg^2\times \V)$. This right-limit belongs to the initial class  $\initialclass_{\tmod_{m+1}-\tmod_{m}}$ defined in \eqref{main:eq-initial-class} and is gauge-equivalent to $(A(\tmod_{m}),\higgs(\tmod_{m}))$, i.e.,  there exists a $g_{m+1}\in \CG_{1-\kappa}$ such that $(A(\tmod_{m}+),\higgs(\tmod_{m}+))=g_{m+1}\act (A(\tmod_{m}),\higgs(\tmod_{m}))$.  
\item\label{main:item-global-decomposition} (Decomposition) For all $m\geq 1$, we have the decompositions $A(\tmod_m)=\Asmoothpartmod(\tmod_m)+\Aroughpartmod(\tmod_m)$ and $\higgs(\tmod_m)=\phismoothpartmod(\tmod_m) + \phiroughpartmod(\tmod_m)$ into smooth and rough parts. 
\item\label{main:item-global-smooth-part} (Uniform-in-time estimates for the smooth part) For all $p\geq 1$, it holds that 
\begin{align*}
\sup_{m\geq 1} \E& \Bigl[\covariantenergy( \Asmoothpartmod(\tmod_{m}),\phismoothpartmod(\tmod_{m}))^p \Bigr] 
\leq C(A_0,\higgs_0,\kappa,p), \\  
\limsup_{m\rightarrow \infty} \E& \Bigl[\covariantenergy( \Asmoothpartmod(\tmod_{m}),\phismoothpartmod(\tmod_{m}))^p \Bigr] 
\leq C(\kappa,p).
\end{align*}
We emphasise that, in contrast to the bound on the supremum, the bound on the limit superior does not depend on the initial data  $(A_0,\higgs_0)$.  
\item\label{main:item-global-rough-part} (Uniform-in-time estimates for the rough part) For all $p\geq 1$, it holds that 
\begin{equation*}
\sup_{m\geq 1} \E \Bigl[ \| \Aroughpartmod(\tmod_m) \|_{\covHolder{\Asmoothpartmod(\tmod_m),F}^{-\kappa}}^p \Bigr] +  \sup_{m\geq 1} \E \Bigl[ \| \phiroughpartmod(\tmod_m) \|_{\covHolder{\Asmoothpartmod(\tmod_m)}^{-\kappa}}^p \Bigr] \leq C(\kappa,p),
\end{equation*}
where $\| \place \|_{\covHolder{\Asmoothpartmod(\tmod_m),F}^{-\kappa}}$ and $\|\place \|_{\covHolder{\Asmoothpartmod(\tmod_m)}^{-\kappa}}$ are as in Definition \ref{defn:covariant_Holder_global}.
\end{enumerate}
\end{theorem}

The properties in \ref{main:item-global-SYMH} and \ref{main:item-global-initial} are exactly as in our soft result and the decompositions in \ref{main:item-global-decomposition} follow directly from our definitions. The main aspects of Theorem \ref{thm:global_hard} are therefore contained in \ref{main:item-global-smooth-part} and \ref{main:item-global-rough-part}, which imply that the smooth and rough parts are uniformly bounded. 

\subsection{Markov process and Yang--Mills--Higgs measure}\label{section:Markov-and-Measure}

Recall that one of the main results of \cite{CCHS22,IlyaHaoYM} was to construct a natural state space 
$\fO_\alpha = \Omega_\alpha^1 / \CG_\alpha$ of
``gauge orbits of distributional connections'' as well as a
 \textit{canonical} Markov process on $\fO$ which had a natural interpretation as the
 ``Langevin dynamics'' for the Yang--Mills measure. Here, $\CG_\alpha$ denotes the group
 of little $\alpha$-Hölder continuous gauge transformations and $\Omega_\alpha^1$ is the
 space of distributional connections introduced in \cite[Sec.~3]{CCHS22}. The exponent $\alpha$ 
 should be thought of as imposing bounds of the form $|\int_\ell A| \lesssim |\ell|^\alpha$
  on one-forms $A \in \Omega_\alpha^1$  for short
 segments $\ell$ (but it also controls how this quantity varies when the line segment
 varies) and can be chosen arbitrarily close to $1$.
 
In view of the results of \cite{CCHS22+} which considers the 3-dimensional case, it is unsurprising 
that an analogous Markov process can also be constructed locally in time in the case of 2D
Yang--Mills--Higgs on the space $\fO_{\kappa}^\YMHs = \Omega_{\kappa}^\YMHs / \CG_{1-\kappa}$ for $\kappa > 0$ small enough.
Here, $\Omega_{\kappa}^\YMHs = \Omega_{1-\kappa}^1 \times \cC^{-\kappa}_x$, with $\Omega_{1-\kappa}^1$ as in \cite{CCHS22}.
The main content of Section~\ref{section:langevin} is to show that this process coincides with the
process constructed in Theorem~\ref{thm:global_soft} (provided of course that one quotients out
the action of the gauge group). In order to formulate this more precisely, we note that 
\eqref{intro:eq-SYMH} has the schematic form
\begin{equ}[e:schematicLimiting]
\d_t Y = (\Delta-1) Y + \linTerm Y + \hat \bil(Y,DY) + \tril(Y,Y,Y) + \boldxi\;,
\end{equ}
for $Y = (A,\Phi)$ and $\boldxi = (\xi,\zeta)$ some distributional $\Y = \fg^2\oplus \V$-valued space-time 
processes and the remaining boldfaced symbols denoting 
suitable multilinear maps. One then constructs $Y$ as the limit of solutions 
to regularised and renormalised equations of the form
\begin{equ}[e:schematicRegularised]
\d_t Y_\eps = (\Delta-1) Y_\eps + (\linTerm + \renormY) Y_\eps + \hat \bil(Y_\eps,DY_\eps) + \tril(Y_\eps,Y_\eps,Y_\eps) + \boldxi_\eps\;,
\end{equ}
for a suitable choice of $\renormY \in \End(\Y)$ and $\boldxi_\eps = \moll\star \boldxi$ with $\moll$
a space-time mollifier which we choose of the form given in  \eqref{e:mollif} below. 
We henceforth say that two renormalisations
$\renormY$ and $\renormYb$ are equivalent if $\lim_{\eps \to 0} |\renormY - \renormYb| = 0$.
Given $Y \in \Omega_{\kappa}^\YMHs$,
we write $[Y] \in \fO_{\kappa}^\YMHs$ for the corresponding gauge orbit. We also write 
$\bar{\fO}_{\kappa}^\YMHs = \fO_{\kappa}^\YMHs \sqcup \{\skull\}$ for the space where we
add a ``point at infinity'' (see \cite{CCHS22} for a natural metric on this space in which
$\skull$ is close to gauge orbits that are ``large''; this can be extended to our setting
incorporating the Higgs field, see Theorem~\ref{theo:Polish} below). 
Combining the results of \cite{CCHS22,CCHS22+} (see also \cite{CS26}), it is by now relatively 
straightforward to show the following convergence result where $\$\place\$_{\kappa}$ is the
norm appearing just after \eqref{e:basicState} below. 

\begin{proposition}\label{prop:limitUsualRegStruc}
There exists a choice of renormalisation $\renormY$ with the following properties.
\begin{enumerate}[label=(\roman*),itemsep=0ex]
\item For every $Y_0 \in \Omega_{\kappa}^\YMHs$ there exists a stopping
time $\tau > 0$ and an $\Omega_{\kappa}^\YMHs$-valued process $Y$ such that, 
writing $Y_\eps$ for the solution to \eqref{e:schematicRegularised} with initial
condition $Y_0$, one has
\begin{equ}
\lim_{\eps \to 0} \sup_{t \in [0,\tau]} \$Y_\eps(t) - Y(t)\$_{\kappa} = 0\;,
\end{equ}
in probability.
\item There exists a Markov process $\bar Y$ with values
in $\bar{\fO}_{\kappa}^\YMHs$ such that, for every $Y_0$ as above, one has
$\bar Y(t) = [Y(t)]$ for all $t \in [0,\tau]$, provided that $\bar Y(0) = [Y_0]$.
\end{enumerate}
\end{proposition}

\begin{remark}
Note that all of these stochastic processes are defined on the same probability space
and they are adapted to the natural filtration generated by the noise $\boldxi$ appearing
in \eqref{e:schematicLimiting}. (In the case of \eqref{e:schematicRegularised} this statement
requires the regularisation $\boldxi_\eps$ to be non-anticipative, but this is not actually required 
for Proposition~\ref{prop:limitUsualRegStruc} to hold.)
\end{remark}

\begin{remark}
While the gauge orbit  process $\bar Y$ may blow up in finite time (in which case
it remains at $\skull$), this statement is also compatible with a situation where
$\bar Y$ never blows up but $Y$ blows up in finite time. This is a consequence of the fact that 
gauge orbits are unbounded.
\end{remark}

\begin{remark}
The fact that the process $Y$ obtained in this way is Markovian is by now standard
and remains true for \textit{any} other renormalisation $\renormYb$ such that
$\lim_{\eps \to 0}(\renormY-\renormYb)$ exists and is finite. The fact that the 
quotient process $\bar Y$ is also Markovian however is highly non-trivial and this
highly constrains $\renormY$.

In fact, within the class of renormalisation constants that naturally arise in this
context, for any other choice of renormalisation $\renormYb$ satisfying the properties
above, there exists a $G$-invariant symmetric
element $M \in \End(\V)$ such that $\renormY - \renormYb$ is equivalent to $0 \oplus M$,
under the splitting $\Y = \fg^2 \oplus \V$. In other words, the choice of renormalisation can be 
absorbed into the choice of $\mass$ in Setting~\ref{intro:assumption-setup},
see for example \cite{CCHS22+,IlyaHaoYM}.
\end{remark}

Consider now the process $Z = (A,\Phi)$ constructed in Theorem~\ref{thm:global_soft}.
We then have the following result.

\begin{proposition}\label{main:prop-gauge-equivalence}
The process $Z$ almost surely takes values in $\Omega_\kappa^\YMHs$,
so that $\bar Z(t) = [Z(t)]$ is well defined.
Furthermore, it is possible to realise both $Y$ and $Z$ on the same probability space
(but driven by different white noises) in such a way that, for every $t > 0$, 
$\bar Y(t) = \bar Z(t)$ almost surely.
In particular, we have the following two consequences
\begin{enumerate}[label=(\roman*),itemsep=0ex]
\item The process $\bar Y$ satisfies $\bar Y(t) \neq \skull$ almost surely.
\item The process $\bar Z$ is Markovian.
\end{enumerate}
\end{proposition}

\begin{remark}
In view of Lemma~\ref{main:lem-gauge-covariance}, we know that the process 
$Z$ is gauge covariant in the sense that 
the solutions starting from gauge equivalent starting points are gauge equivalent (modulo
a change of space-time white noise). However, it is certainly not Markovian. In particular,
it is not obvious at all that there is any choice of renormalisation that makes the process
$\bar Z$ Markovian, and this is where we need the comparison to the construction of
\cite{CCHS22,CCHS22+} provided by this proposition.
\end{remark}

Results of this type (symmetry properties of solutions to singular SPDEs) have recently
been obtained in a number of different contexts and their proofs fall into two categories.
The first consists of ``hard'' estimates where one obtains explicit expressions for
the BPHZ renormalisation constants (see \cite{BHZ19,Analyst} for this terminology) and exploits the 
fact that the BPHZ renormalisations of different approximations to the same equation
converge to the same limit. See for example \cite{CCHS22,Quasilinear}. 
The second consists of ``soft'' 
arguments which are more conceptual in nature and were first introduced in \cite{BGHZ22}, but then
also used in the context of Yang--Mills--Higgs in \cite{CCHS22+} where explicit calculations
would have been prohibitive.

In our case, Proposition~\ref{prop:limitUsualRegStruc} could in principle be obtained by the 
same ``soft'' argument used in \cite{CCHS22+}, but it is less clear how to set
up such a soft argument to prove Proposition~\ref{main:prop-gauge-equivalence}. Therefore,
we decided to go down the ``hard'' route and we obtain in Section~\ref{sec:renConstants} precise (up to equivalence)
expressions for all the renormalisation constants appearing in the equation.

This then allows us to compare three different processes, as described in Section~\ref{sec:threeProcesses}.
The first one is simply the process $Y_\eps^{(1)}$ obtained by solving \eqref{e:schematicRegularised} with some 
initial data $Y_0$. However, the driving noise $\boldxi_\eps = \moll \star \boldxi$ is replaced 
by the noise $\moll \star (g \act \boldxi)$, for some adapted $\liegroup$-valued space-time 
process $g_t$ to be described a bit later.

\begin{remark}\label{rem:indep}
The crucial observation here is that the law of the process 
$Y_\eps^{(1)}$ does not depend on the choice of $g$ as a consequence 
of Itô's isometry.
\end{remark}

Regarding the second process, we start from the solution $Y_\eps^{B}$ to
\eqref{main:eq-SYMH-eps} which includes the term coming from the 
``infinitesimal DeTurck gauge'' with background connection $B$ as opposed to the 
regular DeTurck gauge used for the process $Y_\eps^{(1)}$. We then
define a gauge-transformation-valued process $g_t^{(2)}$ by solving the PDE which,
in the absence of noise and renormalisation, would generate the gauge transformation
that switches between the regular and infinitesimal DeTurck gauges, and we 
set $Y_\eps^{(2)}(t) = g_t^{(2)} \act Y_\eps^{B}(t)$.
If it weren't for the presence of renormalisation and the fact that, in
the definition of $Y_\eps^{(1)}$, we act on the noise $\boldxi$ with $g$
\textit{before} hitting it with the mollifier, the equations solved by 
$Y_\eps^{(1)}$ and $Y_\eps^{(2)}$ would be identical (provided that the process $g$
appearing in $Y_\eps^{(1)}$ is chosen in the right way).

To construct the third process, we consider the process $Y_\eps^{g_0}$
obtained by starting \eqref{e:schematicRegularised} (this time driven by the original noise $\boldxi_\eps$) 
with gauge-transformed initial data $g_0 \act Y_0$. 
We then define a time-dependent gauge transformation $g_t^{(3)}$ by solving 
a suitable parabolic PDE (with right-hand side determined by the term added to 
the YM gradient flow by the DeTurck trick) which constructs the gauge transform
turning the DeTurck term for the original solution starting at $Y_0$ into that
for the solution starting at $g_0 \act Y_0$ (again, this would be the case if it 
weren't for the presence of renormalisation). The process $Y_\eps^{(3)}$ is
then defined by setting $Y_\eps^{(3)}(t) = g_t^{(3)} \act Y_\eps^{g_0}(t)$.

At fixed $\eps > 0$, all three processes $Y_\eps^{(i)}$ solve very similar
equations, the difference being that they are driven by slightly different noises
and with slightly different renormalisation terms.
However, their limits
in probability $Y^{(i)} = \lim_{\eps \to 0} Y_\eps^{(i)}$ (ignoring for the moment
technicalities due to the fact that these may in principle blow up in finite time)
turn out to be such that $Y^{(1)} = Y^{(2)}$ for some choice of the process $g$
and a suitable choice of renormalisation constants.
Since $Y \eqlaw Y^{(1)}$ and $Z = (g_t^{(2)})^{-1} \act Y^{(2)}$, this does indeed
provide the claimed coupling between $\bar Y$ and $\bar Z$.

As a reality check, we verify that the \textit{same} choice of renormalisation 
constants is such that one also has $Y^{(1)} = Y^{(3)}$ for some (different in general) 
choice of $g$. Since the law of $Y^{(1)}$ is independent of that choice of $g$,
one concludes as in \cite{CCHS22} that there is a coupling between the process $Y = Y^{(1)}$
with initial condition $Y_0$ and a gauge-transformed version of the same process
with initial condition $g_0 \act Y_0$, thus showing that the 
gauge orbit process $[Y]$ is also Markovian. Since we know from
Theorem \ref{thm:global_soft} that the process $Y^{(2)}$ is global in time (modulo 
regauging it at a suitable
sequence of stopping times which has no effect on its projection onto
the space of gauge orbits), the claim follows. This leads us to our main result,
where we write again $\$\place\$_{\kappa}$ for the natural norm on
$\Omega_{1-\kappa}^1 \times \cC^{-\kappa}_x$.

\begin{theorem}[Yang--Mills--Higgs measure]
The Markov process $[Y]$ admits a unique invariant measure $\mu_{\YMH}$ 
on $\fO_{\kappa}^\YMHs$. Furthermore, under $\mu_{\YMH}$, 
$\inf_{\bar Y \in [Y]} \$\bar Y\$_{\kappa}$ admits finite moments 
of all orders.
\end{theorem}

Existence of $\mu_{\YMH}$ follows from the standard Krylov--Bogoliubov 
argument, combined with the uniform bounds of Theorem \ref{thm:global_soft}. This also
gives the moments of all orders, so it remains to show uniqueness. Since $\fO_{\kappa}^\YMHs$ is Polish,  
this follows from the standard theory of Markov chains
once we show that the process has the strong Feller property and is topologically irreducible \cite[Prop.~4.1.1\ \&\ Thm~4.2.1]{DPZ2}.
The strong Feller property follows in a way that is virtually identical to \cite[Prop.~10.1]{IlyaHaoYM}. Regarding topological
irreducibility, it would be tempting to try to apply the general result of \cite{SupportGeneral}. Unfortunately,
logarithmic divergences play a special role in that work and require us to check a highly non-trivial
condition ``by hand''. Instead, we develop in Section~\ref{sec:irreducible} a much simpler proof technique based on 
the Cameron--Martin theorem that we hope might be useful in the analysis of a much broader class of singular SPDEs. More precise details will be given in 
Corollary~\ref{cor:main} below.

\section{Preliminaries}\label{section:preliminaries}

We first let $\chi\colon [0,\infty)\ra [0,1]$ be a smooth function such that $\operatorname{supp} \chi \subseteq [1/16,16]$ and such that $\sum_{N\in \dyadic} \chi(N^2 t) =1 $ for all $t\in (0,1]$. For all $N\in \dyadic$, we then define 
\begin{equation}\label{prelim:eq-chi-N}
\chi_N(t)\eqdef \chi(N^2 t). 
\end{equation}
We emphasise that the function $\chi$ here and the function $\rho$ from \eqref{main:eq-mcP} can be chosen completely separately. 

We now define parameters $\lossdel,\reg_1,\reg_2,\theta,\theta^\prime,\kappa,\lambda,\lambda_1,\lambda_2\in (0,1)$ and $r\in [1,\infty)$. These parameters primarily, but not exclusively, appear in our covariant local theory. At an informal level, they are chosen such that 
\begin{equation}
\begin{aligned}
\label{e:small_constants}
0< \lossdel \ll \frac{1}{r} \ll \kappa \ll \reg_1 = \reg_2 \ll \theta\ll \theta^\prime \ll  \lambda \ll \lambda_1 \ll \lambda_2 \ll 1.  
\end{aligned}
\end{equation}
The reason for introducing the two parameters $\reg_1$ and $\reg_2$, despite their values being equal, is that they play different roles in Section~\ref{section:covariant_lwp}. The parameter $\reg_1$ will appear in the regularity of the paracontrolled term, and the parameter $\reg_2$ will appear in the regularity of the smooth remainder. To be more precise, we first choose a sufficiently small constant $\varkappa>0$. Then, we define
\begin{equation}\label{prelim:eq-parameters}
\begin{alignedat}{19}
\lambda_2 &\coloneq \varkappa, \quad &
\lambda_1 &\coloneq \varkappa^2, \quad &
\lambda &\coloneq \varkappa^4, \quad &
\theta^\prime&\coloneq \varkappa^8, \quad &  
\theta &\coloneq \varkappa^{16}, \\ 
\reg_1&\coloneq\varkappa^{32}, \quad &
\reg_2&\coloneq\varkappa^{32}, \quad  & 
\kappa&\coloneq \varkappa^{64}, \quad &
r&\coloneq \varkappa^{-128}, \quad &
\text{and} \quad
\lossdel&\coloneq \varkappa^{256}. 
\end{alignedat}
\end{equation}
Note that, as a consequence of our choices, we also have that 
\begin{equation*}
 \frac{\reg_1}{\theta}=\frac{\reg_2}{\theta} \ll \lambda \qquad \text{and} \qquad \lambda \theta^\prime \gg \theta. 
\end{equation*}
For future use, we also choose 
\begin{equation}\label{prelim:eq-T0}
T_0 = T_0(\varkappa) 
\end{equation}
as sufficiently small depending on the parameters appearing in \eqref{prelim:eq-parameters} or, equivalently, depending on $\varkappa$.  Finally, we define the parameters $\fa,\gamma_0,\gamma_1,\gamma_2$, and $\logtimestep$. These parameters appear primarily, but not exclusively, in our covariant global theory. They are chosen such that 
\begin{equation}\label{e:constants_global_iteration}
\logtimestep \ll  \lossdel, \quad \lambda_2 \ll \fa \ll \gamma_1 - 4 \ll \gamma_0 - 2 \ll 1, \quad \text{and} \quad \gamma_2=3.
\end{equation}
Once all these constants are chosen, choose some $\stepupperbound$ small enough so that it satisfies the requirements of all the lemmas and propositions in this paper. 

\subsection{Notations}\label{section:notations}

In this section, we set notation and conventions that will be used throughout this paper. We use $\genV$ to denote a finite-dimensional vector space endowed with an orthogonal representation of $\liegroup$. When specifying multiple such spaces, we may write $\genV_1, \genV_2$, etc.

\begin{notation}[Dyadic integers]\label{notation:dyadic-integers}
We will often use the capital letters $K$, $L$, $M$, and $N$ to denote dyadic integers, i.e., elements of $\dyadic$. Whenever it is clear from context that $N$ is a dyadic integer, we will omit the restriction to dyadic integers from our notation. For example, instead of $\sum_{N\in \dyadic}$, we will often simply write $\sum_{N}$. 
\end{notation}

\begin{notation}[Space-time indices]\label{notation:space-time-indices}
We use Greek letters such as $\alpha, \beta, \gamma, \mu, \nu$ to denote space-time indices, i.e.\ indices with values in $\{0\} \cup [d] = \{0, 1, \ldots, d\}$, and we use Latin letters such as $i, j, k, \ell \in [d]$ to denote spatial indices. In accordance with Einstein's summation convention, whenever Greek letters are contracted, the sum is over $\{0\} \cup [d]$, while whenever Latin letters are contracted, the sum is just over $[d]$.
\end{notation}

\begin{notation}[Higher derivatives]\label{notation:higher-derivatives}
For $k \geq 0$, we write $\covd_A^{(k)} \Phi$ to denote the collection of all covariant derivatives of the form $\covd_{A}^{i_1} \cdots \covd_{A}^{i_k} \Phi$, where $i_1, \ldots, i_k \in [d]$. Given a norm $\|\place\|$, we write
\begin{equ}
\|\covd_A^{(k)} \Phi\| = \sum_{i_1, \ldots, i_k =1}^d \|\covd_{A}^{i_1} \cdots \covd_{A}^{i_k} \Phi\|,
\end{equ}
i.e.\ $\|\covd_A^{(k)} \Phi\|$ is the sum of the norms of all possible $k$th order derivatives.
\end{notation}

\begin{notation}[Space-time derivatives]\label{notation:space-time-derivatives}
Given a multi-index $\mbf{m} \in \N_0^{\{0\} \cup [d]}$, we define the usual parabolic scaling
\begin{equ}
\stscaling{\mbf{m}} \eqdef 2 \mbf{m}_0 + \sum_{i=1}^d \mbf{m}_i.
\end{equ}
For $m \geq 0$, we let $\covd^{(m)}_{A, t, x}$ denote the collection of all covariant derivatives of the form $\covd_A^{\alpha_1} \cdots \covd_A^{\alpha_k}$ for some $k \geq 0$, $\alpha_1, \ldots, \alpha_k \in \{0\} \cup [d]$ such that the parabolic scaling is $m$. That is, letting $\mbf{m}_\alpha \in \N^{\{0\} \cup [d]}$ denote the multi-index obtained by counting the total number of appearances of the indices $0, 1, \ldots, d$ among the $\alpha_1, \ldots, \alpha_k$, we restrict to the collection of $\alpha_1, \ldots, \alpha_k \in \{0\} \cup [d]$ such that $\stscaling{\mbf{m}_\alpha} = m$. Here, we specify that $\covd_A^0 = \ptl_t$ denotes differentiation in time. Similar to Notation~\ref{notation:higher-derivatives}, given a norm $\|\place\|$, we let
\begin{equ}
\|\covd^{(m)}_{A, t, x} \Phi \| = \sum_{\substack{\alpha_1, \ldots, \alpha_k \\ \stscaling{\mbf{m}_\alpha} = m}} \|\covd_A^{\alpha_1} \cdots \covd_A^{\alpha_k} \Phi\|.
\end{equ}
\end{notation}

\begin{notation}[Wick products]
If $X$ is a $k$-linear expression of the Gaussian driving noises $\xi$, $\zeta$,
we write $\Wick{X}$ for its component in the $k$th Wiener chaos.
\end{notation}

\begin{notation}[Standard and schematic Landau symbols]\label{prelim:notation-Landau-standard-schematic}
We write $\bigO(\place)$ for the standard Landau symbol. To be precise, we write $X=\bigO(Y)$ if there exists an absolute constant $C>0$ such that $|X|\leq C|Y|$. We also write $\schO(\place)$ for a schematic Landau symbol, which is used here for purposes similar to those in \cite{Oh2014}. For elements $X_1,\hdots,X_n$ of $\frkg$ or $\V$ from Setting~\ref{intro:assumption-setup}, or a more general inner product space $\genV$ equipped with an orthogonal representation of $\liegroup$, we write $\schO(X_1 \hdots X_n)$ for any multilinear expression built out of our representations of $\liegroup$. For instance, the three product terms in \eqref{intro:eq-linearized} can be written as 
\begin{equation*}
\big[ F_B^{kj},C_j \big]=\schO(F_B C), \quad 
\big[ C_j , 2 \covd_B^j C^k - \covd_B^k C^j \big]=\schO(C \covd_B C), \quad \text{and} \quad 
\big[ C_j, \big[ C^j , C^k \big]\big]=\schO(CCC). 
\end{equation*}
The two most important aspects of the $\schO(\place)$-notation are the pointwise estimate $|\schO(X_1 \hdots X_n)|=\bigO(|X_1|\hdots |X_n|)$ and the covariant product rule
\begin{equation*}
\covd_{B,j} \schO(X_1 \hdots X_n) = \sum_{m=1}^{n} \schO(X_1 \hdots \covd_{B,j} X_m \hdots X_n), 
\end{equation*}
which will be further discussed in Lemma~\ref{prelim:lem-equivariant-product-rule} below. 
\end{notation}

\begin{notation}[Euclidean heat kernels]\label{notation:euclidean-heat-kernels}
We let $\fp^{\Euc}, \fp$ denote the standard heat kernels on $\R^d, \T^d$. Explicitly, we have that
\begin{equs}
\fp^{\Euc}(s, y; t, x) &= \fp^{\Euc}_{t-s}(x-y) = \frac{1}{(4\pi (t-s))^{\frac{d}{2}}} \exp\Bigl(-\frac{|x-y|^2}{4(t-s)}\Bigr), \\
\fp(s, y; t, x) &= \fp_{t-s}(x-y) = \sum_{n \in \Z^d} \fp^{\Euc}_{t-s}(x-y+2\pi n).
\end{equs}
We refer to both $\fp^{\Euc}$ and $\fp$ as Euclidean heat kernels.
\end{notation}

\begin{notation}[Limits of integration]
In order to simplify the notation, we will omit the limits of integration whenever they are clear from context. For example, we will often write 
\begin{equation*}
 \iint \ds \dy \rho_\eps(t-s) \fq_B(s,y;t,x) \xi(s,y) 
 \quad \text{instead of} \quad \int_0^t \ds \int_{\T^2} \dy \rho_\eps(t-s) \fq_B(s,y;t,x) \xi(s,y).
\end{equation*}
More generally, all time variables will be nonnegative, and the second time variable in any heat kernel will always be greater than or equal to the first time variable. 
\end{notation}

\begin{notation}[Norms and inner products of differential forms]\label{notation:norms-forms}
For $\genV$-valued one-forms $C, D \colon \T^d \ra \genV^d$, let
\begin{equs}
\langle C, D \rangle_{\genV} &\eqdef \delta^{ij} \langle C_i, D_j \rangle_{\genV} = \langle C_i, D^i \rangle_{\genV},  \\
|C|_{\genV} &\eqdef \langle C, C \rangle_{\genV}^{1/2}.
\end{equs}
Similarly, given $\genV$-valued two-forms $Q, R \colon \T^d \ra \genV^{d \times d}$, let
\begin{equs}
\langle Q, R \rangle_{\genV} &\eqdef \frac{1}{2} \delta^{ik} \delta^{j \ell} \langle Q_{ij}, R_{k\ell} \rangle_{\genV} = \frac{1}{2} \langle Q_{ij}, R^{ij} \rangle_{\genV}, \\
|Q|_{\genV} &\eqdef \langle Q, Q \rangle_{\genV}^{1/2}.
\end{equs}
We will often omit the subscript $\genV$ from the norms and write $|C|, |Q|$. On the other hand, we will not do this for the inner product, because we now define the $L_x^2$ inner products
\begin{equs}
\langle C, D \rangle &\eqdef \int_{\T^d} \langle C(x), D(x) \rangle_{\genV} \, \dx ,\\
\langle Q, R \rangle &\eqdef \int_{\T^d} \langle Q(x), R(x) \rangle_{\genV} \, \dx. 
\end{equs}
For $\genV$-valued functions $\Phi, \Psi \colon \T^d \ra \genV$, we also write $\langle \Phi, \Psi \rangle$ to denote their $L_x^2$ inner product. To emphasise that this is the $L^2$ inner product, we may also write $\langle \place, \place \rangle_{L^2}$ or $\langle \place, \place \rangle_{L_x^2}$.
\end{notation}

\begin{notation}[Covariant curl and divergence of a one-form]\label{notation:covariant-curl}
Let $A \colon \T^2 \ra \frkg^2$ be a connection one-form and $C \colon \T^2 \ra \genV^2$ be a $\genV$-valued one-form. Define the covariant curl $\curl_A C \colon \T^2 \ra \genV$ as the function:
\begin{equ}
\curl_A C \eqdef \covd_{A, 1} C_2 - \covd_{A, 2} C_1 .
\end{equ}
We remark that this is (up to a possible sign) the Hodge-$\star$ of the covariant exterior differential of $C$. We also define the covariant divergence $\div_A C \colon \T^2 \ra \genV$ as the function:
\begin{equ}
\div_A C \eqdef \covd_A^\ell C_\ell .
\end{equ}
\end{notation}

\begin{notation}[Covariant divergence of a two-form]\label{notation:covariant-divergence-two-form}
Let $A \colon \T^d \ra \frkg^d$ be a connection one-form and $F \colon \T^d \ra \genV^{d \times d}$ be a $\genV$-valued two-form. We define the covariant divergence $\covd_A^* F \colon \T^d \ra \genV^d$ as the one-form
\begin{equ}
(\covd_A^* F)_i \eqdef - \covd_A^j F_{ji}, ~~ i \in [d]. 
\end{equ}
\end{notation}

\begin{notation}[Covariant Hessian]\label{notation:covariant-hessian}
Let $A \colon \T^d \ra \frkg^d$ be a connection one-form and $\Phi \colon \T^d \ra \genV$ be a $\genV$-valued function. Given vector fields $X, Y$ on $\T^d$, we let the covariant derivative of $\Phi$ in the direction $X$ be
\begin{equ}
\covd_{A, X} \Phi \eqdef X^i \covd_{A, i} \Phi,
\end{equ}
and using this, we define the covariant Hessian
\begin{equ}
(\Hess_A \Phi)(X, Y) \eqdef \covd_{A, X} \covd_{A, Y} \Phi - \covd_{A, \ptl_X Y} \Phi.
\end{equ}
Here, $\ptl_X Y$ refers to the usual directional derivative of a vector field.
\end{notation}

\subsection{Geometric identities}\label{section:geometric-identities}

In this section, we state various geometric identities which will be used throughout the paper. These identities are all quite classical and perhaps do not require proofs. However, we still prefer to give some of the slightly more involved proofs, since we do not expect the reader to be intimately familiar with the various standard manipulations of geometric analysis. Throughout this section, suppose that all functions which appear are smooth, or more generally, have enough derivatives so that all expressions which appear are classically well-defined.

\subsubsection{Covariant derivatives, product formulas, and commutator identities}

The following definitions and results will be used to, at least in parts of our argument, treat the connection one-form and Higgs field on equal footing.

\begin{definition}[\protect{$\liegroup$-equivariant maps}]\label{prelim:def-equivariant}
Let $\Repr_1$ and $\Repr_2$ be two representations of $\liegroup$ on vector spaces $\genV$ and~$\boldZ$, respectively. Furthermore, let $\boldU\subseteq \genV$ be an open set that is invariant under $\Repr_1$, i.e., that satisfies $\Repr_1(g) u \in \boldU$ for all $g\in \liegroup$ and $u\in \boldU$. Then, a map $\Qc\colon \boldU \rightarrow \boldZ$ is called $\liegroup$-equivariant if, for all $g\in \liegroup$ and $u\in \boldU$, 
\begin{equation*}
\Qc(\Repr_1(g)u)= \Repr_2(g) \Qc(u).
\end{equation*}
\end{definition}

\begin{remark}
The reason for allowing $\Qc$ to be defined only on an open set $\boldU \subseteq \genV$, rather than all of $\genV$, is that we will later need $\liegroup$-equivariant maps $\Qc$ defined on the general linear group of a vector space, which is an open subset of the space of endomorphisms. 
\end{remark}

\begin{lemma}\label{prelim:lem-equivariant}
The Lie bracket $[\place,\place]\colon \frkg \times \frkg \rightarrow \frkg$, the linearised representation $\rep_\ast \colon \frkg \times \V \rightarrow \V$, and the linear map $\bil\colon \V \otimes \V \rightarrow \frkg$ are $\liegroup$-equivariant.
\end{lemma}

\begin{proof}
Since the action of $\liegroup$ on $\frkg$ is the adjoint action, the $\liegroup$-equivariance of the Lie bracket follows\footnote{It can also be shown using a similar argument as for $(X,v)\mapsto \rep_\ast(X)v$ below.} from the fact that each $\Ad_g$ is a Lie algebra 
automorphism. Regarding $\liegroup$-equivariance of 
$(X,v)\mapsto \rep_\ast(X)v$, we see that for all $g\in \liegroup$, $X\in \frkg$, and $v\in \V$, it holds that 
\begin{equation}\label{prelim:eq-equivariant-p1}
\begin{aligned}
\rep(g) \rep_\ast(X) v 
&= \rep(g) \frac{\mathrm{d}}{\mathrm{d}t}\bigl( \rep(e^{tX})v \bigr)\big|_{t=0}
= \frac{\mathrm{d}}{\mathrm{d}t}\bigl( \rep(g e^{tX})v \bigr)\big|_{t=0} \\
&= \frac{\mathrm{d}}{\mathrm{d}t}\bigl( \rep( g e^{tX} g^{-1}) \rep(g) v \bigr)\big|_{t=0} 
= \rep_{\ast}(\Ad_g X) \rep(g) v\;,
\end{aligned}
\end{equation}
as claimed. To show the $\liegroup$-equivariance of $\bil$, let $g\in \liegroup$, $X\in\frkg$, and $u,v\in \V$. Using the $\Ad$-invariance of the inner product on $\frkg$, the definition of $\bil$, and the orthogonality
of $\rep$, we have
\begin{equs}
\big\langle \Ad_g \bil(u\otimes v), X \big\rangle_\frkg 
&= \big\langle \bil(u\otimes v), \Ad_{g^{-1}} X \big\rangle_\frkg
= \big\langle u, \rep_\ast (\Ad_{g^{-1}}X) v \big\rangle_{\V}\\
&= \big\langle \rep(g)u, \rep(g) \rep_\ast (\Ad_{g^{-1}}X) v \big\rangle_{\V} 
=  \big\langle \rep(g)u,  \rep_\ast (X) \rep(g) v \big\rangle_{\V}\\
&= \big\langle \bil(\rep(g)u\,\otimes \, \rep(g)v), X \big\rangle_{\frkg}\;.
\end{equs}
Since $X$ was arbitrary, it follows that $\Ad_g \bil(u\otimes v)=\bil(\rep(g)u\otimes \rep(g)v)$ as required.
\end{proof}

\begin{remark}\label{remark:inner-product-G-equivariant}
Because we always assume that the vector space $\genV$ is equipped with an inner product $\langle \place, \place \rangle_\genV$ and an orthogonal action $\Repr$ of $\liegroup$ on $\genV$, we may also view the inner product itself as a $\liegroup$-equivariant map $\genV \times \genV \ra \R$. Here, the action of $\liegroup$ on $\R$ is taken to be the trivial action. 
\end{remark}

In the following lemma, we show that $\liegroup$-equivariant maps satisfy chain and product rules. 

\begin{lemma}\label{prelim:lem-equivariant-product-rule}
Let $A\colon \T^d \rightarrow \frkg^d$ be a connection one-form. Let $k\geq 1$ and, for each $1\leq j \leq k$, let $\Repr_j$ be a representation of $\liegroup$ on a vector space $\genV_j$, let $\boldU_j \subseteq \genV_j$ be an open set that is invariant under $\Repr_j$, and let $u_j \colon \T^d \rightarrow \boldU_j$ be a map. Furthermore, let $\Repr$ be a representation of $\liegroup$ on a vector space $\boldZ$. Finally, let $\Qc\colon \boldU_1 \times \hdots \times \boldU_k \rightarrow \boldZ$ be a differentiable $\liegroup$-equivariant map. Then, it holds that 
\begin{equation}\label{prelim:eq-equivariant-chain-rule}
\covd_A \bigl( \Qc(u_1,\hdots, u_k) \bigr)  = \sum_{j=1}^{k}
 \frac{\partial \Qc}{\partial u_j}(u_1,\hdots,u_k) \covd_A u_j. 
\end{equation}
In particular, if $\Qc\colon \genV_1 \times \hdots \times \genV_k \rightarrow \boldZ$ is a multilinear, $\liegroup$-equivariant map and $w_j\colon \T^d \rightarrow \genV_j$, where $1\leq j \leq k$, then 
\begin{equation}\label{prelim:eq-equivariant-product-rule}
\covd_A \bigl( \Qc (w_1,\hdots,w_k) \bigr) = \sum_{j=1}^k \Qc(w_1,\hdots, \covd_A w_j, \hdots, w_k). 
\end{equation}
\end{lemma}

We note that, with a slight abuse of notation, we denote the different covariant derivatives corresponding to the different representations of $\liegroup$ all by $\covd_A$.  

\begin{proof}
Since \eqref{prelim:eq-equivariant-product-rule} is a special case of \eqref{prelim:eq-equivariant-chain-rule}, it suffices to prove \eqref{prelim:eq-equivariant-chain-rule}. Using Definition~\ref{prelim:def-equivariant}, we obtain for all $X\in \frkg$ and $t\in \R$ that 
\begin{equation*}
\Repr(e^{tX}) \Qc(u_1,\hdots,u_k)= \Qc(\Repr_1(e^{tX})u_1,\hdots,\Repr_k(e^{tX})u_k).
\end{equation*}
By taking the derivative at $t=0$, it then follows that 
\begin{equation*}
\repr(X) \Qc(u_1,\hdots,u_k)= \sum_{j=1}^k \frac{\partial \Qc}{\partial u_j}(u_1,\hdots,u_k) (\Repr_j)_{\ast}(X) u_j. 
\end{equation*}
Together with the definition of the covariant derivatives $\covd_A$ and the chain rule for Euclidean derivatives, this implies \eqref{prelim:eq-equivariant-chain-rule}.
\end{proof}

\begin{remark}
In particular, for smooth $\Phi,\Psi\colon \T^d \rightarrow \frkg$, one has the Leibniz rule 
\begin{equ}
\covd_A^j \big[ \Phi, \Psi \big] = \big[ \covd_A^j \Phi, \Psi \big] + \big[ \Phi, \covd_A^j \Psi \big]. \label{eq:product-formula}
\end{equ}
\end{remark} 

\begin{lemma}[Commutator between covariant heat operator and covariant derivatives]\label{lem:cov-heat-operator-covd-commutator}
Let $A \colon [0, T] \times \T^d \ra \frkg^d$ be a smooth time-varying connection one-form, let $\Psi \colon [0, T] \times \T^d \ra \genV$ be a time-varying $\genV$-valued map, and let $j \in [d]$. Then, with $[\place,\place]$ denoting the commutator of two operators, it holds that 
\begin{equs}
~[\Delta_A, \covd_{A, j}] \Psi &\coloneq \bigl( \Delta_A \covd_{A, j} -  \covd_{A, j} \Delta_A\bigr) \Psi = -2(F_A)_{j\ell} \covd_A^\ell \Psi + (\covd^\ell_A(F_A)_{\ell j}) \Psi, \label{eq:cov-laplacian-covd-commutator} \\
~[\ptl_t - \Delta_A, \covd_{A, j}] \Psi &\coloneq \bigl((\ptl_t - \Delta_A) \covd_{A, j} - \covd_{A, j} (\ptl_t - \Delta_A)\bigr)\Psi \\
&= 2(F_A)_{j\ell} \covd_A^\ell \Psi + (\ptl_t A_j - \covd_A^\ell (F_A)_{\ell j})\Psi.\label{eq:cov-heat-operator-covd-commutator}
\end{equs}
\end{lemma}
\begin{remark}\label{remark:A-YM-GF-drop-term}
In particular, if $A$ is a solution to the Yang--Mills gradient flow, i.e., if $\ptl_t A_j = \covd_A^\ell (F_A)_{\ell j}$, the second term on the right-hand side of \eqref{eq:cov-heat-operator-covd-commutator} drops out. In the Abelian case, this was used in \cite{BC24}, see e.g.\ Section~3.3 of that paper. We will also use this extensively (in the non-Abelian case) in Section~\ref{section:covariant-heat-kernel-estimates}.
\end{remark}
\begin{proof}[Proof of Lemma \ref{lem:cov-heat-operator-covd-commutator}]
For the first commutation identity \eqref{eq:cov-laplacian-covd-commutator}, we apply \eqref{intro:eq-curvature-commutator-identity} twice and \eqref{prelim:eq-equivariant-product-rule} once to obtain
\begin{equs}
\covd_{A}^\ell \covd_{A, \ell} \covd_{A, j} \Psi &= \covd_A^\ell \covd_{A, j} \covd_{A, \ell} \Psi + \covd_{A}^{\ell} \bigl((F_{A})_{\ell j} \Psi\bigr) \\
&= \covd_{A, j} \covd_A^\ell \covd_{A, \ell} \Psi + (F_A)^\ell_{~j} \covd_{A, \ell} \Psi+ \covd_{A}^\ell \bigl((F_{A})_{\ell j} \Psi\bigr) \\
&= \covd_{A, j} \covd_A^\ell \covd_{A, \ell} \Psi - 2 (F_A)_{j \ell} \covd_A^\ell \Psi + (\covd_A^\ell (F_A)_{\ell j}) \Psi.
\end{equs}
The second identity \eqref{eq:cov-heat-operator-covd-commutator} follows by combining the first identity \eqref{eq:cov-laplacian-covd-commutator} along with the fact that
\(
[\ptl_t, \covd_{A, j}]\Psi = \repr(\ptl_t A_j) \Psi
\).
\end{proof}

\begin{lemma}[Bochner formulas]
Let $A \colon [0,T]\times \T^d \ra \frkg^d$ be a time-varying connection one-form and let $\Psi \colon [0,T]\times  \T^d \ra \genV$ be a $\genV$-valued map. Then we have the following Bochner formulas, where $\Delta$ denotes the regular Laplacian
acting on real-valued functions:
\begin{equs}\label{eq:bochner}
\frac{1}{2} \Delta \bigl(|\Psi|_\genV^2\bigr) &= \scal{\Psi, \Delta_A \Psi }_\genV + |\covd_A \Psi|_\genV^2, \\
\frac{1}{2} (\ptl_t - \Delta) \bigl(|\Psi|_\genV^2\bigr) &= \langle \Psi, (\ptl_t - \Delta_A) \Psi \rangle_\genV - |\covd_A \Psi|_\genV^2. \label{eq:parabolic-bochner}
\end{equs}
\end{lemma}
\begin{proof}
We prove \eqref{eq:bochner}; the proof for \eqref{eq:parabolic-bochner} is similar. We compute (recall Remark~\ref{remark:inner-product-G-equivariant} and the product formula \eqref{prelim:eq-equivariant-product-rule})
\begin{equ}
\Delta \langle \Psi, \Psi \rangle_\genV = \ptl_\ell \ptl^\ell \langle \Psi, \Psi \rangle_\genV = 2 \ptl_\ell \langle \Psi, \covd_A^\ell \Psi \rangle_\genV = 2 \langle \Psi, \covd_{A, \ell} \covd_A^\ell \Psi\rangle_\genV + 2 \langle \covd_{A, \ell} \Psi, \covd_A^\ell \Psi \rangle_\genV. 
\end{equ}
The desired result now follows.
\end{proof}

\begin{lemma}[Linearisation]\label{prelim:lem-linearization}
Let $A,B\colon \T^d \rightarrow \frkg^d$ be two connection one-forms, let $C\colon \T^d \rightarrow \frkg^d$ be the one-form such that $A=B+C$, let $\Psi\colon \T^d \rightarrow \genV$, and let $j,k \in [d]$. Then, it holds that
\begin{align}
\covd_A^j \Psi &= \covd_B^j \Psi + C^j \Psi, \label{prelim:eq-linearization-1} \\ 
F_A^{jk} &= F_B^{jk} + \covd_B^j C^k - \covd_B^k C^j + \big[ C^j , C^k \big], \label{prelim:eq-linearization-2}\\
\covd_{A,j} F_A^{jk} &= \covd_{B,j} F_B^{jk}+ \Delta_B C^k + 2 \big[ F_B^{kj},C_j \big] - \covd_B^k \covd_{B}^j C_j  \label{prelim:eq-linearization-3}\\
&\hphantom{=\,}+ \big[ \covd_{B}^j C_j, C^k \big] + \big[ C_j , 2 \covd_B^j C^k - \covd_B^k C^j \big]
+ \big[ C_j, \big[ C^j , C^k \big]\big], \notag \\
\covd_A^k \covd_A^j C_j &= \covd_B^k \covd_B^j C_j - \big[ \covd_B^j C_j, C^k \big],  \label{prelim:eq-linearization-4} \\
\Delta_A \Psi - (\covd_A^j C_j) \Psi 
&= \Delta_B \Psi + 2 C_j \covd_B^j \Psi + C^j C_j \Psi. 
\label{prelim:eq-linearization-5}
\end{align}
\end{lemma}

\begin{proof}
The first identity \eqref{prelim:eq-linearization-1} follows directly from the definition of our covariant derivatives. To obtain the second identity \eqref{prelim:eq-linearization-2}, we write 
\begin{align*}
F_{B+C}^{jk} &= \partial^j (B+C)^k - \partial^k (B+C)^j + \big[ (B+C)^j, (B+C)^k \big] \\ 
&= \bigl( \partial^j B^k - \partial^k B^j + \big[ B^j, B^k \big] \bigr) 
+ \bigl( \partial^j C^k + \big[ B^j , C^k \big] - \partial^k C^j - \big[ B^k, C^j \big] \bigr) + \big[ C^j, C^k \big] \\
&= F_B^{jk} + \bigl( \covd_B C\bigr)^{jk} + \big[ C^j, C^k \big].
\end{align*}
To obtain the third identity \eqref{prelim:eq-linearization-3}, we first use \eqref{prelim:eq-linearization-1} and \eqref{prelim:eq-linearization-2}, which yield that 
\begin{align}
\covd_{A,j} F_A^{jk} 
&= \bigl( \covd_{B,j} + \ad(C_j) \bigr) \bigl( F_B^{jk} + (\covd_B C)^{jk} + \big[C^j,C^k \big] \bigr) \notag \\ 
&= \covd_{B,j} F_B^{jk} + \covd_{B,j} (\covd_B C)^{jk} + \big[ C_j , F_B^{jk} \big] \label{prelim:eq-linearization-p1} \\ 
&\quad+  \covd_{B}^j \big[ C_j , C^k \big] + \big[ C_j, (\covd_B C)^{jk} \big] + \big[ C_j , \big[ C^j , C^k \big]\big].\label{prelim:eq-linearization-p2}
\end{align}
Using the commutator identity \eqref{intro:eq-curvature-commutator-identity}, we see that the second and third summands in~\eqref{prelim:eq-linearization-p1} can be written as 
\begin{equation*}
\covd_{B,j} \bigl( \covd_{B}^j C^k - \covd_{B}^k C^j \bigr) + \big[ C_j, F_B^{jk}\big]
= \Delta_B C^k - \covd_B^k \covd_B^j C_j + 2 \big[ F_B^{kj}, C_j \big].
\end{equation*}
Furthermore, using the covariant Leibniz rule \eqref{eq:product-formula}, we see that the first and second summands in \eqref{prelim:eq-linearization-p2} can be written as
\begin{align*}
    \covd_{B}^j \big[ C_j , C^k \big] + \big[ C_j, (\covd_B C)^{jk} \big] 
    &= \big[ \covd_B^j C_j , C^k \big] + \big[ C_j , \covd_B^j C^k\big]
    + \big[ C_j , \covd_B^j C^k\big] - \big[ C_j, \covd_B^k C^j \big] \\
    &= \big[ \covd_B^j C_j, C^k \big] + \big[ C_j, 2 \covd_B^j C^k - \covd_B^k C^j \big].
\end{align*}
Altogether, we therefore obtain \eqref{prelim:eq-linearization-3}. In order to obtain the fourth identity \eqref{prelim:eq-linearization-4}, we first apply \eqref{prelim:eq-linearization-1} to $\covd_A^j$ and then use the skew-symmetry of the Lie bracket. From this, it follows that
\begin{equation*}
\covd_A^k \covd_A^j C_j = \covd_A^k \bigl( \covd_B^j C_j + \big[ C^j , C_j \big] \bigr) = \covd_A^k \covd_B^j C_j.
\end{equation*}
After applying \eqref{prelim:eq-linearization-1} to $\covd_A^k$ and again using the skew-symmetry of the Lie bracket, we obtain \eqref{prelim:eq-linearization-4}. Finally, in order to obtain the last identity \eqref{prelim:eq-linearization-5}, we use Lemma~\ref{prelim:lem-equivariant} and Lemma~\ref{prelim:lem-equivariant-product-rule}, which yield that
\begin{align*}
\Delta_A \Psi - (\covd_A^j C_j) \Psi 
&= (\covd_B^j + C^j) (\covd_{B,j} + C_j ) \Psi 
- (\covd_B^j C_j + [C^j,C_j]) \Psi \\ 
&= \covd_B^j \covd_{B,j} \Psi 
+ C^j \covd_{B,j} \Psi + \covd_B^j (C_j \Psi) + C^j C_j \Psi  - (\covd_B^j C_j) \Psi \\
&=\Delta_B \Psi + 2 C_j \covd_B^j \Psi + C^j C_j \Psi. \qedhere 
\end{align*}
\end{proof}

In the following, recall our definitions of inner products of differential forms (Notation~\ref{notation:norms-forms}), the covariant curl (Notation~\ref{notation:covariant-curl}), and the covariant divergence of a two-form (Notation~\ref{notation:covariant-divergence-two-form}).

\begin{lemma}\label{lemma:curl-ibp-identity}
Let $A \colon \T^2 \ra \frkg^2$ be a connection one-form, $C \colon \T^2\ra \genV^2$ be a $\genV$-valued one-form, and $F \colon \T^2\ra \genV^{2\times2}$ be a $\genV$-valued two-form. Then we have the integration by parts identity
\begin{equ}
\langle \curl_A C,  F_{12} \rangle = \langle C, \covd_A^* F \rangle . 
\end{equ}
\end{lemma}
\begin{proof}
We have that
\begin{equ}
\langle \curl_A C, F_{12} \rangle = \int_{\T^2} \langle \curl_A C(x) , F_{12}(x) \rangle_\genV \, \dx = \int_{\T^2} \langle \covd_A^1 C^2(x) - \covd_A^2 C^1(x), F_{12}(x) \rangle_{\genV}\, \dx.
\end{equ}
To finish, apply integration by parts and use the antisymmetry of the two-form $F$. 
\end{proof}

\begin{lemma}\label{prelim:lem-Hessian-identity}
Let $A\colon \T^d \rightarrow \frkg^d$ be a connection one-form and let $C\colon \T^d \rightarrow \frkg^d$ be a $\frkg$-valued one-form. Then
\begin{equs}
\label{e:prelim:lem-Hessian-identity-1}
	2 \sum_{1 \leq i < j \leq d}  \langle F_A^{ij}, [C_i, C_j] \rangle
	&= \sum_{1 \leq i < j \leq d} \| \covd_{A, i} C_j - \covd_{A, j} C_i \|_{L^2}^2
	+ \| \covd_A^j C_j \|_{L^2}^2 \\
	&\quad - \sum_{i, j = 1}^{d} \| \covd_{A, i} C_j \|_{L^2}^2 \;.
\end{equs}
Also, we have that
\begin{equs}
\label{e:prelim:lem-Hessian-identity-2}
	2 \sum_{1 \leq i < j \leq d} \langle F_A^{ij}, [C_i, C_j] \rangle
	&= -\sum_{1 \leq i < j \leq d} \| \covd_{A, i} C_j - \covd_{A, j} C_i \|_{L^2}^2
	- \| \covd_A^j C_j \|_{L^2}^2 \\
	&\quad - \sum_{i = 1}^{d}  \langle C_i, \Delta_A C_i + 2 [(F_A)_{ij}, C^j] \rangle \;.
\end{equs}
\end{lemma}
\begin{proof}
By expanding the squares, the right-hand side of \eqref{e:prelim:lem-Hessian-identity-1} is equal to
\begin{equ}
	2 \sum_{1 \leq i < j \leq d} \int_{\T^d}
	\Bigl( \bigl\langle  \covd_{A, i} C_i, \covd_{A, j} C_j \bigr\rangle_{\fg}
	- \bigl\langle \covd_{A, i} C_j, \covd_{A, j} C_i \bigr\rangle_{\fg} \Bigr) \, \dx \;.
\end{equ}
For each $i < j$, integration by parts and \eqref{intro:eq-curvature-commutator-identity} give
\begin{equs}
	\,&\int_{\T^d}
	\Bigl( \bigl\langle \covd_{A, i} C_i, \covd_{A, j} C_j \bigr\rangle_{\fg}
	- \bigl\langle \covd_{A, i} C_j, \covd_{A, j} C_i \bigr\rangle_{\fg} \Bigr) \, \dx \\
	&= \int_{\T^d} \bigl\langle
	(\covd_{A, j}\covd_{A, i} - \covd_{A, i}\covd_{A, j})C_j, C_i
	\bigr\rangle_{\fg} \, \dx \\
	&= \int_{\T^d} \langle F_A^{ij}, [C_i, C_j] \rangle_{\fg} \, \dx \;.
\end{equs}
The equality \eqref{e:prelim:lem-Hessian-identity-1} follows. For \eqref{e:prelim:lem-Hessian-identity-2}, note that
\begin{equs}
	- \sum_{i = 1}^{d} \langle C_i, 2 [(F_A)_{ij}, C^j] \rangle_{\fg} = 2 \sum_{i, j = 1}^{d} \langle [C_i, C_j ], (F_A)_{ij} \rangle_{\fg}  = 4 \sum_{1\leq i < j \leq d} \langle [C_i, C_j ], (F_A)_{ij} \rangle_{\fg}  
\end{equs}
and 
\begin{equs}
- \sum_{i = 1}^{d} \langle C_i, \Delta_A C_i \rangle = \sum_{i,j=1}^d \| \covd_{A,j} C_i \|_{L^2}^2. 
\end{equs}
Combining this with \eqref{e:prelim:lem-Hessian-identity-1} proves \eqref{e:prelim:lem-Hessian-identity-2}. 
\end{proof}

\subsubsection{Identities for the Yang--Mills gradient flow}\label{section:identities-YM-gradient-flow}

Next, we discuss a fully covariant formulation of the Yang--Mills gradient flow, which is described in \cite{Oh2014, Oh2015}. We find this formulation to be quite elegant, and moreover extremely useful for computations. To this end, we first recall the Yang--Mills gradient flow:
\begin{equ}[e:YMheat]
\ptl_t A_i = \covd_A^\ell (F_A)_{\ell i}, ~~ i \in [d].
\end{equ}
Here, $A \colon [0, T] \times \T^d \ra \frkg^d$ is a time-varying connection one-form on $\T^d$. The fully covariant formulation proceeds by appending a temporal component $A_0 \colon [0, T] \times \T^d \ra \frkg$, and viewing $(A_\alpha, \alpha \in \{0, \ldots, d\})$ as a connection one-form on $[0, T] \times \T^d$. Then the above equation may be written as
\begin{equs}
(F_A)_{0i} &= \covd_A^\ell (F_A)_{\ell i}, ~~ i \in [d], \label{eq:YM-heat-flow-fully-covariant}\\
A_0 &= 0. \label{eq:fully-covariant-A-0}
\end{equs}
Here, the index $0$ denotes the time direction, so that $\ptl_0$ denotes differentiation in time. Note that
\begin{equ}
(F_A)_{0i} = \ptl_0 A_i - \ptl_i A_0 + [A_0, A_i] = \ptl_0 A_i - \covd_{A, i} A_0,
\end{equ}
and thus since $A_0 = 0$, we see that the above equation is equivalent to the Yang--Mills gradient flow. The point now is that the choice of $A_0$ precisely corresponds to the choice of gauge, since by the above computation, we may rewrite
\begin{equ}
\ptl_0 A_i = \covd_A^\ell (F_A)_{\ell i} + \covd_{A, i} A_0, ~~ i \in [d],
\end{equ}
and we recall (see e.g.\ the discussion around \cite[Equation (1.6)]{CCHS22}) that the extra pure gradient term $\covd_{A, i} A_0$ is always tangent to the gauge orbit of $A$. For instance, the DeTurck Yang--Mills flow is given by taking $A_0 = \ptl^\ell A_\ell$. 

The name ``fully covariant formulation'' arises from the fact that the equation \eqref{eq:YM-heat-flow-fully-covariant} is covariant even under space-time gauge transformations. Under such transformations, the only thing that may change is the gauge condition \eqref{eq:fully-covariant-A-0}. In particular, given a space-time gauge transformation $g \colon [0, T] \times \T^d \ra \liegroup$, we may hit both sides of the equation \eqref{eq:YM-heat-flow-fully-covariant} by $\Ad_g$, from which we obtain (here $\tilde{A} = g \act A$)
\begin{equs}
(F_{\tilde{A}})_{0i} = \Ad_g (F_A)_{0i} &= \Ad_g \bigl( \covd_A^\ell (F_A)_{\ell i}\bigr) =  \covd_{\tilde{A}}^\ell (F_{\tilde{A}})_{\ell i}, ~~ i \in [d], \\
\tilde{A}_0 &= (g \act A)_0 = \Ad_g A_0 - (\ptl_0 g) g^{-1}.
\end{equs}
As a sample computation, suppose that $\bar{A}$ is a solution to the DeTurck Yang--Mills flow:
\begin{equs}
(F_{\bar{A}})_{0i} &= \covd_{\bar{A}}^\ell (F_{\bar{A}})_{\ell i}, ~~ i \in [d] \\
\bar{A}_0 &= \ptl^\ell \bar{A}_\ell.
\end{equs}
Now suppose that we want to solve for a space-time gauge transformation $g$ such that $g \act \bar{A}$ is a solution to the Yang--Mills gradient flow. This corresponds to solving for $g$ such that $(g \act \bar{A})_0 = 0$, which may be written as the equation
\begin{equ}
(\ptl_0 g) g^{-1} = \Ad_g \bar{A}_0 = \Ad_g (\ptl^\ell \bar{A}_\ell).
\end{equ}
We also remark that the fully covariant formulation is helpful when dealing with equations of the form
\begin{equ}
(\ptl_t - \Delta_A) \phi = G,
\end{equ}
where $\phi, G$ are $\genV$-valued, and $A$ is a solution to the Yang--Mills gradient flow. In particular, it is then natural to write
\begin{equ}
(\covd_{A, 0} - \Delta_A) \phi = G,
\end{equ}
and given a space-time gauge transformation $g$, upon hitting both sides above by $g$ (here the representation of $\liegroup$ on $\genV$ is kept implicit) and then commuting the $g$ in on the left-hand side, we may obtain
\begin{equ}
(\covd_{g \act A, 0} - \Delta_{g \act A}) (g \act \phi) = g \act G.
\end{equ}
This may be converted back into the equation
\begin{equ}\label{eq:g-act-higgs-evolution}
(\ptl_t - \Delta_{g \act A})(g \act \phi) = ((\ptl_t g) g^{-1}) (g \act \phi) + g \act G. 
\end{equ}
Next, we discuss Yang--Mills connections, which are the stationary solutions to the Yang--Mills gradient flow.

\begin{definition}\label{def:YM-connection}
Let $A \colon \T^d \ra \frkg^d$ be a connection one-form. We say that $A$ is Yang--Mills (alternatively, is a Yang--Mills connection) if it solves the equations
\begin{equ}
\covd_A^\ell (F_A)_{\ell i} = 0, ~~ i \in [d].
\end{equ}
We refer to the above equations as the Yang--Mills equations.
\end{definition}

\begin{remark}\label{remark:YM-connection-2D}
In dimension $d = 2$, note that (up to sign) the curvature only has a single component $F_A^{12}$, and the Yang--Mills equations are equivalent to
\begin{equ}
\covd_A^1 (F_A)_{12} = 0, ~~ \covd_A^2 (F_A)_{12} = 0.
\end{equ}
In particular, in 2D, a connection one-form $A$ is Yang--Mills if and only if $\covd_A F_A^{12} = 0$.
\end{remark}

Having introduced the fully covariant formulation of the Yang--Mills gradient flow, we next turn to deriving various identities associated to this flow. First, we record the Bianchi identity:
\begin{equ}\label{eq:bianchi}
\covd_{A, \alpha} (F_A)_{\beta \gamma} + \covd_{A, \beta} (F_A)_{\gamma \alpha} + \covd_{A, \gamma} (F_A)_{\alpha \beta} = 0,
\end{equ}
valid for any values $\alpha, \beta, \gamma$ (recall from Notation~\ref{notation:space-time-indices} that Greek letters are assumed to be in $\{0\} \cup [d]$).

\begin{lemma}\label{lemma:YM-gradient-flow-covariant-divergence-free}
Let $A \colon [0, T] \times \T^d \ra \frkg^d$ be a solution to the Yang--Mills gradient flow
\eqref{e:YMheat}.
Then $\ptl_t A$ is covariant divergence free, i.e.\ $\covd_A^j \ptl_t A_j = 0$.
\end{lemma}
\begin{proof}
We have that $\covd_A^j \ptl_t A_j = \covd_A^j \covd_A^\ell (F_A)_{\ell j}$. By swapping the indices $j \leftrightarrow \ell$ and using \eqref{intro:eq-curvature-commutator-identity}, we obtain
\begin{equs}
\covd_A^j \covd_A^\ell (F_A)_{\ell j} = \covd_A^\ell \covd_A^j (F_A)_{j \ell} &= - \covd_A^\ell \covd_A^j (F_A)_{\ell j} \\
&= - \covd_A^j \covd_A^\ell (F_A)_{\ell j} - [F_A^{\ell j}, (F_A)_{\ell j}] = - \covd_A^j \covd_A^\ell (F_A)_{\ell j},
\end{equs}
from which the desired result follows.
\end{proof}

\begin{lemma}[Equation for the curvature]\label{lemma:curvature-equation}
Let $A \colon [0, T] \times \T^d \ra \frkg^d$ be a solution to the Yang--Mills gradient flow
\eqref{e:YMheat}. Then the curvature $F_A$ satisfies the equation:
\begin{equs}
(\ptl_t - \Delta_A) (F_A)_{ij} = 2 [(F_A)_i^{~\ell}, (F_A)_{\ell j}].
\end{equs}
In the special case $d = 2$, we further have that
$(\ptl_t - \Delta_A) F_A = 0$.
\end{lemma}
\begin{proof}
For simplicity, write $F = F_A$ and $\covd = \covd_A$. Recall the fully covariant formulation of the Yang--Mills gradient flow \eqref{eq:YM-heat-flow-fully-covariant}--\eqref{eq:fully-covariant-A-0}.
By several applications of the Bianchi identity \eqref{eq:bianchi} and curvature commutation \eqref{intro:eq-curvature-commutator-identity}, we obtain
\begin{equs}
\covd_0 F_{ij} &= -\covd_i F_{j0} - \covd_j F_{0i} = \covd_i F_{0j} - \covd_j F_{0i} \\
&= \covd_i \covd^\ell F_{\ell j} - \covd_j \covd^\ell F_{\ell i} = \covd^\ell \covd_i F_{\ell j} - \covd_j \covd^\ell F_{\ell i} + [F_i^{~\ell}, F_{\ell j}] \\
&= - \covd^\ell \covd_\ell F_{ji} - \covd^\ell \covd_j F_{i\ell} - \covd_j \covd^\ell F_{\ell i} + [F_i^{~\ell}, F_{\ell j}] \\
&= \Delta_A F_{ij} - \covd_j \covd^\ell F_{i\ell} - \covd_j \covd^\ell F_{\ell i} - [F^\ell_{~j}, F_{i\ell}] + [F_i^{~\ell}, F_{\ell j}]\\
&= \Delta_A F_{ij} + 2 [F_i^{~\ell}, F_{\ell j}],
\end{equs}
as desired (note that $\covd_0 = \ptl_0$, since $A_0 = 0$). In the special case $d = 2$, there are only two spatial indices, and thus when $i \neq j$, at least one of the terms $F_i^{~\ell}, F_{\ell j}$ must be of the form $F_{11}$ or $F_{22}$, and these are both zero.
\end{proof}

\begin{lemma}[Commutators for higher derivatives]\label{lemma:equations-for-higher-derivatives}
Let $A \colon [0, T] \times \T^d \ra \frkg^d$ be a solution to the Yang--Mills gradient flow
\eqref{e:YMheat}.
Let $\Psi \colon [0, T] \times \T^d \ra \genV$ be a $\genV$-valued map.
Then for any $k \geq 1$, we have the following expression for the commutator:
\begin{equ}
~\big[(\ptl_t - \Delta_A), \covd_A^{(k)} \big]\Psi = \sum_{\substack{j_1+j_2 = k \\ j_1 \geq 0, j_2 \geq 1}} \schO\bigl((\covd_A^{(j_1)} F_A)  (\covd_A^{(j_2)} \Psi)\bigr) .
\end{equ}
Here, as in Notation~\ref{notation:higher-derivatives}, $\covd_A^{(k)}$ denotes any term of the form $\covd_A^{\ell_1} \cdots \covd_A^{\ell_k}$, where $\ell_1, \ldots, \ell_k \in [d]$. Similarly, if $\tilde{A} \colon [0, T] \times \T^d \ra \frkg^d$ is a solution to the backward Yang--Mills gradient flow:
\begin{equ}
\ptl_t \tilde{A}_j = -\covd^\ell_{\tilde{A}} (F_{\tilde{A}})_{\ell j}, \qquad j \in [d],
\end{equ}
then for any $k \geq 1$, we have that
\begin{equ}
\big[(\ptl_t - \Delta_{\tilde{A}}), \covd_{\tilde{A}}^{(k)} \big]\Psi = \sum_{\substack{j_1+j_2 = k \\ j_1, j_2 \geq 0}} \schO\bigl((\covd_{\tilde{A}}^{(j_1)} F_{\tilde{A}}) ( \covd_{\tilde{A}}^{(j_2)} \Psi)\bigr).
\end{equ}
\end{lemma}
\begin{remark}\label{remark:higher-derivatives-equation-with-ad-F-term}
In particular, if $\Psi$ solves the equation:
\begin{equs}
(\ptl_t - \Delta_A) \Psi = G,
\end{equs}
then Lemma~\ref{lemma:equations-for-higher-derivatives} implies
\begin{equ}
(\ptl_t - \Delta_A) \covd_A^{(k)} \Psi = \sum_{\substack{j_1+j_2 = k \\ j_1 \geq 0 \\ j_2 \geq 1}} \schO\bigl((\covd_A^{(j_1)} F_A) (\covd_A^{(j_2)} \Psi)\bigr) + \covd_A^{(k)} G .
\end{equ}
If moreover $G$ has the form $G = \schO(F_A \Psi)$, then we further have that
\begin{equs}
(\ptl_t - \Delta_A) \covd_A^{(k)} \Psi = \sum_{\substack{j_1 + j_2 = k \\ j_1, j_2 \geq 0}} \schO\bigl((\covd_A^{(j_1)} F_A) (\covd_A^{(j_2)} \Psi)\bigr).
\end{equs}
This follows by inductively applying the product formula \eqref{prelim:eq-equivariant-product-rule} to show that $\covd_A^{(k)} \schO(F_A \Psi)$ is of the form given by the right-hand side above. In particular, the restrictions on the indices are now given by $j_1, j_2 \geq 0$, rather than $j_1 \geq 0, j_2 \geq 1$. Similarly, if 
\begin{equ}
(\ptl_t - \Delta_{\tilde{A}}) \tilde{\Psi} = \schO(F_{\tilde{A}} \tilde{\Psi}), 
\end{equ}
then we have that
\begin{equ}
(\ptl_t - \Delta_{\tilde{A}}) \covd_{\tilde{A}}^{(k)} \tilde{\Psi} = \sum_{\substack{j_1 + j_2 = k \\ j_1, j_2 \geq 0}} \schO\bigl((\covd_{\tilde{A}}^{(j_1)} F_{\tilde{A}}) (\covd_{\tilde{A}}^{(j_2)} \tilde{\Psi})\bigr).
\end{equ}
\end{remark}
\begin{proof}[Proof of Lemma \ref{lemma:equations-for-higher-derivatives}]
To show the first claim, we argue by induction. The base case $k = 1$ follows by the heat operator commutation identity \eqref{eq:cov-heat-operator-covd-commutator}, along with the fact that $A$ solves the Yang--Mills gradient flow (recall Remark~\ref{remark:A-YM-GF-drop-term}). Next, suppose that the claim is true for some $k \geq 1$. By writing the commutator as a sum of two terms, and then applying \eqref{eq:cov-heat-operator-covd-commutator} as well as the product formula \eqref{prelim:eq-equivariant-product-rule}, we obtain
\begin{equs}
\big[\ptl_t - \Delta_A, \covd_A^{(k+1)} \big]\Psi &= \big[\ptl_t - \Delta_A, \covd_A\big] \covd_A^{(k)} \Psi + \covd_A\Big(\big[ \ptl_t - \Delta_A, \covd_A^{(k)}\big] \Psi\Big) \\
&= \schO\bigl(F_A \covd_A^{(k+1)} \Psi\bigr) + \covd_A \sum_{\substack{j_1+j_2 = k \\ j_1 \geq 0,\, j_2 \geq 1}} \schO\bigl((\covd_A^{(j_1)} F_A) (\covd_A^{(j_2)} \Psi)\bigr) \\
&=  \sum_{\substack{j_1+j_2 = k +1 \\ j_1 \geq 0,\, j_2 \geq 1}} \schO\bigl((\covd_A^{(j_1)} F_A) (\covd_A^{(j_2)} \Psi)\bigr), 
\end{equs}
as desired. The second claim follows by a similar argument.
\end{proof}

In order to control space-time derivatives, it will sometimes be useful to combine
Lemma~\ref{lemma:equations-for-higher-derivatives} and Remark~\ref{remark:higher-derivatives-equation-with-ad-F-term} with the following result. 

\begin{lemma}\label{prelim:lem-space-time-derivatives-covariant-heat}
Let $A \colon [0, T] \times \T^d \ra \frkg^d$ be a solution to the Yang--Mills gradient flow
\eqref{e:YMheat}. Let $\Psi \colon [0, T] \times \T^d \ra \genV$ be a $\genV$-valued map satisfying 
\begin{equation}\label{prelim:eq-space-time-derivatives-covariant-heat-evoleq}
(\partial_t - \Delta_A) \Psi = \schO(F_A \Psi). 
\end{equation}
Then, it holds for all $k\geq 0$ that
\begin{equation}\label{prelim:eq-space-time-derivatives-covariant-heat-identity}
 \covd_{A,t,x}^{(k)} \Psi = \sum_{\ell \geq 0} \sum_{\substack{j_1+\hdots + j_{\ell+1}=k-2\ell \\ j_1,\hdots, j_{\ell+1} \geq 0}} \schO \Bigl( 
  \bigl( \covd_{A}^{(j_1)} F_A \bigr) 
 \hdots 
 \bigl( \covd_{A}^{(j_\ell)} F_A \bigr)
 \bigl( \covd_{A}^{(j_{\ell+1})} \Psi \bigr)  
 \Bigr). 
\end{equation}
Here, as in Notation~\ref{notation:space-time-derivatives}, $\covd_{A,t,x}^{(k)} \Psi$ is any expression of the form $\covd_A^{\alpha_1} \hdots \covd_A^{\alpha_m} \Psi$, where $m\geq 0$, $\alpha_1,\hdots,\alpha_m \in \{ 0 \} \cup [d]$, and the total number of derivatives (with time derivatives counted twice) equals~$k$. 
\end{lemma}

We emphasise that, unlike the left-hand side, the right-hand side of \eqref{prelim:eq-space-time-derivatives-covariant-heat-identity} only contains spatial covariant derivatives.

\begin{proof}
Similar to the proof of \eqref{eq:cov-heat-operator-covd-commutator}, we may compute (in an intermediate step we apply Lemma~\ref{lemma:YM-gradient-flow-covariant-divergence-free} to conclude that $\covd^\ell F_{0\ell} = 0$)
\begin{equ}\label{eq:heat-temporal-commutator}
[(\ptl_t - \Delta_A), \covd_{A, 0}] \Psi = -[\Delta_A, \covd_{A, 0}] \Psi = 2 (F_A)_{0\ell} \covd_A^\ell \Psi = \schO((\covd_A F_A) (\covd_A \Psi)).
\end{equ}
From Lemma~\ref{lemma:curvature-equation}, it follows that $\partial_t F_A = \Delta_A F_A + \schO(F_A F_A)$. Together with 
\eqref{prelim:eq-space-time-derivatives-covariant-heat-evoleq} and \eqref{eq:heat-temporal-commutator}, this allows us to eliminate time derivatives acting on $\Psi$ or $F_A$ by either adding two spatial derivatives or adding a factor of $F_A$. The desired identity \eqref{prelim:eq-space-time-derivatives-covariant-heat-identity} now follows by induction on the number of time derivatives, and we omit the details. 
\end{proof}

\begin{lemma}[Equation for the covariant exterior derivative]\label{lemma:curl-equation}
Let $A \colon [0, T] \times \T^d \ra \frkg^d$ be a solution to the Yang--Mills gradient flow \eqref{e:YMheat}.
Let $\Psi \colon [0, T] \times \T^d \ra \genV^d$ satisfy
\begin{equ}
(\ptl_t - \Delta_A) \Psi^i = 2 F_A^{i\ell} \Psi_\ell, ~~ i \in [d].
\end{equ}
Define the covariant exterior derivative $\mrm{d}_A \Psi$ as the two-form given by $(\mrm{d}_A \Psi)_{ij} \eqdef \covd_{A, i} \Psi_j - \covd_{A, j} \Psi_i$ for $i, j \in [d]$. Then $\mrm{d}_A \Psi$ satisfies the following equation:
\begin{equs}
(\ptl_t - \Delta_A) (\mrm{d}_A \Psi)_{ij} &= 2 \bigl((F_A)_{i\ell} \covd_A^\ell \Psi_j - (F_A)_{j\ell} \covd_A^\ell \Psi_i\bigr) \\
&\quad \quad \quad + 2 \Bigl(\covd_{A, i}\bigl((F_A)_{j\ell} \Psi^\ell\bigr) - \covd_{A, j} \bigl((F_A)_{i\ell} \Psi^\ell \bigr)\Bigr) \\
&= 2 (F_A)_i^{~\ell} (\mrm{d}_A \Psi)_{\ell j} - 2 (F_A)_j^{~\ell} (\mrm{d}_A \Psi)_{\ell i}  - 2 \bigl(\covd_A^\ell (F_A)_{ij} \bigr) \Psi_\ell.
\end{equs}
In the special case $d = 2$, we further have that (recall from Notation~\ref{notation:covariant-curl} that $\curl_A \Psi = \covd_{A, 1} \Psi_2 - \covd_{A, 2} \Psi_1 = (\mrm{d}_A \Psi)_{12}$)
\begin{equ}
(\ptl_t - \Delta_A) \curl_A \Psi =  2 (\ptl_t A_1) \Psi_2 -2 (\ptl_t A_2) \Psi_1. 
\end{equ}
\end{lemma}
\begin{proof}
For simplicity, write $\covd = \covd_A$ and $F = F_A$. First, by using \eqref{eq:cov-heat-operator-covd-commutator}, we have that
\begin{equ}
(\ptl_t - \Delta_A) \covd_i \Psi_j = 2 F_{i \ell} \covd^\ell \Psi_j + 2 \covd_i \bigl(F_{j\ell} \Psi^\ell\bigr).
\end{equ}
The first desired equation for $(\mrm{d}_A \Psi)_{ij}$ now follows upon taking the difference of the equations for $\covd_i \Psi_j$ and $\covd_j \Psi_i$. For the second desired equation, we first distribute
\begin{equ}
\covd_i \bigl(F_{j\ell} \Psi^\ell\bigr) - \covd_j \bigl( F_{i\ell} \Psi^\ell\bigr) = (\covd_i F_{j\ell} - \covd_j F_{i\ell}) \Psi^\ell + F_{j\ell}  \covd_i \Psi^\ell - F_{i\ell} \covd_j \Psi^\ell.
\end{equ}
By the antisymmetry of $F$ and the Bianchi identity \eqref{eq:bianchi}, we have that
\begin{equ}
\covd_{i} F_{j\ell} - \covd_j F_{i\ell} = \covd_i F_{j\ell} + \covd_j F_{\ell i} = - \covd_\ell F_{ij}.
\end{equ}
Combining the previous few identities, it follows that
\begin{equs}
F_{i\ell}& \covd^\ell \Psi_j - F_{j\ell} \covd^\ell \Psi_i +  \covd_i \bigl(F_{j\ell} \Psi^\ell\bigr) - \covd_j \bigl( F_{i\ell} \Psi^\ell \bigr)\\
&=  F_{i\ell} (\covd^\ell \Psi_j - \covd_j \Psi^\ell) - F_{j\ell} (\covd^\ell \Psi_i - \covd_i \Psi^\ell) -(\covd_\ell F_{ij}) \Psi^\ell,
\end{equs}
which is the second desired equation.

In the case $d = 2$, taking $i = 1, j = 2$, note that $(\mrm{d}_A \Psi)_{11} = (\mrm{d}_A \Psi)_{22} = F_{11} = F_{22} = 0$. Thus starting from the second equation, we further obtain
\begin{equ}
(\ptl_t - \Delta_A) \curl_A \Psi = -2 (\covd_\ell F_{12}) \Psi^\ell = -2 (\covd_1 F_{12})\Psi_1 - 2 (\covd_2 F_{12}) \Psi_2 .
\end{equ}
The desired equation now follows upon using that
$\covd_1 F_{12} = \ptl_t A_2$ and $\covd_2 F_{12} = -\ptl_t A_1$.
\end{proof}

Next, we state the following diamagnetic inequalities (which are also variously referred to as Kato's inequality and the Bochner--Weitzenböck inequality). For proofs, see e.g.\ \cite[Lemmas 4.2 and 4.6]{Oh2015}.

\begin{lemma}[Diamagnetic inequalities]\label{prelim:lem-diamagnetic}
Let $A \colon \T^d \ra \frkg^d$ be a connection one-form and $\Psi \colon \T^d \ra \genV$ be a smooth 
$\genV$-valued map. Then the following inequality holds in a distributional sense:
\begin{equ}\label{eq:gradient-diamagnetic}
|\nabla |\Psi|_\genV| \leq |\covd_A \Psi|_\genV.
\end{equ}
Let $B \colon [0, T] \times \T^d \ra \frkg^d$ be a time-varying connection one-form, and let $\Psi, G \colon [0, T] \times \T^d \ra \genV$ be smooth and such that
$
(\ptl_t - \Delta_B) \Psi = G
$.
Then the following inequality holds in a distributional sense:
\begin{equ}\label{eq:parabolic-diamagnetic}
(\ptl_t - \Delta) |\Psi|_\genV \leq |G|_\genV. 
\end{equ}
\end{lemma}

\begin{remark}\label{remark:curvature-subsolution-heat-equation}
Combining Lemma~\ref{lemma:curvature-equation} with the parabolic diamagnetic inequality \eqref{eq:parabolic-diamagnetic}, we see that in dimension $d = 2$, given a solution $A$ to the Yang--Mills gradient flow \eqref{e:YMheat}, the norm of its curvature $|F_A^{12}|$ is a subsolution to the heat equation, i.e.\ $(\ptl_t - \Delta)|F_A^{12}| \leq 0$.
\end{remark}

\subsubsection{Parallel transport}\label{section:parallel-transport}

In this section, we recall and discuss the parallel transport associated to a connection one-form $A$. (This is of course a classical topic in Riemannian geometry, see e.g.\ \cite[Section 1]{BGV2004} for a textbook introduction.) Parallel transport will allow us to obtain refined short-time expansions for various covariant heat kernels that later appear. Such expansions will play important roles in several places, for instance in understanding the high$\times$high interactions in our covariant local theory (see Section~\ref{section:objects-paracontrolled}), as well as in establishing one of the main inputs (see Section~\ref{section:perturbing-around-YM-connection}) behind the decay mechanism that we uncover.

Throughout this section, fix a vector space $\genV$ and a representation $\Repr$ of $\liegroup$ on $\genV$ (as usual, one should think of $\genV = \frkg$ and $\Repr = \Ad$ or $\genV = \V$ and $\Repr = \rep$).

\begin{definition}[Parallel transport]\label{def:parallel-transport}
Let $A \colon [0, T] \times \T^d \ra \frkg^d$ be a time-varying connection one-form. By appending a temporal component $A_0 \eqdef 0$, we may view $A \colon [0, T] \times \T^d \ra \frkg^{d+1}$ as a space-time connection. Fix a base point $(0, y_0)$. Let $B_{1/10}(y_0)$ be the ball of radius $1/10$ around $y_0$ in $\T^d$. The parallel transport $\pt_{\genV, A}(0, y_0; \place, \place) \colon [0, T] \times B_{1/10}(y_0) \ra \Repr(\liegroup)$ is a $\Repr(\liegroup)$-valued map defined by the following equations. First, $\pt_{\genV, A}(0, y_0; 0, y_0) = \id_{\genV}$. Next, for any $(t, x) \in [0, T] \times B_{1/10}(y_0)$, $(t, x) \neq (0, y_0)$, we have that:
\begin{equ}
\covd_{A, (t, x-y_0)} \pt_{\genV, A}(0, y_0; t, x) = 0.
\end{equ}
These equations may be equivalently viewed as transport equations along rays emanating from the base point $(0, y_0)$. In particular, for $v \in \R^{1+d}$, let $\gamma(t) = (0, y_0) + tv$ be the corresponding space-time line. Then along this line (more precisely, for $t$ such that $\gamma(t) \in [0, T] \times B_{1/10}(y_0)$),  $\pt_{\genV, A}$ satisfies the equations:
\begin{equ}\label{eq:pt-equation}
\covd_{A, \dot{\gamma}(t)} \pt_{\genV, A}(0, y_0; \gamma(t)) = 0, \qquad
\pt_{\genV, A}(0, y_0; \gamma(0)) = \Id_\genV.
\end{equ}
More explicitly, the first equation may be written:
\begin{equs}
\covd_A^\alpha \pt_{\genV, A}(0, y_0; \gamma(t)) \dot{\gamma}_\alpha(t) = \covd_A^0 \pt_{\genV, A}(0, y_0; \gamma(t)) \dot{\gamma}_0(t) + \covd_A^\ell \pt_{\genV, A}(0, y_0; \gamma(t)) \dot{\gamma}_\ell(t)  = 0.
\end{equs}
Here, to define the covariant differentiation, we think of $\pt_{\genV, A}$ as taking values in the vector space of endomorphisms $\End(\genV)$, on which the Lie group $\liegroup$ acts by composition of endomorphisms, i.e.\ given $g \in \liegroup$ the action of $g$ sends $M \in \End(\genV)$ to $\Repr(g) M$. 
\end{definition}

\begin{notation}\label{notation:parallel-transport-genV}
For notational brevity, we will often omit the subscript $\genV$ and write $\pt_A$ in place of $\pt_{\genV, A}$.
\end{notation}

\begin{remark}\label{remark:parallel-transport}
We make the following remarks regarding parallel transport.
\begin{enumerate}[label=(\alph*)]
\item\label{remark:pt-geometry} Geometrically, $A$ can be viewed as a connection one-form on a vector bundle $E \ra \T^d$ with fibres modelled on $\genV$. By appending the temporal component $A_0 = 0$ to $A$, we then obtain a connection one-form on the induced space-time bundle $E^{\mrm{st}} \ra [0, T] \times \T^d$, whose fibre $E^{\mrm{st}}_{(t, x)}$ above a given point $(t, x)$ is just $E_x$. Now, one should view $\pt_A(0, y_0; \place, \place) \in \Gamma((E^{\mrm{st}}_{(0, y_0)})^* \otimes E^{\mrm{st}})$, i.e.\ $\pt_A$ is a section of the tensor bundle $(E^{\mrm{st}}_{(0, y_0)})^* \otimes E^{\mrm{st}}$ (or more precisely, the restricted bundle obtained by restricting the spatial point in the base manifold to lie in a neighbourhood of $y_0$). Due to this, the notion of covariant differentiation $\covd_A$ may be extended to $\pt_A$ in the usual manner. (For instance, if we fix $X \in E^{\mrm{st}}_{(0, y_0)}$, then $\pt_A(0, y_0; \place, \place) X \in \Gamma(E^{\mrm{st}})$ which we may covariantly differentiate.)

\item In the following, we will typically take $A$ to be a solution to the Yang--Mills gradient flow, in which case taking $A_0 = 0$ is natural in light of the fully covariant formulation of the Yang--Mills gradient flow \eqref{eq:YM-heat-flow-fully-covariant}--\eqref{eq:fully-covariant-A-0}.

\item\label{remark:A-constant-in-time-ptA-constant-in-time} In the case $A(t, x) = A(x)$ is constant in time, we may further take $\pt_A(0, y_0; t, x) = \pt_A(y_0; x)$ to be constant in time as well. In particular, for any curve $\gamma(t) = y_0 + tv$ for $v \in \R^d$, $\pt_A(y_0; \gamma(t))$ satisfies the equation:
\begin{equ}
\covd_{A, \dot{\gamma}(t)} \pt_A(y_0; \gamma(t)) = \covd_{A}^j \pt_A(y_0; \gamma(t)) \dot{\gamma}_j(t) = 0, \quad \pt_A(y_0; \gamma(0)) = \Id_\genV,
\end{equ}
at least for $t$ small enough so that $\gamma(t) \in B_{1/10}(y_0)$.

\item\label{remark:pt-isometry} By using the defining equation \eqref{eq:pt-equation}, one may readily verify that $\pt_{A}$ is always an isometry of $\genV$. Indeed, for $X \in \genV$, the time derivative of $\langle \pt_{A}(0, y_0; \gamma(t)) X, \pt_{A}(0, y_0; \gamma(t)) X \rangle_\genV$ is seen to be zero by an application of \eqref{eq:pt-equation}.

\item\label{remark:pt-inverse} By~\ref{remark:pt-isometry}, we have that $\pt_A^{-1} = \pt_A^*$. In the setup of~\ref{remark:pt-geometry}, we may view $\pt_A^*$ itself as a section of the bundle $E_{(0, y_0)}^{\mrm{st}} \otimes (E^{\mathrm{st}})^*$. More concretely, we still have that $\pt_A^* \in \End(\genV)$, except that due to the adjoint, the action of $\liegroup$ is slightly different: given $g \in \liegroup$, the action sends $M \in \End(\genV)$ to $M \Repr(g)^*$. This results in covariant differentiation being slightly different as well:
\begin{equ}\label{eq:covd-pt-inverse-def}
\covd_A (\pt_A^{-1}) = \ptl (\pt_A^{-1}) + \pt_A^{-1} \repr(A)^* = \ptl (\pt_A^{-1}) - \pt_A^{-1} \repr(A).
\end{equ}
Now by differentiating the equation $\id_\genV = \pt_A \pt_A^{-1}$ and applying the product formula \eqref{prelim:eq-equivariant-product-rule} (where the $\mbf{Z}$ therein is taken to be $\End(\genV)$, the action of $\liegroup$ on the vector space $\mbf{Z} = \End(\genV)$ therein is taken as action by conjugation), we thus obtain
\begin{equ}
0 = (\covd_A \pt_A) \pt_A^{-1} + \pt_A \covd_A (\pt_A^{-1}),
\end{equ}
which implies that
\begin{equ}\label{eq:covd-pt-inverse}
\covd_A (\pt_A^{-1}) = - \pt_A^{-1} (\covd_A \pt_A) \pt_A^{-1}.
\end{equ}
In particular, $\pt_A^{-1}$ itself is parallel along any line segment $\gamma$ as taken in Definition~\ref{def:parallel-transport}, i.e.
\begin{equ}\label{eq:pt-inverse-parallel}
\covd_{A, \dot{\gamma}(t)} \bigl(\pt_A(0, y_0; \gamma(t))^{-1}\bigr) = 0.
\end{equ}

\end{enumerate}
\end{remark}

The existence of parallel transport follows from standard ODE arguments. In particular, for a time $t_0$ such that the spatial point of $\gamma(t_0)$ is still in the ball $B_{1/10}(y_0)$, we may solve the following ODE for $h \colon [0, t_0] \ra \Repr(\liegroup)$:
\begin{equ}\label{eq:pt-explicit-ode}
\dot{h}(t) = - \repr(A(\gamma(t))) \cdot \dot{\gamma}(t) h(t) = - \repr(A_\alpha(\gamma(t))) \dot{\gamma}^\alpha(t) h(t), ~~ h(0) = \Id_\genV.
\end{equ}
Setting $\pt_A(0, y_0; \gamma(t)) = h(t)$, we claim that this satisfies the defining equation \eqref{eq:pt-equation}. Towards this end, let $t > 0$. We have that
\begin{equ}
(\ptl_{\dot{\gamma}(t)} \pt_A)(0, y_0; \gamma(t)) = \frac{d}{ds} \pt_A(0, y_0; \gamma(t+s)) \Big|_{s = 0} = \frac{d}{ds} h(t+s) \Big|_{s = 0} = \dot{h}(t) = - \repr(A(\gamma(t))) \cdot \dot{\gamma}(t) h(t), 
\end{equ}
from which \eqref{eq:pt-equation} follows upon moving the right-hand side to the left-hand side, and recalling the definition of $\covd_A$.

We remark that $\pt_A(0, y_0; \gamma(t))$ is well-defined. Indeed, if we take a reparametrisation $\tilde{\gamma}(t) = \gamma(ct)$ for some $c > 0$, and then we solve for $\tilde{h}$ as before (except now with $\gamma$ replaced by $\tilde{\gamma}$), we will have that $\tilde{h}(t/c) = h(t)$. 

\begin{remark}[Different conventions for the parallel transport ODE]
Sometimes, the ODE for parallel transport is written
\begin{equ}
\dot{h}(t) = h(t)  \repr(A(\gamma(t))) \cdot \dot{\gamma}(t), \quad h(0) = \Id_\genV.
\end{equ}
For instance, this is the convention taken in \cite{CCHS22, CCHS22+}. To relate this to \eqref{eq:pt-explicit-ode}, one may readily verify that given a solution $h$ to the above ODE, the adjoint $h^*$ is a solution to \eqref{eq:pt-explicit-ode}. Therefore, the two conventions differ by an adjoint (or equivalently by an inverse \dash recall Remark~\ref{remark:parallel-transport}\ref{remark:pt-inverse}).
\end{remark}

In the following, given a space-time line segment $\gamma(t) = (0, y_0) + tv$, and $t \geq 0$ such that $t|v_x| \leq 1/10$ (where $v_x$ denotes the spatial components of $v \in \R^{1+d}$), we write $\pt_A(\gamma(t)) = \pt_A(0, y_0; \gamma(t))$ for brevity. Additionally, given $0 \leq s \leq t$, we let
\begin{equ}
\pt_A(\gamma(s); \gamma(t)) \eqdef \pt_A(\gamma(t)) \pt_A(\gamma(s))^{-1}.
\end{equ}
Note in particular that $\pt_A(0, y_0; \gamma(t)) = \pt_A(\gamma(0); \gamma(t))$.

\begin{lemma}[Transport lemma]\label{lemma:transport-lemma}
Let $A \colon [0, T] \times \T^d \ra \frkg^d$ be a time-varying connection one-form, and
let $\Psi \colon [0, T] \times \T^d \ra \genV$ be a $\genV$-valued map. Let $y_0 \in \T^d$, $v = (v_0, v_x)\in \R^{1+d}$ such that $v_0 \geq 0$, and let $\gamma(t) = (0, y_0) + tv$. Then for any $t \geq 0$ with $t v_0 \leq T$ and $t|v_x| \leq 1/10$, we have that
\begin{equ}
\Psi(\gamma(t)) = \pt_A(0, y_0; \gamma(t)) \Psi(0, y_0) + \int_0^t \ds \,\pt_A(\gamma(s); \gamma(t)) \covd_{A, \dot{\gamma}(s)} \Psi(\gamma(s)).
\end{equ}
\end{lemma}
\begin{remark}
One may view this as a geometric fundamental theorem of calculus, since it expresses the difference between the values of $\Psi$ at the two endpoints of the path $\gamma$ in terms of an integral involving the covariant derivative of $\Psi$ along the path. When making this comparison, it is imperative to parallel transport $\Psi(0, y_0)$, i.e.\ to take $\pt_A(0, y_0; \gamma(t)) \Psi(0, y_0)$, because the expression
\begin{equ}
\Psi(\gamma(t)) - \Psi(0, y_0)
\end{equ}
does not make geometric sense. Geometrically, $\Psi$ is a section of a vector bundle, and there is no way to directly compare the values of the section lying in different fibres of the bundle (indeed, this is the whole point of the connection and its associated parallel transport).
\end{remark}
\begin{proof}
Recall our notation $\pt_A(\gamma(t)) = \pt_A(0, y_0; \gamma(t))$. First, we derive the ODE satisfied by $\pt_A(\gamma(t))^{-1} \Psi(\gamma(t))$. By the product rule \eqref{prelim:eq-equivariant-product-rule} (with the $\mbf{Z}$ therein taken to be $\genV$, and the action of $\liegroup$ on $\mbf{Z} = \genV$ therein taken to be the trivial action) and \eqref{eq:pt-inverse-parallel} (recalling that covariant differentiation of $\pt_A^{-1}$ is given by \eqref{eq:covd-pt-inverse-def}), we have that
\begin{equs}
\frac{d}{dt} \Big(\pt_A(\gamma(t))^{-1} \Psi(\gamma(t))\Big) &= \covd_{A, \dot{\gamma}(t)} \bigl(\pt_A(\gamma(t))^{-1}\bigr) \Psi(\gamma(t)) + \pt_A(\gamma(t))^{-1} \covd_{A, \dot{\gamma}(t)} \Psi(\gamma(t)) \\
&= \pt_A(\gamma(t))^{-1} \covd_{A, \dot{\gamma}(t)} \Psi(\gamma(t)).
\end{equs}
Upon integrating, we obtain
\begin{equs}
\pt_A(\gamma(t))^{-1} \Psi(\gamma(t)) - \Psi(\gamma(0)) = \int_0^t \ds\, \pt_A(\gamma(s))^{-1} \covd_{A, \dot{\gamma}(s)} \Psi(\gamma(s)).
\end{equs}
Rearranging and using that $\pt_A(\gamma(t)) \pt_A(\gamma(s))^{-1} = \pt_A(\gamma(s); \gamma(t))$, we obtain the desired result.
\end{proof}

We will derive more results about parallel transport in Section~\ref{section:parallel-transport-estimates}. We do not prove those results here, as this section is intended to be short and introductory.

\subsection{Covariance and induced linear maps}

We recall that the real vector space $\V$ carries an orthogonal representation $\rep\colon \liegroup \rightarrow \End(\V)$. As already mentioned in the introduction, we use $\rep_\ast$ to denote the induced derivative representation  
\begin{equation*}
\rep_\ast \colon \frkg \rightarrow \End(\V),
\end{equation*} 
which is defined as $\rep_\ast(A)=\frac{\mathrm{d}}{\mathrm{d}\eps}\rep(e^{\eps A})|_{\eps=0}$. In the following, we sometimes omit the representations $\rep$ and $\rep_\ast$ from our notation. That is, for $g\in \liegroup$, $A\in \frkg$, and $\Phi \in \V$, we sometimes write
\begin{equation*}
g \Phi = \rep(g) \Phi \qquad \text{and} \qquad A \Phi = \rep_\ast(A) \Phi. 
\end{equation*}

\begin{definition}\label{prelim:def-Casimir-induced}
We let $(\ONBg_a)_{a}$ and $(\ONBV_m)_{m}$ be orthonormal bases of $\frkg$ and $\V$, respectively. Then, we define the elements $\Cov_\frkg \in \frkg \otimes \frkg$ and $\Cov_\V \in \V \otimes \V$ by
\begin{equation}\label{prelim:eq-Casimir}
\Cov_\frkg \eqdef  \ONBg_a \otimes \ONBg^a \qquad \text{and} \qquad \Cov_\V \eqdef \ONBV_m \otimes \ONBV^m.  
\end{equation}
(Summation over repeated indices is implicit.) For all $A\in \frkg$ and $\Phi \in \V$, we define the linear maps 
\begin{equs}[2]
\Cov_\fg^\fg(A) &\eqdef [\ONBg^a,[\ONBg_a,A]], 
\qquad & \qquad 
\Cov_\V^\fg(A) &\eqdef \bil\bigl(\ONBV_m \otimes \rep_\ast(A)\ONBV^m\bigr), 
\label{prelim:eq-Casimir-induced-1} \\ 
\Cov_\fg^\V(\Phi) &\eqdef \rep_\ast(\ONBg_a)\rep_\ast(\ONBg^a) \Phi, 
\qquad & \qquad 
\Cov_{\V}^{\V} (\Phi) &\eqdef \langle \ONBV_m , \Phi \rangle_{\V} \,  \ONBV^m = \Phi. 
\label{prelim:eq-Casimir-induced-2} 
\end{equs}
\end{definition}
In \eqref{prelim:eq-Casimir-induced-1}--\eqref{prelim:eq-Casimir-induced-2}, the subscript is the space one takes a basis of and the superscript is the space the map acts on.  We note that $\Cov_\fg^\fg$ corresponds to the Killing map used in \cite[Definition 2.23]{BC23}. The reason for not keeping to the notation from \cite{BC23} is that we want to use similar notation for the different terms in \eqref{prelim:eq-Casimir-induced-1}--\eqref{prelim:eq-Casimir-induced-2}. 

In order to apply our gauge covariance result (Lemma \ref{main:lem-gauge-covariance}) with our choice of the renormalisation constants from Definition \ref{main:def-renormalization-constant}, we need the following lemma. 

\begin{lemma}\label{prelim:lem-G-covariance}
The four maps $\Cov_\fg^\fg$, $\Cov_\V^\fg$, $\Cov_\fg^\V$, and $\Cov_\V^\V$ are $\liegroup$-covariant. 
\end{lemma}

\begin{proof}
Since the proofs for the four different maps are quite similar, we only present the details for $\Cov^\V_\frkg$. To this end, let $g\in \liegroup$ and let $\higgs\in \V$. Using our earlier identity \eqref{prelim:eq-equivariant-p1}, we have that 
\begin{equation*}
\Cov^\V_\frkg ( \rep(g) \higgs)= \rep_\ast(\ONBg_a) \rep_\ast(\ONBg^a) \rep(g) \higgs = \rep(g) \rep_\ast(\Ad_{g^{-1}} \ONBg_a) \rep_\ast(\Ad_{g^{-1}} \ONBg^a) \higgs.
\end{equation*}
Since the adjoint representation is orthogonal, $(\Ad_{g^{-1}} \ONBg_a)_a$ is also an orthonormal basis of $\frkg$. By expressing $(\Ad_{g^{-1}} \ONBg_a)_a$ in terms of the orthonormal basis $(\ONBg_a)$, it is easy to see that 
\begin{equation*}
\rep_\ast(\Ad_{g^{-1}} \ONBg_a) \rep_\ast(\Ad_{g^{-1}} \ONBg^a)=
\rep_\ast(\ONBg_a) \rep_\ast(\ONBg^a),
\end{equation*}
which then implies the desired result.
\end{proof}

\begin{example}\label{ex:su_n_case}
Let us consider the case $\liegroup=\SU(N)$, $\frkg=\su(N)$, and $\V=\C^N$, where the representation $\rep$ is given by the standard matrix-vector product. We emphasise that we treat  $\su(N)\subseteq \C^{N\times N}$ as a real Lie algebra and $\C^N$ as a real vector space. In particular, $\dim_{\R}(\su(N))=N^2-1$ and $\dim_{\R}(\C^N)=2N$. The standard inner products on 
$\su(N)$ and $\C^N$ are given by
\begin{equation*}
\langle A, B \rangle_{\su(N)} = \Tr(A^\dagger B) \qquad \text{and} \qquad \langle u , v \rangle_{\C^N} = \Re(u^\dagger v)\;,
\end{equation*}
where $\Tr$ is the standard, unnormalised trace and $\dagger$ denotes the Hermitian conjugate $\C^{N\times M} \to \C^{M\times N}$. 
When equipped with the above inner product, $\C^N$ has the orthonormal basis
$(\ONBV_m) = \{e_j\}_{j \le N} \cup \{\icomplex e_j\}_{j \le N}$,
where the $j$th entry of $e_j \in \C^N$ is equal to one, and all other entries are zero.  
For any $A\in \su(N)$ and $v,w\in \C^N$, it holds that
\begin{equation*}
\langle \bil(v\otimes w), A \rangle_{\su(N)} 
= \langle v, Aw \rangle_{\C^N} 
= \Re( v^\dagger A w) 
= \tfrac{1}{2} (v^\dagger A w - w^\dagger A v) 
= - \Tr\bigl( \tfrac{1}{2} ( v w^\dagger - w v^\dagger) A \bigr). 
\end{equation*}
We now note that $v w^\dagger - w v^\dagger$ is skew-Hermitian and that, although its trace is purely imaginary, it is not necessarily trace-free. After projecting to the space of trace-free matrices, we then obtain that
\begin{equation*}
\bil(v \otimes w) = \tfrac{1}{2} (vw^\dagger - w v^\dagger) - \tfrac{1}{2N} \Tr( vw^\dagger - w v^\dagger ) \Id = P_0 \bigl( \tfrac{1}{2} (vw^\dagger - w v^\dagger) \bigr), 
\end{equation*}
where $P_0 A \eqdef A - \frac{1}{N} \Tr(A) \Id$ is the orthogonal projection onto trace-free matrices. Using the definition of $\Cov_{\frkg}^{\frkg}$, one can show for all  $A,B\in \su(N)$ that $\langle \Cov_{\frkg}^{\frkg} A, B \rangle_{\su(N)}= \Kil(A,B)$, where $\Kil$ is the Killing form of $\su(N)$. Since the Killing form of $\su(N)$ is $\Kil(A,B)=2N \Tr(AB) = -2N \scal{A,B}_{\su(N)}$, this implies that 
\begin{equation*}
\Cov_{\frkg}^{\frkg} A = - 2 N A. 
\end{equation*}
From a direct computation, one sees that 
\begin{equation}\label{prelim:eq-sum-ONB}
\ONBg^a \ONBg_a = - \frac{N^2-1}{N} \cdot \Id \qquad \text{and}  \qquad \ONBV^m \ONBV_m^\dagger = 2 \cdot \Id. 
\end{equation}
Indeed, to obtain the first identity, one first shows that $\ONBg^a \ONBg_a$ commutes with every element $U\in \SU(N)$ and then computes its trace. To obtain the second identity, one simply inserts $(\ONBV_m) = \{e_j\}_{j \le N} \cup \{\icomplex e_j\}_{j \le N}$. We emphasise that the factor of $2$ appears since $e_j e_j^\dagger$ and $(\icomplex e_j) (\icomplex e_j)^\dagger$ make the same contribution. From the first identity in \eqref{prelim:eq-sum-ONB}, we obtain for all $v \in \C^N$ that
\begin{equation*}
\Cov_\frkg^\V v = \ONBg^a \ONBg_a v = - \frac{N^2-1}{N} v. 
\end{equation*}
Similarly, from the second identity in \eqref{prelim:eq-sum-ONB}, we obtain that 
\begin{align*}
\Cov_\V^\frkg& A = \bil(\ONBV_m \otimes A\ONBV^m ) 
= \tfrac{1}{2} P_0 ( \ONBV_m (A\ONBV^m)^\dagger - A\ONBV^m \ONBV_m^\dagger) \\ 
&= \tfrac{1}{2} P_0(\ONBV_m (\ONBV^m)^\dagger A^\dagger - A \ONBV^m \ONBV_m^\dagger)  = -2 P_0 ( A ) = -2 A. 
\end{align*}
In the last identity, we used that since $A\in \su(N)$, it is trace-free. To avoid any confusion with articles in the literature, including \cite{BC24}, we emphasise again that $\V=\C^N$ is treated as a real vector space. If $\V=\C^N$ were treated as a complex vector space, then $\{e_j\}_{j \le N}$ would already be an orthonormal basis, which changes the above calculations. The easiest way to see this is to compare the standard Gaussian vectors on $\V=\C^N$ as a vector space over $\R$ and $\C$. Given a sequence $(X_j)_{j=1}^{2N}$ of independent, standard, real-valued Gaussians, the standard Gaussian vectors are given by 
\begin{equation*}
\sum_{j=1}^{N} X_j e_j + \sum_{j=1}^{N} X_{j+N} (\icomplex e_j)  \qquad \text{and} \qquad 
\sum_{j=1}^{N} \frac{X_j+\icomplex X_{j+N}}{\sqrt{2}} e_j,
\end{equation*}
respectively. Thus, the standard Gaussian vectors on $\V$ as a real or complex vector space differ by a factor of $\sqrt{2}$. 
\end{example}

\begin{lemma}\label{prelim:lem-identity-Casrhog}
For every orthonormal basis $(\ONBV_m)$  of $\V$ and every $\Phi \in \V$, one has
\begin{equation*}
\bil(\ONBV_m \otimes \Phi) \ONBV^m = \PhiCas(\Phi). 
\end{equation*}
\end{lemma}

\begin{proof}
Writing $(\ONBg_a)$ for an orthonormal basis of $\fg$, one has 
\begin{equ}
\bil(\ONBV^k \otimes \Phi) \ONBV_k
= \scal{\bil(\ONBV^k \otimes \Phi),\ONBg^a}_\fg \ONBg_a \ONBV_k
= \ONBg_a  \scal{\ONBV^k , \ONBg^a \Phi}_\V \ONBV_k
= \ONBg_a \ONBg^a \Phi = \PhiCas \Phi\;,
\end{equ}
as claimed. (Here we do not explicitly write $\rep_*$ for the action of $\fg$ on $\V$.)
\end{proof}

\subsection{Analytic preliminaries}\label{section:harmonic-analysis}

In this section, we collect various analytic preliminaries. We begin by setting our convention for Hölder spaces.

\begin{definition}[Hölder spaces]\label{def:Holder-space}
For $\alpha \in \R$, we take $\|\place\|_{\Cs_x^{\alpha}}$ to be the usual (Besov-)Hölder norm, and we define the space $\Cs_x^\alpha$ to be the closure of $C_x^\infty$ with respect to this norm.
\end{definition}

The following gauge covariant Gagliardo--Nirenberg--Sobolev inequality is a direct consequence of the diamagnetic inequality \eqref{eq:gradient-diamagnetic} and the classical Gagliardo--Nirenberg--Sobolev inequality (see e.g.\ \cite[Proposition A.3]{Tao06} for the statement of the inequality on $\R^d$).

\begin{lemma}[Gauge-covariant Gagliardo--Nirenberg--Sobolev inequality]\label{prelim:lem-GNS}
Let $B\colon \T^2 \rightarrow \frkg^2$ be a connection one-form and let $\phi \colon \T^2 \rightarrow \genV$. Furthermore, let $\theta \in (0,1)$ and $r\in [1,\infty)$ satisfy $\theta \geq\frac{2}{r}$.   Then, it holds that 
\begin{equation*}
\| \phi \|_{L_x^\infty} \lesssim_{\theta,r} \| \phi \|_{L_x^r} + \| \phi \|_{L_x^r}^{1-\theta} \| \covd_B \phi \|_{L_x^r}^{\theta},
\end{equation*}
where the implicit constant is uniform in $B$. 
\end{lemma}

We next state the following endpoint parabolic Schauder estimate, whose proof is similar to the proof of \cite[Lemma 2.5]{TaoGRWMIV} (which works on $\R^2$).

\begin{lemma}[Endpoint Schauder estimate]\label{lemma:endpoint-strichartz}
Let $T \in (0, 1]$ and let $u \colon [0, T] \times \T^2 \ra \R$ be a solution to the heat equation:
\begin{equ}
(\ptl_t - \Delta) u = 0.
\end{equ}
Then we have the following endpoint Schauder estimate:
\begin{equ}
\|u\|_{L_t^2 L_x^\infty([0, T])} \lesssim \|u(0)\|_{L_x^2}.
\end{equ}
\end{lemma}
\begin{proof}
We proceed by duality. It suffices to show that for $F \in L_t^2 L_x^1([0, T])$, we have that
\begin{equ}
\langle u, F \rangle_{L_{t,x}^2([0, T])} \lesssim \|u(0)\|_{L_x^2} \|F\|_{L_t^2 L_x^1([0, T])}.
\end{equ}
Towards this end, we may write
\begin{equ}
\langle u, F \rangle_{L_{t,x}^2([0, T])} = \Big\langle u(0), \int_0^T \ds \, e^{s \Delta} F(s) \Big \rangle_{L_x^2} \leq \|u(0)\|_{L_x^2} \bigg\|\int_0^T \ds \, e^{s \Delta} F(s) \bigg\|_{L_x^2}.
\end{equ}
Thus to finish, it suffices to show that
\begin{equ}
\bigg\|\int_0^T \ds \, e^{s \Delta} F(s) \bigg\|_{L_x^2}^2 \lesssim \|F\|_{L_t^2 L_x^1([0, T])}^2.
\end{equ}
Towards this end, observe that (noting a symmetry between $s$ and $s'$ to reduce to $s > s'$)
\begin{equs}
\bigg\|\int_0^T \ds \, e^{s \Delta} &F(s) \bigg\|_{L_x^2}^2 = 2 \int_0^T \ds \int_0^s \ds' \langle F(s), e^{(s + s')\Delta} F(s')\rangle_{L_x^2} \\
&\lesssim \int_0^T \ds \, \|F(s)\|_{L_x^1} s^{-1} \int_0^s \ds' \, \|F(s')\|_{L_x^1}\\
	& = \int_0^{T} \|F(s')\|_{L_x^1} \Big( \int_{s'}^{T} \|F(s)\|_{L_x^1} s^{-1} \ds \Big) \ds' \\
	&\leq \int_0^{T} \|F(s')\|_{L_x^1} \bigg(\int_{s'}^T \|F(s)\|_{L_x^1}^2 s^{-\frac12} \ds\bigg)^{1/2} \cdot \bigg(\int_{s'}^{T} s^{-\frac32} \ds \bigg)^{1/2} \, \ds' \\
	&\lesssim \|F \|_{L^2_t L^1_x([0, T])} \bigg(\int_0^T \Big( \int_{s'}^T \|F(s)\|_{L_x^1}^2 s^{-\frac12} \ds \Big) (s')^{-\frac12} \ds'\bigg)^{1/2} \lesssim \|F \|_{L^2_t L^1_x([0, T])}^2 \;,
\end{equs} 
as desired.
\end{proof}

The following lemma will allow us to control certain norms (see Definition~\ref{def:Fnorm}) of the curvature of a given solution to the Yang--Mills gradient flow. This will allow us to obtain covariant heat kernel and Schauder estimates which behave very similarly to the standard heat kernel and Schauder estimates, which in turn will be crucial in our covariant local theory (Sections~\ref{section:covariant-Ansatz}--\ref{section:covariant_lwp}). In the following, for solutions $u \colon [0, T] \times \T^2 \ra \R$ to the heat equation $(\ptl_t - \Delta) u = 0$, we recall the classical monotonicity formula:
\begin{equ}\label{eq:heat-monotonicity}
\bigl\|u(0)\bigr\|_{L_x^2}^2 - \bigl\|u(T)\bigr\|_{L_x^2}^2 = 2 \bigl\|\nabla u \bigr\|_{L_{t,x}^2([0, T])}^2.
\end{equ}

\begin{lemma}\label{lemma:strichartz-estimate-using-decrease}
Let $T \in (0, 1/2]$, and let $f \colon [0, T] \times \T^2 \ra [0, \infty)$ be a non-negative subsolution to the heat equation, i.e.,
\begin{equ}
(\ptl_t - \Delta) f \leq 0.
\end{equ}
Let $D \colon [0, T] \ra [0, \infty)$ be the decrease of the squared $L_x^2$-norm of $f$: 
\begin{equ}
D(t) \eqdef \|f(0)\|_{L_x^2}^2 - \|f(t)\|_{L_x^2}^2, ~~ t \in [0, T].
\end{equ}
Then we have the following estimate:
\begin{equs}
\|f\|_{L_t^2 L_x^\infty([0, T])} \lesssim T^{1/2} \|f(0)\|_{L_x^2} + \abs{\log T}^{1/2} D(T)^{1/2} .
\end{equs}
\end{lemma}
\begin{proof}
Let $g$ be the solution to the heat equation with initial data $f(0)$. Then we have the pointwise bound $0 \leq f \leq g$, and thus it suffices to estimate $g$ by the right-hand side of the claimed inequality. Towards this end, let $N_0 \eqdef 10^{3} T^{-100}$ (say). We estimate $P_{\leq N_0} g$ and $P_{> N_0} g$ separately. Note that both of these terms are solutions to the heat equation. For the low-frequency part, by Bernstein's inequality (see e.g.\ \cite[Lemma 2.1]{BCD2011}), we have that for any $t \in [0, T]$,
\begin{equs}
\|P_{\leq N_0} g(t)\|_{L_x^\infty} \leq \sum_{N \leq N_0} \|P_N g(t)\|_{L_x^\infty} &\lesssim \|g(t)\|_{L_x^2} + \sum_{1 < N \leq N_0} \|\nabla P_N g(t)\|_{L_x^2} \\
&\lesssim \|g(t)\|_{L_x^2} + (\log N_0)^{1/2} \|\nabla P_{\leq N_0} g(t) \|_{L_x^2},
\end{equs}
where we applied Cauchy--Schwarz in the final inequality (recalling that the sum over $N$ only includes the dyadic numbers!). 
After taking the $L^2_t$-norm, we obtain 
\begin{equs}
\|P_{\leq N_0} g\|_{L_t^2 L_x^\infty([0, T])} \lesssim T^{1/2} \|f(0)\|_{L_x^2} + \abs{\log T}^{1/2} \|\nabla P_{\leq N_0} g\|_{L_{t,x}^2([0, T])},
\end{equs}
where we used \eqref{eq:heat-monotonicity} to bound $\|g(t)\|_{L_x^2} \leq \|g(0)\|_{L_x^2} = \|f(0)\|_{L_x^2}$. Next, observe by \eqref{eq:heat-monotonicity} that
\begin{equs}
\|\nabla P_{\leq N_0} g\|_{L_{t,x}^2([0, T])}^2 \lesssim \|\nabla g\|_{L_t^2 L_x^2([0, T])}^2 &\sim \|g(0)\|_{L_x^2}^2 - \|g(T)\|_{L_x^2}^2 \\
&\leq \|f(0)\|_{L_x^2}^2 - \|f(T)\|_{L_x^2}^2 = D(T),
\end{equs}
where in the second-to-final step, we used the pointwise bound $f \leq g$ and the fact that $f(0) = g(0)$. Combining the previous few displays, it follows that the low-frequency part $P_{\leq N_0} g$ satisfies the required bound. 

For the high-frequency part, by the endpoint Schauder estimate (Lemma~\ref{lemma:endpoint-strichartz}) and the fact that $P_{> N_0} g$ itself is a solution to the heat equation, we have that
\begin{equs}
\|P_{> N_0} g\|_{L_t^2 L_x^\infty([0, T])} \lesssim \|P_{> N_0} g(0)\|_{L_x^2}. 
\end{equs}
By classical heat smoothing estimates (see e.g.\ \cite[Lemma 2.4]{BCD2011}), we also have that
\begin{equs}
\|P_{> N_0} g(T)\|_{L_x^2} \lesssim e^{- c N_0^2 T} \|P_{> N_0} g(0)\|_{L_x^2} \ll \|P_{> N_0} g(0)\|_{L_x^2},
\end{equs}
where the final inequality follows by our choice of $N_0 = 10^{3} T^{-100}$. Combining this with \eqref{eq:heat-monotonicity}, we obtain
\begin{equs}
\|\nabla P_{> N_0} g\|_{L_{t,x}^2([0, T])}^2 \sim \|P_{> N_0} g(0)\|_{L_x^2}^2 - \|P_{> N_0} g(T)\|_{L_x^2}^2 \sim \|P_{> N_0} g(0)\|_{L_x^2}^2,
\end{equs}
and from this we further obtain (again applying \eqref{eq:heat-monotonicity} in the penultimate estimate)
\begin{equs}
\|P_{> N_0} g\|_{L_t^2 L_x^\infty([0, T])}^2 \lesssim \|\nabla P_{> N_0} g\|_{L_{t,x}^2 ([0, T])}^2 &\lesssim \|\nabla g\|_{L_{t,x}^2([0, T])}^2 \\
&\sim \|g(0)\|_{L_x^2}^2 - \|g(T)\|_{L_x^2}^2 \lesssim D(T),
\end{equs}
where in the final inequality, we again used the pointwise bound $f \leq g$ and the fact that $g(0) = f(0)$. The bound for the high-frequency part now follows (in fact it is more than acceptable).
\end{proof}

\begin{lemma}\label{lemma:euclidean-heat-kernel-gradient-torus}
For $t \in (0, 1]$ and $x, y \in \T^d$ such that $|x-y| \leq 10^{-1}$, we have that
\begin{equ}
\bigg| \frac{\nabla \fp(0, y; t, x)}{\fp(0, y; t, x)} + \frac{(x-y)}{2t} \bigg| \lesssim t^{10}.
\end{equ}
\end{lemma}
\begin{proof}
We view $\fp$ as a $2\pi$-periodic function on $\R^d$, and we identify $\T^d \cong (-\pi, \pi]^d$. Letting $\fp^{\Euc}_t(x)$ be the heat kernel on $\R^d$, we have the formula
\begin{equ}
\fp(0, y; t, x) = \sum_{n \in \Z^d} \fp^{\Euc}_t(x-y+2\pi n).
\end{equ}
By translation invariance, we may assume that $y = 0$, and so $x$ is such that $|x| \leq 10^{-1}$. We may then compute
\begin{equs}
\nabla \fp(0, 0; t, x) &= \nabla \fp^{\Euc}_t(x) + \sum_{n \neq 0} \nabla \fp^{\Euc}_t(x+ 2\pi n) \\ 
&= -\frac{x}{2t} \fp^{\Euc}_t(x) - \sum_{n \neq 0} \frac{x+2\pi n}{2t} \fp^{\Euc}_t(x+ 2\pi n) \\
&= - \frac{x}{2t} \fp(0, 0; t, x) - \sum_{n \neq 0}  \frac{\pi n}{t}\fp^{\Euc}_t(x + 2\pi n). 
\end{equs}
Since $|x| \leq 10^{-1}$, we have that for $n \neq 0$,
\begin{equs}
\Big| \frac{\pi n \fp^{\Euc}_t(x + 2\pi n)}{t \fp^{\Euc}_t(x)} \Big| &\lesssim t^{-1} \exp\Bigl(- c \frac{|n|^2}{t}\Bigr) \\
&\leq t^{-1} \exp\Bigl(-\frac{c|n|^2}{2t}\Bigr) \exp\Bigl(-\frac{c|n|^2}{2}\Bigr) \lesssim t^{10} \exp\Bigl(-\frac{c|n|^2}{2}\Bigr) .
\end{equs}
It follows that
\begin{equ}
\sum_{n \neq 0}  \Big| \frac{\pi n \fp^{\Euc}_t(x + 2\pi n)}{t \fp^{\Euc}_t(x)}\Big| \lesssim t^{10},
\end{equ}
from which the desired result follows.
\end{proof}

\begin{lemma}[Gronwall]\label{lemma:gronwall}
Let $f, a, b \colon [0, T] \ra [0, \infty)$ be such that $a\in L^1([0,T])$, $f, b\in L^\infty([0,T])$, and
\begin{equ}
\|f\|_{L^\infty([0,t])} \leq \| a f\|_{L^1([0,t])} + \|b\|_{L^\infty([0,t])}\;,
\end{equ}
uniformly over $t \in [0,T]$. Then we have that
\begin{equ}
\|f\|_{L^\infty([0,t])} \leq  \|b\|_{L^\infty([0,t])} \exp \|a\|_{L^1([0,t])}\;,
\end{equ}
uniformly over $t \in [0,T]$.
\end{lemma}

In the following, we define function spaces of negative temporal regularity. They will be used as part of our modified solution spaces, see Definition~\ref{local:def-modified-solution-space} below. 

\begin{definition}[Temporal regularity spaces]\label{def:temporal-regularity-space}
Let $T > 0$ and $(\mc{B}, \|\place\|_{\mc{B}})$ be a Banach space. For $\alpha \in (0, 1]$, define $\|\place\|_{\Cs_t^\alpha \mc{B}}$ as the usual $\alpha$-Hölder norm for $\mc{B}$-valued functions:
\begin{equ}
\|\phi\|_{\Cs_t^\alpha \mc{B}([0, T])} \eqdef \|\phi\|_{C_t^0 \mc{B}([0, T])} + \sup_{\substack{0 \leq s < t \leq T}} \frac{\|\phi(t) - \phi(s)\|_{\mc{B}}}{|t-s|^\alpha}.
\end{equ}
For $\alpha \in (-1, 0)$, we set
\begin{equ}
\|\phi\|_{\Cs_t^\alpha \mc{B}([0, T])} \eqdef \bigg\| \int_0^t \phi(s) \ds \bigg\|_{\Cs_t^{1+\alpha} \mc{B}([0, T])}.
\end{equ}
In other words, the $\Cs_t^\alpha \mc{B}$ norm of $\phi$ is defined to be the $\Cs_t^{1+\alpha} \mc{B}$ norm of the integral $t \mapsto \int_0^t \phi(s) \ds$.
\end{definition}

\begin{remark}
It follows directly from the definition that for $\alpha \in (-1, 0)$, time integration maps $\phi \in \Cs_t^\alpha \mc{B}$ to $t \mapsto \int_0^t \phi(s) \ds \in \Cs_t^{1+\alpha} \mc{B}$. Conversely, time differentiation maps $\phi \in \Cs_t^{1+\alpha} \mc{B}$ to $\partial_t \phi(t) \in \Cs_t^\alpha \mc{B}$.
\end{remark}

The following version of Young's integration bound \cite{Young} is an immediate consequence of the sewing
lemma \cite[Lem.~4.2]{BookPeter} (just set $\Xi_{s,t} =  (\int_s^t \phi_1(r)\,\dr)\phi_2(s)$).

\begin{lemma}[Product estimate for temporal regularity spaces]\label{prelim:lem-product-temporal}
Let $T \in (0, 1]$ and let $(\mc{B}, \|\place\|_{\mc{B}})$ and $(\mc{B}_i, \|\place\|_{\mc{B}_i})$, where $i = 1, 2$, be Banach spaces such that there is a continuous multiplication map $\mc{B}_1 \times \mc{B}_2 \ra \mc{B}$. Furthermore, let $\alpha \in (-1, 0)$ and $\beta \in (0, 1)$ be such that $\alpha + \beta > 0$. Then, we have for all smooth functions $\phi_i \colon [0, T] \ra \mc{B}_i$, where $i = 1, 2$, that 
\[
\|\phi_1 \phi_2\|_{\Cs_t^\alpha \mc{B}([0, T])} \lesssim \|\phi_1\|_{\Cs_t^\alpha \mc{B}_1([0, T])} \|\phi_2\|_{\Cs_t^\beta \mc{B}_2([0, T])}.\tag*{$\qedsymbol$}
\]
\end{lemma}

The following lemma shows that our temporal regularity spaces are well-behaved under time localisations, which will be used in the proof of our weighted covariant Schauder estimates (see Proposition~\ref{local:prop-weighted-covariant-Schauder}).

\begin{lemma}[Time-localisation]\label{prelim:lem-time-localization}
Let $0<T\leq 1$, let $-1<\alpha \leq \beta \leq 0$, and let $\chi \in C^\infty_c((0,\infty))$ be arbitrary (i.e.\ it is not necessarily the cutoff function fixed at the beginning of this section). Furthermore, let $\mc{B}$ be a Banach space. For all $\phi \colon [0,T]  \rightarrow \mc{B}$ and $N\in \dyadic$, it then holds that 
\begin{equation}\label{prelim:eq-time-localization}
\big\| \chi(N^2 t) \phi \big\|_{\Cs_t^\alpha \mc{B}([0,T])} 
\lesssim N^{2(\alpha-\beta)} \big\| \phi \big\|_{\Cs_t^\beta \mc{B}([0, T])}. 
\end{equation}
\end{lemma}

\begin{proof}
We let $\Phi$ be the anti-derivative of $\phi$ satisfying $\Phi(0)=0$. By definition, it holds that $\| \Phi \|_{\Cs_t^{1+\beta} \mc{B}} \lesssim \| \phi \|_{\Cs_t^\beta \mc{B}}$ and, for all $t\in [0,T]$, it holds that 
\begin{equation}\label{prelim:eq-time-localization-e0}
\big\| \Phi(t) \big\|_{\mc{B}}
= \big\| \Phi(t) - \Phi(0)\big\|_{\mc{B}} \lesssim t^{1+\beta} \| \Phi \|_{\Cs_t^{1+\beta} \mc{B}}. 
\end{equation}
Since 
\begin{equation*}
\int_0^t \chi(N^2 s) \phi(s) \ds = \chi(N^2 t) \Phi(t) - \int_0^t N^2 \chi^\prime(N^2 s) \Phi(s) \ds, 
\end{equation*}
our desired estimate \eqref{prelim:eq-time-localization} can be reduced to the two estimates
\begin{align}
\Big\| \chi(N^2 t) \Phi(t) \Big\|_{\Cs_t^{\alpha^\prime} \mc{B}}
&\lesssim N^{2(\alpha^\prime-\beta^\prime)} \big\| \Phi(t) \big\|_{\Cs_t^{\beta^\prime}\mc{B}}, \label{prelim:eq-time-localization-e1} \\ 
\Big\| N^2 \int_0^t \chi^\prime(N^2 s) \Phi(s) \, \ds \Big\|_{\Cs_t^{\alpha^\prime} \mc{B}}
&\lesssim N^{2(\alpha^\prime-\beta^\prime)} \big\| \Phi(t) \big\|_{\Cs_t^{\beta^\prime}\mc{B}}, \label{prelim:eq-time-localization-e2}
\end{align}
where $\alpha^\prime\coloneq 1+\alpha$ and $\beta^\prime\coloneq 1+\beta$. By interpolation, it suffices to prove \eqref{prelim:eq-time-localization-e1} and \eqref{prelim:eq-time-localization-e2} for $(\alpha^\prime,\beta^\prime)=(0,0),\, (0,1), \, (1,1)$. This follows easily from \eqref{prelim:eq-time-localization-e0} and standard estimates, and we omit the remaining details. 
\end{proof}

\subsection{Supermartingale estimates}\label{section:supermartingale}

In this subsection, we prove uniform-in-time moment bounds for a stochastic process satisfying supermartingale-type conditions. Near the end of the article (see Section \ref{section:global-theory}), they will be used to prove the uniform-in-time bounds in our soft and hard global well-posedness results (Theorems~\ref{thm:global_soft} and~\ref{thm:global_hard}). The reason for presenting our supermartingale-type bounds in a separate subsection, rather than including them in Section \ref{section:global-theory}, is that the arguments are easier to understand without the heavy notation used in Section \ref{section:global-theory}. As previously mentioned in the AI disclosure statement, the proofs in this subsection were suggested to the authors by AI (GPT-6 Astra). The closest arguments in the existing literature, which formed the basis of the AI's response, are in \cite{PR99,DFMS04}.

\begin{proposition}\label{prelim:prop-X-supermartingale}
Let $(\mathcal{G}_m)_{m=0}^\infty$ be a filtration. Let $(X_m)_{m=0}^\infty$, $(Y_m)_{m=0}^\infty$, and $(D_m)_{m=0}^\infty$ be stochastic processes adapted to the filtration and assume that almost surely $X_m\geq 10$ for all $m\geq 0$ and that $Y_0\geq 0$. Let $r\geq 1$, $p>r+\frac{1}{2}$, $0\leq \eps <1$, $c>0$, and $K\geq 1$ be constants. Furthermore, assume the following three properties:
\begin{enumerate}[label=(\roman*)]
\item (Recursive estimate) For all $m\geq 0$, it holds that
\begin{equation}\label{prelim:eq-recursive}
X_{m+1} \leq X_m - D_m + Y_{m+1} +K. 
\end{equation}
\item (Decay) For all $m\geq 0$, it holds that 
\begin{equation}\label{prelim:eq-decay}
D_m \geq c \frac{\sqrt{X_m}}{\log(X_m)}. 
\end{equation}
\item (Fluctuation bounds) For all $m\geq 1$, it holds that
\begin{equation}\label{prelim:eq-fluctuation}
\E \Bigl[ (Y_{m+1}- \eps D_m)_{+}^p \Big| ~ \mathcal{G}_{m-1} \Bigr] \leq K. 
\end{equation}
Furthermore, it holds that 
\begin{equation}\label{prelim:eq-fluctuation-initial}
\E \Bigl[ (Y_{1}- \eps D_0)_{+}^p \Big| ~ \mathcal{G}_{0} \Bigr] \leq Y_0^p. 
\end{equation}
\end{enumerate}
Then, there exists a constant $C=C(\eps,c,p,r,K)$ such that
\begin{equation}\label{prelim:eq-X-conclusion}
\sup_{m\geq 0} \E \bigl[ X_m^r \big| ~ \mathcal{G}_0 \bigr] \leq C \bigl( X_0^r+Y_0^r+1 \bigr) \qquad \text{and} \qquad \limsup_{m\rightarrow \infty} \E \bigl[ X_m^r \big| ~ \mathcal{G}_0 \bigr] \leq C. 
\end{equation}
\end{proposition}

\begin{remark}
Since the conclusions are conditional on $\mathcal{G}_0$, we may always make the qualitative assumption that $X_0$ and $Y_0$ (as well as $Z_0$ below) are bounded.
\end{remark}

\begin{remark}\label{prelim:rem-martingale-global-iteration}
In our global iteration scheme (see Section \ref{section:global-theory}), $X_m$ will be the covariant energy of the smooth part of our solution, $D_m$ will be its expected decay, and $Y_{m+1}$ will capture different stochastic fluctuations. In this context, it will be crucial that the fluctuation bound only involves the conditional moment with respect to $\mathcal{G}_{m-1}$, and not $\mathcal{G}_m$. This allows us to average over the stochastic forcing on the preceding time interval in our iteration, as required by our estimates. 

The reason for treating the initial fluctuation separately, i.e., for separating \eqref{prelim:eq-fluctuation} and \eqref{prelim:eq-fluctuation-initial}, also lies in our global iteration scheme. To be more precise, since the initial data in our global well-posedness theorems (Theorems \ref{thm:global_soft} and \ref{thm:global_hard}) is not required to have any random structure, the estimates in the first step of our global iteration scheme sometimes have to differ from the estimates in later steps.
\end{remark}

Before we turn to the proof of Proposition \ref{prelim:prop-X-supermartingale}, we first state and prove the following simpler lemma. 

\begin{lemma}\label{prelim:lem-Z-supermartingale-bound}
Let $(\mathcal{G}_m)_{m=0}^\infty$ be a filtration and let $(Z_m)_{m=0}^\infty$ be an adapted stochastic process which almost surely satisfies $Z_m\geq 10$ for all $m\geq 0$. Furthermore, let $c>0$, $K\geq 1$, and $p>\frac{3}{2}$. Finally, assume that, for all $m\geq 0$,
\begin{equation}\label{prelim:eq-Z-assumption}
\E \bigl[ Z_{m+1}^p\big|~\mathcal{G}_m \bigr]^{\frac{1}{p}} \leq Z_m - c \frac{\sqrt{Z_m}}{\log(Z_m)} + K.  
\end{equation}
Then, it holds for all $0<r< p - \frac{1}{2}$, $\frac{1}{2}<\alpha<\min(p-r,1)$, and $m\geq 0$ that 
\begin{equation}\label{prelim:eq-Z-conclusion}
 \E\bigl[ Z_m^r \big| ~ \mathcal{G}_0 \bigr] \leq \min\biggl( Z_0^r, \frac{C Z_0^p}{(m+1)^{\frac{p-r}{\alpha}}} \biggr) + C, 
\end{equation}
where $C=C(\alpha,c,p,r,K)$ is a constant.
\end{lemma}

\begin{proof}[Proof of Lemma \ref{prelim:lem-Z-supermartingale-bound}]
For expository purposes, we split the proof into five different steps. 

\proofstep{Step 1: Initial simplification.} After choosing $R=R(\alpha,c,p,K)$ sufficiently large, and then choosing $A=A(\alpha,c,p,K,R)$ sufficiently large, our assumption \eqref{prelim:eq-Z-assumption} implies that 
\begin{alignat}{3}
\E \bigl[ Z_{m+1}^p\big|~\mathcal{G}_m \bigr]  
&\leq Z_m^p- Z_m^{p-\alpha} \qquad &&\text{on } \{ Z_m > R\}, \label{prelim:eq-Z-bound-large-values} \\ 
\E \bigl[ Z_{m+1}^p\big|~\mathcal{G}_m \bigr] 
&\leq A^p \qquad &&\text{on } \{ Z_m \leq R\}.\label{prelim:eq-Z-bound-small-values}
\end{alignat}

\proofstep{Step 2: The supermartingales.} To simplify the notation, set $\beta\eqdef \alpha/p \in (0,1)$. For any $m\geq 0$, we define $H_m \colon [0,\infty) \rightarrow [0,\infty)$ as 
\begin{equation*}
H_m(x) \eqdef \bigl( x^\beta + \beta m \bigr)^{\frac{1}{\beta}}. 
\end{equation*}
From a direct computation of the first and second derivatives of $H_m$, it follows that $H_m$ is both strictly increasing and concave. We further define a stopping time\footnote{The reason for using the symbol $\nu$ for our stopping time, rather than the more conventional choice $\tau$, is that the symbol $\tau$ will be needed in Proposition \ref{prelim:prop-stopping-time} below.} $\nu$ as 
\begin{equation*}
\nu \eqdef \inf \bigl\{ n \geq 0 \colon Z_n \leq R \bigr\},
\end{equation*}
i.e., as the first time at which the process is bounded by $R$. We now claim that  
\begin{equation}\label{prelim:eq-supermartingale}
(H_{\nu \wedge m}(Z_{\nu \wedge m}^p))_{m=0}^{\infty}
\end{equation}
is a supermartingale. To see this, we first split
\begin{align}
\E \bigl[ H_{\nu \wedge (m+1)}(Z_{\nu \wedge (m+1)}^p) \big| ~ \mathcal{G}_m \bigr] 
&= \E \bigl[ \ind\{ \nu\leq m\}  H_{\nu \wedge (m+1)}(Z_{\nu \wedge (m+1)}^p) \big| ~ \mathcal{G}_m \bigr] \label{prelim:eq-supermartingale-first} \\ 
&+ \E \bigl[ \ind\{ \nu> m\}  H_{\nu \wedge (m+1)}(Z_{\nu \wedge (m+1)}^p) \big| ~ \mathcal{G}_m \bigr]. \label{prelim:eq-supermartingale-second}
\end{align}
The first summand \eqref{prelim:eq-supermartingale-first} is clearly equal to $\ind\{ \nu\leq m\} H_{\nu \wedge m}(Z_{\nu \wedge m}^p)$, which is as desired. Using Jensen's inequality and \eqref{prelim:eq-Z-bound-large-values}, the second summand can be estimated by
\begin{align*}
\E \bigl[& \ind\{ \nu> m\}  H_{\nu \wedge (m+1)}(Z_{\nu \wedge (m+1)}^p) \big| ~ \mathcal{G}_m \bigr] 
= \ind\{ \nu> m\} \E \bigl[  H_{m+1}(Z_{m+1}^p) \big| ~ \mathcal{G}_m \bigr] \\
&\leq \ind\{ \nu> m\} H_{m+1} \bigl( \E \bigl[ Z_{m+1}^p \big| ~ \mathcal{G}_m \bigr] \bigr) \leq \ind\{ \nu> m\} H_{m+1}(Z_m^p - Z_m^{p-\alpha}).
\end{align*}
Since $x\in [0,\infty) \mapsto x^\beta$ is concave, it holds that 
\begin{equation*}
(Z_m^p - Z_m^{p-\alpha})^\beta \leq (Z_m^p)^\beta - \beta (Z_m^p)^{\beta-1} Z_m^{p-\alpha}= (Z_m^p)^\beta - \beta, 
\end{equation*}
where we used that $\beta=\frac{\alpha}{p}$.  Together with the definitions of $H_m$ and $H_{m+1}$, this implies that 
\begin{equation*}
H_{m+1}(Z_m^p - Z_m^{p-\alpha}) \leq H_m(Z_m^p). 
\end{equation*}
Thus, the second summand \eqref{prelim:eq-supermartingale-second} can be bounded by $\ind\{ \nu >m \} H_m(Z_m^p)=\ind\{ \nu >m \} H_{\nu \wedge m}(Z_{\nu\wedge m}^p)$, which is as desired. Using simpler arguments than for \eqref{prelim:eq-supermartingale}, one can also show that 
\begin{equation}\label{prelim:eq-supermartingale-simpler}
(Z_{\nu\wedge m}^r)_{m=0}^\infty
\end{equation}
is a supermartingale. Indeed, as in \eqref{prelim:eq-supermartingale-first}--\eqref{prelim:eq-supermartingale-second}, one first decomposes 
\begin{equation}\label{prelim:eq-supermartingale-simpler-decomposition}
\E \bigl[ Z_{\nu \wedge (m+1)}^r \big|~ \mathcal{G}_m \bigr]
= \E \bigl[ \ind\{ \nu \leq m \} Z_{\nu \wedge (m+1)}^r \big|~ \mathcal{G}_m \bigr]
+ \E \bigl[ \ind \{ \nu > m\} Z_{\nu \wedge (m+1)}^r \big|~ \mathcal{G}_m \bigr].
\end{equation}
The first summand in \eqref{prelim:eq-supermartingale-simpler-decomposition} is clearly equal to $\ind\{ \nu \leq m \} Z_{\nu \wedge m}^r$, which is as desired. Using Hölder's inequality and \eqref{prelim:eq-Z-bound-large-values}, the second summand in \eqref{prelim:eq-supermartingale-simpler-decomposition} can be bounded by 
\begin{equation*}
\ind \{ \nu > m\} \E \bigl[  Z_{m+1}^r \big|~ \mathcal{G}_m \bigr] \leq \ind \{ \nu > m\} (Z_m^p - Z_m^{p-\alpha})^{\frac{r}{p}} \leq \ind \{ \nu > m\} Z_m^r = \ind \{ \nu > m\} Z_{\nu \wedge m}^r,
\end{equation*}
which is also as desired. 

\proofstep{Step 3: Estimate for an excursion above $R$.} For all $m\geq 0$, we show that 
\begin{equation}\label{prelim:eq-single-excursion}
\E \bigl[ \ind\{ \nu > m\} Z_m^r \big| ~ \mathcal{G}_0 \bigr] \leq \min\biggl( Z_0^r, \frac{B Z_0^p}{(m+1)^{\frac{p-r}{\alpha}}} \biggr), 
\end{equation}
where $B=B(\alpha,p,r,A)$. The case $m=0$ is trivial, and we may assume that $m \geq 1$. To obtain the bound by $Z_0^r$ in \eqref{prelim:eq-single-excursion}, which is the simpler estimate, we only need to use the supermartingale property for \eqref{prelim:eq-supermartingale-simpler} and that 
$\ind\{\nu>m\}Z_m^r\leq Z_{\nu\wedge m}^r$. The bound by $B Z_0^p/(m+1)^{\frac{p-r}{\alpha}}$ in \eqref{prelim:eq-single-excursion} is more difficult to obtain and relies on our supermartingale from \eqref{prelim:eq-supermartingale}. Indeed, using the supermartingale property, it holds that 
\begin{equation}\label{prelim:eq-supermartingale-application}
\E \bigl[ H_{\nu \wedge m }(Z_{\nu \wedge m}^p) \big| ~ \mathcal{G}_0 \bigr] \leq H_0 (Z_0^p).
\end{equation}
By definition, the right-hand side of \eqref{prelim:eq-supermartingale-application} is equal to $Z_0^p$. The left-hand side of \eqref{prelim:eq-supermartingale-application} can be bounded below by 
\begin{equation*}
\E \bigl[ H_{\nu \wedge m }(Z_{\nu \wedge m}^p) \big| ~ \mathcal{G}_0 \bigr] \geq \E \bigl[ \ind \{\nu>m\} H_{m}(Z_{m}^p) \big| ~ \mathcal{G}_0 \bigr]. 
\end{equation*}
Since $0<r<p$, it holds that 
\begin{equation*}
 H_{m}(Z_{m}^p) =  \bigl((Z_m^p)^{\alpha/p}+\tfrac{\alpha}{p} m\bigr)^{\frac{p}{\alpha}} = \bigl(Z_m^\alpha+\tfrac{\alpha}{p} m\bigr)^{\frac{r}{\alpha}}\bigl(Z_m^\alpha+\tfrac{\alpha}{p} m\bigr)^{\frac{p-r}{\alpha}} \geq Z_m^r \bigl( \tfrac{\alpha}{p} m\bigr)^{\frac{p-r}{\alpha}}.
\end{equation*}
From this, it follows that
\begin{equation*}
\E \bigl[ \ind\{ \nu >m\} Z_m^r \big| ~ \mathcal{G}_0 \bigr] \leq \bigl( \tfrac{\alpha}{p} m\bigr)^{-\frac{p-r}{\alpha}} Z_0^p. 
\end{equation*}
Provided that $B$ is sufficiently large, which also allows us to replace $m$ by $m+1$, this implies the second bound in \eqref{prelim:eq-single-excursion}.

\proofstep{Step 4: Estimate for an excursion starting from below $R$.} In this step, we show for all $0\leq k < m$ that 
\begin{equation}\label{prelim:eq-excursion-starting-below}
\E \bigl[ \ind\{ Z_m > R, \hdots, Z_{k+1} > R , Z_k \leq R \} Z_m^r \big| ~ \mathcal{G}_k \bigr] \leq \frac{A^p B}{(m-k)^\frac{p-r}{\alpha}}. 
\end{equation}
To obtain \eqref{prelim:eq-excursion-starting-below}, we first use the time-shifted version of \eqref{prelim:eq-single-excursion}, which implies that
\begin{equation*}
\E \bigl[ \ind\{ Z_m > R, \hdots, Z_{k+1} > R \} Z_m^r \big| ~ \mathcal{G}_{k+1} \bigr] \leq \frac{B Z_{k+1}^p}{(m-k)^\frac{p-r}{\alpha}}. 
\end{equation*}
Using \eqref{prelim:eq-Z-bound-small-values}, we then obtain that 
\begin{align*}
\E \bigl[& \ind\{ Z_m > R, \hdots, Z_{k+1} > R , Z_k \leq R \} Z_m^r \big| ~ \mathcal{G}_k \bigr] \\
&\leq \E \biggl[  \ind\{ Z_k \leq  R \}  \frac{B Z_{k+1}^p}{(m-k)^\frac{p-r}{\alpha}} \bigg| ~ \mathcal{G}_k \biggr] \leq \frac{A^p B}{(m-k)^\frac{p-r}{\alpha}},
\end{align*}
which yields the desired estimate \eqref{prelim:eq-excursion-starting-below}.

\proofstep{Step 5: Conclusion.} Using \eqref{prelim:eq-single-excursion} and \eqref{prelim:eq-excursion-starting-below}, we obtain that 
\begin{align*}
\E \bigl[ Z_m^r \big| ~ \mathcal{G}_0 \bigr]  
&=  \E \bigl[ \ind\{ Z_m > R, \hdots, Z_0 > R \} Z_m^r \big| ~ \mathcal{G}_0 \bigr]  \\ 
&+ \sum_{k=0}^{m-1} \E \bigl[ \ind\{ Z_m > R, \hdots, Z_{k+1} > R, Z_k \leq R \} Z_m^r \big| ~ \mathcal{G}_0 \bigr]   
+ \E\bigl[ \ind\{ Z_m \leq R\} Z_m^r \big| ~ \mathcal{G}_0 \bigr]  \\
&\leq \min\biggl( Z_0^r, \frac{B Z_0^p}{(m+1)^{\frac{p-r}{\alpha}}} \biggr) + A^p B \sum_{k=0}^{m-1} \frac{1}{(m-k)^{\frac{p-r}{\alpha}}} +R^r \\
&\leq   \min\biggl( Z_0^r, \frac{C Z_0^p}{(m+1)^{\frac{p-r}{\alpha}}} \biggr) + C, 
\end{align*}
where $C=C(\alpha,p,r,A,B,R)$ is a constant. In the last line, we used that $(p-r)/\alpha>1$. From this, we clearly obtain the desired conclusion \eqref{prelim:eq-Z-conclusion}.
\end{proof}

\begin{proof}[Proof of Proposition \ref{prelim:prop-X-supermartingale}]
To simplify the notation, we set $U_{m+1}\eqdef (Y_{m+1}-\eps D_m)_+$. From the recursive estimate \eqref{prelim:eq-recursive} and the definition of $U_{m+1}$, it then follows that 
\begin{equation}\label{prelim:eq-recursive-rewritten}
\begin{aligned}
X_{m+1} &\leq X_m - (1-\eps) D_m + K + U_{m+1} \\ 
 & \leq X_m - (1-\eps) (D_m \wedge X_m) +K +U_{m+1}. 
\end{aligned}
\end{equation}
The reason for introducing the minimum of $D_m$ and $X_m$ in \eqref{prelim:eq-recursive-rewritten} is that $X_m - (1-\eps) (D_m \wedge X_m) +K$ is nonnegative, which will be needed below. We then define the corrected process
\begin{equation}\label{prelim:eq-corrected-process}
Z_m \eqdef X_m + 2 W_m, \qquad \text{where} \qquad W_m \eqdef \E \bigl[ U_{m+1}^p \big|~ \mathcal{G}_m \bigr]^{\frac{1}{p}}. 
\end{equation}
The purpose of the correction $2W_m$ in \eqref{prelim:eq-corrected-process} is to later absorb a contribution resulting from the $U_{m+1}$-term in \eqref{prelim:eq-recursive-rewritten}, while leaving the additional decay term $-W_m$. The stochastic processes $(W_m)_{m=0}^\infty$ and $(Z_m)_{m=0}^\infty$ are adapted to the filtration. From the definition, it follows directly that $X_m \leq Z_m$ for all $m\geq 0$. In order to bound $X_m$, it therefore suffices to bound $Z_m$, for which we want to use 
Lemma \ref{prelim:lem-Z-supermartingale-bound}. For this reason, we now show for all $m\geq 0$ that 
\begin{equation}\label{prelim:eq-Z-reduction-estimate}
\E \bigl[ Z_{m+1}^p \big| \mathcal{G}_m \bigr]^{\frac{1}{p}} \leq Z_m - c^\prime \frac{\sqrt{Z_m}}{\log(Z_m)} + K^\prime,
\end{equation}
where $c^\prime>0$ and $K^\prime\geq 1$ are constants depending on $c,\eps,p,r$, and $K$. From the tower property of conditional expectations and the fluctuation bound \eqref{prelim:eq-fluctuation}, it follows that 
\begin{equation*}
\E \bigl[ W_{m+1}^p \big|~\mathcal{G}_m\bigr] 
= \E \Bigl[ \E \bigl[ U_{m+2}^p \big| ~ \mathcal{G}_{m+1} \bigr] \Big| ~ \mathcal{G}_m \Bigr] 
= \E \bigl[  U_{m+2}^p   \big| ~ \mathcal{G}_m \bigr] \leq K. 
\end{equation*}
Together with our recursive estimate from \eqref{prelim:eq-recursive-rewritten}, it then follows that 
\begin{align*}
\E \bigl[ Z_{m+1}^p \big| ~ \mathcal{G}_m \bigr]^{\frac{1}{p}}
&\leq \E \bigl[ X_{m+1}^p \big| ~ \mathcal{G}_m \bigr]^{\frac{1}{p}}
+ 2 \E \bigl[ W_{m+1}^p \big| ~ \mathcal{G}_m \bigr]^{\frac{1}{p}}  \\
&\leq \E \bigl[ (X_{m}-(1-\eps) (D_m\wedge X_m) +K)^p \big| ~ \mathcal{G}_m \bigr]^{\frac{1}{p}} + \E \bigl[ U_{m+1}^p \big| ~ \mathcal{G}_m \bigr]^{\frac{1}{p}} +2 K^{\frac{1}{p}} \notag \\ 
&= X_m - (1-\eps) (D_m\wedge X_m) + K + W_m + 2 K^{\frac{1}{p}} \\
&= Z_m - (1-\eps) (D_m\wedge X_m) - W_m + K + 2K^{\frac{1}{p}}. 
\end{align*}
In the penultimate identity, we used that $X_m$ and $D_m$ are $\mathcal{G}_m$-measurable and the definition of $W_m$, and in the last identity, we used the definition of $Z_m$. By using our decay condition \eqref{prelim:eq-decay}, as well as the fact that $X_m \geq 10$ for all $m\geq 0$, we easily obtain that 
\begin{equation*}
(1-\eps) (D_m\wedge X_m) + W_m \geq (1-\eps) \min\Bigl( c \frac{\sqrt{X_m}}{\log(X_m)}, X_m \Bigr) + W_m \geq c^\prime \frac{\sqrt{Z_m}}{\log(Z_m)}. 
\end{equation*}
This completes the proof of \eqref{prelim:eq-Z-reduction-estimate}. By using Lemma \ref{prelim:lem-Z-supermartingale-bound} for some $\frac{1}{2}<\alpha<\min(1,p-r)$, we then obtain for all $m\geq 0$ that 
\begin{equation*}
\E \bigl[ Z_m^r \big| ~ \mathcal{G}_0 \bigr] 
\leq  \min\biggl( Z_0^r, \frac{C Z_0^p}{(m+1)^{\frac{p-r}{\alpha}}} \biggr)   + C. 
\end{equation*}
Since $X_m \leq Z_m$ for all $m \geq 0$, and since our fluctuation bound \eqref{prelim:eq-fluctuation-initial} implies that $W_0\leq Y_0$, it then follows that 
\begin{equation*}
\E\bigl[ X_m^r \big| ~ \mathcal{G}_0 \bigr]
\leq \E \bigl[ Z_m^r \big| ~ \mathcal{G}_0 \bigr] 
\leq  \min\biggl( (X_0+2Y_0)^r, \frac{C (X_0+2Y_0)^p}{(m+1)^{\frac{p-r}{\alpha}}} \biggr)   + C. 
\end{equation*}
This readily leads to the desired estimates on the supremum and limit superior in \eqref{prelim:eq-X-conclusion}. 
\end{proof}

In Section \ref{section:global-theory}, Proposition \ref{prelim:prop-X-supermartingale} will be used to prove the estimate \ref{main:item-global-smooth-part} from our hard global well-posedness result (Theorem \ref{thm:global_hard}). While the estimate \ref{main:item-global-smooth-part} is uniform in our iteration index $m\geq 1$, our stopping time $\tmod_m$ is random, and thus \ref{main:item-global-smooth-part} does not directly yield any estimates for the solution at a given deterministic time. To obtain estimates at a deterministic time, like in our soft global well-posedness result (Theorem \ref{thm:global_soft}), we will also need the following proposition.

\begin{proposition}\label{prelim:prop-stopping-time}
Let the same assumptions as in Proposition \ref{prelim:prop-X-supermartingale} be satisfied. In addition, let $q_1>0$, let $0\leq q_2\leq p$,  and let $(\tau_m)_{m=1}^\infty$ be a stochastic process such that, for all $m\geq 0$, $\tau_{m+1}$ is $\mathcal{G}_{m+1}$-measurable and $0<\tau_{m+1}\leq 1$. Furthermore, we assume that, for all $m\geq 1$,
\begin{equation}\label{prelim:eq-time-increment-moments}
\E \bigl[ \tau_{m+1}^{-q_1} \big| ~ \mathcal{G}_{m-1} \bigr] \leq \E\bigl[ X_m^{q_2}\big| ~ \mathcal{G}_{m-1} \bigr] +K. 
\end{equation}
For each $T>0$, let 
\begin{equation*}
M_T \eqdef \inf \Bigl\{ m \geq 1 \colon \sum_{k=0}^{m-1} \tau_{k+1} \geq T \Bigr\}. 
\end{equation*}
Then, $M_T$ is almost surely finite. Furthermore, it holds that 
\begin{equation}\label{prelim:eq-stopping-time-final-estimate}
\sup_{T>0} \E \bigl[ X_{M_T+1}^r \big| ~ \mathcal{G}_0 \bigr] \leq C^\prime \bigl( X_0^r + Y_0^r + 1 \bigr) \qquad \text{and} \qquad \limsup_{T\rightarrow \infty} \E \bigl[ X_{M_T+1}^r \big| ~ \mathcal{G}_0 \bigr] \leq C^\prime, 
\end{equation}
where $C^\prime=C^\prime(\eps,c,p,q_1,q_2,r,K)$ is a constant. 
\end{proposition}

We emphasise two technical aspects of Proposition \ref{prelim:prop-stopping-time}, which will be convenient in the proof of Proposition \ref{prop:modified-iteration-uniform-bounds} below. First, the conditional moment estimate \eqref{prelim:eq-time-increment-moments} is only assumed for $m\geq 1$, but not for $m=0$. Second, \eqref{prelim:eq-stopping-time-final-estimate} controls the process $X=(X_m)_{m=0}^\infty$ not at $M_T$ itself, but one step after~$M_T$. 

Throughout the proof, we will refer to the sequence $(\tau_{m})_{m=1}^\infty$ as the physical time increments and to their partial sums $(t_m)_{m=0}^\infty$ as the physical times. The reason is that, in Section \ref{section:global-theory}, $m\in \mathbb{N}_0$ will correspond to our iteration index, whereas $t_m\in [0,\infty)$ will correspond to the actual time variable in the stochastic Yang--Mills--Higgs model.

\begin{proof}[Proof of Proposition \ref{prelim:prop-stopping-time}] We use the same notation as in the proofs of Proposition \ref{prelim:prop-X-supermartingale} and Lemma~\ref{prelim:lem-Z-supermartingale-bound}. In particular, $Z_m = X_m + 2 W_m$ and $R$, $A$, and $B$ are as in the preceding proofs. For expository purposes, we split the argument into five steps. 

\proofstep{Step 1: Further applications of our supermartingale.} In this step, we prove variants of the supermartingale bounds from the proof of Lemma \ref{prelim:lem-Z-supermartingale-bound}. To this end, we first recall that \eqref{prelim:eq-supermartingale} and \eqref{prelim:eq-supermartingale-simpler} are supermartingales. For any $j\geq 0$, we introduce the stopping time
\begin{equation*}
\nu_j \eqdef \inf\{ m \geq j \colon Z_m \leq R\},
\end{equation*}
which is the first return to the region below $R$ at or after time $j$. We emphasise that $\nu_j - j$ is essentially the same as $\nu$ from the proof of Lemma \ref{prelim:lem-Z-supermartingale-bound}, except that the process $(Z_m)_{m=0}^\infty$ is replaced by the shifted process $(Z_{j+m})_{m=0}^\infty$. Using similar arguments as in the proof of Lemma \ref{prelim:lem-Z-supermartingale-bound}, one can then show that 
\begin{equation}\label{prelim:eq-ST-supermartingales}
(H_{(\nu_j-j) \wedge m}(Z_{j+(\nu_j-j)\wedge m}^p))_{m=0}^\infty \qquad \text{and} \qquad  (Z_{j+(\nu_j-j)\wedge m}^r)_{m=0}^\infty
\end{equation}
are supermartingales with respect to the shifted filtration $(\mathcal{G}_{j+m})_{m=0}^\infty$. For any almost surely finite stopping time $M$ for the original filtration $(\mathcal{G}_m)_{m=0}^\infty$ satisfying $M\geq j$ almost surely, the difference $M-j$ is a stopping time for the shifted filtration $(\mathcal{G}_{j+m})_{m=0}^\infty$. From the optional stopping theorem for nonnegative supermartingales~\cite[Thm 4.8.4]{Durrett19}, it then follows that 
\begin{equation}\label{prelim:eq-ST-supermartingale-bound-1}
\E \bigl[ H_{(\nu_j \wedge M)-j}(Z_{\nu_j \wedge M}^p) \big| ~ \mathcal{G}_j \bigr] = \E \bigl[ H_{(\nu_j-j) \wedge (M-j)}(Z_{j+(\nu_j-j)\wedge (M-j)}^p) \big| ~ \mathcal{G}_j \bigr] \leq H_0(Z_j^p)
\end{equation}
and 
\begin{equation}\label{prelim:eq-ST-supermartingale-bound-1prime}
\E \bigl[ Z_{\nu_j \wedge M}^r \big| ~ \mathcal{G}_j \bigr]
=\E \bigl[ Z_{j+ (\nu_j-j) \wedge (M-j)}^r \big| ~ \mathcal{G}_j \bigr] \leq Z_j^r. 
\end{equation}
By using \eqref{prelim:eq-ST-supermartingale-bound-1}--\eqref{prelim:eq-ST-supermartingale-bound-1prime} and following the same argument as in the derivation of \eqref{prelim:eq-single-excursion}, we further obtain that
\begin{align}
\E \bigl[ \ind\{ \nu_j > M \} Z_M^r \big| ~ \mathcal{G}_j \bigr]  &\leq Z_j^r, \label{prelim:eq-ST-supermartingale-bound-2} \\ 
\E \bigl[ \ind\{ \nu_j > M \} (M-j+1)^{\frac{p-r}{\alpha}} Z_M^r \big| ~ \mathcal{G}_j \bigr]  &\leq B Z_j^p.  \label{prelim:eq-ST-supermartingale-bound-3}  
\end{align}
We note that \eqref{prelim:eq-ST-supermartingale-bound-3} with $M=j+n$, where $n\geq 0$, implies that 
\begin{equation*}
\mathbb{P} \bigl( \nu_j > j+n \big| ~ \mathcal{G}_j \bigr) \leq \frac{BZ_j^p}{(n+1)^{\frac{p-r}{\alpha}}}. 
\end{equation*}
From this, it follows that $\nu_j$ is almost surely finite for all $j\geq 0$. By taking the intersection of the events $\{ \nu_j < \infty\}$ over all $j\geq 0$, we then obtain that 
\begin{equation}\label{prelim:eq-ST-returns-below-R}
\mathbb{P}\bigl( Z_m \leq R ~~\textup{infinitely often} \bigr) =1.
\end{equation}

\proofstep{Step 2: Estimates of physical times.} 
In this step, we examine the physical time increments $(\tau_{m})_{m=1}^\infty$ as well as the physical times $(t_m)_{m=0}^\infty$, the latter of which are defined as 
\begin{equation}\label{prelim:eq-physical-times}
t_0 \eqdef 0 \qquad \text{and} \qquad t_m \eqdef \sum_{k=0}^{m-1} \tau_{k+1} \quad \text{for all } m \geq 1.
\end{equation} 
We emphasise that, since  $\tau_{k+1}$ is $\mathcal{G}_{k+1}$-measurable for all $k\geq 0$, $t_m$ is $\mathcal{G}_m$-measurable. 
We first prove a conditional lower bound on $\tau_{k+2}$ on the event $\{ Z_k \leq R \}$. On this event, our assumption \eqref{prelim:eq-time-increment-moments} implies that $\E [ \tau_{k+2}^{-q_1}|~ \mathcal{G}_k ] 
\leq \E[ X_{k+1}^{q_2}|~\mathcal{G}_k]
+ K$. By using the parameter condition $q_2 \leq p$ and our earlier estimate~\eqref{prelim:eq-Z-bound-small-values}, the remaining conditional moment can be bounded by 
\begin{equation*} 
\E\bigl[ X_{k+1}^{q_2} \big|~\mathcal{G}_k \bigr]
\leq \E\bigl[ Z_{k+1}^{q_2} \big|~\mathcal{G}_k \bigr] 
\leq \E\bigl[ Z_{k+1}^{p} \big|~\mathcal{G}_k \bigr]^{\frac{q_2}{p}}
\leq A^{q_2}. 
\end{equation*} 
For $h\eqdef 2^{-\frac{1}{q_1}} (A^{q_2}+K)^{-\frac{1}{q_1}}$, it then follows from Markov's inequality that, on the same event, 
\begin{equation}\label{prelim:eq-taukp1-lower-bound}
  \mathbb{P}\bigl( \tau_{k+2} \geq  h \big| ~ \mathcal{G}_k \bigr) = 1-  
\mathbb{P}\bigl( \tau_{k+2} < h \big| ~ \mathcal{G}_k \bigr) 
\geq 1-  h^{q_1} (A^{q_2}+K) \geq \tfrac{1}{2}.
\end{equation}
Next, we show that the process $(Z_m)_{m=0}^\infty$ cannot visit the region below $R$ too often while the physical times $(t_m)_{m=0}^\infty$ stay in a fixed interval. To quantify this, we let $I\subseteq [0,\infty)$ be a bounded time interval and set 
\begin{equation}\label{prelim:eq-ST-visits}
N(I) \eqdef \sum_{k\geq 0} \ind\{ Z_k \leq R, \, t_k \in I \}. 
\end{equation}
We note that for each $k\geq 0$, the $k$th summand in \eqref{prelim:eq-ST-visits} is $\mathcal{G}_k$-measurable. Using Tonelli's theorem and~\eqref{prelim:eq-taukp1-lower-bound}, it then holds that
\begin{align*}
\E \bigl[& N(I) \big| ~ \mathcal{G}_0 \bigr] = \sum_{k=0}^\infty \E \bigl[ \ind\{ Z_k \leq R, \, t_k \in I \} \big| ~ \mathcal{G}_0 \bigr] 
\leq  2 \sum_{k=0}^{\infty}  \E \Bigl[ \ind\{ Z_k \leq R, \, t_k \in I \} \mathbb{P} \bigl( \tau_{k+2} \geq  h \big| ~ \mathcal{G}_k \bigr) \Big| ~ \mathcal{G}_0 \Bigr] \\ 
&=  2 \sum_{k=0}^{\infty}  \E \Bigl[ \ind\{ Z_k \leq R, \, t_k \in I, \, \tau_{k+2} \geq  h \} \Big| ~ \mathcal{G}_0 \Bigr] 
= 2   \E \Bigl[ \sum_{k=0}^{\infty} \ind\{ Z_k \leq R, \, t_k \in I, \, \tau_{k+2} \geq  h \} \Big| ~ \mathcal{G}_0 \Bigr]. 
\end{align*}
The intervals $[t_{k+1},t_{k+2})$ corresponding to indices for which $t_k \in I$ and $\tau_{k+2}\geq h$ are pairwise disjoint, have length at least $h$, and are contained in $[\inf(I),\sup(I)+2]$. Consequently, there are at most $(2+\abs{I})/h$ such indices. As a result, we obtain that 
\begin{equation}\label{prelim:eq-ST-visits-bound}
\E \bigl[ N(I) \big| ~ \mathcal{G}_0 \bigr] \leq 2 (2+\abs{I})/h.
\end{equation}
While our main motivation for proving \eqref{prelim:eq-ST-visits-bound} lies in \eqref{prelim:eq-ST-second-summand-p2}   below, it also implies that $t_m\rightarrow \infty$ almost surely as $m\rightarrow \infty$. Indeed, if there were any $b\geq 0$ such that $\lim_{m\rightarrow \infty} t_m \leq b$ with positive probability, then, together with \eqref{prelim:eq-ST-returns-below-R}, this would imply that $N([0,b])=\infty$ with positive probability, which would contradict \eqref{prelim:eq-ST-visits-bound}. 

\proofstep{Step 3: Properties of $M_T$.} Using our physical times $(t_m)_{m=0}^\infty$, we can write 
\begin{equation*}
M_T = \inf\{ m \geq 1 \colon t_m \geq T\}. 
\end{equation*}
Since $t_m$ is $\mathcal{G}_m$-measurable, it follows that $\{ M_T \leq m \}=\{t_m \geq T\} \in \mathcal{G}_m$, and hence $M_T$ is a stopping time. Furthermore, since $t_m \rightarrow \infty$ almost surely as $m\rightarrow \infty$, it follows that $M_T$ is almost surely finite. From the upper bound $\tau_{k+1}\leq 1$, it follows that $t_m \leq m$. Together with our description of $M_T$, this immediately implies that $M_T \geq T$. More generally, on the event $\{ M_T > m\}$, it holds that 
\begin{equation*}
T-t_m \leq t_{M_T}-t_m = \sum_{k=m}^{M_T-1} \tau_{k+1} \leq M_T-m.
\end{equation*}

\proofstep{Step 4: Bounds on excursions.} We first show that 
\begin{equation}\label{prelim:eq-ST-one-excursion}
\E \bigl[ \ind\{ M_T < \nu_0 \} Z_{M_T}^r \big| ~ \mathcal{G}_0 \bigr] 
\leq \min \biggl( Z_0^r, \frac{BZ_0^p}{(1+T)^{\frac{p-r}{\alpha}}} \biggr). 
\end{equation}
The bound by the first argument in \eqref{prelim:eq-ST-one-excursion} follows directly from \eqref{prelim:eq-ST-supermartingale-bound-2} with $M=M_T$ and $j=0$. To obtain the bound by the second argument in \eqref{prelim:eq-ST-one-excursion}, we use \eqref{prelim:eq-ST-supermartingale-bound-3} with $M=M_T$ and $j=0$. Together with $M_T+1\geq 1+T$, this implies that 
\begin{equation*}
(1+T)^{\frac{p-r}{\alpha}} \E \bigl[ \ind\{ M_T < \nu_0 \} Z_{M_T}^r \big| ~ \mathcal{G}_0 \bigr] 
\leq  \E \bigl[ \ind\{ M_T < \nu_0 \} (M_T+1)^{\frac{p-r}{\alpha}} Z_{M_T}^r \big| ~ \mathcal{G}_0 \bigr] 
\leq B Z_0^p.
\end{equation*}
Next, we show for all $k\geq 0$ that 
\begin{equation}\label{prelim:eq-ST-later-excursion}
\begin{aligned}
\E \Bigl[ \ind\{& Z_k \leq R, \, t_k < T \}\, \ind\{ M_T <\nu_{k+1} \} Z_{M_T}^r \Big| ~ \mathcal{G}_k \Bigr] \\ 
&\leq \frac{2^{\frac{p-r}{\alpha}}A^p B}{(1+T-t_k)^{\frac{p-r}{\alpha}}} \, \ind\{ Z_k \leq R, \, t_k < T \}. 
\end{aligned}
\end{equation}
On the event $\{ t_k < T\}$, it holds that $M_T-k\geq (k+1)-k \geq 1$ and $M_T -k \geq T-t_k$, which implies that $M_T - k \geq 2^{-1} (1+T-t_k)$. By combining this with \eqref{prelim:eq-ST-supermartingale-bound-3} for $M=M_T \vee (k+1)$ and $j=k+1$, we then obtain that 
\begin{align*}
&2^{-\frac{p-r}{\alpha}} (1+T-t_k)^{\frac{p-r}{\alpha}} \ind\{ Z_k \leq R, \, t_k < T \} \,  \E \Bigl[  \ind\{ M_T <\nu_{k+1} \} Z_{M_T}^r \Big| ~ \mathcal{G}_{k+1} \Bigr] \\ 
&\quad \leq \ind\{ Z_k \leq R, \, t_k < T \} \E \Bigl[ \ind\{ M_T <\nu_{k+1} \} (M_T-k)^{\frac{p-r}{\alpha}} Z_{M_T}^r \Big| ~ \mathcal{G}_{k+1} \Bigr] \\
&\quad \leq \ind\{ Z_k \leq R, \, t_k < T \} B Z_{k+1}^p. 
\end{align*}
After conditioning on $\mathcal{G}_k$ and using \eqref{prelim:eq-Z-bound-small-values}, we then obtain \eqref{prelim:eq-ST-later-excursion}. 

\proofstep{Step 5: Conclusion.} We first estimate $\E[Z_{M_T}^r|~\mathcal{G}_0]$, and then extend our estimate to $\E[X_{M_T+1}^r|~\mathcal{G}_0]$. We decompose the conditional expectation $\E[Z_{M_T}^r|~\mathcal{G}_0]$ according to whether the process stays above $R$ through time $M_T$, is above $R$ at time $M_T$ but visited the region below $R$ at some time $0\leq k < M_T$, or satisfies $Z_{M_T}\leq R$. From this, it follows that 
\begin{equation}\label{prelim:eq-ST-decomposition}
\begin{aligned}
\E \bigl[ Z_{M_T}^r \big| ~ \mathcal{G}_0 \bigr] 
&= \E \bigl[ \ind\{ M_T < \nu_0 \} Z_{M_T}^r \big| ~ \mathcal{G}_0 \bigr] \\  
&+ \sum_{k=0}^\infty \E \bigl[ \ind\{ Z_k \leq R, t_k < T \} \ind\{ M_T < \nu_{k+1} \} Z_{M_T}^r \big| ~ \mathcal{G}_0 \bigr] \\ 
&+ \E \bigl[ \ind\{ Z_{M_T}\leq R\} Z_{M_T}^r\big| ~ \mathcal{G}_0 \bigr]. 
\end{aligned}
\end{equation}
The first summand in \eqref{prelim:eq-ST-decomposition} can be bounded by using \eqref{prelim:eq-ST-one-excursion}, and the third summand in \eqref{prelim:eq-ST-decomposition} is trivially bounded by $R^r$. In order to bound the second summand in \eqref{prelim:eq-ST-decomposition}, we group the physical times $(t_k)_{k=0}^\infty$ into unit-sized intervals preceding $T$. To be more precise, for each $j\geq 0$, we let 
\begin{equation*}
I_j(T) \eqdef [T-j-1,T-j) \cap [0,\infty).
\end{equation*}
Then, together with \eqref{prelim:eq-ST-later-excursion}, we obtain that 
\begin{align}
&\sum_{k=0}^\infty \E \bigl[ \ind\{ Z_k \leq R, t_k < T \} \ind\{ M_T < \nu_{k+1} \}  Z_{M_T}^r \big| ~ \mathcal{G}_0 \bigr] \allowdisplaybreaks[4] \notag \\  
&=\sum_{k=0}^\infty \E \Bigl[ \E \bigl[ \ind\{ Z_k \leq R, t_k < T \} \ind\{ M_T < \nu_{k+1} \}  Z_{M_T}^r \big| ~ \mathcal{G}_k \bigr] \Big| ~ \mathcal{G}_0 \Bigr] \allowdisplaybreaks[4] \notag \\  
&= \sum_{j=0}^\infty \sum_{k=0}^\infty \E \Bigl[  \E \bigl[ \ind\{ Z_k \leq R, t_k \in I_j(T)\} \ind\{ M_T < \nu_{k+1} \}  Z_{M_T}^r \big| ~ \mathcal{G}_k \bigr] \Big| ~ \mathcal{G}_0\Bigr] \allowdisplaybreaks[4] \notag \\ 
&\leq 2^{\frac{p-r}{\alpha}}A^p B  \sum_{j=0}^\infty \sum_{k=0}^\infty \E \Bigl[ \ind\{ Z_k \leq R, t_k \in I_j(T)\} \frac{1}{(1+T-t_k)^{\frac{p-r}{\alpha}}} \Big| ~ \mathcal{G}_0 \Bigr]. \label{prelim:eq-ST-second-summand-p1}
\end{align}
On the event $\{ t_k \in I_j(T)\}$, it holds that $T-t_k > j$. By combining this with our definition of $N$ from~\eqref{prelim:eq-ST-visits} and our estimate \eqref{prelim:eq-ST-visits-bound}, we then obtain that 
\begin{equation}\label{prelim:eq-ST-second-summand-p2}  
\begin{aligned}
\eqref{prelim:eq-ST-second-summand-p1} 
&\leq 2^{\frac{p-r}{\alpha}}A^p B \sum_{j=0}^{\infty} \sum_{k=0}^\infty  \frac{1}{(1+j)^{\frac{p-r}{\alpha}}} \E \bigl[ \ind\{ Z_k \leq R, t_k \in I_j(T)\} \big| ~ \mathcal{G}_0  \bigr] \\
&\leq  2^{\frac{p-r}{\alpha}}A^p B \sum_{j=0}^{\infty} \frac{1}{(1+j)^{\frac{p-r}{\alpha}}} \E \bigl[ N(I_j(T))  \big| ~ \mathcal{G}_0  \bigr] \\
&\leq 2^{1+\frac{p-r}{\alpha}} (3h^{-1}) A^p B \sum_{j=0}^{\infty} \frac{1}{(1+j)^{\frac{p-r}{\alpha}}}. 
\end{aligned}
\end{equation}
As in Lemma \ref{prelim:lem-Z-supermartingale-bound}, we have that $(p-r)/\alpha>1$, and hence the remaining sum converges. All in all, we therefore obtain that 
\begin{equation*}
 \E \bigl[ Z_{M_T}^r  \big| ~ \mathcal{G}_0  \bigr] \leq \min \biggl( Z_0^r, \frac{BZ_0^p}{(1+T)^{\frac{p-r}{\alpha}}} \biggr) + \widetilde{C},
\end{equation*}
where $\widetilde{C}=\widetilde{C}(\eps,c,p,q_1,q_2,r,K)$ is a suitable constant. This concludes our estimates of $Z_{M_T}^r$, and it now only remains to transfer our estimates to $X_{M_T+1}^r$. In order to replace $M_T$ by $M_T+1$, we first recall \eqref{prelim:eq-Z-bound-large-values} and \eqref{prelim:eq-Z-bound-small-values}. Together with Hölder's inequality, we then obtain that $\E[Z_{m+1}^r |~\mathcal{G}_m]\leq Z_m^r +A^r$ for all $m\geq 0$. Since $M_T$ is a stopping time, it follows that
\begin{align*}
&\E \bigl[ Z_{M_T+1}^r  \big| ~ \mathcal{G}_0  \bigr] 
= \sum_{m=0}^\infty \E \bigl[ \ind\{ M_T =m \} Z_{m+1}^r \big| ~ \mathcal{G}_0 \bigr] 
= \sum_{m=0}^\infty \E \Bigl[ \ind\{ M_T =m \} \E \bigl[  Z_{m+1}^r \big|~ \mathcal{G}_m \bigr] \Big| ~ \mathcal{G}_0\Bigr] \\ 
&\leq \sum_{m=0}^\infty \E \bigl[ \ind\{ M_T =m \} (Z_{m}^r+A^r) \big| ~ \mathcal{G}_0 \bigr]  
= \E \bigl[ Z_{M_T}^r  \big| ~ \mathcal{G}_0  \bigr] + A^r. 
\end{align*}
By combining the above estimates with the fact that $X_{M_T+1}\leq Z_{M_T+1}$ and that, as in the proof of Proposition \ref{prelim:prop-X-supermartingale}, $Z_0 \leq X_0 +2 Y_0$, we obtain  
\begin{equation*}
 \E \bigl[ X_{M_T+1}^r  \big| ~ \mathcal{G}_0  \bigr] \leq \min \biggl( (X_0+2 Y_0)^r, \frac{B(X_0+2 Y_0)^p}{(1+T)^{\frac{p-r}{\alpha}}} \biggr) + \widetilde{C} +A^r. 
\end{equation*}
This readily implies the desired estimates of the supremum and limit superior in \eqref{prelim:eq-stopping-time-final-estimate}. 
\end{proof}

\section{Covariant parabolic estimates}\label{section:covariant-heat-kernel-estimates}

In this section, we state and prove various parabolic estimates that will be heavily used in the following sections. Because we will prove quite a few results, many of which are only tangentially connected, we will split this section into further subsections, such that the main results within each subsection will be closely related. If the reader wishes, they may choose to read only the statements of the main results of each subsection, before moving on to later sections. This way, the reader may more quickly get an idea of the analytic input that will be needed in the coming sections, and then later they can come back to fill in the details.

Before diving into the subsections, we first introduce the main players of the present section, and then prove some preliminary lemmas which will be frequently used in the coming subsections. Throughout this section, we work with a general dimension $d \geq 2$ unless otherwise specified. We also fix a vector space $\genV$ on which $\liegroup$ and $\frkg$ act (via representations $\Repr, \repr$), with the action $\repr$ of $\frkg$ being the linearisation of the action $\Repr$ of $\liegroup$. In keeping with our notational convention, we will typically not explicitly write $\Repr, \repr$. As usual, we assume that $\genV$ is equipped with an inner product such that the action of $\liegroup$ is orthogonal (and therefore that of $\frkg$ is skew-symmetric). One should have in mind that later we will take $\genV = \frkg$ or $\V$ (or direct sums of these vector spaces). All implicit constants may depend on this choice of~$\genV$. 

In the following, let $B \colon [0, T] \times \T^d \ra \frkg^d$ be a smooth time-varying connection one-form. We will later on take $B$ to be a solution to the Yang--Mills gradient flow, but for the current discussion and the initial results, this assumption will not be needed. In contrast to the notation of Section~\ref{section:preliminaries}, we typically no longer explicitly denote the dependence of the covariant derivatives $\covd$ and the curvature $F$ on the underlying connection $B$. We recall the setup discussed in Section~\ref{section:main-heat}. First,  consider two natural linear parabolic equations associated to $B$:
\begin{equs}
(\ptl_t - \Delta_B) \phi &= G, \label{eq:covariant-heat-equation} \\
(\ptl_t - \Delta_B) \psi_i &= 2 F_i^{~\ell} \psi_\ell + H_i, ~~ i \in [d]. \label{eq:covariant-heat-equation-with-curvature}
\end{equs}
In the first equation, $\phi, G$ are $\genV$-valued, while in the second equation, $\psi, H$ are $\genV^d$-valued. Recall also the component-free form of \eqref{eq:covariant-heat-equation-with-curvature}:
\begin{equs}
(\ptl_t - \Delta_B) \psi = 2 F \psi + H,
\end{equs}
and recall from Definition~\ref{kernel:def-kernels} the covariant heat kernels $\fp_B, \fq_B$ associated to \eqref{eq:covariant-heat-equation}--\eqref{eq:covariant-heat-equation-with-curvature}.

\begin{notation}
In the following, when we write $|\fp_B|$ (resp.\ $|\fq_B|$), these should be understood as the operator norm of $\fp_B$ (resp.\ $\fq_B$) as a linear operator on the space $\genV$ (resp.\ $\genV^d$).  
\end{notation}

Next, recall from Definition~\ref{def:Fnorm} the following norm on the curvature $F$, which will appear in many (if not all) of the main results of the coming subsections. Given an interval $I = [s_0, t_0] \sse [0, \infty)$, and $k \geq 0$, we have that 
\begin{equs}
\Fnorm{F}_{I, k} = \sum_{j=0}^k  \Bigl(\big\|(t - s_0)^{(j+1)/2} \covd^{(j)} F\big\|_{L_t^2 L_x^\infty(I)} + \big\|(t-s_0)^{j/2} \covd^{(j)} F \big\|_{L_t^1 L_x^\infty(I)}\Bigr).
\end{equs}
In particular, when $k = 0$, we have that
\begin{equs}
\Fnorm{F}_{I, 0} = \big\|(t - s_0)^{1/2}F\big\|_{L_t^2 L_x^\infty(I)} + \big\|F\big\|_{L_t^1 L_x^\infty(I)}.
\end{equs}
In the case $I = [0, t]$ is an interval starting at $0$, we will often write $\Fnorm{F}_{t, k}$ as a shorthand for $\Fnorm{F}_{[0, t], k}$.

\begin{remark}\label{remark:Fpert-controls-perturbation}
In the following, we will show that $\Fnorm{\place}$ controls the extent to which solutions of the covariant heat equations \eqref{eq:covariant-heat-equation}--\eqref{eq:covariant-heat-equation-with-curvature} behave like solutions to the usual heat equation. In particular, control on $\Fnorm{\place}$ implies that the additional terms involving curvature appearing in these equations can be treated perturbatively, hence the superscript ``pert''.
\end{remark}

We will apply Gronwall (Lemma~\ref{lemma:gronwall}) in many places, which will incur an exponential loss. We next set notation for this loss.

\begin{notation}[Gronwall loss]\label{notation:gronwall-loss}
Given an interval $I \sse \R^+$, let
\begin{equ}
\gronwall_I(F) \eqdef \exp\Big( \big\|F\big\|_{L_t^1 L_x^\infty(I)}\Big).
\end{equ}
When $I = [0, t]$, we may write $\gronwall_t(F)$ instead of $\gronwall_{[0, t]}(F)$. Due to implicit constants associated with the ``$\lesssim$'' symbol, we will often use this notation for $\gronwall_I(CF)$, where $C$ is the appropriate implicit constant. 
\end{notation}

Next, we discuss several preliminary results which will form the backbone for the arguments of the coming subsections. For the reader who just wishes to get an idea of the main results of this section, the ensuing results may be safely skipped.

\begin{proposition}[Covariant monotonicity formulas]\label{prop:monotonicity}
Let $G \colon [0, T] \times \T^d \ra \genV$, and let $\phi$ satisfy
\begin{equ}
(\ptl_t - \Delta_B) \phi = G.
\end{equ}
Furthermore, let $K \colon [0, T]\times \T^d \rightarrow [0, \infty)$ be a non-negative solution to the backward heat equation, i.e.
\begin{equ}
(\ptl_t + \Delta) K = 0.
\end{equ}
Then, it holds that 
\begin{equ}\label{eq:covariant-monotonicity-square}
\ptl_t \Big(K \frac{|\phi|^2}{2}\Big) = -(\Delta K) \frac{|\phi|^2}{2} + K \Delta \Big(\frac{|\phi|^2}{2} \Big) - K |\covd\phi|^2 + K \langle \phi, G\rangle_\genV.
\end{equ}
More generally, for any $p \geq 2$,
\begin{equs}
\ptl_t \Big(K \frac{|\phi|^{p}}{p} \Big) = -(\Delta K) \frac{|\phi|^{p}}{p} &+ K \Delta\Big(\frac{|\phi|^{p}}{p}\Big) \\
&-K \Big(\frac{p-2}{4} \Big| |\phi|^{\frac{p-4}{2}} \nabla (|\phi|^2)\Big|^2 + |\phi|^{p-2} |\covd\phi|^2 - |\phi|^{p-2} \langle \phi, G \rangle_\genV\Big).
\end{equs}
\end{proposition}
\begin{proof}
By the Bochner formula \eqref{eq:bochner}, we have that
\begin{equs}
\ptl_t \Big(\frac{|\phi|^2}{2}\Big) = \langle \phi, \ptl_t \phi \rangle_\genV &= \langle \phi, \Delta_B \phi \rangle_\genV + \langle \phi, G \rangle_\genV = \Delta\Big(\frac{|\phi|^2}{2}\Big) - |\covd \phi|^2 + \langle \phi, G \rangle_\genV,
\end{equs}
and the claim follows from the Leibniz rule (for the calculations for general $p \geq 2$, see e.g.\ the proof of \cite[Proposition 3.1]{BC24}).
\end{proof}

\begin{remark}
By integrating \eqref{eq:covariant-monotonicity-square} in space-time and using integration by parts to cancel the two terms which involve the Laplacian, we have that for any $t_0 \in [0, T]$,
\begin{equ}\label{eq:integrated-covariant-monotonicity-square}
\bigl\|K(t_0)^{1/2} \phi(t_0)\bigr\|_{L_x^2}^2 + 2 \bigl\|K^{1/2} \covd \phi \bigr\|_{L_{t,x}^2([0, t_0])}^2 = \bigl\|K^{1/2}(0) \phi(0)\bigr\|_{L_x^2}^2 + 2 \int_0^{t_0} \dt \int_{\T^d} \dx K \langle \phi, G \rangle_{\genV}.
\end{equ}
Later on, we will need the following special case of the above, where $K \equiv 1$, $G = \covd^k H_k$ for some $H \colon [0, T] \times \T^d \ra \genV^d$. In this case, integrating the last term in \eqref{eq:integrated-covariant-monotonicity-square} by parts, and then applying Hölder and Young, we may obtain
\begin{equ}\label{eq:integrated-covariant-monotonicity-square-ibp}
\bigl\|\phi(t_0)\big\|_{L_x^2}^2 + \bigl\|\covd \phi\bigr\|_{L_{t,x}^2([0, t_0])}^2 \leq \|\phi(0)\|_{L_x^2}^2 + \bigl\| H \bigr\|_{L_{t,x}^2([0, t_0])}^2.
\end{equ}
Here, we absorbed one of the terms arising from Young into the good term $\|\covd \phi\|_{L_{t,x}^2([0, t_0])}^2$ on the left-hand side, which is why this term only appears with a prefactor $1$, instead of the prefactor $2$ as in \eqref{eq:integrated-covariant-monotonicity-square}.
\end{remark}

The covariant monotonicity formula leads to the following estimates, which will be used shortly.

\begin{lemma}\label{lemma:monotonicity-formula-integrated-initial-bounds}
Let $G \colon [0, T] \times \T^d \ra \genV$, $\phi_0 \colon \T^d \ra \genV$, and let $\phi$ be the solution to
\begin{equ}
(\ptl_t - \Delta_B) \phi = G, \qquad
\phi(0) = \phi_0.
\end{equ}
Furthermore, let $K \colon [0, T]\times \T^d \rightarrow (0, \infty)$ be a positive solution to the backward heat equation. We have the bound
\begin{equ}\label{eq:integrated-heat-kernel-weighted-monotonicity}
\|K^{1/2} \phi\|_{L_t^\infty L_x^2([0, T])} +  \|K^{1/2} \covd \phi\|_{L_{t, x}^2([0, T])} \lesssim \|K^{1/2}(0) \phi_0\|_{L_x^2} + \sup_{t_0 \in [0, T]} \biggl| \int_0^{t_0} \int_{\T^d} \dt \dx\, K \langle \phi, G \rangle_\genV \biggr|^{1/2}.
\end{equ}
Additionally, for $\alpha > 0$, we have the following time-weighted estimate:
\begin{equs}\label{eq:integrated-time-and-heat-kernel-weighted-monotonicity}
\big\|t^{\tfrac{\alpha}{2}} K^{1/2} &\phi \big\|_{L_t^\infty L_x^2([0, T])} + \big\| t^{\tfrac{\alpha}{2}} K^{1/2} \covd \phi\big\|_{L_{t,x}^2([0, T])} \\
&\lesssim \sup_{t_0 \in [0, T]} \biggl| \int_0^{t_0} \int_{\T^d} \dt \dx\, t^\alpha K \langle \phi, G \rangle_\genV \biggr|^{1/2} + \alpha^{1/2} \big\|t^{\tfrac{\alpha-1}{2}} K^{1/2} \phi \big\|_{L_{t,x}^2([0, T])}.
\end{equs}
\end{lemma}

\begin{remark}
In most use cases of \eqref{eq:integrated-heat-kernel-weighted-monotonicity}--\eqref{eq:integrated-time-and-heat-kernel-weighted-monotonicity}, we will directly estimate the integral by the $L_{t, x}^1$ norm, e.g.\ in \eqref{eq:integrated-heat-kernel-weighted-monotonicity}, we will estimate
\begin{equ}
\sup_{t_0 \in [0, T]} \biggl| \int_0^{t_0} \int_{\T^d} \dt \dx\, K \langle \phi, G \rangle_\genV \biggr|^{1/2} \leq \bigl\|K \langle \phi, G \rangle_\genV \bigr\|_{L_{t,x}^1([0, T])}^{1/2}.
\end{equ}
We will typically insert this estimate without further comment.
\end{remark}
\begin{proof}
The first estimate follows directly from \eqref{eq:integrated-covariant-monotonicity-square}, upon taking a sup over $t_0 \in [0, T]$. The second estimate follows by integrating the evolution equation for $t^\alpha K |\phi|^2$, which may be computed using the Leibniz rule and \eqref{eq:covariant-monotonicity-square}. 
\end{proof}

Next, we prove the following very general estimate, which will be repeatedly used in the coming subsections. On a first reading, the reader can choose to set $H = 0$, and to ignore the last estimate \eqref{eq:double-time-weighted-estimate}. While the statement of the result may be quite long, the main idea in the proof (especially when $H = 0$) is to combine our previous estimates (Lemma~\ref{lemma:monotonicity-formula-integrated-initial-bounds}) with Hölder and kickback.

\begin{lemma}\label{lemma:higher-regularity-abstract}
Let $\phi, G \colon [0, T] \times \T^d \ra \genV^d$, $H \colon [0, T] \times \T^d \ra \genV$, $\phi_0 \colon \T^d \ra \genV^d$ satisfy
\begin{equ}[e:equationphi]
(\ptl_t - \Delta_B) \phi = \schO(F \phi) + G + \covd H, \qquad
\phi(0) = \phi_0.
\end{equ}
Let $K \colon [0, T] \times \T^d \ra (0, \infty)$ be a positive solution to the backward heat equation. 
For $\alpha > 0$, we have the time-weighted estimate:
\begin{equation}\label{eq:higher-regularity-time-weighted-abstract}
\begin{aligned}
\big\|t^{\alpha/2} K^{1/2} \phi \big\|_{L_t^\infty L_x^2([0, T])} &+ \big\|t^{\alpha/2} K^{1/2} \covd \phi \big\|_{L_{t, x}^2([0, T])} \\
&\lesssim \Big(\bigl\|t^{1/2} F \bigr\|_{L_t^2 L_x^\infty([0, T])} + \alpha^{1/2}\Big) \big\|t^{(\alpha-1)/2} K^{1/2} \phi\big\|_{L_{t, x}^2([0, T])} \\
&\quad + \big\|t^{\alpha/2} K^{1/2} G\big\|_{L_t^1 L_x^2([0, T])} + N_\alpha(H),
\end{aligned}
\end{equation}
where
\begin{equ}
N_\alpha(H) = \big\|t^{\alpha/2} K^{1/2} H\big\|_{L_{t, x}^2([0, T])} + \big\|t^{\alpha/2} (T-t)^{-1/2} K^{1/2} H \big\|_{L_t^1 L_x^2([0, T])} + \big\| t^{\alpha/2} |\Delta K|^{1/2} H \big\|_{L_t^1 L_x^2([0, T])} .
\end{equ}
If $K$ is constant, then one may further take $N_\alpha(H)$ to just be the first term on the right-hand side above. Assuming that $H = 0$, and letting $T_1 > T$, we also have the double time-weighted estimate:
\begin{equs}\label{eq:double-time-weighted-estimate}
\big\|t^{\alpha/2} &(T_1 - t)^{\alpha/2} K^{1/2} \phi\big\|_{L_t^\infty L_x^2([0, T])} + \big\|t^{\alpha/2} (T_1 - t)^{\alpha/2} K^{1/2} \covd \phi\big\|_{L_{t, x}^2([0, T])} \\
&+ \alpha^{1/2} \big\|t^{\alpha/2} (T_1 - t)^{(\alpha-1)/2} K^{1/2} \phi\|_{L_{t, x}^2([0, T])} \\ 
&\quad\lesssim   \alpha^{1/2} \big\|t^{(\alpha-1)/2} (T_1 - t)^{\alpha/2} K^{1/2} \phi \big\|_{L_{t, x}^2([0, T])}  \\
&\qquad+\big\|(T_1 - t)^{1/2} F\big\|_{L_t^2 L_x^\infty([0, T])} \big\|t^{\alpha/2} (T_1 - t)^{(\alpha-1)/2} K^{1/2} \phi\big\|_{L_{t, x}^2([0, T])} \\
&\qquad+   \big\|t^{\alpha/2} (T_1 - t)^{\alpha/2} K^{1/2} G\big\|_{L_t^1 L_x^2([0, T])} .
\end{equs}
\end{lemma}
\begin{proof}
In the following, all space-time norms are on $[0, T] \times \T^d$. By the monotonicity formula \eqref{eq:integrated-time-and-heat-kernel-weighted-monotonicity} (with the $G$ therein taken to be the right-hand side of \eqref{e:equationphi}), we have that
\begin{equs}
\big\|t^{\alpha/2} K^{1/2} \phi \big\|_{L_t^\infty L_x^2} &+ \big\|t^{\alpha/2} K^{1/2} \covd \phi\big\|_{L_{t, x}^2} \lesssim \big\|t^\alpha K |F| |\phi|^2\big\|_{L_{t, x}^1}^{1/2} + \big\|t^\alpha K \langle G, \phi \rangle_\genV \big\|_{L_{t, x}^1}^{1/2} \\
&+ \sup_{t_0 \in [0, T]} \biggl|\int_0^{t_0} \int_{\T^d} \dt \dx\, t^\alpha K \langle \covd H, \phi \rangle_\genV \biggr|^{1/2} + \alpha^{1/2} \big\|t^{(\alpha-1)/2} K^{1/2} \phi\big\|_{L_{t, x}^2}.
\end{equs}
By Hölder and Young, we have that
\begin{equs}
\big\|t^\alpha K |F| |\phi|^2\big\|_{L_{t, x}^1}^{1/2} &\leq \big\|t^{1/2} F\big\|_{L_t^2 L_x^\infty}^{1/2} \big\|t^{\alpha/2} K^{1/2} \phi\big\|_{L_t^\infty L_x^2}^{1/2} \big\|t^{(\alpha-1)/2} K^{1/2} \phi\big\|_{L_{t, x}^2}^{1/2}  \\
&\leq \delta \big\|t^{\alpha/2} K^{1/2} \phi\big\|_{L_t^\infty L_x^2} + C_\delta \big\|t^{1/2} F \big\|_{L_t^2 L_x^\infty} \big\|t^{(\alpha-1)/2} K^{1/2} \phi \big\|_{L_{t, x}^2},
\end{equs}
where $\delta$ is some small constant. Similarly, we have that
\begin{equ}
\big\|t^\alpha K \langle G, \phi \rangle_\genV \big\|_{L_{t, x}^1}^{1/2} \leq \delta \big\|t^{\alpha/2} K^{1/2} \phi \big\|_{L_t^\infty L_x^2} + C_\delta \big\|t^{\alpha/2} K^{1/2} G \big\|_{L_t^1 L_x^2}. 
\end{equ}
Combining, we thus obtain the intermediate estimate that for some constant $C$,
\begin{equs}\label{eq:higher-regularity-time-weighted-abstract-intermediate}
~\hspace{-5mm}\big\|t^{\alpha/2} K^{1/2} \phi &\big\|_{L_t^\infty L_x^2} + \big\|t^{\alpha/2} K^{1/2} \covd \phi \big\|_{L_{t, x}^2} \leq C \Big(\big\|t^{1/2} F \big\|_{L_t^2 L_x^\infty} + \alpha^{1/2}\Big) \big\|t^{(\alpha-1)/2} K^{1/2} \phi \big\|_{L_{t, x}^2} \\
& + C \big\|t^{\alpha/2} K^{1/2} G \big\|_{L_t^1 L_x^2} + C  \sup_{t_0 \in [0, T]} \biggl|\int_0^{t_0} \int_{\T^d} \dt \dx \, t^\alpha K \langle \covd H, \phi \rangle_\genV \biggr|^{1/2}.
\end{equs}
To handle the last term, for $t_0 \in [0, T]$, we first integrate by parts to estimate
\begin{equ}
C \biggl|\int_0^{t_0} \int_{\T^d} \dt \dx \, t^\alpha K \langle \covd H, \phi \rangle_\genV \biggr|^{1/2} \leq C \big\|t^\alpha K \langle H, \covd \phi \rangle_\genV \big\|_{L_{t, x}^1}^{1/2} + C \big\|t^\alpha \nabla K \langle H, \phi \rangle_\genV \big\|_{L_{t, x}^1}^{1/2}.
\end{equ}
For the first term on the right-hand side, again by Hölder and Young, we may estimate
\begin{equ}
C \big\|t^\alpha K \langle H, \covd \phi \rangle_\genV \big\|_{L_{t, x}^1}^{1/2} \leq \delta \big\|t^{\alpha/2} K^{1/2} \covd \phi \big\|_{L_{t, x}^2} + C_\delta \big\|t^{\alpha/2} K^{1/2} H \big\|_{L_{t, x}^2}.
\end{equ}
Taking $\delta \leq 1/10$ (say), the first term above may be absorbed back into the left-hand side of \eqref{eq:higher-regularity-time-weighted-abstract-intermediate}, while the second term is already acceptable. For the other term (the one involving $\nabla K$), first note that if $K$ is constant, then $\nabla K = 0$, and thus the term drops out. Otherwise, in general, note that since $\T^d$ has non-negative Ricci curvature (in fact it has zero Ricci curvature), the function $K$, which is a positive solution to the backward heat equation, satisfies the following Li--Yau Harnack inequality (see \cite[(2.3)]{Ni2008}):
\begin{equ}
\frac{|\nabla K|^2}{K} \leq \Delta K + \frac{d}{2(T - t)} K.
\end{equ}
We may thus estimate
\begin{equs}
C \big\|t^\alpha & \nabla K \langle H, \phi \rangle_\genV \big\|_{L_{t, x}^1}^{1/2} \\
&\leq C \big\|t^{\alpha/2} K^{1/2} \phi\big\|_{L_t^\infty L_x^2}^{1/2} \Big\| t^{\alpha/2} \frac{|\nabla K|}{K^{1/2}} H \Big\|_{L_t^1 L_x^2}^{1/2} \\
&\leq  \delta \big\|t^{\alpha/2} K^{1/2} \phi\big\|_{L_t^\infty L_x^2} +  C_\delta \big\|t^{\alpha/2} |\Delta K|^{1/2} H\big\|_{L_t^1 L_x^2}  +  C_\delta \big\|t^{\alpha/2} (T-t)^{-1/2} K^{1/2} H\big\|_{L_t^1 L_x^2}.
\end{equs}
Taking $\delta \leq 1/10$ (say), the first term above may be absorbed back into the left-hand side of \eqref{eq:higher-regularity-time-weighted-abstract-intermediate}, while the last two terms are already acceptable.
Combining, the first desired estimate \eqref{eq:higher-regularity-time-weighted-abstract} now follows. The second desired estimate \eqref{eq:double-time-weighted-estimate} follows by a similar argument; we omit the proof.
\end{proof}

A small modification of the proof of \eqref{eq:higher-regularity-time-weighted-abstract} in Lemma~\ref{lemma:higher-regularity-abstract} also gives the following result. We avoid combining this with the statement of Lemma~\ref{lemma:higher-regularity-abstract}, because the statement of that lemma is already long enough as it is. In the following, recall our notation (Notation~\ref{notation:gronwall-loss}) for Gronwall loss.

\begin{lemma}\label{lemma:higher-regularity-abstract-small-Lt1-Lxinfty}
There is some constant $C_0$ such that the following holds. Let $\phi, G \colon [0, T] \times \T^d \ra \genV^d$, $H \colon [0, T] \times \T^d \ra \genV$, $\phi_0 \colon \T^d \ra \genV^d$ satisfy \eqref{e:equationphi}.
Let $K \colon [0, T] \times \T^d \ra (0, \infty)$ be a positive solution to the backward heat equation. Then
\begin{equation}\label{eq:higher-regularity-abstract-alternative}
\begin{aligned}
\big\|K^{1/2}& \phi \big\|_{L_t^\infty L_x^2([0, T])} + \big\|K^{1/2} \covd \phi \big\|_{L_{t, x}^2([0, T])} \\
&\lesssim \gronwall_T(C_0 F) \Big(\big\|K^{1/2}(0) \phi_0 \big\|_{L_x^2}   + \big\|K^{1/2} G\big\|_{L_t^1 L_x^2([0, T])} + N(H)\Big),
\end{aligned}
\end{equation}
where
\begin{equ}
N(H) = \big\|K^{1/2} H\big\|_{L_{t, x}^2([0, T])} + \big\|K^{1/2} (T-t)^{-1/2} H \big\|_{L_t^1 L_x^2([0, T])} + \big\| |\Delta K|^{1/2} H \big\|_{L_t^1 L_x^2([0, T])}.
\end{equ} 
\end{lemma}

\begin{remark}
The main difference between \eqref{eq:higher-regularity-abstract-alternative} and \eqref{eq:higher-regularity-time-weighted-abstract} is the lack of the $\|t^{1/2}F\|_{L_t^2 L_x^\infty}$ term on the right-hand side of the former. Instead, this term has been replaced by the Gronwall loss $\gronwall_T(C_0 F)$. The difference arises as follows. After an application of the monotonicity formula \eqref{eq:integrated-heat-kernel-weighted-monotonicity}, we obtain
\begin{equs}
\big\|K^{1/2} \phi\big\|_{L_t^\infty L_x^2}^2 + \big\|K^{1/2} \covd \phi\big\|_{L_{t, x}^2}^2 &\lesssim \big\|K |F| |\phi|^2\big\|_{L_{t, x}^1} + \text{other terms} \\
&\lesssim \Big\| \big\|F(t)\big\|_{L_x^\infty} \big\|K^{1/2} \phi(t)\big\|_{L_x^2}^2\Big\|_{L_t^1} +  \text{other terms}.
\end{equs}
(We previously treated the first term on the right-hand side by a  different application of Hölder.) At this stage, we may apply Gronwall (Lemma~\ref{lemma:gronwall}) and estimate ``other terms'' exactly the same as before, to arrive at \eqref{eq:higher-regularity-abstract-alternative}.
(There is also the fact that \eqref{eq:higher-regularity-abstract-alternative} does not include a time-weight, which is not as essential, but does result in a slight modification of the proof.)
\end{remark}

Recall that $\fp$ denotes the Euclidean heat kernel on $\T^d$.

\begin{lemma}\label{lemma:heat-kernel-weighted-to-L-infty}
For $T > 0$, let $f \colon [0, T] \times \T^d \ra \R$ be continuous. 
For any $(t, x) \in (0, T] \times \T^d$ and $s_1 \in [0, t)$, we have 
\begin{equ}
|f(t, x)| \leq \big\|\fp^{1/2}(s, y; t, x) f(s)\big\|_{L_s^\infty L_y^2([s_1, t])}.
\end{equ}
\end{lemma}
\begin{remark}
In the following, we will often apply this lemma with $f = |\phi|$, or more generally, $f = |\covd^{(k)} \phi|$, for some $\phi \colon [0, T] \times \T^d \ra \genV$ or $\phi \colon [0, T] \times \T^d \ra \genV^d$.
\end{remark}
\begin{proof}
It suffices to note that
$
\int \dy\, \fp(s, y; t, x) |f(s, y)|^2 \to |f(t, x)|^2
$
as $s\uparrow t$.
\end{proof}

In the remainder of the section, unless otherwise specified, we let $B \colon [0, T] \times \T^d \ra \frkg^d$ be a smooth solution to the Yang--Mills gradient flow, i.e.\
\begin{equ}
\ptl_t B_i = \covd^\ell F_{\ell i}, ~~ i \in [d].
\end{equ}
Next, we record several estimates involving $B$ and its curvature that will be useful in several places.

\begin{lemma}
\label{le:DF_bound}
Let $d = 2$, $T > 0$. We have that
\begin{equ}
    \|t^{1/2} \covd F \|_{L_t^\infty L_x^2([0, T])} \lesssim (1 + \|t^{1/2} F\|_{L_t^2 L_x^\infty([0, T])}) \|\covd F\|_{L_t^2 L_x^2([0, T])} \;.
\end{equ}
\end{lemma}
\begin{proof}
By Lemma~\ref{lemma:curvature-equation}, we have that in dimension $d = 2$, the curvature $F$ is a solution to the covariant heat equation. By the identity \eqref{eq:cov-heat-operator-covd-commutator}, we have that $(\ptl_t - \Delta_B - 2F) \covd F = 0$. The desired estimate then follows by \eqref{eq:higher-regularity-time-weighted-abstract} (with the $\alpha$ therein taken to be $1$, $K \equiv 1$, and $G, H \equiv 0$).
\end{proof}

\begin{remark}\label{remark:curvature-pointwise-in-time-bound-Fnorm}
Assume $d = 2$. In this case, recall by Remark~\ref{remark:curvature-subsolution-heat-equation} that $|F^{12}|$ is a subsolution to the heat equation, i.e.\ $(\ptl_t - \Delta) |F^{12}| \leq 0$. It follows that for all $1 \leq p \leq \infty$, the function $t \mapsto \bigl\|F(t)\bigr\|_{L_x^p}$ is non-increasing in $t$. In particular, we have that for any $t > 0$,
\begin{equ}[eq:curvature-pointwise-in-time-bound-Fnorm]
\bigl\|F(t)\bigr\|_{L_x^\infty} \leq \frac{1}{t} \int_0^t \bigl\|F(s)\bigr\|_{L_x^\infty} \, \ds \leq t^{-1} \Fnorm{F}_{t, 0}.
\end{equ}
\end{remark}

\begin{remark}\label{remark:Fpert-less-time-weight}
While $\Fnorm{\place}$ is sufficient to control the extent to which solutions of covariant heat equations behave like solutions to the usual heat equation (as discussed in Remark~\ref{remark:Fpert-controls-perturbation}), we will in fact control a slightly stronger version of this norm with less singular time weights. This will be needed in various places later on, for instance in deriving the heat kernel parametrix (see in particular Corollary~\ref{cor:pB-qB-parametrix}) in Section~\ref{section:heat-kernel-parametrix}. With this in mind, we record here the following estimates which will be used in several places later. Let $0 < s_0 < t_0 \leq T$, and let $\alpha \in [0, 1]$. We have that
\begin{equs}
\bigl\|\ptl_t B \bigr\|_{L_t^1 L_x^\infty([s_0, t_0])} &\lesssim s_0^{-\alpha/2} \bigl\| t^{\alpha/2} \covd F \bigr\|_{L_t^1 L_x^\infty([0, T])} \\
&\leq s_0^{-\alpha/2} T^{-(1-\alpha)/2} \bigl(T^{(1-\alpha)/2} \Fnorm{t^{-(1-\alpha)/2}F}_{T, 1}\bigr) , \label{eq:ptl-t-B-estimate-under-more-stringent-Fpert}\\
\bigl\|F\bigr\|_{L_t^2 L_x^\infty([s_0, t_0])} &\leq s_0^{-\alpha/2} \bigl\|t^{\alpha/2} F \bigr\|_{L_t^2 L_x^\infty([0, T])} \\
&\leq s_0^{-\alpha/2} T^{-(1-\alpha)/2} \bigl(T^{(1-\alpha)/2} \Fnorm{t^{-(1-\alpha)/2}F}_{T, 1}\bigr), \label{eq:F-estimate-under-more-stringent-Fpert} \\
\Fnorm{F}_{[s_0, t_0], 1} &\lesssim  (t_0-s_0)^{(1-\alpha)/2} \Fnorm{t^{-(1-\alpha)/2} F}_{T, 1} \\
&\leq (T^{-1}(t_0-s_0))^{(1-\alpha)/2} \bigl(T^{(1-\alpha)/2} \Fnorm{t^{-(1-\alpha)/2}F}_{T, 1}\bigr). \label{eq:Fpert-estimate-under-more-stringent-Fpert}
\end{equs}
In the first two estimates, we use that $\| t^{\alpha/2} \covd F \|_{L_t^1 L_x^\infty([0, T])}$ and $\|t^{\alpha/2} F \|_{L_t^2 L_x^\infty([0, T])}$ are terms appearing in $\Fnorm{t^{-(1-\alpha)/2}F}_{T, 1}$. We also remark that to estimate the $\bigl\|F\bigr\|_{L_t^1 L_x^\infty([s_0, t_0])}$ norm appearing in the LHS of \eqref{eq:Fpert-estimate-under-more-stringent-Fpert}, we split into cases $s_0 \leq t_0/2$ and $s_0 > t_0/2$. In the former case, we use that 
\begin{equ}
\bigl\|F\bigr\|_{L_t^1 L_x^\infty([s_0, t_0])} \leq t_0^{(1-\alpha)/2} \bigl\|t^{-(1-\alpha)/2} F \bigr\|_{L_t^1 L_x^\infty([s_0, t_0])} \lesssim (t_0 - s_0)^{(1-\alpha)/2} \Fnorm{t^{-(1-\alpha)/2} F}_{T, 1},
\end{equ} 
and in the latter case, we first apply Hölder and then the first inequality in \eqref{eq:F-estimate-under-more-stringent-Fpert} to estimate 
\begin{equ}
\bigl\|F\bigr\|_{L_t^1 L_x^\infty([s_0, t_0])} \leq (t_0 - s_0)^{1/2} \bigl\|F\bigr\|_{L_t^2 L_x^\infty([s_0, t_0])} \lesssim (t_0 - s_0)^{(1-\alpha)/2} \Fnorm{t^{-(1-\alpha)/2} F}_{T, 1}.
\end{equ} 
One should think of the quantity $T^{(1-\alpha)/2} \Fnorm{t^{-(1 - \alpha)/2} F}_{T, 1}$ as a stronger version of $\Fnorm{\place}$, and we will indeed have control on it in our later arguments.
\end{remark}

\subsection{Heat kernel estimates}\label{subsection:covariant-heat-kernel-estimates}

In this subsection, we prove pointwise bounds on the covariant heat kernels $\fp_B, \fq_B$. The main point is that assuming control on $\Fnorm{\place}$ (recall Definition~\ref{def:Fnorm}), the covariant heat kernels behave as if they were the Euclidean heat kernel $\fp$. We proceed to state the main results, the first of which concerns $\fp_B$, and the second of which concerns $\fq_B$ (as well as the difference $\fq_B - \fp_B$).
Recall also that $B$ is assumed to be a solution to the Yang--Mills gradient flow.

\begin{remark}
Throughout this section, we will fix the base point $(s, y) = (0, 0)$ and obtain estimates on  $\fp_B(t, x) = \fp_B(0, 0; t, x)$ and $\fq_B(t, x) = \fq_B(0, 0; t, x)$. The estimates can then be easily translated to the setting of a general base point $(s, y)$. For instance, in the estimates, all appearances of $|x|$ should be replaced by $|x-y|$, all appearances of $t$ should be replaced by $(t-s)$, all appearances of $\Fnorm{F}_{t, k}$ should be replaced by $\Fnorm{F}_{[s, t], k}$, and all appearances of $\gronwall_t$ should be replaced by $\gronwall_{[s, t]}$.
\end{remark}


\begin{proposition}\label{prop:pB-ptwise-bounds}
There is a constant $c > 0$ such that, for any point $(t_0, x_0) \in (0, 1] \times \T^d$, we have the pointwise bounds:
\begin{equs}
|\fp_B(t_0, x_0)| &\leq \fp(t_0, x_0), \label{eq:zeroth-order-pB-ptwise-bd}\\
|\covd^{(k)} \fp_B(t_0, x_0)| &\lesssim_k \bigl(1 + (\Fnorm{F}_{t_0, k-1})^k\bigr) t_0^{-\frac{k+d}2} \exp(-c |x_0|^2/t_0).\qquad \label{eq:higher-order-pB-ptwise-bd}
\end{equs}
The second inequality holds for any $k \geq 1$. In particular, the constant $c$ does not depend on $k$.
\end{proposition}


\begin{proposition}\label{prop:qB-ptwise-bds}
There are constants $C_0, c$ such that the following holds. For any point $(t_0, x_0) \in (0, 1] \times \T^d$ and $k \geq 0$, the following estimate holds:
\begin{equ}
|\covd^{(k)} \fq_B(t_0, x_0)| \lesssim_k \gronwall_{t_0}(C_0 F) \bigl(1 + (\Fnorm{F}_{t_0, k})^{k}\bigr) t_0^{-\frac{k+d}{2}} \exp(-c |x_0|^2/t_0). \label{eq:qB-ptwise-bounds} 
\end{equ}
We also have the pointwise estimates on the difference between $\fq_B$ and $\fp_B$:
\begin{equ}\label{eq:qB-pB-difference-pointwise-bounds}
\bigl|\covd^{(k)} \bigl(\fq_B(t_0,x_0) - \fp_B(t_0,x_0)\bigr) \bigr| \lesssim_k \gronwall_{t_0}(C_0 F)\Fnorm{F}_{t_0, k} \bigl(1 + (\Fnorm{F}_{t_0, k})^{(k-1) \vee 0}\bigr) t_0^{-\frac{k+d}2} \exp(-c |x_0|^2/t_0).
\end{equ}
\end{proposition}

Having stated our main Propositions~\ref{prop:pB-ptwise-bounds} and~\ref{prop:qB-ptwise-bds}, we now proceed to their proofs. 
First, we prove the following intermediate proposition regarding higher regularity estimates for solutions to the covariant heat equation~\eqref{eq:covariant-heat-equation}. Our first main result (Proposition~\ref{prop:pB-ptwise-bounds}) will directly follow from this proposition.

\begin{proposition}[Higher regularity for solutions to the covariant heat equation]\label{prop:covariant-heat-equation-higher-regularity}
Let $\phi, G \colon [0, T] \times \T^d \ra \genV$ be such that
\begin{equ}
(\ptl_t - \Delta_B) \phi = G,
\end{equ}
and let $K \colon [0, T] \times \T^d \ra (0, \infty)$ be a positive solution to the backward heat equation. Then for any $k \geq 0$, we have that
\begin{equation}\label{eq:covariant-heat-equation-higher-regularity}
\begin{aligned}
\big\|&t^{k/2} K^{1/2} \covd^{(k)} \phi\big\|_{L_t^\infty L_x^2([0, T])} + \big\|t^{k/2} K^{1/2} \covd^{(k+1)} \phi \big\|_{L_{t, x}^2([0, T])} \\
&\lesssim_k \Big(1 + (\Fnorm{F}_{T,k-1})^k\Big) \Big(\big\|K^{1/2}(0) \phi(0)\big\|_{L_x^2} + \sum_{j=0}^k \big\|t^{j/2} K^{1/2} \covd^{(j)} G\big\|_{L_t^1 L_x^2([0, T])}\Big).
\end{aligned}
\end{equation}
\end{proposition}
\begin{proof}
In the following, all space-time norms are implicitly taken on $[0, T] \times \T^d$. We proceed by induction. First, note that the base case $k = 0$ follows by \eqref{eq:integrated-heat-kernel-weighted-monotonicity}, and estimating the term
\begin{equ}
\big\|K \langle \phi, G \rangle_\genV \big\|_{L_{t, x}^1}^{1/2} \leq \delta \big\|K^{1/2} \phi \big\|_{L_t^\infty L_x^2} + C_\delta \big\|K^{1/2} G \big\|_{L_t^1 L_x^2},
\end{equ}
where $\delta$ is a small enough constant so that the first term on the right-hand side may be absorbed by kickback. (Note that we do not need to define $\Fnorm{F}_{T,-1}$ since it is being raised to the power $0$.)
Next, suppose that for some $k_0 \geq 0$, the estimate \eqref{eq:covariant-heat-equation-higher-regularity} holds for all $0 \leq k \leq k_0$. We proceed to show that the estimate also holds for $k = k_0 + 1$. Towards this end, recall from Lemma~\ref{lemma:equations-for-higher-derivatives} (and Remark~\ref{remark:higher-derivatives-equation-with-ad-F-term}) that $\covd^{(k_0+1)} \phi$ satisfies the following schematic equation:
\begin{equs}
(\ptl_t - \Delta_B) \covd^{(k_0+1)} \phi &= \sum_{\substack{j_1 + j_2 = k_0 +1\\ j_1 \geq 0 \\ j_2 \geq 1}} \schO((\covd^{(j_1)} F) (\covd^{(j_2)} \phi)) + \covd^{(k_0 + 1)} G  \\
&= \schO(F \covd^{(k_0 + 1)} \phi) + R + \covd^{(k_0 + 1)} G,
\end{equs}
where $R$ consists of the summands corresponding to $1 \leq j_1, j_2 \leq k_0$. Next, by the time-weighted estimate \eqref{eq:higher-regularity-time-weighted-abstract} with $\alpha = k_0 + 1$ (applied with $H = 0$), we have that
\begin{equs}
\big\|&t^{(k_0 + 1)/2} K^{1/2} \covd^{(k_0  +1)} \phi \big\|_{L_t^\infty L_x^2} + \bigl\|t^{(k_0+1)/2} K^{1/2} \covd^{(k_0 + 2)} \phi\bigr\|_{L_{t,x}^2} \\
&\lesssim \Big(\big\|t^{1/2} F\big\|_{L_t^2 L_x^\infty} + 1\Big)\big\|t^{k_0/2} K^{1/2} \covd^{(k_0 + 1)} \phi\big\|_{L_{t, x}^2} + \big\|t^{(k_0 + 1)/2} K^{1/2} R\big\|_{L_t^1 L_x^2} \\
&\quad \quad + \big\|t^{(k_0 + 1)/2} K^{1/2} \covd^{(k_0 + 1)} G\big\|_{L_t^1 L_x^2}.
\end{equs}
The first term on the right-hand side above may be directly handled by our inductive assumption, and the last term is already acceptable.
To estimate the second term, first note that for fixed $1 \leq j_1, j_2 \leq k_0$ such that $j_1 + j_2 = k_0 + 1$, we have by the inductive assumption that
\begin{equs}
\big\|t^{(k_0 + 1)/2} K^{1/2} &(\covd^{(j_1)} F) \covd^{(j_2)} \phi\big\|_{L_t^1 L_x^2} \\
&\lesssim \big\|t^{(j_1+1)/2} \covd^{(j_1)} F\big\|_{L_t^2 L_x^\infty} \big\|t^{(j_2 - 1)/2} K^{1/2} \covd^{(j_2)} \phi \big\|_{L_{t, x}^2} \\
&\lesssim \Big(\Fnorm{F}_{T, k_0} + (\Fnorm{F}_{T, k_0})^{k_0 + 1}\Big)  \Big(\big\|K^{1/2}(0) \phi(0)\big\|_{L_x^2} + \sum_{j=0}^{k_0} \big\|t^{j/2} K^{1/2} \covd^{(j)} G\big\|_{L_t^1 L_x^2}\Big) .
\end{equs}
The estimate of the term involving $R$ now follows, since $R$ consists of a finite sum involving terms of the above form. Now by combining the previous few estimates, we obtain the inductive step, i.e.\ we have shown \eqref{eq:covariant-heat-equation-higher-regularity} for $k = k_0 + 1$. 
\end{proof}

We now use Proposition~\ref{prop:covariant-heat-equation-higher-regularity} to prove our first main result (Proposition~\ref{prop:pB-ptwise-bounds}).

\begin{proof}[Proof of Proposition~\ref{prop:pB-ptwise-bounds}]
The first estimate \eqref{eq:zeroth-order-pB-ptwise-bd} is a consequence of the diamagnetic inequality \eqref{eq:parabolic-diamagnetic}. Next, we show the higher order estimates \eqref{eq:higher-order-pB-ptwise-bd}. Let $k \geq 1$. 
Fix $(t_0, x_0)$ and let $K(t, x) \coloneq \fp(t, x; t_0, x_0)$. Let $I \coloneq [\tfrac{t_0}{2}, t_0]$ and note that $\fp_B$ is a solution to the covariant heat equation on $I \times \T^d$. Applying Proposition~\ref{prop:covariant-heat-equation-higher-regularity} (with the forcing term $G = 0$), we obtain
\begin{equ}\label{eq:pB-ptwise-bounds-p1}
\big\|(t - \tfrac{t_0}{2})^{k/2} K^{1/2} \covd^{(k)} \fp_B\big\|_{L_t^\infty L_x^2(I)} \lesssim \bigl(1 + (\Fnorm{F}_{t_0, k-1})^k\bigr) \big\|(K^{1/2}\fp_B)(\tfrac{t_0}{2})\big\|_{L_x^2}.
\end{equ}
By Lemma~\ref{lemma:heat-kernel-weighted-to-L-infty}, the left-hand side of \eqref{eq:pB-ptwise-bounds-p1} can be bounded below by 
\begin{equation*}
\big\|(t - \tfrac{t_0}{2})^{k/2} K^{1/2} \covd^{(k)} \fp_B\big\|_{L_t^\infty L_x^2(I)} 
\gtrsim t_0^{k/2} \bigl| \covd^{(k)} \fp_B \bigr|(t_0,x_0). 
\end{equation*}
By using the zeroth order estimate \eqref{eq:zeroth-order-pB-ptwise-bd}, the estimate $\fp(t, x) \lesssim t^{-d/2}$, and the semigroup property of the heat kernel, the $K^{1/2}\fp_B$-term on the right-hand side of \eqref{eq:pB-ptwise-bounds-p1} can be bounded above by  
\begin{align*}
\big\|(K^{1/2}\fp_B)(\tfrac{t_0}{2})\big\|_{L_x^2} 
&\lesssim \big\|\fp(\tfrac{t_0}{2},x;t_0,x_0)^{\frac{1}{2}} \fp(\tfrac{t_0}{2},x)\big\|_{L_x^2} \\ 
&\lesssim t_0^{-d/4} \big\|\fp(\tfrac{t_0}{2},x;t_0,x_0)^{\frac{1}{2}} \fp(\tfrac{t_0}{2},x)^{\frac{1}{2}} \big\|_{L_x^2} \\ 
&\lesssim t_0^{-d/4} \fp(t_0,x_0)^{\frac{1}{2}}. 
\end{align*}
Combining, we obtain the desired estimate \eqref{eq:higher-order-pB-ptwise-bd}.
\end{proof}

In dimension $d = 2$, since the curvature $F$ itself is a solution to the covariant heat equation \eqref{eq:covariant-heat-equation} (recall Lemma~\ref{lemma:curvature-equation}), we may also apply Proposition~\ref{prop:covariant-heat-equation-higher-regularity} to obtain higher regularity estimates of $F$. This gives the following corollary, which estimates higher regularity norms of $F$ in terms of the zeroth order norm $\Fnorm{F}_{T, 0}$.

\begin{corollary}\label{cor:curvature-higher-regularity-estimates}
Let $d = 2$. Let $T \leq 1$. For all $k \geq 0$, we have that
\begin{equ}
\Fnorm{F}_{T, k} \lesssim_k \Fnorm{F}_{T, 0} + (\Fnorm{F}_{T, 0})^{(k+1)!}.
\end{equ}
\end{corollary}
\begin{remark}
Since we will only ever need this corollary for a bounded range of $k$ (e.g., $k \leq 10$), we are not worried about the factorial growth in $k$ of the upper bound.
\end{remark}
\begin{proof}
We proceed by induction. The base case $k = 0$ is trivial. Suppose that for some $k_0 \geq 0$, the result is true for all $0 \leq k \leq k_0$ and all $T \in (0, 1]$. We proceed to show that the result is true for $k = k_0 + 1$ and all $T \in (0, 1]$. Towards this end, fix $T \in (0, 1]$. Let $t_0 \in (0, T]$. Since $F$ itself is a solution to the covariant heat equation (Lemma~\ref{lemma:curvature-equation}), we have that
\begin{equ}
\covd^{(k_0 + 1)} F(t_0) = \int_{\T^d} \dy \,\covd^{(k_0 + 1)} \fp_B(t_0/2, y; t_0, \place) F(t_0/2, y).
\end{equ}
Thus by \eqref{eq:higher-order-pB-ptwise-bd}, we have that (using that $\Fnorm{F}_{[t_0/2, t_0], k} \leq \Fnorm{F}_{T, k}$, which follows from the definition of $\Fnorm{\place}$)
\begin{equ}
\bigl\|\covd^{(k_0 + 1)} F(t_0)\bigr\|_{L_x^\infty} \lesssim \bigl(1 + \bigl(\Fnorm{F}_{T, k_0}\bigr)^{k_0+1} \bigr) t_0^{-(k_0 + 1)/2} \bigl\|F(t_0/2)\bigr\|_{L_x^\infty}.
\end{equ}
Since this is true for arbitrary $t_0 \in (0, T]$, we may obtain (using our inductive assumption)
\begin{equs}
\big\|t^{(k_0 + 2)/2} \covd^{(k_0 + 1)}  F\big\|_{L_t^2 L_x^\infty([0, T])} + \,&\big\|t^{(k_0 + 1)/2} \covd^{(k_0 + 1)} F \big\|_{L_t^1 L_x^\infty([0, T])} \\
&\lesssim \bigl(1 + (\Fnorm{F}_{T, k_0})^{k_0+1}\bigr) \Fnorm{F}_{T, 0} \\
&\lesssim \Fnorm{F}_{T, 0} + (\Fnorm{F}_{T, 0})^{(k_0+1)! \times (k_0 + 1) + 1}.
\end{equs}
The inductive step now follows upon noting that $(k_0+1)! \times (k_0 + 1) + 1 \leq (k_0 + 2)!$.
\end{proof}

Next, we proceed towards the proof of the second main result (Proposition~\ref{prop:qB-ptwise-bds}), which gives estimates for $\fq_B$. We begin by observing that the difference $\fq_B - \fp_B$ solves a certain covariant heat equation of the form \eqref{eq:covariant-heat-equation-with-curvature}:
\begin{equ}
(\ptl_t - \Delta_B) (\fq_B - \fp_B) = 2 F \fq_B = 2 F(\fq_B - \fp_B) + 2 F\fp_B.
\end{equ}
Motivated by this, we first prove a higher regularity estimate for solutions to the above equation.

\begin{proposition}[Higher regularity for covariant heat equation with curvature term]\label{prop:covariant-heat-equation-with-curvature-higher-regularity}
There exists a constant $C_0$ such that the following holds. Let functions $\phi, \psi, G \colon [0, T] \times \T^d \ra \genV^d$ be such that:
\begin{equs}[2]
(\ptl_t - \Delta_B) \phi &= G, &\quad \phi(0) &= \phi_0, \\
(\ptl_t - \Delta_B) \psi &= 2 F (\psi +\phi), &\quad \psi(0) &= \psi_0,
\end{equs}
and let $K \colon [0, T] \times \T^d \ra (0, \infty)$ be a positive solution to the backward heat equation. Then
\begin{equs}
~&\big\|K^{1/2} \psi\big\|_{L_t^\infty L_x^2([0, T])} + \big\|K^{1/2} \covd \psi\big\|_{L_{t, x}^2([0, T])} \\
&\lesssim \gronwall_T(C_0 F) \Big( \big\|K^{1/2}(0) \psi_0\big\|_{L_x^2} + \big\|F\big\|_{L_t^1 L_x^\infty([0, T])} \bigl( \big\|K^{1/2}(0) \phi_0\big\|_{L_x^2} + \big\|K^{1/2} G \big\|_{L_t^1 L_x^2([0, T])}\bigr)\Big).
\end{equs}
Moreover, for $k \geq 1$, we have that
\begin{equation}\label{eq:covariant-heat-equation-with-curvature-higher-regularity}
\begin{aligned}
&\big\|t^{k/2} K^{1/2} \covd^{(k)} \psi \big\|_{L_t^\infty L_x^2([0, T])} + \big\|t^{k/2} K^{1/2} \covd^{(k+1)} \psi \big\|_{L_{t, x}^2([0, T])} \lesssim_k \gronwall_T(C_0 F) \Bigl(1 + (\Fnorm{F}_{T, k})^{k-1}\Bigr) \\
&\quad \quad \times \Big(  \Bigl(1 + \Fnorm{F}_{T, k}\Bigr) \big\| K^{1/2}(0) \psi_0 \big\|_{L_x^2} +   \Fnorm{F}_{T, k} \Big(\big\|K^{1/2}(0) \phi_0 \big\|_{L_x^2} + N_k(G)\Big)\Big) ,
\end{aligned}
\end{equation}
where 
\begin{equ}
N_k(G) = \sum_{j=0}^{k-1} \big\|t^{j/2} K^{1/2} \covd^{(j)} G \big\|_{L_t^1 L_x^2([0, T])} .
\end{equ}
\end{proposition}
\begin{proof}
In the following, all space-time norms are implicitly taken on $[0, T] \times \T^d$. The argument will be very similar to the proof of Proposition~\ref{prop:covariant-heat-equation-higher-regularity}, except we incur a Gronwall loss in the first (i.e.\ $k = 0$) step of the estimates. To show the zeroth order estimate, by the estimate \eqref{eq:higher-regularity-abstract-alternative}, Hölder, and Proposition~\ref{prop:covariant-heat-equation-higher-regularity}, we have that
\begin{equation}\label{eq:covariant-heat-equation-curvature-intermediate-difference-zeroth-order-1}
\begin{aligned}
\big\|K^{1/2} &\psi\big\|_{L_t^\infty L_x^2} + \big\|K^{1/2} \covd \psi\big\|_{L_{t, x}^2} \\
&\lesssim  \gronwall_T(C_0 F) \Big(\big\|K^{1/2}(0) \psi_0\big\|_{L_x^2} + \big\|K^{1/2} F \phi\big\|_{L_t^1 L_x^2} \Big) \\
&\lesssim \gronwall_T(C_0 F) \Big(\big\|K^{1/2}(0)\psi_0\big\|_{L_x^2} +  \big\|F\big\|_{L_t^1 L_x^\infty} \Big(\big\|K^{1/2}(0)\phi_0 \big\|_{L_x^2} + \big\|K^{1/2} G\big\|_{L_t^1 L_x^2}\Big)\Big),
\end{aligned}
\end{equation}
which is acceptable. By possibly increasing the constant $C_0$, the above also implies that
\begin{equ}\label{eq:covariant-heat-equation-curvature-intermediate-difference-zeroth-order-2}
\big\|K^{1/2} \psi\big\|_{L_t^\infty L_x^2} + \big\|K^{1/2} \covd \psi\big\|_{L_{t, x}^2} \lesssim  \gronwall_T(C_0 F)\Big(\big\|K^{1/2}(0)\psi_0\big\|_{L_x^2} +  \big\|K^{1/2}(0)\phi_0 \big\|_{L_x^2} + \big\|K^{1/2} G\big\|_{L_t^1 L_x^2}\Big),
\end{equ}
which we will use later.

To prove \eqref{eq:covariant-heat-equation-with-curvature-higher-regularity}, we proceed by induction. For the base case $k = 1$, by Lemma~\ref{lemma:equations-for-higher-derivatives} and Remark~\ref{remark:higher-derivatives-equation-with-ad-F-term}, we have that $\covd \psi$ satisfies the following schematic equation:
\begin{equ}
(\ptl_t - \Delta_B) \covd \psi = \schO(F \covd \psi) + \schO((\covd F) \psi) + \schO(F \covd \phi + (\covd F) \phi). 
\end{equ}
By the time-weighted estimate \eqref{eq:higher-regularity-time-weighted-abstract} with $\alpha = 1$ (and $G$ as the right-hand side of the above equation), we have that
\begin{equs}
\bigl\|t^{1/2} K^{1/2} \covd \psi\bigr\|_{L_t^\infty L_x^2} &+ \bigl\|t^{1/2} K^{1/2} \covd^{(2)} \psi\bigr\|_{L_{t,x}^2} \lesssim \bigl(\bigl\|t^{1/2} F \bigr\|_{L_t^2 L_x^\infty} + 1 \bigr) \bigl\|K^{1/2} \covd \psi\bigr\|_{L_{t,x}^2} \\
&+ \bigl\|t^{1/2} K^{1/2} |\covd F| |\psi|\bigr\|_{L_t^1 L_x^2} + \bigl\|t^{1/2} K^{1/2} \bigl(|\covd F| |\phi| + |F| |\covd \phi|\bigr)\bigr\|_{L_t^1 L_x^2}.
\end{equs}
For the first term on the right-hand side above, we combine \eqref{eq:covariant-heat-equation-curvature-intermediate-difference-zeroth-order-1} and \eqref{eq:covariant-heat-equation-curvature-intermediate-difference-zeroth-order-2} to estimate
\begin{equs}
~&\bigl(\bigl\|t^{1/2} F \bigr\|_{L_t^2 L_x^\infty} + 1 \bigr) \bigl\|K^{1/2} \covd \psi\bigr\|_{L_{t,x}^2}  \\
&\lesssim \gronwall_T(C_0 F) \Big(  \Bigl(1 + \Fnorm{F}_{T, 1}\Bigr) \big\| K^{1/2}(0) \psi_0 \big\|_{L_x^2} +   \Fnorm{F}_{T, 1} \Big(\big\|K^{1/2}(0) \phi_0 \big\|_{L_x^2} + \bigl\|K^{1/2} G \bigr\|_{L_t^1 L_x^2}\Big)\Big),
\end{equs}
which is acceptable. For the second term, we apply Hölder to estimate
\begin{equ}
\bigl\|t^{1/2} K^{1/2} |\covd F| |\psi|\bigr\|_{L_t^1 L_x^2} \leq \bigl\|t^{1/2} |\covd F| \bigr\|_{L_t^1 L_x^\infty} \bigl\|K^{1/2} |\psi|\bigr\|_{L_t^\infty L_x^2} \leq \Fnorm{F}_{T, 1} \bigl\|K^{1/2} |\psi|\bigr\|_{L_t^\infty L_x^2},
\end{equ}
and then apply \eqref{eq:covariant-heat-equation-curvature-intermediate-difference-zeroth-order-2} to obtain an acceptable contribution. The last term which involves $\phi$ and $\covd \phi$ may be similarly estimated by Hölder and the $k = 0$ case of Proposition~\ref{prop:covariant-heat-equation-higher-regularity} to obtain an acceptable contribution. Combining, we thus obtain the base case $k = 1$ of \eqref{eq:covariant-heat-equation-with-curvature-higher-regularity}.

Next, we suppose that for some $k_0 \geq 1$, the desired estimate \eqref{eq:covariant-heat-equation-with-curvature-higher-regularity} is true for all $1 \leq k \leq k_0$. We proceed to show that it is also true for $k = k_0 + 1$. Towards this end, first note by Lemma~\ref{lemma:equations-for-higher-derivatives} and Remark~\ref{remark:higher-derivatives-equation-with-ad-F-term} that $\covd^{(k_0+1)} \psi$ satisfies the following schematic equation:
\begin{equs}
(\ptl_t - \Delta_B) \covd^{(k_0 + 1)} \psi &= \sum_{\substack{\ell + m = k_0 + 1 \\ \ell, m \geq 0}} \schO\bigl((\covd^{(\ell)} F)(\covd^{(m)} \psi + \covd^{(m)} \phi)\bigr) \\
&= \schO\bigl(F \covd^{(k_0 + 1)} \psi\bigr) + \schO\bigl( (\covd^{(k_0 + 1)} F) \psi\bigr) + R_1 + R_2,
\end{equs}
where
\begin{equ}
R_1 \coloneq  \sum_{\substack{\ell + m = k_0 + 1 \\ \ell, m \geq 1}} \schO\bigl( (\covd^{(\ell)} F)(\covd^{(m)} \psi)\bigr),\quad
R_2 \coloneq \sum_{\substack{\ell + m = k_0 + 1 \\ \ell,m \geq 0}} \schO\bigl((\covd^{(\ell)} F) (\covd^{(m)} \phi) \bigr) .
\end{equ}
By the time-weighted estimate \eqref{eq:higher-regularity-time-weighted-abstract} with $\alpha = k_0 + 1$ (and $G = \schO((\covd^{(k_0 + 1)} F)\psi) + R_1 + R_2$, $H = 0$), we have that
\begin{equs}
\big\|t^{(k_0 + 1)/2} K^{1/2} \covd^{(k_0 + 1)} \psi\big\|_{L_t^\infty L_x^2}& + \big\|t^{(k_0 + 1)/2} K^{1/2} \covd^{(k_0 +2)} \psi\big\|_{L_{t, x}^2} \\
&\lesssim \bigl(\big\|t^{1/2} F\big\|_{L_t^2 L_x^\infty} + 1 \bigr) \big\|t^{k_0/2} K^{1/2} \covd^{(k_0+1)} \psi\big\|_{L_{t, x}^2} \\
&\quad + \big\|t^{(k_0 + 1)/2} K^{1/2} |\covd^{(k_0 + 1)} F| |\psi|\big\|_{L_t^1 L_x^2} \\
&\quad + \big\|t^{(k_0 + 1)/2} K^{1/2} |R_1| \big\|_{L_t^1 L_x^2} + \big\|t^{(k_0 + 1)/2} K^{1/2} |R_2| \big\|_{L_t^1 L_x^2} .
\end{equs}
The first term above may be handled by the inductive assumption. To estimate the second term, we apply Hölder to obtain
\begin{equ}
\big\|t^{(k_0 + 1)/2} K^{1/2} |\covd^{(k_0 + 1)} F| |\psi|\big\|_{L_t^1 L_x^2} \leq \big\|t^{(k_0 + 1)/2} \covd^{(k_0 + 1)} F \big\|_{L_t^1 L_x^\infty} \big\|K^{1/2} \psi\big\|_{L_t^\infty L_x^2},
\end{equ}
and by \eqref{eq:covariant-heat-equation-curvature-intermediate-difference-zeroth-order-2}, this is acceptable. To estimate the term involving $R_1$, for $\ell, m \geq 1 $ such that $\ell + m = k_0 + 1$, we may apply Hölder to obtain
\begin{equs}
\big\|t^{(k_0+1)/2} K^{1/2} |\covd^{(\ell)} F | |\covd^{(m)} \psi|\big\|_{L_t^1 L_x^2} &\leq \big\|t^{(\ell + 1)/2} \covd^{(\ell)} F\big\|_{L_t^2 L_x^\infty} \big\|t^{(m-1)/2} K^{1/2}\covd^{(m)} \psi\big\|_{L_{t, x}^2} \\
&\lesssim \Fnorm{F}_{T, k_0+1} \big\|t^{(m-1)/2} K^{1/2} \covd^{(m)} \psi\big\|_{L_{t, x}^2}.
\end{equs}
Inserting the inductive assumption as well as the zeroth order estimate (note $1 \leq m \leq k_0$), we obtain an acceptable contribution. 
Finally, by similar arguments, the term involving $R_2$ may be estimated by Proposition~\ref{prop:covariant-heat-equation-higher-regularity} and Hölder; we omit the details.
\end{proof}

Using Proposition~\ref{prop:covariant-heat-equation-with-curvature-higher-regularity}, we may now prove our second main result (Proposition~\ref{prop:qB-ptwise-bds}).

\begin{proof}[Proof of Proposition~\ref{prop:qB-ptwise-bds}] 
We first show the difference estimate \eqref{eq:qB-pB-difference-pointwise-bounds}. Fix $(t_0, x_0)$ and let $\psi(t, x) = (\fq_B - \fp_B)(t, x)$, so that
\begin{equs}
(\ptl_t - \Delta_B) \psi = 2 F \fq_B = 2F\psi + 2 F \fp_B.
\end{equs}
Thus, $\psi$ satisfies the conditions of Proposition~\ref{prop:covariant-heat-equation-with-curvature-higher-regularity} (where we take $\phi = \fp_B$). We begin by showing the zeroth order estimate
\begin{equ}\label{eq:psi-zeroth-order-estimate-intermediate}
|\psi(t, x)| \lesssim \gronwall_t(C F) \big\|F\big\|_{L_t^1 L_x^\infty([0, t])} \fp(t, x) , \text{ for all $t \in (0, t_0], x \in \T^d$}.
\end{equ}
For $t \in (0, t_0]$, let
\begin{equ}
f(t) \eqdef \sup_{x \in \T^d} \frac{|\psi(t, x)|}{\fp(t, x)}.
\end{equ}
For $x \in \T^d$, by our heat kernel estimate \eqref{eq:zeroth-order-pB-ptwise-bd} on $\fp_B$ and the semigroup property of the heat kernel, we have that
\begin{equs}
|\psi(t, x)| &\lesssim \int \ds  \dy\, \fp(s, y; t, x)  |F(s, y)| |\psi(s, y)| + \int \ds \dy\, \fp(s, y; t, x) |F(s, y)|  \fp(s, y) \\
&\lesssim  \int \ds \dy \, \fp(s, y; t, x) \fp(s, y) f(s) |F(s, y)| + \fp(t, x) \big\|F\big\|_{L_s^1 L_y^\infty([0, t])} \\
&\lesssim \fp(t, x) \bigg(\int_0^t \ds \,  f(s) \big\|F(s)\big\|_{L_y^\infty} + \big\|F\big\|_{L_s^1 L_y^\infty([0, t])}\bigg).
\end{equs}
It follows that
\begin{equ}
f(t) \lesssim \int_0^t \ds f(s) \big\|F(s)\big\|_{L_y^\infty} + \big\|F\big\|_{L_s^1 L_y^\infty([0, t])},
\end{equ}
and the zeroth order estimate \eqref{eq:psi-zeroth-order-estimate-intermediate} now follows by Gronwall (Lemma~\ref{lemma:gronwall}); we remark that $f \in L_t^\infty([0, t_0])$ because of classical short-time heat kernel estimates (recall the background $B$ is assumed to be smooth). Combining \eqref{eq:psi-zeroth-order-estimate-intermediate} with our zeroth order estimate \eqref{eq:zeroth-order-pB-ptwise-bd} on $\fp_B$, we directly obtain the zeroth order estimate on $\fq_B$ (i.e.\ the $k = 0$ case of \eqref{eq:qB-ptwise-bounds}).

To obtain the higher order estimates, let $k \geq 1$, $I = [\tfrac{t_0}{2}, t_0]$, $x_0 \in \T^d$, $K(t, x) = \fp(t, x; t_0, x_0)$. By \eqref{eq:covariant-heat-equation-with-curvature-higher-regularity} from Proposition~\ref{prop:covariant-heat-equation-with-curvature-higher-regularity}, we have that 
\begin{equs}
\big\|&(t - \tfrac{t_0}{2})^{k/2} K^{1/2} \covd^{(k)} \psi\bigr\|_{L_t^\infty L_x^2(I)} \\
&\lesssim  \gronwall_{[t_0/2, t_0]}(C F) \Bigl( \bigl(1 + \bigl(\Fnorm{F}_{t_0, k}\bigr)^k\bigr) \big\|(K^{1/2} \psi)(\tfrac{t_0}{2})\bigr\|_{L_x^2} + \Fnorm{F}_{t_0, k} \bigl(1 + \bigl(\Fnorm{F}_{t_0, k}\bigr)^{k-1}\bigr) \big\|(K^{1/2} \fp_B)(\tfrac{t_0}{2})\bigr\|_{L_x^2} \Bigr) \\
&\lesssim \gronwall_{t_0}(C F) \Bigl( \bigl(1 + \bigl(\Fnorm{F}_{t_0, k}\bigr)^k\bigr) \bigl\|F\bigr\|_{L_t^1 L_x^\infty([0, t_0])} +  \Fnorm{F}_{t_0, k} \bigl(1 + \bigl(\Fnorm{F}_{t_0, k}\bigr)^{k-1}\bigr)\Bigr) \bigl\|K^{1/2} \fp(\tfrac{t_0}{2})\bigr\|_{L_x^2},
\end{equs}
where we used the zeroth order estimate \eqref{eq:psi-zeroth-order-estimate-intermediate} as well as the pointwise estimate \eqref{eq:zeroth-order-pB-ptwise-bd}. To obtain the desired inequality, we split into cases where $\Fnorm{F}_{t_0, k} > 1$ or $\Fnorm{F}_{t_0, k} \leq 1$. In the first case, we simply estimate $\bigl\|F\bigr\|_{L_t^1 L_x^\infty([0, t_0])} \leq \gronwall_{t_0}(F)$ (here we just used that $x \leq e^x$ and the definition of $\gronwall_{t_0}(F)$ (recall Notation~\ref{notation:gronwall-loss})) and then estimate the above by
\begin{equ}
\gronwall_{t_0}((C + 1) F) \Fnorm{F}_{t_0, k} \bigl(1 + \bigl(\Fnorm{F}_{t_0, k}\bigr)^{k-1}\bigr) \bigl\|K^{1/2} \fp(\tfrac{t_0}{2})\bigr\|_{L_x^2},
\end{equ} 
while in the latter case, we may directly obtain the estimate
\begin{equ}
\gronwall_{t_0}(C F) \Fnorm{F}_{t_0, k} \bigl(1 + \bigl(\Fnorm{F}_{t_0, k}\bigr)^{k-1}\bigr) \bigl\|K^{1/2} \fp(\tfrac{t_0}{2})\bigr\|_{L_x^2}.
\end{equ}
Having reached this point, the remainder of the argument is similar to the proof of Proposition~\ref{prop:pB-ptwise-bounds}; we omit the details. We conclude by noting that the difference estimate \eqref{eq:qB-pB-difference-pointwise-bounds} combined with our higher order estimates \eqref{eq:higher-order-pB-ptwise-bd} on $\fp_B$ directly imply the $k \geq 1$ case of the estimate \eqref{eq:qB-ptwise-bounds}.
\end{proof}

\begin{remark}
From the proof of Proposition~\ref{prop:qB-ptwise-bds}, we see that the constant $C_0$ in the proposition is at least one more than the corresponding constant in Proposition~\ref{prop:covariant-heat-equation-with-curvature-higher-regularity}. We take $C_0$ to be such a constant in the following.
\end{remark}

The following result is a direct corollary of Propositions~\ref{prop:covariant-heat-equation-higher-regularity} and~\ref{prop:covariant-heat-equation-with-curvature-higher-regularity}.

\begin{corollary}\label{cor:higher-regularity-with-curvature}
Let $C_0$ be the constant from Proposition~\ref{prop:qB-ptwise-bds}, which is at least one more than the corresponding constant in Proposition~\ref{prop:covariant-heat-equation-with-curvature-higher-regularity}. Let $\phi, G \colon [0, T] \times \T^d \ra \genV^d$ be such that:
\begin{equ}
(\ptl_t - \Delta_B - 2F) \phi = G, \quad \phi(0) = \phi_0.
\end{equ}
Let $K \colon [0, T] \times \T^d \ra (0, \infty)$ be a positive solution to the backward heat equation. Then for any $k \geq 0$, we have that
\begin{equs}
\big\|&t^{k/2} K^{1/2} \covd^{(k)} \phi\big\|_{L_t^\infty L_x^2([0, T])} + \big\|t^{k/2} K^{1/2} \covd^{(k+1)} \phi\big\|_{L_{t, x}^2([0, T])} \\
&\lesssim \gronwall_T(C_0 F) \Big(1 + \bigl(\Fnorm{F}_{T, k}\bigr)^{k}\Big) \Big(\big\|K^{1/2}(0) \phi_0\big\|_{L_x^2} + \sum_{j=0}^k \big\|t^{j/2} K^{1/2} \covd^{(j)} G\big\|_{L_t^1 L_x^2([0, T])} \Big).
\end{equs}
\end{corollary}
\begin{proof}
We split $\phi = \tilde{\phi} + \psi$, where
$(\ptl_t - \Delta_B) \tilde{\phi} = G$ with initial condition $\tilde{\phi}(0) = \phi_0$.
We then estimate $\tilde{\phi}$ by Proposition~\ref{prop:covariant-heat-equation-higher-regularity}, and $\psi$ by Proposition~\ref{prop:covariant-heat-equation-with-curvature-higher-regularity}. We remark that in the case $k = 0$, we obtain the estimate 
\begin{equs}
\bigl\|K^{1/2} &\phi\bigr\|_{L_t^\infty L_x^2([0, T])} + \bigl\|K^{1/2} \covd \phi\bigr\|_{L_{t,x}^2([0, T])} \\
& \lesssim\gronwall_T (C F) \bigl(1 + \bigl\|F\bigr\|_{L_t^1 L_x^\infty([0, T])}\bigr) \bigl(\bigl\|K^{1/2}(0) \phi_0\bigr\|_{L_x^2} + \bigl\|K^{1/2} G \bigr\|_{L_t^1 L_x^2([0, T])} \bigr),
\end{equs}
where $C$ is the constant from Proposition~\ref{prop:covariant-heat-equation-with-curvature-higher-regularity}. In this case, we estimate $\bigl\|F\bigr\|_{L_t^1 L_x^\infty([0, T])} \leq \gronwall_T(F)$ (which follows because $x \leq e^x$), which results in the Gronwall factor $\gronwall_T((C + 1) F)$, as claimed.
\end{proof}

We also record the following corollary of Lemma~\ref{lemma:heat-kernel-weighted-to-L-infty} and Proposition~\ref{prop:covariant-heat-equation-with-curvature-higher-regularity}.

\begin{corollary}\label{cor:heat-equations-comparison}
Let $C_0$ be the constant from Proposition~\ref{prop:covariant-heat-equation-with-curvature-higher-regularity}. Let $\alpha \in [0, 1)$. Let $\phi, \psi \colon [0, T] \times \T^d \ra \genV^d$ be such that $\phi(0) = \psi(0)$, and such that:
\begin{equs}
(\ptl_t - \Delta_B - 2F) \psi &=  0, \\
(\ptl_t - \Delta_B) \phi &= 0.
\end{equs}
Then for any $k \geq 0$, we have that
\begin{equ}
\label{e:heat-equations-comparison}
\bigl\|t^{(k-(1-\alpha))/2} \covd^{(k)} (\psi - \phi)\bigr\|_{L_{t,x}^\infty([0, T])} \lesssim \gronwall_T(C_0 F) \Fnorm{t^{-(1-\alpha)/2} F}_{T, k}  \Bigl(1 + \bigl(\Fnorm{F}_{T, k}\bigr)^{(k-1) \vee 0}\Bigr) \bigl\|\phi(0)\bigr\|_{L_x^\infty}.
\end{equ}
\end{corollary}
\begin{proof}
Let $t_0 \in (0, T]$. By combining Lemma~\ref{lemma:heat-kernel-weighted-to-L-infty} with Proposition~\ref{prop:covariant-heat-equation-with-curvature-higher-regularity} (with the $\psi$ therein taken to be $\psi - \phi$ and the $G$ therein taken to be $0$), we obtain
\begin{equ}
\bigl\|t_0^{k/2} \covd^{(k)}(\psi - \phi)(t_0)\bigr\|_{L_x^\infty} \lesssim \gronwall_T(C_0 F) \Fnorm{F}_{t_0, k} \Bigl(1 + \bigl(\Fnorm{F}_{T, k}\bigr)^{(k-1) \vee 0}\Bigr) \bigl\|\phi(0)\bigr\|_{L_x^\infty}.
\end{equ}
Next, we insert the following estimate into the above:
\begin{equ}
\Fnorm{F}_{t_0, k} \leq t_0^{(1-\alpha)/2} \Fnorm{t^{-(1-\alpha)/2} F}_{t_0, k} \leq t_0^{(1-\alpha)/2} \Fnorm{t^{-(1-\alpha)/2} F}_{T, k},
\end{equ}
and then we finish by taking a sup over $t_0 \in [0, T]$.
\end{proof}

\subsection{Heat kernel estimates with mixed derivatives}\label{section:mixed-derivative-heat-kernel-estimates}

In this section, we extend the covariant heat kernel estimates of Propositions~\ref{prop:pB-ptwise-bounds} and~\ref{prop:qB-ptwise-bds} to the case of mixed derivatives, i.e.\ the case where we differentiate in both spatial variables. Such an expression arises naturally when derivatives fall on the forcing term of a covariant heat equation (see \eqref{eq:dagger-fq}--\eqref{eq:dagger-fp}). We first introduce the following convenient notation.

\begin{notation}\label{notation:dagger}
We introduce the covariant derivatives $\covd_y^{\dagger}$ as follows. For any $\Psi \colon [0,T] \times \T^d \rightarrow \End(\genV)$ and $j \in [d]$, we define 
\begin{equation}\label{eq-cov-dagger}
\covd_y^{\dagger,j} \Psi (s,y) \eqdef - \partial_y^{\, j} \Psi (s,y) + \Psi(s,y)  B^j(s,y). 
\end{equation}
Here, $B^j(s,y)$ is interpreted as an element of $\End(\genV)$ using our representation $\repr$ of $\frkg$. We extend this to $\End(\genV^d)$-valued maps $\Psi$ in the usual componentwise fashion. Similar to Notation~\ref{notation:higher-derivatives}, we use $\covd_y^{\dagger, (k)}$ to denote any $k$th-order derivative of the form $\covd_y^{\dagger, j_1} \cdots \covd_y^{\dagger, j_k}$, for indices $j_1, \ldots, j_k \in [d]$. 
\end{notation}

The reason for introducing the symbol $\dagger$, and for not explicitly referring to $\covd_y^{\dagger}$ as an adjoint, is that we want to avoid any confusion with the Hermitian adjoint (or Hermitian conjugate) taken directly in $\End(\genV)$. 

We remark that $\covd_y^\dagger$ commutes with $\covd_x$ (as well as $\partial_t$) when acting on 
$\End(\genV)$ or $\End(\genV^d)$-valued maps. The reason is that, in their definitions, the Lie-algebra elements act on opposite sides of the endomorphism. Combining this with \eqref{eq-cov-dagger}, it directly follows for all $m, k \geq 0$, $0 \leq s < t$, $G \colon \T^d \rightarrow \genV^d$, and $H \colon \T^d \ra \genV$ that
\begin{equs}
\int \dy \,  \covd_x^{(k)}\fq_{B}(s,y;t,x) (\covd_y^{(m)} G)(y) &= \int \dy \,  (\covd_y^{\dagger, (m)} \covd_x^{(k)} \fq_{B})(s,y;t,x) G(y), \label{eq:dagger-fq}\\
\int \dy \,  \covd_x^{(k)}\fp_{B}(s,y;t,x) (\covd_y^{(m)} H)(y) &= \int  \dy \,  (\covd_y^{\dagger,(m)} \covd_x^{(k)} \fp_{B})(s,y;t,x) H(y). \label{eq:dagger-fp}
\end{equs}

The main result of this section provides estimates for the kernels appearing in the right-hand sides of \eqref{eq:dagger-fq}--\eqref{eq:dagger-fp}. We state it as follows.

\begin{proposition}
\label{prop:qB-mixed-derivatives-ptwise-bounds}
Let $C_0$ be the maximum of the corresponding constants in Lemma~\ref{lemma:higher-regularity-abstract-small-Lt1-Lxinfty} and Proposition~\ref{prop:qB-ptwise-bds}. For any $m, k \geq 0$, $0\leq s_0 < t_0$ such that $t_0 - s_0 \leq 1$, and points $x_0, y_0 \in \T^d$, the following estimates hold:
\begin{equs}
\Bigl|\covd_y^{\dagger, (m)} \covd^{(k)}_x \fp_B(s_0, y_0; t_0, x_0)&\Bigr| \lesssim \Big(1 + (\Fnorm{F}_{t_0, m})^{m}\Big) \Big(1 + (\Fnorm{F}_{[s_0, t_0], k})^{k}\Big) \\
&\times \min(s_0, t_0 - s_0)^{-\frac{m}{2}} (t_0 - s_0)^{-\tfrac{k+d}{2}} \exp\biggl(-c \frac{|x_0 - y_0|^2}{t_0 - s_0}\biggr), \\
\Bigl| \covd_y^{\dagger, (m)} \covd^{(k)}_x \fq_B(s_0, y_0; t_0, x_0)&\Bigr| \lesssim \gronwall_{[s_0, t_0]}(C_0 F)\Big(1 + (\Fnorm{F}_{t_0, m})^{m}\Big) \Big(1 + (\Fnorm{F}_{[s_0, t_0], k})^{k}\Big) \\
&\times \min(s_0, t_0 - s_0)^{-\frac{m}{2}} (t_0 - s_0)^{-\tfrac{k+d}{2}} \exp\biggl(-c \frac{|x_0 - y_0|^2}{t_0 - s_0}\biggr).
\end{equs}
\end{proposition}

In order to prove this proposition, we will need several intermediate results, which we discuss next. The starting point is to derive an evolution equation for $\fp_B, \fq_B$ in $s, y$ instead of $t, x$. Let $\fp_B^*, \fq_B^*$ denote the adjoints of $\fp_B, \fq_B$ as operators on $\genV, \genV^d$. The following lemma shows that the adjoints $\fp_B^*, \fq_B^*$ are in fact the fundamental solutions to the associated \textit{backward} covariant heat equations.

\begin{lemma}[Adjoints and backward heat equations]\label{lemma:adjoints-backwards-heat-equations}
For all $s_0 < t_0$, $x_0, y_0 \in \T^d$, we have that
\begin{equs}
(\ptl_s + \Delta_{B, y}) \fp_B^*(s_0, y_0; t_0, x_0) &= 0, \\
(\ptl_s + \Delta_{B, y} + 2F) \fq_B^*(s_0, y_0; t_0, x_0) &= 0.
\end{equs}
Here, the subscript $y$ in $\Delta_{B, y}$ denotes that the differentiation is in the $y$-variable.
\end{lemma}
\begin{proof}
By standard considerations, given $\psi_{t_0} \colon \T^d \ra \genV^d$, the function 
\begin{equ}
\psi(s, y) \eqdef \int \dx \, \fq_B^*(s, y; t_0, x) \psi_{t_0}(x)
\end{equ}
solves the adjoint equation
\begin{equ}
(\ptl_s + \mc{L}_s^*) \psi = 0,
\end{equ}
where $\mc{L}_s = \Delta_B + 2F$. The covariant Laplacian $\Delta_B$ is self-adjoint, and so is $F$ as an operator on $\genV^d$ (and thus also as an operator on the space of $L^2$ functions $\T^d \ra \genV^d$), since for $X, Y \in \genV^d$, we have that
\begin{equ}
\langle X, F Y \rangle_{\genV^d} = \langle X^i,  F_{ij} Y^j\rangle_\genV = - \langle F_{ij} X^i, Y^j \rangle_{\genV} = \langle F_{ji} X^i, Y^j \rangle_{\genV} = \langle F X, Y \rangle_{\genV^d}.
\end{equ}
Thus $\mc{L}_s^* = \mc{L}_s$, and the equation for $\fq_B^*$ now follows. The equation for $\fp_B^*$ follows similarly (in fact, it is simpler, since there is no $F$ term).
\end{proof}

Next, we relate the kernels estimated in Proposition~\ref{prop:qB-mixed-derivatives-ptwise-bounds} to the adjoints $\fp_B^*, \fq_B^*$.

\begin{lemma}\label{lemma:covd-dagger-adjoint-relation}
For any $m, k \geq 0$, $0 \leq s < t$, and points $x, y  \in \T^d$, we have that
\begin{equs}
\bigl(\covd_y^{\dagger, (m)} \covd_x^{(k)} \fp_B(s, y; t, x)\bigr)^* &= (-1)^m \covd_y^{(m)} \bigl(\covd_x^{(k)} \fp_B(s, y; t, x)\bigr)^*, \\
\bigl(\covd_y^{\dagger, (m)} \covd_x^{(k)} \fq_B(s, y; t, x)\bigr)^* &= (-1)^m \covd_y^{(m)} \bigl(\covd_x^{(k)} \fq_B(s, y; t, x)\bigr)^*.
\end{equs}
\end{lemma}
\begin{proof}
We show the identity for $\fq_B$, as the argument for $\fp_B$ is similar. Let $G, H \colon \T^d \ra \genV^d$. Applying the definition of the adjoint as well as integration by parts, we have that 
\begin{equs}
\int &\dx \int \dy\,  \big\langle H(x), \covd^{(k)}_x \fq_B(s, y; t, x) \covd^{(m)}_y G(y)\big\rangle_{\genV^d}  \\
&= (-1)^m \int \dx \int  \dy \, \big \langle \covd^{(m)}_y (\covd^{(k)}_x \fq_B(s, y; t, x))^* H(x), G(y)\big \rangle_{\genV^d} \\
&= (-1)^m \int  \dx \int  \dy \, \Big\langle H(x), \Big(\covd^{(m)}_y (\covd^{(k)}_x \fq_B(s, y; t, x))^*\Big)^* G(y) \Big\rangle_{\genV^d}.
\end{equs}
Since this is true for arbitrary $H$, the result now follows upon combining this with \eqref{eq:dagger-fq}.
\end{proof}

In summary, so far we have expressed the kernels in Proposition~\ref{prop:qB-mixed-derivatives-ptwise-bounds} in terms of solutions to backward covariant heat equations. We proceed to estimate such solutions. Since the time reversal of a solution to the backward covariant heat equation is a solution to the forward covariant heat equation, a natural approach would be to reduce to the estimates of Section~\ref{subsection:covariant-heat-kernel-estimates}. However, one needs to be careful here, because the norm $\Fnorm{\place}$ (recall Definition~\ref{def:Fnorm}) is not invariant under time reversal. This reflects the important fact that for us, the curvature $F$ may be singular at the initial time, but will become smoother at later times (after all, it solves a covariant heat equation \dash recall Lemma~\ref{lemma:curvature-equation}). However, after time reversal, the curvature will now possibly be singular at the final time, not the initial time. We therefore need to define a new norm which reflects this behaviour.

\begin{definition}[Norms for time reversals]\label{kernel:def-trep}
Let $[s_0, t_0]$ be an interval, and let $t_1 > t_0$. For $k \geq 0$, define the following norm
\begin{equ}
\revFnorm{F}_{[s_0, t_0], t_1, k} \eqdef \sum_{j=0}^k  \Bigl(\big\|(t_1 - t)^{(j+1)/2} \covd^{(j)} F\big\|_{L_t^2 L_x^\infty([s_0, t_0])} + \big\|(t_1 - t)^{j/2} \covd^{(j)} F\big\|_{L_t^1 L_x^\infty([s_0, t_0])} \Bigr).
\end{equ}
This norm will be needed when dealing with time reversals, hence the superscript ``trep'', which is the reversal of ``pert''.
\end{definition}

Next, we state the main intermediate result towards the proof of Proposition~\ref{prop:qB-mixed-derivatives-ptwise-bounds}, which should be considered the analogue of Propositions~\ref{prop:covariant-heat-equation-higher-regularity} and~\ref{prop:covariant-heat-equation-with-curvature-higher-regularity}, except now we use the norm $\revFnorm{\place}$ in place of $\Fnorm{\place}$.

\begin{lemma}\label{lemma:time-reversal-higher-regularity}
Let $C_0$ be the constant from Lemma~\ref{lemma:higher-regularity-abstract-small-Lt1-Lxinfty}. Let $\tilde{B} \colon [0, T] \times \T^d \ra \frkg^d$ be a solution to the backward Yang--Mills gradient flow, i.e.\
\begin{equ}
\ptl_t \tilde{B}_i = -\tilde{\covd}^\ell \tilde{F}_{\ell i}, ~~ i \in [d],
\end{equ}
where $\tilde{\covd} = \covd_{\tilde{B}}$ and $\tilde{F} = F_{\tilde{B}}$. Let $\phi \colon [0, T] \times \T^d \ra \genV^d$ and $\phi_0 \colon \T^d \ra \genV^d$ be such that
\begin{equ}
(\ptl_t - \Delta_{\tilde{B}} - 2\tilde{F}) \phi = 0, \qquad
\phi(0) = \phi_0,
\end{equ}
and let $K \colon [0, T] \times \T^d \ra (0, \infty)$ be a positive solution to the backward heat equation. Then for $T_1 > T$, $k \geq 0$, we have that
\begin{equs}
\big\|&t^{k/2} (T_1 - t)^{k/2} K^{1/2} \tilde{\covd}^{(k)} \phi\big\|_{L_t^\infty L_x^2([0, T])} + \big\|t^{k/2} (T_1 - t)^{k/2} K^{1/2} \tilde{\covd}^{(k+1)} \phi\big\|_{L_{t, x}^2([0, T])} \\
&\lesssim \gronwall_T(C_0 \tilde{F}) T_1^{k/2} \Big(1 + (\revFnorm{\tilde{F}}_{[0, T], T_1, k})^{k}\Big) \big\|K^{1/2}(0) \phi_0\big\|_{L_x^2} . 
\end{equs}
Similarly, if instead $\phi \colon [0, T] \times \T^d \ra \genV$, $\phi_0 \colon \T^d \ra \genV$ are such that
\begin{equ}
(\ptl_t - \Delta_{\tilde{B}}) \phi = 0, \qquad \phi(0) = \phi_0, 
\end{equ}
then for any $K, T_1, k$ as before, we have that
\begin{equs}
\big\|&t^{k/2} (T_1 - t)^{k/2} K^{1/2} \tilde{\covd}^{(k)} \phi\big\|_{L_t^\infty L_x^2([0, T])} + \big\|t^{k/2} (T_1 - t)^{k/2} K^{1/2} \tilde{\covd}^{(k+1)} \phi\big\|_{L_{t, x}^2([0, T])} \\
&\lesssim  T_1^{k/2} \Big(1 + (\revFnorm{\tilde{F}}_{[0, T], T_1, k})^{k}\Big) \big\|K^{1/2}(0) \phi_0\big\|_{L_x^2} . 
\end{equs}
\end{lemma}
\begin{proof}
We focus on the first estimate, as the proof of the second is quite similar. In the following, all space-time norms are implicitly on $[0, T] \times \T^d$. We proceed by induction. The base case $k = 0$ follows by Lemma~\ref{lemma:higher-regularity-abstract-small-Lt1-Lxinfty} (with the $G, H$ therein taken to be zero). Next, suppose that for some $k_0 \geq 0$, the result is true for all $0 \leq k \leq k_0$. We proceed to show that it also holds for $k_0 + 1$. Towards this end, by Lemma~\ref{lemma:equations-for-higher-derivatives} (and Remark~\ref{remark:higher-derivatives-equation-with-ad-F-term}), we have that
\begin{equ}
(\ptl_t - \Delta_{\tilde{B}}) \tilde{\covd}^{(k_0 + 1)} \phi = \schO(\tilde{F} \tilde{\covd}^{(k_0 + 1)} \phi) + \sum_{\substack{j_1 + j_2 = k_0 + 1 \\ j_1 \geq 1 \\ j_2 \geq 0}} \schO((\tilde{\covd}^{(j_1)} \tilde{F}) (\tilde{\covd}^{(j_2)} \phi)).
\end{equ}
For notational convenience, let $w_k(t) = t^{k/2} (T_1 - t)^{k/2}$. By the double time-weighted estimate \eqref{eq:double-time-weighted-estimate}, we have that
\begin{equs}
~&\big\|w_{k_0+1} K^{1/2} \tilde{\covd}^{(k_0 + 1)} \phi\big\|_{L_t^\infty L_x^2} + \big\|w_{k_0 + 1} K^{1/2} \tilde{\covd}^{(k_0 + 2)} \phi\big\|_{L_{t, x}^2} \\
&\lesssim  \big\|(T_1 - t)^{1/2} w_{k_0} K^{1/2} \tilde{\covd}^{(k_0+1)} \phi\big\|_{L_{t, x}^2} + \big\|(T_1 - t)^{1/2} \tilde{F}\big\|_{L_t^2 L_x^\infty} \big\|t^{1/2} w_{k_0} K^{1/2} \tilde{\covd}^{(k_0 + 1)} \phi\big\|_{L_{t, x}^2} \\
&+ \sum_{\substack{j_1 + j_2 = k_0 + 1 \\ j_1 \geq 1 \\ j_2 \geq 0}} \big\|w_{k_0+1} K^{1/2} |\tilde{\covd}^{(j_1)} \tilde{F}| |\tilde{\covd}^{(j_2)} \phi|\big\|_{L_t^1 L_x^2}.
\end{equs}
By the inductive assumption, we may estimate the first term on the right-hand side above by
\begin{equ}
\big\|(T_1 - t)^{1/2} w_{k_0} K^{1/2} \tilde{\covd}^{(k_0+1)} \phi\big\|_{L_{t, x}^2} \lesssim \gronwall_T(C_0 \tilde{F}) T_1^{1/2} T_1^{k_0/2} \Big(1 + (\revFnorm{\tilde{F}}_{[0, T], T_1, k_0})^{k_0}\Big) \big\|K^{1/2}(0) \phi_0\big\|_{L_x^2},
\end{equ}
which is acceptable. For the second term, again applying the inductive assumption, we have that
\begin{equs}
\big\|(T_1 - t)^{1/2}& \tilde{F}\big\|_{L_t^2 L_x^\infty} \big\|t^{1/2} w_{k_0} K^{1/2} \tilde{\covd}^{(k_0 + 1)} \phi\big\|_{L_{t, x}^2} \\
&\lesssim \gronwall_T(C_0 \tilde{F}) T^{1/2} T_1^{k_0/2} \revFnorm{\tilde{F}}_{[0, T], T_1, 0} \Big(1 + (\revFnorm{\tilde{F}}_{[0, T], T_1, k_0})^{k_0}\Big) \big\|K^{1/2}(0) \phi_0\big\|_{L_x^2},
\end{equs}
which is also acceptable. Finally, for the third term, for fixed $j_1 \geq 1$, $j_2 \geq 0$ such that $j_1 + j_2 = k_0 + 1$, we may estimate by Hölder and the inductive assumption
\begin{equs}
\big\|w_{k_0+1} & K^{1/2} |\tilde{\covd}^{(j_1)} \tilde{F}| |\tilde{\covd}^{(j_2)} \phi| \big\|_{L_t^1 L_x^2} \\
&\leq \big\|(T_1 - t)^{j_1/2} \tilde{\covd}^{(j_1)} \tilde{F}\big\|_{L_t^1 L_x^\infty} \big\| t^{j_1/2} w_{j_2} K^{1/2} \tilde{\covd}^{(j_2)} \phi\big\|_{L_t^\infty L_x^2} \\
&\lesssim \gronwall_T(C_0 \tilde{F}) T^{j_1/2} T_1^{j_2/2} \revFnorm{\tilde{F}}_{[0, T], T_1, k_0+1} \Big(1 + (\revFnorm{\tilde{F}}_{[0, T], T_1, k_0})^{k_0}\Big) \big\|K^{1/2}(0) \phi_0\big\|_{L_x^2},
\end{equs}
which is also acceptable. Combining, we have thus proven the inductive step.
\end{proof}

We now have everything in place to prove Proposition~\ref{prop:qB-mixed-derivatives-ptwise-bounds}.

\begin{proof}[Proof of Proposition~\ref{prop:qB-mixed-derivatives-ptwise-bounds}]
It suffices to assume $s_0 > 0$ in the following (since if $m > 0$ then the right-hand side of the claimed inequality is infinity, while if $m = 0$ then the claim reduces to Propositions~\ref{prop:pB-ptwise-bounds} and~\ref{prop:qB-ptwise-bds}). We focus on the estimate involving $\fq_B$, as the one for $\fp_B$ proceeds by a similar argument. By Lemma~\ref{lemma:covd-dagger-adjoint-relation}, it suffices to estimate $\covd_y^{(m)} \bigl(\covd_x^{(k)} \fq_B\bigr)^*$. Let $r_0 = (t_0 + s_0)/2$ be the midpoint between $s_0$ and $t_0$. 
For $s \in [s_0, r_0]$, define
\begin{equ}
\phi(s, y) = \bigl(\covd^{(k)}_x \fq_B(s, y; t_0, x_0)\bigr)^*.
\end{equ}
By the semigroup property of the covariant heat kernel, we have that
\begin{equ}
\phi(s, y)^* = \int \dz\, \bigl(\covd^{(k)}_x \fq_B(r_0, z; t_0, x_0)\bigr) \fq _B(s, y; r_0, z).
\end{equ}
Then by Lemma~\ref{lemma:adjoints-backwards-heat-equations}, we have 
\begin{equ}
(\ptl_s + \Delta_B + 2F) \phi = 0, \qquad
\phi(r_0) = \bigl(\covd^{(k)}_x \fq_B(r_0, \place; t_0, x_0)\bigr)^*.
\end{equ}
Define the (shifted) time-reversal $\tilde{\phi} \colon [0, r_0 - s_0] \times \T^d \ra \End(\genV^d)$ by $\tilde{\phi}(t) = \phi(r_0 - t)$. Similarly, let $\tilde{B}(t) = B(r_0 - t)$, so that $\tilde{F}(t) = F(r_0 - t)$. Since $B$ is a solution to the Yang--Mills gradient flow, $\tilde{B}$ is a solution to the backward Yang--Mills gradient flow. By the chain rule, we have that
\begin{equ}
(\ptl_t - \Delta_{\tilde{B}} - 2\tilde{F}) \tilde{\phi} = 0, \qquad
\tilde{\phi}(0) = \bigl(\covd_x^{(k)} \fq_B(r_0, \place; t_0, x_0)\bigr)^*.
\end{equ}
In the following, write $\tilde{\covd} = \covd_{\tilde{B}}$ and observe that
\begin{equ}
\tilde{\covd}^{(m)} \tilde{\phi}(r_0 - s_0, y_0) = \covd_y^{(m)} \bigl(\covd_x^{(k)} \fq_B(s_0, y_0; t_0, x_0)\bigr)^*,
\end{equ}
which is precisely the quantity we want to estimate. By Lemma~\ref{lemma:time-reversal-higher-regularity} applied with $T = r_0 - s_0$, $T_1 = r_0$ (note $T_1 > T$ since we assumed $s_0 > 0$ at the beginning), we have that for any positive solution $K \colon [0, r_0 - s_0] \times \T^d \ra (0, \infty)$ to the backward heat equation, and any $m \geq 0$,
\begin{equs}
\big\|t^{m/2} (r_0 - t)^{m/2}& K^{1/2} \tilde{\covd}^{(m)} \tilde{\phi} \big\|_{L_t^\infty L_x^2([0, r_0 - s_0])} \\
&\lesssim \gronwall_{r_0 - s_0}(C_0 \tilde{F}) r_0^{m/2} \Big(1 + (\revFnorm{\tilde{F}}_{[0, r_0 - s_0], r_0, m})^{m}\Big)\big\|K^{1/2}(0) \tilde{\phi}(0)\big\|_{L_x^2}.
\end{equs}
Taking $K(t, x) = \fp(t, x; r_0 - s_0, y_0)$, observe by Proposition~\ref{prop:qB-ptwise-bds} that
\begin{equs}
\big\|K^{1/2}(0) \tilde{\phi}(0)\big\|_{L_x^2} &\leq \big\| K^{1/2}(0) |\covd^{(k)}_x \fq_B(r_0, \place; t_0, x_0)| \big\|_{L_x^2} \\
&\lesssim  \gronwall_{[r_0, t_0]}(C_0 F) \Big(1 + (\Fnorm{F}_{[r_0, t_0], k})^{k}\Big) (t_0 - r_0)^{-\tfrac{k+d}{2}} \exp\biggl(- c \frac{|x_0 - y_0|^2}{t_0 - r_0}\biggr).
\end{equs}
Combining this with Lemma~\ref{lemma:heat-kernel-weighted-to-L-infty},
we thus obtain (note that $t_0 - r_0 \sim t_0 - s_0$ by the choice of $r_0$ as the midpoint between $s_0, t_0$, and also that $\gronwall_{r_0 - s_0}(C_0\tilde{F}) \gronwall_{[r_0, t_0]}(C_0F)$ can be estimated by $\gronwall_{[s_0, t_0]}(C_0 F)$)
\begin{equs}
(r_0 - s_0)^{m/2} &s_0^{m/2} |\tilde{\covd}^{(m)} \tilde{\phi}(r_0 - s_0, y_0)| \lesssim \gronwall_{[s_0, t_0]}(C_0 F) r_0^{m/2}\Big(1 +  (\revFnorm{\tilde{F}}_{[0, r_0 - s_0], r_0, m})^{m}\Big) \\
&\times \Big(1 + (\Fnorm{F}_{[s_0, t_0], k})^{k}\Big) (t_0 - s_0)^{-\tfrac{k+d}{2}} \exp\biggl(- c \frac{|x_0 - y_0|^2}{t_0 - s_0}\biggr).
\end{equs}
Next, by the fact that $\tilde{B}(t) = B(r_0 - t), \tilde{F}(t) = F(r_0 - t)$, and a change of variables, we may obtain that for any $\alpha, \beta \in \R$,
\begin{equ}
\int_0^{r_0 - s_0} \dt \, (r_0 - t)^\alpha \big\| \covd_{\tilde{B}(t)}^{(m)} \tilde{F}(t)\big\|_{L_x^\infty}^\beta = \int_{s_0}^{r_0} \dt \,  t^\alpha \big\|\covd^{(m)} F(t)\big\|_{L_x^\infty}^\beta,
\end{equ}
from which we obtain
\begin{equ}
\revFnorm{\tilde{F}}_{[0, r_0 - s_0], r_0, m} \leq \Fnorm{F}_{r_0, m}.
\end{equ}
We thus further obtain
\begin{equs}
|\tilde{\covd}^{(m)} \tilde{\phi}(r_0 - s_0, y_0)| \lesssim \,& \gronwall_{[s_0, t_0]}(C_0 F) (r_0 - s_0)^{-m/2} s_0^{-m/2} r_0^{m/2} \Big(1 + (\Fnorm{F}_{t_0, m})^{m}\Big) \\
&\times \Big(1 + (\Fnorm{F}_{[s_0, t_0], k})^{k}\Big)  (t_0 - s_0)^{-\frac{k+d}{2}} \exp\biggl(- c \frac{|x_0 - y_0|^2}{t_0 - s_0}\biggr).
\end{equs}
Finally, using that
\begin{equ}
(r_0 - s_0)^{-m/2} s_0^{-m/2} r_0^{m/2} \lesssim \min(s_0, t_0 - s_0)^{-m/2},
\end{equ}
the desired result now follows.
\end{proof}

\subsection{Parallel transport estimates}\label{section:parallel-transport-estimates}

In this subsection, we derive identities and estimates involving parallel transport (introduced in Section~\ref{section:parallel-transport}, see Definition~\ref{def:parallel-transport}) and its covariant derivatives. With some exceptions, we will typically work with a general time-varying connection one-form $A$, in particular we do not assume that $A$ is a solution to the Yang--Mills gradient flow \eqref{e:YMHF}. These results will be needed in several places later \dash see e.g.\ Sections~\ref{section:heat-kernel-parametrix},~\ref{section:perturbing-around-YM-connection}, and~\ref{section:objects-paracontrolled}. Because their proofs are all quite direct, we will break from the convention of previous sections, and instead give the proofs directly after stating the results. We begin by setting some notation.

\begin{notation}
Let $A \colon [0, T] \times \T^d \ra \frkg^d$ be a time-varying connection one-form with associated curvature $F_A$. As usual, by appending a temporal component $A_0 = 0$, we may view $A$ as a connection one-form on $[0, T] \times \T^d$. Given vectors $v, w \in \R^{1+d}$, we write
\begin{equ}
F_A(v, w) = (F_A)_{\alpha \beta} v^\alpha w^\beta.
\end{equ}
Note that this is a function of $(t, x)$, and so we may also explicitly write 
\begin{equ}
F_A(v, w)(t, x) = F_A(t, x)(v, w) = (F_A(t, x))_{\alpha \beta} v^\alpha w^\beta.
\end{equ}
However, as this gets quite long, we will often refrain from writing out $(t, x)$. If $v, w \in \R^d$, we will also write $F_A(v, w)$ and $F_A(t, x)(v, w)$, where we now sum only over spatial indices (or equivalently, we append zero components to $v, w$ to view them as vectors in $\R^{1+d}$ and then apply our previous definition). Given an additional vector $u \in \R^{1+d}$, we similarly write
\begin{equ}
(\covd_{A, w} F_A)(v, u) = (\covd_{A, w} (F_A)_{\alpha \beta}) v^\alpha u^\beta = (\covd_{A, \gamma} (F_A)_{\alpha \beta}) v^\alpha u^\beta w^\gamma.
\end{equ}
\end{notation}

We proceed to the first set of main results of this section, which provide identities involving parallel transport. Recall our notation $\pt_A(\gamma(t)) = \pt_A(0, y_0; \gamma(t))$ in the case where $\gamma$ is a space-time line segment such that $\gamma(0) = (0, y_0)$.

\begin{remark}\label{remark:curvature-as-operator-genV}
In the following results (Lemmas~\ref{lemma:covd-parallel-transport} and~\ref{lemma:pt-curvature-commutators}), one should view $F_A(v, w)$ and $(\covd_{A, w}F_A)(v, u)$ as linear operators on $\genV$ (as opposed to Lie algebra elements). In other words, we always implicitly apply the representation $\repr$ to obtain elements of $\End(\genV)$.
\end{remark}

\begin{lemma}[Covariant derivatives of parallel transport]\label{lemma:covd-parallel-transport}
Let $A \colon [0, T] \times \T^d \ra \frkg^d$ be a time-varying connection one-form, and let $y_0 \in \T^d$. Let $z = (z_0, z_x) \in \R^{1+d}$ such that $z_0 \geq 0$, and let $\gamma(t) = (0, y_0) + tz$. For any $t_0 > 0$ such that $t_0 z_0 \leq T$ and $t_0 |z_x| \leq 1/10$ and $v, w\in \R^{1+d}$, we have the formulas
\begin{equs}
(\pt_A^{-1} \covd_{A, v} \pt_A)(\gamma(t_0)) &= \frac{1}{t_0}\int_0^{t_0} \dt \, t  \, (\pt_A^{-1} F_A(z, v) \pt_A)(\gamma(t)) , \\
(\pt_A^{-1} \covd_{A, w} \covd_{A, v} \pt_A) (\gamma(t_0)) &= \frac{1}{t_0^2} \int_0^{t_0} \dt \, t^2 \Big( (\pt_A^{-1} F_A(z, w) \covd_{A, v} \pt_A)(\gamma(t)) + (\pt_A^{-1} F_A(z, v) \covd_{A, w} \pt_A)(\gamma(t))\Big) \\
&+ \frac{1}{t_0^2} \int_0^{t_0} \dt \, t^2 (\pt_A^{-1} (\covd_{A, w} F_A)(z, v) \pt_A)(\gamma(t)) + \frac{1}{t_0^2} \int_0^{t_0} \dt \, t (\pt_A^{-1} F_A(w, v) \pt_A)(\gamma(t)).
\end{equs}
\end{lemma}
\begin{proof}
For notational simplicity, write $\covd = \covd_{A}$ and $F = F_A$. Let $\gamma(t, s) = (0, y_0) + t(z + s v)$, so that in particular $\gamma(t, 0) = \gamma(t)$. We may use $\gamma$ to define a pullback vector bundle over $[0, t_0] \times (-\varep, \varep)$, where $\varep$ is sufficiently small so that $\gamma(t, s)$ is within the domain of $\pt_A$ for all $(t, s) \in [0, t_0] \times (-\varep, \varep)$. Define the vector fields $\dot{\gamma} = \ptl_t \gamma$, and $J = \ptl_s \gamma$. We will denote the induced covariant derivatives on the pullback bundle by $\covd_t = \covd_{\dot{\gamma}}$ and $\covd_s = \covd_J$. By definition of $\pt_A$, we have that $\covd_t \pt_A = 0$. The induced curvature is given by $F(\ptl_t, \ptl_s) = F(\dot{\gamma}, J)$. With all this setup, we may now compute
\begin{equ}
\covd_t \covd_s \pt_A = \covd_s \covd_t \pt_A + F(\ptl_t, \ptl_s) \pt_A = F(\dot{\gamma}, J) \pt_A. 
\end{equ}
Applying the above identity at $s = 0$ (in which case $\dot{\gamma} = z$ and $J = t v$), and then applying Lemma~\ref{lemma:transport-lemma}, we obtain (note that the initial condition is zero since $tv = 0$ at time $t = 0$)
\begin{equ}
(\pt_A^{-1} \covd_{t_0 v} \pt_A)(\gamma(t_0)) = \int_0^{t_0} \dt \, t (\pt_A^{-1} F(z, v) \pt_A)(\gamma(t)).
\end{equ}
Dividing both sides by $t_0$, the first desired identity now follows. 

For the second desired identity, we proceed similarly, except now we define $\gamma(t, s_1, s_2) = (0, y_0) + t(z + s_1 v + s_2 w)$, $\dot{\gamma}(t) = \ptl_t \gamma$, $J_1 = \ptl_{s_1} \gamma$, $J_2 = \ptl_{s_2} \gamma$. We may now use $\gamma$ to define a pullback vector bundle over $[0, t_0] \times (-\varep, \varep)^2$, on which we have a covariant derivative given by $\covd_t = \covd_{\dot{\gamma}}$ and $\covd_{s_i} = \covd_{J_i}$ for $i = 1, 2$. Arguing as before, we may first obtain 
\begin{equ}
\covd_t \covd_{s_1} \pt_A = F(\ptl_t, \ptl_{s_1}) \pt_A. 
\end{equ}
We may then obtain (applying the above identity and the product formula in the second identity)
\begin{equs}
~&\covd_t \covd_{s_2} \covd_{s_1} \pt_A \\
&= F(\ptl_t, \ptl_{s_2}) \covd_{s_1} \pt_A + \covd_{s_2}(\covd_t \covd_{s_1} \pt_A) \\
&= F(\ptl_t, \ptl_{s_2}) \covd_{s_1} \pt_A + F(\ptl_t, \ptl_{s_1}) \covd_{s_2} \pt_A + (\covd_{s_2} F)(\ptl_t, \ptl_{s_1}) \pt_A + \bigl( F(\ptl_{s_2} \ptl_t, \ptl_{s_1}) + F(\ptl_t, \ptl_{s_2} \ptl_{s_1})\bigr) \pt_A \\
&= F(\dot{\gamma}, J_2) \covd_{J_1} \pt_A + F(\dot{\gamma}, J_1) \covd_{J_2} \pt_A + (\covd_{J_2} F)(\dot{\gamma}, J_1) \pt_A + \bigl( F(\ptl_{s_2} \dot{\gamma}, J_1) + F(\dot{\gamma}, \ptl_{s_2} J_1)\bigr) \pt_A. 
\end{equs}
Using that $\ptl_{s_2} \dot{\gamma} = w$ and $\ptl_{s_2} J_1 = 0$, then applying the above identity at $s_1 = s_2 = 0$, in which case $\dot{\gamma} = z$, $J_1  = tv$, $J_2 = tw$, and then applying Lemma~\ref{lemma:transport-lemma}, we may obtain the second desired identity (after multiplying both sides by $\pt_A(\gamma(t_0))^{-1}$ and dividing both sides by $t_0^{2}$).
\end{proof}

\begin{lemma}[Commutator identities for parallel transport and curvature]\label{lemma:pt-curvature-commutators}
Let $A \colon [0, T] \times \T^d \ra \frkg^d$ be a time-varying connection one-form, and let $y_0 \in \T^d$. Let $u = (u_0, u_x) \in \R^{1+d}$ be such that $u_0 \geq 0$, and let $\gamma(t) = (0, y_0) + t u$. For any $t \geq 0$ such that $t u_0 \leq T$ and $t |u_x| \leq 1/10$ and $w, v \in \R^{1+d}$, we have that
\begin{equ}
\big(\pt_A^{-1} F_A(w, v) \pt_A\big)(\gamma(t)) - F_A(\gamma(0))(w, v) = \int_0^t \ds \, \big(\pt_A^{-1} (\covd_{A, \dot{\gamma}(s)} F_A)(w, v)\pt_A \big)(\gamma(s)).
\end{equ}
Consequently, for any $t_0 > 0$ such that $t_0 u_0 \leq T$ and $t_0 |u_x| \leq 1/10$, we have that
\begin{equs}
(\pt_A^{-1} \covd_{A, v} \pt_A)(\gamma(t_0)) ~&- \frac{t_0}{2} F_A(\gamma(0))(\dot{\gamma}(0), v) \\
&= \frac{1}{t_0} \int_0^{t_0} \dt \,  t \int_0^t \ds \, \bigl(\pt_A^{-1} (\covd_{A, \dot{\gamma}(s)} F_A)(\dot{\gamma}(s), v) \pt_A\bigr)(\gamma(s)).
\end{equs}
Additionally, if $d = 2$ and $A$ is a Yang--Mills connection (so in particular it is also constant in time), then we have that for all $t \geq 0$ such that $t |u_x| \leq 1/10$,
\begin{equ}
\big(\pt_A^{-1} F_A(w, v) \pt_A\big) (\gamma(t)) = F_A(\gamma(0))(w, v). 
\end{equ}
\end{lemma}
\begin{proof}
The second identity follows by combining the first identity with Lemma~\ref{lemma:covd-parallel-transport} (note that $\dot{\gamma}(s) = u$ is constant in $s$, so we may apply the first identity with $w = u$). The third identity follows by combining the first identity with the fact that in dimension $d = 2$, $A$ being a Yang--Mills connection implies that all covariant derivatives of $F_A$ are zero (recall Remark~\ref{remark:YM-connection-2D}, which gives that all covariant spatial derivatives are zero, and recall that $A$ is constant in time, and thus the covariant temporal derivative of $F_A$ is also zero). Thus, we focus on proving the first identity. Towards this end, by the product formula and the defining equation \eqref{eq:pt-equation} of parallel transport, we have that
\begin{equ}
\frac{d}{dt} \big(\pt_A^{-1} F_A(w, v) \pt_A\big) (\gamma(t))   = \big( \pt_A^{-1} (\covd_{A, \dot{\gamma}(t)} F_A)(w, v) \pt_A \big)(\gamma(t)),
\end{equ}
whence the first identity follows by integrating both sides.
\end{proof}

Next, we turn to the second set of main results of this section, which provide estimates involving parallel transport.

\begin{lemma}[Estimates on covariant derivatives of parallel transport]\label{lemma:covd-parallel-transport-estimate}
Let $A \colon [0, T] \times \T^d \ra \frkg^d$ be a time-varying connection one-form. For any $(t_0, x_0) \in (0, T] \times \T^d$ such that $|x_0| \leq 10^{-1}$, we have that
\begin{equs}
|\covd_A \pt_A(t_0, x_0)| &\lesssim \bigl( \bigl\|\ptl_t A \bigr\|_{L_t^1 L_x^\infty([0, t_0])} + \bigl\|F_A\bigr\|_{L_t^2 L_x^\infty([0, t_0])}\bigr) \bigl( 1 +  |x_0| t_0^{-1/2}\bigr), \label{eq:covd-parallel-transport-estimate} \\
|\ptl_t \pt_A(t_0, x_0)| &\lesssim  \bigl\|\ptl_t A \bigr\|_{L_t^1 L_x^\infty([0, t_0])} |x_0| t_0^{-1}. \label{eq:time-derivative-parallel-transport-estimate} 
\end{equs}
If, moreover, $A$ is a solution to the Yang--Mills gradient flow \eqref{e:YMHF}, then we have that
\begin{equs}
|\Delta_A \pt_A(t_0, x_0)| &\lesssim \bigl\|\ptl_t A \bigr\|_{L_t^1 L_x^\infty([0, t_0])} |x_0| t_0^{-1}\\
&+\bigl( \bigl\|\ptl_t A \bigr\|_{L_t^1 L_x^\infty([0, t_0])} + \bigl\|F_A\bigr\|_{L_t^2 L_x^\infty([0, t_0])}\bigr)^2 \bigl(1 + |x_0| t_0^{-1/2}\bigr)^2. \label{eq:cov-laplacian-parallel-transport-estimate}
\end{equs}
\end{lemma}
\begin{proof}
In the following, all space-time norms will implicitly be on $[0, t_0] \times \T^d$. Also, for notational simplicity, we write $\covd = \covd_A$ and $F = F_A$. Let $z = (1, x_0 / t_0)$, $\gamma(t) =  t z$. Note that $\gamma(0) = (0, 0)$ and $\gamma(t_0) = (t_0, x_0)$.
Since $\pt_A$ is an isometry (recall Remark~\ref{remark:parallel-transport}\ref{remark:pt-isometry}), it suffices to estimate $\pt_A^{-1} \covd \pt_A$. By Lemma~\ref{lemma:covd-parallel-transport}, we have that for $\ell \in [d]$, $t \in [0, t_0]$,
\begin{equ}\label{eq:covd-parallel-transport-intermediate}
|(\pt_A^{-1} \covd_\ell \pt_A)(\gamma(t))| \lesssim \frac{1}{t} \int_0^{t} \ds \, s |F(z, e_\ell)| \lesssim  \bigl\|\ptl_t A \bigr\|_{L_t^1 L_x^\infty} + \bigl\|F\bigr\|_{L_t^2 L_x^\infty} |x_0| t_0^{-1/2},
\end{equ}
where we used the identity $F(z, e^\ell) = \ptl_t A^\ell + F^{i\ell} (x_0)_i / t_0$, as well as Hölder.
Applying this with $t = t_0$, we obtain the first desired estimate \eqref{eq:covd-parallel-transport-estimate}. The second desired estimate \eqref{eq:time-derivative-parallel-transport-estimate} follows by similar considerations (we replace $e_\ell$ with $e_0$, which is the standard basis vector of $\R^{1+d}$ in the $0$th (i.e.\ temporal) direction). 
Towards \eqref{eq:cov-laplacian-parallel-transport-estimate}, first note that by Lemma~\ref{lemma:covd-parallel-transport}, we have that
\begin{equs}
(\pt_A^{-1}\Delta_A &\pt_A)(\gamma(t_0)) = (\pt_A^{-1}\covd^\ell \covd_\ell \pt_A) (\gamma(t_0)) \\
&= \frac{1}{t_0^2}\int_0^{t_0} \dt \,  t^2 \Big(2  \bigl(\pt_A^{-1} F(z, e^\ell) \covd_\ell \pt_A\bigr)(\gamma(t)) + \bigl(\pt_A^{-1} (\covd^\ell F)(z, e_\ell) \pt_A\bigr) (\gamma(t)) \Big).
\end{equs}
Since $A$ is assumed to be a solution to the Yang--Mills gradient flow \eqref{e:YMHF}, we have that (in the last identity we use that $\ptl_t A$ is covariant divergence free by Lemma~\ref{lemma:YM-gradient-flow-covariant-divergence-free})
\begin{equs}
(\covd^\ell F)(z, e_\ell) = \covd^\ell \ptl_t A_\ell + \covd^\ell F_{i\ell} (x_0^i / t_0) = - \ptl_t A_i (x_0^i / t_0),
\end{equs}
The estimate \eqref{eq:cov-laplacian-parallel-transport-estimate} now follows by inserting this identity and using the estimate \eqref{eq:covd-parallel-transport-intermediate}.
\end{proof}


The following result is a bit unlike the previous results, in that it works with a general base point, and it also assumes control on a slightly different quantity than we have seen before (see \eqref{eq:pt-approximate-semigroup-norm-assumption}). This will be needed in Section~\ref{section:objects-paracontrolled} to compute the resonance of certain high$\times$high-interactions arising in our covariant paracontrolled calculus. In the following, we recall our standing convention in Section~\ref{section:covariant-heat-kernel-estimates} that $B$ denotes a solution to the Yang--Mills gradient flow \eqref{e:YMHF}, and $F = F_B$.

\begin{lemma}[Approximate semigroup property of parallel transport]\label{lemma:parallel-transport-approximate-semigroup} Let $\alpha \in [0, 1]$, and suppose that $L \geq 1$ is such that
\begin{equ}\label{eq:pt-approximate-semigroup-norm-assumption}
T^{(1-\alpha)/2}\Fnorm{t^{-(1-\alpha)/2} F}_{T, 1} \leq L. 
\end{equ}
Let $0 \leq r_0 < s_0 < t_0 \leq T$ and $x_0, y_0, z_0 \in \T^d$ be such that $|x_0 - y_0|, |x_0 - z_0|, |y_0 - z_0| \leq 10^{-1}$. Let $\rho \eqdef \max(|x_0 - z_0|, |y_0 - z_0|)$. We have the following approximate semigroup property:
\begin{equs}
|&\pt_B(r_0, z_0; t_0, x_0) - \pt_B(s_0, y_0; t_0, x_0) \pt_B(r_0, z_0; s_0, y_0)| \\
&\lesssim L T^{-(1-\alpha)/2} r_0^{-\alpha/2} \Bigl( \rho \log \frac{t_0 - r_0}{s_0 - r_0} + |x_0 - y_0|\bigl(1 + \rho (t_0 - s_0)^{-1/2}\bigr) \Bigr).
\end{equs}
\end{lemma}
\begin{proof}
It suffices to show that for any $X \in \genV$ with $|X|_{\genV} = 1$, we have that
\begin{equ}
\Bigl|\Bigl(\pt_B(r_0, z_0; t_0, x_0) - \pt_B(s_0, y_0; t_0, x_0) \pt_B(r_0, z_0; s_0, y_0)  \Bigr)X \Bigr| \lesssim R,
\end{equ}
where $R$ is the claimed upper bound in the lemma statement. Towards this end, fix such an $X$, and define the local section $\Psi(t, x) \eqdef \pt_B(r_0, z_0; t, x) X$ for $|x-z_0| \leq 1/10$ (say). Define also the path $\gamma \colon [0, t_0 - s_0] \ra \R_+ \times \T^d$ from $(s_0, y_0)$ to $(t_0, x_0)$ by 
\begin{equ}
\gamma(t) \eqdef (\gamma_t(t), \gamma_x(t)) \eqdef (s_0, y_0) + t \Big(1, \frac{x_0 - y_0}{t_0 - s_0}\Big).
\end{equ}
By Lemma~\ref{lemma:transport-lemma}, we have that
\begin{equs}
|\Psi(t_0, x_0) - \pt_B(s_0, &y_0; t_0, x_0) \Psi(s_0, y_0)| \\
&\leq \int_0^{t_0 - s_0} \dt \,  |\covd_{\dot{\gamma}(t)} \Psi(\gamma(t))| \\
&\lesssim \int_0^{t_0 - s_0} \dt \, |\ptl_t \pt_B(r_0, z_0; \gamma(t))| + \frac{|x_0 - y_0|}{t_0 - s_0} |\covd \pt_B(r_0, z_0; \gamma(t))|.
\end{equs}
Applying Lemma~\ref{lemma:covd-parallel-transport-estimate} (in particular \eqref{eq:covd-parallel-transport-estimate} and \eqref{eq:time-derivative-parallel-transport-estimate}), we obtain the further upper bound:
\begin{equs}
~&\int_0^{t_0 - s_0} \dt \,\|\ptl_t B\|_{L_t^1 L_x^\infty([r_0, t + s_0])} \frac{|\gamma_x(t) - z_0|}{t + s_0 - r_0} \\
&+ \frac{|x_0 - y_0|}{t_0 - s_0} \int_0^{t_0 - s_0} \dt\, \Big(\|\ptl_t B\|_{L_t^1 L_x^\infty([r_0, t + s_0])} + \|F\|_{L_t^2 L_x^\infty([r_0, t + s_0])} \Bigr) \Bigl(1 +  \frac{|\gamma_x(t) - z_0|}{(t + s_0 - r_0)^{1/2}}\Big).
\end{equs}
To finish, it suffices to show the separate estimates
\begin{equs}
~&\int_0^{t_0 - s_0} \dt \,  \|\ptl_t B\|_{L_t^1 L_x^\infty([r_0, t + s_0])} \frac{|\gamma_x(t) - z_0|}{t + s_0 - r_0} \\
&\quad \quad \quad \quad \quad \quad \lesssim L T^{-\frac{1-\alpha}{2}} r_0^{-\frac{\alpha}{2}} \rho \log \frac{t_0 - r_0}{s_0 - r_0}, \label{eq:parallel-transport-approximate-semigroup-intermediate-1} \\
~&\frac{|x_0 - y_0|}{t_0 - s_0} \int_0^{t_0 - s_0} \dt\, \|\ptl_t B\|_{L_t^1 L_x^\infty([r_0, t + s_0])} \Bigl(1 +  \frac{|\gamma_x(t) - z_0|}{(t + s_0 - r_0)^{1/2}}\Big) \\
&\quad \quad \quad \quad \quad \quad \lesssim L  T^{-\frac{1-\alpha}{2}}  r_0^{-\frac{\alpha}{2}} |x_0 - y_0| (1 + \rho (t_0 - s_0)^{-1/2}), \label{eq:parallel-transport-approximate-semigroup-intermediate-2} \\
~&\frac{|x_0 - y_0|}{t_0 - s_0} \int_0^{t_0 - s_0} \dt \|F\|_{L_t^2 L_x^\infty([r_0, t + s_0])} \Bigl(1 + \frac{|\gamma_x(t) - z_0|}{(t + s_0 - r_0)^{1/2}}\Bigr) \\
&\quad \quad \quad \quad \quad \quad \lesssim L T^{-\frac{1-\alpha}{2}} r_0^{-\frac{\alpha}{2}} |x_0 -y_0| (1 + \rho (t_0 - s_0)^{-1/2}).\qquad \label{eq:parallel-transport-approximate-semigroup-intermediate-3} 
\end{equs}
First, note that for any $t \in [0, t_0 - s_0]$, we have that
\begin{equ}
|\gamma_x(t) - z_0| \leq \max(|x_0 - z_0|, |y_0 - z_0|) = \rho,
\end{equ}
which follows because $(\gamma_x(t), t \in [0, t_0 - s_0])$ traces out the edge $(y_0, x_0)$ of the triangle whose vertices are $x_0, y_0, z_0$. For \eqref{eq:parallel-transport-approximate-semigroup-intermediate-1}, we apply \eqref{eq:ptl-t-B-estimate-under-more-stringent-Fpert} to estimate the left-hand side of \eqref{eq:parallel-transport-approximate-semigroup-intermediate-1} by
\begin{equ}
L T^{-(1-\alpha)/2} r_0^{-\alpha/2} \rho \int_0^{t_0 - s_0} \dt \, (t + s_0 - r_0)^{-1} \lesssim L T^{-(1-\alpha)/2} r_0^{-\alpha/2}  \rho \log \frac{t_0 - r_0}{s_0 - r_0},
\end{equ}
which proves \eqref{eq:parallel-transport-approximate-semigroup-intermediate-1}. Next, for \eqref{eq:parallel-transport-approximate-semigroup-intermediate-2}, we again apply \eqref{eq:ptl-t-B-estimate-under-more-stringent-Fpert} to estimate the left-hand side of \eqref{eq:parallel-transport-approximate-semigroup-intermediate-2} by
\begin{equ}
L T^{-(1-\alpha)/2} r_0^{-\alpha/2} |x_0 - y_0| \bigl(1 + \rho (t_0 - s_0)^{-1/2}\bigr),
\end{equ}
as desired. Finally, \eqref{eq:parallel-transport-approximate-semigroup-intermediate-3} follows by an argument similar to that for \eqref{eq:parallel-transport-approximate-semigroup-intermediate-2}, where we replace the application of \eqref{eq:ptl-t-B-estimate-under-more-stringent-Fpert} by \eqref{eq:F-estimate-under-more-stringent-Fpert}.
\end{proof}

\subsection{Heat kernel parametrix}\label{section:heat-kernel-parametrix}

In this subsection, we obtain more refined short-time expansions for $\fp_B$ and $\fq_B$ (i.e.\ parametrices), which are expressed in terms of the parallel transport $\pt_B$ (which recall was defined in Definition~\ref{def:parallel-transport}). While some versions (and extensions) of these results may be found in the textbook \cite{BGV2004}, we prefer to give self-contained proofs here. We state the main results of this subsection (Proposition~\ref{prop:pB-parametrix} and Corollary~\ref{cor:pB-qB-parametrix}), and then proceed with the proofs. In the following, recall that $\fp$ denotes the Euclidean heat kernel.

\begin{definition}\label{def:ptcutoff}
Let $\ptcutoff \colon \R^d \ra [0, 1]$ be a smooth radial function such that $\ptcutoff(x) = 1$ for all $|x| \leq 1/20$, and $\ptcutoff(x) = 0$ for all $|x| > 1/10$. We will also view $\ptcutoff$ as a function on $\T^d$. For notational convenience, given points $(s, y), (t, x)$, we will write $\ptcutoff(s, y; t, x) = \ptcutoff(x-y)$ as well as $\ptcutoff(t, x) = \ptcutoff(0, 0; t, x) = \ptcutoff(x)$.
\end{definition}

A version of this proposition was already stated as Proposition~\ref{main:prop-parametrix}, but we choose to reproduce it here for the reader's convenience (there is also the slight difference that the following proposition includes the spatial cutoff $\ptcutoff$).

\begin{proposition}\label{prop:pB-parametrix}
There is a constant $c > 0$ such that the following holds. Let $t_0 \in (0, 1]$. Let
\begin{equs}
\mc{N}(B) \eqdef \Big(t_0^{1/2} \big(\|\ptl_t B\|_{L_t^1 L_x^\infty([0, t_0])} + \|F\|_{L_t^2 L_x^\infty([0, t_0])} \big) \Big)^2  & + t_0^{1/2}\|\ptl_t B\|_{L_t^1 L_x^\infty([0, t_0])}  \\
&+ t_0^{10} \bigl(1 + \big\|F\big\|_{L_t^2 L_x^\infty([0, t_0])}\bigr) .
\end{equs}
For any point $x_0 \in \T^d$, we have the pointwise bounds
\begin{equs}
|\fp_B - \fp \ptcutoff \pt_B |(t_0, x_0) &\lesssim \mc{N}(B) t_0^{-\tfrac{d}{2}} \exp(-c|x_0|^2/t_0), \label{eq:pB-parametrix-zeroth-order-estimate} \\
\bigl|\covd(\fp_B - \fp \ptcutoff \pt_B)\bigr|(t_0, x_0) &\lesssim \mc{N}(B) (1 + \Fnorm{F}_{t_0, 0})  t_0^{-\tfrac{d+1}{2}} \exp(- c|x_0|^2/t_0). \label{eq:pB-parametrix-first-order-estimate}
\end{equs}
\end{proposition}
\begin{remark}\label{remark:pB-parametrix-higher-order-derivatives}
One might naturally wonder whether higher order pointwise estimates also hold for the difference $\fp_B - \fp \ptcutoff \pt_B $, similar to what we proved for $\fp_B, \fq_B$, and $\fp_B - \fq_B$. However, we expect that for the second order estimate and higher, one can do no better than estimating each of the terms $\fp_B, \fp \ptcutoff \pt_B$ individually. This is because $\fp \ptcutoff \pt_B$ is only a leading-order approximation for $\fp_B$. In particular, if we begin taking higher derivatives of $\fp_B - \fp \ptcutoff \pt_B$, then these derivatives will fall on the forcing term $\fp \ptcutoff \Delta_B \pt_B$ appearing in the equation \eqref{eq:pB-parametrix-equation} for $\fp_B - \fp \ptcutoff \pt_B$. If we take two or more derivatives, then this will result in a forcing term with a non-integrable singularity, and thus we cannot just integrate this term up to obtain a bound for the difference (which is what we do for the zeroth and first order estimates).
\end{remark}

We state the following corollary of Proposition~\ref{prop:pB-parametrix} which will in fact be the result which is applied later (for instance, in Section~\ref{section:non-resonant-estimates}). Unlike the previous result, this corollary is stated for a general pair of space-time points $(s, y), (t, x)$, and it allows us to gain a sufficiently high power of $t-s$, which is needed for later use. In order to gain this power, we need to assume control on a slightly stronger norm than we have been using thus far \dash see \eqref{eq:Fnorm-more-stringent-assumption} (this assumption is similar to the one appearing in \eqref{eq:pt-approximate-semigroup-norm-assumption}). This norm is stronger because we are placing less of a time weight on $F$ (and $\covd F$).

\begin{corollary}\label{cor:pB-qB-parametrix}
Let $T \leq 1$, $(s, y), (t, x) \in [0, T] \times \T^d$ be such that $s < t$. Let $\alpha \in [0, 1]$, and suppose that for some $L \geq 1$, we have that
\begin{equ}\label{eq:Fnorm-more-stringent-assumption}
\bigl(T^{(1-\alpha)/2}\Fnorm{u^{-(1-\alpha)/2} F}_{T, 1}\bigr)^2 \leq L. 
\end{equ}
For notational convenience, let
\begin{equ}
\mc{N}(B, L, s, t) = L \bigl(T^{-1}(t-s)\bigr)^{1-\alpha} s^{-\alpha} (t-s)^\alpha + \bigl(\big\|\ptl_u B\big\|_{L_u^1 L_x^\infty([s, t])} + \bigl\|F\bigr\|_{L_u^2 L_x^\infty([s, t])} \bigr) (t-s)^{1/2}.
\end{equ}
Then, we have that
\begin{equs}
~\hspace{-5mm}\big|(\fp_B - \fp \ptcutoff \pt_B)(s, y; t, x)\big| &\lesssim \mc{N}(B, L, s, t) (t-s)^{-\tfrac{d}{2}} \exp\Bigl(-c \frac{|x-y|^2}{t-s}\Bigr), \label{eq:pB-parametrix-t-s-power} \\
\hspace{-8mm}\big|(\covd \fp_B - (\nabla \fp) \ptcutoff \pt_B)(s, y; t, x)\big| & \lesssim (1 + \Fnorm{F}_{[s, t], 0})\mc{N}(B, L, s, t) (t-s)^{-\tfrac{d+1}{2}} \exp\Bigl(-c \frac{|x-y|^2}{t-s}\Bigr). \label{eq:covd-pB-parametrix-t-s-power}
\end{equs}
Furthermore, with $C_0$ taken to be the constant from Proposition~\ref{prop:qB-ptwise-bds} we have that, for $0 \leq k \leq 1$,
\begin{equation}
\begin{aligned}
&\quad \bigl| \bigl(\covd^{(k)} \fq_B - \bigl(  (\nabla^{(k)} \fp) \ptcutoff  \pt_B\bigr)\bigr)(s,y;t,x) \bigr| \\  
&\lesssim \gronwall_{[s, t]}(C_0 F) \bigl(\mc{N}(B, L, s, t)(1 + \Fnorm{F}_{[s, t], 0}) + \Fnorm{F}_{[s, t], 1} \bigr) (t-s)^{-\tfrac{d+k}{2}} \exp\Bigl( - c \frac{|x-y|^2}{t-s} \Bigr). \label{eq:qB-parametrix-t-s-power}
\end{aligned}
\end{equation}
\end{corollary}

Next, we proceed to prove the main results. We start with Proposition~\ref{prop:pB-parametrix}. Since the proposition shows that $\fp_B$ is well approximated by $\fp \ptcutoff \pt_B$, it is natural to begin by deriving the equation satisfied by $\fp \ptcutoff \pt_B$. Towards this end, we compute
\begin{equs}
(\ptl_t &- \Delta_B) (\fp \ptcutoff \pt_B) \\
&= \bigl((\ptl_t - \Delta) (\fp \ptcutoff) \bigr) \pt_B + \fp \ptcutoff (\ptl_t - \Delta_B) \pt_B - 2 \ptl^\ell (\fp \ptcutoff)\, \covd_\ell \pt_B \\
&= -2 \ptl^\ell \fp (\ptl_\ell \ptcutoff) \pt_B - \fp (\Delta \ptcutoff) \pt_B + \fp \ptcutoff \Big( \ptl_t \pt_B - 2 (\ptl^\ell \fp / \fp) \covd_\ell \pt_B - \Delta_B \pt_B\Big)  - 2 \fp \ptl^\ell \ptcutoff \covd_\ell \pt_B .
\end{equs}
Now observe that 
\begin{equs}
\ptl_t \pt_B - 2 (\ptl^\ell \fp / \fp) \covd_\ell \pt_B &= \ptl_t \pt_B + (x^\ell /t) \covd_\ell \pt_B -\bigl(2\ptl^\ell \fp / \fp + x^\ell / t\bigr) \covd_\ell \pt_B \\
&= \covd_{\dot{\gamma}(t)} \pt_B(\gamma(t)) - \bigl(2\ptl^\ell \fp / \fp + x^\ell / t\bigr) \covd_\ell \pt_B,
\end{equs}
where for fixed $x, t$, the path $\gamma$ is defined as $\gamma(s) = s (1, x/t)$.
By the definition of parallel transport (Definition~\ref{def:parallel-transport}), the first term on the right-hand side above is zero, at least on the region where $\ptcutoff \neq 0$ (because this ensures that $x$ is close to zero). For ease of notation, let 
\begin{equ}
R \eqdef 2 \nabla \fp / \fp  + x / t.
\end{equ}
We thus obtain
\begin{equs}
(\ptl_t - \Delta_B) (\fp \ptcutoff \pt_B) = -2 \ptl^\ell \fp (\ptl_\ell \ptcutoff) \pt_B - \fp \ptcutoff \Delta_B \pt_B  - \fp \bigl(2 \ptl^\ell \ptcutoff +  \ptcutoff R^\ell \bigr)  \covd_\ell \pt_B - \fp (\Delta \ptcutoff) \pt_B,
\end{equs}
and so the difference $\fp_B - \fp \ptcutoff \pt_B$ solves
\begin{equs}\label{eq:pB-parametrix-equation}
(\ptl_t - \Delta_B)(\fp_B - \fp \ptcutoff \pt_B) &=  2 \ptl^\ell \fp (\ptl_\ell \ptcutoff) \pt_B  + \fp \ptcutoff \Delta_B \pt_B + \fp \bigl(2 \ptl^\ell \ptcutoff +  \ptcutoff R^\ell \bigr) \covd_\ell \pt_B + \fp (\Delta \ptcutoff) \pt_B, \\
(\fp_B - \fp \ptcutoff \pt_B)(0) &= 0.
\end{equs}
Thus in order to estimate the difference $\fp_B - \fp \ptcutoff \pt_B$, we need estimates on $\covd \pt_B, \Delta_B \pt_B$. Lemma~\ref{lemma:covd-parallel-transport-estimate} exactly provides the needed estimates. We now proceed to prove Proposition~\ref{prop:pB-parametrix}.



\begin{proof}[Proof of Proposition~\ref{prop:pB-parametrix}]
By \eqref{eq:pB-parametrix-equation}, we have that
\begin{equ}
(\fp_B - \fp \ptcutoff \pt_B) (t_0, x_0) = I_m + I_r,
\end{equ}
where (here $`m\text{'}$ stands for main, and $`r\text{'}$ stands for remainder)
\begin{equs}
I_m &\eqdef \int \dt \dx \, \fp_B(t, x; t_0, x_0)  (\fp \ptcutoff \Delta_B \pt_B)(t, x), \\
I_r &\eqdef \int \dt \dx \, \fp_B(t, x; t_0, x_0) \bigl(  2 \ptl^\ell \fp (\ptl_\ell \ptcutoff) \pt_B +  \fp \bigl( 2\ptl^\ell \ptcutoff + \ptcutoff R^\ell\bigr)  \covd_\ell \pt_B + \fp (\Delta \ptcutoff) \pt_B\bigr)(t, x).
\end{equs}
Applying the diamagnetic inequality \eqref{eq:zeroth-order-pB-ptwise-bd} and inserting our pointwise estimate of $\Delta_B \pt_B$ \eqref{eq:cov-laplacian-parallel-transport-estimate}, we obtain
\begin{equ}
\abs{I_m} \lesssim \int \dt \dx \, \fp(t, x; t_0, x_0) \fp(t, x) J_m(t, x),
\end{equ}
where 
\begin{equs}
J_m(t, x) \eqdef \bigl( \|\ptl_t B\|_{L_t^1 L_x^\infty([0, t])} + \|F\|_{L_t^2 L_x^\infty([0, t])}\bigl)^2\bigl(1 +  |x| t^{-1/2}\bigr)^2 + \|\ptl_t B\|_{L_t^1 L_x^\infty([0, t])} |x| t^{-1}.
\end{equs}
Combining this with the estimates
\begin{equs}
\int \dt \dx \, \fp(t, x; t_0, x_0) \fp(t, x) |x| t^{-1} &\lesssim t_0^{1/2} t_0^{-\tfrac{d}{2}} \exp(-c|x_0|^2 / t_0), \\
\int \dt \dx \, \fp(t, x; t_0, x_0) \fp(t, x)  \bigl(1 + |x| t^{-1/2} \bigr)^2 &\lesssim t_0 t_0^{-\tfrac{d}{2}} \exp(-c |x_0|^2 / t_0), 
\end{equs}
we may obtain an acceptable bound on $I_m$. The remainder $I_r$ satisfies better estimates, since fewer derivatives fall on the parallel transport $\pt_B$. In particular, by the diamagnetic inequality \eqref{eq:zeroth-order-pB-ptwise-bd} combined with \eqref{eq:covd-parallel-transport-estimate} and Lemma~\ref{lemma:euclidean-heat-kernel-gradient-torus}, we may estimate 
\begin{equ}
\abs{I_r} \lesssim \int \dt \dx\, \fp(t, x; t_0, x_0) \bigl(|\ptl \fp | |\ptl \ptcutoff| + \fp (|\ptl \ptcutoff| + t^{10}) J_r + \fp |\Delta \ptcutoff| \bigr)(t, x),
\end{equ}
where
\begin{equ}
J_r(t, x) \eqdef \big(\big\|\ptl_t B\big\|_{L_t^1 L_x^\infty([0, t])} + \big\|F \big\|_{L_t^2 L_x^\infty([0, t])}\big)\bigl(1 + |x| t^{-1/2}\bigr).
\end{equ}
Note that $\ptl \ptcutoff, \Delta \ptcutoff$ are supported on the region $|x| \sim 1$, and thus 
\begin{equ}
(|\ptl \fp| |\ptl \ptcutoff|) (t, x), ~~ (\fp |\ptl \ptcutoff|)(t, x) ,~~ (\fp |\Delta \ptcutoff|)(t, x) \lesssim t^{10} \exp(-c |x|^2 / t).
\end{equ} 
Combining, we may obtain the further estimate
\begin{equ}
\abs{I_r} \lesssim t_0^{10} \bigl(1 + \big\|F\big\|_{L_t^2 L_x^\infty([0, t_0])}   + \big\|\ptl_t B\big\|_{L_t^1 L_x^\infty([0, t_0])} \bigr) t_0^{-\tfrac{d}{2}} \exp(-c |x_0|^2 / t_0),
\end{equ}
which is acceptable. Combining our estimates of $I_m$ and $I_r$, we obtain the zeroth order estimate \eqref{eq:pB-parametrix-zeroth-order-estimate}.

For the first order estimate \eqref{eq:pB-parametrix-first-order-estimate}, we proceed similarly, where instead of applying the diamagnetic inequality $|\fp_B(t, x; t_0, x_0)| \leq \fp(t, x; t_0, x_0)$, we use the pointwise estimate \eqref{eq:higher-order-pB-ptwise-bd} to bound
\begin{equ}
\big|\covd \fp_B(t, x; t_0, x_0)\big| \lesssim \bigl(1 + \Fnorm{F}_{[t, t_0], 0}\bigr)  (t_0 - t)^{-\tfrac{d+1}{2}} \exp\Bigl(-c\frac{|x_0 - x|^2}{t_0 - t}\Bigr).
\end{equ}
The remainder of the proof proceeds analogously; we omit the details. 
\end{proof}

Having proven Proposition~\ref{prop:pB-parametrix}, we next show how Corollary~\ref{cor:pB-qB-parametrix} follows from Proposition~\ref{prop:pB-parametrix}, Proposition~\ref{prop:qB-ptwise-bds}, and Lemma~\ref{lemma:covd-parallel-transport-estimate}.

\begin{proof}[Proof of Corollary~\ref{cor:pB-qB-parametrix}]
By Proposition~\ref{prop:pB-parametrix}, we have that for $k = 0, 1$,
\begin{equ}
\big|\covd^{(k)} (\fp_B - \fp \ptcutoff \pt_B)(s, y; t, x)\big| \lesssim (1 + \Fnorm{F}_{[s,t], 0})^k \mc{N}(B) \times (t - s)^{-\tfrac{d+k}{2}} \exp\Bigl(-c\frac{|x-y|^2}{t-s}\Bigr),
\end{equ}
where 
\begin{equs}
\mc{N}(B) = &\, \Bigl((t-s)^{1/2} \bigl(\big\|\ptl_t B \big\|_{L_t^1 L_x^\infty([s, t])} + \bigl\|F\bigr\|_{L_t^2 L_x^\infty([s, t])} \bigr)\Bigr)^2 \\
&+ (t - s)^{1/2} \big\|\ptl_t B \big\|_{L_t^1 L_x^\infty([s, t])}  + (t-s)^{10}\bigl(1 + \big\|F\big\|_{L_t^2 L_x^\infty([s, t])}\bigr) .
\end{equs}
Combining \eqref{eq:ptl-t-B-estimate-under-more-stringent-Fpert} and \eqref{eq:F-estimate-under-more-stringent-Fpert}, we may estimate
\begin{equ}
\mc{N}(B) \lesssim  L T^{-(1-\alpha)} (t-s) s^{-\alpha} + \bigl(\big\|\ptl_t B\big\|_{L_t^1 L_x^\infty([s, t])} + \bigl\|F\bigr\|_{L_t^2 L_x^\infty([s, t])}\bigr) (t-s)^{1/2}.
\end{equ}
Combining, we obtain the estimate \eqref{eq:pB-parametrix-t-s-power}.
For \eqref{eq:covd-pB-parametrix-t-s-power}, we additionally estimate (applying \eqref{eq:covd-parallel-transport-estimate}, \eqref{eq:ptl-t-B-estimate-under-more-stringent-Fpert} and \eqref{eq:F-estimate-under-more-stringent-Fpert})
\begin{equs}
\Big|\Big(&\covd (\fp \ptcutoff \pt_B) - (\nabla \fp) \ptcutoff \pt_B\Big)(s, y; t, x)\Big| \\
&\leq | (\fp \ptcutoff \covd \pt_B)(s, y; t, x) | + |(\fp (\nabla \ptcutoff) \pt_B)(s, y; t, x)|\\
&\lesssim \Bigl(\|\ptl_t B\|_{L_t^1 L_x^\infty([s, t])} + \|F\|_{L_t^2 L_x^\infty([s, t])}\Bigr) \Bigl(1 + \frac{|x-y|}{(t-s)^{1/2}}\Bigr) \fp(s, y; t, x) + |\nabla \ptcutoff(x-y)| \fp(s, y; t, x) \\
&\lesssim \Bigl( \|\ptl_t B\|_{L_t^1 L_x^\infty([s, t])} + \|F\|_{L_t^2 L_x^\infty([s, t])}   + (t-s)^{10} \Bigr) (t-s)^{1/2} (t-s)^{-\tfrac{d+1}{2}}\exp\Bigl(-c\frac{|x-y|^2}{t-s}\Bigr),
\end{equs}
where we also used that $\nabla \ptcutoff(x-y)$ is supported on $|x - y| \sim 1$, and for such $x, y$, we have that $\fp(s, y; t, x) \lesssim (t-s)^{10} (t-s)^{-d/2} \exp(-c|x-y|^2/ (t-s))$. The estimate \eqref{eq:covd-pB-parametrix-t-s-power} now follows by combining our previous few estimates. 

Finally, for the estimate \eqref{eq:qB-parametrix-t-s-power}, by \eqref{eq:qB-pB-difference-pointwise-bounds}, we have that
\begin{equ}
\big|\covd^{(k)} (\fq_B - \fp_B)(s, y; t, x)\big| \lesssim \gronwall_{[s, t]}(C_0 F) \Fnorm{F}_{[s, t], 1} (t-s)^{-\tfrac{d+k}{2}} \exp\Bigl(-c\frac{|x-y|^2}{t-s}\Bigr), 
\end{equ}
which is acceptable. Combining this with \eqref{eq:pB-parametrix-t-s-power} and \eqref{eq:covd-pB-parametrix-t-s-power}, 
we obtain the desired result.
\end{proof}


\subsection{Perturbing around a Yang--Mills connection}\label{section:perturbing-around-YM-connection}

In this section, we discuss a different approach towards approximating the covariant heat kernel $\fq_B$ which is particularly useful when $B$ is near a Yang--Mills connection. We will use this refined approximation as one of the main inputs in Section~\ref{section:decay-estimates-in-expectation} towards obtaining decay (and thus, uniform-in-time bounds). The main result of this section is Proposition~\ref{prop:qB-perturb-near-YM}. As this result will not be used until much later, the reader may want to skip ahead first and only come back to this section when they are ready to read Section~\ref{sec:improved_stochastic_estimates}. Before we state the result, we first need several pages of preliminary discussion.

As a starting point, it turns out that when $B$ is a Yang--Mills connection on $\R^2$, there is in fact an exact formula for $\fq_B$ in terms of Mehler kernels. We proceed to discuss the Mehler kernels, and then afterwards discuss how $\fq_B$ can be exactly given in terms of these kernels when $B$ is Yang--Mills.

As in Sections~\ref{subsection:covariant-heat-kernel-estimates} and~\ref{section:heat-kernel-parametrix}, throughout we fix the base point $(s, y) = (0, 0)$, as the results proven in this section can be easily translated to the case of a general base point. It will be convenient to view $\T^2 \cong (-\pi, \pi]^2 \sse \R^2$, since we will often also view $0 \in \R^2$. Later on, we will work in polar coordinates, and so we set the requisite notation now.

\begin{notation}\label{notation:polar-coord}
Define the angular and radial vector fields $e_\theta, e_r$ on $\R^2 \backslash \{0\}$ by
\begin{equ}
e_\theta(x) \eqdef (- x_2, x_1), \qquad e_r(x) \eqdef x,
\end{equ}
as well as the function $r(x) \eqdef |x|$. We will also view $e_\theta, e_r$ as vector fields defined on a neighbourhood of $0$ in $\T^2$, e.g.\ defined for $|x| \leq 1/10$. We will often write $\ptl_\theta$ to mean $\ptl_{e_\theta}$, i.e.\ the directional derivative in the direction $e_\theta$. However, we will refrain from using the common notation $\ptl_r$ to denote $\ptl_{r^{-1} e_r}$ (i.e.\ differentiation in the unit radial direction), since from our notation for $\ptl_\theta$, one might think that $\ptl_r$ should denote $\ptl_{e_r}$.
\end{notation}

For the moment, we will work on $\R^2$, as this will simplify the discussion of the Mehler kernels. Let $\lambda \neq 0$. The Mehler kernel $\mehler_\lambda \colon (0, \infty) \times \R^2 \ra [0, \infty)$ is given by the explicit formula:
\begin{equ}\label{eq:mehler-kernel-formula}
\mehler_\lambda(t, x) = \frac{\lambda}{4\pi \sinh(\lambda t)} \exp\biggl(-\frac{\lambda|x|^2}{4\tanh(\lambda t)}\biggr).
\end{equ}
When $\lambda = 0$, we set $\mehler_0 = \fp^{\Euc}$, the heat kernel on $\R^2$. Note that $\mehler_\lambda \ra \mehler_0$ pointwise as $\lambda \ra 0$.
This kernel solves the equation
\begin{equ}\label{eq:mehler-equation}
(\ptl_t - \Delta + \frac{\lambda^2}{4} |x|^2) \mehler_\lambda = 0, ~~ \mehler_\lambda(0) = \delta_0.
\end{equ}
Because it will be needed later (in the proof of the main result of this section, Proposition~\ref{prop:qB-perturb-near-YM}), for $(s, y) \in [0, \infty) \times \R^2$, we also introduce the following version of the Mehler kernel $\mehler_\lambda^\theta(s, y; \place, \place) \colon (s, \infty) \times \R^2 \ra \C$ which includes an extra phase factor (recall Notation~\ref{notation:polar-coord}):
\begin{equ}
\mehler_\lambda^\theta(s, y; t, x) = \begin{cases}  \frac{\lambda}{4\pi \sinh(\lambda(t-s))} \exp\biggl(-\frac{\lambda|x-y|^2}{4 \tanh(\lambda(t-s))}\biggr) \exp\bigl(\icomplex (\lambda/2)(x_1 y_2 - x_2 y_1)\bigr) & \lambda \neq 0 \\
\fp^{\Euc}(s, y; t, x) & \lambda = 0.
\end{cases}
\end{equ}
When $(s, y) = (0, 0)$, we in fact have that $\mehler_\lambda^\theta(0, 0; t, x) = \mehler_\lambda(t, x)$. However, unlike for $\mehler_\lambda$, we will also need to consider the values of $\mehler_\lambda^\theta(s, y; t, x)$ at general $(s, y)$, basically because $\mehler_\lambda^\theta$ will be the integral kernel corresponding to a certain inhomogeneous linear equation that we will encounter. In that vein, we remark that $\mehler_\lambda^\theta(s, y; \place, \place)$ solves the following equation:
\begin{equ}\label{eq:mehler-theta-equation}
\Bigl(\ptl_t - \Delta + \frac{\lambda^2}{4} |x|^2 - \icomplex \lambda \ptl_\theta\Bigr) \mehler_\lambda^\theta(s, y; \place, \place) = 0, ~~ \mehler_\lambda^\theta(s, y; s, \place) = \delta_y(\place).
\end{equ}
To describe $\fq_B$ in terms of Mehler kernels, we will need to consider vector-valued versions of the equations \eqref{eq:mehler-equation} and \eqref{eq:mehler-theta-equation}. For ease of notation, we make the following definition regarding operators and kernels that will frequently appear.

\begin{definition}\label{def:YMpert-operators}
Let $X \in \frkg$. Define the operators
\begin{equ}
\mc{L}_X \eqdef \Delta + \frac{1}{4} |x|^2 \repr(X)^2 , \qquad
\mc{L}_X^\theta \eqdef \Delta + \frac{1}{4} |x|^2 \repr(X)^2 + \repr(X) \ptl_\theta, 
\end{equ}
acting on functions $\phi \colon \T^2 \ra \genV$ or $\phi \colon \R^2 \ra \genV$. (To be precise, in the former case, the operator acts on $\phi$ in a neighbourhood of $0$, so that the angular and radial vector fields are well-defined.) Next, we define operators on $\genV^{2}$-valued functions as follows. First, let $\maglapF{X} \in \End(\genV^{2})$ be defined as $\maglapF{X}(u, v) = (\repr(X) v, -\repr(X)u)$ for $u, v \in \genV$. (As we will see later, $\maglapF{X}$ is similar to a rotation, thus the name.) We then define the operators
\begin{equs}
\mc{L}_X^F &\eqdef \mc{L}_X + 2 \maglapF{X} = \Delta +  2 \maglapF{X} + \frac{1}{4} |x|^2 \repr(X)^2,\\
\mc{L}_X^{F, \theta} &\eqdef \mc{L}_X^\theta + 2 \maglapF{X} = \Delta + 2 \maglapF{X} + \frac{1}{4} |x|^2 \repr(X) ^2 +  \repr(X)\ptl_\theta .
\end{equs}
In the above, we abuse notation and also view $\repr(X) \in \End(\genV^{2})$ which sends $(u, v) \mapsto (\repr(X) u, \repr(X) v)$. 
\end{definition}

In the following, we will often abuse notation and write $X$ instead of $\repr(X)$. 

\begin{remark}
Note that in two dimensions, the curvature $F_B$ may be viewed as an operator on $\genV^{2}$, and moreover since $F_B^{12} = -F_B^{21}$, we have that $\maglapF{F_B^{12}} = F_B$. Also, for any $X \in \frkg$, $\maglapF{X}$ is a self-adjoint operator on $\genV^{2}$.
\end{remark}

As we will see later on, the action of the operators $\mc{L}_X, \mc{L}_X^\theta$ can be understood via an eigendecomposition of $X$ (viewed as an operator on $\genV$), and on each eigenspace of $X$, they act as the operators appearing on the left-hand sides of \eqref{eq:mehler-equation} and \eqref{eq:mehler-theta-equation}. This will provide the link between $\fq_B$ and the Mehler kernels.

\begin{definition}
\label{defn:s_B_eq}
Let $X \in \frkg$. Define the kernels $\fs_X(\place, \place) \colon [0, \infty) \times \R^2 \ra \End(\genV), \fs_X^{F}(\place, \place) \colon [0, \infty) \times \R^2 \ra \End(\genV^{2})$ by 
\begin{equs}[2]
(\ptl_t - \mc{L}_X) \fs_X &= 0, & \qquad \fs_X(0) &= \delta_0 \Id_{\genV} , \\
(\ptl_t - \mc{L}_X^F) \fs_X^{F} &= 0, & \qquad \fs_X^{F}(0) &= \delta_0 \Id_{\genV^{2}},
\end{equs}
where $\Id_\genV, \Id_{\genV^{2}}$ denote the identity operators on $\genV, \genV^{2}$, respectively. Note that since $\mc{L}_X, \mc{L}_X^F$ are self-adjoint, we have that $\fs_X, \fs_X^F$ are self-adjoint operators on $\genV, \genV^{2}$.
Similarly, define the kernels $\fs_X^\theta(\place, \place; \place, \place) \in \End(\genV_\C), \fs_X^{F, \theta}(\place, \place; \place, \place) \in \End(\genV_\C^{2})$ (here $\genV_\C$ denotes the complexification of $\genV$) by
\begin{equs}[2] 
(\ptl_t - \mc{L}_X^\theta) \fs_X^\theta(s, y; \place, \place) &= 0, &\qquad \fs_X^\theta(s, y; s, \place) &= \delta_y \Id_{\genV_\C},  \\
(\ptl_t - \mc{L}_X^{F, \theta}) \fs_X^{F, \theta}(s, y; \place, \place) &= 0, & \qquad \fs_X^{F, \theta}(s, y; s, \place) &= \delta_y \Id_{\genV_\C^{2}}, \label{eq:fs-B-F-theta}
\end{equs}
for any $s \geq 0$ and $y \in \R^2$.
\end{definition}

We claim that both $\fs_X, \fs_X^F$ can be explicitly described in terms of the Mehler kernel $\mehler_\lambda$. Similarly, $\fs_X^\theta, \fs_X^{F, \theta}$ may be expressed in terms of $\mehler_\lambda^\theta$. Note that $\mehler_\lambda^\theta$ is complex-valued, which is why we defined these latter two kernels as operators on the respective complexifications $\genV_\C$ and $\genV_\C^{2}$.

In the following discussion, we focus primarily on $\fs_X, \fs_X^F$, but the main points in the discussion will also apply to $\fs_X^\theta, \fs_X^{F, \theta}$. First, we recall several elementary facts from linear algebra. Since $X \in \End(\genV)$ is a skew-symmetric operator on a real vector space, it is diagonalisable when viewed as an operator on the complexification $\genV_\C$, with all eigenvalues purely imaginary. Moreover, for all nonzero eigenvalues $\icomplex \lambda$, the complex conjugate $-\icomplex \lambda$ is also an eigenvalue. A consequence of this is that for each nonzero eigenvalue $\icomplex \lambda$, there is in fact a two-dimensional subspace spanned by vectors $u, v \in \genV$ (which themselves are exactly the real and imaginary parts of a complex eigenvector corresponding to $\icomplex \lambda$) such that $X u = - \lambda v$, $X v = \lambda u$ (i.e.\ $X$ acts by rotation and scaling on this subspace). Finally, the endomorphism $X^2$ is also diagonalisable, with eigenvalues $-\lambda^2$, where $\lambda$ is such that $\icomplex \lambda$ is an eigenvalue of $X$.

From this, it follows that on each $-\lambda^2$ eigenspace of $X^2$, the kernel $\fs_X$ acts as multiplication by $\mehler_\lambda$. Put another way, if 
\begin{equ}
X^2 = O D O^T = O \begin{pmatrix} -\lambda_1^2 & 0 & \cdots & 0 \\
0 & - \lambda_2^2 & \cdots & 0 \\
0 & 0 & \cdots & -\lambda_n^2 \end{pmatrix} O^T,
\end{equ}
is a diagonalisation of $X^2$, then we may write
\begin{equ}
\fs_X = O \begin{pmatrix} \mehler_{\lambda_1} & 0 & \cdots & 0 \\
0 & \mehler_{\lambda_2} & \cdots & 0 \\
0 & 0 & \cdots & \mehler_{\lambda_n} \end{pmatrix} O^T.
\end{equ}
Next, recall that $\maglapF{X}$ is defined as an operator on $\genV^{2}$, whose block form is given by
\begin{equ}
\maglapF{X} = \begin{pmatrix} 0 & X \\ -X & 0 \end{pmatrix}.
\end{equ}
To discuss $\fs_X^F$, we first describe the action of $\maglapF{X}$ in terms of the action of $X$. Note that $\maglapF{X}$ is a symmetric operator, and thus all eigenvalues are real. For each vector $v \in \genV$ in the kernel of $X \in \End(\genV)$, we may obtain vectors $(v, 0)$ and $(0, v)$ in the kernel of $\maglapF{X}$. For each pair of vectors $u, v \in \genV$ with $X u = -\lambda v, X v = \lambda u$, we obtain a pair of eigenvectors with eigenvalue $\lambda$, namely $(u,v)$ and $(v, - u)$,
as well as a pair of eigenvectors with eigenvalue $-\lambda$, namely $(v, u)$ and $(u, - v)$.

From the linear independence of $u, v$, one may directly verify that these four eigenvectors are themselves linearly independent. We remark that all eigenspaces of $\maglapF{X}$ may be obtained in this way, which can be seen from dimensional considerations: $\maglapF{X}$ is an operator on a vector space of dimension $2 \dim(\genV)$, while $X$ is an operator on a vector space of dimension $\dim(\genV)$. We have previously discussed how the kernel of $X$ gives rise to a subspace of $\genV^{2}$ of dimension $2\dim(\mrm{ker}(X))$. And each two-dimensional subspace of $\genV$ spanned by vectors $u, v$ on which $X$ acts by rotation and scaling gives rise to a four-dimensional subspace of $\genV^{2}$. Therefore the total dimension of all eigenspaces of $\maglapF{X}$ arising in either of the two cases sums to $2 \dim(\genV)$, which means that we have a complete diagonalisation of $\maglapF{X}$.

From our previous discussion, this diagonalisation of $\maglapF{X} \in \End(\genV^{2})$ has the further property that the direct sum of eigenspaces of $\maglapF{X}$ with nonzero eigenvalues may be decomposed into four-dimensional subspaces, which themselves may be further decomposed into the direct sum of two-dimensional subspaces,
such that one of the two-dimensional subspaces has eigenvalue $\lambda$, while the other two-dimensional subspace has eigenvalue $-\lambda$. In the following, it will be natural to look at the action of $\fs_X^F$ on each of these four-dimensional subspaces.

By a short calculation, one may show that $\fs_X^F = \exp(2t \maglapF{X}) \fs_X$, in other words that
\begin{equ}
(\ptl_t - \mc{L}_X^F) (\exp(2t \maglapF{X}) \fs_X) = 0.
\end{equ}
Here, similar to before, we abuse notation and view $\fs_X \in \End(\genV^{2})$ as a block diagonal matrix, where each diagonal block is given by $\fs_X \in \End(\genV)$. Now, given $\lambda \neq 0$ and a pair of vectors $u, v \in \genV$ with $X u = -\lambda v, X v = \lambda u$ as before, note that $u, v$ are both eigenvectors of $X^2$ with eigenvalue $-\lambda^2$, and thus by the previous discussion, we have that both $u, v$ are eigenvectors of $\fs_X$ with eigenvalue $\mehler_\lambda$. It follows that on the four-dimensional subspace of $\genV^{2}$ obtained out of $u, v$, the action of $\fs_X$ (as an endomorphism on $\genV^{2}$) is given by a multiple $\mehler_\lambda$ of the identity. It therefore follows that $\fs_X^F$ acts as $e^{\pm 2 \lambda t} \mehler_\lambda$ on the two-dimensional $\pm\lambda$-eigenspaces comprising the four-dimensional subspace obtained from $u, v$. Finally, on the kernel of $\maglapF{X}$, both $\fs_X, \fs_X^F$ act as multiplication by the Euclidean heat kernel $\fp$.

To summarise the discussion, we have described explicit formulas for $\fs_X, \fs_X^F$ in terms of the eigendecomposition of $X \in \End(\genV)$ (or more precisely, the eigendecomposition of $X$ when viewed as an operator on the complexification of $\genV$). As mentioned at the beginning of the discussion, these considerations also apply to $\fs_X^\theta, \fs_X^{F, \theta}$. In particular, $\fs_X^{F, \theta}(s, y; t, x) = \exp(2(t-s) \maglapF{X}) \fs_X^\theta(s, y; t, x)$, and for each $\lambda$ such that $i \lambda$ is an eigenvalue of $X$, we have that $\mehler_\lambda^\theta$ is an eigenvalue of $\fs_X^{\theta}$. We record the following explicit estimate:
\begin{equ}\label{eq:mehler-phase-estimate}
\abs{\fs_X^{F, \theta}(s, y; t, x)} \lesssim \exp(C(t-s) \abs{X}) \fp(s, y; t, x),
\end{equ}
which we will use later.

\begin{remark}[Exponential growth]\label{eq:fq-B-exponential-growth}
It follows from this discussion that $\fs_X^F$ has eigenvalues which grow exponentially in time, as long as $X$ is not the zero operator on $\genV$. Indeed, in this case we have that $\fs_X^F$ has an eigenvalue of the form $e^{2\lambda t} \mehler_\lambda$, for some $\lambda > 0$. Recall from the formula \eqref{eq:mehler-kernel-formula} for $\mehler_\lambda$ that there is a $(\sinh(\lambda t))^{-1}$ prefactor in $\mehler_\lambda$. Since $(\sinh(\lambda t))^{-1}$ only decays like $e^{-\lambda t}$, we see that $e^{2\lambda t} \mehler_\lambda$ grows like $e^{\lambda t}$. The decay mechanism that we uncover in Section~\ref{section:decay-estimates-in-expectation} ultimately rests on this exponential growth (combined with the fact that this exponential growth appears with a good sign).
\end{remark}

\begin{remark}\label{remark:fs-X-F-radial-satisfies-L-X-F-equation}
Recall we previously observed that $\mehler_\lambda^\theta(0, 0; t, x) = \mehler_\lambda(t, x)$. Since we have related $\fs_X^F, \fs_X^{F, \theta}$ to $\mehler_\lambda, \mehler_\lambda^\theta$, it also follows that $\fs_X^{F, \theta}(0, 0; t, x) = \fs_X^F(t, x)$. In particular, we have that 
\begin{equ}
(\ptl_t - \mc{L}_X^{F, \theta}) \fs_X^F = 0.
\end{equ}
\end{remark}

We may now return to the claim made at the beginning of the section, which is if $B$ is a Yang--Mills connection, then there is an explicit formula for the kernel $\fq_B$. This is the content of the next result, which says that $\fq_B$ may be expressed in terms of $\fs_X^F$ when $B$ is Yang--Mills. The proof is deferred until later in the section. To be precise, here the kernel $\fq_B$ is defined on the plane, while up to now, we have only worked with the version of $\fq_B$ defined on the torus. We will tolerate this slight clash, because we find the result helpful for exposition purposes, as it points us towards the correct quantity around which to perturb. On the other hand, the result itself is not needed later on (only Proposition~\ref{prop:qB-perturb-near-YM} is needed). 

\begin{lemma}\label{lemma:qB-YM-exact-formula}
If $B \colon \R^2 \ra \frkg^2$ is a Yang--Mills connection, then for $X = F_B^{12}(0)$, we have that 
$\pt_B^{-1} \fq_B = \fs_{X}^F$.
\end{lemma}

\begin{remark}
In the setting of Lemma~\ref{lemma:qB-YM-exact-formula}, it is also true that $\pt_B^{-1} \fp_B = \fs_X$, though we will not prove this, as our main focus is on approximating $\pt_B^{-1} \fq_B$ (because this is the quantity which grows exponentially in time \dash recall Remark~\ref{eq:fq-B-exponential-growth}). Additionally, we choose to focus on approximating $\pt_B^{-1} \fq_B$ rather than $\fq_B$ directly (even though these two tasks should be equivalent) because in some sense the parallel transport $\pt_B$ only serves as a distraction. Geometrically, consider a vector bundle $E \ra \R^2$ with fibres modelled on $\genV$. By post-composing with $\pt_B^{-1}$, the kernel $(\pt_B^{-1} \fq_B)(t, x)$ takes values in the space $\End(E_y)$, while $\fq_B(t, x)$ itself takes values in the space $L(E_y, E_x)$. Therefore to describe $\pt_B^{-1} \fq_B$, it essentially suffices to understand its eigenvalues, and Lemma~\ref{lemma:qB-YM-exact-formula} provides exact formulas for the eigenvalues in the special case where $B$ is Yang--Mills. 
\end{remark}

With this result in mind, it is then natural to try to approximate $\pt_B^{-1} \fq_B$ by $\fs_X^F$ (where we take $X = F_B^{12}(0, 0)$) when $B$ is close to a Yang--Mills connection, which recall is the main goal of the present section. We are now almost ready to state the main result of the section (Proposition~\ref{prop:qB-perturb-near-YM}). First, we define a quantity which captures how far away $B$ is from being a Yang--Mills connection, and which appears in the quantitative estimates of the main result. 

\begin{definition}
Let $I \sse [0, \infty)$ be an interval of finite length and let $B \colon I \times \T^2 \ra \frkg^2$ be a time-varying connection one-form. Define
\begin{equs}
\label{e:YMFnorm}
\YMFnorm{B}_{I} \eqdef 
\abs{I}^{3/2} \big\|\ptl_t^2 B \big\|_{L_{t,x}^\infty(I)} + \abs{I}& \bigl(\bigl\|\ptl_t F_B \bigr\|_{L_{t,x}^\infty(I)} + \bigl\|\covd \ptl_t B \bigr\|_{L_{t,x}^\infty(I)}\bigr) \\
&+ 
\abs{I}^{1/2} \bigl(\bigl\|\covd F_B \bigr\|_{L_{t,x}^\infty(I)} + \bigl\|\ptl_t B \bigr\|_{L_{t,x}^\infty(I)} \bigr).
\end{equs}
When $I = [0, t_0]$, we will write $\YMFnorm{B}_{t_0} = \YMFnorm{B}_{I}$.
\end{definition}

In two dimensions, Yang--Mills connections $B$ satisfy $\covd F_B = 0$, and so $\YMFnorm{B}_{I} = 0$ if and only if $B$ is a constant-in-time Yang--Mills connection. In the following, recall the cutoff function $\ptcutoff$ from Definition~\ref{def:ptcutoff}, which is used to transfer between functions on $\T^2$ and functions on $\R^2$ (and vice versa). The following proposition is the more precise version of Proposition~\ref{main:prop-kernel-YM-connection}.

\begin{proposition}\label{prop:qB-perturb-near-YM}
There is a constant $C_0$ such that the following holds. Let $t_0 \leq 1$ and let $B \colon [0, t_0] \times \T^2 \ra \frkg^2$ be a smooth solution to the Yang--Mills gradient flow. Let $c_0 \geq 0$ be such that
\begin{equ}
\label{e:qB-perturb-near-YM_assumption}
t_0 \bigl(1 + t_0 \bigl\|F_B\bigr\|_{L_{t,x}^\infty([0, t_0])} \bigr)^{5} \YMFnorm{B}_{t_0} \bigl(1 + t_0 \YMFnorm{B}_{t_0}\bigr)^5  \leq c_0.
\end{equ}
Take $X = F_B^{12}(0, 0)$. We have that for all $(t, x) \in (0, t_0] \times \T^2$,
\begin{equ}
\big| \ptcutoff \bigl(\pt_B^{-1} \fq_B - \fs_X^F\bigr)(t, x) \big| \lesssim \bigl(c_0 + t^{10}\bigr) \exp(C_0 t_0 \|F_B\|_{L_{t,x}^\infty([0, t_0])}) t^{-1} \exp(-|x|^2 / (C_0 t)).
\end{equ}
\end{proposition}

\begin{remark}
While reading the arguments in this section, it may be good to keep in mind the following model scenario (see also Lemma~\ref{le:exists_YM_close_interval}): $B$ is a solution to the Yang--Mills gradient flow, and moreover for some $N \gg 1$ and some $c \ll 1$, we have that $t_0 = c N^{-1}\log \log N$, and that
\begin{equ}
\bigl\|F_B\bigr\|_{L_{t,x}^\infty([0, t_0])} \sim N, ~~  \|\covd^{(j)} F_B\|_{L_{t,x}^\infty([0, t_0])} \lesssim \frac{N^{1 + j/2}}{\log N} \text{ for $1 \leq j \leq 3$}.
\end{equ}
One can heuristically think of $\covd F_B$ as being supported on frequencies $\lesssim N^{1/2}$, which is consistent with the second estimate, which says that each additional derivative costs a factor of $N^{1/2}$. It is however important to note the extra $\log N$ gain in the estimate of the derivatives (as compared to $F_B$ itself), which can be interpreted as saying that $B$ is almost a Yang--Mills connection (whose curvature is large).

Continuing, we recall that since $B$ is a solution to the Yang--Mills gradient flow, we have that $\ptl_t F_B = \Delta_B F_B$ and $\ptl_t \covd F_B = \covd^{(3)} F_B + \schO([F_B, \covd F_B])$, which then implies 
$\big\|\ptl_t \covd F_B \big\|_{L_{t,x}^\infty} \lesssim N^{5/2}/\log N$.
It thus follows that 
\begin{equ}
t_0 \YMFnorm{B}_{t_0} \lesssim c \frac{(\log \log N)^{5/2}}{\log N}\ll 1.
\end{equ}
In the context of Proposition~\ref{prop:qB-perturb-near-YM}, we may take $c_0 \sim c (\log \log N)^{100} / \log N$, and the prefactor appearing in Proposition~\ref{prop:qB-perturb-near-YM} may then be estimated (assuming that $c \ll 1$ is sufficiently small depending on $C_0$)
\begin{equ}
\bigl(c_0 + t_0^{10}\bigr) \exp(C_0 t_0 \|F_B\|_{L_{t,x}^\infty([0, t_0])}) \lesssim c \frac{(\log \log N)^{100}}{\log N} \exp(C_0 c \log \log N) \lesssim \frac{1}{(\log N)^{1-c_1}},
\end{equ}
where $c_1 > 0$ is small if $c$ is small. Therefore in this setting, $\pt_B^{-1} \fq_B - \fs_X^F$ is indeed a perturbation.
\end{remark}

We now proceed to prove Lemma~\ref{lemma:qB-YM-exact-formula} and Proposition~\ref{prop:qB-perturb-near-YM}. While the lemma is not strictly needed in the following, its proof is a nice warm-up towards the proof of the proposition. To begin, recall that $\fq_B$ satisfies
\begin{equ}
(\ptl_t - \Delta_B - 2 F_B) \fq_B = 0, \quad \fq_B(0) = \delta_0 \Id_{\genV^{2}}.
\end{equ}
As an intermediate step, we derive the equation satisfied by $\pt_B^{-1} \fq_B$. This will also be useful in the proof of Proposition~\ref{prop:qB-perturb-near-YM}. The following derivation applies both when $\fq_B$ is defined on the torus and when it is defined on the plane.

\begin{lemma}\label{lemma:pt-inverse-qB-equation}
We have that
\begin{equs}
(\ptl_t - \Delta) (\pt_B^{-1} \fq_B) &= -\pt_B^{-1} \bigl((\ptl_t - \Delta_B) \pt_B \bigr) \pt_B^{-1} \fq_B + 2 (\pt_B^{-1} F_B \pt_B) \pt_B^{-1}  \fq_B \\
&\quad + 2 (\pt_B^{-1} \covd_{r^{-1}e_\theta} \pt_B) \ptl_{r^{-1}e_\theta} (\pt_B^{-1} \fq_B) + 2 (\pt_B^{-1} \covd_{r^{-1} e_r} \pt_B) \ptl_{r^{-1} e_r} (\pt_B^{-1} \fq_B).
\end{equs}
\end{lemma}
\begin{proof}
To begin, we compute by the product rule
\begin{equ}
(\ptl_t - \Delta) (\pt_B^{-1} \fq_B) = \bigl((\ptl_t - \Delta_B) (\pt_B^{-1})\bigr) \fq_B + \pt_B^{-1} (\ptl_t - \Delta_B) \fq_B - 2 \covd_\ell (\pt_B^{-1}) \covd^\ell \fq_B.
\end{equ}
We will need the intermediate computations (recall Remark~\ref{remark:parallel-transport}\ref{remark:pt-inverse}, in particular \eqref{eq:covd-pt-inverse}):
\begin{equs}
\ptl_t (\pt_B^{-1}) &= -\pt_B^{-1} (\ptl_t \pt_B )\pt_B^{-1}, \\
\covd_\ell(\pt_B^{-1}) &= -\pt_B^{-1} (\covd_\ell \pt_B) \pt_B^{-1}, \\
\Delta_B (\pt_B^{-1}) &= -\covd^\ell (\pt_B^{-1} (\covd_\ell \pt_B) \pt_B^{-1}) = -\pt_B^{-1} (\Delta_B \pt_B) \pt_B^{-1} + 2 \pt_B^{-1} (\covd_\ell \pt_B) \pt_B^{-1} (\covd^\ell \pt_B) \pt_B^{-1}, \\
\covd^\ell (\fq_B) &= \covd^\ell (\pt_B \pt_B^{-1} \fq_B) = (\covd^\ell \pt_B) \pt_B^{-1} \fq_B + \pt_B \ptl^\ell (\pt_B^{-1} \fq_B).
\end{equs}
Inserting these intermediate computations, we may further obtain   
\begin{equ}
(\ptl_t - \Delta) (\pt_B^{-1} \fq_B) = - \pt_B^{-1} \bigl((\ptl_t - \Delta_B) \pt_B \bigr) \pt_B^{-1} \fq_B + 2 \pt_B^{-1} F_B  \fq_B + 2 \pt_B^{-1} (\covd_\ell \pt_B) \ptl^\ell (\pt_B^{-1} \fq_B).
\end{equ}
To finish, we use that $r^{-1} e_r$, $r^{-1} e_\theta$ are orthonormal, and thus
\begin{equ}
(\covd_\ell \pt_B) \ptl^\ell (\pt_B^{-1} \fq_B) = (\covd_{r^{-1} e_r} \pt_B) \ptl_{r^{-1} e_r} (\pt_B^{-1} \fq_B) + (\covd_{r^{-1} e_\theta} \pt_B) \ptl_{r^{-1} e_\theta}(\pt_B^{-1} \fq_B) .
\end{equ}
The desired result now follows by combining the previous few identities.
\end{proof}

From Lemma~\ref{lemma:pt-inverse-qB-equation}, we see that the Laplacian $\Delta_B \pt_B$ of the
parallel transport map enters into the equation for $\pt_B^{-1} \fq_B$. It will be useful to express the Laplacian in polar coordinates, which is the purpose of the next lemma. This is a standard calculation but we prefer to include it. Recall the angular and radial vector fields ($e_\theta, e_r$) from Notation~\ref{notation:polar-coord}. As before, the following derivation holds both on the plane as well as on a neighbourhood around $0$ in the torus. In the following, recall the covariant Hessian introduced in Notation~\ref{notation:covariant-hessian}.

\begin{lemma}[Polar coordinates]\label{lemma:laplacian-polar-coord}
We have that
\begin{equs}
(\Hess_B \Phi)(r^{-1} e_r, r^{-1} e_r) &= \covd_{r^{-1} e_r} \covd_{r^{-1} e_r} \Phi, \\
(\Hess_B \Phi)(r^{-1} e_\theta, r^{-1} e_\theta) &= \bigl(\covd_{r^{-1} e_\theta} \covd_{r^{-1} e_\theta} + r^{-1} \covd_{r^{-1} e_r} \bigr) \Phi .
\end{equs}
Consequently, we have that
\begin{equ}\label{eq:laplacian-polar-coord}
\Delta_B \Phi = \bigl(\covd_{r^{-1} e_r} \covd_{r^{-1} e_r} + r^{-1} \covd_{r^{-1} e_r} + \covd_{r^{-1} e_\theta} \covd_{r^{-1} e_\theta}\bigr) \Phi.
\end{equ}
\end{lemma}
\begin{proof}
We compute
\begin{equ}
\covd_{r^{-1} e_r} \covd_{r^{-1} e_r} \Phi = (\Hess_B \Phi)(r^{-1} e_r, r^{-1} e_r) + \covd_{\ptl_{r^{-1} e_r} (r^{-1} e_r)} \Phi = (\Hess_B \Phi)(r^{-1} e_r, r^{-1} e_r),
\end{equ}
where we used that $\ptl_{r^{-1} e_r} (r^{-1} e_r) = 0$. Next, we compute
\begin{equ}
\covd_{r^{-1} e_\theta} \covd_{r^{-1} e_\theta} \Phi = (\Hess_B \Phi)(r^{-1} e_\theta, r^{-1} e_\theta) + \covd_{\ptl_{r^{-1} e_\theta} (r^{-1} e_\theta)} \Phi.
\end{equ}
We have that
\begin{equ}\label{eq:derivative-angular-angular}
\ptl_{r^{-1} e_\theta} (r^{-1} e_\theta) = r^{-1} \ptl_{r^{-1} e_\theta} e_\theta = r^{-2} \ptl_{e_\theta} e_\theta = - r^{-2} e_r.
\end{equ}
Since $r^{-1} e_\theta$ and $r^{-1} e_r$ form an orthonormal basis of $\R^2$, we obtain the claimed formula for the covariant Laplacian by tracing over the Hessian over this orthonormal basis and using the formulas for the Hessian that we just established.
\end{proof}

The remaining ingredient in the proof of Lemma~\ref{lemma:qB-YM-exact-formula} is an exact computation of covariant derivatives of parallel transport in the case where $B$ is a Yang--Mills connection. This exact computation captures much of what is special about the geometry associated to a Yang--Mills connection. In the following, recall Remark~\ref{remark:parallel-transport}\ref{remark:A-constant-in-time-ptA-constant-in-time}, which says that if $B$ is constant in time, then so is its parallel transport $\pt_B$.

\begin{lemma}\label{lemma:YM-connection-covd-pt}
Let $B \colon \R^2 \ra \frkg^2$ be a Yang--Mills connection. Then for any point $x\in \R^2$ and vector $v \in \R^2$, we have that
\begin{equs}
(\pt_B^{-1} \covd_{v} \pt_B)(x) &= \frac{1}{2} F_B(0)(x, v), \\
(\pt_B^{-1} \covd_{r^{-1} e_\theta} \pt_B)(x) &= \frac{1}{2} |x| F_B^{12}(0), \\
(\pt_B^{-1}\Delta_B \pt_B)(x) &= \frac{1}{4} |x|^2 (F_B^{12}(0))^2.
\end{equs}
\end{lemma}
\begin{proof}
By Lemmas~\ref{lemma:covd-parallel-transport} and~\ref{lemma:pt-curvature-commutators} (or more precisely, analogues of these lemmas on the plane for constant-in-time connections), we have that (with $\gamma(t) = tx$)
\begin{equs}
(\pt_B^{-1} \covd_{v} \pt_B)(x) &= \int_0^1 \dt \, t \cdot  (\pt_B^{-1} F_B(x, v) \pt_B)(\gamma(t)) \\
&= \int_0^1 \dt \, t \cdot F_B(0)(x, v) = \frac{1}{2} F_B(0)(x, v),
\end{equs}
which gives the first desired identity. 

For the second identity, we apply the first identity to obtain
\begin{equs}
(\pt_B^{-1} \covd_{r^{-1} e_\theta} \pt_B)(x) &= \frac{1}{2} F_B(0)(x, r^{-1} e_\theta)\\
&= \frac{1}{2}\Big( r^{-1} F_B^{12}(0) x_1^2 + r^{-1} F_B^{21}(0) x_2 (-x_2) \Big) = \frac{1}{2}|x| F_B^{12}(0),
\end{equs}
as desired. 

Finally, for the last identity,
we differentiate both sides above by $\ptl_{r^{-1} e_\theta}$ to obtain (it may be helpful to read the first identity from right to left)
\begin{equs}
\pt_B^{-1}  \covd_{r^{-1} e_\theta} \covd_{r^{-1} e_\theta} \pt_B - \pt_B^{-1} (\covd_{r^{-1}e_\theta} \pt_B) \pt_B^{-1} \covd_{r^{-1} e_\theta} \pt_B  &= \ptl_{r^{-1} e_\theta} (\pt_B^{-1} \covd_{r^{-1} e_\theta} \pt_B) \\
&= \frac{1}{2} (\ptl_{r^{-1} e_\theta} r) F_B^{12}(0) = 0,
\end{equs}
where we used that $r = |x|$ is a radial function, so that its angular derivative is zero. Note that by definition, $\covd_{e_r} \pt_B = 0$.
Rearranging and combining this with \eqref{eq:laplacian-polar-coord}, we thus obtain
\begin{equ}
(\pt_B^{-1} \Delta_B \pt_B)(x) = (\pt_B^{-1} \covd_{r^{-1} e_\theta} \pt_B)(x) (\pt_B^{-1} \covd_{r^{-1} e_\theta} \pt_B )(x) = \frac{1}{4} |x|^2 (F_B^{12}(0))^2,
\end{equ}
as desired.
\end{proof}

Because it will be useful later, we record the following calculation from the proof of Lemma~\ref{lemma:YM-connection-covd-pt}:
\begin{equ}\label{eq:curvature-radial-angular-calculation}
F_B(x, e_\theta) = F_B^{12} x_1^2 + F_B^{21} x_2 (-x_2) = F_B^{12} |x|^2.
\end{equ}
We are now ready to prove Lemma~\ref{lemma:qB-YM-exact-formula}.

\begin{proof}[Proof of Lemma~\ref{lemma:qB-YM-exact-formula}]
By Lemma~\ref{lemma:pt-inverse-qB-equation} combined with the facts that when $B$ is Yang--Mills, $\pt_B$ is independent of time and that $\covd_{e_r} \pt_B = 0$ (by definition), we obtain
\begin{equ}
(\ptl_t - \Delta) (\pt_B^{-1} \fq_B) = 2 (\pt_B^{-1} F_B \pt_B) (\pt_B^{-1} \fq_B) + (\pt_B^{-1} \Delta_B \pt_B) (\pt_B^{-1} \fq_B) + 2 (\pt_B^{-1} \covd_{r^{-1} e_\theta} \pt_B) \ptl_{r^{-1} e_\theta}(\pt_B^{-1} \fq_B).
\end{equ}
Next, by Lemmas~\ref{lemma:pt-curvature-commutators} and~\ref{lemma:YM-connection-covd-pt} (the first of which implies that $\pt_B^{-1} F_B \pt_B = F_B(0)$), we further have that 
\begin{equ}
(\ptl_t - \mc{L}_X^F) (\pt_B^{-1} \fq_B) = 2 (\pt_B^{-1} \covd_{r^{-1} e_\theta} \pt_B) \ptl_{r^{-1} e_\theta}(\pt_B^{-1} \fq_B),
\end{equ}
where $\mc{L}_X^F$ is defined in Definition~\ref{def:YMpert-operators}. To finish, observe that $\fs_X^F$ is radial, i.e.\ $\ptl_{e_\theta} \fs_X^F = 0$, and thus it also solves
\begin{equ}
(\ptl_t - \mc{L}_X^F) \fs_X^F =  2 (\pt_B^{-1} \covd_{r^{-1} e_\theta} \pt_B) \ptl_{r^{-1} e_\theta} \fs_X^F.
\end{equ}
Since $\fs_X^F(0) = \delta_0 \Id_{\genV^{2}} = (\pt_B^{-1} \fq_B)(0)$, it follows that $\fs_X^F = \pt_B^{-1} \fq_B$, as desired.
\end{proof}


The proof of Lemma~\ref{lemma:qB-YM-exact-formula} also suggests how to go about proving Proposition~\ref{prop:qB-perturb-near-YM}: we should derive the equation satisfied by $\pt_B^{-1} \fq_B - \fs_X^F$, which should be of the form $(\ptl_t - \mc{L}_X^{F, \theta}) (\pt_B^{-1} \fq_B - \fs_X^F) = G_B$, where $G_B = 0$ if $B$ is Yang--Mills, and more generally $G_B$ should be perturbative if $B$ is almost a Yang--Mills connection. Towards this end, define the remainder kernel (recall that we take $X = F_B^{12}(0, 0)$)
\begin{equ}
\fr_B^F \eqdef \ptcutoff(\pt_B^{-1} \fq_B - \fs_X^F) \colon [0, \infty) \times \R^2 \ra \End(\genV^{2}).
\end{equ}
Here, while $\pt_B^{-1} \fq_B$ is technically only defined on the torus, we use the cutoff function $\ptcutoff$ to view $\ptcutoff \pt_B^{-1} \fq_B$ as a function on the plane. Note that $\fr_B^F$ is precisely the quantity that we need to estimate in order to prove Proposition~\ref{prop:qB-perturb-near-YM}.

Combining Remark~\ref{remark:fs-X-F-radial-satisfies-L-X-F-equation} (which gives that $(\ptl_t - \mc{L}_X^{F, \theta}) \fs_X^F = 0$) with Lemma~\ref{lemma:pt-inverse-qB-equation}, we obtain the following equation for $\fr_B^F$:
\begin{equs}\label{eq:rB-F-equation}
(\ptl_t -  \mc{L}_X^{F, \theta}) \fr_B^F = \,&2 (\pt_B^{-1} F_B \pt_B - \maglapF{X})  (\ptcutoff \pt_B^{-1} \fq_B) \\
&\quad + \Bigl(\pt_B^{-1} \Delta_B \pt_B - \frac{r^2}{4} X^2\Bigr) (\ptcutoff \pt_B^{-1} \fq_B) \\
&\quad + 2 \ptcutoff \Bigl(\pt_B^{-1} \covd_{r^{-1}e_\theta} \pt_B - \frac{r}{2} X\Bigr) \ptl_{r^{-1}e_\theta} (\pt_B^{-1} \fq_B) \\
&\quad - (\pt_B^{-1} \ptl_t \pt_B) (\ptcutoff \pt_B^{-1} \fq_B) + 2 \ptcutoff (\pt_B^{-1} \covd_{r^{-1} e_r} \pt_B) \ptl_{r^{-1} e_r} (\pt_B^{-1} \fq_B) \\
&\quad- r X (\ptl_{r^{-1} e_\theta} \ptcutoff) (\pt_B^{-1} \fq_B - \fs_X^F) \\
&\quad - 2 (\ptl^\ell \ptcutoff) \ptl_\ell (\pt_B^{-1} \fq_B - \fs_X^F) - (\Delta \ptcutoff)(\pt_B^{-1} \fq_B - \fs_X^F).
\end{equs}
The main step towards Proposition~\ref{prop:qB-perturb-near-YM} is to estimate all terms on the right-hand side of \eqref{eq:rB-F-equation} in terms of the quantity $\YMFnorm{B}_{t_0}$ (which we recall captures how close $B$ is to being a Yang--Mills connection). The next lemma is a first step towards the needed estimates.

\begin{lemma}[Commutator estimates of parallel transport and curvature]\label{lemma:pt-curvature-commutator}
Let $B \colon [0, t_0] \times \T^2 \ra \frkg^2$ be a time-varying connection one-form. Let $u, w \in \R^2$, $\gamma(t) = t(1, w)$. For $t \in [0, t_0]$ such that $t|w| \leq 1/10$, we have that for $v \in \R^2$,
\begin{equ}\label{eq:pt-curvature-commutator-gamma-t}
\big| (\pt_B^{-1} F_B(u, v) \pt_B)(\gamma(t)) - F_B(\gamma(0))(u, v) \big| \lesssim t |u| |v| \bigl( \big\|\ptl_t F_B \big\|_{L_{t,x}^\infty([0, t])} +  \big\|\covd F_B \big\|_{L_{t, x}^\infty([0, t])} |w|\bigr).
\end{equ}
If, more generally, $v \in \R^{1+2}$, we also have the estimates
\begin{equs} 
~\hspace{-5mm}&\big| \bigl(\pt_B^{-1} F_B(\dot{\gamma}(t), v) \pt_B\bigr)(\gamma(t)) - F_B(\gamma(0)) (\dot{\gamma}(0), v) \big| \\
&\lesssim t\Bigl( \big\|\ptl_t^2 B\big\|_{L_{t, x}^\infty([0, t])} (|v_x| + |w| |v_0|) + \bigl\|\ptl_t F_B \bigr\|_{L_{t,x}^\infty([0, t])} |w| |v_x| \\
&\qquad \qquad+ \big\|\covd \ptl_t B \big\|_{L_{t, x}^\infty([0, t])} (|w| |v_x| + |w|^2 |v_0|)  + \big\|\covd F_B \big\|_{L_{t, x}^\infty([0, t])} |w|^2 |v_x| \Bigr), \label{eq:pt-F-commutation-estimate}\\
~&\big| (\pt_B^{-1} \covd_v \pt_B)(\gamma(t)) - \frac{t}{2} F_B(\gamma(0))(\dot{\gamma}(0), v)\big| \\
&\lesssim t^2 \Bigl( \big\|\ptl_t^2 B\big\|_{L_{t, x}^\infty([0, t])} (|v_x| + |w| |v_0|) + \bigl\|\ptl_t F_B \bigr\|_{L_{t,x}^\infty([0, t])} |w| |v_x| \\
&\qquad \qquad +\big\|\covd \ptl_t B \big\|_{L_{t, x}^\infty([0, t])} (|w| |v_x| + |w|^2 |v_0|)  + \big\|\covd F_B \big\|_{L_{t, x}^\infty([0, t])} |w|^2 |v_x| \Bigr).\label{eq:covd-pt-estimate}
\end{equs}
\end{lemma}
\begin{proof}
By Lemma~\ref{lemma:pt-curvature-commutators}, we have that
\begin{equ}
\bigl(\pt_B^{-1} F_B(u, v) \pt_B\bigr)(\gamma(t)) - F_B(\gamma(0))(u, v) = \int_0^t \ds \, \bigl(\pt_B^{-1} (\covd_{\dot{\gamma}(s)} F_B)(u, v) \pt_B\bigr)(\gamma(s)).
\end{equ}
The estimate \eqref{eq:pt-curvature-commutator-gamma-t} now follows by integrating the pointwise estimate (note that $\covd_{\dot{\gamma}(s)} F_B = \ptl_t  F_B + \covd_w F_B$)
\begin{equ}
\big|\bigl(\pt_B^{-1} (\covd_{\dot{\gamma}(s)} F_B)(u, v) \pt_B\bigr)(\gamma(s))\big| \lesssim \big\|\ptl_t F_B \big\|_{L_{t,x}^\infty} |u| |v| + \big\|\covd F_B \big\|_{L_{t,x}^\infty} |u| |v| |w|.
\end{equ}
Here, we used the assumption that $v \in \R^2$ is a spatial vector, so that only at most one derivative of $F_B$ appears.

Next, again by Lemma~\ref{lemma:pt-curvature-commutators}, we have that 
\begin{equ}
\bigl(\pt_B^{-1} F_B(\dot{\gamma}(t), v) \pt_B\bigr)(\gamma(t)) - F_B(\gamma(0))(\dot{\gamma}(0), v)  = \int_0^t \ds \, \bigl(\pt_B^{-1} (\covd_{\dot{\gamma}(s)} F_B)(\dot{\gamma}(s), v) \pt_B\bigr)(\gamma(s)).
\end{equ}
Note that $\dot{\gamma}(s) = (1, w)$, and thus we may expand
\begin{equs}
~&(\covd_{\dot{\gamma}(s)} F)(\dot{\gamma}(s), v) \\
&= (\covd_{(1, w)} F)((1, w), v) \\
&= \covd_0 F_{0j} v^j + \covd_0 F_{ij} w^i v^j + \covd_0 F_{i0} w^i v_0 + w^i \covd_i F_{0j} v^j + w^i \covd_i F_{jk} w^j v^k + w^i \covd_i F_{j0} w^j v^0.
\end{equs}
Noting that $F_{0i} = \ptl_t B_i$ and $\covd_0 = \ptl_t$, we may thus obtain the following pointwise estimate:
\begin{equs}
\big| \bigl(\pt_B^{-1} (\covd_{ \dot{\gamma}(s)} F_B)(\dot{\gamma}(s), v) \pt_B\bigr)(\gamma(s))&\big| \lesssim \big\|\ptl_t^2 B\big\|_{L_{t, x}^\infty} (|v_x| + |w| |v_0|)+ \bigl\|\ptl_t F_B \bigr\|_{L_{t,x}^\infty} |w| |v_x| \\
&\qquad  + \big\|\covd \ptl_t B \big\|_{L_{t, x}^\infty} (|w| |v_x| + |w|^2 |v_0|) + \big\|\covd F_B \big\|_{L_{t, x}^\infty} |w|^2 |v_x|.
\end{equs}
Upon integrating, we obtain \eqref{eq:pt-F-commutation-estimate}. For \eqref{eq:covd-pt-estimate}, we integrate once more and again apply Lemma~\ref{lemma:pt-curvature-commutators}.
\end{proof}

\begin{lemma}\label{lemma:F-expansion-computation-errors}
Let $B \colon [0, t_0] \times \T^2 \ra \frkg^2$ be a time-varying connection one-form. Let $x \in \T^2$ be such that $|x| \leq 1/10$, and let $t \in (0, t_0]$. We have that
\begin{equs}
\bigl| (t/2) F_B((1, x/t), r^{-1} e_\theta) - (r/2) F_B^{12}(0, 0)\bigr| &\lesssim  t^{1/2} \YMFnorm{B}_{t_0}, \label{eq:F-radial-angular-expansion}\\
\bigl| t F_B((1, x/t), r^{-1} e_r) \bigr| &\lesssim t^{1/2} \YMFnorm{B}_{t_0}, \label{eq:F-radial-radial-expansion}\\
\bigl|t  F_B((1, x/t), (1, 0)) \bigr| &\lesssim  \YMFnorm{B}_{t_0} (rt^{-1/2}). \label{eq:F-radial-temporal-expansion}
\end{equs}
To be clear, in the above $F_B$ is always evaluated at the space-time point $(0, 0)$.
\end{lemma}
\begin{proof}
To show \eqref{eq:F-radial-angular-expansion}, we first expand (and use \eqref{eq:curvature-radial-angular-calculation})
\begin{equs}
(t/2) F_B((1, x/t), r^{-1} e_\theta) &= (t/2) F_B(x/t, r^{-1} e_\theta) + (t/2) \ptl_t B \cdot (r^{-1} e_\theta) \\
&= (r/2) F_B^{12}(0, 0) + (t/2) \ptl_t B \cdot (r^{-1} e_\theta),
\end{equs}
and to finish, we estimate
\begin{equ}
t \big|\ptl_t B \cdot (r^{-1} e_\theta)| \lesssim t \big\|\ptl_t B \big\|_{L_{t, x}^\infty} \lesssim t^{1/2} \YMFnorm{B}_{t_0}.
\end{equ}
To show \eqref{eq:F-radial-radial-expansion}, we first note that $F_B((1, x/t), r^{-1} e_r) = (\ptl_t B) \cdot (r^{-1} e_r)$ and then apply a similar estimate as above. Finally, to show \eqref{eq:F-radial-temporal-expansion}, we estimate
\[
\big| t F_B((1, x/t), (1, 0))\big| = \big| t (\ptl_t B) \cdot (x/t)\big| \lesssim \YMFnorm{B}_{t_0} |x| t^{-1/2} = \YMFnorm{B}_{t_0} (rt^{-1/2}). \qedhere
\]
\end{proof}

Using Lemmas~\ref{lemma:pt-curvature-commutator} and~\ref{lemma:F-expansion-computation-errors}, we may now estimate all error terms appearing on the right-hand side of \eqref{eq:rB-F-equation}.

\begin{lemma}\label{lemma:YM-perturb-main-intermediate-estimate-clean}
Let $B \colon [0, t_0] \times \T^2 \ra \frkg^2$ be a time-varying connection one-form. Let $x \in \T^2$ be such that $|x| \leq 1/10$, and let $t \in (0, t_0]$. We have that
\begin{equs}
\big| \bigl(\pt_B^{-1} F_B \pt_B\bigr)(t, x) - F_B(0, 0)\big| &\lesssim \YMFnorm{B}_{t_0} \bigl(1 +  (rt^{-1/2})\bigr), \label{eq:curvature-pt-commutation-clean}\\
\big| (\pt_B^{-1} \covd_{r^{-1} e_\theta} \pt_B)(t, x) - \frac{r}{2} F_B^{12}(0, 0) \big| &\lesssim t^{1/2} \YMFnorm{B}_{t_0} \bigl(1 + (rt^{-1/2})^2 \bigr), \label{eq:pt-angular-derivative-estimate-clean}\\
\big| \pt_B^{-1} \covd_{r^{-1} e_r} \pt_B (t, x) \big| &\lesssim t^{1/2} \YMFnorm{B}_{t_0}\bigl(1 + (rt^{-1/2})^2 \bigr), \label{eq:pt-radial-derivative-estimate-clean} \\
\big| (\pt_B^{-1} \ptl_t \pt_B)(t, x) \big| &\lesssim \YMFnorm{B}_{t_0} \bigl(1 + (rt^{-1/2})^2 \bigr), \label{eq:pt-time-derivative-estimate-clean} \\
\hspace{-5mm}\big|(\pt_B^{-1} \Hess_B \pt_B)(r^{-1} e_\theta, r^{-1} e_\theta)(t, x) - \frac{r^2}{4} (F_B^{12}(0, 0))^2 \big| &\lesssim  \bigl(1 + t\big\|F \big\|_{L_{t, x}^\infty([0, t_0])} + t\YMFnorm{B}_{t_0}\bigr) \label{eq:pt-hessian-angular-estimate-clean}\\
&\quad \times \YMFnorm{B}_{t_0} \bigl(1 + (rt^{-1/2})^5 \bigr) ,  \\
\big|(\pt_B^{-1} \Hess_B \pt_B)(r^{-1} e_r, r^{-1} e_r)(t, x)\big| &\lesssim \bigl(1 + t \YMFnorm{B}_{t_0}\bigr) \label{eq:pt-hessian-radial-estimate-clean}\\
&\quad \times \YMFnorm{B}_{t_0} \bigl(1 + (rt^{-1/2})^2 \bigr),  \\
\big| (\pt_B^{-1} \Delta_B \pt_B)(t, x) - \frac{r^2}{4} (F_B^{12}(0, 0))^2 \big| &\lesssim \bigl(1 + t\big\|F \big\|_{L_{t, x}^\infty([0, t_0])} + t\YMFnorm{B}_{t_0}\bigr) \label{eq:pt-laplacian-estimate-clean} \\
&\quad \times \YMFnorm{B}_{t_0}\bigl( 1 + (rt^{-1/2})^5 \bigr). 
\end{equs}
\end{lemma}
\begin{proof}
We proceed with the estimates one-by-one.

\proofstep{Proof of \eqref{eq:curvature-pt-commutation-clean}.} We apply \eqref{eq:pt-curvature-commutator-gamma-t} with $w = x/t$ and arbitrary $u, v \in \R^2, |u| = |v| = 1$ to estimate
\begin{equ}
\big| (\pt_B^{-1} F_B \pt_B)(t, x) - F_B(0, 0) \big| \lesssim \big\|\ptl_t F_B \big\|_{L_{t,x}^\infty} t + \big\|\covd F_B \big\|_{L_{t, x}^\infty} |x| \lesssim \YMFnorm{B}_{t_0}\bigl(1 +  (rt^{-1/2})\bigr).
\end{equ}

\proofstep{Proofs of \eqref{eq:pt-angular-derivative-estimate-clean} and \eqref{eq:pt-radial-derivative-estimate-clean}.}
We apply \eqref{eq:covd-pt-estimate} with $v = (0, r^{-1} e_\theta)$, $w = x /t$, to obtain 
\begin{equs}
\big|(\pt_B^{-1} &\covd_{r^{-1} e_\theta} \pt_B)(t, x) - \frac{t}{2} F_B((1, w), r^{-1} e_\theta) \big| \\
&\lesssim t^2 \bigl( \big\|\ptl_t^2 B \big\|_{L_{t, x}^\infty} + \big\|\ptl_t F_B \big\|_{L_{t, x}^\infty} |x| t^{-1} + \bigl\|\covd \ptl_t B \bigr\|_{L_{t,x}^\infty} |x| t^{-1} + \big\|\covd F_B \big\|_{L_{t, x}^\infty} |x|^2 t^{-2}\bigr) \\
&\lesssim t^{1/2} \YMFnorm{B}_{t_0} \bigl(1 + (rt^{-1/2})^2 \bigr).
\end{equs}
Combining this with \eqref{eq:F-radial-angular-expansion}, we have thus shown \eqref{eq:pt-angular-derivative-estimate-clean}. The proof of \eqref{eq:pt-radial-derivative-estimate-clean} is similar, except we take $v = (0, r^{-1} e_r)$, and we apply \eqref{eq:F-radial-radial-expansion} in place of \eqref{eq:F-radial-angular-expansion}.

\proofstep{Proof of \eqref{eq:pt-time-derivative-estimate-clean}.} We apply \eqref{eq:covd-pt-estimate} with $v = (1, 0)$, $w = x/t$, to obtain
\begin{equs}
\big| (\pt_B^{-1} \ptl_t \pt_B)(t, x) - \frac{t}{2} F_B((1, w), (1, 0)) \big| &\lesssim t^2 \bigl( \big\|\ptl_t^2 B \big\|_{L_{t, x}^\infty} |x| t^{-1} + \big\|\covd \ptl_t B \big\|_{L_{t, x}^\infty} |x|^2 t^{-2} \bigr) \\
&\lesssim \YMFnorm{B}_{t_0} \bigl(1 + (r t^{-1/2})^2 \bigr).
\end{equs}
Combining this with \eqref{eq:F-radial-temporal-expansion}, we may thus obtain \eqref{eq:pt-time-derivative-estimate-clean}.

\proofstep{Proof of \eqref{eq:pt-hessian-angular-estimate-clean}.} We want to apply Lemma~\ref{lemma:covd-parallel-transport}, but we need to be careful because $r^{-1} e_\theta$ is non-constant, in particular $\ptl_{r^{-1} e_\theta}(r^{-1} e_\theta) = - r^{-1} (r^{-1} e_r) \neq 0$ (recall \eqref{eq:derivative-angular-angular}). Towards this end, let $v_1 = v_2 = (r^{-1} e_\theta)(x)$. Note then that $(\Hess_B \pt_B)(r^{-1} e_\theta, r^{-1} e_\theta)(t, x) = (\covd_{v_1} \covd_{v_2} \pt_B)(t, x)$ (the point now is that the vectors $v_1, v_2$ do not vary in space, unlike the vector field $r^{-1} e_\theta$).
Applying Lemma~\ref{lemma:covd-parallel-transport}, we thus obtain (here $\gamma(s) = s (1, x/t)$)
\begin{equs}
(\pt_B^{-1} \Hess_B \pt_B)(r^{-1} &e_\theta, r^{-1} e_\theta)(t, x) \\
&=
(\pt_B^{-1} \covd_{v_1} \covd_{v_2} \pt_B)(t, x) \\
&= \frac{2}{t^2} \int_0^t \ds \, s^2 \bigl(\pt_B^{-1} F_B(\dot{\gamma}(s), r^{-1}e_\theta) \pt_B\bigr)(\gamma(s)) (\pt_B^{-1} \covd_{r^{-1}e_\theta} \pt_B)(\gamma(s)) \\
&\quad \quad + \frac{1}{t^2} \int_0^t \ds \, s^2 \bigl(\pt_B^{-1} (\covd_{r^{-1} e_\theta} F_B)(\dot{\gamma}(s), r^{-1}e_\theta) \pt_B\bigr)(\gamma(s)).
\end{equs}
The second term above is a perturbation, as we may estimate
\begin{equs}
\frac{1}{t^2} \int_0^t \ds \, s^2 \big| \bigl(\pt_B^{-1} (\covd_{r^{-1} e_\theta} F_B)(\dot{\gamma}(s), r^{-1}e_\theta) \pt_B\bigr) (\gamma(s))\big| &\lesssim \int_0^t \ds \bigl(\big\|\covd F_B \big\|_{L_{t, x}^\infty} |x| t^{-1} + \big\|\covd \ptl_t B \big\|_{L_{t, x}^\infty}\bigr) \\
&\lesssim \YMFnorm{B}_{t_0} \bigl(1 + (rt^{-1/2}) \bigr),
\end{equs}
which is acceptable. To handle the leading term, first, for notational brevity, let
\begin{equs}
R_1(s) &= \bigl(\pt_B^{-1} F_B(\dot{\gamma}(s), r^{-1} e_\theta) \pt_B\bigr)(\gamma(s)) -  \frac{r}{t}  F_B^{12}(0, 0), \\
R_2(s) &= (\pt_B^{-1} \covd_{r^{-1}e_\theta} \pt_B)(\gamma(s)) - \frac{s}{2t} r F_B^{12}(0, 0) .
\end{equs}
Note that (recalling \eqref{eq:curvature-radial-angular-calculation}) $F_B(\dot{\gamma}(0), r^{-1} e_\theta) = F_B^{12}(0, 0) r t^{-1} + (\ptl_t B) \cdot r^{-1} e_\theta$. Thus by \eqref{eq:pt-F-commutation-estimate}, we may estimate
\begin{equs}
\big|R_1(s) \big| &\lesssim s \bigl(\big\|\ptl_t^2 B \big\|_{L_{t, x}^\infty} + \bigl( \bigl\|\ptl_t F_B \bigr\|_{L_{t,x}^\infty} +  \big\|\covd \ptl_t B \big\|_{L_{t, x}^\infty}\bigr) |x| t^{-1} + \big\|\covd F_B \big\|_{L_{t, x}^\infty} |x|^2 t^{-2} \bigr) + \bigl\|\ptl_t B\bigr\|_{L_{t,x}^\infty} \\
&\lesssim t^{-1/2} \YMFnorm{B}_{t_0} \bigl(1 + (rt^{-1/2})^2\bigr).
\end{equs}
We also apply \eqref{eq:pt-angular-derivative-estimate-clean} to estimate (here, note that $|x|^{-1} e_\theta(x) = |(s/t) x|^{-1} e_\theta(\gamma(s))$, and thus we apply \eqref{eq:pt-angular-derivative-estimate-clean} with $(t, x)$ replaced by $\gamma(s) = (s, (s/t)x)$)
\begin{equ}
\big|R_2(s) \big| \lesssim s^{1/2} \YMFnorm{B}_{t_0} \bigl(1 + (rt^{-1/2})^2\bigr).
\end{equ}
Inserting these approximations, it follows that
\begin{equs}
~& \frac{1}{t^2}\bigg| \int_0^t \ds \, s^2 \Bigl(2 \bigl(\pt_B^{-1} F_B(\dot{\gamma}(s), r^{-1}e_\theta) \pt_B\bigr)(\gamma(s)) (\pt_B^{-1} \covd_{r^{-1}e_\theta} \pt_B)(\gamma(s)) - \frac{s}{t^2} |x|^2 (F_B^{12}(0, 0))^2\Bigr) \bigg| \\
&\lesssim \frac{1}{t^2}\int_0^t \ds \, s^2 | R_1(s) (\pt_B^{-1} \covd_{r^{-1} e_\theta} \pt_B)(\gamma(s))| + \frac{1}{t^2}\int_0^t \ds \, s^2 \frac{r}{t} \big| F_B^{12}(0, 0) R_2(s) \big| .
\end{equs}
We may combine our estimates on $R_1, R_2$ with the additional estimates
\begin{equs}
r \big|F_B^{12}(0, 0) \big| &\lesssim \big\|F_B \big\|_{L_{t, x}^\infty} r, \\
\big|(\pt_B^{-1} \covd_{r^{-1} e_\theta} \pt_B)(\gamma(s))\big| &\lesssim \bigl| (s/t)r F_B^{12}(0, 0) \bigr| + \bigl|R_2(s)\bigr| \lesssim \big\|F_B \big\|_{L_{t, x}^\infty} r + |R_2(s)|,
\end{equs}
in order to estimate
\begin{equs}
\int_0^t &\ds \, | R_1(s) (\pt_B^{-1} \covd_{r^{-1} e_\theta} \pt_B)(\gamma(s))| \\
&\lesssim t^{1/2} \YMFnorm{B}_{t_0} \bigl(1 + (rt^{-1/2})^2 \bigr) \bigl( \big\|F_B\big\|_{L_{t,x}^\infty} r + t^{1/2} \YMFnorm{B}_{t_0} \bigl(1 + (rt^{-1/2})^2\bigr)\bigr) \\
&\lesssim t \bigl(\big\|F_B \big\|_{L_{t, x}^\infty} + \YMFnorm{B}_{t_0}\bigr) \YMFnorm{B}_{t_0} \bigl(1 + (r t^{-1/2})^5\bigr),
\end{equs}
and similarly for the second term. To finish, we simply compute
\begin{equ}
\frac{1}{t^2} \int_0^t \ds  \frac{s^3}{t^2} |x|^2 (F_B^{12}(0, 0))^2 = \frac{r^2}{4} (F_B^{12}(0, 0))^2. 
\end{equ}
The desired estimate \eqref{eq:pt-hessian-angular-estimate-clean} now follows.

\proofstep{Proof of \eqref{eq:pt-hessian-radial-estimate-clean}.} We start as before, except now since $\ptl_{r^{-1} e_r} (r^{-1} e_r) = 0$, we in fact have that $\covd_{r^{-1} e_r} \covd_{r^{-1} e_r} \pt_B = (\Hess_B \pt_B)(r^{-1} e_r, r^{-1} e_r)$. We then apply Lemma~\ref{lemma:covd-parallel-transport} to obtain (here $\gamma(s) = s(1,  x/t)$)
\begin{equs}
(\pt_B^{-1} & \Hess_B \pt_B)(r^{-1} e_r, r^{-1} e_r) (t, x) \\
&= \frac{2}{t^2} \int_0^t \ds \, s^2 \bigl(\pt_B^{-1} F_B(\dot{\gamma}(s), r^{-1}e_r) \pt_B\bigr)(\gamma(s)) (\pt_B^{-1} \covd_{r^{-1}e_r} \pt_B)(\gamma(s)) \\
&\quad \quad + \frac{1}{t^2} \int_0^t \ds \, s^2 \bigl(\pt_B^{-1} (\covd_{r^{-1} e_r} F_B)(\dot{\gamma}(s), r^{-1}e_r) \pt_B\bigr)(\gamma(s)).
\end{equs}
The second term may be estimated as before:
\begin{equ}
\frac{1}{t^2}\int_0^t \ds \, s^2 \big| \bigl(\pt_B^{-1} (\covd_{r^{-1} e_r} F_B)(\dot{\gamma}(s), r^{-1}e_r) \pt_B\bigr)(\gamma(s))\big| \lesssim  \YMFnorm{B}_{t_0} \bigl(1 + (rt^{-1/2})\bigr).
\end{equ}
The first term is now also a perturbation. To see this, first note that
\begin{equ}
F_B(\dot{\gamma}(s), r^{-1} e_r) = F_B((1, x/t), r^{-1} e_r) = (\ptl_t B) \cdot (r^{-1} e_r).
\end{equ}
Combining this with \eqref{eq:pt-radial-derivative-estimate-clean}, we may estimate
\begin{equs}
\frac{1}{t^2}\bigg| \int_0^t \ds \, s^2 \bigl(\pt_B^{-1} F_B&(\dot{\gamma}(s), r^{-1}e_r) \pt_B\bigr)(\gamma(s)) (\pt_B^{-1} \covd_{r^{-1}e_r} \pt_B)(\gamma(s))\bigg| \\
&\lesssim \int_0^t \ds \, \big\|\ptl_t B \big\|_{L_{t, x}^\infty} s^{1/2} \YMFnorm{B}_{t_0} \bigl(1 + (rt^{-1/2})^2\bigr) \\
&\lesssim t \bigl(\YMFnorm{B}_{t_0}\bigr)^2 \bigl(1 + (rt^{-1/2})^2\bigr),
\end{equs}
which is acceptable. Combining, we may thus obtain \eqref{eq:pt-hessian-radial-estimate-clean}. 

\proofstep{Proof of \eqref{eq:pt-laplacian-estimate-clean}.} We combine \eqref{eq:pt-hessian-angular-estimate-clean} and \eqref{eq:pt-hessian-radial-estimate-clean}.
\end{proof}

\begin{remark}
In the estimates of Lemmas~\ref{lemma:pt-curvature-commutator},~\ref{lemma:F-expansion-computation-errors}, and~\ref{lemma:YM-perturb-main-intermediate-estimate-clean}, we did not assume that $B$ is a solution to the Yang--Mills gradient flow. The only place where we use this assumption in the proof of Proposition~\ref{prop:qB-perturb-near-YM} will be in applying the heat kernel estimates \eqref{eq:qB-ptwise-bounds} of $\fq_B$, as we will shortly see.
\end{remark}

With the estimates from Lemma~\ref{lemma:YM-perturb-main-intermediate-estimate-clean} in hand, we may now prove Proposition~\ref{prop:qB-perturb-near-YM}.

\begin{proof}[Proof of Proposition~\ref{prop:qB-perturb-near-YM}]
Recall that $\fr_B^F$ is precisely the quantity to be estimated in Proposition~\ref{prop:qB-perturb-near-YM}. Towards this end, we begin by slightly rewriting its equation \eqref{eq:rB-F-equation}. We distribute
\begin{equ}
\ptl_{r^{-1} e_\theta} (\pt_B^{-1} \fq_B) = \pt_B^{-1} \covd_{r^{-1}e_\theta} \fq_B - \pt_B^{-1} (\covd_{r^{-1} e_\theta} \pt_B) \pt_B^{-1} \fq_B,
\end{equ}
and similarly for $r^{-1} e_r$. Inserting this into \eqref{eq:rB-F-equation}, we further obtain 
\begin{equ}
(\ptl_t - \mc{L}_X^{F,\theta}) \fr_B^F = R_1 (\ptcutoff \pt_B^{-1} \fq_B) + 2 \ptcutoff R_2\pt_B^{-1} \covd_{r^{-1} e_\theta} \fq_B + 2 \ptcutoff R_3 \pt_B^{-1} \covd_{r^{-1} e_r} \fq_B + R_4,
\end{equ}
where (recall $X = F_B^{12}(0, 0)$)
\begin{equs}
R_1 = \,&2 (\pt_B^{-1} F_B \pt_B - \maglapF{X}) + \Bigl(\pt_B^{-1} \Delta_B \pt_B - \frac{r^2}{4} X^2\Bigr) \\
&- 2 \Bigl(\pt_B^{-1} \covd_{r^{-1}e_\theta} \pt_B - \frac{r}{2}X\Bigr) \pt_B^{-1} \covd_{r^{-1} e_\theta} \pt_B \\
&- \pt_B^{-1} \ptl_t \pt_B - 2 (\pt_B^{-1} \covd_{r^{-1} e_r} \pt_B) (\pt_B^{-1} \covd_{r^{-1} e_r} \pt_B), \\
R_2 = \,&   \Big(\pt_B^{-1} \covd_{r^{-1} e_\theta} \pt_B - \frac{r}{2} X\Big), \\
R_3 = \,&   \pt_B^{-1} \covd_{r^{-1} e_r} \pt_B \\
R_4 = \,& - r X (\ptl_{r^{-1} e_\theta} \ptcutoff) (\pt_B^{-1} \fq_B - \fs_X^F) - 2(\ptl^\ell \ptcutoff) \ptl_\ell(\pt_B^{-1} \fq_B - \fs_X^F) - (\Delta \ptcutoff)(\pt_B^{-1} \fq_B - \fs_X^F).
\end{equs}
By Lemma~\ref{lemma:YM-perturb-main-intermediate-estimate-clean}, we have the following pointwise bounds on $R_1, R_2, R_3$ on the support of $\ptcutoff$ (i.e.\ for $x$ such that $|x| \leq 1/10$):
\begin{equs}
|R_1(t, x)| &\lesssim \bigl(1 + t\big\|F_B \big\|_{L_{t, x}^\infty} + t\YMFnorm{B}_{t_0}\bigr) \YMFnorm{B}_{t_0} \bigl( 1 + (rt^{-1/2})^5 \bigr),\\
|R_2(t, x)| + |R_3(t, x)| & \lesssim t^{1/2} \YMFnorm{B}_{t_0} \bigl(1 + (rt^{-1/2})^2\bigr).
\end{equs}
To estimate $R_4$, note that it is composed of all the cases where derivatives fall on the cutoff function $\ptcutoff$, and thus $R_4$ is only nonzero for values of $x$ with $|x| \gtrsim 1$. Combining this with our pointwise estimates of $\fq_B$ (Proposition~\ref{prop:qB-ptwise-bds}) and our exact formula for $\fs_X^F$ discussed at the beginning of this section, we may estimate
\begin{equ}
|R_4(t, x)| \lesssim t^{10} \bigl(1 + \bigl\|F\bigr\|_{L_{t,x}^\infty([0, t_0])} + \YMFnorm{B}_{t_0} \bigr) \exp\bigl(C t\bigl\|F\bigr\|_{L_{t,x}^\infty([0, t_0])}\bigr) t^{-1} \exp(-c |x|^2 / t).
\end{equ}
Here, the power $10$ could have been any number. The point is that $R_4$ is much smaller than anything else, and thus will result in a (much more than) acceptable contribution.

From the defining equation \eqref{eq:fs-B-F-theta} for $\fs_X^{F, \theta}$, we have that $\fs_X^{F, \theta}$ is the integral kernel of $(\ptl_t - \mc{L}_X^{F, \theta})^{-1}$. Thus upon integrating and combining the previous pointwise bounds with our pointwise bounds on $\fs_X^{F, \theta}$ (recall \eqref{eq:mehler-phase-estimate}) as well as $\fq_B$ and $\covd \fq_B$ from Proposition~\ref{prop:qB-ptwise-bds}, we have that for all $(t_1, x_1) \in (0, t_0] \times \T^2$,
\begin{equs}
\big|\fr_B^F(t_1, x_1)\big| &\lesssim \int \dt \dx \, \big| \fs_X^{F, \theta}(t, x; t_1, x_1)\big| \Big( \big|R_1(t, x)\big| |\fq_B(t, x)\big| \\
&\quad \quad + \bigl(|R_2(t, x)| + |R_3(t, x)|\bigr) |\covd \fq_B(t, x)\big| + |R_4(t, x)| \Big) \\
&\lesssim \exp(C_0 t_0 \|F_B\|_{L_{t,x}^\infty}) \int \dt \dx \,  \fp(t, x; t_1, x_1) t^{-1} \exp(-c |x|^2 / t) \\
&\times \Bigl(t^{10} \bigl(1 + \YMFnorm{B}_{t_0}\bigr) + \bigl(1 + t\big\|F_B \big\|_{L_{t, x}^\infty} + t\YMFnorm{B}_{t_0}\bigr) \YMFnorm{B}_{t_0} \bigl( 1 + (rt^{-1/2})^5 \bigr)\Bigr) \\
&\lesssim \bigl(c_0 + t_1^{10}\bigr)\exp(C_0 t_0 \|F_B\|_{L_{t,x}^\infty}) t_1^{-1} \exp(- c |x_1|^2 / t_1),
\end{equs}
as desired (in the last step we applied the smallness assumption \eqref{e:qB-perturb-near-YM_assumption}).
\end{proof}

\subsection{Schauder estimates}\label{section:schauder-and-commutator-estimates}

In this section, we collect several Schauder estimates which will be needed when carrying out perturbative arguments in the coming sections. These results will be obtained as consequences of the material in Sections~\ref{subsection:covariant-heat-kernel-estimates}--\ref{section:mixed-derivative-heat-kernel-estimates}. We begin by setting notation regarding covariant Duhamel operators. 

\begin{definition}\label{def:covariant-duhamel}
For $H \colon [0, T] \times \T^d \ra \genV$, we define
\begin{equ}
\Duh_B(H)(t, x) \eqdef \int_0^t \ds \int_{\T^d} \dy\, \fp_B(s, y; t, x) H(s, y), \quad (t, x) \in [0, T] \times \T^d.
\end{equ}
Similarly, for $G \colon [0, T] \times \T^d \ra \genV^{d}$, we define
\begin{equ}
\Duh_B^F(G)(t, x) \eqdef \int_0^t \ds \int_{\T^d} \dy \, \fq_B(s, y; t, x) G(s, y), \quad (t, x) \in [0, T] \times \T^d. 
\end{equ}
Note that $\Duh_B(H), \Duh_B^F(G)$ satisfy
\begin{equs}[2]
(\ptl_t - \Delta_B) \Duh_B (H) &= H, \qquad \Duh_B(H)(0) &= 0, \\
(\ptl_t - \Delta_B - 2F) \Duh_B^F(G) &= G, \qquad \Duh_B^F(G)(0) &= 0.
\end{equs}
\end{definition}


We proceed with the main results of this section.

\begin{lemma}[Covariant Schauder estimates]\label{lemma:covariant-schauder-estimates}
Let $C_0$ be the maximum of the constants from Lemma~\ref{lemma:higher-regularity-abstract-small-Lt1-Lxinfty} and Proposition~\ref{prop:qB-ptwise-bds}. Let $T \leq 1$, $G \colon [0, T] \times \T^d \ra \genV^{d}$. For any $1 \leq p \leq \infty$, we have that
\begin{equs}
\bigl\|\Duh_{B}^F(G)\bigr\|_{L_t^\infty L_x^p([0, T])} &\lesssim \gronwall_T(C_0 F) \bigl\| G \bigr\|_{L_t^1 L_x^p([0, T])}, \\
\bigl\|\covd_B\Duh_{B}^F(G)\bigr\|_{L_t^\infty L_x^p([0, T])} &\lesssim \gronwall_T(C_0 F)  \bigl(1 + \Fnorm{F}_{T, 1}\bigr) \sup_{t_0 \in (0, T]} \bigl \| (t_0 - t)^{-1/2} G \bigr\|_{L_t^1 L_x^p([0, t_0])}, \\
\bigl\|\Duh_{B}^F(\covd_B G)\bigr\|_{L_t^\infty L_x^p([0, T])} &\lesssim \gronwall_T(C_0 F) \bigl(1 + \Fnorm{F}_{T, 1}\bigr) \sup_{t_0 \in (0, T]} \bigl\|\min(t, (t_0 - t))^{-1/2} G\bigr\|_{L_t^1 L_x^p([0, t_0])}.
\end{equs}
Similarly, for $H \colon [0, T] \times \T^d \ra \genV$, we have that
\begin{equs}
\bigl\|\Duh_{B}(H)\bigr\|_{L_t^\infty L_x^p([0, T])} &\lesssim \bigl\| H \bigr\|_{L_t^1 L_x^p([0, T])}, \\
\bigl\|\covd_B\Duh_{B}(H)\bigr\|_{L_t^\infty L_x^p([0, T])} &\lesssim  \bigl(1 + \Fnorm{F}_{T, 0}\bigr) \sup_{t_0 \in (0, T]} \bigl \| (t_0 - t)^{-1/2} H \bigr\|_{L_t^1 L_x^p([0, t_0])}, \\
\bigl\|\Duh_B(\covd_B H)\bigr\|_{L_t^\infty L_x^p([0, T])} &\lesssim  \bigl(1 + \Fnorm{F}_{T, 1}\bigr)\sup_{t_0 \in (0, T]} \bigl\|\min(t, (t_0 - t))^{-1/2} H\bigr\|_{L_t^1 L_x^p([0, t_0])}
\end{equs}
\end{lemma}
\begin{proof}
We focus on the estimates for $\Duh_B^F$, as the argument for $\Duh_B$ is essentially the same (with Proposition~\ref{prop:pB-ptwise-bounds} applied in place of Proposition~\ref{prop:qB-ptwise-bds}). Let $t_0 \in [0, T]$. By Proposition~\ref{prop:qB-ptwise-bds} (in particular \eqref{eq:qB-ptwise-bounds} with $k = 0$) and Young's convolution inequality, we have that
\begin{equs}
\bigl\|\Duh_{B}^F(G)(t_0)\bigr\|_{L_x^p} &\leq \int \ds \, \biggl\|\int \dy \, \fq_B(s, y; t_0, x) G(s, y) \biggr\|_{L_x^p} \\
&\lesssim \gronwall_T(C_0 F)\int \ds \,  \biggl\| \int \dy \, (t_0 - s)^{-\tfrac{d}{2}} \exp\bigg(-c\frac{|x - y|^2}{t_0 - s}\bigg) |G(s, y)| \biggr\|_{L_x^p} \\
&\lesssim \gronwall_T(C_0 F) \int \ds \, \bigl\|G(s)\bigr\|_{L_x^p} \lesssim \gronwall_T(C_0 F) \bigl\|G\bigr\|_{L_t^1 L_x^p([0, T])},
\end{equs}
which shows the first estimate. The second estimate may be proven similarly, where we apply \eqref{eq:qB-ptwise-bounds} with $k = 1$. The third estimate also follows similarly, where we additionally apply \eqref{eq:dagger-fq} to first move the derivative off of $G$ and onto $\fq_B$, and then apply Proposition~\ref{prop:qB-mixed-derivatives-ptwise-bounds} in place of the estimate \eqref{eq:qB-ptwise-bounds}.
\end{proof}

\begin{lemma}\label{lemma:higher-regularity-covd-ibp}
Let $C_0$ be the constant from Proposition~\ref{prop:qB-ptwise-bds}. Let $G \colon [0, T] \times \T^d \ra \genV^d$.
For any $k \geq 0$, we have that
\begin{equ}
\big\|t^{k/2} \covd^{(k)} \Duh_B^F(G) \big\|_{L_{t, x}^\infty([0, T])} \lesssim \gronwall_T(C_0 F) \bigl(1 + (\Fnorm{F}_{T, k})^k\bigr) N_k(G),
\end{equ}
where 
\begin{equs}
N_k(G)& = \big\|G\big\|_{L_t^1 L_x^\infty([0, T])} \\
&+ \sum_{j=0}^{k-1} \Bigl( \big\|t^{(j+1)/2} \covd^{(j)} G \big\|_{L_t^2 L_x^\infty([0, T])} + \sup_{u \in [0, T]} \big\|t^{(j+1)/2} (u - t)^{-1/2} \covd^{(j)} G \big\|_{L_t^1 L_x^\infty([0, u])}\Bigr).
\end{equs}
In particular, when $k = 0$, the sum is empty, and $N_0(G) = \big\|G\big\|_{L_t^1 L_x^\infty([0, T])}$. Similarly, for $H \colon [0, T] \times \T^d \ra \genV$, we have that for any $k \geq 0$,
\begin{equ}
\bigl\|t^{k/2} \covd^{(k)} \Duh_B(H)\bigr\|_{L_{t,x}^\infty([0, T])} \lesssim \bigl(1 + \bigl(\Fnorm{F}_{T, {k-1}}\bigr)^k \bigr) N_k(H).
\end{equ}
\end{lemma}
\begin{proof}
We focus on the estimate of $\Duh_B^F(G)$, as the estimate of $\Duh_B(H)$ proceeds almost exactly the same (one just replaces the application of Corollary~\ref{cor:higher-regularity-with-curvature} by Proposition~\ref{prop:covariant-heat-equation-higher-regularity} in the base case of the ensuing inductive argument). Let $\phi = \Duh_B^F(G)$. Let $(t_0, x_0) \in [0, T] \times \T^d$, and let $K(t, x) = \fp(t, x; t_0, x_0)$. In the following, all space-time norms will be taken on $[0, t_0] \times \T^d$, unless otherwise specified. By Lemma~\ref{lemma:heat-kernel-weighted-to-L-infty}, it suffices to show that for all $k \geq 0$, we have that
\begin{equ}\label{eq:higher-regularity-covd-ibp-inductive-assumption}
\big\|t^{k/2} K^{1/2} \covd^{(k)} \phi\big\|_{L_t^\infty L_x^2} + \big\|t^{k/2} K^{1/2} \covd^{(k+1)} \phi\big\|_{L_{t, x}^2} \lesssim \gronwall_T(C_0 F) \bigl(1 + (\Fnorm{F}_{T, k})^k\bigr) N_k(G).
\end{equ}
We proceed by induction. For the base case $k = 0$, by Corollary~\ref{cor:higher-regularity-with-curvature}, we have that
\begin{equ}\label{eq:higher-regularity-covd-ibp-intermediate}
\big\|K^{1/2} \phi\big\|_{L_t^\infty L_x^2} + \big\|K^{1/2} \covd \phi\big\|_{L_{t, x}^2} \lesssim \gronwall_T(C_0 F) \big\|K^{1/2} G \big\|_{L_t^1 L_x^2} \lesssim \gronwall_T(C_0 F)\big\|G\big\|_{L_t^1 L_x^\infty},
\end{equ}
where in the last inequality, we used Hölder and the fact that $\big\|K^{1/2} \big\|_{L_t^\infty L_x^2} = 1$.

Next, suppose that for some $k_0 \geq 0$, the estimate \eqref{eq:higher-regularity-covd-ibp-inductive-assumption} has been proven for all $0 \leq k \leq k_0$. We proceed to prove the estimate for $k_0 + 1$. Towards this end, first note by Lemma~\ref{lemma:equations-for-higher-derivatives} and Remark~\ref{remark:higher-derivatives-equation-with-ad-F-term} that $\covd^{(k_0+1)} \phi$ satisfies the following schematic equation:
\begin{equs}
(\ptl_t - \Delta_B) \covd^{(k_0 + 1)} \phi &= \sum_{\substack{j_1 + j_2 = k_0 + 1 \\ j_1, j_2 \geq 0}} \schO\bigl((\covd^{(j_1)} F)( \covd^{(j_2)} \phi)\bigr) + \covd^{(k_0 + 1)} G \\
&= \schO(F \covd^{(k_0 + 1)} \phi) + R + \covd^{(k_0 + 1)} G,
\end{equs}
where
\begin{equ}
R :=  \sum_{\substack{j_1 + j_2 = k_0 + 1 \\ j_1 \geq 1, j_2 \geq 0}} \schO\bigl((\covd^{(j_1)} F) (\covd^{(j_2)} \phi) \bigr).
\end{equ}
By the time-weighted estimate \eqref{eq:higher-regularity-time-weighted-abstract} with $\alpha = k_0 + 1$ (and the $G, H$ therein taken to respectively be $R$ and $\covd^{(k_0)} G$), we have that
\begin{equs}
\big\|&t^{(k_0 + 1)/2} K^{1/2} \covd^{(k_0 + 1)} \phi\big\|_{L_t^\infty L_x^2} + \big\|t^{(k_0 + 1)/2} K^{1/2} \covd^{(k_0 + 2)} \phi\big\|_{L_{t, x}^2} \\
&\lesssim \Big(1 + \Fnorm{F}_{t_0, 0}\Big) \big\|t^{k_0/2} K^{1/2} \covd^{(k_0 + 1)} \phi\big\|_{L_{t, x}^2} \\
&\quad + \big\|t^{(k_0 + 1)/2} K^{1/2} R \big\|_{L_t^1 L_x^2} + \big\|t^{(k_0 + 1)/2} K^{1/2} \covd^{(k_0)} G \big\|_{L_{t, x}^2} \\
&\quad + \big\|t^{(k_0 + 1)/2} (t_0 - t)^{-1/2} K^{1/2}  \covd^{(k_0)} G \big\|_{L_t^1 L_x^2} + \big\|t^{(k_0 + 1)/2} |\Delta K|^{1/2} \covd^{(k_0)} G \big\|_{L_t^1 L_x^2}.
\end{equs}
By our inductive assumption, the first term on the right-hand side is acceptable. To estimate the last three terms, we apply Hölder to obtain
\begin{equs}
\big\|t^{(k_0 + 1)/2} K^{1/2} \covd^{(k_0)} G \big\|_{L_{t, x}^2} &\leq \big\|K^{1/2} \big\|_{L_t^\infty L_x^2} \big\|t^{(k_0 + 1)/2} \covd^{(k_0)} G \big\|_{L_t^2 L_x^\infty} \\
&\leq \big\|t^{(k_0 + 1)/2} \covd^{(k_0)} G \big\|_{L_t^2 L_x^\infty}, \\
\big\|t^{(k_0 + 1)/2} K^{1/2} (t_0 - t)^{-1/2} \covd^{(k_0)} G \big\|_{L_t^1 L_x^2} &\leq \big\|t^{(k_0 + 1)/2} (t_0 - t)^{-1/2} \covd^{(k_0)} G \big\|_{L_t^1 L_x^\infty} \\
\big\|t^{(k_0 + 1)/2} |\Delta K|^{1/2} \covd^{(k_0)} G \big\|_{L_t^1 L_x^2} &\leq \big\|(t_0 - t)^{1/2} |\Delta K|^{1/2} \big\|_{L_t^\infty L_x^2} \big\| t^{(k_0 + 1)/2} (t_0 - t)^{-1/2} \covd^{(k_0)} G \big\|_{L_t^1 L_x^\infty} \\
&\lesssim \big\| t^{(k_0 + 1)/2} (t_0 - t)^{-1/2} \covd^{(k_0)} G \big\|_{L_t^1 L_x^\infty}.
\end{equs}
By the definition of $N_k(G)$, each of the terms above is acceptable. 
To handle the term involving $R$, recall that $R$ is a sum of terms of the form $\schO((\covd^{(j_1)} F) (\covd^{(j_2)} \phi))$, where $j_1 + j_2 = k_0 + 1$, $j_1 \geq 1$, $j_2 \geq 0$. Note in particular this implies that $j_2 \leq k_0$. Thus for each such term, we may estimate by Hölder and our inductive assumption
\begin{equs}
\big\|t^{(k_0 + 1)/2} K^{1/2} (\covd^{(j_1)} F)(\covd^{(j_2)} \phi)\big\|_{L_t^1 L_x^2} &\leq \big\|t^{j_1/2} \covd^{(j_1)} F \big\|_{L_t^1 L_x^\infty} \big\|t^{j_2/2} K^{1/2} \covd^{(j_2)} \phi\big\|_{L_t^\infty L_x^2} \\
&\lesssim \gronwall_T(C_0 F)\Fnorm{F}_{t_0, k_0 + 1} \bigl(1 + (\Fnorm{F}_{T, k_0})^{k_0}\bigr) N_{k_0}(G),
\end{equs}
which is acceptable. 
\end{proof}

\section{Analysis of deterministic flows}\label{section:deterministic-rough}

In this section, we collect many of the deterministic estimates that will be needed later on. In Sections~\ref{sec:rough_YMHF} and~\ref{sec:higgsdet_estimates} we estimate the various deterministic pieces that will appear in our decomposition of the solution $(A, \higgs)$ to \eqref{intro:eq-SYMH} (recall that various aspects of the decomposition were already informally discussed in Section~\ref{section:intro}). Many of these estimates involve analysing (small) rough perturbations of (large) smooth backgrounds. They will be needed in both the local theory (see e.g.\ Section~\ref{section:sdet-higgsdet-sol}) as well as in the global theory. For the latter, in Section~\ref{section:global-iteration-monotone-quantities-deterministic}, we consider a single step of our global iteration scheme and treat all terms which may be estimated deterministically (by applying the estimates from Sections~\ref{sec:rough_YMHF} and~\ref{sec:higgsdet_estimates}).


We begin with the following lemma, which estimates solutions to the Yang--Mills gradient flow on an interval $[0, T]$ by a combination of gauge covariant quantities (given in terms of the curvature and its covariant derivatives) along with the size of the connection at the end time $T$. In particular, the result shows that solutions to the Yang--Mills gradient flow satisfy an approximate Coulomb gauge condition, in that the divergence is (almost) one derivative better than expected. This will be needed soon in Section~\ref{sec:rough_YMHF} as well as much later in Section~\ref{sec:bounds_from_covariant_norms}. We also remark that the final estimate \eqref{eq:ymgf-der-div-B-estimate} is quite crude, and will only be needed in Section~\ref{sec:bounds_from_covariant_norms} to establish a soft statement, thus the reader may choose to skip its proof on a first reading.

\begin{lemma}\label{lemma:ymgf-approximate-coulomb-gauge}
Let $0 < T \leq 1$, let $B \in C_t^0 \cC_x^{-\kappa}([0, T] \times \T^2) \cap C_{t,x}^\infty((0, T] \times \T^2)$ be a classical solution to the Yang--Mills gradient flow at positive times, and let
\begin{equ}
M_0 \eqdef \bigl\|t^{\frac{1+\kappa}{2}} F_B\bigr\|_{L_{t,x}^\infty([0, T])},\qquad
M_1 \eqdef \bigl\|t^{1+\frac{\kappa}{2}} \covd_B F_B \bigr\|_{L_{t,x}^\infty([0, T])}.
\end{equ}
(These quantities may of course be infinite, in which case some of 
the following estimates are trivially true.) Then we have the bounds
\begin{equs}
\bigl\| t^{\frac{\kappa}{2}} B\bigr\|_{L_{t,x}^\infty([0, T])} &\lesssim \bigl\|B(T)\bigr\|_{L_x^\infty} + M_1,  \label{eq:ymgf-B-Linfty-estimate}\\
\bigl\|B\bigr\|_{C_t^0 \cC_x^{-\kappa}([0, T])} & \lesssim \bigl(M_0 + 1\bigr) \bigl(\bigl\|B(T)\bigr\|_{L_x^\infty} + M_1 + 1\bigr) , \label{eq:ymgf-B-negative-regularity-estimate}\\
\bigl\|t^{\kappa} \ptl^\ell B_\ell \bigr\|_{L_{t,x}^\infty([0, T])} &\lesssim \bigl\|\ptl^\ell B_\ell(T)\bigr\|_{L_x^\infty} + \bigl(\bigl\|B(T)\bigr\|_{L_x^\infty} + M_1\bigr) M_1, \label{eq:ymgf-div-B-estimate} \\
\bigl\|t^{\frac{1}{2} + \kappa} \ptl (\ptl^\ell B_\ell) \bigr\|_{L_{t,x}^\infty([0, T])} &\lesssim \bigl\|\ptl(\ptl^\ell B_\ell)(T)\bigr\|_{L_x^\infty} +\bigl\|\ptl B(T)\bigr\|_{L_x^\infty} \bigl(M_1 + 1 \bigr)\\
 &\quad + \bigl(\bigl\|B(T)\bigr\|_{L_x^\infty} + M_1 + 1\bigr)^2 (M_0 + M_1 + 1)^5 .\label{eq:ymgf-der-div-B-estimate}
\end{equs}
\end{lemma}

\begin{proof}
Since $B$ is a solution to the Yang--Mills gradient flow \eqref{e:YMHF}, we have that for $t_0 \in (0, T]$ and $i \in [2]$,
\begin{equ}\label{eq:B-from-covd-FB}
B_i(t_0) = B_i(T) - \int_{t_0}^T \dt \, \covd_B^\ell (F_B)_{\ell i}(t). 
\end{equ}
The estimate \eqref{eq:ymgf-B-Linfty-estimate} now directly follows from the definition of $M_1$. 

For the estimate \eqref{eq:ymgf-B-negative-regularity-estimate}, let $\phi_x^\lambda$ be a ($\frkg$-valued) $L_x^1$-normalised test function localised to some point $x \in \T^2$ at scale $\lambda \ll 1$. Since $B \in C_t^0 \cC_x^{-\kappa}([0, T])$, it suffices to show that (uniformly in $x \in \T^2$ and $\lambda \ll 1$)
\begin{equ}\label{eq:ymgf-lower-regularity-intermediate}
\max_{i \in [2]} \sup_{t_0 \in (0, T]} \bigl| \langle B_i(t_0), \phi_x^\lambda \rangle\bigr| \lesssim \lambda^{-\kappa} \bigl(M_0 + 1\bigr) \bigl(\bigl\|B(T)\bigr\|_{L_x^\infty} + M_1 + 1\bigr).
\end{equ}
Towards this end, if $t_0 \in (0, T]$ is such that $t_0 \geq \lambda^2$, then by \eqref{eq:ymgf-B-Linfty-estimate} we have that
\begin{equ}
\bigl|\langle B_i(t_0), \phi_x^\lambda \rangle \bigr| \lesssim t_0^{-\frac{\kappa}{2}} \bigl(\bigl\|B(T)\bigr\|_{L_x^\infty} + M_1\bigr) \lesssim \lambda^{-\kappa} \bigl(\bigl\|B(T)\bigr\|_{L_x^\infty} + M_1\bigr),
\end{equ}
which is acceptable. Therefore, we may assume that $t_0 \in (0, T]$ is such that $t_0 \leq \lambda^2$. In this case, integrating the equation $\ptl_t B_i = \covd_B^\ell (F_B)_{\ell i}$ between $t_0$ and $\lambda^2 \wedge T$, we obtain
\begin{equ}
\langle B_i(t_0), \phi_x^\lambda \rangle = \langle B_i(\lambda^2 \wedge T), \phi_x^\lambda \rangle - \int_{t_0}^{\lambda^2 \wedge T} \dt \langle \covd_B^\ell (F_B)_{\ell i}(t), \phi_x^\lambda \rangle. 
\end{equ}
If $\lambda^2 \geq T$, then the first term on the right may be directly estimated by $\|B(T)\|_{L_x^\infty}$, while if $\lambda^2 < T$, the term may be estimated as in the previous case. For the second term, we expand $\covd_B^\ell = \ptl^\ell + \ad(B^\ell)$, and then integrate the derivative term by parts to obtain
\begin{equ}
\int_{t_0}^{\lambda^2 \wedge T} \dt\, \langle \covd_B^\ell (F_B)_{\ell i}(t), \phi_x^\lambda \rangle =- \int_{t_0}^{\lambda^2 \wedge T} \dt \, \langle (F_B)_{\ell i}(t), \ptl^\ell \phi_x^\lambda \rangle + \int_{t_0}^{\lambda^2 \wedge T} \dt\, \langle [B^\ell, (F_B)_{\ell i}], \phi_x^\lambda \rangle.
\end{equ}
Combining the definition of $M_0$ with the fact that $\bigl\|\ptl^\ell \phi_x^\lambda\bigr\|_{L_x^1} \lesssim \lambda^{-1}$ and the estimate \eqref{eq:ymgf-B-Linfty-estimate}, we may estimate both terms on the right-hand side above by acceptable contributions. Combining, we thus obtain \eqref{eq:ymgf-lower-regularity-intermediate}.

For the estimate \eqref{eq:ymgf-div-B-estimate}, recall from Lemma~\ref{lemma:YM-gradient-flow-covariant-divergence-free} that $\ptl_t B$ is covariant divergence free, i.e.\ $\covd_B^i \ptl_t B_i = 0$. Since $\ptl_t B_i = \covd_B^\ell (F_B)_{\ell i}$, and $\covd_B^i = \ptl^i + [B^i, \place]$, we may rewrite this as
\begin{equ}\label{eq:div-B-evolution}
\ptl_t \ptl^i B_i = - [B^i, \covd_B^\ell (F_B)_{\ell i}]. 
\end{equ}
The estimate \eqref{eq:ymgf-div-B-estimate} now follows by integrating this identity back from the end time $T$, and applying the first estimate on $B$ along with the definition of $M_1$ (and directly estimating the boundary term by $\|\ptl^\ell B_\ell(T)\|_{L_x^\infty}$).

We move on to the final estimate \eqref{eq:ymgf-der-div-B-estimate}. We will be quite wasteful in the following bounds (as this will be sufficient for our purposes). To start, we hit both sides of \eqref{eq:div-B-evolution} by a derivative to obtain
\begin{equation}\label{eq:der-div-evolution}
\begin{aligned}
\ptl_t \bigl(\ptl (\ptl^i B_i)\bigr) &= - \ptl [B^i, \covd_B^\ell (F_B)_{\ell i}] \\
&= - [B^i, \covd_B \covd_B^\ell (F_B)_{\ell i}] + [B^i, [B, \covd_B^\ell (F_B)_{\ell i}]] - [\ptl B^i, \covd_B^\ell (F_B)_{\ell i}].
\end{aligned}
\end{equation}
We thus see that we need additional estimates on $\covd_B^{(2)} F_B$ and $\ptl B$. For the former, first observe that by the definition of $M_0, M_1$ and $\Fnorm{F_B}_{T, 1}$ (recall Definition~\ref{def:Fnorm}), we may estimate $\Fnorm{F_B}_{T, 1} \lesssim M_0 + M_1$. Since $F_B$ solves a covariant heat equation (recall Lemma~\ref{lemma:curvature-equation}), we obtain by our heat kernel estimate \eqref{eq:higher-order-pB-ptwise-bd} that
\begin{equ}\label{eq:ymgf-intermediate-M2-estimate}
\bigl\|t^{\frac{3+\kappa}{2}} \covd_B^{(2)} F_B \bigr\|_{L_{t,x}^\infty([0, T])} \lesssim (1 + M_0 + M_1)^2 M_0. 
\end{equ}
For an estimate on $\ptl B$, we hit both sides of \eqref{eq:B-from-covd-FB} by a derivative to obtain
\begin{equ}
\ptl B_i(t_0) = \ptl B_i(T) - \int_{t_0}^T \dt \, \bigl( \covd_B \covd_B^\ell (F_B)_{\ell i} (t) - [B, \covd_B^\ell (F_B)_{\ell i}]\bigr),
\end{equ}
and thus (applying \eqref{eq:ymgf-intermediate-M2-estimate} to estimate the first term in the integral, and \eqref{eq:ymgf-B-Linfty-estimate} to estimate the second)
\begin{equ}
t_0^{\frac{1+\kappa}{2}}\bigl\|\ptl B (t_0)\bigr\|_{L_x^\infty} \lesssim \bigl\|\ptl B(T)\bigr\|_{L_x^\infty} + (1 + M_0 + M_1)^2 M_0 + \bigl(\bigl\|B(T)\bigr\|_{L_x^\infty} + M_1\bigr) M_1.
\end{equ}
Combining our estimates with the evolution \eqref{eq:der-div-evolution} for $\ptl (\ptl^\ell B_\ell)$, we may obtain \eqref{eq:ymgf-der-div-B-estimate}.
\end{proof}

Next, we define gauge invariant Hölder norms which will appear in the various well-posedness results of the coming sections. First, we set some notation which will be used in the definition, as well as in many places later on.

\begin{notation}\label{notation:linear-propagation}
	Let $A \in \cC_x^{-\kappa}$ be a spatial connection one-form, and let $T > 0$. Suppose that there is a unique Yang--Mills gradient flow $B \in C_{t,x}^\infty((0, T] \times \T^2)$ such that $B(t) \ra A$ in $\cC_x^{-\kappa}$ as $t \downarrow 0$. Given a (possibly distributional) $\genV$-valued one-form $u_0$, suppose that there is a unique solution $u$ to the following equation:
	\begin{equ}
		(\partial_t - \Delta_{B} - 2 F_{B}) \, u = 0 \;, \quad u(0) = u_0 \;.
	\end{equ}
    Then for $s \in (0, T]$, we let $e^{s\Delta^F_{A}} u_0 \eqdef u(s)$. 
	
    Similarly, given a $\genV$-valued distribution $v_0$, suppose that there is a unique solution to the following equation:
	\begin{equ}
		(\partial_t - \Delta_{B}) \, v = 0 \;, \quad v(0) = v_0 \;.
	\end{equ}
    Then for $s \in (0, T]$, we let $e^{s\Delta_{A}} v_0 \eqdef v(s)$. 
\end{notation}
\begin{remark}
    If $A$ is smooth, then $B$ exists uniquely and is smooth on $[0, \infty) \times \T^2$, and by standard linear parabolic theory, for any $\alpha \in \R$, $u_0, v_0 \in \cC_x^{\alpha}$, $s \geq 0$, the terms $e^{s \Delta_A^F} u_0, e^{s \Delta_A} v_0$ are well-defined. If, on the other hand, $A \in \cC_x^{-\kappa}$ is genuinely distributional, then one needs to be a bit more careful. In this paper, it will always be the case that for such $A$, the Yang--Mills gradient flow $B$ is such that the quantities $M_0, M_1$ from Lemma~\ref{lemma:ymgf-approximate-coulomb-gauge} are finite, and thus in particular by that lemma, $B$ satisfies the time-weighted $L_{t,x}^\infty$ estimate \eqref{eq:ymgf-B-Linfty-estimate}. Using this, one may show using standard parabolic theory that if $u_0, v_0 \in \cC_x^{-\kappa}$ (this is certainly not sharp, but it is sufficient for our purposes), then the terms $e^{s \Delta_A^F} u_0, e^{s \Delta_A} v_0$ are well-defined.
\end{remark}

\begin{definition}[Gauge invariant Hölder norms]
\label{defn:covariant_Holder}
    Let $A, B$ be as in Notation~\ref{notation:linear-propagation}. Let $\alpha < 1$, $T > 0$. Given a (possibly distributional) $\genV$-valued one-form $u_0$ such that $e^{s \Delta_A^F} u_0$ is well-defined for all $s \in (0, T]$, we let
	\begin{equ}
		\| u_0 \|_{\covHolder{A, F, T}^{\alpha}} \eqdef \max_{0 \leq m \leq 1 } \sup_{s \in (0, T] } \big( s^{\frac{-\alpha + m}{2}} \wedge 1 \big) \|\covd^{(m)}_{B(s)} \, e^{s\Delta^F_{A}} u_0 \|_{L_x^\infty} \;. 
	\end{equ}
	Similarly, given a $\genV$-valued distribution $v_0$ such that $e^{s \Delta_A} v_0$ is well-defined for all $s \in (0, T]$, we let
	\begin{equ}
		\| v_0\|_{\covHolder{A, T}^{\alpha}} \eqdef \max_{0 \leq m \leq 1} \sup_{s \in (0, T] } \big( s^{\frac{-\alpha + m}{2}} \wedge 1 \big) \|\covd^{(m)}_{B(s)} \, e^{s\Delta_{A}} v_0 \|_{L_x^\infty} \;.
	\end{equ}
    More generally, for $k \geq 1$, we let
    \begin{equ}
    \| u_0 \|_{\covHolder{A, F, T, k}^{\alpha}} \eqdef \max_{0 \leq m \leq k} \sup_{s \in (0, T] } \big( s^{\frac{-\alpha + m}{2}} \wedge 1 \big) \|\covd^{(m)}_{B(s)} \, e^{s\Delta^F_{A}} u_0 \|_{L_x^\infty} \;,
    \end{equ}
    and similarly for $\|v_0\|_{\covHolder{A, T, k}^\alpha}$.
\end{definition}

\begin{remark}
We refer to these as gauge invariant Hölder norms because given a smooth gauge transformation $g \in \CG$, we have that
\begin{equ}
\|g \act u_0 \|_{\covHolder{g \act A, F, T, k}^\alpha} = \|u_0\|_{\covHolder{A, F, T, k}^\alpha}, 
\end{equ}
and similarly for $\|\place\|_{\covHolder{A, T, k}^\alpha}$. 
\end{remark}

Before proceeding, we set some notation that will simplify many of the ensuing estimates.

\begin{notation}\label{notation:gronwall-F-pert}
Let $T > 0$, and let $B$ be a smooth solution to the Yang--Mills gradient flow on $(0, T] \times \T^2$. Let $C_0$ be the maximum of the corresponding constants in Lemma~\ref{lemma:higher-regularity-abstract-small-Lt1-Lxinfty} and Propositions~\ref{prop:qB-ptwise-bds} and~\ref{prop:covariant-heat-equation-with-curvature-higher-regularity}. For $k \geq 0$, we let
\begin{equ}\label{eq:gronwall-F-pert}
M_T(B, k) \eqdef \gronwall_T(C_0 F_{B}) \bigl(1 + \Fnorm{F_{B}}_{T, 2}\bigr)^k \;.
\end{equ}
\end{notation}

\subsection{YM gradient flow with rough initial data}
\label{sec:rough_YMHF}

In this section, we prove estimates on solutions to the Yang--Mills gradient flow with rough initial data. We will later apply these estimates in a setting where the initial data is made of two pieces: a large smooth piece and a smaller rough piece. Then in order to obtain estimates, it is natural to treat the smooth piece non-perturbatively and the rough piece perturbatively. We accomplish this by expanding around the background arising from the smooth piece. This is the general method of proof behind the following proposition, which shows that solutions to the Yang--Mills gradient flow exist starting from any $\cC_x^{-\kappa}$ initial data, at least after an initial gauge transformation.

\begin{proposition}[Yang--Mills gradient flow with $\cC_x^{-\kappa}$ initial data]\label{pr:admissible_timescale_general}
There is a constant $c_0 \leq 1$ such that the following holds. Suppose that $B_0^{\sharp} = \YMHFsmooth_0 + \YMHFrough_0$, where $\YMHFsmooth_0$ is smooth and $\YMHFrough_0 \in \cC_x^{-\kappa}$. Let $\YMHFsmooth$ be the solution to the Yang--Mills gradient flow with $\YMHFsmooth(0) = \YMHFsmooth_0$. 
Let $0 < T \leq 1$ and $K \geq 1$ be such that 
\begin{equs}
M_T(\YMHFsmooth, 2) T^{\frac{1-\kappa}{2}} K &\leq c_0\;, \label{e:admissible_timescale_smallness}\\
\label{e:Aroughpart_0_smallness_general}
\|\YMHFrough_0\|_{\covHolder{\YMHFsmooth_0, F, T, 2}^{-\kappa}}  &\leq \roughpartsize \;. 
\end{equs}
Then there exist a gauge transformation $g_0 \in \CG_{1-\kappa}$ and a solution $B \in C_{t,x}^\infty((0, T] \times \T^2)$ to the Yang--Mills gradient flow \eqref{e:YMHF} such that $B(t) \ra g_0 \act B_0^{\sharp}$ in $\cC_x^{-\kappa}$ as $t \downarrow 0$. Moreover, there exists $g \in C_t^0 \CG_{1-\kappa}([0, T])$ such that $g$ is smooth at positive times, $g(0) = g_0$, and $g^{-1} \act B = \YMHFsmooth + \YMHFrough$, where (as before) $\YMHFsmooth$ is the solution to the Yang--Mills gradient flow with initial data $\YMHFsmooth(0) = \YMHFsmooth_0$, and $\YMHFrough$ satisfies the following estimate. First, let $\Blin^-(t) = e^{t \Delta_{\smash{\YMHFsmooth_0}}^F} \YMHFrough_0$. For all $t_0 \in (0, T]$, we have that
\begin{equ}\label{e:admissible_timescale_Brem_estimate}
\max_{0 \leq j \leq 2} \bigl\|t^{\frac{j+\kappa}{2}} \covd_{\YMHFsmooth}^{(j)} (\YMHFrough - \Blin^-)\bigr\|_{L_{t,x}^\infty([0, t_0])}  \lesssim M_T(\YMHFsmooth, 2) K^2 t_0^{\frac{1-\kappa}{2}} \;. 
\end{equ}
Moreover, we have the following $L_x^\infty$ and $L_x^2$ curvature estimates:
\begin{equs}\label{e:Fnorm_bound_roughB_general}
\bigl\|\covd_B^{(j)} F_B(t)\bigr\|_{L_x^\infty} \leq \bigl\|\covd_{\YMHFsmooth}^{(j)} F_{\YMHFsmooth}(t)\bigr\|_{L_x^\infty} + \bigO \bigl(t^{-\frac{\kappa}{2}} \bigl(t^{-\frac{1+j}{2}} + \bigl\|F_{\YMHFsmooth}(t)\bigr\|_{L_x^\infty}^j \bigr) \roughpartsize \bigr)\;, ~~ j = 0, 1, \\
\label{e:Fnorm_bound_roughB_general_L2}
\bigl\|\covd_B^{(j)} F_B(t)\bigr\|_{L_x^2} \leq \bigl\|\covd_{\YMHFsmooth}^{(j)} F_{\YMHFsmooth}(t)\bigr\|_{L_x^2} + \bigO \bigl(t^{-\frac{\kappa}{2}} \bigl(t^{-\frac{1+j}{2}} + \bigl\|F_{\YMHFsmooth}(t)\bigr\|_{L_x^2}^j \bigr) \roughpartsize \bigr) \;, ~~ j = 0, 1,
\end{equs}
and for all $t \in (0, T]$, we have that
\begin{equs}
~&\|F_B(t)\|_{L^2}^2 = \|F_{\YMHFsmooth}(t)\|_{L^2}^2  +  \bigl \| \curl_{\YMHFsmooth(t)} \Blin^-(t) \bigr\|_{L_x^2}^2 + 2 \, \bigl \langle F_{\YMHFsmooth}^{12}(t), \bigl(\curl_{\YMHFsmooth} \Blin^- + [\Blin^{-, 1}, \Blin^{-,2}] \bigr)(t)  \bigr\rangle \\
& + M_T(\YMHFsmooth, 3) \bigO\Bigl(\roughpartsize^2 t^{-\kappa} \bigl\| \covd_{\YMHFsmooth} F_{\YMHFsmooth} \bigr\|_{L_{t,x}^2([0, t])}  +\roughpartsize^3 t^{-\frac12 - \frac32 \kappa} + \bigl\|F_{\YMHFsmooth}(t)\bigr\|_{L_x^2} \roughpartsize^3 t^{\frac{1}{2}-\frac{3\kappa}{2}}\Bigr) \;.\label{e:prop_perturbed_curvature_average_general}
\end{equs}
\end{proposition}

\begin{remark}
In principle, one could always take $\YMHFsmooth_0 = 0$ and set $K = \|B_0^\sharp\|_{\covHolder{0, F, T, 2}^{-\kappa}}$.
However, the statement is particularly useful in situations where the rough part is a perturbation of the smooth
part and one has $K \ll \|B_0^\sharp\|_{\covHolder{0, F, T, 2}^{-\kappa}}$. This allows us to take the existence time-scale $T$ to be much longer than otherwise possible.
\end{remark}

It is natural to wonder in what sense the solution to the Yang--Mills gradient flow $B$ obtained in Proposition~\ref{pr:admissible_timescale_general} is unique. The point here is that the Yang--Mills gradient flow is not fully parabolic, so that low-regularity well-posedness is not as direct as for the DeTurck Yang--Mills flow. We record the following general lemma which clarifies the uniqueness of $B$. If the reader wishes, the statement and its proof (which appears a bit later) may be skipped on a first reading.

\begin{lemma}[Low-regularity uniqueness of the Yang--Mills gradient flow]\label{lemma:ymgf-low-regularity-uniqueness}
Let $A_0 \in \cC_x^{-\kappa}$. Let $T > 0$, and let $B, \hat{B} \in C_t^0 \cC_x^{-\kappa}([0, T] \times \T^2) \cap C_{t,x}^\infty((0, T] \times \T^2)$ be such that $B(t), \hat{B}(t) \ra A_0$ in $\cC_x^{-\kappa}$ as $t \downarrow 0$, and such that $B, \hat{B}$ are solutions to the Yang--Mills gradient flow at positive times. Moreover, suppose that
\begin{equ}
\bigl\|t^{1+\frac{\kappa}{2}} \covd_B F_B \bigr\|_{L_{t,x}^\infty([0, T])}, \bigl\|t^{1+\frac{\kappa}{2}} \covd_{\hat{B}} F_{\hat{B}} \bigr\|_{L_{t,x}^\infty([0, T])} < \infty.
\end{equ}
Then $B = \hat{B}$.
\end{lemma}

\begin{remark}
A more intrinsic formulation of Lemma~\ref{lemma:ymgf-low-regularity-uniqueness} would be to prove a low-regularity well-posedness result for $\cC_x^{-\kappa}$ initial data $A_0$ satisfying an ``approximate Coulomb gauge condition'', in that the divergence $\ptl^\ell (A_0)_\ell$ lies in a higher regularity space than expected (say, $\cC_x^{-1+\varep}$ for some $\varep > 0$). In contrast, Lemma~\ref{lemma:ymgf-low-regularity-uniqueness}
 assumes that we already have a solution to the Yang--Mills gradient flow with low regularity initial data, and then proceeds to show that the solution is unique.

We expect that the more intrinsic result may be obtained by analysing a system for the divergence and curl of $B$. Such a ``div-curl system'' has for instance been studied in \cite[Section 5.4]{OT22}, albeit at much higher regularity. Nonetheless, we expect that one should proceed similarly as in the cited reference, i.e.\ one must use that the curl satisfies a parabolic equation, while the divergence only satisfies an ODE (in particular, there is no smoothing). The main question then reduces to finding appropriate function spaces for the curl and divergence which are compatible with the low-regularity initial data, and such that the nonlinear estimates for the div-curl system may be closed. While this should be doable, we decided to state a more practical version of uniqueness in the present paper, as the proof relies on fewer new ideas and thus is much shorter, and in any case, the uniqueness result will not actually be needed in the eventual global iteration. The more intrinsic version is left as a side question which we find intriguing in its own right.
\end{remark}

Before proving Proposition~\ref{pr:admissible_timescale_general}, we state the following lemma which compares solutions to covariant heat equations with operator $\Delta_B$ and $\Delta_{\YMHFsmooth}$, where $B, \YMHFsmooth$ are as in Proposition~\ref{pr:admissible_timescale_general}. This will be needed later in Section~\ref{section:global-iteration-monotone-quantities-deterministic}. We defer the proof of the lemma until after the proof of the proposition.

\begin{lemma}\label{le:B-YMHFsmooth-heat-comparison}
Let all notation and assumptions be as in Proposition~\ref{pr:admissible_timescale_general}. Upon possibly taking the constant $c_0$ appearing in \eqref{e:admissible_timescale_smallness} smaller, the following holds. Let $g \colon [0, T] \times \T^2 \ra \liegroup$ be given as in the statement of Proposition~\ref{pr:admissible_timescale_general}. Let $\phi_1, \phi_2 \colon [0, T] \times \T^2 \ra \genV$ be such that $(g^{-1} \act \phi_1)(0) = \phi_2(0)$, and moreover suppose that $\phi_1, \phi_2$ satisfy:
\begin{equ}
(\ptl_t - \Delta_B) \phi_1 = 0, \qquad
(\ptl_t - \Delta_{\YMHFsmooth}) \phi_2 = 0.
\end{equ}
Then we have that for any $t_0 \in (0, T]$, $\alpha \in [0, 1/2]$, $1 \leq p \leq \infty$,
\begin{equ}
\bigl\| t^{\frac{\alpha}{2}} \bigl( g^{-1} \act \phi_1 - \phi_2\bigr)  \bigr\|_{L_t^\infty L_x^p([0, t_0])} \lesssim \bigl(1 + \Fnorm{F_{\YMHFsmooth}}_{t_0, 1}\bigr) t_0^{\frac{1-\kappa}{2}} K \bigl\|t^{\frac{\alpha}{2}} \phi_2\bigr\|_{L_t^\infty L_x^p([0, t_0])}.
\end{equ}
Additionally, we have that
\begin{equ}
\bigl\|t^{\frac{\alpha}{2}} \phi_1\bigr\|_{L_t^\infty L_x^p([0, t_0])} + \bigl\|t^{\frac{1+\alpha}{2}} \covd_B \phi_1\bigr\|_{L_t^\infty L_x^p([0, t_0])} \lesssim \bigl(1 + \Fnorm{F_{\YMHFsmooth}}_{t_0, 1}\bigr)^2 \bigl(1 +  t_0^{\frac{1-\kappa}{2}} K \bigr) \bigl\|t^{\frac{\alpha}{2}} \phi_2\bigr\|_{L_t^\infty L_x^p([0, t_0])}.
\end{equ}
Similarly, if $\psi_1, \psi_2 \colon [0, T] \times \T^2 \ra \genV^2$ are such that $(g^{-1} \act \psi_1)(0) = \psi_2(0)$, and in addition they satisfy the equations:
\begin{equ}
(\ptl_t - \Delta_B - 2F_B) \psi_1 = 0, \quad (\ptl_t - \Delta_{\YMHFsmooth} - 2F_{\YMHFsmooth}) \psi_2 = 0,
\end{equ}
then for any $t_0 \in (0, T]$, $\alpha \in [0, 1/2]$ and $1 \leq p \leq \infty$, we have that
\begin{equ}
\bigl\| t^{\frac{\alpha}{2}} \bigl( g^{-1} \act \psi_1 - \psi_2\bigr)  \bigr\|_{L_t^\infty L_x^p([0, t_0])} \lesssim M_T(\YMHFsmooth, 1) t_0^{\frac{1-\kappa}{2}} K \bigl\|t^{\frac{\alpha}{2}} \psi_2\bigr\|_{L_t^\infty L_x^p([0, t_0])}.
\end{equ}    
\end{lemma}

We now begin to prove Proposition~\ref{pr:admissible_timescale_general}. Towards this end, we will implement the classical DeTurck trick. In particular, we first prove local well-posedness in the DeTurck gauge (with appropriate choice of background connection).

\begin{lemma}[Local solution theory for Yang--Mills heat flow in the infinitesimal DeTurck gauge]\label{le:perturbed_YMHF-general}
Under the assumptions of Proposition~\ref{pr:admissible_timescale_general}, there exists a unique $B^{\sharp} \in C_t^0 \cC_x^{-\kappa}([0, T] \times \T^2)$ such that on $(0, T] \times \T^2$, $B^{\sharp}$ is a smooth classical solution to the equation
\begin{equ}\label{eq:B-sharp}
\partial_t B^{\sharp}_k = \covd_{B^{\sharp}}^j (F_{B^{\sharp}})_{j k} + \covd_{B^{\sharp}, k} \covd_{\YMHFsmooth}^j(B^{\sharp}_j - \YMHFsmooth_j), ~~ k \in [2],
\end{equ}
and moreover $B^{\sharp}(t) \ra  B_0^\sharp$ in $\cC_x^{-\kappa}$ as $t \downarrow 0$. Additionally, letting $\Blin^-(t) = e^{t \Delta^F_{\smash{\YMHFsmooth_{0}}}} \YMHFrough_{0}$, we have that for all $t_0 \in (0, T]$ (recall \eqref{eq:gronwall-F-pert} for the definition of $M_T(\YMHFsmooth, 2)$)
\begin{equ}\label{e:perturbed_YMHF_general}
\max_{0 \leq j \leq 2} \bigl\|t^{\frac{j +\kappa}{2}}\covd_{\YMHFsmooth}^{(j)} (B^{\sharp} - \YMHFsmooth - \Blin^-)\bigr\|_{L_{t,x}^\infty([0, t_0])} \lesssim M_T(\YMHFsmooth, 2) \roughpartsize^2 t_0^{\frac{1-\kappa}{2}}.
\end{equ}
Finally, there exists a gauge transformation 
\begin{equ}\label{eq:gauge-transformation-regularity}
g \in  C_t^0 \CG_{1-\kappa}([0, T]) \cap C^{\infty}_{t,x}((0, T] \times \T^2 ),
\end{equ}
such that for each $x \in \T^2$, $g(\place, x)$ is a solution to the ODE
\begin{equs}
\bigl(\ptl_t g g^{-1}\bigr)(t, x) &= \Ad_g \bigl(\covd_{\YMHFsmooth}^\ell (B^{\sharp}_\ell - \YMHFsmooth_\ell)\bigr) (t, x), ~~ t \in (0, T], \label{eq:perturbed-YMHF-g-eqution}\\
g(T, x) &= \groupid,
\end{equs}
and such that $B = g \act B^\sharp$ is a smooth classical solution to the Yang--Mills gradient flow on $(0, T] \times \T^2$, which moreover satisfies $B(t) \ra g(0) \act B_0^\sharp$ in $\cC_x^{-\kappa}$ as $t \downarrow 0$.
\end{lemma}
\begin{proof}
To solve \eqref{eq:B-sharp}, we make the ansatz $B^{\sharp} = \YMHFsmooth + \YMHFrough$, where recall $\YMHFsmooth$ is a solution to the Yang--Mills gradient flow with initial data $\YMHFsmooth_0$. By expanding out the equation \eqref{eq:B-sharp} for $B^{\sharp}$ (using the identity \eqref{prelim:eq-linearization-2}), it suffices to solve for $\YMHFrough$ which satisfies the following equation:
\begin{equ}\label{eq:B-minus-equation}
(\ptl_t - \Delta_{\YMHFsmooth} - 2 \ad(F_{\YMHFsmooth})) \YMHFrough = \schO(\YMHFrough \covd_{\YMHFsmooth} \YMHFrough) + \schO((\YMHFrough)^3), \qquad \YMHFrough(0) = \YMHFrough_0.
\end{equ}
This equation is a non-degenerate semilinear parabolic PDE, so local well-posedness in $\cC_x^{-\kappa}$
and smoothness at positive times are standard. (A 
simple power counting argument shows that one merely needs to take $\kappa < \frac12$.) 
We therefore assume henceforth that $\YMHFrough$ exists and we focus on deriving a priori bounds. 
Towards this end, we set $\Brem = \YMHFrough - \Blin^-$ and we make the soft assumption that 
$\YMHFrough_0$ is smooth (and thus also $\YMHFrough$), which may easily be removed by a limiting argument. The remainder $\Brem$ solves 
\begin{equ}
(\ptl_t - \Delta_{\YMHFsmooth} - 2 \ad(F_{\YMHFsmooth})) \Brem = \schO(\YMHFrough \covd_{\YMHFsmooth} \YMHFrough) + \schO((\YMHFrough)^3), \qquad \Brem(0) = 0.
\end{equ}
Let $G$ be the right-hand side of the above equation. For $t_0 \in (0, T]$, define
\begin{equ}
\bigl\| \YMHFrough \bigr\|_{\perturblemnorm(t_0)} \eqdef \max_{0 \leq j \leq 2} \bigl\|t^{\frac{j+\kappa}{2}} \covd_{\YMHFsmooth}^{(j)} \YMHFrough\bigr\|_{L_{t,x}^\infty([0, t_0])}.
\end{equ}
For each $t_0 \in (0, T]$, it follows from Lemma~\ref{lemma:higher-regularity-covd-ibp} (with $k \in \{0, 1, 2\}$) that
\begin{equ}
\max_{0 \leq j \leq 2} t_0^{\frac{j+ \kappa}{2}}\bigl\|\covd_{\YMHFsmooth}^{(j)} \Brem(t_0)\bigr\|_{L_x^\infty}   \lesssim t_0^{\frac{\kappa}{2}} M_T(\YMHFsmooth, 2) N_2(G),
\end{equ}
where
\begin{equs}
N_2(G)& = \big\|G\big\|_{L_t^1 L_x^\infty([0, t_0])} \\
&+ \sum_{j=0}^{1} \Bigl( \big\|t^{(j+1)/2} \covd_{\YMHFsmooth}^{(j)} G \big\|_{L_t^2 L_x^\infty([0, t_0])} + \sup_{u \in [0, t_0]} \big\|t^{(j+1)/2} (u - t)^{-1/2} \covd_{\YMHFsmooth}^{(j)} G \big\|_{L_t^1 L_x^\infty([0, u])}\Bigr).
\end{equs}
We now estimate the norms on $G$ appearing on the right-hand side above. Towards this end, from the definition of $G$, we may estimate
\begin{equs}
\sup_{u \in (0, t_0]} \bigl\|t^{1/2} (u - t)^{-1/2} G \bigr\|_{L_t^1 L_x^\infty([0, u])} &\lesssim t_0^{\frac{1}{2}-\kappa} \bigl\| \YMHFrough \bigr\|_{\perturblemnorm(t_0)}^2 + t_0^{1-\frac{3\kappa}{2}}\bigl\| \YMHFrough \bigr\|_{\perturblemnorm(t_0)}^3, \\
\sup_{u \in (0, t_0]} \bigl\|t (u - t)^{-1/2} \covd_{\YMHFsmooth} G \bigr\|_{L_t^1 L_x^\infty([0, u])} &\lesssim t_0^{\frac{1}{2}-\kappa} \bigl\| \YMHFrough \bigr\|_{\perturblemnorm(t_0)}^2 + t_0^{1-\frac{3\kappa}{2}}\bigl\| \YMHFrough \bigr\|_{\perturblemnorm(t_0)}^3.
\end{equs}
All the other norms appearing in $N_2(G)$ may be estimated by similar means. Combining the previous few estimates,
we thus obtain a constant $C$ such that for all $t_0 \in (0, T]$ (here we also use \eqref{e:Aroughpart_0_smallness_general} which implies that $\bigl\|\Blin^-\bigr\|_{\perturblemnorm(t_0)} \leq \roughpartsize$, and the fact that $\YMHFrough = \Blin^- + \Brem$)
\begin{equs}
\bigl\|\Brem\bigr\|_{\perturblemnorm(t_0)} &\leq C M_T(\YMHFsmooth, 2)  \bigl( t_0^{\frac{1-\kappa}{2}} \bigl\| \YMHFrough \bigr\|_{\perturblemnorm(t_0)}^2 + t_0^{1-\kappa}\bigl\| \YMHFrough \bigr\|_{\perturblemnorm(t_0)}^3\bigr), \label{eq:Brem-continuity-estimate}\\
\bigl\| \YMHFrough \bigr\|_{\perturblemnorm(t_0)} &\leq \roughpartsize + C M_T(\YMHFsmooth, 2) \bigl( t_0^{\frac{1-\kappa}{2}} \bigl\| \YMHFrough \bigr\|_{\perturblemnorm(t_0)}^2 + t_0^{1-\kappa}\bigl\| \YMHFrough \bigr\|_{\perturblemnorm(t_0)}^3\bigr) \label{eq:YMHFrough-continuity-estimate}. 
\end{equs}
We now take the constant $c_0$ in our smallness assumption \eqref{e:admissible_timescale_smallness} to be small enough such that $10^4 C c_0 \leq 1$, where $C$ is as in \eqref{eq:YMHFrough-continuity-estimate}. We proceed by a standard continuity argument. Let 
\begin{equ}
I = \{t_0 \in (0, T] : \bigl\|\YMHFrough\bigr\|_{\perturblemnorm(t_0)} \leq 10 \roughpartsize\},
\end{equ}
where $K$ is as in the lemma statement. By \eqref{eq:YMHFrough-continuity-estimate}, $I$ contains all sufficiently small $t_0 > 0$. Moreover, since we made the soft assumption that $\YMHFrough_0$ is smooth, we have that $t \mapsto \bigl\|\YMHFrough\bigr\|_{\perturblemnorm(t)}$ is continuous in $t$, and thus $I$ is a closed subset of $(0, T]$. To finish, we show that $I$ is an open subset of $(0, T]$. First, note that by the smallness assumption \eqref{e:Aroughpart_0_smallness_general}, our choice of $c_0$, and the definition of $I$, we have that for all $t_0 \in I$,
\begin{equs}
C M_T(\YMHFsmooth, 2)  t_0^{\frac{1-\kappa}{2}} \bigl\| \YMHFrough\bigr\|_{\perturblemnorm(t_0)} &\leq 10^2 C M_T(\YMHFsmooth, 2) t_0^{\frac{1-\kappa}{2}} K \\
&\leq 10^2 C c_0 \leq \frac{1}{100}.
\end{equs}
Thus by \eqref{eq:YMHFrough-continuity-estimate}, it follows that for all $t_0 \in I$,
\begin{equ}
\bigl\|\YMHFrough\bigr\|_{\perturblemnorm(t_0)} \leq \roughpartsize + \frac{1}{2} \bigl\|\YMHFrough\bigr\|_{\perturblemnorm(t_0)},
\end{equ}
from which we obtain
\begin{equ}
\bigl\|\YMHFrough\bigr\|_{\perturblemnorm(t_0)}  < 10 \roughpartsize,
\end{equ}
and thus $I$ is an open subset of $(0, T]$. It follows that $I = (0, T]$, as desired. 
The bound \eqref{e:perturbed_YMHF_general} is then a direct consequence of the smallness condition \eqref{e:admissible_timescale_smallness}, the estimate \eqref{eq:Brem-continuity-estimate}, and the estimate $\bigl\|\YMHFrough\bigr\|_{\perturblemnorm(t_0)} \leq 10K$ for all $t_0 \in (0, T]$.

It remains to prove the claimed results about the gauge transformation $g$. First, assuming that \eqref{eq:gauge-transformation-regularity} holds, one may show by direct calculation (using for example
\eqref{eq:perturbed-YMHF-g-eqution} and the fully covariant formulation introduced in Section~\ref{section:identities-YM-gradient-flow}) that $B = g \act B^\sharp$ is indeed a classical solution to the Yang--Mills gradient flow on $(0, T] \times \T^2$, which moreover satisfies $B(t) \ra g(0) \act B_0^\sharp$ in $\cC_x^{-\kappa}$ as $t \downarrow 0$. Thus, it remains to show that the solution $g$ to \eqref{eq:perturbed-YMHF-g-eqution} indeed exists and satisfies the regularity conditions specified by \eqref{eq:gauge-transformation-regularity}.

Since $\YMHFsmooth$, $B^{\sharp}$ are smooth on $(0, T] \times \T^2$,  there is a unique smooth solution $g$ to \eqref{eq:perturbed-YMHF-g-eqution} on $(0, T] \times \T^2$ by standard ODE theory. Moreover, by the estimate $\bigl\|\YMHFrough\bigr\|_{\perturblemnorm(T)} \leq 10 \roughpartsize$, we have that 
\begin{equ}\label{eq:YMHF-perturbed-ptl-t-g}
\bigl\|t^{\frac{1+\kappa}{2}} \ptl_t g g^{-1} \bigr\|_{L_{t,x}^\infty([0, T])} \lesssim \roughpartsize.
\end{equ}
From this, it follows that $g$ extends continuously to $[0, T] \times \T^2$, and moreover 
\begin{equ}
g(0, x) = \groupid - \int_0^T \dt \, (\ptl_t g g^{-1})(t, x) g(t, x), \quad x \in \T^2.
\end{equ}
It remains to show that $g \in C_t^0 \CG_{1-\kappa}([0, T])$, i.e.\ we need to estimate the $C_t^0 \Cs_x^{1-\kappa}([0, T] \times \T^2)$ norm of $g$. We will first derive an estimate in $L_t^\infty \Cs_x^{1-\kappa}((0, T] \times \T^2)$, and then discuss how to obtain the continuity at $t = 0$. Towards this end, we will first need an estimate of $(dg g^{-1})(t_0)$ at positive times $t_0 > 0$. In particular, we claim that
\begin{equ}\label{eq:perturbed-YMHF-dg-bound}
\bigl\|t^{\frac{\kappa}{2}} d g g^{-1}\bigr\|_{L_{t,x}^\infty([0, T])} \lesssim 1,
\end{equ}
where the implicit constant may depend on $\YMHFsmooth$ and $\roughpartsize$. Given this claim, for $t \in [0, T]$, $x, y \in \T^2$ such that $|x-y|^2 \leq T$, we let $t_0 = |x-y|^2 \vee t$, and then estimate
\begin{equs}
|g(t, x) - g(t, y)| &\leq \ind(t < t_0) \bigl( |g(t, x) - g(t_0, x)| + |g(t_0, y) - g(t, y)|\bigr) + |g(t_0, x) - g(t_0, y)|  \\
&\lesssim \ind(t < t_0) t_0^{\frac{1-\kappa}{2}} + |x-y| t_0^{-\frac{\kappa}{2}} \lesssim |x-y|^{1-\kappa}.
\end{equs}
Here, we estimated the first two terms (in the case $t < t_0$, otherwise $t = t_0$ and these terms are zero) by $t_0^{\frac{1-\kappa}{2}}$ (using \eqref{eq:YMHF-perturbed-ptl-t-g}), and we estimated the last term by $|x-y| t_0^{-\frac{\kappa}{2}}$ (using \eqref{eq:perturbed-YMHF-dg-bound}). Since the implicit constants are uniform in $t, x, y$, we obtain that
$g \in L_t^\infty \Cs_x^{1-\kappa}((0, T] \times \T^2)$, as desired. Moreover, we obtain that $g(0)$ is $(1-\kappa)$-H\"{o}lder. At the end, we will show that $g(0)$ is a limit of smooth functions in the $\Cs_x^{1-\kappa}$ norm, and thus $g(0) \in \CG_{1-\kappa}$.

To complete the proof of the $L_t^\infty \Cs_x^{1-\kappa}((0, T]\times \T^2)$ estimate, it remains to show \eqref{eq:perturbed-YMHF-dg-bound}. We let $h_\alpha = - (\ptl_\alpha g) g^{-1}$ (here $\alpha \in \{0\} \cup [2]$, recall Notation~\ref{notation:space-time-indices}), so that $h = g \act 0$ (here both sides are viewed as space-time connections). Since $h$ is a pure gauge, its curvature is zero, and thus for each $i \in [2]$, we have that $0 = (F_h)_{0i} = \ptl_0 h_i - \covd_{h, i} h_0$, i.e.\
\begin{equ}\label{eq:DeTurck-YM-lwp-h-intermediate}
\ptl_t h_i = \covd_{g \act 0, i} (-\ptl_t g g^{-1}) = - \Ad_g \ptl_i\bigl( \Ad_g^{-1} (\ptl_t g g^{-1})\bigr).
\end{equ}
Integrating this on $[t_0, T]$ and using that $g(T) \equiv \groupid$ so that $(dg g^{-1})(T) = 0$, we obtain
\begin{equs}
(\ptl_i g g^{-1})(t_0) &= - \int_{t_0}^T \dt \, \bigl(\Ad_g \ptl_i \covd_{\YMHFsmooth}^j \YMHFrough_j\bigr)(t) \\
&= -\int_{t_0}^T \dt \, \bigl(\Ad_g \covd_{\YMHFsmooth, i} \covd_{\YMHFsmooth}^j \YMHFrough_j\bigr)(t) + \int_{t_0}^T \dt \, \bigl(\Ad_g [\YMHFsmooth_i, \covd_{\YMHFsmooth}^j \YMHFrough_j ] \bigr)(t).
\end{equs}
Since $\YMHFsmooth$ is smooth and $\bigl\| \YMHFrough\bigr\|_{\perturblemnorm(T)} \lesssim \roughpartsize$, the second term may be directly estimated by an acceptable contribution, and thus we focus on the first term. For this term, we may obtain an acceptable contribution upon applying the estimate
\begin{equ}
\bigl\| t^{1 + \frac{\kappa}{2}} \covd_{\YMHFsmooth}^{(2)} \YMHFrough\bigr\|_{L_{t,x}^\infty([0, T])} \lesssim 1,
\end{equ}
where the implicit constant may depend on $\YMHFsmooth$ and $\roughpartsize$. This estimate follows by recalling that $\YMHFrough = \Blin^- + \Brem$, and applying the assumption \eqref{e:Aroughpart_0_smallness_general} and the $j=2$ case of the estimate \eqref{e:perturbed_YMHF_general}.

Finally, we discuss the continuity of $g(t)$ at $t = 0$ in the $\Cs_x^{1-\kappa}$ norm (which in particular implies that $g(0) \in \CG_{1-\kappa}$). Towards this end, it suffices to show that as $t \downarrow 0$, $((g^{-1} dg)(t))_{t > 0}$ forms a Cauchy net in $\cC_x^{-\kappa}$. For notational simplicity, let $\tilde{h} = g^{-1} dg = -\Ad_g^{-1} h$. The reason for considering $\tilde{h}$ (as opposed to $h$) is that its time evolution takes a slightly more convenient form. In particular, starting from \eqref{eq:DeTurck-YM-lwp-h-intermediate}, we may obtain
\begin{equ}
\ptl_t \tilde{h}_i =  \ptl_i \covd_{\YMHFsmooth}^\ell \YMHFrough_\ell + \Ad_g^{-1}[\ptl_t g g^{-1}, h_i].
\end{equ}
Thus for $0 < s < t$, we have that
\begin{equ}
\tilde{h}_i(t) - \tilde{h}_i(s) =  \int_s^t \du \, \bigl(\ptl_i \covd_{\YMHFsmooth}^\ell \YMHFrough_\ell + \Ad_g^{-1}[\ptl_t g g^{-1}, h_i]\bigr).
\end{equ}
By the estimates \eqref{eq:YMHF-perturbed-ptl-t-g} and \eqref{eq:perturbed-YMHF-dg-bound}, the second term above may be directly estimated in $L_x^\infty$ (since the resulting time singularity is integrable). Similarly, expanding $\YMHFrough = \Blin^- + \Brem$ and expanding $\Blin^-$ around $e^{t \Delta} \YMHFrough_0$, we see that all terms may be directly estimated in $L_x^\infty$ using \eqref{e:Aroughpart_0_smallness_general} and \eqref{e:perturbed_YMHF_general}, except for the most singular term (which has a non-integrable time singularity in $L_x^\infty$):
\begin{equ}
\int_s^t \du \, \ptl_i \ptl^\ell e^{u \Delta} \YMHFrough_{0, \ell} =  (-\Delta)^{-1} \ptl_i \ptl^\ell (e^{s \Delta} \YMHFrough_{0, \ell} - e^{t \Delta} \YMHFrough_{0, \ell}).
\end{equ}
To finish, we use that
\begin{equ}
\bigl\|(-\Delta)^{-1} \ptl_i \ptl^\ell (e^{s \Delta} \YMHFrough_{0, \ell} - e^{t \Delta} \YMHFrough_{0, \ell})\bigr\|_{\cC_x^{-\kappa}} \lesssim \bigl\|e^{s \Delta} \YMHFrough_{0} - e^{t \Delta} \YMHFrough_{0}\bigr\|_{\cC_x^{-\kappa}},
\end{equ}
combined with the fact that $e^{t \Delta} \YMHFrough_0$ converges in $\cC_x^{-\kappa}$ as $t \downarrow 0$ (to $\YMHFrough_0$ \dash here we use that $\cC_x^{-\kappa}$ was defined as the closure of $C_x^\infty$; recall Definition~\ref{def:Holder-space}).
\end{proof}

Having shown local well-posedness of the DeTurck Yang--Mills heat flow, we next prove the following curvature estimates. These estimates will be needed to show the curvature estimates claimed in Proposition~\ref{pr:admissible_timescale_general}.

\begin{lemma}\label{le:pertrubted_YMHF_curvature}
In the setting of Lemma~\ref{le:perturbed_YMHF-general}, we have for any $t \in (0, T]$ and $x \in \T^2$ the following additional pointwise curvature estimate:
\begin{equ}
|\covd_{B^\sharp}^{(j)} F_{B^{\sharp}}(t, x)| \leq |\covd_{\YMHFsmooth}^{(j)} F_{\YMHFsmooth}(t, x)| + \bigO \bigl(t^{-\frac{\kappa}{2}} \bigl(t^{-\frac{1+j}{2}} +  \bigl|F_{\YMHFsmooth}(t, x)\bigr|^j \bigr) \roughpartsize \bigr), ~~ j \in \{0, 1\},
\label{e:perturbed_curvature_general}
\end{equ} 
as well as the $L_x^2$ curvature estimate \eqref{e:prop_perturbed_curvature_average_general}, except with $F_B$ 
 replaced by $F_{B^{\sharp}}$ on the left-hand side.
\end{lemma}
\begin{proof}
Let $\Brem = B^{\sharp} - \YMHFsmooth - \Blin^-$. By \eqref{e:admissible_timescale_smallness} and \eqref{e:perturbed_YMHF_general}, we have that for all $t \in (0, T]$ (recall~\eqref{eq:gronwall-F-pert} for the definition of $M_T(\YMHFsmooth, 2)$),
\begin{equs}
\bigl\|\Brem(t)\bigr\|_{L_x^\infty} &\lesssim M_T(\YMHFsmooth, 2) \roughpartsize^2 t^{\frac{1}{2}-\kappa} \lesssim c_0 t^{-\frac{\kappa}{2}} \roughpartsize, \label{eq:Brem-estimate}\\
\bigl\| (\covd_{\YMHFsmooth} \Brem)(t)\bigr\|_{L_x^\infty} &\lesssim M_T(\YMHFsmooth, 2) \roughpartsize^2 t^{-\kappa} \lesssim c_0 t^{-\frac{1+\kappa}{2}} \roughpartsize, \label{eq:curl-Brem-estimate}
\end{equs}
so that $\Brem$ is no larger than $\Blin^-$. On the other hand, we may choose to use the intermediate 
estimates above at several places in the following. To begin, we expand out (by \eqref{prelim:eq-linearization-2})
\begin{equ}
F_{B^{\sharp}}^{12} = F_{\YMHFsmooth}^{12} + \curl_{\YMHFsmooth}(\Blin^- + \Brem) + [(\Blin^{-} + \Brem)^1, (\Blin^- + \Brem)^2].
\end{equ}
The pointwise estimate \eqref{e:perturbed_curvature_general} in the case $j = 0$ follows by combining the triangle inequality with the assumed bound on $\Blin^-$ \eqref{e:Aroughpart_0_smallness_general} and our previous estimates on $\Brem$ (and using that $t^{-\kappa} \roughpartsize^2 \ll t^{-\frac{1+\kappa}{2}} \roughpartsize$ by our smallness condition \eqref{e:admissible_timescale_smallness}). The case $j = 1$ follows by applying $\covd_{B^{\sharp}}$ to both sides of the above identity, expanding (recall $\YMHFrough = \Blin^- + \Brem$) $\covd_{B^{\sharp}} = \covd_{\YMHFsmooth} + [\YMHFrough, \place]$ as well as $F_{B^{\sharp}} = F_{\YMHFsmooth} + \mc{O}(\covd_{\YMHFsmooth} \YMHFrough + (\YMHFrough)^2)$ on the right-hand side to obtain
\begin{equ}
\covd_{B^\sharp} F_{B^{\sharp}}^{12} = \covd_{\YMHFsmooth} F_{\YMHFsmooth}^{12} + [\YMHFrough, F_{\YMHFsmooth}] + \mc{O}(\covd_{\YMHFsmooth}^{(2)} \YMHFrough + \YMHFrough \covd_{\YMHFsmooth} \YMHFrough + (\YMHFrough)^3).
\end{equ}
The second and third terms on the right-hand side may be directly estimated by combining the assumed bound on $\Blin^-$ \eqref{e:Aroughpart_0_smallness_general} with the estimates \eqref{e:perturbed_YMHF_general} on $\Brem$ and the smallness condition \eqref{e:admissible_timescale_smallness}.

To show the $L_x^2$ estimate, we expand out bilinearly to obtain
\begin{equs}
\bigl\|F_{B^{\sharp}}(t)\bigr\|_{L_x^2}^2 &=\|F_{\YMHFsmooth}(t)\|_{L^2}^2  +  \bigl \| \curl_{\YMHFsmooth(t)} \Blin^-(t) \bigr\|_{L_x^2}^2 \\
&\quad + 2 \, \bigl \langle F_{\YMHFsmooth}^{12}(t), \bigl(\curl_{\YMHFsmooth} \Blin^-\bigr)(t) + [\Blin^{-, 1}, \Blin^{-, 2}](t)\bigr\rangle \\ 
&\quad + 2 \langle \curl_{\YMHFsmooth} \Blin^{-}(t) , [\Blin^{-, 1}, \Blin^{-, 2}](t)\rangle + \bigl\|[\Blin^{-, 1}, \Blin^{-, 2}](t)\bigr\|_{L_x^2}^2 \\
&\quad + \text{perturbative terms},
\end{equs}
where ``$\text{perturbative terms}$'' collects all terms with at least one factor of $\Brem$. We first estimate the terms on the second-to-last line, and then proceed to estimate the perturbative terms. By the assumption \eqref{e:Aroughpart_0_smallness_general} on the size of $\YMHFrough_0$, we may directly estimate 
\begin{equ}
2 \abs{\langle \curl_{\YMHFsmooth} \Blin^-(t) , [\Blin^{-, 1}, \Blin^{-, 2}](t)\rangle} + \bigl\|[\Blin^{-, 1}, \Blin^{-, 2}](t)\bigr\|_{L_x^2}^2 \lesssim t^{-\frac{1+3\kappa}{2}} K^3 + t^{-\frac{4\kappa}{2}} K^4 \lesssim t^{-\frac{1+3\kappa}{2}} K^3,
\end{equ}
where we used that $t^{\frac{1-\kappa}{2}} K \leq T^{\frac{1-\kappa}{2}} K \leq 1$ (recall the smallness assumption \eqref{e:admissible_timescale_smallness}). This contribution is acceptable. 

Next, we estimate all perturbative terms by an acceptable contribution. For the most delicate perturbative term (which requires that we first integrate by parts), we first apply Lemma~\ref{le:DF_bound} (and the definition of $\Fnorm{F_{\YMHFsmooth}}_{t, 0}$; recall Definition~\ref{def:Fnorm}) to estimate
\begin{equ}\label{eq:covd-FB-time-weighted-estimate-intermediate}
\bigl\|t^{1/2} \covd_{\YMHFsmooth} F_{\YMHFsmooth}\bigr\|_{L_t^\infty L_x^2([0, t])} \lesssim \bigl(1 + \Fnorm{F_{\YMHFsmooth}}_{t, 0}\bigr) \bigl\|\covd_{\YMHFsmooth} F_{\YMHFsmooth}\bigr\|_{L_{t,x}^2([0, t])}.
\end{equ}
We now estimate the most delicate perturbative term by (here we apply integration by parts (Lemma~\ref{lemma:curl-ibp-identity}) in the first step and \eqref{eq:covd-FB-time-weighted-estimate-intermediate} in the final step):
\begin{equs}
\bigl| \langle F_{\YMHFsmooth}^{12}(t), (\curl_{\YMHFsmooth} \Brem)(t) \rangle\bigl|  &= \bigl| \langle \covd_{\YMHFsmooth}^* F_{\YMHFsmooth}(t), \Brem(t) \rangle \bigr| \\
&\lesssim M_T(\YMHFsmooth, 2) \bigl\| \covd_{\YMHFsmooth} F_{\YMHFsmooth}(t)\bigr\|_{L_x^2} \roughpartsize^2 t^{\frac{1}{2}-\kappa} \\
&\lesssim M_T(\YMHFsmooth, 3)  \bigl\| \covd_{\YMHFsmooth} F_{\YMHFsmooth}\bigr\|_{L_{t,x}^2([0, t])} \roughpartsize^2 t^{-\kappa} .
\end{equs}
All other perturbative terms are less delicate, as they may be directly estimated (i.e.\ without first integrating by parts). For instance, one of the perturbative terms is (here we apply \eqref{eq:Brem-estimate})
\begin{equs}
\bigl| \langle F_{\YMHFsmooth}^{12}(t), [\Brem^1, (\Blin^- + \Brem)^2](t)\rangle \bigr| &\lesssim \bigl\|F_{\YMHFsmooth}(t) \bigr\|_{L_x^2} \bigl\|\Brem\bigr\|_{L_x^\infty} \bigl\|\Blin^- + \Brem\bigr\|_{L_x^\infty} \\
&\lesssim M_T(\YMHFsmooth, 2) \bigl\|F_{\YMHFsmooth}(t)\bigr\|_{L_x^2} \roughpartsize^3 t^{\frac{1}{2}-\frac{3\kappa}{2}} .
\end{equs}
Another perturbative term is (here we apply \eqref{eq:curl-Brem-estimate})
\begin{equ}
\bigl| \langle (\curl_{\YMHFsmooth} \Brem)(t), (\curl_{\YMHFsmooth} (\Blin^- + \Brem))(t) \rangle \bigr| \lesssim M_T(\YMHFsmooth, 2) \roughpartsize^3 t^{-\frac{1}{2}-\frac{3\kappa}{2}}.
\end{equ}
All other terms may be estimated by similar arguments.
\end{proof}

Combining the previous few intermediate results, we may now prove Proposition~\ref{pr:admissible_timescale_general}.

\begin{proof}[Proof of Proposition~\ref{pr:admissible_timescale_general}]
We take the constant $c_0$ to be the constant in Lemma~\ref{le:perturbed_YMHF-general}. Then under the assumptions of the proposition, we may apply Lemma~\ref{le:perturbed_YMHF-general} to obtain $B^{\sharp}, g$ with the properties as stated in Lemmas~\ref{le:perturbed_YMHF-general}--\ref{le:pertrubted_YMHF_curvature}. In particular, $B = g \act B^{\sharp}$ is a Yang--Mills gradient flow with the claimed properties.

Next, note that $g^{-1} \act B = \YMHFsmooth + \YMHFrough$, and the claimed estimates \eqref{e:admissible_timescale_Brem_estimate} on $\YMHFrough$ follow directly from \eqref{e:perturbed_YMHF_general}. Finally, for the curvature estimates, we use that $\covd_B^{(j)} F_B = \Ad_g \covd_{B^{\sharp}}^{(j)} F_{B^{\sharp}}$. From this, \eqref{e:perturbed_curvature_general} directly implies \eqref{e:Fnorm_bound_roughB_general} and \eqref{e:Fnorm_bound_roughB_general_L2},
and the second estimate \eqref{e:prop_perturbed_curvature_average_general} follows directly from the corresponding estimate given by Lemma~\ref{le:pertrubted_YMHF_curvature}.
\end{proof}

Next, we prove the uniqueness result Lemma~\ref{lemma:ymgf-low-regularity-uniqueness}. We reiterate that this may be skipped on a first reading. 

\begin{proof}[Proof of Lemma~\ref{lemma:ymgf-low-regularity-uniqueness}]
Without loss of generality, we may take $T \leq 1$ (since $B$ is smooth at positive times and uniqueness of the Yang--Mills gradient flow from smooth initial data is classical \cite{Rad92}). We proceed by solving for a gauge transformation $h$ that puts $B$ into the DeTurck gauge, and then resorting to a more standard well-posedness statement for the DeTurck Yang--Mills heat flow with $\cC_x^{-\kappa}$ initial data. Towards this end, note that if $h \in C_t^0 \cC_x^{2-2\kappa}([0, T] \times \T^2) \cap C_{t,x}^\infty((0, T] \times \T^2)$ is a solution to:
\begin{equ}
(\ptl_t h) h^{-1} = \ptl^\ell ((\ptl_\ell h) h^{-1}) - \ptl^\ell \bigl(\Ad_h B_\ell\bigr), ~~ h(0) = \groupid, 
\end{equ}
then $\tilde{B} = h \act B$ is a classical solution to the DeTurck Yang--Mills heat flow at positive times, which moreover satisfies $\tilde{B}(t) \ra A_0$ in $\cC_x^{-\kappa}$. (To see this, it is helpful to use the fully covariant formulation \dash recall Section~\ref{section:identities-YM-gradient-flow}). To see why a solution $h$ with the claimed regularity exists, note that the equation for $h$ may be written as a nonlinear heat equation with forcing:
\begin{equ}
(\ptl_t h) h^{-1} = (\Delta h) h^{-1} - (\ptl_\ell h) h^{-1} (\ptl^\ell h) h^{-1} - \ad((\ptl^\ell h) h^{-1}) \Ad_h B_\ell - \Ad_h (\ptl^\ell B_\ell).
\end{equ}
Recall by Lemma~\ref{lemma:ymgf-approximate-coulomb-gauge} (here we use the assumed bound on $\covd_B F_B$, as well as our reduction to $T \leq 1$) that 
\begin{equ}
\bigl\|t^{\frac{\kappa}{2}} B \bigr\|_{L_{t,x}^\infty([0, T])}, \bigl\|t^{\kappa} \ptl^\ell B_\ell \bigr\|_{L_{t,x}^\infty([0, T])} < \infty.
\end{equ}
Therefore by a standard Picard iteration argument, for some possibly smaller time $0 < T_1 \leq T$, we may solve for $h \in C_t^0 \cC_x^{2-2\kappa}([0, T_1] \times \T^2)$. Since $B$ is smooth at positive times, we may moreover obtain that $h \in C_{t,x}^\infty((0, T_1] \times \T^2)$. Thus by the previous discussion, we have that $\tilde{B} = h \act B$ is a smooth classical solution to the DeTurck Yang--Mills heat flow on the time interval $(0, T_1]$, and moreover $\tilde{B}(t) \ra A_0$ in $\cC_x^{-\kappa}$ as $t \downarrow 0$. Since the DeTurck Yang--Mills heat flow is parabolic, we may show that $\tilde{B}$ is unique (once the initial data $A_0 \in \cC_x^{-\kappa}$ has been fixed) by a standard parabolic well-posedness argument. Moreover, the equation for the gauge transformation $h$ may be written in terms of $\tilde{B}$ as the collection of ODEs (one for each spatial point $x \in \T^2$):
\begin{equ}
(\ptl_t h) h^{-1} = - \ptl^\ell \tilde{B}_\ell, ~~ h(0) = \groupid.
\end{equ}
(We remark that these ODEs may be integrated up from time zero, since by standard parabolic theory, one may show that $\bigl\|t^{\frac{1+\kappa}{2}} \ptl^\ell \tilde{B}_\ell \bigr\|_{L_{t,x}^\infty([0, T_1])} < \infty$, and so the resulting time singularity is integrable.) Thus since $\tilde{B}$ is uniquely determined, so is the gauge transformation $h$, and therefore the original $B = h^{-1} \act \tilde{B}$ is uniquely determined as well, at least on $[0, T_1]$. The extension to $[0, T]$ follows by uniqueness for smooth initial data \cite{Rad92}.
\end{proof}

Using the same ideas as in the proof of Lemma~\ref{lemma:ymgf-low-regularity-uniqueness}, we can show that the $\invnorm{\place}{-\kappa}$ norm from Definition~\ref{main:def-invnorm} is well-defined. We quickly do this next.

\begin{lemma}\label{lemma:invnorm-well-def}
Let $A_0, \higgs_0$ be as in Definition~\ref{main:def-invnorm}, i.e.\ let $A_0, \higgs_0 \in \Cs_x^{-\kappa}$. Then there exists a time $T > 0$ and a gauge transformation $g \in \CG_{1-\kappa}$ such that $g \act (A_0, \higgs_0) \in \initialclass_T$ (recall \eqref{main:eq-initial-class}), as required in the definition. Moreover, if $T_1, T_2 > 0$ are times and $g_1, g_2 \in \CG_{1-\kappa}$ are gauge transformations such that $g_1 \act (A_0, \higgs_0) \in \initialclass_{T_1}$ and $g_2 \act (A_0, \higgs_0) \in \initialclass_{T_2}$, then $g_1 g_2^{-1}$ is smooth, and the resulting Yang--Mills gradient flows $B_1, B_2$ are related by $B_1 = (g_1 g_2^{-1}) \act B_2$. The same is true for the resulting covariant heat flows $\Psi_1, \Psi_2$. In particular, the quantity defined in \eqref{main:eq-invnorm} is the same for $(B_1, \Psi_1)$ as for $(B_2, \Psi_2)$, and thus the $\invnorm{\place}{-\kappa}$ norm is well-defined.
\end{lemma}
\begin{proof}
The existence of $T > 0$, $g \in \CG_{1-\kappa}$, and the Yang--Mills gradient flow $B$ starting from $g \act A_0$ follows by Proposition~\ref{pr:admissible_timescale_general} (taking the $\YMHFsmooth_0$ therein to be zero, the $\YMHFrough_0$ therein to be $A_0$, and $T > 0$ small enough and $K \geq 1$ large enough so that \eqref{e:admissible_timescale_smallness} and \eqref{e:Aroughpart_0_smallness_general} are satisfied). This proposition also provides the required bounds on $F_B$ and $\covd_B^* F_B$ stated in \eqref{main:eq-initial-class} (from which it follows by Lemma~\ref{lemma:ymgf-low-regularity-uniqueness} that $B$ is unique). To see that the covariant Hölder norm of $g \act \higgs_0$ is finite, we use these bounds and Lemma~\ref{lemma:ymgf-approximate-coulomb-gauge} to obtain bounds on $B$ and $\ptl^\ell B_\ell$, and then a standard perturbative argument to obtain the solution to the covariant heat flow starting from $g \act \higgs_0$.

Next, let $T_1, T_2 > 0$, $g_1, g_2 \in \CG_{1-\kappa}$, $B_1, B_2$, and $\Psi_1, \Psi_2$ be as given in the lemma statement. As in the proof of Lemma~\ref{lemma:ymgf-low-regularity-uniqueness}, there exist $0 < T' \leq T_1 \wedge T_2$ and gauge transformations $h_1, h_2 \in C_t^0 \Cs_x^{2-2\kappa}([0, T']) \cap C_{t,x}^\infty((0, T'])$ such that $h_1 \act B_1, h_2 \act B_2$ are DeTurck Yang--Mills flows which are smooth at positive times, and such that $h_1(0) = h_2(0) \equiv \groupid$.

Having reached this point, we finish by a statement of smooth gauge equivalence at positive times between DeTurck Yang--Mills flows whose initial data are related by a $\CG_{1-\kappa}$ gauge transformation. In particular, by a similar argument as in the proof of Lemma~\ref{lemma:ymgf-low-regularity-uniqueness} (except now using weighted spaces), it is possible to solve for a gauge transformation $h$ which is smooth at positive times and which takes the DeTurck Yang--Mills flow $h_2 \act B_2$ to the DeTurck Yang--Mills flow $h_1 \act B_1$. Recalling the fully covariant formulation of the Yang--Mills gradient flow and DeTurck Yang--Mills flow discussed at the beginning of Section~\ref{section:identities-YM-gradient-flow}, the equation specifying $h$ is given by $h \act (h_2 \act B_2) = h_1 \act B_1$ (as space-time connections) with initial data $h(0) = g_1 g_2^{-1}$. 

From this, letting $\tilde{h} = h_1^{-1} h h_2$ (which is smooth at positive times), it follows that $B_1 = \tilde{h} \act B_2$, at least for small enough positive times. We claim that $\tilde{h}$ is constant in time. To see this, we compute
\begin{equ}
((h_1^{-1} h h_2) \act B_2)_0 = (h_1^{-1} \act ((h h_2) \act B_2))_0 = (h_1^{-1} \act (h_1 \act B_1))_0 = (B_1)_0.
\end{equ}
Since $B_1, B_2$ are Yang--Mills gradient flows, we have that $(B_1)_0 = (B_2)_0 = 0$, and thus it follows that $(\ptl_t \tilde{h}) \tilde{h}^{-1} = 0$, and thus $\tilde{h}$ is constant in time. Note that since $\tilde{h}(0) = h_1^{-1}(0) h(0) h_2(0) = g_1 g_2^{-1}$, it follows that $g_1 g_2^{-1} = \tilde{h}$ is smooth and $B_1 = (g_1 g_2^{-1}) \act B_2$, as desired. To finish, note that (by \eqref{eq:g-act-higgs-evolution}) $\tilde{h} \act \Psi_2$ is also a solution to a covariant heat equation with background connection $\tilde{h} \act B_2 = B_1$ and initial data $(\tilde{h} \act \Psi_2)(0) = (g_1 g_2^{-1}) \act (g_2 \act \higgs_0) = g_1 \act \higgs_0$, and thus (by a standard parabolic well-posedness argument) we have that $\tilde{h} \act \Psi_2 = \Psi_1$. The desired result now follows.
\end{proof}

To finish off this section, we prove Lemma~\ref{le:B-YMHFsmooth-heat-comparison}.

\begin{proof}[Proof of Lemma~\ref{le:B-YMHFsmooth-heat-comparison}]
Let $\tilde{\phi}_1 = g^{-1} \act \phi_1$. Recall from Proposition~\ref{pr:admissible_timescale_general} that $g^{-1} \act B = B^{\sharp} = \YMHFsmooth + \YMHFrough$. Moreover, from the proof of the proposition, the gauge transformation $g$ satisfies \eqref{eq:perturbed-YMHF-g-eqution}. Then we have that (recall the fully covariant discussion in Section~\ref{section:identities-YM-gradient-flow})
\begin{equ}
(\ptl_t - \Delta_{B^{\sharp}}) \tilde{\phi}_1 = -(g^{-1} \ptl_t g) \tilde{\phi}_1 = -(\covd_{\YMHFsmooth}^j \YMHFrough_j) \tilde{\phi}_1,
\end{equ}
Let $\psi = \tilde{\phi}_1 - \phi_2$. We then have that (applying \eqref{prelim:eq-linearization-5} to expand $\Delta_{B^{\sharp}}$ around $\Delta_{\YMHFsmooth}$)
\begin{equs}
(\ptl_t - \Delta_{\YMHFsmooth}) \psi &= 2 \covd_{\YMHFsmooth}^j(\YMHFrough_j (\psi + \phi_2)) - 2 (\covd_{\YMHFsmooth}^j \YMHFrough_j) (\psi + \phi_2) + (\YMHFrough)^j \YMHFrough_j (\psi + \phi_2), \\
\psi(0) &=0 .
\end{equs}
By \eqref{e:admissible_timescale_smallness}, \eqref{e:Aroughpart_0_smallness_general}, and \eqref{e:admissible_timescale_Brem_estimate}, we have that
\begin{equ}
\bigl\|t^{\frac{\kappa}{2}} \YMHFrough\bigr\|_{L_{t,x}^\infty([0, T])} + \bigl\|t^{\frac{1+\kappa}{2}} \covd_{\YMHFsmooth} \YMHFrough\bigr\|_{L_{t,x}^\infty([0, T])} \lesssim K. 
\end{equ}
Combining this with the Schauder estimates of Lemma~\ref{lemma:covariant-schauder-estimates}, we may obtain that for any $t_0 \in (0, T]$,
\begin{equ}
t_0^{\frac{\alpha}{2}} \bigl\|\psi(t_0)\bigr\|_{L_x^p} \lesssim \bigl(1 + \Fnorm{F_{\YMHFsmooth}}_{t_0, 1}\bigr)\bigl( t_0^{\frac{1-\kappa}{2}} K + \bigl(t_0^{\frac{1-\kappa}{2}} K \bigr)^2 \bigr)\bigl(\bigl\|t^{\frac{\alpha}{2}} \psi\bigr\|_{L_t^\infty L_x^p([0, t_0])} + \bigl\|t^{\frac{\alpha}{2}} \phi_2\bigr\|_{L_t^\infty L_x^p([0, t_0])}\bigr).
\end{equ}
(Here, we used that $\alpha \in [0, 1/2]$, so that $\alpha + \kappa < 1$, and thus we only see integrable singularities.) Taking a sup and applying \eqref{e:admissible_timescale_smallness}, we thus obtain that for some constant $C$, for all $t_0 \in (0, T]$,
\begin{equ}
\bigl\|t^{\frac{\alpha}{2}} \psi\bigr\|_{L_t^\infty L_x^p([0, t_0])} \leq C (c_0 + c_0^2) \bigl\|t^{\frac{\alpha}{2}} \psi\bigr\|_{L_t^\infty L_x^p([0, t_0])} +  C \bigl(1 + \Fnorm{F_{\YMHFsmooth}}_{t_0, 1}\bigr) t_0^{\frac{1-\kappa}{2}} K \bigl\|t^{\frac{\alpha}{2}} \phi_2\bigr\|_{L_t^\infty L_x^p([0, t_0])}.
\end{equ}
Thus upon possibly taking $c_0$ smaller, we may absorb the first term on the right-hand side back into the left-hand side, to obtain
\begin{equ}
\bigl\|t^{\frac{\alpha}{2}} \psi\bigr\|_{L_t^\infty L_x^p([0, t_0])} \lesssim \bigl(1 + \Fnorm{F_{\YMHFsmooth}}_{t_0, 1}\bigr) t_0^{\frac{1-\kappa}{2}} K \bigl\|t^{\frac{\alpha}{2}} \phi_2\bigr\|_{L_t^\infty L_x^p([0, t_0])},
\end{equ}
which is the first desired estimate.

For the second estimate, first note that $|\phi_1| = |g^{-1} \act \phi_1|$, and thus combining this with the first estimate, we obtain (for some constant $C$)
\begin{equs}
\bigl\|t^{\frac{\alpha}{2}} \phi_1\bigr\|_{L_t^\infty L_x^p([0, t_0])} &\leq \bigl\| t^{\frac{\alpha}{2}} \bigl( g^{-1} \act \phi_1 - \phi_2\bigr)  \bigr\|_{L_t^\infty L_x^p([0, t_0])} + \bigl\|t^{\frac{\alpha}{2}} \phi_2\bigr\|_{L_t^\infty L_x^p([0, t_0])} \\
&\leq \bigl(1 + C  \bigl(1 + \Fnorm{F_{\YMHFsmooth}}_{t_0, 1}\bigr) t_0^{\frac{1-\kappa}{2}} K \bigr) \bigl\|t^{\frac{\alpha}{2}} \phi_2\bigr\|_{L_t^\infty L_x^p([0, t_0])}.\label{e:B-YMHFsmooth-heat-comparison_intermediate}
\end{equs}
To finish, we estimate $\covd_B \phi_1$. By the pointwise estimates for $\fp_B$ (Proposition~\ref{prop:pB-ptwise-bounds}) and Young's convolution inequality, we have that for any $t \in (0, t_0]$,
\begin{equ}
\bigl\|\covd_B \phi_1(t)\bigr\|_{L_x^p} \lesssim \bigl(1 + \Fnorm{F_B}_{t_0, 0}\bigr) t^{-1/2} \bigl\|\phi_1(t/2)\bigr\|_{L_x^p}.
\end{equ}
Recalling the definition of $\Fnorm{\place}$ (Definition~\ref{def:Fnorm}), we see that the $j = 0$ case of the curvature estimate \eqref{e:Fnorm_bound_roughB_general} directly implies that
\begin{equ}
\Fnorm{F_B}_{t_0, 0} \leq \Fnorm{F_{\YMHFsmooth}}_{t_0, 0} + C t_0^{\frac{1-\kappa}{2}} K \lesssim \Fnorm{F_{\YMHFsmooth}}_{t_0, 0}  + 1.
\end{equ}
Inserting this into the previous math display and then multiplying both sides of that display by $t^{\frac{1+\alpha}{2}}$ and taking a sup over $t \in (0, t_0]$, we see that the estimate for $\covd_B \phi_1$ follows from \eqref{e:B-YMHFsmooth-heat-comparison_intermediate}.

Finally, for the third estimate, we proceed similarly to the proof for the first estimate. We now let $\psi = g^{-1} \act \psi_1 - \psi_2$, and compute (applying \eqref{prelim:eq-linearization-5} as well as the curvature expansion \eqref{prelim:eq-linearization-2})
\begin{equ}
(\ptl_t - \Delta_{\YMHFsmooth} - 2F_{\YMHFsmooth}) \psi = \mc{O}\bigl(\covd_{\YMHFsmooth} \bigl(\YMHFrough (\psi + \psi_2) \bigr) \bigr) + \mc{O}\bigl((\covd_{\YMHFsmooth} \YMHFrough + (\YMHFrough)^2) (\psi + \psi_2)\bigr).
\end{equ}
Note also that $\psi(0) = 0$. Applying the Schauder estimates of Lemma~\ref{lemma:covariant-schauder-estimates} as before, we may obtain that for any $t_0 \in (0, T]$,
\begin{equs}
~&t_0^{\frac{\alpha}{2}} \bigl\|\psi(t_0)\bigr\|_{L_x^p} \\
&\lesssim \gronwall_T(C_0 F_{\YMHFsmooth}) \bigl(1 + \Fnorm{F_{\YMHFsmooth}}_{t_0, 1} \bigr) \bigl(t_0^{\frac{1-\kappa}{2}} K + \bigl(t_0^{\frac{1-\kappa}{2}} K \bigr)^2\bigr) \bigl(\bigl\|t^{\frac{\alpha}{2}} \psi\bigr\|_{L_t^\infty L_x^p([0, t_0])} + \bigl\|t^{\frac{\alpha}{2}} \psi_2\bigr\|_{L_t^\infty L_x^p([0, t_0])}\bigr).
\end{equs}
The rest of the argument proceeds similarly to before; we omit the details.
\end{proof}

\subsection{Estimates of deterministic pieces arising from the Higgs field}
\label{sec:higgsdet_estimates}

In this section, we derive estimates on the pair $(\Sdet, \higgsdet)$ which is introduced to better track the deterministic evolution of the Yang--Mills--Higgs gradient flow. These estimates will be needed in both the local theory (see e.g.\ Section~\ref{section:sdet-higgsdet-sol}) as well as in the global theory (see e.g.\ Section~\ref{section:global-iteration-monotone-quantities-deterministic}). Throughout this section, we take $B \in C_t^0 \Cs_x^{-\kappa}([0, T] \times \T^2) \cap C_{t,x}^\infty((0, T]\times \T^2)$ to be a smooth solution to the Yang--Mills gradient flow at positive times. To avoid any pathologies, we will moreover assume that the quantities $M_0, M_1$ appearing in Lemma~\ref{lemma:ymgf-approximate-coulomb-gauge} which are defined out of $B$ are finite. One should imagine (and this will indeed be the case) that $B$ is obtained by an application of Proposition~\ref{pr:admissible_timescale_general}.

Recall that the pair $(\Sdet, \higgsdet)$ is defined as the solution to the system: 
\begin{equs}
(\ptl_t - \Delta_B - 2\ad(F_B)) \Sdet &= [\Sdet^j, 2\covd_{B, j} \Sdet - \covd_{B} \Sdet_j] + [\Sdet^j, [\Sdet_j, \Sdet]] + \bil(\higgsdet \otimes \covd_{B+\Sdet}\higgsdet ) \label{eq:Sdet-equation}\\ 
\Sdet(0) &= 0, \\
(\ptl_t - \Delta_{B+\Sdet}) \higgsdet &= - (\covd_B^j \Sdet_j) \higgsdet -|\higgsdet|^2 \higgsdet, \label{eq:higgsdet-equation}\\
\higgsdet(0) &= \shiggsdet_0 + \rhiggsdet_0.
\end{equs}
To estimate $\higgsdet$, we decompose $\higgsdet = \rlinhiggsdet + \reghiggsdet$, where
\begin{equs}[2]
(\ptl_t - \Delta_B) \rlinhiggsdet &= 0, &\quad \rlinhiggsdet(0) &= \rhiggsdet_0, \label{eq:rlinhiggsdet-equation} \\
(\ptl_t - \Delta_{B+\Sdet}) \reghiggsdet &= 2\Sdet_j \covd_B^j \rlinhiggsdet + \Sdet_j \Sdet^j \rlinhiggsdet - (\covd_B^j \Sdet_j) \reghiggsdet - |\higgsdet|^2 \higgsdet, &\quad \reghiggsdet(0) &= \shiggsdet_0. \label{eq:reghiggsdet-equation}
\end{equs}
Here, we used that (recall \eqref{prelim:eq-linearization-5})
\begin{equ}
(\Delta_{B+\Sdet} - \Delta_B) \rlinhiggsdet - (\covd_B^j \Sdet_j) \rlinhiggsdet = 2 \Sdet_j \covd_B^j \rlinhiggsdet + \Sdet_j \Sdet^j \rlinhiggsdet.
\end{equ}

We state the first main result of this section. Recall the covariant Hölder norm introduced in Definition~\ref{defn:covariant_Holder} as well as the definition of $M_T(B, k)$ (Notation~\ref{notation:gronwall-F-pert}).

\begin{proposition}[Estimates of $\Sdet, \higgsdet$]\label{prop:Sdet-higgsdet-estimates}
Let $C_0$ be the maximum of the corresponding constants in Lemma~\ref{lemma:higher-regularity-abstract-small-Lt1-Lxinfty} and Propositions~\ref{prop:qB-ptwise-bds} and~\ref{prop:covariant-heat-equation-with-curvature-higher-regularity}. There is a constant $c_1 \leq 1$ such that the following holds. Suppose that $T \in (0, 1]$ satisfies
\begin{equ}\label{eq:Sdet-higgsdet-main-prop-smallness-assumption}
M_T(B, 10) \bigl( T^{\frac{1}{2} - \frac{1}{r}} \big\|\shiggsdet_0\big\|_{L_x^r} + T^{\frac{1-\kappa}{2}} \bigl\|\rhiggsdet_0\bigr\|_{\covHolder{B(0), T}^{-\kappa}}\bigr) \leq c_1.
\end{equ}
Then $\Sdet, \higgsdet$ exist on $[0, T]$ and moreover satisfy the following estimates:
\begin{equs}
\max_{0 \leq j \leq 2} \bigl\|t^{\frac{j}{2}} \covd_B^{(j)} \Sdet\bigr\|_{L_{t,x}^\infty([0, T])} &\lesssim M_T(B, 10) \bigl(T^{\frac{1}{2}-\frac{2}{r}} \big\|\shiggsdet_0\big\|_{L_x^r}^2 + T^{\frac{1}{2}-\kappa} \bigl\|\rhiggsdet_0\bigr\|_{\covHolder{B(0), T}^{-\kappa}}^2\bigr), \label{eq:Sdet-L-infty-estimates} \\ 
\max_{0 \leq j \leq 2} \bigl\|t^{\frac{j+\kappa}{2}} \covd_B^{(j)} \higgsdet\bigr\|_{L_{t,x}^\infty([0, T])} &\lesssim \bigl(1 + \Fnorm{F_B}_{T, 1} \bigr)^2\bigl(T^{\frac{\kappa}{2} - \frac{1}{r}} \big\|\shiggsdet_0\big\|_{L_x^r} + \bigl\|\rhiggsdet_0\bigr\|_{\covHolder{B(0), T}^{-\kappa}}\bigr),\quad \label{eq:higgsdet-L-infty-estimates}\\
\max_{0 \leq j \leq 1} \bigl\|t^{\frac{j}{2}} \covd_B^{(j)}\Sdet\bigr\|_{L_t^\infty L_x^2([0, T])} &\lesssim \gronwall_{T}(C_0 F_B) \bigl(1 + \Fnorm{F_B}_{T, 0} \bigr) \\
& \hspace{-15mm}\times  \bigl( \big\|\bil(\higgsdet \otimes \covd_{B+\Sdet}\higgsdet )\big\|_{L_t^1 L_x^2([0, T])} + \bigl\|t^{1/2}\bil(\higgsdet \otimes \covd_{B+\Sdet}\higgsdet )\bigr\|_{L_{t,x}^2([0, T])}\bigr). \label{eq:Sdet-L2-estimates}
\end{equs}
We remark that unlike the first two estimates, the last estimate only concerns zeroth and first order derivatives.
\end{proposition}

\begin{remark}[Soft smoothness assumption]\label{remark:soft-smoothness-assumption}
In the course of the proof of Proposition~\ref{prop:Sdet-higgsdet-estimates}, we will assume that $\shiggsdet_0, \rhiggsdet_0$ are smooth, and thus also $\higgsdet(0)$ is smooth. This soft smoothness assumption will be convenient when applying monotonicity formulas at various intermediate steps. Since the final estimates only depend on $\|\shiggsdet_0\|_{L_x^r}$, $\bigl\|\rhiggsdet_0\bigr\|_{\covHolder{B(0), T}^{-\kappa}}$, they may be extended to the case of non-smooth initial data by a standard limiting argument.
\end{remark}

The proof of Proposition~\ref{prop:Sdet-higgsdet-estimates} is deferred until later. We proceed to state several more results regarding $\higgsdet$ that will be needed later (e.g.\ in our global iteration scheme \dash Section~\ref{section:global-iteration-monotone-quantities-deterministic}). To state these results, we first further decompose $\reghiggsdet$ (recall \eqref{eq:reghiggsdet-equation}) into $\reghiggsdet = \shiggsdet + \remhiggsdet$, where
\begin{equs}
\label{e:shiggsdet_eq}
(\ptl_t - \Delta_{B+\Sdet}) \shiggsdet &= -(\covd_B^j \Sdet_j)\shiggsdet - |\shiggsdet|^2 \shiggsdet, \\
\shiggsdet(0) &= \shiggsdet_0, \\
\label{e:remhiggsdet_eq}
(\ptl_t - \Delta_{B+\Sdet}) \remhiggsdet &= 2 \Sdet_j \covd_B^j \rlinhiggsdet + \Sdet_j \Sdet^j \rlinhiggsdet - (\covd_B^j \Sdet_j) \remhiggsdet - (|\higgsdet|^2 \higgsdet - |\shiggsdet|^2 \shiggsdet), \\
\remhiggsdet(0) &= 0.
\end{equs}
We record the following identity which will be useful in the coming proofs: by the Bochner formula \eqref{eq:parabolic-bochner}, we have that
\begin{equ}\label{eq:shiggsdet-parabolic-bochner}
\frac{1}{2} (\ptl_t - \Delta) \bigl(|\shiggsdet|^2\bigr) + |\covd_{B + \Sdet} \shiggsdet |^2 =  -\langle \shiggsdet, (\covd_B^j \Sdet_j) \shiggsdet \rangle -  |\shiggsdet|^4 = -  |\shiggsdet|^4,
\end{equ}
where we used that $\langle u, X u \rangle = 0$ for all $u \in \V, X \in \frkg$, since the action of $\frkg$ on
$\V$ is skew-symmetric. 

Our next result gives a lower bound on the decrease of $\shiggsdet$. In Section~\ref{section:global-iteration-monotone-quantities-deterministic}, this decrease will be used to absorb various bad terms in our estimates. (One of the main bad terms is the potential increase in the $L_x^2$ norm of the curvature upon adding $\Sdet$ to $B$ to obtain $B + \Sdet$.)

\begin{proposition}[Decrease estimate for $\shiggsdet$]
\label{pr:shiggsdet_decrease}
Let $T$ be as in Proposition~\ref{prop:Sdet-higgsdet-estimates} and let $\gamma > 0$. For any $t \in [0, T]$, $p \in [2, r]$ and initial data such that $\|\shiggsdet_0\big\|_{L_x^r} \geq 1$, we have that
\begin{equ}
\label{e:shiggsdet_decrease}
\big\|\shiggsdet_0\big\|_{L_x^p}^\gamma - \big\|\shiggsdet(t)\big\|_{L_x^p}^\gamma \gtrsim_\gamma  \big\|\shiggsdet_0\big\|_{L_x^p}^{\gamma + 2} t.
\end{equ}
Moreover, for any $\gamma \geq 4$, we have 
\begin{equ}
\label{e:shiggedet_decrease_L4}
	\big\|\shiggsdet_0\big\|_{L_x^4}^\gamma - \big\|\shiggsdet(t)\big\|_{L_x^4}^\gamma \geq \gamma \|\shiggsdet (t) \|_{L_x^4}^{\gamma - 4} \| |\covd_{B + \Sdet} \shiggsdet | |\shiggsdet| \|_{L^2_{t, x}([0, t])}^2 \;.
\end{equ}
\end{proposition}

One should think of $\shiggsdet$ as the large part of $\higgsdet = \shiggsdet + \rlinhiggsdet + \remhiggsdet$, as $\rlinhiggsdet, \remhiggsdet$ will be much smaller. The next main result shows that $\remhiggsdet$ is indeed a perturbation. (We will take care of $\rlinhiggsdet$ much later, via stochastic estimates \dash see Section~\ref{sec:improved_stochastic_estimates}).

\begin{lemma}[Smallness of $\remhiggsdet$]
\label{le:remhiggsdet_bound}
Let $T \leq 1$. Recall Notation~\ref{notation:gronwall-F-pert}. Let
\begin{equ}
c(T) \eqdef M_T(B, 10) \bigl( T^{\frac{1}{2} - \frac{1}{r}} \big\|\shiggsdet_0\big\|_{L_x^r} + T^{\frac{1-\kappa}{2}} \bigl\|\rhiggsdet_0\bigr\|_{\covHolder{B(0), T}^{-\kappa}}\bigr).
\end{equ}
There is some constant $c_1 \leq 1$ such that if $c(T) \leq c_1$, then 
\begin{equ}
\big\|\remhiggsdet\big\|_{L_{t,x}^\infty([0, T])} + \bigl\|t^{1/2} \covd_B \remhiggsdet\bigr\|_{L_{t,x}^\infty([0, T])} \lesssim c(T)^2 T^{-\frac{\kappa}{2}} \bigl\|\rhiggsdet_0\bigr\|_{\covHolder{B(0), T}^{-\kappa}}.
\end{equ}
\end{lemma}

From the estimate \eqref{eq:Sdet-L2-estimates} in Proposition~\ref{prop:Sdet-higgsdet-estimates}, we see that the $L_x^2$ norms of $\Sdet, \covd_B \Sdet$ are controlled by the sum of the $L_t^1 L_x^2$ norm of $\bil(\higgsdet \otimes \covd_{B+\Sdet} \higgsdet)$ and the corresponding time-weighted $L_{t,x}^2$ norm. The next main result shows that up to smaller errors, these norms are given by the same quantities, except with $\higgsdet$ replaced by the large piece $\shiggsdet$.

\begin{lemma}\label{lemma:bil-covd-higgsdet-higgsdet-replace-shiggsdet}
There is some constant $C$ such that the following holds. Let $T > 0$ satisfy the smallness assumption \eqref{eq:Sdet-higgsdet-main-prop-smallness-assumption} and suppose $K \geq 0$ is such that 
\begin{equ}\label{eq:bil-covd-higgsdet-higgsdet-replace-shiggsdet-remhiggsdet-assumption}
\bigl\|\remhiggsdet\bigr\|_{L_{t,x}^\infty([0, T])} + \bigl\|t^{1/2} \covd_B \remhiggsdet\bigr\|_{L_{t,x}^\infty([0, T])} \leq K.
\end{equ}
Let 
\begin{equ}
K' = C\bigl(1 + \bigl(\Fnorm{F_B}_{T, 1}\bigr)^2\bigr) \bigl(T^{-\frac{\kappa}{2}} \bigl\|\rhiggsdet_0\bigr\|_{\covHolder{B(0), T}^{-\kappa}} + K\bigr) \bigl(T^{\frac{1}{2}-\frac{1}{r}} \bigl\|\shiggsdet_0\bigr\|_{L_x^r} + T^{\frac{1-\kappa}{2}} \bigl\|\rhiggsdet_0\bigr\|_{\covHolder{B(0), T}^{-\kappa}}\bigr).
\end{equ}
Then
\begin{equs}
\bigl\|\bil(\higgsdet \otimes \covd_{B+\Sdet} \higgsdet )\bigr\|_{L_t^1 L_x^2([0, T])} &\leq \bigl\|\bil(\shiggsdet \otimes \covd_{B+\Sdet} \shiggsdet)\bigr\|_{L_t^1 L_x^2([0, T])} + K', \\
\bigl\|t^{1/2}\bil(\higgsdet \otimes \covd_{B+\Sdet} \higgsdet )\bigr\|_{L_{t,x}^2([0, T])} &\leq \bigl\|t^{1/2}\bil(\shiggsdet \otimes \covd_{B+\Sdet} \shiggsdet)\bigr\|_{L_{t,x}^2([0, T])} + K'.
\end{equs}
\end{lemma}

The next main result approximates $\shiggsdet$ by the corresponding linear solution $\slinhiggsdet$, which we define to be the solution to the following equation:
\begin{equ}\label{eq:slinhiggsdet-eq}
(\ptl_t - \Delta_B) \slinhiggsdet = 0, ~~ \slinhiggsdet(0) = \shiggsdet_0.
\end{equ}

\begin{lemma}\label{lemma:shiggsdet-linear-approximation}
Let $C_0, T$ be as in Proposition~\ref{prop:Sdet-higgsdet-estimates}. Then we have that
\begin{equ}
\bigl\|\shiggsdet - \slinhiggsdet\bigr\|_{L_t^\infty L_x^4([0, T])} \lesssim \gronwall_{T}(C_0 F_B) \bigl(1 + \Fnorm{F_B}_{T, 2}\bigr)^{15} \bigl(T^{\frac{1}{2} - \frac{1}{r}} \bigl\|\shiggsdet_0\bigr\|_{L_x^r} + T^{\frac{1-\kappa}{2}} \bigl\|\rhiggsdet_0\bigr\|_{\covHolder{B(0), T}^{-\kappa}}\bigr)^2 \bigl\|\shiggsdet_0\bigr\|_{L_x^4}.
\end{equ}
\end{lemma}

We proceed to prove all the main results that we have stated. We start with Proposition~\ref{prop:Sdet-higgsdet-estimates}. For notational convenience, we make the following definition.

\begin{definition}\label{def:sdet-higgsdet-continuity-norm}
For $t_0 \geq 0$ such that $(\Sdet, \higgsdet)$ exist on $[0, t_0]$, let
\begin{equs}\label{eq:N-t-def}
N_\Sdet(t_0) &\eqdef \max_{0 \leq j \leq 2} \bigl\|t^{\frac{j}{2}} \covd_B^{(j)} \Sdet\bigr\|_{L_{t,x}^\infty([0, t_0])}, \\
N_{\higgsdet}(t_0;\kappa) &\eqdef \max_{0 \leq j \leq 2} \bigl\|t^{\frac{j+\kappa}{2}} \covd_B^{(j)} \higgsdet\bigr\|_{L_{t,x}^\infty([0, t_0])} .
\end{equs}
We will also use $N_{\psi}(t_0; \kappa)$ for other choices of $\psi$, for instance $\psi = \remhiggsdet$ or $\psi = \rlinhiggsdet$.
\end{definition}

By the definition of $\bigl\|\place\bigr\|_{\covHolder{B(0), t_0}^{-\kappa}}$ (Definition~\ref{defn:covariant_Holder}) and the heat kernel estimate \eqref{eq:higher-order-pB-ptwise-bd}, we have that
\begin{equ}\label{eq:rlinhiggsdet-estimate}
N_{\rlinhiggsdet}(t_0; \kappa) \lesssim \bigl(1 + \bigl(\Fnorm{F_B}_{t_0, 1}\bigr)^2\bigr) \bigl\|\rhiggsdet_0\bigr\|_{\covHolder{B(0), t_0}^{-\kappa}}.
\end{equ}
(The heat kernel estimate is used to estimate the second derivatives which appear in Definition~\ref{def:sdet-higgsdet-continuity-norm}.) It also follows from Definition~\ref{def:sdet-higgsdet-continuity-norm} and the above that there is some constant $C > 0$ such that
\begin{equ}\label{eq:N-t-higgsdet-estimate}
N_\higgsdet(t_0; \kappa) \leq N_{\reghiggsdet}(t_0; \kappa) + C \bigl(1 + \bigl(\Fnorm{F_B}_{t_0, 1}\bigr)^2\bigr) \bigl\|\rhiggsdet_0\bigr\|_{\covHolder{B(0), t_0}^{-\kappa}} .
\end{equ}
Next, we obtain the following perturbative estimates on $\Sdet$ (Lemma~\ref{lemma:Sdet-perturbative-estimates}) and $\reghiggsdet$ (Lemma~\ref{lemma:reghiggsdet-perturbative-estimates}), which will serve as the main steps towards Proposition~\ref{prop:Sdet-higgsdet-estimates}. In the following, to help track which terms are the main terms and which terms are perturbations, it may help to think of $t_0^{1/2} N_\Sdet(t_0) \ll 1$ and $t_0^{(1-\kappa)/2} N_{\higgsdet}(t_0; \kappa) \ll 1$, as we will show these estimates in the course of proving the proposition.

\begin{lemma}[Perturbative estimates of $\Sdet$]\label{lemma:Sdet-perturbative-estimates}
Let $t_0 \in (0, 1]$ be such that $(\Sdet, \higgsdet)$ exist on $[0, t_0]$. We have that
\begin{equ}
N_\Sdet(t_0) \lesssim M_{t_0}(B, 2) \Big(t_0^{\frac{1}{2}-\kappa} N_{\higgsdet}(t_0;\kappa)^2 \bigl(1 + t_0^{\frac{1}{2}} N_\Sdet(t_0)\bigr) + t_0^{\frac{1}{2}} N_\Sdet(t_0)^2  + t_0 N_\Sdet(t_0)^3\Big).
\end{equ}
\end{lemma}
\begin{proof}
Let $G = \schO(\Sdet \covd_B \Sdet) + \schO(\Sdet^3) + \bil(\higgsdet \otimes \covd_{B+\Sdet} \higgsdet )$, so that (recall \eqref{eq:Sdet-equation})
\begin{equ}
(\ptl_t - \Delta_B - 2\ad(F_B)) \Sdet = G.
\end{equ}
In the following, all space-time norms are implicitly on $[0, t_0] \times \T^d$, unless explicitly specified. By Lemma~\ref{lemma:higher-regularity-covd-ibp} (with the $k$ therein taken to be $0, 1, 2$), we obtain
\begin{equ}
N_\Sdet(t_0) \lesssim  M_{t_0}(B, 2) N_2(G),
\end{equ}
where
\begin{equ}
N_2(G) = \bigl\| G \bigr\|_{L_t^1 L_x^\infty} + \sum_{j=0}^1 \Bigl(\bigl\|t^{(j+1)/2} \covd_B^{(j)} G \bigr\|_{L_t^2 L_x^\infty} + \sup_{u \in [0, t_0]} \bigl\|t^{(j+1)/2} (u-t)^{-1/2} \covd_B^{(j)} G \bigr\|_{L_t^1 L_x^\infty([0, u])} \Bigr).
\end{equ}
We proceed to estimate the terms on the right-hand side above. We will just give the argument for the term $\bigl\|G\bigr\|_{L_t^1 L_x^\infty}$, as the estimates for the other terms are quite similar. Towards this end, we estimate by Hölder (and inserting $\covd_{B+\Sdet} \higgsdet = \covd_B \higgsdet + \Sdet \higgsdet$)
\begin{equs}
\bigl\|G\bigr\|_{L_t^1 L_x^\infty} &\lesssim \bigl\| \schO(\Sdet \covd_B \Sdet)\bigr\|_{L_t^1 L_x^\infty} + \bigl\|\schO(\Sdet^3)\bigr\|_{L_t^1 L_x^\infty} + \bigl\|\covd_{B+\Sdet} \higgsdet \otimes \higgsdet\bigr\|_{L_t^1 L_x^\infty} \\
&\lesssim \bigl\|t^{-1/2}\bigr\|_{L_t^1} N_\Sdet(t_0)^2  + \bigl\|1\bigr\|_{L_t^1} N_\Sdet(t_0)^3 + \bigl\|\covd_B \higgsdet \otimes \higgsdet\big\|_{L_t^1 L_x^\infty} + \bigl\|\Sdet \higgsdet \otimes \higgsdet\bigr\|_{L_t^1 L_x^\infty} \\
&\lesssim t_0^{1/2} N_\Sdet(t_0)^2 + t_0 N_\Sdet(t_0)^3 + t_0^{\frac{1}{2}-\kappa} N_\higgsdet(t_0; \kappa)^2 + t_0^{1-\kappa} N_\Sdet(t_0) N_\higgsdet(t_0; \kappa)^2,
\end{equs}
which is acceptable. The desired result now follows.
\end{proof}

\begin{lemma}[Perturbative estimates of $\reghiggsdet$]\label{lemma:reghiggsdet-perturbative-estimates}
Let $t_0 \in (0, 1]$ be such that $(\Sdet, \higgsdet)$ exist on $[0, t_0]$. We have that
\begin{equs}
~&N_{\reghiggsdet}(t_0; \kappa) \lesssim \bigl(1 + \Fnorm{F_B}_{t_0, 1} \bigr)^2 \Big( t_0^{\frac{\kappa}{2}-\frac{1}{r}} \big\|\shiggsdet_0\big\|_{L_x^r} + N_{\higgsdet}(t_0;\kappa) t_0^{\frac{1}{2}} N_\Sdet(t_0) \bigl(1 + t_0^{\frac{1}{2}} N_\Sdet(t_0)\bigr) +  t_0^{1-\kappa} N_{\higgsdet}(t_0; \kappa)^3\Big).
\end{equs}
\end{lemma}
\begin{proof}
In the following, all space-time norms are implicitly on $[0, t_0] \times \T^d$.
We rewrite the equation \eqref{eq:reghiggsdet-equation} in the following form (by using \eqref{prelim:eq-linearization-5} to expand $\Delta_{B+\Sdet}$ around $\Delta_B$):
\begin{equ}\label{eq:reghiggsdet-perturbative-estimates-intermediate-equation}
(\ptl_t - \Delta_B) \reghiggsdet = 2 \Sdet_j \covd_B^j \higgsdet + \Sdet^j \Sdet_j \higgsdet - |\higgsdet|^2 \higgsdet \eqdef H, \quad \reghiggsdet(0) = \shiggsdet_0.
\end{equ}
We split $\reghiggsdet = \slinhiggsdet + \tilde{\higgsdet}$, where $\slinhiggsdet$ is defined in \eqref{eq:slinhiggsdet-eq} (i.e.\ it solves the linear equation with initial data $\shiggsdet_0$), and as a consequence the remainder $\tilde{\higgsdet}$ solves
\begin{equ}
(\ptl_t - \Delta_B) \tilde{\higgsdet} = H , \quad \tilde{\higgsdet}(0) = 0.
\end{equ}
Since $\reghiggsdet = \slinhiggsdet + \tilde{\higgsdet}$, it suffices to estimate $N_{\slinhiggsdet}(t_0; \kappa)$ and $N_{\tilde{\higgsdet}}(t_0; \kappa)$. By the kernel estimates \eqref{eq:zeroth-order-pB-ptwise-bd}--\eqref{eq:higher-order-pB-ptwise-bd}, it follows that
\begin{equ}
N_{\slinhiggsdet}(t_0; \kappa) \lesssim \bigl(1 + \Fnorm{F_B}_{t_0, 1}\bigr)^2 t_0^{\frac{\kappa}{2}-\frac{1}{r}} \bigl\|\shiggsdet_0\bigr\|_{L_x^r},
\end{equ}
which is acceptable. Next, by Lemma~\ref{lemma:higher-regularity-covd-ibp} (applied with $k = 0, 1, 2$), we have that for all $t_1 \in (0, t_0]$,
\begin{equ}\label{eq:reghiggsdet-perturbative-estimate-intermediate}
\max_{0 \leq j \leq 2} t_1^{\frac{j}{2}} \bigl\|\covd_B^{(j)} \tilde{\higgsdet}(t_1) \bigr\|_{L_x^\infty} \lesssim \bigl(1 + \Fnorm{F_B}_{t_0, 1}\bigr)^2 N_2(H),
\end{equ}
where $H$ is as defined in \eqref{eq:reghiggsdet-perturbative-estimates-intermediate-equation}, and (in the following, all space-time norms are on $[0, t_1] \times \T^2$)
\begin{equs}
N_2(H) = \big\|H\big\|_{L_t^1 L_x^\infty} + \sum_{j=0}^{1} \Bigl( \big\|t^{(j+1)/2} \covd_B^{(j)} H \big\|_{L_t^2 L_x^\infty} + \sup_{u \in [0, t_1]} \big\|t^{(j+1)/2} (u - t)^{-1/2} \covd_B^{(j)} H \big\|_{L_t^1 L_x^\infty([0, u])} \Bigr).
\end{equs}
We will focus on estimating the $L_t^1 L_x^\infty$ norm, as the estimates for the other norms are similar. Towards this end, we estimate
\begin{equs}
\bigl\|H\bigr\|_{L_t^1 L_x^\infty} &\lesssim \bigl\|\Sdet_j \covd_B^j \higgsdet\bigr\|_{L_t^1 L_x^\infty} + \bigl\|\Sdet^j \Sdet_j \higgsdet\bigr\|_{L_t^1 L_x^\infty} + \bigl\||\higgsdet|^2 \higgsdet\bigr\|_{L_t^1 L_x^\infty}  \\
&\lesssim t_0^{\frac{1-\kappa}{2}} N_\Sdet(t_0) N_{\higgsdet}(t_0; \kappa) + t_0^{1 - \frac{\kappa}{2}} N_\Sdet(t_0)^2 N_{\higgsdet}(t_0; \kappa) + t_0^{1-\frac{3\kappa}{2}} N_{\higgsdet}(t_0; \kappa)^3.  
\end{equs}
Multiplying both sides of \eqref{eq:reghiggsdet-perturbative-estimate-intermediate} by $t_1^{\frac{\kappa}{2}}$ and then taking a sup over $t_1 \in [0, t_0]$, we then obtain an acceptable estimate.
\end{proof}

We are now able to show Proposition~\ref{prop:Sdet-higgsdet-estimates}.

\begin{proof}[Proof of Proposition~\ref{prop:Sdet-higgsdet-estimates}]
First, we prove \eqref{eq:Sdet-L-infty-estimates} and \eqref{eq:higgsdet-L-infty-estimates}. We focus on deriving an a priori bound, as the local well-posedness on the time interval $[0, T]$ may be proven by small adaptations of the same methods. Thus, we assume that $(\Sdet, \higgsdet)$ already exist on $[0, T]$, and are smooth (recall Remark~\ref{remark:soft-smoothness-assumption}). We proceed by a standard continuity argument. We have that $N_\Sdet, N_{\higgsdet}(\place; \kappa), N_{\reghiggsdet}(\place; \kappa)$ are continuous in $t$.

By combining Lemmas~\ref{lemma:Sdet-perturbative-estimates} and~\ref{lemma:reghiggsdet-perturbative-estimates} with \eqref{eq:N-t-higgsdet-estimate} there is some constant $C \geq 1$ such that for $T$ as in the proposition statement and all $t \in [0, T]$, we have that
\begin{equs}
N_\Sdet(t) &\leq C M_T(B, 2) \Big(t^{\frac{1}{2}-\kappa} N_{\higgsdet}(t; \kappa)^2 \bigl(1 + t^{\frac{1}{2}} N_{\Sdet}(t)\bigr) \label{eq:Nt-Sdet-a-priori} + t^{\frac{1}{2}} N_\Sdet(t)^2 + t N_\Sdet(t)^3\Big), \\
N_{\reghiggsdet}(t; \kappa) &\leq C \bigl(1 + \Fnorm{F_B}_{T, 1}\bigr)^2 \Big( t^{\frac{\kappa}{2}-\frac{1}{r}} \big\|\shiggsdet_0\big\|_{L_x^r} \label{eq:Nt-shiggsdet-a-priori} \\
&+ \bigl(C \bigl(1 + \Fnorm{F_B}_{T, 1}\bigr)^2 \bigl\|\rhiggsdet_0\bigr\|_{\covHolder{B(0), T}^{-\kappa}} + N_{\reghiggsdet}(t; \kappa)\bigr) t^{\frac{1}{2}}  N_\Sdet(t) \bigl(1 + t^{\frac{1}{2}} N_\Sdet(t)\bigr)  + t^{1-\kappa} N_{\higgsdet}(t; \kappa)^3\Big), \\
N_{\higgsdet}(t; \kappa) &\leq N_{\reghiggsdet}(t; \kappa) + C \bigl(1 + \bigl(\Fnorm{F_B}_{T, 1}\bigr)^2\bigr)  \bigl\|\rhiggsdet_0\bigr\|_{\covHolder{B(0), T}^{-\kappa}}.
\end{equs}
In the following, let the constant $c_1$ appearing in \eqref{eq:Sdet-higgsdet-main-prop-smallness-assumption} be chosen so that
\begin{equ}\label{eq:small-constant-c1-choice-continuity}
c_1 \eqdef (10 C)^{-100}.
\end{equ}
Now let $I \sse [0, T]$ be the set of all $t \in [0, T]$ such that
\begin{equs}
N_\Sdet(t) &\leq (10C)^{10} M_T(B, 10) \Big(t^{\frac{1}{2}-\frac{2}{r}} \big\|\shiggsdet_0\big\|_{L_x^r}^2 +  t^{\frac{1}{2}-\kappa} \bigl\|\rhiggsdet_0\bigr\|_{\covHolder{B(0), T}^{-\kappa}}^2 \Big), \label{eq:continuity-Sdet-condition}\\
N_{\reghiggsdet}(t; \kappa) &\leq 10 C (1 + \Fnorm{F_B}_{T, 1}\bigr)^2 \bigl(t^{\frac{\kappa}{2} - \frac{1}{r}} \big\|\shiggsdet_0\big\|_{L_x^r} + \bigl\|\rhiggsdet_0\bigr\|_{\covHolder{B(0), T}^{-\kappa}}\bigr). \label{eq:continuity-shiggsdet-condition}
\end{equs}
Note that $0 \in I$, so that $I$ is nonempty. In the following, let $c_2 = (10 C)^{10} M_T(B, 8)^{-1} c_1$. For all $t \in I$, we have by the smallness assumption \eqref{eq:Sdet-higgsdet-main-prop-smallness-assumption} that
\begin{equs}
\bigl(1 + \Fnorm{F_B}_{T, 1}\bigr)^2 t^{\frac{1}{2}} N_\Sdet(t) &\leq (10C)^{10} M_T(B, 12) \Big(\bigl(t^{\frac{1}{2}-\frac{1}{r}} \big\|\shiggsdet_0\big\|_{L_x^r}\bigr)^2 + \bigl(t^{\frac{1-\kappa}{2}} \bigl\|\rhiggsdet_0\bigr\|_{\covHolder{B(0), T}^{-\kappa}}\bigr)^2\Big) \\
&\leq (10 C)^{10} M_T(B, 8)^{-1} c_1^2 \leq c_2, \label{eq:continuity-Sdet-smallness}\\
t^{\frac{1-\kappa}{2}} N_{\reghiggsdet}(t;\kappa) &\leq 10 C M_T(B, 8)^{-1} c_1 \leq c_2, \label{eq:continuity-shiggsdet-smallness} \\
t^{\frac{1-\kappa}{2}} N_{\higgsdet}(t; \kappa) &\leq t^{\frac{1-\kappa}{2}} N_{\reghiggsdet}(t; \kappa) + C \bigl(1 + \bigl(\Fnorm{F_B}_{T, 1}\bigr)^2\bigr) t^{\frac{1-\kappa}{2}} \bigl\|\rhiggsdet_0\bigr\|_{\covHolder{B(0), T}^{-\kappa}} \leq c_2.
\end{equs}
Inserting these estimates as well as \eqref{eq:continuity-Sdet-condition}, \eqref{eq:continuity-shiggsdet-condition} back into \eqref{eq:Nt-Sdet-a-priori}, we obtain
\begin{equs}
~&N_S(t) \leq C M_T(B, 2)  \\
&\times \Big(t^{\frac{1}{2}-\kappa} \Big(11 C (1 + \Fnorm{F_B}_{T, 1}\bigr)^2 \bigl(t^{\frac{\kappa}{2} - \frac{1}{r}} \big\|\shiggsdet_0\big\|_{L_x^r} + \bigl\|\rhiggsdet_0\bigr\|_{\covHolder{B(0), T}^{-\kappa}}\bigr)\Big)^2 (1 + c_2) + c_2 N_S(t) + c_2^2 N_S(t)\Big) \\
&\leq (10C)^5 M_T(B, 6) \Bigl(t^{\frac{1}{2}-\frac{2}{r}} \big\|\shiggsdet_0\big\|_{L_x^r}^2 + t^{\frac{1}{2}-\kappa} \bigl\|\rhiggsdet_0\bigr\|_{\covHolder{B(0), T}^{-\kappa}}^2\Bigr) + \frac{1}{2} N_\Sdet(t),
\end{equs}
where we used the definition of $c_2$ as well as our choice \eqref{eq:small-constant-c1-choice-continuity} of $c_1$ small enough depending on $C$. From this, it follows that for all $t \in I$, we have that
\begin{equ}
N_\Sdet(t) \leq 2 \cdot (10 C)^5 M_T(B, 10) \Big(t^{\frac{1}{2}-\frac{2}{r}} \big\|\shiggsdet_0\big\|_{L_x^r}^2 + t^{\frac{1}{2}-\kappa} \bigl\|\rhiggsdet_0\bigr\|_{\covHolder{B(0), T}^{-\kappa}}^2\Big),
\end{equ}
which is strictly better than \eqref{eq:continuity-Sdet-condition}. To finish, we estimate $N_{\reghiggsdet}$ by inserting our previous estimates back into \eqref{eq:Nt-shiggsdet-a-priori} to obtain
\begin{equs}
N_{\reghiggsdet}(t; \kappa) \leq C& \bigl(1 + \Fnorm{F_B}_{T, 1}\bigr)^2\Big(t^{\frac{\kappa}{2}-\frac{1}{r}} \big\|\shiggsdet_0\big\|_{L_x^r}\\
& + \bigl(C c_2(1 + c_2) + C \bigl(1 + \Fnorm{F_B}_{T, 1}\bigr)^2 c_2^2\bigr) \bigl\|\rhiggsdet_0\bigr\|_{\covHolder{B(0), T}^{-\kappa}} \Big) + \frac{1}{2} N_{\reghiggsdet}(t; \kappa),
\end{equs}
from which we further obtain
\begin{equ}
N_{\reghiggsdet}(t; \kappa) \leq 2 C \bigl(1 + \Fnorm{F_B}_{T, 1}\bigr)^2\Big(t^{\frac{\kappa}{2}-\frac{1}{r}} \big\|\shiggsdet_0\big\|_{L_x^r} + \bigl\|\rhiggsdet_0\bigr\|_{\covHolder{B(0), T}^{-\kappa}}\Big),
\end{equ}
which is strictly better than \eqref{eq:continuity-shiggsdet-condition}. Combining our previous estimates (and using the fact that $N$ is continuous), it follows that $I$ is an open subset of $[0, T]$. Since $I$ is also closed, we conclude that $I = [0, T]$, thus proving \eqref{eq:Sdet-L-infty-estimates} and \eqref{eq:higgsdet-L-infty-estimates} (for the latter, we also use \eqref{eq:N-t-higgsdet-estimate}).

Next, we prove the $L_x^2$ estimate \eqref{eq:Sdet-L2-estimates}. By Corollary~\ref{cor:higher-regularity-with-curvature} (applied with $K \equiv 1$), we have that
\begin{equs}\label{eq:Sdet-L2-estimate-intermediate}
\big\|\Sdet &\big\|_{L_t^\infty L_x^2} + \big\|\covd_B \Sdet\big\|_{L_{t, x}^2} \\
&\lesssim \gronwall_{T} (C_0 F_B) \Big( \big\|\bil(\higgsdet \otimes \covd_{B+\Sdet} \higgsdet) \big\|_{L_t^1 L_x^2} + \big\|\schO(\Sdet\covd_B \Sdet)\big\|_{L_t^1 L_x^2} + \big\|\schO(\Sdet^3)\big\|_{L_t^1 L_x^2}\Big)  \\
&\lesssim \gronwall_{T}(C_0 F_B) \Big(\big\|\bil(\higgsdet \otimes \covd_{B+\Sdet} \higgsdet) \big\|_{L_t^1 L_x^2} + \bigl(T^{1/2} N_S(T) + T N_S(T)^2\bigr) \big\|\Sdet\big\|_{L_t^\infty L_x^2} \Big).
\end{equs}
Upon possibly taking the constant $c_1$ in the smallness assumption \eqref{eq:Sdet-higgsdet-main-prop-smallness-assumption} on $T$ to be even smaller, we may absorb (recall \eqref{eq:continuity-Sdet-smallness}) the second term on the right-hand side back into the left-hand side, and thus we obtain
\begin{equ}\label{eq:Sdet-L2-estimate-regularity-zero-intermediate}
\big\|\Sdet\big\|_{L_t^\infty L_x^2} + \big\|\covd_B \Sdet\big\|_{L_{t, x}^2} \lesssim \gronwall_{T}(C_0 F_B) \big\|\bil(\higgsdet\otimes \covd_{B+\Sdet} \higgsdet ) \big\|_{L_t^1 L_x^2},
\end{equ}
which gives an acceptable estimate for the first term on the left-hand side of \eqref{eq:Sdet-L2-estimates}. To complete the proof of \eqref{eq:Sdet-L2-estimates}, we proceed to estimate $\covd_B \Sdet$. Towards this end, note that
\begin{equ}
(\ptl_t - \Delta_B) \covd_B \Sdet = \schO(\ad(F_B) \covd_B \Sdet) + \covd_B H, 
\end{equ}
where $H = 2\ad(F_B)\Sdet + \schO(\Sdet \covd_B \Sdet) + \schO(\Sdet^3) + \bil(\higgsdet \otimes \covd_{B+\Sdet} \higgsdet)$. Thus by \eqref{eq:higher-regularity-time-weighted-abstract} (with $K \equiv 1$, $\alpha = 1$, $G = 0$) and integration by parts, we have that
\begin{equ}
\big\|t^{1/2} \covd_B \Sdet\big\|_{L_t^\infty L_x^2} + \big\|t^{1/2} \covd_B^{(2)} \Sdet\big\|_{L_{t,x}^2} \lesssim \bigl(1 + \Fnorm{F_B}_{T, 0}\bigr) \big\|\covd_B \Sdet\big\|_{L_{t, x}^2} + \big\|t^{1/2} H \big\|_{L_{t,x}^2}.
\end{equ}
The first term on the right-hand side may be directly estimated by \eqref{eq:Sdet-L2-estimate-regularity-zero-intermediate}. For the second term, we estimate
\begin{equ}
\big\|t^{1/2} H \big\|_{L_{t,x}^2} \lesssim \Big( \Fnorm{F_B}_{T, 0} + T^{1/2} N_\Sdet(T) \bigl(1 + T^{1/2} N_\Sdet(T)\bigr) \Big)\big\|\Sdet\big\|_{L_t^\infty L_x^2} + \bigl\|t^{1/2} \bil(\higgsdet \otimes \covd_{B + \Sdet} \higgsdet) \bigr\|_{L_{t,x}^2}.
\end{equ}
To finish, we combine \eqref{eq:continuity-Sdet-smallness} which gives that $T^{1/2} N_\Sdet(T) \ll 1$ with our previous estimate \eqref{eq:Sdet-L2-estimate-regularity-zero-intermediate}. The desired estimate \eqref{eq:Sdet-L2-estimates} now follows.
\end{proof}

Next, we show Proposition~\ref{pr:shiggsdet_decrease}.

\begin{proof}[Proof of Proposition~\ref{pr:shiggsdet_decrease}]
By \eqref{eq:shiggsdet-parabolic-bochner} and the Leibniz rule, it follows that for some constant $c > 0$, we have that
\begin{equ}
\frac{1}{p} \ptl_t \big\|\shiggsdet\big\|_{L_x^p}^p \leq - c \big\|\shiggsdet\big\|_{L_x^p}^{p+2}.
\end{equ}
In particular, the $L_x^p$ norm of $\shiggsdet$ is non-increasing.  Now fix $t_1 \in [0, T]$. If
\begin{equ}
\big\|\shiggsdet(t_1)\big\|_{L_x^p} \leq \frac{1}{2} \big\|\shiggsdet_0\big\|_{L_x^p},
\end{equ}
then the decrease
\begin{equ}
\big\|\shiggsdet_0\big\|_{L_x^p}^\gamma - \big\|\shiggsdet(t_1)\big\|_{L_x^p}^\gamma \gtrsim \big\|\shiggsdet_0\big\|_{L_x^p}^\gamma \gtrsim \big\|\shiggsdet_0\big\|_{L_x^p}^{\gamma + 2} t_1,
\end{equ}
where the final inequality follows because $t_1 \leq T$, and $T \leq \big\|\shiggsdet_0\big\|_{L_x^r}^{-2} \leq \bigl\|\shiggsdet_0\bigr\|_{L_x^p}^{-2}$ (this follows by the smallness condition \eqref{eq:Sdet-higgsdet-main-prop-smallness-assumption} in Proposition~\ref{prop:Sdet-higgsdet-estimates} and the assumption that $p \leq r$). It remains to show the decrease lower bound in the case
\begin{equ}\label{eq:shiggsdet-decrease-intermediate-assumption}
\big\|\shiggsdet(t_1)\big\|_{L_x^p} \geq \frac{1}{2} \big\|\shiggsdet_0\big\|_{L_x^p}.
\end{equ} 
Assume that this is the case. By the chain rule, we may obtain
\begin{equ}
\frac{1}{\gamma} \ptl_t \big\|\shiggsdet\big\|_{L_x^p}^\gamma \lesssim - \big\|\shiggsdet\big\|_{L_x^p}^{\gamma+2}.
\end{equ}
Integrating this between $t \in [0, t_1]$, using that the $L_x^p$ norm of $\shiggsdet$ is non-increasing, and applying \eqref{eq:shiggsdet-decrease-intermediate-assumption}, we obtain
\begin{equ}
\big\|\shiggsdet_0\big\|_{L_x^p}^\gamma - \big\|\shiggsdet(t_1)\big\|_{L_x^p}^\gamma \gtrsim \big\|\shiggsdet(t_1)\big\|_{L_x^p}^{\gamma + 2} t_1 \gtrsim  \big\|\shiggsdet_0\big\|_{L_x^p}^{\gamma + 2} t_1,
\end{equ}
as desired.

For \eqref{e:shiggedet_decrease_L4}, we use that $\| \shiggsdet(t) \|_{L_x^4}$ is non-increasing in $t$ and the Leibniz rule to obtain
\begin{equ}
	\frac{1}{\gamma} \partial_t \| \shiggsdet \|_{L_x^4}^{\gamma} = \frac{1}{4} \| \shiggsdet \|_{L_x^4}^{\gamma - 4} \bigl( \partial_t \| \shiggsdet \|_{L_x^4}^{4}\bigr) \leq - \| \shiggsdet \|_{L_x^4}^{\gamma - 4} \| |\covd_{B + \Sdet} \shiggsdet | |\shiggsdet| \|_{L_x^2}^2 \;. 
\end{equ}
The bound then follows by integrating in time. 
\end{proof}

Next, we proceed to prove Lemma~\ref{le:remhiggsdet_bound}. Towards this end, we first need the following additional estimates on $\shiggsdet$.

\begin{lemma}[A priori estimates of $\shiggsdet$]\label{lemma:shiggsdet-L-infty}
For any $t_0 \in (0, 1]$ such that $(\Sdet, \shiggsdet)$ exists on $[0, t_0]$, we have that for any $1 \leq p \leq \infty$,
\begin{equ}
\bigl\|\shiggsdet\bigr\|_{L_t^\infty L_x^p([0, t_0])} \leq \bigl\|\shiggsdet_0\bigr\|_{L_x^p}, 
\end{equ}
as well as the estimates
\begin{equ}
\bigl\|t^{1/r} \shiggsdet\bigr\|_{L_{t, x}^\infty([0, t_0])} \lesssim \bigl\|\shiggsdet_0\bigr\|_{L_x^r}, \qquad
\bigl\|\covd_{B+\Sdet}\shiggsdet\bigr\|_{L_{t,x}^2([0, t_0])} \lesssim \bigl\|\shiggsdet_0\bigr\|_{L_x^2}.
\end{equ}
\end{lemma}
\begin{proof}
From \eqref{eq:shiggsdet-parabolic-bochner} and the Leibniz rule, it follows that $|\shiggsdet|$ is a sub-solution to the heat equation, i.e.\ $(\ptl_t - \Delta)|\shiggsdet| \leq 0$. The first and second desired estimates then follow from classical heat equation estimates. The last desired estimate follows directly by integrating \eqref{eq:shiggsdet-parabolic-bochner} in space and time.
\end{proof}

\begin{proof}[Proof of Lemma~\ref{le:remhiggsdet_bound}]
We may assume that $c_1$ is smaller than the constant $c_1$ appearing in Proposition~\ref{prop:Sdet-higgsdet-estimates}, so that the condition $c(T) \leq c_1$ implies that the smallness assumption \eqref{eq:Sdet-higgsdet-main-prop-smallness-assumption} is satisfied. (We will take $c_1$ even smaller later in the proof.) Then by that proposition, we have that for all $t \in [0, T]$,
\begin{equ}
t^{1/2} N_\Sdet(t) \lesssim c(T)^2 \bigl(1 + \Fnorm{F_B}_{T, 1}\bigr)^{-10}. \label{eq:remhiggsdet-Sdet-estimate-intermediate}
\end{equ}
Fix $t_1 \in [0, T]$. In the following, all space-time norms are implicitly on $[0, t_1] \times \T^2$ unless otherwise specified. We may rewrite the equation \eqref{e:remhiggsdet_eq} in the form:
\begin{equs}
(\ptl_t - \Delta_B) \remhiggsdet &= 2 \Sdet_j \covd_B^j (\rlinhiggsdet + \remhiggsdet) + \Sdet_j \Sdet^j (\rlinhiggsdet + \remhiggsdet) - (|\higgsdet|^2 \higgsdet - |\shiggsdet|^2 \shiggsdet) \eqdef G, \\
\remhiggsdet(0) &= 0.
\end{equs}
By Lemma~\ref{lemma:higher-regularity-covd-ibp} (with the $k$ therein taken to be $0, 1$), we have that
\begin{equs}\label{eq:remhiggsdet-perturbative-estimate-intermediate}
~&\bigl\|\remhiggsdet\bigr\|_{L_{t,x}^\infty} + \bigl\|t^{1/2}\covd_B \remhiggsdet\bigr\|_{L_{t,x}^\infty} \\
&\lesssim \bigl(1 + \Fnorm{F_B}_{T, 1}\bigr) \Bigl( \bigl\|G\bigr\|_{L_t^1 L_x^\infty} + \bigl\|t^{1/2} G\bigr\|_{L_t^2 L_x^\infty} + \sup_{u \in [0, t_1]} \bigl\|t^{1/2} (u - t)^{-1/2} G \bigr\|_{L_t^1 L_x^\infty([0, u])}\Bigr).
\end{equs}
We proceed to estimate the $L_t^1 L_x^\infty$ norm of $G$. By \eqref{eq:remhiggsdet-Sdet-estimate-intermediate} and the definition of the covariant Hölder norm (Definition~\ref{defn:covariant_Holder}), we may estimate
\begin{equs}
\big\|\Sdet_j \covd_B^j \rlinhiggsdet\big\|_{L_t^1 L_x^\infty} &\lesssim T^{\frac{1-\kappa}{2}} N_\Sdet(T) \bigl\|\rhiggsdet_0\bigr\|_{\covHolder{B(0), T}^{-\kappa}} \lesssim \bigl(1 + \Fnorm{F_B}_{T, 1}\bigr)^{-10} c(T)^2 T^{-\frac{\kappa}{2}} \bigl\|\rhiggsdet_0\bigr\|_{\covHolder{B(0), T}^{-\kappa}}, \\
\big\|\Sdet_j \Sdet^j \rlinhiggsdet\big\|_{L_t^1 L_x^\infty} &\lesssim T^{1-\frac{\kappa}{2}} N_\Sdet(T)^2 \bigl\|\rhiggsdet_0\bigr\|_{\covHolder{B(0), T}^{-\kappa}} \lesssim \bigl(1 + \Fnorm{F_B}_{T, 1}\bigr)^{-10} c(T)^2 T^{-\frac{\kappa}{2}}\bigl\|\rhiggsdet_0\bigr\|_{\covHolder{B(0), T}^{-\kappa}},
\end{equs}
where in the final inequality we used that $c(T) \leq c_1 \leq 1$, so that $c(T)^4 \leq c(T)^2$. Similarly, we may estimate
\begin{equs}
\bigl\|\Sdet_j \covd_B^j \remhiggsdet\bigr\|_{L_t^1 L_x^\infty} + \bigl\|\Sdet_j \Sdet^j \remhiggsdet\bigr\|_{L_t^1 L_x^\infty} &\lesssim \bigl(T^{1/2} N_\Sdet(T) + \bigl(T^{1/2} N_\Sdet(T)\bigr)^2\bigr) \bigl( \bigl\|\remhiggsdet\bigr\|_{L_{t,x}^\infty} + \bigl\|t^{1/2}\covd_B \remhiggsdet\bigr\|_{L_{t,x}^\infty}\bigr) \\
&\leq \frac{1}{10} \bigl(1 + \Fnorm{F_B}_{T, 1}\bigr)^{-10} \bigl( \bigl\|\remhiggsdet\bigr\|_{L_{t,x}^\infty} + \bigl\|t^{1/2}\covd_B \remhiggsdet\bigr\|_{L_{t,x}^\infty}\bigr),
\end{equs}
where the final inequality follows by \eqref{eq:remhiggsdet-Sdet-estimate-intermediate} and the fact that $c(T) \leq c_1$, where we may take $c_1 \ll 1$. To estimate the final term in $G$, we may first estimate
\begin{equ}
\bigl\| |\remhiggsdet|^2 \remhiggsdet\bigr\|_{L_t^1 L_x^\infty} \leq t_1 \bigl\|\remhiggsdet\bigr\|_{L_{t,x}^\infty}^3,
\end{equ}
and then we expand the polynomial around this cubic term:
\begin{equs}
\abs{|\higgsdet|^2 \higgsdet - |\shiggsdet|^2 \shiggsdet - |\remhiggsdet|^2 \remhiggsdet} \lesssim \bigl(|\shiggsdet|^2 |\rlinhiggsdet| &+ |\shiggsdet| |\rlinhiggsdet|^2 + |\rlinhiggsdet|^3\bigr) \\
&+ \bigl(|\shiggsdet|^2 + |\rlinhiggsdet|^2\bigr) |\remhiggsdet| \\
&+ \bigl(|\shiggsdet| + |\rlinhiggsdet|\bigr) |\remhiggsdet|^2,
\end{equs}
where the terms are grouped by the number of appearances of $\remhiggsdet$. For notational convenience, let
\begin{equ}\label{eq:remhiggsdet-upper-bd-notation}
A = T \Bigl( \bigl(T^{-\frac{1}{r}} \big\|\shiggsdet_0\big\|_{L_x^r}\bigr)^2 T^{-\frac{\kappa}{2}} \bigl\|\rhiggsdet_0\bigr\|_{\covHolder{B(0), T}^{-\kappa}} + T^{-\frac{1}{r}} \big\|\shiggsdet_0\big\|_{L_x^r} \bigl(T^{-\frac{\kappa}{2}} \bigl\|\rhiggsdet_0\bigr\|_{\covHolder{B(0), T}^{-\kappa}}\bigr)^2 + \bigl(T^{-\frac{\kappa}{2}} \bigl\|\rhiggsdet_0\bigr\|_{\covHolder{B(0), T}^{-\kappa}}\bigr)^3\Bigr).
\end{equ}
We may estimate (recall Lemma~\ref{lemma:shiggsdet-L-infty})
\begin{equs}
~&\big\||\shiggsdet|^2 |\rlinhiggsdet| + |\shiggsdet| |\rlinhiggsdet|^2 + |\rlinhiggsdet|^3\big\|_{L_t^1 L_x^\infty} \\
&\lesssim t_1 \Big( \bigl(t_1^{-\frac{1}{r}} \big\|\shiggsdet_0\big\|_{L_x^r}\bigr)^2 t_1^{-\frac{\kappa}{2}} \bigl\|\rhiggsdet_0\bigr\|_{\covHolder{B(0), T}^{-\kappa}} + t_1^{-\frac{1}{r}} \big\|\shiggsdet_0\big\|_{L_x^r} \bigl(t_1^{-\frac{\kappa}{2}} \bigl\|\rhiggsdet_0\bigr\|_{\covHolder{B(0), T}^{-\kappa}}\bigr)^2 + \bigl(t_1^{-\frac{\kappa}{2}} \bigl\|\rhiggsdet_0\bigr\|_{\covHolder{B(0), T}^{-\kappa}}\bigr)^3\Big)\leq A.
\end{equs}
Similarly, we may estimate
\begin{equs}
\big\|\bigl(|\shiggsdet|^2 + |\rlinhiggsdet|^2\bigr) |\remhiggsdet|\big\|_{L_t^1 L_x^\infty} &\lesssim t_1 \Big( \bigl(t_1^{-\frac{1}{r}} \big\|\shiggsdet_0\big\|_{L_x^r}\bigr)^2 + \bigl(t_1^{-\frac{\kappa}{2}} \bigl\|\rhiggsdet_0\bigr\|_{\covHolder{B(0), T}^{-\kappa}}\bigr)^2\Big) \big\|\remhiggsdet\big\|_{L_{t,x}^\infty}, \\
\big\|\bigl(|\shiggsdet| + |\rlinhiggsdet|\bigr) |\remhiggsdet|^2\big\|_{L_t^1 L_x^\infty} &\lesssim t_1 \Big(t_1^{-\frac{1}{r}} \big\|\shiggsdet_0\big\|_{L_x^r} + t_1^{-\frac{\kappa}{2}} \bigl\|\rhiggsdet_0\bigr\|_{\covHolder{B(0), T}^{-\kappa}}\Big) \big\|\remhiggsdet\big\|_{L_{t,x}^\infty}^2.
\end{equs}
Combining, we have estimated the $L_t^1 L_x^\infty$ norm of $G$. The other norms on $G$ appearing in the right-hand side of \eqref{eq:remhiggsdet-perturbative-estimate-intermediate} may be estimated by the same quantities by similar arguments. There is therefore some constant $C$ such that for all $t_1 \in [0, T]$,
\begin{equs}\label{eq:remhiggsdet-intermediate-estimate}
\big\|&\remhiggsdet\big\|_{L_{t,x}^\infty([0, t_1])} + \bigl\|t^{1/2} \covd_B \remhiggsdet\bigr\|_{L_{t, x}^\infty([0, t_1])} \\
&\leq \frac{1}{10}\bigl(\big\|\remhiggsdet\big\|_{L_{t,x}^\infty([0, t_1])} + \bigl\|t^{1/2} \covd_B \remhiggsdet\bigr\|_{L_{t, x}^\infty([0, t_1])}\bigr) + C c(T)^2 T^{-\frac{\kappa}{2}}\bigl\|\rhiggsdet_0\bigr\|_{\covHolder{B(0), T}^{-\kappa}}\\
&+ C \bigl(1 + \Fnorm{F_B}_{T, 1}\bigr) \Bigl( A +   t_1 \Big( \bigl(t_1^{-\frac{1}{r}} \big\|\shiggsdet_0\big\|_{L_x^r}\bigr)^2 + \bigl(t_1^{-\frac{\kappa}{2}} \bigl\|\rhiggsdet_0\bigr\|_{\covHolder{B(0), T}^{-\kappa}}\bigr)^2\Big) \big\|\remhiggsdet\big\|_{L_{t,x}^\infty([0, t_1])} \\
& + t_1 \bigl(t_1^{-\frac{1}{r}} \big\|\shiggsdet_0\big\|_{L_x^r} + t_1^{-\frac{\kappa}{2}} \bigl\|\rhiggsdet_0\bigr\|_{\covHolder{B(0), T}^{-\kappa}}\bigr) \big\|\remhiggsdet\big\|_{L_{t,x}^\infty([0, t_1])}^2 + t_1 \bigl\|\remhiggsdet\bigr\|_{L_{t,x}^\infty([0, t_1])}^3\Bigr).
\end{equs}
To finish, we may proceed by a continuity argument similar to that in the proof of Proposition~\ref{prop:Sdet-higgsdet-estimates}, to obtain that for all $t \in [0, T]$,
\begin{equ}
\big\|\remhiggsdet\big\|_{L_{t,x}^\infty([0, t])} +  \bigl\|t^{1/2} \covd_B \remhiggsdet\bigr\|_{L_{t, x}^\infty([0, t])}  \leq 10 C \Bigl( c(T)^2 T^{-\frac{\kappa}{2}} \bigl\|\rhiggsdet_0\bigr\|_{\covHolder{B(0), T}^{-\kappa}} + \bigl(1 + \Fnorm{F_B}_{T, 1}\bigr) A \Bigr). 
\end{equ}
Finally, we simply note that from the definition of $A$ (see \eqref{eq:remhiggsdet-upper-bd-notation}) and $c(T)$ (in the statement of Lemma~\ref{le:remhiggsdet_bound}), we have that
\begin{equation*}
\bigl(1 + \Fnorm{F_B}_{T, 1}\bigr) A \lesssim c(T)^2 T^{-\frac{\kappa}{2}} \bigl\|\rhiggsdet_0\bigr\|_{\covHolder{B(0), T}^{-\kappa}}. \qedhere
\end{equation*}
\end{proof}

Next, we prove Lemma~\ref{lemma:bil-covd-higgsdet-higgsdet-replace-shiggsdet}, which basically follows from the triangle inequality and our previous estimates.

\begin{proof}[Proof of Lemma~\ref{lemma:bil-covd-higgsdet-higgsdet-replace-shiggsdet}]
We may expand
\begin{equ}
\bil(\higgsdet \otimes \covd_{B+\Sdet} \higgsdet) = \bil(\shiggsdet \otimes \covd_{B+\Sdet} \shiggsdet) + \bil( (\rlinhiggsdet + \remhiggsdet)\otimes \covd_{B+\Sdet} \shiggsdet) + \bil(\higgsdet \otimes \covd_{B+\Sdet} (\rlinhiggsdet + \remhiggsdet)).
\end{equ}
We first estimate the $L_t^1 L_x^2$ norm. By the triangle inequality, it suffices to estimate the $L_t^1 L_x^2$ norms of the last two terms on the right-hand side above. Towards this end, we estimate (applying \eqref{eq:shiggsdet-parabolic-bochner}, the definition of $\covHolder{B(0), T}^{-\kappa}$ (Definition~\ref{defn:covariant_Holder}), and \eqref{eq:bil-covd-higgsdet-higgsdet-replace-shiggsdet-remhiggsdet-assumption})
\begin{equs}
\bigl\|\bil((\rlinhiggsdet + \remhiggsdet)  \otimes \covd_{B+\Sdet} \shiggsdet)\bigr\|_{L_t^1 L_x^2} &\lesssim \bigl\|\covd_{B+\Sdet} \shiggsdet\bigr\|_{L_{t,x}^2} \bigl( \bigl\|\rlinhiggsdet\bigr\|_{L_t^2 L_x^\infty} + \bigl\|\remhiggsdet\bigr\|_{L_t^2 L_x^\infty} \bigr) \\
&\lesssim   \bigl\|\shiggsdet_0\bigr\|_{L_x^2} \bigl(T^{\frac{1-\kappa}{2}}  \bigl\|\rhiggsdet_0\bigr\|_{\covHolder{B(0), T}^{-\kappa}} + T^{1/2} K\bigr),
\end{equs}
which is acceptable. Next, we estimate
\begin{equs}
\bigl\|\bil(\higgsdet & \otimes \covd_{B+\Sdet}(\rlinhiggsdet + \remhiggsdet))\bigr\|_{L_t^1 L_x^2} \\
&\leq \bigl\|\bil( \higgsdet \otimes \covd_{B} (\rlinhiggsdet + \remhiggsdet) )\bigr\|_{L_t^1 L_x^2} + \bigl\|\bil(\higgsdet \otimes \Sdet(\rlinhiggsdet + \remhiggsdet) )\bigr\|_{L_t^1 L_x^2} \\
&\lesssim \Bigl(T^{\frac{1}{2}-\kappa} \bigl\|\rhiggsdet_0\bigr\|_{\covHolder{B(0), T}^{-\kappa}} + T^{\frac{1-\kappa}{2}}K + T^{1-\kappa} N_\Sdet(T)\bigl\|\rhiggsdet_0\bigr\|_{\covHolder{B(0), T}^{-\kappa}}+ T^{1-\frac{\kappa}{2}} N_\Sdet(T) K  \Bigr)\bigl\|t^{\frac{\kappa}{2}} \higgsdet\bigr\|_{L_{t,x}^\infty}.
\end{equs}
To finish, we apply \eqref{eq:higgsdet-L-infty-estimates} and use that Proposition~\ref{prop:Sdet-higgsdet-estimates} implies that $T^{1/2} N_\Sdet(T) \ll 1$. The estimate for the time-weighted $L_{t,x}^2$ norm proceeds similarly; we omit the details.
\end{proof}

Next, we quickly show Lemma~\ref{lemma:shiggsdet-linear-approximation}.

\begin{proof}[Proof of Lemma~\ref{lemma:shiggsdet-linear-approximation}]
From the equations \eqref{e:shiggsdet_eq}, \eqref{eq:slinhiggsdet-eq} for $\shiggsdet$ and $\slinhiggsdet$, we may compute (using \eqref{prelim:eq-linearization-5} to expand $\Delta_{B+\Sdet}$ around $\Delta_B$)
\begin{equ}
(\ptl_t - \Delta_B) (\shiggsdet - \slinhiggsdet) = 2 \covd_B^j (\Sdet_j \shiggsdet) - 2 (\covd_B^j \Sdet_j) \shiggsdet + \Sdet_j \Sdet^j \shiggsdet - |\shiggsdet|^2 \shiggsdet, ~~ (\shiggsdet - \slinhiggsdet)(0) = 0.
\end{equ}
Recalling the definition (Definition~\ref{def:covariant-duhamel}) of $\Duh_B$, it follows that
\begin{equ}
\shiggsdet - \slinhiggsdet = \Duh_B\bigl(2\covd_B^j (\Sdet_j \shiggsdet) - 2 (\covd_B^j \Sdet_j) \shiggsdet + \Sdet_j \Sdet^j \shiggsdet - |\shiggsdet|^2 \shiggsdet\bigr) .
\end{equ}
In the following, all space-time norms are on $[0, T] \times \T^d$ unless otherwise specified. By Lemma~\ref{lemma:covariant-schauder-estimates}, we may estimate
\begin{equs}
\bigl\|\shiggsdet - \slinhiggsdet\bigr\|_{L_t^\infty L_x^4} &\lesssim \bigl(1 + \Fnorm{F_B}_{T, 1}\bigr)\sup_{u \in (0, T]} \bigl\|\min(t, (u-t))^{-1/2} \Sdet \shiggsdet\bigr\|_{L_t^1 L_x^4([0, u])} \\
&+ \bigl\|(\covd_B^j \Sdet_j)\shiggsdet\bigr\|_{L_t^1 L_x^4} + \bigl\|\Sdet_j \Sdet^j \shiggsdet\bigr\|_{L_t^1 L_x^4} + \bigl\||\shiggsdet|^2\shiggsdet\bigr\|_{L_t^1 L_x^4}.
\end{equs}
By Lemma~\ref{lemma:shiggsdet-L-infty}, we may further estimate the right-hand side above by
\begin{equ}
\bigl(1 + \Fnorm{F_B}_{T, 1}\bigr) \bigl( T^{1/2} N_\Sdet(T) + T N_\Sdet(T)^2  + \bigl(T^{\frac{1}{2}-\frac{1}{r}} \bigl\|\shiggsdet_0\bigr\|_{L_x^r}\bigr)^2\bigr) \bigl\|\shiggsdet_0\bigr\|_{L_x^4}.
\end{equ}
By Proposition~\ref{prop:Sdet-higgsdet-estimates}, the terms $T^{1/2} N_\Sdet(T)$ and $T N_\Sdet(T)^2$ may be estimated by an acceptable contribution. The desired result now follows.
\end{proof}

We finish this section off with two results of a more auxiliary nature. The first result is a soft statement that will be needed in the proof of Theorem~\ref{main:thm-covariant-local} at the end of Section~\ref{section:local-proof}.

\begin{lemma}\label{lemma:higgsdet-non-covariant-space-soft}
In the setting of Proposition~\ref{prop:Sdet-higgsdet-estimates}, suppose we further have that $\rhiggsdet_0 \in \Cs_x^{-\kappa}$ and 
\begin{equ}\label{eq:higgsdet-non-covariant-space-soft}
\bigl\|t^{\frac{\kappa}{2}} B \bigr\|_{L_{t,x}^\infty([0, T])} + \bigl\|t^{\kappa} \ptl^\ell B_\ell \bigr\|_{L_{t,x}^\infty([0, T])} < \infty.
\end{equ}
Then $\higgsdet \in C_t^0 \Cs_x^{-\kappa}([0, T])$.
\end{lemma}
\begin{proof}
Recall the equation \eqref{eq:higgsdet-equation} for $\higgsdet$. By expanding the covariant Laplacian $\Delta_{B+\Sdet}$ around the non-covariant Laplacian $\Delta$, we may compute
\begin{equs}
(\ptl_t - \Delta) \higgsdet &= 2 \ptl^j ((B + \Sdet)_j \higgsdet) - (\ptl^j (B + \Sdet)_j) \higgsdet + (B + \Sdet)^j (B + \Sdet)_j \higgsdet - (\covd_B^j \Sdet_j) \higgsdet - |\higgsdet|^2 \higgsdet, \\
\higgsdet(0) &= \shiggsdet_0 +  \rhiggsdet_0.
\end{equs}
Let $G$ be the right-hand side of the $\higgsdet$ equation above. By combining the estimates \eqref{eq:Sdet-L-infty-estimates}, \eqref{eq:higgsdet-L-infty-estimates} with the assumed estimates on $B$, we may obtain by classical Schauder estimates that $\Duh(G) \in C_t^0 \Cs_x^{-\kappa}([0, T])$ (here $\Duh$ is the usual Duhamel operator associated to the heat operator $(\ptl_t - \Delta)$). Since $\higgsdet = e^{t \Delta} \shiggsdet_0 + e^{t \Delta} \rhiggsdet_0 + \Duh(G)$, it remains to show that the linear propagations $e^{t \Delta} \shiggsdet_0, e^{t \Delta} \rhiggsdet_0 \in C_t^0 \Cs_x^{-\kappa}([0, T])$. For $e^{t \Delta} \shiggsdet_0$, this directly follows from the fact that $\shiggsdet_0 \in L_x^r$ and the fact that $\frac{1}{r} \ll \kappa$ (recall \eqref{e:small_constants}), combined with classical heat smoothing estimates. For $e^{t \Delta} \rhiggsdet_0$, this directly follows from classical heat smoothing estimates and the assumption that $\rhiggsdet_0 \in \Cs_x^{-\kappa}$ (for the continuity at $t = 0$, we use that $\Cs_x^{-\kappa}$ is defined as the closure of smooth functions, recall Definition~\ref{def:Holder-space}). 
The desired result now follows.
\end{proof}

Next, using our bounds for $\Sdet$, we obtain the following perturbative estimate on the effect of adding $\Sdet$ into the background connection used in the covariant Hölder norm. This will be needed in the global iteration (Section~\ref{section:decay-estimates-in-expectation}). 
\begin{lemma}
\label{le:B_Sdet_heat_comparison}
	Recall the constants specified in \eqref{prelim:eq-parameters} and \eqref{e:constants_global_iteration}. 
	Let $T_0$ be sufficiently small. Suppose that $T \in (0,T_0]$ and $B \in C_t^0 \Cs_x^{-\kappa}([0, 2T] \times \T^2) \cap C_{t,x}^\infty((0, 2T] \times \T^2)$ is a smooth solution to the Yang--Mills gradient flow at positive times such that
	\begin{equ}
		M_{2T}(B,100) \leq T^{-\lossdel} \;.
	\end{equ}
	Let $\Sdet,\higgsdet$ solve \eqref{eq:Sdet-equation}, \eqref{eq:higgsdet-equation} with initial data satisfying
	\begin{equ}
		\|\shiggsdet_0\|_{L_x^r} \leq T^{-1/\gamma_0} \;, \qquad \|\rhiggsdet_0\|_{\covHolder{B(0),T}^{-\kappa}} \leq T^{-2\fa} \;.
	\end{equ}
	Then, for any $t, t_0 \in (0, T]$ and any distributional one-form $u_0$ we have 
	\begin{equ}
	\label{e:B_Sdet_heat_comparison}
		\|u_0\|_{\covHolder{(B+\Sdet)(t),F, t_0}^{-\kappa}} \lesssim T^{-\lossdel}\|u_0\|_{\covHolder{B(t),F,t_0}^{-\kappa}} \;.
	\end{equ}
	The analogous result also holds for the $\covHolder{(B+\Sdet)(t),t_0}^{-\kappa}$ norm. 
	Moreover, if $T' \in (0,1]$ satisfies
        \begin{equ}
        \label{e:B_Sdet_heat_comparison_long_assumption}
      	  	\gronwall_{T'}\bigl(C_0F_{\Psi_{\place}((B+\Sdet)(t))}\bigr)\bigl(1+\Fnorm{F_{\Psi_{\place}((B+\Sdet)(t))}}_{T',1}\bigr) \leq (T')^{-\lossdel} \;,
        \end{equ}
	where $C_0$ is as in Notation~\ref{notation:gronwall-F-pert}, and $\Psi_{\place}((B+\Sdet)(t))$ denotes the Yang--Mills gradient flow starting from $(B + \Sdet)(t)$, then
	\begin{equs}
	\label{e:B_Sdet_heat_comparison_long}
		\|u_0\|_{\covHolder{(B+\Sdet)(t),F,T'}^{-\kappa}} &\lesssim T^{-\frac{\kappa}{2}-\lossdel}\|u_0\|_{\covHolder{B(t),F,T}^{-\kappa}} \;. 
	\end{equs}
	The analogous result also holds for the $\covHolder{(B+\Sdet)(t),T'}^{-\kappa}$ norm. 
\end{lemma}
\begin{proof}
	First, applying \eqref{eq:Sdet-L-infty-estimates} on $[0,t]$, and using the assumptions on the initial data, gives $\|\Sdet(t)\|_{L^\infty} \lesssim T^{\frac12-\frac2r-\frac2{\gamma_0}-\lossdel}$. 
	Consequently, the heat kernel estimate \eqref{eq:qB-ptwise-bounds} yields
	\begin{equs}
		\|\Sdet(t)\|_{\covHolder{B(t),F,T,2}^{-\kappa}} \lesssim T^{-\lossdel}\|\Sdet(t)\|_{L_x^\infty} 
		\lesssim T^{\frac12-\frac2r-\frac2{\gamma_0}-2\lossdel} \;.
	\end{equs}
	This estimate is enough to treat $\Sdet$ as a valid perturbation in Proposition~\ref{pr:admissible_timescale_general}, where we choose $K = T^{\frac12-\frac2r-\frac2{\gamma_0}-2\lossdel}$, $\YMHFsmooth_0 = B(t)$, and $\YMHFrough_0 = \Sdet(t)$. Note that due to our parameter choices in \eqref{prelim:eq-parameters} and \eqref{e:constants_global_iteration}, we indeed have that $M_T(\Psi_{\place}(B(t)), 2) T^{\frac{1-\kappa}{2}} K \ll 1$, at least if $T_0$ is sufficiently small, so that the smallness assumption \eqref{e:admissible_timescale_smallness} in the proposition is satisfied. Therefore, by Proposition~\ref{pr:admissible_timescale_general} and Lemma~\ref{le:B-YMHFsmooth-heat-comparison} with $\alpha = \kappa $, $p = \infty$, for any $t_0 \in (0, T]$ we have
	\begin{equ}
	\label{e:B_Sdet_heat_comparison_mid}
		\sup_{0<s\leq t_0}s^{\frac{\kappa}{2}}\|e^{s\Delta^F_{(B+\Sdet)(t)}}u_0\|_{L_x^\infty}
		\lesssim \sup_{0<s\leq t_0}s^{\frac{\kappa}{2}}\|e^{s\Delta^F_{B(t)}}u_0\|_{L_x^\infty} \leq \|u_0\|_{\covHolder{B(t),F,t_0}^{-\kappa}} \;.
	\end{equ}
	Here, we are actually applying Lemma~\ref{le:B-YMHFsmooth-heat-comparison} to $g_0 \act u_0$, with background connection $\Psi_{\place}\bigl(g_0 \act (B+\Sdet)(t) \bigr)$, where $g_0$ is the gauge transformation obtained when applying Proposition~\ref{pr:admissible_timescale_general}. Then, we use the fact that 
	\begin{equ}
		\|e^{s\Delta^F_{ g_0 \act (B+\Sdet)(t)}} \bigl( g_0 \act u_0 \bigr) \|_{L_x^\infty} = \|e^{s\Delta^F_{(B+\Sdet)(t)}}u_0\|_{L_x^\infty} 
	\end{equ}
	to get \eqref{e:B_Sdet_heat_comparison_mid}. The desired bounds then follow from \eqref{e:B_Sdet_heat_comparison_mid} and \eqref{eq:qB-ptwise-bounds}. We remark that when proving \eqref{e:B_Sdet_heat_comparison_long}, a direct application would give a factor $(T')^{-\lossdel}$ on the right-hand side, which gives the desired bound if $T' > T$. When $T' \leq T$, \eqref{e:B_Sdet_heat_comparison_long} is the same as \eqref{e:B_Sdet_heat_comparison}. 
\end{proof}

\subsection{Monotone quantity in global iteration}\label{section:global-iteration-monotone-quantities-deterministic}

Recall the energy \eqref{main:eq-covariant-energy} that is used to measure the size of the smooth part of the solution.
In Section~\ref{section:global-theory}, we will perform a careful analysis on the iteration of $\covariantenergy$ to obtain a uniform-in-time bound. As a first step, in this section we prove the following proposition which identifies the terms that dominate $\covariantenergy$ when adding a singular perturbation to the initial data, based on the estimates obtained in Sections~\ref{sec:rough_YMHF} and \ref{sec:higgsdet_estimates}. 

\begin{proposition}
\label{pr:iteration_detarministic_bounds}
	Let $T_0$ be as in \eqref{prelim:eq-T0}, chosen sufficiently small. Recall the constants specified as in \eqref{e:constants_global_iteration}. 
	There exists some constant $C$ such that the following holds. Recall Notation~\ref{notation:gronwall-F-pert} and Definition~\ref{defn:covariant_Holder}.
	Let $B_0^{\sharp} = \YMHFsmooth_0 + \YMHFrough_0$ and $T \in (0, T_0]$ be such that 
\begin{equ}
\label{e:iteration_detarministic_bounds_connection_assumption} 
	M_T(\YMHFsmooth, 100) \leq T^{-\fa} \;, \qquad 
	\bigl\|\YMHFrough_0\bigr\|_{\covHolder{\YMHFsmooth_0, F, T, 2}^{-\kappa}} \leq T^{-2\fa} \;. 
\end{equ}
In the following, we will use $e^{t \Delta^F_{+}}$, $e^{t \Delta_{+}}$ as shorthand for $e^{t \Delta^F_{\smash{\YMHFsmooth_0}}}$, $e^{t \Delta_{\smash{\YMHFsmooth_0}}}$. Let $B, \YMHFsmooth, g$ be given by Proposition~\ref{pr:admissible_timescale_general} (applied with $K$ therein taken to be $T^{-2 \fa}$). Additionally, let $\higgsdet_0 = \shiggsdet_0 + \rhiggsdet_0$ be such that 
	\begin{equ}
	\label{e:phi_0_smallness_global}
		T^{\frac{1}{2}} \big\|\shiggsdet_0\big\|_{L_x^r}  \leq T^{\frac{\gamma_0 - 2}{100}} \;, \quad  \| t^{\frac{\kappa}{2} } e^{\tdet\Delta_{+} } (g_0^{-1} \act \rhiggsdet_0) \|_{L^\infty_t L^4_x([0, T])} \leq T^{-\fa} \;, \quad \bigl\|\rhiggsdet_0\bigr\|_{\covHolder{B(0), T}^{-\kappa}} \leq T^{- 2 \fa} \;.
	\end{equ}
    Let $\Sdet, \higgsdet$ be the solution of \eqref{eq:Sdet-equation}--\eqref{eq:higgsdet-equation}. Then, for any $t \in (0, T]$ we have 
	\minilab{e:iteration_deterministic_bounds}
	\begin{equs}
	\label{e:iteration_deterministic_bounds1}
		\covariantenergy ((B + \Sdet)(t), \higgsdet(t)) 
		&\leq \covariantenergy(\YMHFsmooth_0, \shiggsdet_0) - \frac12 \Dec_{t}(\YMHFsmooth_0, \shiggsdet_0 ) + C t^{-1 + \kappa} \\
	\label{e:iteration_deterministic_bounds2}
		&\quad+  \Big \| \curl_{\YMHFsmooth(t)} e^{t \Delta^F_{+}} \YMHFrough_0 \Big\|_{L^2}^2 \\
	\label{e:iteration_deterministic_bounds3}
		&\quad+ 2 \int_{\T^2} \Big \langle F^{12}_{\YMHFsmooth}(t), \curl_{\YMHFsmooth(t)} e^{t \Delta^F_{+}} \YMHFrough_0 \Big \rangle_{\fg} \, \dx \\
	\label{e:iteration_deterministic_bounds4}
		&\quad + 2 \int_{\T^2} \Big \langle F^{12}_{\YMHFsmooth}(t), \Big[( e^{t \Delta^F_{+}} \YMHFrough_0)^1, ( e^{t \Delta^F_{+}} \YMHFrough_0)^2 \Big]  \Big \rangle_{\fg} \, \dx \\
	\label{e:iteration_deterministic_bounds5}
		&\quad+ C \| \shiggsdet_0 \|_{L^4}^{\gamma_1 - 4} \Big| \int_{\T^2} | e^{t \Delta_{+}} (g_0^{-1} \act \shiggsdet_0) |^2 \\
		& \qquad \qquad \qquad \qquad \qquad \langle e^{t \Delta_{+}} (g_0^{-1} \act \shiggsdet_0), e^{t \Delta_{+}} (g_0^{-1} \act \rhiggsdet_0) \rangle_{\V} \, \dx \Big| \;,
	\end{equs}
	where we used the notation
	\begin{equs}
	\label{e:dec_defn}
		\Dec_{t}(\YMHFsmooth_0, \shiggsdet_0 ) &= 2 \|\covd_{\YMHFsmooth}^* F_{\YMHFsmooth} \|_{L^{2}_{t, x}([0, t])}^2 + C_{\eqref{e:shiggsdet_decrease}, \gamma_2, r} \one_{\| \shiggsdet_0 \|_{L^r} \geq 1 }  \| \shiggsdet_0 \|_{L^r}^{\gamma_2 + 2} t  \\
		&+   \one_{\| \shiggsdet_0 \|_{L^4} \geq 1 }  \Bigl( \bigl( \gamma_1 \|\shiggsdet (t) \|_{L_x^4}^{\gamma_1 - 4} \| |\covd_{B + \Sdet} \shiggsdet | |\shiggsdet| \|_{L^2_{t, x}([0, t])}^2 \bigr) \vee \bigl( C_{\eqref{e:shiggsdet_decrease}, \gamma_1, 4} \| \shiggsdet_0 \|_{L^4}^{\gamma_1 + 2} t \bigr) \Bigr) \;. 
	\end{equs}
	This also gives the following crude bound, which sometimes will be convenient. For $t \in (0, T]$, we have 
	\begin{equ}
	\label{e:iteration_deterministic_bounds_rmk}
		\covariantenergy \bigl( (B + \Sdet)(t), \higgsdet(t) \bigr) 
		\leq \covariantenergy(\YMHFsmooth_0, \shiggsdet_0) + C t^{-\frac74} \;.
	\end{equ}
\end{proposition}
\begin{remark}
	By Propositions~\ref{prop:monotonicity} and~\ref{pr:shiggsdet_decrease} we have 
	\begin{equ}
	\label{e:Dec_lower_bound}
		\Dec_{t}(\YMHFsmooth_0, \shiggsdet_0) \leq  \covariantenergy(\YMHFsmooth(0), \shiggsdet(0)) - \covariantenergy(\YMHFsmooth(t), \shiggsdet(t)) \;. 
	\end{equ}
	This is the reason for using the notation $\Dec$, which stands for ``decrease''. Additionally, recalling the discussion on remnants of exact solvability near the end of Section~\ref{section:main-ideas-brief}, we have that in two dimensions, the covariant divergence $\covd_{\YMHFsmooth}^* F_{\YMHFsmooth}$ controls $\covd_{\YMHFsmooth} F_{\YMHFsmooth}$. This allows us to use the first term in \eqref{e:dec_defn} to absorb terms involving $\covd_{\YMHFsmooth} F_{\YMHFsmooth}$, which we will do without further comment in the following.
\end{remark}

\begin{proof}
    We begin by making several preliminary observations. Note that when \eqref{e:iteration_detarministic_bounds_connection_assumption} holds, by \eqref{e:Fnorm_bound_roughB_general} we have that $\Fnorm{F_B}_{T, 0} \leq \Fnorm{F_{\YMHFsmooth}}_{T, 0} + \bigO(1)$, and thus (in the second inequality we also use Corollary~\ref{cor:curvature-higher-regularity-estimates})
	\begin{equ}
	\label{e:Fpert_B_bounds}
		\gronwall_{T}(C_0 F_B ) \lesssim \gronwall_{T}(C_0 F_{\YMHFsmooth}) \leq T^{-\fa} \;, \quad \Fnorm{F_B}_{T, 2} \lesssim \big(1 + \Fnorm{F_{\YMHFsmooth}}_{T, 0} \big)^{6} \leq T^{-\frac{\fa}{10}} \;.
	\end{equ}
	By Remark~\ref{remark:curvature-pointwise-in-time-bound-Fnorm}, this implies that under \eqref{e:iteration_detarministic_bounds_connection_assumption}, for any $t > 0$ we have 
	\begin{equ}
	\label{e:F_bounds}
		\| F_{\YMHFsmooth}(t) \|_{L^\infty} + \| F_B(t) \|_{L^\infty}  \lesssim t^{-1 - \frac{\fa}{10}} \;. 
	\end{equ}
	Also, combining \eqref{e:Fpert_B_bounds} with \eqref{e:phi_0_smallness_global}, for any $t \leq T$ we have 
	\begin{equ}\label{eq:Sdet-size-estimate}
		M_t(B, 100) \bigl( t^{\frac{1}{2} - \frac{1}{r}} \big\| \shiggsdet_0\big\|_{L_x^r} + t^{\frac{1-\kappa}{2}} \bigl\| \rhiggsdet_0\bigr\|_{\covHolder{B(0), T}^{-\kappa}}\bigr) \lesssim t^{\gamma_1 - 4 } \;.
	\end{equ}
	Here we used that $\gamma_1 - 4 \ll \gamma_0 - 2$.
	In particular, the assumption \eqref{eq:Sdet-higgsdet-main-prop-smallness-assumption} is satisfied when $T$ is small enough, so that $\Sdet, \higgsdet$ are well-defined up to time $T$. 
	
	In the following, recall $\rlinhiggsdet$ (defined by \eqref{eq:rlinhiggsdet-equation}) and $\remhiggsdet$ (defined implicitly via the decomposition $\higgsdet = \shiggsdet + \rlinhiggsdet + \remhiggsdet$, where $\shiggsdet$ is as in \eqref{e:shiggsdet_eq}). By \eqref{e:phi_0_smallness_global} and Lemma~\ref{le:remhiggsdet_bound} we get that for any $t \leq T$, 
	\begin{equ}
		\| \rlinhiggsdet(t) \|_{L^\infty_x} \leq t^{-2 \fa - \frac{\kappa}{2}}  \;, \quad 
		\| \remhiggsdet(t) \|_{L^\infty_{x} } + \| t^{\frac12} \covd_B \remhiggsdet (t) \|_{L^\infty_{x}}  \lesssim t^{\gamma_1 - 4} \;.
	\end{equ}
    We divide the proof of \eqref{e:iteration_deterministic_bounds} into Lemmas~\ref{le:F_BS_expansion}--\ref{le:iteration_deterministic_bounds_mid5} below.
	Based on the discussion above, one can verify all assumptions in Lemmas~\ref{le:F_BS_expansion}--\ref{le:iteration_deterministic_bounds_mid5}.
    
    For \eqref{e:iteration_deterministic_bounds_rmk}, plugging \eqref{e:iteration_detarministic_bounds_connection_assumption} and \eqref{e:phi_0_smallness_global} into the right-hand side of \eqref{e:iteration_deterministic_bounds} and applying \eqref{e:F_bounds}, it follows that for any $t \leq T$,
		\begin{equs}
		\covariantenergy ((B + \Sdet)(t), \higgsdet(t)) &\leq \covariantenergy(\YMHFsmooth_0, \shiggsdet_0)  + C t^{-1 + \kappa} + C \Big( t^{-1 - \fa} T^{- 4 \fa} + t^{-1 - \fa} \cdot t^{-\frac{1+\kappa}{2} } T^{-2\fa} \\
		&\qquad \qquad \qquad \qquad \qquad + t^{-1 - \fa} \cdot t^{ -\kappa } T^{-4\fa} + T^{-\bigl(\frac12-\frac{\gamma_0-2}{100}\bigr)(\gamma_1 - 1)} t^{-\frac{\kappa}{2}} T^{-2\fa} \Big) \\
		&\leq \covariantenergy(\YMHFsmooth_0, \shiggsdet_0) + C t^{-\frac74} \;.
	\end{equs}
	This concludes the proof.
\end{proof}

We record several preliminary estimates which hold in the setting of Proposition~\ref{pr:iteration_detarministic_bounds}, and which will be used several times in the following. Let $T$ be as in the proposition. For any $t \in (0, T]$, by Proposition~\ref{prop:Sdet-higgsdet-estimates} (applied on $[0, t]$), the estimates \eqref{e:Fpert_B_bounds}, and \eqref{e:phi_0_smallness_global}, we have that
	\begin{equs}
	\max_{0 \leq m \leq 1} t^{\frac{m}{2}}\bigl\|\covd_B^{(m)} \Sdet(t) \bigr\|_{L_x^2} &\lesssim t^{-(1.1) \fa} \bigl(\big\|\bil(\higgsdet \otimes \covd_{B+\Sdet}\higgsdet)\big\|_{L^1_{t} L^2_{x}([0, \tdet])} \\
	&\quad \quad \quad \quad \quad + \bigl\|t^{1/2} \bil(\higgsdet \otimes \covd_{B+\Sdet} \higgsdet)\bigr\|_{L_{t,x}^2([0, \tdet])}\bigr) \;, \label{e:iteration_deterministic_bounds_mid1_intermed_Sdet_L2}\\
	\max_{0 \leq m \leq 2} t^{\frac{m}{2}}\bigl\| \covd_B^{(m)}\Sdet(t)\bigr\|_{L_x^\infty} &\lesssim t^{-2\fa} \bigl(t^{- \frac{1}{2} + \frac{\gamma_0-2}{50} - \frac{2}{r}} + t^{\frac{1}{2} - \kappa - 4\fa}\bigr) \label{e:iteration_deterministic_bounds_mid1_intermed_Sdet_Linfty}\;.
	\end{equs}
	Additionally, recall that in the setting of Proposition~\ref{pr:iteration_detarministic_bounds}, we apply Proposition~\ref{pr:admissible_timescale_general} with the choice $K = T^{-2\fa}$. Thus by \eqref{e:Aroughpart_0_smallness_general}, \eqref{e:admissible_timescale_Brem_estimate}, we have that
	\begin{equ}
	\max_{0 \leq m \leq 1} t^{\frac{m}{2}}\bigl\|\covd_{\YMHFsmooth}^{(m)} \YMHFrough(t) \bigr\|_{L_x^\infty} \lesssim t^{-2\fa - \frac{\kappa}{2}} \label{e:iteration_deterministic_bounds_mid1_intermed_rough_part}\;.
	\end{equ} 
We proceed to prove the various intermediate estimates on which the proof of Proposition~\ref{pr:iteration_detarministic_bounds} relies.

\begin{lemma}
\label{le:F_BS_expansion}
	There exists some constant $C$ such that the following holds. In the same setting as Proposition~\ref{pr:iteration_detarministic_bounds}, for any $t \in (0, T]$ we have 
	\begin{equs}
	\label{e:iteration_deterministic_bounds_mid1}
		\| F_{B + \Sdet}&(\tdet) \|_{L^2}^2 -  \| F_{\YMHFsmooth(\tdet)} \|_{L^2}^2  \leq \frac{1}{10}\|\covd_{\YMHFsmooth} F_{\YMHFsmooth}\|_{L^2_{t, x}([0, \tdet])}^2 + C \tdet^{-1+ \kappa}  \\
		&\quad+ C \tdet^{ -1 - 10 \fa } \bigl(\big\|\bil(\higgsdet \otimes \covd_{B+\Sdet}\higgsdet)\big\|_{L^1_{t} L^2_{x}([0, \tdet])}^{2} + \bigl\|t^{1/2} \bil(\higgsdet \otimes \covd_{B+\Sdet} \higgsdet)\bigr\|_{L_{t,x}^2([0, \tdet])}^2\bigr) \\
		&\quad+  \Bigl \| \curl_{\YMHFsmooth(\tdet)} e^{\tdet \Delta^F_{+}} \YMHFrough_0 \Bigr\|_{L^2}^2 \\
		&\quad+ 2 \int_{\T^2} \Big \langle F^{12}_{\YMHFsmooth}(\tdet), \curl_{\YMHFsmooth(\tdet)} e^{\tdet \Delta^F_{+}} \YMHFrough_0 \Big \rangle_{\fg} \, \dx \\
		&\quad + 2 \int_{\T^2} \Big \langle F^{12}_{\YMHFsmooth}(\tdet), \Big[( e^{\tdet \Delta^F_{+}} \YMHFrough_0)^1, ( e^{\tdet \Delta^F_{+}} \YMHFrough_0)^2\Big]  \Big \rangle_{\fg} \, \dx \;. 
	\end{equs}
\end{lemma}
\begin{remark}
	The right-hand side is already in the desired form, except for the term on the second line. This term will be absorbed by the dissipation coming from the decrease of the $L^4$ norm of $\higgsdet$ (see Lemma~\ref{le:iteration_deterministic_bounds_mid4}). This is precisely the reason for including the $L^4$ norm of $\higgsdet$ in the energy.
\end{remark}

\begin{proof}
	We first expand 
	\begin{equ}
	\label{e:F_BS_expansion}
		\|F_{B+\Sdet}\|_{L^2}^2 = \|F_B\|_{L^2}^2 + 2 \, \langle F^{12}_B, \curl_{B} \Sdet \rangle + 2 \langle F^{12}_B, [\Sdet^1, \Sdet^2] \rangle + \lVert \curl_{B} \Sdet + [\Sdet^1, \Sdet^2] \rVert_{L^2}^2 \;.
	\end{equ}	
	For the first term in \eqref{e:F_BS_expansion}, by \eqref{e:prop_perturbed_curvature_average_general}, we have (here we also apply \eqref{e:F_bounds} to estimate the final term in \eqref{e:prop_perturbed_curvature_average_general})
	\begin{equs}
		\|F_B(\tdet)\|_{L^2}^2 &= \|F_{\YMHFsmooth}(\tdet)\|_{L^2}^2 +  \Bigl \| \curl_{\YMHFsmooth(\tdet)} e^{\tdet \Delta^F_{+}} \YMHFrough_0 \Bigr\|_{L^2}^2 \\
		&\quad+ 2 \, \int_{\T^2} \Big \langle F^{12}_{\YMHFsmooth}(\tdet), \curl_{\YMHFsmooth(\tdet)} e^{\tdet \Delta^F_{+}} \YMHFrough_0 + \Big[( e^{\tdet \Delta^F_{+}} \YMHFrough_0)^1, ( e^{\tdet \Delta^F_{+}} \YMHFrough_0)^2 \Big]  \Big \rangle_{\fg} \, \dx  \\
		& \quad + \bigO \Bigl( t^{-10 \fa - \kappa} \bigl\| \covd_{\YMHFsmooth} F_{\YMHFsmooth} \bigr\|_{L_{t,x}^2([0, t])} + t^{-1 + \kappa} \Bigr) \;.
	\end{equs} 
	Applying Young's inequality to the term $t^{-10 \fa - \kappa} \bigl\| \covd_{\YMHFsmooth} F_{\YMHFsmooth} \bigr\|_{L_{t,x}^2([0, t])} $ gives the desired bound. 
	
	For the second term in \eqref{e:F_BS_expansion}, first recall that in Proposition~\ref{pr:admissible_timescale_general}, $g$ and $\YMHFrough$ are defined such that $g^{-1} \act B = \YMHFsmooth + \YMHFrough$. Then we have
	\begin{equs}
		\langle F^{12}_B, \curl_{B} \Sdet \rangle &= \Big \langle \Ad_{g^{-1}} F^{12}_B, \Ad_{g^{-1}} \curl_{B} \Sdet \Big \rangle \\
		\label{e:iteration_deterministic_bounds_mid1_mid1}
		&= \Big \langle F^{12}_{\YMHFsmooth}, \curl_{\YMHFsmooth + \YMHFrough} \Ad_{g^{-1}} \Sdet \Big \rangle  \\
		\label{e:iteration_deterministic_bounds_mid1_mid2}
		&\quad + \Big \langle \curl_{\YMHFsmooth} \YMHFrough + [\YMHFrough_1, \YMHFrough_2], \Ad_{g^{-1}} \curl_B \Sdet \Big \rangle \;. 
	\end{equs}
	For notational brevity, let
	\begin{equ}
	J = \big\|\bil(\higgsdet \otimes \covd_{B+\Sdet}\higgsdet)\big\|_{L^1_{t} L^2_{x}([0, \tdet])} + \bigl\|t^{1/2} \bil(\higgsdet \otimes \covd_{B+\Sdet} \higgsdet)\bigr\|_{L_{t,x}^2([0, \tdet])} \;.
	\end{equ}
	Now by integration by parts, Lemma~\ref{le:DF_bound}, \eqref{e:Fpert_B_bounds}, \eqref{e:iteration_deterministic_bounds_mid1_intermed_Sdet_L2}, and \eqref{e:iteration_deterministic_bounds_mid1_intermed_rough_part}, the term
	\eqref{e:iteration_deterministic_bounds_mid1_mid1} is bounded (up to constants) by
	\begin{equs}
		\| \covd_{\YMHFsmooth} F_{\YMHFsmooth} (t) \|_{L^2} &\| \Sdet (t) \|_{L^2} + \|F_{\YMHFsmooth} (t) \|_{L^2} \|\YMHFrough (t) \|_{L^\infty} \| \Sdet (t) \|_{L^2} \\
		&\lesssim \tdet^{-\frac12 - \fa}  \|\covd_{\YMHFsmooth} F_{\YMHFsmooth} \|_{L^2_{t, x}([0, \tdet])} \cdot \tdet^{-(1.1) \fa} J  + \tdet^{-1 - \frac{\fa}{10}} \cdot t^{-2\fa} \tdet^{-\frac{\kappa}{2}} \cdot \tdet^{-(1.1)\fa} J \;,
	\end{equs}
	which gives the desired bound by Young's inequality. For \eqref{e:iteration_deterministic_bounds_mid1_mid2}, we apply \eqref{e:iteration_deterministic_bounds_mid1_intermed_Sdet_L2} and \eqref{e:iteration_deterministic_bounds_mid1_intermed_rough_part} to estimate this term by
	\begin{equs}
		\bigl( t^{-2\fa} \tdet^{-\frac12 - \frac{\kappa}{2}} + t^{-4\fa} \tdet^{-\kappa} \bigr) \tdet^{-(1.1)\fa - \frac12} J \;,
	\end{equs}
	which gives the desired bound by Young's inequality. 
	
	For the third term in \eqref{e:F_BS_expansion}, we use H\"older's inequality, \eqref{e:iteration_deterministic_bounds_mid1_intermed_Sdet_L2}, \eqref{e:iteration_deterministic_bounds_mid1_intermed_Sdet_Linfty}, and \eqref{e:F_bounds} to estimate
\begin{equ}	
	\langle F^{12}_B(t), [\Sdet^1(t), \Sdet^2(t)] \rangle
	\leq 2 \| F_B(t) \|_{L^\infty_x} \|\Sdet(t) \|_{L^2_x}^2
	\leq 2 t^{-1 - \frac{\fa}{10}} \cdot t^{-(2.2)\fa} J^2 \;,
\end{equ}
	which gives the desired bound. 
	
	For the last term in \eqref{e:F_BS_expansion}, we apply \eqref{e:iteration_deterministic_bounds_mid1_intermed_Sdet_L2} and \eqref{e:iteration_deterministic_bounds_mid1_intermed_Sdet_Linfty} to estimate
	\begin{equs}
		~& 2\|\covd_B \Sdet(t) \|_{L^2}^2 + 2 \| \Sdet(t) \|_{L^\infty}^2 \| \Sdet(t) \|_{L^2}^2
		\\
		&\lesssim \Big( \tdet^{-(1.1)\fa - \frac12} J  + \bigl(\tdet^{-2 \fa} \tdet^{- \frac{1}{2} + \frac{\gamma_0 - 2}{50} - \frac{2}{r}} + \tdet^{-2\fa} \tdet^{\frac12 - \kappa} \tdet^{-4\fa} \bigr) \cdot \tdet^{-(1.1)\fa} J \Big)^2 \lesssim \tdet^{-1 - 8\fa} J^2 \;.
	\end{equs}
	The desired estimate follows. 
\end{proof}

\begin{lemma}
\label{le:iteration_deterministic_bounds_mid2}
	There exists some constant $C$ such that the following holds. Let $T \in (0, 1]$, let $B \in C_t^0 \Cs_x^{-\kappa}([0, T] \times \T^2) \cap C_{t,x}^\infty((0, T] \times \T^2)$ be a smooth solution to the Yang--Mills gradient flow at positive times, and suppose that $T$ satisfies the smallness assumption \eqref{eq:Sdet-higgsdet-main-prop-smallness-assumption} of Proposition~\ref{prop:Sdet-higgsdet-estimates}. Let $\Sdet, \higgsdet$ be the solution to \eqref{eq:Sdet-equation}--\eqref{eq:higgsdet-equation} (which is provided by Proposition~\ref{prop:Sdet-higgsdet-estimates}). Let $\shiggsdet$ solve \eqref{e:shiggsdet_eq}. If $\higgsdet$ has decomposition $\higgsdet = \shiggsdet + \rhiggsdet = \shiggsdet + \rlinhiggsdet + \remhiggsdet$ with bounds 
	\begin{equ}
	\label{e:phi_different_parts_bounds}
		\| t^{\frac{\kappa}{2} + 2 \fa} \rlinhiggsdet(t) \|_{L^\infty_{t, x}([0, T])} \leq 1  \;, \quad 
		\| t^{-(\gamma_1 - 4)} \remhiggsdet(t) \|_{L^\infty_{t, x}([0, T])} \lesssim 1 \;,
	\end{equ}
	then for any $t \in (0, T]$ we have
	\begin{equs}
	\label{e:iteration_deterministic_bounds_mid2}
		\| \higgsdet(\tdet) \|_{L^4}^{\gamma_1} &\leq \| \shiggsdet_0 \|_{L^4}^{\gamma_1} - \frac34 \bigl(\| \shiggsdet_0 \|_{L^4}^{\gamma_1} - \| \shiggsdet(\tdet) \|_{L^4}^{\gamma_1} \bigr) + C \tdet^{-1 + \kappa} \\
	\label{e:iteration_deterministic_bounds_mid22}
		&\quad + C \| \shiggsdet_0 \|_{L^4}^{\gamma_1 - 4} \Big| \int_{\T^2} |\shiggsdet(\tdet)|^2 \langle \shiggsdet(\tdet), \rlinhiggsdet(\tdet) \rangle_{\V} \, \dx \Big|  \;.
	\end{equs}
\end{lemma}
\begin{remark}
	The term \eqref{e:iteration_deterministic_bounds_mid22} will be bounded in Lemma~\ref{le:iteration_deterministic_bounds_mid5}. 
\end{remark}

\begin{proof}
	First we expand the $L^4$ norm at time $\tdet$ to get (here we apply \eqref{e:phi_different_parts_bounds})
	\begin{equ}
		\| \shiggsdet + \higgsroughpart \|_{L^4}^4 \leq \| \shiggsdet \|_{L^4}^4 + 4 \int_{\T^2} |\shiggsdet|^2 \langle \shiggsdet, \rhiggsdet \rangle_{\V} \, \dx + C \bigl(1 + \| \shiggsdet \|_{L^4}^2 \bigr) \cdot t^{-8 \fa} \tdet^{-2\kappa}  \;.
	\end{equ}
	If $\| \shiggsdet(\tdet) \|_{L^4} \leq C_\star \tdet^{-\frac15 }$ with any large constant $C_\star$ (to be chosen later) then the bound \eqref{e:iteration_deterministic_bounds_mid2} directly follows from this expansion and the fact that $\gamma_1 < 5$. From now on we assume that the opposite holds. 
	By choosing $C_\star$ large enough, we would get 
	\begin{equ}
		\| \shiggsdet(\tdet) \|_{L^4} \geq C_\star \tdet^{-\frac15} \geq  \| \rhiggsdet(\tdet) \|_{L^4} \;.
	\end{equ}
	Since from this we get $\|\shiggsdet(0) \|_{L^4} \geq \| \shiggsdet(\tdet) \|_{L^4} \geq 1$, by Proposition~\ref{pr:shiggsdet_decrease} we get
	\begin{equ}
	\label{e:Dec_2_lower_bound}
		\| \shiggsdet_0 \|_{L^4}^{\gamma_1} - \| \shiggsdet(\tdet) \|_{L^4}^{\gamma_1} \geq  C_{\eqref{e:shiggsdet_decrease}, \gamma_1, 4} \| \shiggsdet_0 \|_{L^4}^{\gamma_1 + 2} \tdet \;.
	\end{equ}
	For $\gamma_1 \geq 4$, using that $|(1 + x)^{\gamma_1/4} - 1| \leq C_{\gamma_1} |x|$ for $|x| \leq 1$, we get at time $t$, 
	\begin{equs}
		\,& \| \shiggsdet + \higgsroughpart \|_{L^4}^{\gamma_1} \\
		&= \| \shiggsdet \|_{L^4}^{\gamma_1} \Bigl(1 +  \| \shiggsdet \|_{L^4}^{-4} \Bigl(4 \int_{\T^2} |\shiggsdet|^2 \langle \shiggsdet, \rhiggsdet \rangle_{\V} \, \dx + \bigO \big(\| \shiggsdet \|_{L^4}^{2} \cdot t^{-8\fa} \tdet^{-2\kappa} \big) \Bigr) \Bigr)^{\frac{\gamma_1}{4}} \\
		&\leq \| \shiggsdet \|_{L^4}^{\gamma_1} + C \| \shiggsdet_0 \|_{L^4}^{\gamma_1 - 4} \Big| \int_{\T^2} |\shiggsdet|^2 \langle \shiggsdet, \rlinhiggsdet \rangle_{\V} \, \dx \Big| \\
		\label{e:L4_expansion_mid1}
		&\quad+ C \| \shiggsdet \|_{L^4}^{\gamma_1 - 4} \Big| \int_{\T^2} |\shiggsdet|^2 \langle \shiggsdet, \remhiggsdet \rangle_{\V} \, \dx \Big| \\
		\label{e:L4_expansion_mid2}
		&\quad+ C \| \shiggsdet \|_{L^4}^{\gamma_1 - 2} \tdet^{-8\fa} \cdot \tdet^{- 2\kappa} \;,
	\end{equs} 
	where we chose $C_\star$ large enough so that the remainder term in the Taylor expansion is indeed smaller than $1$.
	Now it suffices to bound \eqref{e:L4_expansion_mid1} $+$ \eqref{e:L4_expansion_mid2} by $\frac14 \bigl(\| \shiggsdet_0 \|_{L^4}^{\gamma_1} - \| \shiggsdet(\tdet) \|_{L^4}^{\gamma_1} \bigr) + C\tdet^{-1 + \kappa}$. 
	For \eqref{e:L4_expansion_mid1}, by \eqref{e:phi_different_parts_bounds}, \eqref{e:Dec_2_lower_bound} and Young's inequality, we have
	\begin{equs}
		 \| \shiggsdet_0 \|_{L^4}^{\gamma_1 - 4} \Big| \int_{\T^2} |\shiggsdet(t) |^2 \langle \shiggsdet(t), \remhiggsdet(t) \rangle_{\V} \, \dx \Big| &\lesssim \| \shiggsdet_0 \|_{L^4}^{\gamma_1 - 1} t^{\gamma_1 - 4 }\\
		 &= \bigl( \| \shiggsdet_0 \|_{L^4}^{\gamma_1 + 2} \tdet \bigr)^{\frac{\gamma_1 - 1}{\gamma_1 + 2} } \tdet^{\gamma_1 - 4 - \frac{\gamma_1 - 1}{\gamma_1 + 2} } \\ 
		 \label{e:young_shiggsdet_decrease_L4}
		 & \leq \frac18 C_{\eqref{e:shiggsdet_decrease}, \gamma_1, 4} \| \shiggsdet_0 \|_{L^4}^{\gamma_1 + 2} \tdet + C \tdet^{(\gamma_1 - 4 - \frac{\gamma_1 - 1}{\gamma_1 + 2} ) \frac{\gamma_1 + 2}{3} }\;. \quad 
	\end{equs}
	The exponent of $t$ is 
	\begin{equ}	
		- \frac{\gamma_1 - 1}{3} + (\gamma_1 - 4) \frac{\gamma_1 + 2}{3} = -1 - \frac{\gamma_1 - 4}{3} + (\gamma_1 - 4) \frac{\gamma_1 + 2}{3} > - 1 + \kappa \;,
	\end{equ}
	so this term is acceptable. 
	For \eqref{e:L4_expansion_mid2}, we proceed similarly and get
	\begin{equs}
	\label{e:L4_expansion_mid2_Young}
		\| \shiggsdet_0 \|_{L^4}^{\gamma_1 - 2} \tdet^{- 2\kappa - 8 \fa} \leq \frac18 C_{\eqref{e:shiggsdet_decrease}, \gamma_1, 4} \| \shiggsdet_0 \|_{L^4}^{\gamma_1 + 2} \tdet + C \tdet^{(-2\kappa - 8 \fa - \frac{\gamma_1 - 2}{\gamma_1 + 2}) \frac{\gamma_1 + 2}{4} } \;,
	\end{equs}
	which is acceptable since $\gamma_1 \in (4, 5)$.  
\end{proof}

\begin{lemma}
\label{le:iteration_deterministic_bounds_mid3}
	In the same setting as Lemma~\ref{le:iteration_deterministic_bounds_mid2}, for any $t \in (0, T]$ we have 
	\begin{equ}
	\label{e:iteration_deterministic_bounds_mid3}
		\| \higgsdet(\tdet) \|_{L^r}^{\gamma_2} \leq \| \shiggsdet_0 \|_{L^r}^{\gamma_2} - \frac34 \bigl( \| \shiggsdet_0 \|_{L^r}^{\gamma_2} - \| \shiggsdet(\tdet) \|_{L^r}^{\gamma_2} \bigr) + C \tdet^{-1 + \kappa} \;.
	\end{equ}
\end{lemma}
\begin{proof}
	If $\| \shiggsdet(t) \|_{L^r} \leq 100 \tdet^{-\kappa/2} \tdet^{-2 \fa}$ then the bound \eqref{e:iteration_deterministic_bounds_mid3} immediately holds by the triangle inequality and \eqref{e:phi_different_parts_bounds}. Otherwise by \eqref{e:phi_different_parts_bounds} we have 
	\begin{equ}
		\| \shiggsdet(t) \|_{L^r} \geq 100 \tdet^{-\kappa/2} \tdet^{- 2 \fa} \geq \|\rhiggsdet(t) \|_{L^r} \;. 
	\end{equ}
	Also from this we get $\|\shiggsdet (0) \|_{L^r} \geq \| \shiggsdet (\tdet)\|_{L^r} \geq 1$, and thus
	by Proposition~\ref{pr:shiggsdet_decrease} we have
	\begin{equ}
		 \| \shiggsdet_0 \|_{L^r}^{\gamma_2} - \| \shiggsdet(\tdet) \|_{L^r}^{\gamma_2} \geq  C_{\eqref{e:shiggsdet_decrease}, \gamma_2, r} \| \shiggsdet_0 \|_{L^r}^{\gamma_2 + 2} \tdet \;.
	\end{equ}
	Then, we use the triangle inequality to get
	\begin{equs}
		\| \higgsdet(t) \|_{L^r}^{\gamma_2} &\leq \| \shiggsdet(t) \|_{L^r}^{\gamma_2} \Big(1 + \frac{\|\higgsroughpart(t) \|_{L^r}}{\|\shiggsdet(t) \|_{L^r}} \Big)^{\gamma_2} \\
		&\leq \| \shiggsdet(t) \|_{L^r}^{\gamma_2} + C \| \shiggsdet(t) \|_{L^r}^{\gamma_2 - 1} \|\higgsroughpart(t) \|_{L^r}  \\
		& \leq \| \shiggsdet_0 \|_{L^r}^{\gamma_2} -\bigl( \| \shiggsdet_0 \|_{L^r}^{\gamma_2} - \| \shiggsdet(\tdet) \|_{L^r}^{\gamma_2} \bigr)  + C \| \shiggsdet_0 \|_{L^r}^{\gamma_2 - 1} \tdet^{-\kappa/2 - 2 \fa} \;. 
	\end{equs}
	The last term is treated using Young's inequality as follows: 
	\begin{equ}
		C\| \shiggsdet_0 \|_{L^r}^{\gamma_2 - 1} \tdet^{-\kappa/2 - 2\fa} \leq \frac14 C_{\eqref{e:shiggsdet_decrease}, \gamma_2, r} \| \shiggsdet_0 \|_{L^r}^{\gamma_2 + 2} \tdet + C' \tdet^{(-\frac{\kappa}{2} - 2 \fa - \frac{\gamma_2 - 1}{\gamma_2 + 2}) \frac{\gamma_2 + 2}{3} } \;,
	\end{equ}
	From the choice of our $\gamma_2$ just after \eqref{prelim:eq-T0}, the exponent in $\tdet$ is better than $\tdet^{-1 + \kappa}$. 
\end{proof}

\begin{lemma}
\label{le:iteration_deterministic_bounds_mid4}
	Let $T_0$ be as in Proposition~\ref{pr:iteration_detarministic_bounds}. For any $\eps > 0$, there exists some constant $C = C_{\eps}$ such that the following holds. 
	Let $B$ be a YM gradient flow on $(0, T]$ with $T \in (0, T_0]$ such that $M_T(B, 10) \leq T^{-5\fa} $.
	Let $\higgsdet_0 = \shiggsdet_0 + \rhiggsdet_0$ be such that the first and third estimates in \eqref{e:phi_0_smallness_global} hold. Let $\Sdet, \higgsdet$ be the solution of \eqref{eq:Sdet-equation}--\eqref{eq:higgsdet-equation} and let $\shiggsdet$ solve \eqref{e:shiggsdet_eq}. If $\higgsdet$ has the decomposition $\higgsdet = \shiggsdet + \rlinhiggsdet + \remhiggsdet$ with bounds 
	\begin{equ}
		\| t^{\frac{\kappa}{2} + 2 \fa} \rlinhiggsdet(t) \|_{L^\infty_{t, x}([0, T])} \leq 1  \;, \quad 
		\| t^{-(\gamma_1 - 4)} \remhiggsdet(t) \|_{L^\infty_{t, x}([0, T])} + \| t^{-(\gamma_1 - 4)} t^{\frac12} \covd_B \remhiggsdet(t) \|_{L^\infty_{t, x}([0, T])} \lesssim 1 \;,
	\end{equ}
	then for any $t \in (0, T]$ we have 
	\begin{equs}
	\label{e:iteration_deterministic_bounds_mid4}
		\tdet^{ -1 - 10 \fa } \bigl(\big\|\bil(\higgsdet \otimes \covd_{B+\Sdet}\higgsdet)\big\|_{L^1_t L^2_{x}([0, \tdet])}^2 &+  \big\|t^{1/2}\bil(\higgsdet \otimes \covd_{B+\Sdet}\higgsdet)\big\|_{L_{t,x}^2([0, \tdet])}^2\bigr) \\
		&\leq \eps \bigl( \| \shiggsdet_0 \|_{L^4}^{\gamma_1} - \| \shiggsdet(\tdet) \|_{L^4}^{\gamma_1} \bigr) + C \tdet^{-1 + \kappa} \;.
	\end{equs}
\end{lemma}
\begin{proof}
	By Lemma~\ref{lemma:bil-covd-higgsdet-higgsdet-replace-shiggsdet} and the last estimate of \eqref{e:phi_0_smallness_global}, we have
	\begin{equs}
		\big\|\bil(\higgsdet \otimes \covd_{B+\Sdet}\higgsdet)\big\|_{L^1_t L^2_{x}([0, \tdet])}^2 &\leq 2 \big\|\bil(\shiggsdet \otimes \covd_{B+\Sdet}\shiggsdet)\big\|_{L^1_t L^2_{x}([0, \tdet])}^2 + C t^{\gamma_1 - 4} \;,\\
		\big\|t^{1/2}\bil(\higgsdet \otimes \covd_{B+\Sdet}\higgsdet)\big\|_{L_{t,x}^2([0, \tdet])}^2 &\leq 2 \big\|t^{1/2}\bil(\shiggsdet \otimes \covd_{B+\Sdet}\shiggsdet)\big\|_{L_{t,x}^2([0, \tdet])}^2 + C t^{\gamma_1 - 4} \;,
	\end{equs}
	which means the left-hand side of \eqref{e:iteration_deterministic_bounds_mid4} is bounded by (here we also apply Hölder's inequality in time)
	\begin{equ}
		4 \tdet^{ -10 \fa } \big\|\bil(\shiggsdet \otimes \covd_{B+\Sdet}\shiggsdet)\big\|_{L^2_{t, x}([0, \tdet])}^2 + C \tdet^{-1 + \kappa} \;. 
	\end{equ}
	By \eqref{e:shiggedet_decrease_L4} with $\gamma = 4$ and $\gamma = \gamma_1$, we have
	\begin{equs}
	\label{e:shiggedet_decrease_L4_appli}
		4 \big\|\bil(\shiggsdet \otimes \covd_{B+\Sdet}\shiggsdet)\big\|^2_{L^2_{t, x}([0, \tdet])} \lesssim \frac{\| \shiggsdet_0 \|_{L^4}^{\gamma_1} - \| \shiggsdet(\tdet) \|_{L^4}^{\gamma_1} }{\|\shiggsdet(\tdet) \|_{L^4}^{\gamma_1 - 4} }\ \wedge \| \shiggsdet_0 \|_{L^4}^{4} \;. 
	\end{equs}
	We discuss the following cases, with $C_\star$ to be chosen later: 
	\begin{enumerate}
		\item If $ \| \shiggsdet_0 \|_{L^4}^{4} \leq C_\star \tdet^{-\frac45}$, then the bound \eqref{e:iteration_deterministic_bounds_mid4} directly follows from the above discussion.
		
		\item If $ \| \shiggsdet_0 \|_{L^4}^{4} \geq C_\star \tdet^{-\frac45}$, and $\| \shiggsdet(\tdet) \|_{L^4} \leq \frac12 \| \shiggsdet_0 \|_{L^4}$, then by \eqref{e:shiggedet_decrease_L4_appli} we have 
		\begin{equs}
			 4 \tdet^{-10 \fa }\big\|\bil(\shiggsdet \otimes \covd_{B+\Sdet}\shiggsdet)\big\|_{L^2_{t, x}([0, \tdet])}^2 &\leq C \tdet^{-10 \fa }  \| \shiggsdet_0 \|_{L^4}^{4} \\
		& \leq C C_\star^{-\frac{50}{4} \fa} \| \shiggsdet_0 \|_{L^4}^{\gamma_1} 
		\leq \eps \bigl( \| \shiggsdet_0 \|_{L^4}^{\gamma_1} - \| \shiggsdet(\tdet) \|_{L^4}^{\gamma_1} \bigr) \;,
		\end{equs}
		provided that we choose $C_\star$ large enough. 
		
		\item If $ \| \shiggsdet_0 \|_{L^4}^{4} \geq C_\star \tdet^{-\frac45}$, and $\| \shiggsdet(\tdet) \|_{L^4} \geq \frac12 \| \shiggsdet_0 \|_{L^4}$, then by \eqref{e:shiggedet_decrease_L4_appli},
		\begin{equs}
			4 \tdet^{-10 \fa }\big\|\bil(\shiggsdet \otimes \covd_{B+\Sdet}\shiggsdet)\big\|_{L^2_{t, x}([0, \tdet])}^2 &\leq C \tdet^{-10 \fa} \| \shiggsdet_0 \|_{L^4}^{- (\gamma_1 - 4)} \bigl( \| \shiggsdet_0 \|_{L^4}^{\gamma_1} - \| \shiggsdet(\tdet) \|_{L^4}^{\gamma_1} \bigr) \\
			&\leq C C_\star^{-\frac{(\gamma_1 - 4)}{4} } \tdet^{-10\fa + \frac15 (\gamma_1 - 4)} \bigl( \| \shiggsdet_0 \|_{L^4}^{\gamma_1} - \| \shiggsdet(\tdet) \|_{L^4}^{\gamma_1} \bigr) \\
			&\leq \eps \bigl( \| \shiggsdet_0 \|_{L^4}^{\gamma_1} - \| \shiggsdet(\tdet) \|_{L^4}^{\gamma_1} \bigr) \;,
		\end{equs}
		provided that we choose $C_\star$ large enough.
	\end{enumerate}
	The bound \eqref{e:iteration_deterministic_bounds_mid4} now follows by combining these three cases. 
\end{proof}

\begin{lemma}
\label{le:iteration_deterministic_bounds_mid5}
Let $T_0$ be as in Proposition~\ref{pr:iteration_detarministic_bounds}. For any $\varepsilon > 0$, there exists some constant $C = C_{\varepsilon}$ such that the following holds. Let $B_0^{\sharp} = \YMHFsmooth_0 + \YMHFrough_0$ and $T \in (0, T_0]$ be such that \eqref{e:iteration_detarministic_bounds_connection_assumption} holds. 
Let $B, \YMHFsmooth, g$ be as given by Proposition~\ref{pr:admissible_timescale_general}. 
Let $\higgsdet_0 = \shiggsdet_0 + \rhiggsdet_0$ be such that \eqref{e:phi_0_smallness_global} holds. Let $\Sdet, \higgsdet$ be the solution of \eqref{eq:Sdet-equation}--\eqref{eq:higgsdet-equation} and let $\shiggsdet$ solve \eqref{e:shiggsdet_eq}.
Then for any $t \in (0, T]$ we have
	\begin{equs}
	\label{e:iteration_deterministic_bounds_mid5}
		\| \shiggsdet_0 \|_{L^4}^{\gamma_1 - 4} \Big| \int_{\T^2} |\shiggsdet(t)|^2 & \langle \shiggsdet(t), \rlinhiggsdet(t) \rangle_{\V} \, \dx \Big| \leq  \varepsilon \bigl(\| \shiggsdet_0 \|_{L^4}^{\gamma_1} - \| \shiggsdet(\tdet) \|_{L^4}^{\gamma_1} \bigr) + C \tdet^{-1 + \kappa}  \\
		& + \| \shiggsdet_0 \|_{L^4}^{\gamma_1 - 4} \Big| \int_{\T^2} | e^{\tdet \Delta_{+}} (g_0^{-1} \act \shiggsdet_0) |^2 \langle e^{\tdet \Delta_{+}} (g_0^{-1} \act \shiggsdet_0) , e^{\tdet \Delta_{+}} (g_0^{-1} \act \rhiggsdet_0) \rangle_{\V} \, \dx \Big| \;.
	\end{equs}
\end{lemma}
\begin{proof}
First we replace $\shiggsdet$ on the left-hand side by $\slinhiggsdet = e^{\tdet \Delta_{B(0)}} \shiggsdet_0$. By Lemma~\ref{lemma:shiggsdet-linear-approximation} and \eqref{eq:Sdet-size-estimate}, we get 
\begin{equ}
\bigl\|\shiggsdet - \slinhiggsdet\bigr\|_{L_t^\infty L_x^4([0, t])} \lesssim t^{\gamma_1 - 4 } \| \shiggsdet_0 \|_{L^4} \;.
\end{equ}
	Therefore, the error when replacing $\shiggsdet$ by $\slinhiggsdet$ on the left-hand side of \eqref{e:iteration_deterministic_bounds_mid5} is bounded by
	\begin{equs}
		 \,& \|\shiggsdet_0 \|_{L^4}^{\gamma_1 - 1} t^{\gamma_1 - 4 } \cdot \tdet^{-\frac{\kappa}{2} - 2\fa} \\
		 & \leq  \|\shiggsdet_0 \|_{L^4}^{\gamma_1 - 1} \tdet^{\frac{\gamma_1 - 4}{2} } \leq \varepsilon  C_{\eqref{e:shiggsdet_decrease}, \gamma_1, 4}  \| \shiggsdet_0 \|_{L^4}^{\gamma_1 + 2} \tdet + C \tdet^{(\frac{\gamma_1 - 4}{2} - \frac{\gamma_1 - 1}{\gamma_1 + 2}) \frac{\gamma_1 + 2}{3} } \;,
	\end{equs}
	which is acceptable by the same argument as for \eqref{e:young_shiggsdet_decrease_L4}. 
	
	Now let $g, B, \YMHFsmooth, \YMHFrough$ be as given by Proposition~\ref{pr:admissible_timescale_general}. In particular, we have that $g^{-1} \act B = \YMHFsmooth + \YMHFrough$. Note that  
	\begin{equs}
		\Big| \int_{\T^2} |\slinhiggsdet |^2 \langle \slinhiggsdet, \rlinhiggsdet \rangle_{\V} \, \dx \Big| = 
		\Big| \int_{\T^2} |g^{-1} \act \slinhiggsdet |^2 \langle g^{-1} \act \slinhiggsdet, g^{-1} \act \rlinhiggsdet \rangle_{\V} \, \dx \Big| 
	\end{equs}
	By Lemma~\ref{le:B-YMHFsmooth-heat-comparison} (applied with $t_0$ therein taken to be $t$), we get
	\begin{equs}
		\| (g^{-1} \act \slinhiggsdet)(\tdet) - e^{\tdet\Delta_{+} } (g_0^{-1} \act \shiggsdet_0) \|_{L^4}  &\lesssim \tdet^{\frac13} \|\shiggsdet_0 \|_{L^4} \;, \\ 
		\| t^{\frac{\kappa}{2} } \big( (g^{-1} \act \rlinhiggsdet) (\tdet) - e^{\tdet\Delta_{+} } (g_0^{-1} \act \rhiggsdet_0 ) \big) \|_{L^4}  &\lesssim \tdet^{\frac13} \| t^{\frac{\kappa}{2} } e^{\tdet\Delta_{+} } (g_0^{-1} \act \rhiggsdet_0) \|_{L^\infty_t L^4_x} \;.
	\end{equs}
	Therefore, we get (applying \eqref{e:phi_0_smallness_global})
	\begin{equs}
		\, &\| \shiggsdet_0 \|_{L^4}^{\gamma_1 - 4} \Bigg| \Big| \int_{\T^2} |\slinhiggsdet(t) |^2 \langle \slinhiggsdet (t), \rlinhiggsdet(t) \rangle_{\V} \, \dx \Big | \\
		& - \Big| \int_{\T^2} | e^{\tdet \Delta_{+}} (g_0^{-1} \act \shiggsdet_0) |^2 \langle e^{\tdet \Delta_{+}} (g_0^{-1} \act \shiggsdet_0) , e^{\tdet \Delta_{+}} (g_0^{-1} \act \rhiggsdet_0) \rangle_{\V} \, \dx \Big|  \Bigg|  \\
		&\lesssim t^{\frac13 - \frac{\kappa}{2}} \| \shiggsdet_0 \|_{L^4}^{\gamma_1 - 1} \| t^{\frac{\kappa}{2} } e^{\tdet\Delta_{+} } (g_0^{-1} \act \rhiggsdet_0) \|_{L^\infty_t L^4_x} \lesssim \| \shiggsdet_0 \|_{L^4}^{\gamma_1 - 1} t^{\gamma_1 - 4} \;,
	\end{equs}
	which can then be estimated as in \eqref{e:young_shiggsdet_decrease_L4} (with $\varepsilon$ in place of $1/8$). The bound \eqref{e:iteration_deterministic_bounds_mid5} now follows.
\end{proof}

We also record the following perturbative result on $\covd_{B + \Sdet} F_{B + \Sdet}$, which will also be needed in the global iteration argument.

\begin{lemma}
\label{le:D_BS_F_BS}
	There exists some constant $C$ such that the following holds. In the same setting as in Proposition~\ref{pr:iteration_detarministic_bounds}, for any $t \in (0,T]$ we have
	\begin{equs}
		\bigl\|\covd_{B+\Sdet}F_{B+\Sdet}(t)\bigr\|_{L^2} & \leq C t^{-\frac12} (1 + \|t^{1/2} F_{\YMHFsmooth}\|_{L_t^2 L_x^\infty([0, T])}) \| \covd_{\YMHFsmooth} F_{\YMHFsmooth} \|_{L^{2}_{t, x}([0, T]) } + C  t^{- \frac32 + \frac{\gamma_0 - 2}{50} - 2\fa - \frac2r}  \\
		\label{e:D_BS_F_BS}
		&+ C t^{-\frac12 + \frac{\gamma_0 - 2}{50} - 2\fa - \frac2r} \Big( \| F_{\YMHFsmooth(t)} \|_{L^2} + C \sqrt{ \| \shiggsdet_0 \|_{L^4}^{\gamma_1} - \| \shiggsdet(\tdet) \|_{L^4}^{\gamma_1} } \Big) \;.
	\end{equs}
\end{lemma}

\begin{proof}
	Expanding $\covd_{B + \Sdet} F_{B + \Sdet}$ and using the triangle inequality, we have
	\begin{equ}
		\bigl\|\covd_{B+\Sdet}F_{B+\Sdet}\bigr\|_{L^2}
		\leq \bigl\|\covd_BF_B\bigr\|_{L^2}
		+ \bigl\|\covd_B\bigl(\curl_B\Sdet+[\Sdet^1,\Sdet^2]\bigr)\bigr\|_{L^2} + \|\Sdet\|_{L^\infty}\|F_{B+\Sdet}\|_{L^2} \;.
	\end{equ}
	By Proposition~\ref{pr:admissible_timescale_general} (recall that in the setting of Proposition~\ref{pr:iteration_detarministic_bounds}, we apply this proposition with $K = T^{-2\fa}$), the $L_x^\infty$ estimate \eqref{e:F_bounds} on $F_{\YMHFsmooth}(t)$, and Lemma~\ref{le:DF_bound}, we have
	\begin{equs}
		\|\covd_B F_B (t)\|_{L^2} & \leq \|\covd_{\YMHFsmooth} F_{\YMHFsmooth} (t)\|_{L^2} + C t^{- 1 - \frac{\kappa}{2} - 3 \fa} \\
		&\lesssim t^{-\frac12} (1 + \|t^{1/2} F_{\YMHFsmooth}\|_{L_t^2 L_x^\infty([0, T])}) \| \covd_{\YMHFsmooth} F_{\YMHFsmooth} \|_{L^{2}_{t, x}([0, T]) } + t^{- 1 - \frac{\kappa}{2} - 3 \fa} \;. 
	\end{equs}
	By \eqref{e:iteration_deterministic_bounds_mid1_intermed_Sdet_Linfty}, we have
	\begin{equs}
		\bigl\|\covd_B\bigl(\curl_B\Sdet+[\Sdet^1,\Sdet^2]\bigr)\bigr\|_{L^2} \lesssim t^{- \frac32 + \frac{\gamma_0 - 2}{50} - 2\fa - \frac2r} \;.
	\end{equs}
	By \eqref{e:iteration_deterministic_bounds_mid1_intermed_Sdet_Linfty} and \eqref{e:iteration_deterministic_bounds_mid1}, we have 
	\begin{equs}
		\|\Sdet\|_{L^\infty}\|F_{B+\Sdet}\|_{L^2} & \lesssim t^{-\frac12 + \frac{\gamma_0 - 2}{50} - 2\fa - \frac2r} \cdot \Big( \| F_{\YMHFsmooth(t)} \|_{L^2} + \| \covd_{\YMHFsmooth} F_{\YMHFsmooth} \|_{L^2_{t, x}([0, t]) } + t^{\frac{-1 + \kappa}{2} } \\
		&+ t^{- \frac12 - 5 \fa } \bigl( \big\|\bil(\higgsdet \otimes \covd_{B+\Sdet}\higgsdet)\big\|_{L^1_{t} L^2_{x}([0, \tdet])} + \big\|t^{1/2}\bil(\higgsdet \otimes \covd_{B+\Sdet}\higgsdet)\big\|_{L_{t,x}^2([0, \tdet])} \bigr) \\
		&+ \| \covd_{\YMHFsmooth(t)} e^{t \Delta_+^F} \YMHFrough_0 \|_{L^\infty} + \| e^{t \Delta_+^F} \YMHFrough_0 \|_{L^\infty}^2  \Big) \;.
	\end{equs}
	which is acceptable by the assumption \eqref{e:iteration_detarministic_bounds_connection_assumption} and Lemma~\ref{le:iteration_deterministic_bounds_mid4}. The lemma now follows by combining the previous estimates.
\end{proof}

\section{Covariant ansatz}\label{section:covariant-Ansatz}

We recall from \eqref{intro:eq-SYMH} that we are interested in the stochastic Yang--Mills--Higgs equations in the infinitesimal DeTurck gauge. They are formally given by 
\begin{equation}\label{local:eq-SYM}
\begin{aligned}
\partial_t A &= - \covd_A^\ast F_A - \covd_A \covd_A^\ast (A-B) +  \bil( \higgs \otimes \covd_A \higgs)+\xi, &\qquad A(0)&=A_0,\\
\partial_t \higgs &= \Delta_A \higgs + ( \covd_A^\ast (A-B)) \higgs - |\higgs|^2 \higgs -\mass \higgs+\zeta, & \higgs(0)
&=\higgs_0, 
\end{aligned}
\end{equation}
where $B$ is the solution to the deterministic Yang--Mills gradient flow 
\begin{equation}\label{local:eq-B}
\partial_t B = - \covd_B^\ast F_B, \qquad B(0)=A_0.
\end{equation}
For suitable initial data $A_0$, the local well-posedness of \eqref{local:eq-B} has already been established in Section~\ref{section:deterministic-rough}. 
The local well-posedness of~\eqref{local:eq-SYM}, however, is much more challenging. It can be established using methods similar to those used in \cite{CCHS22,BC23}, which proved the local well-posedness of the stochastic Yang--Mills equations in the original DeTurck gauge. For our purposes, however, the resulting local well-posedness theorem for \eqref{local:eq-SYM} would be insufficient. As already discussed in Remark \ref{main:rem-local-condition}, the reason for this is that the methods in \cite{CCHS22,BC23} are non-covariant. Used directly, the arguments therefore lead to local results for \eqref{local:eq-SYM} that involve gauge dependent norms of the initial data, whereas our global iteration scheme only uses gauge invariant norms\footnote{Strictly speaking, they are not norms since they are typically norms of mildly nonlinear expressions of the state $(A_0,\higgs_0)$.}. While the gauge dependence can be remedied by passing  to gauge-equivalent initial data (and using the gauge covariance of the solution theory developed in \cite{CCHS22,BC23}), the resulting local time of existence would still depend too heavily on the size of the initial data. 

Instead of using the local theory from \cite{CCHS22,BC23}, we develop a covariant local theory of~\eqref{local:eq-SYM}, which will be distributed over Sections~\ref{section:covariant-Ansatz}--\ref{section:objects}. In this section, we derive a covariant paracontrolled ansatz, which is the starting point of our covariant local theory. Some aspects of this were already described in Section \ref{section:main-local}, i.e., in the statement and discussion of our covariant local well-posedness result. However, given the sheer length of this article, we believe that occasional repetitions are helpful to the reader, and we will therefore start the discussion of our covariant local theory anew (and go into much more detail).

\subsection{Smooth approximation and stochastic objects}

In order to discuss \eqref{local:eq-SYM} rigorously, we need to first replace the space-time white noises $\xi$ and $\zeta$ with smooth approximations, and then later take limits. In \cite{CCHS22,BC23}, these smooth approximations are defined using convolutions against  smooth, compactly supported kernels in the spatial or space-time variables, respectively. However, since this would break the gauge covariance of \eqref{local:eq-SYM}, we work with different, covariant smooth approximations. They were previously defined in Definition \ref{main:def-gauge covariant-approximation-space-time-white-noise} and are given by 
\begin{equation}\label{local:eq-xi-zeta-smooth-approximation}
\begin{aligned}
\xi_{\Beps}(t,x) &=   \int_0^t \ds \int_{\T^2} \dy \, \rho_\eps^\prime(t-s) \fq_B(s,y;t,x) \xi(s,y), \\ 
\zeta_{\Beps}(t,x) &=   \int_0^t \ds \int_{\T^2} \dy \, \rho_\eps^\prime(t-s) \fp_B(s,y;t,x) \zeta(s,y). 
\end{aligned}
\end{equation}
Here, $\eps \in \sdyadic$, $\rho_\eps$ is the cutoff function from \eqref{main:eq-rho-eps}, $\rho_\eps^\prime$ is its time-derivative, and $\fp_B$ and $\fq_B$ are the covariant heat kernels from Definition~\ref{kernel:def-kernels}. 

The precise form of the covariant smooth approximations from \eqref{local:eq-xi-zeta-smooth-approximation} will be useful for many different reasons. For now, we only briefly discuss two aspects of \eqref{local:eq-xi-zeta-smooth-approximation}, for which we focus on $\xi_{\Beps}$:  
First, $\xi_{\Beps}$ converges to $\xi$ in the sense of space-time distributions as $\eps\rightarrow 0$ and therefore deserves to be called a smooth approximation of $\xi$. This follows from the fact that $\rho_\eps^\prime(t-s)$ converges to $\delta(t-s)$ as $\eps\rightarrow 0$ and $\fq_B(s,y;t,x)$ converges to $\delta(x-y) \Id_{\frkg\oplus \frkg}$ as $s\rightarrow t$. Second, even before taking the limit $\eps\rightarrow 0$, $\xi_{\Beps}$ is gauge covariant, and therefore deserves to be called a covariant smooth approximation. To see this, let $g\in \CG$ be any smooth, time-independent gauge transformation. From the gauge covariance of the covariant heat kernel~$\fq_B$, it then follows that 
\begin{equation*}
\Ad_{g(x)} \xi_{\Beps}(t,x) = \iint \ds \dy \,  \rho_\eps^\prime(t-s) \fq_{g\act B}(s,y;t,x) \bigl( \Ad_g \xi \bigr)(s,y) \eqlaw \xi_{g\act B,\eps}(t,x)\;,
\end{equation*}
where the equality in law holds as  processes in $(t,x)$.
In order to define our smooth approximation of \eqref{local:eq-SYM}, we need to make use of the renormalisation terms
previously introduced in Definition~\ref{main:def-renormalization-constant}. Recall that we defined
\begin{equation}\label{local:eq-renormalization-operators}
\begin{aligned}
\renormA &= \Crho \bigl( \Cov^{\frkg}_{\frkg} + \Cov^{\frkg}_{\V}\bigr) -\tfrac{1}{8\pi} \Cov^{\frkg}_{\V}, \\ 
\renormphi &= (-\sigma_{\rho,\eps}^2 + \Crho ) \Cov^{\V}_{\frkg} + (\dim(\V)+2) \sigma_{\rho,\eps}^2 \Id_{\V} - \mass,
\end{aligned}
\end{equation}
where $\Cov^{\frkg}_{\frkg}$, $\Cov^{\frkg}_{\V}$, and $\Cov^{\V}_{\frkg}$  are as in Definition~\ref{prelim:def-Casimir-induced}.

The counterterms in \eqref{local:eq-renormalization-operators}, including the non-divergent $\Crho$- and $\frac{1}{8\pi}$-terms, are needed to ensure that our limiting solution does not depend on our choice of $\rho \in \mc{P}$ and is gauge equivalent (in law) to the solution previously constructed in \cite{CCHS22}; see Section~\ref{section:langevin}. The mass $\mass$, which appeared in the $\higgs$-equation in \eqref{local:eq-SYM}, is included in \eqref{local:eq-renormalization-operators} for notational convenience, and could have been kept separate. Equipped with our covariant smooth approximations of $\xi$ and $\zeta$ and our renormalisations, we previously introduced \eqref{main:eq-SYMH-eps}, i.e., our covariant smooth approximation of the stochastic YMH equations \eqref{local:eq-SYM}. In coordinates, it can be written as 
\begin{equation}\label{local:eq-SYM-coordinates-eps}
\begin{aligned}
\partial_t A_\eps^k &= \covd_{A_\eps,j} F_{A_\eps}^{jk} + \covd_{A_\eps}^k \covd_{A_\eps}^j (A_\eps-B)_j + \bil\bigl( \higgs_\eps \otimes \covd_{A_\eps}^k \higgs_\eps \bigr)+  \renormA (A_\eps -B)^k + \xi_{\Beps}^k, \\
\partial_t \higgs_\eps &= \Delta_{A_\eps} \higgs_\eps - (\covd_{A_\eps}^j(A_\eps-B)_j) \higgs_\eps 
-  |\higgs_\eps|^2 \higgs_\eps+ \renormphi (\higgs_\eps) + \zeta_{\Beps}, \\
A_\eps(0)&=A_0, \qquad \higgs_\eps(0)= \higgs_0.   
\end{aligned}
\end{equation}
 We emphasise that $B$  is still as in \eqref{local:eq-B} and, in particular, does not depend on $\eps$.
As a first step towards our ansatz, we rewrite \eqref{local:eq-SYM-coordinates-eps} in terms of the difference $C_\eps\coloneq A_\eps-B$. From Lemma~\ref{prelim:lem-linearization}, it follows that
\begin{equation}\label{local:eq-C}
\begin{aligned}
( \partial_t - \Delta_B ) C_\eps^k - 2 \big[ F_B^{kj},C_{\eps,j} \big] &=  \big[ C_{\eps,j} , 2 \covd_B^j C_\eps^k - \covd_B^k C_\eps^j \big]
+ \big[ C_{\eps,j}, \big[ C_\eps^j , C_\eps^k \big]\big] \\
&\quad +\bil\bigl( \higgs_\eps \otimes \covd_B^k \higgs_\eps \bigr) + \bil \bigl( \higgs_\eps \otimes C_\eps^k \higgs_\eps\bigr) + \renormA ( C_\eps^k ) + \xi_{\Beps}^k, \\ 
(\partial_t - \Delta_B ) \higgs_\eps &= 2 C_{\eps,j} \covd_B^j \higgs_\eps + C_\eps^j C_{\eps,j} \higgs_\eps -   |\higgs_\eps|^2 \higgs_\eps + \renormphi(\higgs_\eps) + \zeta_{\Beps}. 
\end{aligned}
\end{equation}
We note that in our formulation of \eqref{local:eq-C}, we used that there are two cancellations between \eqref{prelim:eq-linearization-3} and \eqref{prelim:eq-linearization-4}.  We also emphasise that the term $2 \big[ F_B^{kj},C_{\eps,j} \big]$ is one of the most important terms in \eqref{local:eq-C}. The reason is that  the dependence of the local time of existence on the size of the curvature tensor $F_B$ is very important for our purposes. Due to this, we treat $\Delta_B C_\eps^k + 2 \big[ F_B^{kj},C_{\eps,j} \big]$ rather than just $\Delta_B C^k_\eps$ non-perturbatively. As suggested by \eqref{local:eq-C}, we will work covariantly with respect to the connection one-form $B$, which leads to the following notation.

\begin{notation}\label{local:notation-omit-B}
Given a solution $B$ to the YM gradient flow, we write 
\begin{equation*}
\covd = \covd_B \qquad \text{and} \qquad F = F_B.
\end{equation*}
That is, whenever the connection one-form is omitted in covariant derivatives or curvature tensors, it is given by $B$. However, we will never omit $B$ in our notation for the covariant Laplacian $\Delta_B$, as it may lead to confusion. Furthermore, we write 
\begin{equation*}
\bigl( \ad(F) C_\eps\bigr)^k = \ad(F^{kj}) C_{\eps,j} = \big[ F^{kj}, C_{\eps,j} \big].
\end{equation*}
That is, we interpret $\ad(F)$ as a map from $\frkg^2$ into $\frkg^2$. Finally, we omit the superscript $k$ from \eqref{local:eq-C} and related equations. In other words, unless the components are contracted, we no longer write out the components of $C_\eps$, $\ad(F)C_\eps$, or $\covd C_\eps^j$, and instead interpret them directly as elements of $\frkg^2$. 
\end{notation}
In our simplified notation, \eqref{local:eq-C} then takes the form
\begin{equation}\label{local:eq-C-simplified}
\begin{aligned}
(\partial_t - \Delta_B - 2 \ad(F) ) C_\eps  &=  \big[ C_{\eps,j} , 2 \covd^j C_\eps - \covd C_\eps^j \big]
+ \big[ C_{\eps,j}, \big[ C_\eps^j , C_\eps \big]\big] \\
&\quad + \bil\bigl( \higgs_\eps \otimes \covd \higgs_\eps \bigr) + \bil \bigl( \higgs_\eps \otimes C_\eps \higgs_\eps\bigr) + \renormA ( C_\eps )  + \xi_{\Beps}, \\
(\partial_t - \Delta_B ) \higgs_\eps &= 2 C_{\eps,j} \covd^j \higgs_\eps + C_\eps^j C_{\eps,j} \higgs_\eps -   |\higgs_\eps|^2 \higgs_\eps  + \renormphi(\higgs_\eps) + \zeta_{\Beps}. 
\end{aligned}
\end{equation}
We now recall \eqref{main:local-eq-Sdet}--\eqref{main:eq-mixed-quadratic-stochastic}, i.e., the definitions of the deterministic and random terms in our covariant ansatz. The deterministic terms $\Sdet$ and $\higgsdet$ were defined as solutions of the system 
\begin{align}
\hspace{-1.5ex}(\partial_t - \Delta_B - 2 \ad(F)) \Sdet &= [ \Sdet_j, 2 \covd^j \Sdet - \covd \Sdet^j] + [ \Sdet_j , [\Sdet^j , \Sdet ]] 
+ \bil\bigl( \higgsdet \otimes \covd \higgsdet\bigr) + \bil \bigl( \higgsdet \otimes \Sdet \higgsdet \bigr), \label{local:eq-Sdet} \\ 
(\partial_t - \Delta_B) \higgsdet &= 2 \Sdet_j \covd^j \higgsdet  + \Sdet^j \Sdet_j \higgsdet - |\higgsdet|^2 \higgsdet, \label{local:eq-higgsdet} \\ 
\Sdet(0)&=0, \qquad \higgsdet(0)=\higgsdet_0. \notag 
\end{align}
Here, in keeping with the notation from Section \ref{section:deterministic-rough}, we denote the initial data of the Higgs field by $\higgsdet_0 \eqdef \higgs_0$. Later on, we will also decompose $\higgsdet_0$ into the positive and negative regularity terms $\shiggsdet_0$ and $\rhiggsdet_0$ appearing in  Section~\ref{section:deterministic-rough}. The covariant linear stochastic objects are the solutions of the inhomogeneous, linear part of~\eqref{local:eq-C-simplified}, i.e., 
\begin{equation}\label{local:eq-linear-stochastic}
\begin{aligned}
(\partial_t - \Delta_B - 2 \ad(F))  \linear[\Beps]&= \xi_{\Beps}, \qquad \linear[\Beps] (0)=0, \\ 
(\partial_t - \Delta_B )  \philinear[\Beps]&= \zeta_{\Beps}, \qquad \philinear[\Beps] (0)=0. 
\end{aligned}
\end{equation}
From the definitions of $\fp_B$, $\fq_B$, $\xi_{\Beps}$, and $\zeta_{\Beps}$, it follows that the covariant linear stochastic objects can be written as 
\begin{equation}\label{local:eq-linear-stochastic-integral-version}
\begin{aligned}
\linear[\Beps](t,x) &= \iint \ds \dy \, \rho_\eps(t-s) \fq_B(s,y;t,x) \xi(s,y), \\ 
\philinear[\Beps](t,x) &= \iint \ds \dy \,  \rho_\eps(t-s) \fp_B(s,y;t,x) \zeta(s,y). 
\end{aligned}
\end{equation}
This can be seen by acting on the right-hand sides of \eqref{local:eq-linear-stochastic-integral-version} with 
 $\partial_t - \Delta_B - 2 \ad(F)$ or $\partial_t - \Delta_B$ and then using the definitions of $\fq_B$ and $\fp_B$. The covariant quadratic stochastic objects are solutions with zero initial data of  
\begin{align}
(\partial_t - \Delta_B - 2 \ad(F))  \quadratic[\Beps] 
    &= \big[ \linear[\Beps][][j], 2 \covd^j \linear[\Beps] - \covd \linear[\Beps][j]\big], \label{local:eq-A-quadratic-stochastic} \\ 
(\partial_t - \Delta_B - 2 \ad(F))  \phiquadratic[\Beps] 
    &= \bil \bigl(  \philinear[\Beps] \otimes \covd \philinear[\Beps]\bigr), \label{local:eq-phi-quadratic-stochastic}\\ 
(\partial_t - \Delta_B)  \mixedquadratic[\Beps] 
&= 2 \linear[\Beps][j] \covd_j \philinear[\Beps]. \label{local:eq-mixed-quadratic-stochastic}
\end{align}
After this recap of \eqref{main:local-eq-Sdet}--\eqref{main:eq-mixed-quadratic-stochastic}, we now introduce nonlinear remainders, which were not defined in Section \ref{section:main-local}. The 
nonlinear remainders are given by 
\begin{equation}\label{local:eq-Z-def}
\begin{aligned}
\Wrem_\eps &\coloneq A_\eps -B - \linear[\Beps], \qquad &\qquad 
Z_\eps &\coloneq A_\eps -  B - \Sdet-  \linear[\Beps]- \quadratic[\Beps] - \phiquadratic[\Beps], \\
\phiw_\eps &\coloneq \higgs_\eps - \philinear[\Beps], \qquad & \qquad \phiz_\eps &\coloneq \higgs_\eps - \higgsdet - \philinear[\Beps] - \mixedquadratic[\Beps]. 
\end{aligned}
\end{equation}
The nonlinear remainders $\Wrem_\eps$ and $\phiw_\eps$ are only used to simplify the terms in our evolution equations below and, due to this, we do not write out their own evolution equations. The evolution equations for $Z_\eps$ and $\phiz_\eps$ are given by
\begin{equs}
(\partial_t - \Delta_B - 2 &\ad(F) )Z_\eps \\
&=\Big[  \Wrem_{\eps,j} , 2 \covd^j \linear[\Beps] - \covd \linear[\Beps][j] \Big]
+ \Big[ \linear[\Beps][][j], 2 \covd^j \Wrem_\eps - \covd \Wrem^j_\eps \Big] \label{local:eq-Z-new-1} \\ 
&\quad+ \Big[ \Wrem_{\eps,j}  , 2 \covd^j \Wrem_\eps - \covd \Wrem_\eps^j \Big] - \Big[ \Sdet_j , 2 \covd^j \Sdet - \covd \Sdet^j \Big] \label{local:eq-Z-new-2} \\ 
&\quad+ \Big[  \linear[\Beps][][j][]\hspace{-0.5ex}+ \Wrem_{\eps,j}, \Big[ \linear[\Beps][j][][] \hspace{-0.5ex}+ \Wrem_\eps^j,  \linear[\Beps] \hspace{-0.5ex}+ \Wrem_\eps \Big] \Big] - \big[ \Sdet_j, \big[ \Sdet^j , \Sdet \big]\big] \label{local:eq-Z-new-3} \\ 
&\quad+  \bil\Big( \phiw_\eps \otimes \covd \philinear[\Beps]\Big)
+ \bil \Big( \philinear[\Beps]\otimes \covd \phiw_\eps \Big)\label{local:eq-Z-new-4}\\
&\quad+ \bil \Big( \phiw_\eps \otimes \covd \phiw_\eps\Big) - \bil \Big( \higgsdet \otimes \covd \higgsdet\Big)\label{local:eq-Z-new-5}\\
&\quad+ \bil \Big( \bigl( \philinear[\Beps]+\phiw_\eps\bigr) \otimes \bigl(\linear[\Beps]+\Wrem_\eps\bigr) \bigl( \philinear[\Beps] + \phiw_\eps\bigr) \Big) - \bil\Big( \higgsdet \otimes \Sdet \higgsdet \Big) \label{local:eq-Z-new-6}\\ 
&\quad+ \renormA \bigl(  \linear[\Beps]\hspace{-0.5ex}+W_\eps\bigr), \label{local:eq-Z-new-7}  
\end{equs}
and 
\begin{equs}
 (\partial_t - \Delta_B ) \phiz_\eps 
&= 2 \Wrem_\eps^j \covd_j \philinear[\Beps]  +2 \linear[\Beps][j] \covd_j \phiw_\eps \label{local:eq-z-new-1}\\ 
&\quad+ 2 \bigl(\Wrem_\eps^j \covd_j \phiw_\eps - \Sdet^j \covd_j \higgsdet\bigr) \label{local:eq-z-new-2}\\
&\quad+ \bigl( \linear[\Beps]+\Wrem_\eps\bigr)^j \bigl( \linear[\Beps]+\Wrem_\eps\bigr)_j \bigl( \philinear[\Beps]+\phiw_\eps)- \Sdet^j \Sdet_j \higgsdet \label{local:eq-z-new-3}\\ 
&\quad-  \big| \philinear[\Beps] + \phiw_\eps\big|^2 \bigl( \philinear[\Beps]+\phiw_\eps\bigr)  + |\higgsdet|^2 \higgsdet + \renormphi ( \philinear[\Beps][][][] +\phiw_\eps). \label{local:eq-z-new-4}
\end{equs}
As in the local theory of \cite{CCHS22,BC23}, the evolution equations for $Z_\eps$ and $\phiz_\eps$ cannot be solved using a direct contraction-mapping argument. One of the reasons is that the terms $\covd \linear[\Beps]$ in \eqref{local:eq-Z-new-1} only have regularity $-1\minus$ so that $Z_\eps$ has regularity at most $1\minus$. But, since $\linear[\Beps]$ has regularity $0\minus$, this is insufficient to even define interactions of the form $[Z_\eps,\covd \linear[\Beps]]$ or $[\covd Z_\eps, \linear[\Beps]]$. In order to solve the evolution equations for $Z_\eps$ and $\phiz_\eps$ and to have the local time of existence depend only on gauge invariant norms of $B$, we make a covariant paracontrolled ansatz. 
For this, we first need a covariant version of Littlewood--Paley decompositions and paraproducts.

\begin{remark}
In \eqref{local:eq-Z-new-2}, \eqref{local:eq-Z-new-3}, \eqref{local:eq-Z-new-5}, \eqref{local:eq-Z-new-6}, 
and \eqref{local:eq-z-new-2}--\eqref{local:eq-z-new-4}, we subtract different terms involving $\Sdet$ and $\higgsdet$. The reason for this is not their regularity, since all of them are rather smooth, but rather their size. Without subtracting them here, the remainders $Z_\eps$ and $\phiz_\eps$ would become too large for our global iteration argument. For more details on this, see Section~\ref{section:global-theory}. 
\end{remark}

\subsection{A Littlewood--Paley decomposition}\label{section:Littlewood--Paley}

In Section \ref{section:main-local}, we previously explained why the standard Littlewood--Paley operators from \eqref{main:eq-standard-Littlewood-Paley} will not play an important role in this article, and we also introduced our covariant Littlewood--Paley blocks in~\eqref{main:eq-covariant-Littlewood-Paley-block}. In this subsection, we now discuss our motivation behind \eqref{main:eq-covariant-Littlewood-Paley-block} in more detail. To this end, we first discuss two alternatives to \eqref{main:eq-covariant-Littlewood-Paley-block}, both of which encounter different issues in our setting. It is tempting to try to define a covariant Littlewood--Paley decomposition using the covariant heat operator $e^{s \Delta_{B(t)}}$, which would be in the spirit of many existing works (see e.g. \cite{KR06}). Given a parameter $s_0>0$, one could decompose 
\begin{equation}\label{local:eq-heuristic-1}
\linear[\Beps](t,x)= - \int_0^{s_0} \ds \, \partial_s \bigl( e^{s\Delta_{B(t)}} \linear[\Beps] \bigr)(t,x) + e^{s_0 \Delta_{B(t)}} \linear[\Beps](t,x).
\end{equation}
The $e^{s_0 \Delta_{B(t)}}$-term in \eqref{local:eq-heuristic-1} should be interpreted as a smooth remainder, and will not be too important in this discussion. One may then interpret 
\begin{equation*}
- \int_0^{s_0} \ds \chi(N^2 s) \partial_s \bigl( e^{s\Delta_{B(t)}} \linear[\Beps] \bigr)(t,x), 
\end{equation*}
where $\chi$ is our cutoff function (see Section~\ref{section:preliminaries}) and $N\in \dyadic$, as a covariant version of the Littlewood--Paley projection to frequencies of order $N$. This interpretation can be motivated using a similar object in the Euclidean setting, where the corresponding symbol is given by  
\begin{equation*}
n \mapsto - \int_0^{\infty} \ds \chi(N^2 s) \partial_s \bigl( e^{-s|n|^2} \bigr). 
\end{equation*}
This symbol behaves like  $N^{-2} |n|^2 e^{-c N^{-2} |n|^2}$, which is maximal for $|n|\sim N$. Equipped with the covariant Littlewood--Paley decomposition from \eqref{local:eq-heuristic-1}, one can easily write down paraproduct decompositions of multilinear interactions. For example, one can decompose
\begin{equation}\label{local:eq-heuristic-2} 
\begin{aligned}
&\quad\Big[ Z_{\eps,j} (t,x), \covd_{B(t)}^j \linear[\Beps](t,x) \Big] \\
&= \iint_{[0,s_0]^2} \hspace{-0.5ex} \ds \ds^\prime  \Big[ \partial_s \bigl( e^{s \Delta_{B(t)}} Z_{\eps,j} \bigr) (t,x), \covd_{B(t)}^j \partial_{s^\prime} \bigl( e^{s^\prime \Delta_{B(t)}} \linear[\Beps] \bigr)(t,x)   \Big] + \Big\{ \hdots \Big\},
\end{aligned}
\end{equation}
where we omitted terms involving $e^{s_0 \Delta_{B(t)}}$-operators. The low$\times$high-, high$\times$high-, and high$\times$low-interactions then correspond to the regions $s\gg s^\prime$, $s\sim s^\prime$, and $ s \ll s^\prime$, respectively. In our setting, however, there is a serious problem with \eqref{local:eq-heuristic-1} and \eqref{local:eq-heuristic-2}: We do not have suitable estimates of the kernel of the heat propagator $e^{s\Delta_{B(t)}}$ and its covariant derivatives. The reason for this is that $B(t)$ enters at the fixed time $t$ and we have no uniform-in-time bounds for the $L_x^\infty$-norm of the curvature tensor $F_B(t)$. To compensate for this, one could try to replace $e^{s\Delta_{B(t)}} \linear[\Beps]$ with the solution $\linear[\Beps](t,x;s)$ of 
\begin{equation}\label{local:eq-heuristic-3}
\bigl( \partial_s - \Delta_{B(t+s)}\bigr) \linear[\Beps](t,x;s) =0, \qquad \linear[\Beps](t,x;0)=\linear[\Beps](t,x). 
\end{equation} 
The solution of \eqref{local:eq-heuristic-3} can be explicitly represented as 
\begin{equation}\label{local:eq-heuristic-4}
\linear[\Beps](t,x;s) = \int \dy \, \fp_B(t,y;t+s,x) \linear[\Beps](t,y). 
\end{equation}
By the results from Section~\ref{section:covariant-heat-kernel-estimates}, the covariant heat kernel appearing in \eqref{local:eq-heuristic-4} obeys several useful estimates. However, there is an issue with \eqref{local:eq-heuristic-4}, which appears once we examine the counterpart of \eqref{local:eq-heuristic-2}. It is now of the form 
\begin{equation}\label{local:eq-heuristic-5} 
\begin{aligned}
&\quad \Big[ Z_{\eps,j} (t,x), \covd_{B(t)}^j \linear[\Beps](t,x) \Big] \\
&= \iint \dy \dy^\prime \iint_{[0,s_0]^2} \hspace{-0.5ex} \ds \ds^\prime  \Big[ \partial_s \bigl( \fp_B(t,y;t+s,x) Z_{\eps,j}(t,y) \bigr),\\
&\hspace{25ex} \covd_{B(t)}^j \partial_{s^\prime} \bigl( \fp_B(t,y^\prime;t+s^\prime,x) \linear[\Beps](t,y^\prime) \bigr)   \Big] +  \Big\{ \hdots \Big\}. 
\end{aligned}
\end{equation}
In \eqref{local:eq-heuristic-5}, the covariant derivative $\covd_{B(t)}$ involves $B$ at time $t$ and the two covariant heat kernels have end times $t+s$ and $t+s^\prime$. Due to this mismatch, it can be difficult to control $\covd_{B(t)} \fp_B(t,y^\prime;t+s^\prime,x)$ or, after integration by parts, $\covd_{B(t)} \fp_B(t,y;t+s,x)$. In principle, it may be possible to address this by controlling $B(t+s^\prime)-B(t)$ and $B(t+s)-B(t)$ by using our bounds on covariant derivatives of $F_B$ and the fact that $B$ solves the Yang--Mills gradient flow, but we do not pursue this here.

Instead, our covariant Littlewood--Paley decomposition 
uses the heat kernels that already appear in the Duhamel integral formulations of \eqref{local:eq-linear-stochastic}, \eqref{local:eq-A-quadratic-stochastic}--\eqref{local:eq-mixed-quadratic-stochastic}, and \eqref{local:eq-Z-new-1}--\eqref{local:eq-z-new-4}. We start with the covariant Littlewood--Paley decomposition of $\linear[\Beps]$, and only later discuss the other terms in our current ansatz. Using the functions $\chi_N$ from~\eqref{prelim:eq-chi-N} and the 
 integral representation of $\linear[\Beps]$ from~\eqref{local:eq-linear-stochastic-integral-version}, we write  
\begin{equ}\label{local:eq-linear-decomposition}
\linear[\Beps] = \sum_{N} \linear[\Beps][][][N]\;,\qquad
\linear[\Beps][][][N](t,x) 
\eqdef  \iint \ds \dy \,  \fq_{\Beps,N}(s,y;t,x) \xi(s,y), 
\end{equ}
where 
\begin{equation}\label{local:eq-qB-eps-N}
\fq_{\Beps,N}(s,y;t,x) \eqdef \rho_\eps(t-s) \chi_N(t-s)\fq_B(s,y;t,x). 
\end{equation}
Here, as discussed in Notation~\ref{notation:dyadic-integers}, we use the convention that sums over $N$ always run over $N \in \dyadic$ and not $N \in \N_0$. For later use, we also define 
\begin{equation}\label{local:eq-qB-N}
\fq_{B,N}(s,y;t,x) \eqdef \chi_N(t-s)\fq_B(s,y;t,x). 
\end{equation}
Due to the support properties of $\rho_\eps$ and $\chi_N$, we see that $\linear[\Beps][][][N]$ is nonzero only for $N\lesssim \eps^{-1}$; in particular, the sum in \eqref{local:eq-linear-decomposition} is finite. 
Since  $\fq_{\Beps,N}$ and $\fq_{B,N}$ are localisations of  $\fq_B$ to times $t-s\sim N^{-2}$, 
parabolic scaling suggests that $\linear[\Beps][][][N]$ can be interpreted as a covariant Littlewood--Paley projection of $\linear[\Beps]$ to frequencies of size $\sim N$. 
This is also consistent with Lemma~\ref{object:lem-linear} below, which shows that covariant spatial derivatives and time derivatives acting on $\linear[\Beps][][][N]$ cost us $\sim N$ and $\sim N^2$, respectively.

Furthermore, also due to the aforementioned support properties, $\linear[\Beps][][][N]$ only depends on the precise choice of $\rho_\eps$ in the case $N\sim \eps^{-1}$. Indeed, if $N \sim \eps^{-1}$ is not satisfied, then the $\rho_\eps(t-s)$-factor in \eqref{local:eq-qB-eps-N} is just equal to zero or one on the support of $\chi_N(t-s)$. This will later be useful in showing that the limit of our solutions of \eqref{local:eq-SYM-coordinates-eps} does not depend on the precise choice of $\rho$. As for $\linear[\Beps]$, we define the covariant Littlewood--Paley blocks of $\philinear[\Beps]$ as 
\begin{equation}\label{local:eq-philinear-N}
\philinear[\Beps][][][N](t,x) \eqdef \iint \ds \dy \, \fp_{\Beps,N}(s,y;t,x) \zeta(s,y), 
\end{equation}
where 
\begin{equation*}
\fp_{\Beps,N}(s,y;t,x) \eqdef  \rho_\eps(t-s) \chi_N(t-s)\fp_B(s,y;t,x). 
\end{equation*}
As for $\fq_{B,N}$, we also write 
\begin{equation}\label{local:eq-pB-N}
\fp_{B,N}(s,y;t,x) \eqdef \chi_N(t-s)\fp_B(s,y;t,x). 
\end{equation}

In order to define the covariant Littlewood--Paley decompositions of the remaining terms, we start with a decomposition of the Duhamel integrals of $\partial_t - \Delta_B - 2 \ad(F)$ and $\partial_t - \Delta_B$. To this end, we first recall from Definition~\ref{def:covariant-duhamel} that
\begin{align*}
\Duh_{B}^F \bigl( \Psi \bigr)(t,x)
&= \iint \ds \dy \,\fq_B(s,y;t,x) \Psi(s,y),  \\
\Duh_{B} \bigl( \Psi \bigr)(t,x)
&= \iint \ds \dy \,\fp_B(s,y;t,x) \Psi(s,y).
\end{align*}
Here, the map $\Psi$ takes values in $\frkg^2$ or $\V$, respectively. The superscript $F$ in $\Duh_{B}^F$ indicates that $\ad(F)$-terms, such as in \eqref{local:eq-C-simplified}, have been accounted for through $\fq_B$. We can then define time-localised Duhamel integrals as
\begin{equation}\label{local:Duh-LWP-block}
\begin{aligned}
\Duh^F_{B,N}\bigl( \Psi \bigr)(t,x)
&\eqdef \iint \ds \dy \,\fq_{B,N}(s,y;t,x) \Psi(s,y), \\  
\Duh_{B,N} \bigl( \Psi \bigr)(t,x)
&\eqdef \iint \ds \dy \,\fp_{B,N}(s,y;t,x) \Psi(s,y),
\end{aligned}
\end{equation}
where $N\in \dyadic$, $\fq_{B,N}$ is as in \eqref{local:eq-qB-N}, and $\fp_{B,N}$ is as in \eqref{local:eq-pB-N}. From now on, we work on the time interval $[0,T]$ and, for notational reasons, assume that $0<T\leq \frac{1}{128}$. If we restrict our Duhamel integrals to our time interval~$[0,T]$ (both for the support of $\Psi$ and the evaluation of the Duhamel integrals), we have the decompositions
\begin{equation}\label{local:eq-Duhamel-N}
\Duh_{B} = \sum_{N \geq  N_{\lo}} \Duh_{B,N} \quad \text{and} \quad \Duh_{B}^F = \sum_{N \geq  N_{\lo}} \Duh^F_{B,N}, 
\end{equation}
where $N_{\lo}$ is the largest dyadic number not exceeding $(32 T)^{-\frac{1}{2}}$, i.e., 
\begin{equation}\label{local:eq-N-lo}
N_{\lo} \coloneq \max \bigl\{ N \in \dyadic \colon N \leq (32T)^{-\frac{1}{2}} \bigr\}. 
\end{equation}
Using the time-localised Duhamel integrals, we now define the covariant Littlewood--Paley blocks of the remaining terms by acting with $\Duh^F_{B,N}$ or $\Duh_{B,N}$ on the right-hand side of their evolution equations. To be more precise, we define the Littlewood--Paley blocks of $\Sdet$ by 
\begin{equation}\label{local:eq-Sdet-N}
\Sdet_N \coloneq \Duh^F_{B,N} \bigl( \mbox{RHS of \eqref{local:eq-Sdet}} \bigr). 
\end{equation}
Similarly, we define the covariant Littlewood--Paley blocks of $\higgsdet$ by 
\begin{equation}\label{local:eq-higgsdet-N}
\higgsdet_N = \higgsdet_{\lin,N} + \Duh_{B,N} \bigl( \mbox{RHS of \eqref{local:eq-higgsdet}} \bigr),
\end{equation}
where the covariant Littlewood--Paley block of the linear evolution is defined as 
\begin{equation*}
\higgsdet_{\lin,N}(t,x) \eqdef \int \dy \,  \fp_{B,N}(0,y;t,x) \higgsdet_0(y). 
\end{equation*}
Finally, we define the Littlewood--Paley blocks of the covariant quadratic stochastic objects by 
\begin{equation}\label{local:eq-quadratic-N}
\begin{aligned}
\quadratic[\Beps][][][N] &\coloneq \Duh^F_{B,N} \bigl( \mbox{RHS of \eqref{local:eq-A-quadratic-stochastic}}\bigr), \qquad &
\phiquadratic[\Beps][][][N] &\coloneq \Duh^F_{B,N} \bigl( \mbox{RHS of \eqref{local:eq-phi-quadratic-stochastic}}\bigr), \\
\mixedquadratic[\Beps][][][N] &\coloneq \Duh_{B,N} \bigl( \mbox{RHS of \eqref{local:eq-mixed-quadratic-stochastic}}\bigr), 
\end{aligned}  
\end{equation}
and, similarly, we define the Littlewood--Paley blocks of the nonlinear remainders by
\begin{equation}\label{local:eq-Z-z-N-combined}
\begin{aligned}
Z_{\eps,N} &\coloneq  \Duh_{B,N}^F \Big( \mbox{RHS of \eqref{local:eq-Z-new-1}--\eqref{local:eq-Z-new-7}} \Big),  \\ 
\phiz_{\eps,N} &\coloneq \Duh_{B,N} \Big(\mbox{RHS of \eqref{local:eq-z-new-1}--\eqref{local:eq-z-new-4}}  \Big).  
\end{aligned} 
\end{equation}

Before we continue with our definition of the covariant paraproducts, note that
our covariant Littlewood--Paley decomposition does not come with a collection of ``projection operators''
$P_N$ such that the decomposition of a function $A$ is given by $A = \sum_N P_N A$. For this reason, we 
will always keep track of the whole sequence of Littlewood--Paley blocks as a single object and, given
a sequence $(A_N)_{N \ge N_\lo}$, we will denote it by $A_\star$ to indicate that the index $N$ is free.
Again, free indices denoted by $\star$ will only ever take values in $2^{\N_0}$. Whenever we start from a sequence $A_\star=(A_N)_{N\geq N_{\lo}}$, we will denote the corresponding single function by $A$, i.e., we write $A=\sum_{N\geq N_{\lo}} A_N$. \label{local:text-Littlewood-Paley-blocks}

Let $\theta \in (0,1)$ be the parameter from \eqref{prelim:eq-parameters}. For two dyadic scales $M,N\in \dyadic$, we write  
\begin{equation}\label{local:eq-ll-theta}
M \ll_\theta N \quad\Longleftrightarrow\quad M \leq N^{1-\theta}.
\end{equation}
Furthermore, we write $M\gg_\theta N$ if $N\ll_\theta M$ and we write $M\sim_\theta N$ if neither $M\ll_\theta N$ nor $M\gg_\theta N$. From this definition, it follows that $M\sim_\theta N$ if and only if $\min(M,N)>\max(M,N)^{1-\theta}$. 

\begin{definition}[Covariant paraproducts]\label{local:def-paraproduct}
Let $0<T\leq \frac{1}{128}$, let $\genV_1,\genV_2$, and $\genV_3$ be vector spaces, and let $\Qc\colon \genV_1 \times \genV_2 \rightarrow \genV_3$ be a bilinear map. Furthermore, let $U_\star=(U_M)_{M\geq N_\lo}$ and let $V_\star=(V_N)_{N\geq N_{\lo}}$, where $U_M\colon [0,T] \times \T^2 \rightarrow \genV_1$ and $V_N\colon [0,T]\times \T^2 \rightarrow \genV_2$. Then, we define the paraproduct
\begin{equ}[e:defparaproduct]
\Qc \bigl( U_\star \parall V_\star \bigr) 
\coloneq \sum_{\substack{M,N\geq N_{\lo}\colon \\ M\ll_\theta N }} \Qc \bigl( U_M, V_N \bigr). 
\end{equ}
Similarly, we define $\Qc ( U_\star \paragg V_\star)$, $\Qc ( U_\star \parasim V_\star)$, and  $\Qc ( U_\star \paransim V_\star)$ by replacing the condition $M\ll_\theta N$ appearing in 
\eqref{e:defparaproduct} by ``$M \gg_\theta N$'', ``$M \sim_\theta N$'', and  ``$M \not\sim_\theta N$'', respectively. In line with \eqref{e:defparaproduct}, we will also write
\begin{equation*}
\Qc \bigl( U_\star,  V_\star \bigr) \coloneq \sum_{M,N\geq N_{\lo}} \Qc \bigl( U_M, V_N \bigr). 
\end{equation*}
\end{definition}

\begin{remark}
While these definitions apply to general $U_\star$ and $V_\star$, we will only use them for sequences of covariant Littlewood--Paley blocks.
\end{remark}

We emphasise that the parameter $\theta$ appearing in Definition~\ref{local:def-paraproduct} is needed 
for many of our arguments below, see e.g.\ Remarks~\ref{local:rem-theta-low-high} and~\ref{local:rem-theta-high-low}. The reason is that our covariant Littlewood--Paley decomposition only approximately satisfies the Littlewood--Paley trichotomy, which we discuss now. Let $f,g\colon \T^2 \rightarrow \R$ be two functions, let $K,M,N\in \dyadic$, and let $(P_L)_{L\in \dyadic}$ be the standard Littlewood--Paley projections. It then follows from frequency-support considerations that
\begin{equation*}
P_K \bigl( P_M f \cdot P_N g \bigr) =0 
\end{equation*}
unless 
\begin{equation}\label{local:eq-trichotomy}
K \sim M \gg N, \qquad K \sim N \gg M, \qquad \mbox{or} \qquad K \lesssim M \sim N. 
\end{equation}
The conditions in \eqref{local:eq-trichotomy} are often called the Littlewood--Paley trichotomy, and the above fact is heavily used in the paracontrolled calculus of \cite{GIP15}. In our covariant setting, this has no direct counterpart. For example, consider one of the terms coming from \eqref{local:eq-Z-new-1} and~\eqref{local:eq-Z-z-N-combined}, which can be decomposed into terms of the form
\begin{equation}\label{local:eq-trichotomy-term}
\Duh_{B,K}^F \bigl(  \bigl[ Z_{\eps,M,j} ,  \covd^j \linear[\Beps][][][N] \bigr] \bigr). 
\end{equation}
Unlike in the non-covariant paracontrolled calculus of \cite{GIP15}, there is no reason for this term to vanish if, for example, $M \ll N$ and $K\ll N$. When $K\sim 1$ and $M\sim N^{1-\theta}$, the low$\times$high-interaction in~
\eqref{local:eq-trichotomy-term} is still only well-defined under certain $\theta$-dependent regularity conditions on $Z_{\eps,\star}$ (see Lemma~\ref{local:lem-high-low}). Furthermore, when $K\sim N^{1-\theta}$ and $M\sim 1$, a gain from the small ratio $K/N \lesssim N^{-\theta}$ can only be obtained using norms at negative temporal regularity (see Definition~\ref{local:def-modified-solution-space}). Consequently, the regularity conditions in our contraction-mapping argument can only be satisfied by exploiting the additional room provided by the parameter $\theta$.

\subsection{Paracontrolled ansatz}

Equipped with our covariant Littlewood--Paley decomposition and paraproducts, we are now able to derive our covariant paracontrolled ansatz. This derivation will closely follow \cite[Section 3]{BC23}, except that we include a Higgs field and replace all derivatives and paraproducts by their covariant counterparts. Before turning to the derivation, however, we introduce the following notation. 

\begin{notation}[Omitting the subscript $B$]\label{local:notation-omitting-subscript}
Throughout this section, as well as Sections~\ref{sec:covariant_function_spaces}--\ref{section:objects}, the covariant stochastic objects are all defined using the same connection one-form $B$. In order to shorten our expressions, we therefore now drop the $B$-dependence from our notation. That is, we write 
\begin{equation*}
\linear[\oeps][][][\ostar]=\linear[\Beps][][][\ostar], \quad 
\philinear[\oeps][][][\ostar]=\philinear[\Beps][][][\ostar], \quad
\quadratic[\oeps][][][\ostar] = \quadratic[\Beps][][][\ostar], \quad \textup{etc.}
\end{equation*}
Since the non-covariant stochastic estimates only appear in Section~\ref{section:langevin}, which does not use the graphical symbols $\linear[\oeps][][][\ostar]$, $\philinear[\oeps][][][\ostar]$, 
and $\quadratic[\oeps][][][\ostar]$, this should not cause any confusion.
\end{notation}

We now make the ansatz
\begin{equation}\label{local:eq-paracontrolled-Ansatz}
Z_{\eps,\star} \coloneq X_{\eps,\star} + Y_{\eps,\star} \qquad \text{and} \qquad 
\phiz_{\eps,\star}\coloneq \phix_{\eps,\star} + \phiy_{\eps,\star}, 
\end{equation}
where $X_{\eps,\star}$ and $\phix_{\eps,\star}$ are paracontrolled components and $Y_{\eps,\star}$ and $\phiy_{\eps,\star}$ are smooth nonlinear remainders. In total, our covariant ansatz then takes the form 
\begin{equation*}
\begin{aligned}
 A_\eps &=  B + \Sdet +  \linear[\oeps] + \quadratic[\oeps] + \phiquadratic[\oeps] + X_\eps + Y_\eps, \\
 \higgs_\eps &=   \higgsdet + \philinear[\oeps] +  \mixedquadratic[\oeps] + \phix_\eps + \phiy_\eps. 
\end{aligned}
\end{equation*}
The terms in our ansatz for $Z_{\eps,\star}$ and $\phiz_{\eps,\star}$ will be chosen as solutions of the following system of equations\footnote{Compared to \eqref{local:eq-Z-z-N-combined}, we denote the frequency scale by $K$ instead of $N$. This will be notationally convenient in Section~\ref{section:local-paraproduct-estimates} below.}: The paracontrolled component $X_{\eps,\star}$ is determined by $X_{\eps,\star}=\Xwhite_{\eps,\star}+\Xblack_{\eps,\star}$, where 
\begin{align}\label{local:eq-Xwhite-K}
\Xwhite_{\eps,K} &=  \Duh_{B,K}^F \Big(  \Big[  \Wrem_{\eps,\star,j} \parall \bigl( 2 \covd^j \linear[\oeps][][][\ostar]- \covd \linear[\oeps][j][][\ostar] \bigr) \Big] \Big), \\
\Xblack_{\eps,K} &= \Duh_{B,K}^F \Big( \bil\bigl( \phiw_{\eps,\star} \parall \covd \philinear[\oeps][][][\ostar]\bigr) \Big).  \label{local:eq-Xblack-K}
\end{align}
The paracontrolled component $\phix_{\eps,\star}$ is determined by
\begin{equation}\label{local:eq-phix-K}
\begin{aligned}
\phix_{\eps,K} &= 2 \Duh_{B,K} \Big( \Wrem_{\eps,\star}^j \parall \covd_j \philinear[\oeps][][][\ostar] \Big). 
\end{aligned}
\end{equation}
The smooth nonlinear remainder $Y_{\eps,\star}$ is determined by\footnote{The equations for $Y_{\eps,\star}$ and $\phiy_{\eps,\star}$ will be written out again in more detail in \eqref{local:eq-Y-final-new-1}--\eqref{local:eq-y-final-new-last} below. For now, we encourage the reader to ignore all terms which are only referenced, but not written out here.} 
\begin{align}
Y_{\eps,K} = \Duh_{B,K}^F \Big(&  \Big[  \Wrem_{\eps,\star,j}, \bigl( 2 \covd^j \linear[\oeps][][][\ostar] - \covd \linear[\oeps][j][][\ostar]\bigr) \Big] - \Big[  \Wrem_{\eps,\star,j} \parall \bigl( 2 \covd^j \linear[\oeps][][][\ostar] - \covd \linear[\oeps][j][][\ostar]\bigr) \Big] \label{local:eq-Y-new-1} \\
&+  \Big[ \linear[\oeps][][j], 2 \covd^j \Wrem_\eps - \covd \Wrem^j_\eps \Big]\label{local:eq-Y-new-1p} \\ 
&+  \bil\Big( \phiw_{\eps,\star} \otimes \covd \philinear[\oeps][][][\ostar]\Big) 
-  \bil\Big( \phiw_{\eps,\star} \parall \covd \philinear[\oeps][][][\ostar]\Big)
+ \bil \Big( \philinear[\oeps]\otimes \covd \phiw_\eps \Big)
\label{local:eq-Y-new-2}\\ 
&+ \eqref{local:eq-Z-new-2} + \eqref{local:eq-Z-new-3}
+ \eqref{local:eq-Z-new-5} + \eqref{local:eq-Z-new-6}
+ \eqref{local:eq-Z-new-7} \Big).  \label{local:eq-Y-new-3} 
\end{align}
And finally, the smooth nonlinear remainder $\phiy_{\eps,\star}$ is determined by
\begin{align}
\phiy_{\eps,K} = \Duh_{B,K} \Big(& 2 \Wrem_{\eps,\star}^j \covd_j \philinear[\oeps][][][\ostar] 
- 2 \Wrem_{\eps,\star}^j \parall \covd_j \philinear[\oeps][][][\ostar]
+ 2 \linear[\oeps][j][][] \covd_j \phiw_\eps  \label{local:eq-y-new-1}\\ 
&+ \eqref{local:eq-z-new-2} + \eqref{local:eq-z-new-3} + \eqref{local:eq-z-new-4} \Big). \label{local:eq-y-new-2} 
\end{align}
In order to use the paracontrolled ansatz in our estimates of $Y_{\eps,K}$ and $\phiy_{\eps,K}$, we still need to rewrite the terms in \eqref{local:eq-Y-new-1}, \eqref{local:eq-Y-new-1p}, \eqref{local:eq-Y-new-2}, and \eqref{local:eq-y-new-1}. For this, we first prove the following lemma.

\begin{lemma}\label{local:lem-rewriting-high-high}
For $i=1,2,3$, let $\genV_i$ be a vector space equipped with an orthogonal representation of $\liegroup$, and let $j=1,2$. Furthermore, let $U_\star =(U_M)_{M \geq N_{\lo}}$ and $V_\star=(V_N)_{N\geq N_{\lo}}$, where, for each $M,N\geq N_{\lo}$, $U_M \colon [0,T]\times \T^2 \rightarrow \genV_1$ and $V_N \colon [0,T] \times \T^2 \rightarrow \genV_2$. If $\Qc\colon \genV_1 \times \genV_2 \rightarrow \genV_3$ is $\liegroup$-equivariant, 
then it holds that  
\begin{equation*}
\Qc \bigl( U_\star, \covd^j V_\star  \bigr)
= - \Qc \bigl( \covd^j  U_\star \parasim V_\star  \bigr)
+ \covd^j \Qc \bigl( U_\star \parasim   V_\star  \bigr)
+ \Qc \bigl( U_\star \paransim \covd^j V_\star  \bigr).
\end{equation*}
\end{lemma}

\begin{proof}
To obtain this identity, we first use the definition of our paraproducts and then use the covariant product rule for $\liegroup$-equivariant bilinear maps 
(see  Lemma~\ref{prelim:lem-equivariant-product-rule}). 
\end{proof}

Equipped with Lemma~\ref{local:lem-rewriting-high-high}, we can now rewrite the terms in \eqref{local:eq-Y-new-1}, \eqref{local:eq-Y-new-1p}, \eqref{local:eq-Y-new-2}, and \eqref{local:eq-y-new-1}. First, using Lemma~\ref{local:lem-rewriting-high-high} and the skew-symmetry of the Lie bracket, we see that \eqref{local:eq-Y-new-1}--\eqref{local:eq-Y-new-1p} can be written as 
\begin{align}
\big[  \Wrem_{\eps,\star,j}&, \bigl( 2 \covd^j \linear[\oeps][][][\ostar] - \covd \linear[\oeps][j][][\ostar]\bigr) \big] - \big[  \Wrem_{\eps,\star,j} \parall \bigl( 2 \covd^j \linear[\oeps][][][\ostar] - \covd \linear[\oeps][j][][\ostar]\bigr) \big]
+ \big[ \linear[\oeps][][j], 2 \covd^j \Wrem_\eps - \covd \Wrem^j_\eps \big]\notag \\
&= \big[  \Wrem_{\eps,\star,j}\parasim \bigl( 2 \covd^j \linear[\oeps][][][\ostar] - \covd \linear[\oeps][j][][\ostar]\bigr) \big]
+ \big[  \Wrem_{\eps,\star,j}\paragg \bigl( 2 \covd^j \linear[\oeps][][][\ostar] - \covd \linear[\oeps][j][][\ostar]\bigr) \big] \notag \\ 
&\quad -  2 \big[ \covd^j \linear[\oeps][][j][\ostar] \parasim \Wrem_{\eps,\star} \big]
+ 2 \covd^j  \big[ \linear[\oeps][][j][\ostar] \parasim \Wrem_{\eps,\star} \big]
+ 2 \big[ \linear[\oeps][][j][\ostar]\paransim \covd^j \Wrem_{\eps,\star} \big] \notag \\
&\quad + \big[ \covd \linear[\oeps][][j][\ostar] \parasim \Wrem_{\eps,\star}^j \big]
-\covd  \big[ \linear[\oeps][][j][\ostar] \parasim \Wrem_{\eps,\star}^j \big]
- \big[ \linear[\oeps][][j][\ostar] \paransim \covd \Wrem_{\eps,\star}^j \big] \notag \\
&= 2  \big[    \Wrem_{\eps,\star,j}  \parasim \bigl(  \covd^j \linear[\oeps][][][\ostar]- \covd \linear[\oeps][j][][\ostar] \bigr) \big]
+ 2 \big[ \Wrem_{\eps,\star} \parasim  \covd^j \linear[\oeps][][j][\ostar] \big] 
\label{local:eq-para-rewriting-1}\\
&\quad - 2 \covd^j  \big[  \Wrem_{\eps,\star} \parasim \linear[\oeps][][j][\ostar]\big] 
+ \covd  \big[  \Wrem_{\eps,\star}^j \parasim \linear[\oeps][][j][\ostar] \big] 
\label{local:eq-para-rewriting-2}\\
&\quad +\big[ \Wrem_{\eps,\star,j} \paragg \bigl( 2 \covd^j \linear[\oeps][][][\ostar] - \covd \linear[\oeps][j][][\ostar]\bigr)\big] + \big[ \linear[\oeps][][j][\ostar] \paransim \bigl( 2 \covd^j \Wrem_{\eps,\star} - \covd \Wrem_{\eps,\star}^j \bigr)\big]. 
\label{local:eq-para-rewriting-3}
\end{align}
The exact expression in \eqref{local:eq-para-rewriting-1} plays an important role in our argument, since it will need to be renormalised. Since it appears repeatedly below, it is convenient to define
\begin{equation}\label{local:eq-Lie-bracket-term}
\Lb\bigl(\Wrem_{\eps,\star},\covd \linear[\oeps][][][\ostar]\bigr)
\coloneq  2  \big[    \Wrem_{\eps,\star,j}  \parasim \bigl(  \covd^j \linear[\oeps][][][\ostar]- \covd \linear[\oeps][j][][\ostar] \bigr) \big]
+ 2 \big[ \Wrem_{\eps,\star} \parasim  \covd^j \linear[\oeps][][j][\ostar] \big]. 
\end{equation}
We now turn to the terms in \eqref{local:eq-Y-new-2}. Using Lemma~\ref{local:lem-rewriting-high-high} and the skew-symmetry of $\bil$, we obtain 
\begin{align}
\bil\bigl( \phiw_{\eps,\star} \otimes \covd &\philinear[\oeps][][][\ostar]\bigr)
- \bil\bigl( \phiw_{\eps,\star} \parall \covd \philinear[\oeps][][][\ostar]\bigr)
+ \bil \bigl( \philinear[\oeps]\otimes \covd \phiw_\eps \bigr) \notag \\ 
&= 2 \bil \bigl( \phiw_{\eps,\star} \parasim \covd \philinear[\oeps][][][\ostar]\bigr) 
\label{local:eq-para-rewriting-4} \\
&\quad - \covd \Big( \bil\bigl( \phiw_{\eps,\star} \parasim \philinear[\oeps][][][\ostar] \bigr)\Big)
+ \bil\bigl( \phiw_{\eps,\star} \paragg \covd \philinear[\oeps][][][\ostar] \bigr)
+ \bil\bigl( \philinear[\oeps][][][\ostar] \paransim  \covd \phiw_{\eps,\star}  \bigr). 
\label{local:eq-para-rewriting-5} 
\end{align}
Finally, we turn to the terms in \eqref{local:eq-y-new-1}. Using Lemma~\ref{local:lem-rewriting-high-high} again, we decompose
\begin{equation}\label{local:eq-para-rewriting-6}
 2 \Wrem_{\eps,\star}^j  \covd_j \philinear[\oeps][][][\ostar]  
 -  2 \Wrem_{\eps,\star}^j \parall \covd_j \philinear[\oeps][][][\ostar]
 =  2 \Wrem_{\eps,\star}^j \parasim \covd_j \philinear[\oeps][][][\ostar]  
 +  2 \Wrem_{\eps,\star}^j \paragg \covd_j \philinear[\oeps][][][\ostar]  
\end{equation}
and 
\begin{equation}\label{local:eq-para-rewriting-7}
2  \linear[\oeps][j][][] \covd_j \phiw_\eps
= - 2 \bigl( \covd_j \linear[\oeps][j][][\ostar] \bigr) \parasim    \phiw_{\eps,\star}
+ 2 \covd_j \bigl( \linear[\oeps][j][][\ostar]  \parasim  \phiw_{\eps,\star} \bigr)
+ 2   \linear[\oeps][j][][\ostar] \paransim  \covd_j \phiw_{\eps,\star}.
\end{equation}
In order to finish our reformulation of \eqref{local:eq-Y-new-1}--\eqref{local:eq-y-new-2}, it only remains to complete two steps: First, we insert our paracontrolled ansatz into the high$\times$high-products in~\eqref{local:eq-para-rewriting-1}, \eqref{local:eq-para-rewriting-4}, the first summand in~\eqref{local:eq-para-rewriting-6}, and the first summand in~\eqref{local:eq-para-rewriting-7}, which is done using the double Duhamel trick. For this, we need to introduce operators corresponding to the right-hand sides of \eqref{local:eq-Xwhite-K}, \eqref{local:eq-Xblack-K}, and \eqref{local:eq-phix-K}, which are defined as 
\begin{align}
\Xwhiteop_{\eps,K}(\Wrem_{\eps,\star}) &\eqdef  \Duh_{B,K}^F \Big(  \big[  \Wrem_{\eps,\star,j} \parall \bigl( 2 \covd^j \linear[\oeps][][][\ostar]- \covd \linear[\oeps][j][][\ostar] \bigr) \big] \Big), 
\label{local:eq-Xwhiteop-K}\\
\Xblackop_{\eps,K}(\phiw_{\eps,\star}) &\eqdef \Duh_{B,K}^F \Big( \bil\bigl( \phiw_{\eps,\star} \parall \covd \philinear[\oeps][][][\ostar]\bigr) \Big),   
\label{local:eq-Xblackop-K}\\
\phixop_{\eps,K}(\Wrem_{\eps,\star}) &\eqdef 2 \Duh_{B,K} \Big( \Wrem_{\eps,\star}^j \parall \covd_j \philinear[\oeps][][][\ostar] \Big). \label{local:eq-phixop-K}
\end{align}
In line with our earlier decomposition of $X_{\eps,K}$ into $\Xwhite_{\eps,K}$ and $\Xblack_{\eps,K}$, we also write
\begin{equation*}
\Xop_{\eps,K}(\Wrem_{\eps,\star},\phiw_{\eps,\star}) \eqdef  
\Xwhiteop_{\eps,K}(\Wrem_{\eps,\star})+ \Xblackop_{\eps,K}(\phiw_{\eps,\star}). 
\end{equation*}
Second, we renormalise the individual terms in the nonlinearity, and the counterterms will involve the constants $\sigma_{\rho,\eps}^2$ and $\Crho$ from Definition \ref{main:def-renormalization-constant} and the induced linear maps from Definition~\ref{prelim:def-Casimir-induced}. 
In total, we then arrive at the following system of equations: As in \eqref{local:eq-Xwhite-K}--\eqref{local:eq-phix-K} (or as on the right-hand sides of \eqref{local:eq-Xwhiteop-K}--\eqref{local:eq-phixop-K}), the equations for the paracontrolled components, which we repeat for convenience, are
\begin{align}
\Xwhite_{\eps,K} &=  \Duh_{B,K}^F \Big(  \big[  \Wrem_{\eps,\star,j} \parall \bigl( 2 \covd^j \linear[\oeps][][][\ostar]- \covd \linear[\oeps][j][][\ostar] \bigr) \big] \Big), 
\label{local:eq-Xwhite-K-final}\\
\Xblack_{\eps,K} &= \Duh_{B,K}^F \Big( \bil\bigl( \phiw_{\eps,\star} \parall \covd \philinear[\oeps][][][\ostar]\bigr) \Big),   
\label{local:eq-Xblack-K-final}\\
\phix_{\eps,K} &= 2 \Duh_{B,K} \Big( \Wrem_{\eps,\star}^j \parall \covd_j \philinear[\oeps][][][\ostar] \Big). \label{local:eq-phix-K-final}
\end{align}
For the smooth nonlinear remainder $Y_{\eps,\star}$, we have that 
\begin{align}
&\quad Y_{\eps,K} = \Duh_{B,K}^F \Big( \notag \\ 
&\quad \Lb\bigl( \quadratic[\oeps][][][\ostar],\covd \linear[\oeps][][][\ostar]\bigr) 
+ \hlmath{( \sigma_{\rho,\eps}^2+\Crho) \Cov_{\frkg}^{\frkg} (\linear[\oeps][][][])}
+ \Lb\bigl( \phiquadratic[\oeps][][][\ostar],\covd \linear[\oeps][][][\ostar]\bigr)
\label{local:eq-Y-final-new-1}\\
&+ 2 \bil \bigl( \mixedquadratic[\oeps][][][\ostar] \parasim \covd \philinear[\oeps][][][\ostar]\bigr) 
+\hlmath{(\sigma_{\rho,\eps}^2+\Crho) \Cov^{\frkg}_{\V}( \linear[\oeps])}
\allowdisplaybreaks[4]\label{local:eq-Y-final-new-1p}\\
&+\Lb\bigl(\Xop_{\eps,\star}(\Wrem_{\eps,\diamond},\phiw_{\eps,\diamond}),\covd \linear[\oeps][][][\ostar]\bigr) +  2 \bil \bigl( \phixop_{\eps,\star}(\Wrem_{\eps,\diamond}) \parasim \covd \philinear[\oeps][][][\ostar]\bigr) +
\hlmath{(\sigma_{\rho,\eps}^2+\Crho) \bigl( \Cov_{\frkg}^{\frkg} + \Cov^{\frkg}_{\V}\bigr)(W_\eps)}  \label{local:eq-Y-final-new-2}\allowdisplaybreaks[4]\\
&+ \Lb\bigl(\Sdet_\star+Y_{\eps,\star},\covd \linear[\oeps][][][\ostar]\bigr) +  2 \bil \bigl( (\higgsdet_\star+\phiy_{\eps,\star}) \parasim \covd \philinear[\oeps][][][\ostar]\bigr) 
\label{local:eq-Y-final-new-3}\allowdisplaybreaks[4]\\
&-  2 \covd^j  \big[  \Wrem_{\eps,\star} \parasim \linear[\oeps][][j][\ostar]\big] 
+  \covd  \big[  \Wrem_{\eps,\star}^j \parasim \linear[\oeps][][j][\ostar] \big] 
-  \covd \bigl( \bil\bigl( \phiw_{\eps,\star} \parasim \philinear[\oeps][][][\ostar] \bigr)\bigr)
\label{local:eq-Y-final-new-4}\allowdisplaybreaks[4]\\
&+  \big[ \Wrem_{\eps,\star,j} \paragg \bigl( 2 \covd^j \linear[\oeps][][][\ostar] - \covd \linear[\oeps][j][][\ostar]\bigr)\big] +\bil\bigl( \phiw_{\eps,\star} \paragg \covd \philinear[\oeps][][][\ostar] \bigr) 
\label{local:eq-Y-final-new-5}\allowdisplaybreaks[4]\\
&+   \big[ \linear[\oeps][][j][\ostar] \paransim \bigl( 2 \covd^j \Wrem_{\eps,\star} - \covd \Wrem_{\eps,\star}^j \bigr)\big] 
+ \bil\bigl( \philinear[\oeps][][][\ostar] \paransim  \covd \phiw_{\eps,\star}  \bigr) 
\label{local:eq-Y-final-new-6}\allowdisplaybreaks[4]\\
&+ \big[ \Wrem_{\eps,j}  , 2 \covd^j \Wrem_\eps - \covd \Wrem_\eps^j \big] - \big[ \Sdet_j , 2 \covd^j \Sdet - \covd \Sdet^j \big] + \bil \bigl( \phiw_\eps \otimes \covd \phiw_\eps\bigr) - \bil \bigl( \higgsdet \otimes \covd \higgsdet\bigr)
\label{local:eq-Y-final-new-7}\allowdisplaybreaks[4]\\
&+  \big[  \linear[\oeps][][j][]\hspace{-0.5ex}+ \Wrem_{\eps,j}, 
\big[ \linear[\oeps][j][][] \hspace{-0.5ex}+ \Wrem_\eps^j,  \linear[\oeps] \hspace{-0.5ex}+ \Wrem_\eps \big] \big] - \big[ \Sdet_j, \big[ \Sdet^j , \Sdet \big]\big] - \hlmath{\sigma_{\rho,\eps}^2 \Cov_\frkg^{\frkg}( \linear[\oeps][][][] + \Wrem_\eps)} 
\label{local:eq-Y-final-new-8}\allowdisplaybreaks[4]\\
&+  \bil \bigl( \bigl( \philinear[\oeps]+\phiw_\eps\bigr) \otimes \bigl(\linear[\oeps]+\Wrem_\eps\bigr) \bigl( \philinear[\oeps] + \phiw_\eps\bigr) \bigr) - \bil\bigl( \higgsdet \otimes \Sdet \higgsdet \bigr) - \hlmath{\sigma_{\rho,\eps}^2 \Cov_{\V}^{\frkg} ( \linear[\oeps][][][]  + \Wrem_\eps)} 
\label{local:eq-Y-final-new-9}\allowdisplaybreaks[4]\\ 
&- \hlmath{\tfrac{1}{8\pi} \Cov_{\V}^{\frkg} \bigl( \linear[\oeps] + W_\eps \bigr)} \Big)\;,
\label{local:eq-Y-final-new-10}
\end{align}
where we have highlighted the counterterms in light grey.
In \eqref{local:eq-Y-final-new-2}, we used the symbols $\diamond$ and $\star$ for the frequency scales of the input and output, 
respectively, of $\Xop$ and $\phixop$. We emphasise that the terms within each line of \eqref{local:eq-Y-final-new-1}--\eqref{local:eq-Y-final-new-9} are similar to each other. Finally, for the smooth remainder $\phiy_{\eps,\star}$, we have that
\begin{align}
 \phiy_{\eps,K} &= \Duh_{B,K} \Big(
2 \bigl( \quadratic[\oeps][][][\ostar] + \phiquadratic[\oeps][][][\ostar]\bigr)^j \parasim \covd_j \philinear[\oeps][][][\ostar]
+ \hlmath{(\sigma_{\rho,\eps}^2 + \Crho) \Cov^{\V}_{\frkg}( \philinear[\oeps])} 
-  2 \bigl( \covd_j \linear[\oeps][j][][\ostar] \bigr) \parasim    \mixedquadratic [\oeps][][][\ostar]
\label{local:eq-y-final-new-1}\allowdisplaybreaks[4]\\
&\quad+ 
2 \Xop_{\eps,\star}^j(\Wrem_{\eps,\diamond},\phiw_{\eps,\diamond}) \parasim \covd_j \philinear[\oeps][][][\ostar]
+ \hlmath{(\sigma_{\rho,\eps}^2 + \Crho) \Cov^{\V}_{\frkg}( \phiw_\eps)}
- 2 \bigl( \covd_j \linear[\oeps][j][][\ostar] \bigr) \parasim    \phixop_{\eps,\star}(\Wrem_{\eps,\diamond})
\label{local:eq-y-final-new-2}\allowdisplaybreaks[4]\\
&\quad+ 2 \bigl( \Sdet_{\star} + Y_{\eps,\star}\bigr)^j 
 \parasim \covd_j \philinear[\oeps][][][\ostar] 
 - 2 \bigl( \covd_j \linear[\oeps][j][][\ostar] \bigr) \parasim    \bigl( \higgsdet_{\star} + \phiy_{\eps,\star} \bigr)
\label{local:eq-y-final-new-3}\allowdisplaybreaks[4]\\
&\quad+ 2 \covd_j \bigl(   \linear[\oeps][j][][\ostar]  \parasim  \phiw_{\eps,\star} \bigr)
+ 2 \Wrem_{\eps,\star}^j \paragg \covd_j \philinear[\oeps][][][\ostar] 
+ 2   \linear[\oeps][j][][\ostar] \paransim  \covd_j \phiw_{\eps,\star}
\label{local:eq-y-final-new-4}\allowdisplaybreaks[4]\\
&\quad+ 2 \bigl(\Wrem_\eps^j \covd_j \phiw_\eps - \Sdet^j \covd_j \higgsdet\bigr) 
\label{local:eq-y-final-new-5}\allowdisplaybreaks[4]\\
&\quad+ \bigl( \linear[\oeps]+\Wrem_\eps\bigr)^j \bigl( \linear[\oeps]+\Wrem_\eps\bigr)_j \bigl( \philinear[\oeps]+\phiw_\eps\bigr)- \Sdet^j \Sdet_j \higgsdet - \hlmath{2 \sigma_{\rho,\eps}^2 \Cov_{\frkg}^{\V}( \philinear[\oeps][][][] + \phiw_\eps )} 
\label{local:eq-y-final-new-6}\allowdisplaybreaks[4]\\
&\quad- \big| \philinear[\oeps] + \phiw_\eps\big|^2 \bigl( \philinear[\oeps]+\phiw_\eps\bigr) + \hlmath{(\dim(\V)+2) \sigma_{\rho,\eps}^2 (\philinear[\oeps]+\phiw_\eps)}  + |\higgsdet|^2 \higgsdet \label{local:eq-y-final-new-7}\\ 
&\quad - \hlmath{\mass \bigl( \philinear[\oeps] + \phiw_\eps\bigr)} \Big). \label{local:eq-y-final-new-last}
\end{align}

We note that, by Definition \ref{main:def-renormalization-constant}, adding up the counterterms in \eqref{local:eq-Y-final-new-1}--\eqref{local:eq-Y-final-new-10} gives exactly the counterterm in~\eqref{local:eq-Z-new-7}, and adding up the counterterms in \eqref{local:eq-y-final-new-1}--\eqref{local:eq-y-final-new-last} gives exactly the counterterm in~\eqref{local:eq-z-new-4}.

\begin{remark}
By writing out the integral equations \eqref{local:eq-Y-final-new-1}--\eqref{local:eq-y-final-new-last}, we arrived at rather long expressions. The precise form of most terms is irrelevant, i.e., the precise indices, prefactors, or signs are not too important. However, in some terms such as \eqref{local:eq-Y-final-new-1}--\eqref{local:eq-Y-final-new-3}, the precise form is crucial for our argument. For this reason, we believe it is best to write out the integral equations for $Y_{\eps,\star}$ and $\phiy_{\eps,\star}$ explicitly once.
\end{remark}

\begin{remark}
In \eqref{local:eq-Y-final-new-10}, we left the $\frac{1}{8\pi}$-counterterm alone. The reason is that it will only be important in Section~\ref{section:langevin}, where it will be needed to prove that our solution is gauge equivalent (in law) to the solution obtained in \cite{CCHS22}. In contrast, we grouped the $\Crho$-counterterms together with their corresponding interactions. For the arguments in Sections~\ref{section:covariant-Ansatz}--\ref{section:objects}, one could leave the $\Crho$-counterterms alone as well, and only pay attention to them in Section~\ref{section:langevin}. However, we would like to illustrate that, even though it is only proven in Section~\ref{section:langevin}, the arguments in Sections~\ref{section:covariant-Ansatz}--\ref{section:objects} could be used to prove that the limiting solution does not depend on $\rho$. And, for this, we need to correctly account for the $\Crho$-counterterms. 
\end{remark}

\section{Solution spaces, Schauder estimates, and paraproduct estimates}
\label{sec:covariant_function_spaces}

In this section, we introduce the covariant functional framework that will be needed in the proof of our covariant local well-posedness result (Theorem~\ref{main:thm-covariant-local}). Throughout the section, we let $0< T \leq \frac{1}{128}$, let $N_{\lo}$ be as in~\eqref{local:eq-N-lo}, and let $B$ be the Yang--Mills gradient flow from~\eqref{local:eq-B}. Furthermore, as in Notation~\ref{local:notation-omit-B}, we simply write $F$ for $F_B$. 

\subsection{Solution spaces}

In this subsection, we introduce the functional framework that will be used to solve \eqref{local:eq-Xwhite-K-final}--\eqref{local:eq-y-final-new-last}. In order to treat the connection one-form and Higgs field simultaneously, we introduce the following notation.

\begin{notation}\label{local:notation-unified} 
First, we introduce $\Y\eqdef \frkg^2 \oplus \V=\frkg \oplus \frkg \oplus \V$. This space is equipped with the natural representation of $\liegroup$, i.e., the representation $\rep_\Y=\Ad \oplus \Ad \oplus \rep$. Occasionally, we will need to refer to the different components of $\Y$. For this reason, it is convenient to write 
\begin{equation*}
\Y^{(1)}=\Y^{(2)} \eqdef \frkg \qquad \text{and} \qquad \Y^{(3)}\eqdef \V. 
\end{equation*}
For each $N\in \dyadic$, we define the $\End(\Y)$-valued kernels
\begin{equation}\label{object:eq-qb-Y}
\begin{aligned}
\fq_{\Y,B}(\place,\place) &\eqdef \fq_{B}(\place,\place;\frkg^2) \oplus \fp_{B}(\place,\place;\V), \\ 
\fq_{\Y,B,\eps,N}(\place,\place) &\eqdef \fq_{B,\eps,N}(\place,\place;\frkg^2) \oplus \fp_{B,\eps,N}(\place,\place;\V),\\
\fq_{\Y,B,N}(\place,\place) &\eqdef \fq_{B,N}(\place,\place;\frkg^2) \oplus \fp_{B,N}(\place,\place;\V), 
\end{aligned}
\end{equation}
where $\fq_B(\place,\place;\frkg^2)$ and $\fp_B(\place,\place;\V)$ are the $\End(\frkg^2)$- and $\End(\V)$-valued kernels from Definition~\ref{kernel:def-kernels} and their decorated versions are as in \eqref{local:eq-qB-eps-N} and \eqref{local:eq-qB-N}. We emphasise that the second kernel in \eqref{object:eq-qb-Y} is a $\fp_B$- rather than $\fq_B$-kernel, and the reason for this is that the $\higgs$-equation in \eqref{local:eq-C} does not contain a curvature term. We will sometimes also use the $\End(\Y)$-valued kernels
\begin{equation}
\begin{aligned}\label{object:eq-pb-Y}
\fp_{\Y,B}(\place,\place) &\eqdef \fp_{B}(\place,\place;\frkg^2)  \oplus \fp_{B}(\place,\place;\V),\\
\fp_{\Y,B,\eps,N}(\place,\place) &\eqdef \fp_{B,\eps,N}(\place,\place;\frkg^2)  \oplus \fp_{B,\eps,N}(\place,\place;\V),\\
\fp_{\Y,B,N}(\place,\place) &\eqdef \fp_{B,N}(\place,\place;\frkg^2)  \oplus \fp_{B,N}(\place,\place;\V),
\end{aligned}
\end{equation}
where the first kernel on each right-hand side is $\End(\frkg^2)$-valued and the second kernel is $\End(\V)$-valued. In line with \eqref{object:eq-qb-Y} and \eqref{object:eq-pb-Y}, we further write
\begin{equation}\label{object:eq-Duh-Y}
\begin{aligned}
\Duh_{\Y,B,N}^{F}(\Psi)(t,x)&\eqdef \iint \ds \dy\,  \fq_{\Y,B,N}(s,y;t,x) \Psi(s,y), \\  
\Duh_{\Y,B,N}(\Psi)(t,x) &\eqdef \iint \ds \dy \, \fp_{\Y,B,N}(s,y;t,x) \Psi(s,y),
\end{aligned}
\end{equation}
where $\Psi\colon [0,T]\times\T^2 \rightarrow \Y$. 
\end{notation}

We note that, while Notation~\ref{local:notation-unified} will already be useful in this section, it will be most useful in our treatment of stochastic objects, see e.g.\ Sections~\ref{section:covariant-data-set} and \ref{section:objects}. We proceed to introduce our solution spaces. In the following, recall Notation~\ref{notation:space-time-derivatives} regarding space-time derivatives.

\begin{definition}[Solution spaces]\label{local:def-solution-space}
Let $\alpha\in \R$, and let $\eta \leq 0$. For all $M\in \dyadic$ satisfying $M\geq N_{\lo}$ and all $\phi_M \colon [0,T]\times \T^2 \rightarrow \Y$, we then define
\begin{equation}\label{local:eq-solution-space}
\| \phi_M \|_{\sol^{\alpha,\eta}_M} 
= \| \phi_M \|_{\sol^{\alpha,\eta}_M([0,T])}
\coloneq \max_{0\leq m \leq 10} M^{\alpha-m}  \big\| t^{-\eta/2}\covd_{t,x}^{(m)} \phi_M \big\|_{ L_{t,x}^\infty([0,T]\times \T^2)}. 
\end{equation}
Furthermore, let $\phi_\star=(\phi_N)_{N\geq N_{\lo}}$ where, for each $N\geq N_{\lo}$, $\phi_N \colon [0,T] \times \T^2 \rightarrow \Y$. Then, we define 
\begin{equation}\label{local:eq-solution-space-sequence}
\| \phi_\star\|_{\sol^{\alpha,\eta}} 
= \| \phi_\star\|_{\sol^{\alpha,\eta}([0,T])}
\coloneq \sup_{N\geq N_{\lo}} \| \phi_N \|_{\sol^{\alpha,\eta}_N([0,T])}. 
\end{equation}
We also write $\| \phi_M \|_{\sol^{\alpha}_M}\coloneq \| \phi_M \|_{\sol^{\alpha,0}_M}$ and $\| \phi_\star \|_{\sol^{\alpha}}\coloneq \| \phi_\star \|_{\sol^{\alpha,0}}$. 
If the maps take values in the component $\Y^{(a)}$ rather than $\Y$ itself, where $a\in [3]$, the norms from \eqref{local:eq-solution-space} and \eqref{local:eq-solution-space-sequence} are defined similarly. 
\end{definition}

\begin{remark}\label{local:rem-solution-space-scaling}
In Definition~\ref{local:def-solution-space}, it is useful to think of $\alpha+\eta$ as a scaling exponent. Indeed, in the case $B\equiv 0$, $\| \phi_\star \|_{\sol^{\alpha,\eta}}$ and $\| \phi(0)\|_{\Cs_x^{\alpha+\eta}}$ scale the same under the parabolic scaling. Keeping this in mind will be helpful in understanding the parameter conditions in Section~\ref{section:local-paraproduct-estimates} and the powers of $T$ in Section~\ref{section:lwp-nonlinear-estimate}. 
\end{remark}

The solution spaces $\sol^{\alpha,\eta}_M$ will be used to control functions which, in the covariant sense described in Section~\ref{section:Littlewood--Paley}, have regularity $\alpha$ and are supported on frequencies 
of order $ M$. For negative values of $\alpha$, we will often need to work with strengthened versions of \eqref{local:eq-solution-space} and~\eqref{local:eq-solution-space-sequence}. The reason is that, when $\alpha$ is negative, the exponent of $M$ in \eqref{local:eq-solution-space} is always negative, and hence~\eqref{local:eq-solution-space} does not allow us to gain powers of $M$. In the following, recall the temporal regularity spaces introduced in Definition~\ref{def:temporal-regularity-space}.

\begin{definition}[Modified solution spaces]\label{local:def-modified-solution-space}
Let $\alpha \in (-1,0)$. For all $M\in \dyadic$ satisfying $M\geq N_{\lo}$ and all $\phi_M \colon [0,T]\times \T^2\rightarrow \Y$, we then define
\begin{equation}\label{local:eq-modified-solution-space}
\begin{aligned}
\big\| \phi_M \big\|_{\modsol^{\alpha}_M} 
&= \big\| \phi_M \big\|_{\modsol^{\alpha}_M([0,T])} \\
&\coloneq \big\| \phi_M \big\|_{\sol^\alpha_M([0,T])}  
+ \sup_{-1\leq \alpha^\prime \leq \alpha} \max_{0\leq m \leq 10} M^{-\alpha^\prime+\alpha-m} \big\| \covd^{(m)}_{t,x} \phi_M \big\|_{\Cs_t^{\alpha^\prime/2} L_x^r}, 
\end{aligned}
\end{equation}
where $r$ is as in \eqref{e:small_constants}. 
Furthermore, let $\phi_\star=(\phi_N)_{N\geq N_{\lo}}$ where, for each $N\geq N_{\lo}$, $\phi_N \colon [0,T] \times \T^2 \rightarrow \Y$. Then, we define 
\begin{equation}\label{local:eq-modified-solution-space-sequence}
\| \phi_\star\|_{\modsol^{\alpha}} 
= \| \phi_\star\|_{\modsol^{\alpha}([0,T])}
\coloneq \sup_{N\geq N_{\lo}} \big\| \phi_N \big\|_{\modsol^{\alpha}_N([0,T])}. 
\end{equation}
If the maps take values in the component $\Y^{(a)}$ rather than $\Y$ itself, where $a\in [3]$, the norms are again defined similarly. 
\end{definition}
We note that the exponent of $M$ in \eqref{local:eq-modified-solution-space} is positive when $\alpha^\prime<\alpha$ and $m=0$, which will often be used below.

\begin{remark}
As mentioned before, the most important difference between our $\sol^{\alpha,\eta}$-norm and our $\modsol^{\alpha}$-norm is that the latter includes terms at negative temporal regularity, but there are also two other differences: First, our $\modsol^{\alpha}$-norm does not include singular weights in time. This could be done but, since the $\modsol^{\alpha}$-norm will only be used for our covariant stochastic objects, it is simply not needed. Second, while the $\modsol^{\alpha}$-norm includes $L_x^r$-terms, the $\sol^{\alpha,\eta}$-norm is entirely $L_x^\infty$-based. This difference is not reflected in our notation for the $\modsol^{\alpha}$-norm since the parameter $r$ is fixed throughout the article. The reason for using the $L_x^r$-norm is only technical and lies in the proof of Lemma~\ref{lemma:time-regularity-kolmogorov} below. 
\end{remark}

In the following, we state and prove two lemmas related to properties of our solution spaces. In the first lemma, we prove a Hölder estimate in time. 

\begin{lemma}[Hölder estimate in time]\label{local:lem-solution-spaces-hoelder}
Let $\alpha\in \R$, let $\eta \in (-2,0]$, let $0\leq \beta\leq 2$, and let $0\leq m \leq 8$. Furthermore, let $M\in \dyadic$ satisfy $M\geq N_{\lo}$ and let $\phi_M \colon [0,T]\times \T^2 \rightarrow \Y$. Finally, assume that $\phi_M(t)=0$ for all $t\ll M^{-2}$. Then, it holds that 
\begin{equation}\label{local:eq-solution-spaces-hoelder}
\big\| t^{-\eta/2} \covd_{t,x}^{(m)} \phi_M \big\|_{\Cs_t^{\frac{\beta}{2}} L_x^\infty} \lesssim M^{-\alpha+\beta+m}\big\| \phi_M \big\|_{\sol^{\alpha,\eta}_M}.
\end{equation}
\end{lemma}

\begin{proof}
By interpolation, it suffices to treat the cases $\beta=0$ and $\beta=2$. The case $\beta=0$ is trivial. In the case $\beta=2$, we estimate
\begin{align*}
  \big\| \partial_t \bigl( t^{-\eta/2} \covd_{t,x}^{(m)} \phi_M\bigr) \big\|_{L_{t,x}^\infty}
&\lesssim \big\| t^{-\eta/2-1} \covd_{t,x}^{(m)} \phi_M  \big\|_{L_{t,x}^\infty} + \big\| t^{-\eta/2} \covd_{t,x}^{(m+2)} \phi_M  \big\|_{L_{t,x}^\infty} \\
&\lesssim M^2  \big\| t^{-\eta/2} \covd_{t,x}^{(m)} \phi_M  \big\|_{L_{t,x}^\infty} + \big\| t^{-\eta/2} \covd_{t,x}^{(m+2)} \phi_M  \big\|_{L_{t,x}^\infty} \\ 
&\lesssim M^{2-\alpha+m} \big\| \phi_M \big\|_{\sol^{\alpha,\eta}_M}\;,
\end{align*}
where we used that $t^{-1} \lesssim M^2$ on the support of $\phi_M$ to obtain the second inequality. 
\end{proof}

In the second lemma, we prove an estimate that relates solution spaces with different regularities and singular weights.

\begin{lemma}\label{local:lem-sol-different-regularities-weights}
Let $\alpha_1,\alpha_2\in \R$, and let $\eta_1,\eta_2\leq 0$. Furthermore, assume that $\alpha_1+\eta_1\leq \alpha_2+\eta_2$ and that $\eta_1\geq \eta_2$. Finally, let $N\geq N_{\lo}$ and let $\phi_N \colon [0,T] \times \T^2 \rightarrow \Y$ be supported on times $t\gtrsim N^{-2}$. Then, it holds that 
\begin{equation*}
\| \phi_N \|_{\sol^{\alpha_1,\eta_1}_N}
\lesssim T^{\frac{\alpha_2+\eta_2-(\alpha_1+\eta_1)}{2}} \| \phi_N \|_{\sol^{\alpha_2,\eta_2}_N}\;,
\end{equation*}
where the implicit constant is independent of $T$ and $N$.
\end{lemma}

Roughly speaking, this shows that, at the cost of increasing the regularity parameter, we can make the weight more singular.

\begin{proof}
On the support of $\phi_N$, we have that
\begin{equation*}
N^{\alpha_1} t^{-\frac{\eta_1}{2}}
= N^{\alpha_1} t^{\frac{\eta_2-\eta_1}{2}} t^{-\frac{\eta_2}{2}}
\lesssim N^{\alpha_1} N^{-(\eta_2-\eta_1)} t^{-\frac{\eta_2}{2}}
= N^{\alpha_2} N^{\alpha_1+\eta_1-\alpha_2-\eta_2} t^{-\frac{\eta_2}{2}}. 
\end{equation*}
Since $\alpha_1+\eta_1-\alpha_2-\eta_2\leq0$ and $N\gtrsim N_{\lo}\sim T^{-\frac{1}{2}}$, this implies the desired result. 
\end{proof}

After this discussion of the solution spaces themselves, we will soon turn to the corresponding Schauder estimates. To simplify the notation below, we first define the following norm, which will be used for the nonlinearities in \eqref{local:eq-Xwhite-K-final}--\eqref{local:eq-y-final-new-last}.

\begin{definition}\label{local:def-nl-norm}
Let $\alpha\in \R$, and let $\nu\leq 0$. Furthermore, let $M\in \dyadic$ satisfy $M\geq N_{\lo}$ and let $\phi_M \colon [0,T]\times \T^2 \rightarrow \Y$. Then, we define
\begin{equ}\label{local:eq-nl-norm-M}
\big\| \phi_M \big\|_{\nl_M^{\alpha,\nu}} 
= \big\| \phi_M \big\|_{\nl_M^{\alpha,\nu}([0,T])} 
\coloneq 
\bigg\| \iint \ds \dy \, \fq_{\Y,B,M}(s,y;t,x) \phi_M(s,y) \bigg\|_{\sol^{\alpha+2,\nu}_M([0,T])}. 
\end{equ}
Furthermore, let $\phi \colon [0,T] \times \T^2 \rightarrow \Y$. Then, we define 
\begin{equation}\label{local:eq-nl-norm} 
 \| \phi \|_{\nl^{\alpha,\nu}} = \| \phi \|_{\nl^{\alpha,\nu}([0,T])} \coloneq  \sup_{N\geq N_{\lo}} \| \phi \|_{\nl_N^{\alpha,\nu}([0,T])}.
\end{equation}
As for our solution spaces, we write $\nl^{\alpha}$ and $\nl^{\alpha}_M$ in the case $\nu=0$. If $\phi_M$ or $\phi$ takes values in $\frkg^2$ or $\V$ rather than $\Y$, we define $\big\| \phi_M \big\|_{\nl_M^{\alpha,\nu}}$ and $\| \phi \|_{\nl^{\alpha,\nu}}$ by first using the canonical embedding $\frkg^2 \hookrightarrow \Y$ or $\V \hookrightarrow \Y$ and then using the norms for $\Y$-valued maps.
\end{definition}
 We also remark that it is natural that the map $\phi$ in \eqref{local:eq-nl-norm} does not depend on $N$, since this is true for the nonlinearities in \eqref{local:eq-Xwhite-K-final}--\eqref{local:eq-y-final-new-last}. 

\subsection{Schauder and commutator estimates}\label{section:local-schauder-commutator-estimates}

In this subsection, we prove (covariant) Schauder estimates in our solution spaces, i.e., in the spaces from Definitions~\ref{local:def-solution-space} and~\ref{local:def-nl-norm}. We also prove commutator estimates (which are stated in $L_{t,x}^\infty$ spaces). The Schauder estimates will be an important ingredient in our nonlinear estimates, which will be the subject of Section~\ref{section:lwp-nonlinear-estimate}, whereas the commutator estimates will be needed in uncovering the structure of the paracontrolled terms \dash see Section~\ref{section:local-paracontrolled}. We first define the following quantity, which will absorb all losses in the size of $F$ arising from our covariant heat kernel estimates, such as Propositions~\ref{prop:pB-ptwise-bounds},~\ref{prop:qB-ptwise-bds}, and~\ref{prop:qB-mixed-derivatives-ptwise-bounds}. In the following, let $C_0$ be the maximum of the constants $C_0$ appearing in the statements of Lemma~\ref{lemma:higher-regularity-abstract-small-Lt1-Lxinfty}, Proposition~\ref{prop:qB-ptwise-bds}, and Proposition~\ref{prop:covariant-heat-equation-with-curvature-higher-regularity}.

\begin{notation}
In analogy with Notation~\ref{notation:gronwall-loss}, we define
\begin{equation}\label{local:eq-lossfactor}
\lossfactor_T(F)\eqdef \exp\bigl( C_0 \| F \|_{L_t^1 L_x^\infty([0,T])}\bigr)
\bigl( 1+  \Fnorm{F}_{T,100} + T^{\frac{1-\roughdataweight}{2}} \Fnorm{t^{-\frac{1-\roughdataweight}{2}} F}_{T, 1} \bigr)^{100}. 
\end{equation}
\end{notation}

In contrast to Notation~\ref{notation:gronwall-loss}, we have also included the constant $\gronwallconstant$ in the exponential, which is notationally convenient. 

In the following lemma, we state our unweighted Schauder estimates. While we will later prove weighted Schauder estimates as well, the unweighted Schauder estimates are much less technical. For expository purposes, we therefore split the unweighted and weighted Schauder estimates into Lemma~\ref{local:lem-covariant-Schauder} and Proposition~\ref{local:prop-weighted-covariant-Schauder}, respectively. 

\begin{lemma}[Covariant Schauder estimates]\label{local:lem-covariant-Schauder}
Let $N \geq N_{\lo}$, let $\alpha\in  \R$,  let $\beta \in (-2,0]$, and let $1 \leq r \leq \infty$. Furthermore, 
let $\Qc\colon \Y \times \Y \rightarrow \Y$ be a $\liegroup$-equivariant bilinear map, let $\phi \colon [0,T]\times \T^2 \rightarrow \Y$, and let $\psi\colon [0,T]\times \T^2 \rightarrow \Y$. Then, we have the estimates 
\begin{align}
\max_{m=0,1} N^{-\alpha-m} \big\|   \covd^{(m)} \phi \big\|_{\nl^\alpha_N}
&\lesssim \lossfactor_T(F)  \| \phi \|_{L_{t,x}^\infty}, \label{local:eq-covariant-Schauder-1}\\  
\max_{\substack{m=0,1\colon \\ m<2+\beta}} N^{-\alpha-m}\big\|  \covd^{(m)} \phi \big\|_{\nl^\alpha_N} 
&\lesssim \lossfactor_T(F) N^{-\beta+\frac{2}{r}} \| \phi \|_{\Cs_t^{\frac{\beta}{2}} L_x^r},\label{local:eq-covariant-Schauder-2}
\\ 
N^{-\alpha} \big\|  \Qc \bigl(   \phi, \covd \psi \bigr)  \big\|_{\nl^\alpha_N}
&\lesssim \lossfactor_T(F)  \| \phi \|_{L_{t,x}^\infty} 
\big\| \covd \psi \big\|_{L_{t,x}^\infty}, \label{local:eq-covariant-Schauder-3}
\\ 
N^{-\alpha} \big\|  \Qc \bigl(   \phi, \covd \psi \bigr)  \big\|_{\nl^\alpha_N}
&\lesssim \lossfactor_T(F)  \bigl( N \| \phi \|_{L_{t,x}^\infty}
+ \big\| \covd \phi \big\|_{L_{t,x}^\infty} \bigr)  \| \psi \|_{L_{t,x}^\infty}. \label{local:eq-covariant-Schauder-4}
\end{align}
\end{lemma}

\begin{remark}
We note that \eqref{local:eq-covariant-Schauder-1} is the special case of \eqref{local:eq-covariant-Schauder-2} in which $(\beta,r)=(0,\infty)$. However, since \eqref{local:eq-covariant-Schauder-1} will be used repeatedly, we state it as a separate estimate.  
\end{remark}

\begin{remark}\label{local:rem-Schauder-Duh-without-F}
We observe that the Schauder estimates \eqref{local:eq-covariant-Schauder-1}--\eqref{local:eq-covariant-Schauder-4} also hold if the operator $\Duh^{F}_{\Y,B,N}$ in Definition \ref{local:def-nl-norm} is replaced by $\Duh_{\Y,B,N}$. The reason is that all estimates of $\fq_{\Y,B,N}$ that are used in the proof of \eqref{local:eq-covariant-Schauder-1}--\eqref{local:eq-covariant-Schauder-4} also hold for $\fp_{\Y,B,N}$.  The same observation also applies to the weighted Schauder estimates in Proposition~\ref{local:prop-weighted-covariant-Schauder} below.
\end{remark}

\begin{remark} 
Let $a,b,c\in[3]$ and let $\Qc\colon \Y^{(a)} \times \Y^{(b)} \rightarrow \Y^{(c)}$ be a $\liegroup$-equivariant bilinear map. Then, we can define a $\liegroup$-equivariant bilinear map $\widetilde{\Qc}\colon \Y \times \Y \rightarrow \Y$ by 
$\widetilde{\Qc}(y_1,y_2) = \iota_c \Qc(\pi_a y_1,\pi_b y_2)$
for all $y_1,y_2 \in \Y$. Here, $\pi_a\colon \Y\rightarrow \Y^{(a)}$ and $\pi_b\colon \Y \rightarrow \Y^{(b)}$ are the canonical projections, and $\iota_c\colon \Y^{(c)} \rightarrow \Y$ is the canonical embedding. Due to this, all estimates in Lemma~\ref{local:lem-covariant-Schauder} also hold for $\frkg$- or $\V$-valued maps, which will be used repeatedly in Section~\ref{section:local-paraproduct-estimates}.
\end{remark}

Before we turn to the proof of Lemma~\ref{local:lem-covariant-Schauder}, we need to make additional preparations. First, we prove the following pointwise estimate of the curvature tensor $F$ in time. 

\begin{lemma}\label{lemma:curvature-pointwise-in-time-bound}
For any $t \in (0, T]$ and $k \geq 0$, we have that 
\begin{equ}
\big\|\covd_x^{(k)} F(t)\big\|_{L_x^\infty} \lesssim \Fnorm{F}_{T, 0} \bigl( 1+ \Fnorm{F}_{T, (k-1) \vee 0}  \bigr)^{k} t^{-(1+\frac{k}{2})}.   
\end{equ}
\end{lemma}
\begin{proof}
By \eqref{eq:curvature-pointwise-in-time-bound-Fnorm}, we have that for any $t\in (0,T]$,
\begin{equ}
\big\|F(t)\big\|_{L_x^\infty} \leq \Fnorm{F}_{T,0} t^{-1}.
\end{equ}
For $k \geq 1$, by the heat kernel estimate \eqref{eq:higher-order-pB-ptwise-bd} and the fact that $F_B^{12}$ satisfies $(\ptl_t - \Delta_B) F_B^{12} = 0$ (recall Lemma~\ref{lemma:curvature-equation}), we have that 
\begin{equ}
t^{k/2} \bigl\| \covd_x^{(k)} F(t) \bigr\|_{L_x^\infty} \lesssim 
\bigl( 1 + \Fnorm{F}_{T,k-1}\bigr)^k \| F(t/2) \|_{L_x^\infty} 
\lesssim \Fnorm{F}_{T, 0}\bigl( 1 + \Fnorm{F}_{T,k-1}\bigr)^{k} t^{-1}, 
\end{equ}
which implies the desired estimate. 
\end{proof}

Recall from Notation~\ref{local:notation-omit-B} that we omit the subscript $B$ from our covariant derivatives. Additionally, recall the covariant derivatives $\covd_y^{\dagger, (m)}$ from Notation~\ref{notation:dagger}, as well as the fact that they commute with the covariant $x$-derivatives $\covd_x^{(k)}$. We also introduce $\Delta_{B,y}^\dagger = \covd_y^{\dagger,j} \covd_{y,j}^\dagger$. A computation shows that 
\begin{equ}\label{eq:dagger-laplacian-adjoint-relation}
(\Delta_{B, y} (\fq_{\Y, B}^*))^* = \Delta_{B, y}^\dagger \fq_{\Y, B}.
\end{equ}

\begin{lemma}\label{lemma:qB-mixed-all-derivatives}
Let $N\geq N_{\lo}$ and let $0\leq s < t \leq T$ satisfy $64^{-1} N^{-2} \leq t-s \leq 64 N^{-2}$. Then, we have for all $k,m \geq 0$ that 
\begin{equation}\label{eq:qB-mixed-no-s} 
\begin{aligned}
\bigl|& \covd_y^{\dagger,(m)} \covd_{t,x}^{(k)} \fq_{\Y,B}(s,y;t,x) \bigr| \\ 
&\lesssim 
\gronwall_T(C_0 F) \bigl( 1+ \Fnorm{F}_{T, k \vee m}  \bigr)^{k+m+2} 
\min(s,N^{-2})^{-\frac{m}{2}} N^{k+2} \exp \bigl( - c N^2 |x-y|^2 \bigr)
\end{aligned}
\end{equation}
and 
\begin{equation}\label{eq:qB-mixed-with-s} 
\begin{aligned}
\bigl|& \partial_s \covd_y^{\dagger,(m)} \covd_{t,x}^{(k)} \fq_{\Y,B}(s,y;t,x) \bigr| \\  
&\lesssim 
\gronwall_T(C_0 F) \bigl( 1+ \Fnorm{F}_{T, k \vee (m+2)}  \bigr)^{k+m+4} 
\min(s,N^{-2})^{-\frac{m+2}{2}} N^{k+2} \exp \bigl( - c N^2 |x-y|^2 \bigr). 
\end{aligned}
\end{equation}
The same estimates also hold with $\fq_{\Y,B}$ replaced by $\fq_{\Y,B,N}$ or $\fq_{\Y, B, \varep,N}$, where $\eps\in \sdyadic$. 
\end{lemma}

\begin{proof} We first prove the estimate \eqref{eq:qB-mixed-no-s}. Using Lemma~\ref{prelim:lem-space-time-derivatives-covariant-heat}, we have that 
\begin{equation*}
\covd_{t,x}^{(k)} \fq_{\Y,B}(s,y;t,x) = \sum_{\ell\geq0} \sum_{\substack{j_1 + \hdots + j_{\ell+1}=k - 2\ell, \\ j_1,\hdots, j_{\ell+1} \geq 0}} \schO \Big( \bigl( \covd_x^{(j_1)} F(t) \bigr) 
\hdots 
\bigl( \covd_x^{(j_\ell)} F(t) \bigr) 
\bigl( \covd_x^{(j_{\ell+1})} \fq_{\Y,B}(s,y;t,x) \bigr) \Big). 
\end{equation*}
By applying $\covd_{y}^{\dagger,(m)}$ to both sides of the identity, it follows that 
\begin{align*}
\covd_{y}^{\dagger,(m)}& \covd_{t,x}^{(k)} \fq_{\Y,B}(s,y;t,x) \\ 
&= \sum_{\ell\geq0} \sum_{\substack{j_1 + \hdots + j_{\ell+1}=k - 2\ell,  \\ j_1,\hdots, j_{\ell+1} \geq 0}} \schO \Big( \bigl( \covd_x^{(j_1)} F(t) \bigr) 
\hdots 
\bigl( \covd_x^{(j_\ell)} F(t) \bigr) 
\bigl( \covd_{y}^{\dagger,(m)} \covd_x^{(j_{\ell+1})} \fq_{\Y,B}(s,y;t,x) \bigr) \Big). 
\end{align*}
Using Proposition~\ref{prop:qB-mixed-derivatives-ptwise-bounds} and Lemma~\ref{lemma:curvature-pointwise-in-time-bound}, it then follows that 
\begin{align*}
\bigl|&\covd_{y}^{\dagger,(m)} \covd_{t,x}^{(k)} \fq_{\Y,B}(s,y;t,x) \bigr|  \\ 
&\lesssim \sum_{\ell\geq0} \sum_{\substack{j_1 + \hdots + j_{\ell+1}=k - 2\ell,  \\ j_1,\hdots, j_{\ell+1} \geq 0}} 
\bigl| \covd_x^{(j_1)} F(t) \bigr| \hdots \bigl|\covd_x^{(j_\ell)} F(t)\bigr|  
\bigl| \covd_{y}^{\dagger,(m)} \covd_x^{(j_{\ell+1})} \fq_{\Y,B}(s,y;t,x)  \bigr| \\
&\lesssim \sum_{\ell\geq0} \sum_{\substack{j_1 + \hdots + j_{\ell+1}= k - 2\ell,  \\ j_1,\hdots, j_{\ell+1} \geq 0}} \Bigl(  \bigl( 1+ \Fnorm{F}_{T, k}  \bigr)^{j_1+\hdots +j_\ell + \ell}
t^{-(\ell + \frac{j_1+\hdots + j_\ell}{2})} \\ 
 & \quad  \times \gronwall_T(C_0 F) \bigl( 1+ \Fnorm{F}_{T, k \vee m}  \bigr)^{m+j_{\ell+1} +2} \min(s,N^{-2})^{-\frac{m}{2}} N^{j_{\ell+1}+2} \exp\bigl( - c N^2 |x-y|^2 \bigr) \Bigr). 
\end{align*}
Since $t \geq t-s \geq 64^{-1} N^{-2}$, it holds that 
\begin{equation*}
t^{-(\ell + \frac{j_1+\hdots + j_\ell}{2})} N^{j_{\ell+1}+2} 
\lesssim N^{j_1+ \hdots + j_{\ell}+j_{\ell+1}+2\ell+2} = N^{k+2}. 
\end{equation*}
Furthermore, using that $j_1+\hdots +j_{\ell+1}+ \ell = k - \ell \leq k$, it also holds that 
\begin{equation*}
\bigl( 1+ \Fnorm{F}_{T, k}  \bigr)^{j_1+\hdots +j_\ell + \ell}
\bigl( 1+ \Fnorm{F}_{T, k \vee m}  \bigr)^{m+j_{\ell+1} +2} 
\lesssim \bigl( 1+ \Fnorm{F}_{T, k \vee m}  \bigr)^{k+m+2}.
\end{equation*}
By combining our previous three estimates, we obtain our desired estimate \eqref{eq:qB-mixed-no-s}. We now turn to the proof of \eqref{eq:qB-mixed-with-s}, i.e., the estimate for the $\partial_s$-derivative. By using Lemma~\ref{lemma:adjoints-backwards-heat-equations} and taking adjoints (and using \eqref{eq:dagger-laplacian-adjoint-relation}), it follows that 
\begin{equation*}
(\partial_s + \Delta_{B,y}^\dagger) \fq_{\Y,B}(s,y;t,x)  
= \schO \bigl(  \fq_{\Y,B}(s,y;t,x) F(s) \bigr). 
\end{equation*} 
Since $\covd_{t,x}^{(k)}$ commutes with both $(\partial_s + \Delta_{B,y}^\dagger)$ and the right multiplication by $F(s)$, applying it to the above identity yields
\begin{equation*}
(\partial_s + \Delta_{B,y}^\dagger) \covd_{t,x}^{(k)} \fq_{\Y,B}(s,y;t,x)  
= \schO \bigl( \bigl( \covd_{t,x}^{(k)}\fq_{\Y,B}(s,y;t,x) \bigr) F(s) \bigr). 
\end{equation*}
By using Remark~\ref{remark:higher-derivatives-equation-with-ad-F-term} or, more precisely, its analogue for $(\partial_s + \Delta_{B,y}^\dagger)$, we obtain that 
\begin{equation*}
(\partial_s + \Delta_{B,y}^\dagger) \covd_y^{\dagger,(m)}\covd_{t,x}^{(k)} \fq_{\Y,B}(s,y;t,x)  
= \sum_{\substack{j_1+j_2=m\\ j_1,j_2 \geq 0}}\schO \Bigl( \bigl( \covd_y^{\dagger,(j_1)}\covd_{t,x}^{(k)}\fq_{\Y,B}(s,y;t,x) \bigr) \bigl( \covd_y^{(j_2)} F(s) \bigr) \Bigr). 
\end{equation*}
The desired estimate \eqref{eq:qB-mixed-with-s} can now be obtained by first bringing the $\Delta_{B,y}^{\dagger}$-term over to the right-hand side and then using both Lemma~\ref{lemma:curvature-pointwise-in-time-bound} and the estimate \eqref{eq:qB-mixed-no-s}. Finally, in order to see that the estimates remain valid after replacing $\fq_{\Y,B}$ with $\fq_{\Y,B,N}$ or $\fq_{\Y,\Beps,N}$, we rely on the product rule, \eqref{local:eq-qB-eps-N}, and \eqref{local:eq-qB-N}. 
\end{proof}

Equipped with Lemma~\ref{lemma:qB-mixed-all-derivatives}, we can now prove our unweighted covariant Schauder estimates. 

\begin{proof}[Proof of Lemma~\ref{local:lem-covariant-Schauder}] 
We first show that \eqref{local:eq-covariant-Schauder-3} and \eqref{local:eq-covariant-Schauder-4} can be derived from \eqref{local:eq-covariant-Schauder-1} and the product formula (Lemma~\ref{prelim:lem-equivariant-product-rule}). Indeed, \eqref{local:eq-covariant-Schauder-3} follows from 
\begin{equation*}
N^{-\alpha} \big\|  \Qc \bigl(   \phi, \covd \psi \bigr)  \big\|_{\nl^\alpha_N}
\lesssim \lossfactor_T(F) \big\|  \Qc \bigl(   \phi, \covd \psi \bigr) \big\|_{L_{t,x}^\infty} \lesssim 
\lossfactor_T(F)  \| \phi \|_{L_{t,x}^\infty} 
\big\| \covd \psi \big\|_{L_{t,x}^\infty}.
\end{equation*}
Similarly, \eqref{local:eq-covariant-Schauder-4} follows from
\begin{align*}
N^{-\alpha} \big\|  \Qc \bigl(   \phi, \covd \psi \bigr)  \big\|_{\nl^\alpha_N}  
&\leq N^{-\alpha} \big\|  \covd  \Qc \bigl(   \phi, \psi \bigr)  \big\|_{\nl^\alpha_N}
+  N^{-\alpha} \big\|  \Qc \bigl(  \covd \phi,  \psi \bigr)  \big\|_{\nl^\alpha_N} \\ 
&\lesssim \lossfactor_T(F) \bigl( N \bigl\| \Qc \bigl( \phi, \psi \bigr) \bigr\|_{L_{t,x}^\infty} + \bigl\| \Qc\bigl( \covd \phi, \psi \bigr) \bigr\|_{L_{t,x}^\infty} \bigr) \\ 
&\lesssim \lossfactor_T(F)  \bigl( N \| \phi \|_{L_{t,x}^\infty}
+ \big\| \covd \phi \big\|_{L_{t,x}^\infty} \bigr)  \| \psi \|_{L_{t,x}^\infty}.
\end{align*}
We now turn to the proof of \eqref{local:eq-covariant-Schauder-1}. In fact, we prove the more general estimate 
\begin{equation}\label{local:eq-variant-CS-1} 
 \big\|   \covd^{(m)} \phi \big\|_{\nl^\alpha_N}
\lesssim \lossfactor_T(F)  N^{\alpha+m+\frac{2}{r}}   \| \phi \|_{L_{t}^\infty L_x^r}, 
\end{equation}
which also covers the case $\beta=0$ in \eqref{local:eq-covariant-Schauder-2}. To this end, let $0\leq k \leq 10$ and let $0 \leq m \leq 1$. We then first write 
\begin{equation}\label{eq:schauder-estimate-duhamel-expression-intermediate}
\begin{aligned}
\Big(\covd^{(k)}_{t,x} \Duh_{\Y,B, N}^F(\covd^{(m)} \phi)\Big)(t, x) &= \iint \ds \dy \,  \covd^{(k)}_{t,x} \fq_{\Y,B,N}(s, y; t, x) \covd_y^{(m)} \phi(s, y) \\
&= \iint \ds \dy \,  \Bigl( \covd_y^{\dagger,(m)} \covd^{(k)}_{t,x} \fq_{\Y,B,N}(s, y; t, x)\Bigr)  \phi(s, y).
\end{aligned}
\end{equation}
Using Lemma~\ref{lemma:qB-mixed-all-derivatives}, we then estimate
\begin{align*}
\bigl| \eqref{eq:schauder-estimate-duhamel-expression-intermediate} \bigr| 
&\lesssim \lossfactor_T(F) N^{k+2} \int \ds \dy \, \ind\{ t-s \sim N^{-2} \} \min(s,N^{-2})^{-\frac{m}{2}} \exp\bigl( -c N^2 |x - y|^2 \bigr) |\phi(s,y)|  \\
&\lesssim \lossfactor_T(F) N^{k+2} \Bigl\| \ind\{ t-s \sim N^{-2} \} \min(s,N^{-2})^{-\frac{m}{2}} \exp\bigl( -c N^2 |x - y|^2 \bigr) \Bigr\|_{L_s^1 L_y^{r^\prime}} 
\| \phi \|_{L_s^\infty L_y^r} \\ 
&\lesssim \lossfactor_T(F) N^{k+m+\frac{2}{r}-2} \| \phi \|_{L_s^\infty L_y^r}. 
\end{align*}
In the last line, we used that $0\leq m \leq 1$, since otherwise $\min(s,N^{-2})^{-\frac{m}{2}}$ would not be integrable. This completes the proof of \eqref{local:eq-variant-CS-1}. 

It remains to prove \eqref{local:eq-covariant-Schauder-2}. By \eqref{local:eq-variant-CS-1}, we may assume that $\beta<0$. By Definition~\ref{def:temporal-regularity-space}, and since Hölder functions can easily be extended from compact intervals to the real line, there exists a $\widetilde{\phi} \colon \R \times \T^2 \rightarrow \Y$ such that $\partial_t \widetilde{\phi}(t)=\phi(t)$ for all $t\in (0,T)$ and 
\begin{equation*}
\| \widetilde{\phi} \|_{\Cs_t^{1+\beta/2}L_x^r(\R\times \T^2)} \lesssim \| \phi \|_{\Cs_t^{\beta/2}L_x^r([0,T]\times \T^2)}. 
\end{equation*}
 For a given $t\in (0,T]$, we choose $s_0\coloneq0$ if $t\leq 64N^{-2}$, and $s_0\coloneq t$ if $t>64N^{-2}$. From this definition, it then follows for all $s\in [0,T]$ in the support of $\chi_N(t-s)$ that $|s-s_0|\lesssim \min(s,t-s)$. We now let $(P_{L})_{L\geq 1}$ be Littlewood--Paley projections in the time variable to frequencies $\sim L^2$. We choose to denote the frequency scales by $L^2$, rather than $L$, in order to respect parabolic scaling. It then holds that 
\begin{equation}\label{local:eq-covariant-Schauder-Phi-bound}
\begin{aligned}
\big\| P_L \widetilde{\phi}(s) - P_L \widetilde{\phi}(s_0) \big\|_{L_x^r} 
&\lesssim \min(L^2 |s-s_0|,1) L^{-2-\beta} \| \widetilde{\phi} \|_{\Cs_t^{1+\beta/2}L_x^r} \\
&\lesssim \min(L^2 s, L^2 (t-s),1) L^{-2-\beta} \| \widetilde{\phi} \|_{\Cs_t^{1+\beta/2}L_x^r}. 
\end{aligned}
\end{equation}
Starting from \eqref{eq:schauder-estimate-duhamel-expression-intermediate}, we then write
\begin{equation}\label{local:eq-covariant-Schauder-intermediate-decomposed}
\begin{aligned}
\Big(&\covd^{(k)}_{t,x} \Duh_{\Y,B, N}^F(\covd^{(m)} \phi)\Big)(t, x)  \\
&= \iint \ds \dy \, \Bigl( \covd_{y}^{\dagger,(m)} \covd^{(k)}_{t,x} \fq_{\Y,B,N}(s, y; t, x) \Bigr)  \partial_s \bigl( \widetilde{\phi}(s,y) - \widetilde{\phi}(s_0,y) \bigr) \\ 
&= - \sum_{L\geq 1}\iint \ds \dy \,  \Bigl( \partial_s \covd_{y}^{\dagger,(m)} \covd^{(k)}_{t,x} \fq_{\Y,B,N}(s, y; t, x) \Bigr) \bigl( P_L \widetilde{\phi}(s,y) - P_L \widetilde{\phi}(s_0,y) \bigr).
\end{aligned}
\end{equation}
We emphasise that there is no boundary term in \eqref{local:eq-covariant-Schauder-intermediate-decomposed} since, if $s=0$ is in the support of $\fq_{\Y,B,N}(s,y;t,x)$, then $s_0=0$. Using \eqref{local:eq-covariant-Schauder-Phi-bound} and Hölder's inequality, each summand in \eqref{local:eq-covariant-Schauder-intermediate-decomposed} can be estimated by 
\begin{equation}
\begin{aligned}
\Big|& \int \ds \dy \,  \Bigl( \partial_s \covd_{y}^{\dagger,(m)} \covd^{(k)}_{t,x} \fq_{\Y,B,N}(s, y; t, x) \Bigr) \bigl( P_L \widetilde{\phi}(s,y) - P_L \widetilde{\phi}(s_0,y) \bigr)\Big| \\
&\lesssim L^{-2-\beta} \Big\| \min\bigl(L^2 s, L^2 (t-s),1 \bigr) 
\Bigl( \partial_s \covd_{y}^{\dagger,(m)} \covd^{(k)}_{t,x} \fq_{\Y,B,N}(s, y; t, x) \Bigr) \Big\|_{L_s^1 L_y^{r^\prime}} \| \widetilde{\phi} \|_{\Cs_t^{1+\beta/2}L_x^r}.
\end{aligned}
\end{equation}
Using Lemma~\ref{lemma:qB-mixed-all-derivatives}, we can estimate
\begin{align*}
\Big\|& \min\bigl(L^2 s, L^2 (t-s),1 \bigr) 
\Bigl( \partial_s \covd_{y}^{\dagger,(m)} \covd^{(k)}_{t,x} \fq_{\Y,B,N}(s, y; t, x) \Bigr) \Big\|_{L_s^1 L_y^{r^\prime}} \\
&\lesssim \lossfactor_T(F) N^{k+2} 
\Big\| \ind\{ t-s\sim N^{-2} \} \min\bigl(L^2 s, L^2 (t-s),1 \bigr) \min\big( s, N^{-2} \big)^{-\frac{m+2}{2}} \exp\bigl( - c N^2 |x-y|^2\bigr) \Big\|_{L_s^1 L_y^{r^\prime}} \\
&\lesssim \lossfactor_T(F) N^{k+\frac{2}{r}} \Big\| \ind\{ t-s\sim N^{-2} \} \min\bigl(L^2 s, L^2 (t-s),1 \bigr) \min\big( s, N^{-2} \big)^{-\frac{m+2}{2}} \Big\|_{L_s^1}.
\end{align*}
In the $L_s^1$-norm above, the last factor is equivalent to the sum of $s^{-\frac{m+2}{2}}$ and $N^{m+2}$. The contribution of the $s^{-\frac{m+2}{2}}$-term can be bounded by 
\begin{align*}
\big\|& \ind\{s\geq L^{-2}\} \ind\{ t-s \sim N^{-2} \} s^{-\frac{m+2}{2}} \big\|_{L_s^1} 
+ L^2 \big\| \ind\{s\leq L^{-2}\} \ind\{ t-s \sim N^{-2} \}  s^{1-\frac{m+2}{2}} \big\|_{L_s^1} \\
&\lesssim \min\bigl( L^{m+2} N^{-2}, L^m \log(2+L/N) \bigr) + \min\bigl( L^2 N^{m-2}, L^m \bigr) \\
&\lesssim \min\bigl( L^2 N^{m-2}, L^m \bigr) \log(2+L/N), 
\end{align*}
where we used that $m=0,1$. Similarly, the contribution of the $N^{m+2}$-term can be bounded by 
\begin{align*}
&N^{m+2} \big\| \ind\{ t-s \sim N^{-2} \} \min\bigl(L^2 s, L^2 (t-s),1 \bigr)  \big\|_{L_s^1} \\
&\lesssim N^{m+2} \min\Bigl( \big\| \ind\{ t-s \sim N^{-2} \}\big\|_{L_s^1}, 
L^2 \big\| \ind\{ t-s \sim N^{-2} \} (t-s)\big\|_{L_s^1} \Bigr)  
\lesssim \min(N^m, N^{m-2} L^2). 
\end{align*}
Since $\beta<0$ and $m<2+\beta$, it holds that
\begin{equation*}
\sum_{L\geq 1} L^{-2-\beta} \Bigl( \min\bigl( L^2 N^{m-2}, L^m \bigr)  \log(2+L/N) + \min(N^m, N^{m-2} L^2) \Bigr) \lesssim N^{m-\beta-2}, 
\end{equation*}
which yields the desired estimate. 
\end{proof}

So far, we have stated and proved the unweighted Schauder estimates in our solution spaces (Lemma~\ref{local:lem-covariant-Schauder}). We now turn to the weighted Schauder estimates, which are the subject of the following proposition. 

\begin{proposition}[Weighted covariant Schauder estimates]\label{local:prop-weighted-covariant-Schauder}
Let $N \geq N_{\lo}$, let $\alpha\in  \R$,  let $\beta \in (-2,0]$, let $-2 < \eta,\eta_1,\eta_2,\widetilde{\eta}_1,\widetilde{\eta}_2,\nu \leq 0$, and let $1 \leq r \leq \infty$. Furthermore, assume that  
\begin{equation*}
\eta_1+\eta_2 >-2 \quad \text{and} \quad \widetilde{\eta}_1+\widetilde{\eta}_2>-1. 
\end{equation*}
To simplify the notation, let $\eta_{12}\coloneq \eta_1 +\eta_2$ and $\widetilde{\eta}_{12}\coloneq \widetilde{\eta}_1+\widetilde{\eta}_2$. 
Let $\Qc\colon \Y \times \Y\rightarrow \Y$ be a $\liegroup$-equivariant bilinear map, let $\phi \colon [0,T]\times \T^2 \rightarrow \Y$, and let $\psi\colon [0,T]\times \T^2 \rightarrow \Y$. Then, we have the estimates  
\begin{align}
\max_{\substack{m=0,1\colon \\ m<2+\eta}} N^{-\alpha-m} \big\|  \covd^{(m)} \phi   \big\|_{\nl^{\alpha,\nu}_N}
&\lesssim \lossfactor_T(F)  \max \bigl( N^{\nu-\eta},N_{\lo}^{\nu-\eta}\bigr)   \| t^{-\eta/2} \phi \|_{L_{t,x}^\infty}, \label{local:eq-weighted-covariant-Schauder-1}\allowdisplaybreaks[3]\\  
\max_{\substack{m=0,1\colon \\ m<2+\beta+\eta}} N^{-\alpha-m} \big\|  \covd^{(m)} \phi   \big\|_{\nl^{\alpha,\nu}_N}
&\lesssim \lossfactor_T(F) \max \bigl( N^{\nu-\eta},N_{\lo}^{\nu-\eta}\bigr) 
 N^{-\beta+\frac{2}{r}} \| t^{-\eta/2}\phi \|_{\Cs_t^{\frac{\beta}{2}} L_x^r},
 \label{local:eq-weighted-covariant-Schauder-2}
\end{align}
as well as the estimates
\begin{align}
N^{-\alpha}  \big\|   \Qc \bigl(  \phi, \covd \psi \bigr) \big\|_{\nl^{\alpha,\nu}_N}
\lesssim \lossfactor_T(F)&\max \bigl( N^{\nu-\eta_{12}},N_{\lo}^{\nu-\eta_{12}}\bigr)  \| t^{-\eta_1/2} \phi \|_{L_{t,x}^\infty} 
\big\| t^{-\eta_2/2}\covd \psi \big\|_{L_{t,x}^\infty}, \allowdisplaybreaks[3] \label{local:eq-weighted-covariant-Schauder-3}\\ 
N^{-\alpha} \big\|   \Qc \bigl(  \phi, \covd \psi \bigr) \big\|_{\nl^{\alpha,\nu}_N}
\lesssim \lossfactor_T(F)& \max \bigl( N^{\nu-\eta_{12}},N_{\lo}^{\nu-\eta_{12}}\bigr)   \big\| t^{-\eta_1/2} \covd \phi \big\|_{L_{t,x}^\infty}\|t^{-\eta_2/2} \psi \|_{L_{t,x}^\infty} \notag \\ 
+ \lossfactor_T(F)& \max\bigl( N^{\nu-\widetilde{\eta}_{12}},N_{\lo}^{\nu-\widetilde{\eta}_{12}}\bigr)
N \big\| t^{-\frac{\widetilde{\eta}_1}{2}} \phi \big\|_{L_{t,x}^\infty}\| t^{-\frac{\widetilde{\eta}_2}{2}}\psi \|_{L_{t,x}^\infty}. 
\label{local:eq-weighted-covariant-Schauder-4}
\end{align}
\end{proposition}

\begin{remark} In applications of Proposition~\ref{local:prop-weighted-covariant-Schauder}, we will mostly take $\nu\geq \eta,\eta_{12},\widetilde{\eta}_{12}$, in which case the $N_{\lo}$-terms in the maxima above can always be omitted. 
\end{remark}

\begin{proof}
Similarly to the proof of Lemma~\ref{local:lem-covariant-Schauder}, \eqref{local:eq-weighted-covariant-Schauder-3} and \eqref{local:eq-weighted-covariant-Schauder-4} can be derived from \eqref{local:eq-weighted-covariant-Schauder-1} and the product formula. Furthermore, \eqref{local:eq-weighted-covariant-Schauder-1} is the special case of 
\eqref{local:eq-weighted-covariant-Schauder-2} in which $(\beta,r)=(0,\infty)$. Thus, it suffices to prove \eqref{local:eq-weighted-covariant-Schauder-2}. 

The idea behind the proof of \eqref{local:eq-weighted-covariant-Schauder-2} is to combine our earlier, unweighted Schauder estimate \eqref{local:eq-covariant-Schauder-2} with a dyadic decomposition in time. For any $0\leq k \leq 10$, we write 
\begin{equation}\label{local:eq-weighted-Schauder-p1}
\begin{aligned}
\covd_{t,x}^{(k)} \Duh_{\Y,B,N}^{F} \big( \covd_{x}^{(m)} \phi \big)
= \sum_{K,M \geq N_{\lo}} \chi_K \covd_{t,x}^{(k)} \Duh_{\Y,B,N}^{F} \big( \chi_M \covd_{x}^{(m)} \phi \big). 
\end{aligned}
\end{equation}
Since the supports of $(\chi_K)_{K\geq N_{\lo}}$ have bounded overlaps, it then holds that 
\begin{equation}\label{local:eq-weighted-Schauder-p2}
\begin{aligned}
\Big\| t^{-\frac{\nu}{2}} \covd_{t,x}^{(k)} \Duh_{\Y,B,N}^{F} \big( \covd_{x}^{(m)} \phi \big) \Big\|_{L_{t,x}^\infty} 
\lesssim \sup_{K \geq N_{\lo}} \sum_{M\geq N_{\lo}} K^{\nu}  
\Big\| \chi_K \covd_{t,x}^{(k)} \Duh_{\Y,B,N}^{F} \big( \chi_M \covd_{x}^{(m)} \phi \big) \Big\|_{L_{t,x}^\infty}. 
\end{aligned}
\end{equation}
We note that, since $\covd_{x}^{(m)}$ only contains covariant spatial derivatives, $\chi_M $ and $\covd_{x}^{(m)}$ commute. The time-localisations in \eqref{local:eq-weighted-Schauder-p1} are such that $t\sim K^{-2}$, $t-s \sim N^{-2}$, and $s\sim M^{-2}$, where $s$ denotes the integration variable. Since $t=t-s+s$, we then must have that $K\sim M \lesssim N$ or $K\sim N \lesssim M$. We now let $\beta^\prime \in (-2,0]$ satisfy $\beta^\prime \leq \beta$ and $\beta^\prime>-2+m$, but it otherwise remains to be chosen. Note that, due to our assumptions $m<2+\beta+\eta$ and $\eta\leq 0$, such a choice can be made. Using Lemma~\ref{prelim:lem-time-localization} (for the cutoff function $\chi_M(t) (M^2t)^{\frac{\eta}{2}}$) and our unweighted Schauder estimate \eqref{local:eq-covariant-Schauder-2}, we obtain that 
\begin{equation}\label{local:eq-weighted-Schauder-p3}
\begin{aligned}
K^{\nu} \Big\|& \chi_K \covd_{t,x}^{(k)} \Duh_{\Y,B,N}^{F} \big( \chi_M \covd_{x}^{(m)} \phi \big) \Big\|_{L_{t,x}^\infty} \\ 
&\lesssim \lossfactor_T(F) K^{\nu} N^{-2+k+m-\beta^\prime+\frac{2}{r}} \big\| \chi_M \phi \big\|_{\Cs_t^{\frac{\beta^\prime}{2}} L_x^r} \\ 
&\lesssim \lossfactor_T(F) K^{\nu} N^{-2+k+m-\beta^\prime+\frac{2}{r}} M^{\beta^\prime-\beta-\eta} \big\| t^{-\frac{\eta}{2}} \phi \big\|_{\Cs_t^{\frac{\beta}{2}} L_x^r}. 
\end{aligned}
\end{equation}
In the case $M \sim K \lesssim N$, we choose $\beta^\prime=\beta$, which leads to the prefactor $K^{\nu-\eta} N^{-2+k+m-\beta+\frac{2}{r}}$ in 
\eqref{local:eq-weighted-Schauder-p3}. Since the sum over $M$ in \eqref{local:eq-weighted-Schauder-p2} then only contains a bounded number of terms, and since $N_{\lo} \leq K \lesssim N$, this leads to an acceptable contribution. In the case $K\sim N \lesssim M$, we choose $\beta^\prime=-2+m+\eps$, where $\eps>0$ is sufficiently small. Then, the prefactor in \eqref{local:eq-weighted-Schauder-p3} becomes
\begin{equation*}
K^{\nu} N^{k+\frac{2}{r}-\eps} M^{-2+m-\beta-\eta+\eps}. 
\end{equation*}
Due to our assumption $m<2+\beta+\eta$, the exponent of $M$ is negative. Thus, the contribution of the case $K\sim N \lesssim M$ to \eqref{local:eq-weighted-Schauder-p2} can be bounded by
\begin{align*}
&\ind \{ K \sim N \} K^{\nu} N^{k+\frac{2}{r}-\eps} \sum_{M\gtrsim N} M^{-2+m-\beta-\eta+\eps} \\  
&\quad \lesssim \ind \{ K \sim N \} K^{\nu} N^{-2+k+m-\beta-\eta+\frac{2}{r}} \lesssim N^{-2+k+m-\beta+\nu-\eta+\frac{2}{r}}, 
\end{align*}
which is acceptable. 
\end{proof}

In addition to the Schauder estimates above, we also control covariant linear evolutions in our solution spaces. This will later be needed in the proof of Proposition~\ref{local:prop-higgsdet-sdet}. 

\begin{lemma}\label{local:lem-linear-evolution-sol}
Let $\higgsdet_0\colon \T^2 \rightarrow \V$. Furthermore, let $\alpha \in \R$ and $\eta\leq 0$ satisfy $\alpha+\eta\leq -\frac{2}{r}$, where $r$ is as in \eqref{e:small_constants}. Then, it holds for all $N\geq N_{\lo}$ that
\begin{equation}\label{local:eq-linear-evolution-sol-1}
\Big\| \int \dy \, \fp_{B,N}(0,y;t,x) \higgsdet_0(y) \Big\|_{\sol^{\alpha,\eta}_N} \lesssim \lossfactor_T(F)  T^{-\frac{\alpha+\eta}{2}-\frac{1}{r}} \big\| \higgsdet_0 \big\|_{L_x^r}. 
\end{equation}
Alternatively, if $\alpha+\eta\leq -\kappa$ is satisfied, then it holds that
\begin{equation}\label{local:eq-linear-evolution-sol-2}
\Big\| \int \dy \, \fp_{B,N}(0,y;t,x) \higgsdet_0(y) \Big\|_{\sol^{\alpha,\eta}_N} \lesssim \lossfactor_T(F)  T^{-\frac{\alpha+\eta}{2}-\frac{\kappa}{2}} \bigl\| \higgsdet_0\bigr\|_{\covHolder{B(0), T}^{-\kappa}}. 
\end{equation}
\end{lemma}

\begin{proof}
We first note that, on the support of the integrals appearing in \eqref{local:eq-linear-evolution-sol-1} and \eqref{local:eq-linear-evolution-sol-2}, it necessarily holds that $\frac{1}{16}N^{-2} \leq t \leq 16 N^{-2}$. As a result, using Definition~\ref{local:def-solution-space}, we obtain
\begin{equation}\label{local:eq-linear-evolution-sol-p1}
\begin{aligned}
\Big\|& \int \dy \, \fp_{B,N}(0,y;t,x) \higgsdet_0(y) \Big\|_{\sol^{\alpha,\eta}_N}  \\
&\sim \max_{0\leq m \leq 10} N^{\alpha+\eta-m} \Big\| \int \dy \, \bigl( \covd_{t,x}^{(m)} \fp_{B,N} \bigr)(0,y;t,x) \higgsdet_0(y) \Big\|_{L_{t,x}^\infty}. 
\end{aligned}
\end{equation}
We now prove \eqref{local:eq-linear-evolution-sol-1} and \eqref{local:eq-linear-evolution-sol-2} separately. To obtain \eqref{local:eq-linear-evolution-sol-1}, we use our pointwise estimate (Lemma~\ref{lemma:qB-mixed-all-derivatives}) and Young's convolution inequality, which imply that
\begin{align*}
&N^{\alpha+\eta-m} \Big\| \int \dy \, \bigl( \covd_{t,x}^{(m)}  \fp_{B,N} \bigr)(0,y;t,x) \higgsdet_0(y) \Big\|_{L_{t,x}^\infty} \\ 
&\lesssim \lossfactor_T(F) N^{\alpha+\eta+2} 
\Big\| \int \exp\Bigl( - c N^2 |x-y|^2 \Bigr) |\higgsdet_0(y)| \dy \Big\|_{L_{t,x}^\infty}
\lesssim \lossfactor_T(F) N^{\alpha+\eta+\frac{2}{r}} \| \higgsdet_0 \|_{L_x^r}. 
\end{align*} 
Due to our condition $\alpha+\eta\leq -\frac{2}{r}$, the exponent of $N$ is non-positive. Together with $N\geq N_{\lo} \sim T^{-\frac{1}{2}}$, this then implies the first estimate \eqref{local:eq-linear-evolution-sol-1}. It remains to prove the second estimate \eqref{local:eq-linear-evolution-sol-2}. To this end, we recall that $\frac{1}{16} N^{-2} \leq t \leq 16 N^{-2}$ and set $t^\prime \coloneq \frac{1}{32} N^{-2}$. Using the semigroup property of the covariant heat kernel, we then have 
\begin{align*}
\int \dy \, \covd_{t,x}^{(m)} \fp_{B,N}(0,y;t,x)\, \higgsdet_0(y)  
=&  \int \dy \, \covd_{t,x}^{(m)} \bigl( \chi_N(t) \fp_{B}(0,y;t,x) \bigr) \higgsdet_0(y)  \\ 
=& \int  \dx^\prime \,  \covd_{t,x}^{(m)} \bigl( \chi_N(t) \fp_{B}(t^\prime,x^\prime;t,x) \bigr) 
\Bigl( \int \dy \, \fp_B(0,y;t^\prime,x^\prime) \higgsdet_0(y)  \Bigr). 
\end{align*}
It therefore follows that 
\begin{align*}
\Bigl\|&  \int \dy \,  \covd_{t,x}^{(m)} \fp_{B,N}(0,y;t,x)\, \higgsdet_0(y)    \Bigr\|_{L_{t,x}^\infty} \\ 
&\lesssim
\bigl \|  \covd_{t,x}^{(m)}\bigl( \chi_N(t) \fp_{B}(t^\prime,x^\prime;t,x) \bigr) \bigr\|_{L_{t,x}^\infty L_{x^\prime}^1} 
\Bigl \|  \int \dy \,  \fp_B(0,y;t^\prime,x^\prime) \higgsdet_0(y)   \Bigr\|_{L_{x^\prime}^\infty}. 
\end{align*}
Using our pointwise estimates (Lemma~\ref{lemma:qB-mixed-all-derivatives}), the first factor can be bounded by $\lossfactor_T(F) N^m$. Using the definition of the $\|\place\|_{\covHolder{B(0), T}^{-\kappa}}$-norm from Definition \ref{defn:covariant_Holder}, the second factor can be bounded by $N^\kappa \bigl\| \higgsdet_0\bigr\|_{\covHolder{B(0), T}^{-\kappa}}$. After inserting this back into~\eqref{local:eq-linear-evolution-sol-p1} and using the condition $\alpha+\eta\leq -\kappa$, this readily implies \eqref{local:eq-linear-evolution-sol-2}.
\end{proof}

Next, we turn to our commutator estimates, which will be needed when uncovering the structure of the paracontrolled term in our local solution ansatz \dash see Section~\ref{section:local-paracontrolled}. 

\begin{lemma}[Commutator estimate]\label{kernel:lem-commutator-integral}
Let $\genV, \genV_1, \genV_2$ be vector spaces on which $\liegroup$ acts. Furthermore, let $\Qc \colon \genV_1 \times \genV_2 \ra \genV$ be a $\liegroup$-equivariant bilinear map,  let $G \colon [0, T] \times \T^2 \ra \genV_1$, and let $H \colon [0, T] \times \T^2 \ra \genV_2$. Let $M \in \dyadic$ and let $-2<\eta \leq 0$. We then have that
\begin{equs}
~&\bigl\|t^{-\frac{\eta}{2}} \bigl( \Duh_{B, M}(\Qc(G, H)) - \Qc(G, \Duh_{B, M}(H)) \bigr)\bigr\|_{L_{t, x}^\infty([0, T])}  \\
& \lesssim   \lossfactor_T(F) \bigl(M^{-4} \big\|t^{-\frac{\eta}{2}} (\ptl_t - \Delta_B) G \big\|_{L_{t,x}^\infty([0, T])} + M^{-3} \big\|t^{-\frac{\eta}{2}}\covd G \big\|_{L_{t,x}^\infty([0, T])}  \bigr) \big\|H \big\|_{L_{t, x}^\infty([0, T])} .
\end{equs}
\end{lemma}
\begin{proof}
Fix $(t_0, x_0) \in [0, T] \times \T^2$. Let $\phi \colon [0, t_0] \times \T^2 \ra \genV$, $\psi \colon [0, t_0] \times \T^2 \ra \genV_2$ be given by
\begin{equ}
\phi = \Duh_{B}\bigl(\chi_M(t_0 - \place) \Qc(G, H)\bigr), ~~ \psi = \Duh_{B}\bigl(\chi_M(t_0 - \place) H \bigr).
\end{equ}
Using this notation, we can write 
\begin{equ}
\Bigl(\Duh_{B, M}(\Qc(G, H)) - \Qc(G, \Duh_{B, M}(H))\Bigr)(t_0, x_0) = \bigl(\phi - \Qc(G, \psi)\bigr)(t_0, x_0). 
\end{equ}
Additionally, observe that $\phi(t)= \psi(t) = 0$ unless $|t - t_0 | \lesssim M^{-2}$. For this reason, we introduce the interval $I\coloneq[0 \vee (t_0-CM^{-2}),t_0]$, which contains the support of $\phi$ and $\psi$. Next, we record several preliminary estimates. For any function $f\colon I \times \T^2 \rightarrow \R$, it holds that 
\begin{equation}\label{kernel:eq-prelim-estimate-weights}
t_0^{-\frac{\eta}{2}} \| f \|_{L_t^1 L_x^\infty(I)} \lesssim M^{-2} \| t^{-\frac{\eta}{2}} f \|_{L_t^\infty L_x^\infty(I)}. 
\end{equation} 
Indeed, this can be obtained easily by distinguishing the cases $t_0 \gg M^{-2}$ and $t_0 \lesssim M^{-2}$. In the first case, one uses that $t^{-\eta/2}\sim t_0^{-\eta/2}$ on the interval $I$, and in the second case, one uses that $\| t^{\eta/2}\|_{L_t^1([0,t_0])}\lesssim t_0^{1+\eta/2}\lesssim M^{-2} t_0^{\eta/2}$. 
Moreover, by the Schauder estimates of Lemma~\ref{lemma:covariant-schauder-estimates}, we have that
\begin{equs}
\bigl\|\psi\bigr\|_{L_{t, x}^\infty([0, t_0])} &\lesssim \bigl\| \chi_M(t_0 - \place) H \bigr\|_{L_t^1 L_x^\infty([0, t_0])} \lesssim  M^{-2} \|H\|_{L_{t, x}^\infty([0, t_0])}, \label{eq:commutator-intermediate-psi-bound}\\
\bigl\|\covd \psi\bigr\|_{L_{t,x}^\infty([0, t_0])} &\lesssim \bigl(1 + \Fnorm{F}_{T, 0} \bigr) \sup_{t \in (0, t_0]} \bigl\|(t - s)^{-1/2} \chi_M(t_0 - s) H\bigr\|_{L_s^1 L_x^\infty([0, t])} \\
&\lesssim \bigl(1 + \Fnorm{F}_{T, 0} \bigr)M^{-1} \bigl\|H\bigr\|_{L_{t,x}^\infty([0, t_0])}. \label{eq:commutator-intermediate-covd-psi-bound}
\end{equs}
Next, we compute
\begin{equs}
&\quad (\ptl_t - \Delta_B) (\phi - \Qc(G, \psi)) \\ 
&= (\ptl_t - \Delta_B) \phi-\Qc((\ptl_t -\Delta_B) G, \psi) + 2 \Qc(\covd^\ell G, \covd_\ell \psi) - \Qc(G,(\partial_t - \Delta_B)\psi) \\ 
&= -\Qc((\ptl_t -\Delta_B) G, \psi) + 2 \Qc(\covd^\ell G, \covd_\ell \psi). 
\end{equs}
Here, we used that the $(\ptl_t - \Delta_B) \phi$- and $\Qc(G,(\partial_t - \Delta_B)\psi)$-terms cancel each other. From this, we then obtain that 
\begin{equ}
\phi - \Qc(G, \psi) = -\Duh_B(\Qc((\ptl_t -\Delta_B) G, \psi)) + 2 \Duh_B(\Qc(\covd^\ell G, \covd_\ell\psi)). 
\end{equ}
Each of the terms above may be estimated by the Schauder estimates of Lemma~\ref{lemma:covariant-schauder-estimates}. In particular, we estimate the first term by
\begin{align*}
t_0^{-\frac{\eta}{2}} \bigl\|\Duh_B(\Qc((\ptl_t - \Delta_B) G, \psi))(t_0)\bigr\|_{L_{x}^\infty}   
 &\lesssim  t_0^{-\frac{\eta}{2}} \bigl\|\Qc((\ptl_t - \Delta_B) G, \psi)\bigr\|_{L_t^1 L_x^\infty(I)} \\
&\lesssim   t_0^{-\frac{\eta}{2}}\bigl\|(\ptl_t - \Delta_B) G \bigr\|_{L_t^1L_x^\infty(I)} 
\bigl\| \psi \bigr\|_{L_{t,x}^\infty([0, t_0])} \\ 
&\lesssim M^{-4} \bigl\| t^{-\frac{\eta}{2}}(\ptl_t - \Delta_B) G \bigr\|_{L_{t,x}^\infty([0, t_0])} \bigl\|H\bigr\|_{L_{t,x}^\infty([0, t_0])}. 
\end{align*}
In the first line, we used that the support of $\psi$ is contained in the interval $I$. In the final line, we used \eqref{kernel:eq-prelim-estimate-weights} and \eqref{eq:commutator-intermediate-psi-bound}. Similarly, applying \eqref{eq:commutator-intermediate-covd-psi-bound} in place of \eqref{eq:commutator-intermediate-psi-bound}, we estimate the second term
\begin{align*}
t_0^{-\frac{\eta}{2}}\bigl\|\Duh_B(\Qc(\covd^\ell G, \covd_\ell \psi))(t_0)\bigr\|_{L_{x}^\infty} 
&\lesssim 
t_0^{-\frac{\eta}{2}} \bigl\|\Qc(\covd^\ell G, \covd_\ell\psi)\bigr\|_{L_t^1 L_x^\infty(I)} \\
&\lesssim  t_0^{-\frac{\eta}{2}} \bigl\| \covd G \bigr\|_{L_t^1 L_x^\infty(I)} \bigl\| \covd \psi \bigr\|_{L_{t,x}^\infty([0, t_0])} \\ 
&\lesssim  \bigl(1 + \Fnorm{F}_{T, 0} \bigr)M^{-3} \bigl\|t^{-\frac{\eta}{2}}\covd G \bigr\|_{L_{t,x}^\infty([0, t_0])} \bigl\|H\bigr\|_{L_{t,x}^\infty([0, t_0])},
\end{align*}
which completes the proof. 
\end{proof}

\begin{lemma}[Derivative commutator estimate]\label{kernel:lem-commutator-derivative}
Let $G \colon [0, T] \times \T^2 \ra \genV$, let $M \geq N_{\lo}$, and let $i \in [2]$. Then, we have that
\begin{equ}
\bigl\|\covd_i \Duh_{B, M}(G) - \Duh_{B, M}(\covd_i  G)\bigr\|_{L_{t, x}^\infty([0, T])} \lesssim \lossfactor_T(F) T^{-(1-\roughdataweight)/2} M^{-(2-\roughdataweight)}  \bigl\|G\bigr\|_{L_{t,x}^\infty([0, T])}.
\end{equ}
\end{lemma}
\begin{proof}
Fix $t_0 \in [0, T]$. We may write
\begin{equ}
\bigl(\covd_i \Duh_{B, M}(G) - \Duh_{B, M}(\covd_i G) \bigr)(t_0) = (\covd_i \phi - \psi)(t_0),
\end{equ}
where $\phi, \psi \colon [0, t_0] \times \T^2 \ra \genV$ respectively satisfy
\begin{equs}[2]
(\ptl_t - \Delta_B) \phi &= \chi_M(t_0 - \place) G, &\quad \phi(0) &= 0, \\
(\ptl_t - \Delta_B) \psi &= \chi_M(t_0 - \place) \covd_i G, &\quad \psi(0) &= 0.
\end{equs}
Note that $\phi(t) = 0$ unless $|t_0 - t| \lesssim M^{-2}$. We may compute (applying \eqref{eq:cov-heat-operator-covd-commutator} and Remark~\ref{remark:A-YM-GF-drop-term})
\begin{equ}
(\ptl_t - \Delta_B) (\covd_i \phi - \psi) = 2 F_i^{~\ell} \covd_\ell \phi. 
\end{equ}
Thus, by Lemma~\ref{lemma:covariant-schauder-estimates} (in the following, all space-time norms are on $[0, t_0] \times \T^2$ unless otherwise specified) we have that
\begin{equs}
\bigl\|\covd_i \phi - \psi\bigr\|_{L_{t,x}^\infty} \lesssim \bigl\|F \covd \phi \bigr\|_{L_t^1 L_x^\infty} &\lesssim \bigl\|t^{\roughdataweight/2} F\bigr\|_{L_t^2 L_x^\infty} \bigl\|t^{-\roughdataweight/2} \covd \phi\bigr\|_{L_t^2 L_x^\infty} \\
&\lesssim M^{-1} \Fnorm{t^{-(1-\roughdataweight)/2} F}_{T, 0} \bigl\|t^{-\roughdataweight/2} \covd \phi\bigr\|_{L_{t,x}^\infty} \\
&\lesssim \lossfactor_T(F)^{1/100} M^{-1}  T^{-(1-\roughdataweight)/2} \bigl\|t^{-\roughdataweight/2} \covd \phi\bigr\|_{L_{t,x}^\infty},
\end{equs}
where in the third inequality, we used that $\phi(t) = 0$ unless $|t_0 - t| \lesssim M^{-2}$, and in the final inequality, we used the definition \eqref{local:eq-lossfactor} of $\lossfactor_T(F)$. Therefore to finish, it suffices to estimate
\begin{equ}\label{eq:duhamel-derivative-commutator-intermediate}
\bigl\|t^{-\roughdataweight/2} \covd \phi\bigr\|_{L_{t,x}^\infty} \lesssim \bigl(1 + \Fnorm{F}_{T, 0}\bigr) M^{-(1-\roughdataweight)} \bigl\|G\bigr\|_{L_{t,x}^\infty}.
\end{equ}
(Again by \eqref{local:eq-lossfactor}, we may absorb $1 + \Fnorm{F}_{T, 0} \leq \lossfactor_T(F)^{1/100}$.) To this end, by the last estimate in Lemma~\ref{lemma:higher-regularity-covd-ibp}, we have that for any $t_1 \in [0, t_0]$,
\begin{equs}
\bigl\|&t^{1/2} \covd \phi\bigr\|_{L_{t,x}^\infty([0, t_1])} \lesssim \bigl(1 + \Fnorm{F}_{T, 0}\bigr)\Bigl( \bigl\|\chi_M(t_0 - \place) G \bigr\|_{L_t^1 L_x^\infty([0, t_1])} \\
&+ \bigl\|t^{1/2} \chi_M(t_0 - \place) G \bigr\|_{L_t^2 L_x^\infty([0, t_1])} + \sup_{u \in [0, t_1]} \bigl\|t^{1/2} (u-t)^{-1/2} \chi_M(t_0 - \place) G \bigr\|_{L_t^1 L_x^\infty([0,u])}\Bigr).
\end{equs}
We may estimate
\begin{equs}
\bigl\|\chi_M(t_0 - \place) G \bigr\|_{L_t^1 L_x^\infty([0, t_1])} &\lesssim \min(t_1, M^{-2}) \bigl\|G\bigr\|_{L_{t,x}^\infty} ,\\
\bigl\|t^{1/2} \chi_M(t_0 - \place) G \bigr\|_{L_t^2 L_x^\infty([0, t_1])} &\lesssim t_1^{1/2} \min(t_1^{1/2}, M^{-1}) \bigl\|G\bigr\|_{L_{t,x}^\infty} , \\
\sup_{u \in [0, t_1]} \bigl\|t^{1/2} (u-t)^{-1/2} \chi_M(t_0 - \place) G \bigr\|_{L_t^1 L_x^\infty([0, u])} &\lesssim t_1^{1/2} \min(t_1^{1/2}, M^{-1}) \bigl\|G\bigr\|_{L_{t,x}^\infty}.
\end{equs}
Combining, we thus obtain that for all $t_1 \in [0, t_0]$,
\begin{equ}
\bigl\|(\covd \phi)(t_1)\bigr\|_{L_x^\infty} \lesssim t_1^{\roughdataweight/2} M^{-(1-\roughdataweight)} \bigl(1 + \Fnorm{F}_{T, 0}\bigr) \bigl\|G\bigr\|_{L_{t,x}^\infty},
\end{equ}
from which the desired estimate \eqref{eq:duhamel-derivative-commutator-intermediate} follows.
\end{proof}

\subsection{Product and paraproduct estimates}\label{section:local-paraproduct-estimates}

In this subsection, we state and prove different covariant product and paraproduct estimates.

\begin{assumption}[Setup]\label{local:assumption-setup}
We make the following assumptions: 
\begin{enumerate}[label=(\roman*),topsep=0.2ex,itemsep=0ex]
\item The indices $a,b,c\in [3]$ are given. 
\item The bilinear map $\Qc\colon \Y^{(a)} \times \Y^{(b)} \rightarrow \Y^{(c)}$ is $\liegroup$-equivariant. 
\item We let  $U_\star=(U_M)_{M\geq N_{\lo}}$ where, for each $M\geq N_{\lo}$, $U_M \colon [0,T]\times\T^2 \rightarrow \Y^{(a)}$ and $U_M$ is supported on times $t\gtrsim M^{-2}$. 
\item We let  $V_\star=(V_N)_{N\geq N_{\lo}}$ where, for each $N\geq N_{\lo}$, $V_N \colon [0,T]\times\T^2 \rightarrow \Y^{(b)}$ and $V_N$ is supported on times $t\gtrsim N^{-2}$. 
\end{enumerate}
\end{assumption}

Our first lemma will later be used to control \eqref{local:eq-Y-final-new-7} and \eqref{local:eq-y-final-new-5}. While these are among the simplest terms in \eqref{local:eq-Xwhite-K-final}--\eqref{local:eq-y-final-new-last}, this nevertheless illustrates many of the ideas used to control other terms in \eqref{local:eq-Xwhite-K-final}--\eqref{local:eq-y-final-new-last}.

\begin{lemma}[Product estimate]\label{local:lem-product-estimate}
Let Assumption~\ref{local:assumption-setup} be satisfied. In addition, let \mbox{$\alpha,\beta,\gamma \in \R$} and 
$\eta_1,\eta_2,\nu \in (-\infty,0]$ satisfy
\begin{equation*}
0<\alpha,\beta<1, \quad \alpha+\beta>1, \quad 
\eta_1,\eta_2>-1, \quad 
\eta_2 \leq \nu, \quad 
\text{and} \quad \gamma +\nu - \eta_1 - \eta_2 < \beta -1. 
\end{equation*}
Then, it holds that 
\begin{equation}
\big\| \Qc( U_\star, \covd V_\star) \big\|_{\nl^{\gamma,\nu}}
\lesssim \lossfactor_T(F) T^{\frac{\alpha+\beta-\gamma-\nu+\eta_1+\eta_2-1}{2}} \big\| U_\star \big\|_{\sol^{\alpha,\eta_1}} \big\| V_\star \big\|_{\sol^{\beta,\eta_2}}.
\end{equation}
\end{lemma}

\begin{proof} We emphasise that, as in Assumption~\ref{local:assumption-setup}, all dyadic frequency scales in the proof will be bounded below by $N_{\lo}\sim T^{-\frac{1}{2}}$. To simplify the notation, set $\eta\coloneq \eta_1+\eta_2$. Using our Schauder estimates 
\eqref{local:eq-weighted-covariant-Schauder-3}
and \eqref{local:eq-weighted-covariant-Schauder-4}, where we choose $(\widetilde{\eta}_1,\widetilde{\eta}_2)=(0,\eta_2)$, we obtain that 
\begin{align*}
&\hspace{-5ex}\big\| \Qc( U_M, \covd V_N )  \big\|_{\nl^{\gamma,\nu}_K} \\ 
\lesssim \lossfactor_T(F) &\min\Big( 
K^{\gamma+\nu-\eta} \| t^{-\eta_1/2} U_M \|_{L_{t,x}^\infty} \| t^{-\eta_2/2} \covd V_N \|_{L_{t,x}^\infty}, \\
&  K^{\gamma+\nu-\eta_2} K \| U_M \|_{L_{t,x}^\infty} \| t^{-\eta_2/2} V_N \|_{L_{t,x}^\infty} 
+ K^{\gamma+\nu-\eta}  \| t^{-\eta_1/2} \covd U_M \|_{L_{t,x}^\infty} \| t^{-\eta_2/2} V_N \|_{L_{t,x}^\infty} \Big). 
\end{align*}
Using that $U_M$ is supported on times $t\gtrsim M^{-2}$, we have
\begin{equation*}
K^{-\eta_2} \| U_M \|_{L^\infty_{t,x}}
\lesssim K^{-\eta_2} M^{-\eta_1} \| t^{-\eta_1/2} U_M \|_{L^\infty_{t,x}} 
= K^{-\eta+\eta_1} M^{-\eta_1} \| t^{-\eta_1/2} U_M \|_{L^\infty_{t,x}}. 
\end{equation*} 
Using the definition of our solution spaces, we therefore obtain that 
\begin{align*}
&\quad \big\| \Qc( U_M, \covd V_N ) \big\|_{\nl^{\gamma,\nu}_K}  \\
&\lesssim \lossfactor_T(F) \min\bigl( N, \max( K^{1+\eta_1} M^{-\eta_1}, M)\bigr)
K^{\gamma+\nu-\eta} M^{-\alpha} N^{-\beta} \big\| U_M \big\|_{\sol^{\alpha,\eta_1}_M} \big\| V_N \big\|_{\sol^{\beta,\eta_2}_N}. 
\end{align*}
Since $\eta_1 \in (-1,0]$, we have that 
\begin{equation*}
\min\bigl( N, \max( K^{1+\eta_1} M^{-\eta_1}, M)\bigr)
\lesssim \min\bigl( N , \max( K,M )\bigr). 
\end{equation*}
Due to this, it remains to prove that
\begin{align}\label{local:eq-product-estimate-p1}
\sup_{K}\sum_{M,N}\min\bigl( N, \max(K,M)  \bigr) K^{\gamma+\nu-\eta} M^{-\alpha} N^{-\beta} \lesssim N_{\lo}^{1+\gamma+\nu-\eta-\alpha-\beta}.
\end{align}
To obtain \eqref{local:eq-product-estimate-p1}, we first use that $\beta \in (0,1)$, which allows us to estimate the sum over $N$ by 
\begin{align*}
\sum_{N}\min\bigl( N, \max(K,M)  \bigr) K^{\gamma+\nu-\eta} M^{-\alpha} N^{-\beta} 
\lesssim \max(K,M)^{1-\beta} K^{\gamma+\nu-\eta} M^{-\alpha}.
\end{align*}
Using $\alpha>0$ and $\alpha+\beta>1$, we can estimate the remaining sum over $M$ by 
\begin{align*}
\sum_M \max(K,M)^{1-\beta} K^{\gamma+\nu-\eta} M^{-\alpha}
&\lesssim K^{1+\gamma+\nu-\eta-\beta} \sum_M M^{-\alpha} + 
K^{\gamma+\nu-\eta} \sum_M M^{1-\alpha-\beta} \\ 
&\lesssim K^{1+\gamma+\nu-\eta-\beta}  N_{\lo}^{-\alpha} + K^{\gamma+\nu-\eta} N_{\lo}^{1-\alpha-\beta}. 
\end{align*}
Since both $1+\gamma+\nu-\eta-\beta$ and $\gamma+\nu-\eta$ are negative, this implies the desired result. 
\end{proof}
          
For the rest of this subsection, $U_\star$ will be at negative regularity and $V_\star$ will be at positive regularity. We encourage the reader to think of $U_\star$ as $\linear[\oeps][][][\ostar]$ and of $V_\star$ as $\Wrem_{\eps,\star}$ or $Y_{\eps,\star}$. Our next lemma, which is another product estimate, will later be used to control \eqref{local:eq-Xwhite-K-final}--\eqref{local:eq-phix-K-final}. 

\begin{lemma}[Product estimate at negative regularity]\label{local:lem-product-negative-regularity}
Let Assumption~\ref{local:assumption-setup} be satisfied. In addition,  let $\alpha\in (-1,0)$, let $\beta,\gamma \in \R$,
and let $\eta,\nu \in (-\infty,0]$. Furthermore, assume that 
\begin{equation}\label{local:lem-product-neg-c1}
-1<\alpha+\eta<0, \quad \beta>0, \quad \alpha+\beta>0, \quad - 1 < \eta \leq \nu, \quad \text{and} \quad \gamma+\nu-\eta<-1+\alpha-\frac{2}{r}. 
\end{equation}
Then, it holds for all $\eps>0$ that  
\begin{equation}\label{local:eq-product-neg-e1}
\big\|  \covd \Qc( U_\star, V_\star ) \big\|_{\nl^{\gamma,\nu}}
\lesssim_\eps \lossfactor_T(F) T^{\frac{2\alpha+\beta-\gamma-\nu+\eta-1}{2}} T^{-\frac{1}{r}-\eps} \big\| U_\star \big\|_{\modsol^{\alpha}} \big\| V_\star \big\|_{\sol^{\beta,\eta}}.
\end{equation}
In addition to \eqref{local:eq-product-neg-e1}, we also have the same estimate but with $\Qc(U_\star,V_\star)$ replaced by $\Qc(U_\star \parall V_\star)$, $\Qc(U_\star \parasim V_\star)$, or $\Qc(U_\star \paragg V_\star)$.
\end{lemma}

\begin{remark}
We emphasise that the estimate \eqref{local:eq-product-neg-e1} cannot be obtained this way if the $\modsol^{\alpha}$-norm is replaced by the $\sol^\alpha$-norm. The reason is that, since $\alpha$ is negative, the $\sol^\alpha$-norm cannot be used to gain powers of the frequency of the $U_\star$-terms.
\end{remark}

\begin{proof}[Proof of Lemma~\ref{local:lem-product-negative-regularity}]
Without loss of generality, we may assume that $\eps$ is sufficiently small depending on all parameters appearing in \eqref{local:lem-product-neg-c1}. Using our weighted Schauder estimate \eqref{local:eq-weighted-covariant-Schauder-2} and Lemma~\ref{local:lem-solution-spaces-hoelder}, we obtain
\begin{align*}
\big\| \covd \Qc( U_M ,V_N ) \big\|_{\nl^{\gamma,\nu}_K}
&\lesssim \lossfactor_T(F) K^{\gamma+\nu-\eta+1-\alpha+\frac{2}{r}+\eps}
\big\| t^{-\eta/2} \Qc( U_M , V_N )  \big\|_{\Cs_t^{\frac{\alpha-\eps}{2}}L_x^r}\\ 
&\lesssim  \lossfactor_T(F) K^{\gamma+\nu-\eta+1-\alpha+\frac{2}{r}+\eps} \big\| U_M \big\|_{\Cs_t^{\frac{\alpha-\eps}{2}}L_x^r}
\big\| t^{-\eta/2} V_N \big\|_{\Cs_t^{-\frac{\alpha-2\eps}{2}} L_x^\infty} \\ 
&\lesssim \lossfactor_T(F) K^{\gamma+\nu-\eta+1-\alpha+\frac{2}{r}+\eps} M^{-\eps} N^{-\alpha-\beta+2\eps} 
\big\| U_M \big\|_{\modsol^{\alpha}_M} 
\big\| V_N \big\|_{\sol^{\beta,\eta}_N}. 
\end{align*}
Due to our conditions from \eqref{local:lem-product-neg-c1}, the exponents of $K$, $M$, and $N$ above are all negative. Since the sum of the exponents equals $-2\alpha-\beta+\gamma+\nu-\eta+1 + \frac{2}{r}+2\eps$, 
this readily implies the desired estimate~\eqref{local:eq-product-neg-e1}. 
\end{proof}

We now turn to low$\times$high-estimates, which will be used to control the terms in \eqref{local:eq-Y-final-new-5}, \eqref{local:eq-Y-final-new-6}, and the second and third summands in \eqref{local:eq-y-final-new-4}. 

\begin{lemma}[Low$\times$high-estimates]\label{local:lem-low-high}
Let Assumption~\ref{local:assumption-setup} be satisfied. In addition,  let $\alpha,\beta,\gamma \in \R$
and let $\eta,\nu \in (-\infty,0]$. Furthermore, assume that the conditions  
\begin{equation}\label{local:eq-low-high-c1}
\alpha<0, \quad 0<\beta<1, \quad \gamma+\nu-\eta<0, \quad -2 < \eta \leq \nu, \quad  \text{and} \quad (1-\theta) \alpha + \beta > (1-\theta) 
\end{equation}
are satisfied.  Then, it holds that 
\begin{equation}\label{local:eq-low-high-e1}
\big\| \Qc ( \covd U_\star \parall  V_\star ) \big\|_{\nl^{\gamma,\nu}} 
\lesssim \lossfactor_T(F) T^{\frac{(1-\theta)\alpha+\beta-\gamma-\nu+\eta+\theta-1}{2}} \big\| U_\star \big\|_{\sol^\alpha} \big\| V_\star \big\|_{\sol^{\beta,\eta}}.
\end{equation}
If, in addition to \eqref{local:eq-low-high-c1}, the conditions   
\begin{equation}\label{local:eq-low-high-c2}
\begin{aligned}
\eta>-1 \quad \text{and} \quad (1-\theta) \alpha + \beta >  1+\gamma +\nu-\eta
\end{aligned}
\end{equation}
are satisfied, it also holds that 
\begin{equation}\label{local:eq-low-high-e2}
\big\| \Qc( U_\star \parall  \covd V_\star )
\big\|_{\nl^{\gamma,\nu}}
\lesssim \lossfactor_T(F)  T^{\frac{(1-\theta)\alpha+\beta-\gamma-\nu+\eta-1}{2}}  \big\| U_\star \big\|_{\sol^\alpha} \big\| V_\star \big\|_{\sol^{\beta,\eta}}.
\end{equation}
\end{lemma}

\begin{remark}\label{local:rem-theta-low-high}
The last condition in \eqref{local:eq-low-high-c1} is one of the reasons for the parameter $\theta$ in our definition of covariant paraproducts (Definition~\ref{local:def-paraproduct}). If $\theta$ were set equal to zero in the last condition in~\eqref{local:eq-low-high-c1}, the condition would become $\alpha+\beta>1$, which would not allow us to control low$\times$high-interactions between $\covd \linear[B]$ and $Z$. We further emphasise that the last  condition in \eqref{local:eq-low-high-c1} has no counterpart in the usual, non-covariant low$\times$high-estimate. The reason is that the non-covariant Littlewood--Paley decomposition obeys the Littlewood--Paley trichotomy (see \eqref{local:eq-trichotomy}). But, as already discussed in Section~\ref{section:Littlewood--Paley}, the Littlewood--Paley trichotomy is not available in our setting.
\end{remark}

\begin{proof} We first prove the estimate \eqref{local:eq-low-high-e1}. Using Definition~\ref{local:def-solution-space} and our weighted Schauder estimate \eqref{local:eq-weighted-covariant-Schauder-3}, we obtain that 
\begin{equation*}
\big\| \Qc( \covd U_M , V_N ) \big\|_{\nl^{\gamma,\nu}_K}
\lesssim \lossfactor_T(F) K^{\gamma+\nu-\eta} M^{1-\alpha} N^{-\beta} \big\| U_M \big\|_{\sol^\alpha_M} \big\| V_N \big\|_{\sol^{\beta,\eta}_N}.
\end{equation*}
It then remains to show that, under our assumption \eqref{local:eq-low-high-c1}, 
\begin{equation}\label{local:eq-low-high-p1}
\sup_K \sum_{\substack{M,N\colon \\ M \leq N^{1-\theta}}} K^{\gamma+\nu-\eta} M^{1-\alpha} N^{-\beta} \lesssim N_{\lo}^{\gamma+\nu-\eta+(1-\theta)(1-\alpha)-\beta}. 
\end{equation}
This can be easily obtained by first summing in $M$ and then summing in $N$, and we omit the computation. Next, we prove the estimate \eqref{local:eq-low-high-e2}. Using Definition~\ref{local:def-solution-space} and our Schauder estimates \eqref{local:eq-weighted-covariant-Schauder-3} and \eqref{local:eq-weighted-covariant-Schauder-4}, we obtain that 
\begin{equation*}
\big\| \Qc(  U_M , \covd V_N ) \big\|_{\nl^{\gamma,\nu}_K}
\lesssim \lossfactor_T(F) \min \bigl( \max(K,M), N \bigr) K^{\gamma+\nu-\eta} M^{-\alpha} N^{-\beta} \big\| U_M \big\|_{\sol^\alpha_M} \big\| V_N \big\|_{\sol^{\beta,\eta}_N}.
\end{equation*}
Since the conditions in \eqref{local:eq-low-high-c1} are still satisfied, the region in which $\min(\max(K,M),N)$ equals $M$ can be treated as in our proof of \eqref{local:eq-low-high-e1}. It then remains to show that, under our assumptions in \eqref{local:eq-low-high-c1} and~\eqref{local:eq-low-high-c2}, 
\begin{equation}\label{local:eq-low-high-p2}
\sup_K \sum_{\substack{M,N\colon \\ M \leq N^{1-\theta}}} \min(K,N) K^{\gamma+\nu-\eta} M^{-\alpha} N^{-\beta} \lesssim N_{\lo}^{1+\gamma+\nu-\eta-(1-\theta)\alpha-\beta}. 
\end{equation}
We note that, since $\theta>0$, the right-hand side in \eqref{local:eq-low-high-p2} is worse than in \eqref{local:eq-low-high-p1}. In order to obtain \eqref{local:eq-low-high-p2}, we estimate 
\begin{align*}
&\quad \sum_{\substack{M,N\colon \\ M \leq N^{1-\theta}}} \min(K,N) K^{\gamma+\nu-\eta} M^{-\alpha} N^{-\beta}  
\lesssim \sum_{N} \min(K,N) K^{\gamma+\nu-\eta} N^{-(1-\theta)\alpha-\beta} \\ 
&\lesssim \sum_{N \leq K} K^{\gamma+\nu-\eta} N^{1-(1-\theta)\alpha-\beta} + \sum_{N>K} K^{1+\gamma+\nu-\eta} N^{-(1-\theta)\alpha-\beta} 
\lesssim  K^{1+\gamma+\nu-\eta-(1-\theta)\alpha-\beta}.
\end{align*}
In the estimate of the first sum, we used that $(1-\theta)\alpha+\beta<1$, which follows from $\alpha<0$ and $\beta<1$. In the estimate of the second sum, we used that $(1-\theta)\alpha+\beta>1-\theta>0$. 
Since the final exponent $1+\gamma+\nu-\eta-(1-\theta)\alpha-\beta$ is negative, this implies \eqref{local:eq-low-high-p2}.
\end{proof}

In our next lemma, we prove different high$\times$high-estimates, which will later be used to control \eqref{local:eq-Y-final-new-3}, \eqref{local:eq-Y-final-new-4}, \eqref{local:eq-y-final-new-3}, and the first summand in \eqref{local:eq-y-final-new-4}.  

\begin{lemma}[High$\times$high-estimates]\label{local:lem-high-high}
Let Assumption~\ref{local:assumption-setup} be satisfied. In addition,  let $\alpha,\beta,\gamma \in \R$
and let $\eta,\nu \in (-\infty,0]$. Furthermore, assume that the conditions  
\begin{equation}\label{local:eq-high-high-c1}
\begin{aligned}
&\alpha<0,\quad \beta>0, \quad \gamma+\nu-\eta<0, \quad  \eta>-2, \quad \nu > -1, \quad \nu \geq \eta,  \\
&\alpha+(1-\theta)\beta>1+\gamma+\nu-\eta, \quad \text{and} \quad \alpha + (1-\theta)(\beta-1)>0
\end{aligned}
\end{equation}
are satisfied.  Then, it holds that 
\begin{equation}\label{local:eq-high-high-e1}
\big\| \Qc( U_\star \parasim \covd V_\star)\big\|_{\nl^{\gamma,\nu}} + 
\big\| \Qc( \covd U_\star \parasim  V_\star)\big\|_{\nl^{\gamma,\nu}} \lesssim \lossfactor_T(F)  T^{\frac{\alpha+(1-\theta)\beta - \gamma -\nu +\eta -1}{2}} \big\| U_\star  \big\|_{\sol^\alpha} \big\| V_\star \big\|_{\sol^{\beta,\eta}}.
\end{equation}
Alternatively, assume that
\begin{equation}\label{local:eq-high-high-c2}
\alpha<0, \quad 0<\beta<1, \quad \nu \geq \eta >-1,  \quad \text{and} \quad \alpha+(1-\theta) \beta > \max(0,\gamma+\nu-\eta+1). 
\end{equation}
Then, it holds that 
\begin{equation}\label{local:eq-high-high-e2}
\big\|  \covd \Qc  ( U_\star \parasim  V_\star ) \big\|_{\nl^{\gamma,\nu}} \lesssim \lossfactor_T(F) T^{\frac{\alpha+(1-\theta)\beta - \gamma - \nu+\eta -1}{2}} \big\| U_\star  \big\|_{\sol^\alpha} \big\| V_\star \big\|_{\sol^{\beta,\eta}}.
\end{equation}
\end{lemma}

\begin{remark} In the proof of Lemma~\ref{local:lem-nl-Y-high-high} below, in which $V_\star$ will be chosen as $\higgsdet_\star$, it is important that \eqref{local:eq-high-high-c1} only imposes the lower bound $\eta>-2$, and not $\eta>-1$. 
\end{remark}

\begin{proof}
We begin with our proof of \eqref{local:eq-high-high-e1}. Using our Schauder estimate \eqref{local:eq-weighted-covariant-Schauder-3}, we obtain that 
\begin{align*}
\big\| \Qc ( U_M , \covd V_N ) \big\|_{\nl^{\gamma,\nu}_K} 
&\lesssim \lossfactor_T(F) K^{\gamma+\nu-\eta} \| U_M \|_{L_{t,x}^\infty} \| t^{-\eta/2} \covd V_N \|_{L_{t,x}^\infty}\\
&\lesssim \lossfactor_T(F) K^{\gamma+\nu-\eta} M^{-\alpha} N^{1-\beta}  \| U_M \|_{\sol^{\alpha}_M} \| V_N \|_{\sol^{\beta,\eta}_N}. 
\end{align*}
Since both $\Qc( U_\star \parasim \covd V_\star)$ and $\Qc( \covd U_\star \parasim V_\star)$ appear in \eqref{local:eq-high-high-e1}, we will also need to estimate $\Qc ( \covd U_M ,  V_N )$. Since $U_M$ has negative regularity, we will want to use our weighted covariant Schauder estimate \eqref{local:eq-weighted-covariant-Schauder-4}, which was proven using integration by parts. However, since our parameter $\eta$ is allowed to be less than or equal to $-1$, we need to be careful. To this end, we choose any parameter $\widetilde{\eta}$ satisfying $\widetilde{\eta}>-1$, $\widetilde{\eta}\geq \eta$, and $\widetilde{\eta}\leq \min(\nu,\eta+1)$, which is possible due to our assumptions on $\eta$ and $\nu$. Since $V_N$ is supported on times $t\gtrsim N^{-2}$, we may estimate 
\begin{equation*}
\| t^{-\widetilde{\eta}/2} V_N \|_{L_{t,x}^\infty} 
\lesssim N^{\widetilde{\eta}-\eta} N^{-\beta} \| V_N \|_{\sol^{\beta,\eta}_N}. 
\end{equation*}
Together with our weighted Schauder estimate \eqref{local:eq-weighted-covariant-Schauder-4}, we then have that
\begin{align*}
 &\quad \big\| \Qc( \covd U_M ,  V_N ) \big\|_{\nl^{\gamma,\nu}_K} \\   
&\lesssim \lossfactor_T(F) K^{\gamma} \min\Big( K^{\nu-\eta} \| \covd U_M \|_{L_{t,x}^\infty} \| t^{-\eta/2} V_N \|_{L_{t,x}^\infty}, \\
&\hspace{15ex} K^{\nu-\eta} \|  U_M \|_{L_{t,x}^\infty} \| t^{-\eta/2} \covd V_N \|_{L_{t,x}^\infty}
+ K^{\nu-\widetilde{\eta}+1} \| U_M \|_{L_{t,x}^\infty} \| t^{-\widetilde{\eta}/2} V_N \|_{L_{t,x}^\infty}
\Big) \\
&\lesssim \lossfactor_T(F) K^{\gamma+\nu-\eta} M^{-\alpha} N^{-\beta} \min\bigl( M , \max \bigl( N, K^{1+\eta-\widetilde{\eta}}N^{\widetilde{\eta}-\eta}\bigr)\bigr)\| U_M \|_{\sol^{\alpha}_M} \| V_N \|_{\sol^{\beta,\eta}_N}.
\end{align*} 
Since  $0\leq \widetilde{\eta}-\eta \leq 1$, the minimum above can be estimated by 
\begin{align*}
	\min\bigl( M , \max \bigl( N, K^{1+\eta-\widetilde{\eta}}N^{\widetilde{\eta}-\eta}\bigr)\bigr)
\leq \min\bigl( M, \max(K,N)\bigr) \leq \min(K,M)+N. 
\end{align*}
After combining our estimates of both $\Qc( U_M , \covd V_N)$ and $\Qc( \covd U_M ,  V_N )$, it then only remains to show that 
\begin{equation}\label{local:eq-high-high-p1}
\begin{aligned}
&\quad K^{\gamma+\nu-\eta} \sum_{\substack{M,N\colon \\ M\sim_\theta N}} \min(K,M)  M^{-\alpha} N^{-\beta}  
+  K^{\gamma+\nu-\eta} \sum_{\substack{M,N\colon \\ M\sim_\theta N}}   M^{-\alpha} N^{1-\beta} \\
&\lesssim  N_{\lo}^{1+\gamma+\nu-\eta-\alpha-(1-\theta)\beta}.
\end{aligned}
\end{equation}
We first estimate the first sum in \eqref{local:eq-high-high-p1}. By using that $\beta>0$, which allows us to sum over $N$, it can be bounded by 
\begin{equation}\label{local:eq-high-high-p1p}
\begin{aligned}
&K^{\gamma+\nu-\eta} \sum_{\substack{M,N\colon \\ M\sim_\theta N}} \min(K,M)  M^{-\alpha} N^{-\beta} \\ 
&\lesssim K^{\gamma+\nu-\eta} \sum_{\substack{M\colon \\ M \lesssim K}}  M^{1-\alpha-(1-\theta)\beta}
+ K^{1+\gamma+\nu-\eta} \sum_{\substack{M\colon \\ M \gtrsim K}}  M^{-\alpha-(1-\theta)\beta}. 
\end{aligned}
\end{equation}
Due to the third and seventh conditions in \eqref{local:eq-high-high-c1}, the first sum in \eqref{local:eq-high-high-p1p} can be bounded by first estimating $K^{\gamma+\nu-\eta}\lesssim M^{\gamma+\nu-\eta}$ and then summing over $M \gtrsim N_{\lo}$. Due to the seventh and eighth conditions in \eqref{local:eq-high-high-c1}, the second sum in \eqref{local:eq-high-high-p1p} can be estimated by first summing over $M$ and then using that $K\gtrsim N_{\lo}$.
We now turn to the second sum in \eqref{local:eq-high-high-p1}. To this end, we first note that the first and last conditions in \eqref{local:eq-high-high-c1} imply that $\beta>1$. By first using that $\beta>1$,  then using that $\alpha+(1-\theta)(\beta-1)>0$, and finally using that $\gamma+\nu-\eta<0$, the second sum in \eqref{local:eq-high-high-p1} can be estimated by 
\begin{align*}
& K^{\gamma+\nu-\eta} \sum_{\substack{M,N\colon \\ M\sim_\theta N}}  M^{-\alpha} N^{1-\beta}   
\lesssim K^{\gamma+\nu-\eta} \sum_{M} M^{-\alpha-(1-\theta)(\beta-1)} \\ 
&\lesssim K^{\gamma+\nu-\eta} N_{\lo}^{-\alpha-(1-\theta)(\beta-1)}
\lesssim N_{\lo}^{\gamma+\nu-\eta -\alpha-(1-\theta)(\beta-1)},
\end{align*}
which is better than the right-hand side of \eqref{local:eq-high-high-p1}. This completes the proof of \eqref{local:eq-high-high-e1}, and it now remains to prove \eqref{local:eq-high-high-e2}. Using Definition~\ref{local:def-solution-space} and our weighted Schauder estimate \eqref{local:eq-weighted-covariant-Schauder-1} with $m=1$, we obtain that 
\begin{equation*}
\big\| \covd \Qc( U_M , V_N ) \big\|_{\nl^{\gamma,\nu}_K}
\lesssim \lossfactor_T(F) K^{1+\gamma+\nu-\eta} M^{-\alpha} N^{-\beta} \big\| U_M \big\|_{\sol^\alpha_M} \big\| V_N \big\|_{\sol^{\beta,\eta}_N}.
\end{equation*}
Alternatively, using Definition~\ref{local:def-solution-space}, our weighted Schauder estimate \eqref{local:eq-weighted-covariant-Schauder-1} with $m=0$, and the product rule 
\begin{equation*}
\covd \Qc( U_M , V_N ) 
=  \Qc( \covd U_M , V_N ) + \Qc( U_M , \covd  V_N ), 
\end{equation*}
we obtain that 
\begin{equation*}
\big\| \covd \Qc( U_M , V_N )  \big\|_{\nl^{\gamma,\nu}_K}
\lesssim \lossfactor_T(F) \max(M,N) K^{\gamma+\nu-\eta} M^{-\alpha} N^{-\beta} \big\| U_M \big\|_{\sol^\alpha_M} \big\| V_N \big\|_{\sol^{\beta,\eta}_N}.
\end{equation*}
It then remains to show that, under our assumption \eqref{local:eq-high-high-c2}, 
\begin{equation}\label{local:eq-high-high-p2}
\sum_{\substack{M,N\colon \\ M\sim_\theta N}} \min\bigl( K, \max(M,N)\bigr) K^{\gamma+\nu-\eta} M^{-\alpha} N^{-\beta} \lesssim N_{\lo}^{1+\gamma+\nu-\eta-\alpha-(1-\theta)\beta}.
\end{equation}
Since $\beta$ is larger than $\alpha$ and since both $M\sim_\theta N$ and $\max(M,N)$ are symmetric in $M$ and $N$, we may insert the restriction $M\geq N$ into 
\eqref{local:eq-high-high-p2}. Then, we estimate
\begin{align*}
&\quad \sum_{\substack{M,N\colon \\ M^{1-\theta} \leq N \leq M}}\hspace{-3ex} \min\bigl( K, M\bigr) 
K^{\gamma+\nu-\eta} M^{-\alpha} N^{-\beta} \\ 
&\lesssim
\sum_{M}\min(K,M) K^{\gamma+\nu-\eta} M^{-\alpha-(1-\theta)\beta}  
\lesssim  K^{1+\gamma+\nu-\eta-\alpha-(1-\theta)\beta}.
\end{align*}
The last inequality above is obtained by splitting the sum into the two regions $M\leq K$ and $M>K$ and using the two conditions  $1-\alpha-(1-\theta)\beta>0$ and $-\alpha-(1-\theta)\beta<0$, which both follow from \eqref{local:eq-high-high-c2}. Furthermore, due to \eqref{local:eq-high-high-c2}, $1+\gamma+\nu-\eta-\alpha-(1-\theta)\beta$ is negative, and we therefore obtain \eqref{local:eq-high-high-p2}. 
\end{proof}

In our last lemma, we prove a high$\times$low-estimate, which will later be used to control \eqref{local:eq-Y-final-new-6} and the last summand in \eqref{local:eq-y-final-new-4}.   

\begin{lemma}[High$\times$low-estimate]\label{local:lem-high-low}
Let Assumption~\ref{local:assumption-setup} be satisfied. In addition,  let $\alpha,\beta,\gamma \in \R$
and let $\eta,\nu \in (-\infty,0]$. Furthermore, assume that 
\begin{alignat}{3}
&-1 -\eta  < \alpha < 0, &  \hspace{5ex}  &0<\beta<1, \hspace{15ex} -1 < \eta \leq \nu, \label{local:eq-high-low-c1a} \\
&\theta (1+\alpha) > (1-\theta) (1-\alpha-\beta),  & \hspace{5ex} &\gamma+\nu-\eta < \alpha - \frac{1-\theta}{\theta} (1-\alpha-\beta) - \frac{2}{r}. \label{local:eq-high-low-c1b}
\end{alignat}
Then, it holds for all $\eps>0$ that 
\begin{equation}\label{local:eq-high-low-e1}
\big\| \Qc( U_\star \paragg \covd V_\star ) \big\|_{\nl^{\gamma,\nu}}
\lesssim_\eps \lossfactor_T(F) T^{\frac{-\gamma-\nu+\eta+\alpha}{2} - \frac{1-\theta}{2\theta} (1-\alpha-\beta)}  T^{-\frac{1}{r}-\eps} \big\| U_\star \big\|_{\modsol^{\alpha}} \big\| V_\star \big\|_{\sol^{\beta,\eta}}. 
\end{equation}
\end{lemma}

\begin{remark}\label{local:rem-theta-high-low}
In order to control the high$\times$low-interactions between $\linear[B]$ and $\covd \Wrem$,  we will later use Lemma~\ref{local:lem-high-low} with $\alpha=-\kappa$, $\beta=1-\reg_1$, and $\gamma=-1+\reg_2$. For these parameters, the two conditions in \eqref{local:eq-high-low-c1b} are satisfied only because $\reg_1/\theta$ is still small. Thus, Lemma~\ref{local:lem-high-low} is another reason for the parameter $\theta$ in our definition of the covariant paraproducts (Definition~\ref{local:def-paraproduct}). 
\end{remark}

\begin{proof} Since our desired estimate \eqref{local:eq-high-low-e1} becomes worse as $\eps$ increases, we may assume that $\eps>0$ is sufficiently small depending on the parameters in \eqref{local:eq-high-low-c1a}--\eqref{local:eq-high-low-c1b}. For technical reasons, we let $\eps^\prime>0$ be sufficiently small depending on both the parameters in  \eqref{local:eq-high-low-c1a}--\eqref{local:eq-high-low-c1b} and $\eps$. Furthermore, we let $\alpha^\prime\in [-1,\alpha]$ remain to be chosen. Using Definition~\ref{local:def-modified-solution-space}, Lemma \ref{local:lem-solution-spaces-hoelder}, and our weighted Schauder estimate \eqref{local:eq-weighted-covariant-Schauder-2}, we first estimate
\begin{align*}
\big\| \Qc( U_M, \covd V_N ) \big\|_{\nl^{\gamma,\nu}_K}
&\lesssim \lossfactor_T(F) K^{\gamma+\nu-\eta-\alpha^\prime+\frac{2}{r}} \big\| t^{-\eta/2}\Qc( U_M, \covd V_N )\big\|_{\Cs_t^{\alpha^\prime/2} L_x^r} \\
&\lesssim \lossfactor_T(F)  K^{\gamma+\nu-\eta-\alpha^\prime+\frac{2}{r}} \big\|  U_M \big\|_{\Cs_t^{\alpha^\prime/2} L_x^r}
\big\| t^{-\eta/2} \covd V_N \big\|_{\Cs_t^{-(\alpha^\prime-\eps^\prime)/2} L_x^\infty} \\
&\lesssim \lossfactor_T(F) K^{\gamma+\nu-\eta-\alpha^\prime+\frac{2}{r}} M^{\alpha^\prime-\alpha} N^{1-\alpha^\prime-\beta+\eps^\prime} \big\| U_\star \big\|_{\modsol^{\alpha}} \big\| V_\star \big\|_{\sol^{\beta,\eta}}.
\end{align*}
Since $1-\beta$, $-\alpha^\prime$, and $\eps^\prime$ are all positive, it holds that 
\begin{equation*}
 M^{\alpha^\prime-\alpha} \sum_{\substack{N\colon \\ N \leq M^{1-\theta}}}  N^{1-\alpha^\prime-\beta+\eps^\prime}
\lesssim M^{\alpha^\prime-\alpha} M^{(1-\theta)(1-\alpha^\prime-\beta+\eps^\prime)} 
= M^{(1-\theta)(1-\alpha-\beta)+(1-\theta) \eps^\prime} M^{-\theta (\alpha-\alpha^\prime)}.
\end{equation*}
In order to obtain our desired estimate \eqref{local:eq-high-low-e1}, it then suffices to show that
\begin{equation}\label{local:eq-high-low-p1}
\begin{aligned}
&\quad\sum_{M} K^{\gamma+\nu-\eta-\alpha+\frac{2}{r}} M^{(1-\theta)(1-\alpha-\beta)+(1-\theta)\eps^\prime} \min_{\alpha^\prime=-1,\alpha} \bigl( K M^{-\theta}\bigr)^{\alpha-\alpha^\prime} \\
&\lesssim  K^{\gamma+\nu-\eta-\alpha+\frac{2}{r}+\frac{1-\theta}{\theta} (1-\alpha-\beta) + 2\eps}.
\end{aligned}
\end{equation}
The reason is that, due to the second condition in \eqref{local:eq-high-low-c1b} and since $\eps>0$ is sufficiently small, the exponent on the right-hand side of \eqref{local:eq-high-low-p1} is negative. In order to prove \eqref{local:eq-high-low-p1}, we first note that the minimum in \eqref{local:eq-high-low-p1} is given by $\min(1,(KM^{-\theta})^{1+\alpha})$. Due to our assumption \eqref{local:eq-high-low-c1a}, the exponent of the first $M$-factor in \eqref{local:eq-high-low-p1} is positive. Furthermore, due to the first condition in \eqref{local:eq-high-low-c1b}, the combined exponent of the first $M$-factor and the $M$-factor in the minimum in \eqref{local:eq-high-low-p1} is negative. From this, it follows that the sum in \eqref{local:eq-high-low-p1} is comparable to its value at $M=K^{1/\theta}$, which is given by 
\begin{equation*}
K^{\gamma+\nu-\eta-\alpha+\frac{2}{r}} \times K^{\frac{1-\theta}{\theta} (1-\alpha-\beta) + \frac{1-\theta}{\theta} \eps^\prime}.
\end{equation*}
Since $\eps^\prime$ is sufficiently small, this completes the proof of \eqref{local:eq-high-low-p1}. 
\end{proof}

\section{Covariant local well-posedness}\label{section:covariant_lwp}

In this section, we prove a covariant local well-posedness result for \eqref{local:eq-SYM-coordinates-eps}. It will be a more detailed version of Theorem~\ref{main:thm-covariant-local}, which is one of the main results of this article. For expository purposes, we first make the following assumption on the time $T$ and the background connection one-form $B$.

\begin{assumption}\label{local:assumption-B}
First, we assume that $0<T\leq T_0$, where $T_0$ is as in \eqref{prelim:eq-T0}. Second, we assume that the background connection one-form $B$ is a solution of \eqref{local:eq-B} on the time interval $[0,T]$ and satisfies
\begin{equation}\label{local:eq-assumption-B-bound}
\lossfactor_T(F_B) \leq T^{-\lossdel}, 
\end{equation}
where $\lossdel$ is as in \eqref{prelim:eq-parameters} and $\lossfactor_T(F)$ is as in \eqref{local:eq-lossfactor}. 
\end{assumption}

The reasons behind these assumptions on $B$ are twofold: First, we will be able to show that they are satisfied by the solution of the Yang--Mills gradient flow \eqref{local:eq-B} on each time interval appearing in our global iteration scheme. Second, they are sufficient to obtain our Schauder estimates (from Section~\ref{sec:covariant_function_spaces}) and our estimates of the covariant stochastic objects (in Section~\ref{section:objects} below). In addition to Assumption \ref{local:assumption-B}, we also make the following assumption on the initial Higgs fields. 

\begin{assumption}\label{local:assumption-higgs-initial}
We assume that the initial Higgs fields $\shiggsdet_0$ and $\rhiggsdet_0$, which previously appeared in Section \ref{sec:higgsdet_estimates}, satisfy the estimates 
\begin{equation}\label{local:eq-assumption-higgs-bound}
\max\bigl( \bigl\| \shiggsdet_0 \bigr\|_{L_x^r}, \bigl\|\rhiggsdet_0\bigr\|_{\covHolder{B(0), T}^{-\kappa}} \bigr) \leq T^{-\frac{1}{\gamma_0}},
\end{equation}
where $\gamma_0$ is as in \eqref{e:constants_global_iteration}.
\end{assumption}

The reason for distributing \eqref{local:eq-assumption-B-bound} and \eqref{local:eq-assumption-higgs-bound} over two different assumptions, rather than collecting them in a single assumption, is that many of our estimates only involve the background connection one-form $B$, and not the initial Higgs fields; see Section \ref{section:objects}. Whenever Assumption \ref{local:assumption-higgs-initial} is made, we also make Assumption \ref{local:assumption-B}, and hence the covariant Hölder norm appearing in \eqref{local:eq-assumption-higgs-bound} is well-defined.

We also mention that Assumptions~\ref{local:assumption-B} and~\ref{local:assumption-higgs-initial} imply the smallness condition 
\eqref{eq:Sdet-higgsdet-main-prop-smallness-assumption} from Proposition~\ref{prop:Sdet-higgsdet-estimates}, i.e., the smallness condition needed for our earlier estimates of $\Sdet$ and $\higgsdet$. Indeed, from the definitions of $M_T(B,10)$ and $\lossfactor_T(F)$ in~\eqref{eq:gronwall-F-pert} and~\eqref{local:eq-lossfactor}, it directly follows that $M_T(B,10)\leq \lossfactor_T(F)$. By using Assumptions~\ref{local:assumption-B} and~\ref{local:assumption-higgs-initial}, the left-hand side of~\eqref{eq:Sdet-higgsdet-main-prop-smallness-assumption} can then be bounded by
\begin{equation*}
M_T(B, 10) \bigl( T^{\frac{1}{2} - \frac{1}{r}} \big\|\shiggsdet_0\big\|_{L_x^r} + T^{\frac{1-\kappa}{2}} \bigl\|\rhiggsdet_0\bigr\|_{\covHolder{B(0), T}^{-\kappa}}\bigr) 
\leq T^{-\lossdel} 
\bigl( T^{\frac{1}{2} - \frac{1}{r}} T^{-\frac{1}{\gamma_0}}+ T^{\frac{1-\kappa}{2}} T^{-\frac{1}{\gamma_0}}\bigr). 
\end{equation*}
Since $\gamma_0-2$ is much larger than $\lossdel$ and $\kappa$, and since $T_0$ from Assumption \ref{local:assumption-B} is sufficiently small, this implies~\eqref{eq:Sdet-higgsdet-main-prop-smallness-assumption}.

Equipped with Assumptions~\ref{local:assumption-B} and~\ref{local:assumption-higgs-initial}, we are now ready to state the main result of this section. 

\begin{proposition}[Covariant local well-posedness]\label{local:prop-main}
Let Assumptions~\ref{local:assumption-B} and~\ref{local:assumption-higgs-initial} be satisfied and let the covariant enhanced data sets (see Definitions~\ref{local:def-enhanced-data-set} and~\ref{local:def-data-set-norm} below) satisfy 
\begin{equation}
\label{e:enhanced_data_assumption}
\sup_{\eps \in \sdyadic} \tripn{ \Xi_{\Beps} } \leq T^{-\lambda_1} \qquad \text{and} \qquad \lim_{\eps,\eps^\prime \rightarrow 0} \tripn{ \Xi_{\Beps}; \Xi_{\Beps^\prime} } =0. 
\end{equation}
Then, for all $\eps \in \sdyadic$, the solutions $(X_{\eps,\star}, \phix_{\eps,\star}, Y_{\eps,\star}, \phiy_{\eps,\star})$ to the paracontrolled system \eqref{local:eq-Xwhite-K-final}--\eqref{local:eq-y-final-new-last} exist and, as $\eps \to 0$, they converge to a limit $(X_{\star}, \phix_{\star}, Y_{\star}, \phiy_{\star})$ in 
\begin{equation*}
\bigl(\sol^{1-\reg_1,0} \times \sol^{1-\reg_1,0} \times \sol^{1+\reg_2,-\lambda} \times \sol^{1+\reg_2,-\lambda}\bigr)([0,T]). 
\end{equation*}
Finally, the remainders 
$R \coloneq \quadratic[B][][][] + \phiquadratic[B][][][] + X + Y$ and  $\phir \coloneq \mixedquadratic[B][][][] + \phix + \phiy$ satisfy the estimates
\begin{equation}\label{local:eq-main-R-estimate}
\| (R,\phir) \|_{L_{t,x}^\infty([0,T])} \lesssim T^{\frac{1-\lambda}{2}-\frac{1}{\gamma_0}-\lambda_2}
\qquad \text{and} \qquad 
\| (R_\star,\phir_\star) \|_{\sol^{1-\lambda,0}([0,T])} \lesssim T^{-\frac{1}{\gamma_0}-\lambda_2}. 
\end{equation}
\end{proposition}

\begin{remark}
From the convergence of the covariant enhanced data set $\Xi_{\Beps}$ and the convergence of $(X_{\eps,\star}, \phix_{\eps,\star}, Y_{\eps,\star}, \phiy_{\eps,\star})$ in Proposition \ref{local:prop-main}, it directly follows that the Littlewood--Paley blocks $((A_\eps-B)_\star,\higgs_{\eps,\star})$ converge in $\sol^{-\kappa}([0,T])$. This is the gauge covariant convergence result alluded to in Remark \ref{main:rem-Cminuskappa}.
\end{remark}

As in Section~\ref{section:covariant-Ansatz}, all covariant stochastic objects used in this section will depend on the same connection one-form $B$ (from Assumption~\ref{local:assumption-B}). For this reason, we follow Notation~\ref{local:notation-omitting-subscript} and simply write $\linear[\oeps][][][]$, $\philinear[\oeps][][][]$, $\quadratic[\oeps][][][]$, etc.\ for our covariant stochastic objects. 

Previously, we set up the covariant ansatz (Section~\ref{section:covariant-Ansatz}) and a covariant functional framework (Section~\ref{sec:covariant_function_spaces}). In order to complete the proof of Proposition~\ref{local:prop-main}, we will now proceed as follows: In Section~\ref{section:covariant-data-set}, we introduce the covariant enhanced data set and its norm, which will be estimated in Section~\ref{section:objects}. In Section~\ref{section:sdet-higgsdet-sol}, we control the deterministic terms $\Sdet$ and $\higgsdet$ in our solution spaces. In Section~\ref{section:lwp-nonlinear-estimate}, we estimate most of the terms in \eqref{local:eq-Xwhite-K-final}--\eqref{local:eq-y-final-new-last}, except for the terms that rely on the paracontrolled structure of $X_{\eps,\star}$ and $\phix_{\eps,\star}$. In Section~\ref{section:local-paracontrolled}, we then uncover the paracontrolled structure of $X_{\eps,\star}$ and $\phix_{\eps,\star}$, and estimate the remaining terms. Finally, in Section~\ref{section:local-proof}, we then finish up the proofs of Proposition~\ref{local:prop-main} and Theorem~\ref{main:thm-covariant-local}.

\subsection{Enhanced data set}\label{section:covariant-data-set}
In this subsection, we discuss the bounds on the covariant stochastic objects that are needed to solve \eqref{local:eq-Xwhite-K-final}--\eqref{local:eq-y-final-new-last}. In order to capture these bounds in a single proposition, we will introduce the covariant enhanced data set and its size. Our definitions will be similar to those in \cite[Section 5]{BC23}, except that all stochastic objects will be covariant and that we account for the Higgs field. 

In order to unify the treatment of the connection one-form and Higgs field, we introduce the following notation, which is in the spirit of Notation~\ref{local:notation-unified}. To this end, recall that $\Y=\frkg^2\oplus \V$. 

\begin{notation}\label{local:notation-Y-valued-objects}
We define the $\Y$-valued space-time white noise and covariant stochastic objects by 
\begin{equation}\label{local:eq-graylinear-boldxi}
\boldxi = \xi \oplus \zeta \qquad \text{and} \qquad 
\graylinear[\oeps][][][\ostar] = \linear[\oeps][][][\ostar] \oplus \philinear[\oeps][][][\ostar]. 
\end{equation}
In addition, we define the covariance $\Cov_{\Y}\in \Y \otimes \Y$ by $\Cov_{\Y}=\delta^{mn}\, \ONBY_m \otimes \ONBY_n$, where $(\ONBY_m)_m$ is an orthonormal basis of $\Y$. 
\end{notation}

In order to state the next definition, we recall the linear maps $\Cov^{\frkg}_{\frkg}$, $\Cov^{\frkg}_{\V}$, and $\Cov^{\V}_{\frkg}$ previously defined in Definition~\ref{prelim:def-Casimir-induced}. We also recall that the renormalisation constants $\sigma_{\rho,\eps}^2$ and $\Crho$ were previously defined in Definition \ref{main:def-renormalization-constant}.  

\begin{definition}[Covariant enhanced data set]\label{local:def-enhanced-data-set}
Let $0<T\leq \frac{1}{128}$, let $B\colon [0,T] \times \T^2 \rightarrow \frkg^2$ be the Yang--Mills gradient flow from \eqref{local:eq-B}, and let $\eps\in \sdyadic$. Then, we define the covariant enhanced data set as 
\begin{equation}\label{local:eq-data-set-combined}
\Xi_{\Beps} \eqdef \bigl( 
\Xi_{\Beps}^{\ans},
\Xi_{\Beps}^{\der},
\Xi_{\Beps}^{\pow},
\Xi_{\Beps}^{\para}\bigr). 
\end{equation}
The superscripts stand for ``ansatz'', ``derivative nonlinearity'', 
``power-type nonlinearity'', 
and ``paracontrolled'', respectively, and the corresponding components of the covariant data set are defined as follows: First, we define 
\begin{align}
\Xi_{\Beps}^{\ans}
&\eqdef \Big( \linear[\oeps][][][\ostar], \philinear[\oeps][][][\ostar], 
\quadratic[\oeps][][][\ostar],
\phiquadratic[\oeps][][][\ostar], 
\mixedquadratic[\oeps][][][\ostar]\Big),
\label{local:eq-data-set-ans}
\end{align}
where the stochastic objects are as in \eqref{local:eq-linear-decomposition}, \eqref{local:eq-philinear-N}, and \eqref{local:eq-quadratic-N}. Second, we define 
\begin{align}
\Xi_{\Beps}^{\der}
\eqdef
\Big(& 
 \Lb\bigl( \quadratic[\oeps][][][\ostar],\covd \linear[\oeps][][][\ostar] \bigr)
+ ( \sigma_{\rho,\eps}^2+\Crho) \Cov_{\frkg}^{\frkg} (\linear[\oeps][][][]),
\qquad \Lb\bigl( \phiquadratic[\oeps][][][\ostar],\covd \linear[\oeps][][][\ostar]\bigr), 
\label{local:eq-data-set-der-1} \\
& 
\bil \bigl( \mixedquadratic[\oeps][][][\ostar] \parasim \covd\philinear[\oeps][][][\ostar]\bigr) 
+\tfrac{1}{2} (\sigma_{\rho,\eps}^2 + \Crho) \Cov^{\frkg}_{\V}( \linear[\oeps][][][]), 
\qquad 
\quadratic[\oeps][\ell][][\ostar] \parasim \covd_\ell \philinear[\oeps][][][\ostar], 
\label{local:eq-data-set-der-2} \\
& \phiquadratic[\oeps][\ell][][\ostar] \parasim \covd_\ell \philinear[\oeps][][][\ostar]
+\tfrac{1}{2} (\sigma_{\rho,\eps}^2 + \Crho) \Cov^{\V}_{\frkg}( \philinear[\oeps][][][]),
\qquad \bigl( \covd_\ell \linear[\oeps][\ell][][\ostar] \bigr) \parasim    \mixedquadratic [\oeps][][][\ostar] \Big),
\label{local:eq-data-set-der-3}
\end{align}
where $\bil$ and $\Lb$ are as in \eqref{intro:eq-bil} and \eqref{local:eq-Lie-bracket-term}. Third, we define\footnote{The conditions on the frequency scales in \eqref{local:eq-data-set-para} are exactly as in our high$\times$high-paraproduct; see Definition~\ref{local:def-paraproduct}. However, to make clear that \eqref{local:eq-data-set-para} involves the tensor product in $\Y$, we refrain from using this notation here.}
\begin{align}\label{local:eq-data-set-para}
\Xi_{\Beps}^{\para}
\eqdef \Big(& \sum_{\substack{K,L \colon \\ K \sim_\theta L}}  \covd_B^k \Duh_{\Y,B,K}\bigl( \graylinear[\oeps][][][]\bigr)
\otimes \covd_B^\ell \graylinear[\oeps][][][L] 
- \frac{1}{4}   \delta^{k\ell} \bigl( \sigma_{\rho,\eps}^2 + \Crho \bigr) \Cov_{\Y}  \Big)_{k,\ell\in [2]}.
\end{align}
Finally, we define 
\begin{align}
\hphantom{}\hspace{-1.5ex}\Xi_{\Beps}^{\pow}
\eqdef \Big(&
  \graylinear[\oeps][][][] \otimes \graylinear[\oeps][][][] 
-\sigma_{\rho,\eps}^2 \Cov_{\Y}, 
\label{local:eq-data-set-pow-one}\\  
&   \graylinear[\oeps][][][] \otimes \graylinear[\oeps][][][] \otimes \graylinear[\oeps][][][]
- \sigma_{\rho,\eps}^2 \bigl(\Cov_{\Y}\otimes \graylinear[\oeps][][][] + \ONBY^m \otimes \graylinear[\oeps][][][] \otimes \ONBY_m + \graylinear[\oeps][][][] \otimes \Cov_{\Y} \bigr)  \Big). 
\label{local:eq-data-set-pow-two}
\end{align}
Here, $(\ONBY_m)_m$ is an orthonormal basis of $\Y$. 
\end{definition}

In addition to the covariant enhanced data set itself, we also define the following ``norm''.

\begin{definition}[Size of covariant enhanced data set]\label{local:def-data-set-norm}
Let $0<T\leq \frac{1}{128}$, let $B\colon [0,T] \times \T^2 \rightarrow \frkg^2$ be the Yang--Mills gradient flow from \eqref{local:eq-B}, and let $\eps\in \sdyadic$. Then, we define
\begin{align}
\tripn{ \Xi_{\Beps} } 
&\eqdef \big\|   \graylinear[\oeps][][][\ostar] 
\big\|_{\modsol^{-\kappa}} 
+ \big\| \quadratic[\oeps][][][\ostar] \big\|_{\sol^{1-2\kappa}}^{\frac{1}{2}} 
+ \big\| \phiquadratic[\oeps][][][\ostar] \big\|_{\sol^{1-2\kappa}}^{\frac{1}{2}}  
+ \big\| \mixedquadratic[\oeps][][][\ostar] \big\|_{\sol^{1-2\kappa}}^{\frac{1}{2}}  \notag \\
&+ \big\| \mbox{Terms in \eqref{local:eq-data-set-der-1}--\eqref{local:eq-data-set-der-3}} \big\|_{\Cs_t^{-\theta} L_x^r}^{\frac{1}{3}} 
+ \big\| \eqref{local:eq-data-set-para} \big\|_{\Cs_t^{-\theta}L_x^r}^{\frac{1}{2}} \label{local:eq-data-set-norm} \\ 
&+ \big\| \eqref{local:eq-data-set-pow-one} \big\|_{\Cs_t^{-\kappa} L_x^r}^{\frac{1}{2}}
+ \big\| \eqref{local:eq-data-set-pow-two} \big\|_{\Cs_t^{-\frac{3\kappa}{2}} L_x^r}^{\frac{1}{3}}. 
 \notag 
\end{align}  
Furthermore, if $\eps^\prime\in \sdyadic$, we define the distance 
\begin{equation*}
\tripn{ \Xi_{\Beps}; \Xi_{\Beps^\prime} } 
\end{equation*}
as in \eqref{local:eq-data-set-norm}, but with all stochastic objects replaced by the corresponding difference of stochastic objects. For example, the first summand in \eqref{local:eq-data-set-norm} is replaced with $
\|  \graylinear[\oeps][][][\ostar] -
\graylinear[\oepsprime][][][\ostar]\|_{\modsol^{-\kappa}}$.
\end{definition}
The estimates of both $\tripn{ \Xi_{\Beps} }$
and $\tripn{\Xi_{\Beps}; \Xi_{\Beps^\prime}} $ are the subject of Section~\ref{section:objects}. For now, we simply state our estimates in terms of  $\tripn{\Xi_{\Beps}}$
and $\tripn{ \Xi_{\Beps}; \Xi_{\Beps^\prime} } $.

\subsection{\protect{Deterministic estimates of \TitleEquation{\Sdet}{S} and \TitleEquation{\higgsdet}{phi}}}\label{section:sdet-higgsdet-sol}

In Proposition~\ref{prop:Sdet-higgsdet-estimates} above, we estimated $\Sdet$ and $\higgsdet$, as well as their covariant derivatives, in weighted $L_{t,x}^\infty$-spaces. In the following, we will estimate them in our weighted solution spaces from Definition~\ref{local:def-solution-space}. This is needed in our covariant local theory, but will not enter directly into our global iteration argument. To this end, we recall that the covariant Littlewood--Paley blocks of $\Sdet$ and $\higgsdet$ were defined in \eqref{local:eq-Sdet-N} and~\eqref{local:eq-higgsdet-N}, respectively. 

\begin{proposition}[\protect{Bounds for $\higgsdet_\star$ and $\Sdet_\star$ in solution spaces}]\label{local:prop-higgsdet-sdet}
Let Assumptions~\ref{local:assumption-B} and \ref{local:assumption-higgs-initial} be satisfied, let $0\leq \alpha \leq 2$, and let $\eta\leq 0$. Under the condition $\alpha+\eta \leq -\kappa$, it then holds that  
\begin{align}\label{local:eq-higgsdet-sol}
\| \higgsdet_{\star} \|_{\sol^{\alpha,\eta}} &\lesssim T^{-\lossdel} T^{-\frac{\alpha+\eta+\kappa}{2}} \big( \higgsdetcombined \big).
\end{align}
Furthermore, under the condition $\alpha+\eta \leq 1-2\kappa$, it holds that 
\begin{align}
\| \Sdet_{\star} \|_{\sol^{\alpha,\eta}} 
&\lesssim  T^{-3\lossdel}
T^{\frac{1-\alpha-\eta-2\kappa}{2}} \big( \higgsdetcombined\big)^2. \label{local:eq-sdet-sol}
\end{align}
\end{proposition}

We note that, due to our upper bound on $T$, the right-hand side of \eqref{local:eq-sdet-sol} is smaller than the right-hand side of \eqref{local:eq-higgsdet-sol}. 

\begin{proof}[Proof of Proposition~\ref{local:prop-higgsdet-sdet}]
The proof does not rely on a contraction-mapping argument. Instead, we just use our earlier estimates of $\Sdet$ and $\higgsdet$ from Proposition~\ref{prop:Sdet-higgsdet-estimates}, and insert them into our Schauder estimates. In particular, we will often use that
\begin{align*}
T^{\frac{1}{2}-\lossdel} \bigl( \| \Sdet \|_{L_{t,x}^\infty} + \| t^{\frac{1}{2}} \covd \Sdet \|_{L_{t,x}^\infty} \bigr)
\lesssim T^{1-2\lossdel-\kappa} \big( \higgsdetcombined \big)^2 
\lesssim T^{1-2\lossdel-\kappa} T^{-\frac{2}{\gamma_0}} \lesssim T^{\lossdel}. 
\end{align*}
In the last estimate, we used that $\gamma_0-2$ is much larger than $\lossdel$ and $\kappa$. For notational simplicity, we assume that $\higgsdetcombined \geq 1$; otherwise, one can follow simplified versions of the arguments below.

\step{Step 1: Estimates of \TitleEquation{\higgsdet_{\star}}{phi_*}.} Using the trivial bound $t\leq T$, that $\higgsdet_{N}$ is supported on times  $t\gtrsim N^{-2}$, and Lemma~\ref{local:lem-sol-different-regularities-weights}, we first estimate 
\begin{equation*}
\big\| \higgsdet_{\star} \big\|_{\sol^{\alpha,\eta}}
\lesssim T^{-\frac{\alpha+\eta+\kappa}{2}} 
\big\| \higgsdet_{\star} \big\|_{\sol^{\alpha,-\alpha-\kappa}} 
\lesssim T^{-\frac{\alpha+\eta+\kappa}{2}} 
\big\| \higgsdet_{\star} \big\|_{\sol^{2,-2-\kappa}}. 
\end{equation*}
Due to this, it suffices to estimate the $\sol^{2,-2-\kappa}$-norm, i.e.\  it suffices to prove \eqref{local:eq-higgsdet-sol} for $(\alpha,\eta)=(2,-2-\kappa)$. We recall from 
\eqref{local:eq-higgsdet-N} that 
\begin{equation}\label{local:eq-higgsdet-sol-p1}
\higgsdet_N=\higgsdet_{\lin,N} + \Duh_{B,N} \bigl( 2 \Sdet_j \covd^j \higgsdet  + \Sdet^j \Sdet_j \higgsdet - |\higgsdet|^2 \higgsdet  \bigr).
\end{equation}
Using Lemma~\ref{local:lem-linear-evolution-sol}, the linear evolution in \eqref{local:eq-higgsdet-sol-p1} can be bounded by 
\begin{equation*}
\bigl\| \higgsdet_{\lin,\star}\bigr\|_{\sol^{2,-2-\kappa}}
\lesssim T^{-\lossdel} \Bigl( T^{\frac{\kappa}{2}-\frac{1}{r}} \| \higgsdet_0^+ \|_{L_x^r} + T^{\frac{\kappa}{2}-\frac{\kappa}{2}} \bigl\|\rhiggsdet_0\bigr\|_{\covHolder{B(0), T}^{-\kappa}} \Bigr) 
\lesssim T^{-\lossdel} \big( \higgsdetcombined \big). 
\end{equation*}
In the last inequality, we used that $1/r$ is much smaller than $\kappa$. For $(\alpha,\eta)=(2,-2-\kappa)$, this agrees with the right-hand side of \eqref{local:eq-higgsdet-sol}. The second term in~\eqref{local:eq-higgsdet-sol-p1} can be estimated by
\begin{equation}\label{local:eq-higgsdet-sol-p2}
\begin{aligned}
\bigl\| \Duh_{B,N} &\bigl(  \Sdet_j \covd^j \higgsdet \bigr) \bigr\|_{\sol^{2,-2-\kappa}_N} 
\lesssim T^{\frac{1}{2}} \bigl\| \Duh_{B,N} \bigl(  \Sdet_j \covd^j \higgsdet \bigr) \bigr\|_{\sol^{2,-1-\kappa}_N} \\
&\lesssim T^{\frac{1}{2}-\lossdel} \| t^{\frac{1}{2}+\frac{\kappa}{2}} \Sdet_j \covd^j \higgsdet  \|_{L_{t,x}^\infty} 
\lesssim T^{\frac{1}{2}-\lossdel} \| \Sdet\|_{L_{t,x}^\infty}  \| t^{\frac{1}{2}+\frac{\kappa}{2}}  \covd \higgsdet  \|_{L_{t,x}^\infty}  \\
&\lesssim T^{\lossdel} \times T^{-\lossdel} \big( \higgsdetcombined \big), 
\end{aligned}
\end{equation}
which is acceptable. The third term in \eqref{local:eq-higgsdet-sol-p1} can be bounded by 
\begin{equation}\label{local:eq-higgsdet-sol-p3}
\begin{aligned}
\bigl\| \Duh_{B,N} &\bigl(\Sdet^j \Sdet_j \higgsdet\bigr) \bigr\|_{\sol^{2,-2-\kappa}_N}
\lesssim T  \bigl\| \Duh_{B,N} \bigl(\Sdet^j \Sdet_j \higgsdet\bigr) \bigr\|_{\sol^{2,-\kappa}_N} 
\lesssim T^{1-\lossdel} \| t^{\frac{\kappa}{2}} \Sdet^j \Sdet_j \higgsdet \|_{L_{t,x}^\infty} \\ 
&\lesssim  T^{1-\lossdel} \| \Sdet \|_{L_{t,x}^\infty}^2 \| t^{\frac{\kappa}{2}} \higgsdet \|_{L_{t,x}^\infty} 
\lesssim T^{2\lossdel} \times T^{-\lossdel} \big( \higgsdetcombined \big),  
\end{aligned}
\end{equation}
which is acceptable. Finally, the fourth term in \eqref{local:eq-higgsdet-sol-p1} can be bounded by 
\begin{equation}\label{local:eq-higgsdet-sol-p4}
\begin{aligned}
\bigl\| \Duh_{B,N} &\bigl(|\higgsdet|^2 \higgsdet \bigr) \bigr\|_{\sol^{2,-2-\kappa}_N}
\lesssim T^{1-\kappa} \bigl\| \Duh_{B,N} \bigl(|\higgsdet|^2 \higgsdet \bigr) \bigr\|_{\sol^{2,-3\kappa}_N}  \\
&\lesssim  T^{1-\lossdel-\kappa} \big\| t^{\frac{3\kappa}{2}} |\higgsdet|^2 \higgsdet \big\|_{L_{t,x}^\infty}  
\lesssim T^{1-4\lossdel-\kappa} \big( \higgsdetcombined \big)^3 \\
&\lesssim \big( \higgsdetcombined \big)^{3-(1-4\lossdel-\kappa)\gamma_0}. 
\end{aligned}
\end{equation}
Since $\gamma_0-2$ is much larger than $\lossdel$ and $\kappa$, this is acceptable. 

\step{Step 2: Estimates of \TitleEquation{\Sdet_\star}{S_*}.} As in our estimates of $\higgsdet_\star$, we first use the trivial bound $t\leq T$ and Lemma \ref{local:lem-sol-different-regularities-weights}. Since our solution space $\sol^{\beta,\nu}$ is only defined for $\nu\leq 0$, we need to distinguish the cases $\alpha\geq 1-2\kappa$ and $\alpha \leq 1-2\kappa$. In the first case, we estimate
\begin{equation*}
\| \Sdet_\star \|_{\sol^{\alpha,\eta}}
\lesssim T^{\frac{1-2\kappa-\alpha-\eta}{2}} \| \Sdet_\star \|_{\sol^{\alpha,1-2\kappa-\alpha}} 
\lesssim T^{\frac{1-2\kappa-\alpha-\eta}{2}} \| \Sdet_\star \|_{\sol^{2,-1-2\kappa}}. 
\end{equation*}
In the second case, we estimate
\begin{equation*}
\| \Sdet_\star \|_{\sol^{\alpha,\eta}}
\lesssim T^{-\frac{\eta}{2}} \| \Sdet_\star \|_{\sol^{\alpha,0}} 
\lesssim T^{\frac{1-2\kappa-\alpha}{2}} T^{-\frac{\eta}{2}} \| \Sdet_\star \|_{\sol^{2,-1-2\kappa}}. 
\end{equation*}
Thus, it suffices to treat the case $(\alpha,\eta)=(2,-1-2\kappa)$. We recall from \eqref{local:eq-Sdet-N} that 
\begin{equation}\label{local:eq-Sdet-sol-p1}
\Sdet_N = \Duh_{B,N}^F \Big( [ \Sdet_j, 2 \covd^j \Sdet - \covd \Sdet^j] + [ \Sdet_j , [\Sdet^j , \Sdet ]] 
+ \bil\bigl( \higgsdet \otimes \covd \higgsdet\bigr) + \bil \bigl( \higgsdet \otimes \Sdet \higgsdet \bigr)\Big). 
\end{equation}
The first term in \eqref{local:eq-Sdet-sol-p1} can be bounded by 
\begin{equation}\label{local:eq-Sdet-sol-p2}
\begin{aligned}
\bigl\| \Duh_{B,N}^F &\bigl( [ \Sdet_j, 2 \covd^j \Sdet - \covd \Sdet^j]\bigr)\bigr\|_{\sol^{2,-1-2\kappa}_N} \\ 
&\lesssim T^{\kappa-\lossdel} \big\| t^{\frac{1}{2}} [ \Sdet_j, 2 \covd^j \Sdet - \covd \Sdet^j]\big\|_{L_{t,x}^\infty} 
\lesssim T^{\kappa-\lossdel} \| \Sdet \|_{L_{t,x}^\infty} \| t^{\frac{1}{2}} \covd \Sdet \|_{L_{t,x}^\infty} \\
&\lesssim T^{1-3\lossdel-\kappa} \big( \higgsdetcombined \big)^{4} 
\lesssim \big( \higgsdetcombined \big)^{4-(1-3\lossdel-\kappa)\gamma_0}.  
\end{aligned}
\end{equation}
 The second term in \eqref{local:eq-Sdet-sol-p1} can be bounded by 
\begin{equation}\label{local:eq-Sdet-sol-p3}
\begin{aligned}
\bigl\| \Duh_{B,N}^F &\bigl( [ \Sdet_j , [\Sdet^j , \Sdet ]] \bigr) 
\bigr\|_{\sol^{2,-1-2\kappa}_N} 
\lesssim T^{\frac{1}{2}+\kappa} \bigl\| \Duh_{B,N}^F \bigl( [ \Sdet_j , [\Sdet^j , \Sdet ]] \bigr) 
\bigr\|_{\sol^{2,0}_N} \\ 
&\lesssim  T^{\frac{1}{2}+\kappa-\lossdel} \big\|  [ \Sdet_j , [\Sdet^j , \Sdet ]] \big\|_{L_{t,x}^\infty}
\lesssim T^{\frac{1}{2}+\kappa-\lossdel} \bigl(  
T^{\frac{1}{2}-\lossdel-\kappa} 
\bigl( \higgsdetcombined \bigr)^2 \bigr)^3 \\
&\lesssim \bigl( \higgsdetcombined \bigr)^{6-(2-4\lossdel-2\kappa)\gamma_0}.
\end{aligned}
\end{equation} 
 The third term in \eqref{local:eq-Sdet-sol-p1} can be bounded by 
\begin{equation}\label{local:eq-Sdet-sol-p4}
\begin{aligned}
\bigl\| \Duh_{B,N}^F &\bigl( \bil\bigl( \higgsdet \otimes \covd \higgsdet\bigr) \bigr) 
\bigr\|_{\sol^{2,-1-2\kappa}_N} 
\lesssim  T^{-\lossdel}
\big\| t^{\frac{1}{2}+\kappa} \bil\bigl( \higgsdet \otimes \covd \higgsdet\bigr)  \big\|_{L_{t,x}^\infty} \\  
&\lesssim T^{-\lossdel} \| t^{\frac{\kappa}{2}} \higgsdet\|_{L_{t,x}^\infty}
\| t^{\frac{1}{2}+\frac{\kappa}{2}} \covd \higgsdet \|_{L_{t,x}^\infty}
\lesssim T^{-3\lossdel} \bigl( \higgsdetcombined \bigr)^2.  
\end{aligned}
\end{equation}
Finally, the fourth term in \eqref{local:eq-Sdet-sol-p1} can be bounded by 
\begin{equation}\label{local:eq-Sdet-sol-p5}
\begin{aligned}
\bigl\| \Duh_{B,N}^F &\bigl( \bil \bigl( \higgsdet \otimes \Sdet \higgsdet \bigr) \bigr) 
\bigr\|_{\sol^{2,-1-2\kappa}_N} 
\lesssim  T^{\frac{1}{2}} \bigl\| \Duh_{B,N}^F \bigl( \bil \bigl( \higgsdet \otimes \Sdet \higgsdet \bigr) \bigr) 
\bigr\|_{\sol^{2,-2\kappa}_N} \\
&\lesssim  T^{\frac{1}{2}-\lossdel} \big\|t^{\kappa} \bil \bigl( \higgsdet \otimes \Sdet \higgsdet \bigr)  
\big\|_{L_{t,x}^\infty} 
\lesssim T^{\frac{1}{2}-\lossdel} \| \Sdet \|_{L_{t,x}^\infty} 
\| t^{\frac{\kappa}{2}} \higgsdet \|_{L_{t,x}^\infty}^2\\
&\lesssim  T^{\frac{1}{2}-\lossdel}
\times T^{\frac{1}{2}-\lossdel-\kappa} \bigl( \higgsdetcombined \bigr)^2 \times T^{-2\lossdel} \bigl( \higgsdetcombined \bigr)^2 \\ 
&\lesssim \bigl( \higgsdetcombined \bigr)^{4-(1-4\lossdel-\kappa)\gamma_0}. 
\end{aligned}
\end{equation}
Since $\gamma_0-2$ is much larger than $\lossdel$ and $\kappa$, all contributions from \eqref{local:eq-Sdet-sol-p2}--\eqref{local:eq-Sdet-sol-p5} are acceptable. 
\end{proof}

\subsection{Nonlinear estimates}\label{section:lwp-nonlinear-estimate}

In this subsection, we estimate all terms in \eqref{local:eq-Xwhite-K-final}--\eqref{local:eq-y-final-new-last}, except the terms that can be controlled using only $\tripn{ \Xi_{\Beps} }$ and the terms whose estimates rely on the paracontrolled structure of $X_{\eps,\star}$ or $\phix_{\eps,\star}$, i.e., \eqref{local:eq-Y-final-new-1}, \eqref{local:eq-Y-final-new-1p}, \eqref{local:eq-Y-final-new-2}, \eqref{local:eq-y-final-new-1}, and \eqref{local:eq-y-final-new-2}. In most cases, the estimates follow directly from our general estimates in Section~\ref{section:local-paraproduct-estimates}, and the only difficulty lies in choosing the correct regularity parameters. 

In the derivation of our ansatz (Section~\ref{section:covariant-Ansatz}), we mostly viewed $\Sdet$ and $\higgsdet$ as parts of $\Wrem_\eps$ and~$\phiw_\eps$, respectively. Since both $\Sdet$ and $\higgsdet$ have high regularity (at least in weighted spaces), this made sense. In this subsection, however, the treatment of $\Sdet$ and $\higgsdet$ requires additional care. The reason is that the size of $\Sdet$ and $\higgsdet$ is relatively large, and we need to show that despite the size of $\Sdet$ and $\higgsdet$, we can still obtain the estimates stated at the end of Proposition~\ref{local:prop-main}. To this end, we introduce the variables
\begin{align}
R_\eps &\coloneq W_\eps - \Sdet = \quadratic[\oeps][][][] + \phiquadratic[\oeps][][][] + X_\eps + Y_\eps, \label{local:eq-Reps} \\ 
\phir_\eps &\coloneq \phiw_\eps - \higgsdet = \mixedquadratic[\oeps][][][] + \phix_\eps + \phiy_\eps.  \label{local:eq-phireps}
\end{align}
Throughout this subsection, we rely on the following assumption.

\begin{assumption}[Setup]\label{local:assumption-nonlinear-setup}
We assume that 
\begin{alignat*}{3}
X_{\eps,\star}&\in \sol^{1-\reg_1,-\lambda}([0,T]\times \T^2 \rightarrow \frkg^2), \qquad & \qquad 
Y_{\eps,\star}&\in \sol^{1+\reg_2,-\lambda}([0,T]\times \T^2 \rightarrow \frkg^2), \\ 
\phix_{\eps,\star}&\in \sol^{1-\reg_1,-\lambda}([0,T]\times \T^2 \rightarrow \V), \qquad &\text{and} \qquad 
\phiy_{\eps,\star}&\in \sol^{1+\reg_2,-\lambda}([0,T]\times \T^2 \rightarrow \V).  
\end{alignat*}
Furthermore, we assume that each covariant Littlewood--Paley block $X_{\eps,N}$, $\phix_{\eps,N}$, $Y_{\eps,N}$, or $\phiy_{\eps,N}$ vanishes for times $0\leq t\leq \frac{1}{16}N^{-2}$. 
That is, in this subsection, we treat $X_{\eps,\star}$, $Y_{\eps,\star}$, $\phix_{\eps,\star}$, and $\phiy_{\eps,\star}$ as arbitrary elements\footnote{Since the paracontrolled estimates will only appear in Section~\ref{section:local-paracontrolled}, $X_{\eps,\star}$ and $\phix_{\eps,\star}$ will only appear through $\Wrem_{\eps,\star}$ and~$\phiw_{\eps,\star}$.} of our solution spaces satisfying the support conditions resulting from our time-localised Duhamel integral operators in \eqref{local:Duh-LWP-block}.
\end{assumption}

Furthermore, throughout this subsection and the next, we repeatedly apply estimates involving $\lossfactor_T(F)$, including the covariant Schauder estimates in Lemma~\ref{local:lem-covariant-Schauder} and Proposition~\ref{local:prop-weighted-covariant-Schauder}. In each such application, we use the bound $\lossfactor_T(F)\leq T^{-\lossdel}$ from Assumption~\ref{local:assumption-B}.

In the first lemma of this subsection, we control the terms appearing in the definitions of the paracontrolled terms 
$\Xwhite_\eps$, $\Xblack_\eps$, and $\phix_\eps$, i.e., the terms in \eqref{local:eq-Xwhite-K-final}, \eqref{local:eq-Xblack-K-final}, and \eqref{local:eq-phix-K-final}. 

\begin{lemma}\label{local:lem-nl-Xterms}
Let Assumptions~\ref{local:assumption-B},~\ref{local:assumption-higgs-initial}, and~\ref{local:assumption-nonlinear-setup} be satisfied. Then, it holds that 
\begin{equation}\label{local:eq-nl-Xterms}
\begin{aligned}
\big\|  \big[  \Wrem_{\eps,\star,j} &\parall \bigl( 2 \covd^j \linear[\oeps][][][\ostar]- \covd \linear[\oeps][j][][\ostar] \bigr) \big] \big\|_{\nl^{-1-\reg_1,0}} 
+ \big\| \bil\bigl( \phiw_{\eps,\star} \parall \covd \philinear[\oeps][][][\ostar]\bigr) \big\|_{\nl^{-1-\reg_1,0}} \\
&+  \big\| \Wrem_{\eps,\star}^j \parall \covd_j \philinear[\oeps][][][\ostar]\big\|_{\nl^{-1-\reg_1,0}} \\ 
&\lesssim  T^{-\lambda_1} \tripn{ \Xi_{\Beps} } 
\bigl(  \higgsdetcombined + T^{\frac{1}{2}} \| (R_{\eps,\star},\phir_{\eps,\star}) \|_{\sol^{1-\reg_1,-\lambda}}\bigr). 
\end{aligned}
\end{equation}
\end{lemma}

In \eqref{local:eq-nl-Xterms}, we use the convenient notation 
\begin{equation*}
\| (R_{\eps,\star},\phir_{\eps,\star}) \|_{\sol^{1-\reg_1,-\lambda}}
\eqdef \| R_{\eps,\star} \|_{\sol^{1-\reg_1,-\lambda}}
+ \| \phir_{\eps,\star} \|_{\sol^{1-\reg_1,-\lambda}},
\end{equation*}
which allows us to shorten several expressions. 

\begin{remark}
We want to use \eqref{local:eq-nl-Xterms} to explain one way of looking at the estimates in this subsection. In~\eqref{local:eq-nl-Xterms}, we encounter a loss of the form $T^{-\lambda_1}$. This loss is acceptable in our global iteration argument, since our parameter $\gamma_0$ has been chosen such that $\gamma_0-2$ is much larger than $\lambda_1$. Due to our parameter conditions \eqref{prelim:eq-parameters}, $\lambda_1$ is much larger than expressions involving other parameters appearing in the local theory, such as $\lossdel$, $\reg_2$, $\lambda$, or $\reg_2/\theta$. Hence, losses involving powers of $T^{-\lossdel}$, $T^{-\reg_2}$, $T^{-\lambda}$, or $T^{-\reg_2/\theta}$ can all easily be bounded by $T^{-\lambda_1}$. As a result of this, it is not too important to keep track of the precise powers of $T$, and one only needs to keep track of the part involving explicit constants (such as~$1/2$ or~$1$). 

As discussed in Remark~\ref{local:rem-solution-space-scaling}, the scaling exponent associated with the $\sol^{\alpha,\eta}$-norm is $\alpha+\eta$, and thus the scaling exponent associated with $\| (R_{\eps,\star},\phir_{\eps,\star}) \|_{\sol^{1-\reg_1,-\lambda}}$ is $1-\reg_1-\lambda \approx 1$. Also, the scaling exponents associated with $\| \higgsdet_0^+\|_{L_x^r}$ and with $ \bigl\|\rhiggsdet_0\bigr\|_{\covHolder{B(0), T}^{-\kappa}}$ is $-\frac{2}{r}\approx 0$ and $-\kappa \approx 0$, respectively. This is the reason why 
$\| \higgsdet_0^+\|_{L_x^r}$ and $ \bigl\|\rhiggsdet_0\bigr\|_{\covHolder{B(0), T}^{-\kappa}}$ come without a factor of~$T$, but $\| (R_{\eps,\star},\phir_{\eps,\star}) \|_{\sol^{1-\reg_1,-\lambda}}$ comes with a factor of $T^{\frac{1}{2}}$. 
\end{remark}

\begin{proof}
We only control the second summand in \eqref{local:eq-nl-Xterms}, since the first and third summands can be treated using similar arguments. The only difference is that, instead of using \eqref{local:eq-higgsdet-sol} as below, one has to use \eqref{local:eq-sdet-sol}. 
Using the $\liegroup$-equivariance of the bilinear map $\bil$ (see Lemma~\ref{prelim:lem-equivariant}) and the product rule (Lemma~\ref{prelim:lem-equivariant-product-rule}), we first write 
\begin{equation}\label{local:eq-nl-Xterms-p1}
\bil\bigl( \phiw_{\eps,\star} \parall \covd \philinear[\oeps][][][\ostar]\bigr)
= \covd \bil\bigl( \phiw_{\eps,\star} \parall  \philinear[\oeps][][][\ostar]\bigr) - \bil\bigl( \covd \phiw_{\eps,\star} \parall  \philinear[\oeps][][][\ostar]\bigr). 
\end{equation}
Using Lemma~\ref{local:lem-product-negative-regularity} with 
\begin{equation}\label{local:eq-nl-Xterms-pc1}
\bigl( \alpha,\beta,\gamma,\nu,\eta \bigr) 
\coloneq \bigl( -\kappa, \, 2\kappa, \, -1 -\reg_1 ,\, 0, \, -3\kappa \bigr), 
\end{equation}
the first summand in \eqref{local:eq-nl-Xterms-p1} can be bounded by
\begin{equation*}
\big\| \covd \bil\bigl( \phiw_{\eps,\star} \parall  \philinear[\oeps][][][\ostar]\bigr) \big\|_{\nl^{-1-\reg_1,0}}
\lesssim T^{-\lossdel} T^{\frac{\reg_1}{2}-\bigO(\kappa)} \big\| \philinear[\oeps][][][\ostar]\big\|_{\modsol^{-\kappa}} \| \phiw_{\eps,\star}\|_{\sol^{2\kappa,-3\kappa}}. 
\end{equation*} 
Using Lemma~\ref{local:lem-sol-different-regularities-weights} and Proposition~\ref{local:prop-higgsdet-sdet}, we can further estimate 
\begin{align*}
\| \phiw_{\eps,\star}\|_{\sol^{2\kappa,-3\kappa}}
&\leq \| \higgsdet_\star \|_{\sol^{2\kappa,-3\kappa}}
+ \| \phir_{\eps,\star} \|_{\sol^{2\kappa,-3\kappa}} \\ 
&\lesssim T^{-\lossdel}  \bigl( \higgsdetcombined \bigr) + T^{\frac{1-\lambda-\reg_1+\kappa}{2}} \| \phir_{\eps,\star} \|_{\sol^{1-\reg_1,-\lambda}}, 
\end{align*}
which then leads to an acceptable contribution to \eqref{local:eq-nl-Xterms}. In order to estimate the second summand in~\eqref{local:eq-nl-Xterms-p1}, we further split 
\begin{equation}\label{local:eq-nl-Xterms-p2}
\bil\bigl( \covd \phiw_{\eps,\star} \parall  \philinear[\oeps][][][\ostar]\bigr)
= 
\bil\bigl( \covd \higgsdet_{\star} \parall  \philinear[\oeps][][][\ostar]\bigr)
+
\bil\bigl( \covd \phir_{\eps,\star} \parall  \philinear[\oeps][][][\ostar]\bigr). 
\end{equation}
Using Lemma~\ref{local:lem-high-low} with 
\begin{equation}\label{local:eq-nl-Xterms-pc2}
\bigl( \alpha,\beta,\gamma,\nu,\eta \bigr) 
\coloneq \bigl(-\kappa, \,  1-3\kappa, \,  -1-\reg_1,\, 0, \, -1+2\kappa\bigr)
\end{equation}
and Proposition~\ref{local:prop-higgsdet-sdet}, we obtain that 
\begin{align*}
\big\| \bil\bigl( \covd \higgsdet_{\star} \parall  \philinear[\oeps][][][\ostar]\bigr) \big\|_{\nl^{-1-\reg_1,0}}
&\lesssim T^{-\lossdel} T^{\frac{\reg_1}{2}-\bigO(\kappa/\theta)} 
\| \philinear[\oeps][][][\ostar] \|_{\modsol^{-\kappa}} 
\| \higgsdet_{\star}\|_{\sol^{1-3\kappa,-1+2\kappa}} \\ 
&\lesssim T^{-\lossdel} T^{\frac{\reg_1}{2}-\bigO(\kappa/\theta)} \| \philinear[\oeps][][][\ostar] \|_{\modsol^{-\kappa}} 
\times T^{-\lossdel} \bigl( \higgsdetcombined \bigr). 
\end{align*}
Due to \eqref{prelim:eq-parameters}, $\kappa/\theta$ is much smaller than $\lambda_1$, and hence this contribution is acceptable. Using Lemma~\ref{local:lem-high-low} with 
\begin{equation*}
\bigl( \alpha,\beta,\gamma,\nu,\eta \bigr) 
\coloneq \bigl(-\kappa,\, 1-\reg_1,\, -1-\reg_1,\, 0, \, -\lambda \bigr), 
\end{equation*}
we also obtain that 
\begin{equation*}
\big\| \bil\bigl( \covd \phir_{\eps,\star} \parall  \philinear[\oeps][][][\ostar]\bigr) \big\|_{\nl^{-1-\reg_1,0}}
\lesssim T^{-\lossdel} T^{\frac{1}{2}-\bigO(\lambda)} 
\| \philinear[\oeps][][][\ostar] \|_{\modsol^{-\kappa}} 
\| \phir_{\eps,\star}\|_{\sol^{1-\reg_1,-\lambda}},
\end{equation*}
which is acceptable as well. 
\end{proof}

\begin{remark}
When verifying that the parameters in \eqref{local:eq-nl-Xterms-pc1} and \eqref{local:eq-nl-Xterms-pc2} satisfy the conditions in Lemma~\ref{local:lem-product-negative-regularity}
and Lemma~\ref{local:lem-high-low}, we strongly advise the reader to keep \eqref{prelim:eq-parameters} in mind. Due to~\eqref{prelim:eq-parameters}, most (but not all) of the parameters among $\kappa,\reg_1,\reg_2,\theta,\lambda$, and $r$ can then be ignored in each of the conditions in Lemma~\ref{local:lem-product-negative-regularity}
and Lemma~\ref{local:lem-high-low}. 
\end{remark}

\begin{remark}
In the proof of Lemma~\ref{local:lem-nl-Xterms}, we always wrote out the $T^{-\lossdel}$-factor explicitly, and never made it part of the $\bigO(\narrowdots)$-notation. This was done for expository reasons, since it allowed us to visually separate the loss from $\lossfactor_T(F)$ and the possible loss from the parameters appearing in 
Lemma~\ref{local:lem-product-negative-regularity},
Lemma~\ref{local:lem-high-low}, and 
Proposition~\ref{local:prop-higgsdet-sdet}. From now on, however, we will always make the parameter $\delta$ part of our $\bigO(\narrowdots)$-terms, since it slightly shortens our expressions.
\end{remark}

In the next lemma, we estimate high$\times$high-interactions between the linear stochastic objects and $\Sdet$, $\higgsdet$, $Y_\eps$, or $\phiy_\eps$, i.e., the terms in \eqref{local:eq-Y-final-new-3} and \eqref{local:eq-y-final-new-3}. 

\begin{lemma}\label{local:lem-nl-Y-high-high}
Let Assumptions~\ref{local:assumption-B},~\ref{local:assumption-higgs-initial}, and~\ref{local:assumption-nonlinear-setup} be satisfied. Then, it holds that 
\begin{equation}\label{local:eq-nl-Y-high-high}
\begin{aligned}
\big\| \Lb&\bigl(\Sdet_\star+Y_{\eps,\star},\covd \linear[\oeps][][][\ostar]\bigr) \big\|_{\nl^{-1+\reg_2,-\lambda}}
+ \big\| \bigl( \Sdet_{\star} + Y_{\eps,\star}\bigr)^j 
 \parasim \covd_j \philinear[\oeps][][][\ostar] \big\|_{\nl^{-1+\reg_2,-\lambda}} \\
 &+ \big\| \bil \bigl( (\higgsdet_\star+\phiy_{\eps,\star}) \parasim \covd \philinear[\oeps][][][\ostar]\bigr)  \big\|_{\nl^{-1+\reg_2,-\lambda}}
+ \big\|  \bigl( \covd_j \linear[\oeps][j][][\ostar] \bigr) \parasim    \bigl( \higgsdet_{\star} + \phiy_{\eps,\star} \bigr) \big\|_{\nl^{-1+\reg_2,-\lambda}} \\ 
 &\lesssim T^{-\lambda_1} 
 \tripn{ \Xi_{\Beps} }  \bigl(  
 \higgsdetcombined
 + T^{\frac{1}{2}} \| (Y_{\eps,\star},\phiy_{\eps,\star}) \|_{\sol^{1+\reg_2,-\lambda}} \bigr). 
\end{aligned}
\end{equation}
\end{lemma}

\begin{proof}
We only estimate the third summand in \eqref{local:eq-nl-Y-high-high}, since the estimates for the remaining three terms are similar. First, we use 
\eqref{local:eq-high-high-e1} from Lemma~\ref{local:lem-high-high} with 
\begin{equation*}
\bigl( \alpha,\beta,\gamma,\nu,\eta \bigr) 
\coloneq \bigl(-\kappa,\, 1+2\kappa, \, -1+\reg_2, \, -\lambda, \, -1-3\kappa\bigr)
\end{equation*}
and Proposition~\ref{local:prop-higgsdet-sdet}, which imply that
\begin{align*}
\big\| \bil \bigl( \higgsdet_\star \parasim \covd \philinear[\oeps][][][\ostar]\bigr)  \big\|_{\nl^{-1+\reg_2,-\lambda}} 
&\lesssim T^{-\bigO(\lambda)} \| \philinear[\oeps][][][\ostar] \|_{\sol^{-\kappa}} \| \higgsdet_\star\|_{\sol^{1+2\kappa,-1-3\kappa}} \\ 
&\lesssim T^{-\bigO(\lambda)} \| \philinear[\oeps][][][\ostar] \|_{\sol^{-\kappa}} \bigl( \higgsdetcombined \bigr).   
\end{align*}
Second, we use \eqref{local:eq-high-high-e1} from Lemma~\ref{local:lem-high-high} with 
\begin{equation*}
\bigl( \alpha,\beta,\gamma,\nu,\eta \bigr) 
\coloneq \bigl(-\kappa,\, 1+\reg_2, \, -1+\reg_2, \, -\lambda, \, -\lambda\bigr), 
\end{equation*}
which implies that 
\begin{equation*}
\big\| \bil \bigl( \phiy_{\eps,\star} \parasim \covd \philinear[\oeps][][][\ostar]\bigr)  \big\|_{\nl^{-1+\reg_2,-\lambda}} 
\lesssim T^{\frac{1}{2}-\bigO(\lambda)} \| \philinear[\oeps][][][\ostar] \|_{\sol^{-\kappa}} \| \phiy_{\eps,\star}\|_{\sol^{1+\reg_2,-\lambda}}. 
\end{equation*}
Thus, both terms make an acceptable contribution to \eqref{local:eq-nl-Y-high-high}.
\end{proof}

\begin{remark}
The estimate of $\bil ( \higgsdet_\star \parasim \covd \philinear[\oeps][][][\ostar])$ is one of the reasons why we place $Y_\eps$ in the weighted space $\sol^{1+\reg_2,-\lambda}$, rather than the unweighted space $\sol^{1+\reg_2,0}$. Indeed, if we pick any parameters $\beta\in (0,2]$ and $\eta\leq 0$, then $\| \higgsdet_\star\|_{\sol^{\beta,\eta}}$ can only be controlled under the condition $\beta+\eta\leq - \kappa$. With $(\alpha,\gamma)=(-\kappa,-1+\reg_2)$, the seventh condition in \eqref{local:eq-high-high-c1} then implies that 
\begin{align*}
\nu < \alpha + (1-\theta) \beta - 1 -\gamma +\eta 
\leq -\kappa + (1-\theta) \beta - \reg_2 - \beta - \kappa 
= - \reg_2 -\kappa- \kappa - \theta \beta.
\end{align*}
In particular, $\nu$ must be negative. We also mention that, rather than looking at 
\eqref{local:eq-high-high-c1}, one  can also see that $\nu$ must be negative\footnote{In scaling considerations, it is best to ignore the parameter $\theta$, which is only included in Definition~\ref{local:def-paraproduct} for technical reasons.} from scaling considerations. 
\end{remark}

In the next lemma, we control the terms involving both a total derivative and high$\times$high-interactions, i.e., the terms in 
\eqref{local:eq-Y-final-new-4} and the first summand in \eqref{local:eq-y-final-new-4}.

\begin{lemma}\label{local:lem-nl-total-derivative}
Let Assumptions~\ref{local:assumption-B},~\ref{local:assumption-higgs-initial}, and~\ref{local:assumption-nonlinear-setup} be satisfied. Then, it holds that 
\begin{equation}\label{local:eq-nl-total-derivative}
\begin{aligned}
\big\| \covd^j  &\big[  \Wrem_{\eps,\star} \parasim \linear[\oeps][][j][\ostar]\big]  \big\|_{\nl^{-1+\reg_2,-\lambda}}
+ \big\| \covd  \big[  \Wrem_{\eps,\star}^j \parasim \linear[\oeps][][j][\ostar] \big] \big\|_{\nl^{-1+\reg_2,-\lambda}} \\
&+  \big\| \covd \bigl( \bil\bigl( \phiw_{\eps,\star} \parasim \philinear[\oeps][][][\ostar] \bigr)\bigr)  \big\|_{\nl^{-1+\reg_2,-\lambda}}
+ \big\| \covd_j \bigl( \linear[\oeps][j][][\ostar]  \parasim  \phiw_{\eps,\star} \bigr)  \big\|_{\nl^{-1+\reg_2,-\lambda}} \\
&\lesssim T^{-\lambda_1} 
 \tripn{ \Xi_{\Beps} } \bigl( 
 \higgsdetcombined 
 + T^{\frac{1}{2}} \| (R_{\eps,\star},\phir_{\eps,\star}) \|_{\sol^{1-\reg_1,-\lambda}} \bigr). 
\end{aligned}
\end{equation}
\end{lemma}

\begin{proof}
We only control the third summand in \eqref{local:eq-nl-total-derivative}, since the remaining three terms can be controlled using similar arguments. Using \eqref{local:eq-high-high-e2} from Lemma~\ref{local:lem-high-high} with
\begin{equation*}
\bigl( \alpha,\beta,\gamma,\nu,\eta \bigr) 
\coloneq \bigl(-\kappa,\, \lambda-\kappa,\,-1+\reg_2,\,-\lambda,\,-\lambda \bigr), 
\end{equation*}
we obtain that
\begin{equation*}
 \big\| \covd \bigl( \bil\bigl( \phiw_{\eps,\star} \parasim \philinear[\oeps][][][\ostar] \bigr)\bigr)  \big\|_{\nl^{-1+\reg_2,-\lambda}} 
 \lesssim T^{-\bigO(\lambda)} \big\| \philinear[\oeps][][][\ostar] \big\|_{\sol^{-\kappa}} \| \phiw_{\eps,\star} \|_{\sol^{\lambda-\kappa,-\lambda}}.
\end{equation*}
Using Proposition~\ref{local:prop-higgsdet-sdet} and the trivial estimate $t\leq T$, we can further estimate
\begin{align*}
\| \phiw_{\eps,\star} \|_{\sol^{\lambda-\kappa,-\lambda}}
&\leq \| \higgsdet_{\star} \|_{\sol^{\lambda-\kappa,-\lambda}} 
+ \| \phir_{\eps,\star} \|_{\sol^{\lambda-\kappa,-\lambda}} \\
&\lesssim T^{-\lossdel} \bigl( \higgsdetcombined \bigr) + T^{\frac{1-\lambda-\reg_1+\kappa}{2}} \| \phir_{\eps,\star} \|_{\sol^{1-\reg_1,-\lambda}}. 
\end{align*}
After combining both estimates, this finishes the proof. 
\end{proof}

In the next lemma, we control terms involving $\covd \linear[\oeps][][][]$ or $\covd \philinear[\oeps][][][]$ at low frequencies, i.e., the terms in~\eqref{local:eq-Y-final-new-5} and the second summand in \eqref{local:eq-y-final-new-4}.

\begin{lemma}\label{local:lem-nl-stochastic-low}
Let Assumptions~\ref{local:assumption-B},~\ref{local:assumption-higgs-initial}, and~\ref{local:assumption-nonlinear-setup} be satisfied. Then, it holds that 
\begin{equation}\label{local:eq-nl-stochastic-low}
\begin{aligned}
\big\| 
\big[ \Wrem_{\eps,\star,j} &\paragg \bigl( 2 \covd^j \linear[\oeps][][][\ostar] - \covd \linear[\oeps][j][][\ostar]\bigr)\big]
\big\|_{\nl^{-1+\reg_2,-\lambda}} 
+ \big\| \Wrem_{\eps,\star}^j \paragg \covd_j \philinear[\oeps][][][\ostar]  \big\|_{\nl^{-1+\reg_2,-\lambda}} \\ 
&+ \big\| \bil\bigl( \phiw_{\eps,\star} \paragg \covd \philinear[\oeps][][][\ostar] \bigr)  \big\|_{\nl^{-1+\reg_2,-\lambda}} \\
&\lesssim  T^{-\lambda_1} 
 \tripn{ \Xi_{\Beps} }  \bigl(  
 \higgsdetcombined
 + T^{\frac{1}{2}} \| (R_{\eps,\star},\phir_{\eps,\star}) \|_{\sol^{1-\reg_1,-\lambda}} \bigr).
\end{aligned}
\end{equation}
\end{lemma}

\begin{proof}
We only estimate the third summand in \eqref{local:eq-nl-stochastic-low}, since the remaining two terms can be treated using similar arguments. We use \eqref{local:eq-low-high-e1} from Lemma~\ref{local:lem-low-high} with
\begin{equation*}
\bigl( \alpha,\beta,\gamma,\nu,\eta \bigr) 
\coloneq \bigl(-\kappa,\, 1-\reg_1, \, -1+\reg_2, \, -\lambda, \, -1+\reg_1-\kappa \bigr), 
\end{equation*}
which implies that
\begin{equation*}
\big\| \bil\bigl( \phiw_{\eps,\star} \paragg \covd \philinear[\oeps][][][\ostar] \bigr)  \big\|_{\nl^{-1+\reg_2,-\lambda}} 
\lesssim T^{-\bigO(\lambda)} \| \philinear[\oeps][][][\ostar] \|_{\sol^{-\kappa}} \| \phiw_{\eps,\star} \|_{\sol^{1-\reg_1,-1+\reg_1-\kappa}}. 
\end{equation*}
Using Proposition~\ref{local:prop-higgsdet-sdet}, we also have that
\begin{align*}
\| \phiw_{\eps,\star} \|_{\sol^{1-\reg_1,-1+\reg_1-\kappa}}
&\leq \| \higgsdet_{\star} \|_{\sol^{1-\reg_1,-1+\reg_1-\kappa}} 
+ \| \phir_{\eps,\star} \|_{\sol^{1-\reg_1,-1+\reg_1-\kappa}} \\
&\lesssim T^{-\lossdel} 
\bigl( \higgsdetcombined \bigr) + T^{\frac{1-\lambda-\reg_1+\kappa}{2}} \| \phir_{\eps,\star} \|_{\sol^{1-\reg_1,-\lambda}}. 
\end{align*}
After combining both estimates, this completes the proof. 
\end{proof}

In the next lemma, we control terms involving our covariant linear stochastic objects, covariant derivatives of $W_\eps$ or $\phiw_{\eps}$, and the non-resonant product $\paransim$. That is, we control the terms in \eqref{local:eq-Y-final-new-6} 
and the last summand in \eqref{local:eq-y-final-new-4}.

\begin{lemma}\label{local:lem-nl-nonres}
Let Assumptions~\ref{local:assumption-B},~\ref{local:assumption-higgs-initial}, and~\ref{local:assumption-nonlinear-setup} be satisfied.  Then, it holds that 
\begin{equation}\label{local:eq-nl-nonres}
\begin{aligned}
 \big\| 
\big[ \linear[\oeps][][j][\ostar] &\paransim \bigl( 2 \covd^j \Wrem_{\eps,\star} - \covd \Wrem_{\eps,\star}^j \bigr)\big] 
\big\|_{\nl^{-1+\reg_2,-\lambda}} 
+ \big\| \bil\bigl( \philinear[\oeps][][][\ostar] \paransim  \covd \phiw_{\eps,\star}  \bigr)  \big\|_{\nl^{-1+\reg_2,-\lambda}} \\ 
&+  \big\| \linear[\oeps][j][][\ostar] \paransim  \covd_j \phiw_{\eps,\star}  \big\|_{\nl^{-1+\reg_2,-\lambda}}  \\ 
&\lesssim  T^{-\lambda_1} 
 \tripn{ \Xi_{\Beps} }  \bigl(  
 \higgsdetcombined
 + T^{\frac{1}{2}} \| (R_{\eps,\star},\phir_{\eps,\star}) \|_{\sol^{1-\reg_1,-\lambda}} \bigr).
\end{aligned}
\end{equation}
\end{lemma}

\begin{proof}
We only control the second summand in \eqref{local:eq-nl-nonres}, since the remaining two terms can be treated using similar arguments. To this end, we first decompose
\begin{equation*}
\bil\bigl(\philinear[\oeps][][][\ostar] \paransim  \covd \phiw_{\eps,\star}  \bigr)
= \bil\bigl( \philinear[\oeps][][][\ostar] \parall \covd \phiw_{\eps,\star}  \bigr) 
+ \bil\bigl( \philinear[\oeps][][][\ostar] \paragg \covd \phiw_{\eps,\star}  \bigr).
\end{equation*} 
In addition, we choose the parameters
\begin{equation*}
\bigl( \alpha,\beta,\gamma,\nu,\eta \bigr) 
\coloneq \bigl(-\kappa, \, 1-\reg_1, \, -1+\reg_2, \, -\lambda, \, -1+\reg_1-\kappa \bigr). 
\end{equation*}
Using \eqref{local:eq-low-high-e2} from Lemma~\ref{local:lem-low-high} and Lemma~\ref{local:lem-high-low}, we then obtain that 
\begin{align*}
 \big\| \bil\bigl(\philinear[\oeps][][][\ostar] &\parall  \covd \phiw_{\eps,\star}  \bigr)  \big\|_{\nl^{-1+\reg_2,-\lambda}} 
 +  \big\| \bil\bigl( \philinear[\oeps][][][\ostar] \paragg  \covd \phiw_{\eps,\star}  \bigr)  \big\|_{\nl^{-1+\reg_2,-\lambda}} \\
 &\lesssim T^{-\bigO(\lambda)} \| \philinear[\oeps][][][\ostar] \|_{\modsol^{-\kappa}} \| \phiw_{\eps,\star} \|_{\sol^{1-\reg_1,-1+\reg_1-\kappa}}. 
\end{align*}
Using Proposition~\ref{local:prop-higgsdet-sdet}, we have that 
\begin{align*}
\| \phiw_{\eps,\star} \|_{\sol^{1-\reg_1,-1+\reg_1-\kappa}}
&\leq \| \higgsdet_{\star} \|_{\sol^{1-\reg_1,-1+\reg_1-\kappa}} 
+ \| \phir_{\eps,\star} \|_{\sol^{1-\reg_1,-1+\reg_1-\kappa}} \\  
&\lesssim T^{-\lossdel} \bigl( \higgsdetcombined \bigr) + T^{\frac{1-\lambda-\reg_1+\kappa}{2}} 
\| \phir_{\eps,\star} \|_{\sol^{1-\reg_1,-\lambda}}. 
\end{align*}
By combining both estimates, we arrive at the desired estimate of the second term in \eqref{local:eq-nl-nonres}. 
\end{proof}

In the next lemma, we control the parts of the derivative nonlinearities containing no covariant linear stochastic objects, i.e., the terms in \eqref{local:eq-Y-final-new-7} and \eqref{local:eq-y-final-new-5}. While these terms have better regularity properties than the terms from Lemmas~\ref{local:lem-nl-Xterms}--\ref{local:lem-nl-nonres} above, we need to be mindful of their size. 

\begin{lemma}\label{local:lem-nl-remainders}
Let Assumptions~\ref{local:assumption-B},~\ref{local:assumption-higgs-initial}, and~\ref{local:assumption-nonlinear-setup} be satisfied.  Then, it holds that
\begin{equation}\label{local:eq-nl-remainders}
\begin{aligned}
\big\| \big[ \Wrem_{\eps,j}  &, 2 \covd^j \Wrem_\eps - \covd \Wrem_\eps^j \big] - \big[ \Sdet_j , 2 \covd^j \Sdet - \covd \Sdet^j \big] \big\|_{\nl^{-1+\reg_2,-\lambda}} \\  
&+ \big\| \bil \bigl( \phiw_\eps \otimes \covd \phiw_\eps\bigr) - \bil \bigl( \higgsdet \otimes \covd \higgsdet\bigr) \big\|_{\nl^{-1+\reg_2,-\lambda}} + \big\| \Wrem_\eps^j \covd_j \phiw_\eps - \Sdet^j \covd_j \higgsdet \big\|_{\nl^{-1+\reg_2,-\lambda}} \\ 
&\lesssim  T^{\frac{1}{2}-\lambda_1} 
\| (R_{\eps,\star},\phir_{\eps,\star}) \|_{\sol^{1-\reg_1,-\lambda}}	  
 \bigl( \higgsdetcombined 
+ T^{\frac{1}{2}} \| (R_{\eps,\star},\phir_{\eps,\star}) \|_{\sol^{1-\reg_1,-\lambda}}     \bigr). 
\end{aligned}
\end{equation}
\end{lemma}

\begin{remark}\label{local:rem-nl-remainders}
In order to obtain the estimates at the end of Proposition~\ref{local:prop-main}, it is important that the right-hand side of \eqref{local:eq-nl-remainders} contains no term of the form $T^{-\lambda_1} (\higgsdetcombined)^2$. This is the reason for subtracting the terms  $[ \Sdet_j , 2 \covd^j \Sdet - \covd \Sdet^j]$, $\bil( \higgsdet \otimes \covd \higgsdet)$, and $\Sdet^j \covd_j \higgsdet$ in \eqref{local:eq-nl-remainders}. 
\end{remark}

\begin{proof}
We only control the second summand in \eqref{local:eq-nl-remainders}, since the first and third summands can be treated using similar arguments. To this end, we decompose
\begin{align}
\bil \bigl( \phiw_\eps &\otimes \covd \phiw_\eps\bigr) - \bil \bigl( \higgsdet \otimes \covd \higgsdet\bigr) \notag\\
&= 
\bil \bigl( \phir_\eps \otimes \covd \higgsdet \bigr)
+ \bil \bigl( \higgsdet \otimes \covd \phir_\eps \bigr)
+ \bil \bigl( \phir_\eps \otimes \covd \phir_\eps \bigr). 
\label{local:eq-nl-remainders-p1}
\end{align}
We now estimate all three terms in \eqref{local:eq-nl-remainders-p1} using Lemma~\ref{local:lem-product-estimate}. We will always choose $(\gamma,\nu)=(-1+\reg_2,-\lambda)$, but the choice of $(\alpha,\beta,\eta_1,\eta_2)$ depends on the term. For the first term, we choose
\begin{equation*}
\bigl( \alpha,\beta,\eta_1,\eta_2\bigr) \coloneq 
\bigl( 1-\reg_1-\lambda,\, 1-\reg_1,\, 0, \, -1+\reg_1-\kappa \bigr). 
\end{equation*}
Together with Lemma~\ref{local:lem-sol-different-regularities-weights} and Proposition~\ref{local:prop-higgsdet-sdet}, we then obtain that
\begin{align*}
\big\| \bil \bigl( \phir_\eps \otimes \covd \higgsdet \bigr) 
\big\|_{\nl^{-1+\reg_2,-\lambda}}
&\lesssim T^{\frac{1}{2}-\bigO(\lambda)} \| \phir_{\eps,\star} \|_{\sol^{1-\reg_1-\lambda,0}} \| \higgsdet_\star \|_{\sol^{1-\reg_1,-1+\reg_1-\kappa}} \\ 
&\lesssim T^{\frac{1}{2}-\bigO(\lambda)}  \| \phir_{\eps,\star} \|_{\sol^{1-\reg_1,-\lambda}} \bigl( \higgsdetcombined \bigr).  
\end{align*} 
For the second term, we choose
\begin{equation*}
\bigl( \alpha,\beta,\eta_1,\eta_2\bigr) \coloneq 
\bigl( \reg_1+\kappa,\, 1-\reg_1,\, -\reg_1-2\kappa,-\lambda\bigr). 
\end{equation*}
Together with Proposition~\ref{local:prop-higgsdet-sdet}, we obtain 
\begin{align*}
\big\| \bil \bigl( \higgsdet  \otimes \covd \phir_\eps \bigr) 
\big\|_{\nl^{-1+\reg_2,-\lambda}}
&\lesssim T^{\frac{1}{2}-\bigO(\lambda)} \| \higgsdet_\star \|_{\sol^{\reg_1+\kappa,-\reg_1-2\kappa}} \| \phir_{\eps,\star} \|_{\sol^{1-\reg_1,-\lambda}}  \\ 
&\lesssim T^{\frac{1}{2}-\bigO(\lambda)} \bigl( \higgsdetcombined \bigr) \| \phir_{\eps,\star} \|_{\sol^{1-\reg_1,-\lambda}}. 
\end{align*} 
Finally, for the third term, we choose 
\begin{equation*}
\bigl( \alpha,\beta,\eta_1,\eta_2\bigr) \coloneq 
\bigl( 1-\reg_1,\, 1-\reg_1,\, -\lambda, \, -\lambda\bigr). 
\end{equation*}
We then obtain that 
\begin{equation*}
\big\| \bil \bigl( \phir_\eps  \otimes \covd \phir_\eps \bigr) 
\big\|_{\nl^{-1+\reg_2,-\lambda}}
\lesssim T^{1-\bigO(\lambda)} \| \phir_{\eps,\star} \|_{\sol^{1-\reg_1,-\lambda}}^2. 
\end{equation*}
By combining our estimates for the first, second, and third terms in \eqref{local:eq-nl-remainders-p1}, we complete the proof. 
\end{proof}

In the next lemma, we estimate 
\eqref{local:eq-Y-final-new-8}, 
\eqref{local:eq-Y-final-new-9}, 
\eqref{local:eq-y-final-new-6}, and 
\eqref{local:eq-y-final-new-7}, 
i.e., the power-type nonlinearity. 

\begin{lemma}\label{local:lem-nl-power}
Let Assumptions~\ref{local:assumption-B},~\ref{local:assumption-higgs-initial}, and~\ref{local:assumption-nonlinear-setup} be satisfied. Then, it holds that
\begin{align}
\big\|  \big[ & \linear[\oeps][][j][]\hspace{-0.5ex}+ \Wrem_{\eps,j}, 
\big[ \linear[\oeps][j][][] \hspace{-0.5ex}+ \Wrem_\eps^j,  \linear[\oeps] \hspace{-0.5ex}+ \Wrem_\eps \big] \big] - \big[ \Sdet_j, \big[ \Sdet^j , \Sdet \big]\big] - \sigma_{\rho,\eps}^2 \Cov_\frkg^{\frkg}( \linear[\oeps][][][] + \Wrem_\eps) \big\|_{\nl^{-1+\reg_2,-\lambda}}
\notag \\ 
&+ \big\|   \bil \bigl( \bigl( \philinear[\oeps]+\phiw_\eps\bigr) \otimes \bigl(\linear[\oeps]+\Wrem_\eps\bigr) \bigl( \philinear[\oeps] + \phiw_\eps\bigr) \bigr) - \bil\bigl( \higgsdet \otimes \Sdet \higgsdet \bigr) - \sigma_{\rho,\eps}^2 \Cov_{\V}^{\frkg} ( \linear[\oeps][][][]  + \Wrem_\eps) \big\|_{\nl^{-1+\reg_2,-\lambda}} 
\notag \\ 
&+ \big\| \bigl( \linear[\oeps]+\Wrem_\eps\bigr)^j \bigl( \linear[\oeps]+\Wrem_\eps\bigr)_j \bigl( \philinear[\oeps]+\phiw_\eps\bigr)- \Sdet^j \Sdet_j \higgsdet - 2 \sigma_{\rho,\eps}^2 \Cov_{\frkg}^{\V}( \philinear[\oeps][][][] + \phiw_\eps ) \big\|_{\nl^{-1+\reg_2,-\lambda}} 
\label{local:eq-nl-power}\\
&+ \big\| -  \big| \philinear[\oeps] + \phiw_\eps\big|^2 \bigl( \philinear[\oeps]+\phiw_\eps\bigr) + (\dim(\V)+2) \sigma_{\rho,\eps}^2 (\philinear[\oeps]+\phiw_\eps)  + |\higgsdet|^2 \higgsdet \big\|_{\nl^{-1+\reg_2,-\lambda}} 
\notag \\ 
&\lesssim T^{\frac{1}{2}-\lambda_1}
\bigl( \tripn{\Xi_{\Beps}} + T^{\frac{1}{2}} \big\| (R_{\eps,\star},\phir_{\eps,\star}) \big\|_{\sol^{1-\reg_1,-\lambda}} \bigr) \notag \\   
&\qquad \times \bigl( \tripn{ \Xi_{\Beps} } + \higgsdetcombined +  T^{\frac{1}{2}} \big\| (R_{\eps,\star},\phir_{\eps,\star}) \big\|_{\sol^{1-\reg_1,-\lambda}} \bigr)^2. \notag 
\end{align}
\end{lemma}

\begin{remark}\label{local:rem-nl-power}
As in Lemma~\ref{local:lem-nl-remainders} and 
Remark~\ref{local:rem-nl-remainders}, it is important that the right-hand side of \eqref{local:eq-nl-power} contains no term of the form $T^{\frac{1}{2}-\lambda_1} (\higgsdetcombined)^3$. This is the reason for subtracting the terms
$[ \Sdet_j, [ \Sdet^j , \Sdet ]]$,
$\bil( \higgsdet \otimes \Sdet \higgsdet )$,
$\Sdet^j \Sdet_j \higgsdet$,
and $-|\higgsdet|^2 \higgsdet$
in  the four different arguments in \eqref{local:eq-nl-power}. 
\end{remark}

\begin{proof}
We only control the first summand in \eqref{local:eq-nl-power}, since the remaining three summands can be controlled using similar arguments. To briefly simplify the notation, set
\begin{equation*}
\mathsf{Pow} \eqdef \big[  \linear[\oeps][][j][]\hspace{-0.5ex}+ \Wrem_{\eps,j}, 
\big[ \linear[\oeps][j][][] \hspace{-0.5ex}+ \Wrem_\eps^j,  \linear[\oeps] \hspace{-0.5ex}+ \Wrem_\eps \big] \big] - \big[ \Sdet_j, \big[ \Sdet^j , \Sdet \big]\big] - \sigma_{\rho,\eps}^2 \Cov_\frkg^{\frkg}( \linear[\oeps][][][] + \Wrem_\eps). 
\end{equation*}
We now recall our weighted Schauder estimate \eqref{local:eq-weighted-covariant-Schauder-2}, which, together with 
Lemma~\ref{local:lem-sol-different-regularities-weights}
and Definition~\ref{local:def-nl-norm}, implies that 
\begin{equation*}
\big\| \mathsf{Pow} \big\|_{\nl^{-1+\reg_2,-\lambda}}
\lesssim
\big\| \mathsf{Pow} \big\|_{\nl^{-1+\reg_2+2\lambda,-3\lambda}}
\lesssim T^{\frac{1}{2}-\bigO(\lambda)}
\big\| t^{\frac{3\lambda}{2}} \mathsf{Pow} \big\|_{\Cs_t^{-\frac{3\kappa}{2}}L_x^r}. 
\end{equation*}
We also note that, due to Proposition~\ref{local:prop-higgsdet-sdet} and Lemma~\ref{local:lem-sol-different-regularities-weights}, 
\begin{align*}
\| \Sdet_{\star} \|_{\sol^{5\kappa,-\lambda}}
&\lesssim T^{-3\lossdel} T^{\frac{1-5\kappa+\lambda}{2}-\kappa} \bigl( \higgsdetcombined \bigr)^2 \lesssim T^{-\bigO(\lambda)} \bigl( \higgsdetcombined \bigr), \\ 
\| R_{\eps,\star} \|_{\sol^{5\kappa,-\lambda}}
&\lesssim T^{\frac{1-5\kappa-\reg_1}{2}} \| R_{\eps,\star} \|_{\sol^{1-\reg_1,-\lambda}} \lesssim T^{\frac{1}{2}-\bigO(\lambda)} 
\| R_{\eps,\star}\|_{\sol^{1-\reg_1,-\lambda}}. 
\end{align*}
By using the above three estimates, as well as the identity $\Wrem_{\eps}=\Sdet + R_{\eps}$, we can reduce the desired estimate of the first summand in \eqref{local:eq-nl-power} to the estimate
\begin{equation}\label{local:eq-power-p1}
\begin{aligned}
&\quad \Bigl\| t^{\frac{3\lambda}{2}} \Bigl( \big[  \linear[\oeps][][j][]\hspace{-0.5ex}+ \Wrem_{\eps,j}, 
\big[ \linear[\oeps][j][][] \hspace{-0.5ex}+ \Wrem_\eps^j,  \linear[\oeps] \hspace{-0.5ex}+ \Wrem_\eps \big] \big] - \big[ \Sdet_j, \big[ \Sdet^j , \Sdet \big]\big] - \sigma_{\rho,\eps}^2 \Cov_\frkg^{\frkg}( \linear[\oeps][][][] + \Wrem_\eps) \Bigr) \Bigr\|_{\Cs_t^{-\frac{3\kappa}{2}}L_x^r} \\
&\lesssim \tripn{\Xi_{\Beps}} \bigl( \tripn{\Xi_{\Beps}} + \| \Wrem_{\eps,\star} \|_{\sol^{5\kappa,-\lambda}} \bigr)^2 
+ \| R_{\eps,\star} \|_{\sol^{5\kappa,-\lambda}} 
\bigl(  \| \Sdet_{\star} \|_{\sol^{5\kappa,-\lambda}} + 
 \| R_{\eps,\star} \|_{\sol^{5\kappa,-\lambda}}\bigr)^2. 
\end{aligned}
\end{equation}
In order to prove \eqref{local:eq-power-p1}, we decompose the argument into 
\begin{align}
\big[  \linear[\oeps][][j][] &+ \Wrem_{\eps,j}, \big[  \linear[\oeps][j][][] + \Wrem_{\eps}^j ,  \linear[\oeps][k][][] + \Wrem_{\eps}^k \big] \big] 
- \sigma_{\rho,\eps}^2 \Cov_\frkg^{\frkg}( \linear[\oeps][k][][] + \Wrem_\eps^k)
- \big[ \Sdet_j, \big[ \Sdet^j , \Sdet^k \big]\big] 
\notag \\
&= \big[  \linear[\oeps][][j][], \big[  \linear[\oeps][j][][],   \linear[\oeps][k][][]  \big] \big] 
- \sigma_{\rho,\eps}^2 \Cov_\frkg^{\frkg} \bigl( \linear[\oeps][k][][] \bigr) 
\label{local:eq-power-type-p2} \allowdisplaybreaks[4]\\ 
&\quad+  \big[  \Wrem_{\eps,j}, \big[  \linear[\oeps][j][][],  \linear[\oeps][k][][] \big] \big] 
\label{local:eq-power-type-p3}\allowdisplaybreaks[4]\\ 
&\quad+ \big[  \linear[\oeps][][j][] , \big[   \Wrem_{\eps}^j ,  \linear[\oeps][k][][] \big] \big] + \sigma_{\rho,\eps}^2 \Cov_\frkg^{\frkg} \bigl( \Wrem_{\eps}^k \bigr) 
\label{local:eq-power-type-p4} \allowdisplaybreaks[4]\\ 
&\quad+ 
\big[  \linear[\oeps][][j][] , \big[   \linear[\oeps][j][][] ,  \Wrem_{\eps}^k\big] \big]  - 2 \sigma_{\rho,\eps}^2 \Cov_\frkg^{\frkg} \bigl( \Wrem_{\eps}^k \bigr)   
\label{local:eq-power-type-p5} \allowdisplaybreaks[4]\\ 
&\quad+ \big[  \linear[\oeps][][j][] , \big[   \Wrem_{\eps}^j , \Wrem_{\eps}^k \big] \big] + \big[ \Wrem_{\eps,j}, \big[  \linear[\oeps][j][][] ,  \Wrem_{\eps}^k \big] \big] + \big[  \Wrem_{\eps,j}, \big[   \Wrem_{\eps}^j ,   \linear[\oeps][k][][]\big] \big] 
\label{local:eq-power-type-p6} \allowdisplaybreaks[4]\\ 
&\quad+ \big[  \Wrem_{\eps,j}, \big[   \Wrem_{\eps}^j ,   \Wrem_{\eps}^k \big] \big]- \big[ \Sdet_j, \big[ \Sdet^j , \Sdet^k \big]\big].
\label{local:eq-power-type-p7} 
\end{align}
To estimate the first term \eqref{local:eq-power-type-p2}, we observe that the argument can be written as the image of 
\begin{equation*}
\linear[\oeps][][j][] \otimes \linear[\oeps][j][][] \otimes \linear[\oeps][k][][]  
- \delta^j_j \sigma_{\rho,\eps}^2 \bigl( \ONBg_a \otimes \ONBg^a \otimes \linear[\oeps][k][][] \bigr) 
- \delta^k_j \sigma_{\rho,\eps}^2\bigl( \ONBg_a \otimes \linear[\oeps][j][][] \otimes \ONBg^a \bigr)
- \delta^{kj} \sigma_{\rho,\eps}^2 \bigl( \linear[\oeps][][j][] \otimes \ONBg_a \otimes \ONBg^a \bigr), 
\end{equation*}
where $(\ONBg_a)_a$ is as in Definition~\ref{prelim:def-Casimir-induced}, under the linear map 
\begin{equation*}
\boldsymbol{L}(X\otimes Y \otimes Z)\eqdef [X,[Y,Z]]. 
\end{equation*} 
Indeed, in order to see that the resulting counterterm is as in \eqref{local:eq-power-type-p2}, we use that 
\begin{align*}
&\quad \boldsymbol{L} \Big( \delta^j_j \bigl( \ONBg_a \otimes \ONBg^a \otimes \linear[\oeps][k][][] \bigr) 
+ \delta^k_j \bigl( \ONBg_a \otimes \linear[\oeps][j][][] \otimes \ONBg^a \bigr)
+ \delta^{kj} \bigl( \linear[\oeps][][j][] \otimes \ONBg_a \otimes \ONBg^a \bigr) \Big) \\
&= 2 \Cov_\frkg^{\frkg}( \linear[\oeps][k][][]) -  \Cov_\frkg^{\frkg}( \linear[\oeps][k][][]) + 0 = \Cov_\frkg^{\frkg}( \linear[\oeps][k][][]). 
\end{align*}
Using Definition~\ref{local:def-data-set-norm}, we then obtain that 
\begin{equation*}
\big\| t^{\frac{3\lambda}{2}} \eqref{local:eq-power-type-p2}\big\|_{\Cs_t^{-\frac{3\kappa}{2}}L_x^r}
\lesssim \big\| \eqref{local:eq-power-type-p2}\big\|_{\Cs_t^{-\frac{3\kappa}{2}}L_x^r}
\lesssim \tripn{ \Xi_{\Beps} }^3. 
\end{equation*}
To estimate the terms in \eqref{local:eq-power-type-p3}, \eqref{local:eq-power-type-p4}, and \eqref{local:eq-power-type-p5}, we note that all three of them can be written as linear combinations of images of 
\begin{equation*}
\bigl(  \linear[\oeps][e][][] \otimes \linear[\oeps][f][][] - \delta^{ef} \sigma_{\rho,\eps}^2 \Cov_{\frkg}\bigr) \otimes \Wrem_\eps^h,
\end{equation*}
where $e,f,h\in [2]$, under the map $\boldsymbol{L}$ and permutations of the tensor factors.  For example, the term in \eqref{local:eq-power-type-p5} can be written as 
\begin{equation*}
\delta_{ej} \delta_f^{j} \delta_{h}^{k} \boldsymbol{L} \Big( \bigl(  \linear[\oeps][e][][] \otimes \linear[\oeps][f][][] - \delta^{ef} \sigma_{\rho,\eps}^2 \Cov_{\frkg} \bigr) \otimes \Wrem^h_\eps\Big). 
\end{equation*}
Using Lemma~\ref{prelim:lem-product-temporal}, Lemma~\ref{local:lem-solution-spaces-hoelder}, and Definition~\ref{local:def-data-set-norm}, we then obtain that 
\begin{align*}
&\quad\big\| t^{\frac{3\lambda}{2}} \eqref{local:eq-power-type-p3} \big\|_{\Cs_t^{-\frac{3\kappa}{2}} L_x^r}
+ \big\| t^{\frac{3\lambda}{2}} \eqref{local:eq-power-type-p4} \big\|_{\Cs_t^{-\frac{3\kappa}{2}} L_x^r}
+ \big\| t^{\frac{3\lambda}{2}} \eqref{local:eq-power-type-p5} \big\|_{\Cs_t^{-\frac{3\kappa}{2}} L_x^r} \\
&\lesssim \max_{e,f\in [2]} \big\|   \linear[\oeps][e][][] \otimes \linear[\oeps][f][][] - \delta^{ef} \sigma_{\rho,\eps}^2 \Cov_{\frkg}\big\|_{\Cs_t^{-\kappa}L_x^r}
\big\| t^{\frac{3\lambda}{2}} \Wrem_{\eps} \big\|_{\Cs_t^{2\kappa}L_x^\infty}
\lesssim \tripn{ \Xi_{\Beps} }^2 \| \Wrem_{\eps,\star} \|_{\sol^{5\kappa,-3\lambda}}.  
\end{align*}
To estimate the terms in~\eqref{local:eq-power-type-p6}, we use Lemma~\ref{prelim:lem-product-temporal} and Lemma~\ref{local:lem-solution-spaces-hoelder}, which yield that 
\begin{align*}
 \big\| t^{\frac{3\lambda}{2}} \eqref{local:eq-power-type-p6} \big\|_{\Cs_t^{-\frac{3\kappa}{2}} L_x^r}
&\lesssim \|  \linear[\oeps][][][] \|_{\Cs_t^{-\frac{3\kappa}{4}} L_x^r}
\big\| t^{\frac{3\lambda}{2}} (\Wrem_\eps \otimes \Wrem_\eps) \big\|_{\Cs_t^{\kappa}L_x^\infty}\\
&\lesssim \|   \linear[\oeps][][][] \|_{\Cs_t^{-\frac{3\kappa}{4}} L_x^r} 
\big\| t^{\frac{3\lambda}{4}} \Wrem_\eps \big\|_{\Cs_t^\kappa L_x^\infty}^2 
\lesssim \tripn{ \Xi_{\Beps} }  \big\| \Wrem_{\eps,\star} \big\|_{\sol^{3\kappa,-\frac{3\lambda}{2}}}^2.  
\end{align*}
Finally, it remains to estimate the terms in~\eqref{local:eq-power-type-p7}. To this end, we insert the identity $\Wrem_{\eps}=\Sdet+R_{\eps}$ into the first summand in 
\eqref{local:eq-power-type-p7}. The cubic term in $\Sdet$ then cancels and each remaining term contains at least one factor of $R_{\eps}$. Using this, we obtain that 
\begin{align*}
 \big\| t^{\frac{3\lambda}{2}} \eqref{local:eq-power-type-p7} \big\|_{\Cs_t^{-\frac{3\kappa}{2}} L_x^r}
&\lesssim \big\| t^{\frac{3\lambda}{2}} \eqref{local:eq-power-type-p7} \big\|_{L_t^\infty L_x^\infty} \\
&\lesssim \big\| t^{\frac{\lambda}{2}} R_{\eps} \big\|_{L_t^\infty L_x^\infty}
\bigl(\big\| t^{\frac{\lambda}{2}} \Sdet \big\|_{L_t^\infty L_x^\infty} + \big\| t^{\frac{\lambda}{2}} R_{\eps} \big\|_{L_t^\infty L_x^\infty} \bigr)^2 \\ 
&\lesssim \big\|  R_{\eps,\star} \big\|_{\sol^{\kappa,-\lambda}}
\bigl(\big\| \Sdet_\star \big\|_{\sol^{\kappa,-\lambda}} + \big\| R_{\eps,\star} \big\|_{\sol^{\kappa,-\lambda}} \bigr)^2. 
\end{align*}
This completes our last estimate of the terms in \eqref{local:eq-power-type-p2}--\eqref{local:eq-power-type-p7}, and hence the proof of \eqref{local:eq-power-p1}. 
\end{proof}

In the last lemma of this subsection, we estimate the $\frac{1}{8\pi}$- and the mass term, i.e., the terms in \eqref{local:eq-Y-final-new-10} and \eqref{local:eq-y-final-new-last}.

\begin{lemma}\label{local:lem-8pi-counterterm}
Let Assumptions~\ref{local:assumption-B},~\ref{local:assumption-higgs-initial}, and~\ref{local:assumption-nonlinear-setup} be satisfied. Then, it holds that
\begin{align*}
\big\| \tfrac{1}{8\pi}& \Cov_{\V}^{\frkg} \bigl( \linear[\oeps] + W_\eps \bigr) \big\|_{\nl^{-1+\reg_2,-\lambda}} + \big\| \mass \bigl( \philinear[\oeps] + \phiw_\eps \bigr) \big\|_{\nl^{-1+\reg_2,-\lambda}}  \\ 
&\lesssim T^{\frac{1}{2}-\lambda_1} \bigl( \tripn{ \Xi_{\Beps} } + \higgsdetcombined + T^{\frac{1}{2}} \| (R_{\eps,\star},\phir_{\eps,\star}) \|_{\sol^{1-\reg_1,-\lambda}} \bigr).
\end{align*}
\end{lemma}

\begin{proof}
Recall that $W_\eps=\Sdet+R_\eps$. Using our weighted Schauder estimates from \eqref{local:eq-weighted-covariant-Schauder-1} and \eqref{local:eq-weighted-covariant-Schauder-2}, we obtain that 
\begin{align*}
&\quad \| \linear[\oeps][][][] \|_{\nl^{-1+\reg_2,-\lambda}} 
+ \| \Sdet \|_{\nl^{-1+\reg_2,-\lambda}}
+ \| R_\eps \|_{\nl^{-1+\reg_2,-\lambda}} \\ 
&\lesssim T^{\frac{1}{2}-\bigO(\lambda)} 
\bigl( \| \linear[\oeps][][][] \|_{\Cs_t^{-\kappa}L_x^r} + \| \Sdet \|_{L_{t,x}^\infty} + \| t^{\frac{\lambda}{2}} R_\eps \|_{L_{t,x}^\infty} \bigr) \\
&\lesssim T^{\frac{1}{2}-\bigO(\lambda)} \bigl( \tripn{ \Xi_{\Beps} } + \higgsdetcombined + T^{\frac{1}{2}}\| R_{\eps,\star} \|_{\sol^{1-\reg_1,-\lambda}} \bigr).
\end{align*}
In the last line, we also used the bound from Assumption~\ref{local:assumption-higgs-initial} and Proposition~\ref{prop:Sdet-higgsdet-estimates}. The estimate of $\mass (\philinear[\oeps] + \phiw_\eps)$ is similar, and we omit the details. 
\end{proof}

\subsection{Paracontrolled estimates}\label{section:local-paracontrolled}

We now estimate the terms in \eqref{local:eq-Y-final-new-2}
and \eqref{local:eq-y-final-new-2}, i.e., the high$\times$high-interactions involving the paracontrolled terms $X_\eps$ and $\phix_{\eps}$. The main estimate of this subsection is the following abstract estimate, which will later be used to control our more specific interactions (see Corollary~\ref{local:cor-paracontrolled-specific}). To state it, it is notationally convenient to extend our $\nl^{\alpha,\nu}$-norm (see Definition \ref{local:def-nl-norm}) for $\Y$-valued maps to a seminorm on $(\Y\otimes \Y)$-valued maps. This can be done by first choosing a basis $(\Qc_i)_i$ of the space of $\liegroup$-equivariant linear maps from $\Y\otimes \Y$ into $\Y$ and then defining 
\begin{equation*}
\| \phi \|_{\nl^{\alpha,\nu}} \eqdef \max_{i} \| \Qc_i ( \phi ) \|_{\nl^{\alpha,\nu}},
\end{equation*}
where $\phi\colon [0,T]\times \T^2 \rightarrow \Y \otimes \Y$.

\begin{proposition}[Paracontrolled estimate]\label{local:prop-paracontrolled}
Let Assumption~\ref{local:assumption-B} be satisfied. Furthermore, let $k,\ell\in [2]$ and let $\Qc \colon \Y \times \Y \rightarrow \Y$ be a $\liegroup$-equivariant bilinear map. Furthermore, let $\mathcal{K} \colon \Y \rightarrow \Y \otimes \Y$ be defined as $\mathcal{K}(W)\coloneq \Qc(W,\ONBY_m)\otimes \ONBY^m$, where $(\ONBY_m)_m$ is an orthonormal basis of $\Y$. Finally, let $W_\star =(W_M)_{M\geq N_{\lo}}$ be a sequence of maps from $[0,T]\times \T^2$ into $\Y$ such that each $W_M$ is supported on times $t\gtrsim M^{-2}$. Then, it holds that\footnote{In the first summand on the left-hand side of \eqref{local:eq-paracontrolled-main-estimate}, $\parasim$ captures the high$\times$high-portion of the tensor product in $\Y\otimes \Y$.}
\begin{align}
&\quad \Big\|  \Duh^{F}_{\Y,B,\star} \bigl( 
\Qc \bigl( W_\diamond \parall \covd^k \graylinear[\oeps][][][\odiamond] \bigr)\bigr) \parasim \covd^\ell \graylinear[\oeps][][][\ostar]   - \tfrac{1}{4}  \delta^{k\ell} (\sigma_{\rho,\eps}^2 + \Crho) \mathcal{K}(W) 
\Big\|_{\nl^{-1+\reg_2,-\lambda}} \notag \\ 
&\lesssim  T^{\frac{1-\lambda_1}{2}} \| W_\star \|_{\sol^{1-\reg_1,-1+\reg_1-\kappa}} \tripn{ \Xi_{\Beps} }^2. \label{local:eq-paracontrolled-main-estimate}
\end{align}
\end{proposition}

\begin{remark}
The reason for using the $\| W_\star \|_{\sol^{1-\reg_1,-1+\reg_1-\kappa}}$-norm in Proposition~\ref{local:prop-paracontrolled} is that it is under control for all terms among $\Sdet_\star$, $\higgsdet_\star$, $X_{\eps,\star}$, $\phix_{\eps,\star}$, $Y_{\eps,\star}$, and $\phiy_{\eps,\star}$. 
\end{remark}

\begin{remark}\label{local:rem-paracontrolled-indices}
A similar statement as in Proposition~\ref{local:prop-paracontrolled} also holds if the tensor product in $\Y$ is replaced by a second $\liegroup$-equivariant bilinear map and each $\Y$ is replaced by one of its components.  To be more precise, let $a,b,c,d,e,f\in[3]$ and let $\Qc_1\colon \Y^{(a)} \times \Y^{(b)} \rightarrow \Y^{(c)}$ and $\Qc_2 \colon \Y^{(d)} \times \Y^{(e)} \rightarrow \Y^{(f)}$ be $\liegroup$-equivariant bilinear maps. Then, 
\begin{align}
&\quad \Big\| \Qc_2 \Big( \Duh^{F,dc}_{\Y,B,\star} \bigl( 
\Qc_1 \bigl( W_\diamond^a \parall \covd^k \graylinear[\oeps][b][][\odiamond] \bigr)\bigr) \parasim \covd^\ell \graylinear[\oeps][e][][\ostar] \Big)  - \tfrac{1}{4} \delta^{dc} \delta^{be} \delta^{k\ell} (\sigma_{\rho,\eps}^2 + \Crho) \widetilde{\mathcal{K}}(W^a) 
\Big\|_{\nl^{-1+\reg_2,-\lambda}} \notag \\ 
&\lesssim  T^{\frac{1-\lambda_1}{2}} \| W_\star \|_{\sol^{1-\reg_1,-1+\reg_1-\kappa}} \tripn{ \Xi_{\Beps} }^2.
 \label{local:eq-paracontrolled-main-estimate-alternative}
\end{align}
Here, $\Duh^{F,dc}_{\Y,B,\star}$ denotes the matrix element of $\Duh^{F}_{\Y,B,\star}$ in the $d$-th row and $c$-th column and $\widetilde{\mathcal{K}}\colon \Y^{(a)} \rightarrow \Y^{(f)}$ is defined as 
\begin{equation*}
\widetilde{\mathcal{K}}(W) \eqdef
\begin{cases}
\begin{aligned}
&\Qc_2 \bigl(  \Qc_1 (W,\ONBY^{(b)}_m), \ONBY^{(b),m}\bigr) \qquad &&\mbox{if $b=e$ and $c=d$,} \\ 
&0 &&\mbox{else,}
\end{aligned}
\end{cases}
\end{equation*}
where $(\ONBY^{(b)}_m)_m$ is an orthonormal basis of $\Y^{(b)}$. Indeed, this estimate can be derived by using Proposition~\ref{local:prop-paracontrolled} and using the canonical embeddings and projections between $\Y$ and its components. The formulation~\eqref{local:eq-paracontrolled-main-estimate} above is somewhat cleaner, but the formulation \eqref{local:eq-paracontrolled-main-estimate-alternative} is closer to the estimates in Corollary~\ref{local:cor-paracontrolled-specific}.
\end{remark}

\begin{remark}
The level of generality in Proposition~\ref{local:prop-paracontrolled}, which allows for different $\liegroup$-equivariant bilinear maps, is needed to address the four different terms in Corollary~\ref{local:cor-paracontrolled-specific}. For the sake of concreteness, we encourage the reader to simply think of
\begin{equation*}
4 \bil \Big( \Duh_{B,\star} \bigl( \Wrem_{\eps,\diamond}^j \parall \covd_j \philinear[\oeps][][][\odiamond]\bigr)   \parasim \covd^k \philinear[\oeps][][] [\ostar] \Big) + (\sigma_{\rho,\eps}^2 + \Crho ) \Cov_{\V}^{\frkg}(W_\eps^k). 
\end{equation*}
\end{remark}

The proof of Proposition~\ref{local:prop-paracontrolled} will be postponed until the end of the subsection. In the next lemma, we first prove an approximate version of the Littlewood--Paley trichotomy, which will be used to better understand the structure of the paracontrolled term $X$. 

\begin{lemma}[Approximate Littlewood--Paley trichotomy]\label{local:lem-trichotomy}
Let Assumption~\ref{local:assumption-B} be satisfied. In addition, let $\mathcal{L}\colon \Y \rightarrow \Y$ be a $\liegroup$-equivariant linear map, and let $\Qc\colon \Y\times \Y \rightarrow \Y$ be a $\liegroup$-equivariant bilinear map.  
Furthermore, let $K,M,N\geq N_{\lo}$ satisfy $M\ll_\theta N$ and $K\not \sim_{\theta^\prime} N$, where $\theta^\prime$ is as in \eqref{prelim:eq-parameters}. Finally, let $W_M \colon [0,T]\times \T^2 \rightarrow \Y$ be supported on times $t\gtrsim M^{-2}$ and let $-1+\reg_1-\kappa \leq \eta_1 \leq 0$.  Then, it holds that 
\begin{align}
\Big\| \Duh_{\Y,B,K}^{F} \Big( \mathcal{L} \bigl( \covd_B \graylinear[\oeps][][][N] \bigr) \Big) \Big\|_{L_{t,x}^\infty} &\lesssim K^{-1-2\theta} N^{-\kappa} \| \graylinear[\oeps][][][N]\|_{\modsol^{-\kappa}_N}, \label{local:eq-trichotomy-e1} \\ 
\Big\| t^{-\frac{\eta_1}{2}} \Duh_{\Y,B,K}^{F} \Big( \Qc \bigl( W_M , \covd_B \graylinear[\oeps][][][N] \bigr) \Big) \Big\|_{L_{t,x}^\infty} &\lesssim 
T^{\frac{1}{2}-\lambda} K^{-1-2\theta} N^{-\kappa} \| \graylinear[\oeps][][][N]\|_{\modsol^{-\kappa}_N} \big\| W_M \big\|_{\sol_M^{1-\reg_1,\eta_1}}.  \label{local:eq-trichotomy-e2} 
\end{align}
The estimates \eqref{local:eq-trichotomy-e1} and \eqref{local:eq-trichotomy-e2} also hold if $\Duh_{\Y,B,K}^{F}$ is replaced with simply $\Duh_{\Y,B,K}$. Furthermore, it also holds that
\begin{equation}
\Big\| \covd_B \Duh_{\Y,B,K} \Big( \mathcal{L} \bigl(  \graylinear[\oeps][][][N] \bigr) \Big) \Big\|_{L_{t,x}^\infty} \lesssim K^{-1-2\theta} N^{-\kappa} \| \graylinear[\oeps][][][N]\|_{\modsol^{-\kappa}_N}. \label{local:eq-trichotomy-e3} 
\end{equation} 
\end{lemma}

\begin{remark}
In our proof of Proposition \ref{local:prop-paracontrolled}, we will only need the case $\eta_1 = - 1 +\reg_1 - \kappa$ in Lemma~\ref{local:lem-trichotomy}. However, in the proof of Proposition \ref{objects:prop-quadratic-dlin} below, the case $\eta_1=0$ will be useful. 

We note that \eqref{local:eq-trichotomy-e3} has been stated only for $\Duh_{\Y,B,K}$ and not for $\Duh_{\Y,B,K}^F$. The technical reason for this is that our commutator estimate (Lemma \ref{kernel:lem-commutator-derivative}) is only available for $\Duh_{\Y,B,K}$. While we believe that \eqref{local:eq-trichotomy-e3}  can be extended to $\Duh_{\Y,B,K}^F$, this extension will not be needed in this article, and we therefore do not pursue it here.
\end{remark}

We emphasise that the $N^{-\kappa}$-factors in \eqref{local:eq-trichotomy-e1} and \eqref{local:eq-trichotomy-e2} are not too important, and are only needed to later sum \eqref{local:eq-trichotomy-e1} and \eqref{local:eq-trichotomy-e2} over all $N$. We also note that, due to the constraint $M\ll_\theta N$, \eqref{local:eq-trichotomy-e2} can easily be summed in $M$. 

\begin{proof}[Proof of Lemma~\ref{local:lem-trichotomy}]
Since $\Duh_{\Y,B,K}$ obeys the same Schauder estimates as $\Duh_{\Y,B,K}^{F}$, we only treat the latter operator in \eqref{local:eq-trichotomy-e1} and \eqref{local:eq-trichotomy-e2}. We now first prove the estimate \eqref{local:eq-trichotomy-e1}. For this, it is sufficient to prove that 
\begin{equation}\label{local:eq-trichotomy-q0}
\Big\| \Duh_{\Y,B,K}^{F} \Big( \mathcal{L} \bigl( \covd_B \graylinear[\oeps][][][N] \bigr) \Big) \Big\|_{L_{t,x}^\infty} 
\lesssim T^{-\lossdel} K^{-1} N^{3\kappa} \min\Bigl( \frac{N}{K}, \frac{K}{N} \Bigr)\| \graylinear[\oeps][][][N]\|_{\modsol^{-\kappa}_N}. 
\end{equation}
Indeed, once \eqref{local:eq-trichotomy-q0} has been proven, the assumption $K\not\sim_{\theta^\prime} N$ allows us to use the $\min$-factor in \eqref{local:eq-trichotomy-q0} to obtain the desired estimate \eqref{local:eq-trichotomy-e1}. To see that \eqref{local:eq-trichotomy-q0} holds, we first 
use our Schauder estimate~\eqref{local:eq-covariant-Schauder-1}, which yields
\begin{equation}\label{local:eq-trichotomy-r4}
\begin{aligned}
&\quad \Big\| \Duh_{\Y,B,K}^{F} \Big( \mathcal{L} \bigl( \covd_B \graylinear[\oeps][][][N] \bigr) \Big) \Big\|_{L_{t,x}^\infty} 
\lesssim T^{-\lossdel} K^{-2} \big\| \covd_B \graylinear[\oeps][][][N]  \big\|_{L_{t,x}^\infty} \\ 
&\lesssim T^{-\lossdel}  K^{-2} N^{1+\kappa} \| \graylinear[\oeps][][][N]\|_{\modsol^{-\kappa}_N} = T^{-\lossdel}  K^{-1} N^{\kappa} \times \frac{N}{K}\times \| \graylinear[\oeps][][][N]\|_{\modsol^{-\kappa}_N}. 
\end{aligned}
\end{equation}
Second, we use the $\liegroup$-equivariance of $\mathcal{L}$ and our Schauder estimate~\eqref{local:eq-covariant-Schauder-2} with the choice of parameters $(\alpha,\beta,m)=(-2,-1+2\kappa,1)$, which yield 
\begin{equation}\label{local:eq-trichotomy-r3}
\begin{aligned}
&\quad \Big\| \Duh_{\Y,B,K}^{F} \Big( \mathcal{L} \bigl( \covd_B \graylinear[\oeps][][][N] \bigr) \Big) \Big\|_{L_{t,x}^\infty} 
=  \Big\| \Duh_{\Y,B,K}^{F} \Big( \covd_B  \mathcal{L} \bigl( \graylinear[\oeps][][][N] \bigr) \Big) \Big\|_{L_{t,x}^\infty} \\ 
&\lesssim T^{-\lossdel} K^{\frac{2}{r}-2\kappa} \big\| \graylinear[\oeps][][][N]\big\|_{\Cs_t^{-\frac{1}{2}+\kappa}L_x^r}  \lesssim T^{-\lossdel} K^{\frac{2}{r}-2\kappa} N^{-1+3\kappa}  \| \graylinear[\oeps][][][N]\|_{\modsol^{-\kappa}_N}  \\ 
&= T^{-\lossdel} K^{-1} K^{\frac{2}{r}-2\kappa} N^{3\kappa} \times \frac{K}{N} \times \| \graylinear[\oeps][][][N]\|_{\modsol^{-\kappa}_N}. 
\end{aligned}
\end{equation}
Since $\frac{1}{r}$ is much smaller than $\kappa$, this completes the proof of \eqref{local:eq-trichotomy-q0}, and hence the proof of \eqref{local:eq-trichotomy-e1}. 

Next, we prove the estimate \eqref{local:eq-trichotomy-e3}. For this, it is sufficient to prove that 
\begin{equation}\label{local:eq-trichotomy-r1}
\Big\| \covd_B \Duh_{\Y,B,K} \Big( \mathcal{L} \bigl( \graylinear[\oeps][][][N] \bigr) \Big) \Big\|_{L_{t,x}^\infty} 
\lesssim T^{-\lossdel} K^{-1} N^{3\kappa} \min\Bigl( \frac{N}{K}  ,  \frac{K}{N} \Bigr)^{1-2\kappa} \| \graylinear[\oeps][][][N]\|_{\modsol^{-\kappa}_N}. 
\end{equation}
In the case $K\gg_{\theta^\prime} N$, we estimate 
\begin{align}
\Big\| \covd_B \Duh_{\Y,B,K} \Big( \mathcal{L} \bigl( \graylinear[\oeps][][][N] \bigr) \Big) \Big\|_{L_{t,x}^\infty} 
&\leq \Big\| \bigl( \covd_B \Duh_{\Y,B,K} - \Duh_{\Y,B,K}  \covd_B \bigr) \Big( \mathcal{L} \bigl( \graylinear[\oeps][][][N] \bigr) \Big) \Big\|_{L_{t,x}^\infty} 
\label{local:eq-trichotomy-r2a} \\ 
&+ \Big\| \Duh_{\Y,B,K} \Big( \covd_B  \mathcal{L} \bigl( \graylinear[\oeps][][][N] \bigr) \Big) \Big\|_{L_{t,x}^\infty}. 
\label{local:eq-trichotomy-r2b}
\end{align}
Using our commutator estimate (Lemma \ref{kernel:lem-commutator-derivative}) and $N\geq N_{\lo} \sim T^{-\frac{1}{2}}$, \eqref{local:eq-trichotomy-r2a} can be bounded by
\begin{equation*}
\eqref{local:eq-trichotomy-r2a} \lesssim T^{-\lossdel} T^{-\frac{1-2\kappa}{2}} K^{-(2-2\kappa)} \bigl\| \graylinear[\oeps][][][N]\bigr\|_{L_{t,x}^\infty}
\lesssim T^{-\lossdel} N^{1-2\kappa} K^{-(2-2\kappa)} N^{\kappa} \| \graylinear[\oeps][][][N] \|_{\sol^{-\kappa}_N}.
\end{equation*}
Using our earlier estimate \eqref{local:eq-trichotomy-r4}, except that $\Duh^F_{\Y,B,K}$ is replaced by $\Duh_{\Y,B,K}$, we can also bound \eqref{local:eq-trichotomy-r2b} by $T^{-\lossdel} K^{-2} N^{1+\kappa} \| \graylinear[\oeps][][][N] \|_{\sol^{-\kappa}_N}$. By combining the two estimates of \eqref{local:eq-trichotomy-r2a} and \eqref{local:eq-trichotomy-r2b}, we obtain~\eqref{local:eq-trichotomy-r1}. In the case $K\ll_{\theta^\prime} N$, we use our covariant Schauder estimate \eqref{local:eq-covariant-Schauder-2} with $(\alpha,\beta,m)=(-1,-1+2\kappa,0)$ and Remark \ref{local:rem-Schauder-Duh-without-F}, which imply that 
\begin{align*}
\quad \Big\| \covd_B \Duh_{\Y,B,K} \Big( \mathcal{L} \bigl(  \graylinear[\oeps][][][N] \bigr) \Big) \Big\|_{L_{t,x}^\infty} 
\leq \Big\|  \Duh_{\Y,B,K} \Big(  \mathcal{L} \bigl( \graylinear[\oeps][][][N] \bigr) \Big) \Big\|_{\sol^1_K} \lesssim T^{-\lossdel} K^{\frac{2}{r}-2\kappa} \big\| \graylinear[\oeps][][][N]\big\|_{\Cs_t^{-\frac{1}{2}+\kappa}L_x^r}.  
\end{align*}
After proceeding exactly as in \eqref{local:eq-trichotomy-r3}, we then obtain~\eqref{local:eq-trichotomy-r1}. 

Finally, we prove the estimate \eqref{local:eq-trichotomy-e2}. For this, we first use a dyadic decomposition to reduce to the case $\eta_1=-1+\reg_1-\kappa$. Indeed, assuming that \eqref{local:eq-trichotomy-e2} holds for $\eta_1=-1+\reg_1-\kappa$, we obtain for all $\eta_1^\prime \in [\eta_1,0]$ that 
\begin{align*}
\Big\|& t^{-\frac{\eta_1^\prime}{2}} \Duh_{\Y,B,K}^{F} \Big( \Qc \bigl( W_M , \covd_B \graylinear[\oeps][][][N] \bigr) \Big) \Big\|_{L_{t,x}^\infty([0,T])}\\ 
&= \sup_{k\geq 0}  \Big\| t^{-\frac{\eta_1^\prime}{2}} \Duh_{\Y,B,K}^{F} \Big( \Qc \bigl( W_M , \covd_B \graylinear[\oeps][][][N] \bigr) \Big) \Big\|_{L_{t,x}^\infty([2^{-(k+1)} T,2^{-k} T])} \\ 
&\lesssim \sup_{k\geq 0}  (2^{-k} T)^{\frac{\eta_1-\eta_1^\prime}{2}} \Big\| t^{-\frac{\eta_1}{2}} \Duh_{\Y,B,K}^{F} \Big( \Qc \bigl( W_M , \covd_B \graylinear[\oeps][][][N] \bigr) \Big) \Big\|_{L_{t,x}^\infty([0,2^{-k} T])}\\ 
&\lesssim \sup_{k\geq 0} \Bigl( (2^{-k} T)^{\frac{\eta_1-\eta_1^\prime}{2}} \times T^{\frac{1}{2}-\lambda}  K^{-1-2\theta} N^{-\kappa} \| \graylinear[\oeps][][][N]\|_{\modsol^{-\kappa}_N} \big\| W_M \big\|_{\sol_M^{1-\reg_1,\eta_1}([0,2^{-k}T])} \Bigr) \\ 
&\lesssim T^{\frac{1}{2}-\lambda}  K^{-1-2\theta} N^{-\kappa} \| \graylinear[\oeps][][][N]\|_{\modsol^{-\kappa}_N} \big\| W_M \big\|_{\sol_M^{1-\reg_1,\eta_1^\prime}([0,T])}.
\end{align*}
In the last line, we used that $\eta_1^\prime \geq \eta_1$ implies that 
\begin{equation*}
\big\| W_M \big\|_{\sol_M^{1-\reg_1,\eta_1}([0,2^{-k}T])}
\lesssim (2^{-k} T)^{\frac{\eta_1^\prime-\eta_1}{2}}\big\| W_M \big\|_{\sol_M^{1-\reg_1,\eta_1^\prime}([0,T])}.
\end{equation*}
For the rest of this proof, we now take $\eta_1=-1+\reg_1-\kappa$. Similarly to the argument for \eqref{local:eq-trichotomy-e1}, \eqref{local:eq-trichotomy-e2} can be reduced to the dyadic estimate 
\begin{equation}\label{local:eq-trichotomy-q1}
\begin{aligned}
\Big\| &t^{-\frac{\eta_1}{2}} \Duh_{\Y,B,K}^{F} \Big( \Qc \bigl( W_M , \covd_B \graylinear[\oeps][][][N] \bigr) \Big) \Big\|_{L_{t,x}^\infty} \\ 
&\lesssim T^{\frac{1}{2}-\lambda} K^{-1} N^{-\kappa} \min\Big( N^{2\kappa} \frac{N}{K}, (KN)^{2\kappa} \Bigl( \frac{K}{N} \Bigr)^\lambda + K^{-\lambda} \Big) 
\| \graylinear[\oeps][][][N]\|_{\modsol^{-\kappa}_N} \big\| W_M \big\|_{\sol^{1-\reg_1,\eta_1}_M}. 
\end{aligned}
\end{equation}
Indeed, if $N\ll_{\theta^\prime} K$, we can use the first argument in the $\min$-factor in \eqref{local:eq-trichotomy-q1} and, if $K\ll_{\theta^\prime} N$, we can use the second argument. In the latter case, it is important that $\lambda \theta^\prime$ is much larger than $\theta$. To obtain \eqref{local:eq-trichotomy-q1}, we first use our
weighted Schauder estimate \eqref{local:eq-weighted-covariant-Schauder-1}, which yields
\begin{align*}
\Big\| t^{-\frac{\eta_1}{2}} \Duh_{\Y,B,K}^{F} \Big( \Qc \bigl( W_M , \covd_B \graylinear[\oeps][][][N] \bigr) \Big) \Big\|_{L_{t,x}^\infty}
&\lesssim T^{-\lossdel} K^{-2} \big\| t^{-\frac{\eta_1}{2}} \Qc \bigl( W_M , \covd_B \graylinear[\oeps][][][N] \bigr) \big\|_{L_{t,x}^\infty} \\
&\lesssim T^{-\lossdel} K^{-2} M^{-1+\reg_1} N^{1+\kappa}  \| \graylinear[\oeps][][][N]\|_{\modsol^{-\kappa}_N}  \big\| W_M \big\|_{\sol^{1-\reg_1,\eta_1}_M}. 
\end{align*}
Since $M\geq N_{\lo}\sim T^{-\frac{1}{2}}$, this yields the first bound in \eqref{local:eq-trichotomy-q1}. In order to get the second bound, we split 
\begin{equation}\label{local:eq-trichotomy-p1}
\begin{aligned}
&\Duh_{\Y,B,K}^{F} \Big( \Qc \bigl( W_M , \covd_B \graylinear[\oeps][][][N] \bigr) \Big) \\ 
= &\Duh_{\Y,B,K}^{F} \Big( \covd_B  \Qc \bigl( W_M , \graylinear[\oeps][][][N] \bigr) \Big)
- \Duh_{\Y,B,K}^{F} \Big( \Qc \bigl( \covd_B W_M ,  \graylinear[\oeps][][][N] \bigr) \Big).
\end{aligned}
\end{equation}
To estimate the first summand in \eqref{local:eq-trichotomy-p1}, we let $\alpha \leq -\kappa$ remain to be chosen. Provided that $\alpha \geq \eta_1 + 2\kappa$, we can first estimate
\begin{align*}
\Big\| &t^{-\frac{\eta_1}{2}} \Duh_{\Y,B,K}^{F} \Big( \covd_B  \Qc \bigl( W_M , \graylinear[\oeps][][][N] \bigr) \Big)
\Big\|_{L_{t,x}^\infty} \\ 
&\leq T^{\frac{\alpha-2\kappa-\eta_1}{2}}
\Big\| t^{-\frac{\alpha-2\kappa}{2}} \Duh_{\Y,B,K}^{F} \Big( \covd_B  \Qc \bigl( W_M , \graylinear[\oeps][][][N] \bigr) \Big)
\Big\|_{L_{t,x}^\infty}.
\end{align*}
Provided that $(\alpha-2\kappa)+\alpha>-1$, which is the reason for not working directly with the weight $t^{-\frac{\eta_1}{2}}$, it then follows from our weighted Schauder estimate \eqref{local:eq-weighted-covariant-Schauder-2} that 
\begin{align*}
T^{\frac{\alpha-2\kappa-\eta_1}{2}}&
\Big\| t^{-\frac{\alpha-2\kappa}{2}} \Duh_{\Y,B,K}^{F} \Big( \covd_B  \Qc \bigl( W_M , \graylinear[\oeps][][][N] \bigr) \Big)
\Big\|_{L_{t,x}^\infty}\\
&\lesssim T^{\frac{\alpha-2\kappa-\eta_1}{2}} T^{-\lossdel} K^{-1-\alpha+\frac{2}{r}} 
\big\|t^{-\frac{\alpha-2\kappa}{2}} \Qc \bigl( W_M , \graylinear[\oeps][][][N] \bigr) \big\|_{\Cs_t^{\frac{\alpha}{2}}L_x^r}. 
\end{align*}
From Lemma~\ref{prelim:lem-product-temporal}, Lemma~\ref{local:lem-solution-spaces-hoelder}, and Lemma~\ref{local:lem-sol-different-regularities-weights}, it then follows that
\begin{align*}
&\big\|t^{-\frac{\alpha-2\kappa}{2}} \Qc \bigl( W_M , \graylinear[\oeps][][][N] \bigr) \big\|_{\Cs_t^{\frac{\alpha}{2}}L_x^r}
\lesssim  \big\| t^{-\frac{\alpha-2\kappa}{2}} W_M \big\|_{\Cs_t^{\frac{-\alpha+\kappa}{2}}L_x^\infty} \big\| \graylinear[\oeps][][][N] \big\|_{\Cs_t^{\frac{\alpha}{2}}L_x^r}\\
&\quad\lesssim N^{\alpha+\kappa} \big\| W_M \big\|_{\sol^{-\alpha+\kappa,\alpha-2\kappa}_M}
\| \graylinear[\oeps][][][N]\|_{\modsol^{-\kappa}_N} 
\lesssim N^{\alpha+\kappa} \big\| W_M \big\|_{\sol^{1-\reg_1,-1+\reg_1-\kappa}_M}
\| \graylinear[\oeps][][][N]\|_{\modsol^{-\kappa}_N}.
\end{align*}
In total, using the definition of $\eta_1$, we obtain an upper bound of the form
\begin{equation*}
T^{-\lossdel} T^{\frac{1+\alpha-\reg_1-\kappa}{2}} K^{-1} N^{-\kappa} \times K^{\frac{2}{r}} N^{2\kappa} \Big( \frac{K}{N} \Big)^{-\alpha} \times  \big\| W_M \big\|_{\sol^{1-\reg_1,-1+\reg_1-\kappa}_M}
\| \graylinear[\oeps][][][N]\|_{\modsol^{-\kappa}_N}. 
\end{equation*}
By choosing $\alpha=-\lambda$, which satisfies all the above conditions, we then obtain the contribution corresponding to the first summand in the second argument of the $\min$-factor in \eqref{local:eq-trichotomy-q1}. To estimate the second summand in \eqref{local:eq-trichotomy-p1}, we use Lemma~\ref{local:lem-high-low} with $(\alpha,\beta,\gamma,\nu,\eta)=(-2\kappa,1-\reg_1,-\frac{2\reg_1}{\theta},\eta_1,\eta_1)$. From this, it follows that\footnote{Actually, the final exponent of $T$ in this estimate is positive, but this does not make a significant difference.}
\begin{align*}
&\quad \Big\| t^{-\frac{\eta_1}{2}} \Duh_{\Y,B,K}^{F} \Big(   \Qc \bigl( \covd_B W_M , \graylinear[\oeps][][][N] \bigr) \Big)
\Big\|_{L_{t,x}^\infty}
\lesssim K^{-2+\frac{2\reg_1}{\theta}} \Big\|\Qc \bigl( \covd_B W_M , \graylinear[\oeps][][][N] \bigr)  \Big\|_{\nl^{-2\reg_1/\theta,\eta_1}_K} \\
&\lesssim 
T^{-\lossdel} T^{\frac{\reg_1}{\theta}-\frac{1-\theta}{2\theta}(\reg_1+2\kappa)} T^{-\frac{1}{r}-2\kappa} K^{-2+\frac{2\reg_1}{\theta}}
\big\| \graylinear[\oeps][][][N]\big\|_{\modsol^{-2\kappa}_N} \big\| W_M \big\|_{\sol^{1-\reg_1,\eta_1}_M}\\
&\lesssim T^{-\bigO(\reg_1/\theta)} K^{-2+\frac{2\reg_1}{\theta}} N^{-\kappa}\| \graylinear[\oeps][][][N]\|_{\modsol^{-\kappa}_N} \big\| W_M \big\|_{\sol^{1-\reg_1,\eta_1}_M}. 
\end{align*}
Since $K\geq N_{\lo} \sim T^{-\frac{1}{2}}$, we can further estimate
\begin{equation*}
K^{-2+\frac{2\reg_1}{\theta}} = K^{-1-\lambda} K^{-1+\lambda+\frac{2\reg_1}{\theta}} 
\lesssim K^{-1-\lambda} T^{\frac{1-\lambda}{2}-\frac{\reg_1}{\theta}}. 
\end{equation*}
Thus, since $\lambda$ is much larger than $\reg_1/\theta$, we obtain the contribution corresponding to the second summand in the second argument of the $\min$-factor in \eqref{local:eq-trichotomy-q1}. 
\end{proof}

In our next lemma, we show that for large values of $K$, $\Duh^{F}_{\Y,B,K}$ can be approximated using $\Duh_{\Y,B,K}$. This will be needed later in order to use our bounds on $\Xi_{\Beps}^{\para}$, which was defined in \eqref{local:eq-data-set-para}.

\begin{lemma}\label{object:lem-duhamel-approximation-no-F}
Let Assumption~\ref{local:assumption-B} be satisfied, let $K\geq N_{\lo}$, and let $-2<\eta\leq 0$. Then, one has for all $\phi \colon [0,T]\times \T^2 \rightarrow \Y$ that 
\begin{equation}
\Big\| t^{-\frac{\eta}{2}} \bigl( \Duh^{F}_{\Y,B,K} -  \Duh_{\Y,B,K} \bigr) (\phi) \Big\|_{L_{t,x}^\infty} \lesssim  T^{-\frac{1-\roughdataweight}{2}-\lossdel} K^{-3+\roughdataweight} \|t^{-\frac{\eta}{2}} \phi \|_{L_{t,x}^\infty}. 
\end{equation}
\end{lemma}

\begin{proof} Let $0\leq s \leq t \leq T$. By using Proposition \ref{prop:qB-ptwise-bds}, Remark \ref{remark:Fpert-less-time-weight}, and the definition of $\lossfactor_T(F)$ from \eqref{local:eq-lossfactor}, we obtain
\begin{align*}
\Big|& \fq_{\Y,B}(s,y;t,x) -   \fp_{\Y,B} (s,y;t,x) \Big| \\
&\lesssim \gronwall_T(C_0 F)  \Fnorm{F}_{[s,t],0} (t-s)^{-1} \exp\Bigl( - c \tfrac{|x-y|^2}{t-s} \Bigr) \\ 
&\lesssim \gronwall_T(C_0 F) \bigl( T^{\frac{1-\roughdataweight}{2}} 
\Fnorm{u^{-\frac{1-\roughdataweight}{2}} F(u)}_{[0,T],1} \bigr) (T^{-1} (t-s))^{\frac{1-\roughdataweight}{2}}(t-s)^{-1} \exp\Bigl( - c \tfrac{|x-y|^2}{t-s} \Bigr) \\ 
&\lesssim \lossfactor_T(F)  (T^{-1} (t-s))^{\frac{1-\roughdataweight}{2}}(t-s)^{-1} \exp\Bigl( - c \tfrac{|x-y|^2}{t-s} \Bigr). 
\end{align*}
By distinguishing the cases $t\gg K^{-2}$ and $t\lesssim K^{-2}$, one then easily obtains that 
\begin{equation*}
\Big\| \ind(|t-s| \sim K^{-2}) t^{-\frac{\eta}{2}} s^{\frac{\eta}{2}} \bigl( \fq_{\Y,B}(s,y;t,x) -  \fp_{\Y,B}(s,y;t,x)\bigr) \Big\|_{L_{s,y}^1} \lesssim \lossfactor_T(F) T^{-\frac{1-\roughdataweight}{2}} K^{-3+\roughdataweight}. 
\end{equation*}
Together with the bound from Assumption~\ref{local:assumption-B}, this implies the desired result.
\end{proof}

Equipped with Lemma~\ref{local:lem-trichotomy} and Lemma~\ref{object:lem-duhamel-approximation-no-F}, we now prove the following lemma, which yields information on the structure of the paracontrolled term $X$. 

\begin{lemma}[Structure of the paracontrolled term]\label{local:lem-X-structure}
Let Assumption~\ref{local:assumption-B} be satisfied. Furthermore, let $\Qc\colon \Y \times \Y \rightarrow \Y$ be a $\liegroup$-equivariant bilinear map, let $k\in [2]$, and let $K\geq N_{\lo}$. Finally, let $W_\star=(W_M)_{M\geq N_{\lo}}$, where $W_M\colon [0,T]\times \T^2 \rightarrow \Y$ is supported on times $t\gtrsim M^{-2}$, and let $-1+\reg_1-\kappa \leq \eta_1 \leq 0$. Then, it holds that 
\begin{equation}\label{local:eq-structure}
\begin{aligned}
&\quad \Big\| t^{-\frac{\eta_1}{2}} \Big( \Duh_{\Y,B,K}^{F} \Big( \Qc \bigl( W_\star \parall \covd^k \graylinear[\oeps][][][\ostar] \bigr) \Big) -  \Qc \Big( W, \covd^k \Duh_{\Y,B,K} \bigl( \graylinear[\oeps][][][] \bigr) \Big) \Big) \Big\|_{L_{t,x}^\infty} \\ 
&\lesssim T^{\frac{1-\lambda_1}{2}} K^{-1-2\theta} \big\| W_\star \big\|_{\sol^{1-\reg_1,\eta_1}} \| \graylinear[\oeps][][][\ostar]\|_{\modsol^{-\kappa}}. 
\end{aligned} 
\end{equation}
In addition, for all $M,N\geq N_{\lo}$ satisfying $M\ll_\theta N$, it holds that 
\begin{equation}\label{local:eq-structure-frequency-localized}
\begin{aligned}
&\quad \Big\| t^{-\frac{\eta_1}{2}} \Big( \Duh_{\Y,B,K}^{F} \Big( \Qc \bigl( W_M, \covd^k \graylinear[\oeps][][][N] \bigr) \Big) -  \Qc \Big( W_M, \covd^k \Duh_{\Y,B,K} \bigl( \graylinear[\oeps][][][N] \bigr) \Big) \Big) \Big\|_{L_{t,x}^\infty} \\ 
&\lesssim T^{\frac{1-\lambda_1}{2}} K^{-1-2\theta} N^{-\kappa} \big\| W_M \big\|_{\sol^{1-\reg_1,\eta_1}_M} \| \graylinear[\oeps][][][N]\|_{\modsol^{-\kappa}_N}. 
\end{aligned} 
\end{equation}
\end{lemma}

\begin{remark}
In this subsection, we will only need the estimate \eqref{local:eq-structure}. However, the frequency-localised estimate \eqref{local:eq-structure-frequency-localized} will be useful in the proof of Lemma \ref{objects:lem-quadratic-dlin-23} below. 
\end{remark}

\begin{proof}
Throughout the proof, when using estimates involving the loss-factor $\lossfactor_T(F)$, we always make use of our assumption $\lossfactor_T(F)\leq T^{-\lossdel}$. We first prove the frequency-localised estimate \eqref{local:eq-structure-frequency-localized}, for which we distinguish the cases $K\sim_{\theta^\prime} N$ and $K\not \sim_{\theta^\prime} N$. Afterwards, we will use \eqref{local:eq-structure-frequency-localized}, as well as some easier estimates, to prove \eqref{local:eq-structure}. 

\proofstep{The case $K\sim_{\theta^\prime} N$ in \eqref{local:eq-structure-frequency-localized}:} In this case, we split
\begin{align}
&\quad \Duh_{\Y,B,K}^{F} \Bigl( \Qc \bigl(  W_M ,  \covd^k \graylinear[\oeps][][][N] \bigr) \Bigr) -  \Qc \Big( W_M, \covd^k \Duh_{\Y,B,K} \bigl( \graylinear[\oeps][][][N]\bigr) \Big) \notag \\
&= \Bigl( \Duh^{F}_{\Y,B,K} -  \Duh_{\Y,B,K} \Bigr) \Big( \Qc \bigl( W_M , \covd^k \graylinear[\oeps][][][N] \bigr) \Big) \label{local:eq-structure-q4} \\
&+   \Duh_{\Y,B,K} \Big( \Qc \bigl(  W_M ,  \covd^k \graylinear[\oeps][][][N] \bigr) \Big)
-   \Qc \Big( W_M, \Duh_{\Y,B,K} \bigl(  \covd^k \graylinear[\oeps][][][N]  \bigr)\Big) \label{local:eq-structure-p4} \\
&+   \Qc \Big( W_M, \Duh_{\Y,B,K} \bigl(   \covd^k \graylinear[\oeps][][][N]  \bigr)\Big)
-   \Qc \Big( W_M, \covd^k \Duh_{\Y,B,K} \bigl( \graylinear[\oeps][][][N]\bigr) \Big).\label{local:eq-structure-p5}
\end{align}
Using Lemma~\ref{object:lem-duhamel-approximation-no-F}, the first term can be estimated by 
\begin{align*}
 \big\| t^{-\frac{\eta_1}{2}} \eqref{local:eq-structure-q4} \big\|_{L_{t,x}^\infty}
&\lesssim T^{-\lossdel} T^{-\frac{1}{2}} K^{-3+\bigO(\kappa)} \| t^{-\frac{\eta_1}{2}} W_M \|_{L_{t,x}^\infty} \big\| \covd \graylinear[\oeps][][][N] \big\|_{L_{t,x}^\infty} \\
&\lesssim T^{-\frac{1}{2}-\lossdel} K^{-3} M^{-1} N (KMN)^{\bigO(\reg_1)} 
\| W_M \|_{\sol^{1-\reg_1,\eta_1}_M}
\| \graylinear[\oeps][][][N]\|_{\modsol^{-\kappa}_N}. 
\end{align*}
Since $M\gtrsim T^{-\frac{1}{2}}$ and $K\sim_{\theta^\prime} N$, the prefactor can be bounded by $K^{-2+\bigO(\theta^\prime)} N^{-\kappa}$. Since $\theta^\prime \ll \lambda_1$ and $K  \gtrsim T^{-\frac{1}{2}}$, this clearly yields an acceptable contribution. Using Lemma~\ref{kernel:lem-commutator-integral} and Definition~\ref{local:def-solution-space}, the second term can be estimated by 
\begin{align*}
 \big\| t^{-\frac{\eta_1}{2}} \eqref{local:eq-structure-p4} \big\|_{L_{t,x}^\infty}
&\lesssim T^{-\lossdel} \bigl( K^{-4} \big\| t^{-\frac{\eta_1}{2}} (\partial_t - \Delta_B ) W_M \big\|_{L_{t,x}^\infty} + K^{-3} \big\| t^{-\frac{\eta_1}{2}} \covd W_M \big\|_{L_{t,x}^\infty}  \bigr)
\big\| \covd \graylinear[\oeps][][][N] \big\|_{L_{t,x}^\infty} \\
&\lesssim  T^{-\lossdel} \bigl( K^{-4} M^{1+\reg_1} + K^{-3} M^{\reg_1} \bigr) N^{1+\kappa} 
\big\| W_M \big\|_{\sol^{1-\reg_1,\eta_1}_M}
\| \graylinear[\oeps][][][N]\|_{\modsol^{-\kappa}_N}. 
\end{align*}
Due to our assumptions $M\ll_\theta N$ and $K\sim_{\theta^\prime} N$, it holds that $M\leq N \lesssim K^{1+\bigO(\theta^\prime)}$. The prefactor can therefore again be estimated by $K^{-2+\bigO(\theta^\prime)}N^{-\kappa}$, which is acceptable. Using Lemma~\ref{kernel:lem-commutator-derivative}, the third term can be estimated by 
\begin{align*}
 \big\| t^{-\frac{\eta_1}{2}} \eqref{local:eq-structure-p5} \big\|_{L_{t,x}^\infty}
&\lesssim \big\| t^{-\frac{\eta_1}{2}} W_M  \big\|_{L_{t,x}^\infty} \big\| \Duh_{\Y,B,K} \bigl(   \covd^k \graylinear[\oeps][][][N]  \bigr) - \covd^k \Duh_{\Y,B,K} \bigl( \graylinear[\oeps][][][N]\bigr) \big\|_{L_{t,x}^\infty} \\ 
&\lesssim  T^{-\lossdel} T^{-\frac{1}{2}} K^{-2+2\kappa} \big\| t^{-\frac{\eta_1}{2}} W_M  \big\|_{L_{t,x}^\infty} \big\| \graylinear[\oeps][][][N] \big\|_{L_{t,x}^\infty} \\
&\lesssim T^{-\lossdel} T^{-\frac{1}{2}} K^{-2+2\kappa} M^{-1+\reg_1} N^\kappa   \big\| W_M \big\|_{\sol^{1-\reg_1,\eta_1}_M} \| \graylinear[\oeps][][][N]\|_{\modsol^{-\kappa}_N}. 
\end{align*}
Since, as before, $M\leq N \leq K^{1+\bigO(\theta^\prime)}$ and $M\gtrsim  T^{-\frac{1}{2}}$, the prefactor can be estimated once more by $K^{-2+\bigO(\theta^\prime)} N^{-\kappa}$, which is acceptable. 

\proofstep{The case $K\not\sim_{\theta^\prime} N$ in \eqref{local:eq-structure-frequency-localized}:}
In this case, we primarily rely on Lemma~\ref{local:lem-trichotomy}. It implies that 
\begin{align*}
\Bigl\| t^{-\frac{\eta_1}{2}} \Duh_{\Y,B,K}^{F} \Big( \Qc \bigl( W_M, \covd^k \graylinear[\oeps][][][N] \bigr) \Big) \Bigr\|_{L_{t,x}^\infty} 
\lesssim T^{\frac{1}{2}-\lambda} K^{-1-2\theta} N^{-\kappa} \| W_M \|_{\sol^{1-\reg_1,\eta_1}_M} \| \graylinear[\oeps][][][N] \|_{\modsol^{-\kappa}_N} 
\end{align*}
as well as 
\begin{equation}\label{local:eq-structure-r0}
\begin{aligned}
 \Bigl \|& t^{-\frac{\eta_1}{2}} \Qc \Big( W_M, \covd^k \Duh_{\Y,B,K} \bigl( \graylinear[\oeps][][][N] \bigr) \Big)  \Bigr\|_{L_{t,x}^\infty} \\
 &\lesssim \| t^{-\frac{\eta_1}{2}} W_M \|_{L_{t,x}^\infty} \bigl\| \covd^k \Duh_{\Y,B,K} \bigl( \graylinear[\oeps][][][N] \bigr) \bigr\|_{L_{t,x}^\infty} \\ 
 &\lesssim M^{-(1-\reg_1)} K^{-1-2\theta} N^{-\kappa} 
 \| W_M \|_{\sol^{1-\reg_1,\eta_1}_M} \| \graylinear[\oeps][][][N] \|_{\modsol^{-\kappa}_N}. 
\end{aligned}
\end{equation}
Since $M\geq N_{\lo} \sim T^{-\frac{1}{2}}$, both contributions are acceptable.

\proofstep{The estimate \eqref{local:eq-structure}:} 
In order to prove \eqref{local:eq-structure}, we decompose 
\begin{align}
\Duh_{\Y,B,K}^{F}& \Big( \Qc \bigl( W_\star \parall \covd^k \graylinear[\oeps][][][\ostar] \bigr) \Big) -  \Qc \Big( W, \covd^k \Duh_{\Y,B,K} \bigl( \graylinear[\oeps][][][] \bigr) \Big) \notag \\ 
&= \sum_{\substack{M,N\colon \\ M\ll_\theta N}} \Bigl( \Duh_{\Y,B,K}^{F} \Big( \Qc \bigl( W_M, \covd^k \graylinear[\oeps][][][N] \bigr) \Big) -  \Qc \Big( W_M, \covd^k \Duh_{\Y,B,K} \bigl( \graylinear[\oeps][][][N] \bigr) \Big) \Bigr) \label{local:eq-structure-r1} \\ 
&-  \sum_{\substack{M,N\colon \\ M\gtrsim_\theta N, \\ K \sim_{\theta^\prime} N}} \Qc \Big( W_M, \covd^k \Duh_{\Y,B,K} \bigl( \graylinear[\oeps][][][N] \bigr) \Big) \label{local:eq-structure-r2} \\
&- \sum_{\substack{M,N\colon \\ M\gtrsim_\theta N, \\ K \not \sim_{\theta^\prime} N}} \Qc \Big( W_M, \covd^k \Duh_{\Y,B,K} \bigl( \graylinear[\oeps][][][N] \bigr) \Big).  \label{local:eq-structure-r3}
\end{align}
We remark that $M\gtrsim_\theta N$, which was not defined in or around \eqref{local:eq-ll-theta}, simply means that $M\gg_\theta N$ or $M\sim_\theta N$.
The first term \eqref{local:eq-structure-r1} can easily be bounded using our frequency-localised estimate \eqref{local:eq-structure-frequency-localized}. Using our Schauder estimate \eqref{local:eq-covariant-Schauder-1}, the summands in the second term \eqref{local:eq-structure-r2} can be bounded by 
\begin{align*}
\Bigl\|& t^{-\frac{\eta_1}{2}} \Qc \Big( W_M, \covd^k \Duh_{\Y,B,K} \bigl( \graylinear[\oeps][][][N] \bigr) \Big) \Bigr\|_{L_{t,x}^\infty}  
\lesssim \| t^{-\frac{\eta_1}{2}} W_M \|_{L_{t,x}^\infty} \bigl\| \covd^k \Duh_{\Y,B,K} \bigl( \graylinear[\oeps][][][N] \bigr) \bigr\|_{L_{t,x}^\infty} \\ 
&\lesssim \| t^{-\frac{\eta_1}{2}} W_M \|_{L_{t,x}^\infty} 
T^{-\lossdel} K^{-1} \| \graylinear[\oeps][][][N] \|_{L_{t,x}^\infty}   
\lesssim T^{-\lossdel} M^{-(1-\reg_1)} K^{-1} N^{\kappa} \| W_\star \|_{\sol^{1-\reg_1,\eta_1}} \| \graylinear[\oeps][][][\ostar] \|_{\modsol^{-\kappa}}. 
\end{align*}
Since $K\geq N_{\lo}\sim T^{-\frac{1}{2}}$ and 
\begin{align*}
\sum_{\substack{M,N\colon \\ M\gtrsim_\theta N, \\ K \sim_{\theta^\prime} N}}
M^{-(1-\reg_1)} K^{-1} N^{\kappa} 
\lesssim 
\sum_{\substack{N\colon \\ K \sim_{\theta^\prime} N}} N^{-(1-\theta) (1-\reg_1)+\kappa} K^{-1} 
\lesssim K^{-(1-\theta^\prime) ((1-\theta)(1-\reg_1)-\kappa)} K^{-1},
\end{align*}
this shows that the contribution of \eqref{local:eq-structure-r2} is of size $\sim K^{-2+\bigO(\theta^\prime)}$, which is acceptable. Finally, using a similar argument as in \eqref{local:eq-structure-r0}, the third summand \eqref{local:eq-structure-r3} can be bounded by 
\begin{align*}
\bigl\| t^{-\frac{\eta_1}{2}} \eqref{local:eq-structure-r3} \bigr\|_{L_{t,x}^\infty}
&\lesssim 
 \sum_{\substack{M,N\colon \\ M\gtrsim_\theta N, \\ K \not \sim_{\theta^\prime} N}}
 M^{-(1-\reg_1)} K^{-1-2\theta} N^{-\kappa} 
 \| W_M \|_{\sol^{1-\reg_1,\eta_1}_M} \| \graylinear[\oeps][][][N] \|_{\modsol^{-\kappa}_N} \\ 
 &\lesssim T^{\frac{1-\reg_1}{2}} K^{-1-2\theta} \| W_\star \|_{\sol^{1-\reg_1,\eta_1}} \| \graylinear[\oeps][][][\ostar] \|_{\modsol^{-\kappa}},
\end{align*}
which is clearly acceptable. 
\end{proof}

Equipped with Lemma~\ref{local:lem-X-structure}, we are now ready to prove the main estimates of this subsection. 

\begin{proof}[Proof of Proposition~\ref{local:prop-paracontrolled}] 
To simplify the notation, we set $\eta_1\coloneq -1+\reg_1-\kappa$. By using the definition of the high$\times$high-product and by inserting a term of the same form as in Lemma~\ref{local:lem-X-structure}, we can rewrite the argument in \eqref{local:eq-paracontrolled-main-estimate} as 
\begin{align}
&\quad  \Duh^{F}_{\Y,B,\star} \bigl( 
\Qc \bigl( W_\diamond \parall \covd^k \graylinear[\oeps][][][\odiamond] \bigr)\bigr) \parasim \covd^\ell \graylinear[\oeps][][][\ostar]  - \tfrac{1}{4}   \delta^{k\ell} (\sigma_{\rho,\eps}^2 + \Crho) \mathcal{K}(W) \notag \\ 
&= \sum_{\substack{K,L\colon \\ K \sim_\theta L}} \hspace{-1ex}
 \Bigl( \Duh^{F}_{\Y,B,K} \bigl( 
\Qc \bigl( W_\diamond \parall \covd^k \graylinear[\oeps][][][\odiamond] \bigr)\bigr) -  \Qc \bigl( W, \covd^k \Duh_{\Y,B,K}(\graylinear[\oeps][][][])\bigr) \Bigr) \otimes \covd^\ell \graylinear[\oeps][][][L]  \label{local:eq-paracontrolled-main-p1} \\ 
&+ \sum_{\substack{K,L\colon \\ K \sim_\theta L}} 
\hspace{-1ex}  \Qc \bigl( W, \covd^k \Duh_{\Y,B,K}(\graylinear[\oeps][][][])\bigr) \otimes \covd^\ell \graylinear[\oeps][][][L]  - \frac{1}{4}   \delta^{k\ell} (\sigma_{\rho,\eps}^2 + \Crho)  \mathcal{K}(W). 
\label{local:eq-paracontrolled-main-p2}
\end{align}
Using Lemma~\ref{local:lem-sol-different-regularities-weights} and our weighted Schauder estimate \eqref{local:eq-weighted-covariant-Schauder-1}, as well as the parameter identity $\reg_1=\reg_2$, we obtain that
\begin{align*}
&\hspace{2.5ex} 
\| \eqref{local:eq-paracontrolled-main-p1} \|_{\nl^{-1+\reg_2,-\lambda}} 
\lesssim \| \eqref{local:eq-paracontrolled-main-p1} \|_{\nl^{0,-1+\reg_1-\lambda}} 
\lesssim T^{\frac{\lambda-\kappa}{2}-\lossdel} \| \eqref{local:eq-paracontrolled-main-p1} \|_{\nl^{0,\eta_1}} 
\\  
&\lesssim T^{\frac{\lambda-\kappa}{2}-\lossdel} \sum_{\substack{K,L\colon \\ K \sim_\theta L }}   \big\| \covd^\ell \graylinear[\oeps][][][L] \big\|_{L_{t,x}^\infty} 
\bigl\| t^{-\frac{\eta_1}{2}} \bigl( \Duh^{F}_{\Y,B,K} \bigl( 
\Qc \bigl( W_\diamond \parall \covd^k \graylinear[\oeps][][][\odiamond] \bigr)\bigr) - \Qc \bigl( W, \covd^k \Duh_{\Y,B,K}(\graylinear[\oeps][][][])\bigr) \bigr) \bigr\|_{L_{t,x}^\infty}.   
\end{align*}
We can now estimate the commutator term using Lemma~\ref{local:lem-X-structure} which, together with the definition of $\tripn{ \Xi_{\Beps} }$, allows us to control the above expression by 
\begin{equation*}
 T^{\frac{\lambda-\kappa}{2}-\lossdel}  T^{\frac{1-\lambda_1}{2}} \tripn{ \Xi_{\Beps} }^{2} \| W_\star \|_{\sol^{1-\reg_1,\eta_1}} \sum_{\substack{K,L\colon \\ K \sim_\theta L }} K^{-1-2\theta} L^{1+\kappa}. 
\end{equation*}
Since the exponent of $K$ is $-1-2\theta$, and not just $-1-\theta$, the dyadic sum can easily be bounded by a constant depending on $\kappa$ and $\theta$. To treat the second term, we first write \eqref{local:eq-paracontrolled-main-p2} as 
\begin{equation*}
\boldsymbol{L} \Big( W \otimes \Big( \sum_{\substack{K,L\colon \\ K \sim_\theta L}} \covd^k \Duh_{\Y,B,K}(\graylinear[\oeps][][][]) \otimes \covd^\ell \graylinear[\oeps][][][L] -  \tfrac{1}{4} \delta^{k\ell} (\sigma_{\rho,\eps}^2 + \Crho) \Cov_{\Y} \Big) \Big).
\end{equation*}
Here, $\boldsymbol{L}\colon \Y \otimes \Y \otimes \Y \rightarrow \Y \otimes \Y$ is the linear map defined by $\boldsymbol{L}(X \otimes Y \otimes Z)\coloneq  \Qc (X,Y) \otimes Z$.  
In order to estimate \eqref{local:eq-paracontrolled-main-p2}, we first use Lemma~\ref{local:lem-sol-different-regularities-weights} and then use our weighted covariant Schauder estimate \eqref{local:eq-weighted-covariant-Schauder-2} with 
\begin{equation*}
(\alpha,\beta,\nu,\eta) \coloneq \big(-3\theta, -2\theta, -1+\reg_2+3\theta-\lambda,\eta_1\big). 
\end{equation*}
Since $\theta$ is much larger than $\frac{1}{r}$, the $N^{\alpha-\beta+\frac{2}{r}}$-factor resulting from \eqref{local:eq-weighted-covariant-Schauder-2} is then $N^{-3\theta+2\theta+\frac{2}{r}}\lesssim 1$. Furthermore, since $\lambda$ is much larger than $\reg_2,\kappa$, and $\theta$, our choice of $\nu$ is smaller than $\eta$.  Hence, the $\max(N_{\lo}^{\nu-\eta},N^{\nu-\eta})$-factor in \eqref{local:eq-weighted-covariant-Schauder-2} contributes $T^{-\frac{(\nu-\eta)}{2}}=T^{\frac{\lambda}{2}-\bigO(\theta)}$. All in all, we therefore obtain 
\begin{align*}
&\quad \| \eqref{local:eq-paracontrolled-main-p2} \|_{\nl^{-1+\reg_2,-\lambda}} 
\lesssim \| \eqref{local:eq-paracontrolled-main-p2} \|_{\nl^{-3\theta,-1+\reg_2+3\theta-\lambda}} 
 \lesssim T^{\frac{\lambda}{2}-\bigO(\theta)} \| t^{-\frac{\eta_1}{2}} \eqref{local:eq-paracontrolled-main-p2}\|_{\Cs_t^{-\theta}L_x^r}. 
 \end{align*}
Using product estimates (Lemma~\ref{prelim:lem-product-temporal}), our Hölder estimate (Lemma~\ref{local:lem-solution-spaces-hoelder}), and our bounds on $\Xi_{\Beps}$, we then obtain that 
\begin{align*}
&T^{\frac{\lambda}{2}-\bigO(\theta)} \| t^{-\frac{\eta_1}{2}} \eqref{local:eq-paracontrolled-main-p2}\|_{\Cs_t^{-\theta}L_x^r} \\
&\lesssim  T^{\frac{\lambda}{2}-\bigO(\theta)} 
\| t^{-\frac{\eta_1}{2}} W \|_{\Cs_t^{2\theta}L_x^\infty} \\
&\qquad \times 
\Big\| \sum_{\substack{K,L\colon \\ K \sim_\theta L}} \covd^k \Duh_{\Y,B,K}(\graylinear[\oeps][][][]) \otimes \covd^\ell \graylinear[\oeps][][][L] -  \tfrac{1}{4} \delta^{k\ell} (\sigma_{\rho,\eps}^2 + \Crho) \Cov_{\Y} \Big\|_{\Cs_t^{-\theta}L_x^r} \\ 
&\lesssim  T^{\frac{\lambda}{2}-\bigO(\theta)} \| W_\star \|_{\sol^{5\theta,\eta_1}} \tripn{ \Xi_{\Beps} }^2 
\lesssim T^{\frac{1+\lambda}{2}-\bigO(\theta)} \| W_\star\|_{\sol^{1-\reg_1,\eta_1}} \tripn{ \Xi_{\Beps} }^2. 
\end{align*}
Thus, the contributions of both \eqref{local:eq-paracontrolled-main-p1} and \eqref{local:eq-paracontrolled-main-p2} to \eqref{local:eq-paracontrolled-main-estimate} are acceptable. 
\end{proof}

In the following corollary, we use Proposition~\ref{local:prop-paracontrolled} to control the specific terms appearing in our evolution equations for $Y_{\eps,\star}$ and $\phiy_{\eps,\star}$. 

\begin{corollary}\label{local:cor-paracontrolled-specific}
Let Assumptions~\ref{local:assumption-B},~\ref{local:assumption-higgs-initial}, and~\ref{local:assumption-nonlinear-setup} be satisfied. Then, it holds that 
\begin{align}
&\Big\| \Lb\bigl(\Xop_{\eps,\star}(\Wrem_{\eps,\diamond},\phiw_{\eps,\diamond}),\covd \linear[\oeps][][][\ostar]\bigr) 
+(\sigma_{\rho,\eps}^2+\Crho)  \Cov_{\frkg}^{\frkg} (W_\eps) \Big\|_{\nl^{-1+\reg_2,-\lambda}} 
\label{local:eq-paracontrolled-specific-1}\\
&+\Big\| 2 \bil \bigl( \phixop_{\eps,\star}(\Wrem_{\eps,\diamond}) \parasim \covd \philinear[\oeps][][][\ostar]\bigr) 
+(\sigma_{\rho,\eps}^2+\Crho)  \Cov^{\frkg}_{\V}(W_\eps) \Big\|_{\nl^{-1+\reg_2,-\lambda}} 
\label{local:eq-paracontrolled-specific-2}\\ 
&+\Big\|  2 \Xop^j_{\eps,\star}(\Wrem_{\eps,\diamond},\phiw_{\eps,\diamond}) \parasim \covd_j \philinear[\oeps][][][\ostar]
+ (\sigma_{\rho,\eps}^2 + \Crho) \Cov^{\V}_{\frkg}( \phiw_\eps) \Big\|_{\nl^{-1+\reg_2,-\lambda}} 
\label{local:eq-paracontrolled-specific-3}\\
&+\Big\| - 2 \bigl( \covd_j \linear[\oeps][j][][\ostar] \bigr) \parasim    \phixop_{\eps,\star}(\Wrem_{\eps,\diamond}) \Big\|_{\nl^{-1+\reg_2,-\lambda}} 
\label{local:eq-paracontrolled-specific-4} \\ 
&\lesssim T^{\frac{1}{2}-\lambda_1}\bigl( \higgsdetcombined + T^{\frac{1}{2}} \| (R_{\eps,\star},\phir_{\eps,\star})\|_{\sol^{1-\reg_1,-\lambda}}\bigr)
\tripn{ \Xi_{\Beps} }^2. 
\end{align}
\end{corollary}

\begin{proof}
We first observe that the terms in \eqref{local:eq-paracontrolled-specific-1}--\eqref{local:eq-paracontrolled-specific-4} are all of the same form as the abstract term in Proposition~\ref{local:prop-paracontrolled}. Using Lemma~\ref{local:lem-sol-different-regularities-weights} and Proposition~\ref{local:prop-higgsdet-sdet}, we have that
\begin{equation*}
\| (\Wrem_{\eps,\star}, \phiw_{\eps,\star}) \|_{\sol^{1-\reg_1,-1+\reg_1-\kappa}} 
\lesssim T^{-\bigO(\lambda)} \bigl(  \higgsdetcombined + T^{\frac{1}{2}} \| (R_{\eps,\star},\phir_{\eps,\star})\|_{\sol^{1-\reg_1,-\lambda}}\bigr). 
\end{equation*}
Thus, due to Proposition~\ref{local:prop-paracontrolled}, it only remains to verify that the renormalisation terms in \eqref{local:eq-paracontrolled-specific-1}--\eqref{local:eq-paracontrolled-specific-4} agree with the renormalisation terms resulting from the corresponding applications of 
\eqref{local:eq-paracontrolled-main-estimate}. For this, we discuss each term separately.

\proofstep{The term in \eqref{local:eq-paracontrolled-specific-1}:}  We first split 
\begin{equation*}
\Lb\bigl(\Xop_{\eps,\star}(\Wrem_{\eps,\diamond},\phiw_{\eps,\diamond}),\covd \linear[\oeps][][][\ostar]\bigr) 
= \Lb\bigl(\Xwhiteop_{\eps,\star}(\Wrem_{\eps,\diamond}),\covd \linear[\oeps][][][\ostar]\bigr) 
+ \Lb\bigl(\Xblackop_{\eps,\star}(\phiw_{\eps,\diamond}),\covd \linear[\oeps][][][\ostar]\bigr). 
\end{equation*}
Since $\Xblackop_{\eps,\star}$ is paracontrolled using $\philinear[\oeps][][][\ostar]$, the second summand above does not produce a renormalisation term. Using \eqref{local:eq-Lie-bracket-term} and \eqref{local:eq-Xwhiteop-K}, we write the first summand above as 
\begin{align*}
& \Lb\bigl(\Xwhiteop_{\eps,\star}(\Wrem_{\eps,\diamond}),\covd \linear[\oeps][][][\ostar]\bigr)^i \\
=&   2  \big[    \Xwhiteop_{\eps,\star,j}(\Wrem_{\eps,\diamond})  \parasim \bigl(  \covd^j \linear[\oeps][i][][\ostar]- \covd^i \linear[\oeps][j][][\ostar] \bigr) \big]
+ 2 \big[ \Xwhiteindexop{i}_{\eps,\star}(\Wrem_{\eps,\diamond})\parasim  \covd^j \linear[\oeps][][j][\ostar] \big] \\  
=&  2 \Big[ \Duh_{B,\star,jk}^F \Big(  \big[  \Wrem_{\eps,\diamond,\ell} \parall \bigl( 2 \covd^\ell \linear[\oeps][k][][\odiamond]- \covd^k \linear[\oeps][\ell][][\odiamond] \bigr) \big] \Big) 
\parasim \bigl(  \covd^j \linear[\oeps][i][][\ostar]- \covd^i \linear[\oeps][j][][\ostar] \bigr) \Big] \\ 
+& 2 \Big[ \Duh_{B,\star}^{F,ik} \Big(  \big[  \Wrem_{\eps,\diamond,\ell} \parall \bigl( 2 \covd^\ell \linear[\oeps][][k][\odiamond]- \covd_k \linear[\oeps][\ell][][\odiamond] \bigr) \big] \Big)
\parasim  \covd^j \linear[\oeps][][j][\ostar] \Big]. 
\end{align*}
Here, as in Remark~\ref{local:rem-paracontrolled-indices}, we are using $\Duh_{B,\star,jk}^F$ to denote the matrix element of $\Duh_{B,\star}^F$ in the $j$-th row and $k$-th column. According to Proposition~\ref{local:prop-paracontrolled}, we need to subtract the renormalisation term
\begin{align*}
&2 \times \tfrac{1}{4} (\sigma_{\rho,\eps}^2+\Crho) 
\Big( \delta_{jk} \bigl(2\delta^{j\ell}\delta^{ik}-2\delta^{\ell i} \delta^{kj} - \delta^{kj}\delta^{\ell i} + \delta^{ik} \delta^{j \ell} \bigr) + \delta^{ik} \bigl( 2 \delta^{\ell j}\delta_{kj} - \delta^j_k \delta^{\ell}_j\bigr)\Big) 
\Cov^{\frkg}_{\frkg}(\Wrem_{\eps,\ell}) \\
&= - (\sigma_{\rho,\eps}^2 + \Crho) \Cov^{\frkg}_{\frkg}(\Wrem_{\eps}^i),
\end{align*}
which is consistent with \eqref{local:eq-paracontrolled-specific-1}.

\proofstep{The term in \eqref{local:eq-paracontrolled-specific-2}:} Using \eqref{local:eq-phixop-K}, we have that 
\begin{align*}
2 \bil \bigl( \phixop_{\eps,\star}(\Wrem_{\eps,\diamond}) \parasim \covd^k \philinear[\oeps][][][\ostar]\bigr) 
= 4 \bil \Big( \Duh_{B,\star} \bigl( \Wrem_{\eps,\diamond}^\ell \parall \covd_\ell \philinear[\oeps][][][\odiamond]\bigr) \parasim \covd^k \philinear[\oeps][][][\ostar] \Big).
\end{align*}
Due to Proposition~\ref{local:prop-paracontrolled}, we need to subtract the renormalisation term
\begin{equation*}
4 \times \tfrac{1}{4} \delta^k_\ell (\sigma_{\rho,\eps}^2+\Crho) \bil(W_\eps^\ell \ONBV_m \otimes \ONBV^m),
\end{equation*}
where $(\ONBV_m)_m$ is an orthonormal basis of $\V$. Due to the skew-symmetry of $\bil$ and the definition of $\Cov^{\frkg}_{\V}$, this is consistent with \eqref{local:eq-paracontrolled-specific-2}.

\proofstep{The term in \eqref{local:eq-paracontrolled-specific-3}:} We first decompose 
\begin{equation*}
 2 \Xop^{j}_{\eps,\star}(\Wrem_{\eps,\diamond},\phiw_{\eps,\diamond}) \parasim \covd_j \philinear[\oeps][][][\ostar]
 =  2 \Xwhiteindexop{j}_{\eps,\star}(\Wrem_{\eps,\diamond}) \parasim \covd_j \philinear[\oeps][][][\ostar]
 +  2 \Xblackindexop{j}_{\eps,\star}(\phiw_{\eps,\diamond}) \parasim \covd_j \philinear[\oeps][][][\ostar].
\end{equation*}
Since $\Xwhiteop_{\eps,\star}$ is paracontrolled using $\linear[\oeps][][][\ostar]$, the first summand above does not produce a renormalisation term. Using \eqref{local:eq-Xblackop-K}, the second summand above can be written as 
\begin{equation*}
2 \Big( \Duh_{B,\star}^F \big( \bil \bigl( \phiw_{\eps,\diamond} \parall \covd \philinear[\oeps][][][\odiamond]\bigr)\big) \Big)^j \parasim \covd_j \philinear[\oeps][][][\ostar].
\end{equation*}
Due to Proposition~\ref{local:prop-paracontrolled} (see also Remark \ref{local:rem-paracontrolled-indices}), we need to subtract the renormalisation term
\begin{equation*}
2 \times \tfrac{1}{4} \delta^j_j (\sigma_{\rho,\eps}^2+ \Crho) \bil(\phiw_\eps \otimes \ONBV_m) \ONBV^m. 
\end{equation*}
After using Lemma~\ref{prelim:lem-identity-Casrhog}, we see that this is consistent with \eqref{local:eq-paracontrolled-specific-3}.

\proofstep{The term in \eqref{local:eq-paracontrolled-specific-4}:} Since $\phixop_{\eps,\star}$ is paracontrolled using only $\philinear[\oeps][][][\ostar]$, the term $ - 2 ( \covd_j \linear[\oeps][j][][\ostar]) \parasim    \phixop_{\eps,\star}(\Wrem_{\eps,\diamond})$ does not produce a renormalisation term. 
\end{proof}

\subsection{Proof of covariant local well-posedness}\label{section:local-proof}

Equipped with the estimates from Sections~\ref{section:sdet-higgsdet-sol}--\ref{section:local-paracontrolled}, we are now ready to first prove the main result of this section (Proposition \ref{local:prop-main}), and then our covariant local well-posedness theorem (Theorem \ref{main:thm-covariant-local}). 

\begin{proof}[Proof of Proposition~\ref{local:prop-main}]
The proof is based on a contraction-mapping argument. For this, we first define
\begin{equation}\label{local:eq-ball}
\begin{aligned}
\mathbb{B} \coloneq \Big\{ \bigl(& \Xwhite_{\eps,\star}, \Xblack_{\eps,\star}, \phix_{\eps,\star}, Y_{\eps,\star}, \phiy_{\eps,\star} \bigr)\colon 
\max\Big(
\| \Xwhite_{\eps,\star} \|_{\sol^{1-\reg_1,0}},
\| \Xblack_{\eps,\star} \|_{\sol^{1-\reg_1,0}}, \\ 
&\| \phix_{\eps,\star} \|_{\sol^{1-\reg_1,0}}, 
\| Y_{\eps,\star} \|_{\sol^{1+\reg_2,-\lambda}}, 
\| \phiy_{\eps,\star} \|_{\sol^{1+\reg_2,-\lambda}} 
\Big) \leq T^{-\frac{1}{\gamma_0}-\lambda_2} \\ 
&\textup{and}\,\, \bigl( \Xwhite_{\eps,N}, \Xblack_{\eps,N}, \phix_{\eps,N}, Y_{\eps,N}, \phiy_{\eps,N} \bigr)(t)=0 \quad \text{for all}\, N\geq N_{\lo} \,\, \text{and} \,\,0\leq t\leq \tfrac{1}{16} N^{-2} \Big\}. 
\end{aligned}
\end{equation}
For notational convenience, we also write 
\begin{equation*}
\mathcal{X}_{\eps,\star} \eqdef \bigl( \Xwhite_{\eps,\star}, \Xblack_{\eps,\star}, \phix_{\eps,\star} \bigr) 
\quad \text{and} \quad \mathcal{Y}_{\eps,\star} \eqdef \bigl( Y_{\eps,\star}, \phiy_{\eps,\star} \bigr). 
\end{equation*}
We further define a map  
$\Gamma_\eps ( \mathcal{X}_{\eps,\star}, \mathcal{Y}_{\eps,\star}) 
= ( \Gamma_{\eps,\mathcal{X}}, \Gamma_{\eps,\mathcal{Y}})( \mathcal{X}_{\eps,\star}, \mathcal{Y}_{\eps,\star})$ by 
\begin{align*}
\Gamma_{\eps,\mathcal{X}} \eqdef \Big(  \eqref{local:eq-Xwhite-K-final}, \, \eqref{local:eq-Xblack-K-final},\,  
\eqref{local:eq-phix-K-final} \Big) \quad \text{and} \quad 
\Gamma_{\eps,\mathcal{Y}} \eqdef \Big( \eqref{local:eq-Y-final-new-1}--\eqref{local:eq-Y-final-new-10}, \, \eqref{local:eq-y-final-new-1}--\eqref{local:eq-y-final-new-last} \Big).
\end{align*}
In the individual terms in \eqref{local:eq-Xwhite-K-final}--\eqref{local:eq-y-final-new-last} and the estimates from Sections~\ref{section:sdet-higgsdet-sol}--\ref{section:local-paracontrolled}, which will be used momentarily, we use that 
\begin{alignat*}{3}
\Wrem_{\eps,\star} &= \Sdet_\star + \quadratic[\oeps][][][\ostar] + \phiquadratic[\oeps][][][\ostar] + \Xwhite_{\eps,\star} + \Xblack_{\eps,\star} + Y_{\eps,\star}, \quad \quad  & 
\phiw_{\eps,\star} &= \higgsdet_\star+ \mixedquadratic[\oeps][][][\ostar] + \phix_{\eps,\star} + \phiy_{\eps,\star}, \\
R_{\eps,\star} &= \quadratic[\oeps][][][\ostar] + \phiquadratic[\oeps][][][\ostar] + \Xwhite_{\eps,\star} + \Xblack_{\eps,\star} + Y_{\eps,\star}, \quad \quad & \phir_{\eps,\star} &= \mixedquadratic[\oeps][][][\ostar] + \phix_{\eps,\star} + \phiy_{\eps,\star}. 
\end{alignat*}
That is, we view every term as a map depending on the given objects $\Sdet$, $\higgsdet$, and $\Xi_{\Beps}$ as well as the five unknowns from \eqref{local:eq-ball}. As above, we write $\mathcal{R}_{\eps,\star}=(R_{\eps,\star},\phir_{\eps,\star})$. 
In the following, we only prove that $\Gamma_\eps$ is a self-map, i.e., that $\Gamma_\eps$ maps $\mathbb{B}$ into $\mathbb{B}$. 
Using standard (but tedious) modifications of our arguments, it is possible to show that $\Gamma_\eps$ is a contraction and that, as long as the constraint $\tripn{\Xi_{\Beps}}\leq T^{-\lambda_1}$ from \eqref{e:enhanced_data_assumption} is satisfied, $\Gamma_{\eps}$ depends Lipschitz-continuously on $\Xi_{\Beps}$ (with respect to the distance from Definition~\ref{local:def-data-set-norm}). To see that $\Gamma_\eps$ is a self-map, it suffices to prove that
\begin{align}
&\quad \max\Big( 
\| \Gamma_{\eps,\mathcal{X}} \|_{\sol^{1-\reg_1,0}},
\| \Gamma_{\eps,\mathcal{Y}}  \|_{\sol^{1+\reg_2,-\lambda}} 
\Big) \notag \\ 
&\lesssim T^{\frac{1}{2}-\lambda_1} \tripn{ \Xi_{\Beps} }^3 \notag \\ 
&+ T^{-\lambda_1}\tripn{ \Xi_{\Beps} }  
\bigl( \higgsdetcombined + T^{\frac{1}{2}} \| \mathcal{R}_{\eps,\star} \|_{\sol^{1-\reg_1,-\lambda}} +T^{\frac{1}{2}} \| \mathcal{Y}_{\eps,\star} \|_{\sol^{1+\reg_2,-\lambda}} \bigr) \notag \\ 
&+T^{\frac{1}{2}-\lambda_1}\| \mathcal{R}_{\eps,\star} \|_{\sol^{1-\reg_1,-\lambda}}  
\bigl( \higgsdetcombined + T^{\frac{1}{2}} \| \mathcal{R}_{\eps,\star} \|_{\sol^{1-\reg_1,-\lambda}} \bigr) \notag \\ 
&+ T^{\frac{1}{2}-\lambda_1} \bigl( \tripn{ \Xi_{\Beps} } + T^{\frac{1}{2}} \| \mathcal{R}_{\eps,\star} \|_{\sol^{1-\reg_1,-\lambda}}  \bigr) 
\bigl( \higgsdetcombined + \tripn{ \Xi_{\Beps} } + T^{\frac{1}{2}} \| \mathcal{R}_{\eps,\star} \|_{\sol^{1-\reg_1,-\lambda}} \bigr)^2 \notag \\
&+ T^{\frac{1}{2}-\lambda_1} \bigl( \higgsdetcombined + \tripn{ \Xi_{\Beps} } + T^{\frac{1}{2}} \| \mathcal{R}_{\eps,\star} \|_{\sol^{1-\reg_1,-\lambda}} \bigr) \notag \\ 
&+ T^{\frac{1}{2}-\lambda_1} \tripn{ \Xi_{\Beps} }^2 \bigl( \higgsdetcombined + T^{\frac{1}{2}} \| \mathcal{R}_{\eps,\star} \|_{\sol^{1-\reg_1,-\lambda}} \bigr). \label{local:eq-main-self-map}
\end{align}
Indeed, once \eqref{local:eq-main-self-map} has been proven, the fact that $\Gamma_\eps$ maps $\mathbb{B}$ into $\mathbb{B}$ follows from the fact that $\gamma_0-2$ is much larger than $\lambda_2$ together with the estimates $\| \higgsdet_0^+\|_{L_x^r},\bigl\|\rhiggsdet_0\bigr\|_{\covHolder{B(0), T}^{-\kappa}}\leq T^{-\frac{1}{\gamma_0}}$, $\tripn{ \Xi_{\Beps} } \leq T^{-\lambda_1}$, and 
\begin{align*} 
\| \mathcal{R}_{\eps,\star} \|_{\sol^{1-\reg_1,-\lambda}}
\lesssim  \tripn{ \Xi_{\Beps} }^2 + 
\| \mathcal{X}_{\eps,\star} \|_{\sol^{1-\reg_1,0}}
+ \| \mathcal{Y}_{\eps,\star} \|_{\sol^{1+\reg_2,-\lambda}} 
\lesssim T^{-2\lambda_1} + T^{-\frac{1}{\gamma_0}-\lambda_2} \sim T^{-\frac{1}{\gamma_0}-\lambda_2}. 
\end{align*}
Note that, in each bracket in \eqref{local:eq-main-self-map} involving the $L_x^r$-norm of $\higgsdet_0^+$, it is the worst term (by approximately a factor of $T^{\frac{1}{2}}$).  In order to obtain \eqref{local:eq-main-self-map}, we estimate the terms in \eqref{local:eq-Xwhite-K-final}--\eqref{local:eq-y-final-new-last} as follows:
\begin{enumerate}[label=(\roman*)]
\item \eqref{local:eq-Xwhite-K-final}, \eqref{local:eq-Xblack-K-final}, and \eqref{local:eq-phix-K-final} are estimated using Lemma~\ref{local:lem-nl-Xterms}.
\item \eqref{local:eq-Y-final-new-1}, \eqref{local:eq-Y-final-new-1p}, and \eqref{local:eq-y-final-new-1} are estimated using our weighted Schauder estimate \eqref{local:eq-weighted-covariant-Schauder-2} and $\tripn{ \Xi_{\Beps} }$. 
\item \eqref{local:eq-Y-final-new-2} and \eqref{local:eq-y-final-new-2} are estimated using Corollary~\ref{local:cor-paracontrolled-specific}.
\item \eqref{local:eq-Y-final-new-3} and \eqref{local:eq-y-final-new-3} are estimated using Lemma~\ref{local:lem-nl-Y-high-high}. 
\item \eqref{local:eq-Y-final-new-4} and the first term in \eqref{local:eq-y-final-new-4} are estimated using Lemma~\ref{local:lem-nl-total-derivative}.
\item \eqref{local:eq-Y-final-new-5} and the second term in \eqref{local:eq-y-final-new-4} are estimated using Lemma~\ref{local:lem-nl-stochastic-low}. 
\item \eqref{local:eq-Y-final-new-6} and the third term in \eqref{local:eq-y-final-new-4} are estimated using Lemma~\ref{local:lem-nl-nonres}. 
\item \eqref{local:eq-Y-final-new-7} and \eqref{local:eq-y-final-new-5} are estimated using Lemma~\ref{local:lem-nl-remainders}. 
\item \eqref{local:eq-Y-final-new-8}, \eqref{local:eq-Y-final-new-9}, \eqref{local:eq-y-final-new-6}, and \eqref{local:eq-y-final-new-7} are estimated using Lemma~\ref{local:lem-nl-power}.
\item \eqref{local:eq-Y-final-new-10} and \eqref{local:eq-y-final-new-last} are estimated using Lemma~\ref{local:lem-8pi-counterterm}. 
\end{enumerate}
Thus, all terms in the definition of $\Gamma_\eps$ make acceptable contributions to \eqref{local:eq-main-self-map}. It now only remains to prove the estimates in \eqref{local:eq-main-R-estimate}. Using the support properties of our time-localised Duhamel integral operators together with Lemma~\ref{local:lem-sol-different-regularities-weights}, as well as the definition of $\tripn{\Xi_{\Beps}}$, we have that 
\begin{align*}
\| \mathcal{R}_{\eps,\star} \|_{\sol^{1-\lambda,0}} 
&\lesssim \bigl\| \bigl( \quadratic[\oeps][][][\ostar] , \phiquadratic[\oeps][][][\ostar], \mixedquadratic[\oeps][][][\ostar] \bigr) \bigr\|_{\sol^{1-\lambda,0}} 
+ \| \mathcal{X}_{\eps,\star} \|_{\sol^{1-\lambda,0}} 
+ \| \mathcal{Y}_{\eps,\star} \|_{\sol^{1-\lambda,0}} \\ 
&\lesssim \tripn{\Xi_{\Beps}}^2 
+  \| \mathcal{X}_{\eps,\star} \|_{\sol^{1-\reg_1,0}}
+ \| \mathcal{Y}_{\eps,\star} \|_{\sol^{1+\reg_2,-\lambda}} \\ 
&\lesssim (T^{-\lambda_1})^2 + T^{-\frac{1}{\gamma_0}-\lambda_2} 
\lesssim T^{-\frac{1}{\gamma_0}-\lambda_2}. 
\end{align*}
Furthermore, using the definition of the $\sol^{1-\lambda,0}$-norm, we also have that 
\begin{equation*}
\| \mathcal{R}_{\eps} \|_{L_{t,x}^\infty} 
\lesssim \sum_{N\geq N_{\lo}} N^{-(1-\lambda)} \| \mathcal{R}_{\eps,N} \|_{\sol^{1-\lambda,0}_N} 
\lesssim T^{\frac{1-\lambda}{2}}  T^{-\frac{1}{\gamma_0}-\lambda_2}. 
\end{equation*}
By taking the limit $\eps\rightarrow 0$, we then obtain the desired estimates \eqref{local:eq-main-R-estimate}. 
\end{proof}

\begin{proof}[Proof of Theorem \ref{main:thm-covariant-local}]
Let $(A_0,\higgs_0)$ be as in the theorem statement and let $B$ be the Yang--Mills gradient flow with initial data $A_0$. We also let $\shiggsdet_0 \coloneq 0$ and $\rhiggsdet_0 \coloneq \higgs_0$. From our conditions in \eqref{main:eq-covariant-local-condition-1}--\eqref{main:eq-covariant-local-condition-2}, it directly follows that Assumptions~\ref{local:assumption-B} and~\ref{local:assumption-higgs-initial} are satisfied. Furthermore, from the definition of $\initialclass_T$ in~\eqref{main:eq-initial-class}, it follows that $B\in C_t^0 \Cs_x^{-\kappa}([0,T]) \cap C_{t,x}^\infty((0,T])$ and that\footnote{A slightly weaker version of \eqref{local:eq-localthm-p1}, in which $\kappa$ is replaced by $2\kappa$, can also be derived from Assumption \ref{local:assumption-B}.} 
\begin{equation}\label{local:eq-localthm-p1}
\bigl\| t^{\frac{1+\kappa}{2}} F_B \bigr\|_{L_{t,x}^\infty([0,T])}, 
\bigl\| t^{1+\frac{\kappa}{2}} \covd_B F_B \bigr\|_{L_{t,x}^\infty([0,T])} < \infty,
\end{equation}
where we also used the equivalence of norms of $\covd_B F_B$ and $\covd_B^\ast F_B$ in two dimensions. From Lemma \ref{lemma:ymgf-approximate-coulomb-gauge}, we then obtain the non-covariant, soft estimate
\begin{equation}\label{local:eq-localthm-p2}
\bigl\| t^{\frac{\kappa}{2}} B \bigr\|_{L_{t,x}^\infty([0,T])}, 
\bigl\| t^{\kappa} \partial^\ell B_\ell \bigr\|_{L_{t,x}^\infty([0,T])} < \infty,
\end{equation}
which will be needed momentarily. We now define the event appearing in Theorem \ref{main:thm-covariant-local} as
\begin{align*}
E(T,A_0,\higgs_0) \eqdef \Bigl\{&\sup_{\eps\in \sdyadic} \tripn{ \Xi_{\Beps} } \leq T^{-\lambda_1} \Bigr\} \medcap \Bigl\{ \lim_{\eps,\eps^\prime\rightarrow 0} \tripn{\Xi_{\Beps}; \Xi_{\Bepsprime}} = 0 \Bigr\} \\ 
&\medcap \Bigl\{ \lim_{\eps,\eps^\prime\rightarrow 0} \bigl\| \graylinear[\Beps][][][] - \graylinear[\Bepsprime][][][] \bigr\|_{C_t^0 \Cs_x^{-\kappa}}=0 \Bigr\}. 
\end{align*}
From Proposition \ref{object:prop-data-set} with $R=T^{-\lambda_1/2}$ and Lemma~\ref{lemma:linear-object-convergence-soft} below, it then follows that there exists a constant $c>0$ such that 
\begin{equation*}
\mathbb{P}\bigl( E(T,A_0,\higgs_0) \bigr) \geq 1 - c^{-1} \exp\bigl(- c T^{-\lambda_1}\bigr).
\end{equation*}
After further decreasing the value of $c$, if necessary, this implies the desired probabilistic estimate of $E(T,A_0,\higgs_0)$. For the remainder of this proof, we restrict ourselves to the event $E(T,A_0,\higgs_0)$. After using the trivial embedding $L_x^\infty \hookrightarrow \Cs_x^{-\kappa}$ and Lemma \ref{local:lem-solution-spaces-hoelder}, the desired convergence of $(A_\eps,\higgs_\eps)$ can be deduced from the following six properties:  
\begin{enumerate}[label=(\alph*)]
\item $B \in C_t^0 \Cs_x^{-\kappa}([0,T])$, which follows from the definition of $\initialclass_T$ in 
\eqref{main:eq-initial-class}.  
\item $\Sdet_\star\in\sol^{1-2\kappa}([0,T])$, which follows from Proposition \ref{local:prop-higgsdet-sdet}.  
\item $\higgsdet \in C_t^0 \Cs_x^{-\kappa}([0,T])$, which follows from Lemma~\ref{lemma:higgsdet-non-covariant-space-soft}. The assumptions $\rhiggsdet_0\in \Cs_x^{-\kappa}$ and~\eqref{eq:higgsdet-non-covariant-space-soft} in Lemma~\ref{lemma:higgsdet-non-covariant-space-soft} are satisfied due to the definition of $\initialclass_T$ in \eqref{main:eq-initial-class} and our soft estimate~\eqref{local:eq-localthm-p2}, respectively.
\item The convergence of $\linear[\Beps][][][]$ and $\philinear[\Beps][][][]$ in $C_t^0 \Cs_x^{-\kappa}([0,T])$, which is part of the definition of $E(T,A_0,\higgs_0)$. 
\item The convergence of $\quadratic[\Beps][][][\ostar]$, $\phiquadratic[\Beps][][][\ostar]$, and $\mixedquadratic[\Beps][][][\ostar]$ in $\sol^{1-2\kappa}([0,T])$, which follows from the convergence of the covariant enhanced data set. 
\item The convergence of $(X_{\eps,\star}, \phix_{\eps,\star}, Y_{\eps,\star}, \phiy_{\eps,\star})$ in $\sol^{1-\reg_1-\lambda_1}([0,T])$, which follows from Lemma \ref{local:lem-sol-different-regularities-weights} and Proposition \ref{local:prop-main}. 
\end{enumerate}
\end{proof}

\section{Stochastic estimates}\label{section:objects}

In this section, we prove that the covariant enhanced data set (see Definition~\ref{local:def-enhanced-data-set}) satisfies the estimates needed in our local theory (see Proposition~\ref{local:prop-main}) with high probability. The main ingredients in the proof will be the covariant heat kernel estimates from Section \ref{section:covariant-heat-kernel-estimates} and Wiener-chaos decompositions. In Section \ref{sec:improved_stochastic_estimates} below, we will also prove more refined stochastic estimates, which will later be needed in our global theory (see Section \ref{section:global-theory}). Roughly speaking, the difference between the ``unrefined'' and ``refined'' stochastic estimates is that $\log$-losses are tolerated in the unrefined estimates, but not in the refined estimates. 

In the following proposition, we collect the main result of this section. 

\begin{proposition}[Control of covariant enhanced data set]\label{object:prop-data-set}
Let Assumption~\ref{local:assumption-B} be satisfied and let $c>0$ be a sufficiently small constant. For all $R\geq 1$, it then holds that
\begin{equation*}
\bP \Big( \sup_{\eps \in \sdyadic} \tripn{ \Xi_{\Beps} }  \geq T^{-\lambda_1/2} R \Big) \leq c^{-1} \exp\bigl( - c R^2 \bigr).
\end{equation*}
Furthermore, it holds that 
\begin{equation*}
\lim_{\eps,\eps^\prime\ra 0} \tripn{ \Xi_{\Beps}; \Xi_{\Bepsprime} }  =0 \qquad \mbox{$\mathbb{P}$-a.s.}
\end{equation*}
\end{proposition}

As in Sections~\ref{section:covariant-Ansatz} and~\ref{section:covariant_lwp}, all covariant stochastic objects will depend on the same connection one-form $B$. We therefore follow Notation~\ref{local:notation-omitting-subscript} again and drop the dependence of the covariant stochastic objects on $B$ from our notation. 

\subsection{Non-resonant estimates}\label{section:non-resonant-estimates}

In this subsection, we estimate moments of the non-resonant parts of the objects in the covariant enhanced data set (Definition~\ref{local:def-enhanced-data-set}). (We also estimate an object which does not appear in the data set \dash see Lemma~\ref{object:lem-linear-approximation}.) As usual, due to various standard reductions (which we defer to Section~\ref{section:estimate-enhanced-dataset-final}), the estimates of the stochastic objects in the appropriate spaces essentially reduce to these moment estimates. We begin with an abstract lemma which will be used for many of the moment estimates.

\begin{lemma}\label{object:lem-time-integration-gain}
Let $C_0 > 1$. There exists a constant $C_1$ depending on $C_0$ such that the following holds. Let $T\in [0,1]$, let $m\in \mathbb{N}$, let $N_1,\hdots,N_m\in \dyadic$ , and let $\mathscr{k}\colon ([0,T]\times \T^2)^{m+1}\rightarrow [0,\infty)$ be a kernel. Furthermore, assume that $\mathscr{k}$ is supported on the set
\begin{equation*}
\Big\{ (s_1,y_1,\hdots, s_m,y_m,t,x) \in ([0,T]\times \T^2)^{m+1} :  |t-s_j| \in [C_0^{-1} N_j^{-2}, C_0 N_j^{-2}] \quad \text{for all } j=1,\hdots,m\Big\}.
\end{equation*}
Let $N_{\max} = \max(N_1, \ldots, N_m)$. Then, it holds for all $t_a,t_b\in [0,T]$ satisfying $t_a<t_b$ and $x\in \T^2$ that 
\begin{align*}
&\quad \int \ds_1 \dy_1 \hdots \ds_m \dy_m \Big( \int_{t_a}^{t_b}  \dt \, \mathscr{k}(s_1,y_1,\hdots,s_m,y_m,t,x) \Big)^2 \\
&\leq C_1 |t_b-t_a| \min \Big( |t_b-t_a|, N_{\max}^{-2} \Big) \sup_{t\in [0,T]}\int \ds_1 \dy_1 \hdots \ds_m \dy_m \,  \mathscr{k}(s_1,y_1,\hdots,s_m,y_m,t,x)^2.
\end{align*}
\end{lemma}

\begin{proof}
This follows directly from the Cauchy--Schwarz inequality and Tonelli's theorem.
\end{proof}

Next, we introduce some convenient notation that will appear in intermediate steps of many of the estimates.

\begin{notation}\label{objects:notation-moralfp}
In the following, let $\moralfp(s, y; t, x)$ be some non-negative function which satisfies the bound
\begin{equ}
\moralfp(s, y; t, x) \lesssim (t-s)^{-1} \exp\bigg(-c\frac{|x-y|^2}{t-s}\bigg),
\end{equ}
for some constant $c$. When $\moralfp$ appears in successive lines, the constant $c$ may change from line to line. Similarly, we write $\moralfp_N(s, y; t, x)$ for a time-localised version of $\moralfp$, i.e.\ some non-negative function which satisfies a bound like $\moralfp_N(s, y; t, x) \lesssim \ind(|t-s| \sim N^{-2}) \moralfp(s, y; t, x)$.
\end{notation}

\begin{remark}[Localised heat kernels as parabolic indicator functions]\label{remark:parabolic-indicator}
The following viewpoint will be extremely convenient in the ensuing stochastic estimates. First, one should think of $\moralfp$ as arising from the covariant heat kernel estimates of Proposition~\ref{prop:qB-ptwise-bds}, which, we recall, estimates $\fq_{\Y,B}$ (and its higher derivatives) by $\moralfp$ (with appropriate additional powers of $(t-s)^{-1}$ to account for higher derivatives). Now if we add in a localisation $N \in \dyadic$ and consider $\fq_{\Y,B, N}$, then we may estimate (by Lemma~\ref{lemma:qB-mixed-all-derivatives}) 
\begin{equ}
|\covd^{(m)}_{t,x} \fq_{\Y,B, N}(s, y; t, x)| \lesssim \lossfactor_T(F) N^{m}\ind(|t-s| \sim N^{-2}) \moralfp(s, y; t, x),
\end{equ}
and one should always keep in mind that morally, the expression
\begin{equ}
\ind(|t-s| \sim N^{-2}) \moralfp(s, y; t, x)
\end{equ}
is a rescaled indicator function on a parabolic ball with the schematic form
\begin{equ}
\ind(|t-s| \sim N^{-2}) \moralfp(s, y; t, x) \approx N^2 \ind(|t-s| \sim N^{-2}) \moralind(|x-y| \lesssim N^{-1}),
\end{equ}
where $\moralind$ is morally a compactly supported bump function, except in reality it is a function with exponential tails. For instance, the expression $\moralind(|x-y| \lesssim N^{-1})$ really refers to a function of the form $\exp\bigl(-c N^2 |x-y|^2\bigr)$, which for all intents and purposes is essentially an indicator function $\ind(|x-y| \lesssim N^{-1})$. As an example of this notation, we have that
\begin{equs}
\int &\du \dz \,\ind(|u-s| \sim M^{-2}) \moralfp(s, y; u, z) \ind(|t-u| \sim N^{-2}) \moralfp(u, z; t, x) \\
&\lesssim \max(M, N)^{-2} \ind(|t-s| \sim \min(M, N)^{-2}) \moralfp(s, y; t, x) \\
&\approx \ind(|t-s| \sim \min(M, N)^{-2}) \max(M, N)^{-2} \min(M, N)^2 \moralind(|x-y| \lesssim \min(M, N)^{-1}).
\end{equs}
We will see further examples in the coming estimates.
\end{remark}

We now begin with the moment estimates, starting with the linear objects. In the following, recall from \eqref{local:eq-lossfactor} the quantity $\lossfactor_T(F)$ defined by
\begin{equ}
\lossfactor_T(F) =\exp\bigl( C_0 \| F \|_{L_t^1 L_x^\infty([0,T])}\bigr)
\bigl( 1+  \Fnorm{F}_{T,100} + T^{\frac{1-\roughdataweight}{2}} \Fnorm{t^{-\frac{1-\roughdataweight}{2}} F}_{T, 1} \bigr)^{100},
\end{equ}
which captures all losses appearing in our heat kernel estimates from Propositions~\ref{prop:pB-ptwise-bounds},~\ref{prop:qB-ptwise-bds}, and Corollary~\ref{cor:pB-qB-parametrix}.

\begin{lemma}[Covariant linear stochastic object]\label{object:lem-linear}
For all $\varep \in \sdyadic$, $0 \leq m \leq 10$, $(t, x) \in [0, T] \times \T^2$, $0 \leq t_0 \leq t_1 \leq T$, $p\geq 2$, and $N \in \dyadic$, it holds that
\begin{equation}\label{object:eq-linear-1}
\Big\| \covd^{(m)}_{t, x} \graylinear[\oeps][][N](t,x) \Big\|_{L^p(\Omega)} \lesssim \lossfactor_T(F) p^{1/2} N^{m},
\end{equation}
as well as
\begin{equation}\label{object:eq-linear-2}
\bigg\| \int_{t_0}^{t_1} \dt\, \covd^{(m)}_{t, x} \graylinear[\oeps][][N](t, x)\bigg\|_{L^p(\Omega)} \lesssim \lossfactor_T(F) p^{1/2} N^{m} |t_0 - t_1|^{1/2} \min\bigl(|t_0 - t_1|^{1/2}, N^{-1}\bigr).
\end{equation}
\end{lemma}

\begin{remark}
We emphasise that, for $m=0$, the right-hand side of \eqref{object:eq-linear-2} decays in $N$. This will later be crucial when we place $\graylinear[\oeps][][N]$ in the modified solution spaces from Definition \ref{local:def-modified-solution-space}. 
\end{remark}

\begin{proof}[Proof of Lemma~\ref{object:lem-linear}] 
By Gaussian hypercontractivity, it suffices to prove the estimates for $p = 2$. For \eqref{object:eq-linear-1}, we apply Lemma~\ref{lemma:qB-mixed-all-derivatives} to obtain 
\begin{equs}
\Big\|\covd^{(m)}_{t, x} \graylinear[\oeps][][N](t, x)\Big\|_{L^2(\Omega)}^2 &= \int \ds \dy \big|\covd^{(m)}_{t, x} \fq_{\Y, B, \varep, N}(s, y; t, x) \big|^2 \\
&\lesssim \lossfactor_T(F)^2 N^{2m} \int \ds \dy \ind(|t-s| \sim N^{-2}) \moralfp(s, y; t, x)^2 \lesssim \lossfactor_T(F)^2 N^{2m},
\end{equs}
as desired. Next, for \eqref{object:eq-linear-2}, by Lemma~\ref{object:lem-time-integration-gain} and our previous estimate, we have that
\begin{equs}
\bigg\| \int_{t_0}^{t_1} \dt\, \covd^{(m)}_{t, x} &\graylinear[\oeps][][N](t, x)\bigg\|_{L^2(\Omega)}^2 \\
&= \int \ds \dy \bigg| \int_{t_0}^{t_1} \dt \, \covd^{(m)}_{t, x} \fq_{\Y, B, \varep, N}(s, y; t, x)\bigg|^2 \\
&\lesssim  |t_0 - t_1| \min\bigl(|t_0 - t_1|, N^{-2}\bigr) \sup_{t \in [0, T]} \int \ds \dy |\covd^{(m)}_{t, x} \fq_{\Y, B, \varep, N}(s, y; t, x)|^2 \\
&\lesssim \lossfactor_T(F)^2 |t_0 - t_1| \min\bigl(|t_0 - t_1|, N^{-2}\bigr) N^{2m},
\end{equs}
as desired.
\end{proof}

We also record the following estimate on the gauge invariant H\"older norm of $\linear_B$, $\philinear_B$. It is not needed in the covariant local theory, but will be used in the global iteration argument. Recall the gauge invariant Hölder norms $\|\place\|_{\covHolder{A, F, T}^\alpha}, \|\place\|_{\covHolder{A, T}^\alpha}$ from Definition~\ref{defn:covariant_Holder}.

\begin{lemma}
\label{le:lollipop_negative_holder}
Let $T, S \in (0, 1]$, and suppose that $B$ is a YM gradient flow on $(0, T)$ and $(T, T + S]$, where a jump from $B(T -)$ to $B(T +)$ is allowed at time $T$. Then, for all $p \geq 2$, $x \in \T^2$, and $s \in (0, S]$, it holds that
\begin{equ}
\label{e:lollipop_negative_holder1}
	\Big\|\bigl(e^{s\Delta^F_{B(T +)}}\linear[B](T)\bigr)(x)\Big\|_{L^p(\Omega)}
	\lesssim \gronwall_{T + s}(C_0 F_B) p^{1/2}
	\Bigl(\frac{T}{s} \wedge \bigl(1 + \abs{\log s} \bigr) \Bigr)^{1/2} \;.
\end{equ}
The same estimate holds for $\philinear[B]$, with $\Delta_{B(T +)}$ in place of $\Delta^F_{B(T +)}$.

Moreover, suppose that
\begin{equ}
	\gronwall_{T}(C_0 F_B) \leq T^{-\lossdel} \;, \quad
	\gronwall_{[T, T + S]}(C_0 F_B)\bigl(1 + \Fnorm{F_B}_{[T, T+S], 1}\bigr) \leq S^{-\lossdel} \;.
\end{equ}
Then, for all $p \geq 2$, it holds that
\begin{equ}
\label{e:lollipop_negative_holder2}
	\Big\|\|\linear[B](T)\|_{\covHolder{B(T +), F, S}^{-\kappa}}\Big\|_{L^p(\Omega)}
	\lesssim p^{1/2} S^{\frac{\kappa}{3}} \;.
\end{equ}
The same estimate holds for $\philinear[B]$, with
$\covHolder{B(T +), S}^{-\kappa}$ in place of
$\covHolder{B(T +), F, S}^{-\kappa}$.
\end{lemma}

\begin{proof}
By the semigroup property at time $T$ and \eqref{eq:qB-ptwise-bounds}, we have that
\begin{equ}
	|\fq_B(u, y; T + s, x)| \lesssim \gronwall_{[u, T + s]}(C_0 F_B) (T + s - u)^{-1} \exp\Bigl(-c\frac{|x-y|^2}{T + s - u}\Bigr) \;.
\end{equ}
Therefore, we get
\begin{equs}
	\Big\|\bigl(e^{s\Delta^F_{B(T +)}}\linear[B](T)\bigr)(x)\Big\|_{L^2(\Omega)}^2
	&= \int_0^T \int_{\T^2} |\fq_B(u, y; T + s, x)|^2 \du \dy \\
	&\lesssim \gronwall_{T + s}(2C_0 F_B) \int_0^T \frac{1}{T + s - u} \du \\
	&\lesssim \gronwall_{T + s}(2C_0 F_B) \Bigl(\frac{T}{s} \wedge \bigl(1 + \abs{\log s}\bigr)\Bigr) \;.
\end{equs}
The first estimate \eqref{e:lollipop_negative_holder1} then follows from Gaussian hypercontractivity. The estimate for $\philinear[B]$ follows in the same way.

We now prove \eqref{e:lollipop_negative_holder2}. By \eqref{eq:qB-ptwise-bounds}, for $0 \leq m \leq 1$ and $s \in (0, S]$, we have
\begin{equs}
	s^{\frac{\kappa+m}{2}}\big\|\covd_{B(T + s)}^{(m)} e^{s\Delta^F_{B(T +)}}\linear[B](T)\big\|_{L_x^\infty}
	&\lesssim S^{-\lossdel}s^{\frac{\kappa+m}{2}}\frac{1}{s}\int_0^s (s - u)^{-\frac{1}{r}-\frac{m}{2}}
	\big\|e^{u\Delta^F_{B(T +)}}\linear[B](T)\big\|_{L_x^r}\du \\
	&\lesssim S^{-\lossdel}s^{\frac{\kappa}{2}-\frac{2}{r}}
	\Bigl(\int_0^s\big\|e^{u\Delta^F_{B(T +)}}\linear[B](T)\big\|_{L_x^r}^r\du\Bigr)^{\frac{1}{r}} \;.
\end{equs}
It follows that
\begin{equ}
	\|\linear[B](T)\|_{\covHolder{B(T +), F, S}^{-\kappa}}
	\lesssim S^{\frac{\kappa}{2}-\frac{2}{r}-\lossdel}
	\Bigl(\int_0^S\big\|e^{u\Delta^F_{B(T +)}}\linear[B](T)\big\|_{L_x^r}^r\du\Bigr)^{\frac{1}{r}} \;.
\end{equ}
Let $p \geq r$. By Minkowski's inequality,
\begin{equ}
	\Big\|\Bigl(\int_0^S\big\|e^{u\Delta^F_{B(T +)}}\linear[B](T)\big\|_{L_x^r}^r\du\Bigr)^{\frac{1}{r}}\Big\|_{L^p(\Omega)}
	\leq \Bigl(\int_0^S\Big\|\big\|e^{u\Delta^F_{B(T +)}}\linear[B](T)\big\|_{L_x^r}\Big\|_{L^p(\Omega)}^r\du\Bigr)^{\frac{1}{r}} \;.
\end{equ}
For the integrand, Jensen's inequality on $\T^2$, Gaussian hypercontractivity, and \eqref{e:lollipop_negative_holder1} give
\begin{equs}
	\Big\|\big\|e^{u\Delta^F_{B(T +)}}\linear[B](T)\big\|_{L_x^r}\Big\|_{L^p(\Omega)}^p
	&\leq \int_{\T^2}\E\Big|\bigl(e^{u\Delta^F_{B(T +)}}\linear[B](T)\bigr)(x)\Big|^p\dx \\
	&\lesssim p^{\frac{p}{2}}\Bigl(T^{-2\lossdel}S^{-2\lossdel}
	\Bigl(\frac{T}{u} \wedge \bigl(1 + \abs{\log u}\bigr)\Bigr)\Bigr)^{\frac{p}{2}} \;.
\end{equs}
Consequently, we obtain
\begin{equs}
	\Big\|\Bigl(\int_0^S\big\|e^{u\Delta^F_{B(T +)}}\linear[B](T)\big\|_{L_x^r}^r\du\Bigr)^{\frac{1}{r}}\Big\|_{L^p(\Omega)}
	&\lesssim p^{1/2} T^{-\lossdel} S^{-\lossdel}
	\Bigl(\int_0^S\Bigl(\frac{T}{u} \wedge \bigl(1 + \abs{\log u}\bigr) \Bigr)^{\frac{r}{2}}\du\Bigr)^{\frac{1}{r}} \\
	&\lesssim_r p^{1/2} T^{\frac{1}{2r}-\lossdel} S^{-\lossdel}
	\lesssim p^{1/2} S^{-\lossdel} \;.
\end{equs}
Combining the previous estimates, we conclude that
\begin{equ}
	\Big\|\|\linear[B](T)\|_{\covHolder{B(T +), F, S}^{-\kappa}}\Big\|_{L^p(\Omega)}
	\lesssim_r p^{1/2} S^{\frac{\kappa}{2}-\frac{2}{r}-2\lossdel}
	\lesssim p^{1/2} S^{\frac{\kappa}{3}} \;.
\end{equ}
This proves the desired estimate for $p \geq r$. The case $2 \leq p < r$ follows from Jensen's inequality. The proof for $\philinear[B]$ is the same.
\end{proof}

Before moving on to the other objects in the data set, we state and prove a result which gives a ``parametrix approximation'' of the linear object, which will be useful when computing resonances in Sections~\ref{section:non-paracontrolled-resonant} and~\ref{section:objects-paracontrolled}. First, we define the approximation. In the following, $\fp_N$ denotes the localised version of the usual heat kernel, i.e.\ the kernel $\fp_{B, N}$ when the background $B \equiv 0$. Similarly, $\fp_{\eps, N}$ denotes $\fp_{\Beps, N}$ when $B \equiv 0$ (this notation will be used later).

\begin{definition}\label{def:parametrix-approximate-objects}
For $N \in \dyadic$, we define $\pmgraylinear[\oeps][][N]\colon [0,T]\times \T^2 \rightarrow \Y$ by
\begin{equ}
\pmgraylinear[\oeps][][N](t,x) \eqdef \int \ds \dy \, \rho_\varep(t-s) (\fp_N \ptcutoff \pt_B)(s, y; t, x) \boldxi(s, y). 
\end{equ}
To be clear, here $\pt_B = \tauYB$, but we stick to our usual convention of keeping the underlying vector space implicit. More generally, given $m=0,1$, let
\begin{equ}
\lsup{(m)}\pmgraylinear[\oeps][][N](t, x) \eqdef \int \ds \dy \, \rho_\varep(t-s) ( (\nabla^{(m)} \fp_N) \ptcutoff \pt_B)(s, y; t, x) \boldxi(s, y),
\end{equ}
where the $\nabla^{(m)}$ appearing in $\nabla^{(m)} \fp_N$ acts on the $x$-variable. In particular, $\lsup{(0)} \pmgraylinear = \pmgraylinear$.
\end{definition}

The following lemma shows that the difference between the linear object and its parametrix approximation may be estimated with a gain of almost one whole power of the frequency.

\begin{lemma}[Approximation of covariant linear stochastic object]\label{object:lem-linear-approximation}
For all $0\leq m \leq 1$, $p\geq 2$, $(t, x) \in [0, T] \times \T^2$, and $N\in \dyadic$, it holds that 
\begin{equ}\label{object:eq-linear-approximation-L2}
\Bigl\|\covd^{(m)} \graylinear[\oeps][][][N](t,x)  - \lsup{(m)}\pmgraylinear[\oeps][][N] (t, x) \Bigr\|_{L^p(\Omega)} \lesssim \sqrt{p} \, \lossfactor_T(F)\bigl( T N^2 \bigr)^{-(1-\roughdataweight)/2} N^{m}.
\end{equ}
\end{lemma}
We emphasise the restriction to $0\leq m \leq 1$ in Lemma~\ref{object:lem-linear-approximation}, which is due to Remark~\ref{remark:pB-parametrix-higher-order-derivatives}.

\begin{proof}[Proof of Lemma~\ref{object:lem-linear-approximation}]
By Gaussian hypercontractivity, it suffices to take $p = 2$. Fix $(t,x) \in [0,T] \times \T^2$. Then, it holds that
\begin{align*}
&\E \bigg[ \Big| \covd^{(m)} \graylinear[\oeps][][][N](t,x)  - \lsup{(m)}\pmgraylinear[\oeps][][][N](t, x)\Big|^2 \bigg]  \\
=& \E\bigg[ \Big| \iint \ds \dy \rho_\varep(t-s) \Big( \bigl(\covd^{(m)} \fq_{\Y, B,N} - \bigl(  (\nabla^{(m)} \fp_N) \ptcutoff \tauYB\bigr)\bigr)(s,y;t,x) \Big) \boldxi(s,y) \Big|^2   \bigg] \\
\lesssim& \iint \ds \dy \Big|  \Bigl(\covd^{(m)} \fq_{\Y, B,N} - \bigl(  (\nabla^{(m)} \fp_N) \ptcutoff \tauYB\bigr)\Bigr)(s,y;t,x) \Big|^2. 
\end{align*}
By the estimate \eqref{eq:qB-parametrix-t-s-power} from Corollary~\ref{cor:pB-qB-parametrix} with the $\alpha$ therein taken to be $\roughdataweight$ (and recalling the definition \eqref{local:eq-lossfactor} of $\lossfactor_T(F)$), as well as the estimates \eqref{eq:ptl-t-B-estimate-under-more-stringent-Fpert}, \eqref{eq:F-estimate-under-more-stringent-Fpert}, and \eqref{eq:Fpert-estimate-under-more-stringent-Fpert}, we have the pointwise estimate
\begin{equs}
\Bigl| \Bigl(\covd^{(m)} &\fq_{\Y, B} - (\nabla^{(m)} \fp) \ptcutoff \pt_{\Y, B}\Bigr) (s, y; t, x)\Bigr| \lesssim \lossfactor_T(F) \Bigl( (T^{-1}(t-s))^{1-\roughdataweight}s^{-\roughdataweight} (t-s)^{\roughdataweight} \\
&+ (T^{-1}(t-s))^{(1-\roughdataweight)/2} \max(s^{-\roughdataweight/2} (t-s)^{\roughdataweight/2}, 1)\Bigr) (t-s)^{-\frac{m}{2}} \moralfp(s,y; t, x). 
\end{equs}
Inserting this estimate, we obtain
\begin{align*}
&\quad \iint \ds \dy \Big| \Bigl( \covd^{(m)} \fq_{\Y, B,N}- \bigl(  (\nabla^{(m)} \fp_N) \ptcutoff  \tauYB \bigr)\Bigr)(s,y;t,x) \Big|^2 \\
&\lesssim \lossfactor_T(F)^2 \iint \ds \dy \,  \chi_N(t-s)^2 \bigl(T^{-1}(t-s)\bigr)^{(1-\roughdataweight)} \bigl(s^{-\roughdataweightmultiple[2]} (t-s)^{\roughdataweightmultiple[2]} \\
&\quad \quad \quad \quad \quad \quad \quad \quad + \max(s^{-\roughdataweight} (t-s)^{\roughdataweight}, 1)\bigr) (t-s)^{-m} \moralfp(s, y; t, x)^2 \\
&\lesssim \lossfactor_T(F)^2   (T N^2)^{-(1-\roughdataweight)} N^{2m}.
\end{align*}
The desired estimate now follows.
\end{proof}

We proceed to estimate moments of the other objects appearing in the covariant enhanced data set.

\begin{lemma}[Quadratic stochastic object with derivative]\label{lemma:quadratic-object-non-resonant-bound}
Let $\Qc \colon \Y \times \Y \ra \Y$ be a $\liegroup$-equivariant bilinear map, $k \in [2]$, $\varep \in \sdyadic$, $m \geq 0$, $(t, x) \in [0, T] \times \T^2$, $p \geq 1$, and $M, N_1, N_2 \in \dyadic$. We have that
\begin{equ}
\biggl\| \Wick{\covd^{(m)}_{t,x} \Duh_{\Y, B, M}^{F}\Bigl(\Qc\bigl(\graylinear[\oeps][][N_1], \covd^k \graylinear[\oeps][][N_2] \bigr)\Bigr)} (t, x) \biggr\|_{L^p(\Omega)} \lesssim \lossfactor_T(F)^3 p M^{m} \max(M, N_1, N_2)^{-1}.   
\end{equ}
\end{lemma}
\begin{proof}
By Gaussian hypercontractivity, it suffices to prove the estimate for $p = 2$. Towards this end, we have that
\begin{equs}
X &\eqdef \Wick{\covd^{(m)}_{t,x} \Duh_{\Y, B, M}^{F}\Bigl(\Qc\bigl(\graylinear[\oeps][][N_1], \covd^k \graylinear[\oeps][][N_2] \bigr)\Bigr)} (t, x) \\
&= \Wick{\int \ds_1 \ds_2 \dy_1 \dy_2 \int \ds \dy \, \covd^{(m)}_{t,x} \fq_{\Y, B, M}(s, y; t, x) \Qc\bigl( \fq_{\Y, B, \varep, N_1} (s_1, y_1; s, y) \boldxi(s_1, y_1), \\
&\quad \quad \quad \quad \quad \quad \quad \quad \quad \quad \covd^k \fq_{\Y, B, \varep, N_2}(s_2, y_2; s, y) \boldxi(s_2, y_2)\bigr)}.
\end{equs}
Thus by the kernel estimates of Propositions~\ref{prop:pB-ptwise-bounds},~\ref{prop:qB-ptwise-bds} and Lemma~\ref{lemma:qB-mixed-all-derivatives}, we have that
\begin{equ}
\big\|X\big\|_{L^2(\Omega)}^2 \lesssim \lossfactor_T(F)^6 M^{2m} \int \ds_1 \ds_2 \dy_1 \dy_2 \kernel(s_1, y_1, s_2, y_2)^2,
\end{equ}
where
\begin{equs}
\kernel(s_1, y_1, s_2, y_2) = \int \ds \dy \, &\ind(|t-s | \sim M^{-2}) \ind(|s-s_1| \sim N_1^{-2}) \ind(|s-s_2| \sim N_2^{-2})   \\
& \times N_2 \moralfp(s, y; t, x) \moralfp(s_1, y_1; s, y) \moralfp(s_2, y_2; s, y).
\end{equs}
We focus on the case $\max(M, N_1, N_2) = N_2$, which is the most dangerous, because of the extra factor of $N_2$ (which arises from the covariant derivative falling on $\graylinear[B][][N_2]$). In this case, we put $\moralfp(s, y; t, x)$ and $\moralfp(s_1, y_1; s, y)$ in $L_y^\infty$ and integrate $\moralfp(s_2, y_2; s, y)$, thus giving the estimate (the function $L(y_1, y_2)$ appearing below is defined just after the following display)
\begin{equs}
\kernel(&s_1, y_1, s_2, y_2) \\
&\lesssim  M^2 N_1^2 N_2 L(y_1, y_2) \int \ds \ind(|t-s | \sim M^{-2}) \ind(|s-s_1| \sim N_1^{-2}) \ind(|s-s_2| \sim N_2^{-2})  \\
&\lesssim M^2 N_1^2 N_2^{-1} L(y_1, y_2) \ind(|t-s_1| \sim \min(M, N_1)^{-2}) \ind(|t-s_2| \sim M^{-2}) \ind(|s_1 - s_2| \lesssim N_1^{-2}),
\end{equs}
where $L(y_1, y_2)$ is essentially a product of indicator functions, in that it satisfies the estimate
\begin{equ}
L(y_1, y_2) \lesssim \moralind(|x-y_1| \lesssim \min(M, N_1)^{-1}) \moralind(|x - y_2| \lesssim M^{-1}) \moralind(|y_1 - y_2| \lesssim N_1^{-1}).
\end{equ}
From this, we may obtain
\begin{equ}\label{eq:derivative-quadratic-object-non-resonant-kernel-bd}
\int \ds_1 \ds_2 \dy_1 \dy_2 \, \kernel(s_1, y_1, s_2, y_2)^2 \lesssim M^4 N_1^4 N_2^{-2} M^{-4} N_1^{-4} \lesssim \max(M, N_1, N_2)^{-2},
\end{equ}
and thus we further obtain
\begin{equs}
\big\|X\big\|_{L^2(\Omega)}^2 &\lesssim \lossfactor_T(F)^6 M^{2m} \int \ds_1 \ds_2 \dy_1 \dy_2 \, \kernel(s_1, y_1, s_2, y_2)^2 \\
&\lesssim \lossfactor_T(F)^6 M^{2 m} \max(M, N_1, N_2)^{-2}, 
\end{equs}
from which the desired second moment estimate follows.
\end{proof}

\begin{lemma}[Power-type quadratic and cubic objects]\label{lemma:no-derivative-objects}
For any $\varep \in \sdyadic$, $(t, x) \in [0, T] \times \T^2$, $0 \leq t_0 \leq t_1 \leq T$, $p \geq 1$, and $N_1, N_2 \in \dyadic$, we have that
\begin{equs}
\Big\| \Wick{\graylinear[\oeps][][N_1] \otimes \graylinear[\oeps][][N_2]}(t, x) \Big\|_{L^p(\Omega)} &\lesssim \lossfactor_T(F)^2 p, \\
\bigg\|\int_{t_0}^{t_1} \dt \, \Wick{\graylinear[\oeps][][N_1] \otimes \graylinear[\oeps][][N_2]}(t, x)\bigg\|_{L^p(\Omega)} &\lesssim \lossfactor_T(F)^2 p |t_0 - t_1|^{1/2} \min(|t_0 - t_1|^{1/2}, \max(N_1, N_2)^{-1}).
\end{equs}
Additionally, for any $N_3 \in \dyadic$, we have that
\begin{equs}
\Big\| \Wick{\graylinear[\oeps][][N_1] \otimes \graylinear[\oeps][][N_2] \otimes \graylinear[\oeps][][N_3]}(t, x) \Big\|_{L^p(\Omega)} &\lesssim \lossfactor_T(F)^3 p^{3/2}, \\
\bigg\|\int_{t_0}^{t_1} \dt \,\Wick{\graylinear[\oeps][][N_1] \otimes \graylinear[\oeps][][N_2] \otimes \graylinear[\oeps][][N_3]}(t, x)\bigg\|_{L^p(\Omega)} &\lesssim \lossfactor_T(F)^3 p^{3/2} |t_0 - t_1|^{1/2} \\
&\quad \quad \times \min(|t_0 - t_1|^{1/2}, \max(N_1, N_2, N_3)^{-1}).
\end{equs}
\end{lemma}
\begin{proof}
We focus on the cubic case, as the quadratic case follows by a simpler argument. 
By Gaussian hypercontractivity, it suffices to consider $p = 2$. By the covariant heat kernel estimates of Proposition~\ref{prop:qB-ptwise-bds}, we have that
\begin{equs}
\bigg\|& \Wick{\graylinear[\oeps][][N_1] \otimes \graylinear[\oeps][][N_2] \otimes \graylinear[\oeps][][N_3]}(t, x) \bigg\|_{L^2(\Omega)}^2 \\
&\lesssim \lossfactor_T(F)^6 \int \ds_1 \ds_2 \ds_3 \dy_1 \dy_2 \dy_3 \ind(|t-s_1| \sim N_1^{-2}) \ind(|t-s_2| \sim N_2^{-2}) \ind(|t-s_3| \sim N_3^{-2}) \\
& \quad \quad \quad \quad \times\moralfp(s_1, y_1; t, x)^2 \moralfp(s_2, y_2; t, x)^2 \moralfp(s_3, y_3; t, x)^2 \\
&\lesssim \lossfactor_T(F)^6.
\end{equs}
The first estimate for the cubic object now follows. For the second estimate, we apply Lemma~\ref{object:lem-time-integration-gain} and then proceed by a similar argument.
\end{proof}

\begin{lemma}[Quadratic object with two derivatives]\label{lemma:xi-para-non-resonant-estimates}
Let $k, \ell \in [2]$, $\varep \in \sdyadic$, $(t, x) \in [0, T] \times \T^2$, $0 \leq t_0 \leq t_1 \leq T$, $p \geq 1$, and $M, N_1, N_2 \in \dyadic$. We have that
\begin{equs}
\Big\| \Wick{\covd^k \Duh_{\Y, B, M}\bigl(\graylinear[\oeps][][N_1]\bigr) \otimes \covd^\ell \graylinear[\oeps][][N_2]}(t, x)\Big\|_{L^p(\Omega)} &\lesssim \lossfactor_T(F)^3 p \max(M, N_1)^{-1} N_2, \\
\bigg\| \int_{t_0}^{t_1} \dt\, \Wick{\covd^k \Duh_{\Y, B, M}\bigl( \graylinear[\oeps][][N_1]\bigr) \otimes \covd^\ell \graylinear[\oeps][][N_2]}(t, x)\bigg\|_{L^p(\Omega)}&\lesssim \lossfactor_T(F)^3 p |t_0 -t_1|^{1/2} \\
&\hspace{-45mm}\times \min\bigl(|t_0 - t_1|^{1/2}, \max(\min(M, N_1), N_2)^{-1} \bigr) \max(M, N_1)^{-1} N_2.
\end{equs}
\end{lemma}
\begin{proof}
For the first estimate, we may write
\begin{equs}
X &\eqdef  \Wick{\covd^k  \Duh_{\Y, B, M}\bigl(\graylinear[\oeps][][N_1]\bigr) \otimes \covd^\ell \graylinear[\oeps][][N_2]}(t, x) \\
&= \Wick{\int \ds_1 \dy_1 \ds_2 \dy_2\, \kernel(s_1, y_1, s_2, y_2; t, x) (\boldxi(s_1, y_1) \otimes \boldxi(s_2, y_2))},
\end{equs}
where
\begin{equ}
\kernel(s_1, y_1, s_2, y_2; t, x) = \bigg(\int \ds \dy\, \covd^k \fp_{\Y, B, M} (s, y; t, x)  \fq_{\Y, B, \varep, N_1}(s_1, y_1; s, y)\bigg) \otimes \covd^\ell \fq_{\Y, B, \varep, N_2}(s_2, y_2; t, x).
\end{equ}
By the covariant heat kernel estimates of Propositions~\ref{prop:pB-ptwise-bounds} and~\ref{prop:qB-ptwise-bds}, we have that
\begin{equs}
\big\|X\big\|_{L^2(\Omega)}^2 &\lesssim \lossfactor_T(F)^6 \int \ds_1 \dy_1 \ds_2 \dy_2 \bigg|M \int \ds \dy \, \moralfp_M(s, y; t, x) \moralfp_{N_1}(s_1, y_1; s, y)\bigg|^2 N_2^2 \moralfp_{N_2}(s_2, y_2; t, x)^2 \\
&\lesssim \lossfactor_T(F)^6 N_2^2 \int \ds_1 \dy_1 \, \bigg|M \moralfp(s_1, y_1; t, x) \int\ds \,\ind(|t-s| \sim M^{-2}) \ind(|s - s_1| \sim N_1^{-2})\bigg|^2 \\
&\lesssim \lossfactor_T(F)^6 N_2^2 \max(M, N_1)^{-2} \int \ds_1 \dy_1 \,\moralfp(s_1, y_1; t, x)^2 \ind(|t-s_1| \sim \min(M, N_1)^{-2}) \\
&\lesssim \lossfactor_T(F)^6 N_2^2 \max(M, N_1)^{-2}.
\end{equs}
The desired second moment estimate now follows. For the estimate of the time integral, we apply Lemma~\ref{object:lem-time-integration-gain} and proceed similarly; we note here that the kernel $\kernel(s_1, y_1, s_2, y_2; t, x)$ is supported on times $|t - s_1| \sim \min(M, N_1)^{-2}$ and $|t - s_2| \sim N_2^{-2}$, which is why (upon applying Lemma~\ref{object:lem-time-integration-gain}) we obtain the factor $\max(\min(M, N_1), N_2)^{-1}$ in the final estimate.
\end{proof}

\subsection{Resonant estimates (non-paracontrolled)}\label{section:non-paracontrolled-resonant}

In this subsection, we estimate the resonant parts of the various objects which appear in the enhanced data set, excluding the resonances involving the high$\times$high-products in our paracontrolled calculus, which are handled in Section~\ref{section:objects-paracontrolled}.
We state the main results of this section (Lemmas~\ref{lemma:res-linear-linear} and~\ref{lemma:qres-estimate-sol-space}), and then afterwards turn towards the proofs. First, we make the following definition. 

\begin{definition}\label{objects:def-localized-resonance}
For $(t, x) \in [0, T] \times \T^2$, $M, N \geq N_\lo$, define
\begin{equ}
\sigma_{\varep, M, N}^2(t, x) \eqdef \int \ds \dy \, \rho_\varep(t-s)^2 \chi_N(t-s) \chi_M(t-s) \fp(s, y; t, x)^2.
\end{equ}
\end{definition}

We make several immediate remarks about these resonances. First, by translation invariance, $\sigma_{\varep, M, N}^2$ does not depend on the spatial point $x$, and thus we may write $\sigma_{\varep, M, N}^2(t)$. Note also that $\sigma_{\varep, M, N}^2 \neq 0$ only if $M \sim N$, and that (here we use translation invariance to set $x = 0$)
\begin{equ}\label{eq:sigma-varep-localized-sum}
\sum_{M, N \geq N_{\lo}} \sigma_{\varep, M, N}^2 = \int_0^t \ds \dy \, \rho_\varep(t-s)^2 \fp(s, y; t, 0)^2.
\end{equ}
In order to more explicitly compute the right-hand side above, we will need to shift from $\T^2$ to $\R^2$. To do so, we identify $\T^2$ with $(-\pi, \pi]^2$, and recall (Notation~\ref{notation:euclidean-heat-kernels}) that the heat kernel $\fp$ on $\T^2$ is the periodisation of the heat kernel $\fp^{\Euc}$ on $\R^2$. Thus for $y \in (-\pi, \pi]^2$, we may obtain
\begin{equ}
\fp(s, y; t, 0) = \fp_{t-s}^{\Euc}(y) + \bigO((t-s)^{20}),
\end{equ}
and thus we may compute (recalling that $\fp_{t}^{\Euc}(x) = (4\pi t)^{-1} \exp(-|x|^2 / 4t)$)
\begin{equation}\label{eq:resonance-torus-to-plane}
\begin{aligned}
\int_{\T^2} \dy \, \fp(s, y; t, 0)^2 &= \int_{(-\pi, \pi]^2}\dy \, \fp_{t-s}^{\Euc}(y)^2 + \bigO((t-s)^{10}) \\
&= \int_{\R^2} \dy \, \fp_{t-s}^{\Euc}(y)^2 + \bigO((t-s)^{10}) = \frac{1}{8\pi(t-s)} + \bigO((t-s)^{10}).
\end{aligned}
\end{equation}
Recalling the definition \eqref{main:eq-sigma-eps-Crho} of $\sigma_{\rho, \varep}^2$, we thus have that
\begin{equ}\label{eq:sigma-varep-approximation}
\biggl| \sigma_{\rho, \varep}^2 - \int_0^t \ds \dy \rho_\varep(t-s)^2 \fp(s, y; t, x)^2 \biggr| \lesssim \abs{\log t} \wedge \abs{\log \varep}.
\end{equ}
For later use, we also record for all $j\in[2]$ the estimate 
\begin{equation}\label{eq:resonance-torus-to-plane-derivative}
\int_{\T^2} \dy \, (\partial_x^j \fp)(s, y; t, 0)^2 =\frac{1}{32\pi(t-s)^2} + \bigO((t-s)^{10}),
\end{equation}
which can be proven using similar arguments. 

The first main result of the section estimates the resonances of the power-type quadratic objects. 

\begin{lemma}[Resonances of power-type quadratic objects]\label{lemma:res-linear-linear}
For any $\varep \in \sdyadic$, $(t, x) \in [0, T] \times \T^2$, and $M, N \geq N_\lo$, we have that
\begin{equ}\label{eq:localized-quadratic-resonance-estimate}
\Bigl| \E\Bigl[  \graylinear[\oeps][][M] \otimes \graylinear[\oeps][][N]\Bigr](t, x) - \sigma_{\varep, M, N}^2(t) \Cov_{\Y} \Bigr| \lesssim \lossfactor_T(F)^2  \ind(M \sim N) (T N^2)^{-(1-\roughdataweight)/2} .
\end{equ}
For notational convenience, let
\[ Q_\varep(t, x) = \E\Bigl[ \graylinear[\oeps][][] \otimes \graylinear[\oeps][][]\Bigr](t, x) -  \sigma_{\rho, \varep}^2 \Cov_{\Y}. \]
Then we also have that
\begin{equ}\label{eq:summed-quadratic-resonance-estimate}
\bigl| Q_\varep(t, x) \bigr| \lesssim \lossfactor_T(F)^2 + \abs{\log t}.
\end{equ}
Finally, for any $1 \leq q < \infty$, we have that $Q_\varep$ converges as $\varep \downarrow 0$ in $L_t^q L_x^\infty([0, T])$ to a quantity which does not depend on the mollifier $\rho \in \mc{P}$.
\end{lemma}

The next main result estimates the resonances of the quadratic objects with one derivative. For notational purposes, given $(t, x) \in [0, T] \times \T^2$, $\varep \in \sdyadic$, and $N \geq N_{\lo}$, let
\begin{equ}\label{objects:eq-resops}
 \qres_{\varep, N}^A(t, x) \eqdef \E\Big[\quadratic[\oeps][][N]\Big](t, x) , \quad 
\qres_{\varep, N}^\higgs(t, x) \eqdef \E\Big[\phiquadratic[\oeps][][N]\Big](t, x).
\end{equ}

\begin{lemma}\label{lemma:qres-estimate-sol-space}
Let $0 < T \leq \frac{1}{128}$. For any $\varep \in \sdyadic$, we have that
\begin{equs}
\bigl\|\qres_{\varep, \star}^A \bigr\|_{\sol^{1-2\kappa}}  &\lesssim \lossfactor_T(F)^3 , \\
\bigl\|\qres_{\varep, \star}^\higgs\bigr\|_{\sol^{1-2\kappa}} &\lesssim \lossfactor_T(F)^4.
\end{equs}
Additionally, $\qres_{\varep, \star}^A, \qres_{\varep, \star}^\higgs$ converge in $\sol^{1-2\kappa}$ to limits which are independent of the mollifier $\rho \in \mc{P}$. 
\end{lemma}

Having stated the main results, we proceed with the proofs. We begin with Lemma~\ref{lemma:res-linear-linear}.

\begin{proof}[Proof of Lemma~\ref{lemma:res-linear-linear}]
The second estimate \eqref{eq:summed-quadratic-resonance-estimate} follows from the first estimate \eqref{eq:localized-quadratic-resonance-estimate} by summing in $M, N$ and applying \eqref{eq:sigma-varep-localized-sum} and \eqref{eq:sigma-varep-approximation}. Thus, we focus on the first estimate. First observe that the left-hand side of the first estimate is zero unless $M \sim N$, due to the definition of space-time white noise as well as $\sigma_{\varep, M, N}^2$. Thus, we may assume that $M \sim N$. Next, for notational convenience, let
\begin{equ}
X_M = \graylinear[\oeps][][M] - \pmgraylinear[\oeps][][M].
\end{equ}
We may then expand
\begin{equ}\label{eq:res-linear-linear-intermediate-expansion}
\graylinear[\oeps][][M] \otimes \graylinear[\oeps][][N] = \pmgraylinear[\oeps][][M] \otimes \pmgraylinear[\oeps][][N] + X_M \otimes \graylinear[\oeps][][N] + \pmgraylinear[\oeps][][M] \otimes X_N.
\end{equ}
From Definition~\ref{def:parametrix-approximate-objects} and the fact that parallel transport $\pt_B$ is an isometry (Remark~\ref{remark:parallel-transport}\ref{remark:pt-isometry}), we have that
\begin{equ}
\E\Bigl[\pmgraylinear[\oeps][][M] \otimes \pmgraylinear[\oeps][][N]\Bigr] =   \int \ds \dy \, \rho_\varep(t-s)^2 \chi_M(t-s) \chi_N(t-s) \ptcutoff(x-y)^2 \fp(s, y; t, x)^2 \Cov_{\Y}. 
\end{equ}
Here, we recall that $\ptcutoff$ is the spatial localisation from Definition~\ref{def:ptcutoff}, which is distinct from the quantity $\sigma_{\varep, M, N}^2$ which will also appear soon. Since $\ptcutoff(x-y)$ enforces the localisation $|x-y| \lesssim 1$, we have that $1 - \ptcutoff(x-y)^2$ is supported on regions where $|x-y| \gtrsim 1$, and thus
\begin{equ}
\Bigl| \int \ds \dy \, \rho_\varep(t-s)^2 \chi_M(t-s) \chi_N(t-s) \ptcutoff(x-y)^2 \fp(s, y; t, x)^2 - \sigma_{\varep, M, N}^2(t, x)\Bigr| \lesssim \ind(M \sim N) N^{-10},
\end{equ}
which is more than acceptable. Thus to finish, it suffices to show that the second and third terms on the right-hand side of \eqref{eq:res-linear-linear-intermediate-expansion} may be estimated by acceptable contributions. By Lemmas~\ref{object:lem-linear} and~\ref{object:lem-linear-approximation}, we may estimate the second term by
\begin{equ}
\Bigl| \E\bigl[ X_M (t, x) \otimes \graylinear[\oeps][][N] (t, x) \bigr]\Bigr| \leq \E\bigl[ \bigl| X_M (t, x)\bigr|^2\bigr]^{1/2} \E \bigl[ \bigl|\graylinear[\oeps][][N] (t, x)\bigr|^2\bigr]^{1/2} \lesssim \lossfactor_T(F)^2 (TM^2)^{-(1-\roughdataweight)/2},
\end{equ}
which is acceptable. The third term on the right-hand side of \eqref{eq:res-linear-linear-intermediate-expansion} may be similarly estimated. Combining, the desired estimate \eqref{eq:localized-quadratic-resonance-estimate} follows.

To finish, we discuss convergence of $Q_\varep$ as $\varep \downarrow 0$. By \eqref{eq:sigma-varep-localized-sum}, we may decompose
\begin{equs}
Q_\varep = \,& \Bigl( \int_0^t \ds \dy \rho_\varep(t-s)^2 \fp(s, y; t, x)^2 - \sigma_{\rho, \varep}^2 \Bigr) \Cov_{\Y} \\
&+   \sum_{M, N \geq N_{\lo}} \Bigl(\E\Bigl[  \graylinear[\oeps][][M] \otimes \graylinear[\oeps][][N]\Bigr](t, x) -  \sigma_{\varep, M, N}^2(t) \Cov_{\Y}\Bigr).
\end{equs}
By \eqref{eq:localized-quadratic-resonance-estimate}, each summand in the second term on the right-hand side above is estimated uniformly in $\varep$. Moreover, by direct verification, the limit as $\varep \downarrow 0$ of each summand exists and is independent of $\rho$. Thus by dominated convergence, we obtain that the limit as $\varep \downarrow 0$ of the second term on the right-hand side above exists and is independent of $\rho$. Due to the estimate \eqref{eq:localized-quadratic-resonance-estimate}, the convergence may be taken in $L_{t,x}^\infty([0, T])$, which is stronger than the claimed convergence in $L_t^q L_x^\infty([0, T])$. 

Next, we discuss the convergence of the first term. By \eqref{eq:resonance-torus-to-plane} and the definition \eqref{main:eq-sigma-eps-Crho} of $\sigma_{\rho, \varep}^2$, we have that
\begin{equ}
\int_0^t \ds \dy \rho_\varep(t-s)^2 \fp(s, y; t, x)^2 - \sigma_{\rho, \varep}^2 = -\int_t^1 \ds \, \frac{\rho_\varep(s)^2}{8\pi s} + \int_0^t \ds \, \rho_\varep(s)^2 \biggl(\int_{\T^2} \dy \fp_s(y)^2 - \frac{1}{8\pi s} \biggr). 
\end{equ}
By the estimate \eqref{eq:resonance-torus-to-plane}, the second term above converges in $L_{t,x}^\infty([0, T])$, which again is stronger than the claimed convergence. Next, we have that pointwise in $t > 0$,
\begin{equ}
\lim_{\varep \downarrow 0} \int_t^1 \ds \, \frac{\rho_\varep(s)^2}{8 \pi s}  = \frac{1}{8\pi} \abs{\log t}.
\end{equ}
Combining this with the fact that 
\begin{equ}
\biggl|\int_t^1 \ds \, \frac{\rho_\varep(s)^2}{8 \pi s} \biggr| \lesssim \abs{\log t},
\end{equ}
we may obtain (by dominated convergence) the desired convergence in $L_t^q L_x^\infty([0, T])$ for $1 \leq q < \infty$.
\end{proof}

We proceed towards Lemma~\ref{lemma:qres-estimate-sol-space}. Before proving this lemma, we will need several preliminary results, which we first discuss. The following result shows that $\qres_{\varep, N}^A$ has a total derivative structure. The claimed estimate of $\qres_{\varep, \star}^A$ in Lemma~\ref{lemma:qres-estimate-sol-space} will follow from this result combined with Schauder estimates and our previous result (Lemma~\ref{lemma:res-linear-linear}) estimating the resonance of the power-type quadratic object.

\begin{lemma}\label{lemma:qres-total-derivative}
For $\varep \in \sdyadic$, $N \in \dyadic$, and $(t, x) \in [0,T] \times \T^2$, we have that
\begin{equs}
(\qres_{\varep, N}^A)(t, x) &= \Duh_{B, N}^F\Big(\covd^\ell \E\big[ \linear[\oeps][][\ell], \linear[\oeps][][] \big]\Big)(t, x) \\
&= \int \ds \dy \, \fq_{B, N}(s, y; t, x) \covd^\ell \E\big[ \linear[\oeps][][\ell], \linear[\oeps][][] \big] (s, y).
\end{equs}
\end{lemma}

Next, we prove that $\linear[B][][]$ satisfies a sort of rotation invariance property in law. This will be the main input for the proof of Lemma~\ref{lemma:qres-total-derivative}.

\begin{lemma}\label{lemma:linear-object-rotation-invariance}
For any $\varep \in \sdyadic$, $\theta \in [0, 2\pi]$, we have that
\begin{equ}
\bigl(\linear[\oeps][][1], \linear[\oeps][][2]\bigr) \stackrel{d}{=} \bigl((\cos \theta) \linear[\oeps][][1] - (\sin \theta) \linear[\oeps][][2], (\sin \theta) \linear[\oeps][][1] + (\cos \theta) \linear[\oeps][][2]\bigr).
\end{equ}
\end{lemma}
\begin{proof}
For $\theta \in [0, 2\pi]$, define $R_\theta \in \End(\frkg^2)$ by
\begin{equ}
R_\theta(X_1, X_2) = ((\cos \theta) X_1 - (\sin \theta) X_2, (\sin \theta) X_1 + (\cos \theta) X_2),
\end{equ}
which in some sense acts by ``rotating'' the coordinates of $(X_1, X_2)$. Using $R_\theta$, define
\begin{equ}
\xi^\theta = R_\theta \xi = (\xi_1^\theta, \xi_2^\theta) = \bigl((\cos \theta) \xi_1 - (\sin \theta) \xi_2, (\sin \theta) \xi_1 + (\cos \theta) \xi_2\bigr).
\end{equ}
By a covariance calculation, one may check that $\xi^\theta \stackrel{d}{=} \xi$. Define $\linear[\oeps][\theta][]$ by
\begin{equ}
\linear[\oeps][\theta][](t, x) \eqdef \int \ds \dy \rho_\varep(t-s) \fq_B(s, y; t, x) \xi^\theta(s, y).
\end{equ}
We thus have that $\linear[\oeps][\theta][] \stackrel{d}{=} \linear[\oeps][][]$. On the other hand, one may check that (here we view $F \in \End(\frkg^2)$) $F R_\theta = R_\theta F$, i.e.\ the action of $F$ on $\frkg^2$ commutes with $R_\theta$. As a consequence, it also follows that $R_\theta$ commutes with the heat operator $(\ptl_t - \Delta_B - 2F)$, and thus we obtain
\begin{equ}
R_\theta \linear[\oeps] = \linear[\oeps][\theta][].
\end{equ}
By combining this with the fact that $\linear[\oeps][\theta][] \stackrel{d}{=} \linear[\oeps][][]$, we obtain the desired result.
\end{proof}

\begin{proof}[Proof of Lemma~\ref{lemma:qres-total-derivative}]
Applying Lemma~\ref{lemma:linear-object-rotation-invariance} in the case $\theta = \pi/2$, we obtain
\begin{equ}\label{eq:linear-object-90-degree-rotation-invariance}
\bigl(\linear[\oeps][][1], \linear[\oeps][][2]\bigr)  \stackrel{d}{=} \bigl(-\linear[\oeps][][2], \linear[\oeps][][1]\bigr).
\end{equ}
From the definition \eqref{local:eq-quadratic-N} of the localised quadratic object $\quadratic[\oeps][][N]$, to show the desired result, it suffices to show that for $i \in [2]$,
\begin{equ}
\E\Big[\linear[\oeps][\ell][], 2\covd_\ell \linear[\oeps][][i] - \covd_i \linear[\oeps][][\ell] \Big] = \covd^\ell \E\Big[\linear[\oeps][][\ell], \linear[\oeps][][i]\Big]. 
\end{equ}
Without loss of generality, take $i = 1$. We then have that
\begin{equ}
\Big[\linear[\oeps][\ell][], 2\covd_\ell \linear[\oeps][][1] - \covd_1\linear[\oeps][][\ell] \Big] = 2\Big[\linear[\oeps][][1], \covd_1 \linear[\oeps][][1] \Big] + 2\Big[\linear[\oeps][][2], \covd_2 \linear[\oeps][][1] \Big] - \Big[\linear[\oeps][][1], \covd_1 \linear[\oeps][][1]\Big] - \Big[\linear[\oeps][][2], \covd_1 \linear[\oeps][][2]\Big].
\end{equ}
Upon taking expectations, and using \eqref{eq:linear-object-90-degree-rotation-invariance}, we have that
\begin{equ}
\E\Big[\linear[\oeps][][1], \covd_1 \linear[\oeps][][1]\Big] + \E \Big[\linear[\oeps][][2], \covd_1 \linear[\oeps][][2]\Big] = 2 \E\Big[\linear[\oeps][][1], \covd_1 \linear[\oeps][][1]\Big],
\end{equ}
and thus observing a cancellation, we obtain
\begin{equ}
\E \Big[\linear[\oeps][\ell][], 2\covd_\ell \linear[\oeps][][1] - \covd_1\linear[\oeps][][\ell] \Big] = 2\E\Big[\linear[\oeps][][2], \covd_2 \linear[\oeps][][1] \Big].
\end{equ}
Thus to finish, it suffices to show that
\begin{equ}
2\E\Big[\linear[\oeps][][2], \covd_2 \linear[\oeps][][1] \Big]  = \covd_2\E\Big[\linear[\oeps][][2], \linear[\oeps][][1]\Big] = \covd^\ell \E\Big[\linear[\oeps][][\ell], \linear[\oeps][][1]\Big].
\end{equ}
(Here, the second identity is immediate, so the main point is to show the first identity.) Towards this end, we write
\begin{equ}
2\E\Big[\linear[\oeps][][2], \covd_2 \linear[\oeps][][1] \Big] = 2 \covd_2 \E\Big[\linear[\oeps][][2], \linear[\oeps][][1]\Big] - 2\E\Big[\covd_2 \linear[\oeps][][2], \linear[\oeps][][1]\Big].
\end{equ}
Thus it further suffices to show that
\begin{equ}
\covd_2 \E\Big[\linear[\oeps][][2], \linear[\oeps][][1]\Big] - 2\E\Big[\covd_2 \linear[\oeps][][2], \linear[\oeps][][1]\Big] = 0.
\end{equ}
Towards this end, we have that (applying \eqref{eq:linear-object-90-degree-rotation-invariance} in the second identity below)
\begin{equs}
\covd_2 \E\Big[\linear[\oeps][][2], \linear[\oeps][][1]\Big]  &= \E\Big[\covd_2 \linear[\oeps][][2], \linear[\oeps][][1]\Big] + \E\Big[\linear[\oeps][][2], \covd_2 \linear[\oeps][][1]\Big] \\
&= \E\Big[\covd_2 \linear[\oeps][][2], \linear[\oeps][][1]\Big] - \E\Big[\linear[\oeps][][1], \covd_2 \linear[\oeps][][2]\Big] \\
&= 2 \E\Big[\covd_2 \linear[\oeps][][2], \linear[\oeps][][1]\Big] ,
\end{equs}
as desired.
\end{proof}

Next, we prove the following estimate which serves as the key step towards the claimed estimate of $\qres_{\varep, \star}^\higgs$ in Lemma~\ref{lemma:qres-estimate-sol-space}. We prove a more general version of the needed estimate, which will also be needed later in Section~\ref{section:objects-paracontrolled}.

\begin{lemma}\label{lemma:covd-philinear-times-philinear-resonance}
Let $0 < T \leq \frac{1}{128}$. Let $\varep \in \sdyadic$ and $\Qc \colon \Y \times \Y \ra \Y$ be an equivariant bilinear map such that $\delta^{mn} \Qc(\ONBY_m, \ONBY_n) = 0$, where $(\ONBY_m)$ is any orthonormal basis of $\Y$. Then we have that
\begin{equ}
\big\| \abs{\log t}^{-1} t^{\roughdataweight/2} \E\big[ \Qc \bigl(\graylinear[\oeps] , \covd_B \graylinear[\oeps]\bigr)\big] \big\|_{L_{t, x}^\infty([0, T])} \lesssim \lossfactor_T(F)^3 T^{-(1-\roughdataweight)/2}. \label{eq:expectation-blinear-covd-philinear-philinear}
\end{equ}
Moreover, for any $1 \leq q < \infty$, we have that as $\varep \downarrow 0$, the quantity $\E\bigl[ \Qc \bigl(\graylinear[\oeps], \covd_B \graylinear[\oeps]\bigr)\bigr]$ converges in the weighted norm $\bigl\|t^{\roughdataweight/2} \place \bigr\|_{L_t^q L_x^\infty([0, T])}$ to a quantity which does not depend on the mollifier $\rho \in \mc{P}$.
\end{lemma}
\begin{proof}
Let $(t, x) \in [0, T] \times \T^2$, and let $(\ONBY_m)$ be an orthonormal basis of $\Y$. We have that
\begin{equ}
\E\big[ \Qc \bigl(\graylinear[\oeps], \covd_B \graylinear[\oeps]\bigr)\big](t, x) = \delta^{mn} \int \ds \dy \rho_\varep(t-s)^2 \Qc(\fq_{\Y, B} \ONBY_m , \covd_B \fq_{\Y, B}\ONBY_n)(s, y; t, x).
\end{equ}
Let $\eta = \fq_{\Y, B} - \fp \ptcutoff \pt_B$, $\omega = \covd_B \fq_{\Y, B} - (\nabla \fp) \ptcutoff \pt_B$. We may decompose the above into $I_1^\varep + I_2^\varep + I_3^\varep$, where
\begin{equs}
I_1^\varep &= \delta^{mn} \int \ds \dy \rho_\varep(t-s)^2 ((\nabla \fp) \fp \ptcutoff^2)(s, y; t, x) \Qc(\pt_B \ONBY_m , \pt_B \ONBY_n)(s, y; t, x), \\
I_2^\varep &= \delta^{mn} \int \ds \dy \rho_\varep(t-s)^2 \Qc( \fq_{\Y, B} \ONBY_m , \omega \ONBY_n)(s, y; t, x), \\
I_3^\varep &= \delta^{mn} \int \ds \dy \rho_\varep(t-s)^2 \Qc( \eta \ONBY_m, ((\nabla \fp)\ptcutoff \pt_B)\ONBY_n)(s, y; t, x).
\end{equs}
Due to our assumption that $\delta^{mn} \Qc(\ONBY_m, \ONBY_n) = 0$, we have that $I_1^\eps= 0$. To finish, it remains to estimate $I_2^\eps$ and $I_3^\eps$. We will focus on $I_2^\eps$; the estimate for $I_3^\eps$ is essentially the same. By the estimate \eqref{eq:qB-parametrix-t-s-power} from Corollary~\ref{cor:pB-qB-parametrix} (taking the $\alpha$ therein to be $\roughdataweight$, and the $L$ therein to be $\lossfactor_T(F)$), we have that
\begin{equs}
|\omega(s, y; t, x)| &\lesssim \lossfactor_T(F)\Big( \bigl(T^{-1}(t-s)\bigr)^{1-\roughdataweight} s^{-\roughdataweight} (t-s)^{\roughdataweight} + \Fnorm{F}_{[s, t],1} \\
&+ \bigl( \big\|\ptl_t B\big\|_{L_t^1 L_x^\infty([s, t])} + \big\|F\big\|_{L_t^2 L_x^\infty([s, t])} \bigr) (t-s)^{1/2}\Big) \times (t-s)^{-1/2} \moralfp(s, y; t, x).
\end{equs}
Inserting this bound into $I_2^\varep$, we obtain
\begin{equ}
\abs{I_2^\varep} \lesssim \lossfactor_T(F)^2 \int \ds \dy\, \bigl(T^{-1}(t-s)\bigr)^{1-\roughdataweight} s^{-\roughdataweight} (t-s)^{\roughdataweight} (t-s)^{-1/2} \moralfp(s, y; t, x)^2 + J^\varep,
\end{equ}
where
\begin{equ}
J^\varep =  \lossfactor_T(F)^2 \int \ds \dy \,\bigl(\big\|\ptl_t B\big\|_{L_t^1 L_x^\infty([s, t])} + \big\|F\big\|_{L_t^2 L_x^\infty([s, t])} + (t-s)^{-1/2} \Fnorm{F}_{[s, t], 1} \bigr) \moralfp(s, y; t, x)^2.
\end{equ}
By direct integration, the first term in $I_2^\varep$ may be estimated up to constants by (here we also use that $t \leq T$)
\begin{equ}
\lossfactor_T(F)^2 t^{\frac{1}{2} - \roughdataweight} T^{-(1-\roughdataweight)} \leq \lossfactor_T(F)^2 t^{-\roughdataweight/2} T^{-(1-\roughdataweight)/2},
\end{equ}
which is acceptable. Thus, it remains to estimate $J^\varep$. We split the integral into two regions: $s \in [0, t - t^{100}]$ and $s \in [t - t^{100}, t]$. In the first region, we will need the estimates (which follow by \eqref{eq:ptl-t-B-estimate-under-more-stringent-Fpert} and \eqref{eq:F-estimate-under-more-stringent-Fpert} with the $\alpha$ therein taken to be $\roughdataweight$, combined with the definition \eqref{local:eq-lossfactor} of $\lossfactor_T(F)$)
\begin{equs}
\bigl\|F\bigr\|_{L_t^2 L_x^\infty([s, t])} &\lesssim s^{-\roughdataweight/2} T^{-(1-\roughdataweight)/2} \lossfactor_T(F), \\
\bigl\|\ptl_t B\bigr\|_{L_t^1 L_x^\infty([s, t])} &\lesssim s^{-\roughdataweight/2} T^{-(1-\roughdataweight)/2} \lossfactor_T(F).
\end{equs}
We will also need an additional estimate, which follows by the definitions of $\Fnorm{F}$ (Definition~\ref{def:Fnorm}) and $\lossfactor_T(F)$ \eqref{local:eq-lossfactor}:
\begin{equ}
(t-s)^{-1/2} \Fnorm{F}_{[s, t], 1} \lesssim s^{-\roughdataweight/2} T^{-(1-\roughdataweight)/2} \lossfactor_T(F).
\end{equ}
Inserting these estimates, we may bound
\begin{equs}
\lossfactor_T(F)^2& \int_0^{t - t^{100}} \ds \int_{\T^2} \dy \,\bigl(\big\|\ptl_t B\big\|_{L_t^1 L_x^\infty([s, t])} + \big\|F\big\|_{L_t^2 L_x^\infty([s, t])} + (t-s)^{-1/2} \Fnorm{F}_{[s, t], 1} \bigr) \moralfp(s, y; t, x)^2 \\
&\lesssim \lossfactor_T(F)^3 T^{-(1-\roughdataweight)/2} \int_0^{t - t^{100}} \ds \, s^{-\roughdataweight/2} (t-s)^{-1} \lesssim \lossfactor_T(F)^3 T^{-(1-\roughdataweight)/2} t^{-\roughdataweight/2} \abs{\log t},
\end{equs}
which is acceptable. For the region $s \in [t - t^{100}, t]$, we will need the following estimates. First, by Remark~\ref{remark:curvature-pointwise-in-time-bound-Fnorm}, we have that for all $u \in [s, t]$,
\begin{equ}
\bigl\|F(u)\bigr\|_{L_x^\infty} \leq u^{-1} \Fnorm{F}_{u, 0} \leq s^{-1} \Fnorm{F}_{T, 0}.
\end{equ}
Next, since $F$ is a solution to the covariant heat equation (by Lemma~\ref{lemma:curvature-equation}), by the kernel estimate \eqref{eq:higher-order-pB-ptwise-bd}, and the same argument as above, we have that for all $u \in [s, t]$,
\begin{equ}
\bigl\|\covd F (u)\bigr\|_{L_x^\infty} \lesssim \bigl(1 + \Fnorm{F}_{T, 0}\bigr) u^{-1/2} \bigl\|F(u/2)\bigr\|_{L_x^\infty} \lesssim \bigl(1 +  \Fnorm{F}_{T, 0}\bigr) \Fnorm{F}_{T, 0} s^{-3/2}.
\end{equ}
From these estimates and the definition of $\Fnorm{F}_{[s, t], 1}$ (Definition~\ref{def:Fnorm}), we may also estimate
\begin{equ}
(t-s)^{-1/2} \Fnorm{F}_{[s, t], 1} \lesssim \bigl(1 +  \Fnorm{F}_{T, 0}\bigr) \Fnorm{F}_{T, 0} s^{-3/2} (t-s) + s^{-1} \Fnorm{F}_{T, 0} (t-s)^{1/2}.
\end{equ}
Inserting these estimates, we may bound (recall that $\ptl_t B = \schO(\covd F)$)
\begin{equs}
\lossfactor_T(F)^2 & \int_{t-t^{100}}^t \ds \int_{\T^2} \dy \,\bigl(\big\|\ptl_t B\big\|_{L_t^1 L_x^\infty([s, t])} + \big\|F\big\|_{L_t^2 L_x^\infty([s, t])} + (t-s)^{-1/2} \Fnorm{F}_{[s, t], 1} \bigr) \moralfp(s, y; t, x)^2  \\
&\lesssim \lossfactor_T(F)^2 \int_{t-t^{100}}^t \ds\, \bigl(s^{-3/2}  + s^{-1} (t-s)^{-1/2}\bigr)  \bigl(1 +  \Fnorm{F}_{T, 0} \bigr) \Fnorm{F}_{T, 0} \\
&\lesssim \lossfactor_T(F)^2 \bigl( t^{-3/2} t^{100} + t^{-1} t^{50} \bigr)  \bigl(1 +  \Fnorm{F}_{T, 0} \bigr) \Fnorm{F}_{T, 0}  \lesssim \lossfactor_T(F)^3,
\end{equs}
where we used that $t \leq T \leq 1/128$, so that $s \geq t - t^{100} \gtrsim t$.
This completes the proof of the estimate \eqref{eq:expectation-blinear-covd-philinear-philinear}.

For the convergence, note that for any $(t, x) \in (0, T] \times \T^2$, we have that $I_2^\varep \ra I_2^0$ and $I_3^\varep \ra I_3^0$, where $I_2^0, I_3^0$ are defined by replacing the $\rho_\varep(t-s)^2$ term in the definitions by $1$. Also note that $I_2^0, I_3^0$ manifestly do not depend on $\rho \in \mc{P}$. To obtain the convergence in the desired norm, we need to estimate the differences $I_2^\varep - I_2^0$ and $I_3^\varep - I_3^0$. Focusing on $I_2^\varep - I_2^0$, the argument proceeds similarly as for our estimate of $I_2^\varep$, except now the difference imposes that $t - s \lesssim \varep^2$. By repeating the same estimates as before, except now taking into account this restriction $t -s \lesssim \varep^2$, we may gain factors of $\varep$ except possibly in the case $t^{100} \leq \varep^2$ (say). Since this region shrinks as $\varep \downarrow 0$, we may obtain convergence in the desired weighted $L_t^q L_x^\infty$ norm. We omit the details.
\end{proof}

\begin{proof}[Proof of Lemma~\ref{lemma:qres-estimate-sol-space}]
We first discuss $\qres_{\varep, \star}^A$, which is simpler due to the total derivative structure provided by Lemma~\ref{lemma:qres-total-derivative}. Combining this lemma with the Schauder estimate \eqref{local:eq-covariant-Schauder-2} (applied with $\alpha = -1-2\kappa$, $\beta = -\kappa$, $r = \infty$) and the estimate \eqref{eq:summed-quadratic-resonance-estimate}, we obtain that for any $N \geq N_{\lo}$,
\begin{equs}
N^{2\kappa} \bigl\| \Duh_{B, N}^F(\covd^\ell \E [\linear[\oeps][][\ell], \linear[\oeps][][]])\bigr\|_{\sol_N^{1-2\kappa}} &\lesssim N^{2\kappa} \bigl\| \covd \E [\linear[\varep][][1], \linear[\varep][][2]]\bigr\|_{\nl_N^{-1-2\kappa}} \\
&\lesssim \lossfactor_T(F) N^{\kappa} \bigl\|\E [\linear[\varep][][1], \linear[\varep][][2]]\bigr\|_{\Cs_t^{-\frac{\kappa}{2}} L_x^\infty} \\
&\lesssim \lossfactor_T(F) N^{\kappa} \bigl\| \abs{\log t}^{-1} \E\bigl[ \linear[\varep][][1] \otimes \linear[\varep][][2]\bigr] \bigr\|_{L_{t,x}^\infty([0, T])} \\
&\lesssim  \lossfactor_T(F)^3 N^{\kappa}.
\end{equs}
Taking the factor of $N^{2\kappa}$ over to the right-hand side, and then taking a sup over $N \geq N_{\lo}$, we obtain an acceptable estimate of $\qres_{\varep, \star}^A$. The convergence follows similarly by combining a Schauder estimate with the convergence given by Lemma~\ref{lemma:res-linear-linear}.

For $\qres_{\varep, \star}^\higgs$, we again apply the Schauder estimate \eqref{local:eq-covariant-Schauder-2} (this time with $\alpha = -1 - 2\kappa$, $\beta = -3\kappa$, $r = \infty$) along with Lemma~\ref{lemma:covd-philinear-times-philinear-resonance} (with the $\Qc$ therein taken to be the bilinear map $\Qc(Y, \tilde{Y}) = \bil(Y^{(3)} \otimes \tilde{Y}^{(3)})$) to estimate, for any $N \geq N_{\lo}$,
\begin{equs}
N^{1+2\kappa} \bigl\|\Duh_{B, N}^F \bigl(\E\big[ \bil \bigl(\philinear[\oeps] \otimes &\covd_B \philinear[\oeps]\bigr)\big] \bigr)\bigr\|_{\sol_N^{1-2\kappa}} \\
&= N^{1+2\kappa} \bigl\| \E\big[ \bil \bigl(\philinear[\oeps] \otimes \covd_B \philinear[\oeps]\bigr)\big]\bigr\|_{\nl_N^{-1-2\kappa}} \\
&\lesssim \lossfactor_T(F) N^{3\kappa} \bigl\| \E\big[ \bil \bigl(\philinear[\oeps] \otimes \covd_B \philinear[\oeps]\bigr)\big]\bigr\|_{\Cs_t^{-\frac{3\kappa}{2}} L_x^\infty} \\
&\lesssim \lossfactor_T(F) N^{3\kappa} \bigl\|\abs{\log t}^{-1} t^{\roughdataweight/2} \E\big[ \bil \bigl(\philinear[\oeps] \otimes \covd_B \philinear[\oeps]\bigr)\big]\bigr\|_{L_{t,x}^\infty([0, T])} \\
&\lesssim \lossfactor_T(F)^4 N^{3\kappa} T^{-(1-\roughdataweight)/2}.
\end{equs}
Moving the $N^{1+2\kappa}$ factor from the left-hand side to the right, and using that $N \geq N_{\lo} \gtrsim T^{-1/2}$, we obtain an acceptable estimate of $\qres_{\varep, \star}^\higgs$. The convergence follows similarly by combining a Schauder estimate with the convergence given by Lemma~\ref{lemma:covd-philinear-times-philinear-resonance}.
\end{proof}

\subsection{Estimate of the enhanced data set}\label{section:estimate-enhanced-dataset-final}

In this section, we prove (most of) Proposition~\ref{object:prop-data-set}, i.e.\ we show that the size of the enhanced data set $\Xi_{B, \varep}$ defined in Definition~\ref{local:def-enhanced-data-set} satisfies a sub-Gaussian tail estimate. The most delicate estimates consist of the high$\times$high objects which are needed in the paracontrolled calculus; we defer those estimates to Section~\ref{section:objects-paracontrolled}. We prefer to structure the argument this way, so that we start from the simplest cases (e.g.\ linear object, power-type quadratic objects, etc.) and build towards the more complicated objects. For all but the high$\times$high objects needed in the paracontrolled calculus, most of the work has already been done in Sections~\ref{section:non-resonant-estimates} and~\ref{section:non-paracontrolled-resonant}, so this section mostly consists of putting all the arguments together. We begin with several abstract results, which reduce the problem of estimating the various objects in the appropriate spaces to the moment estimates of Section~\ref{section:non-resonant-estimates}.

\begin{lemma}\label{lemma:object-sobolev-abstract}
Let $T \in (0, 1]$ and let $X \colon [0, T] \times \T^2 \ra \genV$ be a stochastic process with smooth sample paths such that, for some constants $C_j, j = 0, \ldots, 3$ and some integer $k \geq 0$, we have that for $0 \leq m \leq 3$, all $1\leq p <\infty$, and all $(t, x) \in [0, T] \times \T^2$,
\begin{equ}
\big\| \covd^{(m)}_{t, x} X(t, x)\big\|_{L^p(\Omega)} \leq C_{m} p^{k/2}.
\end{equ}
Assume for simplicity that the constants $C_j$ are non-decreasing in $j$. Then for any $1 \leq q < \infty$, we have that
\begin{equ}
\Big\|\big\|X\big\|_{C_t^0 L_x^q([0, T])}\Big\|_{L^p(\Omega)} \lesssim_q C_0^{1-\tfrac{1}{q}} C_2^{\tfrac{1}{q}} p^{k/2},
\end{equ}
and moreover
\begin{equ}
\Big\| \big\|X\big\|_{C_t^0 L_x^\infty([0, T])} \Big\|_{L^p(\Omega)} \lesssim_q C_0^{1-\tfrac{3}{q}} C_3^{\tfrac{3}{q}} p^{k/2}.
\end{equ}
\end{lemma}
\begin{proof}
In the following, all space-time norms are implicitly on $[0, T] \times \T^2$. By Minkowski's integral inequality, for $p \geq q$, we have that for $0 \leq m \leq 3$,
\begin{equ}
\Big\|\big\|\covd^{(m)}_{t, x} X\big\|_{L_{t, x}^q} \Big\|_{L^p(\Omega)} \leq \Big\| \big\| \covd^{(m)}_{t, x} X(t, x)\big\|_{L^p(\Omega)} \Big\|_{L_{t, x}^q} \leq C_{m} p^{k/2} T^{\frac{1}{q}}.
\end{equ}
While if $p \leq q$, we may simply estimate
\begin{equ}
\Big\|\big\|\covd^{(m)}_{t, x} X\big\|_{L_{t, x}^q} \Big\|_{L^p(\Omega)} \leq \Big\|\big\|\covd^{(m)}_{t, x} X\big\|_{L_{t, x}^q} \Big\|_{L^q(\Omega)} \leq C_{m} q^{k/2} T^{\frac{1}{q}}.
\end{equ}
Combining, it follows that for all $1 \leq p < \infty$,
\begin{equ}\label{eq:process-sobolev-Ltxr-intermediate}
\Big\|\big\|\covd^{(m)}_{t, x} X\big\|_{L_{t, x}^q} \Big\|_{L^p(\Omega)} \lesssim_q C_{m} p^{k/2} T^{\frac{1}{q}}.
\end{equ}
By Sobolev embedding in time, we have that
\begin{equ}
\big\|\covd^{(m)}_{t, x} X\big\|_{C_t^0 L_x^q} \lesssim T^{-\frac{1}{q}}\big\|\covd^{(m)}_{t, x} X\big\|_{L_{t, x}^q}  + \big\|\covd^{(m)}_{t, x} X\big\|_{L_{t, x}^q}^{1-\tfrac{1}{q}} \big\|\ptl_t \covd^{(m)}_{t, x} X\big\|_{L_{t, x}^q}^{\tfrac{1}{q}}.
\end{equ}
Thus for $0 \leq m \leq 1$, we have by \eqref{eq:process-sobolev-Ltxr-intermediate} that (using that the $C_j$ are non-decreasing)
\begin{equ}
\Big\| \big\|\covd^{(m)}_{t, x} X\big\|_{C_t^0 L_x^q} \Big\|_{L^p(\Omega)} \lesssim \bigl(C_{m}^{1-\tfrac{1}{q}}\bigr) \bigl( C_{m+2}^{\tfrac{1}{q}}\bigr) p^{k/2}.
\end{equ}
When $m = 0$, this gives the first desired estimate.

For the second desired estimate it suffices to take (say) $q > 3$. We then apply the gauge covariant Gagliardo--Nirenberg--Sobolev inequality (Lemma~\ref{prelim:lem-GNS}) to obtain
\begin{equ}
\big\|X\big\|_{C_t^0 L_x^\infty} \lesssim \big\|X\big\|_{C_t^0 L_x^q} + \big\|X\big\|_{C_t^0 L_x^q}^{1-\tfrac{2}{q}} \big\|\covd X\big\|_{C_t^0 L_x^q}^{\tfrac{2}{q}}.
\end{equ}
Combining this with the previous estimate, we obtain
\begin{equ}
\Big\| \big\|X \big\|_{C_t^0 L_x^\infty} \Big\|_{L^p(\Omega)} \lesssim \Big( C_{0}^{1-\tfrac{1}{q}} C_{2}^{\tfrac{1}{q}}\Big)^{1-\tfrac{2}{q}} \Big( C_1^{1-\tfrac{1}{q}} C_3^{\tfrac{1}{q}} \Big)^{\tfrac{2}{q}} p^{k/2} \lesssim C_0^{1 - \tfrac{3}{q}} C_3^{\tfrac{3}{q}} p^{k/2},
\end{equ}
giving the second desired estimate.
\end{proof}

The following lemma will be needed for the objects which are estimated in the negative temporal regularity space (from Definition~\ref{def:temporal-regularity-space}).

\begin{lemma}\label{lemma:time-regularity-kolmogorov}
Let $T \in (0, 1]$, and let $X \colon [0, T] \times \T^2 \ra \genV$ be a continuous stochastic process such that there exists $\beta \in (0, 1]$, an integer $k \geq 0$, and a constant $C_0 \geq 0$ such that for all $x \in \T^2$, $0 \leq t_1 < t_2 \leq T$, and $1\leq p <\infty$, we have that
\begin{equs}
\big\| X(t_1, x) - X(t_2, x)\big\|_{L^p(\Omega)} &\leq C_0 |t_1 - t_2|^\beta p^{k/2}, \\
\big\|X(0, x)\big\|_{L^p(\Omega)} &\leq C_0 p^{k/2}.
\end{equs}
Then for any $1 \leq q < \infty$ and $\alpha \in (0, \beta)$, we have that for all $1 \leq p <\infty$,
\begin{equ}
\Big\|\big\|X \big\|_{\Cs_t^{\alpha} L_x^q([0, T])}\Big\|_{L^p(\Omega)} \lesssim_{\alpha, q} C_0 p^{k/2}.
\end{equ}
\end{lemma}
\begin{proof}
We first obtain that for fixed $1 \leq q < \infty$ (as in the proof of Lemma~\ref{lemma:object-sobolev-abstract}),
\begin{equs}
\Big\|\big\|X(t_1) - X(t_2)\big\|_{L_x^q}\Big\|_{L^p(\Omega)} &\lesssim_q C_0 |t_1 - t_2|^\beta p^{k/2}, \\
\Big\| \big\|X(0)\big\|_{L_x^q}\Big\|_{L^p(\Omega)} &\lesssim_q C_0 p^{k/2}.
\end{equs}
We conclude by Kolmogorov's continuity theorem (see e.g.\ \cite[Theorem 4.3.2]{Stroock2ed}).
\end{proof}

We now proceed to prove Proposition~\ref{object:prop-data-set}. In the following, recall from Assumption~\ref{local:assumption-B} that
\begin{equ}
\lossfactor_T(F) \leq T^{-\lossdel}.
\end{equ}
This will be used without comment at many points in the following.

\begin{proof}[Proof of Proposition~\ref{object:prop-data-set}]
We estimate the individual objects comprising the enhanced data set one-by-one. We will mainly focus on deriving estimates which are uniform in $\varep$. Along the way, we will also discuss convergence of the resonant parts, but we will defer the discussion of convergence of the non-resonant parts to the very end, because this may be done in a unified manner. In the following, we will obtain $p$th moment estimates with appropriately quantitative dependence on $p$. By standard results, this will imply the Gaussian tail bound claimed in the statement of Proposition~\ref{object:prop-data-set}.

\proofstep{Estimate of linear objects.} 
We need to estimate $\big\|\graylinear[\oeps][][][\ostar]\big\|_{\modsol^{-\kappa}}$, which recall from Definition~\ref{local:def-modified-solution-space} is equal to
\begin{equ}
\sup_{N \geq N_\lo} \Big( \big\|\graylinear[\oeps][][N]\big\|_{\sol^{-\kappa}_N([0, T])} + \sup_{-1 \leq \alpha' \leq -\kappa} \max_{0 \leq m \leq 10} N^{-\alpha' -\kappa - m} \big\|\covd^{(m)}_{t, x} \graylinear[\oeps][][N]\big\|_{\Cs_t^{\alpha'/2} L_x^r}\Big).
\end{equ}
Recall (from Definition~\ref{local:def-solution-space}) that
\begin{equ}
\big\|\graylinear[\oeps][][N]\big\|_{\sol^{-\kappa}_N([0, T])} = \max_{0 \leq m \leq 10} N^{-\kappa - m} \big\|\covd^{(m)}_{t, x} \graylinear[\oeps][][N]\big\|_{L_{t, x}^\infty([0, T])}.
\end{equ}
We first estimate this norm. By combining Lemmas~\ref{lemma:object-sobolev-abstract} and~\ref{object:lem-linear}, we have that for any $N \geq N_{\lo}$, $p \geq 1$,
\begin{equ}
\sup_{\varep \in \sdyadic} \Big \|\big\|\covd^{(m)}_{t, x} \graylinear[\oeps][][N]\big\|_{C_t^0 L_x^\infty}\Big\|_{L^p(\Omega)} \lesssim T^{-\lossdel} N^{m + \kappa_1} p^{1/2},
\end{equ}
where $\kappa_1$ is sufficiently small depending on $\kappa$ and $m$. 
Moreover, note that for $\varep \ll N^{-1}$, we have that $\graylinear[\oeps][][N]$ is constant in $\varep$. Therefore we may give up an additional $\log N$ factor to obtain
\begin{equs}
\Big\| \sup_{\varep \in \sdyadic} \sup_{N \geq N_{\lo}} \max_{0 \leq m\leq 10} &N^{-\kappa - m} \big\|\covd^{(m)}_{t, x} \graylinear[\oeps][][N] \big\|_{L_{t, x}^\infty} \Big\|_{L^p(\Omega)} \\
&= \Big\| \sup_{N \geq N_{\lo}} \sup_{\varep \in \sdyadic} \max_{0 \leq m \leq 10} N^{-\kappa - m} \big\|\covd^{(m)}_{t, x} \graylinear[\oeps][][N] \big\|_{L_{t, x}^\infty} \Big\|_{L^p(\Omega)} \\
&\lesssim T^{-\lossdel} p^{1/2} \sum_{N \geq N_{\lo}} (\log N) N^{\kappa_1 - \kappa} \lesssim T^{-\lossdel} p^{1/2},
\end{equs}
thus giving the desired estimate for the $\sol^{-\kappa}_N$ piece of the $\modsol^{-\kappa}$-norm. For the other piece, first let $X(t_0) = \int_0^{t_0} \dt\, \covd^{(m)}_{t, x} \graylinear[\oeps][][N]$ for $t_0 \in [0, T]$, and fix $\alpha' \in [-1, -\kappa]$. Let $\beta > 1 + \alpha'/2$. By Lemma~\ref{object:lem-linear}, we may obtain 
\begin{equs}
\big\|X(t_1, x) - X(t_2, x)\big\|_{L^p(\Omega)} &\lesssim T^{-\lossdel} p^{1/2} |t_1 - t_2|^{1/2} \min(|t_1- t_2|^{1/2},  N^{-1}) N^{m} \\
&\lesssim T^{-\lossdel} p^{1/2} N^{m - 2(1-\beta)} |t_1 - t_2|^\beta.
\end{equs}
Note also that $X(0) = 0$. Thus by taking $\beta = 1 + (\alpha' + \kappa_1)/2$, where $\kappa_1$ is sufficiently small compared to $\kappa$ so that $\beta < 1$, we may obtain by Lemma~\ref{lemma:time-regularity-kolmogorov} and the definition of $\Cs_t^{\alpha'/2} L_x^r$ (recall Definition~\ref{def:temporal-regularity-space}) that
\begin{equ}
\Big\| \big\|\covd^{(m)}_{t, x} \graylinear[\oeps][][N]\big\|_{\Cs_t^{\alpha'/2} L_x^r} \Big\|_{L^p(\Omega)} =  \Big\| \big\|X\big\|_{\Cs_t^{1+\frac{\alpha'}{2}} L_x^r} \Big\|_{L^p(\Omega)} \lesssim T^{-\lossdel} p^{1/2} N^{m + \alpha' + \kappa_1}. 
\end{equ}
Applying this with the endpoint values $\alpha' = -1$ and $\alpha' = -\kappa$, then using interpolation to estimate the intermediate values, and applying the same observation as before to take the sup in $\varep$ by only giving up an additional $\log N$ factor, we may estimate
\begin{equs}
\Big\| \sup_{\varep \in \sdyadic} \sup_{N \geq N_{\lo}}& \sup_{-1 \leq \alpha' \leq -\kappa} \max_{0 \leq m\leq 10} N^{-\alpha' -\kappa - m} \big\|\covd^{(m)}_{t, x} \graylinear[\oeps][][N]\big\|_{\Cs_t^{\alpha'/2} L_x^r}\Big\|_{L^p(\Omega)} \\
&\lesssim T^{-\lossdel} p^{1/2} \sum_{N \geq N_{\lo}} (\log N) N^{\kappa_1 - \kappa} \lesssim T^{-\lossdel} p^{1/2}.
\end{equs}
Combining, we thus obtain the desired estimate of $\graylinear[\oeps][][]$. We record the estimate here as it will be needed later in Section~\ref{section:objects-paracontrolled}:
\begin{equ}\label{eq:linear-object-solution-space-estimate}
\Bigl\| \sup_{\varep \in \sdyadic} \bigl\| \graylinear[\oeps][][][\ostar]\bigr\|_{\modsol^{-\kappa}} \Bigr\|_{L^p(\Omega)} \lesssim T^{-\lossdel} p^{1/2}.
\end{equ}
Having written out all the details for the estimate of the linear object, we will be a bit brief from here on, as many of the arguments are similar. In particular, we will always be able to obtain estimates that are 
uniform in $\varep$ by giving up a log factor of the maximum frequency. 

\proofstep{Estimate of power-type quadratic objects.} We want to estimate
\begin{equ}
\big\| \graylinear[\oeps][] \otimes \graylinear[\oeps][] -  \sigma_{\rho,\varep}^2 \Cov_{\Y} \big\|_{\Cs_t^{-\kappa} L_x^r([0, T])}.
\end{equ}
First, by Lemma~\ref{lemma:res-linear-linear}, we have that
\begin{equ}
\graylinear[\oeps][] \otimes \graylinear[\oeps][] -  \sigma_{\rho,\varep}^2 \Cov_{\Y} = \Wick{\graylinear[\oeps][] \otimes \graylinear[\oeps][]} + \bigO(T^{-2\lossdel} + \abs{\log t}),
\end{equ}
where $\bigO(T^{-2\lossdel} + \abs{\log t})$ denotes a deterministic quantity which is bounded by $\lesssim T^{-2\lossdel} + \abs{\log t}$ uniformly in $\varep$. This is already acceptable, and thus it suffices to estimate the non-resonant piece, i.e.\ the first term on the right-hand side above. We may further split
\begin{equ}
\Wick{\graylinear[\oeps][] \otimes \graylinear[\oeps][]}  = \sum_{N, M \lesssim \varep^{-1}} \Wick{\graylinear[\oeps][][N] \otimes \graylinear[\oeps][][M]} \;.
\end{equ}
Applying Lemmas~\ref{lemma:no-derivative-objects} and~\ref{lemma:time-regularity-kolmogorov} (and arguing similarly as for the linear objects), we obtain (taking $\kappa_1 \ll \kappa$)
\begin{equ}
\Big\| \big\|\Wick{\graylinear[\oeps][][N] \otimes \graylinear[\oeps][][M]}\big\|_{\Cs_t^{-\kappa} L_x^r} \Big\|_{L^p(\Omega)} \lesssim T^{-2\lossdel} p \max(N, M)^{2(\kappa_1 - \kappa)}.
\end{equ}
Using the fact that the quantity inside $\big\|\place\big\|_{L^p(\Omega)}$ on the left-hand side becomes 
constant in $\varep$ once $\varep \ll \max(N, M)^{-1}$, we may estimate
\begin{equs}
\Big\| \sup_{\varep \in \sdyadic} \big\| \Wick{\graylinear[\oeps][] \otimes \graylinear[\oeps][]}\big\|_{\Cs_t^{-\kappa} L_x^r} \Big\|_{L^p(\Omega)} &\leq \sum_{N, M \geq N_{\lo}} \Big\| \sup_{\varep \in \sdyadic}\big\| \Wick{\graylinear[\oeps][][N] \otimes \graylinear[\oeps][][M]} \big\|_{\Cs_t^{-\kappa} L_x^r} \Big\|_{L^p(\Omega)} \\
&\lesssim T^{-2\lossdel} p \sum_{N, M \geq N_{\lo}} (\log \max(N, M)) \max(N, M)^{2(\kappa_1 - \kappa)} \\
&\lesssim T^{-2\lossdel} p.
\end{equs}
This completes the estimate of the quadratic object. Convergence of the expectations of the power-type quadratic objects follows directly from the convergence statement of Lemma~\ref{lemma:res-linear-linear}. (As previously mentioned, we defer the discussion of the convergence of the non-resonant pieces until the very end.)

\proofstep{Estimate of power-type cubic objects.} The argument is quite similar to the one for the quadratic objects without derivatives, and thus the details are omitted.

\proofstep{Estimate of quadratic objects with derivative.} We focus on the object $\quadratic[\oeps]$, as the other quadratic objects with one derivative may be estimated quite similarly. We decompose
\begin{equ}
\quadratic[\oeps][][M] = \Wick{\quadratic[\oeps][][M]} + \E\Big[ \quadratic[\oeps][][M] \Big].
\end{equ}
The estimate of the resonant part directly follows from Lemma~\ref{lemma:qres-estimate-sol-space} (note that this piece is deterministic, and so we do not need to take moments). Convergence of the resonant part also follows from this lemma. For the non-resonant part, for fixed $M, N_1, N_2$, by combining Lemmas~\ref{lemma:quadratic-object-non-resonant-bound} and~\ref{lemma:object-sobolev-abstract}, we may first obtain (arguing similarly as for the linear object)
\begin{equs}
\sup_{\varep \in \sdyadic} \Big\| \Big\|\covd^{(m)}_{t, x} \Wick{\Duh_{B, M}^F\Big( \big[\linear[\oeps][j][N_1], 2 \covd_j \linear[\oeps][][N_2] - \covd &\linear[\oeps][][j, N_2]\big]\Big)} \Big\|_{C_t^0 L_x^\infty}\Big\|_{L^p(\Omega)} \\
&\lesssim T^{-3\lossdel} M^{m + \kappa_1} \max(M, N_1, N_2)^{-1} p,
\end{equs}
for some $\kappa_1 \ll \kappa$. Then arguing similarly as for the power-type quadratic objects, we may move the sup over $\varep \in \sdyadic$ inside the $L^p(\Omega)$ norm and then finally obtain the desired estimate of the non-resonant part in $\sol^{1-2\kappa}$. 
\proofstep{Estimate of objects in $\Xi_{\Beps}^{\der}$ (i.e.\ all terms in \eqref{local:eq-data-set-der-1}--\eqref{local:eq-data-set-der-3}).}
Corollary~\ref{objects:cor-xi-der} provides the needed estimate as well as convergence. As discussed at the beginning of this section, the corollary and its proof are deferred to Section~\ref{section:objects-paracontrolled}.

\proofstep{Estimate of objects in $\Xi_{\Beps}^{\para}$ (i.e.\ all terms in \eqref{local:eq-data-set-para}).} For the resonant part, Proposition~\ref{object:prop-dduhlin-dlin-main} provides the needed estimate as well as convergence (we remark that $L_t^q L_x^\infty$ embeds into $\Cs_t^{-\theta} L_x^r$ for large enough $q$, thus the $L_t^q L_x^\infty$ estimates provided by the proposition are sufficient). As discussed at the beginning of this section, the proposition and its proof are deferred to Section~\ref{section:objects-paracontrolled}. In the following, we estimate the non-resonant part, which for fixed $\varep \in \sdyadic$ and a fixed choice of indices $k, \ell \in [2]$, may be written as the sum:
\begin{equ}
\sum_{\substack{M, N, L \\ M \sim_\theta L}} Z_{M, N, L}^\varep,
\end{equ}
where 
\begin{equ}
Z_{M, N, L}^\varep(t, x) \eqdef \Wick{\covd^k \Duh_{\Y, B, M}\bigl(\graylinear[\oeps][][N]\bigr) \otimes \covd^\ell \graylinear[\oeps][][L]}(t, x).
\end{equ}
Fix $M, N, L \geq N_{\lo}$ for now. In order to estimate the $\Cs_t^{-\theta} L_x^r$ norm of each term, we form
\begin{equ}
X_{M, N, L}^\varep(t, x) \eqdef \int_0^t \ds\, Z_{M, N, L}^\varep(s, x).
\end{equ}
By Lemma~\ref{lemma:xi-para-non-resonant-estimates} and interpolation, for any $t_1, t_2 \in [0, T]$ and $x \in \T^2$, we have that
\begin{equs}
\bigl\|&X_{M, N, L}^\varep (t_1, x) - X_{M, N, L}^\varep (t_2, x)\bigr\|_{L^p(\Omega)} \\
&\lesssim T^{-3\lossdel} p |t_1 - t_2|^{1/2} \min\bigl(|t_1 - t_2|^{1/2}, \max(\min(M, N), L)^{-1} \bigr)\max(M, N)^{-1} L \\
&\lesssim T^{-3\lossdel} p |t_1 - t_2|^{1/2} |t_1 - t_2|^{(1/2)(1-3\theta/2)} \max(\min(M, N), L)^{-3\theta/2} \max(M, N)^{-1} L \\
&\lesssim T^{-3\lossdel} p |t_1 - t_2|^{1-3\theta/4} L^{1-3\theta/2} \max(M, N)^{-1}.
\end{equs}
From this, we further obtain
\begin{equ}
\Bigl\| \bigl\|X_{M, N, L}^\varep (t_1) - X_{M, N, L}^\varep(t_2)\bigr\|_{L_x^r} \Bigr\|_{L^p(\Omega)} \lesssim T^{-3\lossdel} p |t_1 - t_2|^{1-3\theta/4} L^{1-3\theta/2} \max(M, N)^{-1},
\end{equ}
and thus by the definition of the $\Cs_t^{-\theta} L_x^r$ norm (Definition~\ref{def:temporal-regularity-space}) and Lemma~\ref{lemma:time-regularity-kolmogorov}, we obtain
\begin{equ}
\Bigl\|\bigl\| Z_{M, N, L}^\varep \bigr\|_{\Cs_t^{-\theta} L_x^r}\Bigr\|_{L^p(\Omega)} = \Bigl\|\bigl\|X_{M, N, L}^\varep \bigr\|_{\Cs_t^{1-\theta} L_x^r} \Bigr\|_{L^p(\Omega)} \lesssim T^{-3\lossdel} p L^{1-3\theta/2} \max(M, N)^{-1}.
\end{equ}
Now by summing over $M, N, L \geq N_{\lo}$ subject to the constraint $M \sim_\theta L$ (which implies that $M \geq L^{1-\theta}$), we may obtain
\begin{equs}
\sum_{\substack{M, N, L \geq N_{\lo} \\ M \sim_\theta L}} \sup_{\varep \in \sdyadic}  \Bigl\| \bigl\|Z_{M, N, L}^\varep\bigr\|_{\Cs_t^{-\theta}L_x^r} \Bigr\|_{L^p(\Omega)} &\lesssim \sum_{\substack{M, N, L \geq N_{\lo} \\ M \sim_\theta L}} T^{-3\lossdel} p L^{1-3\theta/2} \max(M, N)^{-1} \\
&\lesssim T^{-3\lossdel} p.
\end{equs}
Finally, by similar arguments as before, we may move the sup over $\varep \in \sdyadic$ inside the $L^p(\Omega)$ norm, which results in the desired estimate of the non-resonant part.

\proofstep{Convergence of all non-resonant pieces.} Finally, we discuss convergence of all non-resonant pieces. Abstractly, each non-resonant piece is of the following general form. There is some collection $Z^\varep(N_1, \ldots, N_k)$ of functions depending on $\varep \in \sdyadic$ as well as $k$ frequency scales $N_1, \ldots, N_k \in \dyadic$, such that the non-resonant piece for $\varep \in \sdyadic$ may be written
\begin{equ}
\sum_{N_1, \ldots, N_k \geq N_{\lo}} Z^\varep(N_1, \ldots, N_k).
\end{equ}
In the estimate of each non-resonant piece, we have shown that, taking $\|\place\|$ to be the appropriate norm for each piece,
\begin{equ}
\sum_{N_1, \ldots, N_k \geq N_{\lo}} \sup_{\varep \in \sdyadic} \bigl\| Z^\varep(N_1, \ldots, N_k)\bigr\| \stackrel{a.s.}{<} \infty.
\end{equ}
(Here we apply our moment estimates with $p = 1$, and we use that if a non-negative random variable has finite first moment, then it is finite a.s.) Moreover, by direct verification, we may obtain that for each choice of $N_1, \ldots, N_k \geq N_{\lo}$, the limit $\lim_{\varep \downarrow 0} Z^\varep(N_1, \ldots, N_k)$ exists (with respect to the norm $\|\place\|$). Thus, convergence of the non-resonant pieces as $\varep \downarrow 0$ now follows by dominated convergence.
\end{proof}

Next, we quickly show the soft statement that the linear object $\graylinear[\oeps][][]$ almost surely converges as $\varep \downarrow 0$ in $C_t^0 \Cs_x^{-\kappa}([0, T])$. This was needed in the proof of Theorem~\ref{main:thm-covariant-local} at the end of Section~\ref{section:local-proof}. Since the main focus of this paper lies in the more delicate quantitative estimates (which for the enhanced data set are provided by Proposition~\ref{object:prop-data-set}), we will be a bit brief in showing the following soft result.

\begin{lemma}\label{lemma:linear-object-convergence-soft}
In the setting of Proposition~\ref{object:prop-data-set}, as $\varep \downarrow 0$, we have that $\graylinear[\oeps][][]$ almost surely converges in $C_t^0 \Cs_x^{-\kappa}([0, T])$.
\end{lemma}
\begin{proof}
As in the proof of Proposition~\ref{object:prop-data-set}, it suffices to show that
\begin{equ}
\sum_{N \in \dyadic} \Bigl\| \sup_{\varep \in \sdyadic} \bigl\| \graylinear[\oeps][][][N]\bigr\|_{C_t^0 \Cs_x^{-\kappa}([0, T])} \Bigr\|_{L^1(\Omega)} < \infty. 
\end{equ}
Again as in the proof of Proposition~\ref{object:prop-data-set}, it suffices to fix $N \in \dyadic$ and prove an estimate such as
\begin{equ}
\sup_{\varep \in \sdyadic} \Bigl\| \bigl\| \graylinear[\oeps][][][N]\bigr\|_{C_t^0 \Cs_x^{-\kappa}([0, T])} \Bigr\|_{L^1(\Omega)}\lesssim N^{-\kappa'},
\end{equ}
for some $\kappa' > 0$. Here, the implicit constant may now include $\lossfactor_T(F)$ factors, which we do not need to track since we are just trying to show a soft finiteness result. Towards this end, let $K \in \dyadic$ be another frequency scale, and let $P_K$ be the corresponding Littlewood--Paley projection. It suffices to show the following estimates for $(t, x) \in [0, T] \times \T^2$:
\begin{equs}
\bigl\| P_K \graylinear[\oeps][][][N] (t, x)\bigr\|_{L^2(\Omega)} &\lesssim \min(1, K N^{-1}), \\
\bigl\| \ptl_t P_K \graylinear[\oeps][][][N] (t, x) \bigr\|_{L^2(\Omega)} &\lesssim N^2 \min(1, K N^{-1}).
\end{equs}
Indeed, Gaussian hypercontractivity would then give $L_{t,x}^p$ estimates for $p$ arbitrarily large, and then Bernstein's inequality would allow us to obtain $L_t^p L_x^\infty$ estimates by giving up a tiny $K^{2/p}$ factor (which would be acceptable by taking $p$ large enough depending on $\kappa$). Finally, we may interpolate the two estimates to obtain regularity in time, i.e.\ the desired $C_t^0 L_x^\infty$ estimate. Both inequalities follow directly from our kernel estimates (Lemma~\ref{lemma:qB-mixed-all-derivatives}) and standard techniques (i.e.\ Young's convolution inequality).
\end{proof}

We end this subsection with a uniform-in-time version of Lemma~\ref{le:lollipop_negative_holder}. This estimate will be used in the global iteration argument. We put it here as the proof proceeds by a Kolmogorov continuity argument similar to what we have recently seen. 

\begin{lemma}
\label{le:lollipop_negative_holder_uniform}
Let $C_B \geq 1$ and $T, S \in (0, 1]$, and suppose that $B$ is a YM gradient flow on $[0, T + S]$ satisfying $\lossfactor_{T+S}(F_B) \leq C_B(T+S)^{-\lossdel}$, with some constant $\lossdel$ small enough. Then, for all $p \geq 2$, it holds that
\begin{equ}
\label{e:lollipop_negative_holder_uniform}
	\Big\|S^{-\frac{\kappa}{5}}\sup_{0 < t \leq T}
	\|\linear[B](t)\|_{\covHolder{B(t), F, S}^{-\kappa}}\Big\|_{L^p(\Omega)}
	\lesssim p^{1/2} \;,
\end{equ}
where the implicit constant may depend on $C_B$. The same estimate holds for $\philinear[B]$, with $\covHolder{B(t), S}^{-\kappa}$ in place of $\covHolder{B(t), F, S}^{-\kappa}$.
\end{lemma}
\begin{proof}
	We only discuss $\linear[B]$, as the proof for $\philinear[B]$ is the same. We first reduce the estimate of the covariant H\"older norm to the zeroth order estimate. Indeed, by \eqref{eq:qB-ptwise-bounds}, for $0 < t \leq T$ and $0 < s \leq S$,
	\begin{equ}
		s^{\frac{\kappa+1}{2}}\|\covd_{B(t+s)}e^{s\Delta^F_{B(t)}}\linear[B](t)\|_{L_x^\infty}
		\lesssim S^{-\lossdel} s^{\frac{\kappa}{2}}\|e^{\frac{s}{2}\Delta^F_{B(t)}}\linear[B](t)\|_{L_x^\infty} \;.
	\end{equ}
	Thus the first order estimate follows from the zeroth order estimate and the fact that $\lossdel \ll \kappa$. Without loss of generality, let us assume $S \in (0, 1] \cap \sdyadic$. It remains to show that for any dyadic $R \in (0, S/2] \cap \sdyadic$, and all sufficiently large $p \geq 2$,
	\begin{equ}
	\label{e:lollipop_negative_holder_uniform_mid}
		\Big\|\sup_{\substack{0<t\leq T\\ R\leq s\leq 2R}}
		s^{\frac{\kappa}{2}}\|e^{s\Delta^F_{B(t)}}\linear[B](t)\|_{L_x^\infty}\Big\|_{L^p(\Omega)}
		\lesssim p^{1/2} R^{\frac{\kappa}{100}}S^{\frac{\kappa}{4}} \;.
	\end{equ}
	Indeed, after dyadically summing over $R$, \eqref{e:lollipop_negative_holder_uniform_mid} gives \eqref{e:lollipop_negative_holder_uniform}.

	We prove \eqref{e:lollipop_negative_holder_uniform_mid} by applying Kolmogorov's continuity theorem to the process
	\begin{equ}
		\fX(s, t) = \big\| s^{\frac{\kappa}{2}}e^{s\Delta^F_{B(t)}}\linear[B](t)\big\|_{L_x^\infty}
	\end{equ}
	on the region $\{ (s, t): s \in [R, 2R], t \in [0, T] \}$. It suffices to show
	\begin{equ}
	\label{e:lollipop_negative_holder_uniform_Kolmogorov1}
		\sup_{s \in [R, 2R], 0 \leq t \leq T} \Big\| \fX(s, t) \Big\|_{L^p(\Omega)} \lesssim p^{1/2} R^{\frac{\kappa}{100}}S^{\frac{\kappa}{3}} \;,
	\end{equ}
	and
	\begin{equs}
	\label{e:lollipop_negative_holder_uniform_Kolmogorov2}
		\sup_{R \leq s_1 < s_2 \leq 2R, 0 \leq t \leq T} \Big\| \fX(s_1, t) - \fX(s_2, t) \Big\|_{L^p(\Omega)}
		&\lesssim p^{1/2} R^{\frac{\kappa}{100}}S^{\frac{\kappa}{3} - \lossdel} \Bigl(\frac{|s_2-s_1|}{R}\wedge 1\Bigr)^{1/2} \;, \\
	\label{e:lollipop_negative_holder_uniform_Kolmogorov3}
		\sup_{\substack{s_i \in [R, 2R],\ 0 \leq t_i \leq T\\ s_1+t_1=s_2+t_2}} \Big\| \fX(s_1, t_1) - \fX(s_2, t_2) \Big\|_{L^p(\Omega)}
		&\lesssim p^{1/2} R^{\frac{\kappa}{100}}S^{\frac{\kappa}{3}}\Bigl(\frac{|s_2-s_1|}{R}\wedge 1\Bigr)^{1/2} \;.
	\end{equs}

	The bound \eqref{e:lollipop_negative_holder_uniform_Kolmogorov1} follows from a similar argument as in Lemma~\ref{le:lollipop_negative_holder}. We omit the details here.
	For the bound \eqref{e:lollipop_negative_holder_uniform_Kolmogorov2} on the increment in $s$, first note that, before inserting the weights $s_i^{\frac{\kappa}{2}}$, we have
	\begin{equ}
		e^{s_2\Delta^F_{B(t)}}\linear[B](t) - e^{s_1\Delta^F_{B(t)}}\linear[B](t)
		= \int_{s_1}^{s_2} \bigl( \Delta_{B(t+u)} + 2 \ad(F_{B(t+u)}) \bigr) e^{u\Delta^F_{B(t)}}\linear[B](t) \, \du \;.
	\end{equ}
	Since $u \in [R,2R]$, the monotonicity of the curvature gives
	\begin{equ}
		\|F_{B(t+u)}\|_{L_x^\infty} \lesssim u^{-1}\Fnorm{F_B}_{[t,t+u],0}
		\lesssim S^{-\lossdel}R^{-1} \;.
	\end{equ}
	Combining this bound with \eqref{eq:qB-ptwise-bounds} and \eqref{e:lollipop_negative_holder_uniform_Kolmogorov1}, we get
	\begin{equs}
		\,&\Big\| s_1^{\frac{\kappa}{2}}\big\| e^{s_2\Delta^F_{B(t)}}\linear[B](t) - e^{s_1\Delta^F_{B(t)}}\linear[B](t) \big\|_{L_x^\infty}\Big\|_{L^p(\Omega)} \\
		&\lesssim S^{-\lossdel}R^{-1}s_1^{\frac{\kappa}{2}}
		\int_{s_1}^{s_2}\Big\|\big\| e^{\frac{u}{2}\Delta^F_{B(t)}}\linear[B](t)\big\|_{L_x^\infty}
		+ \big\| e^{u\Delta^F_{B(t)}}\linear[B](t)\big\|_{L_x^\infty}\Big\|_{L^p(\Omega)}\du \\
		&\lesssim p^{1/2} R^{\frac{\kappa}{100}}S^{\frac{\kappa}{3}-\lossdel}\frac{|s_2-s_1|}{R} \;.
	\end{equs}
	Finally, since $s_1,s_2 \in [R,2R]$, we have $|s_2^{\frac{\kappa}{2}} - s_1^{\frac{\kappa}{2}}|\lesssim R^{-1+\frac{\kappa}{2}} |s_2-s_1|$.
	This concludes the proof of \eqref{e:lollipop_negative_holder_uniform_Kolmogorov2}.

	Now we turn to \eqref{e:lollipop_negative_holder_uniform_Kolmogorov3}. Assume for instance that $t_1<t_2$. Since the final time $s_1+t_1=s_2+t_2$ is the same, the definition of $\linear[B]$ gives
	\begin{equ}
		e^{s_2\Delta^F_{B(t_2)}}\linear[B](t_2) - e^{s_1\Delta^F_{B(t_1)}}\linear[B](t_1)
		= \int_{t_1}^{t_2}\int_{\T^2}\fq_B(u,y;s_1+t_1,\place)\xi(u,y)\du\dy \;.
	\end{equ}
	By It\^o's isometry and \eqref{eq:qB-ptwise-bounds}, for every $x \in \T^2$,
	\begin{equs}
		\E \Big| s_2^{\frac{\kappa}{2}} \int_{t_1}^{t_2}\int_{\T^2}\fq_B(u,y;s_1+t_1,x)\xi(u,y)\du\dy \Big|^2
		&\lesssim S^{-2\lossdel}R^\kappa \int_{t_1}^{t_2}(s_1+t_1-u)^{-1}\du \\
		&\lesssim S^{-2\lossdel}R^\kappa\Bigl(\frac{|t_2-t_1|}{R}\wedge 1\Bigr) \;.
	\end{equs}
	The same heat kernel smoothing argument as in the proof of Lemma~\ref{le:lollipop_negative_holder}, together with Gaussian hypercontractivity, gives
	\begin{equ}
		\Big\| s_2^{\frac{\kappa}{2}}\Big\|\int_{t_1}^{t_2}\int_{\T^2}\fq_B(u,y;s_1+t_1,\place)\xi(u,y)\du\dy\Big\|_{L_x^\infty}\Big\|_{L^p(\Omega)}
		\lesssim p^{1/2} R^{\frac{\kappa}{100}}S^{\frac{\kappa}{3}}\Bigl(\frac{|t_2-t_1|}{R}\wedge 1\Bigr)^{1/2} \;.
	\end{equ}
	It remains to compare the weights, which is the same argument as for \eqref{e:lollipop_negative_holder_uniform_Kolmogorov2}. This proves \eqref{e:lollipop_negative_holder_uniform_mid} for all sufficiently large $p$. The remaining values of $p$ follow from Jensen's inequality.
\end{proof}

\subsection{\protect{Estimates of high$\times$high-products in paracontrolled calculus}}\label{section:objects-paracontrolled}

In this subsection, we state and prove stochastic estimates needed in our analysis of the paracontrolled terms (see Section \ref{section:local-paracontrolled}), as well as related estimates for interactions between quadratic objects and covariant derivatives of linear objects. These estimates will be split over Propositions \ref{object:prop-dduhlin-dlin-main} and \ref{objects:prop-quadratic-dlin}. 

\begin{proposition}\label{object:prop-dduhlin-dlin-main}
Let Assumption~\ref{local:assumption-B} be satisfied and let $k,\ell\in [2]$. Then, it holds for all $1< q <\infty$, $\eps \in \sdyadic$, and $\rho\in \mc{P}$ that
\begin{equation}\label{object:eq-dduhlin-dlin-main}
\begin{aligned}
\Big\|&\sum_{\substack{K,L \colon \\ K \sim_\theta L}} \E \Big[ \covd_B^k \Duh_{\Y,B,K}\bigl( \graylinear[\oeps][][][]\bigr)(t,x)
\otimes \covd_B^\ell \graylinear[\oeps][][][L](t,x)\Big] 
- \tfrac{1}{4}   \delta^{k\ell} \bigl( \sigma_{\rho,\eps}^2 + \Crho \bigr) \Cov_{\Y} \Big\|_{L_t^q L_x^\infty} \\  &\lesssim_q \lossfactor_T(F)^3.
\end{aligned}
\end{equation}
Furthermore, as $\eps\rightarrow 0$, the argument in \eqref{object:eq-dduhlin-dlin-main} converges to a limit in $L_t^q L_x^\infty$ that does not depend on $\rho$, and the distance to the limit can be bounded by $\lossfactor_T(F)^3 \eps^{2\nu}$, where $0<\nu<\min(1/q,\theta)$. 
\end{proposition}

While the convergence rate in Proposition \ref{object:prop-dduhlin-dlin-main} will be useful in the proof of Lemma \ref{objects:lem-quadratic-dlin-23}, we encourage the reader to ignore it on first reading. We emphasise that the non-resonant part of $\covd^k \Duh_{\Y,B,K}\bigl( \graylinear[\oeps][][][]\bigr) \otimes \covd^\ell \graylinear[\oeps][][][L]$, which is easier to control,  has already been bounded in Lemma \ref{lemma:xi-para-non-resonant-estimates}. 

\begin{remark}
In the estimates from Sections \ref{section:non-resonant-estimates} and \ref{section:non-paracontrolled-resonant}, we always saw a gain in the highest frequency scale. Due to this, it is clear that the limit as $\eps\rightarrow 0$ does not depend on the cutoff function $\rho\in \mc{P}$, which was used in the definition of our covariant, smooth approximation (see Definition~\ref{main:def-gauge covariant-approximation-space-time-white-noise}). In this subsection, however, some estimates will not involve a gain in the highest frequency scale (see the proof of Lemma~\ref{object:lem-explicit-integrals}). Because of this, we will be careful and show that the $\rho$-dependence will only appear through $\Crho$ (as well as $\sigma_{\rho,\eps}^2$). Using this, one can show that the limiting local solution obtained in Theorem~\ref{main:thm-covariant-local} does not depend on $\rho$. Strictly speaking, it is not necessary to worry about the $\rho$-dependence here, since the fact that the limiting solution does not depend on $\rho$ also follows from our arguments in Section~\ref{section:langevin}. However, we think it is good to know that this can also be obtained as part of our covariant estimates. 
\end{remark}

\begin{proposition}\label{objects:prop-quadratic-dlin}
Let Assumption \ref{local:assumption-B} be satisfied. Furthermore, let $k,\ell \in [2]$,  let $\Qc \colon \Y \times \Y \rightarrow \Y$ be a $\liegroup$-equivariant bilinear map satisfying $\delta^{mn} \Qc(\ONBY_m, \ONBY_n) = 0$, where $(\ONBY_m)_m$ is an orthonormal basis of $\Y$, and let $\mathcal{K}$ be as in Proposition \ref{local:prop-paracontrolled}. Then, it holds for $p\geq 2$ that
\begin{equation}\label{objects:eq-quadratic-dlin-main}
\Bigl\| \sup_{\eps \in \sdyadic} \,  \Bigl \| \Duh^F_{\Y,B,\star} \bigl(  \Qc\bigl(\graylinear[\oeps][][][] , \covd^k \graylinear[\oeps][][][] \bigr)\bigr) \parasim \covd^\ell \graylinear[\oeps] - \tfrac{1}{4} \delta^{k\ell} \bigl( \sigma_{\rho,\eps}^2 + \Crho \bigr) \mathcal{K}(\graylinear[\oeps][][][]) \Bigr\|_{\Cs_t^{-\theta}L_x^r} 
\Bigr\|_{L^p(\Omega)} 
\lesssim \lossfactor_T(F)^4 p^{\frac{3}{2}}. 
\end{equation}
Furthermore, the argument in \eqref{objects:eq-quadratic-dlin-main} has an a.s.\ limit in $\Cs_t^{-\theta}L_x^r$ as $\eps \rightarrow 0$. 
\end{proposition}

The reason Proposition \ref{objects:prop-quadratic-dlin} is placed in the same subsection as Proposition \ref{object:prop-dduhlin-dlin-main} is that the most difficult term appearing in \eqref{objects:eq-quadratic-dlin-main} is the resonant interaction between the second and third linear stochastic objects, which is strongly related to the expectation in \eqref{object:eq-dduhlin-dlin-main}.  Proposition \ref{objects:prop-quadratic-dlin} will later be used to control (renormalised versions of) the objects
\begin{equation*}
\Lb\bigl( \quadratic[\oeps][][][\ostar]+\phiquadratic[\oeps][][][\ostar],\covd \linear[\oeps][][][\ostar]\bigr), 
\quad \bil \bigl( \mixedquadratic[\oeps][][][\ostar] \parasim \covd \philinear[\oeps][][][\ostar]\bigr), 
\quad \bigl( \quadratic[\oeps][][][\ostar] + \phiquadratic[\oeps][][][\ostar]\bigr)^j \parasim \covd_j \philinear[\oeps][][][\ostar], 
\quad \text{and} \quad \bigl( \covd_j \linear[\oeps][j][][\ostar] \bigr) \parasim    \mixedquadratic [\oeps][][][\ostar], 
\end{equation*}
which appear in \eqref{local:eq-Y-final-new-1}--\eqref{local:eq-Y-final-new-1p} and \eqref{local:eq-y-final-new-1}. The rest of the subsection is devoted to the proofs of Propositions \ref{object:prop-dduhlin-dlin-main} and \ref{objects:prop-quadratic-dlin}, and we begin with the former. Its proof is split over the following two lemmas. In the first lemma, we compute the expectation in \eqref{object:eq-dduhlin-dlin-main} and then approximate all resulting $\fq_{\Y,B}$-factors with the heat kernel parametrix. 

\begin{lemma}[Using the heat kernel parametrix]\label{object:lem-dduhlin-dlin}
Let the same assumptions as in Proposition~\ref{object:prop-dduhlin-dlin-main} be satisfied and, in addition, let $K,L,N\in \dyadic$ satisfy $K,L,N\geq N_{\lo}$. Then, it holds for all $(t,x)\in [0,T]\times \T^2$ that
\begin{equation}\label{object:eq-dduhlin-lin}
\begin{aligned}
&\quad \bigg|  \E \Big[ \covd_B^k \Duh_{\Y,B,K}\bigl( \graylinear[\oeps][][][N]\bigr)(t,x)
\otimes \covd_B^\ell \graylinear[\oeps][][][L](t,x)\Big] \\
&\qquad -  \Big( \iint \ds \dy \iint \dr \dz (\partial^k \fp_K)(r,z; t,x) \fp_{\eps,N}(s,y; r,z) 
(\partial^\ell \fp_{\eps,L})(s,y; t,x) \Big) \Cov_{\Y} \bigg| \\
&\lesssim \lossfactor_T(F)^3 T^{-\frac{1-\roughdataweightmultiple[4]}{2}} \max(K,L,N)^{-1+\roughdataweightmultiple[4]}.  
\end{aligned}
\end{equation}
\end{lemma}

We mention that the frequency scale $N$ in $\graylinear[\oeps][][][N]$ in \eqref{object:eq-dduhlin-lin}, which did not appear in \eqref{object:eq-dduhlin-dlin-main}, will be useful when establishing convergence as $\eps \rightarrow 0$. 

\begin{proof}[Proof of Lemma~\ref{object:lem-dduhlin-dlin}]
 We recall the integral representation of the covariant linear stochastic object from~\eqref{local:eq-linear-stochastic-integral-version} and \eqref{local:eq-graylinear-boldxi}. Using the definition of space-time white noise, the expectation in \eqref{object:eq-dduhlin-lin} can then be written as 
\begin{align}
&\quad   \E \Big[ \covd_B^k \Duh_{\Y,B,K}\bigl( \graylinear[\oeps][][][N]\bigr)(t,x)
\otimes \covd_B^\ell \graylinear[\oeps][][][L](t,x)\Big] \label{object:eq-dduhlin-lin-p1}\\ 
&=  \iint \ds \dy \iint \dr \dz \, \Big(  \bigl( \covd_B^k \fp_{\Y,B,K}(r,z;t,x) \fq_{\Y,B,\eps,N}(s,y;r,z) \bigr) \otimes \covd_B^\ell \fq_{\Y,B,\eps,L}(s,y;t,x) \Big) \Cov_{\Y}. \notag 
\end{align}
We note that \eqref{object:eq-dduhlin-lin-p1}, just like the second summand in \eqref{object:eq-dduhlin-lin}, can only be nonzero provided that $K\gtrsim L$, which we now assume for the rest of the proof. We now claim that \eqref{object:eq-dduhlin-lin-p1} can be approximated by 
\begin{align}
     &\iint \ds \dy \iint \dr \dz \, \bigg(  \widetilde{\ptcutoff}(x,y,z)   \bigl( \partial^k \fp_{K}\bigr)(r,z;t,x) \fp_{\eps,N}(s,y;r,z)  \bigl( \partial^\ell  \fp_{\eps,L}\bigr) (s,y;t,x) \notag \\
   &\times \Big( \tauYB(r,z;t,x) \tauYB(s,y;r,z) \otimes \tauYB(s,y;t,x) \Big) \bigg) 
   \Cov_{\Y},  \label{object:eq-dduhlin-lin-p2}
\end{align}
where  $\widetilde{\ptcutoff}(x,y,z) \eqdef \ptcutoff(x-z) \ptcutoff(z-y) \ptcutoff(x-y)$ and $\ptcutoff$ is as in Definition \ref{def:ptcutoff}. To see this, we use Proposition~\ref{prop:pB-ptwise-bounds} and Proposition~\ref{prop:qB-ptwise-bds}, which imply that $\fp_{\Y,B}$ and $\fq_{\Y,B}$ satisfy similar pointwise estimates as $\fp$ (except for the loss-factor). We also use Corollary~\ref{cor:pB-qB-parametrix} (with $\alpha=\roughdataweight$), which shows that $\fp_{\Y,B}$ and $\fq_{\Y,B}$ can be approximated using $\fp \tauYB$ (see also Remark \ref{remark:Fpert-less-time-weight}). By using this approximation for all three kernels in \eqref{object:eq-dduhlin-lin-p1}, we gain a factor of $\lossfactor_T(F)$ times either
\begin{align*}
T^{-\frac{1-\roughdataweightmultiple[4]}{2}}  r^{-\frac{\roughdataweightmultiple[2]}{2}} (t-r)^{\frac{1-\roughdataweightmultiple[2]}{2}}, \quad 
T^{-\frac{1-\roughdataweightmultiple[4]}{2}}  s^{-\frac{\roughdataweightmultiple[2]}{2}} (r-s)^{\frac{1-\roughdataweightmultiple[2]}{2}}, \quad \text{or} \quad 
T^{-\frac{1-\roughdataweightmultiple[4]}{2}}  s^{-\frac{\roughdataweightmultiple[2]}{2}} (t-s)^{\frac{1-\roughdataweightmultiple[2]}{2}}. 
\end{align*}
Since $r\geq s$ and $r-s,t-r \leq t-s$, the last of the three terms above is always the largest term. All in all, we therefore obtain that  
\begin{align}
\big|& \eqref{object:eq-dduhlin-lin-p1} - \eqref{object:eq-dduhlin-lin-p2}\big|
\lesssim \lossfactor_T(F)^3  T^{-\frac{1-\roughdataweightmultiple[4]}{2}} \iint \ds \dy \iint \dr \dz \bigg( \chi_K(t-r) \chi_N(r-s) \chi_L(t-s)   \notag \\
&\times s^{-\frac{\roughdataweightmultiple[2]}{2}} (t-s)^{\frac{1-\roughdataweightmultiple[2]}{2}}  (t-r)^{-\frac{3}{2}} (r-s)^{-1} (t-s)^{-\frac{3}{2}} \exp\Bigl( - c \Bigl( \tfrac{|x-z|^2}{t-r} + \tfrac{|z-y|^2}{r-s}+\tfrac{|x-y|^2}{t-s} \Bigr) \Bigr)\bigg). \label{object:eq-dduhlin-lin-p3}
\end{align}
To bound the integral in \eqref{object:eq-dduhlin-lin-p3}, one can first compute the $y$-integral and then the $z$-integral, which gain one power of $r-s$ and $t-r$, respectively. The remaining $(s,r)$-integral can then be bounded by
\begin{align*} 
&\lossfactor_T(F)^3 T^{-\frac{1-\roughdataweightmultiple[4]}{2}} \iint \ds \dr \chi_K(t-r) \chi_N(r-s) \chi_L(t-s) (t-r)^{-\frac{1}{2}} \times s^{-\roughdataweightmultiple[1]} (t-s)^{-1-\roughdataweightmultiple[1]} \\ 
&\lesssim \lossfactor_T(F)^3 T^{-\frac{1-\roughdataweightmultiple[4]}{2}}  \Bigl( \int \ds \chi_L(t-s) s^{-\roughdataweightmultiple[1]} (t-s)^{-1-\roughdataweightmultiple[1]} \Bigr) \sup_{s\in [0,T]}\Bigl( \int \dr \chi_K(t-r) \chi_N(r-s) (t-r)^{-\frac{1}{2}} \Bigr) \\ 
&\lesssim  \lossfactor_T(F)^3 T^{-\frac{1-\roughdataweightmultiple[4]}{2}}
 \times L^{\roughdataweightmultiple[4]} \times K \max(K,N)^{-2}. 
\end{align*}
Due to this, the error resulting from our approximation of \eqref{object:eq-dduhlin-lin-p1} by \eqref{object:eq-dduhlin-lin-p2} makes an acceptable contribution to \eqref{object:eq-dduhlin-lin}. We now claim that \eqref{object:eq-dduhlin-lin-p2} can be approximated by 
\begin{equation}\label{object:eq-dduhlin-lin-p4}
\begin{aligned}
   &\iint \ds \dy \iint \dr \dz \, \bigg(  \widetilde{\ptcutoff}(x,y,z)  \bigl( \partial^k \fp_{K}\bigr)(r,z;t,x) \fp_{\eps,N}(s,y;r,z)  \bigl( \partial^\ell  \fp_{\eps,L}\bigr) (s,y;t,x) \\
   &\times \Big( \tauYB(s,y;t,x)  \otimes \tauYB(s,y;t,x) \Big) \bigg) \Cov_{\Y}.  
\end{aligned}
\end{equation}
Going from \eqref{object:eq-dduhlin-lin-p2} to \eqref{object:eq-dduhlin-lin-p4}, we replaced $\tauYB(r,z;t,x) \tauYB(s,y;r,z)$ with $\tauYB(s,y;t,x)$. To obtain this approximation, we let $\eta>0$ remain to be chosen and use Lemma~\ref{lemma:parallel-transport-approximate-semigroup}, which leads to
\begin{align*}
&\quad \big| \tauYB(r,z;t,x) \tauYB(s,y;r,z)-\tauYB(s,y;t,x)\big| \\
&\lesssim \lossfactor_T(F) T^{-\frac{1-\roughdataweight}{2}} s^{-\frac{\roughdataweight}{2}} \Big( 1+\log\bigl(\tfrac{t-s}{r-s}\bigr) + \tfrac{|x-z|}{\sqrt{t-r}} \Big) \max\bigl(|x-y|,|y-z|\bigr) \\
&\lesssim_\eta \lossfactor_T(F) T^{-\frac{1-\roughdataweight}{2}} s^{-\frac{\roughdataweight}{2}} (t-s)^{\frac{1}{2}} \bigl( \tfrac{t-s}{r-s}\bigr)^{\eta}  \exp\Big( \eta \Big( \tfrac{|x-z|^2}{t-r} + \tfrac{|z-y|^2}{r-s}+\tfrac{|x-y|^2}{t-s} \Big)\Big). 
\end{align*}
In the last inequality, we also used that $t-r$ and $r-s$ are smaller than $t-s$. Together with the pointwise estimate for the Euclidean heat kernel on $\T^2$, it then follows that
\begin{align}
\big|& \eqref{object:eq-dduhlin-lin-p2} - \eqref{object:eq-dduhlin-lin-p4}\big|
\lesssim \lossfactor_T(F) T^{-\frac{1-\roughdataweight}{2}} \iint \ds \dy \iint \dr \dz \bigg( \chi_K(t-r) \chi_N(r-s) \chi_L(t-s)  s^{-\frac{\roughdataweight}{2}} (t-s)^{\frac{1}{2}}  \notag \\
&\times   \bigl( \tfrac{t-s}{r-s}\bigr)^{\eta} (t-r)^{-\frac{3}{2}} (r-s)^{-1} (t-s)^{-\frac{3}{2}} \exp\Bigl( - \bigl( \tfrac{1}{4}-\eta\bigr) \Big( \tfrac{|x-z|^2}{t-r} + \tfrac{|z-y|^2}{r-s}+\tfrac{|x-y|^2}{t-s} \Big) \Bigr)\bigg). \label{object:eq-dduhlin-lin-p5}
\end{align}
As for \eqref{object:eq-dduhlin-lin-p3}, bounding the $y$-integral and $z$-integral in \eqref{object:eq-dduhlin-lin-p5} gains us one power of $r-s$ and $t-r$. Provided that $\eta$ is smaller than $\kappa$, the remaining $(s,r)$-integral can then be bounded by
\begin{align*}
T^{-\frac{1-\roughdataweight}{2}}& \iint \ds \dr \chi_K(t-r) \chi_N(r-s) \chi_L(t-s) s^{-\frac{\roughdataweight}{2}} (t-r)^{-\frac{1}{2}} (r-s)^{-\eta} (t-s)^{-1+\eta} \\ 
&\lesssim T^{-\frac{1-\roughdataweight}{2}} K N^{2\eta} L^{2-2\eta} \times \max(K,N)^{-2} L^{-2+2\kappa} 
\lesssim T^{-\frac{1-\roughdataweight}{2}} \max(K,N)^{-1+2\kappa}. 
\end{align*}
Since $K \gtrsim T^{-\frac{1}{2}}$, the error resulting from our approximation of \eqref{object:eq-dduhlin-lin-p2} by \eqref{object:eq-dduhlin-lin-p4} makes an acceptable contribution to \eqref{object:eq-dduhlin-lin}. Next, we note that \eqref{object:eq-dduhlin-lin-p4} is equal to 
\begin{equation}\label{object:eq-dduhlin-lin-p6}
 \bigg(  \iint \ds \dy \iint \dr \dz \,  \widetilde{\ptcutoff}(x,y,z)  \bigl( \partial^k \fp_{K}\bigr)(r,z;t,x) \fp_{\eps,N}(s,y;r,z)  \bigl( \partial^\ell  \fp_{\eps,L}\bigr) (s,y;t,x) \bigg) \Cov_{\Y}.  
\end{equation}
For this, we now only need to use that $\tauYB(s,y;t,x)$ is an isometry, which implies that 
\begin{equation*}
\Big( \tauYB(s,y;t,x)  \otimes \tauYB(s,y;t,x) \Big)  \Cov_{\Y}= \Cov_{\Y}. 
\end{equation*}
In order to complete the proof, it now only remains to show that \eqref{object:eq-dduhlin-lin-p6} can be approximated by
\begin{equation}\label{object:eq-dduhlin-lin-p7}
 \bigg(  \iint \ds \dy \iint \dr \dz \, \bigl( \partial^k \fp_{K}\bigr)(r,z;t,x) \fp_{\eps,N}(s,y;r,z)  \bigl( \partial^\ell  \fp_{\eps,L}\bigr) (s,y;t,x) \bigg) \Cov_{\Y}.  
\end{equation}
In going from \eqref{object:eq-dduhlin-lin-p6} to \eqref{object:eq-dduhlin-lin-p7}, the only change is that the cutoff $\widetilde{\ptcutoff}(x,y,z)$ has been removed. For this, we will use that if $\widetilde{\ptcutoff}(x,y,z) \neq 1$, then $\max(|x-z|,|z-y|,|x-y|)\gtrsim 1$. We will also use the elementary estimate  
\begin{align*}
\ind \big\{& \max(|x-z|,|z-y|,|x-y|)\gtrsim 1 \big\} 
\exp\Bigl( -  \tfrac{1}{4} \Big( \tfrac{|x-z|^2}{t-r} + \tfrac{|z-y|^2}{r-s}+\tfrac{|x-y|^2}{t-s} \Big) \Bigr)  \\ 
&\lesssim (t-s)^{50} \exp\Bigl( -  \tfrac{1}{8} \Big( \tfrac{|x-z|^2}{t-r} + \tfrac{|z-y|^2}{r-s}+\tfrac{|x-y|^2}{t-s} \Big) \Bigr), 
\end{align*}
which holds for all $0\leq s\leq r \leq t$ and $x,y,z\in \R^2$. Using both, we obtain 
\begin{align*}
\big|& \eqref{object:eq-dduhlin-lin-p6} - \eqref{object:eq-dduhlin-lin-p7}\big| \\ 
&\lesssim  \iint \ds \dy \iint \dr \dz \bigg( \chi_K(t-r)  \chi_N(r-s) \chi_L(t-s) (t-r)^{-\frac{3}{2}} (r-s)^{-1} (t-s)^{-\frac{3}{2}} \\ 
&\quad \times  (t-s)^{50} \exp\Bigl( -  \tfrac{1}{8} \Big( \tfrac{|x-z|^2}{t-r} + \tfrac{|z-y|^2}{r-s}+\tfrac{|x-y|^2}{t-s} \Big) \Bigr)\bigg)  \\ 
&\lesssim \iint \ds \dr  \chi_K(t-r) \chi_N(r-s) \chi_L(t-s)(t-r)^{-\frac{1}{2}} (t-s)^{50-\frac{3}{2}}  \lesssim K \max(K,N)^{-2} L^{-99}, 
\end{align*}
which is clearly acceptable. 
\end{proof}

In the second lemma needed for the proof of Proposition~\ref{object:prop-dduhlin-dlin-main}, we study the sum of the explicit integrals from \eqref{object:eq-dduhlin-lin}. 

\begin{lemma}[Explicit integrals]\label{object:lem-explicit-integrals}
Let $0< T \leq \frac{1}{128}$ and let $\Crho$ be as in \eqref{main:eq-sigma-eps-Crho}. Then, there exists a time-dependent constant $\Czero \in L_t^\infty([0,\infty))$, which does not depend on $\rho$, such that the following holds: For all $\eps\in \sdyadic$, $k,\ell \in [2]$, $1< q <\infty$, and $0<\nu<\min(1/q,\theta)$,
\begin{equation}\label{object:eq-explicit-integrals}
\begin{aligned}
\sup_{\rho\in \mc{P}} \bigg\| \sum_{\substack{K,L\geq N_{\lo}: \\ K\sim_\theta L}} &
\iint \ds \dy \iint \dr \dz (\partial^k \fp_K)(r,z;t,x) \fp_\eps(s,y;r,z) 
(\partial^\ell \fp_{\eps,L})(s,y;t,x) \\
&- \tfrac{1}{4} \delta^{k\ell} \bigl( \sigma_{\rho,\eps}^2 + \Crho + \Czero + \tfrac{1}{8\pi} \log(t) \bigr) \bigg\|_{L_t^q L_x^\infty([0,T])} \lesssim \eps^{2\nu}. 
\end{aligned}
\end{equation}
\end{lemma}

We emphasise that the $\log$-term in \eqref{object:eq-explicit-integrals} is only needed because we prefer our renormalisation constant $\sigma_{\rho,\eps}^2$ to be time-independent, and it will not be important in our proof of Proposition \ref{object:prop-dduhlin-dlin-main}.

\begin{proof}
Throughout the proof, all space-time norms will be over the region $[0,T]\times \T^2$, and we omit this from our notation. Using the semigroup property of the heat kernel, we have 
\begin{align*}
\iint &\ds \dy \iint \dr \dz (\partial^k \fp_K)(r,z;t,x) \fp_\eps(s,y;r,z) 
(\partial^\ell \fp_{\eps,L})(s,y;t,x) \\
=& \iint \ds \dy \iint \dr \dz \chi_K(t-r) \rho_\eps(r-s) \bigl( \rho_\eps \chi_L \bigr)(t-s) (\partial^k \fp)(r,z;t,x) \fp(s,y;r,z) 
(\partial^\ell \fp)(s,y;t,x) \\ 
=&  \iint \ds \dy \int_s^t \dr \, \chi_K(t-r) \rho_\eps(r-s) \bigl( \rho_\eps \chi_L \bigr)(t-s) (\partial^k \fp)(s,y;t,x) (\partial^\ell \fp)(s,y;t,x) \\
=&  \iint \ds \dy  \bigl( \rho_\eps \ast \chi_K\bigr) (t-s) \bigl( \rho_\eps \chi_L \bigr)(t-s) (\partial^k \fp)(s,y;t,x) (\partial^\ell \fp)(s,y;t,x). 
\end{align*}
Thus, \eqref{object:eq-explicit-integrals} can be rewritten as 
\begin{equation}\label{object:eq-explicit-integrals-p0}
\begin{aligned}
 \sup_{\rho\in \mc{P}} \bigg\| \sum_{\substack{K,L\geq N_{\lo} \colon \\ K\sim_\theta L}} &\int_0^t \ds \int_{\T^2} \dy \bigl(  \rho_\eps \ast \chi_K \bigr)(t-s) \bigl(\rho_\eps \chi_L\bigr)(t-s) \partial_k \fp(s,y;t,x) \, \partial_\ell \fp(s,y;t,x) \\
&- \tfrac{1}{4} \delta_{k\ell} \bigl( \sigma_{\rho,\eps}^2 + \Crho + \Czero + \tfrac{1}{8\pi} \log(t) \bigr) \bigg\|_{L_t^q L_x^\infty} \lesssim \eps^{2\nu}.
\end{aligned}
\end{equation}
In the case $k \neq \ell$, a parity argument directly implies that the integral in \eqref{object:eq-explicit-integrals-p0} vanishes, which makes the estimate \eqref{object:eq-explicit-integrals-p0} trivial. Together with a symmetry argument, it therefore suffices to treat the case $k=\ell=1$. To this end, we recall the estimate
\begin{equation*}
\qquad \int_{\T^2} \dy\,  \bigl(\partial^1 \fp(s,y;t,x)\bigr)^2 = \frac{1}{32 \pi (t-s)^2} + \bigO((t-s)^{10}),
\end{equation*}
which follows from \eqref{eq:resonance-torus-to-plane-derivative} and the translation invariance of the heat kernel. By using this estimate together with the change of variables $s\ra t-s$ and the definition of $\sigma_{\rho,\eps}^2$ from \eqref{main:eq-sigma-eps-Crho}, we can rewrite the argument in \eqref{object:eq-explicit-integrals-p0} as 
\begin{align}
 &\tfrac{1}{32\pi} \sum_{\substack{K,L\geq N_{\lo}\colon \\ K\sim_\theta L}} \int_0^t \ds \bigl( \rho_\eps \ast \chi_K\bigr)(s) \bigl( \rho_\eps \chi_L\bigr)(s) s^{-2} 
-  \tfrac{1}{32\pi} \int_0^1  \ds \rho_\eps^2(s) s^{-1}  \label{object:eq-explicit-integrals-p0a}\\  
&-\tfrac{1}{4}  \bigl( \Crho + \Czero + \tfrac{1}{8\pi} \log(t) \bigr) + \tfrac{1}{32\pi} \sum_{\substack{K,L\geq N_{\lo}\colon \\ K\sim_\theta L}} \int_0^t \ds \bigl( \rho_\eps \ast \chi_K\bigr)(s) \bigl( \rho_\eps \chi_L\bigr)(s) \bigO(s^{10}).  \label{object:eq-explicit-integrals-p0b}
\end{align}
Due to the $\bigO(s^{10})$-term, the last summand in \eqref{object:eq-explicit-integrals-p0b} is harmless and is easily seen to converge to a constant not depending on $\rho$. After using that $\sum_{N\geq N_{\lo}} \chi_N(t)=1$ for all $t\in (0,T]$ and that $\ind_{>0} \ast \rho_\eps =I\rho_\eps$, where $I$ denotes the time integral starting from $0$, it then suffices to prove that
\begin{align}
 \sup_{\rho \in \mc{P}} \bigg\| \int_0^t \ds (I\rho_\eps)(s) \rho_\eps(s) s^{-2} - \int_0^t \ds \rho_\eps^2(s) s^{-1} 
- 8\pi (\Crho + \Cone) \bigg\|_{L_t^q} &\lesssim \eps^{\frac{2}{q}},  \label{object:eq-explicit-integrals-p1} \\ 
\sup_{\rho \in \mc{P}} \bigg\|  \sum_{\substack{K,L\geq N_{\lo}\colon \\ K\not\sim_\theta L}} \int_0^t \ds \bigl( \rho_\eps \ast \chi_K\bigr)(s) \bigl( \rho_\eps \chi_L\bigr)(s) s^{-2} - 8 \pi \Ctwo \bigg\|_{L_t^\infty} &\lesssim \eps^{2\theta}, \label{object:eq-explicit-integrals-p2} \\
\sup_{\rho \in \mc{P}} \bigg\|  \int_t^1 \ds \rho_\eps^2(s) s^{-1} + \log(t)  \bigg\|_{L_t^q} &\lesssim \bigl( 1+ \abs{\log(\eps)}\bigr) \eps^{\frac{2}{q}}, \label{object:eq-explicit-integrals-r1}
\end{align}
where $\Cone,\Ctwo \in L_t^\infty([0,\infty))$ are (possibly time-dependent) constants not depending on $\rho$. We first prove \eqref{object:eq-explicit-integrals-p1}. Using $\partial_s (s^{-1})=-s^{-2}$ and integration by parts, we can write 
\begin{align}
&\quad \int_0^t \ds (I\rho_\eps)(s) \rho_\eps(s) s^{-2} - \int_0^t \ds \rho_\eps^2(s) s^{-1} - 8\pi (\Crho+\Cone) \notag\\
&= \bigg( - (I\rho_\eps)(s) \rho_\eps(s) s^{-1}\Big|_{s=0}^{s=t}
+ \int_0^t \ds \rho_\eps^2(s) s^{-1} + \int_0^t \ds (I\rho_\eps)(s) \rho^\prime_\eps(s) s^{-1} \bigg) \\
&\qquad - \int_0^t \ds \rho_\eps^2(s) s^{-1} - 8\pi (\Crho+ \Cone) \notag \\
&= \Bigl( - (I\rho_\eps)(t) \rho_\eps(t) t^{-1} - 8 \pi \Cone \Bigr)  
 + \Big( \int_0^t \ds (I\rho_\eps)(s) \rho^\prime_\eps(s) s^{-1} - 8 \pi \Crho \Big). \label{local:eq-explicit-integrals-q1}
\end{align}
We now discuss the first and second terms in \eqref{local:eq-explicit-integrals-q1} separately. To treat the first term, we note that, since $\rho_\eps$ is bounded, 
\begin{equation}\label{local:eq-explicit-integrals-q2}
\big| (I\rho_\eps)(t) \rho_\eps(t) t^{-1}\big| \lesssim t \times 1 \times t^{-1} \lesssim 1. 
\end{equation}
Furthermore, since $\rho_\eps|_{[2\eps^2,\infty)}\equiv 1$, we also have that $|(I\rho_\eps)(t)-t|\lesssim \eps^2$. As a result, we have that
\begin{equation}\label{local:eq-explicit-integrals-q3}
\big| (I\rho_\eps)(t) \rho_\eps(t) t^{-1} - 1 \big| \lesssim \big|  \rho_\eps(t)  - 1 \big| + \eps^2 t^{-1} \lesssim\eps^2 t^{-1}. 
\end{equation}
By combining \eqref{local:eq-explicit-integrals-q2} and \eqref{local:eq-explicit-integrals-q3}, we then obtain that 
\begin{equation*}
\| (I\rho_\eps)(t) \rho_\eps(t) t^{-1} - 1 \|_{L_t^q} \lesssim \| \min(1,\eps^2 t^{-1}) \|_{L_t^q} \lesssim \eps^{\frac{2}{q}}. 
\end{equation*}
Thus, by choosing $\Cone=-(8\pi)^{-1}$, we see that the first term in \eqref{local:eq-explicit-integrals-q1} converges at the desired rate.  The second term in \eqref{local:eq-explicit-integrals-q1} can be written as 
\begin{equation*}
\int_0^t \ds (I\rho_\eps)(s) \rho^\prime_\eps(s) s^{-1} - 8 \pi \Crho 
= - \int_t^\infty \ds (I\rho_\eps)(s) \rho^\prime_\eps(s) s^{-1}, 
\end{equation*}
where we used the definition of $\Crho$ from \eqref{main:eq-sigma-eps-Crho}. Since $\rho_\eps^\prime$ is supported on times $s \sim \eps^{2}$, for which $I\rho_\eps(s)\lesssim \eps^{2}$ and $\rho_\eps^\prime(s) \lesssim \eps^{-2}$, this integral can be estimated by 
\begin{equation*}
\Big| \int_t^\infty \ds   (I\rho_\eps)(s) \rho^\prime_\eps(s) s^{-1}\Big| \lesssim \eps^{2-2} \ind\{ t \lesssim \eps^2 \}  \int_0^\infty \ds \ind\{s\sim \eps^2\} s^{-1} 
\lesssim  \ind\{ t \lesssim \eps^2 \}.
\end{equation*}
Since this can be bounded in $L_t^q$ by $\eps^{\frac{2}{q}}$, this completes the proof of \eqref{object:eq-explicit-integrals-p1}. We now turn to the proof of \eqref{object:eq-explicit-integrals-p2}. To this end, we set\footnote{The convergence of the dyadic sum in the definition of $\Ctwo$, as well as the fact that $\Ctwo\in L_t^\infty([0,\infty))$, follows easily from our estimates below.}
\begin{equation}\label{object:eq-explicit-integrals-p3}
\Ctwo(t) \eqdef \tfrac{1}{8\pi} \sum_{\substack{K,L\geq N_{\lo} \colon \\  K \not \sim_\theta L}} \int_0^t \ds (I\chi_K)(s) \chi_L(s) s^{-2}. 
\end{equation}
We first note that the dyadic summands in both \eqref{object:eq-explicit-integrals-p2} and \eqref{object:eq-explicit-integrals-p3} are only nonzero if $K\gtrsim L$. The reason is that $\rho_\eps \ast \chi_K$ and $I\chi_K$ are supported on times $\gtrsim K^{-2}$ and that $\chi_L$ is supported on times $\sim L^{-2}$. Together with our existing constraint $K\not\sim_\theta L$, this allows us to restrict to $K\gg_\theta L$, i.e., $L\leq K^{1-\theta}$. We also record the elementary estimates
\begin{equation*}
\| \rho_\eps \ast \chi_K \|_{L^\infty}, \| I \chi_K \|_{L^\infty} \lesssim K^{-2} \qquad \text{and} \qquad \| \rho_\eps \ast \chi_K - I \chi_K \|_{L^\infty} \lesssim \min(\eps^2,K^{-2}), 
\end{equation*}
which will be needed momentarily. We now split 
\begin{align}
 &\sum_{\substack{K,L\geq N_{\lo}\colon \\ K\gg_\theta L}} \int_0^t \ds \bigl( \rho_\eps \ast \chi_K\bigr)(s) \bigl( \rho_\eps \chi_L\bigr)(s) s^{-2} 
 -  \sum_{\substack{K,L\geq N_{\lo} \colon \\  K \gg_\theta L}} \int_0^t \ds (I\chi_K)(s) \chi_L(s) s^{-2} \notag \\ 
&\quad = \sum_{\substack{K,L\geq N_{\lo}\colon \\ K\gg_\theta L}}  \int_0^t \ds \bigl( \rho_\eps \ast \chi_K - I\chi_K \bigr)(s)  (\rho_\eps \chi_L)(s) s^{-2} 
\label{object:eq-explicit-integrals-p4} \\ 
&\quad + \sum_{\substack{K,L\geq N_{\lo}\colon \\ K\gg_\theta L}}  \int_0^t \ds  ( I \chi_K)(s)  (\rho_\eps-1)(s) \chi_L(s) s^{-2}. 
\label{object:eq-explicit-integrals-p5}
\end{align}
The first term \eqref{object:eq-explicit-integrals-p4} can be bounded by 
\begin{align*}
\|& \eqref{object:eq-explicit-integrals-p4} \|_{L_t^\infty}
\lesssim \sum_{\substack{K,L\geq N_{\lo}\colon \\ K\gg_\theta L}}   \| \rho_\eps \ast \chi_K - I \chi_K \|_{L^\infty_s} \| \chi_L s^{-2} \|_{L^1_s} \\ 
&\lesssim \sum_{\substack{K,L\geq N_{\lo}\colon \\ K\gg_\theta L}}   \min(\eps^2,K^{-2}) L^2 
\lesssim \sum_{K} \min(\eps^2,K^{-2}) K^{2(1-\theta)}
\lesssim \eps^{2\theta}. 
\end{align*}
The second term \eqref{object:eq-explicit-integrals-p5} can be bounded by 
\begin{align*}
\|& \eqref{object:eq-explicit-integrals-p5} \|_{L_t^\infty}
\lesssim \sum_{\substack{K,L\geq N_{\lo}\colon \\ K\gg_\theta L}}   \| I \chi_K \|_{L^\infty_s} \| (\rho_\eps-1) \chi_L s^{-2} \|_{L^1_s} \\ 
&\lesssim \sum_{\substack{K,L\geq N_{\lo}\colon \\ K\gg_\theta L}} \ind\{L \gtrsim \eps^{-1}\}  K^{-2} L^2 
\lesssim \sum_{K\geq N_{\lo}} \ind\{K \gtrsim \eps^{-1}\}  K^{-2\theta}
\lesssim \eps^{2\theta}.
\end{align*}
This completes the proof of \eqref{object:eq-explicit-integrals-p2}, and it therefore only remains to prove \eqref{object:eq-explicit-integrals-r1}. However, since the argument in \eqref{object:eq-explicit-integrals-r1} can easily be bounded pointwise by $\ind\{t\lesssim \eps^2\} \abs{\log(t)}$, this estimate is trivial. 
\end{proof}

Equipped with Lemma~\ref{object:lem-dduhlin-dlin} and Lemma~\ref{object:lem-explicit-integrals}, we can now prove the first main result of this subsection. 

\begin{proof}[Proof of Proposition~\ref{object:prop-dduhlin-dlin-main}]
The estimate in \eqref{object:eq-dduhlin-dlin-main} follows directly from Lemma~\ref{object:lem-dduhlin-dlin}, Lemma~\ref{object:lem-explicit-integrals}, and 
\begin{align*}
T^{-\frac{1-\roughdataweightmultiple[4]}{2}} \sum_{K,L,N\geq N_{\lo}} \max(K,L,N)^{-1+\roughdataweightmultiple[4]}  \|  1 \|_{L_t^q([0,T])}
&\lesssim T^{\frac{1}{q}} T^{-\frac{1-\roughdataweightmultiple[4]}{2}} \sum_{K \geq N_{\lo}} \log(K)^2 K^{-1+\roughdataweightmultiple[4]}   \\
&\lesssim \abs{\log(T)}^2 T^{\frac{1}{q}}. 
\end{align*}
The convergence rate and the fact that the limit does not depend on $\rho$ follow from the fact that the argument in \eqref{object:eq-dduhlin-lin} only depends on $\eps$ or $\rho$ if $L\gtrsim \eps^{-1}$ or $N\gtrsim \eps^{-1}$, as well as the corresponding statements in Lemma \ref{object:lem-explicit-integrals}.
\end{proof}

We now turn our attention towards the proof of Proposition \ref{objects:prop-quadratic-dlin}. For later use, it is convenient to first record the following crude estimate, which follows from arguments similar to those in the proof of Lemma \ref{object:lem-dduhlin-dlin}.

\begin{lemma}[Crude estimate]\label{object:lem-dduhlin-dlin-crude}
Let Assumption~\ref{local:assumption-B} be satisfied and let $k,\ell\in [2]$. Furthermore, let $\eps \in \sdyadic$, $\rho\in \mc{P}$, and $K,L,N\geq N_{\lo}$. Then, it holds that 
\begin{equ}\label{object:eq-dduhlin-dlin-crude}
\Big\| \E \Big[ \covd_B^k \Duh_{\Y,B,K}\bigl( \graylinear[\oeps][][][N]\bigr)(t,x)
\otimes \covd_B^\ell \graylinear[\oeps][][][L](t,x)\Big] 
\Big\|_{L_{t,x}^\infty} \lesssim  \lossfactor_T(F)^3 L \max(K,L,N)^{-1}. 
\end{equ}
\end{lemma}

We emphasise that the right-hand side of \eqref{object:eq-dduhlin-dlin-crude} is not summable, and hence \eqref{object:eq-dduhlin-dlin-crude} certainly does not imply Proposition \ref{object:prop-dduhlin-dlin-main}. However, \eqref{object:eq-dduhlin-dlin-crude} will still be useful, since it does not involve a difference between two terms. 

\begin{proof}
Using a similar identity as in \eqref{object:eq-dduhlin-lin-p2}, as well as Propositions \ref{prop:pB-ptwise-bounds} and \ref{prop:qB-ptwise-bds}, we have that 
\begin{align}
\Bigl|& \E \Big[ \covd_B^k \Duh_{\Y,B,K}\bigl( \graylinear[\oeps][][][N]\bigr)(t,x)
\otimes \covd_B^\ell \graylinear[\oeps][][][L](t,x)\Big]  \Bigr|  \notag \\ 
&\lesssim \iint \ds \dy \iint \dr \dz \abs{\covd^k \fp_{\Y,B,K}(r,z;t,x)} \abs{\fq_{\Y,B,\eps,N}(s,y;r,z)} \abs{\covd^\ell \fq_{\Y,B,\eps,L}(s,y;t,x)}  \notag  \\ 
&\lesssim \lossfactor_T(F)^3 \iint \ds \dy \iint \dr \dz \Bigl( \chi_K(t-r) \chi_N(r-s) \chi_L(t-s) 
(t-r)^{-\frac{3}{2}} (r-s)^{-1} (t-s)^{-\frac{3}{2}}   \notag  \\ 
&\qquad \times  \exp\Bigl( - c \Bigl( \tfrac{\abs{x-z}^2}{t-r}+\tfrac{\abs{z-y}^2}{r-s} +\tfrac{\abs{x-y}^2}{t-s} \Bigr) \Bigr) \Bigr). \label{object:eq-dduhlin-dlin-crude-p1}
\end{align}
Due to the $\chi$-factors, this expression is only nonzero if $K,N\gtrsim L$. By first integrating in $(y,z)$ and then integrating in $(s,r)$, we obtain 
\begin{align*}
\eqref{object:eq-dduhlin-dlin-crude-p1}
&\lesssim \lossfactor_T(F)^3 \iint \ds \dr \chi_K(t-r) \chi_N(r-s) \chi_L(t-s) (t-r)^{-\frac{1}{2}} (t-s)^{-\frac{3}{2}} \\ 
&\lesssim \lossfactor_T(F)^3 KL^3 \times K^{-2} \max(L,N)^{-2} 
\lesssim \lossfactor_T(F)^3 L \max(K,L,N)^{-1}. \qedhere
\end{align*}
\end{proof}

In preparation for the proof of Proposition \ref{objects:prop-quadratic-dlin}, we also record the following variant of Lemma \ref{object:lem-linear}, i.e., our earlier estimate of the covariant stochastic object $\graylinear[\oeps][][][]$.

\begin{lemma}\label{objects:lem-deterministic-and-linear}
Let Assumption~\ref{local:assumption-B} be satisfied and let $H\colon [0,T]\times \T^2 \rightarrow \End(\Y)$ be deterministic. Furthermore, let $\eps \in \sdyadic$, $N \geq N_{\lo}$, $p\geq 2$, and $\frac{1}{\theta-\kappa} < q \leq \infty$. Then, it holds that 
\begin{equation}\label{objects:eq-deterministic-and-linear}
\begin{aligned}
\bigl\| \bigl\| H \, \graylinear[\oeps][][][N] \bigr\|_{\Cs_t^{-\theta}L_x^r} \bigr\|_{L^p(\Omega)} 
\lesssim \sqrt{p} \lossfactor_T(F) N^{-2(\theta-\kappa-\frac{1}{q})} \| H \|_{L_t^q L_x^\infty([0,T])}. 
\end{aligned}
\end{equation}
\end{lemma}

\begin{remark}
The point of Lemma \ref{objects:lem-deterministic-and-linear} is that no regularity in time is needed for $H$. If we require that $H$ has regularity greater than $\theta$ in time, then the product $H \, \graylinear[\oeps][][][]$ can be estimated more easily using Lemma \ref{prelim:lem-product-temporal}.
\end{remark}

\begin{proof}
Due to Lemma \ref{lemma:time-regularity-kolmogorov}	and Gaussian hypercontractivity, it suffices to prove for all $0\leq t_0 \leq t_1 \leq T$ and $x\in \T^2$ that 
\begin{equation}\label{objects:eq-deterministic-and-linear-p1}
\Bigl\| \int_{t_0}^{t_1} \dt H(t,x)\, \graylinear[\oeps][][][N](t,x) \Bigr\|_{L^2(\Omega)} 
\lesssim \lossfactor_T(F) (t_1-t_0)^{1-\theta+\kappa} N^{-2(\theta-\kappa-\frac{1}{q})} \| H(t,x)\|_{L_t^q}. 
\end{equation}
To this end, we first use It\^{o}'s isometry and Hölder's inequality, which easily yield that  
\begin{align*}
\Bigl\|& \int_{t_0}^{t_1} \dt H(t,x)\, \graylinear[\oeps][][][N](t,x) \Bigr\|_{L^2(\Omega)}^2\\ 
&\lesssim \Bigl\|  \int_{t_0}^{t_1} \dt H(t,x)\, \fq_{\Y,B,\eps,N}(s,y;t,x) \Bigr\|_{L_{s,y}^1} \Bigl\|  \int_{t_0}^{t_1} \dt H(t,x)\, \fq_{\Y,B,\eps,N}(s,y;t,x) \Bigr\|_{L_{s,y}^\infty}.  
\end{align*}
Using Propositions \ref{prop:pB-ptwise-bounds} and \ref{prop:qB-ptwise-bds}, the $L_{s,y}^1$-factor can be bounded by 
\begin{align*}
\Bigl\|&  \int_{t_0}^{t_1} \dt H(t,x)\, \fq_{\Y,B,\eps,N}(s,y;t,x) \Bigr\|_{L_{s,y}^1} \\
&\lesssim \lossfactor_T(F) \int_{t_0}^{t_1} \dt \abs{H(t,x)} \iint \ds \dy \chi_N(t-s) (t-s)^{-1} \exp\Bigl( - c \tfrac{|x-y|^2}{t-s} \Bigr) \\ 
&\lesssim \lossfactor_T(F) N^{-2} (t_1-t_0)^{1-\frac{1}{q}} \| H(t,x) \|_{L_t^q}. 
\end{align*}
Similarly, the $L_{s,y}^\infty$-factor can be bounded by 
\begin{align*}
\Bigl\|&  \int_{t_0}^{t_1} \dt H(t,x)\, \fq_{\Y,B,\eps,N}(s,y;t,x) \Bigr\|_{L_{s,y}^\infty}\\ 
&\lesssim \lossfactor_T(F) \sup_{s\in [0,T]} \int_{t_0}^{t_1} \dt \abs{H(t,x)} \chi_N(t-s) (t-s)^{-1} \\  
&\lesssim\lossfactor_T(F) N^2 \| H(t,x) \|_{L_t^q} \sup_{s \in [0,T]} \| \chi_N(t-s) \|_{L_t^{q^\prime}([t_0,t_1])} \\  
&\lesssim \lossfactor_T(F) N^2 \min(t_1-t_0,N^{-2})^{1-\frac{1}{q}} \| H(t,x) \|_{L_t^q}. 
\end{align*}
Since $\frac{1}{q}<\theta-\kappa$, the desired estimate now easily follows by bounding the $\min$-factor using a geometric mean.
\end{proof}

Equipped with Lemmas \ref{object:lem-dduhlin-dlin-crude} and \ref{objects:lem-deterministic-and-linear}, whose usefulness will become more clear below, we now examine \eqref{objects:eq-quadratic-dlin-main} more directly. To this end, recall that we are looking for an estimate of 
\begin{equation*}
\Duh^F_{\Y,B,\star} \bigl(  \Qc\bigl(\graylinear[\oeps][][][] , \covd^k \graylinear[\oeps][][][] \bigr)\bigr) \parasim \covd^\ell \graylinear[\oeps] - \tfrac{1}{4} \delta^{k\ell} \bigl( \sigma_{\rho,\eps}^2 + \Crho \bigr) \mathcal{K}(\graylinear[\oeps][][][]). 
\end{equation*}
We note that the first summand is an interaction of three linear objects. Using a Wiener chaos decomposition, it can therefore be split into four terms, one of which is non-resonant and three of which have one resonance. To write out the four terms concisely, we let $\redlinear[\oeps][][][\ostar]$ be an independent copy of $\graylinear[\oeps][][][\ostar]$. To be more precise, we first make an independent copy of the $\Y$-valued space-time white noise $\boldxi$ and then, as suggested by \eqref{local:eq-linear-decomposition}, \eqref{local:eq-philinear-N}, and \eqref{local:eq-graylinear-boldxi}, we define $\redlinear[\oeps][][][N]$ as the integral of $\fq_{\Y,B,\eps,N}(s,y;t,x)$ against this independent copy of $\boldxi$. Using the independent copy $\redlinear[\oeps][][][\ostar]$ of $\graylinear[\oeps][][][\ostar]$, as well as the aforementioned Wiener chaos decomposition, we then obtain the identity
\begin{align}
\Duh^F_{\Y,B,\star}& \bigl(  \Qc\bigl(\graylinear[\oeps][][][] , \covd^k \graylinear[\oeps][][][] \bigr)\bigr) \parasim \covd^\ell \graylinear[\oeps] - \tfrac{1}{4} \delta^{k\ell} \bigl( \sigma_{\rho,\eps}^2 + \Crho \bigr) \mathcal{K}(\graylinear[\oeps][][][]) \label{object:eq-quadratic-dlin-decomposition-original} \\ 
&= \Wick{\Duh^F_{\Y,B,\star} \bigl(  \Qc\bigl(\graylinear[\oeps][][][] , \covd^k \graylinear[\oeps][][][] \bigr)\bigr) \parasim \covd^\ell \graylinear[\oeps]} \label{object:eq-quadratic-dlin-decomposition-nonres} \\ 
&+ \redE \Bigl[ \Duh^F_{\Y,B,\star} \bigl(  \Qc\bigl(\redlinear[\oeps][][][] , \covd^k \redlinear[\oeps][][][] \bigr)\bigr) \parasim \covd^\ell \graylinear[\oeps] \Bigr] \label{object:eq-quadratic-dlin-decomposition-12} \\ 
&+ \redE \Bigl[ \Duh^F_{\Y,B,\star} \bigl(  \Qc\bigl(\redlinear[\oeps][][][] , \covd^k \graylinear[\oeps][][][] \bigr)\bigr) \parasim \covd^\ell \redlinear[\oeps]\Bigr] \label{object:eq-quadratic-dlin-decomposition-13} \\ 
&+ \redE \Bigl[ \Duh^F_{\Y,B,\star} \bigl(  \Qc\bigl(\graylinear[\oeps][][][] , \covd^k \redlinear[\oeps][][][] \bigr)\bigr) \parasim \covd^\ell \redlinear[\oeps] \Bigr]
 - \tfrac{1}{4} \delta^{k\ell} \bigl( \sigma_{\rho,\eps}^2 + \Crho \bigr) \mathcal{K}(\graylinear[\oeps][][][]). \label{object:eq-quadratic-dlin-decomposition-23} 
\end{align} 
Here, $\redE$ denotes the expectation over the space-time white noise used in the definition of $\redlinear[\oeps][][][\ostar]$. We emphasise that each term in \eqref{object:eq-quadratic-dlin-decomposition-12}--\eqref{object:eq-quadratic-dlin-decomposition-23} can be written out as an integral of a deterministic kernel against the space-time white noise $\boldxi$. For example, 
\begin{align*}
\redE\Bigl[ \Duh^F_{\Y,B,\star} \bigl(  \Qc\bigl(\redlinear[\oeps][][][] , \covd^k \graylinear[\oeps][][][] \bigr)\bigr) \parasim \covd^\ell \redlinear[\oeps]\Bigr]  
= \iint \ds_2 \dy_2 \kernel(s_2,y_2;t,x) \boldxi(s_2,y_2), 
\end{align*}
where the kernel is defined as 
\begin{align*}
\kernel(&s_2,y_2;t,x) ( \place) 
= \sum_{\substack{K,L\colon \\ K \sim_\theta L}} \iint \ds_1 \dy_1 \iint \dr \dz \Bigl(\\ 
&\fq_{\Y,B,K}(r,z;t,x) \Qc \bigl( \fq_{\Y,B,\eps}(s_1,y_1;r,z) \ONBY_m, \covd^k \fq_{\Y,B,\eps}(s_2,y_2;r,z) (\place) \bigr) \otimes \covd^\ell \fq_{\Y,B,\eps,L}(s_1,y_1;t,x) \ONBY^m \Bigr)
\end{align*}
and where $(\ONBY_m)_m$ is an orthonormal basis of $\Y$. Compared to this integral, however, the notation in \eqref{object:eq-quadratic-dlin-decomposition-12}--\eqref{object:eq-quadratic-dlin-decomposition-23} is more economical.  

\begin{notation}
Since $\redE$ appears in every expression in \eqref{object:eq-quadratic-dlin-decomposition-12}--\eqref{object:eq-quadratic-dlin-decomposition-23}, it may seem like it would have been better to omit $\redE$ from our notation and implicitly take expectations whenever $\redlinear[\oeps][][][\ostar]$ appears. In the proof of Proposition \ref{objects:prop-quadratic-dlin}, however, we will sometimes need to write expressions such as $\redE [\|\redlinear[\oeps][][][\ostar] \|_{\modsol^{-\kappa}}^2 ]$. If we omitted $\redE$ from our notation, this could easily be confused with $\| \redE \redlinear[\oeps][][][\ostar] \|_{\modsol^{-\kappa}}^2$, which is clearly not the same. 
\end{notation}

In the following four lemmas, we separately bound the terms in \eqref{object:eq-quadratic-dlin-decomposition-nonres}, \eqref{object:eq-quadratic-dlin-decomposition-12}, \eqref{object:eq-quadratic-dlin-decomposition-13}, and \eqref{object:eq-quadratic-dlin-decomposition-23}. We start with the non-resonant term \eqref{object:eq-quadratic-dlin-decomposition-nonres}, which is the easiest term. 

\begin{lemma}[Non-resonant term]\label{objects:lem-quadratic-dlin-nonres}
Let Assumption~\ref{local:assumption-B} be satisfied. Furthermore, let $k,\ell \in [2]$ and let $\Qc \colon \Y \times \Y \rightarrow \Y$ be a $\liegroup$-equivariant bilinear map. Finally, let $\eps\in \sdyadic$ and let $K,L,M,N\geq N_{\lo}$ satisfy $K\sim_\theta L$. Then, it holds for all $p\geq 2$ that 
\begin{equation}\label{objects:eq-quadratic-dlin-nonres}
\begin{aligned}
\Bigl\|& \Bigl\| \Wick{  \Duh^F_{\Y,B,K} \bigl( \Qc\bigl(\graylinear[\oeps][][][M] , \covd^k \graylinear[\oeps][][][N] \bigr) \bigr) \otimes  \covd^\ell \graylinear[\oeps][][][L] } \Bigr\|_{\Cs_t^{-\theta}L_x^r} \Bigr\|_{L^p(\Omega)}  \\ 
&\lesssim \lossfactor_T(F)^4 p^{\frac{3}{2}} \max(K,L,M,N)^{-\theta+2\kappa}. 
\end{aligned}
\end{equation}
\end{lemma}

\begin{proof}
We first prove for all $(t,x)\in [0,T]\times \T^2$, $0\leq t_0 \leq t_1 \leq T$, and $K,L,M,N \geq N_{\lo}$ that 
\begin{equation}\label{objects:lem-quadratic-dlin-nonres-p1}
\Big\|\Wick{\Duh_{\Y, B, K}^{F}\Bigl(\Qc\bigl(\graylinear[\oeps][][M], \covd^k \graylinear[\oeps][][N] \bigr)\Bigr) \otimes \covd^\ell \graylinear[\oeps][][L]} (t, x)\Big\|_{L^2(\Omega)} 
\lesssim \lossfactor_T(F)^4 \max(K,M,N)^{-1} L
\end{equation}
as well as
\begin{equation}\label{objects:lem-quadratic-dlin-nonres-p2}
\begin{aligned}
\bigg\|& \int_{t_0}^{t_1} \dt\, \Wick{\Duh_{\Y, B, K}^{F}\Bigl(\Qc\bigl(\graylinear[\oeps][][M], \covd^k \graylinear[\oeps][][N] \bigr)\Bigr) \otimes \covd^\ell \graylinear[\oeps][][L]} (t, x)\bigg\|_{L^2(\Omega)} \\
&\lesssim \lossfactor_T(F)^4  |t_1-t_0|^{1-\theta+\kappa} \max(K,M,N)^{-1} L^{1-2(\theta-\kappa)}. 
\end{aligned}
\end{equation} 
To this end, we first write 
\begin{align*}
\Wick{&\Duh_{\Y, B, K}^{F}\Bigl(\Qc\bigl(\graylinear[\oeps][][M], \covd^k \graylinear[\oeps][][N] \bigr)\Bigr) \otimes \covd^\ell \graylinear[\oeps][][L]} (t, x) \\
&= \Wick{\int \ds_1 \dy_1 \ds_2 \dy_2 \ds_3 \dy_3\, \bigl(\kernel_{12}(s_1, y_1, s_2, y_2; t, x) \otimes \covd^\ell \fq_{\Y, B, \varep, L}(s_3, y_3; t, x)\bigr) \\
&\quad \quad \quad \quad \quad \quad(\boldxi(s_1, y_1) \otimes \boldxi(s_2, y_2) \otimes \boldxi(s_3, y_3))},
\end{align*}
where $\kernel_{12}$ is the kernel such that
\begin{equs}
\Wick{&\Duh_{\Y, B, K}^{F}\Bigl(\Qc\bigl(\graylinear[\oeps][][M], \covd^k \graylinear[\oeps][][N] \bigr)\Bigr)}(t, x) \\
&= \Wick{\int\ds_1 \dy_1 \ds_2 \dy_2 \kernel_{12}(s_1, y_1, s_2, y_2; t, x) (\boldxi(s_1, y_1) \otimes \boldxi(s_2, y_2))}.
\end{equs}
We note that the kernels $\kernel_{12}$ and $\covd^\ell \fq_{\Y, B, \varep, L}$ have already been estimated in the proofs of Lemma \ref{lemma:quadratic-object-non-resonant-bound} 
and Lemma \ref{object:lem-linear}, respectively. There, it was shown that 
\begin{align*}
\bigl\| \kernel_{12}(s_1, y_1, s_2, y_2; t, x)\bigr\|_{L^2_{s_1,s_2,y_1,y_2}} &\lesssim \lossfactor_T(F)^3 \max(K,M,N)^{-1}, \\ 
\bigl\| \covd^\ell \fq_{\Y, B, \varep, L}(s_3, y_3; t, x)\bigr\|_{L^2_{s_3,y_3}} &\lesssim \lossfactor_T(F) L. 
\end{align*}
Together with It\^{o}'s isometry, this readily implies \eqref{objects:lem-quadratic-dlin-nonres-p1}. In order to obtain \eqref{objects:lem-quadratic-dlin-nonres-p2}, we can use the same kernel estimates, the fact that $\kernel_{12}\otimes \covd^\ell \fq_{\Y,B,\eps,L}$ is supported on times satisfying $t-s_3 \sim L^{-2}$, and Lemma \ref{object:lem-time-integration-gain}.
Equipped with \eqref{objects:lem-quadratic-dlin-nonres-p2}, we can then use Lemma \ref{lemma:time-regularity-kolmogorov} and Gaussian hypercontractivity, which imply that  
\begin{align*}
\Bigl\|& \Bigl\| \Wick{  \Duh^F_{\Y,B,K} \bigl( \Qc\bigl(\graylinear[\oeps][][][M] , \covd^k \graylinear[\oeps][][][N] \bigr) \bigr) \otimes  \covd^\ell \graylinear[\oeps][][][L] } \Bigr\|_{\Cs_t^{-\theta}L_x^r} \Bigr\|_{L^p(\Omega)}  \\ 
&\lesssim \lossfactor_T(F)^4 p^{\frac{3}{2}}  \max(K,M,N)^{-1} L^{1-2(\theta-\kappa)}. 
\end{align*}
In order to complete the proof of \eqref{objects:eq-quadratic-dlin-nonres}, it then only remains to note that, since $K\sim_\theta L$, 
\begin{equation*}
\max(K,M,N)^{-1} L^{1-2(\theta-\kappa)} 
\lesssim \max(K,L,M,N)^{-(1-\theta)}  L^{1-2(\theta-\kappa)}  \lesssim \max(K,L,M,N)^{-\theta+2\kappa}. \qedhere 
\end{equation*}
\end{proof}

In the following lemma, we control the term involving a resonance between the first and second linear stochastic objects, i.e., \eqref{object:eq-quadratic-dlin-decomposition-12}. 

\begin{lemma}[Resonance between the first and second objects]\label{objects:lem-quadratic-dlin-12-resonance}
Let Assumption~\ref{local:assumption-B} be satisfied. Furthermore, let $k,\ell \in [2]$ and let $\Qc \colon \Y \times \Y \rightarrow \Y$ be a $\liegroup$-equivariant bilinear map satisfying $\delta^{mn} \Qc(\ONBY_m, \ONBY_n) = 0$, where $(\ONBY_m)_m$ is an orthonormal basis of $\Y$. Then, it holds for all $p\geq 2$ that
\begin{equation}\label{objects:eq-quadratic-dlin-12-resonance}
\begin{aligned}
\Bigl\|& \sup_{\eps \in \sdyadic} \Bigl\| \redE \Bigl[  \Duh^F_{\Y,B,\star} \bigl( \Qc\bigl(\redlinear[\oeps][][][] , \covd^k \redlinear[\oeps][][][] \bigr) \bigr) \parasim \covd^\ell \graylinear[\oeps][][][\ostar] \Bigr] \Bigr\|_{\Cs_t^{-\theta}L_x^r} \Bigr\|_{L^p(\Omega)}  \\ 
&\lesssim \lossfactor_T(F)^4 p^{\frac{1}{2}}. 
\end{aligned}
\end{equation}
Furthermore, the argument in \eqref{objects:eq-quadratic-dlin-12-resonance} has an a.s.\ limit in $\Cs_t^{-\theta} L_x^r([0,T])$ as $\eps \rightarrow 0$.  
\end{lemma}

\begin{proof} The main ingredient in the proof is Lemma \ref{lemma:covd-philinear-times-philinear-resonance}, i.e., our existing resonant estimate for $\redE [\Qc(\redlinear[\oeps][][][] , \covd^k \redlinear[\oeps][][][] )]$. Together with our weighted covariant Schauder estimate \eqref{local:eq-weighted-covariant-Schauder-1} and the fact that $\redE$ commutes with $\Duh^F_{\Y,B,K}$, it implies that
\begin{align*}
\Bigl\|& \redE \Bigl[  \Duh^F_{\Y,B,\star} \bigl( \Qc\bigl(\redlinear[\oeps][][][] , \covd^k \redlinear[\oeps][][][] \bigr) \bigr) \Bigr] \Bigr\|_{\sol^{2-4\kappa}} \\  
&\lesssim \lossfactor_T(F) \Bigl\| t^{\frac{\roughdataweightmultiple[2]}{2}} \redE \Bigl[  \Qc\bigl(\redlinear[\oeps][][][] , \covd^k \redlinear[\oeps][][][] \bigr) \Bigr] \Bigr\|_{L_{t,x}^\infty}
\lesssim \lossfactor_T(F)^4 T^{-\frac{1-2\kappa}{2}}. 
\end{align*}
Using the definition of the $\Cs_t^{-\theta}L_x^r$-norm, it then follows that 
\begin{align*}
\Bigl\|& \redE \Bigl[  \Duh^F_{\Y,B,\star} \bigl( \Qc\bigl(\redlinear[\oeps][][][] , \covd^k \redlinear[\oeps][][][] \bigr) \bigr) \Bigr] \parasim \covd^\ell \graylinear[\oeps][][][\ostar]  \Bigr\|_{\Cs_t^{-\theta}L_x^r} \\ 
&\lesssim T^{\theta} \sum_{\substack{K,L\geq N_{\lo}\colon \\ K\sim_\theta L}} 
\Bigl\| \redE \Bigl[  \Duh^F_{\Y,B,K} \bigl( \Qc\bigl(\redlinear[\oeps][][][] , \covd^k \redlinear[\oeps][][][] \bigr) \bigr) \Bigr] \otimes \covd^\ell \graylinear[\oeps][][][L]  \Bigr\|_{L_{t,x}^\infty} \\  
&\lesssim T^{-\frac{1-2\kappa}{2}+\theta} \lossfactor_T(F)^4 
\sum_{\substack{K,L \geq N_{\lo}\colon \\ K \sim_\theta L}} 
K^{-(2-4\kappa)} L^{1+\kappa} \bigl\| \graylinear[\eps][][][\ostar]\bigr\|_{\sol^{-\kappa}}.  
\end{align*}
By using the constraint $K\sim_\theta L$, the dyadic sum can be bounded by 
\begin{align*}
T^{-\frac{1-2\kappa}{2}+\theta}  \sum_{\substack{K,L\geq N_{\lo}\colon \\ K \sim_\theta L}} 
K^{-(2-4\kappa)} L^{1+\kappa} 
\lesssim T^{-\frac{1-2\kappa}{2}+\theta} 
\sum_{K\geq N_{\lo}} K^{-(2-4\kappa)+\frac{1+\kappa}{1-\theta}} 
\lesssim T^{-\frac{1-2\kappa}{2}+\theta+\frac{(2-4\kappa)}{2}-\frac{1+\kappa}{2(1-\theta)}}.  
\end{align*}
Since $\kappa$ is much smaller than $\theta$, the exponent of $T$ above is positive. Together with our existing estimate of $\graylinear[\oeps][][][\ostar]$ from \eqref{eq:linear-object-solution-space-estimate}, this then implies \eqref{objects:eq-quadratic-dlin-12-resonance}. Since Lemma \ref{lemma:covd-philinear-times-philinear-resonance} also contains a convergence statement, the a.s.\ convergence as $\eps \rightarrow 0$ can be obtained similarly, and we omit the details. 
\end{proof}

In the following lemma, we control the term involving a resonance between the first and third linear stochastic objects, i.e., \eqref{object:eq-quadratic-dlin-decomposition-13}.

\begin{lemma}[Resonance between the first and third objects]\label{objects:lem-quadratic-dlin-13-resonance}
Let Assumption~\ref{local:assumption-B} be satisfied. Furthermore, let $k,\ell \in [2]$ and let $\Qc \colon \Y \times \Y \rightarrow \Y$ be a $\liegroup$-equivariant bilinear map. Finally, let $\eps\in \sdyadic$ and let $K,L,M,N \geq N_{\lo}$. Then, it holds for 
all $p\geq 2$ that 
\begin{align*}
\Bigl\|& \Bigl\| \redE \Bigl[  \Duh^F_{\Y,B,K} \bigl( \Qc\bigl(\redlinear[\oeps][][][M] , \covd^k \graylinear[\oeps][][][N] \bigr) \bigr) \otimes \covd^\ell \redlinear[\oeps][][][L] \Bigr] \Bigr\|_{\Cs_t^{-\theta}L_x^r} \Bigr\|_{L^p(\Omega)}  \\ 
&\lesssim \lossfactor_T(F)^4 p^{\frac{1}{2}} \max(K,L,M,N)^{-2(\theta-\kappa)}. 
\end{align*}
\end{lemma}

\begin{proof}
Due to Lemma \ref{lemma:time-regularity-kolmogorov} and Gaussian hypercontractivity, it suffices to prove for all $(t,x) \in [0,T] \times \T^2$ and $ 0 \leq t_0 \leq t_1 \leq T$ that 
\begin{equation}\label{objects:eq-quadratic-dlin-13-resonance-1}
\Bigl \| \redE \Bigl[ \bigl( \Duh_{\Y,B,K}^{F} \bigl( \Qc \bigl( \redlinear[\oeps][][][M], \covd^k \graylinear[\oeps][][][N] \bigr) \bigr) \otimes   \covd^\ell \redlinear[\oeps][][][L]  \bigr) \Bigr] (t,x) \Bigr \|_{L^2(\Omega)} 
\lesssim \lossfactor_T(F)^4 
\end{equation}
as well as 
\begin{equation}\label{objects:eq-quadratic-dlin-13-resonance-2}
\begin{aligned}
\Bigl \|& \int_{t_0}^{t_1} \dt \redE \Bigl[ \bigl( \Duh_{\Y,B,K}^{F} \bigl( \Qc \bigl( \redlinear[\oeps][][][M], \covd^k \graylinear[\oeps][][][N] \bigr) \bigr) \otimes   \covd^\ell \redlinear[\oeps][][][L] \bigr) \Bigr]  (t,x)  \Bigr \|_{L^2(\Omega)} \\  
&\lesssim \lossfactor_T(F)^4  (t_1 - t_0)^{1-\theta+\kappa} \max(K,L,M,N)^{-2(\theta-\kappa)}.
\end{aligned} 
\end{equation}
By definition, it holds that 
\begin{align*}
\redE \Bigl[ \bigl( \Duh_{\Y,B,K}^{F} \bigl( \Qc \bigl( \redlinear[\oeps][][][M], \covd^k \graylinear[\oeps][][][N] \bigr) \bigr) \otimes   \covd^\ell \redlinear[\oeps][][][L] \bigr) \Bigr]   (t,x)  
= \iint \ds_2 \dy_2 \kernel(s_2,y_2;t,x) \boldxi(s_2,y_2), 
\end{align*}
where $\kernel(s_2,y_2;t,x) \colon \Y \rightarrow \Y \otimes \Y$ is defined as
\begin{align*}
\kernel(&s_2,y_2;t,x) (\place) 
\eqdef \iint \ds_1 \dy_1 \iint \dr \dz \Bigl( \fq_{\Y,B,K}(r,z;t,x)  \\ 
&\Qc\bigl(\fq_{\Y,B,\eps,M}(s_1,y_1;r,z)\ONBY_m, \covd^k \fq_{\Y,B,\eps,N}(s_2,y_2;r,z) (\place) \bigr) \otimes \covd^\ell \fq_{\Y,B,\eps,L}(s_1,y_1;t,x) \ONBY^{m} \Bigr). 
\end{align*}
Since $t-s_1 = t-r+r-s_1\geq t-r$, the kernel $\kernel(s_2,y_2;t,x)$ can only be nonzero provided that $K\gtrsim L$, which we now assume. Using It\^{o}'s isometry, we then obtain that 
\begin{align*}
\Bigl \| \redE \Bigl[ \bigl( \Duh_{\Y,B,K}^{F} \bigl( \Qc \bigl( \redlinear[\oeps][][][M], \covd^k \graylinear[\oeps][][][N] \bigr) \bigr) \otimes   \covd^\ell \redlinear[\oeps][][][L]  \bigr) \Bigr] (t,x) \Bigr \|_{L^2(\Omega)}^2 
\lesssim \iint \ds_2 \dy_2 \abs{\kernel(s_2,y_2;t,x)}^2.  
\end{align*}
In order to prove the first estimate \eqref{objects:eq-quadratic-dlin-13-resonance-1}, it therefore suffices to prove that 
\begin{align}
\| \kernel(s_2,y_2;t,x) \|_{L_{s_2,y_2}^\infty} &\lesssim \lossfactor_T(F)^4 K^2 L^3 N^3 \max(M,L)^{-2} \max(K,N)^{-4},\label{objects:eq-13-resonance-p1} \\ 
\| \kernel(s_2,y_2;t,x) \|_{L_{s_2,y_2}^1} &\lesssim \lossfactor_T(F)^4 K^{-2} L^3 N^{-1} \max(M,L)^{-2}.\label{objects:eq-13-resonance-p2}
\end{align} 
Indeed, once \eqref{objects:eq-13-resonance-p1} and \eqref{objects:eq-13-resonance-p2} have been shown, we easily obtain that 
\begin{align*}
\|& \kernel(s_2,y_2;t,x) \|_{L_{s_2,y_2}^2} 
\leq \| \kernel(s_2,y_2;t,x) \|_{L_{s_2,y_2}^1}^{\frac{1}{2}} \| \kernel(s_2,y_2;t,x) \|_{L_{s_2,y_2}^\infty}^{\frac{1}{2}} \\ 
&\lesssim \lossfactor_T(F)^4 L^{3} N \max(M,L)^{-2} \max(K,N)^{-2} 
\lesssim \lossfactor_T(F)^4 \min(L,N) \max(K,L,M,N)^{-1}.  
\end{align*}
In the last inequality, we first bounded $N$ by $\max(K,N)$ and then used that $L$ can be bounded by either $\max(M,L)$ or $\max(K,N)$. Clearly, this is sufficient for \eqref{objects:eq-quadratic-dlin-13-resonance-1}, and will later also be useful for~\eqref{objects:eq-quadratic-dlin-13-resonance-2}. In order to obtain \eqref{objects:eq-13-resonance-p1}, we use Propositions \ref{prop:pB-ptwise-bounds} and \ref{prop:qB-ptwise-bds}, which yield that 
\begin{align*}
& \abs{\kernel(s_2,y_2;t,x)} \\ 
&\lesssim \lossfactor_T(F)^4 \iint \ds_1 \dy_1 \iint \dr \dz \biggl( \chi_K(t-r) \chi_M(r-s_1) \chi_N(r-s_2) \chi_L(t-s_1) \\ 
&\quad\times (t-r)^{-1} (r-s_1)^{-1} (r-s_2)^{-\frac{3}{2}} (t-s_1)^{-\frac{3}{2}} \\ 
&\quad \times \exp\Bigl( - c \Bigl( \tfrac{\abs{x-z}^2}{t-r}+\tfrac{\abs{z-y_1}^2}{r-s_1} + \tfrac{\abs{z-y_2}^2}{r-s_2} + \tfrac{\abs{x-y_1}^2}{t-s_1} \Bigr) \Bigr) \biggr). 
\end{align*}
To estimate the integral, we first use the localising $\chi$-factors to pull the inverse powers of the time differences out of the integral, which contributes $K^2 M^2 N^3 L^3$. Then, we integrate over $y_1$ and $z$, which leads to a gain of $(r-s_1) \min(t-r,r-s_2) \sim M^{-2} \max(K,N)^{-2}$. Finally, we integrate over $s_1$ and $r$, which contributes $\max(M,L)^{-2} \max(K,N)^{-2}$. By combining these contributions, we arrive at \eqref{objects:eq-13-resonance-p1}.   In order to obtain \eqref{objects:eq-13-resonance-p2}, we integrate our estimate of $\abs{\kernel(s_2,y_2;t,x)}$ over $(s_2,y_2)$, which yields that
\begin{align*}
\|& \kernel(s_2,y_2;t,x) \|_{L_{s_2,y_2}^1} \\ 
&\lesssim \lossfactor_T(F)^4 \iint \ds_1 \dy_1 \iint \ds_2 \dy_2 \iint \dr \dz \biggl( \chi_K(t-r) \chi_M(r-s_1) \chi_N(r-s_2) \chi_L(t-s_1) \\ 
&\quad\times (t-r)^{-1} (r-s_1)^{-1} (r-s_2)^{-\frac{3}{2}} (t-s_1)^{-\frac{3}{2}} \\ 
& \quad\times \exp\Bigl( - c \Bigl( \tfrac{\abs{x-z}^2}{t-r}+\tfrac{\abs{z-y_1}^2}{r-s_1} + \tfrac{\abs{z-y_2}^2}{r-s_2} + \tfrac{\abs{x-y_1}^2}{t-s_1} \Bigr) \Bigr) \biggr). 
\end{align*}
As before, the inverse powers of the time differences contribute $K^2 M^2 N^3 L^3$. Then, we integrate over $y_1$, $y_2$, and $z$, which leads to a gain of $(r-s_1)(r-s_2)(t-r) \sim (MNK)^{-2}$. Finally, we integrate over $s_1$, $s_2$, and $r$, which contributes $(\max(M,L)NK)^{-2}$. By combining these contributions, we arrive at \eqref{objects:eq-13-resonance-p2}. 
It now only remains to prove the estimate \eqref{objects:eq-quadratic-dlin-13-resonance-2}. To this end, note that $\kernel(s_2,y_2;t,x)$ is supported on times satisfying $t-s_2 \sim K^{-2} + N^{-2} \sim \min(K,N)^{-2}$. By combining our existing estimates of $\kernel$ with Lemma \ref{object:lem-time-integration-gain}, we therefore obtain that 
\begin{align*}
\Bigl \|& \int_{t_0}^{t_1} \dt \redE \Bigl[ \bigl( \Duh_{\Y,B,K}^{F} \bigl( \Qc \bigl( \redlinear[\oeps][][][M], \covd^k \graylinear[\oeps][][][N] \bigr) \bigr) \otimes   \covd^\ell \redlinear[\oeps][][][L] \bigr) \Bigr]   (t,x)  \Bigr \|_{L^2(\Omega)}  \\ 
&\lesssim \lossfactor_T(F)^4 (t_1-t_0)^{\frac{1}{2}} \min\bigl( (t_1-t_0)^{\frac{1}{2}}, \min(K,N)^{-1}\bigr) \min(L,N) \max(K,L,M,N)^{-1}. 
\end{align*}
After estimating the minimum using a geometric mean and using that $K\gtrsim L$, we obtain \eqref{objects:eq-quadratic-dlin-13-resonance-2}. 
\end{proof}

In the following lemma, we control the term involving a resonance between the second and third linear stochastic objects, i.e., \eqref{object:eq-quadratic-dlin-decomposition-23}.

\begin{lemma}[Resonance between the second and third objects]\label{objects:lem-quadratic-dlin-23}
Let Assumption~\ref{local:assumption-B} be satisfied. Furthermore, let $k,\ell \in [2]$,  let $\Qc \colon \Y \times \Y \rightarrow \Y$ be a $\liegroup$-equivariant bilinear map, and let $\mathcal{K}$ be as in Proposition \ref{local:prop-paracontrolled}. Then, it holds for all $p\geq 2$ that
\begin{equation}\label{objects:eq-quadratic-dlin-23}
\begin{aligned}
\Bigl\|& \sup_{\eps \in \sdyadic} \,  \Bigl \| \redE \Bigl[ \Duh^F_{\Y,B,\star} \bigl(  \Qc\bigl(\graylinear[\oeps][][][] , \covd^k \redlinear[\oeps][][][] \bigr)\bigr) \parasim \covd^\ell \redlinear[\oeps] \Bigr] - \tfrac{1}{4} \delta^{k\ell} \bigl( \sigma_{\rho,\eps}^2 + \Crho \bigr) \mathcal{K}(\graylinear[\oeps][][][]) \Bigr\|_{\Cs_t^{-\theta}L_x^r} 
\Bigr\|_{L^p(\Omega)} \\ 
&\lesssim \lossfactor_T(F)^4 p^{\frac{1}{2}}. 
\end{aligned}
\end{equation}
Furthermore, the argument in \eqref{objects:eq-quadratic-dlin-23} has an a.s.\ limit in $\Cs_t^{-\theta}L_x^r$ as $\eps \rightarrow 0$. 
\end{lemma}

\begin{proof}
We begin with a decomposition of the argument in \eqref{objects:eq-quadratic-dlin-23}. By first adding and subtracting the term $\redE [ \Qc(\graylinear[\oeps][][][] , \covd^k \Duh_{\Y,B,\star} ( \redlinear[\oeps][][][] ))  \parasim \covd^\ell \redlinear[\oeps] ]$ and then using our Littlewood--Paley decompositions, we obtain that 
\begin{align}
&\quad \quad  \redE \Bigl[ \Duh^F_{\Y,B,\star} \bigl( \Qc\bigl(\graylinear[\oeps][][][] , \covd^k \redlinear[\oeps][][][] \bigr) \bigr) \parasim \covd^\ell \redlinear[\oeps] \Bigr]  -  \tfrac{1}{4} \delta^{k\ell} \bigl( \sigma_{\rho,\eps}^2 + \Crho \bigr) \mathcal{K}(\graylinear[\oeps][][][]) \notag \allowdisplaybreaks[4]\\[1ex]
& = \sum_{\substack{K,L,M,N\colon \\ K\sim_\theta L,K \gtrsim L, \\ M > \min(K,N)^{\reg_1} }} 
\redE \Bigl[ \Duh^F_{\Y,B,K} \bigl( \Qc\bigl(\graylinear[\oeps][][][M] , \covd^k \redlinear[\oeps][][][N] \bigr) \bigr) \otimes  \covd^\ell \redlinear[\oeps][][][L] \Bigr] \label{objects:eq-quadratic-dlin-q1} \allowdisplaybreaks[4]\\ 
&- \sum_{\substack{K,L,M,N\colon \\ K\sim_\theta L,K \gtrsim L, \\ M > \min(K,N)^{\reg_1} }} 
  \redE \Bigl[ \Qc\bigl(\graylinear[\oeps][][][M] , \covd^k \Duh_{\Y,B,K} \bigl( \redlinear[\oeps][][][N] \bigr) \bigr) \otimes  \covd^\ell \redlinear[\oeps][][][L] \Bigr] \label{objects:eq-quadratic-dlin-q2} \allowdisplaybreaks[4]\\
&+  \hspace{-3ex}\sum_{\substack{K,L,M,N\colon \\ K\sim_\theta L,K \gtrsim L, \\ M \leq  \min(K,N)^{\reg_1} }} \hspace{-2ex} \redE \Bigl[ \Bigl(  \Duh^F_{\Y,B,K} \bigl( \Qc\bigl(\graylinear[\oeps][][][M] , \covd^k \redlinear[\oeps][][][N] \bigr) \bigr)  
- \Qc\bigl(\graylinear[\oeps][][][M] , \covd^k \Duh_{\Y,B,K} \bigl( \redlinear[\oeps][][][N] \bigr) \bigr) \Bigr)  \otimes \covd^\ell \redlinear[\oeps][][][L] \Bigr] \label{objects:eq-quadratic-dlin-q3} \allowdisplaybreaks[4]\\
&+ \sum_{\substack{K,L\colon \\ K\sim_\theta L,K \gtrsim L }} \redE \Bigl[ \Qc\bigl(\graylinear[\oeps][][][] , \covd^k \Duh_{\Y,B,K} \bigl( \redlinear[\oeps][][][] \bigr) \bigr)  \otimes \covd^\ell \redlinear[\oeps][][][L]\Bigr] -  \tfrac{1}{4} \delta^{k\ell} \bigl( \sigma_{\rho,\eps}^2 + \Crho \bigr) \mathcal{K}(\graylinear[\oeps][][][]). 
\label{objects:eq-quadratic-dlin-q4}
\end{align}
In the following, we estimate the four terms in \eqref{objects:eq-quadratic-dlin-q1}--\eqref{objects:eq-quadratic-dlin-q4} separately. 

\proofstep{Estimate of \eqref{objects:eq-quadratic-dlin-q1}:}
In order to control  \eqref{objects:eq-quadratic-dlin-q1}, it suffices to prove for all $K,L,M,N\geq N_{\lo}$ that
\begin{equation}\label{objects:eq-quadratic-dlin-q1a}
\begin{aligned}
\bigl\| \bigl\| &\redE \bigl[ \Duh^F_{\Y,B,K} \bigl( \Qc\bigl(\graylinear[\oeps][][][M] , \covd^k \redlinear[\oeps][][][N] \bigr) \bigr) \otimes  \covd^\ell \redlinear[\oeps][][][L] \bigr] \bigr\|_{\Cs_t^{-\theta} L_x^r} \bigr\|_{L^p(\Omega)} \\ 
&\lesssim \lossfactor_T(F)^4 \sqrt{p} \min(K,M)^{-2(\theta-\kappa)} \min(L,N) \max(K,L,M,N)^{-1}. 
\end{aligned}
\end{equation}
Using the restrictions on the dyadic scales in \eqref{objects:eq-quadratic-dlin-q1}, we can then estimate $\min(K,M)\gtrsim \min(L,N)^{\reg_1}$, which allows us to bound the dyadic factors on the right-hand side of \eqref{objects:eq-quadratic-dlin-q1a} by $\max(K,L,M,N)^{-2\reg_1 (\theta-\kappa)}$. By using Lemma \ref{object:lem-time-integration-gain}, Lemma \ref{lemma:time-regularity-kolmogorov}, and Gaussian hypercontractivity, the estimate \eqref{objects:eq-quadratic-dlin-q1a} can be further reduced to 
\begin{equation}\label{objects:eq-quadratic-dlin-q1b}
\begin{aligned}
\Bigl \|& \redE \Bigl[ \bigl( \Duh_{\Y,B,K}^{F} \bigl( \Qc \bigl( \graylinear[\oeps][][][M], \covd^k \redlinear[\oeps][][][N] \bigr) \bigr) \otimes   \covd^\ell \redlinear[\oeps][][][L]  \bigr) (t,x) \Bigr] \Bigr \|_{L^2(\Omega)} \\ 
&\lesssim \lossfactor_T(F)^4 \min(L,N) \max(K,L,M,N)^{-1}. 
\end{aligned}
\end{equation}
The proof of \eqref{objects:eq-quadratic-dlin-q1b}, which we now address, is similar to the proof of Lemma \ref{objects:lem-quadratic-dlin-13-resonance}. We first write 
\begin{align*}
\redE \Bigl[ \bigl( \Duh_{\Y,B,K}^{F} \bigl( \Qc \bigl( \graylinear[\oeps][][][M], \covd^k \redlinear[\oeps][][][N] \bigr) \bigr) \otimes   \covd^\ell \redlinear[\oeps][][][L]  \bigr) \Bigr] (t,x)
= \iint \ds_1 \dy_1 \kernel(s_1,y_1;t,x) \boldxi(s_1,y_1),
\end{align*}
where $\kernel(s_1,y_1;t,x) \colon \Y \rightarrow \Y \otimes \Y$ is defined as 
\begin{align*}
\kernel(&s_1,y_1;t,x) (\place) 
\eqdef \iint \ds_2 \dy_2 \iint \dr \dz \Bigl( \fq_{\Y,B,K}(r,z;t,x)  \\ 
&\Qc\bigl(\fq_{\Y,B,\eps,M}(s_1,y_1;r,z)(\place), \covd^k \fq_{\Y,B,\eps,N}(s_2,y_2;r,z) \ONBY_m \bigr) \otimes \covd^\ell \fq_{\Y,B,\eps,L}(s_2,y_2;t,x) \ONBY^{m} \Bigr). 
\end{align*}
From the definition, it directly follows that the kernel is only nonzero in the case $K\gtrsim L$. We henceforth restrict to this case. As before, we control the kernel in $L_{s_1,y_1}^\infty$ and $L_{s_1,y_1}^1$, i.e., we prove the two estimates 
\begin{align}
\| \kernel(s_1,y_1;t,x) \|_{L_{s_1,y_1}^\infty} &\lesssim \lossfactor_T(F)^4 K^2 L^3 M^2 N \max(K,M)^{-4} \max(L,N)^{-2}, \label{objects:eq-23-resonance-q1c} \\ 
\| \kernel(s_1,y_1;t,x) \|_{L_{s_1,y_1}^1} &\lesssim \lossfactor_T(F)^4 K^{-2} L^3 M^{-2} N \max(L,N)^{-2}. \label{objects:eq-23-resonance-q1d}
\end{align} 
Once \eqref{objects:eq-23-resonance-q1c} and \eqref{objects:eq-23-resonance-q1d} have been proven, the desired estimate \eqref{objects:eq-quadratic-dlin-q1b} follows from 
\begin{align*}
K^2& L^3 M^2 N \max(K,M)^{-4} \max(L,N)^{-2} \times K^{-2} L^3 M^{-2} N \max(L,N)^{-2} \\
&= L^6 N^2 \max(L,N)^{-4} \max(K,M)^{-4} \lesssim \min(L,N)^2 \max(K,L,M,N)^{-2}. 
\end{align*}
In order to prove \eqref{objects:eq-23-resonance-q1c}, we use Propositions \ref{prop:pB-ptwise-bounds} and \ref{prop:qB-ptwise-bds}, which imply that 
\begin{align*}
& \abs{\kernel(s_1,y_1;t,x)} \\ 
&\lesssim \lossfactor_T(F)^4 \iint \ds_2 \dy_2 \iint \dr \dz \biggl( \chi_K(t-r) \chi_M(r-s_1) \chi_N(r-s_2) \chi_L(t-s_2) \\ 
&\times (t-r)^{-1} (r-s_1)^{-1} (r-s_2)^{-\frac{3}{2}} (t-s_2)^{-\frac{3}{2}} \\ 
& \times \exp\Bigl( - c \Bigl( \tfrac{\abs{x-z}^2}{t-r}+\tfrac{\abs{z-y_1}^2}{r-s_1} + \tfrac{\abs{z-y_2}^2}{r-s_2} + \tfrac{\abs{x-y_2}^2}{t-s_2} \Bigr) \Bigr) \biggr). 
\end{align*}
To estimate the integral, we first use the localising $\chi$-factors to pull the inverse powers of the time differences out of the integral, which contributes $K^2 M^2 N^3 L^3$. Then, we integrate over $y_2$ and $z$, which leads to a gain of $(r-s_2) \min(t-r,r-s_1) \sim N^{-2} \max(K,M)^{-2}$. Finally, we integrate over $s_2$ and $r$, which contributes $\max(N,L)^{-2} \max(K,M)^{-2}$. By combining these contributions, we arrive at \eqref{objects:eq-23-resonance-q1c}.   In order to obtain \eqref{objects:eq-23-resonance-q1d}, we again use Propositions \ref{prop:pB-ptwise-bounds} and \ref{prop:qB-ptwise-bds}, which yield that 
\begin{align*}
& \iint \ds_1 \dy_1 \abs{\kernel(s_1,y_1;t,x)} \\ 
&\lesssim \lossfactor_T(F)^4 \iint \ds_1 \dy_1 \iint \ds_2 \dy_2 \iint \dr \dz \biggl( \chi_K(t-r) \chi_M(r-s_1) \chi_N(r-s_2) \chi_L(t-s_2) \\ 
&\times (t-r)^{-1} (r-s_1)^{-1} (r-s_2)^{-\frac{3}{2}} (t-s_2)^{-\frac{3}{2}} \\ 
& \times \exp\Bigl( - c \Bigl( \tfrac{\abs{x-z}^2}{t-r}+\tfrac{\abs{z-y_1}^2}{r-s_1} + \tfrac{\abs{z-y_2}^2}{r-s_2} + \tfrac{\abs{x-y_2}^2}{t-s_2} \Bigr) \Bigr) \biggr). 
\end{align*}
As before, the inverse powers of the time differences contribute $K^2 M^2 N^3 L^3$. Then, we integrate over $y_1$, $y_2$, and $z$, which leads to a gain of $(r-s_1)(r-s_2)(t-r) \sim (MNK)^{-2}$. Finally, we integrate over $s_2$, $s_1$, and $r$, which contributes $(\max(L,N)MK)^{-2}$. By combining the contributions, we arrive at~\eqref{objects:eq-23-resonance-q1d}.

\proofstep{Estimate of \eqref{objects:eq-quadratic-dlin-q2}:} We first note that the summands in \eqref{objects:eq-quadratic-dlin-q2} can only be non-zero if $K\gtrsim L$ and $N\gtrsim L$, which follows from similar arguments as those below \eqref{object:eq-dduhlin-lin-p1}. Since \eqref{objects:eq-quadratic-dlin-q2} already contains the constraint $M>\min(K,N)^{\reg_1}$, we can therefore include the constraint $M\gtrsim L^{\reg_1}$. By using Lemma \ref{object:lem-dduhlin-dlin-crude} and Lemma \ref{objects:lem-deterministic-and-linear} (with $q=\kappa^{-1}$), we obtain that 
\begin{align*}
\bigl\| \bigl\| \eqref{objects:eq-quadratic-dlin-q2} \bigr\|_{\Cs_t^{-\theta} L_x^r} \bigr\|_{L^p(\Omega)} 
\lesssim \lossfactor_T(F)^4 \sqrt{p} \sum_{\substack{K,L,M,N\colon \\ M \gtrsim L^{\reg_1} }} M^{-2(\theta-2\kappa)} \times L \max(K,L,N)^{-1}. 
\end{align*}
Due to the restriction $M\gtrsim L^{\reg_1}$, we can estimate $M^{-2(\theta-2\kappa)} L$ by $L^{1-2\reg_1(\theta-2\kappa)}$, which then clearly leads to a summable contribution. 

\proofstep{Estimate of \eqref{objects:eq-quadratic-dlin-q3}:}
In this dyadic sum, we have that $M\leq N^{\reg_1}$. This clearly implies $M\leq N^{1-\theta}$, and hence $M\ll_\theta N$. Thus, we can apply Lemma \ref{local:lem-X-structure} (with $\eta_1=0$), which yields that
\begin{align*}
\Bigl\|& \redE \Bigl[ \Bigl(  \Duh^F_{\Y,B,K} \bigl( \Qc\bigl(\graylinear[\oeps][][][M] , \covd^k \redlinear[\oeps][][][N] \bigr) \bigr)  
- \Qc\bigl(\graylinear[\oeps][][][M] , \covd^k \Duh_{\Y,B,K} \bigl( \redlinear[\oeps][][][N] \bigr) \bigr) \Bigr) \otimes \covd^\ell \redlinear[\oeps][][][L]  \Bigr] \Bigl\|_{L_{t,x}^\infty} \\ 
&\leq \redE \Bigl[ \Bigl\| \Bigl(  \Duh^F_{\Y,B,K} \bigl( \Qc\bigl(\graylinear[\oeps][][][M] , \covd^k \redlinear[\oeps][][][N] \bigr) \bigr)  
- \Qc\bigl(\graylinear[\oeps][][][M] , \covd^k \Duh_{\Y,B,K} \bigl( \redlinear[\oeps][][][N] \bigr) \bigr) \Bigr) \otimes \covd^\ell \redlinear[\oeps][][][L]  \Bigr\|_{L_{t,x}^\infty} \Bigr] \\  
&\lesssim  \redE \Bigl[ \Bigl\| \Duh^F_{\Y,B,K} \bigl( \Qc\bigl(\graylinear[\oeps][][][M] , \covd^k \redlinear[\oeps][][][N] \bigr) \bigr)  
- \Qc\bigl(\graylinear[\oeps][][][M] , \covd^k \Duh_{\Y,B,K} \bigl( \redlinear[\oeps][][][N] \bigr) \bigr)\Bigr\|_{L_{t,x}^\infty}
\bigl\| \covd^\ell \redlinear[\oeps][][][L]\bigr\|_{L_{t,x}^\infty} \Bigr] \\ 
&\lesssim T^{\frac{1-\lambda_1}{2}} K^{-1-2\theta} L^{1+\kappa} N^{-\kappa} \bigl\| \graylinear[\oeps][][][M] \bigr\|_{\sol^{1-\reg_1}_M} \redE \bigl[ \bigl\| \redlinear[\oeps][][][\ostar] \bigr\|_{\modsol^{-\kappa}}^2 \bigr] \\ 
&\lesssim  T^{\frac{1-\lambda_1}{2}} K^{-1-2\theta} L^{1+\kappa} M^{1-\reg_1+\kappa} N^{-\kappa} \bigl\| \graylinear[\oeps][][][\ostar]\bigr\|_{\modsol^{-\kappa}}
 \redE \bigl[ \bigl\| \redlinear[\oeps][][][\ostar] \bigr\|_{\modsol^{-\kappa}}^2 \bigr]. 
\end{align*}
Since $K \gtrsim L$ and $M\leq K^{\reg_1}$, we have that 
\begin{equation*}
K^{-1-2\theta} L^{1+\kappa} M^{1-\reg_1+\kappa} N^{-\kappa}
\lesssim K^{-1-2\theta+1+\kappa+\reg_1 (1-\reg_1+\kappa)} N^{-\kappa} 
\lesssim \max(K,L,M,N)^{-\kappa},
\end{equation*}
where we used that $\kappa$ and $\reg_1$ are much smaller than $\theta$. Since this is summable, and since the $\modsol^{-\kappa}$-norm of the linear stochastic object has already been controlled in \eqref{eq:linear-object-solution-space-estimate}, this yields an acceptable contribution. 

\proofstep{Estimate of \eqref{objects:eq-quadratic-dlin-q4}:}
We first rewrite \eqref{objects:eq-quadratic-dlin-q4} as 
\begin{equation}\label{objects:eq-quadratic-dlin-p8}
\boldsymbol{L} \Bigl( \graylinear[\oeps][][][] \otimes \Bigl( \sum_{\substack{K,L\colon \\ K\sim_\theta L}} \redE \bigl[ \covd^k \Duh_{\Y,B,K} \bigl( \redlinear[\oeps] \bigr) \otimes \covd^\ell \redlinear[\oeps][][][L] \bigr] - \tfrac{1}{4} \delta^{k\ell} \bigl( \sigma_{\rho,\eps}^2 + \Crho \bigr) \Cov_{\Y} \Bigr) \Bigr), 
\end{equation}
where $\boldsymbol{L}\colon \Y \otimes \Y \otimes \Y \rightarrow \Y \otimes \Y$ is defined as $\boldsymbol{L}(X\otimes Y \otimes Z) \eqdef \Qc\bigl( X, Y \bigr) \otimes Z$. The term in \eqref{objects:eq-quadratic-dlin-p8} is the only term in this proof for which we need to be careful about the convergence as $\eps\rightarrow 0$, as there will be no decay in the highest frequency scale. For this reason, and to simplify the notation, we therefore set\footnote{This definition will only be used in this proof, and the term should not be confused with the remainders in the statement and proof of Proposition \ref{local:prop-main}.}  
\begin{equation*}
R_{\eps}\coloneq \sum_{\substack{K,L\colon \\ K\sim_\theta L}} \redE \bigl[\covd^k \Duh_{\Y,B,K} \bigl( \redlinear[\oeps] \bigr) \otimes \covd^\ell \redlinear[\oeps][][][L] \bigr] - \tfrac{1}{4} \delta^{k\ell} \bigl( \sigma_{\rho,\eps}^2 + \Crho \bigr) \Cov_{\Y}. 
\end{equation*}
Let $\eps,\eps^\prime \in \sdyadic$ satisfy $\eps^\prime < \eps$ and recall that $\graylinear[\oeps][][][M]=\graylinear[\oepsprime][][][M]$ for all $M \lesssim \eps^{-1}$. Furthermore, set $q\coloneq \kappa^{-1}$. Using Proposition \ref{object:prop-dduhlin-dlin-main} and Lemma \ref{objects:lem-deterministic-and-linear}, we then obtain that 
\begin{equation*}
\bigl\| \bigl\| \graylinear[\oeps] \otimes R_{\eps} \bigr\|_{\Cs_t^{-\theta} L_x^r} \bigr\|_{L^p(\Omega)} \lesssim \lossfactor_T(F) \sqrt{p} \| R_{\eps} \|_{L_t^q L_x^\infty} \lesssim \lossfactor_T(F)^4 \sqrt{p}   
\end{equation*}
as well as
\begin{align*}
\bigl\|& \bigl\| \graylinear[\oeps] \otimes R_{\eps} - \graylinear[\oepsprime] \otimes R_{\eps^\prime} \bigr\|_{\Cs_t^{-\theta} L_x^r} \bigr\|_{L^p(\Omega)} \\ 
&\lesssim \bigl\| \bigl\| \bigl( \graylinear[\oeps] - \graylinear[\oepsprime]\bigr) \otimes R_{\eps}  \bigr\|_{\Cs_t^{-\theta} L_x^r} \bigr\|_{L^p(\Omega)}
+ \bigl\|  \bigl\|  \graylinear[\oepsprime] \otimes \bigl( R_{\eps} - R_{\eps^\prime} \bigr)  \bigr\|_{\Cs_t^{-\theta} L_x^r} \bigr\|_{L^p(\Omega)} \\ 
&\lesssim \sqrt{p} \lossfactor_T(F) \eps^{2(\theta-2\kappa)} \| R_{\eps} \|_{L_t^q L_x^\infty} + \sqrt{p} \lossfactor_T(F) \| R_{\eps} - R_{\eps^\prime}\|_{L_t^q L_x^\infty} \\ 
&\lesssim \sqrt{p} \lossfactor_T(F)^4 \eps^{\kappa}. 
\end{align*}
Thus, \eqref{objects:eq-quadratic-dlin-q4} clearly yields an acceptable contribution to \eqref{objects:eq-quadratic-dlin-23} and furthermore, since $\eps$ has been restricted to $\sdyadic$, it almost surely converges as $\eps\rightarrow 0$.   
\end{proof}

Equipped with Lemmas \ref{objects:lem-quadratic-dlin-nonres}-\ref{objects:lem-quadratic-dlin-23}, we can present the (now very short) proof of the second main proposition of this subsection. 

\begin{proof}[Proof of Proposition \ref{objects:prop-quadratic-dlin}]
After using our Wiener chaos decomposition, it suffices to control the four terms in 
\eqref{object:eq-quadratic-dlin-decomposition-nonres}, 
\eqref{object:eq-quadratic-dlin-decomposition-12},
\eqref{object:eq-quadratic-dlin-decomposition-13},
and \eqref{object:eq-quadratic-dlin-decomposition-23}. 
This is done using Lemmas \ref{objects:lem-quadratic-dlin-nonres}, \ref{objects:lem-quadratic-dlin-12-resonance}, \ref{objects:lem-quadratic-dlin-13-resonance}, and \ref{objects:lem-quadratic-dlin-23}, respectively. The almost sure convergence as $\eps \rightarrow 0$ follows from the corresponding convergence statements in Lemmas \ref{objects:lem-quadratic-dlin-12-resonance} and \ref{objects:lem-quadratic-dlin-23}, as well as the decay in the highest frequency scale in Lemmas \ref{objects:lem-quadratic-dlin-nonres} and~\ref{objects:lem-quadratic-dlin-13-resonance}. 
\end{proof}

To finish this section, we show that Proposition \ref{objects:prop-quadratic-dlin} implies the needed stochastic estimates for $\Xi_{\Beps}^{\der}$ from \eqref{local:eq-data-set-der-1}--\eqref{local:eq-data-set-der-3}. 

\begin{corollary}\label{objects:cor-xi-der}
Let Assumption~\ref{local:assumption-B} be satisfied. Then, it holds for all $p\geq 2$ that 
\begin{align}
\Bigl\|  \sup_{\eps \in \sdyadic}\,  \Bigl\| \Bigl(& 
\Lb\bigl( \quadratic[\oeps][][][\ostar],\covd \linear[\oeps][][][\ostar] \bigr)
+ ( \sigma_{\rho,\eps}^2+\Crho) \Cov_{\frkg}^{\frkg} (\linear[\oeps][][][]), 
\quad  
\Lb\bigl( \phiquadratic[\oeps][][][\ostar],\covd \linear[\oeps][][][\ostar]\bigr), \notag \\ 
&\bil \bigl( \mixedquadratic[\oeps][][][\ostar] \parasim \covd\philinear[\oeps][][][\ostar]\bigr) 
+\tfrac{1}{2} (\sigma_{\rho,\eps}^2 + \Crho) \Cov^{\frkg}_{\V}( \linear[\oeps][][][]),  
\quad \quadratic[\oeps][\ell][][\ostar] \parasim \covd_\ell \philinear[\oeps][][][\ostar], \label{objects:eq-xi-der} \\ 
&\phiquadratic[\oeps][\ell][][\ostar] \parasim \covd_\ell \philinear[\oeps][][][\ostar]
+\tfrac{1}{2} (\sigma_{\rho,\eps}^2 + \Crho) \Cov^{\V}_{\frkg}( \philinear[\oeps][][][]),
\quad 
\bigl( \covd_\ell \linear[\oeps][\ell][][\ostar] \bigr) \parasim    \mixedquadratic [\oeps][][][\ostar] \Bigr)
\Bigr\|_{\Cs_t^{-\theta} L_x^r}
\Bigr\|_{L^p(\Omega)} \lesssim T^{-4\lossdel} p^{\frac{3}{2}}.  \notag 
\end{align}
Furthermore, the argument in \eqref{objects:eq-xi-der} has an a.s.\ limit in $\Cs_t^{-\theta}L_x^r$ as $\eps\rightarrow 0$. 
\end{corollary}

\begin{proof}
For all six terms appearing in \eqref{objects:eq-xi-der}, the estimate and convergence follow by first using Proposition \ref{objects:prop-quadratic-dlin} and then using the same calculations as in the proof of Corollary \ref{local:cor-paracontrolled-specific} to verify that \eqref{objects:eq-xi-der} contains the correct renormalisation terms. To illustrate this, we prove the estimate for the (rather representative) third term in \eqref{objects:eq-xi-der}, but we will omit the details for the remaining terms. We first define $\Qc_j \colon \Y \times \Y \rightarrow \Y$, where $j\in [2]$, and $\widetilde{\Qc}\colon \Y \otimes \Y \rightarrow \frkg$ by
\begin{align*}
\Qc_j\bigl( (X_1,\, X_2,\, \Phi) , (Y_1,\, Y_2,\, \Psi) \bigr) 
&\eqdef \bigl( 0, \, 0, \, 2 X_j \Psi\bigr), \\ 
\widetilde{\Qc} \bigl( (X_1,\, X_2,\, \Phi) \otimes (Y_1,\, Y_2,\, \Psi) \bigr) 
&\eqdef \bil( \Phi \otimes \Psi). 
\end{align*}
Using this notation, we can write 
\begin{equation*}
\bil ( \mixedquadratic[\oeps][][][\ostar] \parasim \covd^k \philinear[\oeps][][][\ostar]) 
= \widetilde{\Qc} \Bigl( \Duh^F_{\Y,B,\star} \bigl( \Qc_j \bigl( \graylinear[\oeps][][][], \covd^j \graylinear[\oeps][][][] \bigr) \bigr) \parasim \covd^k \graylinear[\oeps][][][\ostar] \Bigr). 
\end{equation*}
To see that\footnote{The easiest way to convince oneself of this is to remember that $\mixedquadratic[\oeps][][][\ostar]$ is defined using the two independent objects $\linear[\oeps][][][]$ and $\philinear[\oeps][][]$, and hence the quadratic object itself should not produce a renormalisation term.} the condition $\Qc(\ONBY_m,\ONBY^m)=0$ from Proposition \ref{objects:prop-quadratic-dlin} is satisfied for all $j\in [2]$, we first note that this expression is independent of the choice of the orthonormal basis. Since $\Y = \frkg \oplus \frkg \oplus \V$, we may then choose $(\ONBY_m)_m$ as the orthonormal basis consisting of $((\ONBg_a,0,0))_a$, $((0,\ONBg_a,0))_a$, and $((0,0,\ONBV_n))_n$, where $(\ONBg_a)_a$ and $(\ONBV_n)_n$ are orthonormal bases of $\frkg$ and $\V$. Since the definition of $\Qc_j$ involves different components of the two arguments, one then readily sees that $\Qc_j(\ONBY_m,\ONBY^m)=0$.  We mention that, in order to complete the same step for some of the other terms in \eqref{objects:eq-xi-der}, we also need the skew-symmetry of the Lie-bracket and $\bil$. Using Proposition \ref{objects:prop-quadratic-dlin}, we obtain for all $p\geq 2$ that 
\begin{align*}
\Bigl\| \sup_{\eps \in \sdyadic} \Bigl\| \widetilde{\Qc} \Bigl(& \Duh^F_{\Y,B,\star} \bigl( \Qc_j \bigl( \graylinear[\oeps][][][], \covd^j \graylinear[\oeps][][][] \bigr) \bigr) \parasim \covd^k \graylinear[\oeps][][][\ostar] \Bigr)  \\ 
&- \tfrac{1}{4} \delta^{jk} (\sigma_{\rho,\eps}^2 + \Crho) \widetilde{\Qc} \Bigl( \Qc_j \bigl( \graylinear[\oeps], \ONBY_m \bigr) \otimes \ONBY^m \Bigr) \Bigr\|_{\Cs_t^{-\theta} L_x^r} \Bigr\|_{L^p(\Omega)} \lesssim T^{-4\lossdel} p^{\frac{3}{2}}. 
\end{align*}
As in the proof of Corollary \ref{local:cor-paracontrolled-specific}, we obtain that 
\begin{equation*}
- \tfrac{1}{4} \delta^{jk}  \widetilde{\Qc} \Bigl( \Qc_j \bigl( \graylinear[\oeps], \ONBY_m \bigr) \otimes \ONBY^m \Bigr)
= - \tfrac{1}{2} \delta^{jk} \bil \bigl( \linear[\oeps][][j][] \ONBV_n \otimes \ONBV^n \bigr)  
= \tfrac{1}{2} \bil \bigl( \ONBV_n \otimes \linear[\oeps][k][][] \ONBV^n \bigr) 
= \tfrac{1}{2} \Cov^{\frkg}_{\V}( \linear[\oeps][k][][]), 
\end{equation*}
which yields the correct renormalisation term. 
\end{proof}

\section{Refined stochastic estimates}
\label{sec:improved_stochastic_estimates}
In this section, we obtain refined estimates, in particular 
bounds that are stronger than those suggested by naïve regularity counting,
 for several expressions involving expectations of stochastic objects that arise later on (see \eqref{e:iteration_deterministic_bounds_iteration1}--\eqref{e:iteration_deterministic_bounds_iteration4} in Proposition~\ref{pr:expansion_iteration}). 
The main ideas leading to the improved estimates are:
\begin{enumerate}
	\item We prefer to estimate objects in terms of $\covd_B F_B$ rather than $F_B$, since we have the improved bound due to our stopping criterion (see \eqref{e:iteration_F_endtime_lowerbound} in Lemma~\ref{le:iteration_basic_bounds}), which makes $\| \covd_B F_B \|_{L^2} $ satisfy the same estimate as $\|F_B \|_{L^2} $, thus gaining us one derivative. 
	This observation gives better bounds for \eqref{e:iteration_deterministic_bounds_iteration1} when combined with Lemma~\ref{lemma:curl-equation} (see Lemma~\ref{le:curl_square_expectation_smoothB}). 
	\item Spatial integration of the objects allows for better bounds if we estimate moments instead of applying deterministic product estimates. This gives better bounds for \eqref{e:iteration_deterministic_bounds_iteration2} and \eqref{e:iteration_deterministic_bounds_iteration4}.
	\item The term \eqref{e:iteration_deterministic_bounds_iteration3} turns out to be very negative when $B$ is almost a YM connection, which can be captured by the smallness of $\covd_B F_B$. Our intuition for this is that $\linear[B]$ has a certain probability of aligning with the unstable directions of the Hessian of the Yang--Mills action evaluated at $B$. Moreover, the larger $F_B$ is, the more unstable directions the Hessian will have, and thus the more negative the term \eqref{e:iteration_deterministic_bounds_iteration3} will become. This provides the decay mechanism needed to obtain the existence of an invariant measure for the stochastic Yang--Mills--Higgs process. The actual argument is more involved than the others, and we believe this observation is one of the most important inputs in this article, so we will put this in a separate subsection. 
\end{enumerate}

\begin{notation}\label{notation:tstop}
Let $\tstop \in (0, 1]$. 
In this section, unless otherwise specified, we always assume that $B$ is a YM gradient flow on $(0, \tstop)$ and $(\tstop, 2]$, where a jump is allowed at time $\tstop$, from $B(\tstop -)$ to $B(\tstop +)$. 
\end{notation}

Even with a jump, the kernels $\fp_B$, $\fq_B$ are still well-defined, and satisfy the bounds 
\begin{equ}
\label{e:jump_heat_kernel_bound}
	|\fp_B(s, y; t, x)| \leq \fp(s, y; t, x)\;, \quad |\fq_B(s, y; t, x)| \lesssim \gronwall_{[s,t]}(C_0 F_B)\, (t - s)^{-1} \exp\Bigl(-c \frac{|x - y|^2}{t - s} \Bigr) 
\end{equ}
for any $0\leq s < t \leq 2$, $x, y \in \T^2$. If $s < \tstop < t$, this follows by using the semigroup property of the heat kernels to introduce the intermediate time $\tstop$, and then applying \eqref{eq:zeroth-order-pB-ptwise-bd} and \eqref{eq:qB-ptwise-bounds} on each of the intervals $[s, \tstop]$, $[\tstop, t]$.
For higher order derivatives, we use the semigroup property at $\tstop$, \eqref{eq:qB-ptwise-bounds} with $k$ derivatives on $[\tstop, t]$ and \eqref{e:jump_heat_kernel_bound} on $[s, \tstop]$. For any $k \geq 1$ and $0 \leq s < \tstop < t \leq 2$, this argument gives
	\begin{equs}
	\label{e:jump_heat_kernel_bound_derivatives}
		\bigl|(\covd_{B(t)}^{(k)} \fq_B)(s, y; t, x)\bigr|
		&= \Bigl|\int_{\T^2} (\covd_{B(t)}^{(k)} \fq_B)(\tstop, z; t, x) \fq_B(s, y; \tstop, z) \,\dz\Bigr| \\
		&\lesssim_k \gronwall_t(C_0 F_B) \bigl(1 + (\Fnorm{F_B}_{[\tstop, t], k})^k\bigr) \cdot (t - \tstop)^{-\frac{k}{2}}(t - s)^{-1} \exp\Bigl(-c\frac{|x-y|^2}{t-s}\Bigr) \;,
	\end{equs}
where all covariant derivatives act in the $x$-variable.

\subsection{Expectations}

We first discuss the expectation of the squared curl and the squared divergence of the linear stochastic object. We start from the simpler case where $B$ is smooth throughout the full interval $[0, 2]$ as a warm up.
In our later applications $F_B$ will have a singularity at $t = 0$ and a jump at $t = \tstop$. We will address this more involved case in Lemma~\ref{le:curl_square_expectation}.

\begin{lemma}
\label{le:curl_square_expectation_smoothB}
There exists some constant $C$ such that the following holds.
	Let $B$ be a smooth solution to the YM gradient flow on the full interval $[0, 2]$ (in particular, we assume that it does not have a jump at $\tstop$). Then for any $t \in (\tstop, 2]$ and $x \in \T^2$, we have
	\begin{equs}
	\label{e:curl_square_expectation_smoothB}
		\,&\E \Big| \Big( \curl_{B(t)} \bigl( e^{(t - \tstop) \Delta^{F}_{B(\tstop +)}} \linear[B] (\tstop) \bigr) \Big)(t, x) \Big|^2 + \E \Big| \Big( \div_{B(t)} \bigl( e^{(t - \tstop) \Delta^{F}_{B(\tstop +)}} \linear[B] (\tstop) \bigr) \Big)(t, x) \Big|^2 \\
		&\leq C(t - \tstop)^{-1} \\
		&\quad+ C \gronwall_{\tstop}(2C_0 F_B)  \bigl(1 + \Fnorm{F_B}_{\tstop, 0}\bigr)^2 \cdot \tstop^{\frac43} (t - \tstop)^{-\frac13} \| \fp^{1/2}(\place, \place; t, x) \covd_{B} F_{B} \|_{L^2_{s, y}([0, \tstop])}^2 \\
		&\quad+ C \gronwall_{t}(2C_0 F_B)  \bigl(1 + \Fnorm{F_B}_{t, 0}\bigr)^2 \cdot \tstop^{\frac12} (t - \tstop)^{\frac12} \| \fp^{1/2}(\place, \place; t, x) \covd_{B} F_{B} \|_{L^2_{s,y} ([\tstop, t])}^2 \;. 
	\end{equs}
\end{lemma}
\begin{remark}\label{remark:curl_square_expectation_smoothB}
	The key points are that the constant $C$ before $(t - \tstop)^{-1}$ doesn't depend on the Gronwall constant, and that
	\begin{equ}
		\int_{\T^2} 2\| \fp^{1/2}(\place, \place; t, x) \covd_B F_B \|_{L^2_{s,y} ([0, t])}^2 \, \dx = \| F_B(0) \|^2_{L^2}  - \| F_B(\tstop -) \|_{L^2}^2 + \| F_B(\tstop +) \|_{L^2}^2 - \| F_B(t) \|^2_{L^2}  \;,
	\end{equ}
	so this term can be absorbed into the decay of $F_B$. Moreover, the second term on the right-hand side gains the factor $\tstop^{\frac43}$ at the cost of the factor $(t - \tstop)^{-\frac13}$. This will be needed in the later applications, where we wish to have $\tstop^{4/3} \| \covd_{\YMHFsmooth} F_{\YMHFsmooth} \|_{L^2_{s, y}([0, \tstop] )} \ll 1$.  
\end{remark}

\begin{proof}
	We first note that by Lemma~\ref{lemma:linear-object-rotation-invariance} and the rotational invariance of the kernel $\fq_B$, we have 
	\begin{equs}
	\Bigl( \covd_{B(t)}^1 &\bigl(e^{(t - \tstop) \Delta^{F}_{B(\tstop +)}} \linear[B](\tstop)\bigr)^1, \covd_{B(t)}^2 \bigl(e^{(t - \tstop) \Delta^{F}_{B(\tstop +)}} \linear[B](\tstop)\bigr)^2 \Bigr) \\
	&\eqlaw \Bigl( - \covd_{B(t)}^1 \bigl(e^{(t - \tstop) \Delta^{F}_{B(\tstop +)}} \linear[B](\tstop)\bigr)^2, \covd_{B(t)}^2 \bigl(e^{(t - \tstop) \Delta^{F}_{B(\tstop +)}} \linear[B](\tstop)\bigr)^1 \Bigr) \;. 
	\end{equs}
	Therefore, we have 
	\begin{equ}
		\E \Big| \Big( \curl_{B(t)} \bigl( e^{(t - \tstop) \Delta^{F}_{B(\tstop +)}} \linear[B] (\tstop) \bigr) \Big)(t, x) \Big|^2 = \E \Big| \Big( \div_{B(t)} \bigl( e^{(t - \tstop) \Delta^{F}_{B(\tstop +)}} \linear[B] (\tstop) \bigr) \Big)(t, x) \Big|^2 \;,
	\end{equ}
	and in the subsequent proof we only focus on estimating the $\curl_{B(t)}$ term.
	By the duality characterisation of the $L^2$ norm we have 
	\begin{equs}
		\,&\E \Big| \Big( \curl_{B(t)} \bigl( e^{(t - \tstop) \Delta^{F}_{B(\tstop +)}} \linear[B] (\tstop) \bigr) \Big)(t, x) \Big|^2 \\
		&= \int_0^{\tstop} \int_{\T^2} \Big| \int_{\T^2} (\curl_{B(t)} \fq_{B})(\tstop, \xstop; t, x) \fq_{B}(s, y; \tstop, \xstop) \, \md \xstop \Big|^2 \ds \dy \\
		& \lesssim \sup_{\|\h \|_{L^2_{t, x}([0, \tstop] \times \T^2)} = 1} |\curl_{B} u_{\h}(t, x)|^2 \;,
	\end{equs}
	where $\h$ is a one-form, and $u_{\h}: [0, 2] \times \T^2 \rightarrow \frkg^2$ solves
	\begin{equ}
	\label{e:uh_equ}
		(\ptl_t - \Delta_{B}) u_{\h} = 2 \ad(F_{B}) (u_{\h}) + \one_{\place \leq \tstop} \h \; , \qquad u_{\h}(0) = 0 \;.
	\end{equ}
	Some basic estimates for $u_\h$ are collected in Lemma~\ref{le:uh_bound}. Let $v_{\h} \eqdef \curl_{B} u_{\h}$. By Lemma~\ref{lemma:curl-equation}, we have that $v_{\h}$ solves
	\begin{equs}
		\partial_t v_{\h} = \Delta_{B} v_{\h} + (\curl_{B} \h) \one_{\place \leq \tstop} - 2 \ad(\covd_{B}^{2} F_B^{12} ) u_{\h}^2 - 2\ad(\covd_{B}^{1} F_B^{12} ) u_{\h}^1 \;.
	\end{equs}
	To estimate $v_{\h}$, by Duhamel's principle we have 
	\begin{equ}
	\label{e:vh_duhamel}
		v_{\h}(t, x) = \int_{0}^{t} \int_{\T^2} \fp_{B}(s, y; t, x) \Big( (\curl_{B} \h) \one_{\place \leq \tstop} - 2 \ad(\covd_{B}^{2} F_B^{12} ) u_{\h}^2 - 2\ad(\covd_{B}^{1} F_B^{12}) u_{\h}^1 \Big)(s, y) \ds \dy \;.
	\end{equ}
	For the part with $(\curl_{B} \h) \one_{\place \leq \tstop}$, applying \eqref{eq:integrated-covariant-monotonicity-square-ibp}, we get 
	\begin{equ}
		\Big \| \int_{0}^{\tstop} \int_{\T^2} \fp_{B}(s, y; \tstop, \place) \bigl( (\curl_{B} \h) \one_{\place \leq \tstop} \bigr) (s, y) \ds \dy \Big \|_{L^2_x} \leq 1 \;.
	\end{equ}
	By the diamagnetic inequality and the semigroup property of $\fp_B$, this implies that for any $t > \tstop$, 
	\begin{equ}
		\Big \| \int_{0}^{t} \int_{\T^2} \fp_{B}(s, y; t, \place) \bigl( (\curl_{B} \h) \one_{\place \leq \tstop} \bigr) (s, y) \ds \dy \Big \|_{L^\infty_x} \leq (t - \tstop)^{-\frac12} \;.
	\end{equ}
	For the remaining portion of \eqref{e:vh_duhamel}, we break the integral into two parts separated by $\tstop$. For any interval $I \subset [0, t]$, H\"older's inequality gives
	\begin{equs}
	 \Big| \int_{I} &\int_{\T^2} \fp_{B}(s, y; t, x) \Big( 2 \ad(\covd_{B}^{2} F_B^{12} ) u_{\h}^2 + 2\ad(\covd_{B}^{1} F_B^{12}) u_{\h}^1 \Big)(s, y) \ds \dy \Big| \\
	\label{e:curl_square_expectation_smoothB_mid}
		& \lesssim
		\| \fp^{1/2} (\place, \place; t, x) \covd_B F_B \|_{L^2_{s,y}(I)} \| \fp^{1/2} (\place, \place; t, x) \|_{L^2_s L^q_y(I)} \|u_\h \|_{L^\infty_s L^p_y(I)}\;.
	\end{equs}
	where $\frac1p + \frac1q = \frac12$. 
	We apply \eqref{e:curl_square_expectation_smoothB_mid} with $(p, q) = (6, 3)$ on $I = [0, \tstop]$, and with $(p, q) = (4, 4)$ on $I = [\tstop, t]$. The result then follows by Lemma~\ref{le:uh_bound}. 
\end{proof}

\begin{lemma}
\label{le:uh_bound}
	Suppose that $u_{\h}$ satisfies \eqref{e:uh_equ} with some $\h$ such that $\|\h\|_{L^2_{t, x}([0, \tstop] \times \T^2)} = 1$. Then for any $p \geq 2$ and $t \in (0, 2]$ we have 
	\begin{equ}
	\label{e:uh_estimate_lem}
		\sup_{s \in (0, t]} \bigl((s \wedge \tstop)^{-\frac1p} \|u_{\h}(s) \|_{L^p}\bigr) \lesssim \gronwall_{t}(C_0 F_B)  \bigl(1 + \Fnorm{F_B}_{t, 0}\bigr) \;. 
	\end{equ}
	Moreover, for any $t \in (\tstop, 2]$ and $0 < \alpha < 1, \beta > 0, \lambda \in (0, 1/2]$
	 we have 
	\begin{equs}
	\label{e:uh_derivative_estimate1}
		\,& \int_{0}^{\tstop}  s^{-\alpha} (t - s)^{\beta} \|\fp(s, \place; t, x)^{1/2} \covd_{B} u_{\h}  \|^2_{L_y^2} \ds \lesssim   \gronwall_{\tstop}(2C_0 F_B)  \bigl(1 + \Fnorm{F_B}_{\tstop, 0}\bigr)^2 \\
		&\qquad \qquad \qquad \qquad \cdot \int_{0}^{\tstop} s^{-\alpha + \lambda} (t - s)^{\beta - \lambda} \bigl( \|F_B(s) \|_{L_y^\infty} + (t - s)^{-\frac12} \|\h(s) \|_{L^2_y} \bigr) \ds \;. 
	\end{equs}
	In particular, we have 
	\begin{equ}
	\label{e:uh_derivative_estimate2}
		 \| s^{-\frac18} (t - s)^{\frac14} \fp(s, \place; t, x)^{\frac12} \covd_{B} u_{\h} \|_{L^2_{s, y}([0, \tstop])}^2 
		\lesssim  \gronwall_{\tstop}(2C_0 F_B)  \bigl(1 + \Fnorm{F_B}_{\tstop, 0}\bigr)^3 \cdot \tstop^{\frac18} t^{\frac18} \;.
	\end{equ}
\end{lemma}
\begin{proof}
	For the $L^p$ bound, by \eqref{e:jump_heat_kernel_bound}, Young's convolution inequality, Hölder, and the assumption that $\|\h\|_{L_{t,x}^2} = 1$, we get
	\begin{equs}
		\| u_{\h}(s) \|_{L^p} &\lesssim \gronwall_{s}(C_0 F_B) \bigl(1 + \Fnorm{F_B}_{s, 0}\bigr) \int_0^{s \wedge \tstop } (s - s')^{\frac1p - \frac12} \|\h(s') \|_{L^2} \ds'  \\
		&\lesssim  \gronwall_{s}(C_0 F_B) \bigl(1 + \Fnorm{F_B}_{s, 0}\bigr) (s \wedge \tstop)^{\frac1p}  \;. 
	\end{equs}
	For \eqref{e:uh_derivative_estimate1}, by \eqref{eq:covariant-monotonicity-square} we get 
	\begin{equ}
		\partial_t  \|K^{1/2} u_{\h}  \|^2_{L_y^2} +  \|K^{1/2} \covd_{B} u_{\h}  \|^2_{L_y^2} \lesssim  \|F_B  \|_{L_y^\infty}  \|K^{1/2} u_{\h} \|^2_{L_y^2} +  \|\h \|_{L_y^2} \|K u_{\h}\|_{L^2_y} \;,
	\end{equ}
	where in this proof we choose $K$ to be $\fp(\place, \place; t, x)$. From this we get that for any $s \in (0, \tstop)$,
	\begin{equs}
	\partial_s \bigl(s^{-\alpha}& (t - s)^{\beta} \|K^{1/2} u_{\h}  \|^2_{L_y^2} \bigr) + s^{-\alpha} (t - s)^{\beta} \|K^{1/2} \covd_{B} u_{\h}  \|^2_{L_y^2} \\
	\label{e:uh_bound_mid1}
	&\lesssim s^{-\alpha} (t - s)^{\beta} \bigl( \|F_B \|_{L_y^\infty}  \|K^{1/2} u_{\h} \|^2_{L_y^2} +  \|\h \|_{L_y^2} \|K u_{\h}\|_{L^2_y} \bigr) \;.
\end{equs}
	
	Note that by \eqref{e:uh_estimate_lem} we have $\lim_{s \searrow 0} s^{-\alpha} \|K^{1/2}(s) u_{\h}(s)  \|^2_{L_y^2} = 0$ for $\alpha < 1$. Integrating \eqref{e:uh_bound_mid1} and applying \eqref{e:uh_estimate_lem} twice (first with $p = \frac{2}{\lambda}$, then with $p = \frac{1}{\lambda}$), we get 
	\begin{equs}
		\int_{0}^{\tstop}  s^{-\alpha} (t - s)^{\beta} \|K^{1/2} \covd_{B} u_{\h}  \|^2_{L_y^2} \ds &\lesssim  \gronwall_{\tstop}(2C_0 F_B)  \bigl(1 + \Fnorm{F_B}_{\tstop, 0}\bigr)^2 \\
		& \cdot \int_{0}^{\tstop} s^{-\alpha + \lambda} (t - s)^{\beta - \lambda} \bigl( \|F_B(s) \|_{L_y^\infty} + (t - s)^{-\frac12} \|\h(s) \|_{L^2_y} \bigr) \ds \;. 
	\end{equs}
	The result follows.
	The bound \eqref{e:uh_derivative_estimate2} follows by taking $\alpha = \frac14$, $\beta = \frac12$ and $\lambda = \frac38$ and applying H\"older's inequality.
\end{proof}

Now we deal with the case where $B$ is a smooth YM gradient flow plus a rough perturbation. Furthermore, $B$ is allowed to jump at some time $\tstop \in (0, 2]$. 

\begin{lemma}
\label{le:curl_square_expectation}
    There exists a constant $C$ such that the following holds.
	Suppose $B$ is such that the restriction $B \restr (0, \tstop]$ is a YM gradient flow with a decomposition of the form
	\begin{equ}
		B = g \act \YMHFsmooth + \Ad_g \YMHFrough \;,
	\end{equ}
	where $\YMHFsmooth: [0, \tstop] \times \T^2 \rightarrow \fg^2$ is a smooth solution to the YM gradient flow, $g: [0, \tstop] \times \T^2 \rightarrow \liegroup$ is a continuous gauge transform, and $\YMHFrough$ satisfies
	\begin{equ}
		\| t^{\frac{\kappa}{2}} \YMHFrough \|_{L^\infty_{t, x}([0, \tstop])} + \| t^{\frac{1+\kappa}{2}} \curl_{\YMHFsmooth} \YMHFrough \|_{L^\infty_{t, x}([0, \tstop])} < \infty \;.
	\end{equ}
	Furthermore, suppose $B \restr [\tstop, 2]$ is a YM gradient flow from initial data $B(\tstop+)$. Also define $B(\tstop-) \eqdef \lim_{ t \nearrow \tstop} B(t)$. In particular, a jump at $\tstop$ from $B(\tstop-)$ to $B(\tstop+)$ is allowed. 
	Then for any $t \in (\tstop, 2]$ and $x \in \T^2$, we have
	\begin{equs}
		\E \Big| &\Big( \curl_{B(t)} \bigl( e^{(t - \tstop) \Delta^{F}_{B(\tstop +)}} \linear[B] (\tstop) \bigr) \Big)(t, x) \Big|^2 + \E \Big| \Big( \div_{B(t)} \bigl( e^{(t - \tstop) \Delta^{F}_{B(\tstop +)}} \linear[B] (\tstop) \bigr) \Big)(t, x) \Big|^2 \\
		&\leq C(t - \tstop)^{-1} \\
		&\quad+ C \gronwall_{\tstop}(2C_0 F_B)  \bigl(1 + \Fnorm{F_B}_{\tstop, 0}\bigr)^2 \cdot \tstop^{\frac43} (t - \tstop)^{-\frac13} \| \fp^{1/2}(\place, \place; t, x) \covd_{\YMHFsmooth} F_{\YMHFsmooth} \|_{L^2_s L^2_y([0, \tstop])}^2 \\
		&\quad+ C \gronwall_{t}(2C_0 F_B)  \bigl(1 + \Fnorm{F_B}_{t, 0}\bigr)^2 \cdot \tstop^{\frac12} (t - \tstop)^{\frac12} \| \fp^{1/2}(\place, \place; t, x) \covd_{B} F_{B} \|_{L^2_s L^2_y([\tstop, t])}^2 \\
		&\quad + C \gronwall_{\tstop}(2C_0 F_B)\bigl(1 + \Fnorm{F_B}_{\tstop, 0}\bigr)^2
	(t - \tstop)^{-1 + \kappa}\tstop^{1 - \kappa} \|B(\tstop+) - B(\tstop-)\|_{L^\infty_x}^2 \\
		&\quad + C \gronwall_{\tstop}(2C_0 F_B)  \bigl(1 + \Fnorm{F_B}_{\tstop, 1}\bigr)^4 \tstop^{\frac14} (t - \tstop)^{ - \frac38 - \kappa} t^{\frac18}  \\
		& \qquad \cdot \Big( \| s^{\frac{\kappa}{2} } \YMHFrough \|_{L^\infty_{s, y}([0, \tstop])}^2 (\Fnorm{F_{\YMHFsmooth}}_{\tstop, 0})^2 +  \| s^{\frac{1 + \kappa}{2} } \bigl( \curl_{\YMHFsmooth} \YMHFrough + [\YMHFrough_1, \YMHFrough_2] \bigr) \|_{L^\infty_{s, y}([0, \tstop])}^2 \Big) \;.
	\end{equs}
\end{lemma}
\begin{proof}
	Define $u_\h$ as in the proof of Lemma~\ref{le:curl_square_expectation_smoothB}. Let $v_\h = \curl_B u_{\h}$. 
	It suffices to bound $|v_\h(t, x)|^2$ when $\| \h \|_{L^2_{t, x}([0, \tstop] \times \T^2)} = 1$.
	Since now $B$ has a jump at $\tstop$, $v_\h$ also has a jump. Therefore, \eqref{e:vh_duhamel} acquires an additional boundary term due to that jump. Actually, we have 
	\begin{equs}
	v_{\h}(\tstop-, x) &= \int_{0}^{\tstop} \int_{\T^2} \fp_{B}(s, y; \tstop, x) \Big( \curl_{B} \h - 2 \ad(\covd_{B}^{2} F_B^{12} ) u_{\h}^2 - 2\ad(\covd_{B}^{1} F_B^{12}) u_{\h}^1 \Big)(s, y) \ds \dy \;, \\
	J_{\tstop} &\eqdef v_{\h}(\tstop+, x) - v_{\h}(\tstop-, x) \\
	&= \bigl(B_1(\tstop+) - B_1(\tstop-)\bigr)u_{\h}^2(\tstop, x) - \bigl(B_2(\tstop+) - B_2(\tstop-)\bigr)u_{\h}^1(\tstop, x) \;.
\end{equs}
	For $t > \tstop$, we have 
	\begin{equs}
		v_{\h}(t, x) &= e^{(t - \tstop) \Delta_B } v_{\h}(\tstop+) \\
		\label{e:vh_mod_last_term}
		&\quad- \int_{\tstop}^t \int_{\T^2} \fp_{B}(s, y; t, x) \Big( 2 \ad(\covd_{B}^{2} F_B^{12} ) u_{\h}^2 + 2\ad(\covd_{B}^{1} F_B^{12}) u_{\h}^1 \Big)(s, y) \ds \dy \;. 
	\end{equs}
	
	By exactly the same argument as in  \eqref{e:curl_square_expectation_smoothB_mid}, the square of the $L^\infty$ norm of the term \eqref{e:vh_mod_last_term} is bounded by the third term on the right-hand side of the lemma. We therefore focus on the part 
	\begin{equs}
		\label{e:vh_mod_1}
		e^{(t - \tstop) \Delta_B } v_{\h}(\tstop+) &= \int_{0}^{\tstop} \int_{\T^2} \fp_{B}(s, y; t, \place) (\curl_{B} \h)(s, y) \ds \dy \\
		&\quad- \int_{0}^{\tstop} \int_{\T^2} \fp_{B}(s, y; t, \place) \Big( 2 \ad(\covd_{B}^{2} F_B^{12} ) u_{\h}^2 + 2\ad(\covd_{B}^{1} F_B^{12}) u_{\h}^1 \Big)(s, y) \ds \dy \\
		\label{e:vh_mod_2}\\
		\label{e:vh_mod_3}
		 &\quad+ e^{(t - \tstop) \Delta_B } J_{\tstop} \;.
	\end{equs}
	For \eqref{e:vh_mod_1}, the square of its $L^\infty$ norm is bounded by $C(t - \tstop)^{-1}$, using exactly the same argument as in Lemma~\ref{le:curl_square_expectation_smoothB}. For \eqref{e:vh_mod_3}, let $p_\kappa \eqdef \frac{2}{1-\kappa}$. By the diamagnetic inequality and Lemma~\ref{le:uh_bound} with $p = p_\kappa$, we have
\begin{equs}
	\bigl\|e^{(t - \tstop)\Delta_B} J_{\tstop}\bigr\|_{L^\infty_x}^2
	&\lesssim (t - \tstop)^{-\frac{2}{p_\kappa}}\|J_{\tstop}\|_{L^{p_\kappa}_x}^2 \\
	&\lesssim (t - \tstop)^{-1+\kappa}\|B(\tstop+) - B(\tstop-)\|_{L^\infty_x}^2
	\|u_{\h}(\tstop)\|_{L^{p_\kappa}_x}^2 \\
	&\lesssim \gronwall_{\tstop}(2C_0 F_B)\bigl(1 + \Fnorm{F_B}_{\tstop, 0}\bigr)^2
	(t - \tstop)^{-1+\kappa}\tstop^{1-\kappa}\|B(\tstop+) - B(\tstop-)\|_{L^\infty_x}^2 \;.
\end{equs}
	
	It remains to bound \eqref{e:vh_mod_2}. We first decompose
\begin{equs}
	\covd_B F^{12}_B &= \covd_B \Ad_g F^{12}_{\YMHFsmooth} + \covd_B \Ad_g \bigl( \curl_{\YMHFsmooth} \YMHFrough + [\YMHFrough_1, \YMHFrough_2]\bigr) \\
	\label{e:DB_FB_expansion}
	&= \Ad_g \bigl( \covd_{\YMHFsmooth} F^{12}_{\YMHFsmooth} + [\YMHFrough, F^{12}_{\YMHFsmooth}] \bigr) + \covd_B \Ad_g \bigl( \curl_{\YMHFsmooth} \YMHFrough + [\YMHFrough_1, \YMHFrough_2]\bigr) \;.
\end{equs}
	Then it suffices to show that for any $x \in \T^2$,
	\begin{equs}
	 \label{e:vh_mod_2_1}
		 \Big| &\int_0^{\tstop} \int_{\T^2} \fp_{B}(s, y; t, x)  [\Ad_g\covd_{\YMHFsmooth}  F^{12}_{\YMHFsmooth}, u_{\h}  ](s, y) \ds \dy \Big|  \\
		& \leq C \gronwall_{\tstop}(C_0 F_B)  \bigl(1 + \Fnorm{F_B}_{\tstop, 0}\bigr) \cdot \tstop^{\frac23} (t - \tstop)^{-\frac16} \| \fp^{1/2}(\place, \place; t, x) \covd_{\YMHFsmooth} F_{\YMHFsmooth} \|_{L^2_s L^2_y([0, \tstop])} \;,\\
		\label{e:vh_mod_2_2}
		  \Big| &\int_0^{\tstop} \int_{\T^2} \fp_{B}(s, y; t, x)  [\Ad_g [\YMHFrough, F^{12}_{\YMHFsmooth} ], u_{\h}  ] \ds \dy \Big|  \\
		& \leq C \gronwall_{\tstop}(C_0 F_B)  \bigl(1 + \Fnorm{F_B}_{\tstop, 0}\bigr) \cdot \tstop^{\frac18} (t - \tstop)^{ - \frac18 - \kappa} \| s^{\frac{\kappa}{2} } \YMHFrough \|_{L^\infty_{s, y}([0, \tstop])} \Fnorm{F_{\YMHFsmooth}}_{\tstop, 0} \;, \\
		\label{e:vh_mod_2_3}
		\Big|& \int_0^{\tstop} \int_{\T^2} \fp_{B}(s, y; t, x) [\covd_{B} \Ad_g  (\curl_{\YMHFsmooth} \YMHFrough + [\YMHFrough_1, \YMHFrough_2]), u_{\h}  ] \ds \dy \Big|\\
		 &\leq C \gronwall_{\tstop}(C_0 F_B)  \bigl(1 + \Fnorm{F_B}_{\tstop, 1}\bigr)^2 \tstop^{\frac18} (t - \tstop)^{ - \frac{3}{16} - \kappa} t^{\frac{1}{16}} \| s^{\frac{1 + \kappa}{2} } \bigl( \curl_{\YMHFsmooth} \YMHFrough + [\YMHFrough_1, \YMHFrough_2] \bigr) \|_{L^\infty_{s, y}} \;.
	\end{equs}
	The proof of \eqref{e:vh_mod_2_1} is the same as in the proof of Lemma~\ref{le:curl_square_expectation_smoothB} since $|\Ad_g \covd_{\YMHFsmooth} F^{12}_{\YMHFsmooth}| = |\covd_{\YMHFsmooth} F^{12}_{\YMHFsmooth}|$. For \eqref{e:vh_mod_2_2}, by H\"older's inequality, we get that the left-hand side is bounded by
\begin{equs}
	 \| \fp(\place, \place; t, x) s^{- \frac{\kappa}{2} - \frac14} \|_{L^2_s L^{4/3}_y} \| s^{\frac{\kappa}{2} } \YMHFrough \|_{L^\infty_{s, y}} \|s^{\frac12} F_{\YMHFsmooth} \|_{L^2_s L^\infty_y} \|s^{-\frac14} u_{\h} \|_{L^\infty_s L^4_y} \;.
\end{equs}
The bound \eqref{e:vh_mod_2_2} follows by invoking Lemma~\ref{le:uh_bound}. For \eqref{e:vh_mod_2_3}, by integration by parts (specifically, the identity \eqref{eq:dagger-fp}) and H\"older's inequality, the left-hand side is bounded by 
	\begin{equs}
		\,& \|  s^{-\frac14 -\frac{\kappa}{2}} \covd_{B}^{\dagger} \fp_{B}(\place, \place; t, x) \|_{L^1_s L^{4/3}_y}   \| s^{-\frac14} u_{\h} \|_{L^\infty_s L^4_y}  \| s^{\frac12 + \frac{\kappa}{2}} \bigl( \curl_{\YMHFsmooth} \YMHFrough + [\YMHFrough_1, \YMHFrough_2] \bigr) \|_{L^\infty_{s, y}} \label{e:vh_mod_2_31} \\
		&+ \| s^{-\frac38 -\frac{\kappa}{2}} (t - s)^{-\frac14} \fp^{1/2}(\place, \place; t, x) \|_{L^2_{s, y}} \|  s^{-\frac18} (t - s)^{\frac14} \fp(\place, \place; t, x)^{1/2} \covd_{B} u_{\h} \|_{L^2_{s, y}}  \\
		& \cdot \|  s^{\frac{1+\kappa}{2}} \bigl( \curl_{\YMHFsmooth} \YMHFrough + [\YMHFrough_1, \YMHFrough_2] \bigr) \|_{L^\infty_{s, y}} \;. \label{e:vh_mod_2_32}
	\end{equs}
	Applying Proposition~\ref{prop:qB-mixed-derivatives-ptwise-bounds} and Lemma~\ref{le:uh_bound}, we get that \eqref{e:vh_mod_2_31} is bounded by 
	\begin{equs}
		\bigl(1 &+ \Fnorm{F_B}_{\tstop, 1}\bigr) \tstop^{\frac18}(t - \tstop)^{-\frac18 -\frac{\kappa}{2}} \cdot \gronwall_{\tstop}(C_0 F_B)  \bigl(1 + \Fnorm{F_B}_{\tstop, 0}\bigr) \\
		&\cdot \| s^{\frac12 + \frac{\kappa}{2}} \bigl( \curl_{\YMHFsmooth} \YMHFrough + [\YMHFrough_1, \YMHFrough_2] \bigr) \|_{L^\infty_{s, y}} \;,
 	\end{equs}
	and \eqref{e:vh_mod_2_32} is bounded by 
	\begin{equ}
		\tstop^{\frac{1}{16}}(t - \tstop)^{-\frac{3}{16} - \frac{\kappa}{2}} \cdot \gronwall_{\tstop}(C_0 F_B)  \bigl(1 + \Fnorm{F_B}_{\tstop, 0}\bigr)^{3/2} \tstop^{\frac{1}{16}} t^{\frac{1}{16}} \cdot \| s^{\frac12 + \frac{\kappa}{2}} \bigl( \curl_{\YMHFsmooth} \YMHFrough + [\YMHFrough_1, \YMHFrough_2] \bigr) \|_{L^\infty_{s, y}} \;.
	\end{equ}
    Combining, we thus obtain \eqref{e:vh_mod_2_3}.
\end{proof}

\subsection{Decay mechanism}\label{section:refined_decay_mechanism}
This section is devoted to proving Proposition~\ref{main:prop-YM-connection-decay-mechanism}.
To this end, we will first start with the following result, which considers semisimple Lie algebras $\fg$ (so that $|\ad(F)|_{\op} \sim |F|$) under the slightly stronger assumptions \eqref{e:decay_mechanism_YM_closedness_assumption} and \eqref{e:decay_mechanism_curvature_assumption}. Note that \eqref{e:decay_mechanism_YM_closedness_assumption} is stronger than \eqref{e:close_to_YM_average} due to the $\abs{\log T}^{1/2}$ factor in the denominator, and \eqref{e:decay_mechanism_curvature_assumption} is stronger than Assumption~\ref{local:assumption-B} due to the $\log \abs{\log T}$ factor. 
We will later improve this result to Proposition~\ref{main:prop-YM-connection-decay-mechanism} by first finding some interval $[t_0, t_1] \subset [0, T]$ such that the stronger assumption holds, then proving a monotonicity result that allows us to extend the endpoint from $t_1$ to $T$, and finally using the results from Lemma~\ref{le:curl_square_expectation_smoothB} to show that the contributions from the other parts are perturbative. 

To treat general Lie algebras $\fg$, we will break $\fg$ into the semisimple part and the central part. The previous argument already treats the semisimple part. For the central part, we use that the curvature of this part satisfies the usual linear heat equation, and thus it will be relatively easy to control. The proof of Proposition~\ref{main:prop-YM-connection-decay-mechanism} is presented at the end of this section. 

\begin{proposition}
\label{pr:decay_mechanism_expectation}
	For any constant $C_1$, and small enough constant $c_\star$ (depending only on $C_0$ in Proposition~\ref{prop:qB-perturb-near-YM}), there exists some constant $C$ such that the following holds for any sufficiently small $T$.  
	Let $B$ solve the YM gradient flow on $[0, T]$ (in particular, without any jumps). Suppose $\tstop \leq T/2$, and that we have
	\begin{equ}
	\label{e:decay_mechanism_YM_closedness_assumption}
		\YMFnorm{B}_{T} \leq C_1 \frac{\bigl(\log \abs{\log T} \bigr)^{100}}{\abs{\log T}^{1/2} \cdot T } \;, 
	\end{equ}
	where $\YMFnorm{\place}$ is defined in \eqref{e:YMFnorm}, and for any $t \in [0, T], x \in \T^2$, we have 
	\begin{equ}
	\label{e:decay_mechanism_curvature_assumption}
		|F_B(t, x)| \leq c_{\star} \frac{\log \abs{\log T}}{T} \;.
	\end{equ}
	Then we have
	\begin{equs}
		\,&\E \int_{\T^2} \Big\langle F^{12}_{B(T)} , \Big[\bigl(e^{(T - \tstop) \Delta^F_{B(\tstop)} } \linear[B](\tstop) \bigr)^1 , \bigl(e^{(T - \tstop) \Delta^F_{B(\tstop)} } \linear[B](\tstop) \bigr)^2 \Big] \Big\rangle_{\fg} \, \dx \\
		&\leq - \frac{1}{4\pi} \int_0^{\tstop} \int_{\T^2} |\ad(F_B^{12}(s, y))|_{\op}^2 e^{ 2 |\ad(F_B^{12}(s, y))|_{\op}(T - s)} \ds \dy + C \frac{\bigl(\log \abs{\log T}\bigr)^2}{T} \;,
	\end{equs}
	where $|\ad(F_B^{12}(s, y))|_{\op}$ denotes the operator norm of the operator $\ad(F_B^{12}(s, y)) \in \End(\fg)$. 
\end{proposition}
\begin{proof}
	For any $\fg$-valued one-form $C$, we can rewrite 
	\begin{equ}
		\langle F_B^{12}, [C_1, C_2] \rangle_{\frkg} = -(1/2) \langle C, F_B C \rangle_{\frkg^{2}} \;,
	\end{equ}
	where we keep following the convention $(F_B C)_i = (F_B)_i^{~j} C_j$ that is used throughout the paper. 
	Let  $\{e_{\alpha}\}$ be any orthonormal basis for $\fg^2$. Recall the spatial cutoff function $\ptcutoff$ from Definition~\ref{def:ptcutoff}. By \eqref{eq:qB-ptwise-bounds}, we have that for $0 \leq s < t \leq T$, $x, y \in \T^2$,
	\begin{equ}
	\abs{(1 - \ptcutoff(x-y))\fq_B(s, y; t, x)} \lesssim \gronwall_T(C_0 F_B) \exp(-c / (t-s)) \leq \exp(- c / (t-s)),
	\end{equ} 
	where the constant $c$ may have changed in the last inequality (and we used \eqref{e:decay_mechanism_curvature_assumption} to estimate the Gronwall factor). We may thus introduce $\ptcutoff$ at the expense of a tiny error, allowing us to write the expectation as
	\begin{equs}
		\,& - \frac{1}{2} \int_{\T^2 \times \T^2} \int_0^{\tstop} \Big\langle \fq_B(s, y; T, x) e_{\alpha}, F_B(T, x) \fq_B(s, y; T, x) e_{\beta} \Big\rangle_{\frkg^{2}} \delta^{\alpha \beta} \ds \dy \dx \\
	\label{e:decay_mechanism_expectation_mid0}
		&= - \frac{1}{2} \int_{\T^2 \times \T^2} \int_0^{\tstop} \Big\langle \bigl(\ptcutoff \tau_B^{-1}\fq_B\bigr)(s, y; T, x)e_{\alpha}, F_B(s, y)\bigl(\ptcutoff \tau_B^{-1}\fq_B\bigr)(s, y; T, x)e_{\beta} \Big\rangle_{\frkg^{2}} \delta^{\alpha \beta} \ds \dy \dx \qquad \quad  \\
		&\quad - \frac{1}{2} \int_{\T^2 \times \T^2} \int_0^{\tstop} \Big\langle \bigl(\ptcutoff \tau_B^{-1} \fq_B\bigr) (s, y; T, x)e_{\alpha}, \ptcutoff\bigl(\Ad_{\tau_B^{-1}(s, y; T, x) } F_B(T, x) - F_B(s, y)\bigr) \\
	\label{e:decay_mechanism_expectation_mid_commutator}
		&\qquad\qquad\qquad\qquad\qquad\qquad\qquad\qquad\qquad\qquad \cdot \bigl(\tau_B^{-1}\fq_B\bigr)(s, y; T, x)e_{\beta} \Big\rangle_{\frkg^{2}} \delta^{\alpha \beta} \ds \dy \dx \\
		&\quad + \bigO(e^{-c / T}) \;.
	\end{equs}
	 Note that by \eqref{eq:curvature-pt-commutation-clean}, we have 
            	\begin{equ}
            	\label{e:decay_mechanism_commutator}
            		\ptcutoff |\Ad_{\pt_B^{-1}(s,y; T, x)} F_B(T, x) - F_B(s, y) | \lesssim \YMFnorm{B}_{[s, T]} \bigl(1 + \frac{|x - y|}{(T - s)^{\frac12}} \bigr) \;. 
            	\end{equ}
		Combining with \eqref{eq:qB-ptwise-bounds} and our assumptions, we obtain that the absolute value of \eqref{e:decay_mechanism_expectation_mid_commutator} is bounded by
		\begin{equs}
			\, & \int_{\T^2 \times \T^2} \int_0^{\tstop} \frac{\big(\log \abs{\log T}\big)^{100}}{\abs{\log T}^{1/2} \cdot T}  \Bigl(1 + \frac{\abs{x - y}}{(T - s)^{\frac12}}\Bigr) \\
			& \qquad \qquad \qquad \qquad \qquad \qquad \cdot \gronwall_T(2 C_0 F_B) \Big( \frac{1}{T-s} \exp(-c \frac{\abs{x - y}^2}{T - s}) \Big)^2 \ds \dy \dx \\
			&\lesssim \frac{\big(\log \abs{\log T}\big)^{100}}{\abs{\log T}^{1/2} \cdot T} e^{2 C_0 c_{\star} \log \abs{\log T} } \log \Bigl(\frac{T}{T - \tstop} \Bigr) \;. 
		\end{equs}
	Therefore, if $c_{\star}$ is sufficiently small, this term is much smaller than $T^{-1}$ (recall that $\tstop \leq T/2$ by assumption).

		We now focus on \eqref{e:decay_mechanism_expectation_mid0}. 
	From our assumptions \eqref{e:decay_mechanism_YM_closedness_assumption}--\eqref{e:decay_mechanism_curvature_assumption}, we have that
	\begin{equ}
		c_0 = T \bigl(1 + T \bigl\|F_B\bigr\|_{L_{t,x}^\infty([0, T])} \bigr)^{5} \YMFnorm{B}_{T} \bigl(1 + T \YMFnorm{B}_{T}\bigr)^5 \lesssim \bigl(\log \abs{\log T} \bigr)^5 \cdot \frac{\bigl(\log \abs{\log T}\bigr)^{100}}{\abs{\log T}^{1/2}} \ll 1
	\end{equ} 
	when $T$ is small enough. By Proposition~\ref{prop:qB-perturb-near-YM}, if we replace $\bigl(\ptcutoff \tau_B^{-1} \fq_B\bigr)(s, y; T, x)$ in \eqref{e:decay_mechanism_expectation_mid0} by $\bigl( \ptcutoff \fs^F_{F_B^{12}(s, y)} \bigr)(s, y; T, x)$, the error is bounded by
	\begin{equs}
		\,& c_0 \gronwall_T(2 C_0 F_B) \int_{\T^2 \times \T^2} \int_0^{\tstop} |F_B(s, y)| \frac{1}{(T - s)^2 } \exp\Bigl( - \frac{2|x - y|^2}{C_0(T - s)} \Bigr) \ds \dy \dx \\ 
		&\lesssim \frac{\bigl(\log \abs{\log T} \bigr)^{105}}{\abs{\log T}^{1/2}} \cdot e^{2C_0 c_{\star} \log \abs{\log T}} \cdot \frac{\log \abs{\log T}}{T} \log\Bigl(\frac{T}{T-\tstop}\Bigr) \;. 
	\end{equs}
	Therefore, if $c_{\star}$ is sufficiently small, this term is much smaller than $T^{-1}$. Thus, it suffices to look at 
	\begin{equ}
	\label{e:decay_mechanism_expectation_mid}
		- \frac{1}{2} \int_{\T^2 \times \T^2} \int_0^{\tstop} \Big\langle \bigl(\ptcutoff \fs^F_{F_B^{12}(s, y)}\bigr)(s, y; T, x)e_{\alpha}, F_B(s, y)\bigl(\ptcutoff \fs^F_{F_B^{12}(s, y)}\bigr)(s, y; T, x)e_{\beta} \Big\rangle_{\frkg^{2}} \delta^{\alpha \beta} \ds \dy \dx \;.
	\end{equ}
	
	We want to show that \eqref{e:decay_mechanism_expectation_mid} is very negative.
	To this end, we would like to first fix $s, y$ and integrate out $x$. In the following, it may help to recall the linear algebra discussion just after Definition~\ref{defn:s_B_eq}. Let $\pm \lambda_1, \cdots, \pm \lambda_k$ be the nonzero eigenvalues of $\icomplex \ad(F_B^{12}(s, y))$ and we take the convention that $\lambda_i > 0$. Here we don't write out the dependence of $\lambda_i$ on $s, y$ for notational simplicity. Note that we have 
	\begin{equ}
	\label{e:lambda_i_F_relation}
		\max_{1 \leq i \leq k} \lambda_i = |\ad (F_B^{12}(s, y))|_{\op} \lesssim c_{\star} \frac{\log \abs{\log T}}{T} \;.
	\end{equ} 
	Let $\{u_i , v_i\}_{1 \leq i \leq k} \in \fg$ be an orthonormal set such that 
	\begin{equ}
		F_B^{12}(s, y) u_i = -\lambda_i v_i \;, \quad F_B^{12}(s, y) v_i = \lambda_i u_i \;.
	\end{equ}
	Then, by the discussions after Definition~\ref{defn:s_B_eq}, we know that $\frac{1}{\sqrt{2}} (u_i, v_i),\frac{1}{\sqrt{2}} (v_i, -u_i)$ (with eigenvalue $\lambda_i$) and $\frac{1}{\sqrt{2}}(v_i, u_i), \frac{1}{\sqrt{2}}(u_i, -v_i)$ (with eigenvalue $-\lambda_i$) make up all the normalised eigenvectors of $F_B \in \End(\frkg^{2})$ which have nonzero eigenvalues. 
	Note that the choice of orthonormal basis $e_\alpha$ appearing in \eqref{e:decay_mechanism_expectation_mid} is arbitrary and is in particular allowed to
	depend on the integration variables $s, y$. In the following discussions, it will
	be convenient to choose the $e_{\alpha}$ to be the eigenvectors of the operator $F_B(s, y)$. Now we discuss the following three cases. 
	
	\proofstep{Case 1:} $e_\alpha = e_{\beta}$ is of the form $\frac{1}{\sqrt{2}}(u_i, v_i), \frac{1}{\sqrt{2}}(v_i, -u_i)$. 
	In this case, the integrand in \eqref{e:decay_mechanism_expectation_mid} is equal to 
		\begin{equs}
			\,& \frac12 \lambda_i \bigl(\ptcutoff(x - y) e^{2 \lambda_i (T - s)} \mehler_{\lambda_i}(T-s, x-y) \bigr)^2 \\
		\label{e:decay_mechanism_expectation_mid_case1}
			&= \frac12 \lambda_i e^{4 \lambda_i (T - s)} \ptcutoff(x-y)^2 \Big(\frac{\lambda_i}{4\pi \sinh(\lambda_i (T - s))} \Big)^2 \exp\bigg(-\frac{\lambda_i |x - y|^2}{2\tanh(\lambda_i (T - s))}\bigg) \;, 
		\end{equs}
		where $\mehler_{\lambda_i}$ is the Mehler kernel defined in \eqref{eq:mehler-kernel-formula}.
		Using the fact that $\tanh(x) \lesssim x \lesssim \sinh(x)$ for $x > 0$, we get that for $s \in (0, \tstop) \subset (0, T/2 )$, 
		\begin{equs}
			\Big(\frac{\lambda_i}{4\pi \sinh(\lambda_i (T - s))} \Big)^2 \cdot  \exp\bigg(-\frac{\lambda_i |x - y|^2}{2\tanh(\lambda_i (T - s))}\bigg) \Bigl(1 - \ptcutoff(x - y)^2 \Bigr) \leq e^{-\frac{c}{T}} 
		\end{equs}
		for some constant $c$. 
		Therefore, the $\ptcutoff^2$ term in \eqref{e:decay_mechanism_expectation_mid_case1} can be replaced by $1$ up to an error of order $e^{-c/T}$. Integrating out $x$, for $e_\alpha \in \{(u_i, v_i), (v_i, -u_i) \}$, we then obtain
		\begin{equs}	
			\,& - \frac{1}{2} \int_{\T^2 \times \T^2} \int_0^{\tstop} \Big\langle \bigl(\ptcutoff \fs^F_{F_B^{12}(s, y)}\bigr)(s, y; T, x)e_{\alpha}, F_B(s, y)\bigl(\ptcutoff \fs^F_{F_B^{12}(s, y)}\bigr)(s, y; T, x)e_{\alpha} \Big\rangle_{\frkg^{2}} \ds \dy \dx \\
			&= - \frac12 \int_0^{\tstop} \int_{\T^2} \lambda_i e^{4\lambda_i(T-s)} \Bigl(\frac{\lambda_i}{4\pi \sinh(\lambda_i(T-s))}\Bigr)^2 \frac{2\pi \tanh(\lambda_i(T-s))}{\lambda_i} \ds \dy + \bigO(e^{-c/T}) \\
			&= - \frac12 \int_0^{\tstop} \int_{\T^2} \lambda_i e^{4\lambda_i (T-s)} \cdot \frac{\lambda_i}{4\pi \sinh(2\lambda_i (T-s))} \ds \dy + \bigO(e^{-c/T}) \\
			&\leq - \frac{1}{4\pi} \int_0^{\tstop} \int_{\T^2} \lambda^2_i e^{2\lambda_i (T-s)} \ds \dy + \bigO(e^{-c/T}) \;,
		\end{equs} 
		where we used $\sinh(x) \cosh (x) = (1/2) \sinh(2x)$ from the second line to the third line. Therefore, using \eqref{e:lambda_i_F_relation} and only picking out the contribution from the largest $\lambda_i$ we get 
	\begin{equs}
		\,& - \sum_{e_\alpha \in \cup_i \{(u_i, v_i), (v_i, -u_i) \}} \frac{1}{2} \int_{\T^2 \times \T^2} \int_0^{\tstop}  \\
		&\qquad \qquad \qquad \qquad \qquad 
		\Big\langle \bigl(\ptcutoff \fs^F_{F_B^{12}(s, y)}\bigr)(s, y; T, x)e_{\alpha}, F_B(s, y)\bigl(\ptcutoff \fs^F_{F_B^{12}(s, y)}\bigr)(s, y; T, x)e_{\alpha} \Big\rangle_{\frkg^{2}} \ds \dy \dx \\
		& \leq - \frac{1}{4\pi} \int_0^{\tstop} \int_{\T^2} |\ad(F_B^{12}(s, y))|_{\op}^2 e^{ 2 |\ad(F_B^{12}(s, y))|_{\op}(T - s)} \ds \dy + \bigO(e^{-c/T}) \;.
	\end{equs}

	\proofstep{Case 2:} $e_\alpha = e_\beta$ is of the form $\frac{1}{\sqrt{2}}(v_i, u_i), \frac{1}{\sqrt{2}} (u_i, -v_i)$. In this case, the integrand in \eqref{e:decay_mechanism_expectation_mid} is equal to
		\begin{equ}
			- \frac12 \lambda_i \bigl(\ptcutoff(x - y) e^{-2 \lambda_i (T - s)} \mehler_{\lambda_i}(T-s, x-y) \bigr)^2 \;.
		\end{equ}
		After integrating out $x$, the contribution is bounded in absolute value by  
		\begin{equs}
			\frac12 \int_0^{\tstop} \int_{\T^2} \lambda_i e^{-4\lambda_i (T-s)} \cdot \frac{\lambda_i}{4\pi \sinh(2\lambda_i (T-s))} \ds \dy &\leq
			\frac12 \int_0^{\tstop} \int_{\T^2} \frac{\lambda_i^2 }{4\pi \sinh(2\lambda_i (T-s))} \ds \dy  \\
			&\lesssim \int_0^{\tstop} \int_{\T^2} \Bigl(\lambda_i^2 + \frac{\lambda_i}{T - s}\Bigr) \ds \dy \;,
		\end{equs}
		where we used $\frac{1}{\sinh(x)} \lesssim \frac{1}{x \wedge 1} \leq 1 + \frac{1}{x}$. With \eqref{e:lambda_i_F_relation}, this is further bounded by $\bigl(\log \abs{\log T}\bigr)^2/T$. 
		
	\proofstep{Case 3:} $e_{\alpha} = e_\beta$ is in the kernel of $F_B(s, y)$. In this case, the integrand in \eqref{e:decay_mechanism_expectation_mid} is $0$. 
	
	\proofstep{}The proposition now follows by combining the discussions above.
\end{proof}

Now we select the time interval $[t_0, t_1]$ as promised at the beginning of this subsection, in which $B$ is close enough to some YM connection. We emphasise that $t_1 - t_0$ has to be chosen long enough so that $\|F_B \|_{L^1_t L^\infty_x([t_0, t_1])} \gg 1 $. Otherwise the amplitude in the unstable directions in the linear stochastic object will not be enough to decrease the energy. This is reflected in the $\log \abs{\log T}$ factor in \eqref{e:exists_YM_close_interval_a0}. 
\begin{lemma}
\label{le:exists_YM_close_interval}
	For any constants $C_1, C_2$ depending only on the constants in \eqref{prelim:eq-parameters}, there are constants $C, a_0$ such that the following holds. Suppose $B \in C_{t,x}^\infty((0, T] \times \T^2)$ is a YM gradient flow on $(0, T]$ such that Assumption~\ref{local:assumption-B} and \eqref{e:close_to_YM_average} hold. 
	If $T$ is sufficiently small depending only on $C_1$, $C_2$, and the constants in \eqref{prelim:eq-parameters}, there exist some $ t_{0} < t_{1/2} < t_1 \in [T/8, 3T/8]$ such that
	\begin{equ}
	\label{e:exists_YM_close_interval_a0}
		t_{1/2} - t_0 = a_0 \frac{T}{\abs{\log T}} \;, \qquad t_1 - t_0 = a_0 \frac{T \log \abs{\log T}}{\abs{\log T}} \;, 
	\end{equ}
	and (recall \eqref{e:YMFnorm})
	\begin{equs}
	\label{e:exists_YM_close_interval1}
		\|\covd_B F_B \|_{L^\infty_t L^2_x([t_0, t_1])} &\leq C \frac{\abs{\log T}^{\frac{1}{100}} }{T} \;, \\
	\label{e:exists_YM_close_interval2}
		(t_1 - t_0)\YMFnorm{B}_{[t_0, t_1]} &\leq C \frac{\bigl(\log \abs{\log T} \bigr)^{100}}{\abs{\log T}^{1/2}} \;.
	\end{equs}
	Furthermore, if \eqref{e:decay_mechanism_F_lower_bound} holds, then for any $t \in [t_0, t_1]$ and $x \in \T^2$ we have 
	\begin{equ}
	\label{e:YM_close_interval_pointwise_lower_bound} 
		|F_B(t, x)| \geq C_2 \frac{\abs{\log T}}{T} - C \frac{\abs{\log T}^{1/2}}{T} \;. 
	\end{equ}
\end{lemma}
\begin{proof}
	By \eqref{e:close_to_YM_average}, there exists some time $t_{\min} \in [\frac18 T, \frac14 T ]$ such that 
	\begin{equ}
		\| \covd_B F_B (t_{\min}) \|_{L^2} \lesssim \frac1T \;.
	\end{equ}
	By Remark~\ref{remark:curvature-pointwise-in-time-bound-Fnorm}, for any $r \in (0, T]$, we have that 
	\begin{equ}
	\label{e:remark-curvature-pointwise-in-time-bound-Fnorm}
		\|F_B(r) \|_{L^\infty} \leq \frac{1}{r} \Fnorm{F_B}_{[0, T]} \;. 
	\end{equ}
	By Lemma~\ref{lemma:strichartz-estimate-using-decrease} and \eqref{e:close_to_YM_average}, for any $r_0, r_1 \in [\frac18 T, T]$ such that $0 \leq r_1 - r_0 \leq \frac{T}{\abs{\log T}}$, we have
	\begin{equs}
		\Fnorm{F_B}_{[r_0, r_1], 0} &\lesssim (r_1 - r_0)^{1/2} \bigl\|F_B \bigr\|_{L_t^2L_x^\infty} \\ 
		&\lesssim_{C_1} (r_1 - r_0) \|F_{B(r_0)} \|_{L^2} + (r_1 - r_0)^{\frac12} \abs{\log (r_1 - r_0)}^{\frac12} T^{-\frac12}  \\
	\label{e:exists_YM_close_interval_mid}
		&\lesssim (r_1 - r_0) r_0^{-1} \Fnorm{F_B}_{(0, r_0], 0} + 1 \lesssim 1 \;.
	\end{equs}
	Let 
	\begin{equ}
		t_0 = t_{\min} + \frac{T}{\abs{\log T}}\, ,  \quad 
		t_{1/2} = t_0 +  a_0 \frac{T}{\abs{\log T}} \;,
		   \quad t_1 = t_{0} + a_0 \frac{T\log \abs{\log T} }{\abs{\log T}} \;,
	\end{equ}
	where $a_0$ is some small constant to be chosen later.
	Similar to \eqref{e:exists_YM_close_interval_mid}, we have $\Fnorm{F_B}_{[t_{\min}, t_1], 0} \lesssim a_0 \log \abs{\log T}$. 
	Note that by Lemma~\ref{lemma:curvature-equation} and \eqref{eq:cov-heat-operator-covd-commutator}, $\covd_B F_B^{12}$ solves
	\begin{equ}
	\label{e:covd_B_F_B_evolution}
		\ptl_t \covd_B F_B^{12} = \Delta_B \covd_B F_B^{12} + 2 \ad(F_B) \covd_B F_B^{12} \;.
	\end{equ}
	
	Therefore, from the equation of $\covd_B F_B$, \eqref{eq:qB-ptwise-bounds} and Corollary~\ref{cor:curvature-higher-regularity-estimates}, for any $1 \leq k \leq 3$ we have
	\begin{equs}
		\|\covd_B F_B(t) \|_{L^\infty_t L_x^2([t_0, t_1]) } &\lesssim e^{C a_0 \log \abs{\log T}} \| \covd_B F_B(t_{\min}) \|_{L^2} \;, \\
		\|\covd^{(k)}_B F_B(t) \|_{L_{t,x}^\infty([t_0, t_1]) } &\lesssim \Bigl( \frac{T}{\abs{\log T}} \Bigr)^{-\frac{k}{2}} e^{C a_0 \log \abs{\log T}} \bigl(\log \abs{\log T} \bigr)^{2 \cdot 3!} \| \covd_B F_B(t_{\min}) \|_{L^2} \;.
	\end{equs}
	Choosing $a_0$ small enough, we get the bound \eqref{e:exists_YM_close_interval1}, as well as the following estimate for any $1 \leq k \leq 3$,
	\begin{equ}
		(t_1 - t_0)^{\frac{k}{2} + 1} \| \covd^{(k)}_B F_B \|_{L_{t, x}^\infty([t_0, t_1])} \lesssim \frac{\bigl(\log \abs{\log T} \bigr)^{1 + \frac{k}{2}}}{\abs{\log T}} \cdot \bigl(\abs{\log T}\bigr)^{\frac{1}{100} } \bigl(\log \abs{\log T} \bigr)^{2 \cdot 3!} \;.
	\end{equ}
	Furthermore, from the equation of $\covd_B F_B$ we have
	\begin{equs}
		\,&(t_1 - t_0)^{\frac52} \| \partial_t \covd_B F_B \|_{L^\infty_{t, x}([t_0, t_1])} \\
		&\leq (t_1 - t_0)^{\frac52} \bigl(\| \covd_B^{(3)} F_B\|_{L^\infty_{t, x}([t_0, t_1])} + \| \covd_B F_B \|_{L^\infty_{t, x}([t_0, t_1])} \| F_B \|_{L^\infty_{t, x}([t_0, t_1])} \bigr) \lesssim \frac{\bigl( \log \abs{\log T} \bigr)^{20}}{\abs{\log T}^{1/2}} \;. 
	\end{equs}
	The bound \eqref{e:exists_YM_close_interval2} follows. Now we turn to \eqref{e:YM_close_interval_pointwise_lower_bound}. First, we observe that $|F_B|^2$ solves 
	\begin{equ}
		\partial_t |F_B|^2 = \Delta | F_B|^2 - 2|\covd_B F_B|^2 \;.
	\end{equ}
	Therefore, for any $t > t_{\star}$ ($t_{\star}$ to be determined later), we can write
	\begin{equs}
		|F_B(t)|^2 &= \|F_B(t_{\star})\|_{L^2}^2 \\
		\label{e:F_B_YM_fluctuation}
				&\quad + e^{(t - t_{\star})\Delta } \bigl(|F_B(t_{\star})|^2 - \|F_B(t_{\star})\|_{L^2}^2 \bigr) - 2\int_{t_{\star}}^t e^{(t - s)\Delta} \bigl(|\covd_B F_B(s)|^2 \bigr) \ds \;.
	\end{equs} 
	Using the assumptions \eqref{e:decay_mechanism_F_lower_bound} to lower bound the first line, it remains to bound \eqref{e:F_B_YM_fluctuation} in order to obtain \eqref{e:YM_close_interval_pointwise_lower_bound}. For the first term in \eqref{e:F_B_YM_fluctuation}, let $P_K$ be the Littlewood--Paley projection onto frequencies of order $K$ for $K \in \dyadic$ as in \eqref{main:eq-standard-Littlewood-Paley} (see also \cite[Section~2]{BCD2011}). Combining the heat flow estimate, Bernstein's inequality and the diamagnetic inequality \eqref{eq:gradient-diamagnetic}, it follows that for any $p \geq 2$, the $L^\infty$ norm of the first term is bounded by
	\begin{equs}
		C &(t - t_{\star})^{-\frac1p} \| \sum_{K \in \dyadic} P_K \bigl(|F_B(t_{\star})|^2 \bigr) \|_{L^p} 
		\lesssim (t - t_{\star})^{-\frac1p} \sum_{K \in \dyadic} K^{- \frac{2}{p}} \| \nabla P_K \bigl(|F_B(t_{\star})|^2 \bigr) \|_{L^2} \\
		&\lesssim (t - t_{\star})^{-\frac1p} \sqrt{p} \| \nabla \bigl( |F_B(t_{\star})|^2 \bigr) \|_{L^2} 
		\label{e:F_B_YM_fluctuation_mid1}
		\leq  (t - t_{\star})^{-\frac1p} \sqrt{p} \| \covd_B F_B(t_{\star}) \|_{L^2} \| F_B(t_{\star}) \|_{L^\infty} \;.\quad
	\end{equs}
	Now we take $t_{\star} = \frac{t_{\min} + t_0}{2}$.
	By \eqref{e:exists_YM_close_interval_mid}, the equation of $\covd_B F_B$ and \eqref{eq:qB-ptwise-bounds} we get 
	\begin{equ}
		\|\covd_B F_B(t_{\star}) \|_{L^2} \lesssim \|\covd_B F_B(t_{\min}) \|_{L^2} \lesssim \frac1T \;. 
	\end{equ}
	Taking $p$ in \eqref{e:F_B_YM_fluctuation_mid1} to be $\abs{\log T}$, and combining with the above discussions, for any $t \in (t_0, T]$ we have 
	\begin{equ}
        		\bigl\|e^{(t-t_{\star})\Delta}\bigl(|F_B(t_{\star})|^2-\|F_B(t_{\star})\|_{L^2}^2\bigr)\bigr\|_{L^\infty}
        	\lesssim (t-t_{\star})^{-\frac1{\abs{\log T}}}\abs{\log T}^{\frac12}\frac1T\frac{\abs{\log T}}T
        	\lesssim \abs{\log T}^{\frac32}T^{-2}\;.
        \end{equ}
	Now we turn to the second term in \eqref{e:F_B_YM_fluctuation}. By the heat flow estimate and the estimate \eqref{eq:qB-ptwise-bounds} for $\fq_B$, for any $t \in [t_0, t_1]$ we get
	\begin{equs}
		\Bigl\| \int_{t_{\star}}^t e^{(t - s)\Delta} \bigl(|\covd_B F_B(s)|^2 \bigr) \ds \Bigr\|_{L^\infty} &\lesssim \int_{t_{\star}}^t  \| \covd_B F_B(s) \|_{L^\infty}^2 \ds \\
		&\lesssim \gronwall_{[t_{\min},t_1]}(2 C_0 F_B) \int_{t_{\star}}^t (s - t_{\min})^{-1} \| \covd_B F_B(t_{\min}) \|_{L^2}^2 \ds \\
		&\lesssim e^{C a_0 \log \abs{\log T}} \frac{t_1 - t_{\star}}{t_\star - t_{\min}} \frac{1}{T^2} \lesssim \frac{\abs{\log T}^{1/2}}{T^2} \;,
	\end{equs}
	if we choose $a_0$ small enough. 
	
	From the above discussions, we get that for any $t \in [t_0, t_1]$, the $L^\infty$ norm of \eqref{e:F_B_YM_fluctuation} is bounded by $\abs{\log T}^{3/2}/T^2$, and the bound \eqref{e:YM_close_interval_pointwise_lower_bound} follows. 
\end{proof}
\begin{remark}
	If we apply Lemma~\ref{le:DF_bound} to estimate $\| \covd_B F_B (t_1) \|_{L^2}$ for $t_1 \in [T/8, 3T/8]$, we lose an extra $\abs{\log T}$ coming from the $ \|t^{1/2} F_B\|_{L_t^2 L_x^\infty([0, T])}$ term. The resulting loss is too much to prove that $(t_1 - t_0)\YMFnorm{B}_{[t_0, t_1]}$ is small. This is the reason why we have to first pick out the time $t_{\min}$. 
\end{remark}
\begin{remark}
	Actually, the pointwise lower bound for $|F_B(t, x)|$ in \eqref{e:YM_close_interval_pointwise_lower_bound} is not necessary in the ensuing argument (a lower bound in $L^2$ would suffice). However, we still establish the stronger pointwise estimate because it makes the structure of $|F_B|$ more transparent, namely as a large constant field with lower-order perturbations. 
\end{remark}

Now we prove the following monotonicity result. 
\begin{lemma}
\label{le:decrease_expectataion_mono}
	There exists some constant $C$ such that the following holds. 
	Let $B$ solve the YM heat flow on $[0, T]$ such that Assumption~\ref{local:assumption-B} and \eqref{e:close_to_YM_average} hold. 
	For any $t_{1/2} \in (0, 2T/3]$ and $t_1 \in (t_{1/2} \vee \frac{T}{8}, 2T/3] $, we have 
	\begin{equs}
		\,& \E \int_{\T^2} \Big \langle e^{(T - t_{1/2}) \Delta^F_{B(t_{1/2})} } \linear[B](t_{1/2}) , \bigl(-\Delta_{B(T)} - 2\ad(F_{B(T)})\bigr) e^{(T - t_{1/2}) \Delta^F_{B(t_{1/2})} } \linear[B](t_{1/2}) \Big \rangle_{\fg } \, \dx \\ 
		& \leq  \E \int_{\T^2}  \Big \langle e^{(t_1 - t_{1/2}) \Delta^F_{B(t_{1/2})} } \linear[B](t_{1/2}) , \bigl(-\Delta_{B(t_1)} - 2\ad(F_{B(t_1)})\bigr) e^{(t_1 - t_{1/2}) \Delta^F_{B(t_{1/2})} } \linear[B](t_{1/2}) \Big \rangle_{\fg } \, \dx \\
		&\qquad \qquad \qquad \qquad \qquad \qquad \qquad \qquad \qquad \qquad \qquad \qquad \qquad \qquad \qquad + C T^{- 1/2 - 2 \kappa } \;.
	\end{equs}
	
\end{lemma}

\begin{proof}
	To simplify the notation, for $t > t_{1/2}$ let 
	\begin{equ}
		v(t) = e^{(t - t_{1/2}) \Delta^F_{B(t_{1/2})} } \linear[B](t_{1/2}) \;.
	\end{equ} 
	Recall from Notation~\ref{notation:norms-forms} that for $\fg$-valued functions (or one-forms) $f, g$, we use $\langle f, g \rangle_{L^2}$ to denote $\int_{\T^2} \langle f, g \rangle_{\fg} \, \dx $.  
	Note that for $t > t_{1/2}$, 
	\begin{equs}
	\partial_t \Big\langle v(t), \bigl(-\Delta_{B(t)} - 2\ad(F_{B(t)})\bigr)v(t) \Big\rangle_{L^2}
	&= -2 \Big\|\bigl(\Delta_{B(t)} + 2\ad(F_{B(t)})\bigr)v(t)\Big\|_{L^2}^2 \\
	&\quad - 2 \Big\langle v_k, \bigl[\covd_{B,\ell}F_B^{\ell j}, \covd_{B,j}v^k\bigr]\Big\rangle_{L^2}
	- 2 \Big\langle v_k, \bigl[\Delta_B F_B^{kj}, v_j\bigr]\Big\rangle_{L^2} \\
	&\leq - 2 \Big\langle v_k, \bigl[\covd_{B,\ell}F_B^{\ell j}, \covd_{B,j}v^k\bigr]\Big\rangle_{L^2}
	- 2 \Big\langle v_k, \bigl[\Delta_B F_B^{kj}, v_j\bigr]\Big\rangle_{L^2} \;.
	\end{equs}
	Integration by parts gives
	\begin{equ}
	\Big\langle v_k,\bigl[\Delta_BF_B^{kj},v_j\bigr]\Big\rangle_{L^2} =-\Big\langle \covd_{B,\ell}v_k,
		\bigl[\covd_B^\ell F_B^{kj},v_j\bigr]\Big\rangle_{L^2} -\Big\langle v_k,
		\bigl[\covd_B^\ell F_B^{kj},\covd_{B,\ell}v_j\bigr]\Big\rangle_{L^2} \;.
	\end{equ}
	Consequently, H\"older's inequality implies that 
	\begin{equs}
		\E \Big\langle v(T), \bigl(-\Delta_{B(T)}-2\ad(F_{B(T)})\bigr)v(T) \Big\rangle_{L^2} 
		&\leq  \E \Big\langle v(t_1), \bigl(-\Delta_{B(t_1)}-2\ad(F_{B(t_1)})\bigr)v(t_1) \Big\rangle_{L^2} \\
		&+ 6 \E \int_{t_1}^{T} \|\covd_BF_B(t)\|_{L^2} \|v(t)\|_{L^4} \|\covd_{B(t)}v(t)\|_{L^4}  \dt \;.
	\end{equs}
	By Lemma~\ref{le:DF_bound}, \eqref{e:close_to_YM_average}, \eqref{e:lollipop_negative_holder2}, and the fact that $t_1 > T/8$, we get that 
	\begin{equs}
	\,& \E \int_{t_1}^{T} \|\covd_BF_B(t)\|_{L^2} \|v(t)\|_{L^4} \|\covd_{B(t)}v(t)\|_{L^4}  \dt \\
	 &\lesssim \int_{t_1}^{T} t^{-\frac12} \abs{\log T} \|\covd_B F_B \|_{L^2_{t, x}([0, T])} (t - t_{1/2})^{-\frac12 - \kappa} \cdot \E \| \linear[B](t_{1/2}) \|_{\covHolder{B(t_{1/2}), F, T-t_{1/2}}^{-\kappa}}^2 \dt \\
	&\lesssim T^{-\frac12 - \kappa} \abs{\log T} \cdot \E \| \linear[B](t_{1/2}) \|_{\covHolder{B(t_{1/2}), F, T-t_{1/2}}^{-\kappa}}^2 \lesssim T^{-\frac12 - 2\kappa} \;.
	\end{equs}
	The desired result is proved. 
\end{proof}

\begin{corollary}
\label{cor:decrease_expectataion_mono}
	There exists some constant $C$ such that the following holds for all $T$ sufficiently small. 
	Let $B$ solve the YM heat flow on $[0, T]$ such that Assumption~\ref{local:assumption-B} and \eqref{e:close_to_YM_average} hold. For any $T/8 \leq t_0 < t_0 + T/8 < t_{1/2} < t_{1/2} + T^{1 + \kappa} \leq t_1 \leq 3T/8$, we have 
	\begin{equs}
\label{e:decrease_expectataion_mono}
		\,&\E \int_{\T^2} \Big\langle F^{12}_{B(T)} , \Big[\bigl(e^{\frac{T}{2} \Delta^F_{B(T/2)} } \linear[B](T/2 ) \bigr)^1 , \bigl(e^{ \frac{T}{2} \Delta^F_{B(T/2 )} } \linear[B](T/2 ) \bigr)^2 \Big] \Big\rangle_{\fg} \, \dx \\
		&\leq \E \int_{\T^2} \Big\langle F^{12}_{B(t_1)} , \Big[\bigl(e^{(t_1 - t_{1/2}) \Delta^F_{B(t_{1/2})} } \linear[B \restr [t_0, t_{1/2} ]](t_{1/2}) \bigr)^1 , \bigl(e^{(t_1 - t_{1/2}) \Delta^F_{B(t_{1/2})} } \linear[B \restr [t_0, t_{1/2} ]](t_{1/2}) \bigr)^2 \Big] \Big\rangle_{\fg} \, \dx \\
		&+ C (t_1 - t_{1/2})^{-1} \;,
	\end{equs}
	where
	\begin{equ}
		\linear[B \restr [t_0, t_{1/2} ]](t_{1/2}) \eqdef \int_{t_0}^{t_{1/2}} \int_{\T^2} \fq_B(u, y; t_{1/2}, \place) \xi(u, y) \du \dy \;.
	\end{equ}
\end{corollary}
\begin{proof}
	For $0 \leq s_0 < s_1 < t \leq T$ and $x \in \T^2$, we use the shorthand notation
	\begin{equ}
		v(s_0, s_1; t, x) \eqdef \int_{s_0}^{s_1} \int_{\T^2} \fq_B(u, y; t, x) \xi(u, y) \du \dy \;.
	\end{equ}
	We write $v(s_0, s_1; t)$ for the corresponding spatial function. Then we have
	\begin{equ}
		v(0, T/2; T) = v(0, t_0; T) + v(t_0, t_{1/2}; T) + v(t_{1/2}, T/2; T) \;,
	\end{equ}
	with these three terms being independent and mean zero. Therefore, the left-hand side of \eqref{e:decrease_expectataion_mono} is equal to
	\begin{equs}
		\,& \E \int_{\T^2} \Big\langle F_B^{12}(T), \bigl[v^1(0, t_0; T), v^2(0, t_0; T)\bigr] \Big\rangle_{\fg} \,\dx \\
		&+ \E \int_{\T^2} \Big\langle F_B^{12}(T), \bigl[v^1(t_0, t_{1/2}; T), v^2(t_0, t_{1/2}; T)\bigr] \Big\rangle_{\fg} \,\dx \\
		&+ \E \int_{\T^2} \Big\langle F_B^{12}(T), \bigl[v^1(t_{1/2}, T/2; T), v^2(t_{1/2}, T/2; T)\bigr] \Big\rangle_{\fg} \,\dx \;.
	\end{equs}
	By \eqref{e:prelim:lem-Hessian-identity-1} and Lemma~\ref{le:curl_square_expectation_smoothB}, the first and third terms above are bounded by
	\begin{equs}
		\,& \frac12 \E \Bigl(\|\curl_{B(T)}v(0, t_0; T)\|_{L^2}^2 + \|\div_{B(T)}v(0, t_0; T)\|_{L^2}^2\Bigr) \\
		\,& + \frac12 \E \Bigl(\|\curl_{B(T)}v(t_{1/2}, T/2 ; T)\|_{L^2}^2 + \|\div_{B(T)}v(t_{1/2}, T/2; T)\|_{L^2}^2\Bigr) \lesssim T^{-1}  \;.
	\end{equs}
	where we used that $t_0, t_{1/2} \in [T/8, 3T/8]$. 

	For the middle term, 	we first apply \eqref{e:prelim:lem-Hessian-identity-2}, and then use Lemma~\ref{le:decrease_expectataion_mono} with the noise shifted to starting at $t_0$ to get that
	\begin{equs}
		\,& \E \int_{\T^2} \Big\langle F_B^{12}(T), \bigl[v^1(t_0, t_{1/2}; T), v^2(t_0, t_{1/2}; T)\bigr] \Big\rangle_{\fg} \,\dx \\ 
		&\leq \frac12 \E \int_{\T^2} \Big\langle v(t_0, t_{1/2}; T), \bigl(-\Delta_{B(T)} - 2\ad(F_{B(T)})\bigr)v(t_0, t_{1/2}; T) \Big\rangle_{\fg} \, \dx \\
		&\leq \frac12 \E \int_{\T^2} \Big\langle v(t_0, t_{1/2}; t_1), \bigl(-\Delta_{B(t_1)} - 2\ad(F_{B(t_1)})\bigr)v(t_0, t_{1/2}; t_1) \Big\rangle_{\fg} \, \dx + C T^{-\frac12 - 2\kappa} \;. 
	\end{equs}
	Combining the above discussions, we obtain
	\begin{equs}
		\,& \E \int_{\T^2} \Big\langle F_B^{12}(T), \bigl[v^1(0, T/2; T), v^2(0, T/2; T)\bigr] \Big\rangle_{\fg} \,\dx \\
		&\leq \frac12 \E \int_{\T^2} \Big\langle v(t_0, t_{1/2}; t_1), \bigl(-\Delta_{B(t_1)} - 2\ad(F_{B(t_1)})\bigr)v(t_0, t_{1/2}; t_1) \Big\rangle_{\fg} \, \dx + C T^{-1} \;.
	\end{equs}
	Applying \eqref{e:prelim:lem-Hessian-identity-2} again, together with Lemma~\ref{le:curl_square_expectation_smoothB}, Remark~\ref{remark:curl_square_expectation_smoothB}, and the fact that $t_1 - t_{1/2} \geq T^{1+ \kappa}$, we get 
	\begin{equs}
		\,& \frac12 \E \int_{\T^2} \Big\langle v(t_0, t_{1/2}; t_1), \bigl(-\Delta_{B(t_1)} - 2\ad(F_{B(t_1)})\bigr)v(t_0, t_{1/2}; t_1) \Big\rangle_{\fg} \, \dx \\
		&\leq \E \int_{\T^2} \Big\langle F_B^{12}(t_1), \bigl[v^1(t_0, t_{1/2}; t_1), v^2(t_0, t_{1/2}; t_1)\bigr] \Big\rangle_{\fg} \,\dx + C(t_1 - t_{1/2})^{-1} \;. 
	\end{equs}
	The desired result follows. 
\end{proof}
	
The following lemma rules out the case in which under the assumption of Proposition~\ref{main:prop-YM-connection-decay-mechanism}, most of $F_B$ lies in the centre of $\fg$. 
In words, the proof of the lemma says that since the central part of $F_B$ solves the usual heat equation, and is moreover mean zero, there must be linear decay for this piece. Thus if the central part is large, then we would see an amount of decay which is incompatible with the assumption \eqref{e:close_to_YM_average} in Proposition~\ref{main:prop-YM-connection-decay-mechanism} (which places an upper bound on the decay of $F_B$). 

\begin{lemma}
\label{le:semisimple_part}
	There exists some constant $c_\fg > 0$ only depending on the Lie algebra $\fg$ such that the following holds.
	Let $B \in C_{t,x}^\infty((0, T] \times \T^2)$ be a YM gradient flow on $(0, T]$ such that \eqref{e:close_to_YM_average} and \eqref{e:decay_mechanism_F_lower_bound} hold. If $T$ is sufficiently small, we have 
	\begin{equ}
	\label{e:semisimple_part}
		\inf_{t \in [\frac18 T, \frac23 T]} \| |(\ad F_B^{12})(t)|_{\op} \|_{L^2} \geq \frac{C_2}{c_{\fg}} \frac{\abs{\log T}}{T} \;. 
	\end{equ}
	Furthermore, let $[t_0, t_1] \subset  [\frac18 T, \frac23 T]$ be such that \eqref{e:exists_YM_close_interval1} holds on $[t_0, t_1]$. Then the $L^2$ norm in the previous bound can be replaced by the $L^1$ norm, i.e.\ we have 
	\begin{equ}
	\label{e:semisimple_part_L1}
		\inf_{t \in [t_0, t_1]} \| |(\ad F^{12}_B)(t)|_{\op} \|_{L^1} \geq \frac{C_2}{c_{\fg}} \frac{\abs{\log T}}{T} \;. 
	\end{equ}
\end{lemma}
\begin{proof}
	By \cite[Theorems~6.10, 6.24]{Kir08}, we know that $\fg = \fg_\centre \oplus \fg_{\ss}$, where $\fg_\centre$ is the centre of $\fg$ (i.e.\ $\fg_\centre = \{x \in \fg: \ad(x) = 0 \}$)  and $\fg_{\ss}$ is semisimple. Therefore, we can decompose $B$ into $B_{\centre} + B_{\ss}$, where $B_{\centre} \in \fg_\centre$ and $B_{\ss} \in \fg_{\ss}$. We also have the decomposition $F_B = F_{B_{\centre}} + F_{B_{\ss}}$. 
	As a result, we have $|\ad(F_B^{12})|_{\op} \sim_{\fg} |F_{B_{\ss}}|$ and to show \eqref{e:semisimple_part}, it suffices to show that
	\begin{equ}
	\label{e:centre_part}
		\sup_{t \in [\frac18 T, \frac23 T]} \|F_{B_{\centre}}(t) \|_{L^2} \leq \frac{C_2}{100} \frac{\abs{\log T}}{T} \;. 
	\end{equ}
	Note that $\covd_B F_B = \covd_{B_{\centre}} F_{B_{\centre}} + \covd_{B_{\ss}} F_{B_{\ss }}$, and also both $B_{\centre}$, $B_{\ss}$ are solutions to the YM gradient flow. Therefore, $F_{B_{\centre}}$ solves the equation
	\begin{equ}
		\partial_t F_{B_{\centre}} = \Delta F_{B_{\centre}} \;.
	\end{equ}
	Since $F^{12}_{B_{\centre} } = \curl \bigl(B_{\centre} \bigr)$ has mean zero, we know that $\| F_{B_{\centre}} \|_{L^2}$ decays exponentially fast. 
	If there exists some $t \in [\frac18 T, \frac23 T]$ such that $\| F_{B_{\centre}}(t) \|_{L^2} > \frac{C_2}{100} \cdot \frac{\abs{\log T}}{T}$, then we have
	\begin{equs}
		\| \covd_{B_{\centre}} F_{B_{\centre}} \|_{L^{2}_{t, x}([t, t + \frac14 T])}^2 &= \frac12 \Bigl(\|F_{B_{\centre}}(t) \|_{L^2}^2 - \| F_{B_{\centre}}(t + \frac14 T) \|_{L^2}^2 \Bigr) \\
		&\geq \frac12 \bigl( 1 - e^{-cT} \bigr) \|F_{B_{\centre}}(t) \|_{L^2}^2 \gtrsim \frac{\abs{\log T}^2}{T} \;,
	\end{equs}
	which contradicts \eqref{e:close_to_YM_average}. The bound \eqref{e:centre_part} then follows.
	
	For \eqref{e:semisimple_part_L1}, by the Gagliardo--Nirenberg inequality, we have
    \begin{equ}
    \|F_{B_{\ss}}\|_{L^2} \lesssim \|F_{B_{\ss}}\|_{L^1} + \|F_{B_{\ss}}\|_{L^1}^{1/2} \|\covd_{B_{\ss}} F_{B_{\ss}} \|_{L^2}^{1/2}.
    \end{equ}
    The bound \eqref{e:semisimple_part_L1} follows by combining this with \eqref{e:exists_YM_close_interval1} and \eqref{e:centre_part}.
\end{proof}

We are now ready to conclude the proof of Proposition~\ref{main:prop-YM-connection-decay-mechanism}. 
\begin{proof}[Proof of Proposition~\ref{main:prop-YM-connection-decay-mechanism}]
	Let $t_0, t_{1/2}, t_1, a_0$ be as provided by Lemma~\ref{le:exists_YM_close_interval}, with $a_0$ sufficiently small. 
	The conditions \eqref{e:exists_YM_close_interval_a0} and Assumption~\ref{local:assumption-B} give
	\begin{equ}
		(t_1 - t_0)\|F_B\|_{L^\infty_{t,x}([t_0, t_1])} \lesssim a_0 \log\abs{\log T} \leq c_\star \log\abs{\log(t_1 - t_0)} \;,
	\end{equ}
	provided $a_0$ is sufficiently small. Together with \eqref{e:exists_YM_close_interval2}, this verifies the assumptions of Proposition~\ref{pr:decay_mechanism_expectation} on $[t_0, t_1]$, with the $T, \tstop$ therein taken to be $t_1 - t_0, t_{1/2} - t_0$.

	Applying Corollary~\ref{cor:decrease_expectataion_mono} and Proposition~\ref{pr:decay_mechanism_expectation}, we obtain
	\begin{equs}
		\,&\E \int_{\T^2} \Big\langle F_B^{12}(T), \Big[\bigl(e^{(T/2)\Delta^F_{B(T/2)}}\linear[B](T/2)\bigr)^1, \bigl(e^{(T/2)\Delta^F_{B(T/2)}}\linear[B](T/2)\bigr)^2\Big]\Big\rangle_{\fg} \,\dx \\
		&\leq -\frac{1}{4\pi} \int_{t_0}^{t_{1/2}} \int_{\T^2} |\ad(F_B^{12}(s, y))|_{\op}^2 e^{2|\ad(F_B^{12}(s, y))|_{\op}(t_1 - s)} \,\ds \,\dy \\
		&\qquad + C \frac{\bigl(\log\abs{\log(t_1 - t_0)}\bigr)^2}{t_1 - t_0} + C(t_1 - t_{1/2})^{-1} \;.
	\end{equs}
	By \eqref{e:exists_YM_close_interval_a0}, we can bound the last line by $ \frac{\abs{\log T}\log\abs{\log T}}{T}$. For the first term, Jensen's inequality and \eqref{e:semisimple_part_L1} give
	\begin{equs}
		\int_{t_0}^{t_{1/2}} &\int_{\T^2} |\ad(F_B^{12}(s, y))|_{\op}^2 e^{2|\ad(F_B^{12}(s, y))|_{\op}(t_1 - s)} \,\ds \,\dy \\
		&\geq \int_{t_0}^{t_{1/2}} \bigl\||\ad(F_B^{12}(s))|_{\op}\bigr\|_{L^1}^2 e^{2(t_1 - s)\bigl\||\ad(F_B^{12}(s))|_{\op}\bigr\|_{L^1}} \,\ds \\
		&\geq (t_{1/2} - t_0)\frac{C_2^2}{c_\fg^2}\frac{\abs{\log T}^2}{T^2} \exp\Bigl(2(t_1 - t_{1/2})\frac{C_2}{c_\fg}\frac{\abs{\log T}}{T}\Bigr) \\
		&= \frac{a_0 C_2^2}{c_\fg^2}\frac{\abs{\log T}}{T} \exp\Bigl(\frac{2a_0 C_2}{c_\fg}\bigl(\log\abs{\log T} - 1\bigr)\Bigr) 
		\gtrsim \frac{\abs{\log T}^{1 + 2c}}{T} \;,
	\end{equs}
	where $c = a_0 C_2/c_\fg > 0$. 
	This proves the bound \eqref{e:decay_mechanism_interval} for sufficiently small $T$.
\end{proof}

At the end of this subsection, we state the following estimate that compares smoothened linear stochastic objects with different background connections. This will also be useful to establish the decay mechanism in the iteration scheme in Section~\ref{section:decay-estimates-in-expectation}. 
The estimate is not sharp but will be enough for later use.

\begin{lemma}
\label{le:lollipop_difference} 
	Let $\YMHFsmooth$, $B$, $g, T, K$ be as in Proposition~\ref{pr:admissible_timescale_general}, with $c_0$ chosen small enough such that Lemma~\ref{le:B-YMHFsmooth-heat-comparison} holds. Then we have 
	\begin{equs}
	\E \Bigl\| &\int_0^{T/2} \int_{\T^2}
		\Bigl( \Ad_{g(T,x)^{-1}} \fq_B(s,y;T,x) \Ad_{g(s,y)} \xi(s,y) - \fq_{\YMHFsmooth}(s,y;T,x) \xi(s,y) \Bigr) \,\dy \,\ds
		\Bigr\|_{L_x^2}^2 \\
		&\lesssim M_T(\YMHFsmooth,1)^2 \bigl(M_T(B,1) + M_T(\YMHFsmooth,1)\bigr) T^{\frac{1-\kappa}{2}} K  \;.
	\end{equs}

\end{lemma}

\begin{proof}
	Let $\{e_\alpha\}$ be an orthonormal basis of $\fg \oplus \fg$.
	We first write the expectation as
	\begin{equ}
		\sum_\alpha \int_0^{T/2} \int_{\T^2 \times \T^2} \bigl| \Ad_{g(T,x)^{-1}} \fq_B(s,y;T,x) \Ad_{g(s,y)} e_\alpha - \fq_{\YMHFsmooth}(s,y;T,x) e_\alpha \bigr|^2 \dx \dy \ds \;.
	\end{equ}
	For any fixed $s \in (0,T/2)$, $y \in \T^2$, we apply Lemma~\ref{le:B-YMHFsmooth-heat-comparison}, with the equations shifted to start at time $s$ and initial data $\psi_1(s,x) = \Ad_{g(s,x)} e_\alpha \fp_{\eps^2}(x-y)$, $\psi_2(s,x) = e_\alpha \fp_{\eps^2}(x-y)$ for $e_\alpha \in \fg \oplus \fg$ and $\eps \in (0,1)$. From the lemma applied with $p = 1$ and the parameter $\alpha = 0$, combined with the fact that $\|\psi_2\|_{L_t^\infty L_x^1([s,T])} \lesssim M_T(\YMHFsmooth,1)$ from \eqref{eq:qB-ptwise-bounds}, we get
	\begin{equs}
		\int_{\T^2} \Big| &\int_{\T^2} \Big( \Ad_{g(T,x)^{-1}} \fq_B(s,z;T,x) \Ad_{g(s,z)} \fp_{\eps^2}(z-y) e_\alpha - \fq_{\YMHFsmooth}(s,z;T,x) \fp_{\eps^2}(z-y) e_\alpha \Big) \dz \Big| \,\dx \\
		&\lesssim M_T(\YMHFsmooth,1)^2 T^{\frac{1-\kappa}{2}} K \;.
	\end{equs}
	Letting $\eps \rightarrow 0$ and using \eqref{eq:qB-ptwise-bounds} and H\"older's inequality, we get that for fixed $s \in (0,T/2)$, $y \in \T^2$,
	\begin{equs}
		 \int_{\T^2} \bigl| &\Ad_{g(T,x)^{-1}} \fq_B(s,y;T,x) \Ad_{g(s,y)} e_\alpha - \fq_{\YMHFsmooth}(s,y;T,x) e_\alpha \bigr|^2 \dx \\
		 &\lesssim M_T(\YMHFsmooth,1)^2 \bigl(M_T(B,1) + M_T(\YMHFsmooth,1)\bigr) T^{\frac{1-\kappa}{2}} K (T-s)^{-1} \;.
	\end{equs}
	The desired result then follows by integrating over $s \in (0,T/2)$ and $y \in \T^2$.
\end{proof}

\subsection{Fluctuations}
\label{sec:improved_stochastic_estimates_fluctuations}
In this subsection we work under the convention of Notation~\ref{notation:tstop}, where a jump at $\tstop$ is allowed for the YM gradient flow $B$. We will estimate the variance of several expressions involving the linear stochastic objects that will appear in the global iteration scheme.
\begin{lemma}
\label{le:cross_term_expansion}
	There exists some constant $C$ such that for any $p \geq 2$, $k \in [2]$, and $t \in (\tstop, 2]$, we have
	\begin{equ}
	\label{e:cross_term_expansion}
		\Bigl\|\int_{\T^2} \langle F_B(t, x), \bigl(\covd_B^k e^{(t - \tstop)\Delta^F_{B(\tstop +)}} \linear[B](\tstop)\bigr)(x)\rangle_{\fg} \,\dx\Bigr\|_{L^p(\Omega)} 
		\leq C p^{\frac12} \gronwall_t(C_0 F_B) \cdot \tstop^{\frac12} \|\covd_B F_B(t)\|_{L_x^2} \;.
	\end{equ}
\end{lemma}
\begin{proof} 
	By Gaussian hypercontractivity, it suffices to show 
	\begin{equs}
		\E \, \Big| &\int_{\T^2} \langle F_{B}(t, x), \bigl(\covd_{B}^k e^{(t - \tstop) \Delta^F_{B(\tstop +)}} \linear[B](\tstop) \bigr)(x)\rangle_{\fg} \, \dx \Big|^2 \\
		&\leq C \gronwall_{t}(2C_0 F_B)  \cdot \tstop \|\covd_B F_B(t)\|_{L_x^2}^2 \;. 
	\end{equs}
	By integration by parts and \eqref{e:jump_heat_kernel_bound}, the left-hand side is bounded by
	\begin{equs}
		\int_{0}^{\tstop} \ds &\iiint \dx \dx' \dy \abs{\covd_{B}^k F_B(t, x)} \abs{\covd_{B}^k F_B(t, x')} \abs{\fq_B(s, y; t, x)} \abs{\fq_B(s, y; t, x')}  \\
&\lesssim \gronwall_{t}(2C_0 F_B) \cdot \int_0^{\tstop} \ds \iint \dx \dx'\abs{\covd_{B}^k F_B(t, x)} \abs{\covd_{B}^k F_B(t, x')} (t - s)^{-1} \exp\Bigl(-c\frac{|x - x'|^2}{t - s}\Bigr) \;.
	\end{equs}
	For each $s \in [0, \tstop]$, the integral over $x, x'$ is bounded by $\|\covd_B F_B(t)\|_{L_x^2}^2$ by H\"older's inequality and Young's convolution inequality. This concludes the proof.
\end{proof}

\begin{lemma}
\label{le:phi_cross_term_estimate}
	There exists some constant $C$ such that for any $p \geq 2$, $t \in (\tstop, 2]$, and smooth function $\phismoothpart(\tstop) \colon \T^2 \ra \V$, we have
	\begin{equs}
		\Bigl\|&\int_{\T^2} \Big(|e^{(t - \tstop)\Delta_{B(\tstop +)}} \phismoothpart(\tstop)|^2 \cdot \langle e^{(t - \tstop)\Delta_{B(\tstop +)}} \phismoothpart(\tstop), e^{(t - \tstop)\Delta_{B(\tstop +)}} \philinear[B](\tstop)\rangle_{\V}\Big) \,\dx\Bigr\|_{L^p(\Omega)} \\
	\label{e:phi_cross_term_expectation}
		&\leq C p^{\frac12} \tstop^{\frac12} \|e^{(t - \tstop)\Delta_{B(\tstop +)}} \phismoothpart(\tstop)\|_{L^6}^3 \;.
	\end{equs}
\end{lemma}

\begin{proof}
	For notational simplicity in the proof we write $\higgs(t, x)$ as a shorthand for $\bigl(e^{(t - \tstop) \Delta_{B(\tstop +)}} \phismoothpart(\tstop)\bigr)(x)$. 	By Gaussian hypercontractivity, it suffices to show
	\begin{equs}
		 \E \Big| &\int_{\T^2} \Big(|e^{(t - \tstop) \Delta_{B(\tstop +)}} \phismoothpart(\tstop) |^2  \cdot \langle e^{(t - \tstop) \Delta_{B(\tstop +)}} \phismoothpart(\tstop), e^{(t - \tstop) \Delta_{B(\tstop +)}} \philinear[B](\tstop) \rangle_{\V} \Big) \dx \Big|^2 \\
		& \leq C \tstop  \|e^{(t - \tstop) \Delta_{B(\tstop +)}} \phismoothpart(\tstop) \|_{L^6}^6 \;.
	\end{equs}
	By the diamagnetic inequality (i.e.\ the first inequality in \eqref{e:jump_heat_kernel_bound}), the left-hand side is bounded by
	\begin{equs}
		\iint \dx \dx' &|\higgs(t, x)|^2 |\higgs(t, x')|^2 \int_0^{\tstop} \ds \int \dy \,  \fp(s, y; t, x)  \fp(s, y; t, x') \abs{ \higgs(t, x)} \abs{ \higgs(t, x') }  \\
		& \lesssim \iint \dx \dx' |\higgs(t, x)|^3 |\higgs(t, x')|^3 \int_0^{\tstop} \ds \, \fp(0, x; 2t - 2s, x')  \\
		& \lesssim \tstop \| \higgs(t) \|_{L^6}^6 \;,
	\end{equs}
where in the last step we used H\"older's inequality and Young's convolution inequality as in the proof of Lemma~\ref{le:cross_term_expansion}. 
\end{proof}

\begin{lemma}
\label{le:decay_mechanism_fluctuations}
	There exists some constant $C$ such that for any $p \geq 2$ and $t \in (\tstop, 2]$, and any $\fg$-valued function $F \colon \T^2 \ra \frkg$, we have
\begin{equ}
	 \Bigl\| \Wick{\Big\langle F, \Big[\bigl(e^{(t - \tstop)\Delta^F_{B(\tstop +)}} \linear[B](\tstop)\bigr)^1, \bigl(e^{(t - \tstop)\Delta^F_{B(\tstop +)}} \linear[B](\tstop)\bigr)^2\Big] \Big\rangle_{L_x^2}}\Bigr\|_{L^p(\Omega)} 
	\leq C p \gronwall_t(2C_0 F_B) \cdot \tstop^{\frac12} \|F\|_{L_x^\infty} \;.
\end{equ}
\end{lemma}
\begin{proof}
We have
\begin{equs}
	\,& \Wick{\Big\langle F, \Big[\bigl(e^{(t - \tstop)\Delta^F_{B(\tstop +)}} \linear[B](\tstop)\bigr)^1, \bigl(e^{(t - \tstop)\Delta^F_{B(\tstop +)}} \linear[B](\tstop)\bigr)^2\Big] \Big\rangle_{L_x^2}} \\
	&\qquad = \int_{\bigl([0, \tstop] \times \T^2 \bigr)^2 } \ds_1 \dy_1 \ds_2 \dy_2 \, \kernel(s_1, y_1, s_2, y_2) \, \Wick{\xi(s_1, y_1) \otimes \xi(s_2, y_2)} \;,
\end{equs}
where the kernel $\kernel$ is given by the map which sends $v \otimes w \in (\frkg^{2})^{\otimes 2}$ to
\begin{equ}
	\int_{\T^2} \dx \Big\langle F(x), \Big[(\fq_B(s_1, y_1; t, x) v)^1, (\fq_B(s_2, y_2; t, x) w)^2\Big] \Big\rangle_{\fg} \;.
\end{equ}
By Gaussian hypercontractivity, it suffices to estimate
\begin{equ}
	\bigl\|\kernel\bigr\|_{L_{s_1,y_1}^2 L_{s_2,y_2}^2([0, \tstop]^2)} \lesssim \gronwall_t(2C_0 F_B) \cdot \tstop^{\frac12} \|F\|_{L_x^\infty} \;.
\end{equ}
	Towards this end, applying \eqref{e:jump_heat_kernel_bound} and integrating in $x$, we obtain
	\begin{equs}
		|\kernel(s_1, y_1, s_2, y_2)| &\lesssim \gronwall_t(2C_0 F_B) \\
		&\quad \cdot (2t - s_1 - s_2)^{-1} \moralind(|y_1 - y_2| \lesssim (2t - s_1 - s_2)^{1/2})  \|F\|_{L_x^\infty}  \;,
	\end{equs}
	where $\moralind$ is as in Remark~\ref{remark:parabolic-indicator}.
	Squaring this estimate and integrating in $y_1, y_2$, we obtain
	\begin{equ}
		\bigl\|\kernel\bigr\|_{L_{s_1,y_1}^2 L_{s_2,y_2}^2([0, \tstop]^2)}^2
		\lesssim \gronwall_t(4C_0 F_B) \int_0^{\tstop} \int_0^{\tstop} \frac{\ds_1 \ds_2}{2t - s_1 - s_2}  \|F\|_{L_x^\infty}^2 \;.
	\end{equ}
	Since $t > \tstop$, we have
	\begin{equ}[e:boundIntegral]
		\int_0^{\tstop} \int_0^{\tstop} \frac{\ds_1 \ds_2}{2t - s_1 - s_2} \leq \int_0^{\tstop} \int_0^{\tstop} \frac{\ds_1 \ds_2}{2\tstop - s_1 - s_2}  \lesssim \tstop \;.
	\end{equ}
	This gives the desired result.
\end{proof}

\begin{lemma}
\label{le:curl_square_fluctuations}
	There exists some constant $C$ such that for any $p \geq 2$, $k \in [2]$ and $t \in (\tstop, 2]$, we have
	\begin{equ}
		\Bigl\|\Wick{\bigl\|\covd^k_{B(t)} e^{(t - \tstop)\Delta^F_{B(\tstop +)}} \linear[B](\tstop)\bigr\|_{L_x^2}^2}\Bigr\|_{L^p(\Omega)} 
		\leq C p \gronwall_t(2C_0 F_B) \bigl(1 + \Fnorm{F_B}_{[\tstop, t], 1}\bigr)^2 \cdot \tstop^{\frac12}(t - \tstop)^{-1} \;.
	\end{equ}
\end{lemma}
\begin{proof}
	We have that
	\begin{equs}
		\,&\Wick{\bigl\|\covd^k_{B(t)} e^{(t - \tstop)\Delta^F_{B(\tstop +)}} \linear[B](\tstop)\bigr\|_{L_x^2}^2} \\
		&= \int_{\bigl([0, \tstop] \times \T^2\bigr)^2} \ds_1 \dy_1 \ds_2 \dy_2 \, \kernel(s_1, y_1, s_2, y_2) \, \Wick{\xi(s_1, y_1) \otimes \xi(s_2, y_2)} \;,
	\end{equs}
	where the kernel $\kernel$ is given by the map which sends $v \otimes w \in (\frkg^{ 2})^{\otimes 2}$ to
	\begin{equ}
		\int_{\T^2} \dx \, \Big\langle (\covd^k_{B(t)} \fq_B)(s_1, y_1; t, x)v, (\covd^k_{B(t)} \fq_B)(s_2, y_2; t, x)w \Big\rangle_{\fg} \;.
	\end{equ}
	Here $(s_1, y_1), (s_2, y_2) \in [0, \tstop] \times \T^2$. By Gaussian hypercontractivity, it suffices to show
	\begin{equ}
		\bigl\|\kernel\bigr\|_{L_{s_1,y_1}^2 L_{s_2,y_2}^2([0, \tstop]^2)} \lesssim \gronwall_t(2C_0 F_B) \bigl(1 + \Fnorm{F_B}_{[\tstop, t], 1}\bigr)^2 \cdot \tstop^{\frac12}(t - \tstop)^{-1} \;.
	\end{equ}
	Towards this end, applying \eqref{e:jump_heat_kernel_bound_derivatives} and integrating in $x$, we obtain
	\begin{equs}
		|\kernel(s_1, y_1, s_2, y_2)| &\lesssim \gronwall_t(2C_0 F_B) \bigl(1 + \Fnorm{F_B}_{[\tstop, t], 1}\bigr)^2 \\
		&\quad \cdot (t - \tstop)^{-1}(2t - s_1 - s_2)^{-1} \moralind(|y_1 - y_2| \lesssim (2t - s_1 - s_2)^{1/2}) \;,
	\end{equs}
	where $\moralind$ is as in Remark~\ref{remark:parabolic-indicator}.
	Squaring this estimate and integrating in $y_1, y_2$, we obtain
	\begin{equ}
		\bigl\|\kernel\bigr\|_{L_{s_1,y_1}^2 L_{s_2,y_2}^2([0, \tstop]^2)}^2
		\lesssim \gronwall_t(4C_0 F_B) \bigl(1 + \Fnorm{F_B}_{[\tstop, t], 1}\bigr)^4 \cdot (t - \tstop)^{-2} \int_0^{\tstop} \int_0^{\tstop} \frac{\ds_1 \ds_2}{2t - s_1 - s_2} \;,
	\end{equ}
and we conclude by \eqref{e:boundIntegral} as before.
\end{proof}

\section{Global theory}\label{section:global-theory}

\subsection{Choice of time scale}
In this subsection, we first define the time scale used in a single iteration step. To make the covariant local theory work, we choose this time scale so that Assumptions~\ref{local:assumption-B} and~\ref{local:assumption-higgs-initial} hold. Before introducing the time scale, we set notation regarding the YM gradient flow.

\begin{notation}
	In this section, 
we use $\Psi_t(A)$ to denote the time $t$ slice of the YM gradient flow starting from $A$.\label{text-Psit}
\end{notation}

We first recall the following definition of the time scale that we will work with, which has already been defined in Definition~\ref{main:def-time-scales}. 
\begin{definition}
\label{defn:tau}
	Let $\logtimestep \leq 10^{-10}$ and $\stepupperbound \leq T_0$ be small constants that will be determined later. Here $T_0$ is the constant from Assumption~\ref{local:assumption-B}, which is specified in \eqref{prelim:eq-T0}. For any smooth connection one-form $A$ and smooth $\V$-valued function $\higgs$, we define
	\begin{equ}
		T^{(1)}(A, \higgs) \eqdef \min \{ \logtimestep \|F_A\|_{L^2}^{-1} \log \bigl(\|F_A \|_{L^2} \vee 2\bigr) , \| \higgs \|_{L^r}^{-\gamma_0}, \stepupperbound, \logtimestep \} \;,
	\end{equ}
	with $\gamma_0 > 2$ as specified in \eqref{e:constants_global_iteration}, and 
	\begin{equ}
		T^{(2)}(A, \higgs) \eqdef T^{(2)}(A) \eqdef \inf\Bigl\{t > 0:  \|F_A \|_{L^2}^2 - \|F(\Psi_t(A)) \|_{L^2}^2 \geq \logtimestep t^{-1} \log t^{-1}  \Bigr\} \;.
	\end{equ}
	More generally, we denote $T^{(i)}(A) \eqdef T^{(i)}(A, 0)$ for $i = 1, 2$. 
\end{definition}

\begin{remark}
\label{rmk:tau_lower_bound}
	$T^{(2)}(A)$ is well-defined, since the left-hand side is increasing in $t$ and starts from $0$ at 
	$t = 0$ (thus it is non-negative) and the right-hand side is decreasing in $t$, starts from $\infty$ at $t = 0$, and turns negative at $t=1$. Moreover, from the definition of $T^{(2)}$ we get 
	\begin{equ}
		 \|F_A \|_{L^2}^2 \geq \logtimestep T^{(2)}(A)^{-1} \abs{\log T^{(2)}(A) }  \;,
	\end{equ}
	which implies that 
	\begin{equ}
	\label{e:timestep2_lower_bound}
		T^{(2)}(A) \geq \logtimestep \bigl(\|F_A \|_{L^2} \vee 2 \bigr)^{-2}  \log \bigl(\|F_A \|_{L^2} \vee 2 \bigr) \;. 
	\end{equ}
\end{remark}

\begin{lemma}
\label{le:tau_smallness}
	There exists some constant $C$ such that the following holds.
	Let $\YMHFsmooth(t) = \Psi_t(\YMHFsmooth_0)$ be the YM gradient flow starting from some smooth $\YMHFsmooth_0$. Then
	\begin{equ}
	\label{e:tau_smallness}
		\bigl((T^{(1)} \wedge T^{(2)})(\YMHFsmooth_0)\bigr)^{\frac12}\|F_{\YMHFsmooth}\|_{L_t^2L_x^\infty([0,(T^{(1)} \wedge T^{(2)})(\YMHFsmooth_0)])} \leq C\sqrt{\logtimestep}\abs{\log\bigl((T^{(1)} \wedge T^{(2)})(\YMHFsmooth_0)\bigr)} \;.
	\end{equ}
	As a result, we can choose $\logtimestep$ small enough such that Assumption~\ref{local:assumption-B} is satisfied for $\YMHFsmooth$ on $[0, (T^{(1)} \wedge T^{(2)})(\YMHFsmooth_0) ]$.
\end{lemma}

\begin{proof}
	For notational simplicity, in the proof we use $T$ to denote $(T^{(1)} \wedge T^{(2)})(\YMHFsmooth_0)$. Recall from Remark~\ref{remark:curvature-subsolution-heat-equation} that $|F_{\YMHFsmooth}|$ is a non-negative subsolution to the heat equation. 
	Applying H\"older's inequality and Lemma~\ref{lemma:strichartz-estimate-using-decrease} yields
	\begin{equ}
		T^{\frac12}\|F_{\YMHFsmooth}\|_{L_t^2L_x^\infty([0,T])} \lesssim T\|F_{\YMHFsmooth_0}\|_{L^2}+T^{\frac12}\abs{\log T}^{\frac12}\bigl(\|F_{\YMHFsmooth_0}\|_{L^2}^2-\|F_{\YMHFsmooth}(T)\|_{L^2}^2\bigr)^{\frac12} \;.
	\end{equ}
	For the first term, the desired bound follows from the fact that
	\begin{equ}
		T \leq T^{(1)}(\YMHFsmooth_0) 
		\leq \logtimestep \|F_{\YMHFsmooth_0}\|_{L^2}^{-1} \log \big( \|F_{\YMHFsmooth_0} \|_{L^2} \vee 2 \big) \;.
	\end{equ}
	For the second term, by the definition of $T^{(2)}$, we get
	\begin{equ}
	\label{e:T2_decrease_consequence}
		0  \leq \|F_{\YMHFsmooth_0}\|_{L^2}^2 - \| F_{\YMHFsmooth(T)}\|_{L^2}^2 
		\leq \logtimestep T^{-1}  \abs{\log T} \;,
	\end{equ}
	which gives \eqref{e:tau_smallness}. In the following, recall the definition (Definition~\ref{def:Fnorm}) of $\Fnorm{\place}_{T, k}$. Note that \eqref{e:tau_smallness} also gives 
	\begin{equ}
		\Fnorm{F_{\YMHFsmooth}}_{T,0} \lesssim T^{\frac12}\|F_{\YMHFsmooth}\|_{L_t^2L_x^\infty([0,T])} \lesssim \sqrt{\logtimestep}\abs{\log T} \;.
	\end{equ}
	Therefore, $\Fnorm{F_{\YMHFsmooth}}_{T,100}$ can be controlled by Corollary~\ref{cor:curvature-higher-regularity-estimates}. 
	
	Now we proceed to bound the term $T^{\frac{1-\roughdataweight}{2}}\Fnorm{t^{-\frac{1-\roughdataweight}{2}}F_{\YMHFsmooth}}_{T,1}$ in $\lossfactor_T(F_{\YMHFsmooth})$.
	Since $F_{\YMHFsmooth}$ solves the covariant heat equation (by Lemma~\ref{lemma:curvature-equation}), Proposition~\ref{prop:pB-ptwise-bounds} and Young's inequality imply that
	\begin{equ}
		t^{\frac12} \|\covd_{\YMHFsmooth}F_{\YMHFsmooth}(t)\|_{L_x^\infty}
		\lesssim \Bigl(1+\Fnorm{F_{\YMHFsmooth}}_{T,0}\Bigr)
		\|F_{\YMHFsmooth}(t/2)\|_{L_x^\infty} \;.
	\end{equ}
	Expanding Definition~\ref{def:Fnorm} and using the preceding estimate, we obtain
        \begin{equs}
        \,& T^{\frac{1-\roughdataweight}{2}}\Fnorm{t^{-\frac{1-\roughdataweight}{2}}F_{\YMHFsmooth}}_{[0,T],1} \\
        	&=T^{\frac{1-\roughdataweight}{2}}\Bigl(\|t^{\frac{\roughdataweight}{2}}F_{\YMHFsmooth}\|_{L_t^2L_x^\infty([0,T])}+\|t^{-\frac{1-\roughdataweight}{2}}F_{\YMHFsmooth}\|_{L_t^1L_x^\infty([0,T])} \\
        	&\qquad+\|t^{\frac{1+\roughdataweight}{2}}\covd_{\YMHFsmooth}F_{\YMHFsmooth}\|_{L_t^2L_x^\infty([0,T])}+\|t^{\frac{\roughdataweight}{2}}\covd_{\YMHFsmooth}F_{\YMHFsmooth}\|_{L_t^1L_x^\infty([0,T])}\Bigr) \\
        	&\lesssim T^{\frac{1-\roughdataweight}{2}}\Bigl(1+\Fnorm{F_{\YMHFsmooth}}_{[0,T],0}\Bigr)\Bigl(\|t^{\frac{\roughdataweight}{2}}F_{\YMHFsmooth}\|_{L_t^2L_x^\infty([0,T])}+\|t^{-\frac{1-\roughdataweight}{2}}F_{\YMHFsmooth}\|_{L_t^1L_x^\infty([0,T])} \\
        	&\qquad+\|t^{\frac{\roughdataweight}{2}}F_{\YMHFsmooth}(t/2)\|_{L_t^2L_x^\infty([0,T])}+\|t^{-\frac{1-\roughdataweight}{2}}F_{\YMHFsmooth}(t/2)\|_{L_t^1L_x^\infty([0,T])}\Bigr) \\
        	&\lesssim \Bigl(1+\Fnorm{F_{\YMHFsmooth}}_{[0,T],0}\Bigr)T^{\frac12}\|F_{\YMHFsmooth}\|_{L_t^2L_x^\infty([0,T])} \;,
        \end{equs}	
	where the last inequality follows from H\"older's inequality and the change of variables $t\mapsto t/2$.

	Combining the preceding estimates, we obtain
	\begin{equ}
		\lossfactor_T(F_{\YMHFsmooth})\leq\exp\bigl(C\sqrt{\logtimestep}\abs{\log T}\bigr)\Bigl(1+C\sqrt{\logtimestep}\abs{\log T}+C\logtimestep\abs{\log T}^2+C\bigl(\sqrt{\logtimestep}\abs{\log T}\bigr)^{101!}\Bigr)^{100} \;.
	\end{equ}
	This is smaller than $T^{-\lossdel}$ if $\logtimestep$ is sufficiently small. This finishes the proof.
\end{proof}

In practice, our initial data will not be smooth. We will therefore decompose it into a large smooth component and a small rough remainder. The smallness of the remainder will be measured in the following gauge invariant Hölder norm, which is a special case of the norm from Definition~\ref{defn:covariant_Holder} since the time is always taken to be $T^{(1)} \wedge T^{(2)}$.

\begin{definition}
\label{defn:covariant_Holder_global}
	Let $\Asmoothpart$ be some smooth connection. Recall Notation~\ref{notation:linear-propagation}. For any $\alpha < 1$, and $\genV$-valued distributional one-form $\Aroughpart$ we define
	\begin{equ}
		\| \Aroughpart \|_{\covHolder{\Asmoothpart, F}^{\alpha}} \eqdef \sup_{0 \leq m \leq 1 } \sup_{s \in (0, (T^{(1)} \wedge T^{(2)})(\Asmoothpart)] } \big( s^{\frac{-\alpha + m}{2}} \wedge 1 \big) \|\covd^{(m)}_{\Psi_s(\Asmoothpart)} \, e^{s\Delta^F_{\Asmoothpart}} \Aroughpart \|_{L^\infty} \;. 
	\end{equ}
	Also, for any $\genV$-valued distribution $\phiroughpart$ we define
	\begin{equ}
		\| \phiroughpart \|_{\covHolder{\Asmoothpart}^{\alpha}} \eqdef \sup_{0 \leq m \leq 1} \sup_{s \in (0, (T^{(1)} \wedge T^{(2)})(\Asmoothpart)] } \big( s^{\frac{-\alpha + m}{2}} \wedge 1 \big) \|\covd^{(m)}_{\Psi_s(\Asmoothpart)} \, e^{s\Delta_{\Asmoothpart}} \phiroughpart \|_{L^\infty} \;.
	\end{equ}
\end{definition}

\begin{definition}
\label{defn:L_admissible}
	Let $L \geq 1$, $(\Asmoothpart, \Aroughpart)$ be (possibly distributional) $\fg$-valued $1$-forms and $(\phismoothpart, \phiroughpart)$ be (possibly distributional) $\V$-valued functions. We say that the decomposition $A = \Asmoothpart + \Aroughpart$, $\higgs = \phismoothpart +\phiroughpart$ is $L$-admissible if $\Asmoothpart, \phismoothpart$ are smooth, and $\Aroughpart, \phiroughpart$ satisfy
	\begin{equ}
	\label{e:roughpart_smallness}
		\| \Aroughpart \|_{\covHolder{\Asmoothpart, F}^{-\kappa}}  \leq L^{\fa} \;, \quad 
		\| \phiroughpart \|_{\covHolder{\Asmoothpart}^{-\kappa}}  \leq L^{\fa} \;,
	\end{equ}
	where $\fa$ is specified as in \eqref{e:constants_global_iteration}. 
	For $L$-admissible $A = \Asmoothpart + \Aroughpart$, $\higgs = \phismoothpart + \phiroughpart$,
	we define 
	\begin{equ}
		T^{(L)}(\Asmoothpart, \Aroughpart, \phismoothpart, \phiroughpart) \eqdef (T^{(1)} \wedge T^{(2)})(\Asmoothpart, \phismoothpart)  \wedge L^{-1} \;.
	\end{equ} 
\end{definition}

\begin{proposition}
\label{pr:admissible_timescale}
	There exists a choice of $\logtimestep$ independent of $\stepupperbound$ in Definition~\ref{defn:tau} such that the following holds. 
	Let $A_0 = \Asmoothpart_0 + \Aroughpart_0$, $\higgs_0 = \phismoothpart_0 + \phiroughpart_0$ be $L$-admissible. Set $T^{(L)} = T^{(L)}(\Asmoothpart_0, \Aroughpart_0, \phismoothpart_0, \phiroughpart_0)$.
	Then, there exists some YM gradient flow $B \in C^{\infty} ((0, T^{(L)}] \times \T^2)$ and gauge transform $g \in C_t^0 \CG_{1-\kappa}([0, T^{(L)}])$ such that $g(t)$ is smooth for any $t \in (0, T^{(L)}]$ and such that $g^{-1} \act B$ converges to $A_0$ in $\Cs_x^{-\kappa}$ as $t\rightarrow 0$. Moreover, Assumptions~\ref{local:assumption-B} and~\ref{local:assumption-higgs-initial} are satisfied with $T=T^{(L)}$, $\shiggsdet_0=g_0 \act \phismoothpart_0$, and $\rhiggsdet_0=g_0 \act \phiroughpart_0$. In particular, the smallness condition~\eqref{eq:Sdet-higgsdet-main-prop-smallness-assumption} holds. Finally, it holds that $g_0\act (A_0,\higgs_0) \in \initialclass_{T^{(L)}}$, where the initial class is as in \eqref{main:eq-initial-class}.
\end{proposition}

\begin{proof}
	Let $\YMHFsmooth$ be the Yang--Mills gradient flow starting from $\Asmoothpart_0$. Recall the notation $M_T(\place, \place)$ from Notation~\ref{notation:gronwall-F-pert}. Note that this quantity is manifestly increasing in the parameter $T$.
	We first verify the assumptions of Proposition~\ref{pr:admissible_timescale_general}. 
	By choosing $\logtimestep$ sufficiently small, Lemma~\ref{le:tau_smallness} and Corollary~\ref{cor:curvature-higher-regularity-estimates} imply	\begin{equ}
		M_{T^{(L)}}(\YMHFsmooth,2) \leq M_{(T^{(1)} \wedge T^{(2)})(\Asmoothpart_0, \phismoothpart_0)} (\YMHFsmooth, 2) \leq \bigl((T^{(1)} \wedge T^{(2)})(\Asmoothpart_0, \phismoothpart_0)\bigr)^{-\fa} \leq \bigl(T^{(L)}\bigr)^{-\fa} \;.
	\end{equ}
	By the covariant heat kernel estimate \eqref{eq:qB-ptwise-bounds} and \eqref{e:roughpart_smallness},
	\begin{equ}
		\|\Aroughpart_0\|_{\covHolder{\Asmoothpart_0,F,T^{(L)},2}^{-\kappa}}
		\lesssim M_{T^{(L)}}(\YMHFsmooth,2)
		\|\Aroughpart_0\|_{\covHolder{\Asmoothpart_0,F}^{-\kappa}} 
		\leq \bigl(T^{(L)}\bigr)^{-2\fa} \;,
	\end{equ}
	Applying Proposition~\ref{pr:admissible_timescale_general} with $K =\bigl(T^{(L)}\bigr)^{-2\fa}$ gives the required Yang--Mills gradient flow $B$ and gauge transformation $g$.

	Now we verify Assumption~\ref{local:assumption-B}. Note that by Lemma~\ref{le:tau_smallness}, the assumption is already true for $\YMHFsmooth$. Using \eqref{e:Fnorm_bound_roughB_general} to bound the error when replacing $\YMHFsmooth$ by $B$, it is straightforward to see that 	
	\begin{equs}
		\Fnorm{F_B}_{T^{(L)},0}
		&\lesssim \logtimestep^{\frac{1 - \kappa}{2} - 2 \fa } +\sqrt{\logtimestep}\abs{\log T^{(L)}} \;, \\
		\bigl(T^{(L)}\bigr)^{\frac{1-\roughdataweight}{2}}
		\Fnorm{t^{-\frac{1-\roughdataweight}{2}}F_B}_{T^{(L)},1}
		&\lesssim
		\logtimestep^{\frac{1 - \kappa}{2} - 2 \fa } +\sqrt{\logtimestep}\abs{\log T^{(L)}}
		+\logtimestep\abs{\log T^{(L)}}^2 \;.
	\end{equs}
	With the same argument as in Lemma~\ref{le:tau_smallness}, it follows that Assumption~\ref{local:assumption-B} holds if $\logtimestep$ is sufficiently small. 
	
	Next, we verify Assumption~\ref{local:assumption-higgs-initial}.
	From the definition of $T^{(1)}$ and the gauge invariance of the $L_x^r$-norm, it follows that 
	\begin{equation*}
	\bigl\| g_0 \act \phismoothpart_0 \bigr\|_{L_x^r} = \bigl\|  \phismoothpart_0 \bigr\|_{L_x^r} \leq (T^{(L)})^{-\frac{1}{\gamma_0}}. 
	\end{equation*}
	Furthermore, Lemmas~\ref{le:B-YMHFsmooth-heat-comparison} and~\ref{le:tau_smallness},
    together with the $L$-admissibility of the initial data, give
	\begin{equs}
		\|g_0 \act \phiroughpart_0 \|_{\covHolder{B(0), T^{(L)} }^{-\kappa}}
		&\lesssim
		\bigl(1+\Fnorm{F_{\YMHFsmooth}}_{T^{(L)},1}\bigr)^2
		\bigl(1+ \bigl( T^{(L)} \bigr)^{\frac{1-\kappa}{2} - 2 \fa} \bigr)
		\|\phiroughpart_0\|_{\covHolder{\Asmoothpart_0, T^{(L)}}^{-\kappa}} \lesssim \bigl( T^{(L)} \bigr)^{-\lossdel - \fa} \;.
	\end{equs}
	Using that $\gamma_0>2$ and that $\fa$ and $\lossdel$ are small, together with the fact that $T^{(L)} \leq \logtimestep$ is sufficiently small, this can clearly be bounded by $(T^{(L)})^{-\frac{1}{\gamma_0}}$.  

Finally, the smallness condition \eqref{eq:Sdet-higgsdet-main-prop-smallness-assumption}
follows from Assumptions~\ref{local:assumption-B}
and~\ref{local:assumption-higgs-initial},
as explained immediately after their statements, and the fact that $g_0\act (A_0,\higgs_0) \in \initialclass_{T^{(L)}}$ follows directly from Proposition~\ref{pr:admissible_timescale_general} and our estimates of $g_0 \act \phismoothpart_0$ and $g_0 \act \phiroughpart_0$ above.
\end{proof}
\begin{remark}\label{remark:logtimestep-stepupperbound}
From now on, we fix the constant $\logtimestep$ in Definition~\ref{defn:tau} so that the conclusions of Proposition~\ref{pr:admissible_timescale} hold. As for the constant $\stepupperbound$, we can always choose $\stepupperbound$ as small as we want. The main places where we need it to be small enough are Proposition~\ref{prop:modified-iteration-uniform-bounds}, Lemma~\ref{le:linear_mod_conditional_moments}, and Proposition~\ref{prop:det-time-energy-bound}. 
\end{remark}

\begin{definition}
\label{defn:tau_L}
	Let $A_0 = \Asmoothpart_0 + \Aroughpart_0$, $\higgs_0 = \phismoothpart_0 + \phiroughpart_0$ be $L$-admissible. Let $B$ be as in Proposition~\ref{pr:admissible_timescale}. We define the stopping time
	\begin{equs}
		\tau^{(L)} = \sup \Bigl\{ T \in (0, T^{(L)} ] :\;
		 & E(T, g_0\act A_0, g_0 \act \higgs_0)  \text{ holds, and}  \\
		 &\sup_{0 < t \leq T} \|\linear[B](t)\|_{\covHolder{B(t), F, T^{(L)}}^{-\kappa}} \vee \sup_{0 < t \leq T} \|\philinear[B](t)\|_{\covHolder{B(t), T^{(L)}}^{-\kappa} } 
		\leq \frac12 T^{-\lambda_1}  \Bigr\} \;,
	\end{equs} 
	where $E(T, g_0\act A_0, g_0 \act \higgs_0)$ is the event in Theorem~\ref{main:thm-covariant-local}, and $g_0$ is the gauge transform defined in Proposition~\ref{pr:admissible_timescale}.
\end{definition}
\begin{remark}
Here, we also added bounds on the $\covHolder{B(t), F, T^{(L)}}^{-\kappa}$ and $\covHolder{B(t), T^{(L)}}^{-\kappa}$-norm of $\linear[B](t)$ and $\philinear[B](t)$ to the definition of the stopping time. This will not be needed until Section~\ref{section:decay-estimates-in-expectation}. 
\end{remark}

	The convention here is that $T$ denotes a time step that is deterministic (or at least measurable with respect to the initial data), whereas $\tau$ denotes a stopping time depending on the realisation of the noise. This distinction becomes relevant when applying the strong Markov property at certain steps. However, on typical events the values of $T$ and $\tau$ coincide, which is captured by the following lemma.

\begin{lemma}
\label{le:tau_tail_bound}
	There exists some constant $c > 0$ such that the following holds. Let $A_0 = \Asmoothpart_0 + \Aroughpart_0$, $\higgs_0 = \phismoothpart_0 + \phiroughpart_0$ be $L$-admissible. Let $T^{(L)}, \tau^{(L)}$ be defined as in Definitions~\ref{defn:L_admissible} and~\ref{defn:tau_L}. For any $K \geq \big( T^{(L)} \big)^{-1}$, we have 
	\begin{equ}
		\P \bigl( \tau^{(L)} < K^{-1} \bigr) \leq c^{-1} \exp\bigl(- c K^{\lambda_1} \bigr) \;,
	\end{equ}
	where $\lambda_1$ is specified as in \eqref{prelim:eq-parameters}. 
\end{lemma}
\begin{proof}
	For any $K \geq (T^{(L)})^{-1}$, the definition of $\tau^{(L)}$ gives
	\begin{equs}
		\{\tau^{(L)}<K^{-1}\}
		&\subset E(K^{-1}, g_0\act A_0, g_0\act \higgs_0)^c  \\
		&\quad\cup \Bigl\{\sup_{0<t\leq K^{-1}}\|\linear[B](t)\|_{\covHolder{B(t),F,T^{(L)}}^{-\kappa}} \vee \sup_{0 < t \leq K^{-1}} \|\philinear[B](t)\|_{\covHolder{B(t), T^{(L)}}^{-\kappa} } 
 \geq \frac12K^{\lambda_1}\Bigr\} \;.
	\end{equs}
	From Proposition~\ref{pr:admissible_timescale}, we know that $B$ satisfies Assumption~\ref{local:assumption-B} on $[0,T^{(L)}]$. With the monotonicity of the curvature along the Yang--Mills gradient flow, $B$ satisfies Assumption~\ref{local:assumption-B} on $[0,2T^{(L)}]$, with $\lossdel$ replaced by a slightly larger exponent $\lossdel'$. Since $\lossdel'$ is still sufficiently small, Lemma~\ref{le:lollipop_negative_holder_uniform} still applies under this assumption. Using the lemma applied with $T=K^{-1}$ and $S=T^{(L)}$ gives
	\begin{equs}
		\P\Bigl(\sup_{0<t\leq K^{-1}}\|\linear[B](t)\|_{\covHolder{B(t),F,T^{(L)}}^{-\kappa}} \vee \sup_{0 < t \leq K^{-1}} \|\philinear[B](t)\|_{\covHolder{B(t), T^{(L)}}^{-\kappa} } 
 \geq \frac12K^{\lambda_1}\Bigr) & \lesssim \exp\Bigl(-cK^{2\lambda_1}(T^{(L)})^{-\frac{2\kappa}{5}}\Bigr) \\
 &\lesssim \exp(-cK^{\lambda_1}) \;.
	\end{equs}
	The first event is bounded by the same right-hand side by Theorem~\ref{main:thm-covariant-local}. Although the exponent $\lambda_1$ does not appear in the statement of the theorem, this bound is established in its proof at the end of Section~\ref{section:local-proof}. The lemma is proved.
\end{proof}

Later, we will iterate on the collection of data $(\Asmoothpart, \Aroughpart, \phismoothpart, \phiroughpart)$ rather than on $(A, \Phi)$ alone. This will be useful for the global theory, as it includes the information on how to break the solution $(A, \Phi)$ into a large smooth piece and a small, rough remainder.
We first define the evolution of this collection up to time $\tau^{(L)}$.

\begin{definition}
\label{defn:one_step_iteration_scheme}
	Let $A_0 = \Asmoothpart_0 + \Aroughpart_0$, $\higgs_0 = \phismoothpart_0 + \phiroughpart_0$ be $L$-admissible.
	We define $\flow_t \bigl(\Asmoothpart_0, \Aroughpart_0, \phismoothpart_0, \phiroughpart_0; \xi, \zeta \bigr) = \bigl(\Asmoothpart(t), \Aroughpart(t), \phismoothpart(t), \phiroughpart(t)\bigr)$ up to time $\tau^{(L)}$ as follows. 
	
	First we define the deterministic evolution up to time $T^{(L)}$. Let $\YMHFsmooth$ be the YM gradient flow starting from $\Asmoothpart_0$. 	
	Let the gauge transform $g \in C_t^0 \CG_{1-\kappa}([0, T^{(L)}])$ and the YM gradient flow $B(t)$ started from $g_0 \act A_0$ be as obtained in Proposition~\ref{pr:admissible_timescale}. Then let
$(\Sdet, \higgsdet)$ solve \eqref{eq:Sdet-equation}--\eqref{eq:higgsdet-equation} with $\higgsdet_0 = \shiggsdet_0 + \rhiggsdet_0 = g_0 \act \phismoothpart_0 + g_0 \act \phiroughpart_0$. 
For $t \in \big( 0, T^{(L)} \big] $ we set 
	\begin{equ}[e:smooth_part_evolution] 
		\Asmoothpart (t) \eqdef B(t) + \Sdet(t)\;, \qquad  \phismoothpart(t) \eqdef \higgsdet(t) \;. 
	\end{equ}

	Finally, let $(A(t), \higgs(t))$ be the stochastic evolution defined as in Theorem~\ref{main:thm-covariant-local} up to time $\tau^{(L)}$, with gauge choice $B$ and initial data $g_0 \act (A_0, \higgs_0)$. 
	For $t \in \big( 0, \tau^{(L)} \big] $ we set
	\begin{equ}[e:rough_part_evolution] 
		\Aroughpart (t) \eqdef A(t) - \Asmoothpart(t) \;, \qquad  \phiroughpart(t) \eqdef \higgs(t) - \phismoothpart(t) \;. 
	\end{equ}
\end{definition}

\begin{notation}
\label{notation:YMHF}
	Throughout the rest of this section, whenever $A_0 = \Asmoothpart_0 + \Aroughpart_0$ is $L$-admissible, we use $B$ to denote the YM gradient flow started from the rough initial data $g_0 \act A_0$, as in Proposition~\ref{pr:admissible_timescale}. We use $\YMHFsmooth$ to denote the YM gradient flow started from $\Asmoothpart_0$ (without applying the gauge transform $g_0$!). The same convention applies to other instances of $B$ and $\YMHFsmooth$, possibly with decorations (e.g.\ $B_m$, $\YMHFsmooth_m$ in Definition~\ref{defn:time_steps}).
\end{notation}

We can furthermore single out the most irregular part in $(\Aroughpart, \phiroughpart)$ as follows. This is the only part of the argument where the estimates obtained in the covariant local well-posedness theory enter the global theory.
\begin{proposition}
\label{pr:remainder_estimate}
	Let $A_0 = \Asmoothpart_0 + \Aroughpart_0$, $\higgs_0 = \phismoothpart_0 + \phiroughpart_0$ be $L$-admissible. Let $(\Aroughpart, \phiroughpart)$ be as defined in Definition~\ref{defn:one_step_iteration_scheme}. We can further decompose
	\begin{equ}
	\label{e:solution_decom}
		\Aroughpart = \linear[B] + \Arem \;, \qquad \phiroughpart = \philinear[B] + \phirem \;,
	\end{equ}
	where $\linear[B], \philinear[B]$ solve \eqref{local:eq-linear-stochastic}, but with $\xi_{\Beps}, \zeta_{\Beps}$ replaced by $\xi, \zeta$. 
	For any $2 \lossdel \leq \beta < 1- \lambda$ and $t \in (0, \tau^{(L)}]$, we have the bounds
	\begin{equs}
		\| \Arem(t) \|_{\covHolder{\Asmoothpart(t), F}^{\beta - 2 \lossdel}} & \lesssim \big(\tau^{(L)} \big)^{\frac{1 - \lambda}{2} - \frac{1}{\gamma_0} - \lambda_2} \Bigl( \bigl(\tau^{(L)} \bigr)^{-\frac{\beta}{2}} + \bigl(T^{(1)} \wedge T^{(2)} \bigr)(\Asmoothpart(t))^{-\frac{\beta}{2}} \Bigr) \;, \\
		\| \phirem(t) \|_{\covHolder{\Asmoothpart(t)}^{\beta - 2 \lossdel}} & \lesssim \big(\tau^{(L)} \big)^{\frac{1 - \lambda - \beta}{2} - \frac{1}{\gamma_0} - \lambda_2} \;.
	\end{equs}
\end{proposition}

\begin{proof}
	We only prove the bound for $\Arem$. The bound for $\phirem$ follows similarly. By Proposition~\ref{local:prop-main},  $\Arem(t)$ has a decomposition $\Arem(t) = \sum_{\substack{N \geq N_{\lo}\\ N \in \dyadic} } \Arem[N](t)$ such that for $0 \leq m \leq 1$ and $0 \leq t \leq \tau^{(L)}$,
	\begin{equ}
		\| \covd_{B(t)}^{(m)} \Arem[N](t) \|_{L^\infty} \lesssim N^{m - (1-\lambda) } \big(\tau^{(L)} \big)^{-\frac{1}{\gamma_0} - \lambda_2} \;,
	\end{equ}
	with $N_{\lo} \sim \bigl(\tau^{(L)} \bigr)^{-\frac12}$ (recall \eqref{local:eq-N-lo}). 
	Combining this with Proposition~\ref{prop:Sdet-higgsdet-estimates} we have 
	\begin{equs}
		 \| \covd_{\Asmoothpart(t)} \Arem[N](t) \|_{L^\infty} 
		 &\leq \| \covd_{B(t)} \Arem[N](t) \|_{L^\infty} + \| \Sdet(t) \|_{L^\infty} \| \Arem[N](t) \|_{L^\infty} \\
		 &\lesssim N^{\lambda} \big(\tau^{(L)} \big)^{-\frac{1}{\gamma_0} - \lambda_2} + \big(\tau^{(L)} \big)^{- \lossdel + \frac12 - \frac2r - \frac{2}{\gamma_0}} \cdot N^{-(1 - \lambda)} \big(\tau^{(L)} \big)^{-\frac{1}{\gamma_0} - \lambda_2 } \\
		 &\lesssim N^{\lambda } \big(\tau^{(L)} \big)^{-\frac{1}{\gamma_0} - \lambda_2} \;,
	\end{equs}
	where in the last step we used $N^{-1} \leq N_{\lo}^{-1} \lesssim \bigl( \tau^{(L)} \bigr)^{\frac{1}{2}}$.
	The desired result then follows from Lemma~\ref{le:heat_flow_norm_sol_norm} (applied with $M = N_{\lo}$ and $\nu = 1 - \lambda$), together with the fact that $N_{\lo} \sim \bigl(\tau^{(L)} \bigr)^{-\frac12}$. 
	\end{proof}

\begin{lemma}
\label{le:heat_flow_norm_sol_norm}
	Let $\Asmoothpart$ be some smooth connection, $0 \leq \nu \leq 1$, and $M \in \dyadic$ a dyadic integer. Recall Definition~\ref{defn:covariant_Holder_global}. Suppose that for some $\genV$-valued distributional one-form $\Aroughpart$ we have the decomposition $\Aroughpart = \sum_{\substack{N \geq M \\ N \in \dyadic} } \Aroughpart_N$. Then for any $2\lossdel \leq \beta < \nu \leq 1$ we have
	\begin{equ}
	\label{e:Arough_heat_flow_norm_bound}
		\| \Aroughpart \|_{\covHolder{\Asmoothpart, F}^{\beta - 2\lossdel}}  \lesssim M^{-\nu} \bigl( M^{\beta} + \bigl(T^{(1)} \wedge T^{(2)} \bigr)(\Asmoothpart)^{-\frac{\beta}{2} } \bigr) \sup_{m \in \{0, 1 \}} \sup_{\substack{N \geq M\\ N  \in \dyadic} } N^{\nu - m} \| \covd_{\Asmoothpart}^{(m)} \Aroughpart_N \|_{L^{\infty}} \;.
	\end{equ}
	Similarly, for any $\phiroughpart = \sum_{\substack{N \geq M\\ N \in \dyadic} } \phiroughpart_N$ we have
	\begin{equ}
	\label{e:phirough_heat_flow_norm_bound}
		\| \phiroughpart \|_{\covHolder{\Asmoothpart}^{\beta - 2 \lossdel}} \lesssim M^{-\nu + \beta} \sup_{m \in \{0, 1 \}} \sup_{\substack{N \geq M\\ N \in \dyadic} } N^{\nu - m} \| \covd_{\Asmoothpart}^{(m)} \phiroughpart_N \|_{L^{\infty}}  \;.
	\end{equ}
\end{lemma}
\begin{proof}
	We will first prove \eqref{e:phirough_heat_flow_norm_bound}, and then reduce \eqref{e:Arough_heat_flow_norm_bound} to this case. 
	For \eqref{e:phirough_heat_flow_norm_bound}, it suffices to bound $\| \covd_{\Psi_s(\Asmoothpart)}^{(m)} e^{s \Delta_{\Asmoothpart}} \phiroughpart \|_{L^\infty}$ for $0 \leq m \leq 1$ and $s \leq (T^{(1)} \wedge T^{(2)})(\Asmoothpart)$. We first observe that, by our choice of $\logtimestep$ in the definition of $T^{(1)}$, Assumption~\ref{local:assumption-B} is satisfied. This fact will be used frequently in the sequel without further mention.
	The case of $m = 0$ directly follows from Proposition~\ref{prop:qB-ptwise-bds}. Now we focus on the case of $m = 1$. In the following proof we denote 
	\begin{equ}
		Q_{\Aroughpart} = \sup_{m \in \{0, 1 \}} \sup_{\substack{N \geq M\\ N  \in \dyadic} } N^{\nu - m} \| \covd_{\Asmoothpart}^{(m)} \Aroughpart_N \|_{L^{\infty}} \;, \quad
		Q_{\phiroughpart} = \sup_{m \in \{0, 1 \}} \sup_{\substack{N \geq M\\ N  \in \dyadic} } N^{\nu - m} \| \covd_{\Asmoothpart}^{(m)} \phiroughpart_N \|_{L^{\infty}}  \;.
	\end{equ}
	
	Note that by \eqref{eq:cov-heat-operator-covd-commutator} and Proposition~\ref{prop:qB-ptwise-bds}, we get that for $s \leq (T^{(1)} \wedge T^{(2)})(\Asmoothpart) $,
	\begin{equs}
		\| \covd_{\Psi_s(\Asmoothpart)} e^{s \Delta_{\Asmoothpart}} \phiroughpart_N \|_{L^\infty} &= \| e^{s \Delta_{\Asmoothpart}^F} \covd_{\Asmoothpart} \phiroughpart_N \|_{L^\infty} \\
		& \lesssim \bigl((T^{(1)} \wedge T^{(2)})(\Asmoothpart) \bigr)^{-\lossdel} \| \covd_{\Asmoothpart} \phiroughpart_N \|_{L^\infty} \;.
	\end{equs}
	On the other hand, by Proposition~\ref{prop:pB-ptwise-bounds}, we get 
	\begin{equs}
		\| \covd_{\Psi_s(\Asmoothpart)} e^{s \Delta_{\Asmoothpart}} \phiroughpart_N \|_{L^\infty} &\lesssim \big((T^{(1)} \wedge T^{(2)})(\Asmoothpart) \big)^{-\lossdel} s^{-\frac{1}{2}} \| \phiroughpart_N \|_{L^\infty}  \;.
	\end{equs}
	Interpolating between these two bounds and summing over $N$, we get that for $s \leq (T^{(1)} \wedge T^{(2)})(\Asmoothpart)  $, 
	\begin{equs}
		s^{\lossdel + \frac{1 - \beta}{2}} \| \covd_{\Psi_s(\Asmoothpart)} e^{s \Delta_{\Asmoothpart}} \phiroughpart \|_{L^\infty} &\lesssim \sum_{\substack{N \geq M\\ N  \in \dyadic} } \| \covd_{\Asmoothpart} \phiroughpart_N \|_{L^\infty}^{\beta} \| \phiroughpart_N \|_{L^\infty}^{1 - \beta}  \\ 
		&\leq Q_{\phiroughpart} \sum_{\substack{N \geq M\\ N  \in \dyadic} }  N^{- \nu} \cdot N^{\beta} \lesssim Q_{\phiroughpart} M^{-(\nu - \beta)}  \;,
	\end{equs}
	and the desired estimate follows. 
	
	Now we reduce \eqref{e:Arough_heat_flow_norm_bound} to the previous case. It suffices to bound the error when replacing $e^{s\Delta_{\Asmoothpart}^F}$ by $e^{s\Delta_{\Asmoothpart}}$. 
	By Corollary~\ref{cor:heat-equations-comparison} with $\alpha = 2 \kappa $, for any $0 \leq m \leq 1$ and $s \leq (T^{(1)} \wedge T^{(2)})(\Asmoothpart)$ we have
	\begin{equ}
		\| \covd_{\Psi_s(\Asmoothpart)}^{(m)} \bigl(e^{s \Delta_{\Asmoothpart}^{F}} - e^{s \Delta_{\Asmoothpart}}\bigr)\Aroughpart_N \|_{L^\infty} \lesssim \big((T^{(1)} \wedge T^{(2)})(\Asmoothpart) \big)^{-\lossdel - \frac{1 - 2 \kappa}{2} } s^{\frac{1 - m - 2 \kappa} {2}}  \| \Aroughpart_N \|_{L^\infty} \;. 
	\end{equ}
	Here we used the fact that Assumption~\ref{local:assumption-B} holds on $[0, (T^{(1)} \wedge T^{(2)})(\Asmoothpart)]$ to bound the right-hand side of \eqref{e:heat-equations-comparison}.
	Then, we get that 
	\begin{equs}
		s^{\lossdel + \frac{1 - \beta}{2}} \| \covd_{\Psi_s(\Asmoothpart)} \bigl(e^{s \Delta_{\Asmoothpart}^{F}} - e^{s \Delta_{\Asmoothpart}}\bigr) \Aroughpart \|_{L^\infty} &\lesssim \big((T^{(1)} \wedge T^{(2)})(\Asmoothpart) \big)^{-\frac{\beta}{2} } Q_{\Aroughpart} \sum_{\substack{N \geq M\\ N \in \dyadic} }  N^{-\nu} \\
		& \lesssim \big((T^{(1)} \wedge T^{(2)})(\Asmoothpart) \big)^{-\frac{\beta}{2} } Q_{\Aroughpart} M^{- \nu }  \;,
	\end{equs}
	which gives the desired bound.
\end{proof}

\subsection{Iteration}
\label{sec:iteration}

Now we state our iteration scheme corresponding to the times $T^{(L)}$, $\tau^{(L)}$ defined in Definitions~\ref{defn:L_admissible} and~\ref{defn:tau_L}. We remark that this is not the final iteration scheme that will be used to prove the uniform-in-time bound; however, this is already enough to prove global well-posedness. The final iteration scheme used to obtain the uniform-in-time bound will be introduced in Section~\ref{section:decay-estimates-in-expectation}. We choose to organise the presentation in this way because the iteration scheme introduced here is the natural one to consider first. The final iteration scheme differs from it by an additional tweak. The motivation for this will become clear after the calculations in this subsection.

\begin{definition}
\label{defn:time_steps}
	Suppose at time $0$ we are given $L_0$ and $(\Asmoothpart_0, \Aroughpart_0, \phismoothpart_0, \phiroughpart_0)$ such that it is $L_0$ admissible. In the following, we will inductively define a sequence of stopping times $0 = t_0 < t_1 < \cdots < t_m < \cdots$ along with the time differences $\tau_m = t_m - t_{m-1}$, as well as their modified versions $t_m^\star$, $\tau_m^\star$, and data $L_m \in \N^+$. Moreover, on $(t_{m-1}, t_m]$, we will define functions
	\begin{equ}
		g_m \;, \quad B_m \;, \quad \YMHFsmooth_m \;, \quad \Sdet_m \;, \quad \higgsdet_m \;, \quad \linear[B_m] \;, \quad \philinear[B_m] \;,
	\end{equ}
	and decompositions
	\begin{equ}
		A = \Asmoothpart + \Aroughpart \;, \qquad \higgs = \phismoothpart + \phiroughpart \;.
	\end{equ}
	On $(t_{m-1}, t_m^\star]$, we will define their starred versions
	\begin{equ}
	\label{e:iteration_star_objects}
		g_m^\star\;, \quad B_m^\star \;, \quad B_m^{+, \star} \;, \quad \Sdet_m^\star \;, \quad \higgsdet_m^\star \;, \quad \linear[B^\star_m] \;, \quad \philinear[B^\star_m] \;.
	\end{equ}
	
	To begin, suppose the decompositions
	\begin{equ}
		A(t) = \Asmoothpart(t) + \Aroughpart(t) \;, \qquad
		\higgs(t) = \phismoothpart(t) + \phiroughpart(t)
	\end{equ}
	are defined up to stopping time $t_m$, and are $\fF_{t_m}$-measurable and $L_m$-admissible at time $t_m$. We define
	\begin{equ}
		\tau_{m+1}^{\star} \eqdef T^{(L_m)}\bigl(\Asmoothpart(t_m), \Aroughpart(t_m), \phismoothpart(t_m), \phiroughpart(t_m) \bigr) \;, \quad
		t_{m+1}^{\star} = t_m + \tau_{m+1}^{\star} \;.
	\end{equ}
	Let $g_{m+1}^\star$ be the gauge transformation given by Proposition~\ref{pr:admissible_timescale} for the decomposition
	$A(t_m) = \Asmoothpart(t_m) + \Aroughpart(t_m), \higgs(t_m) = \phismoothpart(t_m) + \phiroughpart(t_m)$, shifted to start at time $t_m$ so that $g_{m+1}^\star$ is defined on $[t_m, t_{m+1}^\star]$. 
	For $t \in (t_m, t_{m+1}^{\star}]$, let $B_{m+1}^{+, \star}(t)$ be the YM gradient flow starting from $B_{m+1}^{+, \star}(t_m) = \Asmoothpart(t_m)$, and $B^\star_{m+1}(t)$ be the YM gradient flow starting from $B^\star_{m+1}(t_m) = g_{m+1}^\star(t_m) \act A(t_m)$. 
	Then, we define $\Sdet^\star_{m+1}, \higgsdet^\star_{m+1}$ to be the solution to \eqref{eq:Sdet-equation}--\eqref{eq:higgsdet-equation}, shifted to start at time $t_m$, with $B$ replaced by $B^\star_{m+1}$ and initial data
	\begin{equ}
		\Sdet^\star_{m+1}(t_m) = 0 \;, \qquad \higgsdet^\star_{m+1}(t_m) = g_{m+1}^\star(t_m) \act \bigl(\phismoothpart(t_m) + \phiroughpart(t_m) \bigr) \;.
	\end{equ}
	We let $\linear[B^{\star}_{m+1}], \philinear[B^{\star}_{m+1}]$ solve
\begin{equs}
	\partial_t \linear[B^{\star}_{m+1}] &= \Delta_{B^{\star}_{m+1}} \linear[B^{\star}_{m+1}] + 2 \ad(F_{B^{\star}_{m+1}})\linear[B^{\star}_{m+1}] + \xi \;,& \quad \linear[B^{\star}_{m+1}](t_{m})& = 0  \;, \\
	\partial_t \philinear[B^{\star}_{m+1}] &= \Delta_{B^{\star}_{m+1}} \philinear[B^{\star}_{m+1}] + \zeta \;, &\quad \philinear[B^{\star}_{m+1}](t_{m})& = 0  \;.
\end{equs}	
	We emphasise that the previously defined quantities with $\star$ are all measurable with respect to $\fF_{t_m}$ (except for $\linear[B^{\star}_{m+1}], \philinear[B^{\star}_{m+1}]$). Next, we define quantities that depend on the noise after time $t_m$. Recall Definitions~\ref{defn:tau_L} and~\ref{defn:one_step_iteration_scheme}. 
	First, we define the stopping times
	\begin{equ}
		\tau_{m+1} \eqdef \tau^{(L_m)}\bigl(\Asmoothpart(t_m), \Aroughpart(t_m), \phismoothpart(t_m), \phiroughpart(t_m); \shift_{t_m}\xi, \shift_{t_m}\zeta\bigr) \;, \quad
		t_{m+1} \eqdef t_m + \tau_{m+1} \;.
	\end{equ}
	Here $\shift_t$ is the shift operator, defined by $\langle \shift_t \xi, f \rangle = \langle \xi, f(\place - t, \place) \rangle$ for any distribution $\xi$ and space-time test function $f$. Let $g_{m+1}, B_{m+1}, \YMHFsmooth_{m+1}, \Sdet_{m+1}, \higgsdet_{m+1}, \linear[B_{m+1}], \philinear[B_{m+1}]$ be the restrictions of the corresponding starred quantities to the time interval $(t_m, t_{m+1}]$.
	For $t \in (t_m, t_{m+1}]$, we define
	\begin{equ}
		\bigl(\Asmoothpart, \Aroughpart, \phismoothpart, \phiroughpart \bigr)(t)
		\eqdef \flow_{t-t_m}\bigl(\Asmoothpart(t_m), \Aroughpart(t_m), \phismoothpart(t_m), \phiroughpart(t_m); \shift_{t_m}\xi, \shift_{t_m}\zeta \bigr) \;.
	\end{equ}
	Finally, choose the smallest integer $L_{m+1}$ such that
	\begin{equ}
	\label{e:Lm_lower_bound}
		L_{m+1} \geq 1 \vee \bigl(\covariantenergy(\Asmoothpart(t_{m+1}), \phismoothpart(t_{m+1}))\bigr)^{\frac1{10}}
	\end{equ}
	and such that the decompositions
	\begin{equ}
		A(t_{m+1}) = \Asmoothpart(t_{m+1}) + \Aroughpart(t_{m+1}) \;, \qquad
		\higgs(t_{m+1}) = \phismoothpart(t_{m+1}) + \phiroughpart(t_{m+1})
	\end{equ}
	are $L_{m+1}$-admissible. 
	
	If one has $\limsup_{m\rightarrow\infty} t_m \neq \infty$, then, for any $t \geq \limsup_{m\rightarrow\infty} t_m$, we define
	$\bigl(\Asmoothpart, \Aroughpart, \phismoothpart, \phiroughpart\bigr)(t) = \skull$.
\end{definition}

\begin{figure}[t]
\centering
\begingroup
\setlength{\tabcolsep}{0pt}

\ifdefined\tikzsetnextfilename
  \tikzsetnextfilename{iteration-scheme}
\fi
\begingroup
\begin{tikzpicture}[
  x=0.76cm,
  y=0.82cm,
  line cap=round,
  line join=round
]
  \colorlet{futurecolor}{blue}
  \colorlet{gaugecolor}{red!65!black}

  \tikzset{
    iterline/.style={line width=0.55pt},
    iterdet/.style={iterline,dash pattern=on 4pt off 3pt},
    iterfuture/.style={iterline,futurecolor,dash pattern=on 4pt off 3pt},
    itertick/.style={line width=0.55pt},
    itergauge/.style={->,line width=0.65pt,gaugecolor},
    iterlab/.style={font=\small}
  }

  \def\xL{1.35}
  \def\xM{6.00}
  \def\xR{11.30}
  \def\xConnMerge{5.70}

  \def\yTime{4.80}
  \def\yConnB{3.62}
  \def\yConnMerge{3.31}
  \def\yConnS{2.96}
  \def\yHiggs{2.1}
  \def\yRoughA{1}
  \def\yRoughPhi{0.25}
  \def\yGaugeTop{-0.20}
  \def\yGaugeLabel{-0.64}
  \def\yGaugeBottom{-1.03}
  \def\yLowerB{-1.45}
  \def\yLowerS{-2.12}
  \def\yLowerPhi{-2.80}
  \def\yLowerRoughA{-3.55}
  \def\yLowerRoughPhi{-4.30}
  
  \def\yLegendDet{3.25}
  \def\yLegendLine{2.65}
  \def\yLegendFuture{2.05}

  \draw[iterline] (\xL,\yTime) -- (\xR,\yTime);

  \foreach \X in {\xL,\xM,\xR}
    \draw[itertick] (\X,\yTime+0.12) -- (\X,\yTime-0.12);
    
  \draw[itertick] (\xR-1.18,\yTime+0.12) -- (\xR-1.18,\yTime-0.12);

  \node[iterlab,anchor=south] at (\xL,\yTime+0.18) {$t_{m-1}$};
  \node[iterlab,anchor=south] at (\xM,\yTime+0.18) {$t_m$ ($t_m^\star$)};
  \node[iterlab,anchor=south] at (\xR-1.18,\yTime+0.18) {$t_{m+1}$};
  \node[iterlab,anchor=south] at (\xR,\yTime+0.18) {$t_{m+1}^{\star}$};

  \draw[iterdet] (\xL,\yConnB) -- (\xConnMerge,\yConnB) -- (\xM,\yConnMerge);
  \draw[iterdet] (\xL,\yConnS) -- (\xConnMerge,\yConnS) -- (\xM,\yConnMerge);
  \draw[iterfuture] (\xM,\yConnMerge) -- (\xR,\yConnMerge);

  \foreach \Y in {\yConnB,\yConnS}
    \draw[itertick] (\xL,\Y+0.12) -- (\xL,\Y-0.12);
  \draw[itertick] (\xM,\yConnMerge+0.12) -- (\xM,\yConnMerge-0.12);
  \draw[itertick,futurecolor] (\xR,\yConnMerge+0.12) -- (\xR,\yConnMerge-0.12);

  \node[iterlab,anchor=east] at (1.12,\yConnB) {$A(t_{m-1})$};
  \node[iterlab,anchor=east] at (1.12,\yConnS) {$0$};
  \node[iterlab,anchor=south] at (3.70,\yConnB) {$B_m$};
  \node[iterlab,anchor=south] at (3.70,\yConnS) {$S_m$};
  \node[iterlab,anchor=south west] at (6,\yConnMerge+0.17) {$A^{+}(t_m)$};

  \node[iterlab,futurecolor,anchor=south] at (8.85,\yConnMerge+0.06) {$B_{m+1}^{+,\star}$}; 

  \draw[iterdet] (\xL,\yHiggs) -- (\xM,\yHiggs);

  \foreach \X in {\xL,\xM}
    \draw[itertick] (\X,\yHiggs+0.12) -- (\X,\yHiggs-0.12);

  \node[iterlab,anchor=east] at (1.12,\yHiggs) {$\Phi(t_{m-1})$};
  \node[iterlab,anchor=south] at (3.70,\yHiggs ) {$\higgsdet_m$};
  \node[iterlab,anchor=south] at (\xM,\yHiggs+0.05) {$\phismoothpart(t_m)$};

  \draw[iterline] (\xL,\yRoughA) -- (\xM,\yRoughA);

  \foreach \X in {\xL,\xM}
    \draw[itertick] (\X,\yRoughA+0.12) -- (\X,\yRoughA-0.12);

  \node[iterlab,anchor=east] at (1.12,\yRoughA) {$0$};
  \node[iterlab,anchor=south] at (\xM,\yRoughA+0.06) {$A^{-}(t_m)$};

  \draw[iterline] (\xL,\yRoughPhi) -- (\xM,\yRoughPhi);

  \foreach \X in {\xL,\xM}
    \draw[itertick] (\X,\yRoughPhi+0.12) -- (\X,\yRoughPhi-0.12);

  \node[iterlab,anchor=east] at (1.12,\yRoughPhi) {$0$};
  \node[iterlab,anchor=south] at (\xM,\yRoughPhi+0.06) {$\phiroughpart(t_m)$};

  \draw[itergauge] (\xM,\yGaugeTop) .. controls (6.6,\yGaugeTop-0.25) and (6.6,\yGaugeBottom+0.25) .. (\xM,\yGaugeBottom);
  \node[iterlab,gaugecolor,anchor=east] at (5.80,\yGaugeLabel) {$g_{m+1}^{\star}(t_m)$};

  \draw[iterdet] (\xM,\yLowerB) -- (\xR,\yLowerB);
  \draw[iterdet] (\xM,\yLowerS) -- (\xR,\yLowerS);

  \foreach \Y in {\yLowerB,\yLowerS}
    \draw[itertick] (\xM,\Y+0.12) -- (\xM,\Y-0.12);
  \foreach \Y in {\yLowerB,\yLowerS}
    \draw[itertick] (\xR,\Y+0.12) -- (\xR,\Y-0.12);

  \node[iterlab,gaugecolor,anchor=east,xshift=-3pt] at (\xM,\yLowerB) {$A(t_m)$};
  \node[iterlab,anchor=east,xshift=-3pt] at (\xM,\yLowerS) {$0$};
  \node[iterlab,anchor=south] at (8.85,\yLowerB - 0.05) {$B_{m+1}^{\star}$};
  \node[iterlab,anchor=south] at (8.85,\yLowerS-0.05) {$S_{m+1}^{\star}$};

  \draw[iterdet] (\xM,\yLowerPhi) -- (\xR,\yLowerPhi);

  \draw[itertick] (\xM,\yLowerPhi+0.12) -- (\xM,\yLowerPhi-0.12);
  \draw[itertick] (\xR,\yLowerPhi+0.12) -- (\xR,\yLowerPhi-0.12);

  \node[iterlab,gaugecolor,anchor=east,xshift=-3pt] at (\xM,\yLowerPhi) {$\Phi(t_m)$};
  \node[iterlab,anchor=south] at (8.85,\yLowerPhi-0.1) {$\higgsdet_{m+1}^{\star}$};

  \draw[iterline] (\xM,\yLowerRoughA) -- (\xR-1.18,\yLowerRoughA);
  \draw[iterline] (\xM,\yLowerRoughPhi) -- (\xR-1.18,\yLowerRoughPhi);

  \foreach \Y in {\yLowerRoughA,\yLowerRoughPhi}
    \draw[itertick] (\xM,\Y+0.12) -- (\xM,\Y-0.12);
  \foreach \Y in {\yLowerRoughA,\yLowerRoughPhi}
    \draw[itertick] (\xR-1.18,\Y+0.12) -- (\xR-1.18,\Y-0.12);

  \node[iterlab,anchor=east,xshift=-3pt] at (\xM,\yLowerRoughA) {$0$};
  \node[iterlab,anchor=east,xshift=-3pt] at (\xM,\yLowerRoughPhi) {$0$};
  \node[iterlab,anchor=south] at (8.85,\yLowerRoughA-0.05) {$\Aroughpart$};
  \node[iterlab,anchor=south] at (8.85,\yLowerRoughPhi-0.05) {$\phiroughpart$};
  
  \draw[iterdet] (12.05,\yLegendDet) -- (13.15,\yLegendDet);
  \node[iterlab,anchor=west] at (13.40,\yLegendDet) {deterministic flow};

  \draw[iterline] (12.05,\yLegendLine) -- (13.15,\yLegendLine);
  \node[iterlab,anchor=west] at (13.40,\yLegendLine) {stochastic / rough part};

  \draw[iterfuture] (12.05,\yLegendFuture) -- (13.15,\yLegendFuture);
  \node[iterlab,anchor=west] at (13.40,\yLegendFuture) {deterministic flow for next step};

\end{tikzpicture}
\endgroup

\par
\endgroup
\caption{Illustration of Definition~\ref{defn:time_steps}. Here for simplicity we assumed $t_m = t_m^\star$.}
\label{fig:iteration_scheme}
\end{figure}
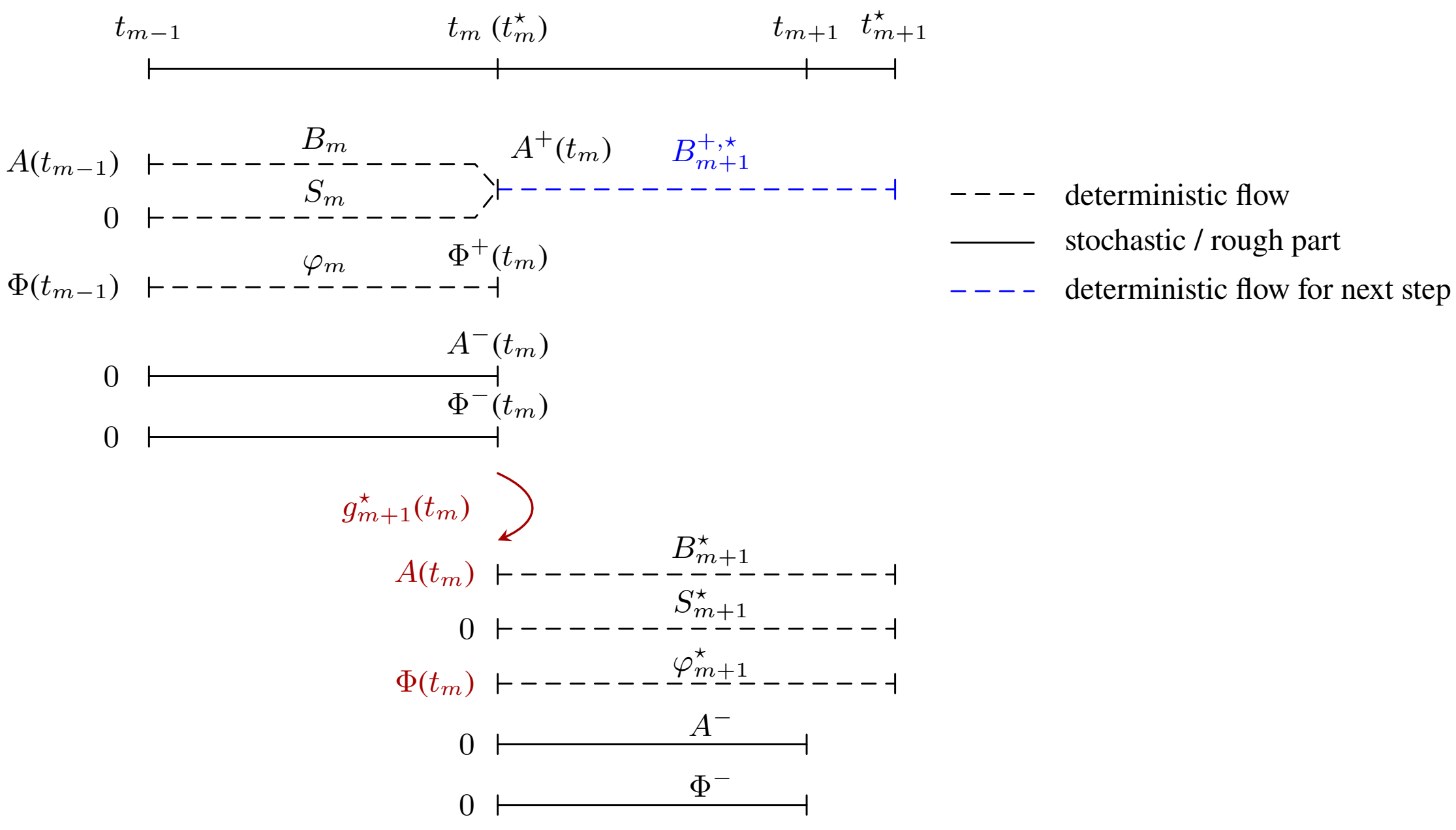


\begin{remark}
\label{rmk:iteration_measurability}
	From the definition it follows that:
	\begin{enumerate}
		\item For all $m\geq 1$, $t_m$ is a stopping time. The starred quantities $g_{m+1}^\star, B_{m+1}^\star, B_{m+1}^{+, \star}, \Sdet_{m+1}^\star, \higgsdet_{m+1}^\star$, $\tau^\star_{m+1}$, $t^\star_{m+1}$, $L^\star_{m+1}$ and $L_m$ are $\fF_{t_{m}}$-measurable, while the non-starred versions are only $\fF_{t_{m+1}}$-measurable.
		\item With high probability, $t_m = t_m^\star$ and the starred quantities and their non-starred versions are the same. Thus on a first reading, the reader may pretend that the equality always holds.
		\item For any $t > 0$, $(\Asmoothpart, \Aroughpart, \phismoothpart, \phiroughpart)(t)$ is measurable with respect to $\fF_{t}$. 
\end{enumerate}
\end{remark}
\begin{remark}
	At first, we will focus only on the sizes of the solution at the stopping times $t_m$ (i.e.\ $\covariantenergy\bigl(\Asmoothpart(t_m), \phismoothpart(t_m)\bigr)$), and we will aim to obtain bounds in expectation that are uniform in $m$. The case of a general deterministic time $t$ will only be addressed at the end of the argument. However, since each time interval $[t_m, t_{m+1}]$ is chosen short enough for the covariant local theory to apply, it is natural to expect that once the bounds at the endpoints are obtained, the corresponding bounds at arbitrary deterministic times follow by perturbative arguments.
\end{remark}

\begin{lemma}
\label{le:iteration_basic_bounds}
	Let $m \geq 0$. In the setting of Definition~\ref{defn:time_steps}, the shifted versions of \eqref{e:iteration_detarministic_bounds_connection_assumption}, \eqref{e:phi_0_smallness_global}, Assumption~\ref{local:assumption-B}, and Assumption~\ref{local:assumption-higgs-initial} hold on the interval $(t_m, t_{m+1}^\star]$. Furthermore, we record the following basic inequalities. For any $m \geq 1$, we have
	\begin{equ}
\label{e:iteration_F_bounds1}
	\sup_{s \in (0, \tau_{m+1}^\star]} s^{1 + \fa} \Bigl( \| F_{B_{m+1}^{+, \star}}(t_m + s) \|_{L^\infty} + \| F_{B_{m+1}^\star}(t_m + s) \|_{L^\infty} \Bigr) \lesssim 1 \;, 
	\end{equ}
	and
	\begin{equs}
	\sup_{s \in [0, \tau_{m+1}^\star]} \| F_{B_{m+1}^{+, \star}}(t_m + s) \|_{L^2}
	&\lesssim \bigl((T^{(1)}\wedge T^{(2)})(\Asmoothpart(t_m), \phismoothpart(t_m))\bigr)^{-1 - \fa}\\
\label{e:iteration_F_bounds2}
	&\leq \bigl( \tau_{m+1}^\star \bigr)^{-1 - \fa} \;, \\
	\| \covd_{B_{m+1}^{+, \star}} F_{B_{m+1}^{+, \star}} \|_{L^2_{t, x}([t_m, t_m + \tau_{m+1}^\star])}
	&\lesssim \bigl((T^{(1)}\wedge T^{(2)})(\Asmoothpart(t_m), \phismoothpart(t_m))\bigr)^{-\frac12 - \lossdel} \wedge \| F_{\Asmoothpart (t_m)} \|_{L^2} \\ 
\label{e:iteration_F_endtime_lowerbound}
	&\leq \bigl(\tau_{m+1}^\star \bigr)^{-\frac12 - \lossdel} \wedge \| F_{\Asmoothpart (t_m)} \|_{L^2} \;, \\
	\| \linear[B_m](t_m) \|_{\covHolder{\Asmoothpart(t_m),F}^{-\kappa}}
	+ \| \Aroughpart(t_m) \|_{\covHolder{\Asmoothpart(t_m),F}^{-\kappa}}
	&\lesssim L_m^{\fa} \;, \quad \| \Aroughpart(t_m) - \linear[B_m](t_m) \|_{\covHolder{\Asmoothpart(t_m),F}^{0}} \lesssim \tau_m^{100 \fa} \;, \\
\label{e:Aroughpart_crude_bound}\\
	\| \philinear[B_m](t_m) \|_{\covHolder{\Asmoothpart(t_m)}^{-\kappa}}
	+ \| \phiroughpart(t_m) \|_{\covHolder{\Asmoothpart(t_m)}^{-\kappa}}
	&\lesssim L_m^{\fa} \;, \quad \| \phiroughpart(t_m) - \philinear[B_m](t_m) \|_{\covHolder{\Asmoothpart(t_m)}^{0}} \lesssim \tau_m^{100 \fa} \;. \\
\label{e:phiroughpart_crude_bound}  
\end{equs}
	Moreover, we have the following relation between the lengths of two time intervals:
	\begin{equs}
	\label{e:two_timesteps_relation}
		\tau_m &\lesssim \stepupperbound \wedge \logtimestep \wedge \| F_{\Asmoothpart(t_m) } \|_{L^2}^{-\frac12 }
		\wedge \| \phismoothpart(t_m) \|_{L^{\infty}}^{- \frac{1}{\frac{1}{\gamma_0} + \frac1r + \lossdel} } \\
		&\lesssim \bigl( (T^{(1)} \wedge T^{(2)})(\Asmoothpart(t_m), \phismoothpart(t_m) )  \bigr)^{\frac{1}{10} } \;.
	\end{equs}
\end{lemma}
\begin{proof}
	The first claim follows from the definition of our stopping time and discussions in the previous subsection (see e.g. Proposition~\ref{pr:admissible_timescale} and Definition~\ref{defn:L_admissible}). The only extra thing to check is that in \eqref{e:iteration_detarministic_bounds_connection_assumption}, two derivatives of $\Aroughpart$ are required, whereas the definition of $\covHolder{\Asmoothpart, F}^{-\kappa}$ norm involves only one derivative. This, however, directly follows from our heat kernel estimate \eqref{eq:qB-ptwise-bounds}. 
	
	The bound \eqref{e:iteration_F_bounds1} 
	follows immediately from the definition of the stopping time (Definition~\ref{defn:tau}), the first claim (which in particular implies that the $L_t^1 L_x^\infty([t_m, t_{m+1}^\star])$ norms of the curvatures are bounded by $(\tau_{m+1}^\star)^{-\lossdel}$) and \eqref{eq:curvature-pointwise-in-time-bound-Fnorm}. 
	
	The bound \eqref{e:iteration_F_bounds2} follows from monotonicity of $\| F_{B_{m+1}^{+, \star}}(t_m + s) \|_{L^2}$ and the definition of $T^{(1)}$. 
	
	The bound \eqref{e:iteration_F_endtime_lowerbound} follows from \eqref{e:T2_decrease_consequence} and \eqref{eq:parabolic-bochner}. 
	For \eqref{e:two_timesteps_relation}, by \eqref{e:iteration_deterministic_bounds_rmk}, we get that 
	\begin{equ}
		 \| F_{\Asmoothpart(t_m) } \|_{L^2}^2 \leq \covariantenergy(\Asmoothpart(t_{m-1}), \phismoothpart(t_{m-1}) ) + C \tau_m^{-\frac74} \;.
	\end{equ}
	From the definition of $T^{(1)}$ and $\covariantenergy$ in Definition~\ref{defn:tau} and \eqref{main:eq-covariant-energy}, it follows that 
	\begin{equ}
		\covariantenergy(\Asmoothpart(t_{m-1}), \phismoothpart(t_{m-1}) ) \lesssim \tau_m^{-4} \;.
	\end{equ} 
	To prove $\tau_m \lesssim  \| \phismoothpart(t_m) \|_{L^{\infty}}^{- \frac{1}{\frac{1}{\gamma_0} + \frac1r + \lossdel} }$, we apply \eqref{eq:higgsdet-L-infty-estimates}, along with Lemma~\ref{le:B-YMHFsmooth-heat-comparison}, and the fact that $L_{m-1} \leq \tau_m^{-1}$ to get 
	\begin{equs}
	\| \phismoothpart(t_m) \|_{L^\infty}
	&\lesssim \tau_m^{-\frac{\kappa}{2} - \lossdel} 
	\Bigl(\tau_m^{\frac{\kappa}{2}-\frac1r}\|\phismoothpart(t_{m-1})\|_{L^r}
	+ \|g_m(t_{m-1})\act\phiroughpart(t_{m-1})\|_{\covHolder{B_m(t_{m-1}),\tau_m}^{-\kappa}}\Bigr) \\
	&\lesssim \tau_m^{-\frac{\kappa}{2}-\lossdel}
	\Bigl(\tau_m^{\frac{\kappa}{2}-\frac1r-\frac1{\gamma_0}}
	+ \tau_m^{-\lossdel}\|\phiroughpart(t_{m-1})\|_{\covHolder{\Asmoothpart(t_{m-1})}^{-\kappa}}\Bigr) \\
	&\lesssim \tau_m^{-\frac{\kappa}{2}-\lossdel}
	\Bigl(\tau_m^{\frac{\kappa}{2}-\frac1r-\frac1{\gamma_0}}
	+ \tau_m^{-\lossdel}L_{m-1}^{\fa}\Bigr)
	\lesssim \tau_m^{-\frac1{\gamma_0}-\frac1r-\lossdel} \;.
	\end{equs}
	
	Combining the discussions above proves the first inequality. 
	The second inequality follows from Remark~\ref{rmk:tau_lower_bound}.
	
	Now we turn to proving \eqref{e:Aroughpart_crude_bound} and \eqref{e:phiroughpart_crude_bound}. 	
	For \eqref{e:Aroughpart_crude_bound}, it follows from Proposition~\ref{pr:remainder_estimate} with $\beta=2\lossdel$ and \eqref{e:two_timesteps_relation} that
	\begin{equs}
		 \| \Aroughpart(t_m) - \linear[B_m](t_m) \|_{\covHolder{\Asmoothpart(t_m),F}^{0}} \lesssim \tau_m^{\frac{1-\lambda}{2}-\frac{1}{\gamma_0}-\lambda_2}
		\Bigl(\tau_m^{-\lossdel}+\bigl((T^{(1)}\wedge T^{(2)})(\Asmoothpart(t_m), \phismoothpart(t_m))\bigr)^{-\lossdel}\Bigr) \lesssim \tau_m^{100 \fa} \;. 
	\end{equs}
	This proves the second inequality in \eqref{e:Aroughpart_crude_bound}.
	The first inequality in \eqref{e:Aroughpart_crude_bound} follows from the second and the definition of $L_m$. 
	The bound \eqref{e:phiroughpart_crude_bound} follows similarly.
\end{proof}

By Lemma~\ref{le:iteration_basic_bounds} and the definition of the time steps, 
the assumptions in Proposition~\ref{pr:iteration_detarministic_bounds} are satisfied on $[t_m, t_{m+1}]$. We now apply that proposition to get the following.
\begin{lemma}
\label{le:iteration_deterministic_bounds_mid_iteration}
	There is some constant $C$ such that in the setting of Definition~\ref{defn:time_steps}, for any $m \geq 1$ we have
	\minilab{e:iteration_deterministic_bounds_mid_iteration}
	\begin{equs}
		\covariantenergy&(\Asmoothpart(t_{m+1}), \phismoothpart(t_{m+1})) \\
		&\leq \covariantenergy(\Asmoothpart(t_{m}), \phismoothpart(t_{m})) - \frac12 \Dec_{m+1} + C\bigl(\tau_{m+1} \bigr)^{-1 + \kappa} \\
		\label{e:iteration_deterministic_bounds_mid_iteration1}
		&\quad+  \Big \| \curl_{\YMHFsmooth_{m+1}(t_{m+1})} e^{\tau_{m+1} \Delta^F_{\Asmoothpart(t_{m})}} \Aroughpart(t_m) \Big\|_{L^2}^2 \\
		\label{e:iteration_deterministic_bounds_mid_iteration2}
		&\quad+ 2 \int_{\T^2} \Big \langle F_{\YMHFsmooth_{m+1}}^{12}(t_{m+1}), \curl_{\YMHFsmooth_{m+1}(t_{m+1})} e^{\tau_{m+1} \Delta^F_{\Asmoothpart(t_{m})}} \Aroughpart(t_m) \Big \rangle_{\fg} \, \dx \\
		\label{e:iteration_deterministic_bounds_mid_iteration3}
		&\quad + 2 \int_{\T^2} \Big \langle F_{\YMHFsmooth_{m+1}}^{12}(t_{m+1}), \Big[( e^{\tau_{m+1} \Delta^F_{\Asmoothpart(t_{m})}} \Aroughpart(t_m))^1, ( e^{\tau_{m+1} \Delta^F_{\Asmoothpart(t_{m})}} \Aroughpart(t_m))^2 \Big]  \Big \rangle_{\fg} \, \dx\qquad  \\
		&\quad+ C \| \phismoothpart(t_m) \|_{L^4}^{\gamma_1 - 4} \Big| \int_{\T^2} |e^{\tau_{m+1} \Delta_{\Asmoothpart(t_m)}} \phismoothpart(t_m) |^2 \\
		\label{e:iteration_deterministic_bounds_mid_iteration4}
		& \qquad \qquad \qquad \qquad \qquad \cdot \langle e^{\tau_{m+1} \Delta_{\Asmoothpart(t_m)}} \phismoothpart(t_m), e^{\tau_{m+1} \Delta_{\Asmoothpart(t_m)}} \phiroughpart(t_m) \rangle_{\V} \,  \dx \Big| \;,
	\end{equs}
	where
	\begin{equs}
	\label{e:iteration_deterministic_bounds_iteration_Dec1}
	\Dec_{m+1} &\eqdef \Dec_{\tau_{m+1}}(\Asmoothpart(t_m), \phismoothpart(t_m))  \\
	&= \Big(\|F_{\YMHFsmooth_{m+1}(t_m)} \|_{L^2}^2 - \| F_{\YMHFsmooth_{m+1}(t_{m+1})} \|_{L^2}^2 \Big) + C_{\eqref{e:shiggsdet_decrease}, \gamma_2, r} \one_{\| \phismoothpart(t_m) \|_{L^r} \geq 1} \| \phismoothpart(t_m) \|_{L^r}^{\gamma_2 + 2} \tau_{m+1} \\
	& + \one_{\| \phismoothpart(t_m) \|_{L^4} \geq 1} \Bigl( \gamma_1 \| \shiggsdet_{m+1}(t_{m+1}) \|_{L_x^4}^{\gamma_1 - 4} \| |\covd_{B_{m+1} + \Sdet_{m+1}} \shiggsdet_{m+1}| |\shiggsdet_{m+1}| \|_{L^2_{t,x}((t_m,t_{m+1}])}^2 \\
	\label{e:iteration_deterministic_bounds_iteration_Dec2}
	&\qquad \vee C_{\eqref{e:shiggsdet_decrease}, \gamma_1, 4} \| \phismoothpart(t_m) \|_{L^4}^{\gamma_1 + 2} \tau_{m+1} \Bigr) \;. 
	\end{equs}
\end{lemma}

	As indicated in Proposition~\ref{pr:remainder_estimate}, $\Aroughpart(t_m)$ and $\phiroughpart(t_m)$ are essentially given by $\linear[B_m](t_m)$ and $\philinear[B_m](t_m)$. In particular, upon replacing $\Aroughpart(t_m)$ and $\phiroughpart(t_m)$ by $\linear[B_m](t_m)$ and $\philinear[B_m](t_m)$ in \eqref{e:iteration_deterministic_bounds_mid_iteration1}--\eqref{e:iteration_deterministic_bounds_mid_iteration4}, the resulting expectations can be computed explicitly. 
	The following proposition bounds the errors in the replacements.
\begin{proposition}
\label{pr:expansion_iteration}
	There is some constant $C$, such that when replacing $\Aroughpart(t_m)$ and $\phiroughpart(t_m)$ by $\linear[B_m](t_m)$ and $\philinear[B_m](t_m)$ in \eqref{e:iteration_deterministic_bounds_mid_iteration1}--\eqref{e:iteration_deterministic_bounds_mid_iteration4}, the total errors can be bounded by $\frac14 \Dec_{m+1} + C\bigl(\tau_{m+1} \bigr)^{-1 + \kappa}$. 
In particular, one has the bound  
	\minilab{e:iteration_deterministic_bounds_iteration}
	\begin{equs}
		\covariantenergy&(\Asmoothpart(t_{m+1}), \phismoothpart(t_{m+1})) \\
		&\leq \covariantenergy(\Asmoothpart(t_{m}), \phismoothpart(t_{m})) - \frac14 \Dec_{m+1} + C\bigl(\tau_{m+1} \bigr)^{-1 + \kappa} \\
		\label{e:iteration_deterministic_bounds_iteration1}
		&\quad+  \Big \| \curl_{\YMHFsmooth_{m+1}(t_{m+1})} e^{\tau_{m+1} \Delta^F_{\Asmoothpart(t_{m})}} \linear[B_m](t_m) \Big\|_{L^2}^2 \\
		\label{e:iteration_deterministic_bounds_iteration2}
		&\quad+ 2 \int_{\T^2} \Big \langle F_{\YMHFsmooth_{m+1}}^{12}(t_{m+1}), \curl_{\YMHFsmooth_{m+1}(t_{m+1})} e^{\tau_{m+1} \Delta^F_{\Asmoothpart(t_{m})}} \linear[B_m](t_m) \Big \rangle_{\fg} \, \dx \\
		\label{e:iteration_deterministic_bounds_iteration3}
		&\quad + 2 \int_{\T^2} \Big \langle F_{\YMHFsmooth_{m+1}}^{12}(t_{m+1}), \Big[( e^{\tau_{m+1} \Delta^F_{\Asmoothpart(t_{m})}} \linear[B_m](t_m))^1, ( e^{\tau_{m+1} \Delta^F_{\Asmoothpart(t_{m})}} \linear[B_m](t_m))^2 \Big]  \Big \rangle_{\fg} \, \dx\qquad \\
		&\quad+ C \| \phismoothpart(t_m) \|_{L^4}^{\gamma_1 - 4} \Big| \int_{\T^2} |e^{\tau_{m+1} \Delta_{\Asmoothpart(t_m)}} \phismoothpart(t_m) |^2 \\
		\label{e:iteration_deterministic_bounds_iteration4}
		& \qquad \qquad \qquad \qquad \qquad \cdot \langle e^{\tau_{m+1} \Delta_{\Asmoothpart(t_m)}} \phismoothpart(t_m), e^{\tau_{m+1} \Delta_{\Asmoothpart(t_m)}} \philinear[B_m](t_m) \rangle_{\V} \,  \dx \Big| \;,
	\end{equs}
	We will let $\itobjects_{m+1}$ denote the sum of the expressions \eqref{e:iteration_deterministic_bounds_iteration1}--\eqref{e:iteration_deterministic_bounds_iteration4}. 
\end{proposition}
\begin{remark}
	In Proposition~\ref{pr:expectation_one_step} later, we will estimate the expectations of \eqref{e:iteration_deterministic_bounds_iteration1}--\eqref{e:iteration_deterministic_bounds_iteration4} essentially by quantities of order $\tau_{m+1}^{-1}$, which is only slightly larger than the $\tau_{m+1}^{-1 + \kappa}$ term on the right-hand side of the inequality in Proposition~\ref{pr:expansion_iteration}. It is therefore quite believable that the $\tau_{m+1}^{-1+\kappa}$ term can absorb all errors when replacing $\Aroughpart(t_m), \phiroughpart(t_m)$ by the corresponding linear objects (ultimately because the errors arising when performing the replacements can be estimated by tiny powers of the time interval). The one slight caveat is that the extra $\frac14 \Dec_{m+1}$ term is needed to handle the error during the replacement of \eqref{e:iteration_deterministic_bounds_iteration4}.
\end{remark}

\begin{proof}
	We compare \eqref{e:iteration_deterministic_bounds_mid_iteration1}--\eqref{e:iteration_deterministic_bounds_mid_iteration4} with \eqref{e:iteration_deterministic_bounds_iteration1}--\eqref{e:iteration_deterministic_bounds_iteration4} one by one. 
		
	\proofstep{Error when replacing \eqref{e:iteration_deterministic_bounds_mid_iteration1} by \eqref{e:iteration_deterministic_bounds_iteration1}}: 
	Applying Proposition~\ref{pr:remainder_estimate} with $\beta = \frac{1 - \lambda - \frac{2}{\gamma_0} - 2\lambda_2}{10} $, the error is bounded by 
	\begin{equs}
\label{e:iteration_deterministic_bounds_mid_iteration1_mid}
	\,& C \tau_{m+1}^{- \frac{\kappa + 1}{2} } \Bigl( \| \linear[B_m] (t_m) \|_{\covHolder{\Asmoothpart(t_m), F}^{- \kappa}} + \| \Aroughpart (t_m) \|_{\covHolder{\Asmoothpart(t_m), F}^{- \kappa}} \Bigr)  \\
	&\qquad \qquad \qquad \cdot \tau_{m+1}^{ \frac{\beta - 1}{2} - \lossdel}  \tau_m^{\frac{1 - \lambda}{2} - \frac{1}{\gamma_0} - \lambda_2} \Bigl(\tau_m^{- \frac{\beta}{2} } + \Bigl( (T^{(1)}\wedge T^{(2)}) \bigl(\Asmoothpart(t_m), \phismoothpart(t_m) \bigr)  \Bigr)^{- \frac{\beta}{2} } \Bigr) \;.
\end{equs}
	Using \eqref{e:Aroughpart_crude_bound} and \eqref{e:two_timesteps_relation}, we get that \eqref{e:iteration_deterministic_bounds_mid_iteration1_mid} is bounded by
	\begin{equ}
		C \tau_{m+1}^{- \fa - \frac{\kappa}{2} - 1 - \lossdel + \frac{\beta}{2} } \bigl(\tau_m^{5\beta - \frac{\beta}{2} } + \tau_m^{5\beta - 5 \beta}  \bigr) \lesssim \tau_{m+1}^{- \fa - \frac{\kappa}{2} - 1 - \lossdel + \frac{\beta}{2} } \leq \tau_{m+1}^{-1 + \kappa}\;,
	\end{equ}
	since $\gamma_0 - 2 \gg \fa, \kappa, \lambda, \lambda_2, \lossdel$.
	
	\proofstep{Error when replacing \eqref{e:iteration_deterministic_bounds_mid_iteration2} by \eqref{e:iteration_deterministic_bounds_iteration2}}: By integration by parts combined with Lemmas~\ref{le:DF_bound} and~\ref{le:iteration_basic_bounds}, the error is bounded by
	\begin{equs}
		C \| \bigl(\covd_{\YMHFsmooth_{m+1}} F_{\YMHFsmooth_{m+1}} \bigr) (t_{m+1}) \|_{L^2} \tau_m^{ 100 \fa} 
		&\lesssim \tau_{m+1}^{-\lossdel - \frac12} \| \covd_{\YMHFsmooth_{m+1}} F_{\YMHFsmooth_{m+1}} \|_{L^2_{t, x}([t_m, t_{m+1}]) } \tau_m^{100 \fa } \\ 
		&\lesssim 
		\tau_{m+1}^{-\lossdel - \frac12} \cdot \Bigl( (T^{(1)}\wedge T^{(2)}) \bigl(\Asmoothpart(t_m), \phismoothpart(t_m) \bigr)  \Bigr)^{-\frac12 - \lossdel + 10 \fa } \leq \tau_{m+1}^{-1 + \kappa} \;,
	\end{equs}
	which gives the desired bound.
		
	\proofstep{Error when replacing \eqref{e:iteration_deterministic_bounds_mid_iteration3} by \eqref{e:iteration_deterministic_bounds_iteration3}}: Similar to above, by Lemma~\ref{le:iteration_basic_bounds}, the error is bounded by (using that $\fa \gg \kappa$)
	\begin{equs}
		\,& C \|F_{\YMHFsmooth_{m+1}}(t_{m+1})\|_{L^2}  \cdot \tau_{m+1}^{- \frac{\kappa}{2} } \Bigl( \| \linear[B_m] (t_m) \|_{\covHolder{\Asmoothpart(t_m), F}^{- \kappa}} + \| \Aroughpart (t_m) \|_{\covHolder{\Asmoothpart(t_m), F}^{- \kappa}} \Bigr) \cdot \tau_m^{100 \fa} \\
		& \lesssim \Bigl(  (T^{(1)}\wedge T^{(2)}) \bigl(\Asmoothpart(t_m), \phismoothpart(t_m) \bigr)   \Bigr)^{-1 - \fa} \cdot \tau_{m+1}^{- \frac{\kappa}{2} - \fa } \cdot \Bigl(  (T^{(1)}\wedge T^{(2)}) \bigl(\Asmoothpart(t_m), \phismoothpart(t_m) \bigr)   \Bigr)^{10 \fa }  \\
		&\leq \tau_{m+1}^{-1 + \kappa} \;.
	\end{equs}
	
	\proofstep{Error when replacing \eqref{e:iteration_deterministic_bounds_mid_iteration4} by \eqref{e:iteration_deterministic_bounds_iteration4}}: 
	By Lemma~\ref{le:iteration_basic_bounds} and Young's inequality, the error is bounded by 
		\begin{equs}
		\,&C \| \phismoothpart (t_m) \|_{L^4}^{\gamma_1 - 1 } \tau_m^{\frac{1 - \lambda}{2} - \frac{1}{\gamma_0} - \lambda_2 - \lossdel } \\
		&\lesssim C \| \phismoothpart(t_m) \|_{L^4}^{\gamma_1 - 1 - \frac{\frac{1-\lambda}{2} - \frac{1}{\gamma_0} - \lambda_2 - \lossdel}{\frac{1}{\gamma_0} + \frac1r + \lossdel}} \\
		&= C \| \phismoothpart(t_m) \|_{L^4}^{\gamma_1 - 1 - c_\star(\gamma_0 - 2) + \bigO(\lambda_2)} \\
		&= C \bigl(\| \phismoothpart(t_m) \|_{L^4}^{\gamma_1 + 2} \tau_{m+1} \bigr)^{\frac{\gamma_1 - 1 - c_\star(\gamma_0 - 2) + \bigO(\lambda_2)}{\gamma_1 + 2}}
			\tau_{m+1}^{-\frac{\gamma_1 - 1 - c_\star(\gamma_0 - 2) + \bigO(\lambda_2)}{\gamma_1 + 2}} \\
			&\leq \frac{1}{4} C_{\eqref{e:shiggsdet_decrease}, \gamma_1, 4} \| \phismoothpart(t_m) \|_{L^4}^{\gamma_1 + 2} \tau_{m+1}
			+ C \tau_{m+1}^{\eta^\star} \;.
		\end{equs}
		Here, $c_\star = \frac{1}{2 \gamma_0 (\frac{1}{\gamma_0} + \frac1r + \lossdel)} \approx \frac12 $, and
		\begin{equs}
			\eta^\star &= -\frac{\gamma_1 - 1 - c_\star(\gamma_0 - 2) + \bigO(\lambda_2)}{\gamma_1 + 2 - \bigl(\gamma_1 - 1 - c_\star(\gamma_0 - 2) + \bigO(\lambda_2) \bigr)}  \\
			&= -1 + \frac{2c_\star(\gamma_0 - 2) - (\gamma_1 - 4) + \bigO(\lambda_2)}{3 + c_\star(\gamma_0 - 2) + \bigO(\lambda_2)} 
			> -1 + \kappa \;,
		\end{equs}
		where we used $\lambda_2 \ll \gamma_1 - 4 \ll \gamma_0 - 2$. This gives the desired result by Proposition~\ref{pr:shiggsdet_decrease}. 
\end{proof}

Now we plug in the bounds obtained in Section~\ref{sec:improved_stochastic_estimates} into Proposition~\ref{pr:expansion_iteration} to estimate the expectation of $\itobjects_{m+1}$. The specific structure of $\itobjects_{m+1}$ as well as the independence between the space-time white noise and the deterministic evolutions lead to improved stochastic estimates, which are discussed in more detail at the beginning of Section~\ref{sec:improved_stochastic_estimates}. 
When applying the results in Section~\ref{sec:improved_stochastic_estimates}, one has to be a bit careful, as the bounds in Section~\ref{sec:improved_stochastic_estimates} are only valid when $\tstop, t$ are deterministic times, while $t_{m}, t_{m+1}$ are random times even when conditioning on $\fF_{t_{m-1}}$. 
We therefore first restrict to a typical event on which $t_m = t_m^\star$ and $t_{m+1}$ coincides with a corresponding $\fF_{t_{m-1}}$-measurable time, introduced in the following definition. 

\begin{definition}
\label{defn:star_objects}
	Recall the starred quantities defined in Definition~\ref{defn:time_steps}. For any $m \geq 1$, choose the smallest integer $L_m^\star$ such that
	\begin{equ}
		L_m^\star \geq 1 \vee \bigl(\covariantenergy((B_m^\star+\Sdet_m^\star)(t_m^\star), \higgsdet_m^\star(t_m^\star))\bigr)^{\frac1{10}} \;,
	\end{equ}
	where $\covariantenergy$ is as defined in \eqref{main:eq-covariant-energy}. We then define the two-star quantities, which are deterministic next time step quantities starting from $t_m^{\star}$, as follows. First, we define
	\begin{equ}
	\tau_{m+1}^{\star \star}
	\eqdef \bigl((T^{(1)}\wedge T^{(2)})\bigl((B_m^\star + \Sdet_m^\star)(t_m^\star), \higgsdet_m^\star(t_m^\star)\bigr)\bigr)
	\wedge (L_m^\star)^{-1} \;, \quad
	t_{m+1}^{\star \star} = t_m^\star + \tau_{m+1}^{\star \star} \;.
\end{equ}
	Define the good event
	\begin{equ}
		\event_m \eqdef \{ \tau_m = \tau_m^\star \}
		\cap \{ L_m = L_m^\star \}
		\cap \{ \tau_{m+1} = \tau_{m+1}^{\star \star} \} \;.
	\end{equ}
	For $t \in (t_m^\star, t_{m+1}^{\star \star}]$, define
	\begin{equ}
		B_{m+1}^{+, \star \star}(t)
		\eqdef \Psi_{t - t_m^\star}\bigl((B_m^\star + \Sdet_m^\star)(t_m^\star)\bigr) \;,
	\end{equ}
	where $\Psi$ denotes the YM gradient flow. 
	We also define
	\begin{equs}
		\Dec_{m+1}^{\star}
		&\eqdef \Big(\|F_{B_{m+1}^{+, \star \star}(t_m^\star)} \|_{L^2}^2
		- \| F_{B_{m+1}^{+, \star \star}(t_{m+1}^{\star \star})} \|_{L^2}^2 \Big) \\
		&\quad+ C_{\eqref{e:shiggsdet_decrease}, \gamma_1, 4}
		\one_{\| \higgsdet_m^\star(t_m^\star) \|_{L^4} \geq 1 }
		\| \higgsdet_m^\star(t_m^\star) \|_{L^4}^{\gamma_1 + 2} \tau_{m+1}^{\star \star} \\
		\label{e:dec_star}
		&\quad+ C_{\eqref{e:shiggsdet_decrease}, \gamma_2, r}
		\one_{\| \higgsdet_m^\star(t_m^\star) \|_{L^r} \geq 1 }
		\| \higgsdet_m^\star(t_m^\star) \|_{L^r}^{\gamma_2 + 2} \tau_{m+1}^{\star \star} \;.
	\end{equs}
\end{definition}

\begin{remark}
	All the quantities defined above are measurable with respect to $\fF_{t_{m-1}}$. On the typical event $\event_m$, all the quantities will be equal to the corresponding quantities without stars (e.g. \eqref{e:iteration_deterministic_bounds_iteration_star}), and $(\Asmoothpart(t_m), \phismoothpart(t_m)) = ((B_m^\star + \Sdet_m^\star)(t_m^\star), \higgsdet_m^\star(t_m^\star)) $. 
\end{remark}

One can show that $\event_m$ is indeed a typical event. We do not prove this here, since the current iteration scheme is merely a simplified illustration of the final scheme in Section~\ref{section:decay-estimates-in-expectation} that yields the uniform-in-time bounds. For the moment, we therefore work only on the good event $\event_m$. The treatment of the corresponding bad event for the final iteration scheme is postponed to Section~\ref{section:decay-estimates-in-expectation}, where that event is introduced and analysed in full. The behaviour on the typical event is described by the following proposition.

\begin{proposition}
\label{pr:expectation_one_step}
	Recall $\itobjects_{m+1}$ that is defined in Proposition~\ref{pr:expansion_iteration}. There exists some constant $C$ independent of $\stepupperbound$, such that for any $m \geq 1$ we have 
	\begin{equ}
	\label{e:one_step_change}
		\E \bigl( \itobjects_{m+1} \, \one_{\event_m} | \, \fF_{t_{m-1}} \bigr) \leq C \bigl(\tau_{m+1}^{\star \star} \bigr)^{-1} \;. 
	\end{equ}
\end{proposition}
\begin{proof}
	First note that on the event $\event_m$, we can replace all terms in $\itobjects_{m+1}$ by their corresponding starred versions. More precisely, we have 
	\minilab{e:iteration_deterministic_bounds_iteration_star}
	\begin{equs}
		\itobjects_{m+1} \, \one_{\event_m}
		\label{e:iteration_deterministic_bounds_iteration_star1}
		&= \bigg( \Big \| \curl_{B_{m+1}^{+, \star \star}(t_{m+1}^{\star \star})} e^{\tau_{m+1}^{\star \star} \Delta^F_{(B_m^\star + \Sdet_m^\star)(t_{m}^{\star})}} \linear[B_m^\star](t_m^\star) \Big\|_{L^2}^2 \\
		\label{e:iteration_deterministic_bounds_iteration_star2}
		&+ 2 \int_{\T^2} \Big \langle F^{12}_{B_{m+1}^{+, \star \star}}(t_{m+1}^{\star \star}), \curl_{B_{m+1}^{+, \star \star}(t_{m+1}^{\star \star})} e^{\tau_{m+1}^{\star \star} \Delta^F_{(B_m^\star + \Sdet_m^\star)(t_{m}^{\star})}} \linear[B_m^\star](t_m^\star) \Big \rangle_{\fg} \, \dx \\
		& + 2 \int_{\T^2} \Big \langle F^{12}_{B_{m+1}^{+, \star \star}}(t_{m+1}^{\star \star}), \\
		\label{e:iteration_deterministic_bounds_iteration_star3}
		&\qquad  \Big[( e^{\tau_{m+1}^{\star \star} \Delta^F_{(B_m^\star + \Sdet_m^\star)(t_{m}^{\star})}} \linear[B_m^\star](t_m^\star))^1, ( e^{\tau_{m+1}^{\star \star} \Delta^F_{(B_m^\star + \Sdet_m^\star)(t_{m}^{\star})}} \linear[B_m^\star](t_m^\star))^2 \Big]  \Big \rangle_{\fg} \, \dx \qquad \quad \\
		&+ C \| \higgsdet_m^\star(t_m^\star) \|_{L^4}^{\gamma_1 - 4} \Big| \int_{\T^2} |e^{\tau_{m+1}^{\star \star} \Delta_{(B_m^\star + \Sdet_m^\star)(t_m^\star)}} \higgsdet_m^\star(t_m^\star) |^2 \\
		\label{e:iteration_deterministic_bounds_iteration_star4}
		& \qquad  \cdot \langle e^{\tau_{m+1}^{\star \star} \Delta_{(B_m^\star + \Sdet_m^\star)(t_m^\star)}} \higgsdet_m^\star(t_m^\star), e^{\tau_{m+1}^{\star \star} \Delta_{(B_m^\star + \Sdet_m^\star)(t_m^\star)}} \philinear[B_m^\star](t_m^\star) \rangle_{\V} \,  \dx \Big| \bigg) \one_{\event_m} \;. \qquad  
	\end{equs}
After this replacement, we take conditional expectation on $\fF_{t_{m-1}}$ on both sides, and upper bound the conditional expectations of \eqref{e:iteration_deterministic_bounds_iteration_star1}--\eqref{e:iteration_deterministic_bounds_iteration_star4} by applying the results of Section~\ref{sec:improved_stochastic_estimates}. 
Note that technically speaking, this replacement is important, since the estimates obtained in Section~\ref{sec:improved_stochastic_estimates} require independence between the times $\tstop, t$ and the noise, and we only gain the independence after this replacement. We will apply the results of Section~\ref{sec:improved_stochastic_estimates} (after shifting time by $t_{m-1}$) with the times $\tstop = \tau_m^\star$, $t = \tstop + \tau_{m+1}^{\star \star}$ and the YM gradient flows taken to be $B = B^{\star}_{m}$ on the first interval and $B = B^{+, \star \star}_{m+1}$ on the second interval. The reader should keep in mind that due to our parameter choices, the power of $\tau_m^\star$ appearing in the final estimates of each term in \eqref{e:iteration_deterministic_bounds_iteration_star} will always be positive, and thus we may ultimately throw this away.

\proofstep{Term \eqref{e:iteration_deterministic_bounds_iteration_star1}}: 
	Since this term is always positive, we can throw away the indicator $\one_{\event_m}$. 
	By Lemma~\ref{le:curl_square_expectation} we have
	\begin{equs}
		\E \Big( \eqref{e:iteration_deterministic_bounds_iteration_star1} \one_{\event_m} \Big| \fF_{t_{m-1}} \Big) &\lesssim \bigl( \tau_{m+1}^{\star \star} \bigr)^{-1} + \big( \tau_m^\star \big)^{-2 \lossdel}  
		\cdot \big( \tau_m^\star \big)^{\frac43} \bigl( \tau_{m+1}^{\star \star} \bigr)^{-\frac13} \|\covd_{B^{+, \star}_{m}} F_{B^{+, \star}_{m}} \|_{L^2_{s,y}([t_{m-1}, t^\star_m])}^2  \\
		&+ \big( \tau_m^\star \big)^{-2 \lossdel } \bigl( \tau_{m+1}^{\star \star} \bigr)^{-2\lossdel} \cdot \big( \tau_m^\star \big)^{\frac12} \bigl( \tau_{m+1}^{\star \star} \bigr)^{\frac12} \|\covd_{B^{+, \star \star}_{m+1}} F_{B^{+, \star \star}_{m+1}} \|_{L^2_{s,y}([t^{\star}_m, t^{\star}_m + \tau_{m+1}^{\star \star}])}^2 \\
		& + \big( \tau_m^\star \big)^{-2 \lossdel } \big( \tau_m^\star \big)^{1 - \kappa} \bigl( \tau_{m+1}^{\star \star} \bigr)^{-1 + \kappa} \| \Sdet^\star_{m} (t^\star_m) \|_{L^\infty_x}^2   \\
		& + \big( \tau_m^\star \big)^{-2\lossdel} \big( \tau_m^\star \big)^{\frac14} \bigl( \tau_{m+1}^{\star \star} \bigr)^{-\frac38 - 2\kappa} (\tau_m^\star + \tau_{m+1}^{\star \star})^{\frac18} \bigl(L_{m-1}^{2\fa} \big( \tau_m^\star \big)^{-\lossdel} + L_{m-1}^{4\fa} \bigr) \;.
	\end{equs}
	In the above estimate, we used the facts listed in Lemma~\ref{le:iteration_basic_bounds} and \eqref{eq:Sdet-L-infty-estimates}. In particular, we have 
	\begin{equs}
		\| \Sdet^\star_{m} (t^\star_m) \|_{L^\infty_x}^2 &\lesssim \bigl( \tau_m^\star \bigr)^{-2 \lossdel} \Bigl( \bigl(\tau_m^\star \bigr)^{1 - \frac4r - \frac{4}{\gamma_0} } + \bigl( \tau_m^\star \bigr)^{1 - 2 \kappa} L_{m-1}^{4 \fa} \Bigr) \;, \\
		 \|\covd_{B^{+, \star}_{m}} F_{B^{+, \star}_{m}} \|_{L^2_{s,y}([t_{m-1}, t^\star_m])}^2  &\lesssim \bigl(\tau_m^\star \bigr)^{-1 - 2\lossdel} \;, \\
		  \|\covd_{B^{+, \star \star}_{m+1}} F_{B^{+, \star \star}_{m+1}} \|_{L^2_{s,y}([t^{\star}_m, t^{\star}_m + \tau_{m+1}^{\star \star}])}^2 & \lesssim \bigl(\tau_{m+1}^{\star \star} \bigr)^{-1 - 2\lossdel}  \;. 
	\end{equs}
	Combining this with the fact that $\tau_m^\star \leq L_{m-1}^{-1}$,
	we get that 
	\begin{equ}
		\E \Big( \eqref{e:iteration_deterministic_bounds_iteration_star1} \one_{\event_m} \Big| \fF_{t_{m-1}} \Big) \lesssim C \bigl( \tau_{m+1}^{\star \star} \bigr)^{-1} \;. 
	\end{equ}

\proofstep{Term \eqref{e:iteration_deterministic_bounds_iteration_star2}}: 
	For this term, we put absolute value on \eqref{e:iteration_deterministic_bounds_iteration_star2} and throw away the indicator $\one_{\event_m}$. 
	By Lemmas~\ref{le:cross_term_expansion},~\ref{le:DF_bound} and \eqref{e:iteration_F_endtime_lowerbound}, we get that
	\begin{equs}
		\E \Big( | \eqref{e:iteration_deterministic_bounds_iteration_star2} | \one_{\event_m}  \Big| \fF_{t_{m-1}}  \Big) & \leq \big( \tau_m^\star \big)^{- \lossdel} \bigl( \tau_{m+1}^{\star \star} \bigr)^{- \lossdel} \big( \tau_m^\star \big)^{\frac12 } \| \covd_{B^{+, \star \star}_{m+1}(t^{\star}_m + \tau^{\star \star}_{m+1})} F_{B^{+, \star \star}_{m+1}(t^{\star}_m + \tau^{\star \star}_{m+1})} \|_{L^2}  \\
		& \lesssim \big( \tau_m^\star \big)^{- \lossdel} \bigl( \tau_{m+1}^{\star \star} \bigr)^{- \lossdel} \big( \tau_m^\star \big)^{\frac12 } \\
		&\qquad \qquad \cdot \bigl( \tau_{m+1}^{\star \star} \bigr)^{-\lossdel - \frac12} \bigl( (T^{(1)}\wedge T^{(2)})((B_m^\star + \Sdet_m^\star)(t_m^\star), \higgsdet_m^\star (t_m^\star)) \bigr)^{- \frac12 - \lossdel}  \;.
	\end{equs}
	By Lemma~\ref{le:iteration_basic_bounds}, we have
	\begin{equ}
		\big( \tau_m^\star \big)^{\frac14} \bigl( (T^{(1)}\wedge T^{(2)})((B_m^\star + \Sdet_m^\star )(t_m^\star), \higgsdet_m^\star(t_m^\star)) \bigr)^{-\frac{1}{40} } \lesssim 1 \;,
	\end{equ}
	which implies that the right-hand side is bounded by $C\bigl(\tau_{m+1}^{\star\star}\bigr)^{-1}$.
	 
\proofstep{Term \eqref{e:iteration_deterministic_bounds_iteration_star3}}: 
	We cannot give an upper bound for the absolute value of this term. However, we can obtain a one-sided upper bound for it thanks to the identity in Lemma~\ref{prelim:lem-Hessian-identity}. 
	From the lemma we get 
	\begin{equs}
		\E \Big( \eqref{e:iteration_deterministic_bounds_iteration_star3} \one_{\event_m}  \Big| \fF_{t_{m-1}}  \Big)  &\leq \E \Bigl( \Big \| \curl_{B^{+, \star \star}_{m+1}(t^{\star}_{m} + \tau^{\star \star}_{m+1} )} e^{\tau^{\star \star}_{m+1} \Delta^F_{\Asmoothpart(t^{\star}_m)}} \linear[B_m^\star](t^{\star}_m) \Big\|_{L^2}^2 \\
		&\quad + \Big \| \div_{B^{+, \star \star}_{m+1}(t^{\star}_{m} + \tau^{\star \star}_{m+1} )} e^{\tau^{\star \star}_{m+1} \Delta^F_{\Asmoothpart(t^{\star}_m)}} \linear[B_m^\star](t^{\star}_m) \Big\|_{L^2}^2 \Big| \fF_{t_{m-1}}   \Bigr) \;.
	\end{equs}
	This can be treated in the same way as \eqref{e:iteration_deterministic_bounds_iteration_star1}.

\proofstep{Term \eqref{e:iteration_deterministic_bounds_iteration_star4}}: 
	This term is also positive so we can throw away the indicator $\one_{\event_m}$. By Lemmas~\ref{le:phi_cross_term_estimate} and~\ref{le:iteration_basic_bounds}, we have
	\begin{equs}
		\E \Big(\eqref{e:iteration_deterministic_bounds_iteration_star4} \one_{\event_m} \Big| \fF_{t_{m-1}}  \Big)   &\lesssim \|\higgsdet_m^\star(t^{\star}_m) \|_{L^4}^{\gamma_1 - 4} \big( \tau_m^\star \big)^{\frac12} \|\higgsdet_m^\star(t^{\star}_m) \|_{L^6}^3  \\
		&\lesssim \|\higgsdet_m^\star(t^{\star}_m) \|_{L^4}^{\gamma_1 - 4 + \frac{2 - 12 r^{-1}}{1 - 4r^{-1} } } \big( \tau_m^\star \big)^{\frac12} \|\higgsdet_m^\star(t^{\star}_m) \|_{L^r}^{ \frac{1}{1 - 4 r^{-1}} }  \\
		&\lesssim  \|\higgsdet_m^\star(t^{\star}_m) \|_{L^4}^{\gamma_1 - 4 + \frac{2 - 12 r^{-1}}{1 - 4r^{-1} } } \lesssim \bigl(\tau_{m+1}^{\star \star} \bigr)^{-1} \;,
	\end{equs}
	where we used that $\gamma_1 - 4 \ll \gamma_0 - 2$.
	
	~\\
	The proposition follows from combining the above steps.
\end{proof}

\begin{remark}
	Actually, it turns out that when $B_m$ is close to a YM connection, the expectation of \eqref{e:iteration_deterministic_bounds_iteration_star3} is very negative, with magnitude much larger than $\bigl( \tau_{m+1}^{\star \star} \bigr)^{-1}$. This is consistent with the fact that only an upper bound (not a two-sided bound) can be obtained there, and this observation is exploited in Proposition~\ref{main:prop-YM-connection-decay-mechanism} to obtain the decay estimate.
\end{remark}

From Propositions~\ref{pr:expansion_iteration} and \ref{pr:expectation_one_step}, we can see that if one wants to get a uniform in time bound, the main task is to compare the terms $ C_{\eqref{e:one_step_change}} \bigl(\tau^{\star}_{m+1} \bigr)^{-1}$ and $\frac14 \Dec_{m+1}$. Here we changed $\tau^{\star \star}_{m+1}$ in Proposition~\ref{pr:expectation_one_step} to $\tau^{\star}_{m+1}$ since the error can be ignored (see Lemma~\ref{le:taumod_rare_event}). This replacement will make both quantities $\fF_{t_m}$ measurable so it is conceptually more convenient. 
The following lemma discusses what happens when the decrease term is smaller.

\begin{lemma}
\label{le:close_YM_connection}
	Let $A_0 = \Asmoothpart_0 + \Aroughpart_0$, $\higgs_0 = \phismoothpart_0 + \phiroughpart_0$ be $L$-admissible. Recall $T^{(L)}$ and $\Dec$ from Definitions~\ref{defn:L_admissible} and~\eqref{e:dec_defn}.
	Suppose that   
	\begin{equ}
	\label{e:close_YM_connection_assumption}
		\Dec_{T^{(L)}}(\Asmoothpart_0,\phismoothpart_0) \leq \bigl(16 C_{\eqref{e:one_step_change}} + 1\bigr)\bigl(T^{(L)}\bigr)^{-1} \;, 
	\end{equ}
	and $T^{(L)} < L^{-1}$, then $T^{(L)} = T^{(1)}(\Asmoothpart_0)$. If furthermore we have $T^{(L)} < \stepupperbound \wedge \logtimestep$, then there exist some constants $c, \delta_0$ such that if $\stepupperbound < \delta_0$, then 
	\begin{equ}
	\label{e:close_YM_connection}
		\| F_{\YMHFsmooth(T^{(L)})} \|_{L^2} \geq c \frac{\abs{\log T^{(L)}}}{T^{(L)}} \;.
	\end{equ}
\end{lemma}
\begin{proof}
	To show that $T^{(L)} = T^{(1)}(\Asmoothpart_0)$, recalling Definition~\ref{defn:tau}, we need to rule out the following two cases.
	\begin{enumerate}
		\item If $T^{(L)} = T^{(2)}(\Asmoothpart_0,\phismoothpart_0)$, then by definition of $T^{(2)}$, we have 
		\begin{equ}
			\Dec_{T^{(L)}}(\Asmoothpart_0,\phismoothpart_0) \geq \logtimestep \bigl(T^{(L)}\bigr)^{-1} \abs{\log T^{(L)}} \;. 
		\end{equ}
		This contradicts \eqref{e:close_YM_connection_assumption} if $T^{(L)}$ is small enough depending on $\logtimestep$ and $C_{\eqref{e:one_step_change}}$. 
		
		\item If $T^{(L)} = \| \phismoothpart_0 \|_{L^r}^{-\gamma_0}$, then we have 
		\begin{equs}
			\Dec_{T^{(L)}}(\Asmoothpart_0,\phismoothpart_0) &\geq C_{\eqref{e:shiggsdet_decrease}, \gamma_2, r} \| \phismoothpart_0 \|_{L^r}^{\gamma_2 + 2} T^{(L)} \\
			& \geq C_{\eqref{e:shiggsdet_decrease}, \gamma_2, r} \bigl(T^{(L)}\bigr)^{-\frac{\gamma_2 + 2}{\gamma_0} + 1} \geq C_{\eqref{e:shiggsdet_decrease}, \gamma_2, r} \bigl(T^{(L)}\bigr)^{-1.1}\;. 
		\end{equs}
		This contradicts \eqref{e:close_YM_connection_assumption} if $T^{(L)}$ is small enough depending on $C_{\eqref{e:shiggsdet_decrease}, \gamma_2, r}, C_{\eqref{e:one_step_change}}$.
	\end{enumerate}
	To show the bound in \eqref{e:close_YM_connection}, first note that by the definition of $T^{(1)}$, if $T^{(L)} < \stepupperbound \wedge \logtimestep$, we have 
	\begin{equ}
		T^{(L)} = \logtimestep \|F_{\YMHFsmooth(0)} \|_{L^2}^{-1} \log \bigl(\|F_{\YMHFsmooth(0)} \|_{L^2} \vee 2\bigr) \;. 
	\end{equ}
	Then, from the definition of $\Dec_{T^{(L)}}(\Asmoothpart_0,\phismoothpart_0)$, \eqref{e:timestep2_lower_bound} and \eqref{e:close_YM_connection_assumption}, we have 
	\begin{equs}
		 \| F_{\YMHFsmooth(T^{(L)})} \|_{L^2}^2 &\geq \|F_{\YMHFsmooth(0)} \|_{L^2}^2 - \Dec_{T^{(L)}}(\Asmoothpart_0,\phismoothpart_0) \\
		 & \geq\Bigl( \logtimestep \frac{\log \abs{T^{(L)}}}{T^{(L)}} \Bigr)^2 - \bigl(16 C_{\eqref{e:one_step_change}} + 1\bigr)\bigl(T^{(L)}\bigr)^{-1} 
		 \geq \Bigl( c \frac{\abs{\log T^{(L)}}}{T^{(L)}} \Bigr)^2 
	\end{equs}
	if $T^{(L)} < \delta_0$ is small enough. The bound in \eqref{e:close_YM_connection} is proved. 
\end{proof}

Combining Propositions~\ref{pr:expansion_iteration} and \ref{pr:expectation_one_step} and Lemma~\ref{le:close_YM_connection}, and ignoring the bad events, we arrive at 
\begin{equ}
	\E \Big( \covariantenergy(\Asmoothpart(t_{m+1}), \phismoothpart(t_{m+1})) \Big) \leq \E \Big( \covariantenergy(\Asmoothpart(t_{m}), \phismoothpart(t_{m})) \Big) + C_{\eqref{e:one_step_change}} \E  \Big[ \Big( T^{(1)} \bigl(\Asmoothpart(t_m) \bigr) \Big)^{-1} \Big] \;. 
\end{equ}
This essentially gives global well-posedness, up to some additional technical details. Heuristically, let us ignore all the expectations here and assume that everything is deterministic. Suppose that we start from initial data with energy $\covariantenergy_0 \geq 1$, and let $M \geq 1$ be the smallest integer such that $\covariantenergy_M \geq 2 \covariantenergy_0$. Let $N$ denote the number of indices $0 \leq m < M$ for which $\covariantenergy_{m+1} > \covariantenergy_m$. By Lemma~\ref{le:close_YM_connection}, this can only happen for a given index $m$ if the following estimate holds:
\begin{equ}
	\tau_{m+1} = T^{(1)}\bigl(\Asmoothpart(t_m), 0 \bigr) = \min \Bigl\{
		\logtimestep \|F_{\Asmoothpart(t_m)}\|_{L^2}^{-1}
		\log\bigl(\|F_{\Asmoothpart(t_m)}\|_{L^2} \vee 2\bigr),
		\stepupperbound \Bigr\}
	 \geq \frac{c}{\sqrt{\covariantenergy_m}} \geq \frac{c}{2 \sqrt{\covariantenergy_0}} \;.
\end{equ}
Since the total energy increase $\covariantenergy_M - \covariantenergy_0$ is at least $\covariantenergy_0$, it follows that
\begin{equ}
	N \geq \frac{c \covariantenergy_0}{ 2C_{\eqref{e:one_step_change}} \sqrt{\covariantenergy_0}}  \;.
\end{equ}
Consequently, since $\covariantenergy_0 \geq1$, we get
\begin{equ}
	t_M \geq \sum_{\substack{0 \leq m < M\\ \covariantenergy_{m+1}>\covariantenergy_m}}\tau_{m+1} \geq N \cdot \frac{c}{2 \sqrt{\covariantenergy_0} }
	\gtrsim 1 \;.
\end{equ}
Thus, the time required for the energy to double is bounded below by a positive constant independent of its initial size. Iterating this argument yields global well-posedness.

\subsection{Modified iteration scheme}\label{section:decay-estimates-in-expectation}

To obtain uniform in time bounds, we need some decay mechanism when the scenario in Lemma~\ref{le:close_YM_connection} happens. The idea, as already discussed in the introduction, is that in this case $B$ is close to some unstable YM connection, so that when perturbing $B$ along its unstable directions, the energy will decrease. Since the linear stochastic object always has probability of aligning with such directions, one would hope that adding the (smoothed out) linear stochastic object back into $\Asmoothpart$ would decrease the energy in expectation. 
Proposition~\ref{main:prop-YM-connection-decay-mechanism} has made this heuristic rigorous. We therefore will modify our definition of $\Asmoothpart$, and therefore the definitions of other related quantities in the iteration scheme as follows. 
\begin{remark}
\label{rmk:decay_mechanism_m_step_application}
	We will apply Proposition~\ref{main:prop-YM-connection-decay-mechanism} when the scenario in Lemma~\ref{le:close_YM_connection} happens. In that case, in the statement of Proposition~\ref{main:prop-YM-connection-decay-mechanism} we can choose  
	\begin{equ}
		C_1 = \sqrt{16 C_{\eqref{e:one_step_change}} + 1} \;, \quad C_2 = c_{\eqref{e:close_YM_connection}}. 
	\end{equ}
	In the following, we will always choose $\stepupperbound$ small enough such that Proposition~\ref{main:prop-YM-connection-decay-mechanism} applies. 
\end{remark}

\begin{definition}
\label{defn:final_timestep}
	Let $A_0 = \Asmoothpart_0 + \Aroughpart_0$, $\higgs_0 = \phismoothpart_0 + \phiroughpart_0$ be $L$-admissible.
	Recall Definition~\ref{defn:one_step_iteration_scheme}. Let $B$ and $\YMHFsmooth$ be defined as in Definition~\ref{defn:one_step_iteration_scheme}.
	We define $T^{(L)}_{1/2}(\Asmoothpart_0, \Aroughpart_0, \phismoothpart_0, \phiroughpart_0)$ as follows: 
	\begin{enumerate}
		\item Recall the definition \eqref{e:dec_defn} of $\Dec_t$. If we have  
		\begin{equ}
			\Dec_{T^{(L)}}(\Asmoothpart_0, \phismoothpart_0 ) \leq \big(16 C_{\eqref{e:one_step_change}} + 1 \big) \bigl( T^{(L)} \bigr)^{-1}  \quad \text{ and } T^{(L)} < \stepupperbound \wedge \logtimestep \;,
		\end{equ}
		then we define
		\begin{equ}
			T^{(L)}_{1/2}(\Asmoothpart_0, \Aroughpart_0, \phismoothpart_0, \phiroughpart_0) \eqdef \frac12 
		T^{(L)}(\Asmoothpart_0, \Aroughpart_0, \phismoothpart_0, \phiroughpart_0) \;, \quad \tau^{(L)}_{1/2} \eqdef T^{(L)}_{1/2} \wedge \tau^{(L)} \;. 
		\end{equ}
		Moreover, for any $t \in (0, \tau^{(L)}]$ let $\flow_{t}(\Asmoothpart_0, \Aroughpart_0, \phismoothpart_0, \phiroughpart_0; \xi, \zeta) = (\Asmoothpart(t), \Aroughpart(t), \phismoothpart(t), \phiroughpart(t))$ be defined as in Definition~\ref{defn:one_step_iteration_scheme}. 
		On the event $\{ \tau^{(L)} = T^{(L)} \}$, we define 
		\begin{equ}
			\wt{\flow}_{T^{(L)}}(\Asmoothpart_0, \Aroughpart_0, \phismoothpart_0, \phiroughpart_0; \xi, \zeta) = (\Asmoothpartmod(T^{(L)}), \Aroughpartmod(T^{(L)}), \phismoothpartmod(T^{(L)}), \phiroughpartmod(T^{(L)}))
		\end{equ}
        as follows:
		\begin{equs}
			\Asmoothpartmod(T^{(L)}) &\eqdef \Asmoothpart(T^{(L)}) + e^{(T^{(L)} - T^{(L)}_{1/2})\Delta^F_{B(T^{(L)}_{1/2})}} \linear[B](T^{(L)}_{1/2}) \;, &\quad \Aroughpartmod(T^{(L)}) &\eqdef A(T^{(L)}) - \Asmoothpartmod(T^{(L)}) \;, \\
			\phismoothpartmod(T^{(L)}) &\eqdef \phismoothpart(T^{(L)}) \;, &\quad \phiroughpartmod(T^{(L)}) &\eqdef \phiroughpart(T^{(L)}) \;. 
		\end{equs}
		If $t \in (0, \tau^{(L)})$ or $t = \tau^{(L)}$ but $\tau^{(L)} \neq T^{(L)}$, we define $\wt{\flow}_t = \flow_t$ as in Definition~\ref{defn:one_step_iteration_scheme}. 
		
		\item Otherwise, everything remains the same as before, i.e.\ $\wt{\flow}_t = \flow_t$. We define $\tau^{(L)}_{1/2} = T^{(L)}_{1/2} = 0$. 
	\end{enumerate}
\end{definition}

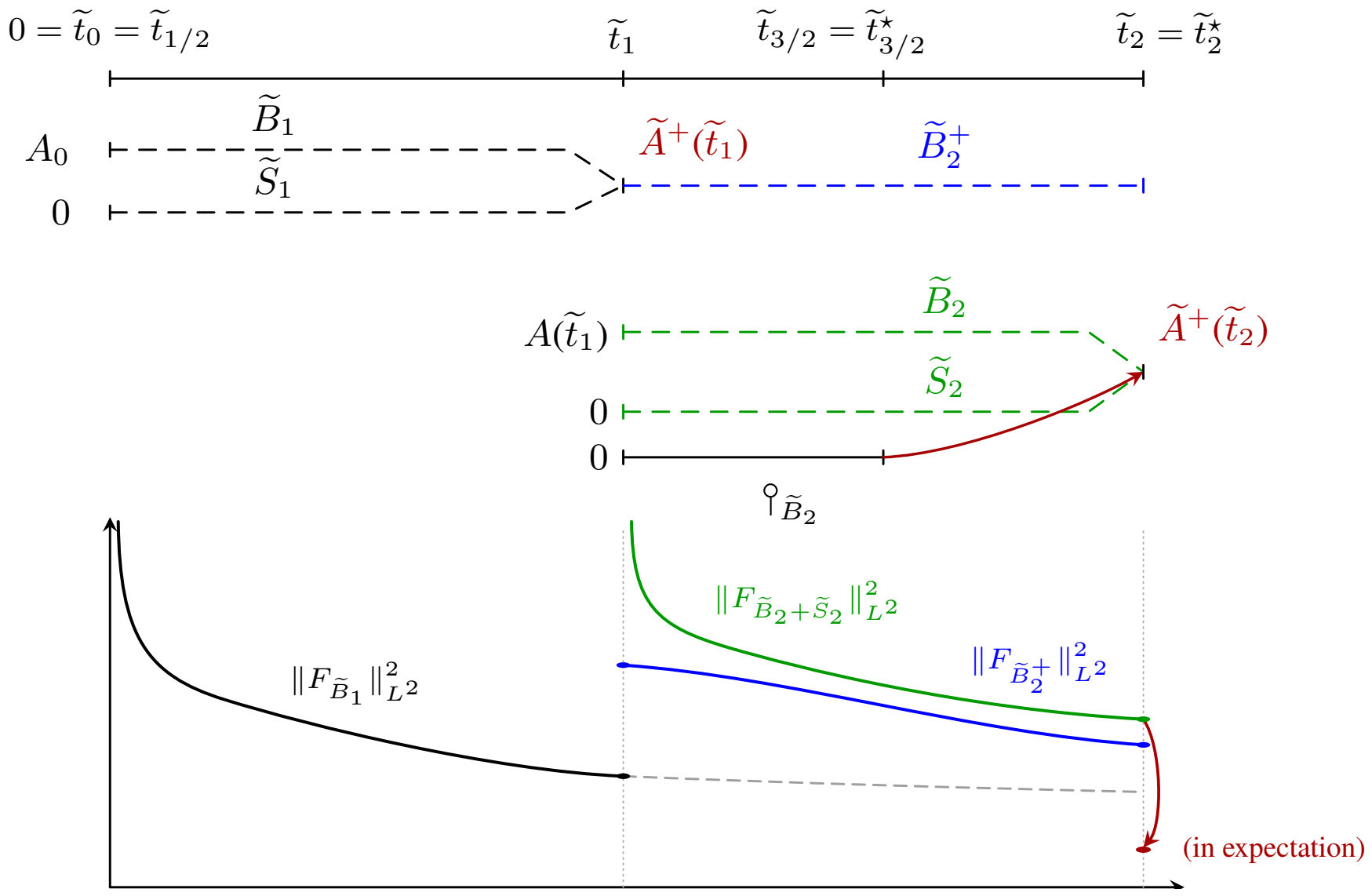
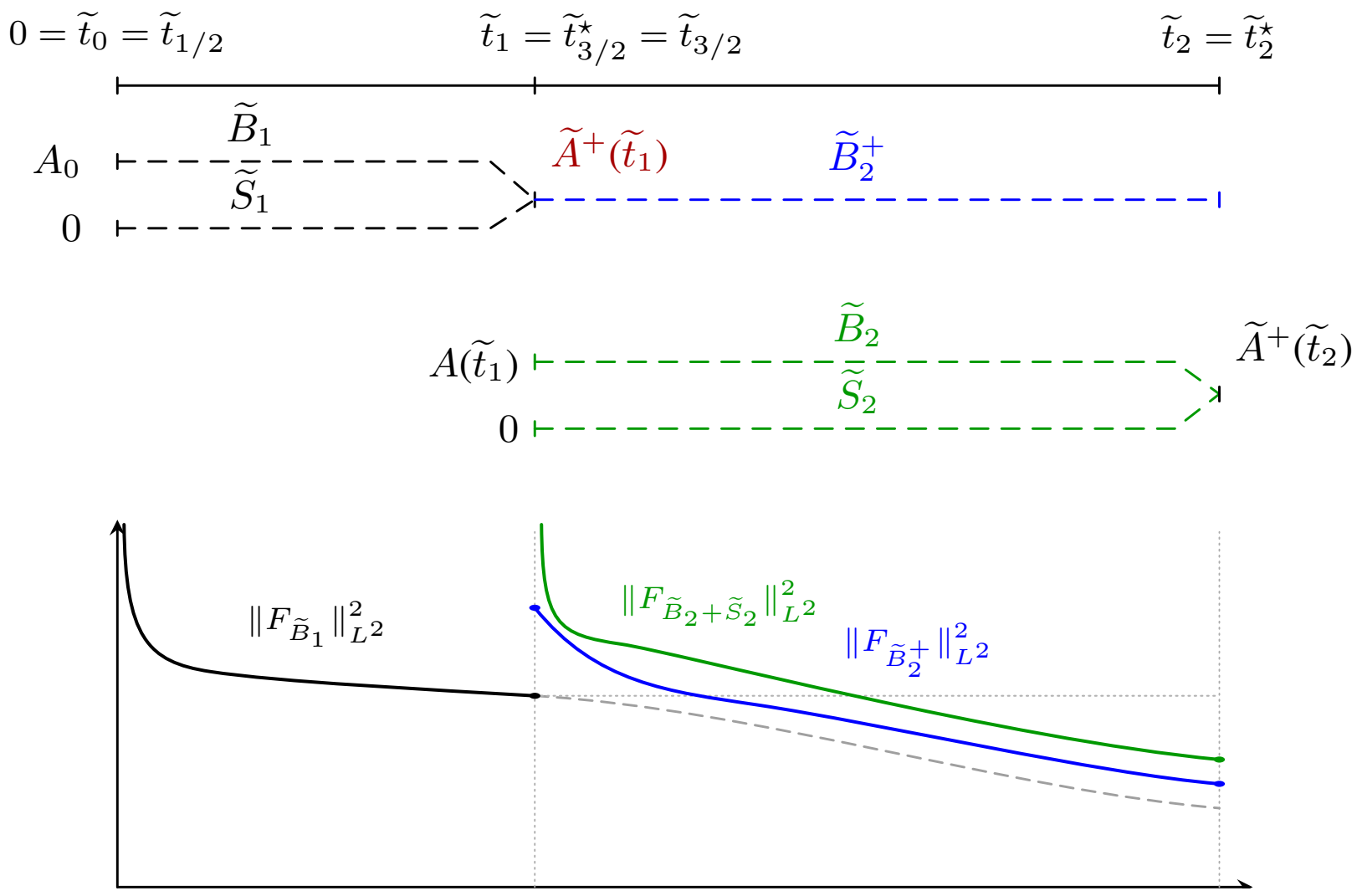
\begin{figure}[p]
\centering
\begin{subfigure}{\textwidth}
  \centering
  \ifdefined\tikzsetnextfilename
    \tikzsetnextfilename{iteration-scheme-mod-small-decrease}
  \fi
  \begingroup
\NewDocumentCommand{\iterlinear}{m}{%
  \tikzsetnextfilename{iteration-linear-\pdfmdfivesum{\noexpand#1}}%
  \tikz[baseline=\linearbaseline,x=1cm,y=1cm,
    every node/.style={inner sep=0.45pt},scale=\linearscale,thin]{
    \coordinate (A) at (0.08,0);
    \coordinate (Bcenter) at (0.08,0.55);
    \draw[white,line width=0] (-0.25,0) -- (-0.25,0);
    \draw[white,line width=0] (0.25,0) -- (0.25,0);
    \draw[heat] ($(A)+(-0.09,0)$) -- ($(Bcenter)+(-0.09,0)$);
    \node[blackdot] at (Bcenter) {};
    \node[right] at (\frequencyxoffset,\frequencyyoffset)
      {$\scriptstyle #1$};
  }%
}
\begin{tikzpicture}[
  x=1.154cm,
  y=0.48cm,
  line cap=round,
  line join=round
]
  \colorlet{lollipopcolor}{red!65!black}
  \colorlet{futurecolor}{blue}
  \colorlet{restartcolor}{green!60!black}
  \renewcommand{\linearscale}{0.45}

  \tikzset{
    iterline/.style={line width=0.55pt},
    iterdet/.style={iterline,dash pattern=on 4pt off 3pt},
    iterfuture/.style={iterline,futurecolor,dash pattern=on 4pt off 3pt},
    iterrestart/.style={iterline,restartcolor,dash pattern=on 4pt off 3pt},
    itertick/.style={line width=0.55pt},
    iterlollipop/.style={->,line width=0.65pt,lollipopcolor},
    energyfirst/.style={line width=0.75pt},
    energyplus/.style={line width=0.75pt,futurecolor},
    energyrestart/.style={line width=0.75pt,restartcolor},
    energyold/.style={line width=0.55pt,gray!75,dash pattern=on 3pt off 2pt},
    energydrop/.style={->,line width=0.70pt,lollipopcolor},
    iterguide/.style={line width=0.35pt,gray!65,densely dotted},
    itertitle/.style={font=\small},
    iterlab/.style={font=\footnotesize},
    iterflowlab/.style={iterlab,font=\small},
    iterenergylab/.style={font=\scriptsize}
  }

  \def\xZero{1.35}
  \def\xHalfOne{4.40}
  \def\xOne{5.10}
  \def\xThreeHalf{7.00}
  \def\xTwo{8.90}
  \def\xMergeOne{4.70}
  \def\xMergeTwo{8.50}

  \def\yTitle{9.25}
  \def\yTime{8.05}
  \def\yBOne{6.80}
  \def\yMergeOne{6.17}
  \def\ySOne{5.70}
  \def\yLinOne{4.78}
  \def\yBTwo{3.60}
  \def\yMergeTwo{2.90}
  \def\ySTwo{2.20}
  \def\yLinTwo{1.40}
  \def\yPhi{0.65}
  \def\yEnergyTop{0.15}
  \def\yEnergyAxis{-6.15}


  \draw[iterline] (\xZero,\yTime) -- (\xTwo,\yTime);
  \foreach \X in {\xZero,\xOne,\xThreeHalf,\xTwo}
    \draw[itertick] (\X,\yTime+0.13) -- (\X,\yTime-0.13);

  \node[iterlab,anchor=south] at (\xZero,\yTime+0.18) {$0=\tmod_0 = \tmod_{1/2}$};
  \node[iterlab,anchor=south] at (\xOne,\yTime+0.18) {$\tmod_1$};
  \node[iterlab,anchor=south] at (\xThreeHalf-0.3,\yTime+0.08) {$\tmod_{3/2} = \tmod_{3/2}^\star$};
  \node[iterlab,anchor=south] at (\xTwo+0.2,\yTime+0.18)
    {$\tmod_2 = \tmod_2^\star $};

  \draw[iterdet] (\xZero,\yBOne) -- (\xMergeOne,\yBOne) -- (\xOne,\yMergeOne);
  \draw[iterdet] (\xZero,\ySOne) -- (\xMergeOne,\ySOne) -- (\xOne,\yMergeOne);
  \foreach \Y in {\yBOne,\ySOne}
    \draw[itertick] (\xZero,\Y+0.12) -- (\xZero,\Y-0.12);
  \draw[itertick] (\xOne,\yMergeOne+0.12) -- (\xOne,\yMergeOne-0.12);

  \node[iterflowlab,anchor=east] at (\xZero-0.18,\yBOne) {$A_0$};
  \node[iterflowlab,anchor=east] at (\xZero-0.18,\ySOne) {$0$};
  \node[iterflowlab,anchor=south] at (2.55,\yBOne) {$\Bmod_1$};
  \node[iterflowlab,anchor=south] at (2.55,\ySOne) {$\Sdetmod_1$};

  \node[iterflowlab,lollipopcolor,anchor=south west]
    at (\xOne,\yMergeOne+0.17)
    {$\Asmoothpartmod(\tmod_1)$};

  \draw[iterfuture] (\xOne,\yMergeOne) -- (\xTwo,\yMergeOne);
  \draw[itertick,futurecolor] (\xTwo,\yMergeOne+0.12) -- (\xTwo,\yMergeOne-0.12);
  \node[iterflowlab,futurecolor,anchor=south]
    at (7.45,\yMergeOne+0.06) {$\widetilde B_2^+$};

  \draw[iterrestart] (\xOne,\yBTwo) -- (\xMergeTwo,\yBTwo) -- (\xTwo,\yMergeTwo);
  \draw[iterrestart] (\xOne,\ySTwo) -- (\xMergeTwo,\ySTwo) -- (\xTwo,\yMergeTwo);
  \draw[itertick,restartcolor] (\xOne,\yBTwo+0.12) -- (\xOne,\yBTwo-0.12);
  \draw[itertick,restartcolor] (\xOne,\ySTwo+0.12) -- (\xOne,\ySTwo-0.12);
  \draw[itertick] (\xTwo,\yMergeTwo+0.12) -- (\xTwo,\yMergeTwo-0.12);

  \node[iterflowlab,anchor=east] at (\xOne,\yBTwo) {$A(\tmod_1)$};
  \node[iterflowlab,anchor=east] at (\xOne,\ySTwo) {$0$};
  \node[iterflowlab,restartcolor,anchor=south]
    at (7.45,\yBTwo) {$\Bmod_2$};
  \node[iterflowlab,restartcolor,anchor=south]
    at (7.45,\ySTwo) {$\Sdetmod_2$};

  \draw[iterline] (\xOne,\yLinTwo) -- (\xThreeHalf,\yLinTwo);
  \draw[itertick] (\xOne,\yLinTwo+0.12) -- (\xOne,\yLinTwo-0.12);
  \draw[itertick] (\xThreeHalf,\yLinTwo+0.12) -- (\xThreeHalf,\yLinTwo-0.12);
  \node[iterflowlab,anchor=east] at (\xOne,\yLinTwo) {$0$};
  \node[iterflowlab,anchor=north east]
    at (6.65, \yLinTwo-0.20)
    {$\iterlinear{\Bmod_2}$};
  \draw[iterlollipop] (\xThreeHalf,\yLinTwo)
    .. controls (7.60,1.46) and (8.50,2.35) .. (\xTwo,\yMergeTwo);
  \node[iterflowlab,lollipopcolor,anchor=south west]
    at (\xTwo,\yMergeTwo+0.17)
    {$\Asmoothpartmod(\tmod_2)$};


  \draw[iterline,->] (\xZero,\yEnergyAxis) -- (\xTwo+0.30,\yEnergyAxis);
  \draw[iterline,->] (\xZero,\yEnergyAxis) -- (\xZero,\yEnergyTop+0.20);

  \foreach \X in {\xOne,\xTwo}
    \draw[iterguide] (\X,\yEnergyAxis) -- (\X,\yEnergyTop);

  \draw[energyfirst]
    (\xZero+0.06,\yEnergyTop+0.12)
    .. controls (\xZero+0.07,-1.75) and (1.62,-2.430) .. (2.18,-2.85)
    .. controls (2.80,-3.315) and (4.15,-4.105) .. (\xOne,-4.20);
  \draw[energyold]
    (\xOne,-4.20)
    .. controls (6.32,-4.322) and (7.63,-4.407) .. (\xTwo,-4.477);
  \fill (\xOne,-4.20) circle (0.050);
  \node[iterenergylab,anchor=south]
    at (3.15,-3.20) {$\|F_{\widetilde B_1}\|_{L^2}^2$};

  \draw[energyplus]
    (\xOne,-2.25)
    .. controls (6.30,-2.45) and (7.65,-3.450) .. (\xTwo,-3.65);
  \fill[futurecolor] (\xOne,-2.25) circle (0.050);
  \fill[futurecolor] (\xTwo,-3.65) circle (0.050);
  \node[iterenergylab,futurecolor,anchor=south]
    at (8.15,-3) {$\|F_{\widetilde B_2^+}\|_{L^2}^2$};

  \draw[energyrestart]
    (\xOne+0.06,\yEnergyTop+0.12)
    .. controls (\xOne+0.07,-1.20) and (5.34,-1.572) .. (5.88,-1.95)
    .. controls (6.35,-2.279) and (7.55,-3.018) .. (\xTwo,-3.20);
  \node[iterenergylab,restartcolor,anchor=south]
    at (6.45,-1.78) {$\|F_{\widetilde B_2+\widetilde S_2}\|_{L^2}^2$};

  \draw[energydrop]
    (\xTwo,-3.20)
    .. controls (\xTwo+0.14,-3.70) and (\xTwo+0.14,-5.05) .. (\xTwo,-5.45);
  \fill[restartcolor] (\xTwo,-3.20) circle (0.050);
  \fill[lollipopcolor] (\xTwo,-5.49) circle (0.055);
  \node[iterenergylab,lollipopcolor,anchor=west]
    at (\xTwo+0.18,-5.49) {(in expectation)};

  \path ([xshift=-1.25cm]\xZero,\yEnergyAxis)
    -- ([xshift=2.85cm]\xTwo,\yEnergyAxis);

\end{tikzpicture}
\endgroup
  \caption{When $\|F_{\widetilde B_2^+}\|_{L^2}^2$ does not decay sufficiently, adding the linear stochastic object from
$[\tmod_1,\tmod_{3/2}]$ produces a further decrease of the energy (in expectation), where $\tmod_{3/2}$ is the midpoint of $[\tmod_1, \tmod_2]$. The
corresponding effect is shown in red in the lower graph.}
  \label{fig:iteration_scheme_mod_small_decrease}
\end{subfigure}

\vspace{1em}

\begin{subfigure}{\textwidth}
  \centering
  \ifdefined\tikzsetnextfilename
    \tikzsetnextfilename{iteration-scheme-mod-large-decrease}
  \fi
  \begingroup
\NewDocumentCommand{\iterlinear}{m}{%
  \tikzsetnextfilename{iteration-linear-\pdfmdfivesum{\noexpand#1}}%
  \tikz[baseline=\linearbaseline,x=1cm,y=1cm,
    every node/.style={inner sep=0.45pt},scale=\linearscale,thin]{
    \coordinate (A) at (0.08,0);
    \coordinate (Bcenter) at (0.08,0.55);
    \draw[white,line width=0] (-0.25,0) -- (-0.25,0);
    \draw[white,line width=0] (0.25,0) -- (0.25,0);
    \draw[heat] ($(A)+(-0.09,0)$) -- ($(Bcenter)+(-0.09,0)$);
    \node[blackdot] at (Bcenter) {};
    \node[right] at (\frequencyxoffset,\frequencyyoffset)
      {$\scriptstyle #1$};
  }%
}
\begin{tikzpicture}[
  x=0.88cm,
  y=0.48cm,
  line cap=round,
  line join=round
]
  \colorlet{lollipopcolor}{red!65!black}
  \colorlet{futurecolor}{blue}
  \colorlet{restartcolor}{green!60!black}
  \renewcommand{\linearscale}{0.45}

  \tikzset{
    iterline/.style={line width=0.55pt},
    iterdet/.style={iterline,dash pattern=on 4pt off 3pt},
    iterfuture/.style={iterline,futurecolor,dash pattern=on 4pt off 3pt},
    iterrestart/.style={iterline,restartcolor,dash pattern=on 4pt off 3pt},
    itertick/.style={line width=0.55pt},
    iterlollipop/.style={->,line width=0.65pt,lollipopcolor},
    energyfirst/.style={line width=0.75pt},
    energyplus/.style={line width=0.75pt,futurecolor},
    energyrestart/.style={line width=0.75pt,restartcolor},
    energyold/.style={line width=0.55pt,gray!75,dash pattern=on 3pt off 2pt},
    energythreshold/.style={line width=0.40pt,gray!60,densely dotted},
    iterguide/.style={line width=0.35pt,gray!65,densely dotted},
    itertitle/.style={font=\small},
    iterlab/.style={font=\footnotesize},
    iterflowlab/.style={iterlab,font=\small},
    iterenergylab/.style={font=\scriptsize}
  }

  \def\xZero{1.35}
  \def\xHalfOne{4.40}
  \def\xOne{5.10}
  \def\xTwo{11.25}
  \def\xMergeOne{4.70}
  \def\xMergeTwo{10.85}

  \def\yTitle{9.25}
  \def\yTime{8.05}
  \def\yBOne{6.80}
  \def\yMergeOne{6.17}
  \def\ySOne{5.70}
  \def\yLinOne{4.78}
  \def\yBTwo{3.50}
  \def\yMergeTwo{2.97}
  \def\ySTwo{2.40}
  \def\yPhi{1.65}
  \def\yEnergyTop{0.70}
  \def\yEnergyAxis{-5.15}


  \draw[iterline] (\xZero,\yTime) -- (\xTwo,\yTime);
  \foreach \X in {\xZero,\xOne,\xTwo}
    \draw[itertick] (\X,\yTime+0.13) -- (\X,\yTime-0.13);
  \node[iterlab,anchor=south] at (\xZero,\yTime+0.18) {$0=\tmod_0=\tmod_{1/2}$};
  \node[iterlab,anchor=south] at (\xOne + 0.7,\yTime+0.06) {$\tmod_1 = \tmod_{3/2}^\star = \tmod_{3/2}$};
  \node[iterlab,anchor=south] at (\xTwo,\yTime+0.18)
    {$\tmod_2=\tmod_2^\star $};

  \draw[iterdet] (\xZero,\yBOne) -- (\xMergeOne,\yBOne) -- (\xOne,\yMergeOne);
  \draw[iterdet] (\xZero,\ySOne) -- (\xMergeOne,\ySOne) -- (\xOne,\yMergeOne);
  \foreach \Y in {\yBOne,\ySOne}
    \draw[itertick] (\xZero,\Y+0.12) -- (\xZero,\Y-0.12);
  \draw[itertick] (\xOne,\yMergeOne+0.12) -- (\xOne,\yMergeOne-0.12);

  \node[iterflowlab,anchor=east] at (\xZero-0.18,\yBOne) {$A_0$};
  \node[iterflowlab,anchor=east] at (\xZero-0.18,\ySOne) {$0$};
  \node[iterflowlab,anchor=south] at (2.55,\yBOne) {$\Bmod_1$};
  \node[iterflowlab,anchor=south] at (2.55,\ySOne) {$\Sdetmod_1$};

  \node[iterflowlab,lollipopcolor,anchor=south west]
    at (\xOne,\yMergeOne+0.17)
    {$\Asmoothpartmod(\tmod_1)$};

  \draw[iterfuture] (\xOne,\yMergeOne) -- (\xTwo,\yMergeOne);
  \draw[itertick,futurecolor] (\xTwo,\yMergeOne+0.12) -- (\xTwo,\yMergeOne-0.12);
  \node[iterflowlab,futurecolor,anchor=south] at (8.00,\yMergeOne+0.06) {$\widetilde B_2^+$};

  \draw[iterrestart] (\xOne,\yBTwo) -- (\xMergeTwo,\yBTwo) -- (\xTwo,\yMergeTwo);
  \draw[iterrestart] (\xOne,\ySTwo) -- (\xMergeTwo,\ySTwo) -- (\xTwo,\yMergeTwo);
  \draw[itertick,restartcolor] (\xOne,\yBTwo+0.12) -- (\xOne,\yBTwo-0.12);
  \draw[itertick,restartcolor] (\xOne,\ySTwo+0.12) -- (\xOne,\ySTwo-0.12);
  \draw[itertick] (\xTwo,\yMergeTwo+0.12) -- (\xTwo,\yMergeTwo-0.12);

  \node[iterflowlab,anchor=east] at (\xOne,\yBTwo) {$A(\tmod_1)$};
  \node[iterflowlab,anchor=east] at (\xOne,\ySTwo) {$0$};
  \node[iterflowlab,restartcolor,anchor=south]
    at (8.00,\yBTwo) {$\Bmod_2$};
  \node[iterflowlab,restartcolor,anchor=south]
    at (8.00,\ySTwo) {$\Sdetmod_2$};
  \node[iterflowlab,anchor=south west] at (\xTwo,\yMergeTwo+0.17)
    {$\Asmoothpartmod(\tmod_2)$};


  \draw[iterline,->] (\xZero,\yEnergyAxis) -- (\xTwo+0.30,\yEnergyAxis);
  \draw[iterline,->] (\xZero,\yEnergyAxis) -- (\xZero,\yEnergyTop+0.20);
  \foreach \X in {\xOne,\xTwo}
    \draw[iterguide] (\X,\yEnergyAxis) -- (\X,\yEnergyTop);

  \draw[energyfirst]
    (\xZero+0.06,\yEnergyTop+0.12)
    .. controls (\xZero+0.07,-1.05) and (1.62,-1.443) .. (2.18,-1.60)
    .. controls (2.80,-1.774) and (4.15,-1.895) .. (\xOne,-2.00);
  \draw[energyold]
    (\xOne,-2.00)
    .. controls (7.10,-2.220) and (9.35,-3.55) .. (\xTwo,-3.85);
  \draw[energythreshold] (\xOne,-2.00) -- (\xTwo,-2.00);
  \fill (\xOne,-2.00) circle (0.050);
  \node[iterenergylab,anchor=south]
    at (3.15,-1.45) {$\|F_{\widetilde B_1}\|_{L^2}^2$};

  \draw[energyplus]
    (\xOne,-0.55)
    .. controls (5.65,-1.80) and (6.30,-1.930) .. (7.10,-2.15)
    .. controls (8.20,-2.453) and (10.15,-3.28) .. (\xTwo,-3.45);
  \fill[futurecolor] (\xOne,-0.55) circle (0.050);
  \fill[futurecolor] (\xTwo,-3.45) circle (0.050);
  \node[iterenergylab,futurecolor,anchor=south]
    at (8.55,-1.95) {$\|F_{\widetilde B_2^+}\|_{L^2}^2$};

  \draw[energyrestart]
    (\xOne+0.06,\yEnergyTop+0.12)
    .. controls (\xOne+0.07,-0.85) and (5.30,-0.993) .. (5.88,-1.15)
    .. controls (6.35,-1.277) and (9.35,-2.72) .. (\xTwo,-3.05);
  \fill[restartcolor] (\xTwo,-3.05) circle (0.050);
  \node[iterenergylab,restartcolor,anchor=south west]
    at (5.72,-1.1) {$\|F_{\widetilde B_2+\widetilde S_2}\|_{L^2}^2$};

  \path ([xshift=-1.25cm]\xZero,\yEnergyAxis)
    -- ([xshift=2.85cm]\xTwo,\yEnergyAxis);

\end{tikzpicture}
\endgroup
\caption{When $\|F_{\widetilde B_2^+}\|_{L^2}^2$ decays sufficiently, no modification
is needed. We have $\tmod_{3/2}=\tmod_1$, and this step coincides with the iteration step
of Section~\ref{sec:iteration}.}
  \label{fig:iteration_scheme_mod_large_decrease}
\end{subfigure}
\caption{Illustration of Definition~\ref{defn:time_steps_mod}. For simplicity, we assume all the starred times are the same as the non-starred times.  
In each panel, the upper
graph shows the evolution of the connection one-forms in the modified iteration,
while the lower graph shows the corresponding evolution of their squared $L^2$
curvature norms. We only focus on the evolution of the connection one-forms as it is the only part that is not treated well enough in Section~\ref{sec:iteration}. }
\label{fig:iteration_scheme_mod}
\end{figure}

Now we discuss how to modify our iteration scheme from Definition~\ref{defn:time_steps}. 
The main difference from Definition~\ref{defn:time_steps} is that we also introduce the intermediate time labeled by $1/2$ in the subscript and the corresponding quantities. Moreover, we also modify the split $A(T^{(L)}) = \Asmoothpartmod(T^{(L)}) + \Aroughpartmod(T^{(L)})$, as in Definition~\ref{defn:final_timestep}. See Figure~\ref{fig:iteration_scheme_mod} for an illustration. This modification will change all the quantities appearing in the previous iteration scheme. We add a tilde to every quantity to address this. 

\begin{definition}
\label{defn:time_steps_mod}
	Suppose at time $0$ we are given $\Lmod_0$ and $(\Asmoothpart_0, \Aroughpart_0, \phismoothpart_0, \phiroughpart_0)$, such that the latter is $\Lmod_0$ admissible. In the following, we will define a sequence of stopping times
	\begin{equ}
		 0 = \tmod_0 \leq \tmod_{1/2} \leq \tmod_1 \leq \tmod_{3/2} \leq \tmod_2 \leq \cdots \leq \tmod_m \leq \tmod_{m+1/2} \leq \tmod_{m+1} \leq \cdots
	\end{equ}
	as well as the time differences $\taumod_m = \tmod_m - \tmod_{m-1}$, their deterministic comparison versions $\tmod_m^\star$, $\tmod_{m-1/2}^\star$, $\taumod_m^\star$, and data $\Lmod_m \in \N^+$ for any $m \geq 1$. Moreover, on $(\tmod_{m-1}, \tmod_m]$, we will define functions
	\begin{equ}
		\widetilde{g}_m \;, \quad \Bmod_m \;, \quad \YMHFsmoothmod_m \;, \quad \Sdetmod_m \;, \quad \higgsdetmod_m \;, \quad
		\linear[\Bmod_m] \;, \quad \linear[\Bmod_m, 1/2] \;, \quad \philinear[\Bmod_m] \;,
	\end{equ}
	and decompositions
	\begin{equ}
		A = \Asmoothpartmod + \Aroughpartmod \;, \qquad \higgs = \phismoothpartmod + \phiroughpartmod \;.
	\end{equ}
	Finally, on $(\tmod_{m-1}, \tmod_m^\star]$, we will define their starred versions
	\begin{equ}
		\widetilde{g}_m^\star \;, \quad \Bmod_m^\star \;, \quad \Bmod_m^{+, \star} \;, \quad \Sdetmod_m^\star \;, \quad \higgsdetmod_m^\star \;, \quad
		\linear[\Bmod_m^\star] \;, \quad \linear[\Bmod_m^\star, 1/2] \;, \quad \philinear[\Bmod_m^\star] \;.
	\end{equ}
	
	To begin, suppose the decompositions
	\begin{equ}
		A(t) = \Asmoothpartmod(t) + \Aroughpartmod(t) \;, \qquad
		\higgs(t) = \phismoothpartmod(t) + \phiroughpartmod(t)
	\end{equ}
	are defined up to the stopping time $\tmod_m$, and are $\fF_{\tmod_m}$-measurable and $\Lmod_m$-admissible at time $\tmod_m$. We define
	\begin{equs}
		\taumod_{m+1}^{\star} &\eqdef T^{(\Lmod_m)}\bigl(\Asmoothpartmod(\tmod_m), \Aroughpartmod(\tmod_m), \phismoothpartmod(\tmod_m), \phiroughpartmod(\tmod_m)\bigr) \;, \quad
		\tmod_{m+1}^{\star} = \tmod_m + \taumod_{m+1}^{\star} \;, \\
		\taumod_{m+1/2}^{\star} &\eqdef T^{(\Lmod_m)}_{1/2}\bigl(\Asmoothpartmod(\tmod_m), \Aroughpartmod(\tmod_m), \phismoothpartmod(\tmod_m), \phiroughpartmod(\tmod_m)\bigr) \;, \quad
		\tmod_{m+1/2}^{\star} = \tmod_m + \taumod_{m+1/2}^{\star} \;.
	\end{equs}
	Let $\widetilde{g}_{m+1}^\star$ be the gauge transform in Proposition~\ref{pr:admissible_timescale}, shifted to start at time $\tmod_m$ for the decomposition
	$A(\tmod_m) = \Asmoothpartmod(\tmod_m) + \Aroughpartmod(\tmod_m), \higgs(\tmod_m) = \phismoothpartmod(\tmod_m) + \phiroughpartmod(\tmod_m)$, defined on $[\tmod_m, \tmod_{m+1}^\star]$. 
	For $t \in (\tmod_m, \tmod_{m+1}^{\star}]$, let $\Bmod_{m+1}^{+, \star}(t)$ be the YM gradient flow starting from $\Bmod_{m+1}^{+, \star}(\tmod_m) = \Asmoothpartmod(\tmod_m)$, and $\Bmod^\star_{m+1}(t)$ be the YM gradient flow starting from $\Bmod^\star_{m+1}(\tmod_m) = \widetilde{g}_{m+1}^\star(\tmod_m) \act A(\tmod_m)$. 
	Then, we define $\Sdetmod^\star_{m+1}, \higgsdetmod^\star_{m+1}$ to be the solution to \eqref{eq:Sdet-equation}--\eqref{eq:higgsdet-equation}, shifted to start at time $\tmod_m$, with $B$ replaced by $\Bmod^\star_{m+1}$ and initial data
	\begin{equ}
		\Sdetmod^\star_{m+1}(\tmod_m) = 0 \;, \qquad
		\higgsdetmod^\star_{m+1}(\tmod_m) = \widetilde{g}_{m+1}^\star(\tmod_m) \act \bigl(\phismoothpartmod(\tmod_m) + \phiroughpartmod(\tmod_m) \bigr) \;.
	\end{equ}
        	We let $\linear[\Bmod_{m+1}^{\star}]$, $\linear[\Bmod_{m+1}^{\star}, 1/2]$ and $\philinear[\Bmod_{m+1}^{\star}]$, starting from $\tmod_m$, $\tmod_{m+1/2}^{\star}$, $\tmod_m$ respectively, 
	solve
        \begin{equs}
        	\partial_t \linear[\Bmod_{m+1}^{\star}]
        	&= \Delta_{\Bmod_{m+1}^{\star}} \linear[\Bmod_{m+1}^{\star}]
        	+ 2 \ad(F_{\Bmod_{m+1}^{\star}})\linear[\Bmod_{m+1}^{\star}] + \xi \;,
        	&\quad \linear[\Bmod_{m+1}^{\star}](\tmod_m)&=0 \;, \\
        	\partial_t \linear[\Bmod_{m+1}^{\star}, 1/2]
        	&= \Delta_{\Bmod_{m+1}^{\star}} \linear[\Bmod_{m+1}^{\star}, 1/2]
        	+ 2 \ad(F_{\Bmod_{m+1}^{\star}})\linear[\Bmod_{m+1}^{\star}, 1/2] + \xi \;,
        	&\quad \linear[\Bmod_{m+1}^{\star}, 1/2](\tmod_{m+1/2}^{\star})&=0 \;, \\
        	\partial_t \philinear[\Bmod_{m+1}^{\star}]
        	&= \Delta_{\Bmod_{m+1}^{\star}} \philinear[\Bmod_{m+1}^{\star}] + \zeta \;,
        	&\quad \philinear[\Bmod_{m+1}^{\star}](\tmod_m)&=0 \;.
        \end{equs}
	Note that the previously defined quantities with $\star$ are all measurable with respect to $\fF_{\tmod_m}$, except for $\linear[\Bmod_{m+1}^{\star}]$, $\linear[\Bmod_{m+1}^{\star}, 1/2]$ and $\philinear[\Bmod_{m+1}^{\star}]$. Now we define more quantities that depend on the noise after time $\tmod_m$. 
	First, we define the stopping times
	\begin{equs}
		\tmod_{m+1/2} &= \tmod_m + \tau^{(\Lmod_m)}_{1/2}\bigl(\Asmoothpartmod(\tmod_m), \Aroughpartmod(\tmod_m), \phismoothpartmod(\tmod_m), \phiroughpartmod(\tmod_m); \shift_{\tmod_m}\xi, \shift_{\tmod_m}\zeta\bigr) \;, \\
		\tmod_{m+1} &= \tmod_m + \tau^{(\Lmod_m)}\bigl(\Asmoothpartmod(\tmod_m), \Aroughpartmod(\tmod_m), \phismoothpartmod(\tmod_m), \phiroughpartmod(\tmod_m); \shift_{\tmod_m}\xi, \shift_{\tmod_m}\zeta\bigr) \;.
	\end{equs}
	We set $\taumod_{m+1} \eqdef \tmod_{m+1} - \tmod_m$. 
	Let $\widetilde{g}_{m+1}, \Bmod_{m+1}, \YMHFsmoothmod_{m+1}, \Sdetmod_{m+1}, \higgsdetmod_{m+1}, \linear[\Bmod_{m+1}], \linear[\Bmod_{m+1}, 1/2]$ and $\philinear[\Bmod_{m+1}]$ be the restrictions of the corresponding starred quantities to the time interval $(\tmod_m, \tmod_{m+1}]$.
	For $t \in (\tmod_m, \tmod_{m+1}]$, we define
	\begin{equ}
		\bigl(\Asmoothpartmod(t), \Aroughpartmod(t), \phismoothpartmod(t), \phiroughpartmod(t)\bigr)
		\eqdef \wt{\flow}_{t-\tmod_m}\bigl(\Asmoothpartmod(\tmod_m), \Aroughpartmod(\tmod_m), \phismoothpartmod(\tmod_m), \phiroughpartmod(\tmod_m); \shift_{\tmod_m}\xi, \shift_{\tmod_m}\zeta \bigr) \;.
	\end{equ}
	We also set
	\begin{equ}
	\label{e:Alinmod_defn}
		\Alinmod(\tmod_{m+1}) =
		\begin{cases}
			\linear[\Bmod_{m+1}, 1/2](\tmod_{m+1}) \;, &\quad \text{if } \tmod_{m+1/2} \neq \tmod_m \text{ and } \taumod_{m+1} = \taumod_{m+1}^\star \;, \\
			\linear[\Bmod_{m+1}](\tmod_{m+1}) \;, &\quad \text{otherwise} \;, 
		\end{cases}
	\end{equ}
	and define
	\begin{equ}
		\Aremmod(\tmod_{m+1}) \eqdef \Aroughpartmod(\tmod_{m+1}) - \Alinmod(\tmod_{m+1}) \;, \qquad
		\phiremmod(\tmod_{m+1}) \eqdef \phiroughpartmod(\tmod_{m+1}) - \philinear[\Bmod_{m+1}](\tmod_{m+1}) \;.
	\end{equ}
	Finally, choose the smallest integer $\Lmod_{m+1}$ such that
	\begin{equ}
	\label{e:Lmod_m_lower_bound}
		\Lmod_{m+1} \geq 1 \vee \bigl(\covariantenergy(\Asmoothpartmod(\tmod_{m+1}), \phismoothpartmod(\tmod_{m+1})) \bigr)^{\frac{1}{10}} \;,
	\end{equ}
	and such that $\bigl(\Asmoothpartmod(\tmod_{m+1}), \Aroughpartmod(\tmod_{m+1}), \phismoothpartmod(\tmod_{m+1}), \phiroughpartmod(\tmod_{m+1}) \bigr)$ is $\Lmod_{m+1}$-admissible. 
	If $\limsup_{m \rightarrow \infty} \tmod_m \neq \infty$, then, for any $t \geq \limsup_{m \rightarrow \infty} \tmod_m$, we define $\bigl(\Asmoothpartmod(t), \Aroughpartmod(t), \phismoothpartmod(t), \phiroughpartmod(t) \bigr) = \skull$.
\end{definition}

\afterpage{\clearpage}

We now state the main result of this section, which says that we have uniform-in-time control of the iteration scheme that we just introduced.

\begin{proposition}\label{prop:modified-iteration-uniform-bounds}
For $\stepupperbound$ sufficiently small (possibly depending on the parameters in \eqref{e:small_constants} and \eqref{e:constants_global_iteration}), the following holds. Let $p \geq 1$. There exists a constant $C_p$ which depends on $p$ and $\stepupperbound$ such that the following holds. Let $A_0^+ \in C_x^\infty(\T^2, \frkg^2), \higgs_0^+ \in C_x^\infty(\T^2, \V)$, $A_0^- \in \Cs_x^{-\kappa}(\T^2, \frkg^2)$, $\higgs_0^- \in \Cs_x^{-\kappa}(\T^2, \V)$. Let $\Lmod_0 \geq 1$ be large enough so that $(A_0^+, A_0^-, \higgs_0^+, \higgs_0^-)$ is $\Lmod_0$-admissible. Consider the iteration scheme from Definition~\ref{defn:time_steps_mod}. We have that
\begin{equs}
\sup_{m \geq 0} \E\bigl[ \covariantenergy(\Asmoothpartmod(\tmod_m), \phismoothpartmod(\tmod_m))^p ~\big|~ \mc{F}_0\bigr] &\leq C_p \bigl(\covariantenergy(A_0^+, \higgs_0^+)^{4 p} + \Lmod_0^{4 p}\bigr) \;, \\
\limsup_{m \toinf} \E\bigl[ \covariantenergy(\Asmoothpartmod(\tmod_m), \phismoothpartmod(\tmod_m))^{p} ~\big|~ \mc{F}_0\bigr] &\leq C_p \;.
\end{equs}
Moreover, for $t > 0$, let $M_t = \inf\{m \geq 0 : \tmod_m \geq t\}$. Then $M_t < \infty$ a.s.\ for each $t > 0$, and we have that
\begin{equs}
\sup_{t > 0} \E\bigl[ \covariantenergy(\Asmoothpartmod(\tmod_{M_t+1}), \phismoothpartmod(\tmod_{M_t+1}))^p ~\big|~ \mc{F}_0\bigr] & \leq C_p \bigl(\covariantenergy(A_0^+, \higgs_0^+)^{4 p} + \Lmod_0^{4 p}\bigr) \;, \\
\limsup_{t \toinf}  \E\bigl[ \covariantenergy(\Asmoothpartmod(\tmod_{M_t+1}), \phismoothpartmod(\tmod_{M_t+1}))^p ~\big|~ \mc{F}_0\bigr] & \leq C_p \;.
\end{equs}
\end{proposition}

The main parts of our hard global well-posedness result Theorem~\ref{thm:global_hard} will follow directly from Proposition~\ref{prop:modified-iteration-uniform-bounds}. The estimates at the stopping times will be needed in Section~\ref{sec:conclusions} to prove our soft global well-posedness result Theorem~\ref{thm:global_soft}. Both theorems will be proven together in Section~\ref{sec:conclusions}. The remainder of this section is devoted to the proof of Proposition~\ref{prop:modified-iteration-uniform-bounds}. We next state the main intermediate step towards this proposition.

\begin{proposition}
\label{pr:main_iteration}
	There exists some constant $K$ such that the following holds. Consider the iteration scheme from Definition~\ref{defn:time_steps_mod}. For any $m \geq 1$, 
	\begin{equ}
		\covariantenergy\bigl(\widetilde{A}^{+}(\tmod_{m+1}), \phismoothpartmod(\tmod_{m+1})\bigr) \leq \covariantenergy\bigl(\widetilde{A}^{+}(\tmod_m), \phismoothpartmod(\tmod_m)\bigr) - 3C_{\eqref{e:one_step_change}}  \bigl( \taumod_{m+1}^\star \bigr)^{-1} + \fluctmod_{m+1} + K \;,
	\end{equ}
	where $\fluctmod_{m+1}$ is some $\fF_{\tmod_{m+1}}$ measurable random variable, such that for any $p \geq 1$, there exists some $C_p$ independent of $m$ such that
	\begin{equ}
	\label{e:main_iteration_fluct}
	\E \Bigl( \Bigl( \fluctmod_{m+1} - 2 C_{\eqref{e:one_step_change}} \bigl( \taumod_{m+1}^\star \bigr)^{-1} \Bigr)_+^p \Big| \fF_{\tmod_{m-1}} \Bigr) \leq C_p \;.
	\end{equ}
	Moreover, for the initial step, we have that 
	\begin{equ}
	\label{e:main_iteration_fluct_init}
		\E \Bigl[ \Bigl( \covariantenergy\bigl(\Asmoothpartmod(\tmod_1), \phismoothpartmod(\tmod_1)\bigr) \Bigr)^p \Bigm| \fF_0 \Bigr] 
		\leq C_p \bigl(\covariantenergy(A_0^+, \higgs_0^+)^{4p} + \Lmod_0^{4p}\bigr) \;.
	\end{equ}
\end{proposition}

Proposition~\ref{prop:modified-iteration-uniform-bounds} essentially follows from Proposition~\ref{pr:main_iteration} and the abstract supermartingale estimates of Section~\ref{section:supermartingale}. However, in order to verify some of the conditions in the abstract supermartingale estimates, we will need several minor lemmas which are only introduced later on. Thus, we choose to defer the proof of Proposition~\ref{prop:modified-iteration-uniform-bounds} to the end of this section. We now work towards the proof of Proposition~\ref{pr:main_iteration}. We first collect several basic facts about the iteration scheme. 

\begin{lemma}\label{lemma:time-step-cond-exp-bd}
Let $p \geq 1$. There is some constant $C_p > 0$ such that the following holds. For all $m \geq 0$, we have that
\begin{equ}
\E[ (\taumod_{m+1})^{-p} ~|~ \mc{F}_{\tmod_m}] \leq (\taumod_{m+1}^\star)^{-p} + C_p \;.
\end{equ}
\end{lemma}

\begin{remark}
From the proof of Lemma~\ref{lemma:time-step-cond-exp-bd}, we will see that we may take the constant $C_p$ as
\begin{equ}
C_p = \frac{p}{\lambda_1 c^{1 + p / \lambda_1}} \Gamma(p / \lambda_1) \;,
\end{equ}
where $\Gamma$ is the usual Gamma function, and the constant $c$ is from Lemma~\ref{le:tau_tail_bound}.
\end{remark}

\begin{proof}
For notational convenience, let $T = \taumod_{m+1}^\star$ and $\mc{F}_m = \mc{F}_{\tmod_m}$. By Lemma~\ref{le:tau_tail_bound}, we have that
\begin{equs}
\E[ (\taumod_{m+1})^{-p} ~|~ \mc{F}_m] &= \int_0^\infty p x^{p-1} \P( (\taumod_{m+1})^{-1} > x ~|~ \mc{F}_m) \,\dx \\
&\leq T^{-p} + c^{-1} \int_{T^{-1}}^\infty p x^{p-1} \exp(-c x^{\lambda_1}) \, \dx 
\leq T^{-p} + C_p \;,
\end{equs}
as desired.
\end{proof}

Next, we show that inverse powers of the starred versions $\taumod_m^\star$ of our iteration time steps may be estimated in terms of the covariant energy.

\begin{lemma}\label{lemma:taumod-star-det-bound}
Consider the iteration scheme from Definition~\ref{defn:time_steps_mod}. For any $m \geq 0$,
we have that
\begin{equ}
(\taumod_{m+1}^\star)^{-1} \leq \frac{2}{\logtimestep} \covariantenergy(\Asmoothpartmod(\tmod_m), \phismoothpartmod(\tmod_m)) + \stepupperbound^{-1} + \logtimestep^{-1} + \Lmod_m \;.
\end{equ}
\end{lemma}
\begin{proof}
Due to the inductive definition of the iteration scheme, it suffices to take $m = 0$, in which case $\tmod_m = \tmod_0 = 0$, and $\Asmoothpartmod(0) = A_0^+, \phismoothpartmod(0) = \higgs_0^+$, $\Aroughpartmod(0) = A_0^-, \phiroughpartmod(0) = \higgs_0^-$. By definition, 
\begin{equ}
\taumod_1^\star = T^{(\Lmod_0)}(A_0^+, A_0^-, \higgs_0^+, \higgs_0^-) = (T^{(1)} \wedge T^{(2)})(A_0^+, \higgs_0^+) \wedge \Lmod_0^{-1} \;. 
\end{equ}
By Definition~\ref{defn:tau} and \eqref{e:timestep2_lower_bound}, we have that
\begin{equ}
\bigl((T^{(1)} \wedge T^{(2)})(A_0^+, \higgs_0^+) \bigr)^{-1} \leq \frac{2}{\logtimestep} \covariantenergy(A_0^+, \Phi_0^+) + \stepupperbound^{-1} + \logtimestep^{-1},
\end{equ}
where we used that $\log 2 > 1/2$ and that $\gamma_0 < \gamma_2$. The desired result now follows.
\end{proof}

Note that now $\Asmoothpartmod(\tmod_{m})$ might include the extra term 
\begin{equ}
	e^{(\tmod^{\star}_{m} - \tmod^{\star}_{m - 1/2})\Delta^F_{\Bmod^{\star}_m(\tmod^{\star}_{m - 1/2})}} \linear[\Bmod^{\star}_m](\tmod^{\star}_{m-1/2})  \;. 
\end{equ}
Therefore, the remainder bound Proposition~\ref{pr:remainder_estimate} needs to be modified correspondingly. 
\begin{lemma}
\label{le:remainder_estimate_mod}
	In the setting of Definition~\ref{defn:time_steps_mod}, for any $m \geq 1$ and $2\lossdel \leq \beta
	< 1-\lambda$, we have
	\begin{equs}
		\| \Aremmod(\tmod_m) \|_{\covHolder{\Asmoothpartmod(\tmod_m),F}^{\beta-2\lossdel}}
		&\lesssim (\taumod_m)^{\frac{1-\lambda}{2}-\frac{1}{\gamma_0}-\lambda_2}
		\Bigl((\taumod_m)^{-\frac{\beta}{2}}
		+\bigl((T^{(1)}\wedge T^{(2)})(\Asmoothpartmod(\tmod_m),
		\phismoothpartmod(\tmod_m))\bigr)^{-\frac{\beta}{2}}\Bigr) \;, \\
		\| \phiremmod(\tmod_m) \|_{\covHolder{\Asmoothpartmod(\tmod_m)}^{\beta-2\lossdel}}
		&\lesssim (\taumod_m)^{\frac{1-\lambda-\beta}{2}-\frac{1}{\gamma_0}-\lambda_2} \;.
	\end{equs}
\end{lemma}

\begin{proof}
	If $\tmod_{m - 1/2} = \tmod_{m-1}$ or $\taumod_m \neq \taumod_m^\star$, then the splitting at $\tmod_m$ is not modified, and the result follows directly from Proposition~\ref{pr:remainder_estimate}. Therefore, we only discuss the case when $\tmod_{m - 1/2} \neq \tmod_{m-1}$ and $\taumod_m = \taumod_m^\star$, and in this case we have
	\begin{equ}
		\Asmoothpartmod(\tmod_m) = \Bmod_m(\tmod_m) + \Sdetmod_m(\tmod_m) + 
	e^{(\tmod^{\star}_{m} - \tmod^{\star}_{m - 1/2})\Delta^F_{\Bmod^{\star}_m(\tmod^{\star}_{m - 1/2})}} \linear[\Bmod^{\star}_m](\tmod^{\star}_{m-1/2})  \;. 
	\end{equ}
	Following the proof of
	Proposition~\ref{pr:remainder_estimate}, we have a decomposition $ \Aremmod(\tmod_m) = \sum_{\substack{N \geq N_{\lo}\\ N \in \dyadic}} \Aremmod[N](\tmod_m)$
	such that, for $0 \leq j \leq 1$,
	\begin{equ}
		\bigl\|\covd_{\Bmod_m(\tmod_m)}^{(j)}
		\Aremmod[N](\tmod_m)\bigr\|_{L^\infty}
		\lesssim N^{j-(1-\lambda)}
		\bigl(\taumod_m\bigr)^{-\frac{1}{\gamma_0}-\lambda_2} \;,
	\end{equ}
	where $N_{\lo}\sim\bigl(\taumod_m\bigr)^{-1/2}$. Moreover, by the definition of our stopping time (see Definition~\ref{defn:tau_L}), we have 
	\begin{equs}
		\|e^{(\tmod^{\star}_{m} - \tmod^{\star}_{m - 1/2})\Delta^F_{\Bmod^{\star}_m(\tmod^{\star}_{m - 1/2})}} \linear[\Bmod^{\star}_m](\tmod^{\star}_{m-1/2})   \|_{L^\infty}
		&\lesssim (\tmod^{\star}_m-\tmod^{\star}_{m-1/2})^{-\frac{\kappa}{2}}\|\linear[\Bmod^{\star}_m](\tmod^{\star}_{m-1/2})\|_{\covHolder{\Bmod^{\star}_m(\tmod^{\star}_{m-1/2}),F,\taumod_m^{\star}}^{-\kappa}} \\
		\label{e:smoothed_lollipop_bound0}
		&\lesssim (\taumod^{\star}_m)^{- \frac{\kappa}{2} - \lambda_1 } \;.
	\end{equs}
	Here from the first line to the second line we used the fact that we are in the case of $\tmod_m = \tmod_m^{\star}$ and $\tmod_{m - 1/2} = \tmod_{m - 1/2}^{\star}$, and by definition (Definition~\ref{defn:time_steps_mod}) we have $\tmod_{m}^{\star} - \tmod_{m - 1/2}^{\star} = \frac12 \taumod_{m}$. 
	
	Therefore, an argument similar to that in Proposition~\ref{pr:remainder_estimate} shows 
	\begin{equs}
		\,& \bigl\|\covd_{\Asmoothpartmod(\tmod_m)} \Aremmod[N](\tmod_m)\bigr\|_{L^\infty} \\
		&\leq \bigl\|\covd_{\Bmod_m(\tmod_m)} \Aremmod[N](\tmod_m)\bigr\|_{L^\infty} \\
		&\quad+\Bigl(\|\Sdetmod_m(\tmod_m)\|_{L^\infty}
		+\|e^{(\tmod^{\star}_{m} - \tmod^{\star}_{m - 1/2})\Delta^F_{\Bmod^{\star}_m(\tmod^{\star}_{m - 1/2})}} \linear[\Bmod^{\star}_m](\tmod^{\star}_{m-1/2})  \|_{L^\infty}\Bigr) \|\Aremmod[N](\tmod_m)\|_{L^\infty} \\
		&\lesssim N^\lambda \bigl(\taumod_m\bigr)^{-\frac{1}{\gamma_0}-\lambda_2} + \Bigl(
		\bigl(\taumod_m\bigr)^{-\lossdel+\frac12-\frac2r-\frac{2}{\gamma_0}} + \bigl(\taumod_m\bigr)^{-\frac{\kappa}{2} - \lambda_1} \Bigr) N^{-(1-\lambda)}
		\bigl(\taumod_m\bigr)^{-\frac{1}{\gamma_0}-\lambda_2} \\
		&\lesssim N^\lambda \bigl(\taumod_m\bigr)^{-\frac{1}{\gamma_0}-\lambda_2} \;. 
	\end{equs}
	The desired result then follows from Lemma~\ref{le:heat_flow_norm_sol_norm}, together with $N_{\lo}\sim\bigl(\taumod_m\bigr)^{-1/2}$. 
	The bound for $\phiremmod(\tmod_m)$ follows in the same way.
\end{proof}

	With Lemma~\ref{le:remainder_estimate_mod}, the modified version of Lemma~\ref{le:iteration_basic_bounds} follows, where one simply replaces all quantities by their tilded versions.
\begin{lemma}
\label{le:iteration_basic_bounds_mod}
	In the setting of Definition~\ref{defn:time_steps_mod}, for every $m \geq 0$, the shifted versions of \eqref{e:iteration_detarministic_bounds_connection_assumption}, \eqref{e:phi_0_smallness_global} and Assumptions~\ref{local:assumption-B} and \ref{local:assumption-higgs-initial} hold on $(\tmod_m, \tmod_{m+1}^\star]$. Furthermore, for every $m \geq 0$, we have
	\begin{equ}
	\label{e:iteration_F_bounds1_mod}
		\sup_{s \in (0, \taumod_{m+1}^\star]} s^{1 + \fa} \Bigl(\|F_{\Bmod_{m+1}^{+, \star}}(\tmod_m + s)\|_{L^\infty} + \|F_{\Bmod_{m+1}^\star}(\tmod_m + s)\|_{L^\infty}\Bigr) \lesssim 1 \;,
	\end{equ}
	and
	\begin{equs}
		\sup_{s \in [0, \taumod_{m+1}^\star]} \|F_{\Bmod_{m+1}^{+, \star}}(\tmod_m + s)\|_{L^2}
		&\lesssim \bigl((T^{(1)} \wedge T^{(2)})(\Asmoothpartmod(\tmod_m), \phismoothpartmod(\tmod_m))\bigr)^{-1 - \fa} \\
	\label{e:iteration_F_bounds2_mod}
		&\leq \bigl(\taumod_{m+1}^\star\bigr)^{-1 - \fa} \;, \\
		\|\covd_{\Bmod_{m+1}^{+, \star}} F_{\Bmod_{m+1}^{+, \star}}\|_{L^2_{t,x}([\tmod_m, \tmod_{m+1}^\star])}
		&\lesssim \bigl((T^{(1)} \wedge T^{(2)})(\Asmoothpartmod(\tmod_m), \phismoothpartmod(\tmod_m))\bigr)^{-\frac12 - \lossdel} \wedge \|F_{\Asmoothpartmod(\tmod_m)}\|_{L^2} \\
	\label{e:iteration_F_endtime_lowerbound_mod}
		&\leq \bigl(\taumod_{m+1}^\star\bigr)^{-\frac12 - \lossdel} \wedge \|F_{\Asmoothpartmod(\tmod_m)}\|_{L^2} \;, 
	\end{equs}
	For the remainder parts, for every $m \geq 1$ we have
	\begin{equs}
	\label{e:Aroughpart_crude_bound_mod}
		\|\Alinmod(\tmod_m)\|_{\covHolder{\Asmoothpartmod(\tmod_m),F}^{-\kappa}} + \|\Aroughpartmod(\tmod_m)\|_{\covHolder{\Asmoothpartmod(\tmod_m),F}^{-\kappa}}
		&\lesssim \Lmod_m^\fa \;, \qquad \|\Aremmod(\tmod_m)\|_{\covHolder{\Asmoothpartmod(\tmod_m),F}^{0}} \lesssim \taumod_m^{100 \fa} \;, \\
	\label{e:phiroughpart_crude_bound_mod}
		\|\philinear[\Bmod_m](\tmod_m)\|_{\covHolder{\Asmoothpartmod(\tmod_m)}^{-\kappa}} + \|\phiroughpartmod(\tmod_m)\|_{\covHolder{\Asmoothpartmod(\tmod_m)}^{-\kappa}}
		&\lesssim \Lmod_m^\fa \;, \qquad \|\phiremmod(\tmod_m)\|_{\covHolder{\Asmoothpartmod(\tmod_m)}^{0}} \lesssim \taumod_m^{100 \fa} \;.
	\end{equs}
	Moreover, for every $m \geq 1$ we have the following relation between the lengths of two time intervals:
	\begin{equs}
	\label{e:two_timesteps_relation_mod}
		\taumod_m
		&\lesssim \stepupperbound \wedge \logtimestep \wedge \|F_{\Asmoothpartmod(\tmod_m)}\|_{L^2}^{-\frac12} \wedge \|\phismoothpartmod(\tmod_m)\|_{L^\infty}^{-\frac{1}{\frac{1}{\gamma_0} + \frac1r + \lossdel}} \\
		&\lesssim \bigl((T^{(1)} \wedge T^{(2)})(\Asmoothpartmod(\tmod_m), \phismoothpartmod(\tmod_m))\bigr)^{\frac{1}{10}} \;.
	\end{equs}
\end{lemma}
\begin{proof}
	The only part whose proof differs from that of Lemma~\ref{le:iteration_basic_bounds} is the fact that $ \|F_{\Asmoothpartmod(\tmod_m)}\|_{L^2} \lesssim \taumod_m^{-2 }$ in \eqref{e:two_timesteps_relation_mod}, since now the extra noise term $e^{(\tmod_m^\star - \tmod_{m-1/2}^\star)\Delta^F_{\smash{\Bmod_m^\star(\tmod_{m-1/2}^\star)}}} \linear[\Bmod_m^\star](\tmod_{m-1/2}^\star)$ can be added to $\Asmoothpartmod(\tmod_m)$ when $\tmod^\star_{m-1/2} \neq \tmod_{m-1}$ and $\taumod_m = \taumod_m^\star$. In this case, this bound will be proved in Lemma~\ref{le:smoothed_lollipop_addition_expansion} that will appear soon. 
\end{proof}

Now we state the parallel version of Proposition~\ref{pr:expansion_iteration}. 
\begin{proposition}
\label{pr:expansion_iteration_mod}
	There exists some constant $C$ such that in the setting of Definition~\ref{defn:time_steps_mod}, for any $m \geq 1$ we have
	\minilab{e:iteration_deterministic_bounds_iteration_mod}
\begin{equs}
	\covariantenergy&\big((\Bmod_{m+1} + \Sdetmod_{m+1})(\tmod_{m+1}), \higgsdetmod_{m+1}(\tmod_{m+1}) \big) \\
	&\leq \covariantenergy(\Asmoothpartmod(\tmod_{m}), \phismoothpartmod(\tmod_{m})) - \frac14 \Decmod_{m+1} + C\bigl(\taumod_{m+1} \bigr)^{-1 + \kappa} \\
	\label{e:iteration_deterministic_bounds_iteration1_mod}
	&\quad+  \Big \| \curl_{\YMHFsmoothmod_{m+1}(\tmod_{m+1})} e^{\taumod_{m+1} \Delta^F_{\Asmoothpartmod(\tmod_{m})}} \Alinmod(\tmod_m) \Big\|_{L^2}^2 \\
	\label{e:iteration_deterministic_bounds_iteration2_mod}
	&\quad+ 2 \int_{\T^2} \Big \langle F_{\YMHFsmoothmod_{m+1}}^{12}(\tmod_{m+1}), \curl_{\YMHFsmoothmod_{m+1}(\tmod_{m+1})} e^{\taumod_{m+1} \Delta^F_{\Asmoothpartmod(\tmod_{m})}} \Alinmod(\tmod_m) \Big \rangle_{\fg} \, \dx \\
	\label{e:iteration_deterministic_bounds_iteration3_mod}
	&\quad + 2 \int_{\T^2} \Big \langle F_{\YMHFsmoothmod_{m+1}}^{12}(\tmod_{m+1}), \Big[( e^{\taumod_{m+1} \Delta^F_{\Asmoothpartmod(\tmod_{m})}} \Alinmod(\tmod_m))^1, ( e^{\taumod_{m+1} \Delta^F_{\Asmoothpartmod(\tmod_{m})}} \Alinmod(\tmod_m))^2 \Big]  \Big \rangle_{\fg} \, \dx\qquad \\
	&\quad+ C \| \phismoothpartmod(\tmod_m) \|_{L^4}^{\gamma_1 - 4} \Big| \int_{\T^2} |e^{\taumod_{m+1} \Delta_{\Asmoothpartmod(\tmod_m)}} \phismoothpartmod(\tmod_m) |^2 \\
	\label{e:iteration_deterministic_bounds_iteration4_mod}
	& \qquad \qquad \qquad \qquad \qquad \cdot \langle e^{\taumod_{m+1} \Delta_{\Asmoothpartmod(\tmod_m)}} \phismoothpartmod(\tmod_m), e^{\taumod_{m+1} \Delta_{\Asmoothpartmod(\tmod_m)}} \philinear[\Bmod_m](\tmod_m) \rangle_{\V} \,  \dx \Big| \;, 
\end{equs}
where 
\begin{equ}
\Decmod_{m+1} \eqdef \Dec_{\taumod_{m+1}}(\Asmoothpartmod(\tmod_{m}), \phismoothpartmod(\tmod_{m})) \;.
\end{equ}
We use $\itobjectsmod^{(1)}_{m+1} $ to denote the sum of \eqref{e:iteration_deterministic_bounds_iteration1_mod}--\eqref{e:iteration_deterministic_bounds_iteration4_mod}. 

When $m = 0$, we record the crude bound 
\begin{equ}
\label{e:iteration_deterministic_bounds_iteration_mod_init}
	\covariantenergy \big((\Bmod_1 + \Sdetmod_1)(\tmod_1), \higgsdetmod_1(\tmod_1) \big) 
	\leq \covariantenergy(\Asmoothpart_0, \phismoothpart_0) - \frac14 \Decmod_1 + C\bigl(\taumod_1\bigr)^{-\frac74} \;.
\end{equ}
\end{proposition}
\begin{proof}
We note that Proposition~\ref{pr:iteration_detarministic_bounds} still applies to the modified iteration scheme. Then it suffices to show a parallel version of Proposition~\ref{pr:expansion_iteration}, in which we bound the errors when replacing $\Aroughpartmod(\tmod_m)$ by $\Alinmod(\tmod_m)$. Actually, the two main ingredients for this replacement are Proposition~\ref{pr:remainder_estimate} and Lemma~\ref{le:iteration_basic_bounds}, which are guaranteed now by Lemma~\ref{le:remainder_estimate_mod} and Lemma~\ref{le:iteration_basic_bounds_mod}. 
This concludes the proof.
\end{proof}

Note that in the modified scheme, when $\tmod_{m + 1/2}^{\star} \neq \tmod_{m}$ and $\tmod_{m+1} = \tmod_{m+1}^{\star}$, $\Asmoothpartmod(\tmod^{\star}_{m+1})$ is not $(\Bmod_{m+1} + \Sdetmod_{m+1})(\tmod_{m+1})$, since we need to add part of the linear stochastic object to it to move it away from the near-YM connection $\Bmod_{m+1}$. 
The following lemma performs another round of perturbative expansion, to single out the extra singular terms in the $ \|F_{\Asmoothpartmod(\tmod^{\star}_{m+1})}\|_{L^2}^2$ due to this modification. 
\begin{lemma}
\label{le:smoothed_lollipop_addition_expansion}
	There exists some constant $C$ such that the following holds. In the setting of Definition~\ref{defn:time_steps_mod},
	for any $m \geq 0$, if $\tmod_{m + 1/2}^{\star} \neq \tmod_{m}$ and $\tmod_{m+1} = \tmod_{m+1}^{\star}$, then 
	\minilab{e:smoothed_lollipop_addition_expansion}
	\begin{equs}
	\,& \|F_{\widetilde{A}^{+, \star}(\tmod^{\star}_{m+1})}\|_{L^2}^2 \\
	&\leq \|F_{(\Bmod^{\star}_{m+1}+\Sdetmod^{\star}_{m+1})(\tmod^{\star}_{m+1})}\|_{L^2}^2 
	+C \bigl( \taumod_{m+1}^{\star} \bigr)^{-\frac12 - \frac{3 \kappa}{2} - 3 \lambda_1 } \\
	\label{e:smoothed_lollipop_addition_expansion1}
	&\quad+2\int_{\T^2}\Big\langle F^{12}_{(\Bmod^{\star}_{m+1}+\Sdetmod^{\star}_{m+1})(\tmod^{\star}_{m+1})},
	\Big[
	\Bigl(e^{(\tmod^{\star}_{m+1}-\tmod^{\star}_{m+1/2})\Delta^F_{\Bmod^{\star}_{m+1}(\tmod^{\star}_{m+1/2})}}\linear[\Bmod^{\star}_{m+1}](\tmod^{\star}_{m+1/2})\Bigr)^1, \\
	&\qquad \qquad \qquad \qquad \qquad \qquad \qquad \qquad
	\Bigl(e^{(\tmod^{\star}_{m+1}-\tmod^{\star}_{m+1/2})\Delta^F_{\Bmod^{\star}_{m+1}(\tmod^{\star}_{m+1/2})}}\linear[\Bmod^{\star}_{m+1}](\tmod^{\star}_{m+1/2})\Bigr)^2
	\Big]\Big\rangle_{\fg}\,\dx \\
	\label{e:smoothed_lollipop_addition_expansion2}
	&\quad+2\int_{\T^2}\Big\langle F^{12}_{(\Bmod^{\star}_{m+1}+\Sdetmod^{\star}_{m+1})(\tmod^{\star}_{m+1})},
	\\
	&\qquad \qquad \qquad \qquad \quad \curl_{(\Bmod^{\star}_{m+1}+\Sdetmod^{\star}_{m+1})(\tmod^{\star}_{m+1})}
	e^{(\tmod^{\star}_{m+1}-\tmod^{\star}_{m+1/2})\Delta^F_{\Bmod^{\star}_{m+1}(\tmod^{\star}_{m+1/2})}}\linear[\Bmod^{\star}_{m+1}](\tmod^{\star}_{m+1/2})\Big\rangle_{\fg}\,\dx \\
	&\quad+\Bigl\|\curl_{(\Bmod^{\star}_{m+1}+\Sdetmod^{\star}_{m+1})(\tmod^{\star}_{m+1})}
	e^{(\tmod^{\star}_{m+1}-\tmod^{\star}_{m+1/2})\Delta^F_{\Bmod^{\star}_{m+1}(\tmod^{\star}_{m+1/2})}}\linear[\Bmod^{\star}_{m+1}](\tmod^{\star}_{m+1/2})\Bigr\|_{L^2}^2
	\label{e:smoothed_lollipop_addition_expansion3} \;.
\end{equs}

We use $\itobjectsmod^{(2)}_{m+1} $ to denote the sum of \eqref{e:smoothed_lollipop_addition_expansion1}--\eqref{e:smoothed_lollipop_addition_expansion3}.
	Moreover, we have
	\begin{equ}
	\label{e:smoothed_lollipop_addition_expansion_loose}
		\|F_{\widetilde{A}^{+, \star}(\tmod_{m+1}^{\star})}\|_{L^2}^2
		\leq C \big(\taumod_{m+1}^{\star} \big)^{-4} \;.
	\end{equ}
\end{lemma}

\begin{proof}
	For simplicity in this proof we use the shorthand notation 
	\begin{equ}
	\label{e:smoothed_out_lollipop}
		v_{\Bmod_{m+1}^{\star} } = e^{(\tmod^{\star}_{m+1}-\tmod^{\star}_{m+1/2})\Delta^F_{\Bmod^{\star}_{m+1}(\tmod^{\star}_{m+1/2})}} \linear[\Bmod^{\star}_{m+1}](\tmod^{\star}_{m+1/2}) \;.
	\end{equ}
	In the case of $\tmod_{m + 1/2}^{\star} \neq \tmod_{m}$ and $\tmod_{m+1} = \tmod_{m+1}^{\star}$, we have
	\begin{equ}
		\Asmoothpartmod(\tmod_{m+1}) = (\Bmod_{m+1} + \Sdetmod_{m+1})(\tmod_{m+1}) + v_{\Bmod_{m+1}^{\star} } \;.
	\end{equ}
	By the expansion of the curvature (see \eqref{prelim:eq-linearization-2}), we get
	\begin{equs}
		\,& \|F_{\widetilde{A}^{+, \star}(\tmod^{\star}_{m+1})}\|_{L^2}^2 \\
		&= \|F_{(\Bmod^{\star}_{m+1}+\Sdetmod^{\star}_{m+1})(\tmod^{\star}_{m+1})}\|_{L^2}^2 \\
		&\quad+2\int_{\T^2}\Big\langle F^{12}_{(\Bmod^{\star}_{m+1}+\Sdetmod^{\star}_{m+1})(\tmod^{\star}_{m+1})},
		\big[\bigl(v_{\Bmod_{m+1}^{\star} }\bigr)^1, \bigl(v_{\Bmod_{m+1}^{\star} }\bigr)^2\big]\Big\rangle_{\fg}\,\dx \\
		&\quad+2\int_{\T^2}\Big\langle F^{12}_{(\Bmod^{\star}_{m+1}+\Sdetmod^{\star}_{m+1})(\tmod^{\star}_{m+1})},
		\curl_{(\Bmod^{\star}_{m+1}+\Sdetmod^{\star}_{m+1})(\tmod^{\star}_{m+1})}v_{\Bmod_{m+1}^{\star} }\Big\rangle_{\fg}\,\dx \\
		&\quad+\bigl\|\curl_{(\Bmod^{\star}_{m+1}+\Sdetmod^{\star}_{m+1})(\tmod^{\star}_{m+1})}v_{\Bmod_{m+1}^{\star} }\bigr\|_{L^2}^2 \\
		&\quad+2\int_{\T^2}\Big\langle
		\curl_{(\Bmod^{\star}_{m+1}+\Sdetmod^{\star}_{m+1})(\tmod^{\star}_{m+1})}v_{\Bmod_{m+1}^{\star} }, \big[\bigl(v_{\Bmod_{m+1}^{\star} }\bigr)^1, \bigl(v_{\Bmod_{m+1}^{\star} } \bigr)^2\big]\Big\rangle_{\fg}\,\dx \\
		&\quad+\Big\|\big[\bigl(v_{\Bmod_{m+1}^{\star} }\bigr)^1, \bigl(v_{\Bmod_{m+1}^{\star} }\bigr)^2\big]\Big\|_{L^2}^2 \;.
	\end{equs}
	It suffices to bound
	\begin{equ}
		2\int_{\T^2}\Big\langle \curl_{(\Bmod^{\star}_{m+1}+\Sdetmod^{\star}_{m+1})(\tmod^{\star}_{m+1})}v_{\Bmod_{m+1}^{\star} }, \big[\bigl(v_{\Bmod_{m+1}^{\star} } \bigr)^1, \bigl(v_{\Bmod_{m+1}^{\star} } \bigr)^2\big]\Big\rangle_{\fg}\,\dx 
		+ \Big\|\big[ \bigl(v_{\Bmod_{m+1}^{\star} }\bigr)^1, \bigl(v_{\Bmod_{m+1}^{\star} }\bigr)^2\big]\Big\|_{L^2}^2 \;.
	\end{equ}
	Similar to the argument in \eqref{e:smoothed_lollipop_bound0}, we have
	\begin{equs}
	\label{e:smoothed_lollipop_bound}
		\|v_{\Bmod_{m+1}^{\star} }\|_{L^\infty} &\lesssim (\taumod^{\star}_{m+1})^{-\frac{\kappa}{2} - \lambda_1} \;, \\
		\|\covd_{\Bmod_{m+1}(\tmod^{\star}_{m+1})}v_{\Bmod_{m+1}^{\star} }\|_{L^\infty} &\lesssim (\taumod^{\star}_{m+1})^{-\frac12 - \frac{\kappa}{2} -\lambda_1} \;.
	\end{equs}
	The bounds for the above two terms then follow from \eqref{e:smoothed_lollipop_bound} and the bound for $\Sdetmod_{m+1}$ from \eqref{eq:Sdet-L-infty-estimates}. The second bound in the statement follows from \eqref{e:smoothed_lollipop_bound} and \eqref{e:iteration_deterministic_bounds_rmk}.
\end{proof}

As in the previous subsection, we wish to estimate the conditional expectation of $\itobjectsmod_{m+1}^{(1)} $ and $\itobjectsmod^{(2)}_{m+1} $. As before, we need to first restrict to the good event to apply the results in Section~\ref{sec:improved_stochastic_estimates}. 
The following is the definition parallel to Definition~\ref{defn:star_objects}. 
\begin{definition}
\label{defn:star_objects_mod}
	Recall the starred objects defined in Definition~\ref{defn:time_steps_mod}.
	For any $m \geq 1$, we first define
	\begin{equs}
		\widetilde{A}^{+,\star}(\tmod_m^\star)
		&\eqdef (\Bmod_m^\star + \Sdetmod_m^\star)(\tmod_m^\star) \\
		&\quad+ \one_{\{\tmod_{m-1/2}^\star \neq \tmod_{m-1}\}}
		e^{(\tmod_m^\star-\tmod_{m-1/2}^\star)\Delta^F_{\Bmod_m^\star(\tmod_{m-1/2}^\star)}}
		\linear[\Bmod_m^\star](\tmod_{m-1/2}^\star) \;.
	\end{equs}
	Define the sigma algebra
	\begin{equ}
		\fF_{\tmod_{m-1/2}^{\star}}^{\,\xi} \eqdef \fF_{\tmod_{m-1}}\vee\sigma\Bigl(\xi\restr\bigl((\tmod_{m-1},\tmod_{m-1/2}^\star]\times\T^2\bigr)\Bigr) \;.
	\end{equ}
	Choose the smallest integer $\Lmod_m^\star$ such that
	\begin{equ}
		\Lmod_m^\star \geq 1\vee\covariantenergy\bigl(\widetilde A^{+,\star}(\tmod_m^\star), \higgsdetmod_m^\star(\tmod_m^\star)\bigr)^{\frac1{10}} \;.
	\end{equ}
	We then define the two-star quantities, which are the next time step objects starting from $\tmod_m^\star$. First, define
	\begin{equ}
		\taumod_{m+1}^{\star \star}
		\eqdef T^{(\Lmod_m^\star)}\bigl(\widetilde{A}^{+,\star}(\tmod_m^\star),
		\higgsdetmod_m^\star(\tmod_m^\star)\bigr) \;, \quad
		\tmod_{m+1}^{\star \star} \eqdef \tmod_m^\star + \taumod_{m+1}^{\star \star} \;.
	\end{equ}
	For $t \in (\tmod_m^\star, \tmod_{m+1}^{\star \star}]$, define
	\begin{equ}
		\Bmod_{m+1}^{+, \star \star}(t)
		\eqdef \Psi_{t - \tmod_m^\star}\bigl(\widetilde{A}^{+,\star}(\tmod_m^\star)\bigr) \;,
	\end{equ}
	where $\Psi$ denotes the YM gradient flow.
	Note that now the previous objects are all $\fF_{\tmod_{m-1/2}^{\star}}^{\,\xi}$ measurable.
	
	Finally, we define the good event
	\begin{equ}
		\eventmod_m \eqdef \{ \taumod_m = \taumod_m^\star \}
		\cap \{ \Lmod_m = \Lmod_m^\star \}
		\cap \{ \taumod_{m+1} = \taumod_{m+1}^{\star \star} \} \;,
	\end{equ}
	which depends on the noise from time $\tmod_{m-1}$ to $\tmod_{m+1}$. 
\end{definition}

For the bad event $\eventmod_m^c$, we have the decomposition 
	\begin{equs}
	\label{e:rare_event_decom_mod}
		\eventmod_m^c
		&\subset \{\taumod_m\neq\taumod_m^\star\}
		\cup\bigl(\{\taumod_m=\taumod_m^\star\}
		\cap\{\Lmod_m>\Lmod_m^\star\}\bigr) \\
		&\quad\cup\bigl(\{\taumod_m=\taumod_m^\star\}
		\cap\{\Lmod_m=\Lmod_m^\star\}
		\cap\{\taumod_{m+1}\neq\taumod_{m+1}^{\star}\}\bigr) \;.
	\end{equs}
	The following lemma shows that each event in this decomposition is a rare event. 
\begin{lemma}
\label{le:rare_event_mod}
	There exist constants $c, C > 0$, independent of $\stepupperbound$, such that for every $m \geq 1$,
	\begin{equ}
	\label{e:tau_neq_T_unlikely1_mod}
		\P\Bigl(\taumod_m\neq\taumod_m^\star\Bigm|\fF_{\tmod_{m-1}}\Bigr)
		\leq C\exp\Bigl(-c\bigl(\taumod_m^\star\bigr)^{-\lambda_1}\Bigr) \;.
	\end{equ}
	Moreover,
	\begin{equs}
	\label{e:tau_neq_T_unlikely2_mod}
		\, &\P\Bigl(\{\taumod_m=\taumod_m^\star\}\cap\{\Lmod_m=\Lmod_m^\star\}
		\cap\{\taumod_{m+1}\neq\taumod_{m+1}^{\star\star}\}
		\Bigm|\fF_{\tmod_{m-1}}\Bigr) \\
		&\leq C\E\Bigl(\exp\Bigl(-c\covariantenergy\bigl(\widetilde{A}^{+,\star}(\tmod_m^\star),
		\higgsdetmod_m^\star(\tmod_m^\star)\bigr)^{\frac{\lambda_1}{10}}\Bigr)\Bigm|\fF_{\tmod_{m-1}}\Bigr) \;,
	\end{equs}
	and
	\begin{equs}
	\label{e:tau_neq_L_unlikely_mod}
	\, &\P\Bigl(\{\taumod_m=\taumod_m^\star\}\cap\{\Lmod_m>\Lmod_m^\star\}
		\Bigm|\fF_{\tmod_{m-1}}\Bigr) \\
		&\leq C\E\Bigl(\exp\Bigl(-c\covariantenergy\bigl(\widetilde{A}^{+,\star}(\tmod_m^\star),
		\higgsdetmod_m^\star(\tmod_m^\star)\bigr)^{\frac{\fa}{5}}\Bigr)\Bigm|\fF_{\tmod_{m-1}}\Bigr) \;.
	\end{equs}
\end{lemma}

\begin{proof}
	Applying Lemma~\ref{le:tau_tail_bound} to the conditional expectation after shifting the initial time to $\tmod_{m-1}$ proves \eqref{e:tau_neq_T_unlikely1_mod}. 
	
	To prove \eqref{e:tau_neq_T_unlikely2_mod}, we first condition on $\fF_{\tmod_m}$ and apply \eqref{e:tau_neq_T_unlikely1_mod}. Then, the desired bound follows from the definition of $\taumod_{m+1}^{\star \star}$, since we have 
	\begin{equ}
		\bigl( \taumod_{m+1}^{\star \star} \bigr)^{-1} \geq \Lmod_m^\star \geq \covariantenergy\big( \widetilde{A}^{+,\star}(\tmod_m^\star), \higgsdetmod_m^\star(\tmod_m^\star) \big)^{\frac{1}{10} } \;. 
	\end{equ}
	
	We now turn to \eqref{e:tau_neq_L_unlikely_mod}. 
	Recall the definitions (Definitions~\ref{defn:time_steps_mod} and \ref{defn:star_objects_mod}) of $\Lmod_m$ and $\Lmod_m^\star$. 
	In particular, since $\Lmod_m$ is chosen to be the smallest integer satisfying \eqref{e:Lmod_m_lower_bound} and such that $\bigl(\Asmoothpartmod(\tmod_m), \Aroughpartmod(\tmod_m), \phismoothpartmod(\tmod_m), \phiroughpartmod(\tmod_m) \bigr) $ is $\Lmod_m$-admissible (recall Definition~\ref{defn:L_admissible}), we get that either
	\begin{equ}
		\| \Aroughpartmod(\tmod_m) \|_{\covHolder{\Asmoothpart(\tmod_m),F}^{-\kappa}} \geq \bigl(\Lmod_m^\star \bigr)^{\fa} 
		\geq  \bigl( \covariantenergy\bigl(\Asmoothpartmod(\tmod_m), \phismoothpartmod(\tmod_m)\bigr)^{\frac{1}{10}} \vee 1 \bigr)^{\fa} \;, 
	\end{equ}
	or
	\begin{equ}
		\| \phiroughpartmod(\tmod_m) \|_{\covHolder{\Asmoothpart(\tmod_m)}^{-\kappa}} \geq \bigl(\Lmod_m^\star \bigr)^{\fa} \geq \bigl( \covariantenergy\bigl(\Asmoothpartmod(\tmod_m), \phismoothpartmod(\tmod_m)\bigr)^{\frac{1}{10}} \vee 1 \bigr)^{\fa} \;. 
	\end{equ}
	Suppose that the first case where $\Aroughpartmod(\tmod_m)$ is large happens. Recall the splitting of $\Aroughpartmod$ at the end of Definition~\ref{defn:time_steps_mod}. 
	Using Lemma~\ref{le:remainder_estimate_mod} and this splitting, we get that 
	\begin{equ}
		\| \Alinmod(\tmod_m) \|_{\covHolder{\Asmoothpart(\tmod_m),F}^{-\kappa}} \geq \bigl( \covariantenergy\bigl(\Asmoothpartmod(\tmod_m), \phismoothpartmod(\tmod_m)\bigr)^{\frac{1}{10}} \vee 1 \bigr)^{\fa} - C  \;.
	\end{equ}
	Conditioning on $\fF_{\tmod_{m-1/2}^\star}^{\,\xi}$, the result then follows from Lemma~\ref{le:lollipop_negative_holder} after shifting the initial time to $\tmod_{m-1/2}^\star$.
	Note that on the event $\{ \taumod_m = \taumod_m^\star \}$, $\tmod_m$ is deterministic after conditioning on $\fF_{\tmod_{m-1/2}^\star}^{\,\xi}$, and the definitions of $\taumod_m^\star$ and $(T^{(1)} \wedge T^{(2)})(\Asmoothpartmod(\tmod_m)) $ ensure that the assumption in Lemma~\ref{le:lollipop_negative_holder} holds.
	If we are in the case where $\phiroughpartmod(\tmod_m)$ is large, the argument is the same, except that	we should shift the initial time to $\tmod_{m-1}$ when applying Lemma~\ref{le:lollipop_negative_holder}. 
\end{proof}

\begin{lemma}
\label{le:linear_mod_conditional_moments}
	Let $\stepupperbound$ be sufficiently small. For any $p \geq 1$, there exists some constant $C_p$ such that for every $m \geq 1$,
	\begin{equ}
	\label{e:linear_mod_conditional_moments}
		\E\Bigl(\|\Alinmod(\tmod_m)\|_{\covHolder{\Asmoothpartmod(\tmod_m),F}^{-\kappa}}^p + \|\philinear[\Bmod_m](\tmod_m)\|_{\covHolder{\Asmoothpartmod(\tmod_m)}^{-\kappa}}^p \Bigm| \fF_{\tmod_{m-1}}\Bigr) \leq C_p \;.
	\end{equ}
\end{lemma}
\begin{proof}
	On the event $\{ \taumod_m=\taumod_m^\star \} $, the desired bound follows directly from Lemma~\ref{le:lollipop_negative_holder}. We note that when applying this lemma to treat $\Alinmod(\tmod_m)$, which is $\linear[\Bmod_m^\star,1/2](\tmod_m^\star) $ on this event, we will first condition on $\fF_{\tmod_{m-1/2}^{\star}}^{\,\xi}$.  

	On the event $\{ \taumod_m \neq \taumod_m^\star \} $, we have $\Alinmod(\tmod_m) = \linear[\Bmod_m^\star](\tmod_m) $. 
	Lemma~\ref{le:B_Sdet_heat_comparison} implies that 
	\begin{equ}
	\label{e:B_Sdet_heat_comparison_application}
		\|\Alinmod(\tmod_m)\|_{\covHolder{\Asmoothpartmod(\tmod_m),F}^{-\kappa}} \lesssim (\taumod_m^\star)^{-\kappa/2 - \lossdel}\sup_{\tmod_{m-1}<t\leq\tmod_m^\star}\|\linear[\Bmod_m^\star](t)\|_{\covHolder{\Bmod_m^\star(t),F,\taumod_m^\star}^{-\kappa}} \;,
	\end{equ}
	and similarly for the object $\philinear[\Bmod_m](\tmod_m)$. 
	By Lemma~\ref{le:lollipop_negative_holder_uniform}, we have
	\begin{equ}
		\E\Bigl(\sup_{\tmod_{m-1}<t\leq\tmod_m^\star}\|\linear[\Bmod_m^\star](t)\|_{\covHolder{\Bmod_m^\star(t),F,\taumod_m^\star}^{-\kappa}}^{2p} \Bigm| \fF_{\tmod_{m-1}}\Bigr) \leq C_p \;.
	\end{equ}
	Since $\taumod_m^\star$ is $\fF_{\tmod_{m-1}}$-measurable, H\"older's inequality, the previous discussions and \eqref{e:tau_neq_T_unlikely1_mod} imply
	\begin{equs}
		\E\Bigl(\one_{\{\taumod_m<\taumod_m^\star\}}\|\Alinmod(\tmod_m)\|_{\covHolder{\Asmoothpartmod(\tmod_m),F}^{-\kappa}}^p \Bigm| \fF_{\tmod_{m-1}}\Bigr) \lesssim (\taumod_m^\star)^{-p(\kappa/2 + \lossdel) }\P\bigl(\taumod_m<\taumod_m^\star\bigm|\fF_{\tmod_{m-1}}\bigr)^{1/2} 
		\lesssim 1 \;.
	\end{equs}
	This proves the desired bound for $\Alinmod$. The bound for $\philinear[\Bmod_m]$ follows similarly. 
\end{proof}

The following proposition is the main step towards proving the decay when the flow becomes close to the Yang--Mills connection. 
The estimate for the main term yielding decay has already been abstracted out in Section~\ref{section:refined_decay_mechanism}. Here the main task is to show that all the other terms are of lower order than this main term. 
\begin{proposition}
	\label{pr:itobjectsmod_expectation}
	There exist some constants $C_1, C_2$ such that $C_1$ is independent of $\stepupperbound$, and for any $m \geq 0$, if $\tmod_{m+1/2}^\star \neq \tmod_m$ and $\taumod_{m+1}^\star < \Lmod_m^{-1}$, then we have 
	\begin{equ}
		\label{e:itobjectsmod_2_expectation}
		\E \bigl( \itobjectsmod^{(2)}_{m+1} \one_{\tmod_{m+1} = \tmod_{m+1}^{\star} } | \, \fF_{\tmod_m }  \bigr) \leq - C_1 \frac{\abs{\log \taumod_{m+1}^{\star} }^{1 + c a_0} }{\taumod_{m+1}^{\star}} + C_2 \;.
	\end{equ}
	\end{proposition}
\begin{remark}
	The inequality \eqref{e:itobjectsmod_2_expectation} provides the decay mechanism when the dynamics is close to some Yang--Mills connection. 
\end{remark}

	\begin{proof}
	We begin with some preliminary discussions. Recall Lemma~\ref{le:close_YM_connection}. First note that under the assumption $\tmod_{m + 1/2}^{\star} \neq \tmod_m$ and $\taumod_{m+1}^\star < \Lmod_m^{-1}$, we have 
	\begin{equ}
	\label{e:expectation_decay_near_YM_Phicondition}
		\taumod^{\star}_{m+1} = T^{(1)}(\Asmoothpartmod(\tmod_m))  = \logtimestep \|F_{\Asmoothpartmod(\tmod_m)}\|_{L^2}^{-1} \log \bigl(\|F_{\Asmoothpartmod(\tmod_m)} \|_{L^2} \vee 2\bigr) \wedge \stepupperbound  \wedge \logtimestep \leq \| \phismoothpartmod(\tmod_m) \|_{L^r}^{-\gamma_0} \;,
	\end{equ}
	which means 
	\begin{equs}
		\| F_{\Asmoothpartmod(\tmod_{m})} \|_{L^2}^2  \leq \covariantenergy\big(\Asmoothpartmod(\tmod_{m}), \phismoothpartmod(\tmod_{m})\big) &\leq \|F_{\Asmoothpartmod(\tmod_{m})} \|_{L^2}^2 + \Bigl(\frac{ \|F_{\Asmoothpartmod(\tmod_m)}\|_{L^2}}{\logtimestep \log \bigl( \|F_{\Asmoothpartmod(\tmod_m)}\|_{L^2} \vee 2 \bigr)} \Bigr)^{\frac{\gamma_1}{\gamma_0}} + C \\
		&\lesssim \|F_{\Asmoothpartmod(\tmod_{m})} \|_{L^2}^2 + 1 \;,
	\end{equs}
	since we have chosen $\gamma_1 - 4 \ll \gamma_0 - 2$. Moreover, the assumptions imply that the decrease $\| F_{\Bmod_{m+1}^{+, \star}(\tmod_{m})} \|_{L^2}^2 - \| F_{\Bmod_{m+1}^{+, \star}(\tmod^\star_{m+1})} \|_{L^2}^2 $ is small. Together with Lemma~\ref{le:close_YM_connection}, this implies that the assumptions \eqref{e:close_to_YM_average} and \eqref{e:decay_mechanism_F_lower_bound} in Proposition~\ref{main:prop-YM-connection-decay-mechanism} hold with $B = \Bmod_{m+1}^{+, \star}$ and $[0,T]$ replaced by $[\tmod_m,\tmod_{m+1}^\star]$.
	Finally, from the definitions of $\tmod_{m+1/2}^{\star}, \tmod_{m+1}^{\star}$ (recall Definition~\ref{defn:time_steps_mod}), we have
	\begin{equ}
		\tmod_{m+1/2}^{\star}-\tmod_m  = \tmod_{m+1}^{\star}-\tmod_{m+1/2}^{\star} = \frac12 \taumod_{m+1}^\star \;.
	\end{equ}
	Now we start to prove \eqref{e:itobjectsmod_2_expectation}. We will show that if $ \tmod_{m + 1/2}^{\star} \neq \tmod_{m} $, we have 
	\begin{equs}
	\,& \E \Bigl( \int_{\T^2} \Big\langle F^{12}_{(\Bmod^{\star}_{m+1} + \Sdetmod^{\star}_{m+1})(\tmod^{\star}_{m+1})}, \Big[v^1_{\Bmod_{m+1}^{\star}}, v^2_{\Bmod_{m+1}^{\star}}\Big] \Big\rangle_{\fg} \, \dx \\
	& \qquad \cdot \one_{\tmod_{m + 1} = \tmod_{m + 1}^{\star}} \Big| \fF_{\tmod_m} \Bigr) \leq - C_1 \frac{\abs{\log \taumod_{m+1}^{\star}}^{1 + c a_0}}{\taumod_{m+1}^{\star}} + C_2 \;,
	\label{e:unstablepart_curvature1}\\
	\,& \E \Bigl( \Big| \int_{\T^2} \Big\langle F^{12}_{(\Bmod^{\star}_{m+1} + \Sdetmod^{\star}_{m+1})(\tmod^{\star}_{m+1})}, \curl_{(\Bmod^{\star}_{m+1} + \Sdetmod^{\star}_{m+1})(\tmod^{\star}_{m+1})} v_{\Bmod_{m+1}^{\star}} \Big\rangle_{\fg} \, \dx \Big| \\
	& \qquad \cdot \one_{\tmod_{m + 1} = \tmod_{m + 1}^{\star}} \Big| \fF_{\tmod_m} \Bigr) \leq C \bigl( \taumod^{\star}_{m+1} \bigr)^{-1 + \kappa} \;,
	\label{e:unstablepart_curvature2}\\
	\,& \E \Bigl( \bigl\|\curl_{(\Bmod^{\star}_{m+1} + \Sdetmod^{\star}_{m+1})(\tmod^{\star}_{m+1})} v_{\Bmod_{m+1}^{\star}}\bigr\|_{L^2}^2 \\
	& \qquad \cdot \one_{\tmod_{m + 1} = \tmod_{m + 1}^{\star}} \Big| \fF_{\tmod_m} \Bigr) \leq C \bigl( \taumod_{m+1}^{\star} \bigr)^{-1} \;.
	\label{e:unstablepart_curvature3}
\end{equs}
Here we used the same shorthand notation as in the proof of Lemma~\ref{le:smoothed_lollipop_addition_expansion}.

	To show \eqref{e:unstablepart_curvature1}, we first define $\linear[\Bmod^{+, \star}_{m+1}]$ as the solution of 
	\begin{equ}
		\partial_t \linear[\Bmod^{+, \star}_{m+1}] = \Delta_{\Bmod^{+, \star}_{m+1}} \linear[\Bmod^{+, \star}_{m+1}] + 2 \ad(F_{\Bmod^{+, \star}_{m+1}})\linear[\Bmod^{+, \star}_{m+1}] + \xi \;, \qquad \linear[\Bmod^{+, \star}_{m+1}](\tmod_m) = 0 \;,
	\end{equ}
	and we also use the shorthand
	\begin{equ}
		v_{\Bmod_{m+1}^{+, \star}} = e^{(\tmod^{\star}_{m+1} - \tmod^{\star}_{m+1/2}) \Delta^{F}_{\Bmod^{+, \star}_{m+1}(\tmod^{\star}_{m+1/2})}} \linear[\Bmod^{+, \star}_{m+1}](\tmod^{\star}_{m+1/2}) \;. 
	\end{equ}
	By Proposition~\ref{main:prop-YM-connection-decay-mechanism}, we have 
	\begin{equ}
		\E  \Bigl( \Big \langle F^{12}_{\Bmod^{+, \star}_{m+1}}(\tmod_{m+1}^{\star}), \Big [\bigl( v_{\Bmod_{m+1}^{+, \star}} \bigr)^1, \bigl(v_{\Bmod_{m+1}^{+, \star}} \bigr)^2 \Big] \Big \rangle \Big| \fF_{\tmod_{m}}  \Bigr) 
		\leq -C \frac{\abs{\log \taumod_{m+1}^{\star} }^{1 + c a_0} }{\taumod_{m+1}^{\star}} \;.
	\end{equ}
	Using Lemma~\ref{le:lollipop_difference} and Lemma~\ref{le:lollipop_negative_holder} to replace $v_{\Bmod_{m+1}^{+, \star}}$ by $v_{\Bmod_{m+1}^{\star}} $, we get 
	\begin{equ}
		\E  \Bigl( \Big \langle \Ad_{\widetilde{g}_{m+1}^{\star}(\tmod_{m+1}^{\star}) }  F^{12}_{\Bmod^{+, \star}_{m+1}}(\tmod_{m+1}^{\star}), 
		\Big [\bigl( v_{\Bmod_{m+1}^{\star}} \bigr)^1, \bigl(v_{\Bmod_{m+1}^{\star}} \bigr)^2 \Big] \Big \rangle  \Big| \fF_{\tmod_{m}}  \Bigr)  \leq - C_1  \frac{\abs{\log \taumod_{m+1}^{\star} }^{1 + c a_0} }{\taumod_{m+1}^{\star}}  + C_2  \;.
	\end{equ}
	Here, we used Lemma~\ref{le:lollipop_difference} to get 
	\begin{equs}
		\,&\E \Bigl( \Bigl\| \int_{\tmod_{m}}^{\tmod_{m+1/2}^\star} \int_{\T^2}
		\Bigl( \Ad_{(\widetilde{g}_{m+1}^\star(\tmod_{m+1}^\star, x) )^{-1} } \fq_{\Bmod_{m+1}^\star}(s,y;\tmod_{m+1}^\star,x) \Ad_{\widetilde{g}_{m+1}^\star(s,y) } \xi(s,y) \\
		&\qquad \qquad \qquad \qquad \qquad - \fq_{\Bmod_{m+1}^{+,\star}}(s,y;\tmod_{m+1}^\star,x) \xi(s,y) \Bigr)  \,\dy \,\ds
		\Bigr\|_{L_x^2}^2 \,\Big|\, \fF_{\tmod_{m}} \Bigr) \lesssim (\taumod_{m+1}^\star)^{\frac14} \;,
	\end{equs}
	and then used the fact that when conditioned on $\fF_{\tmod_m}$, $\Ad_{\widetilde{g}_{m+1}^{\star}(s, y) } \xi(s, y)$ has the same law as white noise on $[\tmod_m, \tmod_{m+1}^\star] \times \T^2$. 
		
	By \eqref{e:smoothed_lollipop_bound0}, Lemma~\ref{le:iteration_basic_bounds_mod} and Lemma~\ref{le:rare_event_mod}, we can add the indicator $\one_{\taumod_{m+1} = \taumod_{m+1}^\star }$ in the above inequality
	\begin{equs}
		 \,&\E  \Bigl(  \one_{\taumod_{m+1} = \taumod_{m+1}^\star } \Big \langle  \Ad_{\widetilde{g}_{m+1}^{\star}(\tmod_{m+1}^{\star}) } F^{12}_{\Bmod^{+, \star}_{m+1}}(\tmod_{m+1}^{\star}), \Big [\bigl( v_{\Bmod_{m+1}^{\star}} \bigr)^1, \bigl(v_{\Bmod_{m+1}^{\star}} \bigr)^2 \Big] \Big \rangle  \Big| \fF_{\tmod_{m}}  \Bigr) 
		\\
		&\leq - C_1  \frac{\abs{\log \taumod_{m+1}^{\star} }^{1 + c a_0} }{\taumod_{m+1}^{\star}}  + C_2 \;.
	\end{equs}
	It remains to replace $\Ad_{\widetilde{g}_{m+1}^{\star}(\tmod_{m+1}^{\star}) } F^{12}_{\Bmod^{+, \star}_{m+1}}$ by $F^{12}_{\Bmod^{\star}_{m+1} + \Sdetmod^\star_{m+1}}$. From the proof of Lemma~\ref{le:pertrubted_YMHF_curvature}, we get that 
	\begin{equ}
	\bigl\|F^{12}_{\Bmod^{\star}_{m+1}}(\tmod_{m+1}^{\star}) - \Ad_{\widetilde{g}_{m+1}^{\star}(\tmod_{m+1}^{\star})} F^{12}_{\Bmod^{+, \star}_{m+1}}(\tmod_{m+1}^{\star})\bigr\|_{L^\infty} \lesssim \bigl(\taumod_{m+1}^{\star} \bigr)^{-\frac12 - \frac{\kappa}{2} - 2\fa} \;.
\end{equ}
	Moreover, applying \eqref{eq:Sdet-L-infty-estimates} gives 
\begin{equ}
	\|\Sdetmod_{m+1}^{\star}(\tmod_{m+1}^{\star})\|_{L^\infty} + (\taumod_{m+1}^{\star})^{\frac12} \|(\covd_{\Bmod_{m+1}^{\star}}\Sdetmod_{m+1}^{\star})(\tmod_{m+1}^{\star})\|_{L^\infty} \lesssim \big(\taumod_{m+1}^{\star} \big)^{- \lossdel + \frac12 - \frac2r - \frac{2}{\gamma_0}} 
	\;.
\end{equ}
The error in the replacement can be bounded straightforwardly by the above and \eqref{e:smoothed_lollipop_bound0}, together with the fact that $\gamma_0 - 2$ is larger than all the other small exponents. 
	This concludes the proof of \eqref{e:unstablepart_curvature1}.

	For \eqref{e:unstablepart_curvature2}, since we are taking the absolute value, we can throw away the indicator. With the same proof as in Lemma~\ref{le:cross_term_expansion}, the left-hand side is bounded by
\begin{equ}
	C \bigl( \taumod^{\star}_{m+1} \bigr)^{\frac12 - \lossdel} \| \covd_{\Bmod^{\star}_{m+1} + \Sdetmod^{\star}_{m+1}} F_{\Bmod^{\star}_{m+1} + \Sdetmod^{\star}_{m+1}} (\tmod^{\star}_{m+1}) \|_{L^2} \;.
\end{equ}
By Lemma~\ref{le:D_BS_F_BS} and Lemma~\ref{le:iteration_basic_bounds_mod}
 this term is further bounded by 
\begin{equs}
	\,& C \bigl(\taumod^{\star}_{m+1}\bigr)^{-\frac12- 3\lossdel} 
	+ C\bigl(\taumod^{\star}_{m+1}\bigr)^{- \lossdel - 1 + \frac{\gamma_0 - 2}{50} - 2\fa - \frac2r} 
	+ C \bigl(\taumod^{\star}_{m+1}\bigr)^{-\lossdel 
	+ \frac{\gamma_0 - 2}{50} - 2\fa-\frac2r} \| \phismoothpartmod(\tmod_m) \|_{L^4}^{\gamma_1 / 2 } \;. 
\end{equs}
Note that \eqref{e:expectation_decay_near_YM_Phicondition} implies that the last term is bounded by
\begin{equ}
	\bigl( \taumod^{\star}_{m+1} \bigr)^{-\lossdel + \frac{\gamma_0 - 2}{50} - 2\fa-\frac2r - \frac{\gamma_1}{2\gamma_0}} \leq \bigl( \taumod^{\star}_{m+1} \bigr)^{-1 + \kappa} \;,
\end{equ}
where we used that $\lossdel, \fa, \frac1r \ll \gamma_1 - 4 \ll \gamma_0 - 2 $. 
Using the fact that $\gamma_0 - 2$ is much larger than all the other constants, the other two terms are also bounded by $\bigl( \taumod^{\star}_{m+1} \bigr)^{-1 + \kappa}$. This concludes the proof of \eqref{e:unstablepart_curvature2}.
	
	For \eqref{e:unstablepart_curvature3}, we first expand the curl and use the triangle inequality to get that the left-hand side is bounded by 
	\begin{equ}
	\label{e:unstablepart_curvature3_mid}
		\E \Bigl( 2 \Bigl(\|\curl_{\Bmod^{\star}_{m+1}} v_{\Bmod_{m+1}^{\star}}  \|_{L^2}^2 + \|\Sdetmod^{\star}_{m+1}\|_{L^\infty}^2\|v_{\Bmod_{m+1}^{\star}}  \|_{L^2}^2\Bigr) \Big| \fF_{\tmod_m} \Bigr) \;.
	\end{equ}
	By Lemma~\ref{le:curl_square_expectation} and Lemma~\ref{le:iteration_basic_bounds_mod}, the first term is bounded by
	\begin{equs}
		\,& C\bigl(\tmod_{m+1}^{\star}-\tmod_{m+1/2}^{\star}\bigr)^{-1} \\
		&+C\bigl(\taumod_{m+1}^{\star}\bigr)^{-\lossdel}\bigl(\tmod_{m+1/2}^{\star}-\tmod_m\bigr)^{\frac43}\bigl(\tmod_{m+1}^{\star}-\tmod_{m+1/2}^{\star}\bigr)^{-\frac13}\bigl(\taumod_{m+1}^{\star}\bigr)^{-1-2\lossdel} \\
		&+C\bigl(\taumod_{m+1}^{\star}\bigr)^{-\lossdel}\bigl(\tmod_{m+1/2}^{\star}-\tmod_m\bigr)^{\frac12}\bigl(\tmod_{m+1}^{\star}-\tmod_{m+1/2}^{\star}\bigr)^{\frac12}\bigl\|\covd_{\Bmod_{m+1}^{\star}}F_{\Bmod_{m+1}^{\star}}\bigr\|_{L_{s,y}^2([\tmod_{m+1/2}^{\star},\tmod_{m+1}^{\star}])}^2 \\
		&+C\bigl(\taumod_{m+1}^{\star}\bigr)^{-\lossdel}\bigl(\tmod_{m+1/2}^{\star}-\tmod_m\bigr)^{\frac14}\bigl(\tmod_{m+1}^{\star}-\tmod_{m+1/2}^{\star}\bigr)^{-\frac38-2\kappa}\bigl(\taumod_{m+1}^{\star}\bigr)^{\frac18}\Bigl(\Lmod_m^{2\fa}\bigl(\taumod_{m+1}^{\star}\bigr)^{-\lossdel}+\Lmod_m^{4\fa}\Bigr) \;.
	\end{equs}
	Applying Lemma~\ref{le:iteration_basic_bounds_mod} to bound the $\bigl\|\covd_{\Bmod_{m+1}^{\star}}F_{\Bmod_{m+1}^{\star}}\bigr\|_{L_{s,y}^2([\tmod_{m+1/2}^{\star},\tmod_{m+1}^{\star}])}^2$ term, we get that all terms except for the first term above are bounded by $C \bigl(\taumod_{m+1}^{\star} \bigr)^{-1 + \kappa}$. Therefore, using \eqref{e:expectation_decay_near_YM_Phicondition} we get that 
	\begin{equ}
		\E \Bigl( \|\curl_{\Bmod^{\star}_{m+1}}v_{\Bmod_{m+1}^{\star}}  \|_{L^2}^2  \Big| \fF_{\tmod_m} \Bigr) \lesssim 
		C \bigl( \tmod_{m+1}^{\star}-\tmod_{m+1/2}^{\star} \bigr)^{-1} \sim \bigl( \taumod_{m+1}^{\star} \bigr)^{-1} \;.
	\end{equ}
 	
	For the second term in \eqref{e:unstablepart_curvature3_mid}, by \eqref{eq:Sdet-L-infty-estimates} and \eqref{e:smoothed_lollipop_bound0}, we get that
	\begin{equs}
		\,&\E \Bigl( \|\Sdetmod^{\star}_{m+1}\|_{L^\infty}^2\|v_{\Bmod_{m+1}^{\star}}  \|_{L^2}^2 \Big| \fF_{\tmod_m} \Bigr) \\
		&\lesssim \bigl(
		\taumod_{m+1}^{\star} \bigr)^{- 2 \lossdel + \kappa - \frac2r - \frac{2}{\gamma_0} } \cdot \bigl(
		\taumod_{m+1}^{\star} \bigr)^{- \kappa - 2 \lambda_1} \leq \bigl(\taumod_{m+1}^{\star} \bigr)^{-1 + \kappa} \;,
	\end{equs}	
	which is acceptable. This concludes the proof of \eqref{e:unstablepart_curvature3}.
	The result is then proved by combining the discussions above. 
\end{proof}

We now prove Proposition~\ref{pr:main_iteration}, conditional on a collection of lemmas which appear after (Lemmas~\ref{le:Decmod_rare_event}, \ref{le:taumod_rare_event}, \ref{le:itobjectsmod_1_expectation_rare}, \ref{le:itobjectsmod_1_expectation}, \ref{le:itobjects1_fluctuations} and \ref{le:itobjects2_fluctuations}). These lemmas will combine to give the fluctuation bound \eqref{e:main_iteration_fluct}. We organise things this way, so that the reader can more easily grasp the structure of the proof of Proposition~\ref{pr:main_iteration}.

\begin{proof}[Proof of Proposition~\ref{pr:main_iteration}]
	First, note that from Lemma~\ref{le:smoothed_lollipop_addition_expansion}, for any $m \geq 1$ we have
\begin{equs}
	\,& \one_{\{\tmod_{m+1} = \tmod_{m+1}^{\star}\}} \covariantenergy\bigl(\widetilde{A}^{+}(\tmod_{m+1}), \phismoothpartmod(\tmod_{m+1})\bigr)  \\
	&\leq \one_{\{\tmod_{m+1} = \tmod_{m+1}^{\star}\}} \bigg[ \covariantenergy\bigl((\Bmod_{m+1} + \Sdetmod_{m+1})(\tmod_{m+1}), \higgsdetmod_{m+1}(\tmod_{m+1})\bigr) \\
	&\qquad \qquad \qquad \qquad \qquad \qquad + \one_{\{\tmod^\star_{m+1/2} \neq \tmod_m\}}
\Bigl(\itobjectsmod^{(2)}_{m+1} + C\bigl(\taumod_{m+1}\bigr)^{-\frac12-\frac{3\kappa}{2}-3\lambda_1}\Bigr) \bigg]
\end{equs}
Combining it with Proposition~\ref{pr:expansion_iteration_mod}, we have 
\begin{equs}
	\,& \covariantenergy\bigl(\widetilde{A}^{+}(\tmod_{m+1}), \phismoothpartmod(\tmod_{m+1})\bigr) \\
	&\leq  \covariantenergy\bigl((\Bmod_{m+1} + \Sdetmod_{m+1})(\tmod_{m+1}), \higgsdetmod_{m+1}(\tmod_{m+1})\bigr) \\
	\label{e:main_iteration_first_step}
	&\qquad \qquad \qquad \qquad + \one_{\{\tmod_{m+1} = \tmod_{m+1}^{\star}\}} \one_{\{\tmod^\star_{m+1/2} \neq \tmod_m\}}
\Bigl(\itobjectsmod^{(2)}_{m+1} + C\bigl(\taumod_{m+1}\bigr)^{-\frac12-\frac{3\kappa}{2}-3\lambda_1}\Bigr) \\
	&\leq  \covariantenergy\bigl(\widetilde{A}^{+}(\tmod_m) , \higgsdetmod_{m} (\tmod_m) \bigr) - \frac14 \Decmod_{m+1} + C \bigl(\taumod_{m+1} \bigr)^{-1 + \kappa} + \itobjectsmod^{(1)}_{m+1} \\
	&\qquad \qquad \qquad \qquad + \one_{\{\tmod_{m+1} = \tmod_{m+1}^{\star}\}} \one_{\{\tmod^\star_{m+1/2} \neq \tmod_m\}}
\Bigl(\itobjectsmod^{(2)}_{m+1} + C\bigl(\taumod_{m+1}\bigr)^{-\frac12-\frac{3\kappa}{2}-3\lambda_1}\Bigr)  \;.
\end{equs}
Using the fact that $\bigl(\taumod_{m+1}\bigr)^{-\frac12-\frac{3\kappa}{2}-3\lambda_1} \leq \eps \bigl(\taumod_{m+1}\bigr)^{-1 + \kappa} + C(\eps)$ for any $\eps \in (0, 1)$, 
we can rewrite the previous inequality as 
\begin{equ}
	\covariantenergy\bigl(\widetilde{A}^{+}(\tmod_{m+1}), \phismoothpartmod(\tmod_{m+1})\bigr) \leq \covariantenergy\bigl(\widetilde{A}^{+}(\tmod_m), \phismoothpartmod(\tmod_m)\bigr) + \driftmod_{m+1} + \fluctmod_{m+1} + K \;,
\end{equ}
where $K$ is some constant, and 
\begin{equs}
	\driftmod_{m+1} 
	&\eqdef - \frac14 \Dec_{\taumod_{m+1}^\star }\bigl(\Asmoothpartmod(\tmod_m), \phismoothpartmod(\tmod_m)\bigr)  + C \bigl(\taumod_{m+1}^\star \bigr)^{-1 + \kappa} - 4 C_\eqref{e:one_step_change} \one_{ \{ \taumod_{m+1}^{\star} = \Lmod_m^{-1} \} } \bigl(\taumod_{m+1}^{\star} \bigr)^{-1} \\
	&+\one_{\{\tmod^\star_{m+1/2} \neq \tmod_m\} \cap \{ \taumod_{m+1}^{\star} < \Lmod_m^{-1} \} } \E \Bigl( \one_{\{\tmod_{m+1} = \tmod_{m+1}^{\star}\}} \itobjectsmod^{(2)}_{m+1} \Big| \, \fF_{\tmod_m} \Bigr) \;, \\ [5pt]
	\fluctmod_{m+1} &\eqdef \Bigl( \frac14 \Dec_{\taumod_{m+1}^\star }\bigl(\Asmoothpartmod(\tmod_m), \phismoothpartmod(\tmod_m)\bigr)  - \frac14 \Dec_{\taumod_{m+1} }\bigl(\Asmoothpartmod(\tmod_m), \phismoothpartmod(\tmod_m)\bigr)  \Bigr) \\
	&+ \Bigl( C \bigl(\taumod_{m+1} \bigr)^{-1 + \kappa} - C \bigl(\taumod_{m+1}^\star \bigr)^{-1 + \kappa} \Bigr) + 4 C_\eqref{e:one_step_change} \one_{ \{ \taumod_{m+1}^{\star} = \Lmod_m^{-1} \} } \bigl(\taumod_{m+1}^{\star} \bigr)^{-1} \\ 
	&+ \itobjectsmod^{(1)}_{m+1} \\
	&+ \one_{\{\tmod^\star_{m+1/2} \neq \tmod_m\} } \bigg( \one_{\{\tmod_{m+1} = \tmod_{m+1}^{\star}\} } \itobjectsmod^{(2)}_{m+1} - \E \Bigl( \one_{\{\tmod_{m+1} = \tmod_{m+1}^{\star}\} \cap \{ \taumod_{m+1}^{\star} < \Lmod_m^{-1} \}  } \itobjectsmod^{(2)}_{m+1} \Big| \, \fF_{\tmod_m} \Bigr) \bigg) \;.
\end{equs}
	It remains to show that 
	\begin{equ}
	\label{e:main_iteration_drift}
		\driftmod_{m+1} \leq -3 C_{\eqref{e:one_step_change}} \bigl(\taumod_{m+1}^\star\bigr)^{-1} + K' \;,
	\end{equ}
	for some constant $K'$
	and that the bound \eqref{e:main_iteration_fluct} holds. To show \eqref{e:main_iteration_drift}, first we apply \eqref{e:itobjectsmod_2_expectation} to get 	\begin{equs}
		\driftmod_{m+1}
		&\leq - \frac14 \Dec_{\taumod_{m+1}^{\star}}\bigl(\Asmoothpartmod(\tmod_m), \phismoothpartmod(\tmod_m)\bigr) + C \bigl(\taumod_{m+1}^{\star}\bigr)^{-1 + \kappa} - 4 C_\eqref{e:one_step_change} \one_{ \{ \taumod_{m+1}^{\star} = \Lmod_m^{-1} \} } \bigl(\taumod_{m+1}^{\star} \bigr)^{-1}  \\
		&\quad - C_{\eqref{e:itobjectsmod_2_expectation},1}\one_{\{\tmod^\star_{m+1/2} \neq \tmod_m\} \cap \{ \taumod_{m+1}^{\star} < \Lmod_m^{-1} \} } \frac{\abs{\log \taumod_{m+1}^{\star}}^{1 + c a_0}}{\taumod_{m+1}^{\star}} + C_{\eqref{e:itobjectsmod_2_expectation},2} \;.
	\end{equs}
	Now we discuss the following cases.
\begin{enumerate}
	\item If $\taumod_{m+1}^{\star} = \Lmod_m^{-1}$, then the bound \eqref{e:main_iteration_drift} immediately follows. 

	\item If $\taumod_{m+1}^{\star} < \Lmod_m^{-1}$ and
	\begin{equ}
		\Dec_{\taumod_{m+1}^{\star}}\bigl(\Asmoothpartmod(\tmod_m), \phismoothpartmod(\tmod_m)\bigr) > \bigl(16 C_{\eqref{e:one_step_change}} + 1\bigr)\bigl(\taumod_{m+1}^{\star}\bigr)^{-1} \;,
	\end{equ} 
	then the bound \eqref{e:main_iteration_drift} also immediately follows. 

	\item If neither of the preceding cases occurs and $\taumod_{m+1}^{\star} < \stepupperbound \wedge \logtimestep$, then $\taumod_{m+1}^{\star} < \Lmod_m^{-1}$ and $\tmod_{m+1/2}^{\star} \neq \tmod_m$. The bound \eqref{e:main_iteration_drift} then follows from the decay in the second line.

	\item If $\taumod_{m+1}^{\star} = \stepupperbound \wedge \logtimestep$, then \eqref{e:main_iteration_drift} follows by choosing a large enough $K'$. 
\end{enumerate}
The bound \eqref{e:main_iteration_drift} follows from the above discussions.

	For \eqref{e:main_iteration_fluct}, it follows from Lemmas~\ref{le:Decmod_rare_event}, \ref{le:taumod_rare_event}, \ref{le:itobjectsmod_1_expectation_rare}, \ref{le:itobjectsmod_1_expectation}, \ref{le:itobjects1_fluctuations} and \ref{le:itobjects2_fluctuations}. We note that although Lemmas~\ref{le:itobjectsmod_1_expectation} and \ref{le:itobjects1_fluctuations} are stated with $\taumod_{m+1}^{\star \star}$, the error from replacing it by $\taumod_{m+1}^\star$ is bounded by Lemma~\ref{le:taumod_rare_event}. 
	
	For \eqref{e:main_iteration_fluct_init}, applying \eqref{e:two_timesteps_relation_mod}, it follows that
	\begin{equ}
		\covariantenergy\bigl(\Asmoothpartmod(\tmod_1), \phismoothpartmod(\tmod_1)\bigr) \lesssim \taumod_1^{-4} + \taumod_1^{-\gamma_1(\frac{1}{\gamma_0} + \frac1r + \lossdel)} + \taumod_1^{-\gamma_2(\frac{1}{\gamma_0} + \frac1r + \lossdel)} + 1 \lesssim \taumod_1^{-4} \;. 
	\end{equ}
	The desired bound then follows from Lemmas~\ref{lemma:time-step-cond-exp-bd} and \ref{lemma:taumod-star-det-bound}. 	

\end{proof}

Next, we prove all the lemmas which were needed for the fluctuation bound \eqref{e:main_iteration_fluct}.

\begin{lemma}
\label{le:Decmod_rare_event}
	For any $p > 0$, there exists some constant $C$ such that for any $m \geq 1$, we have
	\begin{equ}
	\label{e:Decmod_rare_event}
		\E \Bigl( \Big| \Dec_{\taumod_{m+1}^{\star}}\bigl(\Asmoothpartmod(\tmod_m), \phismoothpartmod(\tmod_m)\bigr) - \Dec_{\taumod_{m+1}}\bigl(\Asmoothpartmod(\tmod_m), \phismoothpartmod(\tmod_m)\bigr) \Big|^p \Bigm| \fF_{\tmod_{m-1}} \Bigr) \leq C \;.
	\end{equ}
\end{lemma}

\begin{proof}
	Note that on the event $\{\taumod_{m+1} = \taumod_{m+1}^{\star}\}$, the difference in \eqref{e:Decmod_rare_event} vanishes. Moreover, from the definition \eqref{e:dec_defn}, $\Dec_t$ is nonnegative and increasing in $t$. Since $\taumod_{m+1} \leq \taumod_{m+1}^{\star}$, it suffices to show
	\begin{equ}
		\E \Bigl( \one_{\{\taumod_{m+1} \neq \taumod_{m+1}^{\star}\}} \Dec_{\taumod_{m+1}^{\star}}\bigl(\Asmoothpartmod(\tmod_m), \phismoothpartmod(\tmod_m)\bigr)^p \Bigm| \fF_{\tmod_{m-1}} \Bigr) \leq C \;.
	\end{equ}
	By \eqref{e:Dec_lower_bound}, \eqref{e:Lmod_m_lower_bound} and the definition of $\taumod_{m+1}^{\star}$, we have
	\begin{equ}
		\Dec_{\taumod_{m+1}^{\star}}\bigl(\Asmoothpartmod(\tmod_m), \phismoothpartmod(\tmod_m)\bigr) \leq \covariantenergy\bigl(\Asmoothpartmod(\tmod_m), \phismoothpartmod(\tmod_m)\bigr) \leq \Lmod_m^{10} \leq \bigl(\taumod_{m+1}^{\star}\bigr)^{-10} \;.
	\end{equ}
	The result then follows from \eqref{e:tau_neq_T_unlikely1_mod}. 
\end{proof}

\begin{lemma}
\label{le:taumod_rare_event}
	For any $\alpha > 0 $ and $p > 0$, there exists some constant $C$ such that for any $m \geq 1$, we have   
	\begin{equ}
	\label{e:taumod_rare_event1}
		\E \Bigl( \Big | \bigl( \taumod_{m+1} \bigr)^{- \alpha} - \bigl( \taumod_{m+1}^{\star} \bigr)^{ - \alpha } \Big |^p  \Big | \, \fF_{\tmod_{m-1} } \Bigr)  + \E \Bigl( \Big | \bigl( \taumod_{m+1}^\star \bigr)^{- \alpha} - \bigl( \taumod_{m+1}^{\star \star} \bigr)^{ - \alpha } \Big |^p  \Big | \, \fF_{\tmod_{m-1} } \Bigr) \leq C \;.  
	\end{equ}
	Moreover, we have 
	\begin{equ}
	\label{e:taumod_rare_event2}
		\E \Bigl( \one_{ \{ \taumod_{m+1}^{\star} = \Lmod_m^{-1} \} } \bigl(\taumod_{m+1}^{\star} \bigr)^{-p}  \Big | \, \fF_{\tmod_{m-1} } \Bigr) \leq C \;. 
	\end{equ}
	and
	\begin{equ}
	\label{e:Lmod_expectation_bound}
	\E \Bigl( \Lmod_m^p \Bigm| \fF_{\tmod_{m-1}} \Bigr) \leq C \E \Bigl( \covariantenergy\bigl(\Asmoothpartmod(\tmod_m), \phismoothpartmod(\tmod_m)\bigr)^{\frac{p}{10}} \Bigm| \fF_{\tmod_{m-1}} \Bigr) + C \;.
	\end{equ}
\end{lemma}

\begin{proof}
	The first term in \eqref{e:taumod_rare_event1} follows from Lemma~\ref{lemma:time-step-cond-exp-bd} and \eqref{e:tau_neq_T_unlikely1_mod} immediately. For the second term, 
	note that on the event $\{\taumod_m = \taumod_m^\star \} \cap \{\Lmod_m = \Lmod_m^\star \}$, we have $\taumod_{m+1}^\star = \taumod_{m+1}^{\star \star}$. Therefore, together with the triangle inequality, it suffices to show 
	\begin{equs}
	\label{e:taumod_rare_event_mid1}
		\E \Bigl( \one_{\{\taumod_m \neq \taumod_m^\star\}} \Bigl( \bigl(\taumod_{m+1}^\star\bigr)^{-\alpha p} + \bigl(\taumod_{m+1}^{\star \star}\bigr)^{-\alpha p} \Bigr) \Bigm| \fF_{\tmod_{m-1}} \Bigr) &\leq C \;, \\
	\label{e:taumod_rare_event_mid2}
		\E \Bigl( \one_{\{\taumod_m = \taumod_m^\star\} \cap \{\Lmod_m > \Lmod_m^\star\}} \Bigl( \bigl(\taumod_{m+1}^\star\bigr)^{-\alpha p} + \bigl(\taumod_{m+1}^{\star \star}\bigr)^{-\alpha p} \Bigr) \Bigm| \fF_{\tmod_{m-1}} \Bigr) &\leq C \;.
	\end{equs}
	First we note that by Lemma~\ref{lemma:taumod-star-det-bound}, \eqref{e:two_timesteps_relation_mod} and definition of $\Lmod_m$, we have 
	\begin{equs}
	\bigl(\taumod_{m+1}^{\star}\bigr)^{-1} &\lesssim \covariantenergy\bigl(\Asmoothpartmod(\tmod_m), \phismoothpartmod(\tmod_m)\bigr) + \Lmod_m \\
	&\lesssim \taumod_m^{-10} \wedge \covariantenergy\bigl(\Asmoothpartmod(\tmod_m), \phismoothpartmod(\tmod_m)\bigr)+ \|\Alinmod(\tmod_m)\|_{\covHolder{\Asmoothpartmod(\tmod_m),F}^{-\kappa}}^{\frac1{\fa}} + \|\philinear[\Bmod_m](\tmod_m)\|_{\covHolder{\Asmoothpartmod(\tmod_m)}^{-\kappa}}^{\frac1{\fa}} \;, \\
	\bigl(\taumod_{m+1}^{\star \star}\bigr)^{-1} &\lesssim \covariantenergy\bigl(\widetilde{A}^{+,\star}(\tmod_m^\star), \higgsdetmod_m^\star(\tmod_m^\star)\bigr) + \Lmod_m^\star \lesssim \covariantenergy\bigl(\widetilde{A}^{+,\star}(\tmod_m^\star), \higgsdetmod_m^\star(\tmod_m^\star)\bigr) \;.
	\end{equs}
	With Lemma~\ref{le:rare_event_mod}, it follows that 
	\begin{equs}
	\E \Bigl( \bigl( \one_{\{\taumod_m \neq \taumod_m^\star\}} + \one_{\{\taumod_m = \taumod_m^\star\} \cap \{\Lmod_m > \Lmod_m^\star\}}  \bigr) \cdot \bigl(\taumod_m^{-10} \wedge \covariantenergy\bigl(\Asmoothpartmod(\tmod_m), \phismoothpartmod(\tmod_m)\bigr)  \bigr)^p  \Bigm| \fF_{\tmod_{m-1}} \Bigr) &\leq C \;.
	\end{equs}
	Together with Lemma~\ref{le:linear_mod_conditional_moments}, this gives the desired bound for the $\bigl( \taumod_{m+1}^\star \bigr)^{-1}$ term in \eqref{e:taumod_rare_event_mid1} and \eqref{e:taumod_rare_event_mid2}. 
	
	For the $\bigl( \taumod_{m+1}^{\star \star} \bigr)^{-1}$ term, note that in Lemma~\ref{le:rare_event_mod} we actually proved the stronger bound 
	\begin{equ}
		\P\Bigl(\{\taumod_m=\taumod_m^\star\}\cap\{\Lmod_m>\Lmod_m^\star\}
		\Bigm| \fF_{\tmod_{m-1/2}^\star}^{\,\xi} \Bigr) 
		\leq C \exp\Bigl(-c\covariantenergy\bigl(\widetilde{A}^{+,\star}(\tmod_m^\star),
		\higgsdetmod_m^\star(\tmod_m^\star)\bigr)^{\frac{\fa}{5}}\Bigr) \;.
	\end{equ}
	This gives the desired bound for the $\bigl( \taumod_{m+1}^{\star \star} \bigr)^{-1}$ term in \eqref{e:taumod_rare_event_mid2}. For \eqref{e:taumod_rare_event_mid1}, we first apply Proposition~\ref{pr:expansion_iteration_mod} and Lemma~\ref{le:smoothed_lollipop_addition_expansion} to get that 
	\begin{equ}
		\E\Bigl(\covariantenergy\bigl(\widetilde{A}^{+,\star}(\tmod_m^\star), \higgsdetmod_m^\star(\tmod_m^\star)\bigr)^p \Bigm| \fF_{\tmod_{m-1}}\Bigr) \lesssim \bigl(\taumod_m^\star\bigr)^{-10p} \;,
	\end{equ}
	and the desired bound follows from this and Lemma~\ref{le:rare_event_mod}. 
	
	Now we turn to \eqref{e:taumod_rare_event2}. By \eqref{e:taumod_rare_event_mid1} and \eqref{e:taumod_rare_event_mid2}, it suffices to estimate the contribution on the event $\{\taumod_m = \taumod_m^\star\} \cap \{\Lmod_m = \Lmod_m^\star\}$, on which we have 
	\begin{equ}
		\Lmod_m = \Lmod_m^\star \leq 1 + \covariantenergy\bigl(\Asmoothpartmod(\tmod_m), \phismoothpartmod(\tmod_m)\bigr)^{\frac1{10}} \;.
	\end{equ}
	On the event $\{\taumod_{m+1}^\star = \Lmod_m^{-1}\}$, we have
\begin{equ}
	\covariantenergy\bigl(\Asmoothpartmod(\tmod_m), \phismoothpartmod(\tmod_m)\bigr) \lesssim \Lmod_m^4 \;.
\end{equ}
	Therefore, on the intersection of these events we have 
	\begin{equ}
		\Lmod_m\lesssim 1 + \Lmod_m^{\frac25} \;,
	\end{equ}
	which implies that $\Lmod_m \lesssim 1$. This proves \eqref{e:taumod_rare_event2}. 
	
	Finally, to prove \eqref{e:Lmod_expectation_bound}, note that $\Lmod_m \leq (\taumod_{m+1}^\star)^{-1}$. Thus, \eqref{e:taumod_rare_event_mid1} and \eqref{e:taumod_rare_event_mid2}, with $\alpha = 1$, give
\begin{equ}
	\E \Bigl( \bigl( \one_{\{\taumod_m \neq \taumod_m^\star\}} + \one_{\{\taumod_m = \taumod_m^\star\} \cap \{\Lmod_m > \Lmod_m^\star\}} \bigr) \Lmod_m^p \Bigm| \fF_{\tmod_{m-1}} \Bigr) \leq C \;.
\end{equ}
On the remaining event $\{\taumod_m = \taumod_m^\star\} \cap \{\Lmod_m = \Lmod_m^\star\}$, we have
\begin{equ}
	\Lmod_m^p \leq C \Bigl(1 + \covariantenergy\bigl(\Asmoothpartmod(\tmod_m), \phismoothpartmod(\tmod_m)\bigr)^{\frac{p}{10}}\Bigr) \;.
\end{equ}
	The bound \eqref{e:Lmod_expectation_bound} follows. 
\end{proof}

\begin{lemma}
\label{le:itobjectsmod_1_expectation_rare}
	Recall $\itobjectsmod_{m+1}^{(1)}$ defined in Proposition~\ref{pr:expansion_iteration_mod}. For any $p > 0$, there exists some constant $C$ such that for any $m \geq 1$, we have
	\begin{equ}
	\label{e:itobjectsmod_1_expectation_rare}
		\E \Bigl( \big| \itobjectsmod_{m+1}^{(1)} \big|^p \one_{\eventmod_m^c} \Bigm| \fF_{\tmod_{m-1}} \Bigr) \leq C \;.
	\end{equ}
\end{lemma}

\begin{proof}
	With Lemma~\ref{le:iteration_basic_bounds_mod}, one can get the crude bound 
	\begin{equ}
		\big| \itobjectsmod_{m+1}^{(1)} \big| \lesssim \taumod_{m+1}^{-2} \;.
	\end{equ}
	Therefore, it suffices to bound the conditional expectation of $\taumod_{m+1}^{-2p} \one_{\eventmod_m^c}$. Using the decomposition \eqref{e:rare_event_decom_mod}, it suffices to show
	\begin{equs}
	\,& \E \Bigl( \taumod_{m+1}^{-2p} \one_{\{\taumod_m \neq \taumod_m^\star\}} \Bigm| \fF_{\tmod_{m-1}} \Bigr) + \E \Bigl( \taumod_{m+1}^{-2p} \one_{\{\taumod_m = \taumod_m^\star\} \cap \{\Lmod_m > \Lmod_m^\star\}} \Bigm| \fF_{\tmod_{m-1}} \Bigr) \\
	&+ \E \Bigl( \taumod_{m+1}^{-2p} \one_{\{\taumod_m = \taumod_m^\star\} \cap \{\Lmod_m = \Lmod_m^\star\} \cap \{\taumod_{m+1} \neq \taumod_{m+1}^{ \star}\}} \Bigm| \fF_{\tmod_{m-1}} \Bigr) \leq C
	\end{equs}
	For the first two terms, after conditioning on $\fF_{\tmod_m}$ and applying Lemma~\ref{lemma:time-step-cond-exp-bd}, they are reduced to the estimates \eqref{e:taumod_rare_event_mid1} and \eqref{e:taumod_rare_event_mid2}. For the last term, it directly follows from Lemma~\ref{le:tau_tail_bound} after conditioning on $\fF_{\tmod_m}$. 
\end{proof}

\begin{lemma}
\label{le:itobjectsmod_1_expectation}
	Recall $\itobjectsmod_{m+1}^{(1)}$ defined in Proposition~\ref{pr:expansion_iteration_mod}.
	Let $C_{\eqref{e:one_step_change}}$ be the constant in Proposition~\ref{pr:expectation_one_step}. 
	For any $m \geq 1$ we have 
		\begin{equ}
		\label{e:itobjectsmod_1_expectation}
			\E \bigl( \itobjectsmod_{m+1}^{(1)} \one_{\eventmod_m} | \, \fF_{\tmod_{m-1/2}^{\star}}^{\,\xi}  \bigr) \leq C_{\eqref{e:one_step_change}} \bigl(\taumod_{m+1}^{\star \star} \bigr)^{-1} \;.
		\end{equ}
\end{lemma}
\begin{proof}
	The proof is essentially the same as that of Proposition~\ref{pr:expectation_one_step}. We first use the fact that we are on the event $\eventmod_m$ to change \eqref{e:iteration_deterministic_bounds_iteration1_mod}--\eqref{e:iteration_deterministic_bounds_iteration4_mod} to the corresponding version with stars, and then apply results in Section~\ref{sec:improved_stochastic_estimates}. 
	The subtlety here is that we condition on $\fF_{\tmod_{m-1/2}^{\star}}^{\,\xi}$, since the time $\tmod_{m+1}^{\star \star}$ and the background YM gradient flow $\Bmod_{m+1}^{+, \star \star}$ are only measurable with respect to it rather than $\fF_{\tmod_{m-1}}$. 
	After this conditioning, the estimates in Section~\ref{sec:improved_stochastic_estimates} apply to \eqref{e:iteration_deterministic_bounds_iteration1_mod}--\eqref{e:iteration_deterministic_bounds_iteration4_mod} and yield the same upper bound $C_{\eqref{e:one_step_change}} \bigl(\taumod^{\star \star}_{m+1} \bigr)^{-1}$. We emphasise that here we get the same constant $C_{\eqref{e:one_step_change}}$, as $\Alinmod$ is still a linear stochastic object, except that the forcing space-time white noise may come from a subinterval of the previous linear object. This modification does not affect the constant $C_{\eqref{e:one_step_change}}$ in the estimate. 
\end{proof}

\begin{lemma}
\label{le:itobjects1_fluctuations}
	Recall $\itobjectsmod_{m+1}^{(1)}$ defined in Proposition~\ref{pr:expansion_iteration_mod}.
	For any $p \in \N^*$, there exists some constant $C_p$, such that for any $m \geq 1$ we have
	\begin{equ}
	\label{e:itobjects1_fluctuations1} 
	\E \biggl( \Big| \one_{\eventmod_m} \itobjectsmod^{(1)}_{m+1} - \E \bigl( \one_{\eventmod_m} \itobjectsmod^{(1)}_{m+1} \big| \, \fF_{\tmod_{m-1/2}^{\star}}^{\,\xi} \bigr)  \Big|^p  \Big | \, \fF_{\tmod_{m-1} } \biggr) \leq C_p \E \Bigl( \bigl(\taumod_{m+1}^{\star \star} \bigr)^{p(-1 + 2 \kappa ) } \Big | \, \fF_{\tmod_{m-1} } \Bigr)  \;. 
	\end{equ}
	Moreover, we have 
	\begin{equ}
	\label{e:itobjects1_fluctuations2} 
	\E \biggl( \Bigl( \one_{\eventmod_m} \itobjectsmod^{(1)}_{m+1} - \E \bigl( \one_{\eventmod_m} \itobjectsmod^{(1)}_{m+1} \big| \, \fF_{\tmod_{m-1/2}^{\star}}^{\,\xi} \bigr) - \bigl(\taumod_{m+1}^{\star \star}\bigr)^{-1 + \kappa} \Bigr)_+^p \Big| \, \fF_{\tmod_{m-1}} \biggr) \leq C_p \;.
	\end{equ}
\end{lemma}
\begin{proof}
	From the definition of $\Alinmod$ in \eqref{e:Alinmod_defn}, when $\tmod_{m}^\star = \tmod_{m}$ we have $\Alinmod(\tmod_m) = \linear[\Bmod_m^\star, 1/2](\tmod_m^\star)$. 
	Actually, if $\tmod_{m - 1/2}^\star \neq \tmod_{m-1}$, this is true directly from definition. If $\tmod_{m - 1/2}^\star = \tmod_{m-1}$, then we have $\linear[\Bmod_m^\star, 1/2] = \linear[\Bmod_m^\star]$ from their definitions in Definition~\ref{defn:time_steps_mod}, since now they both start with zero initial data from time $\tmod_{m-1}$. 
	Therefore, on the event $\eventmod_m$, $\itobjectsmod_{m+1}^{(1)}$ is equal to its starred version, defined by
	\minilab{e:itobjectsmod1}
\begin{equs}
	\itobjectsmod_{m+1}^{(1),\star}
	\label{e:itobjectsmod1_1}
	&\eqdef \Big\|\curl_{\Bmod_{m+1}^{+,\star\star}(\tmod_{m+1}^{\star\star})} e^{\taumod_{m+1}^{\star\star}\Delta^F_{\widetilde{A}^{+,\star}(\tmod_m^\star)}} \linear[\Bmod_m^\star,1/2](\tmod_m^\star)\Big\|_{L^2}^2 \\
	\label{e:itobjectsmod1_2}
	&\quad+ 2\int_{\T^2}\Big\langle F_{\Bmod_{m+1}^{+,\star\star}}^{12}(\tmod_{m+1}^{\star\star}), \curl_{\Bmod_{m+1}^{+,\star\star}(\tmod_{m+1}^{\star\star})} e^{\taumod_{m+1}^{\star\star}\Delta^F_{\widetilde{A}^{+,\star}(\tmod_m^\star)}} \linear[\Bmod_m^\star,1/2](\tmod_m^\star)\Big\rangle_{\fg}\,\dx \\
	\label{e:itobjectsmod1_3}
	&\quad+ 2\int_{\T^2}\Big\langle F_{\Bmod_{m+1}^{+,\star\star}}^{12}(\tmod_{m+1}^{\star\star}), \\
	&\qquad\qquad \Big[\bigl(e^{\taumod_{m+1}^{\star\star}\Delta^F_{\widetilde{A}^{+,\star}(\tmod_m^\star)}}\linear[\Bmod_m^\star,1/2](\tmod_m^\star)\bigr)^1, \bigl(e^{\taumod_{m+1}^{\star\star}\Delta^F_{\widetilde{A}^{+,\star}(\tmod_m^\star)}}\linear[\Bmod_m^\star,1/2](\tmod_m^\star)\bigr)^2\Big]\Big\rangle_{\fg}\,\dx \\
	\label{e:itobjectsmod1_4}
	&\quad+ C\|\higgsdetmod_m^\star(\tmod_m^\star)\|_{L^4}^{\gamma_1-4}
	\Big|\int_{\T^2}\big|e^{\taumod_{m+1}^{\star\star}\Delta_{\widetilde{A}^{+,\star}(\tmod_m^\star)}}\higgsdetmod_m^\star(\tmod_m^\star)\big|^2 \\
	&\qquad\qquad \cdot\Big\langle e^{\taumod_{m+1}^{\star\star}\Delta_{\widetilde{A}^{+,\star}(\tmod_m^\star)}}\higgsdetmod_m^\star(\tmod_m^\star), e^{\taumod_{m+1}^{\star\star}\Delta_{\widetilde{A}^{+,\star}(\tmod_m^\star)}}\philinear[\Bmod_m^\star](\tmod_m^\star)\Big\rangle_{\V}\,\dx\Big| \;.
\end{equs}
	We claim that for every $p \geq 1$,
\begin{equ}
\label{e:itobjectsmod_1_rare}
	\E\Bigl(\one_{\eventmod_m^c}\big|\itobjectsmod_{m+1}^{(1),\star}\big|^p \Bigm| \fF_{\tmod_{m-1}}\Bigr) 
	\leq C_p \;. 
\end{equ}
	Actually, from the decomposition \eqref{e:rare_event_decom_mod}, it suffices to show that 
	\begin{equs}
		\,&\E\Bigl(\one_{\{\taumod_m\neq\taumod_m^\star\}}\big|\itobjectsmod_{m+1}^{(1),\star}\big|^p \Bigm| \fF_{\tmod_{m-1}}\Bigr)
+ \E\Bigl(\one_{\{\taumod_m=\taumod_m^\star\}\cap\{\Lmod_m>\Lmod_m^\star\}}\big|\itobjectsmod_{m+1}^{(1),\star}\big|^p \Bigm| \fF_{\tmod_{m-1}}\Bigr) \\
\label{e:itobjectsmod_1_rare_mid}
&+ \E\Bigl(\one_{\{\taumod_m=\taumod_m^\star\}\cap\{\Lmod_m=\Lmod_m^\star\}\cap\{\taumod_{m+1}\neq\taumod_{m+1}^{\star}\}}\big|\itobjectsmod_{m+1}^{(1),\star}\big|^p \Bigm| \fF_{\tmod_{m-1}}\Bigr) \leq C_p \;.
	\end{equs}
	Applying Lemma~\ref{le:iteration_basic_bounds_mod}, Lemma~\ref{lemma:taumod-star-det-bound} and the definition of $\Lmod_m^\star$ would give the following crude estimate
\begin{equ}
	\E\Bigl(\big|\itobjectsmod_{m+1}^{(1),\star}\big|^{2p} \Bigm| \fF_{\tmod_{m-1/2}^\star}^{\,\xi}\Bigr)
	\leq 
	C_p\bigl(\taumod_{m+1}^{\star\star}\bigr)^{-4p} \leq C_p' \covariantenergy\bigl(\widetilde{A}^{+,\star}(\tmod_m^\star),\higgsdetmod_m^\star(\tmod_m^\star)\bigr)^{4p}  \;.
\end{equ}
	For the first term in \eqref{e:itobjectsmod_1_rare_mid}, applying Cauchy--Schwarz and \eqref{e:tau_neq_T_unlikely1_mod} gives
\begin{equs}
	\E\Bigl(\one_{\{\taumod_m\neq\taumod_m^\star\}}\big|\itobjectsmod_{m+1}^{(1),\star}\big|^p \Bigm| \fF_{\tmod_{m-1}}\Bigr)
	&\leq \E\Bigl(\big|\itobjectsmod_{m+1}^{(1),\star}\big|^{2p} \Big| \fF_{\tmod_{m-1}}\Bigr)^{1/2}
	\P\Bigl(\taumod_m\neq\taumod_m^\star \Bigm| \fF_{\tmod_{m-1}}\Bigr)^{1/2} \\
	&\lesssim \E \Bigl( \covariantenergy\bigl(\widetilde{A}^{+,\star}(\tmod_m^\star),\higgsdetmod_m^\star(\tmod_m^\star)\bigr)^{4p} \Big| \fF_{\tmod_{m-1}} \Bigr)^{1/2}\exp\Bigl(-c\bigl(\taumod_m^\star\bigr)^{-\lambda_1}\Bigr)
	\lesssim 1 \;.
\end{equs}
	where in the last step we used \eqref{e:smoothed_lollipop_addition_expansion_loose} and \eqref{e:two_timesteps_relation_mod} to bound the energy by the inverse powers of $\taumod_m^\star$. For the second term in \eqref{e:itobjectsmod_1_rare_mid}, note that in the proof of \eqref{e:tau_neq_L_unlikely_mod}, we actually have the stronger estimate 
\begin{equ}
	\P\Bigl(\{\taumod_m=\taumod_m^\star\}\cap\{\Lmod_m>\Lmod_m^\star\} \Bigm| \fF_{\tmod_{m-1/2}^\star}^{\,\xi}\Bigr) 
	\leq C\exp\Bigl(-c\covariantenergy\bigl(\widetilde{A}^{+,\star}(\tmod_m^\star),\higgsdetmod_m^\star(\tmod_m^\star)\bigr)^{\fa/5}\Bigr) \;.
\end{equ}
	The result then directly follows from this and the previous crude bound on $\itobjectsmod_{m+1}^{(1),\star}$. The bound for the third term in \eqref{e:itobjectsmod_1_rare_mid} is similar. We omit the details. 
	
Now, Jensen's inequality for the conditional expectation gives
\begin{equs}
	\,&\E\Biggl(\Big| \one_{\eventmod_m}\itobjectsmod_{m+1}^{(1)} 
	-\E\bigl(\one_{\eventmod_m}\itobjectsmod_{m+1}^{(1)}\bigm| \fF_{\tmod_{m-1/2}^{\star}}^{\,\xi} \bigr)
	- \itobjectsmod_{m+1}^{(1),\star} 
	+\E\bigl(\itobjectsmod_{m+1}^{(1),\star}\bigm| \fF_{\tmod_{m-1/2}^{\star}}^{\,\xi} \bigr)
	\Big|^p \Bigm| \fF_{\tmod_{m-1}}\Biggr) \\
	&\lesssim \E\Bigl(\one_{\eventmod_m^c}\big|\itobjectsmod_{m+1}^{(1),\star}\big|^p \Bigm| \fF_{\tmod_{m-1}}\Bigr)
	\lesssim 1 \;.
\end{equs}
Therefore, it suffices to show
\begin{equ}
	\label{e:itobjects1_fluctuations1_mid} 
	\E \Biggl( \Big| \itobjectsmod_{m+1}^{(1),\star} - \E\bigl(\itobjectsmod_{m+1}^{(1),\star}\bigm|\fF_{\tmod_{m-1/2}^{\star}}^{\,\xi}\bigr) \Big|^p \Bigm| \fF_{\tmod_{m-1/2}^{\star}}^{\,\xi} \Biggr) \leq C_p \bigl(\taumod_{m+1}^{\star \star} \bigr)^{p(-1 + 2 \kappa) } \;, 
\end{equ}
and 
\begin{equ}
	\label{e:itobjects1_fluctuations2_mid} 
	\E \Biggl( \Bigl( \itobjectsmod_{m+1}^{(1),\star} - \E\bigl(\itobjectsmod_{m+1}^{(1),\star}\bigm|\fF_{\tmod_{m-1/2}^{\star}}^{\,\xi}\bigr) - \bigl(\taumod_{m+1}^{\star \star}\bigr)^{-1 + \kappa} \Bigr)_+^p \Bigm| \fF_{\tmod_{m-1/2}^{\star}}^{\,\xi} \Biggr) \leq C_p \;.
\end{equ}

	We estimate the fluctuation of each term in $\itobjectsmod_{m+1}^{(1), \star}$ using the results in Section~\ref{sec:improved_stochastic_estimates_fluctuations}. 
	
	\proofstep{Fluctuation of \eqref{e:itobjectsmod1_1}}:  
	Applying Lemma~\ref{le:curl_square_fluctuations}, together with Assumption~\ref{local:assumption-B}, it follows that 
	\begin{equ}
	\E\Bigl(\Big|\eqref{e:itobjectsmod1_1} - \E\bigl(\eqref{e:itobjectsmod1_1}\bigm|\fF_{\tmod_{m-1/2}^{\star}}^{\,\xi}\bigr)\Big|^p \Bigm|\fF_{\tmod_{m-1/2}^{\star}}^{\,\xi}\Bigr)^{1/p}
	\lesssim \bigl(\taumod_m^\star\bigr)^{\frac12-2\lossdel}\bigl(\taumod_{m+1}^{\star\star}\bigr)^{-1-2\lossdel} \;.
\end{equ}
	The desired bound now follows from \eqref{e:two_timesteps_relation_mod}.
	 
	\proofstep{Fluctuation of \eqref{e:itobjectsmod1_2}}:    
	This term has zero conditional expectation given $\fF_{\tmod_{m-1/2}^{\star}}^{\,\xi}$. By Lemmas~\ref{le:cross_term_expansion} and~\ref{le:DF_bound}, together with \eqref{e:iteration_F_endtime_lowerbound_mod}, we have
\begin{equs}
	\E\Bigl(\big|\eqref{e:itobjectsmod1_2}\big|^p\Bigm|\fF_{\tmod_{m-1/2}^{\star}}^{\,\xi}\Bigr)^{1/p}
	&\lesssim \bigl(\taumod_m^\star\bigr)^{\frac12-\lossdel}\bigl(\taumod_{m+1}^{\star\star}\bigr)^{-\lossdel}
	\bigl\|(\covd_{\Bmod_{m+1}^{+,\star\star}}F_{\Bmod_{m+1}^{+,\star\star}})(\tmod_{m+1}^{\star\star})\bigr\|_{L^2} \\
	&\lesssim \bigl(\taumod_m^\star\bigr)^{\frac12-\lossdel}\bigl(\taumod_{m+1}^{\star\star}\bigr)^{-1-3\lossdel} \;.
\end{equs}
	The desired bound now follows from \eqref{e:two_timesteps_relation_mod}.

	\proofstep{Fluctuation of \eqref{e:itobjectsmod1_3}}:  Applying Lemma~\ref{le:decay_mechanism_fluctuations}, together with Assumption~\ref{local:assumption-B} and \eqref{eq:curvature-pointwise-in-time-bound-Fnorm}, gives
\begin{equs}
	\E\Bigl(\Big|\eqref{e:itobjectsmod1_3} - \E\bigl(\eqref{e:itobjectsmod1_3}\bigm|\fF_{\tmod_{m-1/2}^{\star}}^{\,\xi}\bigr)\Big|^p\Bigm|\fF_{\tmod_{m-1/2}^{\star}}^{\,\xi}\Bigr)^{1/p} &\lesssim \bigl(\taumod_m^\star\bigr)^{\frac12-2\lossdel}\bigl(\taumod_{m+1}^{\star\star}\bigr)^{-2\lossdel}\bigl\|F_{\Bmod_{m+1}^{+,\star\star}}(\tmod_{m+1}^{\star\star})\bigr\|_{L^\infty} \\
	&\lesssim \bigl(\taumod_m^\star\bigr)^{\frac12-2\lossdel}\bigl(\taumod_{m+1}^{\star\star}\bigr)^{-1-3\lossdel} \;.
\end{equs}
	The desired bound now follows from \eqref{e:two_timesteps_relation_mod}.

	\proofstep{Fluctuation of \eqref{e:itobjectsmod1_4}}: 
	Applying Lemma~\ref{le:phi_cross_term_estimate}, the diamagnetic inequality \eqref{eq:zeroth-order-pB-ptwise-bd}, and the fact that $\|f\|_{L^6}^3 \leq \|f\|_{L^4}^{\frac{2r-12}{r-4}}\|f\|_{L^r}^{\frac{r}{r-4}}$ gives
	\begin{equ}
		\E\Bigl(\Big|\eqref{e:itobjectsmod1_4} - \E\bigl(\eqref{e:itobjectsmod1_4}\bigm|\fF_{\tmod_{m-1/2}^{\star}}^{\,\xi}\bigr)\Big|^p\Bigm|\fF_{\tmod_{m-1/2}^{\star}}^{\,\xi}\Bigr)^{1/p} \lesssim \bigl(\taumod_m^\star\bigr)^{\frac12}\|\higgsdetmod_m^\star(\tmod_m^\star)\|_{L^4}^{\gamma_1-4+\frac{2r-12}{r-4}}\|\higgsdetmod_m^\star(\tmod_m^\star)\|_{L^r}^{\frac{r}{r-4}} \;.
	\end{equ}
	Then, using \eqref{e:two_timesteps_relation_mod} and the definition of $\taumod_{m+1}^{\star \star}$ (in particular, Definition~\ref{defn:tau}), we have
	\begin{equs}
	\,& \bigl(\taumod_m^\star\bigr)^{\frac12}\|\higgsdetmod_m^\star(\tmod_m^\star)\|_{L^4}^{\gamma_1-4+\frac{2r-12}{r-4}}\|\higgsdetmod_m^\star(\tmod_m^\star)\|_{L^r}^{\frac{r}{r-4}} \\
	&\lesssim \|\higgsdetmod_m^\star(\tmod_m^\star)\|_{L^4}^{\gamma_1-4+\frac{2r-12}{r-4}}
	\lesssim \bigl(\taumod_{m+1}^{\star\star}\bigr)^{-\frac1{\gamma_0}(\gamma_1-4+\frac{2r-12}{r-4})}
	\lesssim \bigl(\taumod_{m+1}^{\star\star}\bigr)^{-1+2\kappa} \;,
\end{equs}
	where we also used the fact that $\frac1r , \kappa \ll \gamma_1 - 4 \ll \gamma_0 - 2$. 

	This concludes the proof of \eqref{e:itobjects1_fluctuations1_mid}. For \eqref{e:itobjects1_fluctuations2_mid}, note that for any $a > 0$ and $q > p > 0$, and $x \in \R$, we have the elementary inequality
	\begin{equ}
	(x-a)_+^p \leq |x|^p \one_{\{x>a\}} \leq a^{p-q}|x|^q \;.
	\end{equ}
	Therefore, we get 
\begin{equs}
	\,&\E\Biggl(\Bigl(\itobjectsmod_{m+1}^{(1),\star} - \E\bigl(\itobjectsmod_{m+1}^{(1),\star}\bigm|\fF_{\tmod_{m-1/2}^{\star}}^{\,\xi}\bigr) - \bigl(\taumod_{m+1}^{\star\star}\bigr)^{-1 + \kappa}\Bigr)_+^p \Bigm|\fF_{\tmod_{m-1/2}^{\star}}^{\,\xi}\Biggr) \\
	&\leq \bigl(\taumod_{m+1}^{\star\star}\bigr)^{(q-p)(1 - \kappa)}\E\Biggl(\Big|\itobjectsmod_{m+1}^{(1),\star} - \E\bigl(\itobjectsmod_{m+1}^{(1),\star}\bigm|\fF_{\tmod_{m-1/2}^{\star}}^{\,\xi}\bigr)\Big|^q \Bigm|\fF_{\tmod_{m-1/2}^{\star}}^{\,\xi}\Biggr) \leq C_q\bigl(\taumod_{m+1}^{\star\star}\bigr)^{q\kappa-p(1 - \kappa)} \;.
\end{equs}
	The bound \eqref{e:itobjects1_fluctuations2_mid} then follows by choosing $q$ large enough depending on $p$.
\end{proof}

\begin{lemma}
	\label{le:itobjects2_fluctuations}
	Recall $\itobjectsmod^{(2)}_{m+1} $ defined in Lemma~\ref{le:smoothed_lollipop_addition_expansion}.
	For any $p \in \N^*$, there exists some constant $C_p$, such that for any $m \geq 1$ we have
	\begin{equ}
	\label{e:itobjects2_fluctuations1}
	\E \biggl( \Big| \one_{\taumod_{m+1} = \taumod_{m+1}^\star} \itobjectsmod^{(2)}_{m+1} - \E \bigl( \one_{\taumod_{m+1} = \taumod_{m+1}^\star} \itobjectsmod^{(2)}_{m+1} \big| \, \fF_{\tmod_m} \bigr) \Big|^p \Big| \, \fF_{\tmod_m} \biggr) \leq C_p \bigl(\taumod_{m+1}^\star\bigr)^{p(-1 + 2 \kappa)} \one_{\tmod_{m+1/2}^{\star} \neq \tmod_m} \;.
	\end{equ}
	Moreover, we have
	\begin{equs}
	\label{e:itobjects2_fluctuations2}
	\,& \E \biggl( \Bigl( \one_{\taumod_{m+1} = \taumod_{m+1}^\star} \itobjectsmod^{(2)}_{m+1}
	- \E \bigl( \one_{\{\taumod_{m+1} = \taumod_{m+1}^\star\} \cap \{\taumod_{m+1}^\star < \Lmod_m^{-1}\}} \itobjectsmod^{(2)}_{m+1} \big| \, \fF_{\tmod_m} \bigr) \\
	&\qquad \qquad \qquad \qquad \qquad \qquad \qquad \qquad \qquad - \bigl(\taumod_{m+1}^{\star}\bigr)^{-1 + \kappa} \Bigr)_+^p \Big| \, \fF_{\tmod_{m-1}} \biggr) \leq C_p \;.
	\end{equs}
\end{lemma}

\begin{proof}
	Applying Lemma~\ref{le:iteration_basic_bounds_mod} and Lemma~\ref{le:lollipop_negative_holder}, we obtain the crude estimate
\begin{equ}
	\E\Bigl(\big|\itobjectsmod^{(2)}_{m+1}\big|^{2p} \Bigm| \fF_{\tmod_m}\Bigr) \lesssim \bigl(\taumod_{m+1}^{\star}\bigr)^{-4p} \one_{\tmod_{m+1/2}^{\star} \neq \tmod_m } \;.
\end{equ}
	Combining this and \eqref{e:tau_neq_T_unlikely1_mod}, we get that for every $p \geq 1$,
\begin{equ}
	\E\Bigl(\one_{\taumod_{m+1} \neq \taumod_{m+1}^\star} \big|\itobjectsmod^{(2)}_{m+1}\big|^p \Bigm| \fF_{\tmod_m}\Bigr) \leq C_p \one_{\tmod_{m+1/2}^{\star} \neq \tmod_m } \;.
\end{equ}
	Then to show \eqref{e:itobjects2_fluctuations1}, it suffices to show 
	\begin{equ}
	\label{e:itobjects2_fluctuations1_mid}
	\E \biggl( \Big| \itobjectsmod^{(2)}_{m+1} - \E \bigl(  \itobjectsmod^{(2)}_{m+1} \big| \, \fF_{\tmod_m} \bigr)  \Big|^p  \Big | \, \fF_{\tmod_m } \biggr) \lesssim  \bigl(\taumod_{m+1}^\star \bigr)^{p(-1 + 2\kappa)} \one_{\tmod_{m+1/2}^{\star} \neq \tmod_m }  \;.
	\end{equ}
	Now we use the results in Section~\ref{sec:improved_stochastic_estimates_fluctuations} to bound the fluctuation of each term in \eqref{e:smoothed_lollipop_addition_expansion}. From now on we assume that $\{\tmod_{m+1/2}^{\star} \neq \tmod_m\}$ happens. 
	
\proofstep{Fluctuation of \eqref{e:smoothed_lollipop_addition_expansion1}}:
Applying Lemma~\ref{le:decay_mechanism_fluctuations}, we obtain
\begin{equ}
\E\Bigl(\Big|\eqref{e:smoothed_lollipop_addition_expansion1} - \E\bigl(\eqref{e:smoothed_lollipop_addition_expansion1}\bigm|\fF_{\tmod_m}\bigr)\Big|^p \Bigm|\fF_{\tmod_m}\Bigr)^{1/p} \lesssim \bigl(\taumod_{m+1}^{\star}\bigr)^{\frac12-2\lossdel} \bigl\|F_{(\Bmod_{m+1}^{\star}+\Sdetmod_{m+1}^{\star})(\tmod_{m+1}^{\star})}\bigr\|_{L^\infty} \;.
\end{equ}
It remains to bound $\bigl\|F_{(\Bmod_{m+1}^{\star}+\Sdetmod_{m+1}^{\star})(\tmod_{m+1}^{\star})}\bigr\|_{L^\infty}$. In the analysis of the term \eqref{e:unstablepart_curvature1}, the difference between this term and $\Ad_{\widetilde{g}_{m+1}^{\star}(\tmod_{m+1}^{\star}) } F^{12}_{\Bmod^{+, \star}_{m+1}} $ is discussed. The desired bound then follows from \eqref{e:iteration_F_bounds1_mod}. 
	
\proofstep{Fluctuation of \eqref{e:smoothed_lollipop_addition_expansion2}}:
This term is Gaussian when conditioned on $\fF_{\tmod_m}$. The desired result then directly follows from the proof of \eqref{e:unstablepart_curvature2} and Gaussian hypercontractivity.

\proofstep{Fluctuation of \eqref{e:smoothed_lollipop_addition_expansion3}}:
With the shorthand notation \eqref{e:smoothed_out_lollipop} and the schematic notation from Notation~\ref{prelim:notation-Landau-standard-schematic}, we expand the covariant curl as
\begin{equ}
\curl_{(\Bmod_{m+1}^{\star}+\Sdetmod_{m+1}^{\star})(\tmod_{m+1}^{\star})}v_{\Bmod_{m+1}^{\star}} = \curl_{\Bmod_{m+1}^{\star}(\tmod_{m+1}^{\star})}v_{\Bmod_{m+1}^{\star}} + \schO\Bigl(\Sdetmod_{m+1}^{\star}(\tmod_{m+1}^{\star})\,v_{\Bmod_{m+1}^{\star}}\Bigr) \;.
\end{equ}
By Lemma~\ref{le:curl_square_fluctuations}, we get that 
\begin{equs}
	\,&\E\Biggl(\Big|\bigl\|\curl_{\Bmod_{m+1}^{\star}(\tmod_{m+1}^{\star})}v_{\Bmod_{m+1}^{\star}}\bigr\|_{L^2}^2 - \E\Bigl(\bigl\|\curl_{\Bmod_{m+1}^{\star}(\tmod_{m+1}^{\star})}v_{\Bmod_{m+1}^{\star}}\bigr\|_{L^2}^2 \Bigm|\fF_{\tmod_m}\Bigr)\Big|^p \Bigm|\fF_{\tmod_m}\Biggr)^{1/p} \\
	&\lesssim \bigl(\taumod_{m+1}^{\star}\bigr)^{\frac12-2\lossdel}\bigl(\tmod_{m+1}^{\star}-\tmod_{m+1/2}^{\star}\bigr)^{-1} \lesssim \bigl(\taumod_{m+1}^{\star}\bigr)^{-1+2\kappa} \;.
\end{equs}
For the cross terms in the expansion of the $L^2$ norm, they can be estimated easily using the bound \eqref{eq:Sdet-L-infty-estimates} on $\Sdetmod_{m+1}^{\star}(\tmod_{m+1}^{\star})$. We omit the details here. 

This concludes the proof of \eqref{e:itobjects2_fluctuations1}. For \eqref{e:itobjects2_fluctuations2}, we first remove the indicator on the event $\{ \taumod_{m+1}^\star < \Lmod_m^{-1} \}$. By \eqref{e:taumod_rare_event2}, we have 
\begin{equs}
	\,& \E \Bigl( \one_{\{\taumod_{m+1}^\star = \Lmod_m^{-1}\}} \big|\one_{\taumod_{m+1} = \taumod_{m+1}^\star} \itobjectsmod^{(2)}_{m+1}\big|^p \Bigm| \fF_{\tmod_{m-1}} \Bigr) \\
	&\lesssim \E \Bigl( \one_{\{\taumod_{m+1}^\star = \Lmod_m^{-1}\}} \bigl(\taumod_{m+1}^\star\bigr)^{-2p} \Bigm| \fF_{\tmod_{m-1}} \Bigr) \leq C_p \;.
\end{equs}
After removing this indicator, the rest of the argument is the same as in the end of Lemma~\ref{le:itobjects1_fluctuations}. The lemma is proved. 
\end{proof}

Finally, we prove the main result of this section, Proposition~\ref{prop:modified-iteration-uniform-bounds}.

\begin{proof}[Proof of Proposition~\ref{prop:modified-iteration-uniform-bounds}]
For this proof, we stick to the notation of the abstract result Proposition~\ref{prelim:prop-X-supermartingale}, so we will seek to obtain conditional $r$th moment bounds (so in particular, $r$ is a free parameter, rather than the one fixed in \eqref{e:small_constants}), while the parameter $p$ will play the same role as in Proposition~\ref{prelim:prop-X-supermartingale}. Consider the setting of Proposition~\ref{pr:main_iteration}. To place our iteration scheme within the abstract setup of Propositions~\ref{prelim:prop-X-supermartingale} and~\ref{prelim:prop-stopping-time}, we let $K\geq 1$ be a sufficiently large constant (depending on $r$ and $p$, our parameters from \eqref{e:small_constants} and \eqref{e:constants_global_iteration}, and the earlier estimates in this section). For $m \geq 1$, we let $X_m = \covariantenergy(\Asmoothpartmod(\tmod_m), \phismoothpartmod(\tmod_m))$, $D_m = 3C_{\eqref{e:one_step_change}}  \bigl( \taumod_{m+1}^\star \bigr)^{-1}$, $Y_{m+1} = \fluctmod_{m+1}$, $\mc{G}_m = \mc{F}_{\tmod_m}$, and $\varep = 2/3$. We also let $\mc{G}_0 = \mc{F}_0$. With these choices, Proposition~\ref{pr:main_iteration} implies that the recursive bound~\eqref{prelim:eq-recursive} and decay lower bound~\eqref{prelim:eq-decay} in the assumptions of Proposition~\ref{prelim:prop-X-supermartingale} are satisfied for all $m \geq 1$.  For the decay lower bound~\eqref{prelim:eq-decay}, we remark that this also uses the definition of $T^{(1)}$ (Definition~\ref{defn:tau}) and our choice of parameters~\eqref{e:constants_global_iteration} (in particular $\gamma_1 / 2 < \gamma_0$, and $3/2 < \gamma_0$). Furthermore, Proposition~\ref{pr:main_iteration} also implies that 
\begin{equation*}
\E[(Y_{m+1} - \varep D_m)_+^p ~|~ \mc{F}_{\tmod_{m-1}}] \leq K,
\end{equation*} 
and hence the first fluctuation estimate \eqref{prelim:eq-fluctuation} is satisfied. 

For $m = 0$, we set $X_0 = \covariantenergy(A_0^+, \higgs_0^+)$, $D_0 = X_0^{1/2} / \log X_0$, and $Y_1 = X_1 - X_0 + D_0$. From their definitions, it then directly follows that the conditions \eqref{prelim:eq-recursive} and \eqref{prelim:eq-decay} in Proposition \ref{prelim:prop-X-supermartingale} are satisfied for $m=0$. From \eqref{e:main_iteration_fluct_init}, we obtain that
\begin{equation*}
\E[(Y_{1} - \varep D_0)_+^p ~|~ \mc{F}_0] \leq  \E[X_1^p ~|~ \mc{F}_0] \leq C_p \bigl(X_0^{4p} + \Lmod_0^{4p}\bigr). 
\end{equation*}
Hence, by choosing $Y_0 = C (X_0^4 + \Lmod_0^4 )$ with $C$ large enough, we can ensure that the second fluctuation bound \eqref{prelim:eq-fluctuation-initial} is satisfied.

For any $r \geq 1$, we can now choose any $p > r + \frac{1}{2}$ and apply Proposition~\ref{prelim:prop-X-supermartingale}, which yields that 
\begin{equs}
\sup_{m \geq 0} \E[ X_m^r ~|~ \mc{F}_0] \leq C_p (X_0^{r} + Y_0^{r} +1), \qquad \text{and} \qquad \limsup_{m\rightarrow \infty} \E \bigl[ X_m^r \big| ~ \mc{F}_0 \bigr] \leq C_p \;.
\end{equs}
The first two estimates in Proposition~\ref{prop:modified-iteration-uniform-bounds} now follow. 

For the final two estimates, we first verify the condition \eqref{prelim:eq-time-increment-moments} in Proposition~\ref{prelim:prop-stopping-time}. Towards this end, by Lemmas~\ref{lemma:time-step-cond-exp-bd} and~\ref{lemma:taumod-star-det-bound}, for any $m \geq 1$, $q_1 \geq 1$, we have that
\begin{equ}
\E[\taumod_{m+1}^{-q_1} ~|~ \mc{F}_{\tmod_m}] \leq C_p\bigl( X_m^{q_1} + \Lmod_m^{q_1}+ 1 \bigr) \;.
\end{equ}
Taking an additional conditional expectation over $\mc{F}_{\tmod_{m-1}}$ and applying \eqref{e:Lmod_expectation_bound} yields
\begin{equ}
\E[\taumod_{m+1}^{-q_1} ~|~ \mc{F}_{\tmod_{m-1}}] \leq \E[C_p (X_m^{q_1} + X_m^{q_1/10}) ~|~ \mc{F}_{\tmod_{m-1}}] + C_p \leq \E[X_m^{2q_1} ~|~ \mc{F}_{\tmod_{m-1}}] + C_p \;,
\end{equ}
where we applied Young in the last step (possibly enlarging $C_p$). 
Taking $q_1 = \max(r/2, 1)$ (say), it follows that \eqref{prelim:eq-time-increment-moments} holds with $q_2 = 2q_1$. Since our fluctuation bound \eqref{e:main_iteration_fluct} holds for arbitrarily large $p$, we may take $p > \max(r + \frac{1}{2}, q_2)$ so that the conditions of Proposition~\ref{prelim:prop-stopping-time} are satisfied (after potentially enlarging $K$ by a constant depending only on $p$ and $\stepupperbound$). The last two estimates of Proposition~\ref{prop:modified-iteration-uniform-bounds} (as well as the fact that $M_t < \infty$ a.s.\ for all $t > 0$) now directly follow from Proposition~\ref{prelim:prop-stopping-time}, which provides conditional $r$th moment estimates on $X_{M_t+1}$.
\end{proof}

\subsection{Conclusion}
\label{sec:conclusions}

In this section, we upgrade our estimate of the covariant energy from the random times $\tmod_m$ to deterministic times. The main result is as follows.

\begin{proposition}\label{prop:det-time-energy-bound}
Let $A_0\in \Cs_x^{-\kappa}(\T^2, \frkg^2)$ and $\higgs_0 \in \Cs_x^{-\kappa}(\T^2, \V)$. Let $A_0^+ = 0, A_0^- = A_0, \higgs_0^+ = 0, \higgs_0^- = \higgs_0$. Let 
\begin{equ}
\Lmod_0 = \covariantenergy(A_0^+, \higgs_0^+)^{\frac{1}{10}} +  \bigl\|A_0^-\bigr\|_{\covHolder{A_0^+, F}^{-\kappa}}^{\frac{1}{\fa}} + \bigl\|\higgs_0^- \bigr\|_{\covHolder{A_0^+}^{-\kappa}}^{\frac{1}{\fa}} \;,
\end{equ}
so that $(A_0^+, A_0^-, \higgs_0^+, \higgs_0^-)$ is $\Lmod_0$-admissible. Consider the iteration scheme from Definition~\ref{defn:time_steps_mod}. For all $\stepupperbound$ small enough, we have that for any $p \geq 1$, 
\begin{equs}
\sup_{t \geq 0} \E\bigl[\invnorm{(A(t),\higgs(t))}{-\kappa}^p ~\big|~ \mc{F}_0 \bigr] &\leq C(A_0, \higgs_0, \stepupperbound, p)  \;, \\
\limsup_{t \rightarrow \infty} \E\bigl[\invnorm{(A(t),\higgs(t))}{-\kappa}^p  ~\big|~ \mc{F}_0 \bigr] &\leq C(\stepupperbound, p) \;.
\end{equs}
\end{proposition}

Before we proceed to prove Proposition~\ref{prop:det-time-energy-bound}, we first quickly show how the main global well-posedness results (Theorems~\ref{thm:global_soft} and~\ref{thm:global_hard}) follow.

\begin{proof}[Proof of Theorems~\ref{thm:global_soft} and~\ref{thm:global_hard}]
Fix $\stepupperbound$ small enough so that Proposition~\ref{prop:modified-iteration-uniform-bounds}, Lemma~\ref{le:linear_mod_conditional_moments}, and Proposition~\ref{prop:det-time-energy-bound} apply (recall Remark~\ref{remark:logtimestep-stepupperbound}). We first prove Theorem~\ref{thm:global_hard}. The first two items~\ref{main:item-global-SYMH} and~\ref{main:item-global-initial} follow by the definition of our iteration scheme and the local well-posedness result Theorem~\ref{main:thm-covariant-local}. To avoid confusion, we note that $\widetilde{g}_{m+1}$ in \ref{main:item-global-initial} of Theorem \ref{thm:global_hard} corresponds to $\widetilde{g}_{m+1}(\tmod_m)$ in the notation of Definition \ref{defn:time_steps_mod}. The third item~\ref{main:item-global-decomposition} also follows by the definition of our iteration scheme. For the fourth item \ref{main:item-global-smooth-part}, we apply Proposition~\ref{prop:modified-iteration-uniform-bounds} with the smooth parts $A_0^+, \higgs_0^+ = 0$, the rough parts $A_0^-, \higgs_0^-$ taken to be the initial data $A_0, \higgs_0$ in the statement of Theorem~\ref{thm:global_hard}, and $\Lmod_0$ taken as in the statement of Proposition~\ref{prop:det-time-energy-bound}. Item \ref{main:item-global-rough-part} follows by (recall $\Aroughpartmod = \Alinmod + \Aremmod$, and similarly for $\phiroughpartmod$) Lemma~\ref{le:linear_mod_conditional_moments} and the remainder bounds in \eqref{e:Aroughpart_crude_bound_mod} and \eqref{e:phiroughpart_crude_bound_mod}. Finally, note that the fact that $\tmod_m \uparrow \infty$ a.s. follows from the fact that $M_t < \infty$ a.s. for all $t > 0$, where $M_t$ is as in Proposition~\ref{prop:modified-iteration-uniform-bounds}. This proves Theorem~\ref{thm:global_hard}.

For Theorem~\ref{thm:global_soft}, note that the first two items are exactly the same as in Theorem~\ref{thm:global_hard}, and thus they have already been shown. The third item follows directly from Proposition~\ref{prop:det-time-energy-bound}.
\end{proof}

In the remainder of this section, we work towards the proof of Proposition~\ref{prop:det-time-energy-bound}. We begin with one of the main intermediate results towards the proposition, which estimates the $\invnorm{\place}{-\kappa}$ norm in terms of inverse powers of the iteration time steps $\taumod_m$. The proof is mostly a tedious exercise, so on an initial reading, the reader may want to read the statement of the next lemma but then skip the proof and directly move on to the other main intermediate results.

\begin{lemma}\label{lemma:inv-norm-estimate-time-step}
There is a constant $C$ such that for $\stepupperbound$ small enough, the following holds. Consider the iteration scheme from Definition~\ref{defn:time_steps_mod}. For all $m \geq 0$, we have that
\begin{equ}
\sup_{t \in (\tmod_m, \tmod_{m+1}]} \invnorm{(A(t),\higgs(t))}{-\kappa} \leq C (\taumod_{m+1})^{-3-\kappa} \abs{\log \taumod_{m+1}}^2   \;.
\end{equ}
\end{lemma}
\begin{proof}
Due to the inductive definition of the iteration scheme, it suffices to take $m = 0$. Fix $t \in (0, \tmod_1]$. From the definition of the $\invnorm{\place}{-\kappa}$ norm (recall Definition~\ref{main:def-invnorm}), we want to find a gauge transformation $g \in \CG_{1-\kappa}$ such that there is a Yang--Mills gradient flow $\hat{B} \in C_t^0 \Cs_x^{-\kappa}([0, 1] \times \T^2) \cap C_{t,x}^\infty((0, 1] \times \T^2)$ with $\hat{B}(0) = g \act A(t)$, such that, with $\hat{\Phi}(s) = e^{s \Delta_{\hat{B}(0)}} (g \act \Phi(t))$ in the following, we have that
\begin{equ}
\bigl\| s^{\frac{1+\kappa}{2}} F_{\hat{B}}(s) \bigr\|_{L_{s,x}^\infty([0, 1])} + \bigl\|s^{\frac{2+\kappa}{2}} \covd_{\hat{B}} F_{\hat{B}}(s)\bigr\|_{L_{s,x}^\infty([0, 1])} + \bigl\|s^{\frac{\kappa}{2}}\hat{\Phi}(s)\bigr\|_{L_{s,x}^\infty([0, 1])} \lesssim (\taumod_{1})^{-3-\kappa} \abs{\log \taumod_1}^2  \;.
\end{equ}
We first make some preliminary reductions. Since $F_{\hat{B}}$ solves a covariant heat equation (Lemma~\ref{lemma:curvature-equation}), by our heat kernel bounds \eqref{eq:higher-order-pB-ptwise-bd}, we have that
\begin{equ}
\bigl\|s^{\frac{2+\kappa}{2}} \covd_{\hat{B}} F_{\hat{B}}(s)\bigr\|_{L_{s,x}^\infty([0, 1])} \lesssim \bigl(1 + \Fnorm{F_{\hat{B}}}_{1, 0}\bigr)\bigl\| s^{\frac{1+\kappa}{2}} F_{\hat{B}}(s) \bigr\|_{L_{s,x}^\infty([0, 1])} \;,
\end{equ}
and moreover from the definition (Definition~\ref{def:Fnorm}) of $\Fnorm{\place}$, we have that
\begin{equ}
\Fnorm{F_{\hat{B}}}_{1, 0} \lesssim \bigl\| s^{\frac{1+\kappa}{2}} F_{\hat{B}}(s) \bigr\|_{L_{s,x}^\infty([0, 1])} \;.
\end{equ}
Therefore, it just suffices to show that
\begin{equ}\label{eq:inv-norm-bd-intermediate}
\bigl\| s^{\frac{1+\kappa}{2}} F_{\hat{B}}(s) \bigr\|_{L_{s,x}^\infty([0, 1])} +  \bigl\|s^{\frac{\kappa}{2}}\hat{\Phi}(s)\bigr\|_{L_{s,x}^\infty([0, 1])} \lesssim (\taumod_{1})^{-\frac{3+\kappa}{2}} \abs{\log \taumod_1} \;.
\end{equ}
We proceed as follows. By construction (recall Definition~\ref{defn:time_steps_mod}), we have that
\begin{equs}
A(t) &= \Bmod(t) + \linear[\Bmod](t) + \Sdetmod(t) + \Arem(t) \;, \\
\higgs(t) &= \higgsdetmod(t) + \philinear[\Bmod](t) + \phirem(t) \;,
\end{equs}
where $\Arem(t) = A(t) - \Bmod(t) - \linear[\Bmod](t) - \Sdetmod(t)$ and $\phirem(t) = \higgs(t) - \higgsdetmod(t) - \philinear[\Bmod](t)$ (we do not put tildes on these remainders, as the remainders with tildes are only defined at $\tmod_{m+1}$ in the iteration scheme \dash recall Definition~\ref{defn:time_steps_mod}). Here, for notational simplicity, since we are only working on the first time interval, we have dropped the subscript $1$. To obtain the gauge transformation $g$, we will seek to apply Proposition~\ref{pr:admissible_timescale_general}, which recall provides local well-posedness of the Yang--Mills gradient flow. We proceed to verify the assumptions of the proposition. Let $T = \taumod_1$ in the following. We will apply the proposition with this choice of $T$. Let $\Bmod_{t}(s) = \Bmod(t + s)$ be the shifted version of $\Bmod$. Note that a priori $\Bmod$ is only defined on $[0, \tmod_1]$, but since $\Bmod$ is smooth at positive times, and since the Yang--Mills gradient flow is globally well-posed for smooth initial data, we may extend $\Bmod$ to $[0, \infty)$. Thus, $\Bmod_t$ is well-defined at any time at which it is evaluated. This will be used implicitly in the following.

We also recall that the initial data $\Bmod(0) = \widetilde{g}(0) \act (A_0^+ + A_0^-)$, and $\higgsdetmod(0) = \widetilde{g}(0) \act (\higgs_0^+ + \higgs_0^-)$. By construction of the iteration scheme, the collection $(A_0^+, A_0^-, \higgs_0^+, \higgs_0^-)$ is $\Lmod_0$-admissible, which by definition implies that
\begin{equ}\label{eq:inv-norm-bd-admissibility}
\bigl\|A_0^- \bigr\|_{\covHolder{A_0^+, F, T}^{-\kappa}} \leq \Lmod_0^{\fa}, \quad \bigl\|\higgs_0^-\bigr\|_{\covHolder{A_0^+, T}^{-\kappa}} \leq \Lmod_0^{\fa}.
\end{equ}
Combining the second estimate with Lemma~\ref{le:B-YMHFsmooth-heat-comparison} (recalling that in the iteration scheme, we obtained $\Bmod, \widetilde{g}(0)$ by applying Proposition~\ref{pr:admissible_timescale_general} to $(A_0^+, A_0^-)$), we obtain
\begin{equ}\label{eq:inv-norm-bd-adjusted-admissibility}
\max_{0 \leq j \leq 1} \bigl\|s^{\frac{j + \kappa}{2}} \covd_{\Bmod(s)}^{(j)} e^{s \Delta_{\Bmod(0)}} (\widetilde{g}(0) \act \higgs_0^-)\bigr\|_{L_{s,x}^\infty([0, T])} \lesssim T^{-\lossdel} \bigl\|s^{\frac{\kappa}{2}} e^{s \Delta_{A_0^+}} \higgs_0^-\bigr\|_{L_{s,x}^\infty([0, T])} \lesssim T^{-2\fa},
\end{equ}
where we used that $T \leq \Lmod_0^{-1}$ by definition (Definition~\ref{defn:L_admissible}), as well as the fact that $\lossdel \ll \fa$ (recall \eqref{e:small_constants} and \eqref{e:constants_global_iteration}) in the final inequality. 

We record some additional preliminary estimates. Due to the monotonicity of the curvature (Remark~\ref{remark:curvature-pointwise-in-time-bound-Fnorm}), we have that
\begin{equ}
\gronwall_{T}(C_0 F_{\Bmod_t}) \bigl(1 + \Fnorm{F_{\Bmod_t}}_{T, 0}\bigr)^{100} \leq \gronwall_T(C_0 F_{\Bmod}) \bigl(1 + \Fnorm{F_{\Bmod}}_{T, 0}\bigr)^{100} \leq T^{-\lossdel}, 
\end{equ}
where the second inequality follows because the iteration scheme is constructed so that Assumption~\ref{local:assumption-B} holds on every interval. Additionally, by Corollary~\ref{cor:curvature-higher-regularity-estimates}, we have that
\begin{equ}
\Fnorm{F_{\Bmod_t}}_{T, 2} \lesssim \Fnorm{F_{\Bmod_t}}_{T, 0} + \bigl(\Fnorm{F_{\Bmod_t}}_{T, 0}\bigr)^6 \;.
\end{equ}
Combining, we thus obtain
\begin{equ}\label{eq:Bmod_t-gronwall-control}
M_T(\Bmod_t, 2) = \gronwall_T(C_0 F_{\Bmod_t})\bigl(1 + \Fnorm{F_{\Bmod_t}}_{T, 2}\bigr)^2 \lesssim T^{-\lossdel}.
\end{equ}
With these preliminary estimates in hand, we now define $B_0^+ = \Bmod_t(0)$, and $B_0^- = A(t) - \Bmod_t(0) = \linear[\Bmod](t) + \Sdetmod(t) + \Arem(t)$. We will apply Proposition~\ref{pr:admissible_timescale_general} with this choice of $B_0^+, B_0^-$. To do so, we need to estimate
\begin{equ}
K = \bigl\|\linear[\Bmod](t)\bigr\|_{\covHolder{\Bmod_t(0), F, T, 2}^{-\kappa}} + \bigl\|\Sdetmod(t)\bigr\|_{\covHolder{\Bmod_t(0), F, T, 2}^{-\kappa}} + \bigl\|\Arem(t)\bigr\|_{\covHolder{\Bmod_t(0), F, T, 2}^{-\kappa}} + 1\;. 
\end{equ}
Due to the definition of $T = \taumod_1$ (recall Definition~\ref{defn:tau_L}), we have that
\begin{equ}
\bigl\|\linear[\Bmod](t)\bigr\|_{\covHolder{\Bmod_t(0), F, T}^{-\kappa}} \leq \frac{1}{2} T^{-\lambda_1}.
\end{equ}
The above norm controls up to one derivative of the object (Definition~\ref{defn:covariant_Holder}), while we require control of up to two derivatives. Thus, we perform one additional estimate: we apply our heat kernel bounds \eqref{eq:qB-ptwise-bounds} in conjunction with \eqref{eq:Bmod_t-gronwall-control} to estimate
\begin{equ}
\bigl\|\linear[\Bmod](t)\bigr\|_{\covHolder{\Bmod_t(0), F, T, 2}^{-\kappa}} \lesssim T^{-\lossdel} \bigl\|\linear[\Bmod](t)\bigr\|_{\covHolder{\Bmod_t(0), F, T}^{-\kappa}}  \lesssim T^{-2\lambda_1},
\end{equ}
where the final inequality follows due to our choice $\lossdel \ll \lambda_1$ (recall \eqref{e:small_constants}). Next, we proceed to estimate $\Sdetmod(t)$. By \eqref{eq:Sdet-L-infty-estimates} from Proposition~\ref{prop:Sdet-higgsdet-estimates}, combined with the fact that $(A_0^+, A_0^-, \higgs_0^+, \higgs_0^-)$ is $\Lmod_0$-admissible (recall in particular \eqref{eq:inv-norm-bd-admissibility} and \eqref{eq:inv-norm-bd-adjusted-admissibility}) and the fact that $T = \taumod_1 \leq \Lmod_0^{-1}$, we have that
\begin{equ}
\bigl\|\Sdetmod(t)\bigr\|_{L_x^\infty} \lesssim T^{-\lossdel}\bigl(T^{\frac{1}{2}-\frac{2}{r}} \bigl\|\higgs_0^+\bigr\|_{L_x^r}^2 + T^{\frac{1}{2}-\kappa} T^{-4\fa}\bigr) \lesssim T^{\frac{1}{2} - \frac{2}{r}-\lossdel} \bigl\|\higgs_0^+\bigr\|_{L_x^r}^2 + T^{\frac{1}{2}-10\fa} \; ,
\end{equ}
where we used that $\fa$ is the largest parameter appearing in the above (recall \eqref{e:small_constants} and \eqref{e:constants_global_iteration}). Combining this with the heat kernel estimates \eqref{eq:qB-ptwise-bounds} and \eqref{eq:Bmod_t-gronwall-control}, we have that
\begin{equ}
\bigl\|\Sdetmod(t)\bigr\|_{\covHolder{\Bmod_t(0), F, T, 2}^{-\kappa}} \lesssim T^{-\lossdel} \bigl\|\Sdetmod(t)\bigr\|_{L_x^\infty} \lesssim T^{\frac{1}{2}-\frac{3}{r}} \bigl\| \higgs_0^+\bigr\|_{L_x^r}^2 + 1\; .
\end{equ}
In the final estimate, we used that $\lossdel \ll 1/r$ (recall \eqref{e:small_constants}). Finally, to estimate $\Arem(t)$, combining \eqref{local:eq-main-R-estimate} from Proposition~\ref{local:prop-main} with the same heat kernel estimates and \eqref{eq:Bmod_t-gronwall-control}, we have that
\begin{equ}
\bigl\|\Arem(t)\bigr\|_{\covHolder{\Bmod_t(0), F, T, 2}^{-\kappa}} \lesssim T^{-\lossdel} \bigl\|\Arem(t)\bigr\|_{L_x^\infty} \lesssim 1 \; .
\end{equ}
Combining our estimates of $\linear[\Bmod](t)$, $\Sdetmod(t)$, and $\Arem(t)$, for some sufficiently large constant $C_1$, we have that
\begin{equ}
K  \leq C_1 \bigl(T^{\frac{1}{2}-\frac{3}{r}} \bigl\|\higgs_0^+\bigr\|_{L_x^r}^2 + T^{-2\lambda_1}\bigr) \;. 
\end{equ}
We now check the assumption \eqref{e:admissible_timescale_smallness}. If $\stepupperbound$ is sufficiently small, we have that
\begin{equ}
M_T(\Bmod_t, 2) T^{\frac{1-\kappa}{2}} K \leq C_1 T^{-\lossdel} \bigl(T^{1 - \frac{3}{r} - \frac{\kappa}{2}} \bigl\|\higgs_0^+\bigr\|_{L_x^r}^2 + T^{\frac{1-\kappa}{2} - 2\lambda_1}\bigr) \ll 1 \;,
\end{equ}
and thus \eqref{e:admissible_timescale_smallness} holds. Here, we used that $T \leq \bigl\|\higgs_0^+\bigr\|_{L_x^r}^{-\gamma_0}$ (recall Definition~\ref{defn:tau}), as well as the fact that $\gamma_0 - 2$ is larger than all other parameters appearing in the above (recall \eqref{prelim:eq-parameters} and \eqref{e:constants_global_iteration}). Applying Proposition~\ref{pr:admissible_timescale_general}, we obtain a gauge transformation $g \in \CG_{1-\kappa}$ and a Yang--Mills gradient flow $\hat{B} \in C_t^0 \Cs_x^{-\kappa}([0, T] \times \T^2) \cap C_{t, x}^\infty((0, T] \times \T^2)$ such that $\hat{B}(0) = g \act A(t)$. Moreover, by the estimate \eqref{e:Fnorm_bound_roughB_general}, we have that
\begin{equ}
\bigl\|s^{\frac{1+\kappa}{2}} F_{\hat{B}}\bigr\|_{L_{s,x}^\infty([0, T])} \leq \bigl\|s^{\frac{1+\kappa}{2}} F_{\Bmod_t}\bigr\|_{L_{s,x}^\infty([0, T])} + \bigO(K) \;.
\end{equ}
Due to the monotonicity of the curvature (Remark~\ref{remark:curvature-pointwise-in-time-bound-Fnorm}), we may estimate (applying the above estimate in the second inequality)
\begin{equ}
\bigl\|s^{\frac{1+\kappa}{2}} F_{\hat{B}} \bigr\|_{L_{s,x}^\infty([0, 1])} \leq T^{-\frac{1+\kappa}{2}} \bigl\|s^{\frac{1+\kappa}{2}} F_{\hat{B}}\bigr\|_{L_{s,x}^\infty([0, T])} \leq  T^{-\frac{1+\kappa}{2}} \bigl(\bigl\|s^{\frac{1+\kappa}{2}} F_{\Bmod_t}\bigr\|_{L_{s,x}^\infty([0, T])} + \bigO(K)\bigr)\;.
\end{equ}
To finish the estimate of $F_{\hat{B}}$, we again use the monotonicity of the curvature (Remark~\ref{remark:curvature-pointwise-in-time-bound-Fnorm}), the construction of $\Bmod$ (which was obtained by applying Proposition~\ref{pr:admissible_timescale_general} to $A_0^+, A_0^-$), the estimate \eqref{e:Fnorm_bound_roughB_general} from Proposition~\ref{pr:admissible_timescale_general}, and the fact that $A_0^-$ is $\Lmod_0$-admissible by construction (recall \eqref{eq:inv-norm-bd-admissibility}) in order to estimate
\begin{equs}
\bigl\|s^{\frac{1+\kappa}{2}} F_{\Bmod_t}\bigr\|_{L_{s,x}^\infty([0, T])} &\leq \bigl\|s^{\frac{1+\kappa}{2}} F_{\Bmod}\bigr\|_{L_{s,x}^\infty([0, T])} \\
&\leq \bigl\|s^{\frac{1+\kappa}{2}} F_{\Bmod^{+, \star}}\bigr\|_{L_{s,x}^\infty([0, T])} + \bigO((\taumod_1^\star)^{-2\fa})  \\
&\lesssim \bigl\|F_{A_0^+}\bigr\|_{L_x^2} + T^{-2\fa} \;,
\end{equs}
where $\Bmod^{+, \star}$ denotes the YM gradient flow starting from $A_0^+$, and in the final step we used the fact that $F_{\Bmod^{+, \star}}$ solves a covariant heat equation (Lemma~\ref{lemma:curvature-equation}) and our heat kernel bounds \eqref{eq:zeroth-order-pB-ptwise-bd}. Combining the previous few estimates and using that $T = \taumod_1 \leq \taumod_1^\star \leq \bigl\|F_{A_0^+}\big\|_{L_x^2}^{-1} \log(\bigl\|F_{A_0^+}\bigr\|_{L_x^2}\vee 2\bigr)$ and $T \leq \bigl\|\higgs_0^+\bigr\|_{L_x^r}^{-\gamma_0}$ (recall Definition~\ref{defn:tau}), we obtain
\begin{equ}
\bigl\|s^{\frac{1+\kappa}{2}} F_{\hat{B}} \bigr\|_{L_{s,x}^\infty([0, 1])}  \lesssim T^{-\frac{1+\kappa}{2}} \bigl( T^{-1} \abs{\log T} + T^{\frac{1}{2}-\frac{3}{r}} \bigl\|\higgs_0^+\bigr\|_{L_x^r}^2\bigr) \lesssim T^{ - \frac{3+\kappa}{2}} \abs{\log T}\;.
\end{equ}
It remains to estimate $\hat{\Phi}$ as in \eqref{eq:inv-norm-bd-intermediate}. Recall that $\hat{\Phi}(s) = e^{s \Delta_{\hat{B}(0)}} (g \act \Phi(t))$. By Lemma~\ref{le:B-YMHFsmooth-heat-comparison}, \eqref{eq:Bmod_t-gronwall-control}, and the fact that $T^{\frac{1-\kappa}{2}} K \ll 1$ (which we previously verified), we have that
\begin{equ}
\bigl\|s^{\frac{\kappa}{2}} \hat{\Phi}(s)\bigr\|_{L_{s,x}^\infty([0, T])} \lesssim T^{-\lossdel} \bigl\|s^{\frac{\kappa}{2}} e^{s \Delta_{\Bmod_t(0)}} \Phi(t)\bigr\|_{L_{s,x}^\infty([0, T])} \;.
\end{equ}
Recalling that $\Phi(t) = \higgsdetmod(t) + \philinear[\Bmod](t) + \phirem(t)$, we proceed to individually estimate these terms. If $t \geq T/2$, then by \eqref{eq:higgsdet-L-infty-estimates}, the fact that $\higgs_0^-$ is $\Lmod_0$-admissible (which recall implies \eqref{eq:inv-norm-bd-adjusted-admissibility}), and our heat kernel bounds \eqref{eq:zeroth-order-pB-ptwise-bd}, we have that 
\begin{equ}
\bigl\|s^{\frac{\kappa}{2}} e^{s \Delta_{\Bmod_t(0)}} \higgsdetmod(t)\bigr\|_{L_{s,x}^\infty([0, T])} \leq T^{\frac{\kappa}{2}} \bigl\|\higgsdetmod(t)\bigr\|_{L_x^\infty} \lesssim T^{-\lossdel} \bigl(T^{\frac{\kappa}{2}- \frac{1}{r}} \bigl\|\higgs_0^+\bigr\|_{L_x^r} + T^{-2\fa}\bigr) \;. 
\end{equ}
If instead $t < T/2$, then we decompose $\higgsdetmod(t) = \shiggsdet(t) + \rlinhiggsdet(t) + \remhiggsdet(t)$, where $\shiggsdet$ solves \eqref{e:shiggsdet_eq} (with the $B, \Sdet$ therein taken to be $\Bmod, \Sdetmod$, and $\widetilde{g}(0) \act \higgs_0^+$ as the initial data), $\rlinhiggsdet(t)$ solves \eqref{eq:rlinhiggsdet-equation} (with the $B$ therein taken to be $\Bmod$ and $\widetilde{g}(0) \act \higgs_0^-$ as the initial data), and $\remhiggsdet = \higgsdetmod - \shiggsdet - \rlinhiggsdet$. Since $t < T/2$, we have that (applying \eqref{eq:inv-norm-bd-adjusted-admissibility} in the second step)
\begin{equ}
\bigl\|s^{\frac{\kappa}{2}} e^{s \Delta_{\Bmod_t(0)}} \rlinhiggsdet(t)\bigr\|_{L_{s,x}^\infty([0, T/2])} \leq \bigl\|s^{\frac{\kappa}{2}} e^{s \Delta_{\Bmod(0)}} (\widetilde{g}(0) \act \higgs_0^-)\bigr\|_{L_{s,x}^\infty([t, T])} \lesssim T^{-2\fa} \;.
\end{equ}
Then by our heat kernel bounds \eqref{eq:zeroth-order-pB-ptwise-bd}, we have that
\begin{equ}
\bigl\|s^{\frac{\kappa}{2}} e^{s \Delta_{\Bmod_t(0)}} \rlinhiggsdet(t)\bigr\|_{L_{s,x}^\infty([0, T])} \lesssim \bigl\|s^{\frac{\kappa}{2}} e^{s \Delta_{\Bmod_t(0)}} \rlinhiggsdet(t)\bigr\|_{L_{s,x}^\infty([0, T/2])} \lesssim T^{-2\fa} \;. 
\end{equ}
For $\shiggsdet$, we use that $\bigl\|\shiggsdet(t)\bigr\|_{L_x^r} \leq \bigl\|\higgs_0^+\bigr\|_{L_x^r}$ (which e.g.\ follows by the Bochner formula \eqref{eq:shiggsdet-parabolic-bochner} and the fact that the $L_x^r$ norm is gauge invariant), our heat kernel bounds \eqref{eq:zeroth-order-pB-ptwise-bd}, and the fact that $\frac{1}{r} \ll \kappa$ (recall \eqref{e:small_constants}) to estimate
\begin{equ}
\bigl\|s^{\frac{\kappa}{2}} e^{s \Delta_{\Bmod_t(0)}} \shiggsdet(t) \bigr\|_{L_{s,x}^\infty([0, T])} \lesssim T^{\frac{\kappa}{2}-\frac{1}{r}} \bigl\|\higgs_0^+\bigr\|_{L_x^r} \;.
\end{equ}
For $\remhiggsdet$, we apply Lemma~\ref{le:remhiggsdet_bound}, which implies that $\bigl\|\remhiggsdet(t)\bigr\|_{L_x^\infty} \ll 1$ (here we also use that $T \leq \bigl\|\higgs_0^+\bigr\|_{L_x^r}^{-\gamma_0}$ (by Definition~\ref{defn:tau}) as well as \eqref{eq:inv-norm-bd-adjusted-admissibility}), as well as our heat kernel bounds \eqref{eq:zeroth-order-pB-ptwise-bd} to estimate
\begin{equ}
\bigl\|s^{\frac{\kappa}{2}} e^{s \Delta_{\Bmod_t(0)}} \remhiggsdet(t)\bigr\|_{L_{s,x}^\infty([0, T])} \lesssim 1 \;.
\end{equ}
Combining, we have thus shown that in both cases $t \geq T/2$ and $t < T/2$, we have that
\begin{equ}
\bigl\|s^{\frac{\kappa}{2}} e^{s \Delta_{\Bmod_t(0)}} \higgsdetmod(t)\bigr\|_{L_{s,x}^\infty([0, T])} \lesssim T^{-\lossdel} \bigl(T^{\frac{\kappa}{2}- \frac{1}{r}} \bigl\|\higgs_0^+\bigr\|_{L_x^r} + T^{-2\fa}\bigr) \;. 
\end{equ}
Next, by Definition~\ref{defn:tau_L}, we have that
\begin{equ}
\bigl\|s^{\frac{\kappa}{2}} e^{s \Delta_{\Bmod_t(0)}} \philinear[\Bmod](t)\bigr\|_{L_{s,x}^\infty([0, T])} \leq \bigl\|\philinear[\Bmod](t)\bigr\|_{\covHolder{\Bmod_t(0), T}^{-\kappa}} \leq \frac{1}{2} T^{-\lambda_1} \;.
\end{equ}
Finally, by the estimate \eqref{local:eq-main-R-estimate} from Proposition~\ref{local:prop-main} combined with our heat kernel bounds \eqref{eq:zeroth-order-pB-ptwise-bd}, we have that
\begin{equ}
\bigl\|s^{\frac{\kappa}{2}} e^{s \Delta_{\Bmod_t(0)}} \phirem(t) \bigr\|_{L_{s,x}^\infty([0, T])} \leq \bigl\|\phirem(t)\bigr\|_{L_x^\infty} \ll 1 \;. 
\end{equ}
Combining, we thus obtain (here we use $\lambda_1 \ll \fa$ to absorb the contributions of $\philinear[\Bmod](t)$ and $\phirem(t)$ into $T^{-2\fa - \lossdel}$)
\begin{equs}
\bigl\|s^{\frac{\kappa}{2}} \hat{\higgs}(s)\bigr\|_{L_{s,x}^\infty([0, T])} &\lesssim T^{-\lossdel} \bigl(T^{\frac{\kappa}{2} - \frac{1}{r} - \lossdel} \bigl\|\higgs_0^+\bigr\|_{L_x^r} + T^{-2\fa - \lossdel}\bigr) \\
&\lesssim \bigl\|\higgs_0^+ \bigr\|_{L_x^r} + T^{-10\fa} \lesssim T^{-\frac{1}{\gamma_0}}\;,
\end{equs}
where we used that $\frac{1}{r}, \lossdel \ll \kappa \ll \fa$ in the second estimate, and the fact that $T \leq \bigl\|\higgs_0^+\bigr\|_{L_x^r}^{-\gamma_0}$ (recall Definition~\ref{defn:tau}) as well as that $\fa \ll 1$ in the final estimate. To finish, we use the fact that $\hat{\higgs}$ solves a covariant heat equation, combined with our heat kernel bounds \eqref{eq:zeroth-order-pB-ptwise-bd}, to estimate
\begin{equ}
\bigl\|s^{\frac{\kappa}{2}} \hat{\higgs}(s)\bigr\|_{L_{s,x}^\infty([0, 1])} \leq T^{-\frac{\kappa}{2}} \bigl\|s^{\frac{\kappa}{2}} \hat{\higgs}(s)\bigr\|_{L_{s,x}^\infty([0, T])} \lesssim T^{-\frac{\kappa}{2} - \frac{1}{\gamma_0}}\;. 
\end{equ}
Since $\gamma_0 > 2$, this completes the proof of the estimate \eqref{eq:inv-norm-bd-intermediate}, and thus as previously discussed, we obtain the desired result.
\end{proof}

Lemma~\ref{lemma:inv-norm-estimate-time-step} reduces the problem of estimating the $\invnorm{\place}{-\kappa}$ norm to estimating inverse powers of our iteration time steps $\taumod_m$. From the definition of $\taumod_m$ (Definition~\ref{defn:tau_L}), we expect that (at least on typical events) inverse powers of $\taumod_m$ are controlled by powers of the covariant energy, which is precisely the quantity that we estimated in the previous section. The following result provides such an estimate. Ultimately, this will be the reason why Proposition~\ref{prop:det-time-energy-bound} follows from our estimates of the covariant energy at the random times $\tmod_m$. 

\begin{lemma}\label{lemma:taumod-uniform-exp-bound}
Consider the iteration scheme from Definition~\ref{defn:time_steps_mod}. For any $p \geq 1$, we have that
\begin{equ}
\sup_{m \geq 1} \E[(\taumod_{m})^{-p} ~|~ \mc{F}_0] \lesssim \sup_{m \geq 0} \E[ \covariantenergy(\Asmoothpartmod(\tmod_{m}), \phismoothpartmod(\tmod_m))^p ~|~ \mc{F}_0] + \Lmod_0^p + \stepupperbound^{-p} \;,
\end{equ}
where the implicit constant may depend on $p$.
\end{lemma}
\begin{proof}
Let $m \geq 0$. By Lemmas~\ref{lemma:time-step-cond-exp-bd} and~\ref{lemma:taumod-star-det-bound}, we have that
\begin{equs}
\E[ (\taumod_{m+1})^{-p} ~|~ \mc{F}_0] \lesssim \stepupperbound^{-p} + \E[ \covariantenergy(\Asmoothpartmod(\tmod_{m}), \phismoothpartmod(\tmod_m))^p ~|~ \mc{F}_0] + \E [\Lmod_m^p ~|~ \mc{F}_0 ] \; .
\end{equs}
The first two terms are already acceptable. For the last term, note that if $m = 0$, then this may be bounded by $\Lmod_0^{p}$. For $m \geq 1$, we apply \eqref{e:Lmod_expectation_bound} to obtain an acceptable contribution.
\end{proof}

Before proving Proposition~\ref{prop:det-time-energy-bound}, we first establish an intermediate result which will cover the case where the deterministic time $t \lesssim \stepupperbound$. This will be needed to establish the full proposition, because our plan will be to wait until the iteration time $\tmod_m$ gets within (say) $3\stepupperbound$ of our fixed deterministic time $t$. Upon restarting at this stopping time, we may then apply the intermediate result, because $t - \tmod_m \lesssim \stepupperbound$ by construction.

\begin{lemma}\label{lemma:short-time-det-uniform-bound}
Consider the iteration scheme from Definition~\ref{defn:time_steps_mod}. For any $p \geq 1$, we have that
\begin{equ}
\E\Bigl[\max_{\substack{m \colon \tmod_m \leq 10\stepupperbound \\ m \geq 1}} \taumod_m^{-p} ~\Big|~ \mc{F}_0 \Bigr] \lesssim \sup_{n \geq 0} \E[ \covariantenergy(\Asmoothpartmod(\tmod_{n}), \phismoothpartmod(\tmod_n))^{4p} ~|~ \mc{F}_0] + \Lmod_0^{4p} + \stepupperbound^{-2p}\;, 
\end{equ}
where the implicit constant may depend on $p$.
\end{lemma}

Before proving Lemma~\ref{lemma:short-time-det-uniform-bound}, we first show that as $m$ gets large, it is very unlikely that $\tmod_m \leq 10 \stepupperbound$. This will then allow us to crudely bound the $\max$ appearing in the lemma by a sum.

\begin{lemma}\label{lemma:inverse-time-tail-bound}
Consider the iteration scheme from Definition~\ref{defn:time_steps_mod}. For any $T > 0$, $p \geq 1$, and $m \geq 1$, we have that
\begin{equ}
\P(\tmod_m \leq T ~|~ \mc{F}_0) \lesssim T^p m^{-p} \Bigl(\sup_{n \geq 0} \E[ \covariantenergy(\Asmoothpartmod(\tmod_{n}), \phismoothpartmod(\tmod_n))^p ~|~ \mc{F}_0] + \Lmod_0^p + \stepupperbound^{-p} \Bigr) \;,
\end{equ}
where the implicit constant may depend on $p$.
\end{lemma}
\begin{proof}
For notational simplicity, let $Q$ denote the quantity
\begin{equ}
Q = \sup_{n \geq 0} \E[ \covariantenergy(\Asmoothpartmod(\tmod_{n}), \phismoothpartmod(\tmod_n))^p ~|~ \mc{F}_0] + \Lmod_0^p + \stepupperbound^{-p} \;.
\end{equ}
This quantity arises from an application of Lemma~\ref{lemma:taumod-uniform-exp-bound}, which implies that
\[ \sup_{n \geq 1} \E[(\taumod_n)^{-p} ~|~ \mc{F}_0] \lesssim Q \;. \]
Since $\tmod_m = \taumod_1 + \cdots + \taumod_m$, we have by Cauchy--Schwarz that
\[ m^2 \leq \tmod_m \sum_{i=1}^m \taumod_i^{-1}, \]
and thus we obtain the estimate
\begin{equ}
\bigl\|\tmod_m^{-1} \bigr\|_{L^p(\place | \mc{F}_0)} \leq m^{-1} \sup_{n \geq 1} \bigl\|\taumod_n^{-1} \bigr\|_{L^p(\place |\mc{F}_0)} \lesssim m^{-1} Q^{1/p} \; .
\end{equ}
We now finish by Markov's inequality.
\end{proof}

\begin{proof}[Proof of Lemma~\ref{lemma:short-time-det-uniform-bound}]
For notational simplicity, let
\[ Q = \sup_{n \geq 0} \E[ \covariantenergy(\Asmoothpartmod(\tmod_{n}), \phismoothpartmod(\tmod_n))^{4p} ~|~ \mc{F}_0] + \Lmod_0^{4p} + \stepupperbound^{-4p} \;.\]
We may estimate
\begin{equ}
\max_{\substack{m \colon \tmod_m \leq 10 \stepupperbound \\ m \geq 1}} \taumod_m^{-1} \leq \sum_{m=1}^\infty \taumod_m^{-1} \one_{\tmod_m \leq 10 \stepupperbound} \;.
\end{equ}
Then by triangle, Hölder, and Lemmas~\ref{lemma:taumod-uniform-exp-bound} and~\ref{lemma:inverse-time-tail-bound} (the former applied with $2p$ in place of $p$, and the latter applied with $4p$ in place of $p$), we have that
\begin{equs}
\Bigl\|\max_{\substack{m \colon \tmod_m \leq 10 \stepupperbound \\ m \geq 1}}  \taumod_m^{-1} \Bigr\|_{L^p(\place | \mc{F}_0)} &\leq \sum_{m=1}^\infty \Bigl\| \taumod_m^{-1} \one_{\tmod_m \leq 10\stepupperbound}\Bigr\|_{L^p(\place | \mc{F}_0)} \\
&\leq \sup_{m \geq 1} \Bigl\| \taumod_m^{-1}\Bigr\|_{L^{2p}(\place | \mc{F}_0)}  \sum_{m=1}^\infty \P(\tmod_m \leq 10 \stepupperbound ~|~ \mc{F}_0)^{\frac{1}{2p}} \\
&\lesssim Q^{\frac{1}{2p}} \sum_{m=1}^\infty \Bigl( (10\stepupperbound)^{4p} m^{-4p} Q\Bigr)^{\frac{1}{2p}} \\
&\lesssim \stepupperbound^2  Q^{\frac{1}{p}} \;,
\end{equs}
which implies the desired inequality. 
\end{proof}

Finally, we use the previous results to prove Proposition~\ref{prop:det-time-energy-bound}.

\begin{proof}[Proof of Proposition~\ref{prop:det-time-energy-bound}]
Let $\stepupperbound$ be small enough so that both Proposition~\ref{prop:modified-iteration-uniform-bounds} and Lemma~\ref{lemma:inv-norm-estimate-time-step} apply. For notational simplicity, for $m \geq 0$, let
\begin{equ}
E_m = \covariantenergy(\Asmoothpartmod(\tmod_m), \phismoothpartmod(\tmod_m)) \;.
\end{equ}
It suffices to assume $t > 0$. Suppose first that $t > 3 \stepupperbound$. Let $M = \inf\{m \colon \tmod_m \geq t - 3\stepupperbound\}$. Note that in the notation of Proposition~\ref{prop:modified-iteration-uniform-bounds}, $M = M_{t - 3 \stepupperbound}$. By the same proposition, we have that $\tmod_m \uparrow \infty$ a.s., and thus $M < \infty$ a.s. Since the time steps are of size at most $\stepupperbound$, we have that $\tmod_M \leq t - 2\stepupperbound$, and thus $\tmod_{M+1} \leq t - \stepupperbound < t$. By the strong Markov property, we may restart at time $\tmod_{M+1}$, and then we have by Lemmas~\ref{lemma:inv-norm-estimate-time-step} and~\ref{lemma:short-time-det-uniform-bound} that
\begin{equs}
\E\Bigl[\invnorm{(A(t),\higgs(t))}{-\kappa}^p ~\Big|~ \mc{F}_{\tmod_{M+1}} \Bigr] &\lesssim \E\Bigl[ \sup_{t' \in (\tmod_{M+1}, \tmod_{M+1} + 5\stepupperbound]} \invnorm{(A(t'),\higgs(t'))}{-\kappa}^p ~\Big|~ \mc{F}_{\tmod_{M+1}}\Bigr] \\
&\lesssim \E\Bigl[ \max_{\substack{m \colon \tmod_m \leq \tmod_{M+1} + 10\stepupperbound \\ m \geq M+2}} \taumod_{m}^{-4p} ~\Big|~ \mc{F}_{\tmod_{M+1}} \Bigr] \\
&\lesssim \sup_{n \geq 0} \E\Bigl[E_{M + 1 +n}^{16p} ~\Big|~ \mc{F}_{\tmod_{M+1}} \Bigr] + \Lmod_{M + 1}^{16p} + \stepupperbound^{-8p} \;.
\end{equs}
By the strong Markov property and Proposition~\ref{prop:modified-iteration-uniform-bounds}, we have that the above is further bounded (up to constants) by
\begin{equ}
E_{M+1}^{64 p} + \Lmod_{M+1}^{64 p} + \stepupperbound^{-8p} \;.
\end{equ}
Taking expectations with respect to $\mc{F}_{0}$, we therefore obtain that
\begin{equs}
\E\Bigl[\invnorm{(A(t),\higgs(t))}{-\kappa}^p ~\Big|~ \mc{F}_{0} \Bigr] &\lesssim C(A_0, \higgs_0, \stepupperbound, p) + 1\;,
\end{equs}
where we applied Proposition~\ref{prop:modified-iteration-uniform-bounds} to estimate the term involving $E_{M+1}$, while for the term involving $\Lmod_{M+1}$, we applied the tower property to introduce an intermediate conditional expectation with respect to $\mc{F}_{\tmod_M}$, and then applied the strong Markov property and \eqref{e:Lmod_expectation_bound} to reduce to an estimate for $E_{M+1}$, which again follows by Proposition~\ref{prop:modified-iteration-uniform-bounds}. The first desired result now follows, at least when $t > 3 \stepupperbound$. If $t \leq 3 \stepupperbound$, we may directly apply Lemmas~\ref{lemma:inv-norm-estimate-time-step} and~\ref{lemma:short-time-det-uniform-bound} starting from time $\tmod_0 = 0$ (combined with Proposition~\ref{prop:modified-iteration-uniform-bounds}) to obtain the desired estimate. For the second desired result, we apply the final estimate in Proposition~\ref{prop:modified-iteration-uniform-bounds} to estimate the limsup as $t \toinf$ of conditional expectations of $E_{M+1}^{64 p}$ and $\Lmod_{M+1}^{64p}$ given $\mc{F}_0$ (note as before, we may reduce the latter estimate to the former) by a constant which does not depend on the initial data $A_0, \higgs_0$.
\end{proof}

\section{Construction and irreducibility of the Markov process}
\label{sec:regStruc}

In this section, we recall the construction of the Yang--Mills--Higgs process in the framework of regularity
structures as in \cite{CCHS22,CCHS22+}. The construction in \cite{CCHS22} yields a very `nice' state space
which, since it is Polish, is in particular suitable for applying the general ergodic theory of Markov processes.
It does not include the Higgs field, unlike the work \cite{CCHS22+} which however uses a much `uglier' state space
due to the fact that it considers the 3D case. Here, we would like to have the best of both worlds. The
modifications to \cite{CCHS22} required to accommodate the Higgs field turn out to be quite minor. We therefore mainly focus
on the differences, on setting up the notations required for the next section, and on showing that the resulting
Markov process has good ergodic properties.
 
For this, we first introduce in Section~\ref{sec:stateSpace} a suitable state space $\fO_{\kappa}^\YMHs$
which essentially consists of equivalence classes of pairs $(A,\higgs)$ related by sufficiently regular
gauge transformations. This is a minor adaptation of the space constructed in \cite{CCHS22} to take the Higgs
field into account. In Section~\ref{sec:Markov} we then apply the general theory of regularity structures
to construct a Markov process with values in $\fO_{\kappa}^\YMHs$
(but possibly finite explosion time). 
We then show in Section~\ref{sec:irreducible} that this Markov process is topologically 
irreducible. 
Finally, we will show in the next section that, provided that we make the correct choice of renormalisation
constants, the process previously constructed in Theorem~\ref{thm:global_soft}
(for which the Markov property is not a priori obvious at all) does in fact coincide with the process
constructed here in Section~\ref{sec:Markov}. 
Together with the strong Feller property, this then provides the final ingredients required 
to complete the proof of our main theorem. 

\subsection{State space}\label{sec:stateSpace}

We start by constructing a natural space of gauge orbits for the 2D Yang--Mills--Higgs process
and we show that in particular this space admits continuous spatially localised smoothened Wilson line 
observables. Recall that $\CG_\alpha$ denotes the space of little-$\alpha$-Hölder continuous\footnote{A function is said to be little-$\alpha$-Hölder if it
is a limit of smooth functions in the $\alpha$-Hölder metric.}
gauge transformations $g \colon \T^2 \to \liegroup$ and consider the space $\Omega^1_\alpha$ of
distributional one-forms defined in \cite[Def.~3.22\,\&\,Eq.~3.2]{CCHS22}. Given 
$\kappa \in (0,\f13)$, we then define the space
\begin{equ}[e:basicState]
\Omega_{\kappa}^\YMHs = \Omega^1_{1-\kappa} \times \cC^{-\kappa}_x(\T^2,\V)\;,
\end{equ}
equipped with its natural norm $\$(A,\higgs)\$_{\kappa} = |A|_{1-\kappa} + \|\higgs\|_{\cC^{-\kappa}_x}$.
Recall that smooth connection one-forms are dense in $\Omega^1_{1-\kappa}$ and that the action of 
the group $\CG$ of smooth gauge transformations extends to a continuous action of 
$\CG_{1-\kappa}$ on $\Omega^1_{1-\kappa}$, see \cite[Thm~3.27]{CCHS22}.

\begin{remark}\label{rem:embedding}
We will also make use of the fact that $\cC^{\frac{1-\kappa}2}_x$ embeds continuously into $\Omega_{\kappa}^\YMHs$,
see \cite[Rem.~3.23]{CCHS22}.
\end{remark}

Since the pointwise product extends continuously to $\cC^{1-\kappa}_x \times \cC^{-\kappa}_x \to \cC^{-\kappa}_x$ 
provided that $\kappa < \f12$, it follows that $\CG_{1-\kappa}$ acts continuously on $\Omega_{\kappa}^\YMHs$
by
\begin{equ}
g \act (A,\higgs) \eqdef (g\act A, \rep(g)\higgs)\;.
\end{equ}
We then define the space of gauge orbits 
\begin{equation}\label{e:gauge-orbits}
\fO_{\kappa}^\YMHs \eqdef \Omega_{\kappa}^\YMHs / \CG_{1-\kappa}
\end{equation}
endowed with the quotient topology. As in \cite[Thm~3.45]{CCHS22}, we have the following fundamental result.

\begin{theorem}\label{theo:Polish}
The space $\fO_{\kappa}^\YMHs$ is Polish.
\end{theorem}

\begin{proof}
The proof is virtually identical to that of \cite[Thm~3.45]{CCHS22}, so we only give a rough sketch. 
Based on the norm $\$\place\$_{\kappa}$, we construct 
the distance $k_{\kappa}$ on $\Omega_{\kappa}^\YMHs$ as in \cite[Def.~3.41]{CCHS22} and we 
write $D_{\kappa}$ for the corresponding 
Hausdorff distance between elements of $\fO_{\kappa}^\YMHs$ (with each gauge orbit viewed as a subset of $\Omega_{\kappa}^\YMHs$).
The main tool on which the proof of \cite[Thm~3.45]{CCHS22} relies is the fact that, given
$A \in \Omega^1_{1-\kappa}$, the gauge transformation $g$ can be
recovered completely (modulo possibly a global gauge) from the knowledge of $g \act A$
and one has a quantitative bound of the form $\|g\|_{\cC^{1-\kappa}_x} \lesssim 1 + |A|_{1-\kappa} + |g\act A|_{1-\kappa}$.
Since the map $(g,\higgs) \mapsto \rep(g)\higgs$ is Lipschitz continuous on bounded sets of 
$\CG_{1-\kappa} \times \cC^{-\kappa}_x$, we can again deduce from this that gauge orbits are closed in
$\Omega_{\kappa}^\YMHs$. The remainder of the proof then follows \textit{mutatis mutandis}.
\end{proof}

The following compactness result will be useful later on to show that 
the transition probabilities of our process are tight.

\begin{lemma}\label{lem:compactEmbedding}
Bounded
sets in $\fO_{\bar\kappa}^\YMHs$ are relatively 
compact in $\fO_{\kappa}^\YMHs$ whenever $\bar \kappa < \kappa$.
\end{lemma}

\begin{proof}
Since we know that $\CC^{-\bar \kappa}_x$ embeds compactly into $\CC^{-\kappa}_x$,
it suffices to show that the embedding of Banach spaces 
$\Omega^1_{1-\bar\kappa} \subset \Omega^1_{1-\kappa}$ is compact.

For this, recall that distributional connections as defined in \cite[Def.~3.1]{CCHS22} are continuous
functions on line segments on the torus, which we identify with elements $\ell = (x,v) \in \T^2 \times \R^2_*$ (with $\R^2_* = \R^2 \setminus \{0\}$),
where $x$ is the starting point of the line segment and $v$ its (oriented) length. As in \cite{CCHS22},
we write $|\ell| = |v|$. When quantifying over $\ell$ we will always restrict ourselves
to $|\ell| \le 1$ to make sure that our Hölder-type norms are finite for constant $A$'s.
The crucial observation is that if we define on $\T^2 \times \R^2_*$ the hyperbolic
Riemannian metric 
\begin{equ}
ds^2 = \frac{dx_1^2 + dx_2^2 + dv_1^2 + dv_2^2}{|v|^2}
\end{equ}
and write $d$ for the corresponding distance
then, for any two line segments $\ell,\bar\ell$ that are ``not far'' and ``form a vee''
in the sense
of \cite[Defs.~3.3\,\&\,3.8]{CCHS22}, one has $\Area(\ell,\bar \ell) \sim (|\ell| + |\bar \ell|)^2 d(\ell,\bar \ell)$
in the sense that the ratio of both quantities is uniformly bounded from above and below.
By \cite[Thm~3.11]{CCHS22}, this then shows that $\Omega_{1-\kappa}^1$ is a closed subset 
(because of the additivity constraint in \cite[Def.~3.1]{CCHS22}) of the space of little
$\f{1-\kappa}2$-Hölder functions on $\T^2 \times \R^2_*$ with weight $|\ell|^{1-\kappa}$ in the sense that 
\begin{equ}
\|A\| \sim \sup_{d(\ell, \bar \ell) \le 1} \frac{|A(\ell)-A(\bar\ell)|}{(|\ell|^{1-\kappa} + |\bar\ell|^{1-\kappa})d(\ell,\bar \ell)^{\f{1-\kappa}2}}
+ \sup_\ell \frac{|A(\ell)|}{|\ell|^{1-\kappa}}\;.
\end{equ}
The weight is ``nice'' in the sense that $|\ell|^{1-\kappa} \sim |\bar\ell|^{1-\kappa}$ uniformly over
pairs of segments with $d(\ell, \bar \ell) \le 1$. The compactness of the embedding
$\Omega^1_{1-\bar\kappa} \subset \Omega^1_{1-\kappa}$ is then a consequence of the compactness of the 
embeddings of weighted Hölder spaces when one varies the exponent and the weight.
\end{proof}

We observe that the state space $\fO_{\kappa}^\YMHs$
contains several ``nice'' observables. 
Recall first that for a smooth function $\gamma \colon [0,1] \to \T^2$
and a smooth one-form $A$, the integral of $A$ along $\gamma$ is given by
\begin{equ}[e:defgammaA]
\gamma_A(t) = \int_0^t A_{\gamma(s)}(\dot \gamma(s))\,ds\;.
\end{equ}
We then define the holonomy $\hol(A,\gamma)$ as the solution at time $1$ to the $\liegroup$-valued
differential equation
\begin{equ}[e:ODEWilson]
dg(t) \,g^{-1}(t) = -d\gamma_A(t)\;,\qquad g(0) =\groupid\;.
\end{equ}
Using Young integration \cite{Young}, the solution to \eqref{e:ODEWilson} extends continuously to all 
Hölder continuous maps $\gamma_A \colon [0,1] \to \fg$ with Hölder exponent strictly greater than
$\f12$. Furthermore, we know by \cite[Cor.~3.19]{CCHS22} that \eqref{e:defgammaA} extends to all $A \in \Omega^1_{1-\kappa}$
with the resulting $\gamma_A$ being $(1-\kappa)$-Hölder continuous, with norm bounded by some multiple of 
$|A|_{1-\kappa}$ that only depends on the $\cC^2$ norm of $\gamma$.

Writing $[g]$ for the conjugacy class of $g \in \liegroup$, one consequence of 
this (see \cite[Thm~2.1]{CCHS22}) is that Wilson loop
observables $(A,\higgs) \mapsto [\hol(A,\gamma)]$ for smooth closed loops $\gamma$
can be realised as continuous maps on it (since the map $\fO_{\kappa}^\YMHs \to \fO_{1-\kappa}$
obtained by forgetting the Higgs field is well posed and continuous), but it turns out that it also 
contains a smoothened out version of Wilson lines. 
Recall that if $\gamma_{x,y}\colon [0,1] \to \T^2$ is a smooth map such that $\gamma_{x,y}(0) = x$ and
$\gamma_{x,y}(1) = y$, then the quantity given for smooth pairs $(A,\higgs)$ by
\begin{equ}[e:WilsonLine]
(A,\higgs) \mapsto \scal{\rep(\hol(A,\gamma_{x,y}))\higgs(x),\higgs(y)}_\V\;,
\end{equ}
is gauge invariant as an immediate consequence of the fact that 
$\hol(g\act A,\gamma_{x,y}) = g(y)\hol(A,\gamma_{x,y})g(x)^{-1}$. Of course, if we fix $x$ and $y$, then there is no
hope that the map \eqref{e:WilsonLine} extends continuously to $\fO_{\kappa}^\YMHs$ since
we cannot evaluate the Higgs field $\higgs$ pointwise. The next best thing one can hope for is that, if we allow $x$ and $y$ to vary,
\eqref{e:WilsonLine} extends to a map $\fO_{\kappa}^\YMHs \to \cC^{-2\kappa}_{x,y}(\T^2 \times \T^2)$,
provided that the collection $(x,y) \mapsto \gamma_{x,y}$ is sufficiently regular.

\begin{remark}\label{rem:jumps}
Since the fundamental group of $\T^2$ is non-trivial, it is of course impossible to find
a continuous map $(x,y) \mapsto \gamma_{x,y}$ with values in smooth functions such that 
$\gamma_{x,y}(0) = x$ and $\gamma_{x,y}(1) = y$ for every $x,y \in \T^2$. It is however possible to guarantee that 
discontinuities take place along finitely many smooth hypersurfaces of codimension $1$.
Alternatively, one can lift the curves to the universal cover and hope to get Wilson
line observables $\fO_{\kappa}^\YMHs \to \cC^{-2\kappa}_{x,y}(\R^2 \times \R^2)$, which is the 
approach that we will follow.
\end{remark}

We have the following regularity result of the holonomies 
as a function of the path $\gamma$ along which they are taken.

\begin{lemma}\label{lem:contHol}
Given $\alpha \in (\frac23,1]$, $\eta \in (0,\frac{(3\alpha-2)(2\alpha-1)}{8\alpha})$, and $A \in \Omega^1_\alpha$, the map
$\gamma \mapsto \hol(A,\gamma)$ is $\eta$-Hölder continuous in the supremum norm, uniformly over balls
in $\cC^2$.
\end{lemma}

\begin{proof}
As in the proof of \cite[Prop.~3.21]{CCHS22}, we first use \cite[Thm~3.18]{CCHS22} to obtain
\begin{equ}
\|\gamma_A - \bar \gamma_A\|_\infty \lesssim \sup_{t \in [0,1]}|A(\ell_\gamma^t)- A(\ell_{\bar \gamma}^t)|
+ |A|_\alpha |\gamma;\bar \gamma|_{\beta}^{\alpha/\beta}\;,
\end{equ}
for $\beta \in [\frac23,\alpha)$, where $\ell_\gamma^t$ denotes the line segment joining $\gamma(0)$ to $\gamma(t)$.

By \cite[Thm~3.11]{CCHS22}, it is straightforward to show that 
one has the bound
\begin{equ}
|A(\ell_\gamma^t)- A(\ell_{\bar \gamma}^t)| \lesssim |A|_\alpha \|\gamma-\bar \gamma\|_\infty^{\alpha/2}\;,
\end{equ}
while \cite[Prop.~3.20]{CCHS22} yields the bound
$|\gamma;\bar \gamma|_{\beta} \lesssim \|\gamma-\bar \gamma\|_\infty^{\frac{3\beta-2}{4}}$,
so that 
\begin{equ}[e:boundLinfty]
\|\gamma_A - \bar \gamma_A\|_\infty 
\lesssim |A|_\alpha \|\gamma-\bar \gamma\|_\infty^{\frac{3\beta-2}{4\beta}\alpha}\;.
\end{equ}
Furthermore, as in the proof of \cite[Cor.~3.19]{CCHS22}, we find that for smooth $\gamma$ one has
\begin{equ}[e:apriorialpha]
\|\gamma_A\|_\alpha \lesssim |A|_\alpha\;.
\end{equ}
As in the proof of \cite[Prop.~3.21]{CCHS22}, we then conclude from interpolating between \eqref{e:boundLinfty}
and \eqref{e:apriorialpha} that, for $\bar \alpha < \alpha$,
\begin{equ}
\|\gamma_A - \bar \gamma_A\|_{\bar \alpha} \lesssim 
|A|_\alpha \|\gamma-\bar \gamma\|_\infty^{\frac{3\beta-2}{4\beta}(\alpha-\bar\alpha)}\;.
\end{equ}
Since the solution to \eqref{e:ODEWilson} depends in a locally 
Lipschitz continuous way on $\gamma_A$
in any $\bar \alpha$-Hölder norm as long as $\bar \alpha > \frac12$
(see for example \cite{TerryYoung} or \cite[Thm~3.14]{GaleatiYoung}), we conclude that 
\begin{equ}
d \bigl(\hol(A,\gamma), \hol(A,\bar \gamma)\bigr)
\lesssim |A|_\alpha \|\gamma-\bar \gamma\|_\infty^{\frac{3\beta-2}{4\beta}(\alpha-\bar\alpha)}\;,
\end{equ}
where $d$ is the bi-invariant distance on $\liegroup$ induced by our given
bi-invariant scalar product on $\fg$.  
Choosing $\beta$ close to $\alpha$ and $\bar \alpha$ close to $\frac12$, the claim follows.
\end{proof}

\begin{remark}
The Hölder exponent $\eta$ appearing in this result is probably suboptimal, but it is sufficient
for our purpose. When $\alpha$ is close to $1$, which is the relevant case for us since
we always take $\alpha = 1-\kappa$, 
then one can make $\eta$ close to $\frac18$.
\end{remark}

Similarly, it is straightforward to show that, for every $\gamma \in \cC^2$,
$\hol(\place,\gamma)$ extends to a continuous (in fact even locally Lipschitz continuous) function from $\Omega^1_\alpha$ to $\liegroup$.
As a consequence, we can show that smoothened out ``Wilson lines'' are well-defined observables
on the space $\Omega_{\kappa}^\YMHs$ provided that $\kappa$ is small enough.

\begin{corollary}\label{cor:WilsonLines}
There exists $\kappa_0 > 0$ such that, for all $\kappa \in (0,\kappa_0]$ the following
is true. Let $\gamma$ be a locally Lipschitz continuous map
\begin{equ}
 \R^2 \times \R^2 \ni (x,y) \mapsto \gamma_{x,y} \in \cC^2([0,1],\R^2)\;,
\end{equ}
such that $\gamma_{x,y}(0) = x$, $\gamma_{x,y}(1) = y$ for every $x,y \in \R^2$. Then the map
\begin{equ}
W_{\gamma} \colon (A,\higgs) \mapsto \scal{\rep(\hol(A,\gamma_{x,y})) \higgs(x),\higgs(y)}_\V\;,
\end{equ}
is continuous as a map $\Omega_{\kappa}^\YMHs \to \cC^{-2\kappa}_{x,y}$. 
It is also invariant under the action
of $\CG_{1-\kappa}$, yielding a continuous map from $\fO_{\kappa}^\YMHs$ to $\cC^{-2\kappa}_{x,y}$.
\end{corollary}

\begin{proof}
It suffices to recall that the map $\higgs \mapsto \higgs \otimes \higgs$ is continuous
from $\cC^{-\kappa}_x(\R^2)$ to $\cC_{x,y}^{-2\kappa}(\R^4)$ and that, if  $\kappa$ is small enough, one has $\kappa < \frac{(3\alpha-2)(2\alpha-1)}{16\alpha}$ with
$\alpha = 1-\kappa$, so that Lemma~\ref{lem:contHol} guarantees
 that the map $(x,y) \mapsto \hol(A,\gamma_{x,y})$ is Hölder continuous
of exponent strictly greater than $2\kappa$.
\end{proof}

\subsection{Construction of the Markov process}\label{sec:Markov}

We now construct a Markov process that is a natural candidate for ``the'' solution to the stochastic YMH
equation in the DeTurck gauge, namely \eqref{intro:eq-SYMH} with $B\equiv 0$ (see also Remark~\ref{rem:DeTurck}
and \cite[Eq.~1.8]{CCHS22+}).
Setting $\Y = \fg^2 \times \V$, we can write this equation schematically as
\begin{equ}[e:schematic]
\d_t X = (\Delta-1) X + \linTerm X + \hat \bil(X,DX) + \tril(X,X,X) + \boldxi\;,
\end{equ}
where $\linTerm \colon \Y \to \Y$, $\hat \bil \colon \Y \times \Y^2 \to \Y$, and $\tril\colon \Y^3 \to \Y$ 
are suitable linear, bilinear and 
trilinear maps (we only construct local solutions in this section, so the specifics of these maps are 
irrelevant since we will use neither geometric information, nor the signs of the terms), see
for example \cite[Eq.~1.12]{CCHS22+}. We've explicitly subtracted the constant $1$ from the Laplacian (which is then 
undone by the definition of $\linTerm$) in order to provide a large-scale cutoff which enters the
definition of the BPHZ renormalisation scheme. We will see later on that even though the BPHZ scheme \textit{does}
depend on that choice of constant, the correct renormalisation which guarantees that our process
coincides with the process considered in the previous sections does not.

We can then fit \eqref{e:schematic} into the framework developed in \cite{BCCH}, in particular the
version for vector-valued SPDEs developed in \cite[Sec.~5]{CCHS22}. For the moment, we 
consider \eqref{e:schematic}
as a one-component vector-valued equation with one driving noise, so we set
\begin{equ}
\fL_+ = \{\fx\}\;,\quad \fL_- = \{\bullet\}\;,
\end{equ}
where $\fx$ represents the solution $X$ and $\bullet$ represents the driving noise $\xi$.

The corresponding target spaces $W_\ft$ and kernel spaces $\CK_\ft$ as in \cite[Sec.~5.6]{CCHS22}
are given by
\begin{equs}[0][e:spaces]
W_{\fx} = W_{\bullet} = \Y\;,\quad  \CK_{\fx} = \R\;.
\end{equs}
We fix some $\kappa > 0$ small enough and we postulate that the integration operators associated to 
$\fx$ have degree $2$, while the
noise $\bullet$ has degree $-2-\kappa$. 
We then take as our ``rule'' \cite[Sec.~5]{BHZ19} the minimal rule allowing to 
formulate the right-hand side of \eqref{e:schematic}, analogously to \cite[Eq.~6.10]{CCHS22+}.

Recall that the construction given in \cite{BHZ19,CCHS22} (which is really just the general procedure
outlined in \cite[Sec.~8.1]{H14} translated to the vector-valued setting)
yields a collection of symbols usually denoted either
by formal expressions obtained from the noise symbol $\Xi_\bullet$
(which we henceforth simply denote by $\Xi$) and the ``Taylor monomials''
$\mathbf{X}^k$ by applications of products and ``integration operators'' $\CI = \CI_\fx$, or
as the corresponding syntactic trees. Each symbol $\tau$ has a natural degree $\deg\tau$ and
is associated to a vector space $\CT[\tau]$ obtained by associating a tensor factor $W_\bullet^* = \Y^*$ to
each copy of $\Xi$ in $\tau$. See \cite[Sec.~5]{CCHS22} for details on how one keeps track 
of which copy is which in a systematic way, how symmetries are being taken into account, and how
to define the canonical lift $\LL$ of a smooth noise in this setting, provided that one associates an
integration kernel to $\CI$. Since we will never need elements of degree greater than $2$, 
we write $\CT = \bigoplus_{\deg \tau \le 2} \CT[\tau]$.\label{text:CT-appearence}

\begin{notation}\label{notationTensor}
Given a symbol $\tau$, assuming that we have a canonical isomorphism $\CT[\tau]\simeq W$ for
some Banach space $W$, we will typically denote elements of $\CT[\tau]$ by $x \otimes \tau$ with
$x\in W$, see the discussion in \cite[Sec.~5.1]{CCHS22}. 
This is purely notational and allows us to distinguish for example elements of
$\CT[\<Xi>]$ from elements of $\CT[\<IXi>]$ even though both spaces are canonically isomorphic to $\Y^*$.
(Note that the convention used in \cite{CCHS22} was rather $\tau \otimes x$, but the convention used here
allows us to naturally denote elements of $\Y \otimes \CT[\tau]$ by  $x \otimes \tau$ with $x \in \Y \otimes W$.)
\end{notation}

\begin{remark}\label{rem:NotationNoise}
The reason why it is the dual spaces that are naturally attached to symbols
is the following. Take for example $\ell \in W_\bullet^*$ and write $\ell \otimes\Xi$ for the corresponding
element in $\CT[\Xi]$ of our regularity structure, then the canonical lift $(\Pi,\Gamma) = \LL(\xi)$ 
of a smooth noise 
$\xi$ is such
that $\bigl(\Pi_x (\ell\otimes\Xi)\bigr)(z) = \scal{\ell, \xi(z)}$, where $\scal{\place,\place}$ denotes the duality pairing
between $W_\bullet^*$ and $W_\bullet$. In other words, in order to describe a $W$-valued noise, the corresponding
subspace of the regularity structure should be identified with $W^*$ rather than $W$.
Accordingly, we will simply write $\<IXi>$ for the element of $\Y \otimes \CT[\<IXi>]$
corresponding to the identity $\id \in \End(\Y) \simeq \Y \otimes \Y^* \simeq\Y \otimes \CT[\<IXi>]$.
\end{remark}

For the remainder of this section, we associate to $\CI$ the integration kernel given by 
space-time convolution with the massive heat kernel
\begin{equ}[e:massiveHeat]
\hat \fp_t(x) = e^{-t}\fp_t^\Euc(x)\;.
\end{equ}
Note that this is a kernel on $\R \times \R^2$ and we will always
consider models on the whole space. However, these models will always be obtained by lifting 
spatially periodic noise functions.

We write $\models_0$\label{text:models} for the 
closure, in the space of all models over $\R\times \R^2$ 
for the regularity structure that we just constructed,
of all the models obtained by taking the canonical lift of a smooth function that is periodic in space 
and compactly supported in time, and acting on it with an
element of the renormalisation group $\RR$\label{text:renormalization-group}. 

\begin{remark}
Since we only consider models that are periodic in space, \cite[Eqns~2.16--2.17]{H14} yields a way of assigning
sizes $\$Z\$_I$ and distances $\$Z;\bar Z\$_I$ to models over any compact time interval $I\subset \R$.
We then constrain the behaviour at infinity by setting for example
$\$Z\$ = \sup_{n\in \Z} (1+|n|)^{-1} \$Z\$_{[n,n+1]}$ and similarly for the distance.
Since we only consider solutions over $[0,1]$, the precise weight we use doesn't matter, but it is
convenient for the model to be defined on the whole space and we want the behaviour at infinity
to be loose enough to allow for $\<Xi>$ to represent space-time white noise.
\end{remark}

The group $\RR$ is abelian,
so just a finite-dimensional vector space, but acts on $\models_0$ in a somewhat non-trivial way, see
\cite{BHZ19,CCHS22} for details. 
We have a continuous map
\begin{equ}
\iota \colon \models_0 \to \CC^{-\kappa}(\R\times \T^2, \Y) \;,\qquad Z = (\Pi,\Gamma) \mapsto \Pi_0 \<IXi>\;.
\end{equ}
We will henceforth also use the convention that if $X$ is a function / distribution with values in
$\Y$, then we write $X_A$ for its component in $\fg^2$ and $X_\higgs$ for its component
in $\V$.

We consider the space $\CX = C(\R,\Omega_{\kappa}^\YMHs)$,
where $\Omega_{\kappa}^\YMHs$ is as defined in \eqref{e:basicState}. 
Writing $\pi \colon \Omega_{1-\kappa}^1 \to \CC^{-\kappa}_x$ for the map given in
\cite[Prop.~3.21]{C19}, there is also a natural map
\begin{equ}
\iota \colon \CX  \to \CC^{-\kappa}(\R\times \T^2,\Y)
\end{equ}
given by 
\begin{equ}[e:defiotalollipop]
\bigl(\iota \<IXib>\bigr)(\psi) = \int \bigl(\<IXib>_\higgs(t)\bigr)(\psi_\higgs(t,\cdot)) + \bigl(\pi \<IXib>_A(t)\bigr)(\psi_A(t,\cdot)) \,dt\;.
\end{equ}
We then have a natural noise space $\CN = \{(Z,\<IXib>) \in \models_0 \times \CX\,:\, \iota \<IXib> = \iota Z\}$.\label{text:CN}
Since the maps $\iota$ are continuous, $\CN$ is a closed subset of the Polish space
$\models_0 \times \CX$
and is therefore Polish itself. The ``correct'' way to interpret this construction is that 
elements of $\CN$ are simply models $Z = (\Pi,\Gamma)$ such that $\<IXib> = \Pi_0 \<IXi>$ belongs
to $\CX$, rather than just $\CC^{-\kappa}$. Therefore, we will from now on denote generic elements of 
$\CN$ simply by $Z$. We also write $M_g \colon \CN \to \CN$ for the action of $g \in \RR$ onto $\CN$.

\begin{remark}\label{rem:randomModel}
By \cite{CH16} (but see also \cite{BHZ19,CCHS22,HS23}, 
in particular regarding the fact that the 
solution to the stochastic heat equation takes values in $\CX$),
we know that the BPHZ renormalisation of the canonical lift of the regularised
noise $\boldxi_\eps$ converges in probability to a non-trivial
random variable $\CZ$ with values in $\CN$.
\end{remark}

Furthermore, we have a good local solution theory for equations driven
by \textit{any} model $Z \in \CN$ which we describe now. This is by now well understood
\cite{H14,BCCH}, so we mainly use this description to fix notations and to emphasise details that are
particular to our setup.
Given $Z \in \CN$, recall that it yields for each
sector $\CV \subset \CT$ a space
$\CD^{\gamma,\eta}$\label{text:CD} of $\CV$-valued modelled distributions (that depends on $Z$), 
see \cite[Def.~6.2]{H14}.
In fact, we want to consider modelled distributions whose reconstruction yields a $\Y$-valued function, 
so our spaces $\CD^{\gamma,\eta}$
actually consist of $(\Y \otimes \CV)$-valued functions. In order to emphasise this, we will sometimes
write $\CD^{\gamma,\eta}(\Y)$ or $\CD^{\gamma,\eta}(I,\Y)$ if we want to furthermore restrict
the time variable to some interval $I \subset \R$. (Regarding the space variable, we will only ever consider 
functions that are defined on all of $\R^2$ but are periodic.)

\begin{remark}
Since $\Y$ is finite-dimensional, there is
no ambiguity in the meaning of the tensor product and all definitions from \cite{H14} simply apply
componentwise.
\end{remark}

The definition given in \cite{H14} depends on a hyperplane $P$ at which the modelled distribution
is allowed to become singular. In our case, we will always take $P = \{(t,x)\,:\, t = s\}$ for some $s \in [0,1]$
and we write $\CD^{\gamma,\eta}_s$ if we want to specify the value of $s$. (If it is not specified, we take $s=0$.)

\begin{remark}
We now fix once and for all the exponents $\kappa, \gamma, \eta$ in such a way that
$\kappa < 1/10$ is sufficiently small so that all the results from the previous sections apply, $\gamma = 3/2$, and $\eta \in (-\frac14, -\kappa)$.
\end{remark}

In order to reformulate \eqref{e:schematic} as a fixed point problem that is stable under the convergence
$\eps \to 0$ when $\boldxi$ is replaced by a mollified noise, we write 
 $\hPicard_{s,t}$\label{text:hPicard} for the space of pairs $(X,Z)$ with $Z \in \CN$ and 
$X \in \CD^{\gamma,\eta}_{s,t} \subset \CD^{\gamma,\eta}_s([s,t],\Y)$. 
Since the space $\CD^{\gamma,\eta}_s$ depends on the model $Z$, 
the space $\hPicard_{s,t}$ isn't a product space, but rather a kind
of vector bundle $\pi \colon \hPicard_{s,t} \to \CN$ over $\CN$.
We say ``kind of'' because even though each fibre is a Banach space, the space 
$\CD^{\gamma,\eta}_s$ does not come equipped with local trivialisations.

In our case, the sector $\CV$ containing the solution itself is of regularity $-\kappa$
and the maps given by $X \mapsto \hat \bil(X,DX)$ and $X \mapsto \tril(X,X,X)$ 
as in \eqref{e:schematic}
can be lifted to maps $(\Y \oplus \Y^2) \otimes \CV \to \Y \otimes \CT$ 
in the usual way, see \cite[Sec.~5]{CCHS22} (here the copy of $\Y$ represents the value of $X$
while the copy of $\Y^2$ represents the value of $DX$). 
As in \cite[Sec.~5]{H14}, we furthermore have an integration operator satisfying, for $\bar \gamma > 0$
and $\bar \eta > -2$, 
$\CP \colon \CD^{\bar \gamma,\bar \eta}_{s,t} \to \CD^{\bar \gamma+2,\bar \eta + 2}_{s,t}$ with the property that 
$\CP F - \CI F$ takes values in $\Y \otimes \bar \CT$ (with $\bar \CT$ denoting the polynomial sector of
the regularity structure) and such that the reconstruction operator satisfies
$\CR \CP F = \hat \fp \star \CR F$, with $\star$ denoting space-time convolution. Here, elements of 
$\CD^{\bar \gamma,\bar \eta}_{s,t}$ are extended to all of $\R \times \R^2$ by setting them to $0$ outside
of the time interval $[s,t]$. 

In particular, this shows that the map $\CM_{s,t}$ given by
\begin{equ}\label{e:fixedpointmap}
\CM_{s,t}(X) = \CP \bigl((\linTerm + \extraC) X  + \hat \bil(X,DX) + \tril(X,X,X)\bigr)\;,
\end{equ}
is locally Lipschitz continuous from each fibre $\CD^{\gamma,\eta}_{s,t}$ to itself.
Here, we added ``by hand'' an additional counterterm $\extraC X$ which will later
be chosen in such a way that the solution to the corresponding fixed point problem
with the BPHZ model of Remark~\ref{rem:randomModel} yields a gauge-covariant 
process. In the splitting $\Y = \fg^2 \oplus \V$, $\extraC$ will always be of the form
$\extraC = \extraC_A \oplus \extraC_\Phi$, with both $\extraC_A$ and $\extraC_\Phi$
being $\liegroup$-invariant.

It also follows as in 
\cite[Thm~7.1]{H14} that, on bounded
sets of $\hPicard_{s,t}$, its Lipschitz constant is of order $|t-s|^{\beta}$ for some (possibly small)
exponent $\beta > 0$.
In particular, for every (large) $R>1$ and every (small) $\eps > 0$, there exists $\delta > 0$
such that, for all $0 \le s < t \le 1$ with $|t-s| \le \delta$, one has
\begin{equ}[e:contract]
\|\CM_{s,t}(X) - \CM_{s,t}(Y)\|_{\Picard_{s,t}} \le \eps \|X-Y\|_{\Picard_{s,t}}\;,\qquad \|\CM_{s,t}(X)\|_{\Picard_{s,t}} \le \eps\;,
\end{equ}
for all $Z \in \models_0$ with $\$Z\$ \le R$ and all $X,Y \in \Picard_{s,t}(Z)$
with $\|X\|_{\Picard_{s,t}} \vee \|Y\|_{\Picard_{s,t}} \le R$.

A final ingredient is that there exists a linear map $\CJ \colon \CC^1(I\times \R^2) \to \CC(I\times \R^2, \CT)$ obtained by taking the first order 
``Taylor jet'', namely
\begin{equ}[e:TaylorJet]
\CJ f (t,x) = f(t,x) \one + \nabla f(t,x) \X\;,
\end{equ}
where $\nabla$ denotes the spatial gradient.
If furthermore $f \in \CC^{\gamma,\eta}(\R_+\times \T^2,\Y)$ (which is the case for example for 
the solution $(t,x) \mapsto \bigl(\hat \fp_t X_0\bigr)(x)$ to the massive heat 
equation with 
periodic initial 
condition $X_0 \in \cC^{\eta}_x$) then, by definition, one has $\CJ f \in \Picard_{s,t}$
for any admissible model $Z$, provided that $\gamma \le 2$. 

With all these ingredients in place, given any time interval $[s,t]$ and 
initial condition $X_s \in \cC^{\eta}_x(\T^2,\Y)$, 
we can reformulate \eqref{e:schematic} as the fixed point problem in $\CD^{\gamma,\eta}_{s,t}$
\begin{equ}[e:FPgen]
\hat X = \CM_{s,t}(\hat X) + \<IXi> + \CJ \hat \fp_{\place-s} (X_s - \<IXib>_s)\;,
\end{equ}
where we wrote $\<IXib>_s \in \Omega_{\kappa}^\YMHs$ for the value of the $\<IXib>$-component of
$Z$ at time $s$. We know that \eqref{e:FPgen} has a unique solution if $|t-s|$ is small enough and,
in general, that it admits at most one solution (since one can show that the restriction
of any solution to \eqref{e:FPgen} to some smaller interval $[u,v] \subset [s,t]$
also solves the same problem on $[u,v]$, with $X_s$ replaced by $(\CR\hat X)(u,\place)$).

\begin{remark}
Recall that, for a given model $Z$, elements of $\CD^{\gamma,\eta}_{s,t}$ consist of maps
$X \colon [s,t] \times \T^2 \to \Y \otimes \CT$.
As already discussed in Remark~\ref{rem:NotationNoise}, $\<IXi>$  should naturally be interpreted as 
an element of $\Y \otimes \CT[\<IXi>] \subset \Y \otimes \CT$ and therefore as a 
constant element of $\CD^{\gamma,\eta}_{s,t}$.
\end{remark}

We now write $\DD_{s,t} \subset \Omega_{\kappa}^\YMHs \times \CN$ for the set of 
pairs $(X_s,Z)$ such that \eqref{e:FPgen} admits a solution 
$\hat X$ with $(\hat X,Z) \in \hat\CD^{\gamma,\eta}_{s,t}$ and we set
\begin{equ}[e:defSst]
\CS_{s,t}(X_s,Z) = \CR\CM_{s,t}(\hat X) + \<IXib> + \hat \fp_{\place-s} (X_s - \<IXib>_s)\;,
\end{equ}
where $\CR$ denotes the reconstruction operator.

\begin{remark}
In light of \eqref{e:FPgen}, we could have just written $\CS_{s,t}(X_s,Z) = \CR \hat X$. Unlike \eqref{e:defSst},
this however would not have made it clear that $\CS_{s,t}$ takes values in $\cC([s,t], \Omega_{\kappa}^\YMHs)$.
\end{remark}

Denote by $\RR$ the renormalisation group (defined as in \cite{H14,BHZ19}) associated
to the regularity structure $\CT$ and let $\RR_0 \subset \RR$ be the subgroup generated 
by all BPHZ renormalisations of stationary Gaussian noises with a covariance that 
is symmetric under the two spatial reflections.
We then collect the properties of this construction in
one place.
\begin{theorem}\label{theo:localFPRegStruc}
The solution map $\CS$ and its domain $\DD$ have the following properties.
\begin{enumerate}[label=(\roman*),itemsep=0ex]
\item\label{pt:wellposedness} For any fixed $s,t$, $\DD_{s,t}$ is open and $\CS_{s,t}$ is continuous as a map from
$\DD_{s,t}$ to $\cC([s,t], \Omega_{\kappa}^\YMHs)$.
\item For every $R \ge 1$ and $T > 0$ there exists $\delta > 0$ such that, for all $(X_s,Z) \in \Omega_{\kappa}^\YMHs\times \CN$ 
with $\$X_s\$_{\kappa} + \$Z\$ \le R$ and all $0\le s<t\le T$ with $t-s \le \delta$, one has $(X_s,Z) \in \DD_{s,t}$. 
\item For $s<t<u$ and $(X_s,Z)$ such that $(X_s,Z) \in \DD_{s,u}$, one has $(X_s,Z) \in \DD_{s,t}$ and
\begin{equ}
\CS_{s,t}(X_s,Z) = \CS_{s,u}(X_s,Z) \restr [s,t]\;.  
\end{equ}
\item\label{pt:iv} Conversely, for $s<t<u$ and $(X_s,Z)$ such that $(X_s,Z) \in \DD_{s,t}$ and 
$(\CS_{s,t}(X_s,Z)(t),Z) \in \DD_{t,u}$, one has $(X_s,Z) \in \DD_{s,u}$ and
\begin{equ}
\CS_{s,u}(X_s,Z) = \CS_{s,t}(X_s,Z) \sqcup \CS_{t,u}(\CS_{s,t}(X_s,Z)(t),Z)\;.
\end{equ}
Here, $\sqcup$ denotes the natural concatenation operator.
\item If $X_s$ is smooth and $Z = M_g \LL(\eta)$ for a smooth $\eta$ and some $g \in \RR_0$,\footnote{Since 
the action of $\RR$ onto $\models_0$ leaves $\Pi_0\<IXi>$ invariant, it naturally lifts to an action on $\CN$.} then 
$\CS_{s,t}(X_s,Z)$ is the unique classical solution to the parabolic PDE
\begin{equ}[e:renormGen]
\d_t X = (\Delta-1) X + \linTerm X + \hat \bil(X,DX) + \tril(X,X,X) + (\extraC + C_g) X + \eta\;,
\end{equ}
with initial condition $X(s) = X_s$, where $C_g\in \End(\Y)$ depends on the value of $g$.
\item\label{pt:Markov} If $\boldxi$ is a $\Y$-valued space-time white noise and $\boldxi \mapsto \CZ(\boldxi)$ denotes the 
BPHZ lift of $\boldxi$ mentioned in Remark~\ref{rem:randomModel}, 
then the map $\boldxi \mapsto \CS_{s,t}(X_s,\CZ(\boldxi))$ is measurable
with respect to the $\sigma$-algebra $\F_{s,t}$ generated by evaluations of $\boldxi$ against test functions
supported in the time interval $[s,t]$. In particular, $\CS$ generates a white noise random dynamical system
(and therefore also a Markov process) on $\Omega_{\kappa}^\YMHs \cup \{\skull\}$.\footnote{Here, $\skull$ 
denotes a cemetery state which we include at this stage since solutions could potentially blow up.}
\end{enumerate}
\end{theorem}

\begin{proof}
These results are by now standard, so we mainly provide relevant pointers to the existing literature.
\begin{enumerate}[label=(\roman*),itemsep=0ex]
\item\label{pt:i} If we choose $\gamma > 1+\kappa$, $\eta \in (-\frac14,-\kappa)$, and $\kappa$ small enough,
it follows from \cite[Prop.~6.12]{H14} that the map $\CM_{s,t}$ is continuous from $\CD^{\gamma,\eta}_{s,t}$ into
$\CD^{\bar \gamma,\bar \eta}_{s,t}$, provided that $\bar \gamma < \gamma + 1-\kappa$ and $\bar \eta \le 2\eta+1$.
Since it furthermore takes values in a sector of regularity $1-2\kappa$ 
(excluding the polynomial terms), it follows from \cite[Prop.~3.28]{H14}
and \cite[Def.~6.2]{H14} that, provided that we choose $\bar \eta > \frac12$ (which we can since we chose $\eta > -\frac14$), 
the map $(\hat X,Z) \mapsto \CR \CM_{s,t}(\hat X)$ is continuous from 
$\CD^{\gamma,\eta}_{s,t}$ into $\CC([s,t],\CC^{1/2}_x)$ and therefore also into $\CC([s,t],\Omega_{\kappa}^\YMHs)$ by Remark~\ref{rem:embedding}.
The fact that $\hat\fp_{\place-s}(X_s - \<IXib>_s) \in \CC([s,t],\Omega_{\kappa}^\YMHs)$ follows from the strong 
continuity of the heat semigroup on $\Omega_{\kappa}^\YMHs$, see \cite[Prop.~4.6]{CCHS22}.
\item This is a consequence of the argument in point~\ref{pt:i}, see also the discussion before \eqref{e:contract}.
\item This follows from the definition of the integration operator $\CP$, the uniqueness of solutions, 
and the fact that the heat kernel is non-anticipative.
\item This is essentially the content of \cite[Prop.~7.11]{H14}.
\item This follows from simple power counting arguments, see for example the general form of the counterterms
given by the renormalisation group determined in \cite{BCCH}. Note that there is no constant term generated by
$\RR_0$ as a consequence of the structure of the nonlinearity, see the discussion preceding \eqref{e:counterTerm6}
below. (But even if we had such a term this wouldn't change our argument.) 
\item This follows from the fact that the approximate driving noise $\boldxi_\eps$ mentioned in Remark~\ref{rem:randomModel} 
is adapted to the filtration generated by $\boldxi$, combined with point~\ref{pt:iv}
and the fact that it only looks into the past by time about $\eps^2$.\qedhere
\end{enumerate}
\end{proof}

We henceforth write $\{\CP_t\}_{t \ge 0}$ for the Markov semigroup on
$\bar \Omega_{\kappa}^\YMHs = \Omega_{\kappa}^\YMHs \cup \{\skull\}$ generated by the process $X_t$
obtained in point~\ref{pt:Markov} of this theorem.
In the remainder of this section, we will show that $\CP_t$ has ``nice'' properties
from the perspective of the general theory of Markov processes. 
The reader may worry about the
fact that different choices of $\extraC$ lead to 
different Markov processes. This is irrelevant in this section since our results will in fact be
independent of that choice. The specific choice that we will make
is the one that allows to identify the process constructed in Theorem~\ref{theo:localFPRegStruc} with
the one studied in the previous sections. This will be detailed in the next section,
but we refer to \cite{CCHS22,CCHS22+,CS26} where it was shown that there is a choice
leading to a process that is still Markov under the projection
$\Omega_{\kappa}^\YMHs \to \fO_{\kappa}^\YMHs$ and that this in fact uniquely
determines the choice of $\extraC_A$.

Recall now the definition of a generative probability measure on
the space of càdlàg functions $\R_+ \to \bar\Omega_{\kappa}^\YMHs$
given in \cite[Def.~2.11]{CCHS22} (with the obvious modifications to account for
the presence of the Higgs fields). Loosely speaking, these are laws of processes
following the solution map constructed in Theorem~\ref{theo:localFPRegStruc},
but allowing for the solution to be re-gauged at a sequence of stopping times
in order to avoid explosion along gauge orbits. We also write 
$E_t \colon D(\R_+,\bar\Omega_{\kappa}^\YMHs) \to \bar\Omega_{\kappa}^\YMHs$ for 
the evaluation map at time $t$ and $\Pi\colon\bar\Omega_{\kappa}^\YMHs \to 
\fO_{\kappa}^\YMHs \cup \{\skull\}$ for the canonical projection (extended 
to map the cemetery state to the cemetery state).

\begin{proposition}\label{prop:semigroupOrbits}
There exists a unique choice of $\extraC_A$ such that, for every $\liegroup$-invariant
choice of $\extraC_\Phi$, there exists 
a unique Markov semigroup $\hat \CP$ on $\bar\fO_{\kappa}^\YMHs$
with the property that, for every $X_0 \in \Omega_{\kappa}^\YMHs$, every $t > 0$
and every generative probability measure $\mu$ with $\mu(X(0) = X_0) = 1$,
one has $\hat \CP_t([X_0], \place) = \Pi^* E_t^* \mu$.
Furthermore, one has $\hat \CP_t([X_0], \place) \ge \Pi^*\CP_t(X_0,\place)$ 
on $\fO_{\kappa}^\YMHs$ and
\begin{equ}[e:TVconvergence]
\lim_{t \to 0} \|\hat \CP_t([X_0], \place) - \Pi^*\CP_t(X_0,\place)\|_\TV = 0\;.
\end{equ}
\end{proposition}

\begin{proof}
This is essentially the content of \cite[Thm~2.13]{CCHS22}, since incorporating 
the Higgs field is straightforward. We will see in Proposition~\ref{prop:processesEqual} 
below what is the ``correct''
choice of $\extraC$ which guarantees that the process constructed here 
coincides with the one constructed in the previous sections. 
\end{proof}

\subsection{Irreducibility of the Markov process}\label{sec:irreducible}

The aim of this section is to show that the Markov process constructed in Section~\ref{sec:Markov}
is topologically irreducible. A version of irreducibility that is sufficient for our purpose
is that, for every initial condition
$X_0 \in \fO_{\kappa}^\YMHs$ there exists $T>0$ such that, for every non-empty
open set $\CA \subset \fO_{\kappa}^\YMHs$, one has
$\P(X_T \in \CA) > 0$, where $X_t$ denotes the process at time $t$ with initial condition $X_0$.
In fact, we will show the following slightly stronger version of this statement which holds on pathspace
and before quotienting by the gauge group.

\begin{theorem}\label{theo:fullSupport}
Given $X_0 \in \Omega_{\kappa}^\YMHs$, let $(X,\tau)$ be the maximal solution to 
\eqref{e:schematic} with (potential) explosion time $\tau> 0$. Then, there exists $T > 0$ such that
$\P(\tau > T) > 0$ and such that the law of $X$ restricted to the event $\tau > T$ has
support in $\CC([0,T],\Omega_{\kappa}^\YMHs)$ given by all functions with given initial condition $X_0$.
\end{theorem}

\begin{remark}
Note that this statement neither uses nor implies the fact that our solutions are global.
In fact, nothing in this article rules out the possibility that the solutions 
to \eqref{e:schematic} explode in finite time. Our global result relies crucially on the fact that
we have the freedom to ``regauge'' them on the fly.
\end{remark}

\begin{remark}
In principle, one would expect that irreducibility follows from the results of
\cite{SupportGeneral}. The problem is
that the elements of ``real'' homogeneity $0$ of the regularity structure 
(of which there are many in our setting) cause some non-trivial technical difficulties in that
paper. In particular, Assumption~6 of \cite{SupportGeneral} fails, so would have to be replaced by Assumptions~7 and~8. 
Assumption~7 is easy to verify, but whether Assumption~8 holds is much less obvious.
Instead, we choose to give a relatively short self-contained proof. We believe that this
proof extends to a large class of singular SPDEs driven by non-degenerate noise, so it should be of interest
beyond the setting of this article.
\end{remark}

In order to compare modelled distributions corresponding to different models, it is also 
convenient to note that, for any given $s < t$, there exists a Banach space $\CB$ and continuous injections 
(that are not bijections!) $\iota \colon \Picard_{s,t}(Z) \to \CB$.
The space $\CB$ consists of continuous functions with values in the 
sector $\CT_{<\gamma}$ endowed with the norm $|F| = \sup_{r \in (s,t]}\sup_{x\in \T^2} (r-s)^{3}|F(r,x)|$.
The reason for mentioning this is just to make the point that there is a natural way of
stating whether elements of $\Picard_{s,t}(Z)$ with different values of $Z$ are equal.

\begin{remark}\label{rem:genFP}
A crucial fact that will be used in the sequel is the following extension of the relation between
\eqref{e:renormGen} and \eqref{e:FPgen}.
Consider a model $Z\in \CN$ of the form $M_g \LL(\eta)$ for a smooth function $\eta$
and some $g \in \RR_0$. 
Then, any solution $\hat X\in \Picard_{s,t}(Z)$ to
\begin{equ}[e:FPgen2]
\hat X = \CM_{s,t}(\hat X) + \<IXi> + \CJ f + \CJ \hat\fp_{\place-s} (X_s-\<IXib>_s)\;,
\end{equ}
with smooth initial data $X_s$, where $f \colon [s,t]\times\T^2 \to \Y$ is any smooth function with $\CJ f \in \Picard_{s,t}(Z)$, 
must be such that $X = \rR \hat X$ solves the fixed point problem
\begin{equ}
X = \hat \fp_{s,\place} \bigl((\linTerm + \extraC) X + \hat \bil(X,DX) + \tril(X,X,X) + C_g X + \eta\bigr) + \hat \fp_{\place-s}X_s + f\;,
\end{equ}
in $\CC^1([s,t] \times \T^2,\Y)$. Here, $\hat \fp_{s,t} F = \int_s^t \hat \fp_{t-r} F_r\,dr$, with $\hat \fp_t$ denoting as usual
convolution with the massive heat kernel at time $t$.
\end{remark}

Given $h \in \CC^{-\kappa}$, we also write $\tau_h \colon \models_0 \to \models_0$ for the shift map with the property
that $\tau_h \bigl(M_g \LL(\eta)\bigr) = M_g \LL(\eta + h)$ for all smooth $\eta$. Note that in the setting 
of \cite[Sec.~4.1]{HS23}, one has $\tau_h = Z^1$ with $Z^\eps$ as in \cite[Prop.~4.13]{HS23}.
In particular, this shows that this map is well-defined and, by \cite[Prop.~4.15]{HS23} (combined with the
continuity of the action of the renormalisation group on $\models_0$ and the continuity of the lift of $h$ 
into the extended regularity structure which holds by \cite[Thm~5.14 \;\&\; Prop.~3.31]{H14}),
jointly continuous as a map $\CC^{-\kappa} \times \models_0 \to \models_0$.\footnote{We can simply take
$\theta = 0$ in the construction of \cite[Sec.~4.1]{HS23}. While this causes $\tilde \Xi$ to have (slightly)
negative degree, this is harmless.} See also \cite{MallGPAM,StrongFeller} for very similar constructions.

Since convolution with the
heat kernel maps $\CC^{-\kappa}$ into $\CC^{2-\kappa}$ which is in turn included in $\CC(\R, \Omega_{\kappa}^\YMHs)$
by Remark~\ref{rem:embedding}, the map $(h,Z) \mapsto \tau_h(Z)$ extends to a continuous 
map $\CC^{-\kappa} \times \CN \to \CN$.
We then have the following result.

\begin{theorem}\label{theo:idenShiftH}
Given $Z \in \CN$, $Y_s \in \Omega_\kappa^\YMHs$, and $h \in \CC([s,t]\times \T^2,\Y)$, if $\hat Y$ solves 
\begin{equ}[e:FPgeninhom]
\hat Y = \CM_{s,t}(\hat Y) + \<IXi> + \CJ \bigl(\hat\fp_{\place-s} (Y_s - \<IXib>_s) + \hat\fp_{s,\place} h\bigr)
\end{equ}
in $\Picard_{s,t}(Z)$, then one has $(Y_s,\tau_hZ) \in \DD_{s,t}$ and $\CS_{s,t}(Y_s,\tau_h Z)$
is given by
\begin{equ}
\CS_{s,t}(Y_s,\tau_h Z) = \CR \hat Y \eqdef \CR\CM_{s,t}(\hat Y) + \<IXib> + \bigl(\hat\fp_{\place-s} (Y_s - \<IXib>_s) + \hat\fp_{s,\place} h\bigr)\;.
\end{equ}
\end{theorem}

\begin{proof}
By the definition of $\CN$, there exist sequences $\eta_n$ and $g_n \in \RR_0$ with $\eta_n$ smooth such that
$Z = \lim_{n \to \infty} Z_n$ with $Z_n =  M_{g_n} \LL(\eta_n)$. In particular, it follows 
from Remark~\ref{rem:genFP} and the continuity of the solution map that $Y = \rR\hat Y$ is given 
by $Y= \lim_{n \to \infty} Y_n$ where $Y_n$ solves
\begin{equ}
\d_t Y_n = (\Delta-1) Y_n + (\linTerm + \extraC) Y_n +  \hat \bil(Y_n,DY_n) + \tril(Y_n,Y_n,Y_n) + C_{g_n} Y_n + \eta_n + h\;.
\end{equ}
On the other hand, it follows from the definition of the translation operator $\tau_h$ that
$\tau_h Z = \lim_{n \to \infty} M_{g_n} \LL (\eta_n+h)$, so that the claim follows from
combining \eqref{e:FPgen} with Remark~\ref{rem:genFP}. 
\end{proof}

We now have all the ingredients in place to provide a proof of Theorem~\ref{theo:fullSupport}.

\begin{proof}[Proof of Theorem~\ref{theo:fullSupport}]
Since $\CN$ is Polish, the law of the random model $\CZ$ from Remark~\ref{rem:randomModel}
has non-empty support, so we
can fix an element $Z \in \supp \Law(\CZ)$. As usual, we write $\<IXib>$ for the $\CX$-component of $Z$.
We also write $C$ for the norm of the bounded linear operator $\CJ\hat \fp_{\place} \colon \Omega_{\kappa}^\YMHs \to \CD^{\gamma,\eta}_{0,1}$.
Fixing $X_0 \in \Omega_{\kappa}^\YMHs$, we write 
$R = 2C(1+\$X_0\$ + \$Z\$)$, and we fix $T>0$ small enough so that \eqref{e:contract} holds with $\eps = 1/2$.

Fix now some $\bar\kappa \in (1/2,1-2\kappa)$ as well as an arbitrary $X \in \CC([0,T],\Omega_{\kappa}^\YMHs)$ with $X(0) = X_0$ and
such that, for some $R_X > 0$, one has $\|X - \<IXib> - \hat\fp_{\place} (X_0-\<IXib>_0)\|_{\CC^{\bar\kappa}} \le R_X$ on the interval $[0,T]$. 
(Here we use $\CC^{\bar\kappa}$ to denote the corresponding space-time Hölder space in the parabolic scaling.) 
By \cite[Rem.~3.23]{CCHS22}, the space of
 $\CC^{\bar \kappa}_x$ functions on the torus embeds densely into $\Omega_{\kappa}^\YMHs$ as soon as $\bar \kappa \ge (1-\kappa)/2$,
 which is the case here. In particular, any element $\tilde X \in \CC([0,T],\Omega_{\kappa}^\YMHs)$ with initial condition $X_0$
can be approximated to arbitrary precision by such an $X$, provided that $R_X$ is taken large enough.
 
We also fix
 an arbitrary $\eps \in (0,1]$, and set $R' = R + R_X$.
We want to show that it is possible to find a bounded $h $ such that $(X_0,\tau_h Z) \in \DD_{0,T}$
 and such that the solution $Y = \CS_{0,T}(X_0,\tau_h Z)$ satisfies 
\begin{equ}[e:desiredTracking]
\sup_{t \in [0,T]} \|Y_t - X_t\|_{\CC^{1/2}_x} \le \eps\;,
\end{equ}
which then implies the desired statement by \cite[Lem.~3.12]{SupportGeneral}.

In order to construct $h$, we subdivide $[0,T]$ into intervals of size $\delta \ll 1$ and we proceed to construct $h$
one interval at a time. Our main workhorse for that construction is the claim
that it is possible to find $\delta$ small enough so that, for
every $s \in [0,T-\delta]$ and every $Y_s \in \Omega_{\kappa}^\YMHs$ with $\|Y_s - X_s\|_{\CC^{1/2}_x} \le \eps /2$, there 
exists $h \in L^\infty([s,s+\delta]\times \T^2,\Y)$  
such that 
\begin{enumerate}[label=(\roman*),itemsep=0ex]
\item\label{claim1} Setting $F_t = X_t - \<IXib>_t - \hat\fp_{t-s} (Y_s - \<IXib>_s)$, one has $\|F_t -  \hat\fp_{s,t}h\|_{\CC^{1/2}_x} \le 3\eps/4$ for all $t \in [s,s+\delta]$
and furthermore this quantity is bounded by 
$\eps / 4$ at time $t = s+\delta$.
\item\label{claim2} The solution $\hat Y$ to \eqref{e:FPgeninhom} with $Y_s$ as just described exists, is unique, and
satisfies the bound $\sup_{t \in [s,s+\delta]}\bigl\| \bigl(\rR\CM_{s,s+\delta}(\hat Y)\bigr)_t\bigr\|_{\CC^{1/2}_x} \le \eps/4$.
\end{enumerate}

The point here is that, by \eqref{e:FPgeninhom}, one has 
$X_t - Y_t = F_t - \hat\fp_{s,t} h - \bigl(\rR\CM_{s,s+\delta}(\hat Y)\bigr)_t$, so that 
these two claims taken together immediately yield the bound
$\|X_t - Y_t\|_{\CC^{1/2}_x} \le \eps$ over $[s,s+\delta]$ for 
$Y_t = \CS_{s,t}(Y_s,\tau_h Z)_t$.

To show the desired claim, fix $s \ge 0$ and rewrite $F_t$ as 
\begin{equ}
F_t = X_t - \<IXib>_t - \hat\fp_t (X_0 - \<IXib>_0) - \hat\fp_{t-s} (Y_s - \<IXib>_s - \hat\fp_s (X_0 - \<IXib>_0) )\;.
\end{equ}
Since $\|X_t - \<IXib>_t - \hat\fp_t( X_0 - \<IXib>_0)\|_{\CC^{\bar\kappa}} \le R'$ by assumption 
and since the remaining term in $F$ solves the massive
heat equation (and $\bar \kappa > \frac12$), there exists an exponent $a >0$ such that
\begin{equ}[e:heatEquFt]
\|F_t - \hat\fp_{t-r} F_r\|_{\CC^{1/2}_x} \lesssim R'|t-r|^{a}\;.
\end{equ}

We also fix
a smooth function $\chi\colon \R_+ \to \R_+$ with $\supp \chi \subset (0,\frac12]$ such that $\int_0^\infty \chi(x)\,dx = 1$ and we write $\chi_{\delta}(t) = \delta^{-1}\chi(t/\delta)$. Given $\delta > 0$, we then write
\begin{equ}[e:defHfinal]
H_t = \int_s^t \chi_{\delta}(t-r) \hat\fp_{t-r} F_r\,dr\;.
\end{equ}
One important feature of this definition is that for $\gamma \le 2$ and
$\eta < 0$ one has $\|H\|_{\CC^{\gamma,\eta}} \lesssim \delta^{|\eta/2|}$, with the
$\CC^{\gamma,\eta}$-norm taken on $[s,s+\delta]\times \T^2$. This is a consequence of the fact that
the supremum norm of $F$ is bounded by some constant independent of $\delta$ by assumption and therefore,
by standard properties of the heat kernel combined with the assumptions on $\chi_\delta$, $H$ is smooth
and satisfies a bound of the type $|\d_t^k D^\alpha H(t,x)| \lesssim \delta^{-k-|\alpha|/2} \lesssim \delta^{|\eta/2|} |t-s|^{ \frac {\eta-2k-|\alpha|}2}$. (Here we used the fact that $\eta < 0$ and $|t-s| \le \delta$ to obtain the second inequality.)

A simple calculation then shows that $(\d_t + 1 - \Delta) H_t = h_t$ with  
\begin{equ}[e:defhfinal]
h_t = \int_s^t \chi'_{\delta}(t-r) \hat\fp_{t-r} F_r\,dr\;,
\end{equ}
so that, since $H_s = 0$, one has $H = \hat\fp_{s,\place} h$.
Setting $\rho_{\delta}(r) = \int_0^r \chi_{\delta}(r')\,dr'$, we then have
\begin{equ}
\|H_t - \rho_{\delta}(t-s) F_t\|_{\CC^{1/2}_x} \le \int_s^t \chi_{\delta}(t-r)  \|\hat\fp_{t-r} F_r - F_t\|_{\CC^{1/2}_x}\,dr
\lesssim \int_s^t \chi_{\delta}(t-r)  (t-r)^{a}\,dr \lesssim \delta^{a}\;.
\end{equ}
Since it furthermore follows from \eqref{e:heatEquFt} and our assumption
that $\|F_s\|_{\CC^{1/2}_x} = \|X_s-Y_s\|_{\CC^{1/2}_x} \le \frac\eps2$ that
\begin{equ}
\|F_t\|_{\CC^{1/2}_x}
\le \|F_t - \hat\fp_{t-s}F_s\|_{\CC^{1/2}_x} + \|\hat\fp_{t-s}F_s\|_{\CC^{1/2}_x}
\le |t-s|^{a} + \frac\eps2\;,
\end{equ}
we conclude that there exists $a>0$ such that one has for $t \in [s,s+\delta]$ the bound
\begin{equ}
\|H_t - F_t\|_{\CC^{1/2}_x}
\le \frac\eps2 \one_{|t-s| \le \delta/2} + C\delta^{a}\;.
\end{equ}
Since $H_t - F_t$ is precisely the quantity appearing in the claim \ref{claim1},
this shows that the first property holds with the choice of $h$ given by \eqref{e:defhfinal},
provided that $\delta$ is chosen small enough.

To show that \ref{claim2} also holds, note first that,
by \eqref{e:defhfinal} combined
with the smoothing properties of the heat kernel,
 there exists a constant $C$
(possibly depending on $\eps$ and $R'$)
such that $\|h\|_{L^\infty} \le C \delta^{-1}$. 
Furthermore, we have already seen in the discussion following \eqref{e:defHfinal} that, for any $\eta < 0$,
\begin{equ}[e:boundJPh]
\|\CJ \hat\fp_{s,\place}h\|_{\Picard_{s,s+\delta}} = \|\CJ H\|_{\Picard_{s,s+\delta}} \lesssim  \delta^{|\eta/2|}
\le 1\;,
\end{equ}
where the last inequality possibly requires choosing $\delta$ a bit smaller.
We now choose $\gamma > 1+\kappa$ and $\eta \in (-\frac14,-\kappa)$ as in the proof of point~\ref{pt:i} of
Theorem~\ref{theo:localFPRegStruc}.
If we choose $\delta$ small enough, we can guarantee that the fixed point problem 
\eqref{e:FPgeninhom} admits a unique solution $\hat Y$ over $[s,s+\delta]$ which is such that
$\|\hat Y\|_{\Picard_{s,s+\delta}}$ is bounded, uniformly over $\delta$ small enough.
As in the proof of point~\ref{pt:i} of Theorem~\ref{theo:localFPRegStruc}, we
conclude that $\CM(\hat Y) \in \CD^{\bar \gamma,\bar \eta}_{s,t}$ with norm
bounded uniformly in $\delta$, provided that $\bar \gamma < \gamma + 1-\kappa$ and $\bar \eta \le 2\eta+1$. Since we can take $\bar \eta > \frac12$ and since 
$\CM(\hat Y)_s$ vanishes up to degree $1-2\kappa$, we
can find an exponent $a>0$ such that
 $\sup_{t \in [s,s+\delta]}\|\rR\CM(\hat Y)_t\|_{\CC^{1/2}_x} \le C \delta^a  \le \eps/4$.

In order to show that \eqref{e:desiredTracking} holds, it now suffices to iterate the fixed
point problem \eqref{e:FPgeninhom} 
on intervals of size $\delta$ and to use Theorem~\ref{theo:idenShiftH} to identify
the resulting process with a solution to the original equation with model $\tau_h Z$.
Assume that $s$ is such that $\|Y_s - X_s\|_{\CC^{1/2}_x} \le \eps/2$ (which is
certainly the case for $s=0$ since we consider the solution starting at $X_0$).
On the time interval $[s,s+\delta]$, the solution $Y = \rR \hat Y$ to our fixed point problem
for the model shifted by $h$ chosen as in \eqref{e:defhfinal} is then given by
\begin{equ}
Y_t = X_t - F_t + H_t + \bigl( \rR \CM(\hat Y)\bigr)_t\;.
\end{equ}
Claim~\ref{claim1} guarantees that on each interval the solution to the linearised equation
tracks the target $X$ with an error $F-H$ of at most $3\eps/4$ (in $\CC^{1/2}_x$), while
claim~\ref{claim2}
ensures that the effect of the nonlinearity is at most $\eps/4$, so that one does indeed
have $\|Y_t - X_t\|_{\CC^{1/2}_x} \le \eps$ for all $t \in [s,s+\delta]$.
Furthermore, and this is crucial, claim~\ref{claim1} guarantees that 
one has $\|Y_{s+\delta} - X_{s+\delta}\|_{\CC^{1/2}_x} \le \eps/2$, which allows us to iterate the
argument up to the final time $T$, so that \eqref{e:desiredTracking} does indeed hold, concluding
the proof.
\end{proof}

%
%
%

\section{Identification of the two dynamics}
\label{section:langevin}

The goal of this section is to show that, when quotiented by the action of the gauge group and for a
suitable choice of renormalisation constant,
the maximal solutions to the dynamic constructed in Theorem~\ref{thm:global_soft} 
coincide with those constructed in
Theorem~\ref{theo:localFPRegStruc}. As already pointed out, this is important for two reasons:
\begin{enumerate}
\item We know from \cite{IlyaHaoYM} that, in the absence of the Higgs field, the invariant measure of the
dynamic constructed in \cite{CCHS22} coincides with the Yang--Mills measure that was previously
constructed in \cite{Sengupta1,L03,C19}. This provides a good reality check showing that the non-Abelian Higgs
model constructed in this article is the ``correct'' one.
\item It is clear that the  dynamic constructed in Theorem~\ref{thm:global_soft} is gauge covariant
(in law) as shown in Lemma~\ref{main:lem-gauge-covariance} already for fixed $\eps > 0$. However, while 
gauge covariance holds for every $\liegroup$-invariant
choice of the renormalisation operators $\renormA$ and $\renormphi$, it is  
simply not the case in general that the limiting process is Markovian. 
On the other hand, we know from \cite{CCHS22} (and in fact from several previous
constructions using regularity structures) that the process constructed in Theorem~\ref{theo:localFPRegStruc} is 
Markovian for every choice of renormalisation (provided of course that it yields convergence of 
the model).
\end{enumerate}

Showing that it is possible to choose the renormalisation in such a way that the two processes 
coincide therefore provides a guarantee both that we
are constructing the ``correct'' measure and that the process we are considering is in fact a Markov
process. Since we already showed that this Markov process is topologically irreducible, we can 
then conclude that it admits a unique invariant measure, provided that we show that 
it has the strong Feller property and is tight.

\subsection{The three main processes}
\label{sec:threeProcesses}

Besides showing that Theorems~\ref{thm:global_soft} and~\ref{theo:localFPRegStruc}
yield the same processes (in law and modulo gauge transformation), we
 also verify that this is the exact same renormalisation which causes the latter
process to be gauge covariant.
While we already know from the argument of \cite{CCHS22+} that there is \textit{some} renormalisation
procedure for which this is the case, it is not clear a priori that this is the same
as the one that causes it to coincide with the process of Theorem~\ref{thm:global_soft}.
We could in principle verify this by comparing with the constants from \cite{CCHS22}, but since that
article doesn't include the Higgs field and since it follows almost immediately from the exact 
same calculations required to show that the processes of Sections~\ref{section:covariant_lwp} and~\ref{sec:regStruc} coincide, we
include this here.

%
%
%
%

Recall that, if we write \eqref{e:schematic} more explicitly, the construction of the
previous section provides a limit in probability as $\eps \to 0$
to the solution to the system of SPDEs
\begin{equs}[e:wantedSPDE]
\d_t \hat A_\eps &= - \covd_{\hat A_{\eps}}^\ast F_{\hat A_{\eps}} - \covd_{\hat A_{\eps}} \covd_{\hat A_{\eps}}^\ast \hat A_\eps + \renormA \hat A_\eps  + \bil(\hat\higgs_\eps \otimes \covd_{\hat A_{\eps}}\hat\higgs_\eps) +\xi_{\eps}\;,\\
\d_t \hat \higgs_\eps &= -\covd_{\hat A_{\eps}}^\ast \covd_{\hat A_{\eps}} \hat \higgs_\eps + 
\bigl(\covd_{\hat A_{\eps}}^\ast \hat A_\eps \bigr) \hat \higgs_\eps +
\renormphi \hat \higgs_\eps - \hat \higgs_\eps |\hat \higgs_\eps|^2 + \zeta_{\eps}\;,
\end{equs}
where $\bil$ is as in \eqref{intro:eq-bil},  the renormalisation
operators $\renormA$ and $\renormphi$ are as in \eqref{main:eq-renormalization-operators},
$\xi_\eps$ is defined by
\begin{equ}[e:mollif]
\xi_{\eps}(t,x) = (\moll \star \xi)(t,x) \eqdef \iint \ds \dy \rho_\eps^\prime(t-s) \hat\fp(s,y;t,x) \xi(s,y)\;,
\end{equ}
where $\star$ denotes space-time convolution, and similarly for $\zeta_\eps$.
In view of \eqref{local:eq-xi-zeta-smooth-approximation}, one essentially 
has $\xi_\eps = \xi_{0,\eps}$, except for the presence of the mass term and the
fact that the time integral in \eqref{e:mollif} runs over all of $\R$, so that
$\xi_\eps$ is space-time stationary.
Note also that $t \mapsto \xi_\eps(t,\place)$ (and again the same for $\zeta_\eps$) 
is a stochastic process that is adapted
to the filtration generated by $(\xi,\zeta)$.

\begin{remark}
When considering the canonical lift of models, we will always use the massive heat 
kernel as our integration operator, which has the advantage of not requiring any
arbitrary large scale cut-off. (Of course the value of the mass is arbitrary
though.)
This is also why we are using the massive heat kernel $\hat\fp$ in \eqref{e:mollif},
as it then yields nicer explicit
expressions for the BPHZ renormalisation constants.
\end{remark}

Given a background connection $B$, we also recall that the process $(A_\eps,\higgs_\eps)$
is given by
\begin{equs}[e:covariant_SPDE]
\d_t  A_\eps &= - \covd_{A_{\eps}}^\ast F_{A_{\eps}} - \covd_{A_{\eps}} \covd_{A_{\eps}}^\ast (A_\eps - B) +\renormA (A_\eps - B)+  \bil(\higgs_\eps\otimes \covd_{A_{\eps}}\higgs_\eps) + \xi_{\Beps}\;,\\
\d_t \higgs_\eps &= -\covd_{A_{\eps}}^\ast \covd_{A_{\eps}} \higgs_\eps + \bigl(\covd_{A_{\eps}}^\ast (A_\eps - B)\bigr)\higgs_\eps +
\renormphi  \higgs_\eps -  \higgs_\eps | \higgs_\eps|^2 + \zeta_{\Beps}\;,
\end{equs}
where $\xi_{\Beps}$ is defined in \eqref{local:eq-xi-zeta-smooth-approximation}.

%
%

Our aim now is to show that it is possible to find time-dependent gauge transformations
$g_\eps$ converging to some limit $g$ as $\eps \to 0$ and such that
$g_\eps \act A_\eps$ and $\hat A_\eps$ have the same limit in law as $\eps \to 0$.
In fact, we will construct a coupling between the noises driving the two equations
such that the limits of their maximal solutions agree.
For this, we make use of the fact that, if $g$ is an adapted process, then one has
$\Ad_g \xi \eqlaw \xi$, provided that the product appearing in this expression is
interpreted via Itô integration. In fact, it turns out in our case that even though $g$
is too irregular to be able to define $\Ad_g \xi$ via Young products, it is regular enough so that
Itô and Stratonovich integration coincide.

An explicit calculation shows that if we define $g_\eps$ as the solution to
\begin{equ}[e:equationG]
\d_t g_\eps^{(2)}\, \bigl(g_\eps^{(2)}\bigr)^{-1} = -d^* \bigl(dg_\eps^{(2)}\, \bigl(g_\eps^{(2)}\bigr)^{-1}\bigr)
+ d^* \bigl(\Ad_{g_\eps^{(2)}} B\bigr) - [g_\eps^{(2)} \act A_\eps, g_\eps^{(2)}\act B]\;,
\end{equ}
then $(A_\eps^{(2)},\higgs_\eps^{(2)}) \eqdef g_\eps^{(2)}\act  (A_\eps,\higgs_\eps)$ solves the equation
\begin{equs}
\d_t A_\eps^{(2)} &= - \covd_{A_\eps^{(2)}}^\ast F_{A_\eps^{(2)}} - \covd_{A_\eps^{(2)}} \covd_{A_\eps^{(2)}}^\ast A_\eps^{(2)} +
\renormA (A_\eps^{(2)} - g_\eps^{(2)} \act B) + \bil(\higgs_\eps^{(2)}\otimes \covd_{A_\eps^{(2)}}\higgs_\eps^{(2)}) +
\Ad_{g_\eps^{(2)}} \xi_{\Beps}\;,\\
\d_t \higgs_\eps^{(2)} &= -\covd_{A_\eps^{(2)}}^\ast \covd_{A_\eps^{(2)}} \higgs_\eps^{(2)} + \bigl(\covd_{A_\eps^{(2)}}^\ast  A_\eps^{(2)}\bigr) \higgs_\eps^{(2)} +
\renormphi \higgs_\eps^{(2)} - \higgs_\eps^{(2)} |\higgs_\eps^{(2)}|^2 + \rep(g_\eps^{(2)}) \zeta_{\Beps}\;,\label{e:second-main-process}
\end{equs}
This is ``almost'' the same as the equation \eqref{e:wantedSPDE}, but with two differences.
First, the noises $\xi_\eps$ and $\zeta_\eps$ are replaced by $\Ad_{g_\eps^{(2)}} \xi_{\Beps}$
and $\rep(g_\eps^{(2)}) \zeta_{\Beps}$. If they were
instead replaced by $\moll * (\Ad_{g_\eps^{(2)}} \xi)$ (and similarly for $\zeta$), then 
this would be the same
in law as $\xi_\eps$ since $\Ad_{g_\eps^{(2)}} \xi \eqlaw \xi$. Second, we have an additional
counterterm $-\renormA g_\eps^{(2)} \act B$ in the equation for $A_\eps^{(2)}$. 

In order to shorten notations, it will sometimes be convenient to combine, as in the 
previous section, the processes $A_\eps^{(2)}$ and $\higgs_\eps^{(2)}$
into a single $\Y$-valued process $Y_\eps^{(2)} = (A_\eps^{(2)},\higgs_\eps^{(2)})$.
We can then rewrite the equation solved by this process schematically similarly to 
\eqref{e:schematic} as
\begin{equ}[e:sys2]
\d_t Y_\eps^{(2)} = \Delta Y_\eps^{(2)} + \tril(Y_\eps^{(2)})
+ \hat\bil(Y_\eps^{(2)}, D Y_\eps^{(2)}) + \renormY \bigl(Y_\eps^{(2)} - g_\eps^{(2)} \act B\bigr) + \Repr(g_\eps^{(2)}) \boldxi_{\Beps}\;,
\end{equ} 
where we set $\renormY = \renormA \oplus \renormphi$ and we view $g_\eps^{(2)} \act B$ as taking
values in $\fg^2 \subset \Y$. 
We furthermore  wrote $\boldxi_{\Beps}$ for the $\Y$-valued noise $(\xi_{\Beps},\zeta_{\Beps})$.
We also used $\Repr = \Ad \oplus \Ad \oplus \rep$ for the representation of
$\liegroup$ on $\Y$ (we will use $\repr$ for its linearisation)
and we extended $B$ to a $\Y$-valued process by setting the component in $\V$ to $0$.
For conciseness, we also wrote $\tril(Y)$ rather than $\tril(Y,Y,Y)$.

We then consider the following system which is clearly equal in law to $(\hat A_\eps,\hat \higgs_\eps)$
(by Itô isometry and the fact that $g_\eps$ is sufficiently regular so that Itô and Stratonovich integrations
coincide) and therefore provides a coupling between
$(\hat A_\eps,\hat \higgs_\eps)$ and the gauge transformed copy of $(A_\eps,\higgs_\eps)$:
\begin{equs}
\d_t A_\eps^{(1)} &= - \covd_{A_{\eps}^{(1)}}^\ast F_{A_{\eps}^{(1)}} - \covd_{A_{\eps}^{(1)}} \covd_{A_{\eps}^{(1)}}^\ast A_\eps^{(1)} + \renormA A_\eps^{(1)} + \bil(\higgs_\eps^{(1)} \otimes \covd_{A_{\eps}^{(1)}}\higgs_\eps^{(1)}) +\moll * (\Ad_{g_\eps^{(1)}} \xi)\;,\\
\d_t \higgs_\eps^{(1)} &= -\covd_{A_{\eps}^{(1)}}^\ast \covd_{A_{\eps}^{(1)}} \higgs_\eps^{(1)} + 
\bigl(\covd_{A_{\eps}^{(1)}}^\ast A_\eps^{(1)} \bigr) \higgs_\eps^{(1)} +
\renormphi \higgs_\eps^{(1)} - \higgs_\eps^{(1)} |\higgs_\eps^{(1)}|^2 + \moll * (\rep(g_\eps^{(1)}) \zeta)\;,\\
\d_t g_\eps^{(1)}\, \bigl(g_\eps^{(1)}\bigr)^{-1} &= -d^* \bigl(d g_\eps^{(1)}\, \bigl(g_\eps^{(1)}\bigr)^{-1}\bigr)
+ d^* \bigl(\Ad_{g_\eps^{(1)}} B\bigr) - [A_\eps^{(1)}, g_\eps^{(1)}\act B]\;.
\end{equs}
Our aim then is to show that the triple $(A_\eps^{(1)},\higgs_\eps^{(1)},g_\eps^{(1)})$ has 
the same limit (in probability, not
just in law) as the triple $(A_\eps^{(2)}, \higgs_\eps^{(2)}, g_\eps^{(2)})$ as $\eps \to 0$.
As before, we will write $Y_\eps^{(1)} = (A_\eps^{(1)},\higgs_\eps^{(1)})$.

Finally, we will also consider the system of equations given by
\begin{equs}
\d_t A_\eps^{(3)} &= - \covd_{A_{\eps}^{(3)}}^\ast F_{A_{\eps}^{(3)}} - \covd_{A_{\eps}^{(3)}} \covd_{A_{\eps}^{(3)}}^\ast A_\eps^{(3)} + \renormA \bigl(A_\eps^{(3)}+d g_\eps^{(3)}\, \bigl(g_\eps^{(3)}\bigr)^{-1} \bigr) + \bil(\higgs_\eps^{(3)}\otimes\covd_{A_{\eps}^{(3)}}\higgs_\eps^{(3)}) +\Ad_{g_\eps^{(3)}} \xi_\eps\;,\\
\d_t \higgs_\eps^{(3)} &= -\covd_{A_{\eps}^{(3)}}^\ast \covd_{A_{\eps}^{(3)}} \higgs_\eps^{(3)} + 
\bigl(\covd_{A_{\eps}^{(3)}}^\ast A_\eps^{(3)} \bigr) \higgs_\eps^{(3)} +
\renormphi \higgs_\eps^{(3)} - \higgs_\eps^{(3)} |\higgs_\eps^{(3)}|^2 + \rep(g_\eps^{(3)}) \zeta_\eps\;,\\
\d_t g_\eps^{(3)}\, \bigl(g_\eps^{(3)}\bigr)^{-1} &= -d^* \bigl(d g_\eps^{(3)}\, \bigl(g_\eps^{(3)}\bigr)^{-1}\bigr)
+ d^* \bigl(\Ad_{g_\eps^{(3)}} B\bigr) - [A_\eps^{(3)}, g_\eps^{(3)}\act B]\;.
\end{equs}
As discussed in \cite[Sec.~2.2]{CCHS22} (see in particular Theorem~2.9 there), 
the gauge equivariance of the process $Y^{(1)} = \lim_{\eps \to 0} Y^{(1)}_\eps$
follows if one can show that $Y^{(1)}_\eps$ and $Y^{(3)}_\eps$ converge to the same limit
in the particular case $B=0$.
In principle, we do not need this since we already know that $Y^{(2)}_\eps$ is gauge equivariant,
so that if we can show that $Y^{(1)}_\eps$ and $Y^{(2)}_\eps$ have the same limit, then gauge equivariance
also follows for the former. However, the calculations needed to compare $Y^{(1)}_\eps$ and $Y^{(2)}_\eps$ 
allow us to also compare these to $Y^{(3)}_\eps$ without any additional effort, so we also include 
the comparison of $Y^{(1)}_\eps$ with $Y^{(3)}_\eps$ in the special case $B=0$ as it is quite
instructive to see how it automatically follows if one imposes that $Y^{(1)}_\eps$ and $Y^{(2)}_\eps$ 
have the same limit.

Most of the remainder of this section is devoted to showing that, at least locally in time, $Y_\eps^{(1)}$
and $Y_\eps^{(2)}$ have the same limit in probability, provided that 
one makes the ``correct'' choice of the renormalisation constant $\renormA$. In particular, we 
will see that the choice made in Definition~\ref{main:def-renormalization-constant} is correct (and this is in fact the only possible choice,
up to terms that vanish as $\eps \to 0$). Once this is done, we show that the resulting Markov process
has nice ergodic properties, which then completes the proof of our main result.

The main trick is to find one single regularity structure and a family of fixed point problems
(where we add a counterterm ``by hand'', depending on two parameters $\discrep$
and $\hat \discrep$) in the space
of modelled distributions over that structure such that the $Y_\eps^{(i)}$ can all three be viewed as
reconstructions of solutions to these fixed point problems, but for three different models. The only
disadvantage of this technique is that it forces us to coerce the system \eqref{e:schematic} into a slightly more
complicated setting than what we did in Section~\ref{sec:Markov}. The required
convergence statement then boils down to the fact that there is a choice of 
$\discrep$ and $\hat \discrep$
such that one can take these three models to all be suitable BPHZ lifts
(but for slightly different choices of noises and of integration kernels), 
for which it then follows immediately from 
existing results \cite{CH16,HS23} that they have the same limit as $\eps \to 0$.

\subsection{A priori bounds on the covariant noise}

In order to compare these processes, our first step is to approximate 
$\boldxi_{\Beps}$ (which we recall is given by $\xi_\Beps \oplus \zeta_\Beps$ with
these processes defined in \eqref{main:eq-xi-zeta-smooth-approximation}) in a way that is more
amenable to analysis via regularity structures. For this, we define an $\R^2 \otimes \Y$-valued process $\hat \boldxi_\eps$ by
setting
\begin{equ}\label{e:hat-boldxi}
\hat \boldxi_\eps (t,x)
 \eqdef 2\iint  (t-s)\rho_\eps^\prime(t-s) \nabla_x\hat\fp(s,y;t,x) \otimes \boldxi(s,y)\,\ds \,\dy\;.
\end{equ}
(Again, the time integral runs over all times here, not just positive times.)
The following bound then turns out to be sufficient to show that we can replace $\boldxi_{\Beps}$
by $\boldxi_\eps + \repr(B)\hat \boldxi_\eps$ without any risk of changing the limiting process.

\begin{proposition}\label{prop:descrNoise}
Assume that $B \in \CC^{1,-\kappa}$ for some $\kappa \in (0,\f12)$ is a solution to the Yang--Mills gradient flow 
and set  $\eta_\eps = \boldxi_{\Beps} - \boldxi_\eps - \repr(B)\hat \boldxi_\eps$
for positive times. (For
negative times, we set $\eta_\eps=0$.)
Then, $\eta_\eps$ is a centred Gaussian process satisfying for every $\nu \in [0,2]$ the bounds
\begin{equ}
\E |\scal{\eta_\eps, \phi_z^\lambda}|^2 \lesssim \eps^\nu (r \wedge 1)^{-\kappa-1} \lambda^{-\nu} + 
 \one_{r \lesssim \lambda^2+\eps^2}  \eps^\nu(\eps \vee \lambda)^{-4-\nu}\;,
\end{equ}
uniformly over $\lambda, \eps \le 1$ and $z = (r,u)$ with $r > 0$.

In particular, one has $\lim_{\eps \to 0} \eta_\eps = 0$ in probability 
in $\CC^{\alpha,\mu} \cap \CC^{\mu}$ for every $\alpha < 0$ and every $\mu < -2-\kappa$.
The same is true if we multiply $\eta_\eps$ by any deterministic (or independent
of the noise) function $U_0$.
\end{proposition}

\begin{remark}
Here, the two components of $B$ are contracted with the $\R^2$ factor of the noise $\hat \boldxi_\eps$, so that $\repr(B)\hat \boldxi_\eps$ takes values in $\Y$.
\end{remark}

\begin{proof}
Throughout this proof, we will use the notation 
$P^{(k)}_t(x) = t^{-\frac k2-1}\exp(-c |x|^2/t)$ for a constant $c$ which we allow ourselves to
take smaller and smaller from one line to the next in our calculations. 
Note that the exponents are such that $|D^k \hat\fp_t| \lesssim P_t^{(k)}$
and that, for any $\alpha > 0$ and $k \in \R$, one has the useful bound
\begin{equ}[e:usefulPkt]
\bigl(|t|^{\frac\alpha2} + |x|^{\alpha}\bigr)P^{(k)}_t(x) 
\lesssim P^{(k-\alpha)}_t(x) \;.
\end{equ}

Fix $t_0 > 0$ and $x_0 \in \R^2$ and write $q_B(t,x) = \fq_B(t_0,x_0;t,x)$
for $t > t_0$ (recall that this is an $\End(\fg^2)$-valued kernel). We also write $B^f = B(t,x) \in \fg^2$.
Recall furthermore that the action of the covariant Laplacian $\Delta_B$ on $\fg$-valued functions $f$ is given by
\begin{equ}
\Delta_B f = (\Delta - 1) f + 2 [B_i, \d_i f] + [\d_i B_i,f] + [B_i,[B_i,f]] + f\;.
\end{equ}

It then follows from the variation of constants formula that
\begin{equ}[e:varConstpB]
q_B(t,\place) = \hat\fp_{t-t_0}(\place-x_0) + 2\int_{t_0}^t \hat\fp_{t-s} [B_{s},\nabla q_B(s,\place)]\,ds + \int_{t_0}^t \hat\fp_{t-s} \bigl(G_{s} q_B(s,\place)\bigr)\,ds\;,
\end{equ}
for some function $G$ satisfying
$\|G_s\|_{L^\infty} \lesssim s^{-\frac{1+\kappa}2}$, and the same is true of 
$p_B(t,\place)$ (just for a different $G$ satisfying the same bounds).
By Proposition~\ref{prop:qB-ptwise-bds} and since $\Fnorm{F}_{[t_0,1],0}$
is bounded, uniformly over $t_0 \in (0,1]$, one has
$q_B \lesssim  P^{(0)}_{t-t_0}(\place-x_0)$,

In particular, since
$s^{-\frac{1+\kappa}2} \le t_0^{-\f{\delta+\kappa}2} |s-t_0|^{-\frac{1-\delta}2}$ for any $\delta \in [0,1]$, the last term in \eqref{e:varConstpB} is bounded by
some multiple of 
\begin{equ}
t_0^{-\f{\delta+\kappa}2}  |t-t_0|^{\frac{1+\delta}2} P^{(0)}_{t-t_0}(\place-x_0)\;.
\end{equ}
Regarding the second term in the right-hand side of \eqref{e:varConstpB},
we write it as
\begin{equs}
2\int_{t_0}^t &\hat\fp_{t-s} [B_{s},\nabla q_B(s,\place)]\,ds \\
&= 2\int_{t_0}^t \hat\fp_{t-s} [B_{s},\nabla \hat\fp_{s-t_0}(\place-x_0)]\,ds
+ 4\int_{t_0}^t \hat\fp_{t-s} \Bigl[B_{s},\nabla \int_{t_0}^s \hat\fp_{s-r} [B_{r},\nabla q_B(r,\place)]\,dr\Bigr]\,ds \\
&\quad + 2\int_{t_0}^t \hat\fp_{t-s} \Bigl[B_{s},\nabla \int_{t_0}^s \hat\fp_{s-r} \bigl(G_r q_B(r,\place)\bigr)\,dr\Bigr]\,ds\\
&=  2[B_t, (t-t_0)\nabla \hat\fp_{t-t_0}(\place-x_0)] + 2\int_{t_0}^t \hat\fp_{t-s} [B_{s}-B_t(x),\nabla \hat\fp_{s-t_0}]\,ds \\
&\quad + 4\int_{t_0}^t \hat\fp_{t-s} \Bigl[B_{s},\nabla \int_{t_0}^s \hat\fp_{s-r} [B_{r},\nabla q_B(r,\place)]\,dr\Bigr]\,ds\\
&\quad + 2\int_{t_0}^t \hat\fp_{t-s} \Bigl[B_{s},\nabla \int_{t_0}^s \hat\fp_{s-r} \bigl(G_r q_B(r,\place)\bigr)\,dr\Bigr]\,ds\\
&\eqdef  2[B_t, (t-t_0)\nabla \hat\fp_{t-t_0}(\place-x_0)] + T_1 + T_2 + T_3\;.\label{e:decompKernel}
\end{equs}

To bound $T_1$ we use the fact that, by assumption,
$|B_{s}(y)-B_t(x)| \lesssim \bigl(|x-y| + |s-t|^{\f12}\bigr) |s|^{-\f{1+\kappa}2}$, so that, as a consequence of \eqref{e:usefulPkt},
\begin{equ}[e:boundT1]
|T_1| \lesssim t_0^{-\f{\delta+\kappa}2} \int_{t_0}^{t} \int_{\R^2} P_{t-s}^{(-1)}(\place-y)  P^{(2-\delta)}_{s-t_0}(y-x_0)\,dy \,ds
\lesssim t_0^{-\f{\delta+\kappa}2} (t-t_0)^{\frac{1+\delta}2} P_{t-t_0}^{(0)}(\place-x_0)\;.
\end{equ}
(Here we need to impose $\delta > 0$ in order to avoid creating a non-integrable 
singularity at $s=t_0$.)
To bound $T_2$, we use the fact that $\nabla \hat\fp_t f = \hat\fp_t \nabla f$ to move the derivative hitting
$q_B$ outside of the integrals. We then use again the pointwise bound 
$|q_B| \lesssim P^{(0)}$ from
Proposition~\ref{prop:qB-ptwise-bds},  yielding
\begin{equs}
|T_2| 
&\lesssim 
\int_{t_0}^t P^{(0)}_{t-s} \|B_{s}\|_{L^\infty}  \int_{t_0}^s P^{(1)}_{s-r} \|D B_{r}\|_{L^\infty} P_{r-t_0}^{(0)}(\place-x_0)\,dr\,ds\\
&\quad + 
\int_{t_0}^t P^{(0)}_{t-s} \|D B_{s}\|_{L^\infty} \int_{t_0}^s P^{(1)}_{s-r} \|B_{r}\|_{L^\infty} P_{r-t_0}^{(0)}(\place-x_0)\,dr\,ds\\
&\quad + 
\int_{t_0}^t P^{(1)}_{t-s} \|B_{s}\|_{L^\infty} \int_{t_0}^s P^{(1)}_{s-r} \|B_{r}\|_{L^\infty} P_{r-t_0}^{(0)}(\place-x_0)\,dr\,ds\;.
\end{equs}
Using the fact that $\|B_{s}\|_{L^\infty} \lesssim s^{-\f\kappa2}$
and $\|DB_{r}\|_{L^\infty} \lesssim r^{-\f{1+\kappa}2} \le t_0^{-\f{\delta+\kappa}2}|r-t_0|^{-\f{1-\delta}2}$, we obtain
\begin{equs}
\int_{t_0}^t P^{(0)}_{t-s} \|B_{s}\|_{L^\infty} & \int_{t_0}^s P^{(1)}_{s-r} \|D B_{r}\|_{L^\infty} P_{r-t_0}^{(0)}(\place-x_0)\,dr\,ds \\
&\lesssim 
  t_0^{-\f{\delta+\kappa}2} P^{(0)}_{t-t_0}(\place-x_0) \int_{t_0}^t  s^{-\f\kappa2} \int_{t_0}^s |s-r|^{-\frac12}|r-t_0|^{-\f{1-\delta}2}\,dr\,ds \\
&\lesssim t_0^{-\f{\delta+\kappa}2} P^{(0)}_{t-t_0}(\place-x_0)\int_{t_0}^t (s-t_0)^{\f{\delta-\kappa}2}\,ds
\lesssim 
t_0^{-\f{\delta+\kappa}2} (t-t_0)^{\f{2-\kappa+\delta}2} P_{t-t_0}^{(0)}(\place-x_0)\;.
\end{equs}
The other two terms can be bounded in a similar way. (In the case 
of the last term we obtain the same bound as \eqref{e:boundT1}.)
Since $\kappa < 1$, this
bound is stronger than the bound \eqref{e:boundT1} on $T_1$.
The term $T_3$ satisfies the same bound as $T_2$ since $\|G_s\|_{L^\infty}$ 
satisfies the same bound as $\|D B_s\|_{L^\infty}$.

The remaining terms are the first term in the right-hand side of \eqref{e:varConstpB}
and the first term in the right-hand side of \eqref{e:decompKernel}. These
two terms are almost equal to $\boldxi_\eps$ and $[B,\hat \boldxi_\eps]$ respectively. 
The only difference is that the time integral runs over $\R_+$ rather than all of 
$\R$. 

We can therefore write
$\eta(t,x) = \int Q(s,y;t,x) \boldxi(s,y)$ for a kernel $Q$ satisfying the bound
\begin{equs}
|Q(s,y;t,x)| &\lesssim |\rho_\eps'(t-s)| P^{(0)}_{t-s}(x-y)  \Bigl(\one_{s>0} s^{-\f{\delta+\kappa}2} (t-s)^{\frac{1+\delta}2}  + \one_{s < 0} \bigl(1 + t^{-\f\kappa2}(t-s)^{\f12}\bigr)\Bigr) \\
&\lesssim \one_{t-s \in [\eps^2,2\eps^2]}\eps^{-2}  P^{(0)}_{\eps^2}(x-y)
\Bigl( \one_{s>0} s^{-\f{\delta+\kappa}2} \eps^{1+\delta} + \one_{s < 0} \bigl(1 + \eps t^{-\f\kappa2}\bigr)\Bigr)
\;.
\end{equs}
Take now a test function $\phi_z^\lambda$ at scale $\lambda$ centred around some point $z \in \R\times \T^2$ with
temporal coordinate $r$ and rewrite $\eta(\phi_z^\lambda)$ as
\begin{equ}
\eta(\phi_z^\lambda) = \boldxi(\hat \phi_z)\;,\qquad
\hat \phi_z(s,y) = \int_{\R_+\times \T^2} Q(s,y;t,x) \phi_z^\lambda(t,x)\,dt\,dx\;.
\end{equ}
The reason why the time integral is restricted to $\R_+$ is because
we did set $\eta_\eps(t,x) = 0$ for $t < 0$.
 
Since $\phi_z^\lambda(t,x) \lesssim \lambda^{-2} P^{(0)}_{\lambda^2}(x-u)\one_{|t-r| \le \lambda^2}$, it follows that one has
\begin{equ}[e:boundhatphi]
|\hat \phi_z(s,y)| \lesssim (\eps\vee \lambda)^{-2} P^{(0)}_{\eps^2 + \lambda^2}(y-u) \one_{|s-r| \le \eps^2+\lambda^2} \Bigl(\one_{s>0} s^{-\f{\delta+\kappa}2}\eps^{1+\delta}  + \one_{-2\eps^2<s<0}\Bigr) \;.
\end{equ}
Since \eqref{e:boundhatphi} is in particular true for both
$\delta = 1$ and for $\delta = \frac12$ (say), we obtain the bound
\begin{equ}
 \|\hat \phi_z\|_{L^2}^2 \lesssim (\eps \vee \lambda)^{-4} \Bigl((r \vee \lambda^2 \vee \eps^2)^{-1-\kappa}\eps^{4}  
 + \one_{r \lesssim \lambda^2+\eps^2}(\eps \vee \lambda)^{-2}\eps^{2} \Bigr)\;.
\end{equ}
As a consequence, we have for every $\nu \in [0,2]$ the bound
\begin{equ}
\E |\eta(\phi_z^\lambda)|^2 = \|\hat \phi_z\|_{L^2}^2 
\lesssim \eps^\nu (r \vee \lambda^2)^{-\kappa-1} \lambda^{-\nu}
+ \one_{r \lesssim \lambda^2+\eps^2}  \eps^\nu(\eps \vee \lambda)^{-4-\nu}\;,
\end{equ}
uniformly over $\eps, \lambda \in (0,1]$. 
Recalling that one has 
equivalence of moments by Gaussianity of $\eta$, a Kolmogorov-type criterion 
yields the last claim by choosing $\nu$ small enough.
\end{proof}

\subsection{Formulation in the vector-valued regularity structures framework}

Just as in \cite{CCHS22}, we first write the equation for $g$ in a ``linearised'' way by
introducing variables $h \in \fg^2$ and $U \in \End(\fg) \oplus \End(\V)$ given by $h = dg\,g^{-1}$ and 
$U = \Ad_g \oplus \rep(g)$. We will again write $\Repr$ for this representation and $\repr$ for its
linearisation.
Note that we will also without further ado view $U$ as an element of $\End(\Y)$ (just duplicate the component
in $\End(\fg)$), $\End(\fg^2)$ (duplicate the component
in $\End(\fg)$ and discard the one in $\End(\V)$), or $\End(\fg)$ (extract the component in $\End(\fg)$).
In particular, we will simply write $g \act B = UB-h$. 
The same calculation as in \cite[Lem.~6.3]{CCHS22+} and \cite[Lem.~7.2]{CCHS22} then shows that 
if $g$ solves an equation of the type
\begin{equ}
\d_t  g\,  g^{-1} = - d^* \bigl(d g\, g^{-1}\bigr) + H\;,
\end{equ}
where $H$ is any smooth $\fg$-valued process, then $h$ and $U$ solve
\begin{equ}[e:evolhU]
\d_t h =  \Delta h - [h,dh] + \covd_{-h} H\;,\qquad 
\d_t U = \Delta U - \repr(h^\ell)\repr(h_\ell)\circ U + \repr(H)\circ U\;,
\end{equ}
where we used the shorthand $[h,dh]_i = [h_\ell,\d_\ell h_i]$ with $i \in \{1,2\}$ and summation over $\ell$ implicit.
Furthermore, for any smooth $\fg$-valued one-form $B$, one has the identity
$[U,d^*]B = [h,UB]$.

Writing $\HH \colon \Y \times \fg^2 \times \End(\Y) \times \fg^2 \times \fg^4 \to \fg$ for the function
\begin{equ}[e:defHeps]
\HH(Y,h,U,B,DB) =  U d^*B - [h,U B] - [Y_A, U B-h]\;,
\end{equ}
as well as $\HH_\eps^{(i)}$ for the process $\HH(Y_\eps^{(i)},h_\eps^{(i)},U_\eps^{(i)},B,DB)$ (with the values $i$ referring
to the three settings described in Section~\ref{sec:threeProcesses}),
we can then rewrite the equations for the processes $g_\eps^{(i)}$ as
\begin{equs}[e:equationUh]
\d_t h_\eps^{(i)} &=  \Delta h_\eps^{(i)} - [h_\eps^{(i)}, dh_\eps^{(i)}] + \covd_{-h_\eps^{(i)}} \HH_\eps^{(i)}\;, \\ 
\d_t U_\eps^{(i)} &= \Delta U_\eps^{(i)} - \repr(h_\eps^{(i)})^2\circ U_\eps^{(i)} + \repr(\HH_\eps^{(i)})\circ U_\eps^{(i)}\;.
\end{equs}
We rewrite the equations for $Y_\eps^{(i)}$, $h_\eps^{(i)}$ and $U_\eps^{(i)}$ in mild form in a slightly unconventional 
way by breaking out the direct effect
of the noise $\boldxi_{\eps}$ into a separate component. This will allow us to treat the cases $i \in \{1,2,3\}$ 
in a unified setting by incorporating the convolution with $\moll$ arising in the noise term for
$Y_\eps^{(1)}$ into the integration operator. 

Ignoring for the moment the presence of the renormalisation terms involving $\renormY = \renormA \oplus \renormphi$ 
which are different
for each equation (these terms are of course absolutely crucial and will reappear later), we note that 
one can rewrite the equations satisfied by the three 
sets of processes under consideration in the following unified way. We set
\begin{equs}[e:dataSystem]
\fp_\eps^{(1)} &= \moll \star \hat\fp\;,&\qquad
\boldxi_\eps^{(1)} &= \boldxi\;,&\qquad
\hat\boldxi_\eps^{(1)} &= 0\;,&\qquad
\eta_\eps^{(1)} &= 0\;,\\
\fp_\eps^{(2)} &= \hat\fp\;,&\qquad
\boldxi_\eps^{(2)} &= \boldxi_\eps\;,&\qquad
\hat\boldxi_\eps^{(2)} &= \hat\boldxi_\eps\;,&\qquad
\eta_\eps^{(2)} &= \eta_\eps\;,\\
\fp_\eps^{(3)} &= \hat\fp\;,&\qquad
\boldxi_\eps^{(3)} &= \boldxi_\eps\;,&\qquad
\hat\boldxi_\eps^{(3)} &= 0\;,&\qquad
\eta_\eps^{(3)} &= 0\;,
\end{equs}
and we write
$Y_\eps^{(\place)} = Y_\eps^{(\place,1)} + Y_\eps^{(\place,2)}$ and similarly for $h_\eps^{(\place)}$.
By the mild formulation, one then has
\begin{equs}[e:bigsystem]
Y_\eps^{(\place,1)} &= \hat\fp \star \one_+ \bigl(Y_\eps^{(\place)} + \tril(Y_\eps^{(\place)}) + \hat \bil(Y_\eps^{(\place)},DY_\eps^{(\place)}) + U_\eps^{(\place)} \repr(B) \hat\boldxi_\eps^{(\place)} + U_\eps^{(\place)} \eta_\eps^{(\place)}\bigr)\\
&\quad  + \hat\fp \bigl(Y_{\eps,0}^{(\place)}- Y_{\eps}^{(\place,2)}(0)\bigr)\;,\\
Y_\eps^{(\place,2)} &= \fp_\eps^{(\place)} \star \bigl(U_\eps^{(\place)} \boldxi_\eps^{(\place)}\bigr)\;,\\
h_\eps^{(\place,1)} &= \nabla\hat\fp \star \one_+  \bigl(\HH(Y_\eps^{(\place)},h_\eps^{(\place)},U_\eps^{(\place)},B,DB)
	\bigr) + \hat\fp h_{\eps,0}^{(\place)}\;,\\
h_\eps^{(\place,2)} &= \hat\fp \star \one_+  \bigl(h_\eps^{(\place)}	- [h_\eps^{(\place)}, dh_\eps^{(\place)} + \HH(Y_\eps^{(\place)},h_\eps^{(\place)},U_\eps^{(\place)},B,DB)]\bigr)\;,\\
U_\eps^{(\place)} &= \hat\fp \star \one_+  \bigl(G_U(h_\eps^{(\place)},U_\eps^{(\place)},B,DB) + T_U(Y_\eps^{(\place)}, U_\eps^{(\place)} B-h_\eps^{(\place)}, U_\eps^{(\place)})\bigr) + \hat\fp U_{\eps,0}^{(\place)}\;,
\end{equs}
where the trilinear map
$T_U$ is the map obtained by inserting the last term of \eqref{e:defHeps} into
the evolution equation \eqref{e:evolhU} for $U$ while $G_U$ contains the remaining terms.
The presence of the additional term $Y_\eps^{(\place)}$ in the right-hand side
for $Y_\eps^{(\place,1)}$ (and similarly for $h_\eps^{(\place)}$ and $U_\eps^{(\place)}$) 
comes from the fact that $\hat\fp$ denotes the massive heat kernel
with unit mass. We denoted the initial conditions for the three processes 
by $Y_{\eps,0}^{(\place)}$, etc and we extend $U_\eps^{(\place)}$ to negative times
by setting it equal to its initial condition $U_{\eps,0}^{(\place)}$ there.

Note that all of the nonlinearities appearing in \eqref{e:bigsystem} are composition operators. In particular,
none of them involve the derivatives of the solution, except for the terms $\hat\bil$ and $[h_\eps, dh_\eps]$, where
this is explicitly indicated in the notation.
The reason for choosing $G_U$ in such a way that the right-hand side
of the equation for $U_\eps$ breaks up in this particular way is that $G_U$ takes values in 
$\cC^\alpha$ for every $\alpha < 1$ (and therefore yields a contribution of regularity
almost $\cC^3$ to $U_\eps$) while the remaining term only takes values in 
$\cC^\alpha$ for $\alpha < 0$ and will therefore be the only term that contributes to the 
renormalisation of the term  $U_\eps \boldxi_{\eps}$. Similarly, all the products appearing
in the right-hand sides for $h_\eps^{(\place,1)}$ and $h_\eps^{(\place,2)}$ are classically well-defined and do not require any renormalisation.

We can fit \eqref{e:bigsystem} into the framework developed in \cite{BCCH}, in particular the
version for vector-valued SPDEs developed in \cite[Sec.~5]{CCHS22}. For this, we are going to
index the components of the system \eqref{e:bigsystem} by
\begin{equ}
\fL_+ = \{\fy_1, \fy_2, \fh_1, \fh_2, \fu, \fb\}\;,\quad \fL_- = \{\bullet,\wedge\}\;,
\end{equ}
where $\fy_i$ represents $Y_\eps^{(\place,i)}$, $\fh_i$ represents $h_\eps^{(\place,i)}$, $\fu$ represents $U_\eps^{(\place)}$, 
$\fb$ represents the fixed background connection $B$, $\bullet$ represents the driving noise $\boldxi_\eps^{(\place)}$
and $\wedge$ represents the noise $\hat \boldxi_\eps^{(\place)}$. The noise $\eta_\eps^{(\place)}$
is sufficiently regular so that it will be treated classically (i.e.\ via its Taylor lift).
The corresponding target spaces $W_\ft$ and kernel spaces $\CK_\ft$ as in \cite[Sec.~5.6]{CCHS22}
are given by
\begin{equs}[0][e:spaces2]
W_{\fy_1} = W_{\fy_2} = W_{\bullet} = \Y\;,\quad  W_{\fh_1} = W_{\fh_2} = W_{\fb} = \fg^2\;,\quad W_\wedge = \R^2 \otimes \Y \;,\\
W_{\fu} = \End(\fg)\oplus \End(\V)\;, \qquad \CK_{\fy_1} = \CK_{\fy_2} = \CK_{\fh_2} = \CK_{\fb} =\CK_{\fu} = \R\;,\qquad  \CK_{\fh_1} = \R^2\;.
\end{equs}
(The reason why the kernel space $\CK_{\fh_1}$ is given by $\R^2$ is because the corresponding
kernel is the gradient of the heat kernel.)
We fix some $\bar\kappa > 0$ small enough ($\bar\kappa < \kappa/2$ is sufficient) 
and we postulate that the integration operators associated to the 
elements of $\fL_+$ have degree $2-\bar\kappa$ except for the one associated to $\fh_1$ which has degree $1-\bar\kappa$, while the
noise $\bullet$ has regularity $-2-\bar\kappa$ and the
noise $\wedge$ has regularity $-1-\bar\kappa$. The reason why we do not simply decree
that the integration operators have degree $2$ is that this gives us the wiggle room necessary to 
compare the corresponding integration kernels.
We take as our rule the minimal
rule allowing to formulate the right-hand side of \eqref{e:bigsystem}, analogously to
\cite[Eq.~6.10]{CCHS22+}.

Recall that the construction there yields a collection of symbols usually denoted either
as formal expressions obtained from the noise symbols $\{\Xi_\bullet, \Xi_\wedge\}$
(which we henceforth denote by $\{\Xi,\hat \Xi\}$ to suggest that they will correspond
to $\boldxi$ and $\hat\boldxi$) and the ``Taylor monomials''
$\mathbf{X}^k$ by applications of products and ``integration operators'' $\CI_\ft$, or
as the corresponding syntactic tree. Each symbol $\tau$ has a natural degree and, as before,
is associated to the vector space $\CT[\tau]$ obtained by associating tensor factors $\CK_\ft^*$ \slash $W_\fl^*$ to
each copy of $\CI_\ft$ \slash $\Xi_\fl$ in $\tau$ \dash see \cite[Sec.~5]{CCHS22} for details on how we keep track 
of which copy is which in a systematic way. 

We then introduce \textit{six} modelled distributions $\CY_1$, $\CY_2$, $\CH_1$, $\CH_2$, $\CU$, and $\CB$ representing 
the components of \eqref{e:bigsystem} and $B$ respectively, we write $\CY = \CY_1 + \CY_2$, $\CH=\CH_1+\CH_2$
 and, for any choices $\eta \in C^\infty(\R_+\times\T^2,\Y)$, $C \in L(\Y\oplus\fg^2\oplus\fg^2,\Y)$, 
 $\massA,\hat\massA \in \End(\fg^2)$ and 
 $\hat\mass \in \End(\V)$, we consider the fixed point problem
\begin{equs}
\CY_1 &= \CP_{\fy_1}\one_+ \bigl(\tril(\CY) + \hat \bil(\CY,D\CY) + \CY - \hat \mass\CY_\higgs + \massA \CY_A  + \hat \massA (\CY_A + \CH)
 + \CU \repr(\CB)\hat \Xi\bigr)\\
 &\quad + \CJ \bigl(\hat\fp\star\one_+((\CR \CU) \eta) + \hat\fp Y_0\bigr) + \CP_{\fy_1}\one_+ \bigl(C(\CY,\CH,\CU\CB) \bigr) \\
\CY_2 &= \CP_{\fy_2} \bigl(\CU \Xi\bigr)\\
\CH_1 &= \CP_{\fh_1}\one_+  \bigl(\HH(\CY,\CH,\CU,\CB,D\CB)
	\bigr) + \CJ \bigl(\hat\fp h_0\bigr) \label{e:fpabstract}\\
\CH_2 &= \CP_{\fh_2}\one_+  \bigl(\CH- [\CH, d\CH + \HH(\CY,\CH,\CU,\CB,D\CB)]\bigr) \\
\CU &= \CP_{\fu}\one_+  \bigl(G_U(\CH,\CU,\CB,D\CB) - \repr([\CY, \CU \CB-\CH])\circ \CU\bigr) + \CJ \bigl(\hat\fp U_0\bigr)\\
\CB &= \CJ B
\end{equs}
where $\CP_\ft$ denotes the integration operator associated to
the component $\ft$, $\CR$ is the reconstruction operator, and the nonlinear functions $\tril$, $\hat\bil$, $\HH$, and $G_U$
have been lifted to the regularity structure generated from the data above in the usual
way (see \cite[Sec.~5.8.2]{CCHS22} for the vector-valued case relevant here).
Similarly to the notation previously used in \eqref{e:defiotalollipop}, 
$\CY_A$ denotes the component of $\CY$ in $\fg^2$ and $\CY_\higgs$
denotes its component in $\V$. 
In view of \eqref{e:bigsystem}, what we call $Y_0$ here is not going to be the initial 
condition for $Y_{\eps}^{(\place)}$, but rather $Y_{\eps,0}^{(\place)}- Y_{\eps}^{(\place,2)}(0)$.

The maps $\hat \mass, \massA, \hat \massA$ will be chosen as suitable counterterms
that are independent of $\eps$ and that give us some additional 
freedom to deal with renormalisation. 
The linear map $C = C_\eps$ will always be chosen such that
\begin{equ}[e:limCzero]
\lim_{\eps \to 0} C_\eps = 0\;,
\end{equ}
and therefore does not affect the limit as $\eps \to 0$.
As such, it should be considered as an unimportant tweak that merely serves
to guarantee exact identification of the solution to this fixed point
problem at finite $\eps > 0$.
As before, $\CJ F$ denotes the Taylor jet up to parabolic order $2$ of a sufficiently 
smooth function $F$, see \eqref{e:TaylorJet}.
We will also henceforth use the notation $\hat \CD^{\gamma,\eta} \subset \CD^{\gamma,\eta}$ for the set of those singular modelled distributions that actually belong
to $\CD^\eta$ and vanish at negative times, see for example \cite[Appendix~A]{CCHS22+}.

We have the following well-posedness result for solutions to \eqref{e:fpabstract}. Here,
for $\beta < \alpha < 0$, we write $\CC^{\alpha,\beta}$ for the
space of distributions $\zeta$ such that
\begin{equ}
|\zeta(\phi_{(t,x)}^\lambda)| \le (\sqrt{|t|}\wedge 1)^{\beta-\alpha} \lambda^{\alpha}\;,
\end{equ}
uniformly over $\lambda \le (\sqrt{|t|}\wedge 1)/2$ and locally uniformly
over $(t,x) \in \R_+\times \T^2$. 
(The convention regarding the power of $|t|$ is
chosen in order to make this consistent with \cite[Def.~6.2]{H14}.) We will make use of the fact that if
$\zeta \in \CC^{\alpha,\beta}$ and either $\beta > -2$ or $\beta < -2$ but one
also has $\zeta \in \CC^{\beta}$, then 
$\fp\star \zeta \in \CC^{\alpha+2,\beta+2}$. (If $\alpha +2 > 0$, then this should
be interpreted in the sense of \cite[Def.~6.2]{H14} for the polynomial regularity
structure.)

\begin{proposition}\label{prop:equalProcesses}
Assume that, uniformly over $\eps \in (0,1]$, one has $Y_0 \in \CC_x^{-\kappa}$, $h_0 \in \CC_x^{1-\kappa}$, $U_0 \in \CC_x^{2-\kappa}$, $B \in \CC^{2-\kappa,-\kappa}$,
$\eta \in \CC^{-\kappa,-2-2\kappa}$,
and $\one_+(U_0 \eta) \in \CC^{-\kappa,-2-2\kappa}\cap \CC^{-2-2\kappa}$.

For any admissible model 
with non-anticipative integration kernels, there exists a time $T>0$ such that
the fixed point problem \eqref{e:fpabstract} admits a unique solution up to time $T$
with $\CY_1,\CY_2 \in \CD^{\gamma,\eta}$, $\CH_1,\CH_2 \in \CD^{\gamma,1+\eta}$
and $\CU \in \CD^{1+\gamma,2+\eta}$ with $\gamma = 1 + 2\kappa$ and $\eta = -2\kappa$. 
That solution furthermore depends continuously on all this data in the topologies 
indicated above.
\end{proposition}

\begin{proof}
It appears at first sight that $\eta$ is too singular in order to get uniform bounds
near the origin. We can however exploit our assumptions 
in the following
way. We rewrite the term $(\CR \CU) \eta$ as
\begin{equ}
(\CR \CU) \eta = \bigl((\CR \CU) - U_0\bigr)\eta + U_0\eta\;,
\end{equ}
and we then exploit the fact that
$\one_+(U_0 \eta) \in \CC^{-\kappa,-2-2\kappa}\cap \CC^{-2-2\kappa}$ by assumption,
so that its convolution with $\hat \fp$ belongs to $\CC^{2-2\kappa,-2\kappa}$.
It remains to note that since we extend $\CU$ to be constant and equal to 
$\CJ U_0$ for negative times and we have continuity at time $0$, the fact 
that $\CU \in \CD^{1+\gamma,2+\eta}$ and $U_0 \in \CC_x^{2-\kappa}$
implies that one actually has
\begin{equ}
\CU - \CJ U_0 \in \hat\CD^{2+\eta,2+\eta}\;,
\end{equ}
in the notations of \cite[Appendix~A]{CCHS22+}. 
This allows us to apply \cite[Lem~A.3]{CCHS22+}, which in turn allows us to gain 
enough integrability near time $0$ so that the rest of the proof proceeds in the 
usual way.
\end{proof}

Our strategy is then to find three $\eps$-dependent random models $\hat Z^{(\eps,i)}$ 
for our regularity structure such that the corresponding solutions to \eqref{e:fpabstract}
coincide precisely (after reconstruction) with the processes $(Y_\eps^{(i)},h_\eps^{(i)}, U_\eps^{(i)})$ defined in Section~\ref{sec:threeProcesses}. (In the case $i=3$,
this turns out to only be the case when $B=0$, which is sufficient for our purpose.)
If one can then show that $\lim_{\eps \to 0} \hat Z^{(\eps,i)} = \hat Z$ for some random model 
$\hat Z$ that is independent of $i$, this shows that the three processes under consideration have indeed
the same limit (at least on short times).
It is important in this argument that our choice of $\hat\mass$, $\massA$, and $\hat\massA$ doesn't depend on $i$, even though $C_\eps$ is allowed to depend on it, provided that \eqref{e:limCzero} holds.

We will consider three natural models for this regularity structure, matching the data
already mentioned in \eqref{e:dataSystem}.

\begin{definition}\label{def:Models}
The models $Z^{(\eps,i)}$ are defined as follows.
\begin{enumerate}
\item The model $Z^{(\eps,1)} = (\Pi^{(\eps,1)},\Gamma^{(\eps,1)})$ is obtained by the canonical lift associating
$\boldxi$ to the noise symbol $\Xi$, $0$ to the noise symbol $\hat \Xi$, and the heat kernel to each of the integration
kernels, except for $\fh_1$ which we associate to the gradient of the heat kernel 
and $\fy_2$ which we associate to the kernel $\moll \star \hat\fp$. One furthermore takes $\eta = 0$.
\item The model $Z^{(\eps,2)} = (\Pi^{(\eps,2)},\Gamma^{(\eps,2)})$ is obtained by the canonical lift associating
$\boldxi_\eps$ to the noise symbol $\Xi$, $\hat \boldxi_\eps$ to the noise symbol $\hat \Xi$, and the heat kernel to each 
of the integration kernels, except for $\fh_1$ which we associate to the gradient of the heat kernel. One furthermore takes $\eta = \eta_\eps$. Note that, as a consequence
of Proposition~\ref{prop:descrNoise}, this choice satisfies the assumptions of Proposition~\ref{prop:equalProcesses} for any 
deterministic choice of initial data $U_0$. 
\item The model $Z^{(\eps,3)} = (\Pi^{(\eps,3)},\Gamma^{(\eps,3)})$ is obtained by the canonical lift associating
$\boldxi_\eps$ to the noise symbol $\Xi$, $0$ to the noise symbol $\hat \Xi$, and the heat kernel to each 
of the integration kernels, except for $\fh_1$ which we associate to the gradient of the heat kernel. One furthermore takes $\eta = 0$.
\end{enumerate}
\end{definition}

With this choice, solutions to \eqref{e:fpabstract} for the model $Z^{(\eps,i)}$ reconstruct to 
the processes that solve ``almost'' the same equations as the tuples 
$( Y_\eps^{(\place)}, h_\eps^{(\place)}, U_\eps^{(\place)})$,
except that the terms involving $\renormY$ are replaced by the ``wrong'' counterterm
\begin{equ}
\bigl( \massA A_\eps^{(\place)}  + \hat\massA (A_\eps^{(\place)} + h_\eps^{(\place)})\bigr) 
\oplus (- \hat\mass \higgs_\eps^{(\place)}) + C_\eps(Y_\eps^{(\place)}, h_\eps^{(\place)}, U_\eps^{(\place)} B) \;.
\end{equ} 

\begin{remark}
Note that even though $\boldxi$ is
quite irregular, it still admits a canonical lift when $\eps > 0$ in the case of the model $Z^{(\eps,1)}$.
The reason is that, thanks to the fact that the integration operator $\CP_{\fy_2}$ is associated to the 
smooth kernel $\fp_\eps^{(1)}$ in this case, $\boldxi$ only ever gets 
multiplied by smooth functions when performing that lift. 
\end{remark}

Since solutions to \eqref{e:fpabstract} are continuous with respect to the model and the additional data
(namely the initial conditions, the process $\eta_\eps$, and the integration kernel $\fp_\eps^{(\place)}$) 
and since $\fp_\eps^{(\place)} \to \hat\fp$ and $\eta_\eps \to 0$ by Proposition~\ref{prop:descrNoise}, 
it suffices to show that the three canonical lifts we just described can be renormalised in 
such a way that they converge to the same limiting model, and that this renormalisation
has the effect of yielding the ``correct'' counterterms proportional to $\renormY$ appearing in
the equations for $Y_\eps^{(\place)}$. It turns out that we have set up everything in 
such a way that the BPHZ renormalisation does have exactly that property.

It therefore remains to compute the BPHZ renormalisation constants for the models $Z^{(\place)}_\eps$,
as well as the counterterms that these give rise to in the right-hand side for $Y$.
An important remark is that if we restrict ourselves to the symbols naturally
arising from the system of equations \eqref{e:fpabstract}, then there are only
finitely many symbols of negative degree. 
In fact, besides the symbols of the type $\X^k\Xi$ with $|k|\le 2$ and using the graphical 
notation of \cite[Sec.~6.2.3]{CCHS22}, 
the symbols of negative degree appearing in this construction are given by
\begin{equ}[e:symbolList]
\<IXi>, \enskip
\<I'Xi>, \enskip
\<I'Xi2>, \enskip
\<I'XiX>,\enskip
\<IXi^2>, \enskip
\<IXi^3>, \enskip
\<IXiI'Xi>, \enskip
\<I'XiIXi2>, \enskip
\<IXiI'Xi2>, \enskip
\<IXiI'XXi>, \enskip
\<IXXiI'Xi>, \enskip
\<I'[IXiI'Xi]>, \enskip
\<I[IXiI'Xi]I'Xi>, \enskip
\<IXiI'[IXiI'Xi]>, \enskip
\<I[I'Xi]I'Xi>, \enskip
\<IXiI'[I'Xi]>, \enskip
\<I'[I'Xi]>, \enskip
\<XiIIXi>\;.
\end{equ}
These should be read as syntactic trees, with round nodes $\<Xi>$ representing $\Xi$, square nodes $\<Xi2>$ representing $\hat \Xi$, thin lines
representing one of the integration operators $\CI_\ft$ with $\ft \neq \fh_1$ (these are associated to
integration with respect to the heat kernel or mollifications thereof), thick lines representing  $\CI_{\fh_1}$
or $\CI_{\ft}'$ (these are associated to
integration with respect to the gradient of the heat kernel or mollifications thereof), and crossed circles 
$\X\Xi$. Each of these trees yields a BPHZ renormalisation constant (these constants are really elements of the vector spaces
$\CT[\tau]^*$,
but we will see that in all cases they either vanish or are described by a single real number) 
and a corresponding counterterm. Section~\ref{sec:renConstants} is devoted to the calculation of these constants
while the counterterms themselves are computed in Section~\ref{sec:renCounterterms}.

\begin{remark}\label{rem:determined}
Since the lines can stand for more than
one possible integration operator, each of these symbols potentially represents
a collection of symbols $\tau$ indexing the subspaces $\CT[\tau]$ of the structure space $\CT$.
There is however much less ambiguity than that and the types of the lines appearing in
the symbols \eqref{e:symbolList} are in fact uniquely determined by the shape of the symbol. 
Since the noise $\<Xi>=\Xi$ only appears on the right-hand
side of the component $\fy_2$ of the system \eqref{e:fpabstract}, every line that terminates with a noise symbol
must correspond to either $\CI_{\fy_2}$ (if it is thin) or $\CI_{\fy_2}'$ (if it is thick). A similar
fact holds for $\<Xi2> = \hat \Xi$ which only appears on the right-hand side of $\fy_1$.
This
already determines the types of the lines appearing in the first eleven symbols of \eqref{e:symbolList}.
Since the derivative of $\CY_2$ only appears on the right-hand side of the component $\fy_1$, this implies that the 
remaining line in the next six symbols is always of type $\CI_{\fy_1}$ or $\CI_{\fy_1}'$, depending on whether it is
thin or thick. This leaves
the lower line of the last symbol. Since the only component that multiplies the noise $\Xi$ in the right-hand side
of \eqref{e:fpabstract} is $\CU$, this line is necessarily of type $\fu$.
\end{remark}

\subsection{Expressions for the renormalisation constants}
\label{sec:renConstants}

Recall the three canonical models $Z^{(\eps,i)}$ defined in 
Definition~\ref{def:Models} and
 write $\PPi^{(\eps,i)}$ for the corresponding uncentered models
(defined as in \cite[Sec.~6.2]{BHZ19}), we furthermore set
$\ell^{(\eps,i)}(\tau) = \E \bigl(\PPi^{(\eps,i)}\tau\bigr)(0)$.

\begin{remark}
Recall that the BPHZ renormalisation of a model is defined in \cite[Thm~6.18]{BHZ19}
as the unique model of the type
\begin{equ}
\hat \PPi^{(\eps,i)}\tau = \bigl(g^{(\eps,i)} \otimes \PPi^{(\eps,i)}\bigr)\Deltam \tau\;,
\end{equ}
such that $\E \bigl(\hat \PPi^{(\eps,i)}\tau\bigr)(z) = 0$ for every $z \in \R^3$.
(In particular the theorem states that such a model exists if $\PPi^{(\eps,i)}$ is the canonical
lift of a stationary noise, which is the case here.) Here, $g^{(\eps,i)}$ is a character on the 
free algebra generated by the regularity structure $\CT$ and $\Deltam$ is the operator appearing in
\cite[Cor.~6.37]{BHZ19}.

Writing $\CA_-$ for the negative twisted antipode
introduced in \cite[Prop.~6.6]{BHZ19}, one finds that for the vast majority of the symbols appearing in 
this article, one has $\CA_- \tau = -\tau$.
The only symbols where this is not the case are those involving a factor $X$, namely
$\<IXiI'XXi>$ and $\<IXXiI'Xi>$, but for these symbols $\CA_- \tau$ differs from $-\tau$ by a 
term that has a factor $X$.
As a consequence, inspecting the identity for $g^{(\eps,i)}$ given in \cite[Thm~6.18]{BHZ19},
we see that the BPHZ model is simply given by setting
\begin{equ}[e:BPHZantipode]
g^{(\eps,i)}(\tau) = - \ell^{(\eps,i)}(\tau)\;,
\end{equ}
for $\tau \in \CT$.
\end{remark}

The goal of this subsection is to derive somewhat explicit formulae for the renormalisation 
constants $\ell^{(\eps,i)}(\tau)$ and in particular to exhibit relations between these constants
(in the limit $\eps \to 0$) and the form of their dependence on the arbitrary cutoff function $\rho$
used in the regularisation \eqref{e:mollifNoise} of our driving noise.
In order to keep notations as compact as possible, we will use throughout this section the following convenient 
convention which is similar to the one previously used in \cite{Etienne} (see the formulation of Theorem~1.2 therein).
Whenever $\tau$ is depicted by
a symbol (drawn in blue), the same symbol drawn in red is used to denote $\ell^{(\eps,2)}(\tau)$, so one
has for example
\begin{equ}\label{e:IXi^2r}
\<IXi^2r> = \E \bigl(\PPi^{(\eps,2)} \<IXi^2>\bigr)(0) = \sscal{\hat\rho_\eps \star \hat\fp,\hat\rho_\eps \star \hat\fp}\;, \end{equ}
where we defined
\begin{equ}
\sscal{f,g} \eqdef \sum_{n \in \Z^2} \int_{\R}\int_{\R^2} f(t,x) g(t,x+2\pi n)\,dx\,dt
= \scal{f,g} +  \sum_{n \neq 0} \int_{\R}\int_{\R^2} f(t,x) g(t,x+2\pi n)\,dx\,dt\;.
\end{equ}
This integral converges around $0$ thanks to the cutoff $\rho_\eps$ and at infinity thanks to the
fact that we added a mass in the definition of $\hat\fp$. The sum over $n$ converges thanks to the 
strong localisation properties of the heat kernel.
If we want to consider $\ell^{(\eps,i)}(\tau)$ for $i \neq 2$ instead, we add the corresponding superscript,
writing for example $\<IXi^2r>^{(1)}$ (which happens to coincide with $\<IXi^2r>$).
The reason why we take model $2$ as our reference is that this has the largest number of nonzero constants
and the renormalisation constants for the other two models can all be written in terms of those of model $2$.

Before we proceed, we would like to address the fact that each of the symbols $\tau$
from \eqref{e:symbolList} really corresponds to a whole subspace $\CT[\tau]$ of the regularity structure $\CT$.
As in \cite[Sec.~5.6]{CCHS22}, these subspaces are suitable finite tensor products of the type $\bigotimes_{a} V_a^*$ 
where $a$ runs over the various nodes \slash edges of $\tau$ and $V_a$ is the corresponding target or
kernel space of \eqref{e:spaces2}.
This shows that one should in principle interpret $\ell^{(\eps,i)}(\tau)$ as an element
of the corresponding tensor product $\bigotimes_{a} V_a$. In order to keep things simple in this section,
we do however want to reduce this to a scalar value.
In order to do this, we only consider the case
$\Y = \fg = \R$ so that the spaces $\CT[\tau]$ are of the type $(\R^2)^{\otimes k}$, where
each appearance of a derivative of the heat kernel (thick line), factor $X$ or noise $\hat \Xi = \<Xi2>$ in $\tau$
contributes one tensor factor $\R^2$. Inspection of \eqref{e:symbolList} shows that in all cases of interest to us 
we have $k \in \{0,1,2\}$. One furthermore has $\ell^{(\eps,i)}(\tau) = 0$ when $k=1$.
In the remaining case $k=2$, viewing $\ell^{(\eps,i)}(\tau) = (\ell^{(\eps,i)}(\tau))_{k\ell}$ with $k,\ell \in [2]$ 
as an element of $\R^2 \otimes \R^2$, it follows from 
basic symmetry considerations that these expressions are always proportional to the identity matrix, so we consider them
again as scalars. 
In Section~\ref{sec:renCounterterms} we will see how these constants are then used
to determine the counterterms, also in the case of arbitrary choices for 
$\V$ and $\fg$ (and therefore also for $\Y$).

There are a number of exact relations between renormalisation constants for our specific regularisation.
In fact we will see that, modulo errors that converge to $0$ as $\eps \to 0$, they can all be expressed 
in terms of $\<IXi^2r> \sim \log \eps$ and $\<XiIIXir> \sim 1$. As a consequence of the Gaussianity 
of the noise as well as the spatial symmetry of its covariance, one has
\begin{equ}[e:vanishSymmetry]
\<IXir> = \<I'Xir> = \<I'Xi2r> = \<I'XiXr> = \<I'[I'Xi]r>
= \<IXi^3r> = \<I[IXiI'Xi]I'Xir> = \<IXiI'[IXiI'Xi]r> = \<I'[IXiI'Xi]r> = \<IXiI'Xir> = 0\;,
\end{equ}
(and similarly for $\<IXir>^{(1)}$, $\<IXir>^{(3)}$, etc), so that we only need to consider the 
renormalisation constants associated to the symbols
\begin{equ}[e:symbolListNontrivial]
\<IXi^2>, \enskip
\<I'XiIXi2>, \enskip
\<IXiI'Xi2>, \enskip
\<IXiI'XXi>, \enskip
\<IXXiI'Xi>, \enskip
\<I[I'Xi]I'Xi>, \enskip
\<IXiI'[I'Xi]>, \enskip
\<XiIIXi>\;.
\end{equ}
In order to identify the limit $\eps \to 0$ of a renormalised equation, it is sufficient
to determine the values of the renormalisation constants up to an error that converges to $0$
as $\eps \to 0$. Given two $\eps$-dependent constants $A_\eps$ and $B_\eps$, we
will therefore write $A_\eps \eqeps B_\eps$ as a shorthand for 
$\lim_{\eps \to 0} (A_\eps - B_\eps) = 0$. 

The main result of this section can then be summarised as follows. 
\begin{proposition}\label{prop:relConstants}
There exist constants $\Crho$ and $\hCrho$ depending only on the choice of $\rho$,
as well as constants $\discrep$ and $\hat \discrep$ depending only on the choice
of mass in the heat kernel and on the size of the torus, such that
the BPHZ renormalisation constants associated to the symbols
\eqref{e:symbolListNontrivial} are given by
\begin{equs}
\<IXi^2r> &\eqeps \frac{1}{4\pi}\abs{\log\eps} + \hCrho + \hat \discrep\;,
\quad&\quad \<IXi^2r>^{(1)} &= \<IXi^2r>\;,
\quad&\quad \<IXi^2r>^{(3)} &= \<IXi^2r>\;, \\
\<XiIIXir> &\eqeps \Crho\;,
\quad&\quad \<XiIIXir>^{(1)} &= 0\;,
\quad&\quad \<XiIIXir>^{(3)} &= \<XiIIXir>\;, \\
\<I[I'Xi]I'Xir> &\eqeps \frac14\<IXi^2r> + \frac14\<XiIIXir> - \frac\discrep4\;,
\quad&\quad \<I[I'Xi]I'Xir>^{(1)} &= \<I[I'Xi]I'Xir>\;, 
\quad&\quad \<I[I'Xi]I'Xir>^{(3)} &= \<I[I'Xi]I'Xir>\;, \\
\<IXiI'[I'Xi]r> &\eqeps -\frac14\<IXi^2r> - \frac14\<XiIIXir> + \frac\discrep4\;, 
\quad&\quad \<IXiI'[I'Xi]r>^{(1)} &= \<IXiI'[I'Xi]r>\;, 
\quad&\quad \<IXiI'[I'Xi]r>^{(3)} &= \<IXiI'[I'Xi]r>\;, \\
\<IXXiI'Xir> &\eqeps \frac12 \<IXi^2r> + \frac12 \<XiIIXir> - \frac{\discrep}2\;,
\quad&\quad \<IXXiI'Xir>^{(1)} &\eqeps \frac12\<IXi^2r> - \frac\discrep2 + \frac{1}{16\pi}\;,
\quad&\quad \<IXXiI'Xir>^{(3)} &= \<IXXiI'Xir>\;, \\
\<IXiI'XXir> &\eqeps \frac12 \<IXi^2r> - \frac12 \<XiIIXir> + \frac{\discrep}2 \;,
\quad&\quad \<IXiI'XXir>^{(1)} &\eqeps \frac12\<IXi^2r> + \frac\discrep2 - \frac{1}{16\pi}\;,
\quad&\quad \<IXiI'XXir>^{(3)} &= \<IXiI'XXir>\;, \\
\<I'XiIXi2r> &\eqeps \frac{1}{16\pi} - \frac12 \<XiIIXir> \;,
\quad&\quad \<I'XiIXi2r>^{(1)} &= 0\;, 
\quad&\quad \<I'XiIXi2r>^{(3)} &= 0\;, \\
\<IXiI'Xi2r> &\eqeps -\frac{1}{16\pi} + \frac12 \<XiIIXir>\;,
\quad&\quad \<IXiI'Xi2r>^{(1)} &= 0\;, 
\quad&\quad \<IXiI'Xi2r>^{(3)} &= 0\;.
\end{equs}
\end{proposition}

\begin{proof}
In order to find relations between these renormalisation constants, we note that they can all be expressed 
in terms of space-time kernels that have spatial Fourier transforms\footnote{We are talking here about functions
defined on all of $\R^2$, not on the torus!} of the type
\begin{equ}[e:defP]
\hat P_{f,g}(t,k) = \one_{t > 0}\exp\bigl(-(1+|k|^2) t\bigr) f(t) g(k)\;.
\end{equ}
In particular, it follows from \eqref{e:mollif} that our space-time noise mollifier
is given by $\hat \rho_\eps = P_{\rho_\eps',1}$.
An explicit calculation shows that the space-time convolution of two such kernels is again
a kernel of the same type, namely
\begin{equ}[e:idenConvol]
P_{f_1,g_1} \star P_{f_2, g_2} = P_{f,g}\;,\qquad f(t) = \int_0^t f_1(s)f_2(t-s)\,ds \;,\quad g(k) = g_1(k)g_2(k)\;.
\end{equ}
We will sometimes write $P_{f(t),g(k)}$ instead of $P_{f,g}$ when $f$ and \slash or $g$ are explicit
functions. For example, $P_{t,1}$ means that $f(t) = t, g(k) = 1$.
We will also use the fact that
\begin{equ}[e:idenMult]
x_j P_{f,g}(t,x) = -2 P_{tf(t),i k_j g(k)}(t,x) + P_{f,i \d_j g}(t,x)\;.
\end{equ}

Writing $\sscal{\place,\place}$ for the periodised $L^2$ scalar product defined above, we then have for example
\begin{equ}
\<I[I'Xi]I'Xir> = \sscal{P_{\rho_\eps',1} \star P_{1,ik_1}\star P_{1,1},P_{\rho_\eps',1} \star P_{1,ik_1}} 
= \sscal{P_{\int \rho_\eps,ik_1}, P_{\rho_\eps,ik_1}}\;.
\end{equ}
At this point we note that one can perform the integration over $k$ explicitly, using the useful 
identities
\begin{equ}[e:idenHeat]
\int_{\R^2} e^{-2|k|^2t}\,dk = \frac{\pi}{2t}\;,\qquad 
\int_{\R^2} |k_1|^2 e^{-2|k|^2t}\,dk = \frac12 \int_{\R^2} |k|^2 e^{-2|k|^2t}\,dk = \frac{\pi}{8t^2}\;.
\end{equ}
Let now $\sscal{\place,\place}_t$ be defined like $\sscal{\place,\place}$, but with $t$ held fixed instead of
being integrated out.
As a consequence of \eqref{e:idenHeat} and the decay properties of the 
massive heat kernel, one then
shows that there exist smooth functions $\chi, \hat \chi\colon \R_+ \to \R$
such that $\sup_{t > 0} (t^m + t^{-m}) (|\chi(t)| + |\hat \chi(t)|) < \infty$ for every $m > 0$ and such that
\begin{equ}
\sscal{P_{1,1},P_{1,1}}_t = \frac{e^{-2t}}{4\pi^2}\frac{\pi}{2t} + \chi(t)\;,\qquad
\sscal{P_{1,ik_1},P_{1,ik_1}}_t = \frac{e^{-2t}}{4\pi^2} \frac{\pi}{8t^2} + \hat \chi(t)\;.
\end{equ}
Here, the additional factor $\frac{1}{4\pi^2}$ comes from the fact that $\|f\|_{L^2}^2 = \frac{1}{4\pi^2} \|\hat f\|_{L^2}^2$ in two dimensions and the factor $e^{-2t}$ comes from the mass term included in \eqref{e:defP}.

Using the fact that
$\rho_\eps \to 1$ in a strong enough topology\footnote{The fact that $\sup_{t>0} (t^k \wedge t^{-k})|\rho_\eps(t) - 1| \to 0$ for some (in fact any) $k > 0$ is good enough for example.}, it then follows that
\begin{equ}[e:constBigSymbol]
\<I[I'Xi]I'Xir> \eqeps  \frac{1}{4\pi^2}\int_0^\infty \frac{\pi e^{-2t}}{8t^2} \rho_\eps(t) \int_0^t \rho_\eps(s)\,ds\,dt
+ \int_0^\infty t \hat \chi(t)\,dt\;.
\end{equ}
Similarly, we obtain from the first identity in \eqref{e:idenHeat}
\begin{equ}[e:XiIIXir]
\<XiIIXir> \eqeps \frac{1}{4\pi^2}\int_0^\infty \frac{\pi e^{-2t}}{2t} \rho_\eps'(t) \int_0^t \rho_\eps(s)\,ds\,dt\;.
\end{equ}
As a consequence of the fact that $\rho'$ is compactly supported, it is straightforward to see that
\begin{equ}[e:exprBasicRenorm]
\<XiIIXir> \eqeps \frac{1}{4\pi^2}\int_0^\infty \frac{\pi}{2t} \rho_\eps'(t) \int_0^t \rho_\eps(s)\,ds\,dt = \Crho\;,
\end{equ}
as claimed. To derive the identities for $\<I[I'Xi]I'Xir>$ and $\<IXiI'[I'Xi]r>$
appearing in the statement of the proposition (and in particular identify an
expression for $\discrep$), we integrate by parts in 
\eqref{e:XiIIXir}, so that
\begin{equ}[e:exprXiIIXi]
\<XiIIXir> \eqeps \int_0^\infty \frac{e^{-2t}}{4\pi t} \rho_\eps(t) \int_0^t \rho_\eps(s)\,ds\,dt - \int_0^\infty \frac{e^{-2t}}{8\pi t} \rho_\eps^2(t)\,dt + 
\int_0^\infty \frac{e^{-2t}}{8\pi t^2} \rho_\eps(t) \int_0^t \rho_\eps(s)\,ds\,dt\;.
\end{equ}
Since $\rho_\eps \uparrow 1$, one can verify that
\begin{equ}[e:appearance8pi]
\int_0^\infty \frac{e^{-2t}}{4\pi t} \rho_\eps(t) \int_0^t \rho_\eps(s)\,ds\,dt \eqeps \frac{1}{8\pi}\;.
\end{equ}
The constant $\<IXi^2r>$ on the other hand is given by
\begin{equs}
\<IXi^2r> &\eqeps \frac{1}{4\pi^2}\int_0^\infty \frac{\pi }{2t} e^{-2t} \rho_\eps^2(t)\,dt + \int_0^\infty \chi(t)\,dt\\
&\eqeps \sigma_{\rho,\eps}^2 + \frac{1}{4\pi^2}\int_0^1 \frac{\pi }{2t} (e^{-2t}-1)\,dt
+ \frac{1}{4\pi^2}\int_1^\infty \frac{\pi }{2t} e^{-2t}\,dt + \int_0^\infty \chi(t)\,dt \\
&\eqcolon \sigma_{\rho,\eps}^2 + \hat\discrep\;,\label{e:constCherry}
\end{equs}
where $\hat \discrep$ is independent of the choice of cutoff $\rho\in \mc{P}$,
but does depend on our (arbitrary) choice of mass in the Laplacian and on the size of the torus we work on.
Using the properties \eqref{main:eq-mcP} of $\rho$, one has
\begin{equ}
\sigma_{\rho,\eps}^2 = \frac{1}{8\pi}\int_0^1 \frac{\rho_\eps^2(t)}{t}\,dt
= \frac{1}{8\pi}\int_{2\eps^2}^1 \frac{dt}{t} + \frac{1}{8\pi}\int_{1}^2 \rho^2(t)\frac{dt}{t}\;,
\end{equ}
(for small enough $\eps$)
which is indeed of the desired form $(4\pi)^{-1} \abs{\log \eps} + \hCrho$ for some constant
$\hCrho$ that only depends on our choice of $\rho$.

Combining this with \eqref{e:constBigSymbol}, \eqref{e:exprXiIIXi} and 
\eqref{e:appearance8pi} yields the claimed expressions for $\<I[I'Xi]I'Xir>$ and
$\<IXiI'[I'Xi]r>$, namely
\begin{equ}[e:idenConst1]
4 \<I[I'Xi]I'Xir> = - 4\<IXiI'[I'Xi]r> \eqeps \<IXi^2r> + \<XiIIXir> - \discrep\;,
\end{equ}
where we defined the constant
\begin{equ}[e:defdiscrep]
\discrep \eqdef \frac{1}{8\pi} - \int_0^\infty \bigl(4t \hat \chi(t) - \chi(t)\bigr)\,dt\;.
\end{equ}
Finally, we note that
\begin{equs}
\<IXXiI'Xir> &= 2 \sscal{P_{\rho_\eps,ik_1} ,P_{\int \rho_\eps,ik_1}} =  2\<I[I'Xi]I'Xir>\;,\\
\<IXiI'XXir> &= \sscal{P_{\rho_\eps,1} ,P_{\rho_\eps,1}-2 P_{\int \rho_\eps,k_1^2}} =  \<IXi^2r> - 2\<I[I'Xi]I'Xir>\;,\\
\<I'XiIXi2r> &=  - \<IXiI'Xi2r> = 2\sscal{P_{\rho_\eps,ik_1} , P_{\int s \rho_\eps'(s),ik_1}} \eqeps \frac{1}{4\pi^2}\int_0^\infty \frac{\pi}{4t^2} 
e^{-2t}\rho_\eps(t)\int_0^t s\rho_\eps'(s)\,ds\,dt \\
&= \frac{1}{4\pi^2}\int_0^\infty \frac{\pi}{4t} e^{-2t} \rho_\eps^2(t)\,dt - \frac{1}{4\pi^2}\int_0^\infty \frac{\pi}{4t^2} 
\rho_\eps(t)e^{-2t}\int_0^t \rho_\eps(s)\,ds\,dt\\
&\eqeps \frac12 \<IXi^2r> - 2 \<I[I'Xi]I'Xir> + \frac12\int_0^\infty \bigl(4t\hat \chi(t) - \chi(t)\bigr)\,dt \eqeps \frac{1}{16\pi} - \frac12 \<XiIIXir> \;.\label{e:idenConst2}
\end{equs}
where we integrated by parts in the inner integral, which concludes the proof of the
expressions in the first column of the statement.

We now note that by \eqref{e:idenConvol} and \eqref{e:idenMult} one has
\begin{equs}
\PPi^{(\eps, 2)} \<IXiX>_j &\eqsim (-x_j P_{1,1}) \star \boldxi_\eps = 2 P_{t, ik_j} \star P_{\rho_\eps',1} \star\boldxi
= 2 P_{\int \rho_\eps,ik_j} \star \boldxi\;,\\
\PPi^{(\eps, 1)} \<IXiX>_j &\eqsim (-x_j (\hat\rho_\eps \star P_{1,1})) \star \boldxi = (-x_j P_{\rho_\eps,1}) \star \boldxi
= 2 P_{t\rho_\eps(t),ik_j}\star\boldxi\;,
\end{equs}
where $\eqsim$ means that the identity holds modulo some
multiple of $\Pi^{(\eps, 2)} \<IXi>$.
On the other hand, we have
\begin{equ}
\PPi^{(\eps, 2)} \<IXi2>_j = 2 P_{1,1} \star P_{t\rho_\eps'(t),ik_j} \star \boldxi\;.
\end{equ}
Since $\int_0^t s\rho_\eps'(s)\,ds = t \rho_\eps(t) - \int_0^t \rho_\eps(s)\,ds$, we conclude that
\begin{equ}
\PPi^{(\eps, 2)} \<IXi2>_j \eqsim \PPi^{(\eps, 1)} \<IXiX>_j - \PPi^{(\eps, 2)} \<IXiX>_j\;.
\end{equ}
Since $\PPi^{(\eps, 2)} \<I'Xi> = \PPi^{(\eps, 1)} \<I'Xi>$ and since $\<IXiI'Xir> = \<IXiI'Xir>^{(1)} = 0$
by \eqref{e:vanishSymmetry}, it follows that
\begin{equ}[e:idenRenorm]
	\<IXXiI'Xir> = - \<I'XiIXi2r> + \<IXXiI'Xir>^{(1)} \;.
\end{equ}
Similarly, we have (this time modulo multiples of $\Pi^{(\eps, 2)} \<I'Xi>$):
\begin{equs}
\PPi^{(\eps, 2)} \<I'XiX>_{jj} &\eqsim (-x_j P_{1,ik_j}) \star \boldxi_\eps = (P_{1,1}-2 P_{t, k_j^2}) \star P_{\rho_\eps',1} \star\boldxi
= (P_{\rho_\eps,1}-2 P_{\int \rho_\eps,k_j^2}) \star \boldxi\;,\\
\PPi^{(\eps, 1)} \<I'XiX>_{jj} &\eqsim (-x_j (\hat\rho_\eps \star P_{1,ik_j})) \star \boldxi = (-x_j P_{\rho_\eps,ik_j}) \star \boldxi
= (P_{\rho_\eps,1}-2 P_{t\rho_\eps(t), k_j^2})\star\boldxi\;,
\end{equs}
whence it follows that $\PPi^{(\eps, 2)} \<I'Xi2>_{jj} \eqsim \PPi^{(\eps, 1)} \<I'XiX>_{jj} - \PPi^{(\eps, 2)} \<I'XiX>_{jj}$ and
therefore
\begin{equ}[e:idenRenorm2]
	\<IXiI'XXir> \eqeps - \<IXiI'Xi2r> + \<IXiI'XXir>^{(1)} \;.
\end{equ}
As a consequence of \eqref{e:idenConst1}, \eqref{e:idenRenorm}, 
\eqref{e:idenRenorm2}, and \eqref{e:idenConst2},
we obtain the claimed expressions
\begin{equs}
\<IXXiI'Xir>^{(1)} &\eqeps 2\<I[I'Xi]I'Xir> + \frac{1}{16\pi} - \frac12 \<XiIIXir>\;,\\
\<IXiI'XXir>^{(1)} &\eqeps \<IXi^2r> - 2\<I[I'Xi]I'Xir> 
 + \frac12 \<XiIIXir> - \frac{1}{16\pi}
 \eqeps 2\<I[I'Xi]I'Xir> + \discrep - \frac{1}{16\pi} - \frac12 \<XiIIXir>\;. 
\end{equs}
In particular, one has the useful identity
\begin{equ}[e:combinX]
\<IXXiI'Xir>^{(1)} - \<IXiI'XXir>^{(1)} \eqeps \frac{1}{8\pi} - \discrep\;.
\end{equ}
The remaining renormalisation constants for the model $Z^{(\eps,1)}$ are the same
as for $Z^{(\eps,2)}$ except that,
as a consequence of the fact that $\Xi$ reconstructs to space-time white noise and that the 
kernel appearing in $\<XiIIXir>^{(1)}$ is non-anticipative,
as well as the fact that $\hat \Xi$ reconstructs to $0$,
one has
\begin{equ}[e:vanishing]
\<XiIIXir>^{(1)} = \<I'XiIXi2r>^{(1)} = \<IXiI'Xi2r>^{(1)} = 0\;.
\end{equ}
Regarding the model $Z^{(\eps,3)}$, its renormalisation constants are exactly the same as for 
$Z^{(\eps,2)}$ with the sole exception that
\begin{equ}[e:vanishZ3]
\<I'XiIXi2r>^{(3)} = \<IXiI'Xi2r>^{(3)} = 0\;,
\end{equ}
since $\hat \Xi$ again reconstructs to $0$.
\end{proof}

\subsection{Expressions for the counterterms}
\label{sec:renCounterterms}

We now compute the counterterms arising from the BPHZ renormalisations of the 
models $Z_\eps^{(\place)}$ used to solve the fixed point problems \eqref{e:fpabstract}.
For this, we first note that as far as the calculation of the renormalisation constants
is concerned, we can put this system back into the more conventional differential form 
\begin{equs}
	\partial_t A_i &= \Delta A_i + \big[A_j,2\partial_j A_i - \partial_i A_j + [A_j,A_i]\big] 
	+ \bil(\higgs \otimes (\partial_{i} \higgs + A_{i}\higgs))+ U_A\xi_i + U_A[B_j,\hat \xi_{ij}] \;,\\
	\partial_t\higgs  &= 
	\Delta \higgs -\mass\higgs + 2 A_{j} \partial_{j}\higgs + A_{j}^{2}\higgs - |\higgs|^2 \higgs + U_\higgs \zeta+ U_\higgs \repr(B_j) \hat \zeta_{j}\;, \label{e:bigSys}\\
\d_t U &= \Delta U + F + \repr([A_j, h_j - (U_A B)_j])U\;.
\end{equs}
(A more detailed explanation of why this consists only in a harmless change of notations
is given in \cite[Sec.~5.3]{CCHS22}.)
Here $U_A$ and $U_\higgs$ denote the projections of $U$ onto
$\End(\fg)$ and $\End(\V)$ respectively.
The term $F$ contains only terms that are sufficiently regular so that they cannot contribute
to the counterterms, so we ignore it for the purpose of this section. We have also ignored
the presence of the process $\eta$, of the finite counterterms, as well as the
whole equation for $h$ for the same reason.

\begin{remark}
Note that the equation for $U$ does not require any renormalisation either since $U$, $B$ and $h$
are sufficiently regular so that the right-hand side is classically well-posed. However, it will
play a role in the renormalisation of the equation for $A$ due to the term $U\xi$
which is classically ill-posed.
\end{remark}

The procedure for computing the counterterms arising in the equations for 
$A$ and $\Phi$ is then as follows. First, note that the nonlinearities 
appearing in the right-hand sides are
local functions of $A \in \fg^2$, $DA \in \fg^4$, $\Phi \in \V$, $D\Phi\in \V^2$, 
$U \in \End(\Y)$, $h \in \fg^2$, $B \in \fg^2$, $\xi \in \fg^2$, $\zeta \in \V$, 
$\hat \xi \in \fg^4$, and $\hat\zeta \in \V^2$. We collect all of these components
into one single vector space $\CA$ and view the nonlinearities as polynomial 
functions $P_A \colon \CA \to \fg^2$, $P_\higgs \colon \CA \to \V$, and
$P_U \colon \CA \to \End(\Y)$.
For every symbol $\tau$, one then defines
inductively (over the ``height'' of $\tau$, starting from the noise symbols) 
maps $\bar\bUpsilon_i[\tau] \colon \CA \to \CT[\tau]\otimes \fg$ (corresponding to
the component $A_i$)
and $\bar\bUpsilon_\higgs[\tau] \colon \CA \to \CT[\tau]\otimes \V$ by the procedure 
described in \cite[Sec.~5.8]{CCHS22} (see also \cite{BCCH}). 

Given a symbol $\tau$, write  now $\CV[\tau]$ for the corresponding
subspace in the regularity structure where we set $\fg = \V = \R$.
If the symbol $\tau$ contains $k$ noises, then $\CT[\tau]$
can be naturally identified with 
\begin{equ}
\CT[\tau] \simeq \CV[\tau]\otimes L_\tau\;,
\end{equ}
where $L_\tau$ denotes the space of $k$-linear forms on
$\genV_1\times\dots\times \genV_k$,
where each $\genV_i$ is a copy of $\fg$ or $\V$ depending
on which noises appear in $\tau$. 
Given a concrete element $\hat\tau \in \CV[\tau]$,
we can therefore interpret $\bar\bUpsilon_i[\hat\tau]$ as a function of $(A,DA,\ldots)$
taking values in a space of $\fg$-valued $k$-linear forms.
We will write this as an expression
with $k$ placeholders `$\place$' for the arguments of the $k$-linear form.
In fact, in all cases of interest to us here, one has $k=2$, so 
$\bar\bUpsilon_i[\hat\tau]$ is of the form
\begin{equ}
\bar\bUpsilon_i[\hat\tau](A,DA,\dots)(\place,\place)\;.
\end{equ}
As already alluded to in the previous section, for all symbols $\tau$ 
and all models $Z^{(\eps,i)}$ of interest to us,
the character $\ell^{(\eps,i)}$ on $\CT[\tau]$ factors into
a character (which we again call $\ell^{(\eps,i)}$ in keeping with the
previous section) on
$\CV[\tau]$ and an element of $\genV_1\otimes \genV_2$.
That element furthermore vanishes when $\genV_1 \neq \genV_2$
and is otherwise equal to either $\Cov_\fg$ or $\Cov_\V$
(see Definition~\ref{prelim:def-Casimir-induced}),
depending on whether the noises appearing in $\tau$ are $\fg$-valued or $\V$-valued.

What we then call the counterterm generated by a concrete
element $\hat\tau \in \CV[\tau]$ is the expression
\begin{equ}
\textcolor{darkred}{\hat\tau} \bar\bUpsilon_i[\hat\tau](A,DA,\dots)(\Cov_\fg)\;,
\end{equ}
(or possibly the same with $\Cov_\V$ instead of $\Cov_\fg$) where 
$\textcolor{darkred}{\hat\tau}$ is just a placeholder for a
concrete renormalisation constant $\ell^{(\eps,i)}(\hat\tau)$.

For example, we have $\bar\bUpsilon_i[\<Xi>_i](A,\ldots) = U_A \place$ which
since $U_A \in \End(\fg)$ is indeed a function of $U$ with values
in linear maps from $\fg$ (the space that the noise $\xi_i$ takes values in)
to $\fg$ (the space that $A_i$ takes values in). The counterterm that this 
generates vanishes since the BPHZ character vanishes on $\<Xi>$.
The symbol \<IXi^2> naturally comes with two indices $j,k$ with $j,k \in \{1,2,\star\}$
corresponding to the components
of the two noises that appear in it (here $\star$ denotes the noise $\zeta$).
In the terminology of the discussion above, this is just another way of
saying that $\CV[\<IXi^2>] \simeq \R^3 \otimes \R^3$ which is consistent with 
\eqref{e:spaces} and the fact that $\Y \simeq \R^3$ in the case $\fg = \V = \R$.

Since the BPHZ character vanishes on 
$\<IXi^2>_{jk}$ when $j\neq k$, it suffices to consider the case $j=k$.
The only terms in the right-hand side of the equation for $A_i$ that yield 
counterterms in $\bar\bUpsilon_i[\<IXi^2>_{jj}]$ are the two cubic terms since
the quadratic terms all involve spatial derivatives. With the notation just 
explained, we then have
\begin{equ}{}
\bar\bUpsilon_i[\<IXi^2>_{jj}](A,\ldots) = (1-\delta_{ij})[U_A\place,[U_A\place,A_i]]\;,
\end{equ}
as well as
\begin{equ}{}
\bar\bUpsilon_i[\<IXi^2>_{\star\star}](A,\ldots) = \bil(U_\higgs\place \otimes A_i  U_\higgs\place)\;.
\end{equ}
Summing over the three possible values for $j$, we conclude
that the counterterm generated by the symbols of type \<IXi^2> for the $A$ equation 
is given by
\minilab{e:counterTerms}
\begin{equ}[e:counterTerm1]
\<IXi^2r> (\Cas + \VCas)A\;.
\end{equ}

\begin{remark}\label{rem:placeholderBPHZ}
As already explained, here we don't interpret $\<IXi^2r>$ as a shortcut for
$\<IXi^2r>^{(2)}$ as we did in Section~\ref{sec:renConstants}, but rather 
as a placeholder
that will later be replaced by $-\<IXi^2r>^{(i)}$ when computing the BPHZ
counterterm arising in the equation for $Y^{(i)}_\eps$, or in other
words coming from the BPHZ renormalisation of the model $Z^{(\eps,i)}$.
(Here the minus sign is the one coming from \eqref{e:BPHZantipode}.)
\end{remark}

Similarly, we find that the corresponding counterterm 
for the $\higgs$ equation is given by
\minilab{e:counterTermsPhi}
\begin{equ}[e:counterTermPhi1]
\<IXi^2r> \bigl(2\PhiCas - (2 + \dim \V)\Id_\V\bigr)\higgs\;,
\end{equ}
where the factor $2$ in front of $\PhiCas$ comes from the summation over $j$ in the term $A_j^2\higgs$.

We now turn to symbols of the type $\<I[I'Xi]I'Xi>$ and $\<IXiI'[I'Xi]>$.
It is natural to assign to each of these two symbols four spatial \slash noise indices $ijk\ell$
which, reading the graphical representation of the symbol counterclockwise from the top, represent the 
component $i \in \{1,2,\star\}$ of the top noise (here the $\star$ component represents $\zeta$), 
the direction $j\in \{1,2\}$ in which the derivative of the `upper' integration kernel is taken,
the direction $k\in \{1,2\}$ in which the derivative of the `lower' integration kernel is taken,
and finally the component $\ell \in \{1,2,\star\}$ of the second noise.

The symbol 
$\<I[I'Xi]I'Xi>_{\ell kk\ell}$ appears in the expansion of $[A_j,2\d_j A_i - \d_i A_j]$ for $A_i$, which yields a coefficient
\begin{equ}
\bigl(2\delta_{i\ell}(2\delta_{k\ell}-1) - \delta_{ik}(2-\delta_{k\ell})\bigr)[[A_i,U\place],U\place]\;.
\end{equ}
Similarly, the symbol $\<IXiI'[I'Xi]>_{\ell kk\ell}$ comes with the coefficient
\begin{equ}
\bigl(2\delta_{\ell k}(2 \delta_{i\ell}-\delta_{ik})-\delta_{ik}(2-\delta_{\ell k})\bigr)[U\place, [A_i,U\place]]\;.
\end{equ}
Summing over $k,\ell$, 
we find that the resulting counterterm is given by
\minilab{e:counterTerms}
\begin{equ}[e:counterTerm2]
-\bigl(3 \<I[I'Xi]I'Xir> - \<IXiI'[I'Xi]r>\bigr) \Cas A_i \;.
\end{equ}

Regarding the counterterms for $\higgs$, the symbol 
$\<I[I'Xi]I'Xi>_{\ell kk\ell}$ does not appear in the expansion of the right-hand side of the
$\higgs$ equation since none of the terms includes a factor $\d A$.
At first sight, the symbol $\<IXiI'[I'Xi]>_{\ell kk\ell}$ might appear in the expansion of
$2 A_j \d_j \higgs$ (at least for $j = \ell$), but since no factor $\d A$ appears in the right-hand side for $\higgs$,
it actually doesn't.

The symbols $\<I[I'Xi]I'Xi>_{\star kk\star}$ and $\<IXiI'[I'Xi]>_{\star kk\star}$ appear in the expansion of $\bil(\d_i\higgs \otimes \higgs)$
with respective coefficients
\begin{equ}
2\delta_{ik}\bil((U\place)\otimes (A_i U\place))
\qquad\text{and}\qquad
2\delta_{ik}\bil((A_i U\place)\otimes ( U\place))\;,
\end{equ}
which (since $\bil$ comes with a minus sign in the equation) yields the counterterm 
\minilab{e:counterTerms}
\begin{equ}[e:counterTerm3]
- 2(\<I[I'Xi]I'Xir> - \<IXiI'[I'Xi]r>) \VCas A_i \;,
\end{equ}
after summing over $k$ and contracting with the covariance of the noise.

Regarding the equation for $\higgs$, the symbol $\<IXiI'[I'Xi]>_{\star kk\star}$ does not 
contribute since the right-hand side of the equation contains no term with a factor $\higgs \d A$ (or $\higgs \d\higgs$).
The symbol $\<I[I'Xi]I'Xi>_{\star kk\star}$ however appears in the expansion for
$2 A_j \d_j \higgs$ with coefficient $-2\delta_{jk} \bil(U\place \otimes \higgs)U\place$. 
By Lemma~\ref{prelim:lem-identity-Casrhog}, this
leads
to the counterterm 
\minilab{e:counterTermsPhi}
\begin{equ}[e:counterTermPhi2]
-4\<I[I'Xi]I'Xir> \PhiCas \higgs
\end{equ}
for the $\higgs$ equation (the additional factor $2$ comes from the fact that the sum of $\delta_{jk}$ over all values of $j$ and $k$
is $d=2$). 

The symbol $\<XiIIXi>$ has coefficient $[[U\place,h - UB],U\place]$ in the equation
for $A$, yielding a counterterm
\minilab{e:counterTerms}
\begin{equ}[e:counterTerm4]
\<XiIIXir> \Cas (UB - h)\;.
\end{equ}
No counterterm is generated by this symbol for the $\higgs$ equation since the right-hand side
of the equation for $U$ does not involve $\higgs$ explicitly.

Regarding the symbols of type \<I'XiIXi2> and \<IXiI'Xi2>, we again
assign indices to them by reading from left to right. This time the third index corresponds 
to the direction in which the derivative appearing in the definition of $\hat \xi_\eps$ 
(represented by the square node) is taken.
The terms giving rise to counterterms for these symbols are the terms including
first order derivatives, namely
\begin{equ}
2[A_j, \d_j A_i] - [A_j, \d_i A_j] - \bil(\d_i \higgs \otimes \higgs)\;,
\end{equ}
regarding the equation satisfied by $A_i$, as well as potentially the expression 
$2 A_j \d_j \higgs$ appearing in the right-hand side of the equation for $\higgs$.

When considering the modelled distribution for the term $2[A_j, \d_j A_i]$, we see that the only
symbols of the type \<I'XiIXi2> and \<IXiI'Xi2> having nonzero coefficients  
are $\<I'XiIXi2>_{ijkj}$ with coefficient $2[U[B_k, \place] , U\place]$ and $\<IXiI'Xi2>_{jjki}$
with coefficient $2[U\place, U[B_k, \place]]$. Since, for any orthonormal basis $\{g_\alpha\}$ of $\fg$
and any $X\in \fg$ we have
\begin{equ}
\Cas X = \sum_\alpha [[X,g_\alpha],g_\alpha]\;,
\end{equ}
since $U$ is unitary, and since by \eqref{e:idenConst2} both symbols yield the same
contribution, this leads to the counterterm $2\delta_{ij}\delta_{kj} (\<I'XiIXi2r> - \<IXiI'Xi2r>) \Cas U B_k$,
eventually leading to 
\begin{equ}[e:counterTermCancel1]
2(\<I'XiIXi2r> - \<IXiI'Xi2r>) \Cas U B_i\;,
\end{equ}
when summing over $j,k$.

Similarly, when considering the modelled distribution for the $-[A_j, \d_i A_j]$ term, we see that 
$\<I'XiIXi2>_{jikj}$ appears with coefficient $-[U[B_k,\place], U\place]$, while 
$\<IXiI'Xi2>_{jikj}$ appears with coefficient $-[U\place, U[B_k, \place]]$. 
Similarly to before, this yields the counterterm 
$\delta_{ik} (\<IXiI'Xi2r> - \<I'XiIXi2r>) \Cas  U B_k$ which leads to the counterterm
\begin{equ}[e:counterTermCancel2]
2(\<IXiI'Xi2r> - \<I'XiIXi2r>) \Cas U B_i\;,
\end{equ}
when summing over $j,k$ (since $j$ can take $2$ values in 
dimension $2$ and this expression is independent of $j$), thus cancelling exactly the previous counterterm
\eqref{e:counterTermCancel1}.

Regarding the term $\bil(\d_i \higgs \otimes \higgs)$, we find that
$\<I'XiIXi2>_{\star ik\star}$ appears with coefficient 
$\bil((U \place) \otimes (U B_k\place))$, while 
$\<IXiI'Xi2>_{\star ik\star}$ appears with coefficient 
$\bil((U B_k \place) \otimes (U\place))$ thus, similarly to before, yielding the counterterm
\minilab{e:counterTerms}
\begin{equ}[e:counterTerm5]
(\<IXiI'Xi2r> - \<I'XiIXi2r>) \VCas UB_i\;.
\end{equ}
This term will of course only be relevant in the covariant case since this is the
only case for which $\hat \xi$ is replaced by a nonzero noise.
Since the symbols $\<I'XiIXi2>$ and $\<IXiI'Xi2>$ can only appear in the expansion of  
terms containing either a factor $A\d A$ or a factor $\higgs \d\higgs$ (mixed terms yield 
vanishing renormalisation constants in all of the models under consideration)
and since such terms do not appear in the right-hand side of the equation
for $\higgs$, they do not lead to any counterterms there.

Regarding the symbols $\<IXiI'XXi>$ and $\<IXXiI'Xi>$, the exact same 
calculation as for \<I'XiIXi2> and \<IXiI'Xi2> shows that the counterterms
generated by the two expressions that are bilinear in $A$ and $\d A$ cancel out exactly,
whatever the values of the renormalisation constants assigned to $\<IXiI'XXir>$ and $\<IXXiI'Xir>$.
The modelled distribution for the term $\bil(\d_i \higgs \otimes \higgs)$ gives 
rise to $\<IXiI'XXi>_{\star ik\star}$
with coefficient $\bil((\d_k U \place) \otimes (U\place))$ and to 
$\<IXXiI'Xi>_{\star ik\star}$ with coefficient $\bil((U \place) \otimes (\d_k U\place))$.
Note further that $(\d_k U) U^* = \repr(h_k)$, so since $\bil$ has a minus sign
we finally end up generating the counterterm given by
\minilab{e:counterTerms}
\begin{equ}[e:counterTerm6]
(\<IXiI'XXir> - \<IXXiI'Xir>) \VCas h_i\;.
\end{equ}
Since the equation for $\higgs$ does not contain any term with a factor of the type 
$A \d A$ or $\higgs \d\higgs$, these symbols do not generate any counterterm for $\higgs$ either.

Since none of our symbols exhibit subdivergences, the BPHZ character is simply given by $-\ell_\eps$
and the counterterms all come with the opposite sign from what we just computed,
see Remark~\ref{rem:placeholderBPHZ}.
Combining all the expressions \eqref{e:counterTerm1}--\eqref{e:counterTerm6}, we thus obtain in general the counterterm
\begin{equ}[e:generalCounter]
 (4\<I[I'Xi]I'Xir>-\<IXi^2r>)( \Cas +\VCas)A + \<XiIIXir> \Cas (h - UB)
 + (\<IXXiI'Xir> - \<IXiI'XXir>)  \VCas h - (\<IXiI'Xi2r> - \<I'XiIXi2r>)\VCas U B\;, 
\end{equ}
for the $A$ equation. Similarly, combining \eqref{e:counterTermPhi1}--\eqref{e:counterTermPhi2}, we obtain the counterterm
\begin{equ}[e:generalCounterPhi]
\<IXi^2r> \bigl((2 + \dim \V)\Id_\V - 2\PhiCas\bigr)\higgs + 4\<I[I'Xi]I'Xir> \PhiCas \higgs
\;,
\end{equ}
for the $\Phi$ equation.
In order to obtain these expressions, the only identity between renormalisation constants that we made use of
is the identity $\<I[I'Xi]I'Xir> = - \<IXiI'[I'Xi]r>$ which holds for every single model.

We now use Proposition~\ref{prop:relConstants} to compare this term for all three choices of
model $Z_\eps^{(i)}$. For ease of reference, we refer to the case $i=1$ as the ``basic case'',
the case $i=2$ as the ``covariant case'' and the case $i=3$ as the ``gauged case''.
In the basic case, we use the renormalisation constants for the model $\PPi^{(\eps,1)}$, which 
are such that $\<XiIIXir>^{(1)} = \<IXiI'Xi2r>^{(1)} = \<I'XiIXi2r>^{(1)} = 0$ (so that the second and fourth terms
of \eqref{e:generalCounter} vanish in that case) and $\<IXXiI'Xir>^{(1)} - \<IXiI'XXir>^{(1)}
\eqeps \frac{1}{8\pi} - \discrep$, see Proposition~\ref{prop:relConstants}. 
Since furthermore $4\<I[I'Xi]I'Xir>-\<IXi^2r> \eqeps \<XiIIXir> - \discrep$,
the counterterm is given (modulo terms vanishing as $\eps \to 0$) by
\begin{equ}[e:counterBasic]
 \bigl(\<XiIIXir> - \discrep \bigr)( \Cas +\VCas)A + \Bigl(\frac{1}{8\pi} - \discrep\Bigr)\VCas h\;.
\end{equ}
Using again Proposition~\ref{prop:relConstants} (in particular that $\<IXXiI'Xir> - \<IXiI'XXir> \eqeps \<XiIIXir> - \discrep$), we see that the counterterm for the gauged case (where we use the model $Z^{(\eps,3)}$ so that
in particular $\<IXiI'Xi2r>^{(3)} = 0$) is given by
\begin{equs}[e:counterGauged]
\bigl(\<XiIIXir> - \discrep\bigr)&( \Cas +\VCas)A + \bigl(\<XiIIXir> \Cas  + (\<IXXiI'Xir> - \<IXiI'XXir>)\VCas \bigr)h - \<XiIIXir> \Cas UB \\
&\eqeps  \bigl(\<XiIIXir> - \discrep \bigr)( \Cas +\VCas)(A+h) + \discrep\Cas h- \<XiIIXir> \Cas UB\;.
\end{equs}

Finally, the counterterm in the covariant case is given by \eqref{e:generalCounter} which,
making use of the identity $2 \<IXiI'Xi2r> \eqeps \<XiIIXir> - \frac{1}{8\pi}$ and
using the fact that $g\act B = UB - h$, can be written as
\begin{equ}[e:counterCovariant]
 (\<XiIIXir> - \discrep)( \Cas +\VCas)\bigl(A - g\act B\bigr)
 - \discrep  \Cas (g\act B)
  - \Bigl(\discrep - \frac{1}{8\pi}\Bigr) \VCas U B\;. 
\end{equ}
Recall furthermore that we chose the abstract fixed point problem \eqref{e:fpabstract}
in such a way that the right-hand side of the equation for $A$ already contains 
an additional counterterm
\begin{equ}[e:extraTerm]
\massA A  + \hat\massA (A+h)\;.
\end{equ}
We would now like to choose the constants $\massA$ and $\hat\massA$ in such a way 
that there exists
a linear map $\renormA$ with the following properties
\begin{enumerate}
\item The counterterm in the basic case is given by $\renormA A$.
\item The counterterm in the covariant case is given by $\renormA (A-g\act B)$ for every choice of $B$.
\item The counterterm in the gauged case is given by $\renormA (A+h)$, provided that one takes $B=0$.
\end{enumerate}

\begin{remark}\label{rem:dependMass}
While the resulting limiting process should of course be allowed to depend on physical 
parameters of the model like the size of the torus, it should
not depend on arbitrary choices like the choice of cutoff function $\rho$ or the choice of mass in $\fp$.
Since the BPHZ model is independent of $\rho$ but does in general depend on the choice of mass 
in $\fp$ (which is the analogue in
our case to the large-scale cutoff used in \cite{BHZ19}), this means that the 
constants $\hat \mass$, $\massA$ and $\hat\massA$ should also be independent of $\rho$, but \textit{are}
in fact allowed to depend on this mass, provided that the resulting linear map $\renormA$
(which appears in a formulation where this mass doesn't appear) doesn't depend on it.
\end{remark}

\begin{proposition}\label{prop:processesEqual}
Writing $\hat Z^{(\eps,i)}$ for the BPHZ renormalisation of
the models $Z^{(\eps,i)}$ given in Definition~\ref{def:Models}, there
exists an admissible model $\hat Z$ such that 
\begin{equ}
\lim_{\eps \to 0} \$\hat Z; \hat Z^{(\eps,i)}\$ = 0\;,
\end{equ}
in probability for all $i$. 

Furthermore, let $\massA = \discrep\Cas$, $\hat\massA = \bigl(\discrep - \frac{1}{8\pi}\bigr)\VCas$, and
\begin{equ}[e:choiceHiggsMass]
\hat\mass
= \mass + \hat\discrep(2+\dim\V)\Id_\V
  - (\hat\discrep+\discrep)\PhiCas\;.
\end{equ}
Then, there exist choices of $C_\eps^{(i)} \in L(\Y\oplus\fg^2\oplus\fg^2,\Y)$ 
tending to $0$ such that
the solutions to \eqref{e:fpabstract} with models $\hat Z^{(\eps,i)}$ 
and with $C_\eps = C_\eps^{(i)}$ coincide 
after reconstruction with
the processes $Y_\eps^{(i)}$ defined in Section~\ref{sec:threeProcesses}.
(In the case $i=3$, this is the case only when $B=0$.)

In particular, there exists a stopping time $\tau>0$ and 
a time-dependent adapted gauge transformation $g$ such that, writing
$Y$ for the solution to \eqref{e:covariant_SPDE} and $\hat Y$ for the solution
to \eqref{e:wantedSPDE} (with the same initial data $Y_0 \in \CC_x^{-\kappa}$ and
the same $B \in \CC^{2-\kappa,-\kappa}$), 
one has the almost sure identity
\begin{equ}[e:processesEqual]
\bigl(g(t)\act Y(t)\bigr)(\boldxi) =  \bigl(\hat Y(t)\bigr)(g \act \boldxi) \;,
\end{equ}
for all $t \le \tau$. Furthermore, for every $k>0$ there 
exist $C>0$ and $n>0$ such that
\begin{equ}[e:boundComparisonTime]
\E (\tau \wedge 1)^{-k} \le C (1 + \|Y_0\|_{\CC_x^{-\kappa}}
+ \|B\|_{\CC^{2-\kappa,-\kappa}})^n\;.
\end{equ}
\end{proposition}

\begin{proof}
The first claim follows from the continuity properties of the BPHZ
model with respect to the noise (see for example \cite{CH16,HS23}), combined 
with its continuity with respect to the choice of integration kernel.
(The latter follows from the same proof as in \cite{HS23}.)

In order to identify the corresponding solutions, it remains to show that
the choices for $\massA$, $\hat\massA$ and $\hat\mass$ given in the statement 
yield the counterterms appearing in the equations for the processes $Y_\eps^{(i)}$.
Turning first to the equations for the $A$ components, if we 
add the term \eqref{e:extraTerm} with these choices of $\massA$ and $\hat\massA$
to each of the three counterterms
we just computed in \eqref{e:counterBasic}--\eqref{e:counterCovariant}, then the 
counterterm in the basic case becomes
\begin{equ}
 \<XiIIXir>( \Cas +\VCas)A - \frac{1}{8\pi}\VCas A\;,
\end{equ}
which does indeed coincide asymptotically with $\renormA A$, where
$\renormA$ is as in Definition~\ref{main:def-renormalization-constant}.
By \eqref{e:exprBasicRenorm}, this does indeed satisfy the property of Remark~\ref{rem:dependMass},
namely it is independent of the choice of mass appearing in $\fp$.

Writing $h = dg\,g^{-1}$ and $U = \Ad_g$ as before,
the counterterm in the covariant case is then
\begin{equ}
 \<XiIIXir>( \Cas +\VCas)\bigl(A +h-U B\bigr)
  -\frac{1}{8\pi} \VCas \bigl(A +h-U B\bigr) \eqeps \renormA (A-g\act B)\;,
\end{equ}
while the counterterm in the gauged case is
\begin{equ}
 \<XiIIXir>( \Cas +\VCas)(A+h) - \frac{1}{8\pi}\VCas (A+h) - \<XiIIXir> \Cas UB \eqeps \renormA(A+h) - \<XiIIXir> \Cas UB\;,
\end{equ}
which is of the desired form when $B=0$.
This does indeed satisfy all three conditions above.

Regarding the equation for $\higgs$, its BPHZ counterterms are
identical in all three cases. Combining \eqref{e:generalCounterPhi}
with the expression \eqref{e:constCherry} for $\<IXi^2r>$, 
the identity \eqref{e:idenConst1} for $\<I[I'Xi]I'Xir>$, and the fact that
$\<XiIIXir> \eqeps \Crho$ by
\eqref{e:exprBasicRenorm}, we see that the choice \eqref{e:choiceHiggsMass}
is such that the sum of the BPHZ counterterm
\eqref{e:generalCounterPhi} and the term $-\hat\mass\higgs$
already present in \eqref{e:fpabstract} agrees with
$\renormphi\higgs$ up to a term that vanishes as $\eps\to0$. 
This remaining discrepancy can be
cancelled exactly by the Higgs component of $C_\eps$.
Note that the expression \eqref{e:choiceHiggsMass} for $\hat\mass$
is clearly independent of the choice of cutoff $\rho$.

The last statement is then an immediate consequence of Proposition~\ref{prop:equalProcesses} and of the way in which we
derived the equations for the $Y_\eps^{(i)}$ in the first place.
In particular, the process $g$ is obtained as the limit as $\eps \to 0$ of
the processes $g_\eps^{(2)}$ solving \eqref{e:equationG} with initial condition $g_\eps^{(2)}(0)=\id$.
Note in particular that, since we always take
$A \in \CC^{-\kappa}$ and $B \in \CC^{2-\kappa,-\kappa}$, \eqref{e:equationG}
is classically well-posed with values in $\CC^{2-\kappa}$ even in the limit $\eps \to 0$
and one can take as the stopping time $\tau$ the first time at which the 
$\CC^{1/2}_x$-norm of $U = \Repr(g)$ exceeds some fixed threshold.

Since $(Y,h,U)$ satisfies a singular SPDE with polynomial right-hand side,
one can show that there exist constants $c>0$ and $a>0$ such that, for 
every $R > 1$, every admissible model $Z$ with $\$Z\$ \le R$ and every initial data
$\|Y_0\|_{\CC_x^{-\kappa}} + \|B\|_{\CC^{2-\kappa,-\kappa}} \le R$
one has $\tau \ge c R^{-a}$. Since our model of interest comes from the
BPHZ lift of a Gaussian model, it admits moments of all orders and
the claim follows at once. 
\end{proof}

%
\subsection{Bounds from gauge invariant norms}
\label{sec:bounds_from_covariant_norms}

In order to prove that the gauge invariant Markovian dynamic from \cite{CCHS22} admits global solutions,
the last ingredient we need is a purely analytical result showing that the a priori bounds obtained
in Section~\ref{sec:conclusions} are sufficiently strong to control the solutions to the process
constructed in Theorem~\ref{theo:localFPRegStruc}. We exploit the fact that we are only trying to show that the two processes
are equal modulo gauge transformation, so we will construct a suitable gauge in which the 
desired control holds.

As usual, we fix a compact Lie group $\liegroup$ with Lie algebra $\fg$, endowed with an
$\Ad$-invariant norm.
We recall from Definition~\ref{main:def-invnorm} 
that the invariant $\|\place\|_{\cov,-\kappa}$-``norms'' of smooth pairs 
$(B_0,\higgs_0)$, where $B_0$ is a 
$\liegroup$-equivariant connection on the trivial principal $\liegroup$-bundle over $\T^2$ and 
$\higgs_0 \colon \T^2 \to \V$ is a smooth function, are given by
\begin{equ}[e:covariantNorm]
\|B_0,\higgs_0\|_{\cov,-\kappa} = \sup_{s \in (0,1]} s^{\frac{\kappa}2}\bigl(\|\higgs(s)\|_{L^\infty} + s^{\frac12} \|F(s)\|_{L^\infty} + s\|\covd F(s)\|_{L^\infty}\bigr)\;.
\end{equ}
Here $F(s) = F_{B(s)}$ denotes the curvature tensor for the solution $B(s)$ at time $s$ 
to the Yang--Mills gradient flow with
initial condition $B_0$ and $\higgs(s)$ denotes the solution to the covariant heat equation
\begin{equ}
\d_s \higgs(s) = \Delta_{B(s)} \higgs(s) \;,\qquad \higgs(0) = \higgs_0\;.
\end{equ}
The quantity given by \eqref{e:covariantNorm} is of course not a genuine norm since the relation 
$(B_0,\higgs_0) \mapsto (F,\higgs)$ is nonlinear. It is indeed gauge invariant since, for any gauge
transformation $g$, setting $(B_0^g, \higgs_0^g) = g \act (B_0,\higgs_0)$, one has 
$\higgs^g(s) = \rep(g) \higgs(s)$ and similarly for $F^g$ and $\covd F^g$.

One then has the following bound which is crucial for our approach since it allows us to 
compare the solutions constructed in this article with those constructed in \cite{CCHS22,CCHS22+}.

\begin{theorem}\label{theo:cov2norm}
Given $\kappa \in (0,\frac12)$, there exist constants $c, C>0$ such that the following is true.
For every smooth pair $(B_0,\higgs_0)$ as above, there exists a gauge transformation $g$ such that
\begin{equ}[e:apriori]
\|g\act (B_0,\higgs_0)\|_{\CC^{-\kappa}_x} + \sup_{s \in (0,1]}s^{\kappa} \bigl(\lVert\div(g\act  B(s))\rVert_{L^\infty} + s^{1\over2}\lVert \d\div(g\act  B(s))\rVert_{L^\infty} \bigr)  \le C (1 + \|B_0,\higgs_0\|_{\cov,-\kappa})^c\;.
\end{equ}
\end{theorem}


\begin{remark}\label{rem:extendCkappa}
By Lemma~\ref{lemma:invnorm-well-def}, one can remove the (qualitative) smoothness
assumption and replace it by the a priori knowledge that $B_0$ and $\higgs_0$
belong to $\CC^{-\kappa}_x$.
\end{remark}

\begin{remark}
In fact, writing $\tilde B = g\act B$, one also has $\tilde B \in \CC^{\gamma,-\kappa}$ for every $\gamma < 2$ (see \cite[Def.~6.2]{H14} for a definition of 
these spaces; the chosen scaling is the parabolic one and one uses the polynomial regularity structure). 
Since the Yang--Mills gradient 
flow is not strictly parabolic (we do not add the DeTurck term), it is not completely obvious at first sight
that this follows from parabolic regularity. However,
the evolution equation for $\tilde B$ can be written as
a strictly parabolic equation driven by a nonlinearity of lower order, 
plus the forcing term
$\d_\ell \div \tilde B$, see \cite[p.10]{Oh2014} or \cite[Eq.~1.5]{CCHS22}. 
If we consider $\div \tilde B$ as \textit{given} and
use the a priori bounds of \eqref{e:apriori}, which one can think of as some kind of ``approximate
Coulomb gauge'', then the claim $\tilde B \in \CC^{\gamma,-\kappa}$ does follow from standard parabolic
regularity estimates.
\end{remark}

\begin{remark}
Theorem~\ref{theo:cov2norm} was part of the ``community experiment'' of 
\href{https://1stproof.org/second-batch.html}{First Proof batch~2}. 
One of the proofs \href{https://icarm.zulipchat.com/user_uploads/77843/TxteSbs9T-cYuRBPFbu3bsK1/Alexander-C3proof.pdf}{generated during that experiment} suggests that \eqref{e:apriori} holds in fact with $c=2$.
\end{remark}

The main workhorse for the proof of Theorem~\ref{theo:cov2norm} is the following 
quantitative form of Uhlenbeck compactness in the case when $\liegroup$ is simply connected. We will later 
show how to reduce ourselves to that case.

\begin{proposition}\label{prop:Uhlenbeck}
In the setting of Theorem~\ref{theo:cov2norm} but for simply connected $\liegroup$, there exist 
constants $c,C>0$ such that,
for every smooth connection $B$ with curvature $F$, 
there exists a gauge transformation $g$ such that $\|g\act B\|_{L^\infty} \le C (1 + \|F\|_{L^\infty})^c$. 
\end{proposition}

\begin{remark}
At first glance, one has the impression that, by choosing $p$ large and using the embeddings
$L^\infty \subset L^p$ (for $F$) and $W^{1,p} \subset L^\infty$ (for $g\act B$), 
this should be an immediate consequence of keeping track of constants
in the proof of Uhlenbeck's weak compactness theorem \cite{U82,U82b}, see
the detailed exposition in \cite[Thm~7.1]{BookCompactness}. This does
unfortunately not seem to be the case for the following reason. When considering a sequence of 
connections, that proof first 
exhibits a candidate limit which is such that, for a covering of the torus
by small patches, one can find on each patch a gauge transformation leading to a good control 
of the supremum norm of that limit. The argument then proceeds to show that, for connections with
curvature close to that of the limit, one can find gauge transformations such that their \textit{differences}
with the limit are controlled in $W^{1,p}$. In particular, one never needs to get any quantitative 
control on the limit in any global gauge! We will therefore use a different patching argument, which
uses in a crucial way the fact that our base space is a two-dimensional manifold. (Our argument
would still work in three dimensions, but it breaks down in dimension four.)
\end{remark}

Before turning to the proof of this result, we state two relatively 
straightforward extension results, the second of which makes crucial use of
the fact that $\liegroup$ is simply connected.

\begin{lemma}\label{lem:interp}
For every Lipschitz continuous $f_0 \colon [0,1] \to \liegroup$ there exists 
$f \colon [0,1]^2 \to \liegroup$ such that $f(0,x) = f_0(x)$ and $f(1,x) = e$ for all $x\in [0,1]$.
Furthermore, there exists a constant $C$ such that $\Lip(f) \le C(1+ \Lip(f_0))$.
\end{lemma}

\begin{proof}
Since $\liegroup$ is compact, the torus theorem \cite[Thm~11.9]{Hall} shows that there exists 
$h \in \fg$ such that $\exp(h) = f_0(0)$. It then suffices to set for example
$f(t,x) = f_0((1-t)x)\exp(-th)$ and to note that since $\liegroup$ is compact there exists a 
constant $C$ such that $h$ can always be chosen such that $|h| \le C$.
\end{proof}

\begin{lemma}\label{lem:interpCircle}
Assume that $\liegroup$ is compact and simply connected. 
Then, for every Lipschitz continuous $f_0 \colon S^1 \to \liegroup$, there exists 
$f \colon [0,1] \times S^1 \to \liegroup$ such that $f(0,x) = f_0(x)$ and $f(1,x) = e$ for all $x\in S^1$.
Furthermore, there exists a constant $C$ such that $\Lip(f) \le C(1+ \Lip(f_0))$.
\end{lemma}

\begin{proof}
Let $K>0$ be such that $\exp \colon \fg \to \liegroup$ is a diffeomorphism on the ball $B(K)$ of radius $K$ 
and such that both $\exp$ and $\exp^{-1}$ have Lipschitz constant less than $2$ on $B(K)$, respectively
$\exp(B(K))$. (In particular
$\exp^{-1}$ is a smooth function on the ball of radius $K/2$ around the neutral element of $\liegroup$.)
Fix a finite set $V \subset \liegroup$ that is sufficiently dense such that every $g \in \liegroup$ is at distance at
most $K/8$ from some element $v(g) \in V$. For every $v \in V$, we also fix $h_v \in \fg$ with $|h_v|\le C$
such that $v = \exp(h_v)$. Finally, for any pair $(u,v) \in V$ with $d(u,v) \le K/2$, we
write $g_{uv} \colon [0,1] \to \liegroup$ for the Lipschitz continuous loop given by
\begin{equ}[e:basicLoop]
g_{uv}(x) = 
\left\{\begin{array}{cl}
	\exp(3xh_u) & \text{if $3x \in [0,1]$,} \\
	\exp\bigl((3x-1)\exp^{-1}(vu^{-1})\bigr)u & \text{if $3x\in [1,2]$,} \\
	\exp((3-3x)h_v) & \text{if $3x\in [2,3]$.}
\end{array}\right.
\end{equ}
Since $\liegroup$ is simply connected, we can extend $g_{uv}$ to a Lipschitz continuous
map $g_{uv} \colon [0,1]^2 \to \liegroup$ such that $g_{uv}(0,x)$ is given by \eqref{e:basicLoop} and such that
$g_{uv}(t,0) = g_{uv}(t,1) = g_{uv}(1,x) = e$. (The fact that a continuous such 
map exists is the definition of simple connectedness and it can be made Lipschitz
continuous by a `Lipschitz' version
of Whitney's approximation theorem, see for example \cite[Chap.~6]{Manifolds}.)

Writing $M = 1+\Lip(f_0)$, partition now $S^1$ into $\bigO(M)$ intervals $I_k = [x_k,x_{k+1}]$ 
of equal length at most $K/(4M)$, 
so that in particular, for every $k$, $d(f_0(x), f_0(y)) \le K/4$ for all $x,y \in I_k$.
Let furthermore $v_k \in V$ be such that $d(v_k, f_0(x_k)) \le K/8$, so that $d(v_k,v_{k+1}) \le K/2$, 
and set $h_k = h_{v_k}$.
We then write $f^{(1)} \colon S^1 \to \liegroup$ for the (unique) continuous map such that 
$f^{(1)}(x_k) = v_k$ and such that its restriction to each $I_k$ is a constant speed geodesic
of length less than $K/2$. In particular, for $y \in [x_k, (x_k+x_{k+1})/2]$ one has 
\begin{equ}[e:distf01]
d(f^{(1)}(y),f_0(y)) \le \f12 d(v_k, v_{k+1}) + d(v_k,f_0(x_k)) + d(f_0(x_k), f_0(y))
\le
\frac K4 + \frac K8 + \frac K8 \le \frac K2\;,
\end{equ}
which holds for all $y$ by possibly swapping the roles of $v_k$ and $v_{k+1}$.
Note that one has again $\Lip(f^{(1)}) = \bigO(M)$. 
For $t \in [0,1/3]$, we then set
\begin{equ}
f(t,y) = \exp \bigl(3t \exp^{-1}(f^{(1)}(y)f_0(y)^{-1})\bigr) f_0(y)\;.
\end{equ}
Note that one has indeed $\Lip(f) = \bigO(M)$ as a consequence of the fact that, for $\liegroup$-valued functions,
one has $\Lip(f g) \lesssim \Lip(f) + \Lip(g)$, as well as the fact that both $\exp$ and $\exp^{-1}$ have
Lipschitz constant $2$ on their domain of application by \eqref{e:distf01} and the definition of $K$.

For $t \in [1/3,2/3]$ and $x \in I_k$, we then set $f(t,x) = f_k^{(2)}(3t-1,(x-x_k)/|I_k|)$,
where $f_k^{(2)}$ is defined by
\begin{equ}
f_k^{(2)}(t,x) = 
\left\{\begin{array}{cl}
	\exp((3x+1-t)h_k) & \text{if $3x \in [0,t]$,} \\
	\exp\bigl(\frac{3x-t}{3-2t}\exp^{-1}(v_{k+1} v_k^{-1})\bigr) v_k & \text{if $3x\in [t,3-t]$,} \\
	\exp((4-3x-t)h_{k+1}) & \text{if $3x \in [3-t,3]$.} 
\end{array}\right.
\end{equ}
Since $f_k^{(2)}(0,\place)$ is the constant speed geodesic from $v_k$ to $v_{k+1}$ and
since the $f_k^{(2)}$ are all Lipschitz continuous with Lipschitz constant of order $1$
(independent of $M$),
this definition of $f$ is still Lipschitz continuous with $\Lip(f) = \bigO(M)$.

Furthermore we note that, for every $k$, one has $f_k^{(2)}(1,\place) = g_{v_k,v_{k+1}}$
defined as in \eqref{e:basicLoop}.
For $t \in [2/3,1]$ and $x \in I_k$, it therefore suffices to set
$f(t,x) =  g_{v_k,v_{k+1}}(3t-2,(x-x_k)/|I_k|)$, which concludes the proof.
\end{proof}

\begin{proof}[Proof of Proposition~\ref{prop:Uhlenbeck}]
We set $M = 1 + \|F\|_{L^\infty}$. 
By the proof of weak Uhlenbeck compactness as exposited in \cite[Chap.~7]{BookCompactness}, there exists an 
exponent $q$ such that, if one divides
the torus into a collection $A$ of $\bigO(M^{2q})$ overlapping squares of sidelength about $M^{-q}$, one can find for each $a \in A$ a gauge transformation $g_a \colon a \to \liegroup$ such that 
$\|g_a \act B \restr a\|_{W^{1,p}}\lesssim 1$. In particular one has $\|g_a \act B \restr a\|_{L^\infty}\lesssim 1$
by choosing $p > 2$. To be more precise regarding the locations of the squares, we place the centres of the
squares on a regular grid 
of grid size $\delta \approx M^{-q} \in \N^{-1}$ and we choose them of sidelength $4\delta/3$, so that neighbouring squares overlap on a strip of size $\frac\delta3 \times \frac{4\delta}3$.

For any overlapping pair of squares $a,b\in A$, set $g_{ab} \eqdef g_a g_b^{-1}$. By 
\cite[Lemma~A.8]{BookCompactness},
the $\fg$-valued one-forms  $h_{ab} = dg_{ab} g_{ab}^{-1}$ satisfy $\|h_{ab}\|_{W^{1,p}}\le 1$.
In particular, it follows that the transition functions $g_{ab}$ are Lipschitz continuous with 
Lipschitz constants of order $1$, independently of $M$.

Assume now that one has $a = I_a\times J$ and $b = I_b \times J$ with the overlap consisting of the
rectangle $I \times J$ with $I = I_a \cap I_b$ (to fix ideas, say that $a$ is to the left of $b$). 
By rescaling space suitably, we can assume that $I = J = [0,1]$
and that $g_{ab}$ is defined on $[0,1]^2$. We then use Lemma~\ref{lem:interp} to construct a Lipschitz
function $f \colon [0,1]^2 \to \liegroup$ which agrees with $g_{ab}$ on $\{0\} \times [0,1]$
and equals the neutral element $e$ on $\{1\} \times [0,1]$ and we set $g_{a\cap b}(x) = f(x) g_{b}(x)$. 
If we now set $\hat a = a \cup b$ and
define $g_{\hat a}$ by $g_a$ on $a \setminus b$, $g_b$ on $b \setminus a$ and $g_{a \cap b}$ on
$a\cap b$, we see that 
\begin{equ}
\|g_{\hat a} \act B \restr \hat a\|_{L^\infty} \le \|g_a \act B \restr a\|_{L^\infty} \vee \bigl(\|g_b \act B \restr b\|_{L^\infty} + \Lip(f)\bigr)\;.
\end{equ}
Since $g_{ab}$ has Lipschitz constant of order $1$, it furthermore follows from the construction of 
Lemma~\ref{lem:interp} that $\Lip(f) \lesssim M^q$. 

We now divide the collection $A$ into horizontal strips and, for each of the strips, 
we sequentially perform the above construction for 
all pairs of neighbouring squares within the strip. 
At this stage, we then have a collection $\hat A$ of strips and, for each
strip $a \in \hat A$, a gauge $g_a$ such that $\|g_{a} \act B \restr a\|_{L^\infty} \lesssim M^q$.
Furthermore, we claim that the transition functions $g_{ab}$ between neighbouring strips are such that 
$\Lip(g_{ab}) \lesssim M^q$. This is because they either consist of the original transition functions
(which have Lipschitz constants of order $1$) or they are of the form $f_1 g_{a_1}g_{a_2}^{-1} f_2^{-1}$
for $f_i$ some of the interpolating functions constructed above. Since one has $\Lip(fg) \lesssim \Lip(f) + \Lip(g)$ and since the interpolating functions have Lipschitz constants of order at most $M^q$, the
claim follows. 

It now remains to patch the strips together. This is done again by the same construction, but
this time we have to use Lemma~\ref{lem:interpCircle} since the boundaries of the strips consist
of circles. We furthermore lose another factor $M^q$ since this is the reciprocal width of the overlaps
of neighbouring strips. Once all the strips have been patched together, this indeed yields a 
gauge $g$ such that $\|g \act B\|_{L^\infty} \lesssim M^{2q}$, as claimed.
\end{proof}

\begin{remark}
With a little bit of additional effort, one could of course get bounds on $g\act B$ in $W^{1,p}$ and
not just in $L^p$, but since we obtain much better bounds by a bootstrapping argument anyway 
we do not worry about this.
\end{remark}

\begin{remark}
This proof relies on the fact that
$\liegroup$ is simply connected, which is a situation that we can reduce ourselves to. 
The same argument should work in arbitrary dimension $d \ge 2$ provided that 
the homotopy group 
$\pi_{d-1}(\liegroup)$ is trivial. While this is the case for \textit{every} compact Lie group 
in dimension $3$ (see for example
\cite[Prop.~V.7.5]{BtD85}), it clearly fails for $\liegroup = SU(2) \simeq S^3$ when $d=4$.
\end{remark}

Before we turn to the proof of Theorem~\ref{theo:cov2norm}, we state two additional results which
are crucial for removing the assumption that $\liegroup$ is simply connected.

\begin{lemma}\label{lem:Uhlenbeck-abelian}
The conclusion of Proposition~\ref{prop:Uhlenbeck} also holds in the special case
$\liegroup = U(1)$.
\end{lemma}

\begin{proof}
Here it will be convenient to assume that our torus $\T^2$ has size $1$ rather than $2\pi$.
Setting $\liegroup = \R / \Z$ and $\fg = \R$, one identifies as usual gauge transformations with 
functions $g \colon \R^2 \to \R$ such that $g(x+y) - g(x) \in \Z$ for $y \in \Z^2$
and we have $g \act B = B - \nabla g$. Writing $k \in \Z^2$ for an element
such that $|k - \bar B| \le 1$ where $\bar B$ is the spatial average of $B$, 
it then suffices to choose 
\begin{equ}
g(x) = \Delta^{-1} \div B(x) + k x\;.
\end{equ} 
The claim then follows from the fact that, by the Helmholtz decomposition and the
identity $F_{12} = \curl B$, one has
$g\act B = (\bar B - k) + \nabla^\perp \Delta^{-1} F_{12}$.
\end{proof}

Given a Lie group of the type $\liegroup = \liegroup_1 \times \liegroup_2$, we write its elements 
as $g = g_1 \oplus g_2$ and similarly for its Lie algebra. We then have the following.
\begin{lemma}\label{lem:Uhlenbeck-product}
Given a connection one-form $B = B_1 \oplus B_2$ and a gauge transformation
$g = g_1 \oplus g_2$, one has
\begin{equ}
g\act B = (g_1 \act B_1) \oplus (g_2 \act B_2)\;.
\end{equ}
\end{lemma}

\begin{proof}
It follows immediately from the fact that $\Ad_g h = (\Ad_{g_1} h_1) \oplus (\Ad_{g_2} h_2)$
and $dg\,g^{-1} = dg_1\,g_1^{-1} \oplus dg_2\,g_2^{-1}$.
\end{proof}

We now have all the ingredients in place for the proof of the main result of this section.

\begin{proof}[Proof of Theorem~\ref{theo:cov2norm}]
In a first step, we remove the assumption that $\liegroup$ is simply connected in Proposition~\ref{prop:Uhlenbeck}.
We may however assume that $\liegroup$ is connected since we can always decide to only consider
gauge transformations taking values in the connected component of the identity element. It then follows from the structure theorem 
for compact Lie groups \cite[Thm~V.8.1]{BtD85} that 
$\liegroup\simeq (T\times K)/Z$, where $T$ is a torus, $K$ is a simply-connected, semisimple, compact Lie group, 
and $Z$ is a finite subgroup of the centre of $T\times K$. We now let $\Hliegroup\coloneq T\times K$ and let $\frkh$ 
be its Lie algebra. Since $Z$ is finite, $\liegroup$ admits a fundamental domain in $\Hliegroup$ with non-empty interior,
so that one has a canonical identification $\frkh \simeq \fg$. In particular, writing $\pi \colon \Hliegroup \to \liegroup$
for the canonical projection and setting $g = \pi\circ h$ for some smooth $h \colon \T^2 \to \Hliegroup$,
one has
\begin{equ}[e:identification]
\Ad_g = \Ad_h\;,\qquad dg\,g^{-1} = dh\,h^{-1}\;.
\end{equ}

By Proposition~\ref{prop:Uhlenbeck}, Lemma~\ref{lem:Uhlenbeck-abelian}, and Lemma~\ref{lem:Uhlenbeck-product}, 
given an $\fh$-valued connection one-form $B$, we can find an $\Hliegroup$-valued gauge transformation $h$ such that 
the bound $\|h\act B\|_{L^\infty} \le C (1 + \|F\|_{L^\infty})^c$ holds. Interpreting $B$ as a $\fg$-valued connection, it then
follows at once from \eqref{e:identification} that the conclusion of Proposition~\ref{prop:Uhlenbeck}
holds with $g = \pi\circ h$.

In a second step, we show that, writing $M = 1+\|B_0,\higgs_0\|_{\cov,-\kappa}$, 
there exists a constant $c$ (independent of $M$) 
and a gauge transformation $g$ such that $\|g \act B_1\|_{\CC^2_x} \lesssim M^c$.

For this, note first that by Proposition~\ref{prop:Uhlenbeck}, we can find gauge transformations
$g_\tau$ such that, uniformly over $\tau \in [1/2,1]$, one has $\|g_\tau \act B_{\tau}\|_{L^\infty} \lesssim M^c$. Setting $\hat B_0 = g_\tau \act B_{\tau}$ for some fixed $\tau$,
consider then the Yang--Mills--DeTurck heat flow $\hat B$ with initial condition $\hat B_0$.
Simple scaling arguments show that this equation is locally well-posed for initial conditions
$\hat B_0 \in \CC^{\alpha}_x$ for every $\alpha > -1/2$ which in particular
includes $\hat B_0 \in L^\infty$.
Since the nonlinearity for this equation is a linear combination of terms of type $B \d B$ and $B^3$,
a standard Picard iteration argument for semilinear parabolic PDEs shows that, after a time
$t$ of order $M^{-4c}$, one has $\|\hat B_t\|_{C^k_x} \lesssim M^{c(1+2k)}$.
By the gauge equivalence of the Yang--Mills--DeTurck heat flow and the Yang--Mills gradient flow (see for
example \cite{D85} or \cite[Sec.~4.1]{Struwe}),
there exist gauge transformations $\hat g_t$ such that $\hat B_t = \hat g_t \act (g_\tau \act B_{\tau+t})$.
It now suffices to choose $\tau = 1 - CM^{-4c}$ for a suitable constant so that one can set
$\tau + t = 1$, thus yielding the desired bound
$\|g \act B_1\|_{C^2_x} \lesssim M^{5c}$ for $g = \hat g_t \cdot g_\tau$.

From now on, we will allow ourselves to change the precise value of the exponent $c$ from one step
to the next in our proofs. We furthermore place ourselves in the (time independent!) 
gauge $g$ we just constructed, so we write $B_0$ and $B(s)$ rather than $g \act B_0$ and $g \act B(s)$ 
from now on, and similarly for $\higgs_0$ and $\higgs(s)$.
We also now take for granted the bounds we just obtained, namely
\begin{equ}[e:baseCase]
\sup \bigl\{\|B(1)\|_{L^\infty},\, \|\d B(1)\|_{L^\infty},\, \|\d^2 B(1)\|_{L^\infty}\bigr\} \lesssim M^c\;.
\end{equ}
(We will sometimes write `time' as an argument of $B$ and sometimes as a subscript. We
hope that this will not cause confusion: time values will always be $s,t$, while spatial
indices will be $j$, $\ell$, etc.)

We now apply Lemma~\ref{lemma:ymgf-approximate-coulomb-gauge}. 
By the definition \eqref{e:covariantNorm} of $\|\place\|_{\cov,-\kappa}$, one has $M_0 + M_1 \lesssim M$,
so that, when combined with \eqref{e:baseCase}, 
\eqref{eq:ymgf-B-negative-regularity-estimate} and \eqref{eq:ymgf-div-B-estimate} do indeed
yield the desired bound \eqref{e:apriori} for $B$. In fact, \eqref{eq:ymgf-B-Linfty-estimate} 
also yields the bounds
\begin{equ}[e:boundBval]
\|B(t)\|_{L^\infty} \lesssim M^c t^{-\kappa/2}\;,\quad
\norm{\div B}_{L^\infty} = \|\covd^\ell B_\ell\|_{L^\infty} \lesssim M^c t^{-\kappa}\;,\quad
\norm{\d\div B}_{L^\infty} \lesssim M^c t^{-\kappa-\f12}\;.
\end{equ}

Finally, in order to bound the $\CC^{-\kappa}_x$ norm of $\higgs_0$, we use the bound
\begin{equ}
\|\higgs_0\|_{\CC^{-\kappa}_x} \lesssim \sup_{t \in (0,1]} t^{\kappa/2}\|\fp_t \higgs_0\|_{L^\infty}\;.
\end{equ}
The proof of this fact is a simple variant on the argument given for example
in \cite[Thm~12.4]{LorenzoFrancesco}. In order to bound $\fp_t \higgs_0$ we use the expression of
$\Delta_B$ in coordinates (see for example \cite[Eq.~1.7]{CCHS22+}) and the variation of constants 
formula to write
\begin{equ}
\fp_t \higgs_0 = \higgs(t) + \int_0^t \fp_{t-s} \bigl((\div B)\higgs - B_j^2\higgs - 2 \d_j (B_j \higgs)\bigr)(s)\,ds\;.
\end{equ}
Since $\|\higgs(s)\|_{L^\infty} \lesssim M s^{-\kappa/2}$ by assumption and since the operator norm
of $\fp_{t-s}\d_j$ from $L^\infty$ to itself is of order $(t-s)^{-1/2}$, it is straightforward 
from \eqref{e:boundBval} to bound the integral appearing in this 
expression by (some multiple of) $M^c (t^{1-\frac{3\kappa}2} + t^{\frac12-\kappa})$.
The required bound $\|\higgs_0\|_{\CC^{-\kappa}_x} \lesssim M^c$ follows at once, completing the proof.
\end{proof}

\subsection{Proof of main result}

This finally leads to the main result of this article.

\begin{corollary}\label{cor:main}
The Yang--Mills--Higgs measure on $\fO_{\kappa}^\YMHs$ exists and is unique.
\end{corollary}

\begin{proof}
Writing $\hat\CP_t$ for the Markov semigroup on $\fO_{\kappa}^\YMHs \sqcup \{\skull\}$ given by Proposition~\ref{prop:semigroupOrbits}, we want  to show that it admits a unique
invariant probability measure $\mu_\YMH$ on $\fO_{\kappa}^\YMHs$, which is then 
the sought-after Yang--Mills--Higgs measure.

We first show that $\hat\CP_t$ is conservative in the
sense that $\hat\CP_t(Y_0, \fO_{\kappa}^\YMHs) = 1$ for all 
$Y_0 \in \fO_{\kappa}^\YMHs$. For this, note first that, by
Theorem~\ref{thm:global_soft}, there exists an increasing sequence 
of stopping times
$t_n \to \infty$ (with $t_0=0$ by convention), a sequence of 
gauge transformations $g_n \in \CG_{1-\kappa}$, and a 
$\CC_x^{-\kappa}$-valued process $Y = (A,\Phi)$ such that
\begin{enumerate}[label=(\roman*),itemsep=0ex,topsep=0.5ex]
\item For every $n \ge 0$, the solution $(B_n, \higgs_n)\eqdef Y_{t_n}$ belongs to the 
space $\initialclass_1$ defined in \eqref{main:eq-initial-class}.
\item On $[t_n, t_{n+1})$, $Y = \lim_{\eps \to 0} Y_\eps$
with $Y_\eps$ solving \eqref{e:covariant_SPDE}, where $B$ is the YM gradient
flow with initial condition $B_n$.
\item One has $Y_{t_n} = g_n \act Y_{t_n^-}$.
\item\label{pt:unifBounds} For every $n \ge 0$, one has $Y \in C([t_n,t_{n+1}), \CC_x^{-\kappa})$
and one has $\sup_{t>0} \E \|Y_t\|_{\cov,-\kappa}^p < \infty$.
\end{enumerate}
(Here we changed convention at the jump times and chose $Y$ to be right-continuous
there.)
By Proposition~\ref{prop:processesEqual}, in particular \eqref{e:processesEqual},
we can then construct a white noise $\hat\boldxi$ with respect to the same filtration as
the noise $\boldxi$ driving $Y$, an $\Omega_{\kappa}^\YMHs$-valued
process $\hat Y$, a sequence of stopping times $\sigma_n \to \infty$,
and an adapted sequence of gauge transformations $\hat g_n \in \CG_{1-\kappa}$ such that
\begin{enumerate}[label=(\roman*),itemsep=0ex,topsep=0.5ex] \setcounter{enumi}{4}
\item On $[\sigma_n, \sigma_{n+1})$, $\hat Y$ is given by the solution
map constructed in Theorem~\ref{theo:localFPRegStruc}. In particular, one has 
$\hat Y = \lim_{\eps \to 0} \hat Y_\eps$ 
with $\hat Y_\eps$ solving \eqref{e:wantedSPDE}.
\item\label{pt:gaugeRS} One has $\hat Y_{\sigma_n} = \hat g_n \act \hat Y_{\sigma_n^-}$.
\item\label{pt:equiv} One has $[\hat Y_t] = [Y_t]$ for all $t \ge 0$, almost surely.
\end{enumerate}
More precisely, the construction of $\hat Y$ goes as follows. 
Write $t_\star = \{t_n\}_{n \ge 0} \subset \R$.

Starting from a given initial 
data $Y_0\in \Omega_{\kappa}^\YMHs$, assume that we are given the process $Y$.
Take then $\sigma_1 = \tau \wedge \inf \{t > 0\,:\, t \in t_\star\}$, where $\tau$ is the stopping time
given by Proposition~\ref{prop:processesEqual} and set $\hat\boldxi = g\act \boldxi$
on $[0,\sigma_1)$ so that \eqref{e:processesEqual} guarantees that one does indeed have
$[\hat Y] = [Y]$ on that interval. We then simply set $\hat Y_{\sigma_1} = Y_{\sigma_1}$, 
which can be made consistent with point~\ref{pt:gaugeRS} by taking $\hat g_1 = g(\sigma_1)^{-1}$ if $\sigma_1 < t_1$ and
$\hat g_1 = g_1\, g(\sigma_1)^{-1}$ otherwise. 
Since one has $g(\sigma_1) \in \CG_{1-\kappa}$ by 
parabolic regularity and since the $g_n$ take values in $\CG_{1-\kappa}$
by Theorem~\ref{thm:global_soft}, 
this does indeed guarantee that one has $\hat Y_{\sigma_1} \in \Omega_{\kappa}^\YMHs$.

One then iterates this construction by taking $\sigma_1$ as the new starting time
(which of course also entails shifting the set $t_\star$ by $\sigma_1$).
The bound \eqref{e:boundComparisonTime}, combined with the fact that we have bounds
on the moments of our model and $Y$ is continuous in time (up to time $t_1$) with values in 
$\CC_x^{-\kappa}$ guarantees that these stopping times cannot accumulate.

This construction makes it clear that the law of the process $\hat Y$ is a 
generative probability measure in the sense of Proposition~\ref{prop:semigroupOrbits}
and that its projection onto the space of gauge orbits $\fO_{\kappa}^\YMHs$
is Markov with semigroup $\hat \CP$. In particular, this shows that 
$\hat \CP$ is conservative.
Using \eqref{e:TVconvergence} combined with parabolic regularisation, it is 
furthermore relatively straightforward to obtain a tightness result.

To show the existence of $\mu_\YMH$, we apply the 
Krylov--Bogoliubov theorem \cite[Thm~3.1.1]{DPZ2} which requires us to show tightness
of $\{\hat \CP_t(Y_0,\place)\}_{t \ge 1}$ on $\fO_{\kappa}^\YMHs$.
By Lemma~\ref{lem:compactEmbedding}, this follows if we can show that, for some
$\bar \kappa < \kappa$ one has
\begin{equ}
P_t(M) \coloneq \P \Bigl(\inf_g \|g \act \hat Y(t)\|_{\CC_x^{\bar\kappa}} \ge M \Bigr) \lesssim M^{-p}\;,
\end{equ}
uniformly over $t \ge 1$. We therefore fix  $\bar \kappa \in (0,\kappa)$ from now on.
A straightforward consequence of the fixed point problem \eqref{e:FPgen} is that
there exist
constants $C>0$ and $a>0$ such that, for every $R > 10$, if $\|X_0\|_{\CC^{-\kappa}_x} \le R$, $\sup_{s \in[0,1]}\$\<IXib>_s\$_{\bar \kappa} \le R$, and $\$Z\$\le R$ (where $Z$ denotes the model used for the fixed point problem), then
for $\delta = R^{-a}$ one has the bound
\begin{equ}
\sup_{t \le \delta} t^{\f14+\f\kappa2}\$X_t\$_{\bar \kappa}
\le CR\;,
\end{equ}
and in particular $\$X_\delta\$_{\bar \kappa} \le CR \delta^{-\f14-\f\kappa2} \le R^b$ for some 
exponent $b > 0$. Here we used the fact that $\CC^{1/2}_x$ embeds into $\Omega_{1-\bar\kappa}^1$ by Remark~\ref{rem:embedding}.

As a consequence of the fact that, for any choice of
measurable selection $S \colon \fO_{\kappa}^\YMHs \to \Omega_\kappa^\YMHs$, one
has  $\hat \CP_t \ge (\Pi^*\CP_\delta) S^*\hat \CP_{t-\delta}$ by Proposition~\ref{prop:semigroupOrbits}, 
it therefore follows that, setting $R = M^{1/b}$, one has
\begin{equ}
P_t(M) \le \P\Bigl(\inf_{g \in \CG_{1-\kappa}} \|g\act\hat Y_{t-\delta}\|_{\CC^{-\kappa}_x} > R\Bigr)
+ \P \Bigl(\sup_{s \le 1}\$\<IXib>_s\$_{\bar \kappa} > R\Bigr)
+ \P \bigl(\$Z\$ > R\bigr)\;.
\end{equ}
The first term converges to $0$ faster than any power of $R$ (and therefore of $M$) by
points \ref{pt:equiv} and \ref{pt:unifBounds} above combined 
with Theorem~\ref{theo:cov2norm}, while the remaining two terms do so by standard probabilistic estimates.

The fact that the 
semigroup $\hat\CP_t$ is Feller
(in the sense that transition probabilities are continuous in the weak convergence topology)
is a consequence of the stronger strong Feller property, which
holds by combining \eqref{e:TVconvergence} and the results of \cite{StrongFeller} 
in the same way as \cite[Prop.~10.1]{IlyaHaoYM}.

We finally turn to the uniqueness of $\mu_\YMH$. It follows immediately from
Theorem~\ref{theo:fullSupport} that $\mu_\YMH$ has full support. However, the strong Feller property
implies that any two ergodic invariant measures for $\hat\CP_t$ have disjoint support and the claim follows. 
\end{proof}

\section{Symbolic index}
\renewcommand{\leftmark}{\textsc{Symbolic index}}
\renewcommand{\rightmark}{\textsc{Symbolic index}}

\begin{center}
\renewcommand{\arraystretch}{1.1}
\begin{longtable}{lll}
\toprule
Symbol & Description & Page\\
\midrule
\endfirsthead
\toprule
Symbol & Description & Page\\
\midrule
\endhead
\bottomrule
\endfoot
\bottomrule
\endlastfoot

$\linear[B]$, $\philinear[B]$ &
Covariant linear stochastic objects & 
\pageref{main:eq-linear-stochastic-A} \\ 

$\linear[\Beps]$, $\philinear[\Beps]$ &
Smooth approximations of covariant linear stochastic objects & 
\pageref{main:eq-linear-stochastic-A} \\ 

$\graylinear[B]$  & 
$\Y$-valued covariant linear stochastic object & 
\pageref{local:notation-Y-valued-objects} \\ 

$\redlinear[B]$  & 
Independent copy of $\graylinear[B]$ & 
\pageref{object:eq-quadratic-dlin-decomposition-original} \\ 

$\pmgraylinear[B][][N]$ & 
Approximation of $\graylinear[B][][][N]$ using heat kernel parametrix & 
\pageref{def:parametrix-approximate-objects}\\

$\linear[B_m]$,
$\philinear[B_m]$ &
Covariant linear stochastic objects in the $m$-th step of the iteration &
\pageref{defn:time_steps} \\ 

$\linear[\Bmod_m,1/2]$, $\linear[\Bmod_m^\star,1/2]$ &
Covariant linear stochastic objects started at the intermediate time &
\pageref{defn:time_steps_mod} \\

$\quadratic[\Beps]$, $\phiquadratic[\Beps]$, 
$\mixedquadratic[\Beps]$ & 
Covariant quadratic stochastic objects & 
\pageref{local:eq-A-quadratic-stochastic} \\  

$\<IXi>,\, \<I'Xi>,\,\<I'Xi2>,\, \<I'XiX>,\,\<IXi^2>, \narrowdots$ & 
Syntactic trees & 
\pageref{e:symbolList} \\ 

$\<IXir>,\, \<I'Xir>,\,\<I'Xi2r>,\, \<I'XiXr>,\,\<IXi^2r>, \narrowdots$ & 
Renormalisation constants of syntactic trees & 
\pageref{e:IXi^2r} \\ 

$\Fnorm{\place}_{T,k}$ &
Curvature norm used in our perturbative arguments & 
\pageref{def:Fnorm} \\ 

$\revFnorm{\place}_{[s_0, t_0], t_1, k}$ & 
Curvature norm used in time reversal arguments & 
\pageref{kernel:def-trep} \\ 

$\YMFnorm{\place}_{I}$ & 
Norm to quantify distance to YM connections & 
\pageref{e:YMFnorm} \\ 

$\invnorm{(\place, \place)}{-\kappa}$ & 
Gauge-invariant Hölder norm without background&
\pageref{main:def-invnorm}\\ 

$\|\place \|_{\covHolder{A, F, T}^{-\kappa}}$, $\|\place \|_{\covHolder{A, T}^{-\kappa}}$  & 
Gauge-invariant Hölder norm& 
\pageref{defn:covariant_Holder} \\ 

$\|\place \|_{\covHolder{A, F}^{-\kappa}}$, $\|\place \|_{\covHolder{A}^{-\kappa}}$  & 
Gauge-invariant Hölder norm on an $A$-dependent time scale & 
\pageref{defn:covariant_Holder_global} \\ 

$\tripn{\Xi_{\Beps}}$ & 
Size of the covariant enhanced data set & 
\pageref{local:def-enhanced-data-set} \\ 	

$\$(\place,\place)\$_{\kappa}$ & 
Norm on space of connection one-forms and Higgs fields & 
\pageref{e:basicState} \\ 

$\stscaling{\mbf{m}}$  & 
Parabolic scaling of multi-index $\mbf{m} \in \N^{\{0\} \cup [d]}$ & 
\pageref{notation:space-time-derivatives}  \\ 

$\parall, \paragg,\narrowdots$ & Paraproduct-type operations & 
\pageref{local:def-paraproduct}\\

$\ll_\theta, \gg_\theta,\narrowdots$ & 
Comparison relations used in the paraproducts & 
\pageref{local:eq-ll-theta} \\ 

$[d]$ & 
Spatial indices & 
\pageref{notation:space-time-indices} \\ 

$\gamma_0$ &
Exponent in $T^{(1)}$ controlling the Higgs part &
\pageref{e:constants_global_iteration} , \pageref{defn:tau} \\

$\gamma_1,\gamma_2$ & 
Parameters appearing in covariant energy & 
\pageref{e:constants_global_iteration} \\

$\gamma_A$ & 
Integral of $A$ along the path $\gamma$ & 
\pageref{e:defgammaA} \\ 

$\Delta_A$ & 
Covariant Laplacian & 
\pageref{intro:eq-YMH-Langevin} \\

$\lossdel$ & 
Parameter &
\pageref{e:small_constants} \\ 

$\reg_1$ & 
Regularity parameter &
\pageref{e:small_constants} \\ 

$\reg_2$ & 
Regularity parameter &
\pageref{e:small_constants} \\ 

 $\stepupperbound$ &
Upper bound of iteration time scales &
\pageref{defn:tau} \\

$\logtimestep$ &
Small parameter controlling $T^{(1)}$ and $T^{(2)}$ &
\pageref{defn:tau} \\

$\zeta$ & 
$\V$-valued space-time white noise & 
\pageref{intro:eq-YMH-Langevin-generic}, \pageref{intro:eq-YMH-Langevin} \\

$\zeta_{\Beps}$ & 
Gauge-covariant smooth approximation of $\zeta$ & 
\pageref{main:def-gauge covariant-approximation-space-time-white-noise} \\ 

$\theta$ & 
Parameter appearing in paraproducts & 
\pageref{e:small_constants}, \pageref{local:def-paraproduct} \\ 

$\theta^\prime$ & 
Parameter larger than $\theta$ & 
\pageref{e:small_constants} \\ 

$\kappa$ & 
Regularity parameter & 
\pageref{e:small_constants} \\ 

$\lambda$ & 
Parameter appearing in weighted solution spaces & 
\pageref{e:small_constants}, \pageref{local:prop-main} \\ 
$\lambda_1$, $\lambda_2$ & 
Parameters larger than $\lambda$ & 
\pageref{e:small_constants}\\ 

$\Xi_{\Beps}$ &
Covariant enhanced data set & 
\pageref{local:def-enhanced-data-set} \\ 

$\boldxi$ & 
$\Y$-valued space-time white noise & 
\pageref{local:notation-Y-valued-objects} \\ 

$\hat \boldxi_{\eps}$ & 
Stochastic forcing appearing in the analysis of $\boldxi_{\Beps}-\boldxi_{\eps}$ & 
\pageref{e:hat-boldxi} \\

$\xi$ &
$\frkg^2$-valued space-time white noise &
\pageref{intro:eq-YMH-Langevin} \\ 

$\xi_{\eps}$ & 
Usual smooth approximation of $\xi$ & 
\pageref{e:mollif} \\ 

$\xi_{\Beps}$ & 
Gauge-covariant smooth approximation of $\xi$ & 
\pageref{main:def-gauge covariant-approximation-space-time-white-noise} \\ 

$\PPi^{(\eps,i)}$ & 
Uncentred models associated with our three main processes & 
\pageref{sec:renConstants} \\ 

$\pi$ & 
Representation of $\liegroup$ on $\Y$ or on a generic space $\genV$ & 
\pageref{e:defDA}, \pageref{e:sys2} \\ 

$\pi_*$ & 
Linearisation of $\pi$ & 
\pageref{e:defDA}, \pageref{e:sys2} \\ 

$\rho$ & 
Cutoff function from the set $\mc{P}$  & 
\pageref{main:eq-mcP} \\ 

$\moll$ & 
Mollifier built from $\rho$  & 
\pageref{e:mollif} \\ 

$\rep$ & 
Orthogonal representation of $\liegroup$ on $\V$ & 
\pageref{intro:assumption-setup} \\ 

$\sigma$ & 
Cutoff function & 
\pageref{def:ptcutoff} \\ 

$\sigma_{\rho,\eps}^2$ & 
Renormalisation constant & 
\pageref{main:def-renormalization-constant} \\ 

$\sigma_{\eps,M,N}^2$ & 
Dyadically localised resonance & 
\pageref{objects:def-localized-resonance}\\

$\tau^{(L)}$ &
Stopping time controlling the size of the stochastic objects &
\pageref{defn:tau_L} \\

$\tau^{(L)}_{1/2}$ &
Minimum of $\tau^{(L)}$ and $T^{(L)}_{1/2}$ &
\pageref{defn:final_timestep} \\

$\tau_m$, $\tau_m^\star$ &
Step lengths in the first global iteration scheme &
\pageref{defn:time_steps} \\

$\taumod_m$, $\taumod_m^\star$ &
Step lengths in the modified iteration scheme &
\pageref{defn:time_steps_mod} \\

$\taumod_{m+1/2}^\star$ &
Increment from $\tmod_m$ to $\tmod_{m+1/2}^\star$ &
\pageref{defn:time_steps_mod} \\

 $\tau_{m+1}^{\star\star}$ &
$\fF_{t_{m-1}}$-measurable time step equal to $\tau_{m+1}$ on good event &
\pageref{defn:star_objects} \\

$\taumod_{m+1}^{\star\star}$ &
$\fF_{\tmod_{m-1/2}^\star}^{\,\xi}$-measurable time step
equal to $\taumod_{m+1}$ on good event &
\pageref{defn:star_objects_mod} \\

$\pt_A$ & 
Parallel transport induced by $A$ & 
\pageref{def:parallel-transport} \\ 

$\itobjects_{m+1}$, $\itobjectsmod_{m+1}^{(1)}$ &
Sum of stochastic terms in the one-step energy expansion &
\pageref{pr:expansion_iteration}, \pageref{pr:expansion_iteration_mod} \\

$\itobjectsmod_{m+1}^{(2)}$ &
Additional stochastic terms from modification of iteration scheme &
\pageref{le:smoothed_lollipop_addition_expansion} \\

$\higgs$ &
Higgs field in SYMH equations & 
\pageref{intro:eq-YMH-energy} \\

$\higgs_\eps$ & 
Solution of smooth approximation of SYMH & 
\pageref{local:eq-SYM-coordinates-eps} \\ 

$\hat \higgs_\eps$ & 
Solution of smooth approximation of SYMH in DeTurck gauge &
\pageref{e:wantedSPDE} \\ 

$\higgs_\eps^{(i)}$ & 
Higgs field in one of the three main processes & 
\pageref{e:second-main-process}\\ 

$\phismoothpart$, $\phiroughpart$ &
Smooth and rough components in the first global iteration scheme &
\pageref{defn:one_step_iteration_scheme} \\

$\phismoothpartmod$, $\phiroughpartmod$ &
Smooth and rough components in modified global iteration scheme &
\pageref{defn:final_timestep} \\

$\higgsdet$ & 
Deterministic term in expansion of $\higgs$ &
\pageref{eq:higgsdet-equation} \\ 

$\higgsdet_m$, $\higgsdet_m^\star$ &
$\higgsdet$ term in the $m$-th iteration step &
\pageref{defn:time_steps} \\

$\higgsdetmod_m$, $\higgsdetmod_m^\star$ &
$\higgsdet$ term in the modified iteration &
\pageref{defn:time_steps_mod} \\

$\rlinhiggsdet$,  $\reghiggsdet$ & 
Decomposition of $\higgsdet$ & 
\pageref{eq:rlinhiggsdet-equation} \\ 

$\shiggsdet$, $\remhiggsdet$ & 
Decomposition of $\reghiggsdet$ & 
\pageref{e:shiggsdet_eq} \\ 

$\chi_N$ & 
Cutoffs used in the dyadic decomposition in time & 
\pageref{prelim:eq-chi-N} \\ 

$\Psi_t(A)$ & 
YM gradient flow with initial data $A$ at time $t$ & 
\pageref{text-Psit} \\ 

$\Omega_{\kappa}^\YMHs$ & 
Space of gauge and Higgs fields & 
\pageref{e:basicState} \\ 

$A$ &
Connection one-form (often in SYMH equations) & 
\pageref{intro:eq-YMH-energy} \\ 

$A_\star$ &
Sequence of Littlewood--Paley blocks corresponding to $A$ & 
\pageref{local:text-Littlewood-Paley-blocks} \\ 

$A_\eps$ & 
Solution of smooth approximation of SYMH & 
\pageref{local:eq-SYM-coordinates-eps} \\ 

$\hat A_\eps$ & 
Solution of smooth approximation of SYMH in DeTurck gauge &
\pageref{e:wantedSPDE} \\ 

$A_\eps^{(i)}$ & 
Connection one-form in one of the three main processes & 
\pageref{e:second-main-process}\\ 

$\Asmoothpart$, $\Aroughpart$ &
Smooth and rough components in the first global iteration scheme &
\pageref{main:eq-covariant-energy}, \pageref{defn:one_step_iteration_scheme} \\

$\Asmoothpartmod$, $\Aroughpartmod$ &
Smooth and rough components in the modified iteration scheme &
\pageref{defn:final_timestep} \\

$\widetilde{A}^{+,\star}(\tmod_m^\star)$ &
$\fF_{\tmod_{m-1/2}^\star}^{\xi}$-measurable connection equal to $\Asmoothpartmod(\tmod_m)$ on good event &
\pageref{defn:star_objects_mod} \\

$\event_m$, $\eventmod_m$ &
Good events where actual and starred quantities coincide &
\pageref{defn:star_objects}, \pageref{defn:star_objects_mod}\\

$\fa$ &
Exponent for $L$-admissibility &
\pageref{e:constants_global_iteration}, \pageref{defn:L_admissible} \\

$\Ad_g$ & 
Adjoint action of $g\in \liegroup$ on $\frkg$  & 
\pageref{intro:assumption-setup} \\ 

$\ad_X$ & 
Adjoint action of $X\in \frkg$ on $\frkg$ & 
\pageref{intro:assumption-setup} \\ 

$\bil$ & 
Linear map $\V \otimes \V \rightarrow \frkg$ & 
\pageref{intro:eq-bil} \\ 

$\hat\bil$ & 
Bilinear map in schematic equation & 
\pageref{e:schematic} \\ 

$B$ & 
Connection one-form (often a YM gradient flow) & 
\pageref{e:YMHF} \\ 

$B^\sharp$ & 
YM gradient flow in infinitesimal DeTurck gauge & 
\pageref{eq:B-sharp}\\ 

$B^+$, $B^-$ & 
Positive and negative regularity terms in $B^\sharp$ & 
\pageref{le:perturbed_YMHF-general}, \pageref{eq:B-minus-equation} \\ 

$B_m$, $B_m^\star$ &
YM gradient flows in the $m$-th iteration step &
\pageref{defn:time_steps} \\

$\Bmod_m$, $\Bmod_m^\star$ &
YM gradient flows in the modified iteration &
\pageref{defn:time_steps_mod} \\

$\YMHFsmooth_m$, $B_m^{+,\star}$ &
YM gradient flow from the smooth part of the initial data &
\pageref{defn:time_steps} \\

$\YMHFsmoothmod_m$, $\Bmod_m^{+,\star}$ &
YM gradient flow from smooth part in modified iteration &
\pageref{defn:time_steps_mod} \\

$B_{m+1}^{+,\star\star}$ &
$\fF_{t_{m-1}}$-measurable Yang--Mills flow equal to $B_{m+1}$ on good event &
\pageref{defn:star_objects} \\ 

$\Bmod_{m+1}^{+,\star\star}$ &
$\fF_{\tmod_{m-1/2}}^{\xi}$-measurable Yang--Mills flow equal to $\Bmod_{m+1}$ on good event &
\pageref{defn:star_objects_mod} \\

$\Cs^\alpha_x$, $\Cs^\alpha_{t,x}$ & 
Spatial and space-time Hölder spaces of regularity $\alpha$ &  
\pageref{def:Holder-space} \\ 

$\Cov_\frkg$, $\Cov_\V$ & 
Covariance operators & 
\pageref{prelim:def-Casimir-induced} \\ 

$\Cov_{\Y}$ & 
Covariance operator in $\Y$ & 
\pageref{local:notation-Y-valued-objects} \\

$\Cov^\frkg_\frkg$, $\Cov^\frkg_\V$, $\narrowdots$ & 
Linear maps induced by covariance operators & 
\pageref{prelim:def-Casimir-induced} \\

$\Crho$ & 
Renormalisation constant & 
\pageref{main:def-renormalization-constant} \\ 

$\renormA$ & 
Renormalisation constant in the $A$-equation & 
\pageref{main:def-renormalization-constant} \\ 

$\renormphi$ & 
Renormalisation constant in the $\higgs$-equation & 
\pageref{main:def-renormalization-constant} \\ 

$C$ & 
Difference between connection one-forms $A$ and $B$ & 
\pageref{local:eq-C} \\ 

$\curl_A$ & 
Covariant curl & 
\pageref{notation:covariant-curl} \\ 

$\CD^{\gamma,\eta}$ & 
Space of modelled distributions & 
\pageref{text:CD} \\ 

$\hPicard_{s,t}$ & 
Space of pairs $(X,Z)$ with $Z \in \CN$ and 
$X \in \CD^{\gamma,\eta}_{s,t}$ &  
\pageref{text:hPicard}\\  

$\covd_A$ & 
Covariant derivative & 
\pageref{e:defDA} \\ 

$\covd_A^{(k)}$ & 
Covariant spatial derivatives of $k$th order & 
\pageref{notation:higher-derivatives} \\ 

$\covd_{A,t,x}^{(m)}$ &  
Covariant space-time derivatives of $m$th order &
\pageref{notation:space-time-derivatives} \\ 

$\covd_{A,y}^{\dagger}$ & 
Covariant derivative after integration by parts in $y$-variable& 
\pageref{notation:dagger}  \\

$\div_A$ & 
Covariant divergence & 
\pageref{notation:covariant-divergence-two-form} \\  

$\Dec_{t}$ & 
Decrease related to covariant energy $\covariantenergy$ & 
\pageref{e:dec_defn}\\

$\Dec_m$, $\Decmod_m$  &
$\Dec_t$ at the $m$-th step in the iteration schemes &
\pageref{e:iteration_deterministic_bounds_iteration_Dec1}, \pageref{pr:expansion_iteration_mod} \\

$\Duh_B$, $\Duh_B^F$ & 
Covariant Duhamel integrals & 
\pageref{def:covariant-duhamel} \\ 

$\Duh_{B,N}$, $\Duh_{B,N}^F$ & 
Littlewood--Paley blocks of covariant Duhamel integrals & 
\pageref{local:Duh-LWP-block} \\ 

$\covariantenergy$ & 
Covariant energy used in global theory & 
\pageref{main:eq-covariant-energy}\\

$e^{s\Delta_{B}}$, $e^{s\Delta^F_{B}}$ & 
Linear propagators for covariant heat equations & 
\pageref{notation:linear-propagation} \\ 

$F_A$ & 
Curvature tensor & 
\pageref{intro:eq-curvature-commutator-identity} \\ 

$\fF_{t_m}$, $\fF_{\tmod_m}$ &
Filtrations up to the stopping times $t_m$, $\tmod_m$ &
\pageref{defn:time_steps}, \pageref{defn:time_steps_mod} \\

$\fF_{\tmod_{m-1/2}^\star}^{\,\xi}$ &
$\fF_{\tmod_{m-1}}$ enlarged by the connection noise up to $\tmod_{m-1/2}^\star$ &
\pageref{defn:star_objects_mod} \\

$\liegroup$ &
Compact Lie group &
\pageref{intro:assumption-setup}\\

$\CG$ &  
Group of smooth gauge transformations & 
\pageref{intro:def-CG} \\ 

$\CG_\alpha$ &  
Closure of smooth gauge transformations in the $\alpha$-Hölder norm & 
\pageref{intro:def-CG} \\

$\fg$ &
Lie algebra of $\liegroup$ &
\pageref{intro:assumption-setup}\\

$g_m$, $g_m^\star$ &
Gauge transformations for the $m$-th iteration step &
\pageref{defn:time_steps} \\

$\widetilde{g}_m$, $\widetilde{g}_m^\star$ &
Gauge transformations in the modified iteration &
\pageref{defn:time_steps_mod} \\

$g\act A, g \act \higgs$ & 
Gauge action on connection one-forms and Higgs fields & 
\pageref{intro:eq-gauge-transform} \\ 

$\gronwall_I(F)$ & 
Gronwall loss & 
\pageref{notation:gronwall-loss} \\ 

$\hol(A,\gamma)$ & 
Holonomy of the connection one-form $A$ along the path $\gamma$ &
\pageref{intro:eq-holonomy-process} \\

$h$ &
Right logarithmic derivative of $g$ & 
\pageref{e:evolhU}\\ 

$\initialclass_T$ & 
Class of initial data & 
\pageref{main:eq-initial-class} \\ 

$\linTerm$ & 
Linear map in schematic equation & 
\pageref{e:schematic} \\ 

$\Lb$ & 
Combination of Lie-bracket terms & 
\pageref{local:eq-Lie-bracket-term} \\

$\mc{L}_X$, $\mc{L}_X^F$, $\narrowdots$ & 
Linear operators in perturbation around YM connections & 
\pageref{def:YMpert-operators} \\ 

$\lossfactor_T(F)$ & 
Loss factor & 
\pageref{local:eq-lossfactor} \\ 

$\LL$ & 
Canonical lift of a smooth noise to a model & 
\pageref{rem:NotationNoise} \\ 

$L$-admissibility &
Parameter controlling the rough components &
\pageref{defn:L_admissible} \\

$L_m$, $L_m^\star$ &
Admissibility parameter in the iteration steps &
\pageref{defn:time_steps} \\

$\Lmod_m$, $\Lmod_m^\star$ &
Admissibility parameter in the modified iteration steps &
\pageref{defn:time_steps_mod}, \pageref{defn:star_objects_mod} \\

$L(X,Y)$ &
Space of all linear maps $X\to Y$ &
\pageref{rem:geometric} \\

$\ell^{(\eps,i)}(\tau)$ & 
Renormalisation constant associated with the symbol $\tau$ & 
\pageref{sec:renConstants} \\ 

$\models_0$ & 
Closure of the renormalised canonical lifts of smooth noises & 
\pageref{text:models} \\ 

$\CM_{s,t}$ & 
Fixed-point map & 
\pageref{e:fixedpointmap} \\ 

$\mass$ & 
Mass of the Higgs field & 
\pageref{intro:eq-YMH-energy}\\

$M_T(B, k)$ & Gronwall and $\Fnorm{\place}$ losses & \pageref{notation:gronwall-F-pert} \\

$\modsol^{\alpha}$, $\modsol^{\alpha}_M$ & 
Modified solution spaces & 
\pageref{local:def-modified-solution-space}\\

$\mehler_\lambda$ & 
Mehler kernel & 
\pageref{eq:mehler-kernel-formula}\\

$\CN$ & 
Natural noise space & 
\pageref{text:CN}\\ 

$\mc{N}(B)$ & 
Quantity controlling the heat kernel parametrix error & 
\pageref{prop:pB-parametrix} \\ 

$N_{\lo}$ & 
Lowest frequency in Littlewood--Paley blocks & 
\pageref{local:eq-N-lo} \\ 

$\nl^{\alpha,\nu}$, $\nl^{\alpha,\nu}_M$ & 
Norms for nonlinearity in local theory & 
\pageref{local:def-nl-norm} \\ 

$\schO$ &
Schematic Landau symbol & 
\pageref{prelim:notation-Landau-standard-schematic} \\ 

$\bigO$ & 
Standard Landau symbol & 
\pageref{prelim:notation-Landau-standard-schematic} \\ 

$\fO_{\kappa}^\YMHs$ & 
Space of gauge orbits for gauge and Higgs fields & 
\pageref{rem:GaugeOrbits}, \pageref{e:gauge-orbits} \\ 

$\mc{P}$  & 
Set of cutoff functions & 
\pageref{main:eq-mcP} \\ 

$\flow_t$, $\wt{\flow}_t$ &
One-step evolution in the global iteration schemes &
\pageref{defn:one_step_iteration_scheme}, \pageref{defn:final_timestep} \\

$\fp$ & 
Euclidean heat kernel on torus & 
\pageref{notation:euclidean-heat-kernels} \\ 

$\fp^{\Euc}$ & 
Euclidean heat kernel on $\R^2$ & 
\pageref{notation:euclidean-heat-kernels} \\ 

$\hat\fp$ & 
Massive Euclidean heat kernel on $\R^2$& 
\pageref{e:massiveHeat} \\ 

$\fp_B$ & 
Covariant heat kernel for $\partial_t - \Delta_B$ & 
\pageref{kernel:def-kernels} \\

$\fp_{B,N}$ & 
Dyadic localisation of $\fp_B$ & 
\pageref{local:eq-pB-N}\\ 

$\fp_{\Y,B,N}$ & 
$\End(\Y)$-valued version of $\fp_{B,N}$ & 
\pageref{local:notation-unified} \\ 

$P_N$, $P_{\le N}$ & 
Littlewood--Paley projections & 
\pageref{main:eq-standard-Littlewood-Paley} \\ 

$\fq_B$ & 
Covariant heat kernel for $\partial_t - \Delta_B - 2 F_B$ & 
\pageref{kernel:def-kernels} \\ 

$\fq_{B,N}$ & 
Dyadic localisation of $\fq_B$ & 
\pageref{local:eq-qB-N}\\ 

$\fq_{\Y,B,N}$ & 
$\End(\Y)$-valued version of $\fq_{B,N}$ & 
\pageref{local:notation-unified} \\ 

$\RR$ & 
Renormalisation group & 
\pageref{text:renormalization-group} \\ 

$\Arem$, $\phirem$, 
$\Aremmod$, $\phiremmod$ &
Remainders after subtracting the linear stochastic objects &
\pageref{e:solution_decom} \\

$\qres_{\eps,N}^A, \qres_{\eps,N}^\higgs$ & 
Resonant parts of covariant quadratic stochastic objects & 
\pageref{objects:eq-resops} \\

$\moralfp$ & 
Nonnegative function satisfying a Gaussian heat kernel upper bound & 
\pageref{objects:notation-moralfp} \\ 

$r$ & 
Parameter & 
\pageref{e:small_constants} \\

$\solutionmap_{T,\rho,\eps}$ & 
Data-to-solution map &
\pageref{main:eq-data-to-solution} \\

$\fs_X$, $\fs_X^{F}$, $\narrowdots$ & 
Kernels appearing in perturbation around YM connections & 
\pageref{defn:s_B_eq} \\ 

$\sol^{\alpha,\eta}$, $\sol^{\alpha,\eta}_M$ & 
Solution spaces & 
\pageref{local:def-solution-space} \\ 

$\CS_{s,t}$ & 
Solution operator & 
\pageref{e:defSst} \\ 

$\Sdet$ & 
Deterministic term appearing in expansion of $A$ & 
\pageref{eq:Sdet-equation} \\ 

$\Sdet_m$, $\Sdet_m^\star$ &
$\Sdet$ term in the $m$-th iteration step &
\pageref{defn:time_steps} \\

$\Sdetmod_m$, $\Sdetmod_m^\star$ &
$\Sdet$ term in the modified iteration &
\pageref{defn:time_steps_mod} \\

$\tril$ & 
Trilinear term in schematic equation & 
\pageref{e:schematic} \\ 

$\CT$, $\CT[\tau]$ &
Model space and its subspace associated with the symbol $\tau$& 
\pageref{text:CT-appearence} \\

$T^{(1)}$, $T^{(2)}$ &
Time scales determined by initial data and curvature decrease &
\pageref{defn:tau} \\

$T^{(L)}$ &
Deterministic time step for $L$-admissible initial data &
\pageref{defn:L_admissible} \\

$T^{(L)}_{1/2}$ &
Intermediate deterministic time in the modified evolution &
\pageref{defn:final_timestep} \\

$t_m$, $t_m^\star$ &
Stopping times in the first global iteration scheme &
\pageref{defn:time_steps} \\

$\tmod_m$, $\tmod_m^\star$ &
Stopping times in the modified iteration scheme &
\pageref{defn:time_steps_mod} \\

$\tmod_{m+1/2}$, $\tmod_{m+1/2}^\star$ &
Intermediate stopping times in the modified iteration &
\pageref{defn:time_steps_mod} \\

$t_{m+1}^{\star\star}$ &
$\fF_{t_{m-1}}$-measurable time equal to $t_{m+1}$ on good event &
\pageref{defn:star_objects} \\

$\tmod_{m+1}^{\star\star}$ &
$\fF_{\tmod_{m-1/2}}^{\xi}$-measurable time equal to $\tmod_{m+1}$ on good event &
\pageref{defn:star_objects_mod} \\

$\tstop$ & 
Fixed time at which the background connection may jump & 
\pageref{notation:tstop} \\ 

$U$ & 
Linear part of the gauge action & 
\pageref{e:evolhU} \\ 

$(\ONBg_a)_a$, $(\ONBY_m)_m$ & 
Orthonormal bases of $\frkg$ and $\Y$ & 
\pageref{prelim:def-Casimir-induced},\pageref{local:notation-Y-valued-objects} \\ 

$\V$ & 
Target space for the Higgs field & 
\pageref{intro:assumption-setup} \\ 

$(\ONBV_m)_m$ & 
Orthonormal basis of $\V$ & 
\pageref{prelim:def-Casimir-induced} \\ 

$\genV$ & 
Generic vector space endowed with an action of $\liegroup$ & 
\pageref{e:defDA} \\ 

$W_\eps$, $\phiw_\eps$ & 
Nonlinear remainders & 
\pageref{local:eq-Z-def} \\

$X_\eps, \phix_\eps$ & 
Paracontrolled terms & 
\pageref{local:eq-paracontrolled-Ansatz} \\ 

$\Xwhite_{\eps}$, $\Xblack_{\eps}$ & 
Decomposition of $X_\eps$ & 
\pageref{local:eq-Xwhite-K} \\ 

$\Y$ & 
$\frkg^2 \oplus \V$ & 
\pageref{local:notation-unified} \\ 

$Y_\eps^{(i)}$ & 
Three main $\Y$-valued processes & 
\pageref{e:sys2}\\ 

$Y_\eps, \phiy_\eps$ & 
Smooth nonlinear remainders  & 
\pageref{local:eq-paracontrolled-Ansatz} \\ 

$Z^{(\eps,i)}$ & 
Canonical models associated with our three main processes& 
\pageref{def:Models} \\ 

$Z_\eps$, $\phiz_\eps$ & 
Nonlinear remainders & 
\pageref{local:eq-Z-def} \\

\end{longtable}
\end{center}

\renewcommand{\leftmark}{\textsc{References}}
\renewcommand{\rightmark}{\textsc{References}}
\bibliography{Library}
\bibliographystyle{Martin}

\end{document}